\documentclass[11pt,leqno]{book}

\usepackage{amsmath,amsfonts,amscd,amssymb,theorem}

\long\def\comment#1\endcomment{}

\makeatletter
\begingroup
\gdef\th@dotted{\normalfont\itshape
  \def\@begintheorem##1##2{%
        \item[\hskip\labelsep \theorem@headerfont ##1\ ##2.]}%
\def\@opargbegintheorem##1##2##3{%
   \item[\hskip\labelsep \theorem@headerfont ##1\ ##2\ (##3).]}}
\endgroup
\makeatother

\theoremstyle{dotted}

\newtheorem{theorem}{Theorem}[subsection]
\newtheorem{lemma}[theorem]{Lemma}

\newtheorem{prop}[theorem]{Proposition}
\newtheorem{corr}[theorem]{Corollary}

\makeatletter
\begingroup
\gdef\th@upshape{\normalfont
  \def\@begintheorem##1##2{%
        \item[\hskip\labelsep \theorem@headerfont ##1\ ##2.]}%
\def\@opargbegintheorem##1##2##3{%
   \item[\hskip\labelsep \theorem@headerfont ##1\ ##2\ (##3).]}}
\endgroup
\makeatother

\theoremstyle{upshape}

\newtheorem{defn}[theorem]{Definition}
\newtheorem{remark}[theorem]{Remark}
\newtheorem{exa}[theorem]{Example}

\makeatletter
\begingroup
\gdef\th@upnonum{\normalfont
  \def\@begintheorem##1##2{%
        \item[\hskip\labelsep \theorem@headerfont ##1.]}%
\def\@opargbegintheorem##1##2##3{%
   \item[\hskip\labelsep \theorem@headerfont ##1.]}}
\endgroup
\makeatother

\theoremstyle{upnonum}

\newtheorem{remarknonum}[theorem]{Remark}

\makeatletter
\renewcommand{\subsection}{\@startsection{subsection}{2}{0pt}{-3ex
plus -1ex minus -0.2ex}{-2mm plus -0pt minus
-2pt}{\normalfont\bfseries}} 
\renewcommand{\subsubsection}{\@startsection{subsubsection}{3}{0pt}{-3ex
plus -1ex minus -0.2ex}{-2mm plus -0pt minus
-2pt}{\normalfont\bfseries}} 
\makeatother

\makeatletter
\@addtoreset{equation}{subsection}
\makeatother
\renewcommand{\theequation}{\thesubsection.\arabic{equation}}

\newcommand{\cntrct}                
{\hspace{2pt}\raisebox{1pt}{\text{$\lrcorner$}}\hspace{2pt}}

\newcommand{\proof}[1][Proof.]{\smallskip\noindent{\em #1}}
\def\endproof{\hfill\ensuremath{\square}\par\medskip}

\renewcommand{\labelenumi}{{\normalfont(\roman{enumi})}}

\def\eqref#1{\thetag{\ref{#1}}}

\let\latexref=\ref
\def\ref#1{{\normalfont{\latexref{#1}}}}

\newcommand{\wt}{\widetilde}

\newcommand{\dg}{\dagger}

\newcommand{\idot}{{\:\raisebox{1pt}{\text{\circle*{1.5}}}}}
\newcommand{\hdot}{{\:\raisebox{3pt}{\text{\circle*{1.5}}}}}
\newcommand{\eps}{\varepsilon}
\newcommand{\ups}{\upsilon}
\renewcommand{\phi}{\varphi}

\newcommand{\F}{\mathcal{F}}

\newcommand{\Hom}{\operatorname{Hom}}
\newcommand{\Hhom}{\operatorname{\mathcal{H}{\it om}}}
\newcommand{\RHom}{\operatorname{RHom}}

\newcommand{\Ex}{\operatorname{Ex}}

\newcommand{\Fun}{\operatorname{Fun}}
\newcommand{\vFun}{\overrightarrow{\Fun}}
\newcommand{\fFun}{\operatorname{{\mathcal F}{\mathit u}{\mathit n}}}

\newcommand{\Sec}{\operatorname{Sec}}
\newcommand{\sSec}{\operatorname{{\mathcal S}{\mathit e}{\mathit c}}}

\newcommand{\Ch}{\operatorname{Ch}}

\newcommand{\id}{\operatorname{\sf id}}
\newcommand{\Id}{\operatorname{\sf Id}}

\newcommand{\Tot}{\natural}
\newcommand{\Toth}{{\Tot}h}

\newcommand{\A}{\mathcal{A}}
\newcommand{\cD}{\mathcal{D}}
\newcommand{\D}{{\sf D}}
\newcommand{\C}{{\cal C}}

\newcommand{\ff}{{\sf f}}
\newcommand{\dd}{{\sf d}}

\newcommand{\hush}{\natural}
\newcommand{\hhush}{{\hush\hush}}
\newcommand{\fhush}{{\flat\hush}}
\newcommand{\hushf}{{\hush\flat}}
\newcommand{\hash}{\sharp}
\newcommand{\tr}{\triangle}

\newcommand{\trl}{\triangleleft}
\newcommand{\trr}{\triangleright}

\newcommand{\Unf}{K}
\newcommand{\oUnf}{\overline{\Unf}}
\newcommand{\unf}{k}

\newcommand{\Sets}{\operatorname{Sets}}
\newcommand{\sSets}{\operatorname{\mathcal{S}{\it ets}}}
\newcommand{\Cat}{\operatorname{Cat}}
\newcommand{\cCat}{\operatorname{\mathcal{C}{\it at}}}
\newcommand{\Maps}{\operatorname{Maps}}
\newcommand{\Iso}{{\star}}
\newcommand{\Aut}{{\operatorname{Aut}}}
\newcommand{\End}{{\operatorname{End}}}

\newcommand{\Env}{\operatorname{Env}}

\newcommand{\iInd}{\operatorname{\mathcal{I}\mathit{n}\mathit{d}}}
\newcommand{\Indsec}{\operatorname{\mathcal{I}\mathit{n}\mathit{d}\mathcal{S}\mathit{e}\mathit{c}}}
\newcommand{\Comp}{\operatorname{\mathcal{C}\mathit{o}\mathit{m}\mathit{p}}}

\newcommand{\Cone}{\operatorname{Cone}}
\newcommand{\cCone}{\operatorname{\mathcal{C}{\mathit o}{\mathit
      n}{\mathit e}}}

\newcommand{\Ind}{\operatorname{\sf Ind}}

\newcommand{\ppt}{{\sf pt}}

\newcommand{\bpi}{\overline{\pi}}

\newcommand{\Ab}{\operatorname{Ab}}

\newcommand{\M}{\operatorname{\sf M}}

\newcommand{\copr}{\sqcup}

\newcommand{\Z}{{\mathbb Z}}
\newcommand{\N}{{\mathbb N}}

\newcommand{\E}{\mathcal{E}}
\newcommand{\Kk}{\mathcal{K}}

\newcommand{\sk}{\operatorname{\text{\rm sk}}}

\newcommand{\Exp}{\operatorname{Exp}}

\newcommand{\Con}{\operatorname{Con}}

\newcommand{\Ho}{\operatorname{Ho}}
\newcommand{\Hho}{\operatorname{\mathcal{H}\mathit{o}}}

\newcommand{\hocolim}{\operatorname{\sf hocolim}}
\newcommand{\holim}{\operatorname{\sf holim}}
\newcommand{\colim}{\operatorname{\sf colim}}
\renewcommand{\lim}{\operatorname{\sf lim}}
\newcommand{\colimh}{\colim^h}
\newcommand{\limh}{\lim^h}

\newcommand{\V}{{\sf V}}
\newcommand{\W}{{\sf W}}
\newcommand{\VV}{\mathbb{V}}

\newcommand{\Y}{{\sf Y}}
\newcommand{\y}{{\sf y}}
\newcommand{\X}{{\sf X}}
\newcommand{\ZZ}{{\sf Z}}

\newcommand{\Xx}{\mathcal{X}}
\newcommand{\Yy}{\mathcal{Y}}
\newcommand{\Zz}{\mathcal{Z}}

\newcommand{\ev}{\operatorname{\sf ev}}

\newcommand{\bC}{\overline{\C}}

\newcommand{\G}{\mathcal{G}}

\newcommand{\tors}{\operatorname{\text{\rm-Tors}}}

\newcommand{\Hh}{\mathcal{H}}

\newcommand{\Top}{\operatorname{Top}}
\newcommand{\CW}{\operatorname{CW}}

\renewcommand{\SS}{{\sf S}}

\newcommand{\Tel}{\operatorname{Tel}}

\newcommand{\Pos}{\operatorname{PoSets}}
\newcommand{\pPos}{\operatorname{\mathcal{P}{\it o}\mathcal{S}{\it ets}}}
    
\newcommand{\Posf}{\operatorname{Pos}}
\newcommand{\Posff}{\operatorname{pos}}

\newcommand{\Fib}{\operatorname{Fib}}

\newcommand{\Rel}{\operatorname{BiPoSets}}
\newcommand{\bRel}{\overline{\Rel}}
\newcommand{\Relf}{\operatorname{BiPos}}
\newcommand{\bRelf}{\overline{\Relf}}

\newcommand{\Seg}{\operatorname{Seg}}
\newcommand{\CSS}{\operatorname{CSS}}

\newcommand{\Cyl}{{\sf C}}

\newcommand{\Red}{\operatorname{Red}}
\newcommand{\cRed}{\operatorname{\mathcal{R}{\mathit e}{\mathit d}}}

\newcommand{\pP}{\mathcal{P}}

\newcommand{\I}{\mathcal{I}}
\newcommand{\II}{\mathbb{I}}

\newcommand{\QQ}{{\sf Q}}
\newcommand{\Nn}{{\sf N}}
\newcommand{\bQQ}{\overline{\QQ}}

\newcommand{\T}{{\sf T}}
\newcommand{\bT}{\overline{\T}}

\newcommand{\bb}{\operatorname{/\!\!/}}
\newcommand{\bbd}{\operatorname{\setminus\!\!\setminus}}
\newcommand{\bbdi}{\bbd_\Iso}
\newcommand{\bbi}{\operatorname{/\!\!/\!_\Iso}}

\newcommand{\mm}{\operatorname{/\!\!/\!_\flat}}
\newcommand{\mc}{\operatorname{/\!\!/\!_\sharp}}
\newcommand{\mmh}{\mm^h}

\newcommand{\sE}{{\sf E}}
\newcommand{\bE}{\overline{\E}}
\newcommand{\bsE}{\overline{\sE}}

\newcommand{\bN}{\overline{\Nn}}

\newcommand{\dm}{\diamond}

\newcommand{\wV}{\wt{\V}}

\newcommand{\Cof}{\operatorname{Cof}}
\newcommand{\Dec}{\operatorname{Dec}}

\newcommand{\Tw}{\operatorname{\sf tw}}
\newcommand{\Ar}{\operatorname{\sf ar}}
\newcommand{\Arr}{\operatorname{\sf Ar}}

\newcommand{\Arc}{\Ar_c}

\newcommand{\J}{\mathcal{J}}

\newcommand{\bsigma}{\overline{\sigma}}

\newcommand{\hht}{\operatorname{ht}}

\newcommand{\edge}{\operatorname{\sf e}}

\newcommand{\Gr}{\operatorname{\sf Gr}}

\newcommand{\ssetminus}{\smallsetminus}

\def\emptyset{\varnothing}

\newcommand{\bB}{{\sf b}}

\title{Enhancement for categories and homotopical algebra}

\author{D. Kaledin}

\begin{document}

\maketitle

\tableofcontents

\chapter*{Introduction.}

\section*{Generalities.}

One can say that algebra and geometry are related by something like
the Heisenberg Uncertainty Principle: a perfect algebraic
formalization of a mathematical theory often obscures its geometric
features, while a clear geometric idea can be devilishly difficult
to put into words and turn into a formal algebraic proof. Algebra
and geometry complement each other but do not mix well. The pair has
been famously compared to an angel and a devil, although which is
which is perhaps debatable.

Linguistically, geometry is about earth --- not too angelic, mind
you, although not necessarily chthonic, either. Philosophically, the
basic underlying notion of geometry is that of a ``space''. This has
been around for millenia, since the ancient Greeks at least, and it
can be made precise in different ways, none of them encompassing the
whole thing. In full generality, all one can venture to say is that
a space has points, and there is also some notion of ``closedness''
and/or ``continuity'' associated to this collection of points. The
latter idea is again very old and ingrained -- to the point that
when Cantor was creating modern set theory, he really had to
struggle to convince people that it makes sense to consider a set as
simply a set, with no continuity in any form anywhere in the
picture.

Algebra is much more recent. The word itself is medieval, but it
used to mean something tautological; it is probably safe to say that
algebra with actual mathematical content starts with Galois
theory. This brings in another fundamental and ancient idea that,
however, was never formalised before Galois --- the idea of
``symmetry'', with the corresponding formal notions of a group and a
group action.

The spirit here is decidedly modern; as befits something born in the
industrial era, algebra is instrumental rather than meditative. One
can say that if at the end of the day, geometry is about pictures,
then algebra is about operations: you do not just observe things but
act on them.  Geometrically, if you want to connect two points, you
imagine a curve passing through them, while algebraically, you
construct an operator sending one to the other.

\medskip

Even for a text of more than 600 pages, one page of philosophy is
perhaps enough, so let us get to the point. The point is homotopy
theory.

\medskip

The origin of homotopy theory is of course geometric --- it first
appeared as ``algebraic topology'', the study of topological spaces
and continuous maps between them considered ``up to homotopy''. Here
a space is simply a topological space, maybe with some conditions
making it ``nice enough'', and two points are ``close'' if they are
connected by a continuous path. Topology becomes ``algebraic'' once
you observe that up to homotopy, paths can be inverted, and based
loops form a group --- the fundamental group.

Algebraically, the same fundamental group emerges through the notion
of a covering, in a well-known beautiful parallel with Galois
theory. A choice of a base point determines a choice of a universal
covering, and then the fundamental group acts transitively on the
fibers of this universal covering. To preserve all the symmetries of
the situation and avoid choices, it is better to consider all points
at once; the relevant structure is then not a group but a groupoid,
the fundamental groupoid of a topological space. In this
description, we no longer consider the set of points as a ``space'',
with some continuity properties; rather, the additional structure is
isomorphisms between points. Following recent fashion, one might say
that the set is ``animated'' by allowing points to have non-trivial
isomorphisms (and in particular, non-trivial automorphisms).

If one is only interested in topology, then it is not so important
whether one defines the fundamental group by looking at based loops
or via the Galois-style approach with coverings; the latter
streamlines the proof of the Van Kampen Theorem, but that's the
extent of it. However, what the Galois description strongly suggests
is that homotopy theory is something much more general and
widespread than mere algebraic topology. Indeed, many mathematical
objects --- starting with sets themselves, actually --- do not
naturally form a set, nor a ``space'', whatever that is. Even if we
ignore the size issues and only consider finite sets, it is still
completely unnatural, and some would say meaningless, to state that
two sets are ``equal''. The natural notion is ``isomorphic'', not
equal, and an isomorphism is an additional structure that one needs
to preserve. We have moved beyond Cantor. It is groupoids that
appear in nature, not sets, and ``animation'' is a crucial part of
it. It would be much better to have a formulation of homotopy theory
that has this fact at its core, since it might then be applicable in
many situations when there are no spaces in sight.

\medskip

The problem is, it is not clear what this might be --- in
particular, just looking at groupoids is definitely not enough.

\medskip

The most obvious deficiency is that groupoids do not describe all
``homotopy types'', or ``spaces up to a homotopy equivalence'': a
simply connected space is not necessarily contractible. What one
needs to describe them all is some kind of ``enhanced Galois
theory'' that would also include all the higher homotopy groups (in
particular, all the homotopy groups of spheres). However, even if
one ignores topological realities and looks for a purely formal
algebraic theory, it has purely formal problems that also suggest a
need for some sort of enhancement.

Ideally, one would like to do something similar to what Grothendieck
did with homological algebra. It started its life as the study of
homology and cohomology of topological spaces, a younger cousin to
homotopy theory. However, very soon it became clear that the actual
story has very little to do with topology. At least since
\cite{toho}, we understand that the real subject of homological
algebra is abelian categories and derived functors. There is a clear
problem --- some functors are only exact on the left or on the
right. Looking for a solution to this problem, one is almost
inevitably led to looking at resolutions, characterizing derived
functors by a universal property, replacing objects by chain
complexes, and eventually, to the derived categories of
\cite{ver.the}.

With homotopy theory, we do not have such a clear-cut single problem
whose solution would allow one to bootstrap the whole story starting
from first principles. However, we do have problems we can use as an
entry point. One set of problems concerns the procedure of {\em
  localization}.

\section*{Localization (and its discontents).}

In a nutshell, localization is a universal procedure for creating
new symmetries. Formally, one observes that many mathematical
objects form not just a groupoid but a category --- not all
morphisms are invertible. For any category $\C$ with a class of
morphisms $W$, say that a functor $\gamma:\C \to \E$ to some
category $\E$ {\em inverts $W$} if $\gamma(w)$ is invertible for any
$w \in W$. Then the {\em localization} --- or the {\em category of
  fractions}, in the original terminology of \cite{GZ} --- of $\C$
in the class $W$ is a category $h^W(\C)$ equipped with a functor
$h:\C \to h^W(\C)$ such that \thetag{i} $h$ inverts $W$, and
\thetag{ii} for any functor $\gamma:\C \to \E$ that inverts $W$,
there exists a functor $\gamma':h^W(\C) \to \E$ equipped with an
isomorphism $a:\gamma \to \gamma' \circ h$, and the pair $\langle
\gamma',a \rangle$ is unique up to a unique isomorphism.

Informally, localizing a category forgets some information deemed
in\-ess\-ential. If we have a set, we can impose an equivalence
relation and force some elements to become equal (e.g.\ forget an
integer and only remember its residue $\mod n$ for some $n$). For a
category, this is not a good idea, since objects in a category are
never equal. Instead, we force some morphisms to become
isomorphisms. In particular, this forces some objects to become
isomorphic, but it also keeps track of the animation. In the extreme
case, we can take $W$ to be the class of all morphisms. Then the
localization $h^W(\C)$ is a groupoid whose isomorphism classes of
objects are simply the connected components of the original category
$\C$, but the automorphism groups can be quite complicated.

Starting from approximately \cite{GZ}, and certainly from
\cite{qui.ho}, it became an accepted wisdom that the real content of
homotopy theory is exactly the localization procedure, and what it
gives as its output. Unfortunately, the procedure itself is highly
non-trivial. If the category $\C$ is large, the localization need
not exists. If $\C$ is small, $h^W(\C)$ always exists, but it is in
general very hard to describe and control.

The motivating example for \cite{GZ} is Verdier localization of
\cite{ver.the}, where the localization can be controlled rather
effectively. However, there is an even simpler example of a
well-behaved localization problem that is so fundamental as to pass
almost unnoticed. Namely, let $\C$ be the category $\Cat$ of small
categories, and let $W$ be the class of equivalences of
categories. Then $h^W(\Cat)$ exists, and coincides with the category
$\Cat^0$ whose objects are again small categories, and whose
morphisms are isomorphism classes of functors (this is an easy and
standard exercise, but just in case, a proof can be found below in
Lemma~\ref{cat.loc.le}).

In real life, categories --- and groupoids, for that matter --- are
only considered up to an equivalence; one would like to work in
$\Cat^0$ rather than $\Cat$. For example, a ``commutative square of
categories'' is usually understood as a square
\begin{equation}\label{sq.0.cat}
\renewcommand{\theequation}{$\star$}
\begin{CD}
\C_{01} @>{\gamma^1_{01}}>> \C_1\\
@V{\gamma_{01}^0}VV @VV{\gamma_1}V\\
\C_0 @>{\gamma_0}>> \C
\end{CD}
\end{equation}
equipped with an isomorphism $\alpha:\gamma_0 \circ \gamma^0_{01}
\to \gamma_1 \circ \gamma^1_{01}$; the isomorphism is usually not
written down but always implied. Even if all the categories are
small, such a square is {\em not} a commutative square in $\Cat$
unless $\alpha$ happens to be an identity map, but it always defines
a commutative square in $\Cat^0$. A square \thetag{$\star$} can be {\em
  cartesian} or {\em cocartesian} if it has the usual universal
properties with respect to squares of the same form (for a precise
definition, see Subsection~\ref{dia.subs} below). Cartesian squares
always exists --- one can take $\C_{01}$ to be the category $\C_0
\times_{\C} \C_1$ of triples $\langle c_0,c_1,\alpha \rangle$ of
objects $c_0 \in \C_0$, $c_1 \in \C_1$ and an isomorphism
$\alpha:\gamma_0(\C_0) \to \gamma_1(\C_1)$. Cocartesian squares do
not always exist if the categories are large, but they are still
useful --- for example, a localization is an example of a
cocartesian square (take $\C_0 = \C$, $\C_1 = W$ considered as
discrete category, and let $\C_{01} = W \times [1]$, where $[1]$ is
the single arrow category with two non-trivial objects $0,1 \in [1]$
and a single non-identity map $0 \to 1$). However, note that
cartesian resp.\ cocartesian squares \thetag{$\star$} {\em are not
  cartesian resp.\ cocartesian} as commutative squares in
$\Cat^0$. What happens is, the corresponding universal property uses
the unwritten isomorphism $\alpha:\gamma_0 \circ \gamma^0_{01} \to
\gamma_1 \circ \gamma^1_{01}$, but simply passing to $\Cat^0$
completely forgets it. In fact, $\Cat^0$ does {\em not} have finite
limits nor colimits, and treating it as simply a category is not a
good idea. While the property of a square \thetag{$\star$} to be cartesian
or cocartesian only depends on the corresponding square in $\Cat^0$,
so it makes sense to speak of ``enhanced-cartesian'' or
``enhanced-cocartesian'' squares in $\Cat^0$, those do not have a
universal property and cannot be described purely in terms of the
category $\Cat^0$. Given $\C_0,\C_1,\C \in \Cat^0$ and isomorphism
classes of functors $\gamma_l:\C_l \to \C$, $l=0,1$, one can
construct an enhanced-cartesian square \thetag{$\star$} in $\Cat^0$, and
it will be unique up to an isomorphism, but the isomorphism is not
unique.

Verdier localization follows the same general pattern. One starts
with an abelian category $\A$, considers the category $C_\idot(\A)$
of chain complexes in $\A$, and localizes with respect to the class
of chain-homotopy equivalences (see Subsection~\ref{add.subs} below
for a formal definition). The resulting category $\Ho(\A)$ can be
equivalently described as the category of chain complexes in $\A$
and chain-homotopy classes of morphisms between them. Then while
$C_\idot(\A)$, being an abelian category, has kernels and cokernels,
these do not survive after passage to $\Ho(\A)$. What survives is
cones. These are unique up to a non-unique isomorphism, and cannot
be recovered from the category $\Ho(\A)$; they consitute an
additional structure on $\Ho(\A)$.

The only thing one can do is to axiomatize this additional
structure, and this is exactly what Verdier did, by introducing the
notion of a ``triangulated category''. A triangulated structure on a
category is an early form of enhancement --- in fact, the earliest
form of enhancement in the literature --- and it does allow one to
prove something: once you have a triangulated category $\cD$, a
general Verdier Localization Theorem allows you to localize $\cD$
with respect to certain classes of morphisms in a very controlled
way. Moreover, the resulting category also inherits a triangulated
structure, so the procedure can be repeated. This is exactly how
derived categories are defined in \cite{ver.the} (what one inverts
is the class of quasiisomorphisms in $\Ho(\A)$).

One obvious drawback of Verdier localization is that it only works
in the additive context --- we are dealing with homology, not
homotopy.  The stated goal of \cite{GZ} is to try to understand what
is going on, and maybe even try to find a non-additive
generalization. In principle, such a generalization is not
impossible. Instead of kernel and cokernels, one would deal with
cartesian and cocartesian squares; these would need to be imposed as
an extra structure on the category in question, and subjected to
some axioms (hopefully, a finite number of those). However, to the
best of our knowledge, this has never been. Partly, this could be
just too difficult --- homotopy is much harder to compute than
homology, so while it might be possible that a workable set of
axioms exists, this is by no means guaranteed. However, it equally
might be just not worth the effort. Already in the additive context,
it was realized pretty soon after \cite{ver.the} that the
enhancement given by a triangulated structure is much too weak.

Various problem with triangulated categories are very well
illuminated in the existing literature (for a partial list written
from the same perspective as here, see
e.g.\ \cite[Introduction]{ka.glue}). Let us only mention one of
them. In general, given an ``enhanced'' triangulated category $\cD$,
one expect to have a spectrum $\cD(A,B)$ of morphisms for any two
objects $A,B \in \cD$; the ``naive'' set of morphisms $A \to B$
should be $\pi_0$ of this spectrum. If one is only interested in a
linear situation --- say, if one works with derived category of an
abelian category linear over a field $k$ --- one expects to at least
have an object $\RHom^\hdot(A,B) \in \cD(k)$ in the derived category
of $k$-vector spaces, with naive morphisms given by degree-$0$
homology classes. Moreover, this should be compatible with the
category structure, so that in particular, $\RHom^\hdot(A,A)$ should
be a DG algebra over $k$, well-defined up to a quasiisomorphism. The
triangulated structure gives you all the homotopy groups of this
hypothetical spectrum, with $\pi_n\cD(A,B) = \pi_0\cD(A,B[-n])$, but
not the spectrum itself. In the $k$-linear situation, a complex of
$k$-vector spaces is quasiisomorphic to its cohomology, so we do
recover $\RHom^\hdot(A,B) \in \cD(k)$, but the construction is not
functorial, and certainly does not give you a DG algebra
$\RHom^\hdot(A,A)$.

In 1980-ies, a great codification of homological algebra was done by
Gelfand and Manin in \cite{GM}, the first textbook on derived and
triangulated categories. The problem of enhancement features quite
prominently there, and a search for a good theory of enhanced
triangulated categories was a very active pursuit, back then. By the
end of the decade, people gave up. Even in a linear context, the
problem was considered hopeless. The consensus became, roughly, that
an enhanced $k$-linear triangulated category is the same thing as a
$k$-linear DG category ``considered up to a quasiequivalence'', and
that's the best one can say. At this point, the story goes full
circle and comes back to homotopy. Indeed, while individual DG
categories are in some sense linear, the category of DG categories
certainly is not, and localizing it with respect to
quasiequivalences is a non-linear localization problem. Thus even
when dealing with homology, one needs the full force of non-additive
localization.

\section*{Categories of models.}

At this point, we go back to 1960-ies and to the perhaps most
important work in the history of homotopical algebra --- to
Quillen's \cite{qui.ho}. This appeared slightly later than
\cite{GZ}, gave the name to the subject, and to a large extent, set
the paradigm that is current even today.

It is not clear whether Quillen was particularly happy with
\cite{qui.ho} --- it seems that he never referred to it in his later
work --- or whether he intended it to become such a foundational
text. His stated purpose is quite modest and
utilitarian. Originally, when the only model for homotopy types was
topological spaces, people were quite happy to work with them ``up
to homotopy'' and leave it at that. However, already in 1950-ies,
especially with the discovery of Kan fibrations, it became clear
that ``the same'' homotopy theory can be successfully modeled by a
completely different gadget --- namely, by simplicial sets (or
alternatively, by simplicial groups). Moreover, the situation does
not only occur for homotopy types --- for example, in rational
homotopy theory, commutative DG algebras over a field $k$ of
characteristic $0$, with some restrictions, have ``the same''
homotopy theory as DG Lie algebras over the same field, again with
some restrictions. The problem is, in order to prove this, one first
has to define the meaning of homotopy theory being ``the same''.

Like Verdier, Quillen starts by axiomatizing the situation. He
considers various ``categories of models'' and finds a short and, in
retrospect, very successful set of additional structures and axioms
one needs to impose in order to have a ``homotopy theory''. Unlike
Verdier's triangulated structure, this really only applies to the
categories of models and not to their localizations --- the very
first axiom says that a model category (or rather, in Quillen's
terminology, a ``closed model category'') has all finite limits and
colimits, and as we saw, localization usually destroys
this. Nevertheless, Quillen's axioms insure that the localization
problem is solvable, and in a very controlled way. In addition to a
class $W$ of ``weak equivalences'', a model category $\C$ is
supposed to have classes $C$ and $F$ of ``cofibrations'' and
``fibrations'', and the whole structure insures that the
localization $h^W(\C)$ exists and admits a rather explicit
description. As an application of this, Quillen shows that under
some conditions, an adjoint pair of functors between model
categories $\C_0$, $\C_1$ gives rise to an adjoint pair of functors
between their localizations $h^{W_0}(\C_0)$, $h^{W_1}(\C_1)$, and
moreover, under further conditions, the functors between
localizations are an inverse pairs of equivalences. This is known as
``Quillen adjunction'' and ``Quillen equivalence'', and it is used
by Quillen to prove that various pairs of model categories have
localizations that are naturally equivalent.

However, there is more --- at least in the current paradigm of
homotopical algebra, there is the following fundamental principle.
\begin{itemize}
\item Quillen-equivalent model categories have the same homotopy
  theory.
\end{itemize}
This is an article of faith and not a theorem, since ``homotopy
theory'' --- or, say, ``enhanced localization'' --- is not defined
at this point. Moreover, it cannot be defined, unless one already
assumes that \thetag{$\bullet$} holds. The problem is exactly in the
``enhanced'' part.

At the very least, one expects that, just as in the additive
situation, an ``enhanced localization'' $\Hh^W(\C)$ of a category
$\C$ with respect to a class of morphisms $W$ includes a homotopy
type $\Hh^W(\C)(c,c')$ of morphisms between any two objects $c,c'
\in \C$ (although in the non-linear situation, it is just an
unstable homotopy type and not a spectrum). The morphisms $c \to c'$
in the ``naive'' localization $h^W(\C)$ form the set
$\pi_0\Hh^W(\C)(c,c')$, but $\Hh^W(\C)(c,c')$ can have higher
homotopy groups. In the basic example $\C = \Cat$,
$\Hh^W(\Cat)(I,I')$ is the homotopy type corresponding to the
groupoid of functors $I \to I'$, so $\pi_n = 0$ for $n \geq 2$, but
$\pi_1$ can be non-trivial.

Now, homotopy groups $\pi_\idot\Hh^W(\C)(c,c')$ are just groups, so
one can meaningfully ask whether one can define them when $\C$ is a
model category, and whether they are preserved by Quillen
equivalences --- one can, they are. This is actually proved in
\cite{qui.ho}, where, moreover, \thetag{$\bullet$} is stated and
proved as a theorem --- but under a very restricted understanding of
``homotopy theory'' that Quillen himself calls unsatisfactory. In
general, however, a homotopy type is much more than a collection of
homotopy groups, and it has to be modeled by something, up to some
notion of equivalence. But as we saw, the whole reason for the
formalism is that --- unlike in the linear case, where the only
practical model is a chain complex --- we have at least two
different models, topological spaces and simplicial sets, and one
needs to be sure that they do model the same thing. We are back to
\thetag{$\bullet$}.

In practice, one not only needs homotopy types $\Hh^W(\C)(-,-)$ but
also composition operations on them, and these have to be ``homotopy
coherent'' in some sense. This leads to a large variety of possible
models for, for lack of better term, ``enhanced categories''. Let us
list some of them.
\begin{enumerate}
\item The most naive thing one can consider is small categories
  enriched in topological spaces, or in simplicial sets (that is, we
  have a topological space or a simplicial set of maps $c \to c'$
  for any two objects $c$, $c'$, and all the operations are
  compatible with this additional structure on the nose). This is
  not very convenient for practical work, since weak equivalences
  are not very explicit, and enriched categories are not always easy
  to construct, but it can be done --- a model structure
  exists. Moreover, the Dwyer-Kan localization procedure
  \cite{dw-ka}, probably the first homotopical localization
  procedure in the literature, actually produces a category enriched
  in simplicial sets, on the nose.
\item To make enriched categories slightly more flexible, one can do
  two things. Firstly, a small category is conveniently encoded by
  its nerve, and conversely, a simplicial set is the nerve of a
  small category iff it satisfies the so-called {\em Segal
    condition} (for formal statements, see
  e.g.\ Subsection~\ref{seg.X.subs} belows). For a small category
  enriched in simplicial sets, the nerve is naturally a bisimplicial
  set. Then one can put a model structure on the category of
  bisimplicial sets, and only impose the Segal condition on the
  homotopical level. Secondly, a further improvement has been
  suggested by Ch. Rezk \cite{rzk}: he introduces the notion of a
  {\em complete Segal space}, a fibrant bisimplicial set satisfying
  an additional condition, and shows that bisimplicial sets admit a
  model structure whose fibrant objects are exactly complete Segal
  spaces, and weak equivalences between those are just pointwise
  weak equivalences in one of the two simplicial directions.
\item As an alternative to Segal spaces, it was discovered by
  A. Joyal that one can actually squeeze the two simplicial
  directions of a bisimplicial set into one, by combining the Segal
  condition and the defining properties of a Kan complex in a clever
  way. The result is a new model structure on the category of
  simplicial sets whose fibrant objects are the so-called {\em
    quasicategories}.
\item In a different vein, as shown by C. Barwick and D. Kan
  \cite{BK}, one can actually start with the general principle that
  every enhanced category is obtained by localization; one defines a
  {\em relative category} as a pair of a small category $\C$ and a
  collection $W$ of morphisms in $\C$, and shows that pairs $\langle
  \C,W \rangle$ also form a model category. Weak equivalences are
  functors that induce weak equivalences of Dwyer-Kan localizations,
  so they are even less explicit than in \thetag{i}, but the
  localization procedure is pleasantly simple.
\end{enumerate}
All of the model categories listed above are Quillen-equivalent, so
in the paradigm of \cite{qui.ho}, either one can be used as
foundations for abstract homotopy theory. Complete Segal spaces are
probably closest to intuition and the most pleasing conceptually, so
it is tempting to conclude that this is it (as has been argued very
persuasively by B. To\"en in \cite{toen-inf}). Quasicategories allow
for shorter proofs; this is probably the reason they have been
chosen by J. Lurie in \cite{lu1} and his subsequent foundational
work. Relative categories of Barwick and Kan are important as a
proof-of-concept; while not much practical work has been done in
this setting, the fact that it exists is very illuminating.

\section*{Re-animating the sets.}

In a sense, as long as one stays in the paradigm of \cite{qui.ho},
one can consider the enhanced localization problem completely
solved. Just as expected, there indeed exists an abstract homotopy
theory that can be applied in a wide variety of contexts, certainly
much more general than just topology, and a choice of a particular
model becomes a matter of convenience and/or personal
taste. However, the situation is paradoxical. Localization is
supposed to get rid of irrelevant information and create new
symmetries. Enhanced localization as we have it does create new
symmetries but only modulo homotopy --- that is, well, after
localization; in reality, it does not. There are no ``animated
sets'', we are still dealing with spaces. We just managed to produce
them in a much larger generality. Nor do we really get rid of
irrelevant information, we just hide it under the rug and use the
machinery of \cite{qui.ho} to make sure it stays there. This is
especially clear in the relative categories model: here ``enhanced
localization'' by definition does exactly nothing. We do not really
solve the problem, we just find a way to say precisely which
problems have ``the same'' solution, and leave it at that.

Of course, modern abstract homotopy theory looks somewhat more
algebraic since ``spaces'' now mean simplicial sets rather than
actual topological spaces. However, we now understand --- see
e.g.\ \cite{dr} or \cite{lein} --- that this is an illusion:
simplicial sets are really very close to topological spaces, and
rely on identical intuition and working procedures. We still work up
to homotopy, not up to an isomorphism. The parallel with Galois
theory is completely lost.

Since approximately late 1990-ies, this situation is simply accepted
as a fact of life, but in the 30-odd years between \cite{qui.ho} and
modern treatments of the subject such as \cite{hovey}, there were
attempts to find something better. The most famous is of course
Grothendieck's ``thousand pages long letter'' to Quillen,
\cite{champ}, that has no theorems but a lot of fantastic ideas.
One such was an idea of an ``$\infty$-groupoid'', a gadget that
would successfully model all homotopy types by not only having
objects and morphisms, but also morphisms between morphisms, and so
on. Unfortunately, this turned out to be something of a red herring
--- while many rigourous definitions were proposed by many people
over the years, none of them were too natural, and they all were
defined up to some elaborate notion of a weak equivalence that had
to be controlled by model category techniques. Thus in practice, it
was more of the same, just more complicated. At the end of the day,
in current usage, an ``$\infty$-groupoid'' just means a Kan-fibrant
simplicial set; ``enhanced Galois theory'' occasionally appears on
the level of slogans but has no mathematical content.

However, there is another fantastic idea of \cite{champ} that has
received much less attention, and this is strange, since this idea
does indeed work. If one does consistently what Grothendieck
suggested, one arrives at a formulation of homotopy theory that is
pretty close to the original Galois-style intuition, and completely
independent of \cite{qui.ho} and the corresponding paradigm.

Within the context of Grothendieck's work, the idea is actually
quite natural. For comparison, recall that a scheme $X$ can be
thought of in two ways. On one hand, we have the formal definition
(a locally ringed space that is locally isomorphic to the Zariski
spectrum of a commutative ring). On the other hand, we have a set of
points of $X$, and moreover, a set $X(A)$ of $A$-points with values
in any commutative ring $A$. Just a set of points remembers nothing
of the algebraic structure; however, taken together, these sets form
a functor of points, and it remembers all there is to know about
$X$. Formally, we have a fully faithful embedding from the category
of schemes to the category of contravariant functors from rings to
sets. If needed, one can characterize the essential image of this
embedding, or in fact enlarge it (e.g.\ by considering algebraic
spaces).

Analogously, no matter what one undestands by an ``enhanced
category'' $\C$, at the very least, it has its ``truncation'', a
usual category $\tau(\C)$ with $\Hom$-sets obtained by taking
$\pi_0$ of the corresponding homotopy types of morphisms. However,
there is more. Any usual small category $I$ can be thought of an
enhanced category, with trivial enhancement --- that is, discrete
homotopy types of morphisms. Moreover, for any enhanced category
$\C$, we should have an enhanced category $\Fun(I,\C)$ of functors
from $I$ to $\C$. Combining the two observations, we obtain a
collection of categories $\tau(\Fun(I,\C))$ indexed by all small
categories $I$. Moreover, for any functor $\gamma:I_0 \to I_1$, we
have a pullback functor $\gamma^*:\tau(\Fun(I_1,\C)) \to
\tau(\Fun(I_0,\C))$, and there are isomorphisms between these
functors for composable pairs. It is natural to ask, then, whether
the family of the unenhanced categories $\tau(\Fun(I,\C))$ and
functors $\gamma^*$ between them can serve as a ``functor of
points'' that allows one to recover the enhanced category $\C$.

Grothendieck's name for such a collection is a ``derivator'',
probably because at the time of \cite{champ} one was mostly
interested in enhancements for derived categories; however, the idea
is completely general. It has been pursued to some extent, see
e.g.\ \cite{grotz} for a very good and relatively recent
overview. However, it seems that the story has never been taken to
its logical conclusion. The goal of the present text is to provide
such a conclusion, of sorts --- or at least, to show that with some
modifications, Grothendieck's idea leads to an axiomatic formulation
of abstract homotopy theory that is just as strong as those present
in the literature, but free from the problems described above.

\section*{A trailer.}

Let us describe very briefly what we do (there are also more
detailed introductions at the beginning of each chapter, and a
general overview available as a companion paper \cite{ka.umn}). So
as not to usurp a term with a well-established meaning, we avoid
using the word ``derivator'' and prefer to call our gadgets simply
``enhanced categories''. The modifications we make to the original
suggestion of \cite{champ} are the following two.
\begin{enumerate}
\item It seems that it is not necessary, and in fact not desirable,
  to index our families on all small categories; it is better to
  restrict to the case when $I$ is a partially ordered set. This is
  enough because, roughly speaking, any enhanced category can be
  obtained as an enhanced localization of a partially ordered set,
  and it works better because unlike $\Cat$, the category $\Pos$ of
  partially ordered sets has no $2$-categorical structure, and does
  not need to be localized (a map between partially ordered sets is
  an equivalence iff it is an isomorphism).
\item One needs a convenient packaging of the notion of a ``family
  of categories''. In \cite{champ}, Grothendieck makes a simplifying
  assumption that the family is ``strict'' --- for any two functors
  $\gamma_0:I_0 \to I_1$, $\gamma_1:I_1 \to I_2$, we have $(\gamma_1
  \circ \gamma_0)^* = \gamma_0^* \circ \gamma_1^*$. This has been
  adopted in the subsequent work. We prefer to get rid of this
  simplifying assumption --- among other things, it is probably
  prudent to assume that functors simply cannot be equal --- and
  package the whole thing using an earlier discovery of Grothendieck,
  the ``Grothendieck construction'' of \cite{SGA} (with a slight
  modification, see Subsection~\ref{fib.subs} below). The data
  defining an enhanced category in our sense is then a category $\C$
  equipped with a Grothendieck fibration $\C \to \Pos$, with fibers
  $\C_J$, $J \in \Pos$, and transition functors $f^*:\C_{J'} \to
  \C_J$ for any map $f:J \to J'$ in $\Pos$. An enhanced functor is a
  functor cartesian over $\Pos$.
\end{enumerate}
Just as in the standard approach based on \cite{qui.ho}, our
enhanced categories are considered up to an equivalence --- this is
unavoidable already for the $1$-truncated homotopy type represented
by usual unenhanced groupoids. However, the equivalence in question
is simply an equivalence of categories, and just as for $\Cat$,
localization with respect to such equivalences does not require
model category techniques. Just to be clear (see
Subsection~\ref{fun.subs}), if we have categories $I$, $\C_0$,
$\C_1$ and functors $\pi_0:\C_0 \to I$, $\pi_1:\C_1 \to I$, then a
functor $\C_0 \to \C_1$ over $I$ means a pair $\langle \gamma,\alpha
\rangle$ of a functor $\gamma:\C_0 \to \C_1$ and an isomorphism
$\alpha:\pi_0 \cong \pi_1 \circ \gamma$, and a morphism $\langle
\gamma,\alpha \rangle \to \langle \gamma',\alpha' \rangle$ over $I$
is a morphism $a:\gamma \to \gamma'$ such that $\pi_1(a) \circ
\alpha = \alpha'$. Then if $I$ is small, small categories $\C$
equipped with a Grothendieck fibration $\C \to I$ form a category
that we denote by $\Cat \mm I$, with morphisms given by cartesian
functors over $I$, and localizing $\Cat \mm I$ with respect to the
class of equivalences produces the category $(\Cat \mm I)^0$ with
the same objects, and morphisms given by cartesian functors over $I$
up to an isomorphism over $I$. In words, the effect of the
localization is that objects stay the same, and morphisms become
isomorphism classes of functors rather than functors themselves. If
$I$ is not small, there may be size problems (there may be too many
functors over $I$). Nevertheless, we define ``enhanced categories''
as Grothendieck fibrations over $\Pos$ satisfying several explicit
axioms, and we prove that for any two such fibrations $\C,\C' \to
\Pos$ with essentially small fibers, isomorphism classes of
cartesian functors $\C \to \C'$ over $\Pos$ form a set and not a
proper class. This gives a well-defined category $\Cat^h$ of ``small
enhanced categories'', with morphisms given by isomorphism classes
of cartesian functors.

The full set of axioms defining an enhanced category is given below
in Section~\ref{enh.sec} (or in \cite[Section 3]{ka.umn}), so let us
just give an overview. Firstly, let $\ppt$ be the one-point set, and
for any element $j \in J$ in a partially ordered set $J$, let
$\eps(j):\ppt \to J$ be the embedding onto $j$. Say that a
Grothendieck fibration $\C \to \Pos$ is {\em non-degenerate} if for
any $J \in \Pos$, the transition functor $\prod_{j \in
  J}\eps(j)^*:\C_J \to \prod_{j \in J}\C_\ppt$ is
conservative. Moreover, say that $\C$ is {\em additive} if for any
disjoint union $J = \copr_{s \in S}J_s$ of some $J_s \in \Pos$
indexed by a set $S$, with embedding maps $\eps(s):J_s \to J$, the
transition functor $\prod_{s \in S}\eps(s)^*:\C_J \to \prod_{s \in
  S}\C_{J_s}$ is an equivalence.  Next, consider the partially
ordered set $[1] = \{0,1\}$ corresponding to the single arrow
category (the order is $0 \leq 1$), and for any $J \in \Pos$, let
$s(J),t(J):J \to J \times [1]$ be the embeddings onto $0 \times J$
resp.\ $1 \times J$, and let $e(J):J \times [1] \to J$ be the
projection. Say that a Grothendieck fibration $\C \to \Pos$ is a
{\em reflexive family} if for any $J$, the isomorphism $t(J)^* \circ
e(J)^* \cong \id$ defines an adjunction between the transition
functors $t(J)^*$ and $e(J)^*$ (with $t(J)^*$ being right-adjoint to
$e(J)^*$). Moreover, for any such reflexive family, applying
$s(J)^*$ to the adjunction map $e(J)^* \circ t(J)^* \to \id$
produces a functor $\nu_J:\C_{J \times [1]} \to \Fun([1],\C_J)$ to
the category of functors $[1] \to \C_J$. Say that a functor is an
{\em epivalence} if it is conservative, essentially surjective and
full --- that is, any object in its target lifts to an object in its
source, uniquely up to a non-unique isomorphism --- and say that a
reflexive family $\C \to \Pos$ is {\em separated} if the
corresponding functor $\nu_J$ is an epivalence for any $J \in \Pos$.

To describe the remaining axioms, note that a commutative square in
a category $I$ is the same thing as a functor $\eps:[1]^2 \to I$,
and a commutative square \thetag{$\star$} is the same thing as a
Grothendieck fibration over $[1]^2$. Say that a square \thetag{$\star$} is
{\em semicartesian} if the corresponding functor $\C_{01} \to \C_0
\times_{\C} \C_1$ is an epivalence, and say that a fibration $\C \to
I$ is {\em semicartesian} over a commutative square $\eps:[1]^2 \to
I$ if so is the induced square $\eps^*\C \to [1]^2$. Then an
enhanced category in our sense is a non-degenerate additive
separated reflexive family $\C \to \Pos$ that satisfies five
additional axioms (``semiexactness'', ``excision'',
``semicontinuity'', the ``cylinder axiom'' and ``bar-invariance'',
see Section~\ref{enh.sec} or \cite[Section 3]{ka.umn}), all of them
requiring that $\C$ is semicartesian over commutative squares in
$\Pos$ of some prescribed form. An enhanced category $\C$ is {\em
  small} if all the fibers $\C_J$ are essentially small.

Let us remark right away that the precise shape of our definition of
an enhanced category is not set in stone. For example, in this
introduction, we consider fibrations over the whole $\Pos$, mostly
because it gives the simplest definition of the opposite enhanced
category --- namely, for any fibration $\C \to I$, we also have a
fibration $\C^o_\perp \to I$ whose fibers and transition functors
are opposite to those of $\C$, and then the {\em enhanced-opposite}
$\C^\iota$ to an enhanced category $\C$ is given by $\C^\iota =
\iota^*\C^o_\perp$, where $\iota:\Pos \to \Pos$ sends a partially
ordered set $J$ to the same set with the opposite order, denoted
$J^o$. However, it suffices to work over the full subcategory
$\Posf^\pm \subset \Pos$ of {\em left-finite sets} --- that is, sets
$J \in \Pos$ such that for any element $j \in J$, the subset $J/j =
\{j' \in J|j' \leq j\} \subset J$ is finite. A suitable fibration
over $\Posf^\pm$ then canonically and uniquely extends to the whole
$\Pos$ (this is Proposition~\ref{lf.prop}). Analogously, one can
vary the defining commutative squares in $\Pos$ --- at the end of
the day, this amounts to choosing a set of generating relations in
an algebraic structure, and it can be done in different ways without
changing the structure.

In practice, one does not expect to ever have to actually check the
axioms. Instead, new enhanced categories and functors are produced
from existing ones by a couple of standard procedures, starting from
basic examples that can be used as black box. Here are these
examples (for more details, see Subsection~\ref{enh.exa.subs} or
\cite[Subsection 3.4]{ka.umn}).
\begin{enumerate}
\item For any category $I$, we have an enhanced category $\Unf(I)$
  with fibers $\Unf(I)_J = \Fun(J^o,I)$. This is a ``trivial''
  example where all the epivalences in the definition are
  equivalences. We call such enhanced categories {\em
    $1$-truncated}. In fact, an enhanced category $\C \to \Pos$ can
  be thought of as an enhancement for its fiber $\C_\ppt$, and then
  $\Unf$ is right-adjoint to the truncation operation $\C \mapsto
  \C_\ppt$ (see Proposition~\ref{trunc.prop} for a precise
  statement).
\item For any model category $\C$ with a class $W$ of weak
  equivalences, we have the ``enhanced localization'' $\Hh^W(\C) \to
  \Pos$ whose fibers are given by $\Hh^W(\C)_J \cong
  h^W(\Fun(J^o,\C))$.
\item For any additive category $\A$, we have the ``enhanced
  homotopy category'' $\Hho(\A) \to \Pos$ giving an enhancement for
  $\Ho(\A)$ (explicitly, the fiber $\Hho(\A)_J$ is the
  chain-homotopy category $\Ho(\Fun(J^o,\A))$ further localized with
  respect to maps $f$ such that $\eps(j)^*(f)$ is a chain-homotopy
  equivalence for any $j \in J$).
\item Finally, for any complete Segal space $X$, we have an enhanced
  category $\Unf(X) \to \Pos$ whose fiber $\Unf(X)_J$ is the
  truncation of the complete Segal space of functors $J^o \to X$.
\end{enumerate}
We then prove a version of Brown Representability Theorem (namely,
Theorem~\ref{enh.thm}) saying that for small enhanced categories,
the last example is universal: if we let $\CSS$ be the homotopy
category of complete Segal spaces, then $\Unf:\CSS \to \Cat^h$ is an
equivalence.

It is useful to explain how our formalism accomodates homotopy
types. There appear as {\em enhanced groupoids} --- enhanced
categories $\C \to \Pos$ whose fibers $\C_J$ are groupoids. In this
situation, some of the axioms hold automatically, so the definition
is simpler (and we do it separately, in Section~\ref{grp.sec}, as
warm-up for the general theory). The resulting gadget is probably as
close to an ``enhanced Galois theory'' as one can get. We do not
have ``higher groupoids'' with ``symmetries between symmetries'',
whatever it means. We just have groupoids in the perfectly ordinary
sense of the word, but a whole bunch of them. The truncation
$\C_\ppt$ is the fundamental groupoid $\pi_{\leq 1}(X)$ of our
homotopy type $X$; but for any partially ordered set $J$, we also
have the groupoid $\C_J$, the fundamental groupoid of the space of
maps from the nerve of $J$ to $X$. If one want to look at higher
homotopy groups $\pi_n(X)$, one takes any $J \in \Pos$ whose nerve
is an $n$-sphere (for example, the set of non-degenerate simplices
of the standard simplicial $n$-sphere). The whole structure is of
course richer; in particular, $J$ need not be finite.

\begin{remarknonum}
We do have to allow infinite $J$. The Brown Representability Theorem
that we generalize is the original theorem of \cite{brown};
unfortunately, the trick of Adams \cite{A} that reduces everything
to finite simplicial sets does not seem to generalize. It is
probably simply not true that a small enhanced groupoid, let alone a
general enhanced category is completely defined by its restriction
to the subcategory of finite partially ordered sets.
\end{remarknonum}

In general, the representability theorem is important for two
reasons. Conceptually, it justifies our claim that our formulation
of abstract homotopy theory is at least as strong as the one based
on complete Segal spaces (and {\em a posteriori}, as all the other
approaches in the paradigm of \cite{qui.ho}). For example, there
used to be some controversy as to whether ``a derivator is enough to
recover algebraic $K$-theory'', see a brief recapitulation at the
end of this Introduction; in our approach, the point it
moot. Technically, while the representing object is only unique up
to a weak equivalence and not explicit at all, its mere existence
already has useful corollaries. One such is that the category
$\Cat^h$ is cartesian-closed --- and more generally, for any
enhanced categories $\C$, $\E$ with $\C$ small, there exists an
enhanced functor category $\fFun^h(\C,\E)$ with the usual universal
property (Definition~\ref{enh.fun.def} and
Corollary~\ref{J.C.corr}). In fact, almost all of the enhanced
category theory we develop in Chapter~\ref{enh.ch} (or \cite[Section
  4]{ka.umn}) is done by standard categorical arguments, with
representability used purely as a black box. In most arguments,
Theorem~\ref{enh.thm} is not even used directly but enters through
one of its corollaries, Proposition~\ref{univ.prop}.

Some of the enhanced categorical notions do not need even that, and
carry over {\em verbatim} from the usual category theory. One such
is the notion of a fully faithful functor, see
Corollary~\ref{enh.ff.corr}, while another is that of an adjoint
functor, Definition~\ref{enh.refl.def} --- an enhanced functor
$\gamma$ is right or left-adjoint to an enhanced functor
$\gamma^\dg$ iff it is right resp.\ left-adjoint in the unenhanced
sense, and the adjunction maps are maps over $\Pos$. This is useful
because uniqueness of adjoints, one of the most fundamental facts of
category theory, also carries over without any changes at
all. Another thing that comes almost for free is an enhancement
$\cCat^h$ for the category $\Cat^h$, and an enhanced version of the
equivalence $\CSS \cong \Cat^h$ (Subsection~\ref{cath.subs}). It is
natural to expect that the fibers $\cCat^h_J$ of such an enhancement
would be the categories of small enhanced categories $\C$ equipped
with an enhanced functor $\gamma:\C \to \Unf(J)$ that is a fibration
in some enhanced sense. In reality, no enhancement is needed ---
$\gamma$ should simply be a Grothendieck fibration.

Representability, in the form of Proposition~\ref{univ.prop},
reappears in the construction of enhanced-cartesian squares. Just as
in the case of $\Cat^0$, these cannot be defined purely in terms of
$\Cat^h$; we really need not only small enhanced categories $\C$,
$\C_0$, $\C_1$ and morphisms $\C_l \to \C$, $l=0,1$ in $\Cat^h$, but
also enhanced functors $\gamma_l:\C_l \to \C$, $l=0,1$ representing
the morphisms. Under special curcumstances --- e.g.\ if $\gamma_0$
is fully faithful, or if $\C$ is $1$-truncated --- the product $\C_0
\times_{\C} \C_1$ is an enhanced category. In general, it is not:
fibered products of semicartesian squares are not necessarily
semicartesian. However, this can be corrected. Firstly --- this is
Lemma~\ref{sq.h.le} --- there is an enhanced category $\C_0
\times^h_{\C} \C_1$ equipped with an epivalence $\gamma:\C_0
\times^h_{\C} \C_1 \to \C_0 \times_{\C} \C_1$, cartesian over $\Pos$
(in other words, we have a semicartesian square \thetag{$\star$} of
enhanced categories and enhanced functors). Secondly --- this is
Corollary~\ref{sq.h.corr} --- for any such epivalence $\gamma$, and
for any other enhanced category $\C_{01}$, any functor $\C_{01} \to
\C_0 \times_{\C} \C_1$ cartesian over $\Pos$ factors through
$\gamma$, uniquely up to an isomorphism.

Unlike the unenhanced universal property for $\C_0 \times_{\C}
\C_1$, the isomorphism here is {\em not} unique, and $\gamma$ is
only an epivalence, not an equivalence --- this is just as it should
be, since it is here that all the higher homotopy groups
appear. However, the universal property is still sufficient to
recover the semicartesian product $\C_0 \times^h_{\C} \C_1$ up to a
unique equivalence, that is, as an object in $\Cat^h$. The
comparison $\CSS \cong \Cat^h$ then identifies it with the homotopy
fibered product in $\CSS$. Moreover, it turns out to be sufficiently
functorial to allow for most of the usual categorical constructions,
see Subsection~\ref{univ.funct.subs}, and this is one of the main
tools used to construct new enhanced categories and functors.

Once we have enhanced functor categories and semicartesian products,
the rest of the theory hardly needs representability at all. We
define enhanced comma-categories in Subsection~\ref{enh.co.subs}.
This leads to the notions of an enhanced fibration $\C \to \E$, see
Subsection~\ref{enh.fib.subs} (if the base $\E$ of the fibration is
$1$-truncated, it reduces to a Grothendieck fibration in the usual
sense, see Lemma~\ref{noenh.fib.le}). An enhanced version of the
Grothendieck construction (Subsection~\ref{enh.groth.subs})
identifies small enhanced fibrations $\C \to \E$ and enhanced
functors $\E^\iota \to \cCat^h$, and this leads to a fully faithful
enhanced Yoneda embedding (Subsection~\ref{enh.yo.subs}). We then
have enhanced limits, colimits and Kan extensions
(Section~\ref{enh.lim.sec}), enhanced-complete and cocomplete
categories, and so on and so forth. For further details, see
\cite[Section 4]{ka.umn}.

One last thing that needs mentioning is what happens with large
categories. Already in the unenhanced setting, these present a
problem: since $\Hom$-sets need to be small, functors between two
large categories need not form a category. In the enhanced setting,
the problem becomes worse: already the semicartesian products $-
\times^h -$ of Lemma~\ref{sq.h.le} need not exist.

In principle, this could be handled by introducing the Universe
Axiom and enlarging the universe. In our approach, this would also
require enlarging the domain of definition of our enhanced
categories --- these need to become fibrations defined over the
category of partially ordered sets in the new universe --- but
theoretically, it probably could be done. However, a better way is
to use the machinery of filtered colimits and accessible categories
of Grothendieck \cite{sga4} and Gabriel and Ulmer \cite{GU}. In
effect, one restricts one's attention to categories and functors
that are $\kappa$-accessible for some regular cardinal $\kappa$ that
is large enough (and can be enlarged if needed). This has one
standard problem --- the opposite to an accessible category is
typically not accessible, so that simply taking the opposite
category is not a valid operation in the accessible world --- but
that's life. Modulo that, most of the large categories one needs in
practice are accessible, and treating them as such more-or-less
solves the problem.

In our enhanced setting, we develop the machinery of filtered
colimits and accessible categories in Section~\ref{enh.large.sec},
and it parallels very closely the unenhanced theory (for comparison,
see e.g.\ \cite{ka.acc}, an overview of the classical theory written
as an unenhanced counterpart to Section~\ref{enh.large.sec}). The
same caveats apply, but to a large extent, the technology works, and
most of our enhanced categorical constructions extend to the
accessible setting.

\section*{Sociology, or why we do it.}

Let me now address the elephant in the room. For almost two decades
now, since \cite{lu1}, Jacob Lurie and his school have been
developing a set of very comprehensive and thorough foundations for
abstract homotopy theory colloquially known as the theory of
``$\infty$-categories''. At the moment, the theory is infinitely
more advanced than whatever was the state of the art at the time of,
say, \cite{rzk}, and is actively and productively used by many
mathematicians around the world. Hundreds if not thousands of
graduate students invested uncounted hours of their valuable time
into studying the texts. What's wrong with Lurie's foundations? ---
how could one possibly justify opening the can of worms yet again,
and looking for another approach?

The short answer is, quite possibly, nothing. The feeling of deep
unease experienced by some when encountering the
$\infty$-categorical formalism might just be a personal problem, or
a mental infirmity of some sort --- or just old age, for that
matter. And it is not that important anyway, since the vast majority
of Lurie's work is actually completely model-independent, at least
in spirit, and would survive a change of engine with very little
pain.

Still, the feeling of unease is there, and it is hard to ignore
it. In effect, at the time of writing, a simple browsing of recent
papers that could have used $\infty$-categories reveals a deep
split. Some papers --- or maybe, papers in some areas --- do use
them as a completely standard tool, without even discussing
this. Other papers, in areas that are mathematically very close,
treat the formalism of $\infty$-categories as non-existent. Of
course, the formal reason for this might be simply that most of the
foundational papers of Lurie are not published in a conventional
way, even as preprints, but there is probably more to it that just
that.

The situation was already the same when I was writing
\cite{ka.glue}. At the time, I was trying to rationalize this
unease, at least to myself, and the explanation that I came up with
was roughly the following. The way Lurie's foundations are set up is
pretty much all-or-nothing: either you move the whole of mathematics
to an $\infty$-categorical setting, or you cannot use it at
all. This even extends to some gratuitous rebranding --- why on
earth would someone swap the meaning of ``final'' and ``cofinal'',
for example? --- but actually goes deeper than such minor quirks. By
its very nature, $\infty$-categorical formalism does not allow for a
clean separation of definitions (that should be simple) and proofs
(that can be as difficult as needed but can be safely used as black
box). Well, the definition of a quasicategory is simple enough, but
already to understand adjunction, you really need quite a bit of
simplicial combinatorics, so that already a definition needs to be
put in a black box. In practice, to write an honest paper using the
formalism, one really needs to study thousands pages of text, and
use precise page references to various relevant facts spread out
over all these pages (one example of such an honest paper that comes
to mind now is \cite{NS}). Thus both the specific model used by
Lurie, and the whole foundational paradigm of \cite{qui.ho} are
hard-coded into the formalism as it exists. Forcing all this unto
the whole of mathematics is bound to cause some resentment.

Coming back to the present, the above paragraph still stands, to
some extent, and maybe even becomes stronger --- as the current text
purports to show, both the simplicial machinery and the model
category paradigm are not only cumbersome but also not really
related to the actual mathematical content of the theory. However,
when looking at how the use of $\infty$-categories has developed, I
can now see another reason for worry: the way things are, it is
simply not safe.

\medskip

The thing is, what people really want is to work with enhanced
categories, and more importantly, think about them as if they were
unenhanced. It is well-known and accepted that this is not possible
all the time, but it should be possible most of the time, with some
well-defined points where one just needs to ``say something'' to
make it all work.

\medskip

By itself, this is not a bad thing. The ability to use something
like Grothendieck's six-functor formalism in a general homotopical
setting is deeply liberating, and in XXI century, one shouldn't have
to worry about choosing a particular resolution, be it additive or
non-additive. However, this creates ample space for errors. The
common sloppiness actually goes back to working with usual
unenhanced categories, when people write down commutative squares
\thetag{$\star$}, or more complicated commutative diagrams, without
mentioning the connecting isomorphism or isomorphisms. It is
implicitly assumed, but rarely mentioned, that if the constructions
are ``natural enough'', these things would take care of
themselves. The situation becomes more problematic when we move to
monoidal categories and functors, $2$-categories and suchlike. Here
it is almost impossible to actually write down {\em all} the higher
constraints that need to be in the picture, for typographical
reasons if nothing else, so usually there is some sort of compromise
(or recourse to a strange and unnatural gadget such as Mac Lane
Strictification Theorem of \cite{mcl}). In the $\infty$-categorical
setting, the problems become infinitely worse. Quite often, people
still think $1$-categorically, with the corresponding
$1$-categorical intuition (defining functors by what they do to
objects and morphisms, using monads and the ``Barr-Beck Theorem'',
and so on). It is still assumed that if push comes to shove,
everything can be ``strictified'', or made to ``commute on the
nose'' by a good choice of a model, so that the whole house of cards
does not collapse (and in the $\infty$-categorical setting, the
number of those cards tends to be infinite). Unfortunately, the
formalism itself does not help --- occasionally, to make things
rigourous, you do need to make a choice of a model. There goes your
liberty.

As a baby example of this, let us consider a very common categorical
procedure --- defining things by pullback (for example, taking
fibers).

On the level of homotopy categories such as $\CSS$ or our $\Cat^h$,
pullbacks are of course not functorial and do not have a universal
property. So, if you have a morphism $\C_1 \to \C$ and some other
morphism $\C_0 \to \C$, you cannot just define the homotopy fibered
product $\C_0 \times^h_{\C} \C_1$ while staying in $\CSS \cong
\Cat^h$. What you do in $\CSS$ or in $\Cat^h$ is similar but
slightly different.
\begin{enumerate}
\item In $\CSS$ --- or in the category of quasicategories, if you
  wish --- you choose complete Segal spaces $X$, $X_0$, $X_1$ and
  maps $X_l \to X$, $l=0,1$ representing $\C$, $\C_0$, $\C_1$ and
  the morphisms. Moreover, you arrange, by possibly changing $X_0$,
  that the map $X_0 \to X$ is a fibration in the model category
  sense. Then you take the strict fiber product $X_0 \times_X X_1$,
  and let $\C_0 \times^h \C_1$ be the corresponding object in
  $\CSS$. If you then have a commutative square \thetag{$\star$} in
  $\CSS$, you also choose a complete Segal space $X_{01}$
  representing $\C_{01}$, and a homotopy between the two maps
  $X_{01} \to X$. Then you use a simple model category argument to
  possibly change $X_{01}$ so that the homotopy becomes an identity,
  and obtain the map $X_{01} \to X_0 \times_X X_1$ by the universal
  property of the cartesian square. Then you descend back to $\CSS$.

  As an alternative, you could say that $\CSS$ is complete ``as an
  $\infty$-category'', and then the homotopy fiber product becomes
  functorial and acquires a universal property. However, then your
  input data --- namely, the diagram $\C_0 \to \C \gets \C_1$ --- needs
  to be lifted to an $\infty$-diagram, and this amounts to exactly
  the same thing as above.
\item In $\Cat^h$, you lift $\C$, $\C_0$, $\C_1$ to enhanced
  categories, and lift your morphisms to enhanced functors. However,
  instead of simply taking the cartesian product, you take the
  semicartesian product $\C_0 \times^h_{\C} \C_1$ provided by a
  machine hidden in a black box (namely, Lemma~\ref{sq.h.le}). The
  universal property and the comparison functor $\C_{01} \to \C_0
  \times^h_{\C} \C_1$ is then provided by another machine,
  Corollary~\ref{sq.h.corr}, once you choose an isomorphism between
  the two enhanced functors $\C_{01} \to \C$.
\end{enumerate}
The difference is somewhat subtle, and it is this. In our approach,
you do not strictify all the way, and {\em do not} obtain the
required map $\gamma:\C_{01} \to \C_0 \times^h_{\C} \C_1$ ``on the
nose''. In fact, the formalism does not even allow you to do
this. This is not a bug but a feature. You only get $\gamma$ up to
an isomorphism, and this is exactly what can be done in an invariant
way. The ability to strictify completely just means that you kept
too much information: you are secretly still remembering those
explicit resolutions from which you wanted to be finally free.

In other words, our formalism is designed to imitate the usual
category theory much closer than anything based on the model
category techniques, but with a caveat: in some well-defined places,
you need to explicitly invoke something from a black box. A category
of models where things just commute, period, simply does not
exist. This might be cumbersome, from time to time, but that is how
life is; thinking otherwise would be a mistake. Our formalism is
designed to make such mistakes impossible.

\section*{Safety features presentation.}

Let us now list some other things that in our approach, you simply
cannot do, by design.
\begin{enumerate}
\item {\em Nothing is ever defined by hand.}
\end{enumerate}
In ordinary category theory, you can in principle define a category
by saying what are its objects, morphisms and compositions, and you
can define a functor by saying what it does to objects and
morphisms. In the enhanced setting, this is not possible, because it
requires an infinite amount of data. In fact, this is not
practically possible already in the formalisms based on models, or
more generally, wherever one needs some kind of ``homotopy
coherence''. For a classical example, while it is true that an
infinite loop space is an algebra over an $E_\infty$-operad, for any
particular fixed choice of such an operad, all the infinite loop
spaces in existence were actually constructed by the Segal
machine. In our enhancement formalism, constructing things by hand
is impossible even theoretically.

There is one exception to this rule. We say that {\em enhanced
  objects} in an enhanced category $\C \to \Pos$ are objects in its
truncation $\C_\ppt$ --- or equivalently, enhanced functors $\ppt^h
\to \C$ from the terminal enhanced category $\ppt^h = \Unf(\ppt)
\cong \Pos$. Enhanced morphisms between enhanced objects are
morphisms in $\C_\ppt$, or equivalently, enhanced functors
$\Unf([1]) \to \C$. However, one can also consider ``generalized
objects'' --- namely, all objects in $\C$, in all the fibers $\C_J
\subset \C$, $J \in \Pos$. Then we do have a general reconstruction
result, Lemma~\ref{enh.exp.le}, saying that to define an enhanced
functor $\C \to \C'$ between enhanced categories, it suffices to say
what it does do all generalized objects $c \in \C$, functorially
with respect to isomorphisms, and functorially with respect to the
transition functors $f^*:\C_{J'} \to \C_J$, for all maps $f:J \to
J'$. You do not have to worry about generalized morphisms
(generalized morphisms in $\C_J \subset \C$ are generalized objects
in $\C_{J \times [1]}$).
\begin{enumerate}
\stepcounter{enumi}
\item {\em Nothing is ever equal, and nothing commutes ``on the
  nose''.}
\end{enumerate}
Objects in a ordinary category can only be isomorphic, but morphisms
between two fixed objects can be equal and form a set. Categories
can only be equivalent, functors can only be isomorphic, but again,
morphisms between functors can be equal. In our enhanced formalism,
the last step disappears. Enhanced morphisms between two enhanced
objects form an enhanced groupoid, not a set, so just as objects,
they can only be isomorphic. Enhanced functors should be treated as
enhanced objects in the functor category, so morphisms between also
form an enhanced groupoid and not a set. Technically, our enhanced
functors are also functors in the usual sense, we do have a
``naive'' set of morphisms, so two morphisms can be equal, but this
is an artefact of the formalism and has no invariant meaning; in
particular, such equalities are not preserved by semicartesian
products. Saying that a morphism is invertible is safe; saying that
two morphisms are equal is not. In fact --- see
Remark~\ref{fun.C.rem} --- if one considers enhanced functors over
some enhanced category $\C$, then the ``naive'' equality is simply
the wrong thing, already on the level of $\pi_0$.

We note that this phenomenon causes surprisingly few problems when
building the theory. In fact, the only place in the standard
category theory where an actual equality of morphisms is used is in
the definition of adjunction (see \eqref{adj.eq} in
Subsection~\ref{adj.subs}). Since for us, enhanced adjunction is the
same thing as an unenhanced one, this could cause problems; however,
this can be handled by an alternative description of adjunction in
terms of enhanced comma-categories, Lemma~\ref{adj.eq.le}.
\begin{enumerate}
\stepcounter{enumi}\stepcounter{enumi}
\item {\em All the commutative diagrams have to be enhanced.}
\end{enumerate}
In a sense, this is a continuation of \thetag{ii}, but we cannot
emphasize this point strongly enough since it is the single major
source of gaps and outright mistakes (especially when constructing
enhanced limits or colimits). If you need a commutative diagram in
an enhanced category $\C$, you should really number its vertices and
arrows, turn it into a partially ordered set $J$, or a category $I$,
or even an enhanced category $\E$ if you so wish, and construct an
enhanced functor $\E \to \C$. If $\C = \cCat^h$ --- that is, if you
want a commutative diagram of enhanced categoies and enhanced
functors --- then you can apply the enhanced Grothendieck
construction and consider a fibration $\E' \to \E^\iota$ rather than
an enhanced functor $\E \to \cCat^h$. If $\E = \Unf(J)$, then this
is the same thing as usual Grothendieck fibration $\E' \to
\Unf(J^o)$, but even this is much more than a fibration $\E' \to
J^o$. If $J$ is very simple, then enhancement can be automatic, to
some extent --- see Corollary~\ref{ccat.dim1.corr} or
\cite[Subsection 4.5]{ka.umn} --- and for commutative squares, it is
enough to choose a connecting isomorphism (in particular, enhanced
fibered products in $\cCat^h$ are given by semicartesian squares,
see Lemma~\ref{ccat.sq.le}). Anything more complicated has to be
split into parts or described explicitly via the Grothendieck
construction. As a rule of thumb, if your commutative diagram of
enhanced categories and functors needs something more than the
standard {\tt amscd} package, you are doing something wrong.

\medskip

All of the above begs the question: how {\em do} you define things?
What are safe procedures for constructing new enhanced categories
and functors from those that already exist?

\medskip

For enhanced functors, the most obvious procedure that turns out to
be quite powerful is taking an adjoint. As we have mentioned,
uniqueness is simply inherited from the usual category theory, and
for existence, there are various criteria (for example,
Lemma~\ref{enh.fib.adj.le} is rather useful). For enhanced
categories, one thing you can do to take a functor category
$\fFun^h(\C,\E)$. For both categories and functors, there are also
semicartesian products $- \times^h -$. In effect, for any small
enhanced category $\E$, one can define an enhanced category $\cCat^h
\bbi^h \E$ of small enhanced categories and functors over $\E$, see
Subsection~\ref{lax.subs}, and then for any enhanced functors
$\gamma:\E \to \E'$ between small enhanced categories, we have an
enhanced functor
\begin{equation}
\renewcommand{\theequation}{$\star\star$}
\gamma^*:\cCat^h \bbi^h \E' \to \cCat^h \bbi^h \E
\end{equation}
sending an enhanced object $\C \to \E'$ to $\gamma^*\C = \E
\times^h_{\E'} \C \to \E$. The formal construction of
\thetag{$\star\star$} is in Subsection~\ref{enh.rel.subs}, where it
is also observed that it has an obvious left-adjoint $\gamma_\trr$
sending $\C \to \E$ to $\C \to \E \to \E'$. Thus semicartesian
products are sufficiently functorial to define a pullback operation
$\gamma^*$ both on enhanced categories and on enhanced functors. On
categories, $\gamma^*$ sends enhanced fibrations to enhanced
fibrations (Lemma~\ref{enh.fib.sq.le}) and on functors, it sends
adjoint pairs to adjoint pairs (Corollary~\ref{enh.adj.corr}) and
cartesian functors to cartesian functors (again
Lemma~\ref{enh.fib.sq.le}).

Moreover, it is shown in Subsection~\ref{enh.rel.subs} that under
favourable circumstances --- for example, when $\E$ and $\E'$ are
small, and $\gamma$ or $\gamma^\iota$ is an enhanced fibration ---
the enhanced functor \thetag{$\star\star$} also has a right-adjoint
functor $\gamma_\trl:\cCat^h \bbi^h \E \to \cCat \bbi^h \E'$. The
{\em relative enhanced functor category} $\fFun^h(\E|\E',\C)$ is
defined as $\gamma_\trr\gamma^*(\C \times^h \E')$, and this
construction also turns out to be very useful. If $\E' = \ppt^h$,
this recovers the usual enhanced functor category $\fFun^h(\E,\C)$.

There are also two constructions that are specific to enhanced
fibrations. Any enhanced fibration $\C \to \E$ has enhanced fibers
$\C_e$ over enhanced objects $e \in \E_\ppt$, and enhanced
transition functors $f^*:\C_{e'} \to \C_e$ for enhanced morphisms
$f:e \to e'$ (Subsection~\ref{enh.fiber.subs}). Then firstly, just
as in the unenhanced setting, there is a general procedure that
produces a fibration $\C^\iota_{h\perp} \to \E$ with enhanced fibers
$\C_e^\iota$ and transition functors $(f^*)^\iota$. Moreover,
$\C^\iota_{h\perp}$ can be characterized by a universal property,
see Corollary~\ref{enh.perp.corr}. Secondly, if all the transition
functors $f^*$ admit enhanced left-adjoints, then the opposite
enhanced functor $\C^\iota \to \E^\iota$ is also an enhanced
fibration (Corollary~\ref{enh.bifib.corr}). Iterating the two
constructions, and throwing in some pullbacks for a good measure,
can actually get you quite far. For an example in an unenhanced
setting, see e.g.\ the construction of the so-called ``Fourier-Mukai
$2$-category'' given in \cite[Subsection 4.4]{ka.adj}.

Finally, let us mention that the enhanced category $\cCat^h$ is
enhanced-complete and enhanced-cocomplete, and in particular, it has
enhanced colimits and enhanced-cocartesian squares. A formal
corollary of this is that for small enhanced categories, we always
have enhanced localization given by an explicit enhanced-cocartesian
square (namely, \cite[(4.38)]{ka.umn}). Unfortunately, just as in
the unenhanced setting, the resulting enhanced category is usually
very hard to describe explicitly, but at least, it definitely
exists.

\section*{Leitfaden, or how we do it.}

Now, having explained both what we do in the text and why we think
it is worth doing, let us explain how we do it --- how the text is
organized, what depends on what, etc.

\medskip

First of all, while the formalism we build here does not depend on
model categories, nor on simplicial sets, and relegates both to the
realm of technical gadgets, it does not mean that we should not use
these gadgets while building the formalism. Quite the contrary in
fact, both are really helpful. For model categories, we stay within
the framework of \cite{qui.ho} amplified by the Reedy model
structures of \cite{reedy}, and avoid embellishments added in
\cite{hovey} and subsequent work (these are geared towards
applications that are foundational rather than technical, and do not
work well in our context). For simplicial sets, we do use standard
things such as a description of the Segal condition in terms of the
horn embeddings; one slightly less standard device we also found
very useful are simplicial replacements of \cite{bou3} (and we
develop an enhanced version of them in
Subsection~\ref{enh.repl.subs}).

\medskip

Next, we should say how we deal with the three main problems of
category theory: atrocious notation, idiosyncratic terminology, and
inherent triviality of the arguments.

\medskip

For the latter, well, we prove what needs to be proved, but we also
need many small things that are too trivial to split into a
statement and a proof; we call these ``examples''. Our convention is
that such ``examples'' are an integral part of the text, and can be
referred to further on (unlike ``remarks'' that can be safely
skipped).

For notation, the best we can do is to choose one set of atrocities
over another one. Our preferred notation for left resp.\ right
comma-categories for a functor $\C_0 \to \C_1$ is $\C_0 / \C_1$
resp.\ $\C_1 \setminus \C_0$ (and not, say, $\C_0 \uparrow \C_1$ or
$\C_0 \downarrow \C_1$). The opposite category to some $\C$ is
denoted $\C^o$, and same for functors. We denote by $\C^>$
resp.\ $\C^<$ a category $\C$ with an added terminal resp.\ initial
object. The category of functors $I \to \C$ is $\Fun(I,\C)$, and we
shorten $\Fun(I^o,\C)$ to $I^o\C$, as in $\Delta^o\Sets$. This
simplifies formulas quite a bit, and hopefully, does not lead to
confusion. The category of elements of a functor $I^o \to \Sets$ is
$IX$. We have found it very convenient to denote by $\Cat \bb I$ the
category of small categories $\C$ equipped with a functor $\C \to
I$, with morphisms given by lax functors over $I$, and we use a
variation of this notation for various categories of similar
kind. An ``$I$-augmentation'' is a functor to $I$, and
``$I$-augmented'' is ``equipped with a functor to $I$''. The rest is
more-or-less standard.

For terminology, we try to use common conventions whenever possible,
but sometimes one has to choose between two opposite conventions
that are equally common. For us, ``cofinal'' means ``commutes with
colimits'' (not limits), and a ``cone'' --- unless we are talking
about chain complexes --- is something that creates colimits (again,
not limits). A ``fibration'' is a Grothendieck fibration, and a
``cofibration'' is a functor opposite to a fibration --- an
alternative is ``opfibration'', but the sheer sound of it is too
horrible to contemplate (Lurie uses ``cartesian''
resp.\ ``cocartesian fibration'' to avoid confusion with the model
category usage, but we are free from that problem). We use
``epivalence'' in the sense explained above --- conservative, full,
essentially surjective; we do not know how standard it is. We also
do not know how standard are ``cylinders'' of
Subsection~\ref{cyl.subs}. We treat ``anodyne'' of \cite{GZ} as an
orphaned term, and take the liberty to borrow it for something
related but different (the original meaning was ``trivial Kan
cofibration''). ``Left'' resp.\ ``right-reflexive'' means ``admits a
left'' resp.\ ``right-adjoint''. All full subcategories are strictly
full. A full subcategory is ``left'' resp.\ ``right-admissible'' if
the embedding functor is left resp.\ right-reflexive (this is
standard terminology only for people who work with derived
categories of coherent sheaves, but it is convenient enough to merit
more general usage). From time to time, we need an adjective for
something that needs to be named but does not deserve a permanent
name; in these cases, we plagiarize algebraic geometry and assign
words such as ``ample'', ``proper'' or ``separated'', more-or-less
randomly.

We apologize to Graeme Segal for effectively turning his name into
an adjective (and as many people before us, we do it with great
admiration and respect for his work).

\medskip

For better or for worse, the text is self-contained --- the only
things we skip are the definition of a model category, and the proof
of Quillen Adjunction and the existence of Reedy model
structures. Partially because of this, the first time something
technically new appears is in
Section~\ref{reedy.sec}. Chapter~\ref{cat.ch} is completely and
utterly standard, and serves the twin purpose of fixing notation and
terminology, and giving a template for its enhanced version
developed in Chapter~\ref{enh.ch}. As usual, we lament that fact
that there is no convenient textbook (it seems that a tradition
established at least since \cite{mcl} is that all category theory
textbooks omit exactly the things they ought to contain). We split
Chapter~\ref{cat.ch} into four parts, with Section~\ref{fib.sec}
specifically devoted to the Grothendieck construction, and
Section~\ref{fun.fib.sec} containing some additional constructions
that are less well-known (but certainly also known, at least to the
experts). Chapter~\ref{ord.ch} collects everything that we need
about partially ordered sets. This is a lot, but again, nothing is
new (except for possibly some small lemmas that are not in the
literature because nobody ever needed them). Our approach is to
treat partially ordered sets as small categories of a special type,
and this turns out to be more productive than one would
expect. Almost all the partially ordered sets we need are
``left-bounded'', in the sense that $J/j$ has finite chain dimension
for any element $j \in J$, so from the point of view of the general
theory, these are rather trivial. The sets themselves are in
Section~\ref{pos.sec}, and Section~\ref{bi.sec} contains firstly, an
$I$-augmented version of them, for a category $I$, and then a
technical gadget that we found to be very useful --- ``biordered
sets'' (partially ordered sets equipped with a factorization
system). Partially ordered sets appear as examples further on, but
they are only really needed starting from approximately
Chapter~\ref{semi.ch}; biordered sets are not needed before
Section~\ref{seg.sp.sec}.

Chapter~\ref{simp.ch} is also standard, and contains everything we
will need about simplicial things --- horns and spheres, nerves, the
Segal condition, simplicial replacements and so on. The chapter is
rather short since we do not need that much. However, we also
include Section~\ref{reedy.sec} that deals with Reedy categories,
and here there are things that are to some extent new. Namely, we
introduce the notion of a ``cellular Reedy category'' by expanding
an idea from \cite{GZ}, and prove a result ---
Proposition~\ref{ano.I.prop} --- saying that ``a cellular Reedy
category has an anodyne resolution''. By itself, it does not pass
the bar for a new result, since it involves notions introduced in
the text, but one of its corollaries does: for any
finite-dimensional simplicial set $X$, there exists a functorial
weak equivalence $Y \to X$ such that firstly, $\dim Y \leq \dim X$
and secondly, $Y$ is the nerve of a partially ordered set. The main
idea here is Thomason's double barycentric subdivision of
\cite{tho}, \cite{FL}, but the packaging is different, and the end
result is thus stronger. The result should have been known before
us, but we could not find the statement in the literature, so it
might be new.

On the contrary, nothing new at all happens in Chapter~\ref{loc.ch}
devoted to localization; some statements here are not just old by
very old (e.g.\ Lemma~\ref{kan.le} is in \cite{GZ}). This chapter
also includes Section~\ref{mod.sec} devoted to model categories, and
this is more-or-less the only place in the text where they
appear. Technically, Proposition~\ref{mod.prop} might be new, but
again, even if nobody proved this before, it is only because nobody
cared.

\medskip

Things that are new start with Chapter~\ref{semi.ch} that is
essentially devoted to Brown Representability for families of
groupoids. First, we prepare the ground by introducting the notion
of a ``CW-category'', a stripped-down version of a model category
that only has the classes $C$ and $W$ but not $F$. Apart from model
categories, main examples are $\Pos$ (with $C$ being the class of
``left-closed embeddings''), and the category of finite-dimensional
partially ordered sets.  We could not resist the temptation to
develop the theory further than strictly necessary, up to a
non-trivial theorem (Proposition~\ref{reg.prop}); it is probably
safe to skip most of Section~\ref{CW.sec} until
needed. Section~\ref{semi.sec} discusses what we call ``semiexact
families of groupoids'' over our CW-categories, and this is what
goes into our version of Brown Representability; the parts that are
really crucial for what follows are Subsection~\ref{fram.subs} and
Subsection~\ref{mv.subs}. Brown Representability itself, in several
versions, is in Section~\ref{repr.sec}; the simplest is
Theorem~\ref{coho.bis.thm}. It appears to pass the bar and be
genuinely new.

Chapter~\ref{nerve.ch} is what can be called the ``technical heart''
of the text, meaning the part with the most unpleasant notation and
tedious proofs. Here we interpolate between families of groupoids
over simplicial and bisimplicial sets (covered by Brown
Representability) and over $\Pos$ (that we need for our
applcations). As a warm-up, we start with Section~\ref{grp.sec} that
deals with simplicial sets and enhanced groupoids. As one can
expect, going from $\Pos$ to $\Delta^o\Sets$ is via the nerve
functor; in the other direction, one uses the functorial anodyne
resolutions of Proposition~\ref{ano.I.prop}. In
Section~\ref{seg.sp.sec}, we move on to bisimplicial sets. Here we
are only interested in families that satisfy the appropriate version
of the Segal condition, and these actually correspond rather
precisely to families over biordered sets of
Section~\ref{bi.sec}. Then in Section~\ref{css.sec}, we show that
the completeness condition of \cite{rzk} --- or again, a version of
it --- allows to go back to families over $\Pos$. The final
Section~\ref{aug.sec} studies what happens when we also have an
augmentation.

\medskip

Finally, in Chapter~\ref{enh.ch}, we reap the rewards: we build the
enhanced category theory. The chapter is long but surprisingly
straighforward, starting from Section~\ref{refl.sec} devoted to
reflexive families, and then onward along the lines pretty parallel
to the usual category theory of Chapter~\ref{cat.ch}. We note that
while Chapter~\ref{enh.ch} crucially depends on the results we prove
before, it is practically independent of the constructions; the
results work as perfect black boxes. Thus a viable strategy of
reading the text is to start directly with Chapter~\ref{enh.ch}
(possibly after a cursory glance at Chapter~\ref{cat.ch} to acquaint
oneself with our terminology and notation).

\medskip

Finally, let me mention three things that are very conspicuously {\em
  not} in the present text.
\begin{enumerate}
\item {\em Stable categories and Verdier localization.}
\end{enumerate}
Roughly, this corresponds to Lurie's \cite{lu.st}, and this is
definitely a must for any usable theory. Moreover, as M. Groth has
shown in his beautiful paper \cite{grotz.st}, things work even
better in the derivator formalism, up to and including what is
probably the shortest definition of a stable enhanced category (``an
enhanced category is stable iff it has all finite limits and
colimits, and they commute''). The reason I skip this material is
sheer exhaustion; I will return to this elsewhere.
\begin{enumerate}
\stepcounter{enumi}
\item {\em Grothendieck topologies.}
\end{enumerate}
This is really important for applications to algebraic geometry, so
it features prominently in \cite{lu1}. Given what we already do in
the text, it should not be too difficult to write down the notion of
an enhanced sheaf on a site and develop the associated
machinery. However, from the general categorical point of view, it
is not at all clear to me that this is all there is to it; there
well might be some enhancement already for the notion of a ``site'',
and in general, a genuinely enhanced version of the topos theory,
starting with an appropriate version of Lawvere's subobject
classifier, see e.g.\ \cite{J}. This needs (and deserves) to be
researched further.
\begin{enumerate}
\stepcounter{enumi}\stepcounter{enumi}
\item {\em Axiomatization.}
\end{enumerate}
In some sense, this is a continuation of \thetag{ii} above, and it
concerns developing homotopy theory from first principles ---
something that we still cannot do. Can one describe the category
$\Cat^h$ axiomatically, similarly to what Lawvere did for
Grothendieck's toposes? At present, I do not know; a rather feeble
attempt at it can be found in \cite[Subsection 5.5]{ka.umn}.

\section*{Standing on shoulders.}

To finish the introduction, let us try to show how the present work
fits into the context of exisiting literature, and indeed why we
feel that it fills a gap.

As mentioned above, the original idea of a ``derivator'' is due to
Grothen\-dieck \cite{champ}; another important early sources are the
pioneering works \cite{F} of J. Franke and \cite{hel} of
A. Heller. Since then, there has been a lot of development and some
false hopes. By now, there seems to be a consensus among the experts
that the theory cannot be used as an independent foundational
tool. Effectively, what we do in the present text cannot be done.

To understand the reasoning behind this consensus, let us mention
one thing we do not do: we most definitely {\em do not} attempt to
show that derivators model enhanced categories in the conventional
paradigm based on \cite{qui.ho}.

Namely, at the end of the day, the theory of derivators provides a
fully faithful embedding from enhanced categories to families of
unenhanced ones. The source $\Cat^h$ of the embedding is the
homotopy category of small enhanced categories, in any of the
equivalent definitions, and the target $(\Cat \mm I)^0$ is the
category of families of usual unenhanced categories indexed by some
$I$, localized with respect to equivalences of categories. Since we
need to localize anyway, we may replace $\Cat \mm I$ with the
functor category $I^o\Cat$, again localized with respect to
pointwise equivalences. In our approach, $I$ is $\Pos$.  More
conventionally, and ignoring set-theoretical issues, one can take
say $I=\Cat$; a functor $\Cat^o \to \Cat$ is what is usually called
a {\em prederivator}. A {\em derivator} is a prederivator satisfying
some additional conditions that insure, roughly speaking, that it
lies in the image of the embedding (in our language, this is very
close to a ``complete family'' of Definition~\ref{compl.fam.def}
that describes enhanced categories that are enhanced-complete).

In the paradigm of \cite{qui.ho}, a natural next step would be to
find a model structure on $I^o\Cat$, or maybe on a full subcategory
in $I^o\Cat$, with weak equivalences given by pointwise equivalences
of categories, and then show that it is Quillen-equivalent to the
existing models for $\Cat^h$. Of course, in the light of everything
we explained above, this might seem a rather pointless exercise,
since it would reduce the theory to merely another model for
enhanced categories, one of plenty. However, a worse problem is that
this is also impossible. Indeed, no matter what model structure one
puts on $I^o\Cat$, as soon as one fixes $W$ to be the class of
pointwise equivalences of categories, the enhanced localization
would coincide with the Dwyer-Kan localization with respect to this
class, and it is easily seen to be $2$-truncated --- that is,
homotopy types of morphisms have $\pi_i=0$ for $i \geq 2$. For
$\Cat^h$, this is certainly not the case.

In our terminology, while $\Cat^h$ has a natural enhancement
$\cCat^h$, this enhancement is {\em not} obtained by
localization. The fully faithful embedding $\Cat^h \to (\Cat \mm
\Pos)^0$ has a enhancement $\cCat^h \to \cCat \mm^h \pPos$, but this
enhancement is very far from being fully faithful. Technically, the
fibration $\cCat^h \to \pPos$ has fibers $\cCat^h_J$ spanned by
small enhanced categories $\C$ equipped with a $J$-augmentation ---
that is, a fibration $\C \to \Unf(J)$ --- while fibers $(\cCat \mmh
\pPos)_J$ are spanned by fibrations $\C \to J \times \Pos$. The two
are very different unless $J$ is very simple (see
Corollary~\ref{ccat.dim1.corr} for a comparison result for $\dim J =
1$).

Historically, the observation above was first made in the context of
derivator $K$-theory that itself has interesting history. While
$K_0$ of a stable enhanced category obviously only depends on the
underlying unenhanced triangulated category, the same is not at all
obvious for higher $K$-groups. After Schlichting \cite{schli} proved
that the statement is actually wrong --- different enhancements on
the same triangulated category can give different $K$-groups ---
people turned to derivators. Maltsiniotis \cite{ma} proposed a
conjectural construction under a very weak notion of a derivator; in
particular, it would work for derivators indexed on finite
categories. However, almost immediately To\"en and Vezzosi
\cite[Section 7]{toe-vez} showed that it cannot possibly recover the
correct $K$-theory. The argument is exactly by functoriality:
Maltsiniotis' construction was too functorial, and would factor
through the Dwyer-Kan quotient of the category of derivators. But
this quotient is $2$-truncated, and the standard $K$-theory provably
does not factor through anything $2$-truncated. The situation was
explored further by Muro and Raptis \cite{mr1}, \cite{mr2} who gave
a pretty exhaustive treatment. One side result of this is a ``folk
theorem'' saying that ``one cannot recover $K$-theory from the
derivator''. This is of course wrong, and Muro and Raptis claim
nothing of the kind, but this is what became stuck in public
perception.

\begin{remarknonum}
For the record, here is how it works in in our approach. To recover
$K$-theory from an enhanced category, one just repeats the standard
$S$-construction, and to recover it with full functoriality, it is
enough to do it in a sufficiently natural way. For example, one can
do it using the Grothendieck construction, along the lines of
\cite[Section 3.1]{ka.K}: given a stable small enhanced categry
$\C$, one constructs a semicartesian square
$$
\begin{CD}
  \wt{S}\C @>>> \fFun^h(\Delta_\idot|\Delta,\C)\\
  @VVV @VV{\xi}V\\
  \ppt^h @>{\eps^h(0)}>> \C,
\end{CD}
$$
where the relative enhanced functor category
$\fFun^h(\Delta_\idot|\Delta,\C) = \Delta^\Tot_h\C$ is the enhanced
simplicial expansion of $\C$ in the sense of \eqref{enh.del.wr},
$\xi$ is the enhanced version of corresponding functor
\eqref{xi.eq}, and the bottom arrow corresponds to the enhanced
object $0 \in \C$. Explicitly, $\Delta^\Tot_h\C$ is equipped with an
enhanced fibration $\Delta^\Tot_n\C \to \Unf(\Delta)$ with enhanced
fibers $\fFun^h(\Unf([n]),\C)$, and $\xi$ sends an enhanced functor
$\Unf([n]) \to \C$ to its value at $n \in [n]$. Then one observes
that the induced projection $\wt{S}\C \to \Unf(\Delta)$ is also an
enhanced fibration, and takes as $S\C = \Hh^\Tot(\wt{S}\C_{h\flat})$
the total localization of the corresponding enhanced family of
groupoids $\wt{S}\C_{h\flat} \to \Unf(\Delta)$ of \eqref{h.flat.sq}.
\end{remarknonum}

One a more positive side, it was suggested in \cite{mr2} and further
elaborated by Raptis in \cite{rapt.K} that one can still stay within
Quillen's paradigm if one agrees to work with a more complicated
version of a weak equivalence. The first result here is the paper
\cite{renaud} by O. Renaudin that actually predates
\cite{mr2}. Renaudin only works with derivators, and one only
expects them to describe enhanced-complete enhanced categories; these
are also obtained as localizations of combinatorial model
categories. Renaudin proves that this gives rise to a fully faithful
embedding --- in effect, one can model combinatorial model
categories by derivators. For arbitrary enhanced categories, one
needs arbitrary prederivators; the corresponding embedding theorem
was proved by Arlin in \cite{Ar}, with further pretty definitive
results obtained by Fuentes-Keuthan, Kedziorek and Rovelli in
\cite{FKR}.

The results of \cite{Ar} and \cite{FKR} are actually quite parallel
to our main results: they prove that the enhanced category of small
enhanced categories embeds fully faithfully into the enhanced
category of prederivators, and characterize the image of the
embedding. However, they prove both more and less. ``More'' because
they work in Quillen's paradigm and actually obtain a Quillen
equivalence between the source and a full subcategory of the
target. ``Less'' because the notion of a weak equivalence between
derivators is necessarily quite complicated, and very far from the
naive notion that works for us. Thus the results in \cite{Ar} and
\cite{FKR} are quite strong and very interesting, but from the
foundational point of view, they lead nowhere: we obtain one more
model for enhanced categories, even more inexplicit than the ones in
the literature, and what is the point, then.

Altogether, one might summarize the situation as follows: there are
two and only two approaches to the foundations of homotopy theory,
either \cite{qui.ho} or derivators. The two are not compatible for
very solid and fundamental reasons that cannot be circumvented.

From this point of view, perhaps the only thing that is genuinely
new in the present text is the following suggestion: if this is how
things are, then by far the most natural and easy way out of the
diffuculty is to ditch \cite{qui.ho}.

\section*{Acknowledgements.}

The work presented here is the end result of a rather long project;
over the years, I had a lot of help.

I was first exposed to the ideas of \cite{champ} back in 1991, when
I was a first-year graduate student, by Vladimir Voevodsky and David
Kazhdan. At the time, the common emphasis was on ``higher
categories'', not on derivators. I tried to do something in
area. David was my Ph.D. advisor, and he certainly would have
supported this, were I able to come up with anything meaningful;
unfortunately, this was very clearly not the case.

In the following years, I was doing other things, with the problem
of enhancements somewhere on the back burner, and I would only come
back to it from time to time (for example, during a memorable visit
to MSRI in Spring 2006). I would also give talks on the subject,
from time to time, and while I still had vague ideas at most, and no
actual results, it is amazing how supportive and understanding many
colleagues were. This certainly includes friends and colleagues in
Moscow, especially Misha Verbitsky, Sasha Kuznetsov, Dima Orlov, and
later, Sasha Efimov, Anton Fonarev and Chris Brav. But in other
places people were just as wonderful. I still remember a talk in
Tokyo, in 2008, after which a prominent algebraic geometer came to
me and said something like this: ``I liked your talk; of course, the
last thing the world needs are new foundations for homological
algebra, but at least, there was a story''. This was one of the best
pieces of advice I ever had: no matter what you do, people will
listen if there is a story.

It is of course impossible to write down the whole list of places
where I told these stories over the years, in their evolving form,
and of people who listened to them and in turn, told me back
theirs. However, there are some I just have to mention. Israel in
general, and Hebrew University of Jerusalem in particular is an
amazing place, and I am very happy that I was able to visit it
several times, at the invitation of David Kazhdan. Talking with
David --- about all sorts of mathematics, but enhancements in
particular --- has been always extremely valuable, I learned a lot
from him and I continue to do so. Equally valuable were long
discussions with Vladimir Hinich; especially in thinking about
enhancements and homotopical algebra, I owe him and David a huge
intellectual debt.

Another really amazing place with amazing people is
Chicago/Evanston, especially with Victor Ginzburg, Sasha Beilinson
and Vladimir Drinfeld at the University of Chicago, and Boris
Tsygan, Dmitry Tamarkin and Ezra Getzler at Northwestern. Over the
years, I had an opportunity to visit both places for extended
periods of time, and this in particular helped a lot to develop the
ideas presented in this work.

A different group of colleagues to whom I owe a lot in my
understanding of homotopical algebra are Carlos Simpson, Bertrant
To\"en and Gabriele Vezzosi, the founders of modern Derived
Algebraic Geometry. Talking to Carlos has always been, and continues
to be very inspirational and illuminating. Bertrand does not even
like fibered categories, and is quite happy to work ``up to
homotopy'', but his expertise in the area and understanding of it is
nothing short of astonishing; it is hard to overestimate how much I
learned from him, Gabriele and Carlos.

Somewhat similarly, almost all of actual homotopy theory that I know
I learned from Lars Hesselholt, and he has been a tremendous
influence on my understanding of the subject, both in its
mathematical content, and ideologically, in his passionate push for
proper and modern foundations for it. Lars probably will not much
like this work either --- he seems quite happy with the
$\infty$-categorical formalism as it already exists --- but well, at
least, there is a story. If it amuzes him, I will be glad.

\medskip

The actual work on this text was started in 2017, after I finished
\cite{ka.glue}, my only attempt to do something with model
categories that ended with much more of a whimper than a bang. The
first few pages were written in the Hausdorff Center in Bonn where I
visited at the kind invitation of Tobias Dyckerhoff. The general
shape of the theory was largely clear by mid-2018, with the bulk of
the work done in Moscow, in Mexico City and in Miami (I am very
grateful to Jacob Mostovoy for an invitation to CINVESTAV, and to
Ludmil Katzarkov for an invitation to the University of Miami). I
first presented the theory at a conference in Miami in January 2018,
to an audience that included Maxim Kontsevich and Mikhail Kapranov;
I do appreciate the opportunity to tell stories to such an audience
without anything resembling a finished text, and for the patience
and toleration of that audience.

Unfortunately, finishing the text took much longer than expected,
with other things getting in the way, and a lot of time wasted on
what turned out to be dead-ends. In the intervening years, I kept
giving talks and promising the text in the future, and I am grateful
to many people for listening to all this, trusting my word, and
often giving me very valuable feedback. At the very least, this
includes people in two Oberwolfach meetings, in 2018 and 2024, at
seminars and conferences in Bonn, Mainz, Paris, Antwerp, London,
Evanston, Miami, Simons Center in Stony Brook and University of
Colorado in Boulder --- and lots of places in Russia, both in Moscow
and outside of it. Of the many colleagues with whom I discussed this
project privately, and who helped a lot, I should at least mention
Tobias Dyckerhoff, Moritz Groth, Thomas Nikolaus, Wendy Lowen and
Boris Shoikhet; it was an honor to also present the work to Vladimir
Drinfeld and Nick Rozenblyum, and hear what they had to say. I am
also very grateful for the attention and feedback from Richard
Thomas and people at the Imperial College in London.

The final part of the project, Chapter~\ref{enh.ch}, was written
largely in Spring 2024 at MSRI --- or rather, SLMath now --- where I
spent three months as a Chern Research Professor. The place is just
as perfect as it was back in 2006, and I am very grateful for its
generosity, hospitality and wonderful working atmosphere.

After the first version of the text was posted, I had further
opportunities to discuss the theory presented here with colleagues;
in particular, it was a great pleasure to talk to Fernando Muro. I
also had an extremely useful and illuminating e-mail exchange with
George Raptis who suggested a lot of very relevant references that I
should have known myself. I am quite grateful to George for
effectively doing my homework for me; the last section of this
Introduction has been added as a result.

\medskip

In addition, I would like to thank anonymous referees for their hard
and diligent work on such a challenging and long text, and making a
lot of very useful suggestions and corrections.

\medskip

Finally, there are two colleagues and friends whom I would like to
thank separately, both for help with this particular project and for
continued help and support throughout my mathematical life. The
first is Sasha Beilinson. I have already mentioned him when speaking
about Chicago, but we of course go back much further than
that. While he never was my teacher or advisor in any formal sense
of the word, I learned an enourmous amount of mathematics from him.
What's more, in some sense, I learned from him how to think about
mathematics --- or at least, I learned as much as I was able, and I
will try to learn more. The second is Tony Pantev, a mathematician
of amazing strength and knowledge, and one of the most modest and
unassuming people that I know. Again, he was never my formal
teacher, but informally, he taught me a lot, ever since my grad
school days, and he was always interested in my work --- be it this
particular project, or any of the other things that I was doing. His
feedback and his constant, consistent interest and encouragement was
and is extremely valuable to me, and I very grateful and honored for
the good fortune of having both.

\chapter{Category theory.}\label{cat.ch}

This chapter contains preliminaries. Nothing whatsoever is
new. Section~\ref{cat.sec} mostly serves to fix notation and
terminology. For terminology, non-standard terms are ``dense
subcategory'', ``closed class'', ``left-saturated class'', and also
possibly an ``epivalence'' and a ``semicartesian square''. For any
category $I$, an ``$I$-augmentation'' for a category $\C$ is a
functor $\C \to I$. For notation, $\C^o$ is the opposite category to
$\C$, and $\C^>$ resp.\ $\C^<$ is $\C$ with a new terminal
resp.\ initial object. We use $|\C|$ resp.\ $\|\C\|$ to denote the
cardinality of a small resp.\ essentially small category; the
precise meaning is below in Subsection~\ref{cat.subs}
resp.\ Subsection~\ref{fun.subs}.  We use $\C /_\pi I$ resp.\ $I
\setminus_\pi \C$ for the left resp.\ right comma-categories of a
functor $\pi:\C \to I$, and we use $\Cat \bb I$ and variations of
this notation for different categories of $I$-augmented small
categories. We denote by $\Fun(I,\C)$ the category of functors from
$I$ to $\C$, and we shorten $\Fun(I^o,\C)$ to
$I^o\C$.

Section~\ref{adj.sec} is equally standard and contains various
useful facts about adjunction, limits and Kan extensions. Left
resp.\ right Kan extension along a functor $\gamma$ are denoted by
$\gamma_!$ resp.\ $\gamma_*$. ``Left'' resp.\ ``right-reflexive''
means ``admits a left'' resp.\ ``right-adjoint''. A full subcategory
is left resp.\ right-admissible if its embedding functor is left
resp.\ right-reflexive.

Section~\ref{fib.sec} deals with the Grothendieck construction and
contains both general facts and lots of examples of Grothendieck
fibrations and cofibrations. We then add some facts about fibrations
whose fibers are groupoids (those are called ``families of
groupoids'').

Finally, Section~\ref{fun.fib.sec} contains some slightly less
standard facts about functoriality for various categories of
augmented categories, and a generalization of the right Kan
extension to fibrations (understood as contravariant pseudofunctors
to the category of categories).

\section{Categories and functors.}\label{cat.sec}

\subsection{Categories.}\label{cat.subs}

We work in the minimal set-theoretical setup where there are sets
and classes, and small and large categories, and that is the end of
it. We denote by $\Sets$ the category of sets. For any subset $S'
\subset S$ in a set $S$, we denote by $S \ssetminus S' \subset S$
the complement (note that this is different from $S \setminus S'$
that we will need for other things). A cardinal is an isomorphism
class of sets, and for any set $S$, we denote the corresponding
cardinal by $|S|$. Cardinals are ordered in the usual way ($|S| \leq
|S'|$ iff there exists an injective map $S \to S'$). Categories are
large unless indicated otherwise. We denote by $\Cat$ the category
of small categories.  For any category $\C$, we write $c \in \C$ as
a shorthand for ``$c$ is an object of $\C$'', and for any two
objects $c,c' \in \C$, we denote by $\C(c,c')$ the set of morphisms
$f:c \to c'$. Note that by definition, even if $\C$ is large,
$\C(c,c')$ is always a set. The composition of morphisms $f:c \to
c'$, $f':c' \to c''$ is denoted by $f' \circ f:c \to c''$.  By abuse
of terminology, we use ``maps'', ``morphisms'' and ``arrows''
interchangeably, so that a map $f:c \to c'$ is the same thing as a
morphism $f:c \to c'$ and the same thing as an arrow $f:c \to
c'$. An ``isomorphism'' is the same thing as an ``invertible map'',
and $c \cong c'$ is shorthand for ``there exists an isomorphism $c
\to c'$'', with $f:c \cong c'$ standing for ``$f:c \to c'$ that is
invertible''.

Objects of the {\em arrow category} $\Ar(\C)$ of a category $\C$ are
arrows $c' \to c$ in $\C$, with morphisms from $c'_0 \to c_0$ to
$c'_1 \to c_1$ given by commutative squares
\begin{equation}\label{ar.sq}
\begin{CD}
c_0' @>>> c_0\\
@VVV @VVV\\
c_1' @>>> c_1.
\end{CD}
\end{equation}
An object $x \in \C$ is a {\em retract} of an object $y \in \C$ if
there exist maps $a:x \to y$, $b:y \to x$ such that $b \circ a =
\id$. In this case, $p=a \circ b:y \to y$ satisfies $p^2=p$ and is
known as a {\em projector} or an {\em idempotent endomorphism} of
the object $y$, and $x$ is called the {\em image} of the projector
$p$. If the image exists, it is unique up to a unique isomorphism
--- for any other $x'$, $a':x' \to y$, $b':y \to x'$, $a' \circ b' =
p$, the maps $b \circ a':x' \to x$ and $b' \circ a:x \to x'$ are
inverse to each other --- and automatically preserved by any
functor. We say that a morphism $f$ in $\C$ is a {\em retract} of a
morphism $g$ if $f$ is a retract of $g$ in the arrow category
$\Ar(\C)$. A category $\C$ is {\em Karoubi-closed} if all projectors
in $\C$ have images, and a subcategory $\C' \subset \C$ in a
category $\C$ is {\em Karoubi-dense} if any object in $\C$ is a
retract of an object in $\C'$.

For any category $\C$, we denote by $\C^o$ the opposite category ---
that is, the category with the same objects as $\C$ and morphism
sets $\C^o(c,c') = \C(c',c)$. We denote by $\ppt$ the point category
(one object, one morphism). For any category $\C$, there exists a
unique {\em tautological projection} $\C \to \ppt$, and for any $c
\in \C$, we denote by
\begin{equation}\label{eps.c}
\eps(c):\ppt \to \C
\end{equation}
a unique embedding from $\ppt$ onto $c$. An object $o \in \C$ is
{\em initial} resp.\ {\em terminal} if for any $c \in \C$, there
exists exactly one morphism $o \to c$ resp.\ $c \to o$. An
initial resp.\ a terminal object is unique up to a unique isomorphism,
if it exists. For any category $\C$, we denote by $\C^>$ the
category obtained by formally adding a terminal object $o$ to $\C$,
and we denote by $\C^<$ the category obtained by adding an initial
object $o$, so that $\C^{<o} \cong \C^{o>}$. A category is {\em
  pointed} if it has an initial object $0$ and a terminal object
$1$, and the unique map $0 \to 1$ is an isomorphism. The point
category $\ppt$ is of course pointed; the smallest non-pointed
category is the {\em single arrow category} $[1]$ that has two
objects $0$, $1$ and the unique non-identity map $0 \to 1$. By
definition, we have $[1] \cong \ppt^< \cong \ppt^>$. Another
category that deserves a special name is
\begin{equation}\label{V.eq}
\V = \{0,1\}^<.
\end{equation}
Explicitly, $\V$ has three objects $o,0,1$ and two non-identity
arrows $o \to 0$, $o \to 1$. We note that we have an isomorphism
$[1]^2 = [1] \times [1] \cong \V^>$ sending $1 \times 0$ resp.\ $0
\times 1$ to $0$ resp.\ $1$, and the initial resp.\ terminal object
in $[1]^2$ to the initial resp.\ terminal object in $\V^> \cong
\{0,1\}^{<>}$.

A category $\C$ is {\em connected} if it is not empty, and any two
objects $c,c' \in \C$ can be connected by a finite zigzag of
morphisms (equivalently, $\C$ does not admit a decomposition $\C
\cong \C_1 \copr \C_2$ with non-empty $\C_1$, $\C_2$). A category
$\C$ is {\em $\Hom$-bounded} by an infinite cardinal $\kappa$ iff
$|\C(c,c')| < \kappa$ for any $c,c' \in \C$, and {\em $\Hom$-finite}
iff it is $\Hom$-bounded by all infinite cardinals (that is,
$\Hom$-sets are finite). For any $n \geq 0$, a {\em chain} of length
$n$ in a category $\C$ is a diagram
\begin{equation}\label{c.idot}
\begin{CD}
c_0 @>>> \dots @>>> c_n,
\end{CD}
\end{equation}
and a chain is {\em non-degenerate} if none of the maps $c_{l-1} \to
c_l$, $1 \leq l \leq n$ in \eqref{c.idot} is an identity map. If
$\C$ is small, then all non-degenerate chains in $\C$ form a set
$\Ch(\C)$, and we denote $|\C| = |\Ch(\C)|$. In particular, chains
of length $0$ are objects in $\C$, and chains of length $1$ are
morphisms. If the subset $\Ch_{\leq 1}(\C) \subset \Ch(\C)$ of
chains of length $\leq 1$ is infinite, then $|\C| = |\Ch_{\leq
  1}(\C)|$, but in general, this need not be true.

\begin{exa}\label{P.exa}
Let $P$ be the category with one object and one non-identity
morphism $p$ such that $p^2=p$. Then $P$ has exactly one
non-degenerate chain of each length $n$, so $|P|$ is infinite, while
$|\Ch_{\leq 1}(P)| = 2$.
\end{exa}

A subcategory $\C' \subset \C$ is {\em full} if for any $c,c' \in
\C' \subset \C$, we have $\C'(c,c')=\C(c,c')$ --- informally, $\C'$
has the same morphisms as $\C$ (but less objects). A full
subcategory is thus completely defined by the collection $S$ of the
objects $c \in \C$ that lie in $\C'$ (we say that $\C'$ is {\em
  spanned} by objects in $S$). By abuse of terminology, all full
subcategories $\C' \subset \C$ are also assumed to be {\em strictly
  full}, in the sense that if $c \cong c'$ and $c' \in \C'$, then $c
\in \C'$. It is also useful to reserve a term for a subcategory $\C'
\subset \C$ that has the same objects as $\C$; we call such
subcategories {\em dense}. We say that a class $v$ of morphisms in a
category $\C$ is {\em closed} if it is closed under compositions and
contains all the identity maps. Then giving a dense subcategory $\C'
\subset \C$ is equivalent to giving a closed class $v$ of morphisms
in $\C$, and we say that $v$ {\em defines} the dense subcategory,
and denote it by $\C_v \subset \C$.

The minimal closed class is the class $\Id$ of all identity maps; a
category $\C$ is {\em discrete} if $\C = \C_{\Id}$. The maximal
class is the class $\Tot$ of all maps in $\C$, $\C_{\Tot} = \C$. We
also have the closed class $\Iso$ of all invertible maps.

\begin{defn}\label{iso.def}
The {\em isomorphism groupoid} of a category $\C$ is the dense
subcategory $\C_{\Iso} \subset \C$, and a category $\C$ is {\em
  rigid} if $\C_{\Iso}$ is discrete.
\end{defn}

\begin{exa}\label{fin.rig.exa}
A finite category --- that is, a small category $\C$ such that $|\C|
< \kappa$ for any infinite cardinal $\kappa$ --- is automatically
rigid (already a single non-trivial endomorphism $f:c \to c$ of an
object creates an infinite number of non-degenerate chains
\eqref{c.idot}, with $c_\idot = c$ and the maps $f:c_\idot \to
c_{\idot+1}$))
\end{exa}

\begin{defn}\label{sat.def}
A class of morphisms $v$ in a category $\C$ is {\em left-saturated}
if it is closed under retracts, contains all the identity maps, and
for any composable pair of maps $f$, $g$ in $\C$ such that $f \in
v$, we have $g \in v$ if and only if $f \circ g \in v$. A class $v$
is {\em saturated} if both $v$ and the opposite class $v^o$ are
left-saturated. The {\em saturation} $s(\C,W)$ of a class $W$ of
morphisms in $\C$ is the smallest saturated closed class such that
$W \subset s(\C,W)$.
\end{defn}

Equivalently, a class $v$ is saturated if it is closed under
retracts and has the ``two-out-of-three'' property: for any
composable pair $f$, $g$, as soon as two of the three maps $f$, $g$,
$f \circ g$ are in $v$, so is the third. The classes $\Id$, $\Tot$,
$\Iso$ in any category $\C$ are saturated. A left-saturated class is
automatically closed. For any functor $\gamma:\C' \to \C$ and a
class of maps $v$ in $\C$, we denote by $\gamma^*v$ the class of
maps $f$ in $\C'$ such that $\gamma(f) \in v$. If $v$ is closed
and/or saturated, then so is $\gamma^*v$, and $\gamma$ induces a
functor $\gamma:\C'_{\gamma^*v} \to \C_v$. We will say that $\gamma$
{\em inverts a map} $f$ iff $f \in \gamma^*(\Iso)$ (that is,
$\gamma(f)$ is invertible in $\C$).

A {\em factorization system} $\langle L,R \rangle$ on a category $C$
is given by two closed classes of morphisms $L$, $R$ such that $L
\cap R = \Iso$, any morphism $f:c \to c'$ in $\C$ decomposes as
\begin{equation}\label{facto.eq}
\begin{CD}
c @>{l}>> \wt{c} @>{r}>> c',
\end{CD}
\end{equation}
with $l \in L$, $r \in R$, and the decomposition is unique up to a
unique isomorphism. It is known that each of the two classes $L$,
$R$ in a factorization system completely determines the other one,
and for any commutative square
\begin{equation}\label{facto.sq}
\begin{CD}
c_0 @>{f}>> c_1\\
@V{l}VV @VV{r}V\\
c_0' @>{f'}>> c_1'
\end{CD}
\end{equation}
with $l \in L$, $r \in R$, there exists a unique map $q:c_0' \to
c_1$ such that $f = q \circ l$ and $f' = r \circ q$. Moreover,
factorizations \eqref{facto.eq} are functorial with respect to $f$
(considered as an object in the arrow category $\Ar(\C)$). For these
facts and other details on factorization systems, we refer the
reader to \cite{bou2}.

\subsection{Functors.}\label{fun.subs}

For any categories $\C_0$, $\C_1$ and functor $\gamma:\C_0 \to
\C_1$, we denote by $\gamma^o:\C_0^o \to \C_1^o$ the opposite
functor between the opposite categories.  We let
\begin{equation}\label{iota.C.eq}
  \iota:\Cat \to \Cat
\end{equation}
be the involution sending $\C$ to $\C^o$ and $\gamma$ to
$\gamma^o$. For any map $a:\gamma \to \gamma'$ between two functors
$\gamma,\gamma':\C \to \C'$, and for any functors $\gamma_0:\C_0 \to
\C$, $\gamma'_0:\C' \to \C'_0$, we denote by $\gamma'_0(a):\gamma'_0
\circ \gamma \to \gamma'_0 \circ \gamma'$ and $(a)\gamma_0:\gamma
\circ \gamma_0 \to \gamma' \circ \gamma_0$ the natural maps induced
by $a$.

For any category $I$, an {\em $I$-augmented category} or simply a
{\em category over $I$} is a category $\C$ equipped with a functor
$\pi:\C \to I$.  For any two such categories $\langle \C,\pi
\rangle$, $\langle \C',\pi' \rangle$, a {\em lax functor from $\C$
  to $\C'$ over $I$} is a functor $\gamma:\C \to \C'$ equipped with
a map $\alpha:\pi \to \pi' \circ \gamma$.  For any two such functors
$\gamma,\gamma':\C \to \C'$ we will say that a {\em map from
  $\gamma$ to $\gamma'$ over $I$} is a map $a:\gamma \to \gamma'$
such that $\pi'(a) \circ \alpha=\alpha'$. A {\em functor from $\C$
  to $\C'$ over $I$} is a lax functor $\langle \gamma,\alpha
\rangle$ such that $\alpha$ is an isomorphism. A {\em section} of a
functor $\pi:\C \to I$ is a functor $\sigma:I \to \C$ over $I$ (or
explicitly, a functor $\sigma:I \to \C$ equipped with an isomorphism
$\pi \circ \sigma \cong \id$). A {\em factorization} of a functor
$\gamma:\C \to \C'$ is a diagram
\begin{equation}\label{fac1.dia}
\begin{CD}
\C @>{\gamma'}>> \C'' @>{\gamma''}>> \C'
\end{CD}
\end{equation}
equipped with an isomorphism $\gamma \cong \gamma'' \circ
\gamma'$. For any category $I$, we denote by $\Cat \bb I$ the
category of $I$-augmented small categories $\C$, with morphisms
given by lax functors over $I$, and we let $\Cat \bbi I \subset \Cat
\bb I$ be the dense subcategory whose morphisms are functors over
$I$.

Dually, a {\em category under $I$} means a category $\C$ equipped with
a functor $\phi:I \to \C$. For any two such categories $\langle
\C,\phi \rangle$, $\langle \C',\phi' \rangle$, a {\em lax functor
  from $\C$ to $\C'$ under $I$} is a functor $\gamma:\C \to \C'$
equipped with a map $\alpha:\gamma \circ \phi \to \phi'$, and a {\em
  functor under $I$} is a lax functor $\langle \gamma,\alpha
\rangle$ under $I$ such that $\alpha$ is an isomorphism. For any lax
functors $\gamma,\gamma':\C \to \C'$ under $I$, a {\em map from
  $\gamma$ to $\gamma'$ under $I$} is a map $a:\gamma \to \gamma'$
such that $\alpha' \circ (a)\phi=\alpha$. If $I$ is small, we
denote by $I \bbd \Cat$ the category of all small categories under
$I$ and lax functors under $I$ between them, and we let $I \bbdi
\Cat \subset I \bbd \Cat$ be the dense subcategory whose morphisms
are functors under $I$.

For any functor $\pi:\C \to I$, and any object $i \in I$, we will
denote by $\C/_\pi i$ resp.\ $i \setminus_\pi \C$ the category of
objects $c \in \C$ equipped with a map $\alpha:\pi(c) \to i$
resp.\ $\alpha:i \to \pi(c)$, and we will drop $\pi$ from notation
when it is clear from the context. In particular, for any $c \in
\C$, we write $\C/c = \C/_{\id} c$ and $c \setminus \C = c
\setminus_{\id} \C$. The categories $\C/_\pi i$, $i \setminus_\pi
\C$ are known as the left and right {\em comma-fibers} of the
functor $\pi$, and we denote by
\begin{equation}\label{p.i}
\sigma(i):\C /_\pi i \to \C, \qquad \tau(i):i \setminus_\pi \C \to \C
\end{equation}
the forgetful functors sending $\langle c,\alpha \rangle$ to
$c$. For any map $f:i \to i'$ in $I$, we also have functors
\begin{equation}\label{f.i}
f_!:\C /_\pi i \to \C /_\pi i', \quad f^*:i' \setminus_\pi
\C \to i \setminus_\pi \C
\end{equation}
sending $\langle c,\alpha \rangle$ to $\langle c,f \circ \alpha
\rangle$ resp.\ $\langle c,\alpha \circ f\rangle$, and we have
natural isomorphisms $\tau(i) \circ f^* \cong \tau(i')$,
$\sigma(i') \circ f_! \cong \sigma(i)$. If we have another $I$-augmented
category $\langle \C',\pi' \rangle$, then for any $i \in I$, a lax
functor $\langle \gamma,\alpha \rangle:\C' \to \C$ over $I$ induces
a functor
\begin{equation}\label{lax.i}
i \setminus_{\pi'} \C' \to i \setminus_\pi \C, \qquad \langle
c,\alpha \rangle \mapsto \langle \gamma(c),\alpha(c) \circ \alpha
\rangle.
\end{equation}
The {\em fiber} $\C^\pi_i$ of the functor $\pi$ is the full
subcategory in $\C/_\pi i$ (or in $i \setminus_\pi \C$) spanned by
$\langle c,\alpha \rangle$ with invertible $\alpha$, and we will
again shorten $\C_i = \C^\pi_i$ when $\pi$ is clear from the
context. The {\em right} resp.\ {\em left comma-categories} $I
\setminus_\pi \C$ resp.\ $\C/_\pi I$ are the categories of triples
$\langle c,i,\alpha \rangle$, $c \in \C$, $i \in I$, $\alpha$ a map
from $i$ to $\pi(c)$ resp.\ from $\pi(c)$ to $i$; as before,
we drop $\pi$ from notation if it is clear from the context. We have
the forgetful functors
\begin{equation}\label{co.st.eq}
\tau:\C /_\pi I \to I, \quad \sigma:I \setminus_\pi \C \to I, \qquad
\langle c,i,\alpha \rangle \mapsto i,
\end{equation}
and the usual fibers of the functors \eqref{co.st.eq} are the
comma-fibers of the functor $\pi$: for any object $i \in I$, we have
identifications $(\C /_\pi I)_i \cong \C /_\pi i$ and $(I
\setminus_\pi \C)_i \cong i \setminus_\pi \C$. If $\C$ is $\Cat$ or
some other category whose objects can themselves be treated as
categories in a natural way, our notation becomes ambiguous: for any
$I \in \C$, $\C / I$ and $I \setminus \C$ are either the
comma-categories with respect to some functor $\C \to I$, or the
comma-fibers of $\C$ over $I$ considered as an object in $\C$. To
avoid ambiguity, we denote the latter by $\C / [I]$, $[I] \setminus
\C$.

\begin{exa}\label{ar.exa}
For any category $\C$ with the arrow category $\Ar(\C)$ and identity
functor $\id:\C \to \C$, we tautologically have $\C /_{\id} \C \cong
\C \setminus_{\id} \C \cong \Ar(\C)$, the projections
$\sigma,\tau:\Ar(\C) \to \C$ of \eqref{co.st.eq} send an arrow to
its source resp.\ its target, and the fibers of these projections
are the left resp.\ right comma-fibers of $\id:\C \to \C$: we have
$\Ar(\C)^\sigma_c \cong c \setminus \C$ and $\Ar(\C)^\tau_c \cong
\C/c$.
\end{exa}

A functor $\gamma:\C' \to \C$ is {\em faithful} resp.\ {\em full} if
for any $c,c' \in \C'$, the map $\gamma:\C'(c,c') \to
\C(\gamma(c),\gamma(c'))$ is injective resp.\ surjective, and
$\gamma$ is {\em fully faithful} if it is both faithful and full. In
particular, the embedding functor of a full subcategory $\C' \subset
\C$ is fully faithful. A functor $\gamma:\C' \to \C$ is {\em
  essentially surjective} if it induces a surjection on isomorphism
classes of objects (that is, for any object $c \in \C$ there exists
an object $c' \in \C'$ such that $\gamma(c') \cong c$). We will say
that $\gamma$ is {\em essentially bijective} if it induces a
bijection on isomorphism classes of objects, or in other words, if
$\gamma$ is essentially surjective, and for any $c'_0,c'_1 \in \C'$,
$\gamma(c'_0) \cong \gamma(c'_1)$ implies $c'_0 \cong c'_1$.

A functor $\gamma$ is an {\em equivalence} if it is essentially
surjective and fully faithful. A category $\C$ is {\em essentially
  small} if there exists an equivalence $\C' \to \C$ with small
$\C'$, and the same applies to any other adjective instead of
``small'' (including but not limited to ``finite'', ``discrete'',
``rigid''). We note that a small category $\C'$ equivalent to an
essentially small category $\C$ can be chosen canonically -- for
example, one can take the unique small category $\C_{red}$ that is
equivalent to $\C$ and has exactly one object in each isomorphism
class. We denote $\|\C\| = |\C_{red}|$. A category $\C$ is {\em
  locally small} if $\C/c$ is essentially small for any object $c
\in \C$. A functor $\gamma:\C' \to \C$ is {\em small} if for any
essentially small full subcategory $\C_0 \subset \C$,
$\gamma^{-1}(\C_0) \subset \C'$ is essentially small (or
equivalently, the fiber $\C'_c$ is essentially small for any $c \in
\C$).

We will also need a relaxation of the notion of equivalence. Namely,
we recall that a functor $\gamma:\C' \to \C$ is {\em conservative}
if it reflects isomorphisms --- that is, $\gamma^*(\Iso) =
\Iso$. Note that a fully faithful functor is automatically
conservative.

\begin{defn}\label{epi.def}
  A functor $\gamma$ is an {\em epivalence} if it is full,
  essentially surjective, and conservative.
\end{defn}

An epivalence $\gamma:\C \to \C'$ is an equivalence if and only if
it is faithful. Even if it is not, it is automatically essentially
bijective, and surjective in the following strong sense: if a
functor $\phi:\C \to \E$ to some category $\E$ is isomorphic to the
composition $\phi' \circ \gamma$ for some $\phi':\C' \to \E$, then
$\phi$ together with the isomorphism $\phi \cong \phi' \circ \gamma$
define $\phi'$ uniquely up to a unique isomorphism. In such a
situation, we will say that $\phi$ {\em factors through} the
epivalence $\gamma$. In particular, if a composition of two
epivalences is an equivalence, then both are actually equivalences.

\begin{lemma}\label{epi.kar.le}
Assume given an epivalence $\gamma:\C \to \C'$ and a Karoubi-dense
full subcategory $\C'_0 \subset \C'$. Then $\C_0=\gamma^{-1}(\C'_0)
\subset \C$ is Karoubi-dense.
\end{lemma}

\proof{} For any $c \in \C$, choose $c' \in \C'_0$ and maps
$a:\gamma(c) \to c'$, $b:c' \to \gamma(c)$ such that $b \circ a
=\id$. Then $c' \cong \gamma(c'')$ for some $c'' \in \C_0$ since
$\gamma$ is essentially surjective, $a$ and $b$ lift to maps $a':c
\to c''$, $b':c'' \to c$ since $\gamma$ is full, and $b' \circ a'$ is
invertible since $\gamma$ is conservative, so we can replace $b'$
with $b'' = (b' \circ a')^{-1} \circ b'$ to insure that $b'' \circ
a' = \id$.
\endproof

\begin{lemma}\label{epi.epi.le}
Assume given a conservative functor $p:\C' \to \C$ and a functor
$q:\C'' \to \C'$ such that $p \circ q:\C'' \to \C$ is full. Then if
$q$ is essentially surjective, $p$ is full, and if $p \circ q$ is an
epivalence, so is $p$. Conversely, if $p \circ q$ is essentially
surjective, and $p$ is full, then $p$ is an epivalence and $q$ is
essentially surjective.
\end{lemma}

\proof{} Assume that $q$ is essentially surjective. Then for any two
objects $c'_0,c'_1 \in \C'$ and a map $f:p(c'_0) \to p(c'_1)$, we
can choose $c_0'',c_1'' \in \C''$ such that $q(c''_l) \cong c'_l$,
$l=0,1$, then choose a map $f'':c_0'' \to c_1''$ such that
$p(q(f''))=f$, and let $f'=q(f'):c_0' \to c_1'$. Then $p(f')=f$, so
that $p$ is full. Conversely, if $p \circ q$ is essentially
surjective, then $p$ is trivially essentially surjective, thus an
epivalence. Moreover, for any object $c' \in \C'$ there exists an
object $c'' \in \C''$ and an isomorphism $p(q(c'') \cong p(c')$, and
since $p$ is full, the isomorphism comes from an isomorphism $q(c'')
\cong c'$, so that $q$ is essentially surjective.
\endproof

\subsection{Commutative diagrams.}\label{dia.subs}

By abuse of notation and terminology, we say that a square of
categories and functors
\begin{equation}\label{cat.dia}
\begin{CD}
\C_{01} @>{\gamma^1_{01}}>> \C_1\\
@V{\gamma_{01}^0}VV @VV{\gamma_1}V\\
\C_0 @>{\gamma_0}>> \C
\end{CD}
\end{equation}
is {\em lax commutative} if it is equipped with a map
$\alpha:\gamma_0 \circ \gamma^0_{01} \to \gamma_1 \circ
\gamma^1_{01}$, and {\em commutative} if the map is an isomorphism
(thus commutativity of a square \eqref{cat.dia} is not a condition
but an extra bit of structure). We will say that \eqref{cat.dia} is
{\em strictly commutative} in the rare case when the categories are
small and the functors commute on the nose, so that \eqref{cat.dia}
is a commutative square in the category $\Cat$ (one case where this
is automatic is when all the categories are small, and $\C$ is
rigid). We denote by $\C_0 \times_{\C} \C_1$ the category of triples
$\langle c_0,c_1,\alpha \rangle$, $c_0 \in \C_0$, $c_1 \in \C_1$,
$\alpha:\gamma_0(c_0) \to \gamma_1(c_1)$ an isomorphism. We will
also denote $\C_0 \times_\C \C_1$ by $\gamma_0^*\C_1$ when we want to
emphasize the dependence on $\gamma_0$, and $\gamma_1$ is clear from
the context (for example, a fiber $\C'_c$ of a functor $\C' \to \C$
over an object $c \in \C$ is given by $\C'_c \cong \eps(c)^*\C'$,
where $\eps(c)$ is the embedding \eqref{eps.c}). A commutative
square \eqref{cat.dia} induces a functor
\begin{equation}\label{c.01.f}
\C_{01} \to \C_0 \times_{\C} \C_1,
\end{equation}
and we say that the square is {\em cartesian} if \eqref{c.01.f} is
an equivalence. Moreover, a square \eqref{cat.dia} is {\em weakly
  semicartesian} resp.\ {\em semicartesian} if \eqref{c.01.f} is
essentially surjective resp.\ an epivalence in the sense of
Definition~\ref{epi.def}. We note that while the actual functor
\eqref{c.01.f} depends on the choice of the isomorphism
$\alpha:\gamma_0 \circ \gamma^0_{01} \cong \gamma_1 \circ
\gamma^1_{01}$ in \eqref{cat.dia}, the property of a square to be
cartesian or semicartesian or weakly semicartesian only depends on
the isomorphisms classes of the functors $\gamma_0$, $\gamma_1$,
$\gamma^0_{01}$, $\gamma^1_{01}$, and does not depend on the choice
of $\alpha$ at all. Indeed, if we replace $\gamma_0$, $\gamma_1$,
$\gamma^0_{01}$, $\gamma^1_{01}$ by isomorphic functors, and modify
$\alpha$, then the new functor \eqref{c.01.f} is isomorphic to the
old one. More generally, a diagram
$$
\begin{CD}
  \C_0'' @>{\gamma_0'}>> \C_0' @>{\gamma_0}>> \C_0\\
  @V{\pi''}VV @VV{\pi'}V @VV{\pi}V\\
  \C_1'' @>{\gamma_1'}>> \C_1' @>{\gamma_1}>> \C_1
\end{CD}
$$
of categories and functors is {\em commutative} if so are both the
square on the left and the square on the right --- in other words,
we are given isomorphisms $\alpha:\gamma_1' \circ \pi' \to \pi \circ
\gamma_0$, $\alpha':\gamma_1' \circ \pi'' \to \pi' \circ \gamma_0'$
--- and in this case, the outer square becomes automatically
commutative once equipped with the isomorphism $\alpha'' = (\alpha)\gamma_0'
\circ \gamma_1(\alpha'):\gamma_1 \circ \gamma_1' \circ \pi'' \to \pi
\circ \gamma_0 \circ \gamma_0'$.

\begin{lemma}\label{epi.loc.le}
Assume given categories $\C$, $\C'$ over some category $I$, and a
functor $\gamma:\C' \to \C$ over $I$. Then $\gamma$ is full,
faithful, conservative, essentially surjective, an epivalence or an
equivalence if and only if for any $\nu:[1] \to I$, the same holds
for the induced functor $\nu^*(\gamma):\nu^*\C' \to \nu^*\C$.
\end{lemma}

\proof{} Clear.\endproof

Abusing terminology even further, we say that a commutative square
\eqref{cat.dia} is {\em cocartesian} if for any commutative square
$$
\begin{CD}
\C_{01} @>{\gamma^1_{01}}>> \C_1\\
@V{\gamma^0_{01}}VV @VV{\gamma'_1}V\\
\C_0 @>{\gamma'_0}>> \C',
\end{CD}
$$
there exists a functor $h:\C \to \C'$ and isomorphisms
$\alpha_l:\gamma'_l \to h \circ \gamma_l$, $l=0,1$, and the triple
$\langle h,\alpha_0,\alpha_1\rangle$ is unique up to a unique
isomorphism. Again, the property of being cocartesian only depends
on the isomorphism classes of the functors in \eqref{cat.dia}.  In
general, cocartesian squares of categories and functors are
difficult to construct, but there are special situations when it is
easy. Here is an example.

\begin{defn}\label{cl.def}
A full subcategory $\C_0 \subset \C$ is {\em left-closed} if for any
morphism $c \to c_0$ with $c_0 \in \C_0$, we have $c \in \C_0$ (that
is, $\C \setminus \C_0 \cong \C_0 \setminus \C_0$). A full
subcategory $\C_0 \subset \C$ is {\em right-closed} if $\C_0^o
\subset \C^o$ is left-closed.
\end{defn}

\begin{exa}\label{cl.exa}
Assume given a category $\C$ and two left-closed full subcategories
$\C_0,\C_1 \subset \C$ such that $\C = \C_0 \cup \C_1$. Then
$\C_{01} = \C_0 \cap \C_1$ is left-closed both in $\C_0$ and in
$\C_1$, and if we let all the functors in \eqref{cat.dia} be the
embedding functors, then the square is both cartesian and
cocartesian.
\end{exa}

\begin{exa}\label{cocart.S.exa}
A commutative square
\begin{equation}\label{S.sq}
\begin{CD}
S_0 @>{e}>> S_0'\\
@VVV @VVV\\
S_1 @>>> S_1'
\end{CD}
\end{equation}
of sets can be treated as a commutative square of the corresponding
discrete categories, and if \eqref{S.sq} is cocartesian in $\Sets$
and $e$ is injective, so that $S_0' \cong S_0 \copr (S_0' \ssetminus
S_0)$ and $S_1' \cong S_1 \copr (S_0' \ssetminus S_0)$, then
\eqref{S.sq} is cartesian and cocartesian as a square of
categories. If $i$ is not injective, \eqref{S.sq} can be cocartesian
as a commutative square of sets but not as a commutative square of
categories (e.g.\ take $S_0 = \{0,1\}$ and $S_0'=S_1=S_1' =\ppt$).
\end{exa}

In general, for any left-closed full subcategory $\C_0 \subset \C$
in a category $\C$, the complement $\C_1 = \C \ssetminus \C_0
\subset \C$ is right-closed, and vice versa. Both $\C_0 \subset \C$
and $\C_1 \subset \C$ are fibers of a unique {\em characteristic
  functor}
\begin{equation}\label{chi.eq}
\C \to [1],
\end{equation}
and conversely, for any functor \eqref{chi.eq}, the embeddings $\C_0
\subset \C$ resp.\ $\C_1 \subset \C$ of its fibers are fully
faithful and left resp.\ right-closed. The simplest example is the
product $\C = \C' \times [1]$ for some $\C'$, with the projection
$\C' \times [1] \to [1]$ and fibers $\C_0 \cong \C_1 \cong \C'$. To
describe how a category $\C$ is glued out of its left-closed
subcategory $\C_0 \subset \C$ and right-closed complement $\C_1 = \C
\ssetminus \C_0 \subset \C$, one can consider the $\Hom$-pairing
\begin{equation}\label{hom.pair.eq}
\Hom:\C^o \times \C \to \Sets, \qquad c \times c' \mapsto \C(c,c'),
\end{equation}
with its restriction $F = \Hom \circ (\eps_0^o \times \eps_1):\C_0^o
\times \C_1 \to \Sets$, where $\eps_l:\C_l \to \C$, $l=0,1$ are the
embedding functors. Conversely, for any categories $\C_0$, $\C_1$
and functor $F:\C_0^o \times \C_1 \to \Sets$, we can define a
category $\C = \C_0 \copr_F \C_1$ as the union $\C_0 \copr \C_1$,
with $\Hom$-sets given by
\begin{equation}\label{cl.glue.eq}
(\C_0 \copr_F \C_1)(c,c') = \begin{cases} \C_l(c,c'), &\quad c,c'
    \in \C_l, l=0,1,\\
    F(c,c'), &\quad c \in \C_0,c' \in \C_1,
\end{cases}
\end{equation}
and no morphisms from object in $\C_1$ to objects in $\C_0$. Then
the two constructions are inverse to each other, so that any
category $\C$ equipped with a functor $\C \to [1]$ is of the form
$\C \cong \C_0 \copr_F \C_1$ for a unique functor $F:\C_0^o \times
\C_1 \to \Sets$.

\subsection{Categories of functors.}\label{cat.fun.subs}

For any category $\E$ and essentially small category $I$, functors
from $I$ to $\E$ form a category that we denote by $\Fun(I,\E)$. We
will shorten $\Fun(I^o,\E)$ to $I^o\E$. A functor $\gamma:I_0 \to
I_1$ between essentially small categories induces a pullback functor
$\gamma^*:\Fun(I_1,\E) \to \Fun(I_0,\E)$, $E \mapsto E \circ
\gamma$, and a commutative square \eqref{cat.dia} of essentially
small categories induces a commutative square
\begin{equation}\label{fun.dia}
\begin{CD}
\Fun(\C_{01},\E) @<{\gamma_{01}^{1*}}<< \Fun(\C_1,\E)\\
@A{\gamma_{01}^{0*}}AA @AA{\gamma^*_1}A\\
\Fun(\C_0,\E) @<{\gamma^*_0}<< \Fun(\C,\E),
\end{CD}
\end{equation}
cartesian if \eqref{cat.dia} is cocartesian. The functors
\eqref{eps.c} provide a natural identification $\E \cong
\Fun(\ppt,\E)$, $e \mapsto \eps(e)$. If $I = [1]$, then
$\Fun([1],\E)$ is the arrow category $\Ar(\E)$; an arrow $f:e \to
e'$ in $\E$ defines a functor
\begin{equation}\label{eps.f}
\eps(f):[1] \to \E
\end{equation}
equipped with isomorphisms $\eps(f) \circ \eps(0) \cong \eps(e)$,
$\eps(f) \circ \eps(1) \cong \eps(e')$. For any essentially small
$I$, we have the {\em evaluation functor}
\begin{equation}\label{ev.eq}
\ev = \ev(I,\E):I \times \Fun(I,\E) \to \E, \qquad i \times E \mapsto E(i),
\end{equation}
and the pair $\langle \Fun(I,\E),\ev \rangle$ can be characterized by a
universal property: for any category $\C$, a functor $\gamma:I \times \C
\to \E$ factors as
\begin{equation}\label{ev.dia}
\begin{CD}
I \times \C @>{\id \times \gamma'}>> I \times \Fun(I,\E) @>{\ev(I,\E)}>> \E
\end{CD}
\end{equation}
for some $\gamma':\C \to \Fun(I,\E)$, uniquely up to a unique
isomorphism. Explicitly, we also have the {\em coevaluation functor}
$\ev^\dg(I,\E):\E \to \Fun(I,I \times \E)$ sending $e \in \E$ to the
functor $- \times e:I \to I \times \E$, and the compositions
$$
\begin{CD}
I \times \E @>{\id \times \ev^\dg(I,\E)}>> I \times \Fun(I,I \times \E)
@>{\ev(I,I \times \E)}>> I \times \E,\\
\Fun(I,\E) @>{\ev^\dg(I,\Fun(I,\E))}>> \Fun(I,I \times \Fun(I,\E))
@>{\Fun(I,\ev(I,\E))}>> \Fun(I,\E),
\end{CD}
$$
are both isomorphic to identity functors; then
$\gamma'$ in \eqref{ev.dia} is given by $\gamma'=\Fun(I,\gamma)
\circ \ev^\dg(I,\C)$. For any $i \in I$, the functor \eqref{ev.eq}
restricts to a functor $\ev_i = \ev \circ (\eps(i) \times
\id):\Fun(I,\E) \to \E$, $E \mapsto E(i)$.

\begin{exa}\label{yo.exa}
For any essentially small category $I$, the $\Hom$-pairing
\eqref{hom.pair.eq} corresponds by \eqref{ev.dia} to the fully
faithful {\em Yoneda embedding}
\begin{equation}\label{yo.yo.eq}
\Y:I \to I^o\Sets,
\end{equation}
where as per our general convention, $I^o\Sets = \Fun(I^o,\Sets)$.
\end{exa}

\begin{exa}\label{yo.cat.exa}
For any small $\C$ equipped with a functor $\pi:\C \to I$ to some
$I$, and for any object $i \in I$, the comma-fiber $i \setminus_\pi
\C$ is small. Therefore we have a well-defined functor $I^o \times
\Cat \bb I \to \Cat$ sending $i \times \langle \C,\pi \rangle$ to $i
\setminus_\pi \C$ that acts on morphisms via \eqref{f.i} and
\eqref{lax.i}. If $I$ is essentially small, this corresponds to a
functor
\begin{equation}\label{yo.cat.eq}
\Cat \bb I \to I^o\Cat
\end{equation}
that we call the {\em extended Yoneda embedding}. It is also fully
faithful.
\end{exa}

\begin{exa}\label{cl.glue.exa}
A functor $F:\C^o_0 \times \C_1 \to \Sets$ with essentially small
$\C_0$ corresponds by \eqref{ev.dia} to a functor $\lambda:\C_1 \to
\C_0^o\Sets$. By abuse of notation, we will denote by
\begin{equation}\label{cl.glue.bis.eq}
\C_0 \copr_\lambda \C_1 = \C_0 \copr_F \C_1
\end{equation}
the corresponding category \eqref{cl.glue.eq} glued out of $\C_0$
and $\C_1$, with its characteristic functor \eqref{chi.eq},
left-closed full embedding $\C_0 \to \C_0 \copr_\lambda \C_1$ and
the complementary right-closed full embedding $\C_1 \to \C_0
\copr_\lambda \C_1$.
\end{exa}
  
If categories $I$ and $\E$ are equipped with functors $\nu:I \to
I'$, $\pi:\E \to I'$ to some category $I'$, we denote by
$\vFun_{I'}(I,\E)$ the category of lax functors from $I$ to $\E$
over $I'$, we let $\Fun_{I'}(I,\E) \subset \vFun_{I'}(I,\E)$ be the
full subcategory of functors over $I'$, and if $I' = \ppt$, we drop
it from notation so that the notation stays consistent. For any
class $v$ of morphisms in $I$, we denote by $\Fun^v(I,\E) \subset
\Fun(I,\E)$ the full subcategory spanned by functors $I \to \E$ that
invert maps in $v$. For any closed class $v$ of maps in $\E$, we
denote by $\Fun(I,\E)_v \subset \Fun(I,\E)$ the dense subcategory
defined by the class of maps that are pointwise in $v$. We also
denote by $\Sec(I,\E) = \Fun_I(I,\E)$ the category of sections of a
functor $\E \to I$, and if we have another functor $\gamma:I'' \to
I$, then we simplify notation by writing $\Sec(I'',\E) =
\Sec(I'',\gamma^*\E)$ when $\gamma$ is clear from the context. With
this convention, we have
\begin{equation}\label{sec.fun.I}
\Fun_I(I'',\E) \cong \Sec(I'',\E).
\end{equation}
By definition, $\Fun_{I'}(I,\E)$ fits into a cartesian square
\begin{equation}\label{fun.I.sq}
\begin{CD}
\Fun_{I'}(I,\E) @>>> \Fun(I,\E)\\
@VVV @VV{\pi \circ -}V\\
\ppt @>{\eps(\nu)}>> \Fun(I,I'),
\end{CD}
\end{equation}
where $\eps(\nu)$ is the functor \eqref{eps.c}, and $\pi \circ -$
is post-composition with $\pi$.

If one fixes the target category $\E$, one can combine all the
functor categories $\Fun(I,\E)$, $I \in \Cat$ together by
considering the category $\Cat\bb\E$ of
Subsection~\ref{fun.subs}. Namely, we have the forgetful functor
\begin{equation}\label{fun.bb}
\Cat\bb\E \to \Cat, \qquad \langle \C,\pi \rangle \mapsto \C,
\end{equation}
and its fibers are $(\Cat\bb\E)_I \cong \Fun(I,\E)$. Note that
$(\Cat\bb\E)_{\ppt}$ is $\E$ itself, so that we have a canonical
full embedding
\begin{equation}\label{ups.eq}
\ups:\E \cong (\Cat\bb\E)_{\ppt} \to \Cat\bb\E,
\end{equation}
and if $\E$ is essentially small, the extended Yoneda embedding
\eqref{yo.cat.eq} of Example~\ref{yo.cat.exa} restricts to the usual
Yoneda embedding \eqref{yo.yo.eq} on the image $\ups(\E) \subset
\Cat\bb\E$; this justifies our terminology. In this case, we also
have the forgetful functor $\E \bbd \Cat \to \Cat$ with fibers $(\E
\bbd \Cat)_{\C} \cong \Fun(\E,\C)$.

For any $\E$, a closed class of maps $v$ in $\E$ defines a closed
class of maps $\Cat(v)$ in $\Cat\bb\E$ consisting of lax functors
$\langle \gamma,\alpha \rangle$ such that $\alpha$ is pointwise in
$v$, and the fibers of the induced forgetul functor
$(\Cat\bb\E)_{\Cat(v)} \to \Cat$ are then the functor categories
$\Fun(I,\E)_v$, $I \in \Cat$. In particular, $\Cat(\Iso)$ is the
class of all functors over $\E$, so that $(\Cat\bb\E)_{\Cat(\Iso)} =
\Cat\bbi\E$. More generally, if we have some category $\C$ equipped
with a functor $\phi:\C \to \Cat$, we will denote
\begin{equation}\label{C.bb}
\C\bb^{\phi}\E = \phi^*(\Cat\bb\E),
\end{equation}
and we will drop $\phi$ from notation when it is clear from the
context. A closed class $v$ of maps in $\E$ then defines a closed
class $\C(v) = \phi^*\Cat(v)$ of maps in $\C\bb^{\phi}\E$, and if
$v=\Iso$, we denote $\C\bbi^{\phi}\E = (\C\bb^{\phi}\E)_{\C(\Iso)}$.

\begin{exa}\label{iota.bb.exa}
Let $\iota:\Cat \to \Cat$ be the involution \eqref{iota.C.eq}. Then
for any category $\E$, objects in $\Cat \bb^\iota \E$ are small
categories $\C$ equipped with a functor $\pi:\C^o \to \E$, and
morphisms $\langle \C',\pi' \rangle \to \langle \C,\pi \rangle$ are
given by pairs $\langle \gamma,\alpha \rangle$ of a functor
$\gamma:\C' \to \C$ and a morphism $\alpha:\pi' \to \pi \circ
\gamma^o$. The subcategory $\Cat \bbi^\iota \E \subset \Cat \bb^\iota
\E$ has the same objects, and morphisms $\langle \gamma,\alpha
\rangle$ such that $\alpha$ is invertible. Note that we have an
equivalence $\Cat \bbi^\iota \E \cong \Cat \bbi \E^o$ sending $\langle
\C,\pi \rangle$ to $\langle \C^o,\pi \rangle$ and $\langle
\gamma,\alpha \rangle$ to $\langle \gamma^o,\alpha^{-1} \rangle$.
\end{exa}

\section{Adjunction and limits.}\label{adj.sec}

\subsection{Generalities on adjunction.}\label{adj.subs}

We recall that an {\em adjoint pair of functors} between categories
$\C$, $\C'$ is by definition a pair of functors $\lambda:\C' \to \C$,
$\rho:\C \to \C'$ equipped with {\em adjunction maps} $a_\dg:\id \to
\rho \circ \lambda$, $a^\dg:\lambda \circ \rho \to \id$ such that
the compositions
\begin{equation}\label{adj.eq}
\rho(a^\dg) \circ (a_\dg)\rho:\rho \to \rho, \qquad (a^\dg)\lambda
\circ \lambda(a_\dg):\lambda \to \lambda
\end{equation}
are both equal to the identity. The functor $\lambda$ is {\em
  left-adjoint} to $\rho$, and the functor $\rho$ is {\em
  right-adjoint} to $\lambda$. The adjunction maps induce
isomorphisms
\begin{equation}\label{adj.iso}
\C(\lambda(c'),c) \cong \C'(c',\rho(c)), \qquad c \in \C,c' \in \C',
\end{equation}
functorial in $c$ and $c'$. Either one of the functors $\lambda$,
$\rho$ determines the other one and the adjunction maps uniquely up
to a unique isomorphism, so that having an adjoint is a condition
and not a structure. For brevity, we will say that a functor
$\phi:\C \to \C'$ is {\em left} resp.\ {\em right-reflexive} if it
admits a left-adjoint $\phi_\dg:\C' \to \C$ resp.\ a right-adjoint
$\phi^\dg:\C' \to \C$. If we have a diagram
\begin{equation}\label{bc.0.dia}
\begin{CD}
\C_0 @>{\phi_0}>> \C_0'\\
@V{\gamma}VV @VV{\gamma'}V\\
\C_1 @>{\phi_1}>> \C_1'
\end{CD}
\end{equation}
of categories and functors with right-reflexive $\phi_0$, $\phi_1$,
then any morphism $\phi_1 \circ \gamma \to \gamma' \circ \phi_0$ induces
by adjunction a morphism
\begin{equation}\label{bc.eq}
\gamma \circ \phi_0^\dg \to \phi_1^\dg \circ \gamma',
\end{equation}
and dually for left-reflexive functors. We will call the maps
\eqref{bc.eq} the {\em base-change maps}, and we say that a square
\eqref{bc.0.dia} is {\em right-reflexive} if $\phi_0$, $\phi_1$ are
right-reflexive, and \eqref{bc.eq} is an isomorphism. Dually,
\eqref{bc.0.dia} is {\em left-reflexive} if the corresponding square
of the opposite categories and functors is right-reflexive. We note
that just as for cartesian squares, the property of a square
\eqref{bc.0.dia} to be left or right-reflexive only depends on the
isomorphism classes of the functors $\phi_0$, $\phi_1$, $\gamma$,
$\gamma'$.

In particular, if we have two categories $\C$, $\C'$ equipped with
functors $\pi:\C \to I$, $\pi':\C' \to I$, and a functor $\gamma:\C
\to \C'$ over $I$ that admits a right-adjoint $\gamma^\dg:\C' \to
\C$, then we have a base change map $\pi \circ \gamma^\dg \to
\pi'$. We say that $\gamma$ is {\em right-reflexive over $I$} if
this map is an isomorphism. In this case, $\gamma^\dg:\C' \to \C$ is
also naturally a functor over $I$. Dually, if $\gamma$ has a
left-adjoint $\gamma_\dg:\C' \to \C$, then it is {\em left-reflexive
  over $I$} if the base change map $\pi' \to \pi \circ \gamma_\dg$
is an isomorphism. In this case, $\gamma_\dg:\C' \to \C$ is again
naturally a functor over $I$. Note that if the left or right adjoint
$\gamma_\dg$ is fully faithful, then $\gamma \circ \gamma_\dg \cong
\id$, so that it is automatically adjoint over $I$. If we have a
functor $\nu:I' \to I$ from some category $I'$, then an adjunction
between $\gamma$ and $\gamma_\dg$ induces an adjunction between
$\nu^*(\gamma)$ and $\nu^*(\gamma_\dg)$.

\begin{exa}\label{fun.adj.exa}
Assume given categories $\C$, $\C'$ and an adjoint pair of functors
$\rho:\C \to \C'$, $\lambda:\C' \to \C$. Then for any essentially
small $I$, postcomposition with $\rho$ and $\lambda$ defines an
adjoint pair of functors $\Fun(I,\C) \to \Fun(I,\C')$, $\Fun(I,\C')
\to \Fun(I,\C)$, with adjunction maps induced by the adjunction maps
\eqref{adj.eq}. Dually, if $\C$ and $\C'$ are essentially small,
then for any $I$, $\rho^*:\Fun(\C',I) \to \Fun(\C,I)$ and
$\lambda^*:\Fun(\C,I) \to \Fun(\C',I)$ are adjoint.
\end{exa}

\begin{lemma}\label{epi.sq.le}
Assume given a right-reflexive square \eqref{bc.0.dia} such that
$\gamma'$ is conservative, and $\phi_0$, $\phi_1$ are fully
faithful. Then the square \eqref{bc.0.dia} is cartesian.
\end{lemma}

\proof{} Let $\C''_0 = \phi_1^*\C'_0$. Then since \eqref{bc.0.dia}
is commutative, $\phi_0$ factors through the fully faithful
embedding $\C''_0 \to \C'_0$, so that the corresponding functor
$\C_0 \to \C''_0$ is automatically fully faithful, and we need to
check that it is essentially surjective. Indeed, an object $c \in
\C'_0$ lies in the essential image of $\phi_0$ if and only if the
adjunction map $c \to \phi_0(\phi_0^\dg(c))$ is an isomorphism, and
since \eqref{bc.eq} is an isomorphism and $\gamma'$ is conservative,
this happens if and only if the adjunction map $\gamma'(c) \to
\phi_1(\phi_1^\dg(\gamma'(c)))$ is an isomorphism. The latter
happens if and only if $\gamma'(c)$ lies in the essential image of
$\phi_1$.
\endproof

Reflexivity is a local property with respect to the target category,
in the sense that a functor $\phi:\C' \to \C$ is right-reflexive if
and only if for any $c \in \C$, the left comma-fiber $\C' /_\phi c$
has a terminal object $\langle c',\alpha \rangle$. If this happens,
$c'$ is functorial in $c$, thus gives the adjoint functor
$\phi^\dg$, $\phi^\dg(c) = c'$, and the maps $\alpha$ give one of
the adjunction maps. The same holds for left-reflexive functors and
right-comma fibers.

\begin{exa}\label{equi.adj.exa}
If $\phi:\C' \to \C$ is an equivalence --- that is, essentially
surjective and fully faithful --- then any object $\langle c',\alpha
\rangle \in \C'_c \subset \C' /_\phi \C$ is terminal, so that $\phi$
is left-reflexive. The adjunction maps in this case are invertible,
so that $\phi^\dg$ is inverse to $\phi$, and it is also an
equivalence by \eqref{adj.iso}. One says that a category $\C$ is
{\em equivalent} to a category $\C'$ if there exists an equivalence
$\C' \to \C$, and we see that the relation is symmetric. One can
also use right comma-fibers, and the left-adjoint $\phi_\dg$ is
isomorphic to $\phi^\dg$. If $\C'$, $\C$ are equipped with functors
$\C',\C \to I$ to some $I$, and an equivalence $\phi$ is a functor
over $I$, then it is automatically left and right-reflexive over $I$
(so that an ``equivalence over I'' is the same thing as a functor
over $I$ that is an equivalence).
\end{exa}

\begin{exa}\label{cat.I.exa}
Assume given a small category $I$, and consider the category
$\Cat\bbi I$ of small categories $\C \to I$ and functors over
$I$. Then unless $I$ is rigid in the sense of
Definition~\ref{iso.def}, $\Cat\bbi I$ is {\em not} the comma-fiber
$\Cat/ [I] = \Cat /_{\id} I$ of the identity functor $\id:\Cat \to
\Cat$. We have an embedding $\Cat / [I] \to \Cat \bbi I$, and it has
a right-adjoint functor $\Cat \bbi I \to \Cat / [I]$, $\C \mapsto \C
\times_I I$, but neither is an equivalence.
\end{exa}

\begin{lemma}\label{3.adj.le}
Assume given categories $I_0$, $I_1$, $I_2$ equipped with functors
$\gamma_0:I_0 \to I_1$, $\gamma_1:I_1 \to I_2$, $\gamma_2:I_2 \to
I_0$ and morphisms $a_0:\gamma_2 \circ \gamma_1 \circ \gamma_0 \to
\id$, $a_1:\id \to \gamma_0 \circ \gamma_2 \circ \gamma_1$,
$a_2:\gamma_1 \circ \gamma_0 \circ \gamma_2 \to \id$ such that
$\gamma_0(a_0) \circ (a_1)\gamma_0=\id$, $(a_1)(\gamma_2 \circ
\gamma_1) \circ (\gamma_2 \circ \gamma_1)(a_0) = \id$, $(\gamma_1
\circ \gamma_0)(a_2) \circ (a_1)(\gamma_1 \circ \gamma_0) = \id$,
and $(a_2)\gamma_1 \circ \gamma_1(a_1) = \id$. Moreover, assume that
$\gamma_1$ is essentially surjective. Then $\gamma_2$ is fully
faithful.
\end{lemma}

\proof{} By \eqref{adj.eq}, the conditions insure that $\gamma_0$ is
right-adjoint to $\gamma_2 \circ \gamma_1$, with the adjunction maps
$a_0$, $a_1$, and $\gamma_1$ is left-adjoint to $\gamma_0 \circ
\gamma_2$, with the adjunction maps $a_1$, $a_2$. Then for any
$i_2,i_2' \in I_2$, we have
$$
I_2(i_2,i_2') \cong I_1(i_1,\gamma_0(\gamma_2(i_2'))) \cong
I_0(\gamma_2(\gamma_1(i_1)),\gamma_2(i_2')) \cong
I_0(\gamma_2(i_2),\gamma_2(i_2')),
$$
where $i_1 \in I_1$ is any object equipped with an isomorphism $i_2
\cong \gamma_1(i_1)$.
\endproof

\subsection{Defining an adjunction.}

If we are given two abstract functors $\rho:\C \to \C'$,
$\lambda:\C' \to \C$, then to specify an adjunction between them, it
suffices to provide only one of the adjunction maps. To axiomatize
the situation, it is convenient to introduce the following.

\begin{defn}\label{adj.map.def}
Assume given two functors $\lambda:\C' \to \C$, $\rho:\C \to \C'$. A
map $a:\id \to \rho \circ \lambda$ {\em defines an adjunction}
between $\rho$ and $\lambda$ if for any $c \in \C$, $c' \in \C'$, any
map $f:c' \to \rho(c)$ factors as
\begin{equation}\label{adj.facto}
\begin{CD}
c' @>{a}>> \rho(\lambda(c')) @>{\rho(f')}>> \rho(c)
\end{CD}
\end{equation}
for a unique map $f':\lambda(c') \to c$. Dually, a map $a:\lambda
\circ \rho \to \id$ defines an adjunction between $\rho$ and
$\lambda$ iff $a^o:\id \to \lambda^o \circ \rho^o$ defines an
adjunction between $\lambda^o$ and $\rho^o$.
\end{defn}

If $a$ defines an adjunction between $\rho$ and $\lambda$ in the
sense of Definition~\ref{adj.map.def}, then by virtue of their
uniqueness, the factorizations \eqref{adj.facto} provide
isomorphisms \eqref{adj.iso}, and then $a$ defines the second
adjunction map, so that $\rho$ and $\lambda$ indeed form an adjoint
pair.

\begin{lemma}\label{adj.pb.le}
Assume given a commutative diagram
\begin{equation}\label{l.r.dia}
\begin{CD}
  \C_0' @>{\lambda'}>> \C_1' @>{\rho'}>> \C_0'\\
  @V{\gamma}VV @VVV @VV{\gamma}V\\
  \C_0 @>{\lambda}>> \C_1 @>{\rho}>> \C_0
\end{CD}
\end{equation}
with cartesian squares such that there exists an isomorphism $\id
\cong \rho \circ \lambda$ that defines an adjunction between
$\lambda$ and $\rho$, or between $\rho$ and $\lambda$. Then there
exists an isomorphism $\id \cong \rho' \circ \lambda'$ that defines
an adjunction between $\lambda'$ and $\rho'$ resp.\ between $\rho'$
and $\lambda'$.
\end{lemma}

\proof{} Clear. \endproof

\begin{lemma}\label{adj.ff.le}
Assume given a fully faithful embedding $\nu:\C_0 \to \C$ and a
functor $\gamma:\C' \to \C_0$ such that $\nu \circ \gamma:\C' \to
\C$ is right resp.\ left-reflexive, with adjoint $\gamma^\dg$. Then
$\gamma$ is right resp.\ left-reflexive, with adjoint $\gamma^\dg
\circ \nu$.
\end{lemma}

\proof{} In the right-reflexive case, the adjunction map $\id \to
\gamma^\dg \circ \nu \circ \gamma$ obviously also defines an
adjunction between $\gamma$ and $\gamma^\dg \circ \nu$, and in the
left-reflexive case, the same holds for the adjunction map
$\gamma^\dg \circ \nu \circ \gamma \to \id$.
\endproof

\begin{lemma}\label{a.b.adj.le}
Assume given two functors $\lambda:\C' \to \C$, $\rho:\C \to \C'$,
and a map $a^\dg:\lambda \circ \rho \to \id$. Moreover, assume that
for any $c' \in \C'$, there exists a map $a(c'):c' \to
\rho(\lambda(c'))$ such that $\rho(a^\dg \circ \lambda(f)) \circ
a(c') = f$ for any $c \in \C$ and any map $f:c' \to \rho(c)$ in
$\C'$. Then $a^\dg$ defines an adjunction between $\rho$ and
$\lambda$, and $a(c')$ is functorial with respect to $c'$ and gives
the second adjunction map.
\end{lemma}

\proof{} Clear. \endproof

\begin{corr}\label{a.b.adj.corr}
Assume given an epivalence $\eps:\C' \to \C$ and a functor
$\lambda:\C \to \E$ such that $\lambda \circ \eps:\C' \to \E$ admits
a fully faithful right-adjoint $\rho:\E \to \C'$. Then the
adjunction isomorphism $\lambda \circ \eps \circ \rho \cong \id$
also defines an adjunction between $\lambda$ and $\eps \circ \rho$.
\end{corr}

\proof{} Since $\eps$ is an epivalence, any object $c \in \C$ is of
the form $c \cong \eps(c')$ for some $c' \in \C'$. Then the
adjunction map $c' \to \rho(\lambda(c))$ induces a map $a:c \to
\eps(\rho(\lambda(c)))$. To check that $a$ satisfies the
universal property required in Lemma~\ref{a.b.adj.le}, note that
since $\eps$ is an epivalence, any map $f:c \to \eps(\rho(e))$, $e
\in \E$ comes from a map $c' \to \rho(e)$ in $\C'$.
\endproof

\begin{defn}\label{adj.rel.def}
Assume given categories $\C$, $\C'$, $\C''$ equipped with functors
$\lambda:\C' \to \C$, $\rho:\C \to \C''$, $\gamma:\C' \to
\C''$. Then a map $a:\gamma \to \rho \circ \lambda$ {\em defines an
  adjunction} between $\rho$ and $\lambda$ {\em over} $\gamma$ if
for any $c \in \C$, $c' \in \C'$, any map $f:\gamma(c') \to \rho(c)$
factors as
\begin{equation}\label{adj.rel.facto}
\begin{CD}
\gamma(c') @>{a}>> \rho(\lambda(c')) @>{\rho(f')}>> \rho(c),
\end{CD}
\end{equation}
for a unique map $f':\lambda(c') \to c$.
\end{defn}

If $\C' = \C''$ and $\gamma=\id$, this reduces to
Definition~\ref{adj.map.def}. However, even in the general case, the
uniqueness of the decompositions \eqref{adj.rel.facto} still
provides isomorphisms
\begin{equation}\label{adj.rel.iso}
\C(\lambda(c'),c) \cong \C''(\gamma(c'),\rho(c)), \qquad c \in \C,c' \in \C',
\end{equation}
functorial in $c$ and $c'$.

\begin{exa}\label{adj.rel.exa}
If for some $\lambda:\C'' \to \C$, $\rho:\C \to \C''$, a map $a:\id
\to \rho \circ \lambda$ defines an adjunction between $\rho$ and
$\lambda$, then for any $\gamma:\C' \to \C''$, the induced map
$\gamma \to \rho \circ \lambda \circ \gamma$ trivially defines an
adjunction between $\rho$ and $\lambda \circ \gamma$. The point of
Definition~\ref{adj.rel.def} is that sometimes, while a functor
$\lambda:\C' \to \C$ does not factor through $\gamma$, the
isomorphisms \eqref{adj.rel.iso} can still exist.
\end{exa}

\subsection{Admissible subcategories.}

Any left resp.\ right-reflexive functor $\phi:\C \to \C'$ is fully
faithful iff the adjunction map $a^\dg:\phi_\dg \circ \phi \to \id$
resp.\ $a_\dg:\id \to \phi^\dg \circ \phi$ is an isomorphism. If
$\phi$ is both left and right-reflexive, then $\phi \circ \phi^\dg$
is right-adjoint to $\phi \circ \phi_\dg$ and $a^\dg:\phi \circ
\phi^\dg \to \id$ is adjoint to $a_\dg:\id \to \phi \circ \phi_\dg$,
so that $\phi^\dg$ is fully faithful iff so is $\phi_\dg$. We will
say that a full subcategory $\C \subset \C'$ is {\em left}
resp.\ {\em right-admissible} iff the embedding functor $\C \to \C'$
is left resp.\ right-reflexive, and {\em admissible} if it is both
left and right-admissible.

\begin{exa}
For any object $c \in \C$ in a category $\C$, the embedding
\eqref{eps.c} is left resp.\ right-reflexive if and only if $c \in
\C$ is an initial resp.\ terminal object in $\C$, and in this case,
the full subcategory spanned by $c$ is full and left
resp.\ right-admissible. It is admissible if and only if $\C$ is
pointed.
\end{exa}

\begin{exa}\label{adj.exa}
For any category $\C$, the projections $\sigma,\tau:\Ar(\C) \to \C$
of Example~\ref{ar.exa} have a common section $\eta:\C \to \Ar(\C)$
sending $c \in \C$ to the identity map $\id:c \to c$. The functor
$\eta$ is fully faithful, its essential image is an admissible
subcategory, and the isomorphisms $\sigma \circ \eta \cong \id$,
$\id \cong \sigma \circ \eta$ define an adjunction between $\sigma$
and $\eta$, resp.\ $\eta$ and $\tau$.
\end{exa}

\begin{exa}\label{facto.exa}
For any class of morphisms $v$ in a category $\C$, denote by
$\Ar^v(\C) \subset \Ar(\C)$ the full subcategory spanned by arrows
in $v$. Then if the category $\C$ is equipped with a factorization
system $\langle L,R \rangle$, the subcategories $\Ar^L(\C),\Ar^R(\C)
\subset \Ar(\C)$ are right resp.\ left-admissible, with the adjoint
functors $\Ar(\C) \to \Ar^L(\C)$, $\Ar(\C) \to \Ar^R(\C)$ sending an
arrow $f:c \to c'$ to the arrows $l$ resp.\ $r$ in its canonical
decomposition \eqref{facto.eq}. Since the decompositions
\eqref{facto.eq} are unique, these adjoint functors actually fit
into a cartesian square
\begin{equation}\label{ar.facto.sq}
\begin{CD}
\Ar(\C) @>>> \Ar^R(\C)\\
@VVV @VV{\sigma}V\\
\Ar^L(\C) @>{\tau}>> \C.
\end{CD}
\end{equation}
Moreover, by abuse of notation, denote by $L$ resp.\ $R$ the classes
of morphisms in $\Ar(\C)$ that are pointwise in $L$ resp.\ $R$; then
the full subcategory $\Ar^R(\C)_L = \Ar^R(\C) \cap \Ar(\C)_L \subset
\Ar(\C)_L$ is also left-admissible, with the same adjoint functor,
and dually for $\Ar^L(\C)_R \subset \Ar(\C)_R$.
\end{exa}

\begin{exa}\label{cl.adj.exa}
For any right-admissible subcategory $\C' \subset \C$, and any
intermediate full subcategory $\C' \subset \C_0 \subset \C$, $\C'$
is tautogically right-admissible in $\C_0$ (with the same adjoint
functor). If a full subcategory $\C_0 \subset \C$ is left-closed in
the sense of Definition~\ref{cl.def}, then $\C' \cap \C_0 \subset
\C_0$ is again right-admissible, with the same adjoint
functor. Indeed, let $\gamma:\C \to \C' \subset \C$ be this adjoint
functor; then for any $c \in \C$, we have the adjunction map
$\gamma(c) \to c$, so $\gamma$ sends $\C_0$ into $\C_0$.
\end{exa}

\begin{defn}\label{proj.def}
A {\em left} resp.\ {\em right projector} on a category $\C$ is a
pair $\langle p,\alpha \rangle$ of a functor $p:\C \to \C$ and a map
$\alpha:\id \to p$ resp.\ $\alpha:p \to \id$ such that $p(\alpha)$
is invertible.
\end{defn}

\begin{lemma}\label{p.le}
\begin{enumerate}
\item For any category $\C$ equipped with a left resp.\ right
  projector $\langle p,\alpha\rangle$, the essential image $\C'
  \subset \C$ of the functor $p$ is left resp.\ right admissible,
  and $\alpha$ defines an adjunction between the embedding
  $\lambda:\C' \to \C$ and the projection $\rho:\C \to
  \C'$. Conversely, for any left resp.\ right admissible full
  embedding $\lambda:\C' \to \C$, with the adjoint functor $\rho:\C
  \to \C'$, the composition $p=\lambda \circ \rho$ equipped with the
  adjunction map is a left resp.\ right projector on $\C$.
\item Assume given categories $\C$, $\C'$ equipped with either left
  or right projectors $\langle p,\alpha \rangle$, $\langle \langle
  p',\alpha' \rangle$, and a functor $\gamma:\C \to \C'$ equipped
  with an isomorphism $\phi:\gamma \circ p \cong p' \circ
  \gamma$. Then $(\alpha')\gamma = \phi \circ \gamma(\alpha)$
  for left projectors, and $(\alpha')\gamma = \gamma(\alpha)
  \circ \phi^{-1}$ for right projectors.
\end{enumerate}
\end{lemma}

\proof{} Clear. \endproof

Note that the adjunction of Example~\ref{adj.exa} has the following
universal property. For any category $\C'$ and admissible full
subcategory $\C \subset \C'$, with the embedding functor $\phi:\C
\to \C'$ and its adjoints $\phi_\dg,\phi^\dg:\C' \to \C$, the
adjunction isomorphisms $\id \cong \phi \circ \phi^\dg$, $\phi_\dg
\circ \phi \cong \id$ together induce a map $a:\phi_\dg \to \phi_\dg
\circ \phi \circ \phi^\dg \to \phi^\dg$. We then have a natural
comparison functor
\begin{equation}\label{adj.fu}
\nu:\C' \to \Ar(\C)
\end{equation}
sending $c \in \C'$ to the arrow $a:\phi_\dg(c) \to \phi^\dg(c)$. By
Lemma~\ref{p.le}~\thetag{i}, the embedding $\C \subset \C'$ induces
a left and a right projector $l'$, $r'$ on the category $\C'$, and
similarly, the embedding $\C \subset \Ar(\C)$ induces projectors
$l$, $r$ on the category $\Ar(\C)$. The functor $\nu$ of
\eqref{adj.fu} comes equipped with isomorphisms $\nu \circ l' \cong
l \circ \nu$, $\nu \circ r' \cong r \circ \nu$, and by
Lemma~\ref{p.le}~\thetag{ii}, any functor $\nu':\C' \to \Ar(\C)$
equipped with such isomorphisms is canonically isomorphic to $\nu$.

\subsection{Cylinders and comma-categories.}\label{cyl.subs}

To reduce the study of arbitrary reflexive functors to the study of
left and right-admissible full embeddings, it is very convenient to
use the following version of the cylinder construction.

\begin{defn}\label{cyl.def}
The {\em cylinder} $\Cyl(\gamma)$ of a functor $\gamma:\C_0 \to
\C_1$ between categories $\C_0$, $\C_1$ is the category
\begin{equation}\label{cyl.def.eq}
\Cyl(\gamma) = \C_0 \copr_{F(\gamma)} \C_1
\end{equation}
of \eqref{cl.glue.eq} for the functor $F(\gamma) = \Hom \circ
(\gamma^o \times \id):\C_0^o \times \C_1 \to \Sets$. The {\em dual
  cylinder} $\Cyl^o(\gamma)$ is given by $\Cyl^o(\gamma) =
\Cyl(\gamma^o)^o$.
\end{defn}

\begin{exa}\label{cat.cyl.exa}
The cylinder $\Cyl(\id)$ of the identity map $\id:\ppt \to \ppt$ is
the single arrow category $[1]$. More generally, for any category
$\C$, the cylinder $\Cyl(\pi)$ of the tautological projection
$\pi:\C \to \ppt$ is the category $\C^>$, and the cylinder
$\Cyl(\id)$ of the identity functor $\id:\C \to \C$ is the product
$\C \times [1]$.
\end{exa}

Informally, the cylinder $\Cyl(\gamma)$ is obtained by taking the
disjoint union $\C_0 \copr \C_1$, and adding a morphism $f:c_0 \to
c_1$, $c_0 \in \C_0 \subset \Cyl(\gamma)$, $c_1 \in \C_1 \subset
\Cyl(\gamma)$ for any morphism $f:\gamma(c_0) \to c_1$. Cylinders
are functorial in the obvious way: a commutative square
\eqref{cat.dia} induces a functor $\Cyl(\gamma^0_{01}) \to
\Cyl(\gamma_1)$ that restricts to $\gamma^1_{01}$ resp.\ $\gamma_0$
on $\C_{01}$ resp.\ $\C_0$. In particular, for any $\gamma:\C_0 \to
\C_1$, the projections to the point induce the functor $\Cyl(\gamma)
\to [1]$ with fibers $\Cyl(\gamma)_l \cong \C_l$, $l = 0,1$. We also
have full embeddings $s:\C_0 \to \Cyl(\gamma)$, $t:\C_1 \to
\Cyl(\gamma)$, so that $\Cyl(\gamma)$ is a category under $\C_0
\copr \C_1$, and the subcategory $\C_1 \subset \Cyl(\gamma)$ is
left-admissible. Conversely, a functor $\C \to [1]$ with fibers
$\C_0$, $\C_1$ is the cylinder of a functor $\gamma:\C_0 \to \C_1$
if and only if $\C_1 \subset \C$ is left-admissible, and $\gamma$
can be recovered as the composition
\begin{equation}\label{cyl.deco}
\begin{CD}
\C_0 @>{s}>> \C @>{t_\dg}>> \C_1,
\end{CD}
\end{equation}
where $t_\dg:\C \to \C_1$ is left-adjoint to $t:\C_1 \to \C$. The
full subcategory $\C_0 \subset \C$ is right-admissible if and only
if $\gamma$ is right-reflexive, and in this case, the right-adjoint
functor $\gamma^\dg:\C_1 \to \C_0$ is given by $\gamma^\dg = s^\dg
\circ t$, where $s^\dg:\C \to \C_0$ is right-adjoint to $s:\C_0 \to
\C$. Passing to the opposite categories, we obtain the same
statements for dual cylinders and left-reflexive functors; in
particular, $\Cyl^o(\gamma)$ is a category under $\C_0 \copr \C_1$,
and $\gamma:\C_0 \to \C_1$ factors as
\begin{equation}\label{cyl.o.deco}
\begin{CD}
\C_0 @>{t}>> \Cyl^o(\gamma) @>{s^\dg}>> \C_1.
\end{CD}
\end{equation}
If $\lambda:\C_0 \to \C_1$ is left-adjoint to $\rho:\C_1 \to \C_0$,
then we have an equivalence $\eps:\Cyl(\lambda) \cong \Cyl^o(\rho)$
under $\C_0 \copr \C_1$, and conversely, such an equivalence defines
an adjunction between $\lambda$ and $\rho$. By functoriality of
cylinders, we have commutative squares
\begin{equation}\label{cyl.o.sq}
\begin{CD}
\C_0 @>{t}>> \C_0 \times [1]\\
@V{\gamma}VV @VVV\\
\C_1 @>{t}>> \Cyl(\gamma)
\end{CD}
\qquad\qquad
\begin{CD}
\C_0 @>{s}>> \C_0 \times [1]\\
@V{\gamma}VV @VVV\\
\C_1 @>{s}>> \Cyl^o(\gamma)
\end{CD}
\end{equation}
for any categories $\C_0$, $\C_1$ and functor $\gamma:\C_0 \to
\C_1$. Note that both squares are cocartesian.

\begin{lemma}\label{cl.cyl.le}
Assume given a full subcategory $\C' \subset \C$, with the embedding
functor $\gamma:\C' \to \C$. Then the functor $t_\dg:\Cyl(\gamma)
\to \C$ of \eqref{cyl.deco} is left-reflexive iff $\C' \subset \C$
is left-closed in the sense of Definition~\ref{cl.def}.
\end{lemma}

\proof{} Assume first that $\C$ has an initial object $o \in \C$,
and note that $\Cyl(\gamma)$ then has an initial object in two
cases: either $\C'$ is empty, or $o \in \C' \subset \C$. Then in
general, $t_\dg$ is left-reflexive iff for any $c \in \C$, the
comma-fiber $c \setminus_{t_\dg} \Cyl(\gamma)$ has an initial
object, and we have $c \setminus_{t_\dg} \Cyl(\gamma) \cong \Cyl(c
\setminus \gamma)$, where $c \setminus \gamma:c \setminus_\gamma \C'
\to c \setminus \C$ is the fully faithful embedding induced by
$\gamma$.
\endproof

The cylinder construction is in some sense dual to the
comma-category construction of Subsection~\ref{fun.subs}. Recall
that the right resp.\ left comma-categories $I \setminus_\gamma\C$
resp.\ $\C/_\gamma I$ of a functor $\gamma:\C \to I$ are the
categories of triples $\langle c,i,\alpha \rangle$, $c \in \C$, $i
\in I$, $\alpha$ a map from $i$ to $\gamma(c)$ resp.\ from
$\gamma(c)$ to $i$. We have $(\C /_\gamma I)^o \cong I^o
\setminus_{\gamma^o} \C^o$, and we have natural cartesian squares
\begin{equation}\label{rc.dia}
\begin{CD}
I \setminus_\gamma \C @>{\tau}>> \C\\
@VVV @VV{\gamma}V\\
\Ar(I) @>{\tau}>> I,
\end{CD}
\qquad\qquad
\begin{CD}
\C /_\gamma I @>{\sigma}>> \C\\
@VVV @VV{\gamma}V\\
\Ar(I) @>{\sigma}>> I,
\end{CD}
\end{equation}
where $\Ar(I)$ and $\sigma$ and $\tau$ in the bottom line are as in
Example~\ref{ar.exa}. The projection $\tau:I \setminus_\gamma \C \to
\C$ resp.\ $\sigma:\C /_\gamma I \to \C$ has a fully faithful right
resp.\ left-adjoint functor $\eta:\C \to I \setminus_\gamma \C$
resp.\ $\eta:\C \to \C /_\gamma I$ , sending $c \in \C$ to $\langle
c,\gamma(c),\id \rangle$, and $\gamma$
factors as
\begin{equation}\label{rc.dec}
\begin{CD}
\C @>{\eta}>> \C /_\gamma I @>{\tau}>> I,
\end{CD}
\qquad
\begin{CD}
\C @>{\eta}>> I \setminus_\gamma \C @>{\sigma}>> I,
\end{CD}
\end{equation}
where $\sigma$ resp.\ $\tau$ are the forgetful functors
\eqref{co.st.eq} whose fibers are the comma-fibers of the functor
$\gamma$. On the fiber $(I \setminus_\gamma \C)_i \cong i
\setminus_\gamma \C$, the projection $\tau$ restricts to the functor
$\tau(i)$ of \eqref{p.i}, and $\eta$ restricts to the full embedding
$\C_i \subset i \setminus_\gamma \C$, and dually for the left
comma-fibers. Note that by \eqref{rc.dec}, $\eta$ is a functor over
$I$, and its left-adjoint $\tau$ is a lax functor over $I$.

\begin{exa}
If we take $\C = I$ and $\gamma = \id$, as in Example~\ref{ar.exa},
then $\C /_{\id} \C \cong \C \setminus_{\id} \C \cong \Ar(\C)$, and
$\eta$ of \eqref{rc.dec} is as in Example~\ref{adj.exa}.
\end{exa}

To see the duality between cylinders and comma-categories, one can
use the functor categories of Subsection~\ref{cat.fun.subs}. In
particular, for any functor $\gamma:I_0 \to I_1$ between essentially
small categories, and any category $\C$, we have natural
equivalences
\begin{equation}\label{cyl.com}
\begin{aligned}
\Fun(\Cyl(\gamma),\C) &\cong \Fun(I_0,\C) \setminus_{\gamma^*} \Fun(I_1,\C),\\
\Fun(\Cyl^o(\gamma),\C) &\cong \Fun(I_1,\C) /_{\gamma^*}
\Fun(I_0,\C),
\end{aligned}
\end{equation}
and the decompositions \eqref{rc.dec} corresponds to the
decompositions \eqref{cyl.deco} (by \eqref{fun.dia}, \eqref{cyl.com}
amounts to observing that both commutative squares in
\eqref{cyl.o.sq} are cocartesian). Another useful observation is
that for any functor $\gamma:\C_0 \to \C_1$, we have natural
identifications
\begin{equation}\label{cyl.co.sec}
  \C_0 /_\gamma \C_1 \cong \Sec([1],\Cyl(\gamma)), \qquad
  \C_1 \setminus_\gamma \C_0 \cong \Sec([1],\Cyl^o(\gamma)),
\end{equation}
where $\Sec([1],-)$ are the categories of sections of the natural
projections $\Cyl(\gamma),\Cyl^o(\gamma) \to [1]$. If we are given
two functors $\lambda:\C_0 \to \C_1$, $\rho:\C_1 \to \C_0$, then by
\eqref{cyl.co.sec}, giving an adjunction between $\lambda$ and
$\rho$ is the same thing as giving an equivalence $\eps:\C_0
/_\lambda \C_1 \cong \C_0 \setminus_\rho \C_1$ over $\C_0 \times
\C_1$. This immediately implies that a functor $\gamma:\C_0 \to
\C_1$ is left resp.\ right-reflexive if and only if $\sigma$
resp.\ $\tau$ in \eqref{rc.dec} is left resp.\ right-reflexive, with
a fully faithful adjoint $\sigma_\dg$ resp.\ $\tau_\dg$, and $\tau
\circ \sigma_\dg$ resp.\ $\sigma \circ \tau_\dg$ is then the adjoint
to $\gamma$.

More generally, assume that categories $\C_0$ and $\C_1$ are
equipped with functors $\pi_l:\C_l \to I$, $l=0,1$ to some category
$I$. Then for any $\pi:\C \to I$, we define the {\em relative arrow
  category} $\Ar(\C|I)$ by the cartesian square
\begin{equation}\label{rel.ar.sq}
\begin{CD}
\Ar(\C|I) @>>> \Ar(\C)\\
@VVV @VV{\Ar(\pi)}V\\
I @>{\eta}>> \Ar(I),
\end{CD}
\end{equation}
and for any functor $\gamma:\C_0 \to \C_1$ over $I$, we define the
{\em relative comma-categories} $\C_0 /_{\gamma,I} \C_1$, $\C_1
\setminus_{\gamma,I} \C_0$ by cartesian squares
\begin{equation}\label{rc.I.dia}
\begin{CD}
\C_1 \setminus_{\gamma,I} \C_0 @>{\tau}>> \C_0\\
@VVV @VV{\gamma}V\\
\Ar(\C_1|I) @>{\tau}>> \C_1,
\end{CD}
\qquad\qquad
\begin{CD}
\C_0 /_{\gamma,I} \C_1 @>{\sigma}>> \C_0\\
@VVV @VV{\gamma}V\\
\Ar(\C_1|I) @>{\sigma}>> \C_1,
\end{CD}
\end{equation}
or equivalently,
\begin{equation}\label{rel.lc.sq}
\begin{CD}
  \C_1 \setminus_{\gamma,I} \C_0 @>>> \C_1 \setminus_\gamma \C_0\\
  @VVV @VV{\pi_1 \setminus \pi_0}V\\
  I @>{\eta}>> \Ar(I),
\end{CD}
\qquad\qquad
\begin{CD}
  \C_0 /_{\gamma,I} \C_1 @>>> \C_0 /_\gamma \C_1\\
  @VVV @VV{\pi_0/\pi_1}V\\
  I @>{\eta}>> \Ar(I).
\end{CD}
\end{equation}
As in the absolute case, $\gamma$ admits natural factorizations
\begin{equation}\label{rel.rlc.facto}
\begin{CD}
\C_0 @>{\eta}>> \C_0 /_{\gamma,I} \C_1 @>{\tau}>> \C_1,
\end{CD}
\quad
\begin{CD}
\C_0 @>{\eta}>> \C_1 \setminus_{\gamma,I} \C_0 @>{\sigma}>> \C_1
\end{CD}
\end{equation}
induced by \eqref{rc.dec}, and $\eta$ is right resp.\ left-adjoint
to the forgetful functor $\sigma:\C_0 /_{\gamma,I} \C_1 \to \C_0$
resp.\ $\tau:\C_1 \setminus_{\gamma,I} \C_0 \to \C_0$ over $I$. Then
if we are given two functors $\lambda:\C_0 \to \C_1$, $\rho:\C_1 \to
\C_0$ over $I$, giving an adjunction between $\lambda$ and $\rho$
over $I$ is the same thing as giving an equivalence $\eps:\C_0
/_{\lambda,I} \C_1 \cong \C_0 \setminus_{\rho,I} \C_1$ over $\C_0 \times
\C_1$, and this immediately implies that a functor $\gamma:\C_0 \to
\C_1$ is left resp.\ right-reflexive over $I$ if and only if
$\sigma$ resp.\ $\tau$ in \eqref{rel.rlc.facto} is left
resp.\ right-reflexive, with a fully faithful adjoint $\sigma_\dg$
resp.\ $\tau_\dg$, and $\tau \circ \sigma_\dg$ resp.\ $\sigma \circ
\tau_\dg$ is then the adjoint to $\gamma$ (recall that since
$\sigma_\dg$ and $\tau_\dg$ are fully faithful, both are
automatically adjoint over $I$).

\begin{exa}\label{ar.ar.exa}
For any category $\C$, if we consider $\Ar(\C)$ as a category over
$\C$ via the projection $\tau$, then $\Ar(\Ar(\C)|\C)$ is the
category of composable pairs $\langle f,g \rangle$ of arrows in
$\C$, with $\sigma,\tau:\Ar(\Ar(\C)|\C) \to \Ar(\C)$ sending
$\langle f,g \rangle$ to $f \circ g$ resp.\ $f$, and we have a
cartesian square
\begin{equation}\label{ar.ar.sq}
\begin{CD}
\Ar(\Ar(\C)|\C) @>{\tau}>> \Ar(\C)\\
@V{\beta}VV @VV{\sigma}V\\
\Ar(\C) @>{\tau}>> \C,
\end{CD}
\end{equation}
where $\beta$ sends $\langle f,g\rangle$ to $g$,
\end{exa}

\begin{lemma}\label{adj.co.le}
Assume given a functor $\lambda:\C_0 \to \C_1$ equipped with a
section $\rho:\C_1 \to \C_0$. Then the isomorphism $\lambda \circ
\rho \cong \id$ defines an adjunction between $\lambda$ and $\rho$
if and only if $\sigma:\C_0 \setminus_{\rho,\C_1} \C_1 \to \C_0$ is
an equivalence.
\end{lemma}

\proof{} The ``if'' part immediately follows from
\eqref{rel.rlc.facto}. For the ``only if'' part, if we do have
an adjunction, then $\C_0 \setminus_{\rho,\C_1} \C_1 \cong \C_0
/_{\lambda,\C_1} \C_1 \cong \C_0$.
\endproof

\begin{corr}\label{adm.sq.corr}
Assume given a cartesian square \eqref{cat.dia} and left
resp.\ right-admissible full subcategories $\C' \subset \C$, $\C'_l
\subset \C_l$, $l=0,1$, with adjoint functors $\lambda:\C \to \C'$,
$\lambda_l:\C_l \to \C'_l$, $l=0,1$, such that $\gamma_l(\C'_l)
\subset \C'$, $l = 0,1$, and \eqref{bc.eq} provides isomorphisms
$\gamma_l \circ \lambda_l \cong \lambda \circ \gamma_l$,
$l=0,1$. Then the full subcategory $\C'_{01} = \C'_0 \times_{\C'}
\C'_l \subset \C_{01}$ is left resp.\ right-admissible, with left
resp.\ right-adjoint $\lambda_{01}:\C_{01} \to \C'_{01}$, and
$\gamma^l_{01} \circ \lambda_{01} \cong \lambda_l \circ
\gamma^l_{01}$, $l=0,1$.
\end{corr}

\proof{} It suffices to consider the left-admissible case (in the
right-admissible case, take the opposite categories). Let $\rho$,
$\rho_0$, $\rho_1$, $\rho_{01}$ be the embedding functors.  Then to
construct $\lambda_{01}$ and an isomorphism $\lambda_{01} \circ
\rho_{01} \cong \id$, take the product of $\lambda_0$ and
$\lambda_1$ over $\lambda$, and to check that the isomorphism
satifies the condition of Lemma~\ref{adj.co.le}, note that the
relative comma-categories for the embedding functors also form a
cartesian square.
\endproof

\subsection{Limits and colimits.}\label{lim.subs}

Basic examples of adjoint functors are given by limits and colimits,
and more generally, left and right Kan extensions. Here is a brief
summary.

\begin{defn}
A {\em cone} of a functor $E:I \to \E$ is a functor $E_>:I^> \to \E$
equipped with an isomorphism $E_> \circ s \cong E$, where $s:I \to
I^>$ is the natural embedding. The {\em vertex} of a cone $E_>$ is
its value $E_>(o) \in \E$ at the terminal object $o \in I^>$. A
cone $E_>$ is {\em universal} if for any other cone $E_>'$, there
exists a unique map of cones $E_> \to E_>'$. The functor $E$ {\em
  has a colimit} if it admits a universal cone $E_>$, and the
colimit is given by its vertex, $\colim_IE = E_>(o)$.
\end{defn}

\begin{exa}\label{cocart.sq.exa}
If $I = \V = \{0,1\}^<$ is the category of \eqref{V.eq}, then we
have $\V^> \cong [1]^2$, and a cone of a functor $\V \to \E$ is a
commutative square in $\E$. If the cone is universal, the square is
called {\em cocartesian}.
\end{exa}

\begin{exa}\label{kar.lim.exa}
Let $P^=$ be the category with two objects $\ppt$, $o$, and
morphisms generated by $a:o \to \ppt$ and $b:\ppt \to o$ subject to
the relation $b \circ a = \id$. Then the full subcategory spanned by
$\ppt$ is the category $P$ of Example~\ref{P.exa}, with the
embedding sending $p$ to $a \circ b$. Removing the map $a$ gives a
dense embedding $P^> \subset P^=$ sending $o$ to $o$, and this is a
universal cone for the full embedding $P \to P^=$. Moreover, a cone
$E_>$ of a functor $E:P \to \E$ to some category $\E$ is universal
if and only if it extends to $P^= \supset P^>$ (and if this
happens, $\colim_E$ is the image of the projector $E(p)$).
\end{exa}

\begin{exa}\label{kron.exa}
Define the {\em Kronecker category} $\II$ as the category with two
objects $0,1 \in \II$ and two non-identity maps $s,t:0 \to 1$. Then
a colimit $\colim_{\II}E$ of a functor $E:\II \to \E$ is known as
the {\em coequalizer} of the two maps $E(s),E(t):E(0) \to E(1)$; it
is given by a map $f:E(1) \to \colim_\II E$ such that $f \circ E(s) =
f \circ E(t)$, and universal with this property.
\end{exa}

By virtue of its universal property, colimit is functorial with
respect to $E$, in that for any $E,E':I \to \E$ such that
$\colim_IE$ and $\colim_IE'$ exist, a map $E \to E'$ induces a map
$\colim_IE \to \colim_IE'$. It is also functorial with respect to
$I$, in that for any functor $\gamma:I' \to I$, we have a natural
map
\begin{equation}\label{colim.eq}
\colim_{I'}E \circ \gamma \to \colim_IE,
\end{equation}
again provided that both sides exist. A category $\E$ {\em has
  colimits} or {\em is cocomplete} if $\colim_I E$ exists for any
essentially small $I$ and functor $E:I \to \E$, and $\E$ is {\em
  $\kappa$-cocomplete} for some infinite cardinal $\kappa$ if this
happens for any $I$ with $\|I\| < \kappa$. If $\kappa$ is the smallest
infinite cardinal, $\E$ is {\em finitely cocomplete}. For any
infinite cardinal $\kappa$, a category $\E$ is $\kappa$-cocomplete
if it is finitely cocomplete and has $\kappa'$-indexed coproducts
for any $\kappa' < \kappa$. For any essentially small category $I$
and cocomplete resp.\ $\kappa$-cocomplete category $\E$, the functor
category $\E^I$ is cocomplete resp.\ $\kappa$-cocomplete.

\begin{remark}\label{coeq.rem}
To insure that a category $\E$ is finitely cocomplete, it suffices
to require that $\E$ has all coequalizers of Example~\ref{kron.exa}
and all finite coproducts; for the latter, it further suffices to
require colimits over discrete categories $\emptyset$ and $\{0,1\}$.
\end{remark}

\begin{exa}\label{sets.S.exa}
The category $\Sets$ is cocomplete. For any infinite cardinal
$\kappa$, let $\Sets_\kappa \subset \Sets$ be the full subcategory
of sets $S$ such that $|S| < \kappa$. Then $\Sets_\kappa$ is
finitely cocomplete, and its $\kappa$-cocomplete if and only if for
any map $S' \to S$ of sets such that $|S| < \kappa$ and $|S'_s| <
\kappa$ for any $s \in S$, we also have $|S'| < \kappa$. By
definition, this means that the cardinal $\kappa$ is {\em
  regular}. Not all infinite cardinals are regular, but regular ones
are cofinal, in the sense that for any infinite $\kappa$, there
exists a regular $\kappa' > \kappa$ (in fact, it suffices to take
the {\em succesor cardinal} $\kappa^+$ defined as the smallest
cardinal $\kappa^+ > \kappa$). Technically, one could also allow
$\kappa$ to be finite, and then $\Sets_\kappa$ is
$\kappa$-cocomplete iff $\kappa=0,1,2$; however, we will not need
this.
\end{exa}

By Example~\ref{cat.cyl.exa}, the category $I^>$ is the cylinder of
the tautological projection $\gamma:I \to \ppt$. More generally,
assume given a functor $\gamma:I \to I'$, and consider its cylinder
$\Cyl(\gamma)$, with the embeddings $s:I \to \Cyl(\gamma)$, $t:I'
\to \Cyl(\gamma)$.

\begin{defn}\label{kan.def}
A {\em relative cone} of a functor $E:I \to \E$ with respect to a
functor $\gamma:I \to I'$ is a functor $E_>:\Cyl(\gamma) \to \E$
equipped with an isomorphism $E_> \circ s \cong E$. A relative cone
$E_>$ is {\em universal} if for any other relative cone $E_>'$,
there exists a unique map of relative cones $E_> \to E_>'$. The
functor $E$ {\em has a left Kan extension} with respect to $\gamma$
if it admits a universal relative cone $E_>$ with respect to
$\gamma$, and the left Kan extension $\gamma_!E$ is given by
$\gamma_!E = E_> \circ t:I' \to \E$.
\end{defn}

Explicitly, a cone $E_>$ of a functor $E:I \to \E$ is given by its
vertex $e = E_>(o) \in \E$ and a morphism $E \to \eps(e) \circ
\gamma$, where $\eps(e) \circ \gamma:I \to \ppt \to \E$ is the
constant functor with value $e$. Analogously, a relative cone with
respect to some $\gamma:I \to I'$ is a pair of a functor $E':I' \to
\E$ and a morphism $E \to E' \circ \gamma$. For any $\gamma:I \to
I'$ and $E:I \to \E$, the left Kan extension is given by
\begin{equation}\label{kan.eq}
\gamma_!E(i) = \colim_{I /_\gamma i}E \circ \sigma(i), \qquad i \in I',
\end{equation}
where $\sigma(i)$ are the forgetful functors \eqref{p.i}, and
$\gamma_!E$ exists if so do the colimits in the right-hand side. For
any morphism $f:i \to i'$ in the category $I'$, the correspoding map
$\gamma_!E(f):\gamma_!E(i) \to \gamma_!E(i')$ is the map
\eqref{colim.eq} for the functor $f_!:I /_\gamma i \to I /_\gamma
i'$ of \eqref{f.i}.

Now, for any functor $\gamma:I_0 \to I_1$ between essentially small
categories, and any category $\E$, we have the pullback functor
$\gamma^*:\Fun(I_1,\E) \to \Fun(I_0,\E)$. Then if the left Kan
extension $\gamma_!E$ exists for any $E \in \Fun(I_0,\E)$, it is
functorial in $E$, and we obtain a functor $\gamma_!:\Fun(I_0,\E)
\to \Fun(I_1,\E)$ left-adjoint to $\gamma^*$. Conversely, if
$\gamma^*$ admits a left-adjoint functor $\gamma_!$, then for any $E
\in \Fun(I_0,\E)$, $\gamma_!E$ with the adjunction map $E \to
\gamma^*\gamma_!E$ is the left Kan extension in the sense of
Definition~\ref{kan.def}. If we have another functor $\gamma':I_1
\to I_2$ to an essentially small $I_2$ such that both $\gamma_!$ and
$\gamma'_!$ exists, then $\gamma'_! \circ \gamma_!  \cong (\gamma'
\circ \gamma)_!$ by adjunction, and if $\gamma:I \to \ppt$ is the
taulogical projection, then $\gamma_! \cong \colim_I$. A commutative
diagram
\begin{equation}\label{bc.sq}
\begin{CD}
I_1' @>{\gamma'}>> I'_0\\
@V{\nu_1}VV @VV{\nu_0}V\\
I_1 @>{\gamma}>> I_0
\end{CD}
\end{equation}
of essentially small categories induces a commutative diagram of
functor categories $\Fun(-,\E)$, and if $\gamma_!$ and $\gamma'_!$
exist, we have a functorial map
\begin{equation}\label{ga.bc.eq}
a(\gamma):\gamma'_! \circ \nu_1^* \to \nu_0^* \circ \gamma_!
\end{equation}
given by the base change map \eqref{bc.eq}.

\begin{exa}\label{kan.adj.exa}
If a functor $\gamma:I_0 \to I_1$ between essentially small
categories admits a right-adjoint $\gamma_\dg:I_1 \to I_0$, then by
Example~\ref{fun.adj.exa}, for any target category $\E$, the left
Kan extension functor $\gamma_!:\Fun(I_0,\E) \to \Fun(I_1,\E)$ is
given by $\gamma_!  \cong \gamma_\dg^*$.
\end{exa}

\begin{exa}\label{kan.coco.exa}
By \eqref{kan.eq}, for any cocomplete target category $\E$, the left
Kan extension functor $\gamma_!:\Fun(I_0,\E) \to \Fun(I_1,\E)$
exists for any functor $\gamma:I_0 \to I_1$ between essentially
small categories, and the same is true if $\E$ is
$\kappa$-cocomplete for some $\kappa$, and $\|I_0\|,\|I_1\| < \kappa$.
\end{exa}

To compute left Kan extensions, it is often convenient to replace
the comma-fibers in \eqref{kan.eq} with smaller subcategories. For
instance, by Example~\ref{kan.adj.exa}, for any left-admissible full
subcategory $I' \subset I$ with the embedding functor $\gamma:I' \to
I$ and its left-adjoint $\gamma^\dg:I \to I'$, and any functor $E:I
\to \E$, we have
\begin{equation}\label{adm.eq}
\colim_{I'}\gamma^*E \cong \colim_{I'}\gamma^\dg_!E \cong \colim_IE,
\end{equation}
so that the left-hand side exists if and only if so does the
right-hand side, and the map \eqref{colim.eq} is an isomorphism.
Here is a more general example of such a situation.

\begin{exa}\label{cofin.exa}
Note that an essentially small category $I$ is connected if and only
if for any category $\E$, the tautological embedding $\E \to \Fun(I,\E)$
is fully faithful. Say that a functor $\gamma:I' \to I$ between
essentially small categories is {\em $0$-cofinal} if for any $i \in
I$, the right comma-fiber $i \setminus_\gamma I'$ is connected. Then
for any functor $E:I \to \E$, the source of the map \eqref{colim.eq}
exists if and only if so does its target, and the map is an
isomorphism.
\end{exa}

\begin{exa}\label{kan.facto.exa}\mbox{}\!\!\!
The prototypical example of a Kan extension that gives rise to the
name occurs when $\gamma:I \to I'$ is the embedding functor of a
full subcategory $I \subset I'$. In this case, for any $i \in I
\subset I'$, the left comma-fiber $I /_\gamma i$ has the terminal
object $\langle i,\id \rangle$, so if the left Kan extension $E' =
\gamma_!E$ exists for some $E:I \to \E$, then the map $E \to
\gamma^*E'$ is an isomorphism by \eqref{kan.eq} and
\eqref{adm.eq}. We then have a factorization
\begin{equation}\label{kan.facto}
\begin{CD}
I @>{\gamma}>> I' @>{E'}>> \E
\end{CD}
\end{equation}
of the functor $E$, and $E'$ indeed extends $E$ to a larger category
$I'$. If $I \subset I'$ is right-admissible, then $E' = \gamma_!E
\cong \gamma_\dg^*E$, hence also the factorization \eqref{kan.facto}
exist for any $E$ by Example~\ref{kan.adj.exa}, and one can
characterize it as follows: this is the unique factorization $E
\cong E' \circ \gamma$ such that $E'$ inverts all maps in the class
$\gamma_\dg^*(\Iso)$. Note that $\gamma_\dg^*(\Iso)$ is the minimal
saturated class that contains all the adjunctions maps
$\gamma(\gamma_\dg(i')) \to i'$, $i' \in I'$. More generally, if $I$
and $I'$ are essentially small, then for any closed class $v$ of
morphisms in $I$, we have an equivalence of categories
\begin{equation}\label{kan.adm.eq}
\gamma_! \cong \gamma_\dg^*:\Fun^v(I,\E) \to
\Fun^{\gamma_\dg^*v}(I',\E),
\end{equation}
with the inverse equivalence given by $\gamma^*$. Moreover, while
the Kan extension $\gamma_{\dg !}E$ might not exist for a general
functor $E:I' \to \E$, it exists for any $E \in
\Fun^{\gamma_\dg^*v}(I',\E)$, and $\gamma_{\dg !}E \cong \gamma^*E$.
\end{exa}

\begin{exa}\label{yo.sets.exa}
As mentioned in Example~\ref{sets.S.exa}, the category $\Sets$ is
cocomplete, and so is the functor category $I^o\Sets$ for any
essentially small $I$. Thus for any two essentially small categories
$I_0$, $I_1$, a functor $x:I_0 \to I_1^o\Sets$ has a left Kan
extension $X = \Y_!(x)$ with respect to the Yoneda embedding
\eqref{yo.yo.eq}. An abstract functor $X:I^o_0\Sets \to I_1^o\Sets$
is of this form if and only if the adjunction map $\Y_!\Y^*X \to X$
is an isomorphism, or equivalently, if $X$ preserves all colimits,
or equivalently, if it is right-reflexive.
\end{exa}

\begin{lemma}\label{cl.copr.le}
Assume given categories $\C_{01}$, $\C_0$, $\C_1$, a left-closed
full embedding $\gamma_0:\C_{01} \to \C_0$, and a functor
$\gamma_1:\C_{01} \to \C_1$ such that the Kan extension
$\gamma^o_{1!}X:\C_1^o \to \Sets$ exists for any $X:\C^o_{01} \to
\Sets$. Then there exists a commutative square of categories
\begin{equation}\label{st.sq}
\begin{CD}
\C_{01} @>{\gamma_0}>> \C_0\\
@V{\gamma_1}VV @VVV\\
\C_1 @>>> \C
\end{CD}
\end{equation}
that is both cartesian and cocartesian in the sense of
Subsection~\ref{dia.subs}.
\end{lemma}

\proof{} For any categories $\C_0$, $\C_1$ and functor $F:\C_0^o
\times \C_1 \to \Sets$, with the corresponding category $\C_0
\copr_F \C_1$ of \eqref{cl.glue.eq}, giving a functor $\C_0 \copr_F
\C_1 \to \E$ to some category $\E$ amounts to giving functors
$E_l:\C_l \to \E$, $l=0,1$, and a map $F \to \Hom \circ (E_0^o
\times E_1)$, where $\Hom:\E^o \times \E \to \Sets$ is the
$\Hom$-pairing \eqref{hom.pair.eq}. Then under the assumptions of
the Lemma, let $\C_0' = \C_0 \ssetminus \C_{01}$, represent $\C_0
\cong \C_{01} \copr_{F_0} \C_0'$ for some $F_0:\C_{01}^o \times
\C_0' \to \Sets$, and note that the universal property of Kan
extensions immediately implies that the Kan extension
$F'=(\gamma_1^o \times \id)_!F_0$ exists, and taking $\C = \C_1
\copr_{F'} \C_0'$ does the job.
\endproof

\begin{exa}\label{st.0.exa}
Since the category $\Sets$ is cocomplete, the assumption on
$\gamma_1$ in Lemma~\ref{cl.copr.le} is always satisfied when
$\C_{01}$ is essentially small, but this is not the only case. For
example, by \eqref{kan.eq}, it is always satisfied when
$\gamma_1:\C_{01} \to \C_1$ is the left-closed embedding. In this
case, we also have $\C_1 \cong \C_{01} \copr_{F_1} \C_1'$ for $\C_1'
= \C_1 \ssetminus \C_{01}$ and some $F_1:\C_{01}^o \times \C_1' \to
\Sets$, and \eqref{kan.eq} then shows that $\C$ in \eqref{st.sq} is
given by $\C \cong \C_{01} \copr_{F_0 \copr F_1} (\C_0' \copr
\C_1')$. This construction is inverse to Example~\ref{cl.exa}.
\end{exa}

Dually, the {\em limit} $\lim_I E$ of a functor $E$ is given by
$\lim_I E = (\colim_{I^o}E^o)^o$, and the {\em right Kan extension}
$\gamma_*E$ is $\gamma_*E = (\gamma^o_!E^o)^o$. A category $\E$ is
{\em complete} or {\em has limits} if $\E^o$ is cocomplete, and {\em
  $\kappa$-complete} if $\E^o$ is $\kappa$-cocomplete. Right Kan
extensions are right-adjoint to pullbacks, they can be expressed by
the dual version of \eqref{kan.eq} with limits instead of colimits
and right comma-fibers instead of left ones, and the rest of the
story including the base change map \eqref{ga.bc.eq},
Example~\ref{kan.adj.exa} and Example~\ref{kan.coco.exa} has an
obvious dual counterpart.

\begin{exa}\label{sets.SS.exa}
For any two infinite sets $S$, $S'$, say that $|S| \ll |S'|$ iff
$|S^S| \leq |S'|$ (or equivalently, $|\{0,1\}^S| \leq |S'|$, where
$\{0,1\}$ is the set with two elements, so that $\{0,1\}^S$ is the
set of all subsets in $S$). Then for any regular cardinal $\kappa$,
the category $\Sets_\kappa$ of Example~\ref{sets.S.exa} is
$\kappa'$-complete iff $\kappa' \ll \kappa$, and the same hold for
$I^o\Sets_\kappa$ for any essentially small $I$.
\end{exa}

\begin{defn}\label{cart.cl.def}
A category $\C$ that has finite products is {\em cartesian-closed}
if for any $c \in \C$, the functor $c \times -:\C \to \C$ is
right-reflexive.
\end{defn}

\begin{exa}\label{cc.sets.exa}
For any essentially small category $I$, the cocomplete category
$I^o\Sets$ of Example~\ref{yo.sets.exa} is also complete. Moreover,
since $X \times -$ for any $X \in I^o\Sets$ obviously preserves
colimits, it is right-reflexive, so that $I^o\Sets$ is
cartesian-closed. Explicitly, the functor $I^o\Sets \to I^o\Sets$
right-adjoint to $X \times -$ sends $Y \in I^o\Sets$ to $\Hhom(X,Y)
\in I^o\Sets$ given by
\begin{equation}\label{cc.XY}
\Hhom(X,Y)(i) = I^o\Sets(X \times \Y(i),Y).
\end{equation}
A useful particular case is $X = \Y(i)$, $i \in I$. In this case,
the dual counterpart of \eqref{kan.eq} provides functorial
identifications
\begin{equation}\label{cc.Yi}
\Hhom(\Y(i),Y) \cong \delta^o_*Y |_{I \times \{i\}}, \qquad i \in I, Y \in
I^o\Sets,
\end{equation}
where $\delta:I \to I \times I$ is the diagonal embedding.
\end{exa}

A category $\I$ is {\em filtered} if any functor $I \to \I$ from a
finite $I$ admits a cone. As in Remark~\ref{coeq.rem}, it actually
suffices to require it for $I=\emptyset,\{0,1\},\II$. A category
$\E$ {\em has filtered colimits} if $\colim_\I E$ exists for any
filtered $\I$ and functor $E:\I \to \E$, and an object $e \in \E$ of
such a category is {\em compact} if $\E(e,-)$ preserves all these
filtered colimits. Filtered colimits in $\Sets$ commute with finite
limits; because of this, many ``algebraically defined'' categories
such as $\Cat$ also have filtered colimits.

\begin{remark}
The fact that filtered colimits in $\Sets$ commute with finite
limits is very standard and fundamental, and although it is not
difficult, it does require a proof. The converse -- a category $I$
such that $\colim_I$ preserves finite limits is filtered -- is also
true and in fact easier. We do not go into any detail here since in
Section~\ref{enh.large.sec}, we will develop the whole story in the
enhanced setting (for an unenhanced version of
Section~\ref{enh.large.sec}, see the overview paper \cite{ka.acc}).
\end{remark}

\begin{remark}
``Filtered'' is actually a mistranslation of the original {\em
    filtrant} of \cite{sga4}; it should be ``filtering''. It is very
  tempting to try to restore the original terminology, but it is
  probably too late.
\end{remark}

\subsection{Additive categories.}\label{add.subs}

For any two objects $A,B \in \A$ in a pointed category $\A$ with the
initial-terminal object $0 \in \A$, we have the unique map $A \to B$
that factors through $0$; it is denoted by $0 \in \A(A,B)$. If $\A$
has finite products and coproducts, then we also have the natural
map
\begin{equation}\label{preadd.eq}
  (\id_A \times 0) \copr (0 \times \id_B):A \copr B \to A \times B,
\end{equation}
and $\A$ is {\em pre-additive} if this map is an isomorphism for any
$A,B \in \A$. Equivalently, a category is preadditive iff it has
finite products and coproducts, and they commute. Once this holds,
$\A(A,B)$ acquires a natural structure of a commutative monoid with
unity element $0$, and $\A$ is {\em additive} if for any $A,B \in
\A$, the monoid $\A(A,B)$ is a group. For an additive category $\A$,
one usually denotes $A \copr B \cong A \times B$ by $A \oplus B$,
and calls it the {\em sum} of the objects $A,B \in \A$.

For any additive category $\A$, the category $C_\idot(\A)$ of chain
complexes in $\A$ is also additive. As usual, for any map $f:A_\idot
\to B_\idot$ between complexes $A_\idot$, $B_\idot$ with
differentials $d_A$, $d_B$, the {\em cone} $\Cone(f)$ is the complex
$B_\idot \oplus A_\idot[1]$ with the differential $d_A + d_B + f$,
and we have a natural termwise-split embedding $B_\idot \to
\Cone(f)$. A complex $A_\idot$ is {\em contractible} if
$\id_A:A_\idot \to A_\idot$ factors through the embedding $A_\idot
\to \Cone(\id_A)$. A morphism $f:A_\idot \to B_\idot$ is a {\em
  chain-homotopy equivalence} if $\Cone(f)$ is contractible. Being
contractible is closed under sums and retracts, and for any two maps
$f:A_\idot \to B_\idot$, $g:B_\idot \to C_\idot$, we have
\begin{equation}\label{cone.oct}
  \Cone(f) \oplus \Cone(g) \cong \Cone(g \circ f) \oplus \Cone(\id_B),
\end{equation}
so that the class of all chain-homotopy equivalences is closed and
saturated in the sense of Definition~\ref{sat.def}. Two morphisms
$f_0,f_1:A_\idot \to B_\idot$ are {\em chain-homotopic} if $f_0
\oplus f_1:A_\idot \oplus A_\idot \to B_\idot$ factors through the
embedding $A_\idot \oplus A_\idot \to \Cone(\delta_A)$, where
$\delta_A:A_\idot \to A_\idot \oplus A_\idot$ is the diagonal
map. Being chain-homotopic is an equivalence relation. The {\em
  chain-homotopy category} $\Ho(\A)$ has the same objects as
$C_\idot(\A)$, and chain-homotopy classes of maps between them as
morphisms.

\section{Fibrations and cofibrations.}\label{fib.sec}

\subsection{The Grothendieck construction.}\label{fib.subs}

Recall that the right resp.\ left comma-categories $I
\setminus_\gamma\C$ resp.\ $\C/_\gamma I$ of a functor $\gamma:\C
\to I$ are defined by the cartesian squares \eqref{rc.dia}, and
$\gamma$ then has natural decompositions \eqref{rc.dec}.

\begin{defn}\label{fib.def}
A functor $\gamma:\C \to I$ is a (Grothen\-dieck) {\em fibration} if
the functor $\eta:\C \to I \setminus_\gamma \C$ of \eqref{rc.dec} is
right-reflexive over $I$. A functor $\gamma$ is a {\em cofibration}
iff $\gamma^o$ is a fibration, and it is a {\em bifibration} if it
is both a fibration and a cofibration. A fibration or a cofibration
$\C \to I$ is {\em small} if it is small as a functor (or
equivalently, has essentially small fibers $\C_i$, $i \in I$).
\end{defn}

\begin{remark}\label{cofib.rem}
One can also describe cofibrations in terms of the left
comma-category $\C /_\gamma I$ and the second decomposition
\eqref{rc.dec}.
\end{remark}

\begin{exa}
The functor $\sigma$ resp.\ $\tau$ in \eqref{rc.dec} is a fibration
resp.\ cofibration.
\end{exa}

\begin{lemma}\label{fib.ar.le}
A functor $\pi:\C \to I$ is a fibration if and only if the functor
\begin{equation}\label{pi.co.eq}
\begin{CD}
\Ar(\C) \cong \C \setminus_{\id} \C @>{\pi \setminus \id}>> I
\setminus_\pi \C
\end{CD}
\end{equation}
admits a fully faithful right-adjoint $(\pi \setminus \id)_\dg$. If
$\pi$ is a fibration, then the relative arrow category $\Ar(\C|I)
\subset \Ar(\C)$ is right-admissible, and we have a diagram
\begin{equation}\label{fib.facto.sq}
\begin{CD}
\Ar(\C|I) @>{\nu}>> \Ar(\C) @>{\nu_\dg}>> \Ar(\C|I)\\
@V{\tau}VV @V{\pi \setminus \id}VV @VV{\tau}V\\
\C @>{\eta}>> I \setminus_\pi \C @>{\eta_\dg}>> \C
\end{CD}
\end{equation}
with cartesian squares, where $\nu$, $\nu_\dg$ are the embedding and
its right-adjoint.
\end{lemma}

\proof{} If the required right-adjoint $\nu:I \setminus_\pi \C \to
\Ar(\C)$ exists, then to obtain the adjoint $\eta_\dg$ of
Definition~\ref{fib.def}, it suffices to take $\eta_\dg = \sigma
\circ (\pi \setminus \id)_\dg$. Conversely, if $\pi$ is a fibration,
then the required right-adjoint $(\pi \setminus \id)_\dg$ is the
functor \eqref{adj.fu} for the admissible full embedding $\eta:\C
\to I \setminus_\pi \C$. Now to simplify notation, let $\C_0 = I
\setminus_\pi \C$; then Lemma~\ref{adj.co.le} for the left-reflexive
full embedding $\rho = (\pi \setminus \id)_\dg:\C_0 \to \Ar(\C)$
provides an equivalence
$$
\Ar(\C) \cong \Ar(\C) \setminus_{\rho,\C_0} \C_0,
$$
and if we consider its target as a full subcategory in
$\Ar(\Ar(\C)|\C)$, the cartesian square \eqref{ar.ar.sq} induces the
cartesian square on the right in \eqref{fib.facto.sq}. By
Lemma~\ref{adj.pb.le}, the adjunction between $\eta$ and $\eta_\dg$
then induces an adjunction between $\nu$ and $\nu_\dg$.
\endproof

The notion of a fibration goes back to \cite{SGA}, but
Definition~\ref{fib.def} is more modern; let us recall the classical
setup and explain the relation between the two. For any functor
$\pi:\C \to I$, a map $f:c' \to c$ in $\C$ defines a map $f':\langle
c',\id \rangle \to \langle c,\pi(f)\rangle$ in the right comma-fiber
$\pi(c') \setminus_{\pi} \C$. The map $f$ is {\em cartesian} over
$I$ if for any $c'' \in \C_{\pi(c')} \subset \pi(c') \setminus_\pi
\C$, any map $g:c'' \to \langle c',\pi(f) \rangle$ factors uniquely
through $f'$. The functor $\pi$ is a classical {\em prefibration} if
for any $c \in \C$ and $g:i \to \pi(c)$, there exists a cartesian
map $f:c' \to c$ such that $\pi(f)=g$ (this is called a {\em
  cartesian lifting} of the map $g$). Equivalently, one can define
the {\em strict fiber} $\C^{st}_i$ for any $i \in I$ as the
subcategory $\C^{st}_i \subset \C$ spanned by objects $c$ with
$\pi(c)=i$ and morphisms $f$ with $\pi(f)=\id_i$; then $\C^{st}_i
\subset \C_i$ is the full subcategory spanned by pairs $\langle
c,\alpha \rangle$ with $\alpha=\id_i$, and $\pi$ is a classical
prefibration if and only if for any $i \in I$, the embedding
$\nu_i:\C^{st}_i \subset \C_i \to i \setminus_\pi \C$ has a
right-adjoint $\nu^\dg_i$. In this case, a map $f:i \to i'$ defines
a functor
\begin{equation}\label{fib.tr}
f^* = \nu_i^\dg \circ f^* \circ \nu_{i'}:\C^{st}_{i'} \to \C^{st}_i,
\end{equation}
called the {\em transition functor} of the prefibration $\pi$, where
$f^*$ in the right-hand side stands for the functor \eqref{f.i}. For
any composable pair of maps $f$, $g$, we have a functorial map
\begin{equation}\label{fib.compo}
g^* \circ f^* \to (f \circ g)^*.
\end{equation}
A classical prefibration is a classical {\em fibration} if a
composition of cartesian maps is cartesian, and this is equivalent
to saying that the maps \eqref{fib.compo} are
isomorphisms. Grothendieck in \cite{SGA} axiomatizes the situation:
he says that a {\em pseudofunctor} from $I^o$ to the category of
categories consists of categories $\C_i$, functors \eqref{fib.tr}
and isomorphisms \eqref{fib.compo}, subject to certain higher
constraints, and proves that a pseudofunctor comes from a fibration,
uniquely in a suitable sense.

Now, a functor $\pi:\C \to I$ is a fibration in the sense of
Definition~\ref{fib.def} if and only if the natural projection
$\pi':I \times_I \C \to I$ is a classical fibration. Indeed, we have
$\C_i \cong (I \times_I \C)^{st}_i$, $i \in I$, the adjoint functors
$\nu^\dg_i:i \setminus_\pi \C \to \C_i$ are induced by the functor
$\nu^\dg:I \setminus_\pi \C \to \C$ right-adjoint to $\nu$, so that
$\pi'$ is a classical prefibration, and it is easy to check that the
maps \eqref{fib.compo} are then isomorphisms. In the other
direction, say that a functor $\pi:\C \to I$ is an {\em
  isofibration} if for any $i \in I$, the embedding $\C^{st}_i
\subset \C_i$ is an equivalence (or equivalently, if all invertible
maps in $I$ admit cartesian liftings). Then $\pi$ is a classical
fibration if and only if it is an isofibration and a fibration in
the sense of Definition~\ref{fib.def}.

The main reason we replace the classical notion of a fibration with
Definition~\ref{fib.def} is that the classical notion is not
invariant under equivalences. In a similar vein, we say that $\pi:\C
\to I$ is a {\em prefibration} if $\pi':I \times_I \C \to I$ is a
classical prefibration (or equivalently, if $\C_i \subset i
\setminus \C$ is right-admissible for any $i \in I$). The notion of
a cartesian map is the same in both approaches; for fibrations, a
map $f:c \to c'$ is cartesian over $I$ iff the corresponding object
in $\Ar(\C)$ lies in the essential image of the fully faithful
right-adjoint of Lemma~\ref{fib.ar.le}.

\begin{exa}\label{triv.fib.exa}
For any category $\C$, the tautological projection $\C \to \ppt$ is
a fibration, and a map $f$ in $\C$ is cartesian if and only if it is
invertible.
\end{exa}

The adjoint functors $\nu^\dg_i$ in \eqref{fib.tr} are unique up to
a unique isomorphism but not strictly unique, so that to obtain a
pseudofunctor, one needs to make a choice. Grothendieck call this a
{\em cleavage} of a fibration. In terms of Definition~\ref{fib.def},
this is codified as follows.

\begin{defn}\label{cleav.def}
A {\em cleavage} of a fibration $\pi:\C \to I$ is a pair of a
functor $\nu^\dg:I \setminus_\pi \C \to \C$ over $I$ and a map $\id
\to \nu^\dg \circ \nu$ that defines an adjunction in the sense of
Definition~\ref{adj.map.def}.
\end{defn}

A cleavage is unique up to a unique isomorphism, and unless
otherwise specified, we will tacitly assume that a cleavage has been
chosen for any fibration in sight, so we can speak about {\em the}
transition functors defined by \eqref{fib.tr}. In particular, we
will say that a fibration is {\em constant} resp.\ {\em
  semiconstant} resp.\ {\em fully faithful} over a morphism $f$ if
the transition functor $f^*$ is an equivalence resp.\ an epivalence
resp.\ fully faithful.

\begin{lemma}\label{sat.C.le}
For any fibration $\C \to I$, the class $W'(\C)$ resp.\ $W(\C)$ of
maps $f$ in $I$ such that $\C$ is fully faithful resp.\ constant
along $f$ is closed and left-saturated resp.\ saturated.
\end{lemma}

\proof{} Both classes $W(\C)$, $W'(\C)$ are obviously closed. Since
equivalences are invertible up to an isomorphism, $W(\C)$ has the
two-out-of-three property. If we have functors $\gamma:\C_i \to
\C_{i'}$, $\gamma':\C_{i'} \to \C_{i''}$ such that $\gamma'$ is
fully faithful, then for any $c,c' \in \C_i$,
$\gamma':\C_{i'}(\gamma(c),\gamma(c')) \to
\C_{i''}(\gamma'(\gamma(c)),\gamma'(\gamma(c'))$ is an isomorphism,
so that $\gamma$ is fully faithful iff so is the composition
$\gamma' \circ \gamma$. To finish the proof, it remains to check
that $W(\C)$ and $W'(\C)$ are closed under retracts. Indeed, assume
given a commutative diagram
\begin{equation}\label{I.retr.dia}
\begin{CD}
i_0 @>{a}>> i @>{b}>> i_0\\
@V{f_0}VV @VV{f}V @VV{f_0}V\\
i'_0 @>{a'}>> i' @>{b'}>> i'_0
\end{CD}
\end{equation}
in $I$ such that $b \circ a = \id$ and $b' \circ a' = \id$. Then
${b'}^*$ is faithful, thus if $f^*$ is fully faithful, so is $f^*
\circ {b'}^* \cong b^* \circ f_0^*$, hence also $f_0^*$. Moreover,
$a^*$ is essentially surjective, thus if $f^*$ is essentially
surjective, so is $f_0^*$. Moreover, for any $c,c' \in \C_{i'_0}$,
$g \in \C_{i_0}(f_0^*(c),f_0^*(c'))$, we have $g = a^*(b^*(g))$, and
if $f^*$ is full, we have $b^*(g) = f^*(g')$ for some $g' \in
\C_{i'}({b'}^*(c),{b'}^*(c'))$, so that $g =
f_0^*(a^*(g'))$. Therefore $f_0^*$ is also full.
\endproof

\begin{remark}\label{cleav.rem}
In the setting of \cite{SGA}, a cleavage of a (classical) fibration
$\C \to I$ is a choice of the adjoint functors $\nu_i^\dg$ that
appear in \eqref{fib.tr}. This defines the transition functors
\eqref{fib.tr} on the nose, and then provides the functorial
isomorphisms \eqref{fib.compo}. A cleavage is {\em strict} if all
the isomorphisms \eqref{fib.compo} are identity maps (this
corresponds to a pseudofunctor that is actually a strict
functor). While it is not true that any classical fibration $\C \to
I$ admits a strict cleavage, it becomes true if $I$ is locally
small, and one is allowed to replace $\C$ with an equivalent
category $\C'$ (see below Example~\ref{cat.exa}).
\end{remark}

The story for cofibrations is dual, and if one wishes, one can also
retell it in terms of left comma-categories and comma-fibers, as in
Remark~\ref{cofib.rem}. In particular, a functor $\pi:\C \to I$ is a
{\em precofibration} --- that is, opposite to a prefibration --- if
and only if $\C_i \subset \C/_\pi i$ is left-admissible for any
$i$. In such a case, we have a transition functor $f_!:\C_i \to
\C_{i'}$ for any map $f:i \to i'$ in $I$, and a natural map
\begin{equation}\label{cofib.compo}
(f \circ g)_! \to f_! \circ g_!
\end{equation}
for any pair $f$, $g$ of composable maps. A precofibration is a
cofibration if and only if all the maps \eqref{cofib.compo} are
isomorphisms.

\begin{lemma}\label{adm.cof.le}
Assume given a fibration resp.\ cofibration $\C \to I$ and a right
resp.\ left-admissible full subcategory $I' \subset I$. Then $\C' =
\C \times_I I' \subset \C$ is right resp.\ left-admissible.
\end{lemma}

\proof{} By duality, it suffice to consider the cofibration case.
Since the full subcategory $I' \subset I$ is left-admissible, the
embedding functor $\gamma:I' \to I$ has a left-adjojnt $\gamma_\dg:I
\to I'$, and we have the adjunction map $a:\id \to \gamma \circ
\gamma_\dg$. Then the left-adjoint $\C' \to \C$ to the full
embedding $\C' \to \C$ sends an object $c \in \C_i$, $i \in I$ to
$a(i)_!(c) \in \C_{\gamma_\dg(i)}$.
\endproof

A useful observation is that if we have a fibration $\pi:\C \to I$,
a map $v:c_0 \to c_1$ in $\pi^*(\Iso)$, and a cartesian map $f:c_1'
\to c_1$, then the fibered product $c_0' = c_1' \times_{c_1} c_0$
exists in $\C$ and is given by $c_0' = \pi(f)^*c_0$. Then one can
define the {\em transpose cofibration} $\pi_\perp:\C_\perp \to I^o$
by letting $\C_\perp$ have the same objects as $\C$, with morphisms
from $c$ to $c'$ given by diagrams
\begin{equation}\label{perp.dia}
\begin{CD}
c @<{f}<< \wt{c} @>{v}>> c'
\end{CD}
\end{equation}
with cartesian $f$ and $\pi(v)=\id$, and with compositions given by
appropriate fibered products. The functor $\pi_\perp$ sends a
diagram \eqref{perp.dia} to $\pi(f)$. The cofibration $\C_\perp/I^o$
has the same fibers as $\C$ and the same transition functors $f^o_!
= f^*:\C_i \to \C_{i'}$, while the maps \eqref{cofib.compo} are
inverse to the maps \eqref{fib.compo}. One can then also consider
the fibration $\C^o_\perp = (\C_\perp)^o \to I$ opposite to the
cofibration $\C_\perp \to I^o$; it has fibers $\C^o_i \cong
(\C_i)^o$ and transition functors opposite to the transition functors
of the original fibration $\C \to I$. Dually, for a cofibration $\C
\to I$, we have the transpose fibration $\C^\perp = (\C^o_\perp)^o
\to I^o$, and the transpose opposite cofibration $\C^o_\perp =
(\C^o)_\perp \to I$ with fibers $\C^o_i$.

\subsection{Cartesian functors.}\label{cart.subs}

Now assume given a commutative square
\begin{equation}\label{rel.cart}
\begin{CD}
  \C' @>{\gamma}>> \C\\
  @V{\pi'}VV @VV{\pi}V\\
  I' @>{\phi}>> I
\end{CD}
\end{equation}
of categories and functors.

\begin{defn}\label{cart.def}
If $\pi$, $\pi'$ in \eqref{rel.cart} are fibrations, with adjoint
functors $\eta_\dg$, $\eta'_\dg$, then $\gamma$ is {\em cartesian
  over $\phi$} if the base change map $\gamma \circ \eta'_\dg \to
\eta_\dg \circ (\phi \setminus \gamma)$ is an isomorphism. If $\pi$
and $\pi'$ are cofibrations, then $\gamma$ is {\em cocartesian over
  $\phi$} if $\gamma^o$ is cartesian over $\phi^o$.
\end{defn}

Equivalently, if we let $(\pi \setminus \id)_\dg$, $(\pi' \setminus
\id)_\dg$ be the fully faithful adjoint functors of
Lemma~\ref{fib.ar.le}, then $\gamma$ in \eqref{rel.cart} is
cartesian over $\phi$ iff the base change map $\Ar(\gamma) \circ
(\pi' \setminus \id)_\dg \to (\pi \setminus \id)_\dg \circ (\id
\setminus \gamma)$ is an isomorphism, and the latter happens iff
$\Ar(\gamma)$ sends the essential image of $(\pi' \setminus
\id)_\dg$ into the essential image of $(\pi \setminus \id)_\dg$. In
other words, $\gamma$ is cartesian over $\phi$ if it sends maps
cartesian over $I'$ to maps cartesian over $I$. Dually, $\gamma$ is
cocartesian if it sends cocartesian maps to cocartesian maps. If
$I'=I$ and $\phi=\id$, so that $\gamma$ is a functor over $I$, then
we say that $\gamma$ is {\em cartesian resp.\ cocartesian over $I$}
if it is cartesian resp.\ cocartesian over $\id:I \to
I$. Explicitly, if we have two fibrations $\C',\C \to I$, then a
functor $\gamma$ from $\C'$ to $\C$ over $I$ is given by a
collection of functors $\gamma_i:\C'_i \to \C_i$, $i \in I$, and a
collection of maps
\begin{equation}\label{fu.fib}
\gamma_f:\gamma_i \circ f^* \to f^* \circ \gamma_{i'},
\end{equation}
one for any morphism $f:i \to i'$ in $I$, satisfying certain
compatibility conditions. We say that $\gamma$ is {\em cartesian
  over $f$} if \eqref{fu.fib} is an isomorphism. In this case,
\eqref{fu.fib} defines a commutative square
\begin{equation}\label{fu.fib.sq}
\begin{CD}
\C'_{i'} @>{\gamma_{i'}}>> \C_{i'}\\
@V{f^*}VV @VV{f^*}V\\
\C'_i @>{\gamma_i}>> \C_i
\end{CD}
\end{equation}
and we further say that $\gamma$ is {\em strongly cartesian} over
$f$ if the square \eqref{fu.fib.sq} is cartesian. Then $\gamma$ is
cartesian over $I$ if and only if it is cartesian over all maps in
$I$. Dually, a functor $\gamma:\C' \to \C$ over $I$ between two
cofibrations $\C,\C' \to I$ is {\em cocartesian over a map} $f$ if
$\gamma^o$ is cartesian over $f^o$, and $\gamma$ is simply {\em
  cocartesian} if it is cocartesian over all maps. A cartesian
functor $\gamma:\C' \to \C$ between fibrations $\C,\C' \to I$
induces a cocartesian functor $\gamma_\perp:\C'_\perp \to \C_\perp$
between transpose cofibrations, and then an opposite cartesian
functor $\gamma^o_\perp$. If the base category $I$ is small, we
denote by
\begin{equation}\label{fib.I}
\Cat \mm I, \Cat \mc I \subset \Cat \bbi I
\end{equation}
the subcategories spanned by small categories fibered
resp.\ cofibered over $I$, and cartesian and resp.\ cocartesian
functors between them.

\begin{exa}
The composition $\phi \circ \gamma:\C' \to I$ of fibrations
$\gamma:\C' \to \C$, $\pi:\C \to I$ is trivially a fibration, and
$\gamma$ is cartesian over $I$.
\end{exa}

\begin{exa}
For any fibration $\pi:\C \to I$ and essentially small category
$I'$, the functor $\Fun(I',\C) \to \Fun(I',I)$ given by
postcomposition with $\pi$ is a fibration, and for any functor
$\gamma:I'' \to I'$ from an essentially small $I''$,
$\gamma^*:\Fun(I',\C) \to \Fun(I'',\C)$ is cartesian over
$\gamma^*:\Fun(I',I) \to \Fun(I'',I)$ (both claims immediately follow
from Example~\ref{fun.adj.exa}).
\end{exa}

\begin{exa}\label{fun.cart.exa}
If a commutative square \eqref{rel.cart} is cartesian, and $\pi$ is
a fibration, then $\pi'$ is also a fibration, and $\gamma$ is
cartesian over $\phi$. For a general square \eqref{rel.cart} with
fibrations $\pi'$, $\pi$, $\gamma$ is cartesian over $\phi$ iff the
induced functor $\C' \to \phi^*\C$ is cartesian over $I'$.
\end{exa}

\begin{lemma}\label{fib.sq.le}
Assume given a cartesian square \eqref{cat.dia}, and a fibration
$\pi:\C \to I$ such that $\pi_l = \pi \circ \gamma_l$ is a fibration
and $\gamma_l$ is cartesian for $l=0,1$. Then $\C_{01} \to I$ is
also a fibration, and $\gamma^0_{01}$, $\gamma^1_{01}$ are cartesian
over $I$.
\end{lemma}

\proof{} Apply Corollary~\ref{adm.sq.corr} to the cartesian square
$$
\begin{CD}
  I \setminus \C_{01} @>>> I \setminus \C_1\\
  @VVV @VVV\\
  I \setminus \C_0 @>>> I \setminus \C
\end{CD}
$$
induced by \eqref{cat.dia}, and the right-admissible full embeddings
$\eta$.
\endproof

\begin{lemma}\label{right.fib.le}
Assume given a fibration $\C \to I$ and a right-admissible full
subcategory $\C' \subset \C$ such that the functor $\gamma:\C \to
\C'$ right-adjoint to the embedding is right-adjoint over $I$. Then
$\C' \to I$ is a fibration, and $\gamma$ is cartesian over $I$.
\end{lemma}

\proof{} The restriction of $\eta:\C \to I \setminus \C$ to $\C'
\subset \C$ is right-reflexive over $I$, with right-adjoint $\gamma
\circ \eta_\dg$, and then by Lemma~\ref{adj.ff.le}, this restricts
to an adjoint $\eta'_\dg:I \setminus \C' \to \C'$. Thus $\C'$ is a
fibration, and $\gamma$ is cartesian by construction.
\endproof

\begin{lemma}\label{cart.fib.le}
Assume given fibrations $\pi_l:\C_l \to I$, $l=0,1$, and a functor
$\gamma:\C_0 \to \C_1$ cartesian over $I$. Then $\sigma$ in
\eqref{rel.rlc.facto} is a fibration, and $\gamma$ is a fibration if
and only if $\gamma \setminus_I \id:\Ar(\C_0|I) \to \C_1
\setminus_{\gamma,I} \C_0$ admits a fully faithful right-adjoint
$(\gamma \setminus_I \id)_\dg$. Moreover, for any two fibrations
$\gamma:\C_0 \to \C_1$, $\gamma':\C_0' \to \C_1$, a functor
$\phi:\C_0 \to \C_0'$ over $\C_1$ is cartesian over $\C_1$ if and
only if \thetag{i} it is cartesian over $I$, and \thetag{ii} the
base change map $(\phi \setminus_I \phi) \circ (\gamma \setminus_I
\id)_\dg \to (\gamma' \setminus_I \id)_\dg \circ (\id \setminus_I
\phi)$ is an isomorphism.
\end{lemma}

\proof{} Since $\gamma$ is cartesian, we have a commutative square
\begin{equation}\label{ga.eta.dia}
\begin{CD}
\C_0 @>{\eta}>> I \setminus_{\pi_0} \C_0 @>{\eta_\dg}>> \C_0\\
@V{\gamma}VV @VV{\id \setminus \gamma}V @VV{\gamma}V\\  
\C_1 @>{\eta}>> I \setminus_{\pi_1} \C_1 @>{\eta_\dg}>> \C_1,
\end{CD}
\end{equation}
where $\eta_\dg$ are as in \eqref{fib.facto.sq}. The diagrams
\eqref{fib.facto.sq} for the given fibrations $\pi_0:\C_0 \to I$ and
$\pi_1:\C_1 \to I$ induce a diagram
\begin{equation}\label{comma.facto.sq}
\begin{CD}
\Ar(\C_0|I) @>{\nu}>> \Ar(\C_0) @>{\nu_\dg}>> \Ar(\C_0|I)\\
@V{\gamma \setminus_I \id}VV @V{\gamma \setminus \id}VV @VV{\gamma
  \setminus_I \id}V\\
\C_1 \setminus_{\gamma,I} \C_0 @>{\nu'}>> \C_1 \setminus_\gamma \C_0
@>{\nu'_\dg}>> \C_1 \setminus_{\gamma,I} \C_0\\
@V{\tau}VV  @V{\pi_1 \setminus \id}VV  @VV{\tau}V\\
\C_0 @>{\eta}>> I \setminus \C_0 @>{\eta_\dg}>> \C_0
\end{CD}
\end{equation}
with cartesian squares, where the middle row is obtained by
pulling back $\tau:\Ar(\C_1|I) \to \C_1$ to the top row in
\eqref{ga.eta.dia}. Then firstly, by Lemma~\ref{adj.pb.le}, an
adjunction between $\eta$ and $\eta_\dg$ induces an adjunction
between $\nu'$ and $\nu'_\dg$, so that $\C_1 \setminus_{\gamma,I} \C_
\subset \C_1 \setminus_\gamma \C$ is right-admissible, and $\sigma$
of \eqref{rel.rlc.facto} is a fibration by
Lemma~\ref{right.fib.le}. And secondly, the fact that $(\gamma
\setminus \id)$ in \eqref{comma.facto.sq} admits a fully faithful
right-adjoint iff so does $(\gamma \setminus_I \id)$ also follows
from Lemma~\ref{adj.pb.le}. Finally, say that a map $v$ in $\C_0$ is
{\em cartesian over $(\C_1|I)$} if $\pi_0(v)$ in invertible, so that
$v$ is a map in some fiber $\C_{0i}$, $i \in I$, and it is then
cartesian over $\C_{1i}$. Then any map $f$ in $\C_0$ cartesian over
$\C_1$ canonically factors as $f = c \circ v$, with $c$ resp.\ $v$
cartesian over $I$ resp.\ $(\C_1|I)$, and the two conditions
\thetag{i} resp.\ \thetag{ii} on the functor $\phi$ mean that it
sends maps cartesian over $I$ resp.\ $(\C_1|I)$ to maps cartesian
over $I$ resp.\ $(\C_1|I)$.
\endproof

\begin{lemma}\label{fib.adj.le}
Assume given fibrations $\C,\C' \to I$, and assume given a functor
$\gamma:\C' \to \C$ cartesian over $I$.
\begin{enumerate}
\item If $\gamma_i:\C'_i \to \C_i$ is left-reflexive for
  any $i \in I$, then $\gamma$ is left-reflexive over $I$, and if the square
  \eqref{fu.fib.sq} is left-reflexive for any $f:i \to i'$, then
  the adjoint functor $\gamma^\dg:\C \to \C'$ is cartesian over $I$.
\item If the square \eqref{fu.fib.sq} is right-reflexive for any $f:i
  \to i'$, then $\gamma$ is right-reflexive over $I$, and the adjoint
  functor $\gamma^\dg$ is cartesian over $I$.
\end{enumerate}
\end{lemma}

\proof{} In \thetag{i}, for reflexivity, Lemma~\ref{cart.fib.le} and
the factorization \eqref{rel.rlc.facto} reduced us to the case $\C =
I$, and then to construct an adjoint, it suffices to observe that a
comma-fiber $i \setminus \C'$ has an initial element iff so does its
right-admissible subcategory $\C'_i \subset i \setminus \C'$. Then
the maps \eqref{fu.fib} for the adjoint functor $\gamma^\dg$ are
exactly the base change maps \eqref{bc.eq} for the squares
\eqref{fu.fib.sq}, so $\gamma^\dg$ is cartesian iff all these
squares are left-reflexive. For \thetag{ii}, consider the
transpose-opposite functor $\gamma^o_\perp$.
\endproof

\begin{lemma}\label{bc.le}
Assume given a commutative diagram 
\begin{equation}\label{bc.dia}
\begin{CD}
I_1' @>{\gamma'}>> I'_0 @>{\pi'}>> I'\\
@V{\nu_1}VV @VV{\nu_0}V @VV{\nu}V\\
I_1 @>{\gamma}>> I_0 @>{\pi}>> I
\end{CD}
\end{equation}
of essentially small categories with cartesian squares, and assume
that $\pi$ and $\pi \circ \gamma$ are cofibrations, and $\gamma$ is
cocartesian over $I$. Then for any target category $\E$ and functor
$E:I_1 \to \E$ such that $\colim_{I_1 /_\gamma i}\sigma(i)^*E$
exists for any $i \in I_0$, both $\gamma_!E$ and $\gamma'_!\nu_1^*E$
exist, and the base change map \eqref{ga.bc.eq} is an
isomorphism.
\end{lemma}

\proof{} The factorization \eqref{rel.rlc.facto} for the functor
$\gamma$ provides an isomorphism $\gamma_!  \cong \tau_! \circ
\sigma^*$, and similarly for $\gamma'$, so that (the dual version
of) Lemma~\ref{cart.fib.le} reduces us to the case when $I_0 = I$
and $\pi=\id$. Moreover, a map of functors is invertible iff it is
invertible at any $i \in I'$, so it further suffices to
consider the case $I' = \ppt$, $\nu=\eps(i)$ for some $i \in
I$. Then the Kan extension $\gamma_!$ is given by \eqref{kan.eq},
and the claim reduces to \eqref{adm.eq} for the left-admissible full
embedding $I_{1i} \subset I_1 /_\gamma i$.
\endproof

If the target category $\E$ in Lemma~\ref{bc.le} is complete and
cocomplete, then the base change isomorphism \eqref{ga.bc.eq}
induces other useful base change maps. Firstly, we have the
adjoint isomorphism
\begin{equation}\label{bc.adj.eq}
  \gamma^* \circ \nu_{0*} \cong \nu_{1*} \circ {\gamma'}^*.
\end{equation}
Secondly, we can apply \eqref{bc.eq} to \eqref{ga.bc.eq} and obtain
a functorial map
\begin{equation}\label{bc.2.eq}
\gamma_! \circ \nu_{1*} \to \nu_{0*} \circ \gamma'_!.
\end{equation}
In general, it is not an isomorphism, but in some situations, it is
(one such appears below in Lemma~\ref{bc.2.le}).

\begin{lemma}\label{epi.le}
Assume given fibrations $\C$, $\overline{\C}$, $\C'$ over a category
$I$, and functors $\rho:\C \to \overline{\C}$, $F:\C \to \C'$
cartesian over $I$. Then $\rho$ is essentially surjective or an
epivalence if and only if so is $\rho(i)$ for any $i \in I$, and if
$\rho$ is an epivalence, then $F$ factors through $\rho$ if and only
if $F(i)$ factors through $\rho(i)$ for any $i \in I$.
\end{lemma}

\proof{} Immediately follows from the description of functors in
terms of the maps \eqref{fu.fib}.
\endproof

\begin{corr}\label{retra.corr}
Assume given a commutative diagram of categories and functors
\begin{equation}\label{retra.dia}
\begin{CD}
\wt{\C}'' @>{\wt{q}}>> \wt{\C}' @>{\wt{p}}>> \wt{\C}\\
@V{s}VV @V{s'}VV @VV{s}V\\
\C'' @>{q}>> \C' @>{p}>> \C
\end{CD}
\end{equation}
such that the square on the left is semicartesian, the vertical
functors are fibrations, $\wt{p}$ is conservative, and $q$ is
essentially surjective. Then if the outer rectangle is cartesian,
both squares are cartesian.
\end{corr}

\proof{} For any object $c \in \C''$, the functors $\wt{q}$ and
$\wt{p}$ induce functors
$$
\begin{CD}
\wt{\C}''_c @>{\wt{q}_c}>> \wt{\C}'_{q(c)} @>{\wt{p}_c}>> \wt{\C}_{p(q(c))}
\end{CD}
$$
between the corresponding fibers of the fibrations $s$, $s'$, the functor
$\wt{q}_c$ is essentially surjective by Lemma~\ref{epi.le}, and
the functor $\wt{p}_c$ is conservative. By Lemma~\ref{epi.le}, the square on the
left is cartesian resp.\ semicartesian iff $\wt{q}_c$ is an
equivalence resp.\ an epivalence for any $c \in \C''$, and since $q$
is essentially surjective, the same holds for the square on the
right and $\wt{p}_c$. It remains to apply Lemma~\ref{epi.epi.le}.
\endproof

If we have fibrations $\C,\C' \to I$ and $\C$ is essentially small,
then for any class $v$ of maps in $I$, we denote by
$\Fun^v_I(\C,\C') \subset \Fun_I(\C,\C')$ the full subcategory
spanned by functors cartesian over all maps $f \in v$, and
analogouly for $\Sec^v(I,\C)$. In particular,
$\Fun_I^\Tot(\C,\C')$ consists of all cartesian functors, and
$\Sec^\Tot(I,\C)$ consists of all cartesian sections. If $I$ is
essentially small, we have a natural equivalence
\begin{equation}\label{sec.perp.eq}
\Sec^\Tot(I,\C)^o \cong \Sec^\Tot(I,\C^o_\perp), \qquad
\gamma \mapsto \gamma^o_\perp,
\end{equation}
for any fibration $\C \to I$ with transpose-opposite fibration
$\C^o_\perp \to I$.

\begin{exa}\label{fun.v.i.exa}
For any essentially small category $I$ and any category $\C$, we
have $\Fun(I,\C) \cong \Sec(I,\C \times I) = \Fun_I(I,\C \times I)$,
where $\C \times I$ is fibered over $I$ via the projection $\C
\times I$. For any class $v$ of maps in $I$, we denote
$$
\Fun^v(I,\C) = \Fun_I^v(I,\C \times I) \subset \Fun_I(I,\C \times I)
\cong \Fun(I,\C).
$$
Then explicitly $\Fun^v(I,\C) \subset \Fun(I,\C)$ is the full
subcategory spanned by functors that invert all maps $f \in v$.
\end{exa}

More generally, if the class $v$ is closed, and $\C,\C' \to I$ are
fibrations only over the dense subcategory $I_v \subset I$, then one
still has the maps \eqref{fu.fib} for any $f \in v$, and it makes
sense to speak about functors and sections cartesian over $v$. We
will use the same notation.

\begin{lemma}\label{sec.adj.le}
Assume given fibrations $\C,\C' \to I$ over an essentially small
category $I$, and a cartesian functor $\gamma:\C \to \C'$ left or
right-reflexive over $I$ such that the adjoint $\gamma^\dg:\C' \to
\C$ is also cartesian. Then the induced functor
$\Sec(\gamma):\Sec^\Tot(I,\C) \to \Sec^\Tot(I,\C')$ is left
resp.\ right-reflexive, with adjoint functor $\Sec(\gamma^\dg)$.
\end{lemma}

\proof{} Clear. \endproof

\begin{lemma}\label{map.le}
Assume given a fibration $\pi:\C \to [1]$, two cartesian sections
$\gamma,\gamma':[1] \to \C$, and two maps $f,f':\gamma \to \gamma'$.
Then $f(1)=f'(1)$ implies $f(0)=f'(0)$, and if $\pi$ is constant,
$f(0)=f'(0)$ implies $f(1)=f'(1)$.
\end{lemma}

\proof{} The first claim immediately follows from the universal
property of a cartesian map. For the second, note that if $\pi$ is a
constant fibration, then so is $\pi^o:\C^o \to [1]^o \cong [1]$.
\endproof

\subsection{Examples.}

Let us give some examples of fibrations and cofibrations less
trivial than Example~\ref{triv.fib.exa}. We start with the following
prototypical example.

\begin{exa}\label{fib.cyl.exa}
A functor $\C \to [1]$ with fibers $\C_0$, $\C_1$ is a fibration if
and only if $\C \cong \Cyl^o(\gamma)$ is the dual cylinder of a
functor $\gamma:\C_1 \to \C_0$ in the sense of
Definition~\ref{cyl.def}, with $\gamma$ being the transition
functor. Dually, a cofibration $\C \to [1]$ with transition functor
$\gamma:\C_0 \to \C_1$ is canonically equivalent to the cylinder
$\Cyl(\gamma)$.
\end{exa}

\begin{exa}\label{fib.lt.exa}
For any fibration $\C \to I$, the induced functor $\C^< \to I^<$ is
a fibration (if $I=\ppt$, then $\C^<$ is the dual cylinder of the
tautological projection $\C \to \ppt$, so this is
Example~\ref{fib.cyl.exa}).
\end{exa}

\begin{exa}\label{cyl.fib.exa}
For any functor $\gamma:I_0 \to I_1$, giving a fibration $\C \to
\Cyl(\gamma)$ is equivalent to giving fibrations $\C_0 \to I_0$,
$\C_1 \to I_1$, and a functor $\gamma^*\C_1 \to \C_0$ cartesian over
$I_0$. In particular, to extend a fibration $\C \to I$ to the
category $I^>$, it suffices to specify the fiber $\C_o$ over the
added terminal object $o \in I^>$, and a functor $r:\C_o \times I
\to \C$ cartesian over $I$. Dually, extending a cofibration $\C \to
I$ to a cofibration over $I^<$ amounts to specifying the fiber
$\C_o$ and a cocartesian functor $l:\C_0 \times I \to \C$.
\end{exa}

\begin{exa}\label{sq.fib.exa}
A fibration $\C' \to [1]^2$ is the same thing as a commutative
diagram \eqref{cat.dia}, with $\C \cong \C'_{0 \times 0}$, $\C_0
\cong \C'_{0 \times 1}$, $\C_1 \cong \C'_{1 \times 0}$, and $\C_{01}
\cong \C'_{1 \times 1}$. If we identify $[1]^2 \cong \V^>$, we have
\begin{equation}\label{sec.sq.eq}
  \Sec^\Tot(\V,\C') \cong \C_0 \times_{\C} \C_1,
\end{equation}
and the functor $r:\C_{01} \times \V \to \C|_{\V}$ of
Example~\ref{cyl.fib.exa} induces \eqref{c.01.f} after taking
$\Sec^\Tot(\V,-)$.
\end{exa}

\begin{defn}\label{sq.cart.def}
Assume given a fibration $\C \to I$ and a commutative square in $I$
represented by a functor $\gamma:[1]^2 \to I$. Then $\C$ is {\em
  cartesian} resp.\ {\em semicartesian} resp.\ {\em weakly
  semicartesian} along $\gamma$ iff so is the commutative diagram
\eqref{cat.dia} corresponding to $\gamma^*\C \to [1]^2$.
\end{defn}

\begin{lemma}\label{car.epi.le}\mbox{}\!\!\!
Assume given two commutative diagrams \eqref{cat.dia} that are
cartesian and represented by fibrations $\C,\C'/[1]^2$, and a
functor $\phi:\C' \to \C$ cartesian over $[1]^2$ that is an
epivalence over $0 \times 0$ and essentially surjective over $1
\times 0$ and $0 \times 1$. Then $\phi$ is essentially
surjective. Moreover, if $\phi$ is an equivalence over $0 \times 0$
and an epivalence over $1 \times 0$ and $0 \times 1$, then it is an
epivalence.
\end{lemma}

\proof{} Clear. \endproof

\begin{lemma}\label{car.rel.le}
Assume given commutative squares \eqref{cat.dia} represented by
squares $\C,\C' \to [1]^2$, and a functor $\phi:\C' \to \C$
cartesian over $[1]^2$. Then if $\C$ is cartesian in the sense of
Definition~\ref{sq.cart.def}, the following conditions are
equivalent:
\begin{enumerate}
\item $\C'$ is cartesian resp.\ semicartesian,
\item $\eps^*\C' \to [1]^2$ is cartesian resp.\ semicartesian for
  any cartesian section $\eps:[1]^2 \to \C$ of the fibration $\C \to
  [1]^2$.
\end{enumerate}
Moreover, if $\C$ is only semicartesian, then the semicartesian part
of \thetag{ii} still implies the semicartesian part of \thetag{i}.
\end{lemma}

\proof{} Clear.
\endproof

\begin{exa}\label{bb.exa}
For any category $\E$, the projection \eqref{fun.bb} with its fibers
$\Fun(I,\E)$, $I \in \Cat$ is a fibration whose transition functors
are given by pullbacks $\gamma^*$ with respect to functors $\gamma:I
\to I'$. For any closed class of maps $v$ in $\E$, the induced
projection $(\Cat\bb\E)_{\Cat(v)} \to \Cat$ with its fibers
$\Fun(I,\E)_v$ is also a fibration, with the same transition
functors. If $v = \Iso$, then the class $\Cat(\Iso)$ that defines
the dense subcategory $\Cat \bbi \E \subset \Cat \bb \E$ consists of
all maps cartesian over $\Cat$.
\end{exa}

\begin{exa}\label{cat.exa}
Dually, for any small category $I$, the forgetful functor $I \bbd
\Cat \to \Cat$ with its fibers $(I \bbd \Cat)_{\C} \cong \Fun(I,\C)$
is a cofibration whose transition functor $\gamma_!:\Fun(I,\C) \to
\Fun(I,\C')$ for some $\gamma:\C \to \C'$ is given by
postcomposition with $\gamma$. In particular, if $I = \ppt$, then
$\Cat_\idot = \ppt \bbd \Cat$ is the category of pairs $\langle I,i
\rangle$ of a small category $I$ and an object $i \in I$, with
morphisms $\langle I,i \rangle \to \langle I',i' \rangle$ given by
pairs $\langle \gamma,f \rangle$ of a functor $\gamma:I \to I'$ and
a morphism $f:\gamma(i) \to i'$. The forgetful functor
$\nu_\idot:\Cat_\idot \to \Cat$, $\langle I,i \rangle \mapsto I$ is
a small cofibration whose fiber over some $I$ is $I$ itself, and
whose transition functor associated to a functor $\gamma:I \to I'$
is $\gamma$. This cofibration is obtained by applying the
Grothendieck construction to the tautological pseudofunctor
$\id:\Cat \to \Cat$. It is universal in the following sense: if the
opposite $I^o$ to the base category $I$ is locally small, then any
small cofibration $\C \to I$ is canonically equivalent to
$\gamma^*\Cat_\idot$, where $\gamma:I \to \Cat$ sends $i \in I$ to
$\Sec^{\Tot}(i \setminus I,\tau(i)^*\C)_{red}$ (this is essentially
a rigidified version of the Grothendieck construction).
\end{exa}

\begin{exa}\label{IX.exa}
A fibration $\C \to I$ is {\em discrete} if all its fibers $\C_i$
are essentially discrete. A discrete fibration is automatically
faithful. Small discrete fibrations correspond by the Grothendieck
construction to functors $X:I^o \to \Sets$ --- for any such $X$, the
category $IX$ of pairs $\langle i,x \rangle$, $i \in I$, $x \in X(i)$
with the forgetful functor $IX \to I$ is a discrete fibration, and
any discrete fibration $\C \to I$ with essentially small fibers is
canonically equivalent to $IX \to I$ for a unique $X$. The category
$IX$ is known as the {\em category of elements} of the functor
$X$. If we extend $X$ to a functor $I^{<o} \to \Sets$ by sending the
terminal object $o \in I^{<o}$ to $\ppt$, then the fibration $I^<X
\cong (IX)^< \to I^<$ is that of Example~\ref{fib.lt.exa}. If $X =
\Y(i)$ is represented by an object $i \in I$, then $IX \cong I/i$
and $I^<X \cong I^</i$. For any $X:I^o \to \Sets$, we have
$\lim_{I^o}X \cong \Sec(I,IX)$, where we note that the category
$\Sec(I,IX)$ is small and discrete, thus a set, and $\colim_{I^o}X
\cong \pi_0(IX)$ is the set of connected components of the small
category $IX$. Dually, a cofibration $\C \to I$ is discrete if so is
the opposite fibration $\C^o \to I^o$, and small discrete
cofibrations correspond to functors $I \to \Sets$. We denote the
discrete cofibration corresponding to some $X:I \to \Sets$ by
$I^\perp X \to I$.
\end{exa}

\begin{exa}\label{comma.exa}
For any category $\C$, the projection $\sigma$ resp.\ $\tau$ of
Example~\ref{ar.exa} is a fibration resp.\ cofibration, and if $\C$
has pullbacks with respect to maps in a closed class $v$, then
$\tau$ is a fibration over $\C_v$. More generally, for any functor
$\gamma:\C \to I$, the projection $\sigma:I \setminus_\gamma \C \to
I$ of \eqref{rc.dec} is a fibration, while the projection $\tau:I
\setminus_\gamma \C \to \C$ is a cofibration, and all the functors
$\tau(i)$ of \eqref{p.i} are discrete cofibrations. The functors
$f^*$ of \eqref{f.i} are transition functors of the fibration
$\sigma$ (so our notation is consistent). The fibration $I
\setminus_\gamma \C \to I$ corresponds by the Grothendieck
construction to the functor $I^o \to \Cat$ obtained by applying the
extended Yoneda embedding \eqref{yo.cat.eq} to $\C \in \Cat \bb
I$. The situation for left comma-category $\C /_\gamma I$ is dual.
\end{exa}

\begin{exa}\label{fun.fib.exa}
For any category $\C$ and essentially small category $I$, the
projection $\ev_o:\Fun(I^>,\C) \to \C$ given by evaluation at $o \in
I^>$ is a cofibration, with fibers $\Fun(I^>,\C)_c \cong
\Fun(I,\C/c)$, $c \in \C$. If we identify $I^> \cong \Cyl(\gamma)$
with the cylinder of the tautological projection $\gamma:I \to
\ppt$, then in terms of \eqref{cyl.com}, $\ev_o$ is the cofibration
$\tau:\Fun(I,\C) \setminus_{\gamma*} \C \to \C$ of
Example~\ref{comma.exa}. Dually, $\ev_o:\Fun(I^<,\C) \to \C$ is a
fibration with fibers $\Fun(I^<,\C)_c \cong \Fun(I,c \setminus \C)$.
\end{exa}

\begin{exa}\label{tw.exa}
For any category $I$, the fibration $\Ar(I)^\perp \to I^o$
transpose to $\tau:\Ar(I) \to I$ of Example~\ref{ar.exa} is the
{\em twisted arrow category} $\Tw(I)$ --- its objects are arrows $i
\to i'$ in $I$, and morphisms from $i_0 \to i'_0$ to $i_1 \to i'_1$
are given by commutative diagrams
$$
\begin{CD}
i_0 @>>> i_0'\\
@VVV @AAA\\
i_1 @>>> i_1'.
\end{CD}
$$
As in Example~\ref{ar.exa}, we denote by $\sigma:\Tw(I) \to I$
resp.\ $\tau:\Tw(I) \to I^o$ the projections sending $i \to i'$ to
$i$ resp.\ $i'$; both are fibrations, and so is the product $\sigma
\times \tau:\Tw(I) \to I^o \times I$. Moreover, the fibration
$\sigma \times \tau$ is small and discrete, and if $I$ is
essentially small, it corresponds to the functor $(I^o \times I)^o
\cong I \times I^o \to \Sets$ sending $i \times i'$ to
$I(i',i)$. Alternatively, we can also consider the fibration
$\sigma:\Ar(I) \to I$, and then $\Tw(I) \cong \Ar(I)^o_\perp$ is
transpose-opposite to this fibration.
\end{exa}

\begin{exa}\label{tw.comma.exa}
For any categories $I$, $I'$, and a functor $\phi:I' \to I^o$,
define the {\em twisted comma-category} $I \setminus_\phi^o I'$ by
the cartesian square
\begin{equation}\label{tw.comma.sq}
\begin{CD}
  I \setminus_\phi^o I' @>>> \Tw(I)\\
  @VVV @VV{\tau}V\\
  I' @>{\phi}>> I^o,
\end{CD}
\end{equation}
where $\Tw(I)$ and $\tau$ are as in Example~\ref{tw.exa}. Then the
left vertical arrow $I \setminus_\phi^o I' \to I'$ in
\eqref{tw.comma.sq} is a fibration that is transpose to the
cofibration $\tau:I \setminus_\phi I' \to I'$, and the discrete
fibration $\Tw(I) \to I^o \times I$ induces a discrete fibration $I
\setminus_\phi^o I' \to I' \times I$.
\end{exa}

\begin{exa}\label{comma.discr.exa}
A discrete fibration $\pi:I' \to I$ is equivalent to the fibration
$\sigma(i):I / i \to I$ for some $i \in I$ iff $I'$ has a terminal
object $o$ (and in this case, we have $i = \pi(o)$).
\end{exa}

\begin{exa}\label{closed.fib.exa}
A fully faithful embedding $\gamma:\C' \to \C$ is a fibration if and
only if the full subcategory $\C' \subset \C$ is left-closed in the
sense of Definition~\ref{cl.def}. In this case, every fiber $\C'_c$,
$c \in \C$ of the fibration $\gamma$ is essentially small, and in
fact, it is either empty, or equivalent to the point category
$\ppt$. The full subcategory in $\Cat$ spanned by $\emptyset$ and
$\ppt$ is equivalent to $[1]$, the characteristic functor
\eqref{chi.eq} for $\gamma$ is provided by the Grothendieck
construction, and $\C' \cong \C_0$ is its fiber over $0 \in [1]$. We
note that the projection $t_\dg:\Cyl(\gamma) \to \C$ is both right
and left-admissible by Lemma~\ref{cl.cyl.le}, and since the
right-adjoint functor $t:\C \to \Cyl(\gamma)$ is fully faithful, so
is the left-adjoint $t_{\dg\dg}:\C \to \Cyl(\gamma)$. The
characteristic functor $\chi$ is obtained by composing $t_{\dg\dg}$
with the projection $\Cyl(\gamma) \to [1]$.
\end{exa}

\begin{exa}\label{1.sets.exa}
As a trivial but useful special case of
Example~\ref{closed.fib.exa}, consider the full embedding $\ppt
\cong \{\emptyset\} \subset \Sets$. Then it is left-closed, and the
characteristic functor $\chi:\Sets \to [1]$ is left-adjoint to the
embedding $[1] \to \Sets$ onto $\{\emptyset,\ppt\} \subset \Sets$
--- it sends $\emptyset$ to $0$ and the rest to $1$.
\end{exa}

\begin{exa}\label{V.closed.exa}
In the situation of Example~\ref{cl.exa}, since $\C = \C_0 \cup
\C_1$, the product $\chi_1 \times \chi_2:\C \to [1]^2$ of the
characteristic functors $\chi_0$, $\chi_1$ for the left-closed
subcategories $\C_0,\C_1 \subset \C$ factors through a joint
characteristic functor $\chi:\C \to \V = [1]^2 \ssetminus \{1 \times
1\} \subset [1]^2$, where $\V$ is as in \eqref{V.eq}. We then have
$\C_0 \cong \C /_{\chi} 0$, $\C_1 \cong \C /_{\chi} 1$, and $\C_{01}
\cong \C /_{\chi} o$. Conversely, for any category $\C$ equipped
with a functor $\chi:\C \to \V$, $\C_o \subset \C/0,\C_1 \subset \C$
are left-closed full embeddings, and $\C_o = (\C/0) \cap (\C/1)$, so
we are in the situation of Example~\ref{cl.exa}.
\end{exa}

\begin{exa}\label{st.exa}
Assume given a category $\C_{01}$ and two left-closed full
embeddings $\gamma_0:\C_{01} \to \C_0$, $\gamma_1:\C_{01} \to \C_1$
to some categories $\C_0$, $\C_1$. Then the cartesian cocartesian
square \eqref{st.sq} of Example~\ref{st.0.exa} can be also recovered
by the Grothendieck construction. Indeed, we have a cofibration $\C'
\to \V$ with fibers $\C'_o = \C_{01}$, $\C'_0 = \C_0$, $\C'_1 =
\C_1$ and transition functors $\gamma_0$, $\gamma_1$. Then as in
Example~\ref{V.closed.exa}, we have $\C' = (\C'/0) \cup (\C'_1)$,
and $(\C'/l) \cong \Cyl(\gamma_l)$, $l=0,1$, so to construct $\C$,
it suffices to take the union $\C \subset \C'$ of the essential
images of the left-adjoint functors $\C_l \to \Cyl(\gamma_l)$
provided by Lemma~\ref{cl.cyl.le}
\end{exa}

\begin{exa}\label{cl.X.exa}
Consider a category $I$ and a functor $X:I^o \to \Sets$, with the
corresponding discrete fibration $\pi:IX \to I$.  A functor
$\gamma:\C \to IX$ is a fibration if and only if so is the
composition $\pi \circ \gamma:\C \to I$, and in this case, $\gamma$
is cartesian over $\pi$. In particular, a map $f:X \to X'$ to
another functor $X':I^o \to \Sets$ induces a cartesian functor
$I(f):IX \to IX'$, and this functor is automatically a fibration. If
$f$ is injective, then $I(f)$ is a left-closed full embedding of
Example~\ref{closed.fib.exa}. In general, we have a canonical
equivalence between $I^o\Sets$ and the full subcategory in $\Cat \mm
I$ spanned by small discrete fibrations. This also provides an
equivalence
\begin{equation}\label{IX.sets}
  (IX)^o\Sets \cong I^o\Sets / X
\end{equation}
for any $X:I^o \to \Sets$. In terms of the equivalence
\eqref{IX.sets}, the forgetful functor $\sigma(X):I^o\Sets / X \to
I^o\Sets$ of \eqref{p.i} is identified with the left Kan extension
$\pi^o_!$. In particular, the constant functor $\ppt:(IX)^o \to
\Sets$ with value $\ppt$ is the terminal object in $(IX)^o\Sets$,
thus identified with the terminal object $\id:X \to X$ in
$(I^o\Sets)/X$, so that $\pi^o_!\ppt \cong X$, and we have
\begin{equation}\label{hom.X}
  \Hom(X,Y) \cong \Hom(\ppt,\pi^{o*}Y) \cong \lim_{\langle i,x
    \rangle \in (IX)^o}Y(i)
\end{equation}
for any $Y:I^o \to \Sets$. More generally, for any complete target
category $\C$ and functor $Y:I^o \to \C$, we can define $\Hom(X,Y)
\in \C$ by \eqref{hom.X}, and even if $\C$ is not complete, we can
still use \eqref{hom.X} as long as the limit in the right-hand side
exists (for example, if $X = \Y(i)$ is represented by an object $i
\in I$, then $\Hom(X,Y) = Y(i)$ is well-defined for any $\C$). For
any functor $X' \in (IX)^o\Sets \cong I^o\Sets/X$, the tautological
map $X' \to \ppt$ induces a map
\begin{equation}\label{IX.spl}
  \lim_{I^o}\pi^o_!X' \to \lim_{I^o}\pi^o_!\ppt \cong \lim_{I^o}X,
\end{equation}
elements in its target correspond bijectively to sections $\sigma:I \to
IX$ of the discrete fibration $\pi$, and for any such section $\sigma$,
\eqref{IX.spl} fits into a cartesian square
\begin{equation}\label{IX.spl.sq}
\begin{CD}
  \lim_{I^o}\sigma^{o*}X' @>>> \ppt\\
  @VVV @VV{\sigma}V\\
  \lim_{I^o}\pi^o_!X' @>>> \lim_{I^o}X,
\end{CD}
\end{equation}
so that \eqref{IX.spl} represents $\lim_{I^o}\pi^o_!X'$ as a
disjoint union of sets $\lim_{I^o}\sigma^{o*}X'$.
\end{exa}

\begin{lemma}\label{bc.2.le}
In the situation of Lemma~\ref{bc.le}, assume that $\E=\Sets$, $I_0
= I$, $\gamma:I_1 \to I$ is a discrete cofibration corresponding to
some functor $X:I \to \Sets$, $\nu$ is fully faithful, and the
adjunction map $X \to \nu^*\nu_*X$ is an isomorphism. Then the base
change map \eqref{bc.2.eq} is an isomorphism.
\end{lemma}

\proof{} Since by assumption, the cofibration $\gamma$ is discrete,
it induces an equivalence $\langle i,x \rangle \setminus I'_1 \cong
i \setminus I'$ for any $\langle i,x \rangle \in I_1 =
(I^oX)^o$. Now to check that the map \eqref{bc.2.eq} becomes
invertible after evaluation at any $X':I'_1 \to \Sets$, compute
$\nu_*\gamma'_!X'$ by (the dual version of) \eqref{kan.eq}, and use
decompositions \eqref{IX.spl} into components \eqref{IX.spl.sq} for
all the comma-fibers $i \setminus I'$.
\endproof

\begin{lemma}\label{cocart.le}
Assume given a category $I$ equipped with three functors $X,X',Y:I^o
\to \Sets$ and two morphisms $e:X \to X'$, $f:X \to Y$, and let $Y'
= Y \copr_X X'$. Then if $e$ is injective, the commutative square
\begin{equation}\label{cocart.sq}
\begin{CD}
IX @>{e}>> IX'\\
@V{f}VV @VVV\\
IY @>>> IY'
\end{CD}
\end{equation}
of the corresponding discrete fibrations is a cartesian cocartesian
square of categories.
\end{lemma}

\proof{} As in Example~\ref{fun.v.i.exa}, we can treat functors $IX
\to \C$ to some category $\C$ as functors $\gamma:IX \to I \times \C$ over
$I$ and then describe them in terms of components $\gamma_i$, $i \in
I$ and maps \eqref{fu.fib}. The same holds for $IX'$, $IY$ and
$IY'$, and this reduced the claim to the case $I=\ppt$, where it
amounts to Example~\ref{cocart.S.exa}.
\endproof

\begin{remark}\label{cl.rem}
If $f$ in Lemma~\ref{cocart.le} is also injective, then
\eqref{cocart.sq} is a square of the type described in
Example~\ref{cl.exa}.
\end{remark}

\begin{defn}\label{lft.v.def}
For any fibration resp.\ cofibration $\pi:\C \to I$, and any closed
class $v$ of maps in $I$, the class $\pi^\flat(v)$
resp.\ $\pi^\sharp(v)$ of maps in $\C$ consists of all cartesian
resp.\ cocartesian liftings of maps $f \in v$.
\end{defn}

\begin{remark}
In the notation of Definition~\ref{lft.v.def}, for any fibration $\E
\to I$, and any fibration $\pi:I' \to I$ of essentially small
categories, \eqref{sec.fun.I} restricts to an equivalence
\begin{equation}\label{sec.fun.v}
\Fun^v_I(I',\E) \cong \Sec^{\pi^\flat v}(I',\E)
\end{equation}
for any closed class $v$ of morphisms in $I$.
\end{remark}

\begin{lemma}\label{facto.fib.le}
\begin{enumerate}
\item Assume given a category $I$ equipped with a factorization
  system $\langle L,R \rangle$, and a functor $\gamma:I' \to I$. If
  $\gamma$ is a fibration, then $\langle \gamma^*(L),\gamma^\flat(R)
  \rangle$ is a factorization system on $I'$. Dually, if $\gamma$ is
  a cofibration, then $\langle \gamma^\sharp(L),\gamma^*(R)\rangle$
  is a factorization system on $I'$.
\item Assume given a category $I$ and a dense subcategory $I_R$
  defined by a class of morphisms $R$. Let $\Ar^R(I)$ and $\sigma$,
  $\tau$, $\eta$ be as in Example~\ref{facto.exa}. Then $R$ fits
  into a factorization system $\langle L,R \rangle$ on $I$ if and
  only if $\tau$ is a cofibration, and in this case, $L =
  \eta^*(\tau^\sharp(\Iso))$ and $\tau^\sharp(\Iso) =
  \sigma^*L$. Dually, $L$ fits into a factorization system $\langle
  L,R \rangle$ if and only if $\sigma:\Ar^L(I) \to I$ is a
  fibration, and then $R=\eta^*(\sigma^\flat(\Iso))$ and
  $\sigma^\flat(\Iso) = \tau^*(R)$.
\end{enumerate}
\end{lemma}

\proof{} Direct inspection. \endproof

\begin{exa}\label{facto.ar.exa}
In the situation of Lemma~\ref{facto.fib.le}~\thetag{ii}, not only
the functor $\tau:\Ar^R(I) \to I$ is a cofibration, but so is the
functor $\Ar(I) \to \Ar^R(I)$ left-adjoint to the embedding
$\Ar^R(I) \subset \Ar(I)$. If for any $i \in I$, we let $I /_Li$
resp.\ $I/_Ri$ be the subcategories spanned by arrows $i' \to i$ in
$L$ resp.\ in $R$, then $\tau:\Ar^R(I) \to I$ has fibers $\Ar^R(I)_i
\cong I /_R i$, $i \in I$, while the fiber of $\Ar(I) \to \Ar^R(I)$
over some arrow $i' \to i$ is $I/_Li'$. The functor $\tau:\Ar^R(I)_L
\to I_L$ is also a cofibration, and the full embedding $\Ar^R(I)_L
\to \Ar^R(I) \times_I I_L$ is cocartesian over $I_L$. The situation
for $\Ar^L(I)$ is dual.
\end{exa}

\begin{remark}\label{facto.rem}
The simplest particular case of Lemma~\ref{facto.fib.le}~\thetag{i}
is the factorization system $\langle \Iso,\Tot \rangle$ on $I$; in
this case, the claim is that we have a factorization system $\langle
\pi^*(\Iso),\pi^\flat(\Tot)\rangle$ on $\C$, where $\pi^\flat(\Tot)$
consists of all maps cartesian over $I$. The adjoint functor
\eqref{pi.co.eq} is then the left-admissible full embedding of
Example~\ref{facto.exa}, and the square \eqref{ar.facto.sq} is
\eqref{fib.facto.sq}.
\end{remark}

\begin{exa}\label{fun.I.exa}
For any essentially small category $\E$, $\Fun(\E,-)$ sends adjoint
functors to adjoint functors by Example~\ref{fun.adj.exa}, and it
commutes with taking comma-categories. Therefore for any fibration
$\pi:\C \to I$, $\Fun(\E,\C) \to \Fun(\E,I)$ is a fibration, and we
can define the {\em relative functor category} $\Fun(\E,\C|I)$ by
the cartesian square
\begin{equation}\label{fun.CI.sq}
\begin{CD}
  \Fun(\E,\C|I) @>>> \Fun(\E,\C)\\
  @VVV @VVV\\
  I @>{\gamma'}>> \Fun(\E,I),
\end{CD}
\end{equation}
where $\gamma'$ is as in \eqref{ev.dia} for the projection $\gamma:I
\times \E \to I$ (that is, $\gamma'(i)$, $i \in I$ is the constant
functor with value $i$). Then the vertical arrow on the left in
\eqref{fun.CI.sq} is also a fibration, with fibers $\Fun(\E,\C_i$,
$i \in I$.  If $\E=[1]$, then $\Fun([1],\C|I) \cong \Ar(\C|I)$ is
the relative arrow category of \eqref{rel.ar.sq}.
\end{exa}

\begin{lemma}\label{fun.i.ar.le}
For any category $\C$ equipped with a factorization system $\langle
L,R \rangle$, and any small $I$, the full subcategory
$\Fun_R(I^>,\C) \subset \Fun(I^>,\C)$ spanned by functors that
factor through $\C_R \subset \C$ is left-admissible, with some
adjoint functor $\gamma:\Fun(I^>,\C) \to \Fun_R(I^>,\C)$, and the
adjunction map $\Id \to \gamma$ is pointwise in $L$.
\end{lemma}

\proof{} If $I=\ppt$, then $I^>=[1]$, $\Fun([1],\C) \cong \Ar(\C)$,
$\Fun_R([1],\C) \cong \Ar^R(\C)$, and the claim is
Example~\ref{facto.exa}. In the general case, consider the
cofibration $\ev_o$ of Example~\ref{fun.fib.exa}, use it to identify
\begin{equation}\label{fun.i.ar.eq}
\Fun(I^>,\C) \cong \Fun(I,\Ar(\C)|\C), \ \Fun_R(I^>,\C) \cong
\Fun(I,\Ar^R(\C)|\C),
\end{equation}
where the relative functor categories of \eqref{fun.CI.sq} are taken
with respect to the projections $\tau$, and use the adjoint functor
$\Ar(\C) \to \Ar^R(\C)$ of Example~\ref{facto.exa} to construct the
adjoint functor $\gamma$.
\endproof

\subsection{Families of groupoids.}

As usual, a category is a {\em groupoid} if all its morphisms are
invertible (in other words, $\C$ coincides with its isomorphism
groupoid $\C_{\Iso}$ of Definition~\ref{iso.def}). A groupoid $\C$
is automatically equivalent to its opposite $\C^o$ --- the
equivalence $\C \cong \C^o$ is identical on objects and inverts all
morphisms. For any group $G$, we denote by $\ppt_G$ the groupoid
with one object with automorphism group $G$. More generally, if $G$
acts on a set $S$, the {\em quotient groupoid} $[S \bb G]$ has
elements $s \in S$ as objects, with morphisms $s \to s'$ given by
elements $g \in G$ such that $gs=s'$. A groupoid $\C$ is {\em
  connected} if it is connected as a category (in particular,
non-empty). A connected groupoid is essentially small. For any
essentially small groupoid $\C$, we denote by $\pi_0(\C)$ the set of
isomorphism classes of objects in $\C$, and we observe that $\C$ is
connected if and only if $\pi_0(\C)$ has exactly one element. For
any object $c \in \C$, we denote by $\pi_1(\C,c)$ the automorphism
group $\Aut(c)$. For any groupoid $\C$ and two objects $c,c' \in
\C$, the group $\pi_1(\C,c)$ acts on the set $\C(c,c')$, and if $\C$
is connected, this defines an equivalence $\C \cong
\pi_1(\C,c)\tors$ between $\C$ and the category of torsors over the
group $\pi_1(\C,c)$. We also have an obvious equivalence
$\ppt_{\pi_1(\C,c)} \cong \C$.

\begin{remark}
Note that even if a group $G$ is finite, $\ppt_G$ is finite as a
category only if $G$ is trivial; otherwise $|\ppt_G|$ is countable
(see Example~\ref{fin.rig.exa}). This is not a bug but a feature.
\end{remark}

A functor $\gamma:\C' \to \C$ between groupoids is an equivalence if
and only if it admits an adjoint (indeed, both adjunction maps are
automatically invertible). Even if $\gamma$ is not an equivalence,
it is tautologically conservative, so it is an epivalence if and
only if it is essentially surjective and full. Even if $\C'$ is not
a groupoid, $\gamma$ is automatically a fibration (although not a
classical fibration). It is essentially surjective if and only if
all its fibers $\C'_c$, $c \in \C$ are non-empty, and an epivalence
if and only if all the fibers $\C'_c$ are connected groupoids.

\begin{exa}\label{wr.exa}
For any group $G$, a functor $\gamma:\C \to \ppt_G$ from some
category $\C$ is automatically a fibration, and it is a good
formalization of the intuitive notion of a ``category with an action
of the group $G$''. For example, if $G$ acts on a set $S$, then for
any category $\C$, it acts on the product $\C^S$ of copies of $\C$
numbered by elements $s \in S$, and this defines a fibration $\C \wr
\ppt_G \to \ppt_G$ with fiber $\C^S$ called the {\em wreath product}
of $G$ and $\C$. If $\C \cong \ppt_H$ is also a connected groupoid,
then $\ppt_H \wr \ppt_G \cong \ppt_{H \wr G}$, where $H \wr G$ is
the wreath product of the groups in the usual sense.
\end{exa}

For any groupoid $\C$ and essentially small category $I$, the
functor category $\C^I$ is a groupoid. The {\em inertia groupoid}
$\I(\C)$ is $\I(\C) = \C^{\ppt_{\Z}}$. The tautological projection
$\tau:\ppt_{\Z} \to \ppt$ defines a fully faithful embedding
$\tau^*:\C \to \I(\C)$.

\begin{lemma}\label{ine.le}
Assume given a functor $\gamma:\C' \to \C$ between groupoids such
that the induced functor $\I(\gamma):\I(\C') \to \I(\C)$ between
inertia groupoids is an epivalence resp.\ full over $\C \subset
\I(\C)$. Then $\gamma$ is an equivalence resp.\ fully faithful.
\end{lemma}

\proof{} It suffices to prove that $\C'_c \cong \ppt$ for any object
$c$ in $\C$ resp.\ in the essential image of $\gamma$, and then
since for any $c \in \C \subset \I(\C)$, we obviously have
$\I(\C')_c \cong \I(\C'_c)$, we may assume right away that $\C \cong
\ppt$. Then by assumption, $\I(\C')$ hence also $\C'=\C'_c$ is
connected, so that $\C'_c \cong \ppt_G$ for some group $G$. But
objects in $\I(\ppt_G)$ are group maps $\Z \to G$, or equivalently,
elements $g \in G$, and two objects are isomorphic iff the group
elements are conjugate. Therefore any $g \in G$ is conjugate to the
identity, so that $G$ is trivial.
\endproof

\begin{lemma}\label{epi.sur.le}
An essentially surjective functor $\gamma:\C' \to \C$ between
group\-oids is an epivalence if and only if the diagonal functor
$\delta:\C' \to \C' \times_{\C} \C'$ is essentially surjective.
\end{lemma}

\proof{} An object in $\C' \times_{\C} \C'$ is a triple $\langle
c,c',a \rangle$, $c,c' \in \C'$, $a:\gamma(c) \to \gamma(c')$ a map,
and such a triple is isomorphic to some $\langle
c'',c'',\id\rangle$, $c'' \in \C'$ if and only if $a=\gamma(a')$ for
some $a':c \to c'$.
\endproof

\begin{defn}
A {\em family of groupoids} over a category $I$ is a fibration $\C
\to I$ whose fiber $\C_i$ is a groupoid for any $i \in I$.
\end{defn}

A ``family of groupoids'' is the best approximation to the original
``cat\'e\-gorie fibr\'ee en groupo\"ides'' that we could find. A
fibration $\pi:\C \to I$ is a family of groupoids if and only if all
maps in $\C$ are cartesian over $I$, or equvalently, if $\pi^*(\Iso)
= \Iso$ (if $I=\ppt$, this is Example~\ref{triv.fib.exa}). There is
also the following criterion that does not require one to know {\em
  a priori} that $\pi$ is a fibration.

\begin{lemma}\label{grp.c.le}
A functor $\pi:\C \to I$ is a family of groupoids if and only if the
corresponding functor \eqref{pi.co.eq} is an equivalence.
\end{lemma}

\proof{} If $\pi$ is a family of groupoids, then \eqref{pi.co.eq}
has a fully faithful right-adjoint by Lemma~\ref{fib.ar.le}, and the
adjunction map is in $(\pi \setminus \id)^*(\Iso)=\Iso$. Conversely,
if \eqref{pi.co.eq} is an equivalence, then $\pi$ is a fibration by
Lemma~\ref{fib.ar.le}, and moreover, for any functor $\gamma:I' \to
I$, with the functor $\gamma':\C = \gamma^*\C \to \C$ and the
pullback fibration $\pi':\C' \to I'$, the diagram
$$
\begin{CD}
\C' \setminus_{\id} \C' @>{\gamma' \setminus \gamma'}>>
\C \setminus_{\id} \C\\
@V{\pi' \setminus \id}VV @VV{\pi \setminus \id}V\\
I' \setminus_{\pi'} \C' @>{\gamma \setminus \gamma'}>> I
\setminus_\pi \C
\end{CD}
$$
is cartesian. Taking as $\gamma$ the embeddings $\eps(i):\ppt \to
I$, $i \in I$ of \eqref{eps.c}, we reduce to the case $I=\ppt$. In
this case, \eqref{pi.co.eq} is just the projection $\tau:\Ar(\C) \to
\C$, and the fact that it is an equivalence means that its fibers
are equivalent to $\ppt$, or explicitly, that that for any maps $f:c
\to c_0$, $f':c' \to c_0$ in $\C$, there exists a unique map $g:c
\to c'$ such that $f = f' \circ g$. This in turn amounts to saying
that every map in $\C$ is invertible.
\endproof

\begin{exa}\label{facto.ar.grp.exa}
The fibration $\sigma:\Ar^L(I)_R \to I_R$ of
Example~\ref{facto.ar.exa} is a family of groupoids; the
corresponding equivalence \eqref{pi.co.eq} is provided by the unique
factorization property of a factorization system.
\end{exa}

For any family of groupoids $\C \to I$, the transpose cofibration
$\C_\perp \to I^o$ is equivalent to the opposite cofibration $\C^o
\to I^o$ over $I^o$. For any two families of groupoids $\C,\C' \to
I$, a functor $\gamma:\C \to \C'$ over $I$ is automatically
conservative and cartesian over $I$, and moreover, it is itself a
family of groupoids. A discrete fibration is trivially a family of
groupoids. A family of groupoids $\gamma:\C \to I$ is {\em
  connected} if it has connected fibers $\C_i$; by
Lemma~\ref{epi.le}, this is equivalent to saying that $\gamma$ is an
epivalence (just treat $\gamma$ as a functor over $I$). For a small
family of groupoids $\gamma:\C \to I$, we denote by
$\pi_0(\gamma):\pi_0(\C |^\gamma I) \to I$ the fibration with
discrete fibers corresponding to the sets $\pi_0(\C_i)$, $i \in I$, we
shorten $\C |^\gamma I$ to $\C | I$ when $\gamma$ is clear from the
context, and we note that $\gamma$ canonically factors as
\begin{equation}\label{grp.eq}
\begin{CD}
\C @>{a}>> \pi_0(\C | I) @>{\pi_0(\gamma)}>> I,
\end{CD}
\end{equation}
where $a$ is an epivalence (thus a connected family of
groupoids). For any section $\sigma:I \to \C$, we denote by
$\pi_1(\C,\sigma)$ the functor from $I^o$ to groups that sends $i
\in I$ to $\pi_1(\C_i,\sigma(i))$. If $\C \to I$ is connected, then
we have a natural equivalence
\begin{equation}\label{grp.trs}
\C \cong \pi_1(\C,\sigma)\tors
\end{equation}
over $I$, where for any functor $G$ from $I^o$ to groups, $G\tors$
is the category of pairs $\langle i,X \rangle$ of an object $i \in
I$ and a $G(i)$-torsor $X$.

\begin{lemma}\label{sq.epi.le}
Assume given a cartesian square \eqref{cat.dia} such that $\gamma_0$ is a
connected family of groupoids, and $\gamma_1$ is essentially
surjective. Then $\gamma^0_{01}$ is a connected family of groupoids, and the
square is cocartesian.
\end{lemma}

\proof{} The first claim is obvious: $\gamma^0_{01}$ is a pullback of a
fibration, thus a fibration, and it has the same fibers as $\gamma_0$. The
second claim then immediately follows from the factorization
property of Lemma~\ref{epi.le} and the identifications $\Fun(I,\E)
\cong \Fun_I(I,\E \times I)$ of Example~\ref{fun.v.i.exa}.
\endproof

\begin{lemma}\label{sq.le}
Assume given two families of groupoids $\C',\C/[1]^2$ and a functor
$\gamma:\C' \to \C$ over $[1]^2$, and assume that $\C$ and $\C'$ are
cartesian over $[1]^2$. Then the middle term $\pi_0(\C'|\C)$ of the
factorization \eqref{grp.eq} is cartesian over $[1]^2$. Moreover,
the same is true is $\C'$ is only semicartesian but $\gamma$ is an
equivalence over $0 \times 0 \in [1]^2$.
\end{lemma}

\proof{} Clear. \endproof

\begin{defn}\label{pres.def}
A {\em presentation} of a functor $\C \to I$ is a pair $\langle
\C',\alpha \rangle$ of a category $\C'$ and an epivalence
$\alpha:\C' \to \C \times_I \C$. For any two categories $\C_0$,
$\C_1$ equipped with functors $\C_0,\C_1 \to I$ and presentations
$\langle \C_0',\alpha_0 \rangle$, $\langle \C_1',\alpha_1
\rangle$, a {\em presentation} of a functor $\gamma:\C_0 \to \C_1$
over $I$ is a functor $\gamma':\C'_0 \to \C_1'$ equipped with an
isomorphism $\alpha_1 \circ \gamma' \cong (\gamma \times \gamma)
\circ \alpha_0$.
\end{defn}

\begin{exa}\label{pres.exa}
For any categories $\C$, $I$, $I'$ and functors $\C \to I' \to I$,
we have a cartesian square
\begin{equation}\label{pres.comp.sq}
\begin{CD}
  \C \times_{I'} \C @>>> \C \times_I \C\\
  @VVV @VVV\\
  I' @>{\delta}>> I' \times_I I,
\end{CD}
\end{equation}
where $\delta$ is the diagonal embedding. Then a presentation $\C'$
of the functor $\C \to I$ induces a presentation $\delta^*\C'$ of
the functor $\C \to I'$.
\end{exa}

\begin{lemma}\label{semi.epi.le}
Assume given families of groupoids $\C,I \to [1]^2$ equipped with an
epivalence $\gamma:\C \to I$ cartesian over $[1]^2$, and assume that
$\gamma$ admits a presentation $\C' \to \C \times_I \C$, also
cartesian over $[1]^2$. Then if $\C$ and $\C'$ are semicartesian, so
is $I$.
\end{lemma}

\proof{} Identify $[1]^2 \cong \V^>$ as in Example~\ref{sq.fib.exa},
with $o \in \V^>$ corresponding to $1 \times 1 \in [1]^2$. Then we
have epivalences $\C_o \to \Sec^\Tot(\V,\C)$, $\C'_o \to
\Sec^\Tot(\V,\C')$, and we need to show that $I_o \to
\Sec^\Tot(\V,I)$ is also an epivalence. Let $\Sec^\Tot(\V,\C)' =
\Sec^{\Tot}(\V,\C) \times_{\Sec^{\Tot}(\V,I)}
\Sec^{\Tot}(\V,\C)$. Then $\Sec^{\Tot}(\V,\C') \to
\Sec^{\Tot}(\V,\C)'$ is essentially surjective by
Lemma~\ref{car.epi.le}, so that $\C_o \times_{I_o} \C_o \to
\Sec^{\Tot}(\V,\C)'$ is also essentially surjective. Since $\C_o \to
I_o$ is an epivalence, $\delta:\C_o \to \C_o \times_{I_o} \C_o$ is
essentially surjective by Lemma~\ref{epi.sur.le}, and then
$\delta:\Sec^{\Tot}(\V,\C) \to \Sec^{\Tot}(\V,\C)'$ is also
essentially surjective. Since $\Sec^{\Tot}(\V,\C) \to
\Sec^{\Tot}(\V,I)$ is essentially surjective by
Lemma~\ref{car.epi.le}, it is then an epivalence by
Lemma~\ref{epi.sur.le}, and then $I_o \to \Sec(\V,I)$ is an
epivalence by Lemma~\ref{epi.epi.le}.
\endproof

In the situation of Definition~\ref{pres.def}, a presentation
$\gamma'$ of a functor $\gamma$ defines a commutative square
\begin{equation}\label{pres.sq}
\begin{CD}
  \C_0' @>{\gamma'}>> \C_1'\\
  @V{\pi_0 \circ \alpha_0}VV @VV{\pi_1 \circ \alpha_1}V\\
  \C_0 @>{\gamma}>> \C_1,
\end{CD}
\end{equation}
where $\pi_l:\C_l \times_I \C_l \to \C_l$, $l=0,1$ is the projection
onto the first factor. Say that the presentation $\gamma'$ is {\em
  strict} if the square \eqref{pres.sq} is cartesian.

\begin{lemma}\label{pres.le}
Assume given connected families of groupoids $\C_0,\C_1 \to I$
equipped with presentations $\C_0'$, $\C_1'$, and a functor
$\gamma:\C_0 \to \C_1$ over $I$ that admits a strict presentation
$\gamma':\C_0' \to \C_1'$. Then $\gamma$ is an epivalence.
\end{lemma}

\proof{} It suffices to check that $\gamma$ is an epivalence over
any object $i \in I$, and then since \eqref{pres.sq} is compatible
with pullbacks with respect to functors $I' \to I$, we may assume
right away that $I=\ppt$. Then $\C_l = \ppt_{G_l}$, $l=0,1$ are
connected groupoids associated to some groups $G_0$, $G_1$. For any
group $G$, a presentation $\C'$ of the projection $\C = \ppt_G \to
\ppt$ is of the form $\C' = \ppt_{G'}$, for some group $G'$ equipped
with a surjective map $G' \to G \times G$, and then the kernel $H
\subset G'$ of the projection $\pi:G' \to G$ onto the first factor
maps surjectively onto the second one. In our situation, we have
groups $G_0'$, $G_1'$ with subgroups $H_l \subset G_l'$, $l=0,1$,
and since \eqref{pres.sq} is cartesian, the map $H_0 \to H_1$
induced by $\gamma'$ is an isomorphism. Therefore the map
$\gamma:G_0 \to G_1$ is surjective.
\endproof

\section{Fibrations and functoriality.}\label{fun.fib.sec}

\subsection{Relative functor categories.}\label{rel.fun.subs}

For any functor $\gamma:I \to I'$, and any category $\C$ equipped
with a functor $\pi:\C \to I$, denote by $\gamma_\trr\C$ the same
$\C$ but considered as a category over $I'$ via the functor $\gamma
\circ \pi:\C \to I'$. For any $\C,\C' \to I$, a lax functor $\phi:\C
\to \C'$ over $I$ induces a lax functor
$\gamma_\trr(\phi):\gamma_\trr\C \to \gamma_\trr\C'$, so that
$\gamma_\trr$ defines a functor
\begin{equation}\label{cat.phi}
  \gamma_\trr:\Cat \bb I \to \Cat \bb I'.
\end{equation}
If we restrict our attention to $\Cat \bbi I$, and more generally,
to all, possibly large categories and functors over $I$, then
$\gamma_\trr$ provides a sort of a $2$-categorical left-adjoint to
the pullback operation $\gamma^*$. Namely, for any $\C' \to I'$, we
have the projection $\alpha:\gamma_\trr\gamma^*\C' \cong I
\times_{I'} \C' \to \C'$, and for any $\C \to I$, any functor
$\phi:\gamma_\trr\C \to \C'$ over $I'$ factors as
\begin{equation}\label{trr.dia}
\begin{CD}
\gamma_\trr\C @>{\gamma_\trr(\phi')}>> \gamma_\trr\gamma^*\C'
@>{\alpha}>> \C',
\end{CD}
\end{equation}
where $\phi':\C \to \gamma^*\C'$ is a functor over $I$; moreover,
the pair of $\phi'$ and the factorization \eqref{trr.dia} is unique
up to a unique isomorphism, as in Definition~\ref{adj.map.def}. This
is just a fancy way of restating the obvious: $\phi$ defines a
commutative square of categories and functors, and $\phi'$ is the
corresponding functor \eqref{c.01.f}.  As a slightly less obvious
application of the Grothendieck construction, let us show that
sometimes, $\gamma^*$ also admits a $2$-categorical right-adjoint.

\begin{defn}\label{fl.def}
A functor $\gamma:I \to I'$ is {\em proper} if for any category $\C
\to I$ over $I$, there exists a category $\gamma_\trl\C \to I'$
over $I'$ equipped with a functor $\ev:\gamma^*\gamma_\trl\C \to
\C$ over $I$ with the following universal property:
\begin{itemize}
\item for any category $\C' \to I'$, any functor $\phi:\gamma^*\C' \to
  \C$ over $I$ factors as
\begin{equation}\label{trl.dia}
\begin{CD}
\gamma^*\C' @>{\gamma^*(\phi')}>> \gamma^*\gamma_{\trl}\C
@>{\ev}>> \C,
\end{CD}
\end{equation}
where $\phi':\C' \to \gamma_\trl\C$ is a functor over $I'$, and the
pair of $\phi'$ and a factorization \eqref{trl.dia} is unique up
to a unique isomorphism.
\end{itemize}
For any proper functor $\gamma:I \to I'$, and any category $\C$, the
{\em relative functor category} $\Fun(I|I',\C)$ is defined by
\begin{equation}\label{rel.fun.eq}
  \Fun(I|I',\C) = \gamma_\trl(I \times \C).
\end{equation}
where $I \times \C$ is considered as a category over $I$ via the
projection $I \times \C \to I$.
\end{defn}

Just as in the usual adjunction situation, if $\C' \to I'$ and $\C
\to I$ are essentially small, the universal factorizations
\eqref{trr.dia}, \eqref{trl.dia} provide equivalences
\begin{equation}\label{fun.2adj.eq}
  \Fun_I(\gamma_\trr\C',\C) \cong \Fun_{I'}(\C',\gamma^*\C), \quad
  \Fun_{I'}(\gamma^*\C,\C') \cong \Fun_I(\C,\gamma_\trl\C'),
\end{equation}
where the second equivalence only exists if $\gamma$ is proper.

\begin{lemma}\label{fun.le}
A small cofibration $\gamma:I \to I'$ is proper in the sense of
Definition~\ref{fl.def}. Moreover, let $I_\sharp \subset I$ be the
dense subcategory of maps cocartesian over $I'$, with the embedding
functor $e:I_\sharp \to I$. Then if $e^*\C \to I_\sharp$ is a
fibration, so is $\gamma_\trl\C \to I'$, and in this case, for any
functor $\phi:\gamma^*\C \to \C'$ over $I$, the corresponding
functor $\phi'$ in \eqref{trl.dia} is cartesian over $I'$ iff $\phi$
is cartesian over $I_\sharp \subset I$.
\end{lemma}

\proof{} We first construct the relative functor categories
\eqref{rel.fun.eq}. To do this, for any category $\C$, we define
$\Fun(I| I',\C) \to I'$ as the fibration with fibers
$\Fun(I|I',\C)_i \cong \Fun(I_i,\C)$, $i \in I'$, and transition
functor $(f_!)^*$ for any map $f:i \to i'$, where $f_!:I_i \to
I_{i'}$ is the transition functor of the cofibration $\gamma:I \to
I'$. Note that a functor $\phi:\C \to \C'$ induces a functor
$\wt{\phi}:\Fun(I|I',\C) \to \Fun(I|I',\C')$, cartesian over $I'$,
and if $\phi$ is a fibration, then by Lemma~\ref{cart.fib.le} and
Example~\ref{fun.I.exa}, so is $\wt{\phi}$. Moreover, for any $i \in
I'$, we have the evaluation functor $\ev_i:I_i \times \Fun(I_i,\C)
\to \C$, and taken together, these define a functor
\begin{equation}\label{ev.I.eq}
\ev:I \times_{I'} \Fun(I|I',\C) \to \C.
\end{equation}
We claim that for any category $\C'$ equipped with a functor
$\pi:\C' \to I'$, a functor $\phi:\gamma^*\C' \to \C$ factors
uniquely as
\begin{equation}\label{fun.facto}
\begin{CD}
\gamma^*\C' = I \times_{I'} \C' @>{\id \times \phi'}>>
I \times_{I'} \Fun(I|I',\C) @>{\ev}>> \C,
\end{CD}
\end{equation}
for a unique functor $\phi':\C' \to \Fun(I|I',\C)$ over
$I'$. Indeed, for $\C'=I'$, we can take the tautological section
$s:I' \to \Fun(I|I',I)$ of the fibration $\Fun(I|I',I) \to I'$
sending $i \in I'$ to the embedding functor $I_i \to I$, and then
\eqref{fun.facto} holds iff $\phi' \cong \phi \circ s$; for a
general $\pi:\C' \to I'$, it then suffices to observe that
$\pi^*\Fun(I|I',\C) \cong \Fun(\gamma^*\C'|\C',\C)$. Now
$\gamma_\trl(I \times \C) = \Fun(I|I',\C)$, with the functor
$\ev$ given by \eqref{ev.I.eq}, satisfies the universal property
of Definition~\ref{fl.def}. But this universal property then implies
that we have more functoriality: for any two categories $\C_0$,
$\C_1$, a functor $\phi:I \times \C_0 \to I \times \C_1$ over $I$
induces a functor $\gamma_\trl(\phi):\gamma_\trl(I \times \C_0) \to
\gamma_\trl(I \times \C_1)$. Then for a general $\pi:\C \to I$, we
can construct $\gamma_\trl(\C)$ by the cartesian square
\begin{equation}\label{trl.sq}
\begin{CD}
\gamma_\trl(\C) @>>> \gamma_\trl(I \times \C)\\
@VVV @VV{\gamma_\trl(\id \times \pi)}V\\
\gamma_\trl(I) @>{\gamma_{\trl}(\delta)}>> \gamma_\trl(I \times I),
\end{CD}
\end{equation}
where $\delta:I \to I \times I$ is the diagonal embedding, and we
note that we have $\gamma_\trl(I) = \Fun(I|I',\ppt) \cong I'$,
$\gamma_\trl(\delta) \cong s$.  Since $\gamma^*$ preserves cartesian
squares, the functor $\ev$ for $\C$ can be induced by the
functors $\ev$ for the other three terms in \eqref{trl.sq}, and
it satisfies the universal property since all three satisfy
it. Finally, if $\pi:\C \to I$ is a fibration over $I_\sharp \subset
I$, then $\gamma_\trl(\id \times \pi) = \wt{\pi}$ is a fibration
over $\Fun(I|I',I_\sharp)$, and the functor $\gamma_\trl(\delta)
\cong s$ in \eqref{trl.sq} factors through the dense subcategory
$\Fun(I|I',I_\sharp) \subset \Fun(I|I',I)$.
\endproof

\begin{exa}\label{cat.bb.exa}
The relative functor category construction allows one to rewrite
Example~\ref{bb.exa} in terms of Example~\ref{cat.exa}: we have
\begin{equation}\label{cat.wr}
\Fun(\Cat_\idot | \Cat,\E) \cong \Cat\bb\E
\end{equation}
for any category $\E$, and this is an equivalence over $\Cat$. If we
denote by $\Cat_\idot \bb \E = \Cat_\idot \times_{\Cat} \Cat\bb\E$
the category of triples $\langle I,\alpha,i \rangle$ such that
$\langle I,\alpha \rangle \in \Cat\bb\E$ and $i \in I$, and let
$\nu_\idot:\Cat_\idot\bb\E \to \Cat\bb\E$ be the cofibration induced
by $\Cat_\idot \to \Cat$, then the functor $\ev:\Cat_\idot\bb\E \to
\E$ of \eqref{ev.I.eq} sends $\langle I,\alpha,i \rangle$ to
$\alpha(i)$.
\end{exa}

\begin{exa}
Let a group $G$ act on a set $S$. Take $I=\ppt_G$, $I' = \ppt$, and
let $I^\hdot = S_G$ be the quotient groupoid $[S \bb G]$. Then for
any category $\C$, $\Fun(S_G | \ppt_G,\C) \cong \C \wr \ppt_G$ is
the wreath product of Example~\ref{wr.exa}.
\end{exa}

\begin{remark}
The term ``proper'' in Definition~\ref{fl.def} is somewhat arbitrary
and chosen by analogy with algebraic geometry (alternatively, one
could say that $I$ is dualizable over $I'$, but it is awkward to
apply ``dualizable'' to a functor). Note that the converse to
Lemma~\ref{fun.le} is not true. At the very least, the opposite to a
proper functor is also obviously proper, so that small fibrations
are proper as well. But properness is not an empty condition. What
happens is, for any proper $I \to I'$, \eqref{fun.2adj.eq} dictates
a canonical description of objects and morphisms of the category
$\Fun(I|I',\C)$. For objects, take $\C' = \ppt$; then objects in
$\Fun(I|I',\C)$ are pairs of an object $i \in I'$ and a functor $I_i
\to \C$. For morphisms, take $\C' = [1]$; then a morphism in
$\Fun(I|I',\C)$ is the pair of a morphism $f:i \to i'$ in $I'$ and a
functor $\eps(f)^*I \to \C$, where $\eps(f):[1] \to I$ is the
functor \eqref{eps.f}. However, while this description makes sense
for any small functor $\gamma:I\to I'$, in general, one cannot
construct compositions. One possible characterization of proper
functors is given below in Remark~\ref{fl.seg.rem}.
\end{remark}

For any proper functor $\gamma:I \to I'$, the relative functor
categories $\Fun(I | I',\C)$ are uniquely defined by the universal
property. If $I' = \ppt$, so that $\gamma:I \to \ppt$ is the
tautological projection, then it is proper as soon as $I$ is
essentially small, and $\Fun(I| \ppt,\C) \cong \Fun(I,\C)$. More
generally, assume given a proper functor $I \to I'$, a category
$\C$, and another category $I''$ and functors $\nu:I\to I''$,
$\pi:\C \to I''$. Then we can modify \eqref{rel.fun.eq} by writing
\begin{equation}\label{rel.I.fun.eq}
  \Fun_{I''}(I|I',\C) = \gamma_\trl\nu^*\C,
\end{equation} 
and as in \eqref{fun.I.sq}, this category fits into a cartesian
square
$$
\begin{CD}
\Fun_{I''}(I|I',\C) @>>> \Fun(I|I',\C)\\
@VVV @VV{\pi \circ -}V\\
I' @>{\nu'}>> \Fun(I|I',I''),
\end{CD}
$$
where $\nu'$ corresponds to $\nu$ under
\eqref{fun.2adj.eq}. Then \eqref{ev.I.eq} induces a relative
evalution functor
\begin{equation}\label{ev.I.rel.eq}
\ev:I \times_{I'} \Fun_{I''}(I | I',\C) \to \C
\end{equation}
over $I''$, and for any $\pi:\C' \to I'$, a functor $\phi:\gamma^*\C' \to
\C$ over $I''$ admits a version of the unique factorization
\eqref{fun.facto}: $\phi$ uniquely decomposes as
\begin{equation}\label{fun.rel.facto}
\begin{CD}
\gamma^*\C' = I \times_{I'} \C' @>{\id \times \phi'}>>
I \times_{I'} \Fun_{I''}(I|I',\C) @>{\ev}>> \C,
\end{CD}
\end{equation}
where $\ev$ is the functor \eqref{ev.I.rel.eq} and $\phi'$ is
a functor over $I$. As a further generalization, we may assume given
a closed class $v$ of maps in $I''$ such that the functors
$I,\C \to I''$ are fibrations over $I''_v$, and consider the full
subcategory $\Fun^v_{I''}(I|I',\C) \subset \Fun_{I''}(I|
I',\C)$ spanned by
\begin{equation}\label{fun.I.v.rel}
\Fun^v_{I''}(I | I',\C)_i = \Fun^v_{I''}(I_i,\C) \subset
\Fun_{I''}(I_i,\C) = \Fun_{I''}(I | I',\C)_i
\end{equation}
for all $i \in I'$. Then $\phi'$ in \eqref{fun.rel.facto} takes
values in $\Fun^v_{I''}(I | I',\C)$ if and only if $\phi$ is
cartesian over all maps in $I''$.

\begin{exa}
If $I'=I''$, $\C \to I'$ is a fibration over all maps in $I'$, and
$I = I' \times I_0$ is a constant cofibration with some essentially
small fiber $I_0$, then $\Fun^\Tot_{I'}(I | I',\C) \cong \Fun(I_0,\C
| I')$ is the relative functor category of \eqref{fun.CI.sq}.
\end{exa}

\subsection{Induction and coinduction.}\label{ind.subs}

For any functor $\pi:\C \to I$, the projection $\sigma:I
\setminus_\pi \C \to I$ of \eqref{rc.dec} is a fibration by
Example~\ref{comma.exa}, and this fibration has two universal
properties. Firstly, for any fibration $\pi':\C' \to I$, a functor
$\gamma:\C \to \C'$ over $I$ factors as
\begin{equation}\label{coind1.eq}
\begin{CD}
\C @>{\eta}>> I \setminus_\pi \C @>{\wt{\gamma}}>> \C',
\end{CD}
\end{equation}
with $\wt{\gamma}$ cartesian over $I$, and the factorization is
unique up to a unique isomorphism. Secondly, for any functor
$\pi':\C' \to I$, a lax functor $\gamma:\C' \to \C$ over $I$ factors
as
\begin{equation}\label{coind2.eq}
\begin{CD}
\C' @>{\wt{\gamma}}>> I \setminus_\pi \C @>{\tau}>> \C,
\end{CD}
\end{equation}
with $\wt{\gamma}$ a functor over $I$, and the factorization is
again unique up to a unique isomorphism. Effectively, if $I$ is
small, we have functors
\begin{equation}\label{cat.tri.dia}
\begin{CD}
  \Cat \mm I @>{a}>> \Cat \bbi I @>{b}>> \Cat \bb I @>{\Y}>> \Cat \mm I,
\end{CD}
\end{equation}
with $a$ and $b$ being the embedding functors, and $\Y$ sending
$\pi:\C \to I$ to $\sigma:I \setminus_\pi \C \to I$. Then
informally, $a$ is $2$-right-adjoint to $\Y \circ b$, and $b$ is
$2$-left-adjoint to $a \circ \Y$, in the appropriate $2$-categorical
sense. The $2$-adjunction maps are functors
\begin{equation}\label{coind.adj}
\eta^\dg:I \setminus_\pi \C \to \C, \quad \eta:\C \to I
\setminus_\pi \C, \quad \tau:I \setminus_\pi \C \to \C,
\end{equation}
where $\eta^\dg$ is right-adjoint to $\eta$, and as in
Lemma~\ref{3.adj.le}, $\eta$ serves as the adjunction map for both
$2$-adjunctions. Explicitly, in \eqref{coind1.eq}, we have
$\wt{\gamma} \cong {\eta'}^\dg \circ (\id \setminus \gamma)$, and in
\eqref{coind2.eq}, $\wt{\gamma} \cong (\id \setminus \gamma) \circ
\eta$, where $\id \setminus \gamma :I \setminus_\pi \C' \to I
\setminus_\pi \C$ is given by the collection of functors
\eqref{lax.i}, for all $i \in I$. By the same argument as in
Lemma~\ref{3.adj.le}, it is a formal consequence of the existence of
the two $2$-adjunctions that $\Y$ is $2$-fully faithful, in the
sense that for any small $I$ and $\C,\C' \in \Cat \bb I$, the
functor
\begin{equation}\label{2.yo.eq}
\Y:\vFun_I(\C,\C') \to \Fun_I^\Tot(I \setminus \C,I \setminus \C')
\end{equation}
is an equivalence. This is not surprising, since $\Y$ corresponds to
the extended Yoneda embedding \eqref{yo.cat.eq} by the Grothendieck
construction.

If one treats the embedding $a$ of \eqref{cat.tri.dia} as a
forgetful $2$-functor that forgets the fact that $\C$ is a
fibration, then its $2$-left-adjoint $\Y \circ b$, $\C \mapsto I
\setminus \C$ is a sort of a coinduction $2$-functor from categories
over $I$ to categories fibered over $I$. If $I$ is essentially
small, one can also construct a $2$-right-adjoint to $a$, a sort of
an induction $2$-functor. This is a corollary of Lemma~\ref{fun.le}
and the canonical decomposition \eqref{fun.rel.facto}.

\begin{corr}\label{ind.corr}
Assume given an essentially small category $I$ and a functor $\C \to
I$ from some category $\C$. Then there exist a fibration $\Ind(\C)
\to I$ and a functor $\alpha:\Ind(\C) \to \C$ over $I$ such that for
any fibration $\C' \to I$, any functor $\gamma:\C' \to \C$ over $I$
factors as
\begin{equation}\label{ind.facto}
\begin{CD}
\C' @>{\gamma'}>> \Ind(\C) @>{\alpha}>> \C,
\end{CD}
\end{equation}
where $\gamma'$ is cartesian, and the factorization is unique up
to a unique isomorphism.
\end{corr}

\proof{} Consider the arrow category $\Ar(I)$, with the functors
$\sigma,\tau:\Ar(I) \to I$ of Example~\ref{ar.exa}, and their common
section $\eta:I \to \Ar(I)$ of Example~\ref{adj.exa}. Let $\Ind(\C)
= \tau_\trl\sigma^*\C$, where $\tau_\trl$ exists by
Lemma~\ref{fun.le} since $\tau$ is a cofibration, and define
$\alpha$ as the composition
$$
\begin{CD}
\Ind(\C) @>{\eta \times \id}>> \Ar(I) \times_I \Ind(\C) =
\tau^*\tau_\trl\sigma^*\C @>{\ev}>> \sigma^*\C @>>> \C,
\end{CD}
$$
where $\ev$ is as in Definition~\ref{fl.def}.  By
Lemma~\ref{facto.fib.le}~\thetag{ii}, morphisms in $\Ar(I)$
cocartesian with respect to $\tau:\Ar(I) \to I$ are exactly the
morphisms inverted by $\sigma:\Ar(I) \to I$ (that is,
$\tau^\sharp(\Tot) = \sigma^*(\Iso)$). This implies that $\Ind(\C)
\to I$ is a fibration, also by Lemma~\ref{fun.le}. Moreover, to
obtain the factorizations \eqref{ind.facto}, note that for any
fibration $\pi:\C' \to I$, we then have
$\tau^\sharp(\pi^\flat(\Tot)) = \eta_\dg^*(\Iso)$, where
$\eta_\dg:\pi^*\Ar(I) = I \setminus_\pi \C' \to \C'$ is
right-adjoint to $\eta:\C' \to I \setminus_\pi \C'$ (or in words,
cocartesian lifting of maps in $\C'$ cartesian over $I$ are exactly
maps in $\pi^*\C'$ inverted by $\eta_\dg$). Then by
Example~\ref{kan.facto.exa}, any functor $\gamma:\C' \to \C$ over
$I$ uniquely factors as $\gamma = \phi \circ \eta$, where $\phi =
\eta_!(\gamma) \cong \eta_\dg^*(\gamma):\pi^*\C' \to \C$ is a
functor over $\Ar(I)$ that inverts all maps in
$\tau_\dg\pi^\dg\Tot$. By Lemma~\ref{fun.le}, the functor $\phi$
then admits a unique factorization \eqref{fun.rel.facto}, with the
functor $\phi'$ cartesian over $I$. Altogether, we have canonical
isomorphisms
$$
\alpha \circ \phi' \cong \ev \circ (\eta \times \id) \circ
\phi' \cong \ev \circ (\eta \times \phi') \cong \ev
\circ (\id \times \phi') \circ \eta \cong \phi \circ \eta
\cong \gamma,
$$
and to obtain \eqref{ind.facto}, it suffices to take $\gamma' =
\phi'$.
\endproof

Note that just as the functor categories, $\Ind(\C)$ in
Corollary~\ref{ind.corr} is uniquely defined by its universal
property. Explicitly, we have
\begin{equation}\label{ind.C.eq}
\Ind(\C)_i \cong \Fun_I(I/i,\C) \cong \Sec(I/i,\sigma(i)^*\C),
\qquad i \in I,
\end{equation}
and the transition functors are the pullbacks $(f_!)^*$ with respect
to the functors \eqref{f.i}. For any fibration $\C' \to I$ with
essentially small $\C'$, \eqref{fun.2adj.eq} provides a canonical
equivalence
\begin{equation}\label{ind.2adj}
\Fun_I(\C',\C) \cong \Fun^\Tot_I(\C',\Ind(\C)),
\end{equation}
the induction property mentioned above. Moreover, assume given a
closed class $v$ of maps in $I$, and assume that $\C$ is a fibration
over $I_v$. Then since $\sigma:\Ar(I) \to I$ is a fibration, we can
consider the full subcategory
\begin{equation}\label{ind.v.eq}
\Ind^v(\C) = \Fun_I^v(\Ar(I)|I,\C)
\subset \Fun_I(\Ar(I)|I,\C) = \Ind(\C)
\end{equation}
given by \eqref{fun.I.v.rel}. The induced projection $\Ind^v(\C) \to
I$ is then also a fibration, and for any fibration $\C' \to I$ and
functor $\gamma:\C' \to \C$ over $I$, the corresponding functor
$\wt{\gamma}$ of \eqref{ind.facto} factors through $\Ind^v(\C)
\subset \Ind(\C)$ if and only if $\gamma$ is cartesian over maps in
$v$. In particular, \eqref{ind.2adj} induces an equivalence
\begin{equation}\label{ind.v.2adj}
\Fun^v_I(\C',\C) \cong \Fun^\Tot_I(\C',\Ind^v(\C))
\end{equation}
for any fibration $\C' \to I$ with $\C'$ essentially small.

\begin{remark}\label{ind.rem}
One can also construct the coinduction operation $\Y \circ b$ of
\eqref{cat.tri.dia} in the same way as in Corollary~\ref{ind.corr}:
for any $\pi:\C \to I$, we have $\Y(b(\C)) = I \setminus_\pi \C =
\sigma_\trr\tau^*\C$, where $\sigma_\trr$ is as in \eqref{cat.phi}.
\end{remark}

\subsection{Kan extensions.}\label{2kan.subs}

Another application of Lemma~\ref{fun.le} is a $2$-catego\-ri\-cal
version of the right Kan extension for fibrations. Assume given a
functor $\gamma:I' \to I$ such that $I'$ is essentially small, and a
fibration $\C \to I'$.

\begin{corr}\label{kan.corr}
There exist a fibration $\gamma_*\C \to I$ equipped with a functor
$\alpha:\gamma^*\gamma_*\C \to \C$, cartesian over $I'$, such that
for any fibration $\C' \to I$, any cartesian functor $F:\gamma^*\C'
\to \C$ factors as
\begin{equation}\label{kan.2facto}
\begin{CD}
\gamma^*\C' @>{\gamma^*F'}>> \gamma^*\gamma_*\C @>{\ev}>> \C
\end{CD}
\end{equation}
with $F':\C' \to \gamma_*\C$ cartesian over $I$, and such a
factorization is unique up to a unique isomorphism.
\end{corr}

\proof{} Take $\gamma_*\C = \Fun_{I'}^\Tot(I' /_\gamma I,\C)$ and
$\alpha = \ev \circ (\eta \times \id)$. To construct
the decomposition \eqref{kan.2facto}, use the same argument as in
Corollary~\ref{ind.corr}.
\endproof

Again, for any small fibrations $\C' \to I$, $\C \to I'$, the
canonical factorization \eqref{kan.2facto} provides an equivalence
\begin{equation}\label{kan.2adj}
\Fun^\Tot_{I'}(\gamma^*\C',\C) \cong
\Fun^\Tot_I(\C',\gamma_*\C), \qquad F \mapsto F',
\end{equation}
an adjunction between $\gamma^*$ and $\gamma_*$. We have cartesian
functors 
\begin{equation}\label{2.adj}
\alpha:\C' \to \gamma_*\gamma^*\C', \qquad
\alpha_\dg:\gamma^*\gamma_*\C \to \C
\end{equation}
that serve as $2$-categorical version of the adjunction maps. For
any fibration $\C \to I'$, both $\gamma_*\C \to I$ and $\alpha$ are
uniquely defined by the universal property, and explicitly, we have
\begin{equation}\label{2kan.eq}
\gamma_*\C_i \cong \Sec^\Tot(I' /_\gamma i,\sigma(i)^*\C), \qquad
i \in I,
\end{equation}
a right Kan extension counterpart of \eqref{kan.eq}, where
categories of cartesian sections serve as $2$-categorical versions
of limits. For practical applications, there is also a counterpart
of \eqref{adm.eq} provided by the following result.

\begin{lemma}\label{2adm.le}
Assume given a fibration $\C \to I$ over an essentially small
category $I$, and a left-admissible full subcategory $I' \subset I$,
with the embedding functor $\gamma:I' \to I$ and its left-adjoint
$\gamma^\dg:I \to I'$. Then the pullback functor
\begin{equation}\label{2adm.eq}
\gamma^*:\Sec(I,\C) \to \Sec(I',\C)
\end{equation}
is right-admissible, its right-adjoint $\gamma_*$ is fully faithful,
and the adjunction map $s \to \gamma_*(\gamma^*(s))$ is an
isomorphism for any section $s \in \Sec(I,\C)$ cartesian over all of
the adjunction maps $i \to \gamma(\gamma^\dg(i))$, $i \in I$.
\end{lemma}

\proof{} The adjoint $\gamma_*$ is given by
\begin{equation}\label{sigma.adj}
\gamma_*(s)(i) = a^*s(\gamma^\dg(i)), \qquad s \in \Sec(I',\C),
\end{equation}
where $a:i \to \gamma(\gamma^\dg(i))$ is the adjunction map. Since
$\gamma^\dg \circ \gamma \cong \id$, we have $\gamma^* \circ \gamma_*
\cong \id$, so that $\gamma_*$ is fully faithful.
\endproof

In particular, Lemma~\ref{2adm.le} implies that \eqref{2adm.eq}
induces an equivalence
\begin{equation}\label{bnd.equi}
\gamma^*:\Sec^\Tot(I,\C) \cong \Sec^\Tot(I',\C),
\end{equation}
between the categories of cartesian sections, with the inverse
equivalence given by \eqref{sigma.adj}. This can be used in
\eqref{2kan.eq}. For example, if $\gamma:I' \to I$ is a cofibration,
one can replace the comma-fibers $I' /_\gamma i$ in \eqref{2kan.eq}
with the usual fibers $I'_i$, and in the situation of
Lemma~\ref{bc.le}, we have a natural equivalence
\begin{equation}\label{2.bc.eq}
\nu_0^*\gamma_*\C \cong \gamma'_*\nu_1^*\C
\end{equation}
for any fibration $\C \to I_1$. As the following example shows,
Lemma~\ref{2adm.le} is also useful in the context of
Corollary~\ref{ind.corr}.

\begin{exa}\label{ab.exa}
Assume given a fibration $\C \to I$ over an essentially small
category $I$. Then by Corollary~\ref{ind.corr}, $\id:\C \to \C$
uniquely factors as $\id \cong \alpha \circ \beta$, for a unique
cartesian functor $\beta:\C \to \Ind(\C)$. It is easy to see that in
this case, the isomorphism $\id \cong \alpha \circ \beta$ defines an
adjunction between $\alpha$ and $\beta$ in the sense of
Definition~\ref{adj.map.def} (so that in particular, $\beta$ is
fully faithful). Both are functors over $I$, and so is the
adjunction. Explicitly, in terms of \eqref{ind.C.eq}, the functor
$\alpha(i):\Ind(\C)_i \to \C_i$ is given by evaluation at the
terminal object $\id:i \to i$ in $I / i$, and its adjoint $\beta(i)$
is provided by Lemma~\ref{2adm.le} and identifies $\C_i$ with the
subcategory in \eqref{ind.C.eq} spanned by cartesian sections.
\end{exa}

\begin{exa}\label{ind.exa}
Assume given a factorization system $\langle L,R \rangle$ on an
essentially small category $I$. Then for any category $\C$, one has
$$
\Ind^L(\C \times I) \cong \Fun(\Ar^R(I)|I,\C),
$$
where the cofibration $\Ar^R(I) \to I$ is $\tau$ of
Lemma~\ref{facto.fib.le}~\thetag{ii}. To see this, combine
Lemma~\ref{2adm.le} and Example~\ref{facto.exa}.
\end{exa}

For another trivial but useful application of the adjunction
\eqref{kan.2adj}, consider the cylinder $I = \Cyl(\gamma)$ of a
functor $\gamma:I_0 \to I_1$ between essentially small categories,
and assume given a functor $\chi:I \to I'$ to some $I'$ with
restrictions $\chi_0 = s^*\chi:I_0 \to I'$, $\chi_1 = t^*\chi:I_1
\to I'$, and a fibration $\C \to I'$. Then $t:I_1 \to I$ has the
left-adjoint $t_\dg$ of \eqref{cyl.deco}, we have $\chi_1^*\C \cong
t^*\chi^*\C \cong t_{\dg*}\chi^*\C$, and then \eqref{2.adj} for
$s:I_0 \to I$ gives a functor
\begin{equation}\label{coli.I.0}
\chi_1^*\C \cong t_{\dg*}\chi^*\C \to t_{\dg*}s_*s^*\chi^*\C \cong
\gamma_*\chi_0^*\C
\end{equation}
cartesian over $I_1$. In particular, if we have a functor
$\gamma:I^\hdot \to I$ between essentially small categories, and a
functor $\chi:I^\hdot \to I'$ that admits a left Kan extension
$\gamma_!\chi:I \to I'$, then the universal relative cone
$\chi_>:\Cyl(\gamma) \to I'$ and \eqref{coli.I.0} provide a functor
\begin{equation}\label{coli.I}
(\gamma_!\chi)^*\C \to \gamma_*\chi^*\C
\end{equation}
cartesian over $I$.

\begin{exa}
For any category $I$, Example~\ref{cat.bb.exa} provides a
cofibration $\nu_\idot:\Cat_\idot\bb I \to \Cat \bb I$ and an
evaluation functor $\ev:\Cat_\idot \bb I \to I$. Then for any family
of groupoids $\C$ over $I$, the target of the functor \eqref{coli.I}
for $\gamma=\nu_\idot$, $\chi = \ev$ is easy to describe: we have a
natural identification
\begin{equation}\label{nu.ev.bb}
\nu_{\idot*}\ev^*\C \cong \Cat\bb\C,
\end{equation}
with the projection to $\Cat\bb I$ induced by the fibration $\pi:\C
\to I$. If $I$ is cocomplete, then the fibration $\Cat\bb I \to
\Cat$ is a bifibration, and $\ev$ admits a left Kan extension
$\nu_{\idot!}\ev:\Cat \bb I \to I$ left-adjoint to the embedding
\eqref{ups.eq}. Then \eqref{coli.I} is an equivalence if and only if
$\C$ is cocomplete and $\pi$ preserves colimits.
\end{exa}

\begin{exa}\label{bb.v.exa}
If $\gamma = \eps(i):\ppt \to I$ is the embedding \eqref{eps.c} for
an object $i \in I$ in a category $I$, then for any category $\C$
considered as a category fibered over $\ppt$, the fibration
$\eps(i)_*\C \to I$ has fibers $(\eps(i)_*\C)_{i'} \cong
\C^{I(i,i')}$, that is, the product of copies of $\C$ numbered by
morphisms $i \to i'$. Using this rather trivial observation, one can
extend Example~\ref{cat.bb.exa} to the situation when we have a
closed class $v$ of maps in $\E$: by definition we have a cartesian
square
\begin{equation}\label{cat.I.sq}
\begin{CD}
(\Cat \bb \E)_{\Cat(v)} @>>> \Cat \bb \E\\
@VVV @VV{a}V\\
\eps(\ppt)_*\E_v @>>> \eps(\ppt)_*\E,
\end{CD}
\end{equation}
where $a$ is adjoint to the equivalence $\eps(\ppt)^*(\Cat\bb\E)
\cong \Fun(\ppt,\E) \cong \E$. This gives another description of
the dense subcategory $(\Cat \bb \E)_{\Cat(v)} \subset \Cat \bb \E$.
\end{exa}

Finally, we note that even if the base $I$ of a fibration $\pi:\C
\to I$ is not small, cartesian sections might still form a
well-defined category $\Sec^{\Tot}(I,\C)$. This happens if for any
two cartesian sections $\sigma,\sigma':I \to \C$ of $\pi$, maps from
$\sigma$ to $\sigma'$ form a set and not a proper class. In this
case, we say that $\C$ is {\em bounded} over $I$. For example, by
\eqref{bnd.equi}, this happens if $I$ admits an essentially small
left-admissible subcategory $I' \subset I$. More generally, this
happens if $I'$ is not small but we do know that $\gamma^*\C \to I'$
is bounded. Analogously, if we have a functor $\gamma:I' \to I$ and
a fibration $\C$ over $I'$, we will say that $\C$ is {\em bounded
  with respect to $\gamma$} if $\sigma(i)^*\C$ is bounded over
$I'/i$ for any $i \in I$. In this case, we have a well-defined
fibration $\gamma_*\C'$ over $I$. Note that if $\gamma$ is a
cofibration, then by \eqref{bnd.equi}, $\C'$ is bounded with respect
to $\gamma$ if and only if it is bounded over all the fibers $I'_i$,
$i \in I$, and if this happens, we have a natural identification
\begin{equation}\label{pd.cof}
(\gamma_*\C)_i \cong \Sec^\Tot(I'_i,\C)
\end{equation}
for any object $i \in I$. If any fibration $\C \to I'$ is bounded
with respect to $\gamma$, we say that $\gamma$ itself is bounded;
this happens, for example, if all the comma-fibers $I' /_\gamma I$
are essentially small, or if $\gamma$ is a cofibration with
essentially small fibers.

\begin{exa}\label{coev.exa}
For any category $I$ with a terminal object $o \in I$, the embedding
$\eps(o):\ppt \to I$ of \eqref{eps.c} is left-admissible, so that
any fibration $\C \to I$ is bounded with respect to the tautological
projection $\pi:I \to \ppt$. In fact, we have $\pi_*\C \cong
\eps(o)^*\C \cong \C_o$, and $\alpha$ of Corollary~\ref{kan.corr} is
then a functor
\begin{equation}\label{coev.eq}
\pi^*\C_o \cong I \times \C_o \to \C
\end{equation}
cartesian over $I$. This coincides with the functor $r$ of
Example~\ref{cyl.fib.exa}.
\end{exa}

\chapter{Orders.}\label{ord.ch}

This chapter is devoted to partially ordered sets and their
generalizations. The chapter is longish since we do need a lot of
small facts, but the facts themselves are easy (and not new). The
sets themselves are in Section~\ref{pos.sec}. Our general attitude
is that a partially ordered set is a special type of a small
category, and this turns out to be quite productive. The two classes
of maps between partially ordered sets that are important for our
applications are both defined categorically --- these are
``left-closed embeddings'' and ``reflexive maps'', finite
compositions of left and right-reflexive maps.

Throughout most of the text, we actually do not need arbitrary
partially ordered sets, and we restrict our attention to those
distinguished by some finiteness condition --- we formalise this by
introducing the notion of an ``ample subcategory'' $\I \subset \Pos$
in the category $\Pos$ of all partially ordered sets. Typical
conditions we need are formulated in terms of ``chain dimension''
and discussed in Subsection~\ref{pos.dim.subs}; the precise notion
of chain dimension that we use is in Definition~\ref{dim.def}. From
the point of view of general theory of partially ordered sets, the
sets we need are rather trivial and small. In fact, the only sets we
need that are not ``left-bounded'' in the sense of
Definition~\ref{pos.fin.def} are $\kappa$-directed partially ordered
sets for various regular cardinals $\kappa$ --- these are used for
the theory of accessible enhanced categories constructed in
Subsection~\ref{enh.large.sec}, and discussed in
Subsection~\ref{pos.size.subs}.

``Barycentric subdivision'' of Subsection~\ref{pos.bary.subs} is a
useful general construction that in particular turns any partially
ordered set into a left-finite one. We also need a bunch of specific
cocartesian squares and colimits in the category $\Pos$, and some
generalities on such squares; this is the subject of
Subsection~\ref{pos.perf.subs} and
Subsection~\ref{pos.cyl.subs}.

In Subsection~\ref{ano.subs}, we introduce another class of maps
between partially ordered sets that we will need, that of ``anodyne
maps''; roughly speaking, these are obtained from reflexive ones by
closing with respect to Quillen's Theorem A (see
Definition~\ref{pos.ano.def} for a precise formulation). All anodyne
maps between partially ordered sets induce weak equivalences between
their nerves. The converse is probably not true, but no matter;
these are the only weak equivalences that we will need. In
particular, while it is not true that the barycentric subdivision
sends reflexive maps to reflexive ones, it does send anodyne maps to
anodyne maps.

Section~\ref{bi.sec} starts with a generalization of some of the
material of Section~\ref{pos.subs} to $I$-augmented partially
ordered sets, for some category $I$. Then we introduce another
generalization that will be very useful for us --- the notion of a
``biordered set'' (roughly speaking, a partially ordered set
equipped with a factorization system). The theory of biordered sets
is the bulk of Section~\ref{bi.sec}; in particular, it includes
biordered versions of reflexive maps and lef-closed embedding, and
also a counterpart of the notion of an anodyne map that we call
``bianodyne''. In the last Subsection~\ref{aug.rel.subs}, we combine
the two stories and consider $I$-augmented biordered sets. We note
that the theory of biordered sets is not needed until approximately
Subsection~\ref{seg.sp.sec}, so it is safe to skip
Section~\ref{bi.sec} until then.

\section{Partially ordered sets.}\label{pos.sec}

\subsection{Generalities.}\label{pos.subs}

A {\em preorder} on a set $S$ is a binary relation $\leq$ that is
reflexive ($s \leq s$ for any $s \in S$) and transitive ($s_0 \leq
s_1$ and $s_1 \leq s_2$ implies $s_0 \leq s_2$), but not necessarily
antisymmetric. A {\em preordered set} is a set equipped with a
preorder. Equivalently, a preordered set is a small category with at
most one morphism between any two objects ($s_0 \leq s_1$ iff there
is a morphism $s_0 \to s_1$). As usual, a {\em partial order} on a
set $S$ is an antisymmetric preorder ($s_0 \leq s_1$ and $s_1 \leq
s_0$ imply $s_0=s_1$). A partially ordered set is a set equipped
with a partial order. A preordered set is partially ordered if and
only if the corresponding small category has exactly one object in
each isomorphism class, and as a category, any preordered set is
equivalent to a partially ordered one, unique up to a unique
isomorphism.

A {\em map} between preordered sets $S_0$, $S_1$ --- and in
particular, a map between partially ordered sets --- is a map $f:S_0
\to S_1$ that preserves the preorder relation, or equivalently, a
functor between the corresponding small categories. We note that
such a functor is automatically faithful. A map is {\em dense} if it
is dense as a functor, or equivalently, if it is a bijection of
sets. For a given set $S$, preorders on $S$ and dense maps between
them themselves form a preorder that in fact is a partial order. The
minimal preorder on $S$ is the discrete partial order ($s_0 \leq
s_1$ iff $s_0 = s_1$); whenever an abstract set $S$ is treated
as a preordered one, it is equipped with this discrete partial
order. The maximal preorder is the one with $s_0 \leq s_1$ for any
$s_0,s_1 \in S$; we denote $S$ with this preorder by $e(S)$ (if $S$
is not empty, $e(S)$ is equivalent to $\ppt$).

For any preordered set $J$, the opposite category $J^o$ is also a
preordered set, and so are the functor category $\Fun(I,J)$ for any
small $I$, and the categories $J^<$ resp.\ $J^>$ (explicitly, they
are obtained by adding a new smallest resp.\ largest element $o$ to
$J$). If $J$ is a partially ordered set, then so are $J^o$,
$\Fun(I,J)$, $J^<$ and $J^>$. Functors between partially ordered
sets are order-preserving maps $f:J' \to J$, and for any such map,
the cylinder $\Cyl(f)$ and the dual cylinder $\Cyl^o(f)$ of
Definition~\ref{cyl.def}, as well as the comma-categories $J' /_f
J$, $J \setminus_f J'$ of \eqref{rc.dia} are partially ordered
sets. A preordered set is discrete as a category iff the preorder is
discrete. Any partially ordered set is rigid in the sense of
Definition~\ref{iso.def}. A {\em full embedding} of partially
ordered sets is a map $f:J \to J'$ that is full as a
functor. Equivalently, $f$ is injective and identifies $J$ with its
image $f(J) \subset J'$ with the induced order (that is, $j \leq j'$
iff $f(j) \leq f(j')$). A subset $J \subset J'$ in a partially
ordered set $J'$ is always equipped with the induced order unless
indicated otherwise. For any map $f:J' \to J$ and element $j \in J$,
the fiber $J'_j$ is the preimage $f^{-1}(j) \subset J'$ (with the
induced order). A map is {\em conservative} if it is conservative as
a functor, or equivalently, if all its fibers are
discrete. Partially ordered sets are equivalent as categories iff
they are isomorphic, and functors between partially ordered sets are
isomorphic iff they are equal. We denote the category of
partially ordered sets by $\Pos$, and as in \eqref{iota.C.eq}, we
let $\iota:\Pos \to \Pos$ be the involution $J \mapsto J^o$. We have
the embedding $\phi:\Pos \to \Cat$ and the tautological cofibration
\begin{equation}\label{nu.pos}
\nu_\idot:\Pos_\idot = \phi^*\Cat_\idot \to \Pos
\end{equation}
whose fiber over some $J$ is $J$ itself. Explicitly, $\Pos_\idot$ is
the category of pairs $\langle J,j \rangle$, $J \in \Pos$, $j \in
J$, with maps $\langle J,j \rangle \to \langle J',j' \rangle$ given
by maps $f:J \to J'$ such that $f(j) \leq j'$. We denote by
$$
\nu^\hdot:\Pos^\hdot = \Pos^o_{\idot\perp} \to \Pos
$$
the transpose opposite cofibration to \eqref{nu.pos}. Then
$\Pos^\hdot \cong \iota^*\Pos_\idot$. Explicitly, $\Pos^\hdot$ is
the category of pairs $\langle J,j \rangle$, with maps $\langle J,j
\rangle \to \langle J',j' \rangle$ given by maps $f:J \to J'$ such
that $f(j) \geq j'$. We also have functors $\nu_>:\Pos_\idot \to
\Pos$ and $\nu_<:\Pos^\hdot \to \Pos$ given by
\begin{equation}\label{nu.io.pos}
\nu_>(\langle J,j \rangle) = J/j, \qquad \nu_<(\langle J,j \rangle)
= j \setminus J.
\end{equation}
Explicitly, for any $j \in J$, the comma-set $J/j \subset J$
consists of elements $j ' \in J$ such that $j' \leq j$. In
particular, $j \in J/j$ is the largest element, and if we denote
\begin{equation}\label{J.j.eq}
  J /' j = (J / j) \ssetminus \{j\},
\end{equation}
then $J/j \cong (J /' j)^>$. Dually, $j \setminus J \cong (j
\setminus' J)^<$, where $j \setminus' J = (j \setminus J) \ssetminus
\{j\}$.

A commutative square \eqref{cat.dia} of partially ordered sets is an
honest commutative square in $\Pos$, and it is cartesian as a
commutative square of categories iff it is a cartesian square in
$\Pos$. In particular, $\Pos$ has finite limits. It is also
cartesian-closed in the sense of Definition~\ref{cart.cl.def}, in
that for any $J,J' \in \Pos$, the partially ordered set $\Pos(J,J')$
of maps $f:J \to J'$ with pointwise order satisfies the correponding
universal property. For any $J \in \Pos$, the comma-category
$\Pos/J$ has finite limits, and we have a pair of adjoint functors
\begin{equation}\label{I.f.eq}
f_!:\Pos/J \to \Pos/J', \qquad f^*:\Pos/J' \to \Pos/J
\end{equation}
for any map $f:J \to J'$ in $\Pos$. For any $J \in \Pos$ and objects
$J_0,J_1 \in \Pos/J$ represented by maps $\chi_l:J_l \to J$,
$l=0,1$, a lax functor from $J_0$ to $J_1$ in the sense of
Subsection~\ref{fun.subs} is the same thing as an order-preserving
map $f:J_0 \to J_1$ such that $\chi_1 \circ f \geq \chi_0$ (since
$J$ is a partially ordered set, the map $\alpha:\chi_0 \to f \circ
\chi_1$ is unique if it exists, so that being lax is a condition and
not a structure). Dually, we will say that $f$ is {\em co-lax over
  $J$} if $\chi_1 \circ f \leq \chi_0$. Note that a lax map $f:J_0
\to J_1$ resp.\ a co-lax map $f':J_0 \to J_1$ induces natural maps
\begin{equation}\label{lax.pos.eq}
j \setminus f:j \setminus J_0 \to j \setminus J_1, \qquad f' / j:J_0
/ j \to J_1 / j
\end{equation}
for any element $j \in J$.

For any integer $n \geq 0$, we denote by $[n]$ the totally ordered
set of integers $\{0,\dots,n\}$ with the standard order
(equivalently, we can say that $[0] = \ppt$ is the point, and $[n+1]
= [n]^>$ for any $n \geq 0$). In particular, $[1]$ is the single
arrow category, so our notation is consistent. We denote by $\N$ the
partially ordered set of all integers $n \geq 0$, with the standard
order. The category $\V$ of \eqref{V.eq} is also a partially ordered
set, and we will treat it as such.

For any preordered set $J$, the $\Hom$-pairing \eqref{hom.pair.eq}
factors through the full subcategory $[1] \subset \Sets$ spanned by
$\{\emptyset\}$ and $\ppt$, as in
Example~\ref{closed.fib.exa}. Conversely, for any small category
$I$, we can define a preordered set $\cRed(I)$ by keeping the same
objects, and composing \eqref{hom.pair.eq} with the adjoint functor
$\Sets \to [1]$ of Example~\ref{1.sets.exa}. Explicitly, $\cRed(I)$
has exactly one map $i \to i'$ for any objects $i,i' \in I$ with
non-empty $I(i,i')$. Then $\cRed(I)$ is equivalent to a unique
partially ordered set $\Red(I)$ that we call the {\em reduction} of
$I$, and the tautological functor $I \to \cRed(I)$ descends to the
{\em reduction functor}
\begin{equation}\label{red.eq}
R:I \to \Red(I).
\end{equation}
Any functor $I \to J$ to a partially ordered set $J$ factors
uniquely through the reduction functor $R$, so that $\Red:\Cat \to
\Pos$ is left-adjoint to the full embedding $\Pos \to \Cat$.

\begin{exa}\label{kron.red.exa}
Let $\II$ be the Kronecker category of Example~\ref{kron.exa}. Then
$\Red(\II) \cong [1]$, with the reduction functor sending $0$ to
$0$, $1$ to $1$ and $s$, $t$ to the unique map $0 \to 1$ in
$[1]$. This example is universal: a small category $\C$ is a
preordered set iff any functor $\II \to \C$ factors through $R:\II
\to [1]$.
\end{exa}

\subsection{Left-closed embeddings.}

A partially ordered set $J$ is {\em directed} if it is not empty,
and for any $j_0,j_1 \in I$, there exists $j \in J$ such that
$j_0,j_1 \leq j$ -- or, equivalently, if $J$ is filtered as a
category.

\begin{exa}\label{U.S.exa}
For any non-empty set $S$, the set $P(S)$ of all subsets $S' \subset
S$ ordered by inclusion is directed. A map $f:S_0 \to S_1$ between
sets $S_0$, $S_1$ induces a map $f^*:P(S_1) \to P(S_0)$, $S' \mapsto
f^{-1}(S')$.
\end{exa}

\begin{defn}
A full embedding $f:J \to J'$ in $\Pos$ is {\em left-closed} iff for
any $j \leq j' \in J'$, $j' \in J$ implies $j \in J$, {\em
  right-closed} iff $f^o:J^o \to {J'}^o$ is left-closed, and {\em
  locally closed} if $f(J) \subset J'$ is the intersection of a
right-closed and a left-closed subset (or explicitly, if for any
$j_0 \leq j \leq j_1$, $j_0,j_1 \in f(J)$ we have $j \in f(J)$).
\end{defn}

\begin{exa}\label{st.cl.exa}
For any integers $n \geq m \geq 0$, the only left-closed embedding
$[m] \to [n]$ is the embedding $s:[m] \to [n]$ onto the initial
segment $\{0,\dots,m\} \subset \{0,\dots,n\} = [n]$ of the totally
ordered set $[n]$, and the only right-closed embedding is the
embedding $t = s^o:[m] \to [n]$ onto the terminal segment
$\{n-m,\dots,n\} \subset [n]$.
\end{exa}

\begin{exa}\label{sN.cl.exa}
For any $m \geq 0$, the only left-closed embedding $[m] \to \N$ is
the embedding $s:[m] \to \N$ onto the initial segment $\{0,\dots,m\}
\subset \N$ that identifies $[m]$ with the comma-set $\N/m$. The
right comma-fiber $m \setminus \N$ is identified with $\N$ itself,
by means of the right-closed embedding $q^m:\N \to \N$, where $q:\N
\to \N$ is the shift map $l \mapsto l+1$.
\end{exa}

\begin{exa}\label{well.exa}
More-or-less by definition, a partially ordered set $J$ is
well-ordered iff any non-empty right-closed subset $J' \subset J$ is
of the form $J' = j \setminus J$ for some $j \in J$ (namely, $j$ is
the minimal element in $J'$). In this case, any left-closed subset
$J' \subset J$ not equal to $J$ is of the form $J' = J /' j$, where
$J /'j$ is as in \eqref{J.j.eq}.
\end{exa}

In the categorical language, a full embedding $f:J \to J'$ is
left-closed iff it is left-closed in the sense of
Definition~\ref{cl.def}, and by Example~\ref{closed.fib.exa}, this
happens iff $f$ is a fibration. Dually, $f$ is right-closed iff it
is a cofibration. In the set-theorical language, for any subset $S
\subset S'$ in a set $S'$, we have a unique characteristic function
$\chi_S:S' \to [1]$ such that $S = \chi_S^{-1}(0)$, and $J \subset
J'$ is left-closed iff its characteristic function $\chi_J:J' \to
[1]$ is order-preserving. In this case, $\chi_J$ coincides with the
characteristic functor \eqref{chi.eq}. If we have left-closed
subsets $J_0 \subset J'_0$, $J_1 \subset J_1'$ with some
characteristic functions $\chi_0:J'_0 \to [1]$, $\chi_1:J_1' \to
[1]$, then an order-preserving map $f:J_0' \to J_1'$ sends $J_0$
into $J_1$ iff it is co-lax over $[1]$ (that is, $\chi_1 \circ f
\leq \chi_0$). For any left-closed $J \subset J'$, the complement
$J' \ssetminus J \subset J'$ is right-closed, and vice versa. For
any map $f:J' \to J$ and element $j \in J$, the left comma-fiber
$J'/j$ is left-closed in $J'$ (explicitly, $J'/j \subset J'$
consists of all $j' \in J'$ such that $f(j') \leq j$). An element $j
\in J$ is {\em maximal} iff $\{j\} \subset J$ is right-closed, and a
left-closed embedding $J_0 \subset J$ is {\em elementary} if $J
\ssetminus J_0$ consists of a single maximal element $j \in J$. More
generally, giving a filtration
\begin{equation}\label{cl.J.eq}
J_0 \subset \dots \subset J_n = J
\end{equation}
of a partially ordered set $J$ by left-closed subsets $J_i \in J$ is
equivalent to giving a map $\chi:J \to [n]$, with $J_i =
J/_{\chi}i$. We say that a filtration \eqref{cl.J.eq} is {\em dense}
if so is the map $\chi$, or equivalently, if all the left-closed
embeddings $J_i \subset J_{i+1}$ are elementary. A partially ordered
set $J$ admits a dense filration iff it is finite. Such a filtration
is in general not unique, and any filtration \eqref{cl.J.eq}
on a finite partially ordered set can be refined to a dense one.

For any partially ordered set $J$, the full embedding $J \to J^>$ is
left-closed and elementary, and it is universal in the following
sense: for any left-closed full embedding $J' \to J''$, a map $f:J'
\to J$ gives rise to a unique map $p:J'' \to J^>$ that fits into a
cartesian square
\begin{equation}\label{p.J.sq}
\begin{CD}
J' @>{f}>> J\\
@VVV @VVV\\
J'' @>{p}>> J^>.
\end{CD}
\end{equation}
If $J = \ppt$ and $J' \to \ppt$ is the tautological map, then $J^> =
\ppt^> = [1]$, and $p$ in \eqref{p.J.sq} is the characteristic map
$\chi:J'' \to [1]$ of the subset $J' \subset J''$. Dually, the
embedding $J \to J^<$ is right-closed, and enjoys a similar
universal property with respect to right-closed embeddings.

\begin{exa}\label{ar.S.exa}
Let $S$ be a discrete partially ordered set, and consider the set
$S^<$ and the arrow set $\Ar(S^<)$, with the projection
$\tau:\Ar(S^<) \to S^<$. Then we can identify $\Ar(S^<) \cong (S
\times [1])^<$, and under this identification, $\tau$ is universal
extension of the right-closed embedding $S \to (S \times [1])^<$
onto $S \times \{1\} \subset S \times [1] \subset (S \times [1])^<$.
\end{exa}

The intersection $J_{01} = J_0 \cap J_1 \subset J$ of two
left-closed subsets $J_0,J_1 \subset J$ is a left-closed subset, and
if $J = J_0 \cup J_1$, then the square \eqref{st.sq} of
Example~\ref{cl.exa} gives a cartesian cocartesian square
\begin{equation}\label{pos.st.sq}
\begin{CD}
J_{01} @>>> J_0\\
@VVV @VVV\\
J_1 @>>> J
\end{CD}
\end{equation}
in $\Pos$. We call such squares {\em standard pushout squares}. By
Example~\ref{V.closed.exa}, they correspond to maps $J \to \V$, and
by Example~\ref{st.exa}, for any partially ordered set $J_{01}$
equipped with left-closed full embeddings $J_{01} \to J_0$, $J_{01}
\to J_1$, there exists a coproduct $J = J_0 \copr_{J_{01}} J_1$ in
$\Pos$ that gives a standard square \eqref{pos.st.sq}. Morphisms
between two standard pushout squares correspond to morphisms co-lax
over $\V$. Alternatively, for any standard pushout square
\eqref{pos.st.sq}, one can consider the commutative square
\begin{equation}\label{pos.2.st.sq}
\begin{CD}
  J_{01} \copr J_{01} @>>> J_0 \copr J_1\\
  @VVV @VVV\\
  J_{01} @>>> J,
\end{CD}
\end{equation}
and this square is also cartesian and cocartesian. 

\begin{exa}\label{elem.exa}
If $J_0 = J \ssetminus \{j\} \subset J$ is an elementary left-closed
full embeding, then we have a standard pushout square
\eqref{pos.st.sq} with $J_1 = J/j$, and $J_{01} = J_0 \cap J/j =
J_0/j = J /'j$.
\end{exa}

As a generalization of standard pushout squares \eqref{pos.st.sq},
assume given partially ordered sets $J$, $J'$, and an
order-preserving map $\phi:J' \to J$. Then we have a functor $J \to
\Pos$, $j \mapsto J'/_\phi j$, and we have an isomorphism
\begin{equation}\label{lc.coli.eq}
  J' \cong \colim_{j \in J}J'/_\phi j.
\end{equation}
We call colimits of this form {\em standard colimits}. Dually, we
can also consider the functor $J^o \to \Pos$, $j \mapsto j
\setminus_\phi J'$, and then we have
\begin{equation}\label{lc.o.coli.eq}
  J' \cong \colim_{j \in J^o}j \setminus_\phi J'.
\end{equation}
These colimits are called {\em co-standard}.

\begin{lemma}\label{hom.I.le}
For any partially ordered set $I$ that has a largest element $o \in
I$, $\Pos(I,-)$ preserves the standard colimits \eqref{lc.coli.eq}.
\end{lemma}

\proof{} If we let $\phi(I) = \phi \circ \ev_o:\Pos(I,J') \to J$,
where $\ev_o$ is evaluation at $o \in I$, then
$\Pos(I,J')/_{\phi(I)}j \cong \Pos(I,J'/_\phi j)$ for any $j \in J$.
\endproof

For any partially ordered set $J$, we denote by $L(J) \subset P(J)$
the subset of left-closed subsets $J' \subset J$. Sending a subset
$J' \in L(J)$ to its characteristic function $\chi_{J'}:J \to [1]$
provides an isomorphism $L(J)^o \cong \Fun(J,[1])$, where we need to
take the opposite order on $L(J)$ since $J'_0 \subset J_1'$ iff
$\chi_{J'_0} \geq \chi_{J'_1}$. Alternatively, $L(J) \cong J^o[1]$,
with the isomorphism sending $J' \subset J$ to the characteristic
function $\chi:J^o \to [1]$ of the complement $(J \ssetminus J')^o
\subset J^o$. For any map $f:J_0 \to J_1$, the corresponding map
$f^*:P(J_1) \to P(J_0)$ of Example~\ref{U.S.exa} sends $L(J_1)
\subset P(J_1)$ into $L(J_0) \subset P(J_0)$; in terms of the
identifications $L(J_l) \cong J_l^o[1]$, $l=0,1$, the induced map
$L(J_1) \to L(J_0)$ is the pullback map $f^{o*}$. The $\Hom$-pairing
\eqref{hom.pair.eq} for any $J \in \Pos$ factors through $[1]
\subset \Sets$, and the Yoneda embedding \eqref{yo.yo.eq} reduces to
a map
\begin{equation}\label{yo.pos.eq}
  \Y:J \to L(J) \cong J^o[1], \qquad j \mapsto J / j.
\end{equation}
The gluing construction of Example~\ref{cl.glue.exa} applies to
partially ordered sets, and in this case, \eqref{cl.glue.bis.eq}
reduces to an isomorphism
\begin{equation}\label{pos.glue.eq}
J \cong J_0 \copr_\lambda J_1,
\end{equation}
where $J$ is a partially ordered set equipped with a map $J \to
[1]$, with left-closed fiber $J_0 \subset J$ and right-closed fiber
$J_1 \subset J$, glued along a map $\lambda:J_1 \to
L(J_0)$. Explicitly, $J \cong J_0 \copr J_1$, with an extra order
relation $j_0 \leq j_1$ for any $j_1 \in J_1$ and $j_0 \in
\lambda(j_1) \subset J_0$.

For any $J \in \Pos$, we denote by $L_\idot(J) \subset J \times
L(J)$ the ``incidence subset'' of elements $j \times J'$ such that
$j \in J'$. Equivalently, we have a tautological functor $\phi:L(J)
\to \Pos$, and $L_\idot(J) = \phi^*\Pos_\idot$, where $\Pos_\idot
\to \Pos$ is the cofibration \eqref{nu.pos}. In terms of the
identification $L(J) \cong J^o[1]$, we have the left-closed subset
$L^\hdot(J) \subset J^o \times J^o[1]$ corresponding to the
evaluation pairing $J^o \times J^o[1] \to [1]$; then $L_\idot(J)$ is
the same subset but considered as a subset in $J \times L(J)$. More
generally, for any $\J \in \Pos$ equipped with a map $e:\J \to
L(J)$, we let $\J_\idot = (e \times \id)^*L_\idot(J) \subset J
\times \J$; explicitly, $\J_\idot$ is the set of pairs $\langle j,j'
\rangle$ such that $j \in e(j')$. We denote by
\begin{equation}\label{s.t.J}
\sigma(J):\J_\idot \to J, \qquad \tau(J):\J_\idot \to \J
\end{equation}
the two projections, and we note that $\tau(J)$ is the cofibration
induced from that of Example~\ref{cat.exa}.

\begin{exa}\label{e.J.exa}
Take $\J = J$, and let $e=\Y$ be the Yoneda embedding
\eqref{yo.pos.eq}. Then $\J_\idot = \Ar(J)$, with the projections
$\sigma(J) = \sigma$, $\tau(J) = \tau$ of Example~\ref{ar.exa}.
\end{exa}

\subsection{Reflexive maps.}\label{pos.exa.subs}

Another general categorical notion that can be usefully applied to
partially ordered sets is that of left and right-reflexive functors
of Subsection~\ref{adj.subs}. Here is the corresponding definition.

\begin{defn}\label{po.refl.def}
An order-preserving map $f:J \to J'$ is {\em left} resp.\ {\em
  right-refl\-exive} if there exists a map $f_\dg:J' \to J$ such
that $f_\dg \circ f \leq \id$, $f \circ f_\dg \geq \id$
resp.\ $f_\dg \circ f \geq \id$, $f \circ f_\dg \leq \id$. A full
embedding is {\em reflexive} if it is a finite composition of left
and right-reflexive full embeddings.
\end{defn}

Categorically, $g$ is left resp.\ right-adjoint to $f$; in
particular, it is unique, so that reflexivity is a condition and not
a structure, and Definition~\ref{po.refl.def} is consistent with our
general categorical usage. If $f$ is a left or right-reflexive full
embedding, then $g \circ f = \id$. For any left or right-reflexive
full embedding $J \to J'$, the embedding $J \subset J''$ into any
$J'' \subset J'$ containing $J$ is tautologically left
resp.\ right-reflexive, with the same adjoint map $g$, and for any
partially ordered set $J_0$, the embedding $\id \times f:J_0 \times
J \to J_0 \times J'$ is left resp right-reflexive, with the adjoint
map $\id \times g$.

\begin{exa}\label{U.adj.exa}
For any map $f:S_0 \to S_1$, the map $f^*:P(S_1) \to P(S_0)$ of
Example~\ref{U.S.exa} is left-reflexive, with the adjoint
$f_!:P(S_0) \to P(S_1)$ sending a subset $S' \subset S_0$ to its
image $f(S') \subset S_1$. In fact, $f^*$ is also right-reflexive,
but the corresponding adjoint $f_*$ does not seem to have an
accepted name.
\end{exa}

\begin{lemma}\label{P.adj.le}
An order-preserving map $f:J' \to J$ is right-reflexive if and only
if $(f^*)^{-1}(\Y(J')) \supset \Y(J)$, where $\Y$ is the Yoneda
embedding \eqref{yo.pos.eq}, and $f^*:L(J) \to L(J')$ is induced by
$f^*:P(J) \to P(J')$ of Example~\ref{U.S.exa}.
\end{lemma}

\proof{} For any $j \in J$, we have $f^*(J/j) = J'/_fj$ and $f$ has
a right-adjoint $g$ iff all these sets have largest elements
(namely, $J'/_fj \cong J'/g(j)$).
\endproof

\begin{exa}
For any $J \in \Pos$, the embedding $L(J) \to P(J)$ is
left-reflexive, with the adjoint
\begin{equation}\label{La.eq}
  \Lambda:P(J) \to L(J), \qquad J' \mapsto \bigcup_{j \in J'}J / j =
  J \setminus J',
\end{equation}
where $J \setminus J'$ on the right is the right comma-set of the
embedding $J' \to J$.
\end{exa}

\begin{exa}\label{lc.exa}
For any map $g:J' \to J$ in $\Pos$, we have the
factorization
\begin{equation}\label{lc.pos}
\begin{CD}
J' @>{\eta}>> J' /_g J @>{\tau}>> J
\end{CD}
\end{equation}
opposite to \eqref{rc.dec}. Then $\tau$ is a cofibration, and $\eta$
is right-reflexive, with the right-adjoint map $\sigma:J' /_g J \to
J'$.
\end{exa}

\begin{exa}\label{n.refl.exa}
For any $n,m \geq 0$, a map $f:[n] \to [m]$ is right-reflexive if
and only if $f(0)=0$, with the adjoint map $f_\dg$ sending an
element $l \in [m]$ to $\max\{l' \in [n]|f(l') \leq l\}$. In
particular, if $n \geq m$, then the left-closed embedding $s:[m] \to
[n]$ of Example~\ref{st.cl.exa} is right-reflexive, with the adjoint
$s_\dg:[n] \to [m]$ given by $s_\dg(l) = l$ if $0 \leq l \leq m$,
and $m$ if $m \leq l \leq n$. Dually, $f:[n] \to [m]$ is
left-reflexive if and only if $f(n)=m$, with the adjoint $f_\dg$
sending $l \in [m]$ to $\min\{l' \in [n]|f(l') \geq l\}$, and in
particular, if $n \geq m$, the right-closed embedding $t:[m] \to
[n]$ is left-reflexive, with the adjoint map $t_\dg$ given by
$t_\dg(l) = n-m$, $0 \leq l \leq n-m$, and $t_\dg(l) = l$
otherwise. Note that $s_\dg$ is a fibration and a cofibration, and
so it $t_\dg$.
\end{exa}

\begin{exa}\label{n.m.refl.exa}
More generally, any surjective map $f:[n] \to [m]$ is both right and
left-reflexive, and it is both a fibration and a cofibration.
\end{exa}

\begin{exa}\label{n.N.refl.exa}
For any $m \geq 0$, the left-closed embedding $s:[m] \to \N$ of
Example~\ref{sN.cl.exa} is right-reflexive, with the adjoint map
$s_\dg:\N \to [m]$ sending $l \in \N$, $l < m$ to itself, and the
rest to $m$. The adjoint map $s_\dg$ is again a cofibration. The
embedding $\N \cong m \setminus \N \to \N$ is left-reflexive, with
the adjoint map sending $n$ to $\max\{0,n-m\}$.
\end{exa}

\begin{exa}\label{Z.refl.exa}
Let $\Z$ be the set of integers with the usual order. Then the
embedding $\eps(0):\ppt \to \Z$ onto $0 \in \Z$ is reflexive ---
indeed, it factors as $\ppt \to \N \to \Z$, the first embedding is
right-reflexive, and the second one is left-reflexive.
\end{exa}

\begin{exa}\label{J.N.refl.exa}
For any $J \in \Pos$ equipped with a map $\chi:J \to \N$, the map
$\id \times \chi:J \to J \times \N$ is reflexive. Indeed, by
\eqref{lc.pos}, it factors as
$$
\begin{CD}
  J @>{\eta}>> J /_\chi \N @>{\sigma \times \tau}>> J \times \N,
\end{CD}
$$
where $\eta$ is right-reflexive, $\sigma:J /_\chi \N \to J$ is a
fibration, and the embedding $\sigma \times \tau$ is cartesian over
$\N$, with fibers given by embeddings $\N \cong \chi(j) \setminus \N
\to \N$, $j \in J$, thus left-reflexive by
Lemma~\ref{fib.adj.le}~\thetag{ii} and Example~\ref{n.N.refl.exa}.
\end{exa}

A left-closed embedding $J \to J'$ is right-reflexive iff its
characteristic functor $\chi_J:J' \to [1]$ is a fibration. For any
map $f:J \to J'$, the cylinder $\Cyl(f)$ and the dual cylinder
$\Cyl^o(f)$ of Definition~\ref{cyl.def} are functorial with respect
to $f \in \Ar(\Pos)$, the projection $\Cyl(f) \to [1]$ is a
cofibration, the projection $\Cyl^o(f) \to [1]$ is a fibration, and
\eqref{cyl.deco} is a diagram of partially ordered sets, again
functorial with respect to $f$. If $f:J \to J'$ is a right-reflexive
full embedding, with the adjoint map $g$, then we have maps
\begin{equation}\label{refl.dia}
\begin{CD}
\Cyl(f) @>{\wt{f}}>> J' \times [1] @>{\wt{g}}>> \Cyl^o(g) \cong \Cyl(f)
\end{CD}
\end{equation}
over $[1]$ given by $f$ and $g$ over $0 \in [1]$ and $\id$ over $1$,
and $\wt{f}$ is a right-reflexive full embedding with the adjoint
map $\wt{g}$. In particular, $\Cyl(f)$ is a retract of $J' \times
[1]$, and the embeddings $s:J \to \Cyl(f)$ and $t:J' \to \Cyl(f)$
are then retracts of the embeddings $\id \times s,\id \times t:J'
\to J' \times [1]$. For a left-reflexive $f$, the picture is the
same, but one has to use $\Cyl^o(f) \cong \Cyl(g)$ and swap $s$ and
$t$.

\begin{lemma}\label{refl.sat.le}
The smallest saturated class of maps in $\Pos$ that contains the
projections $J \times [1] \to J$ for all $J \in \Pos$ also contains
all reflexive full embeddings.
\end{lemma}

\proof{} Since saturated classes are closed with respect to
compositions, it suffices to consider separately left and
right-reflexive full embeddings. For any right-reflexive full
embedding $f:J \to J'$, we have the retraction \eqref{refl.dia}, and
since $\id \times s,\id \times t:J' \to J' \times [1]$ are both
one-sided inverses to the projection $J' \times [1] \to J'$, they
lie in our saturated class, and then so do their retracts $s:J \to
\Cyl(f)$, $t:J' \to \Cyl(f)$. Then moreover, so does the one-sided
inverse $t_\dg:\Cyl(f) \to J'$, and by \eqref{cyl.deco}, so does
$f$. For a left-reflexive $f$ with adjoint $g$, replace $\Cyl(f)$
with $\Cyl(g)$ and swap $s$ and $t$.
\endproof

\begin{exa}\label{cyl.exa}
The cylinder $\Cyl(s)$ of the universal left-closed embedding $s:[0]
\to [1]$ is naturally identified with the ordinal $[2]$, and $s$ and
$t_\dg$ in the corresponding decomposition \eqref{cyl.deco} are as
in Example~\ref{n.refl.exa}. The left-adjoint functor
$t_{\dg\dg}:[1] \to [2] \cong \Cyl(s)$ provided by
Lemma~\ref{cl.cyl.le} is then given by the embedding $m:[1] \to [2]$
onto $\{0,2\} \subset [2]$. For any left-closed full embedding of
categories $\gamma:\C' \to \C$ with the characteristic functor
$\chi:\C \to [1]$, \eqref{cyl.deco} fits into a commutative diagram
\begin{equation}\label{m.t.dia}
\begin{CD}
\C' @>{s}>> \Cyl(\gamma) @>{t_\dg}>> \C\\
@VVV @VVV @VV{\chi}V\\
[0]=\ppt @>{s}>> [2] @>{t_\dg}>> [1]
\end{CD}
\end{equation}
with cartesian squares, and the functor $t_{\dg\dg}:\C \to
\Cyl(\gamma) \cong [2] \times_{[1]} \C$ of Lemma~\ref{cl.cyl.le} is
given by $t_{\dg\dg} = (m \circ \chi) \times \id$.
\end{exa}

\begin{lemma}\label{cyl.le}
Assume given a map $q:J_0 \to J_1$ in $\Pos$, and a left-closed
embedding $p:J_0' \to J_0$. Then the natural map $f:\Cyl(q \circ p)
\to \Cyl(q)$ between the cylinders of the maps $q \circ p:J_0' \to
J_1$, $q:J_0 \to J_1$ induced by $p$ is a left-reflexive full
embedding.
\end{lemma}

\proof{} By definition, the map $f$ is a full embedding, and the
complement to its image is naturally identified with $J_0 \ssetminus
p(J_0')$ by the map $s:J_0 \to \Cyl(q)$. Then the map $g:\Cyl(q) \to
\Cyl(q \circ p)$ adjoint to $f$ is equal to $\id$ on the image of
$f$, and it sends $s(j)$ for some $j \in J_0 \ssetminus p(J_0')$ to
$t(q(j))$.
\endproof

\begin{exa}\label{cl.refl.exa}
For any $J$ with the tautological map $q:J \to \ppt$, we have
$\Cyl(q) \cong J^>$. Then for any left-closed embedding $J_0 \to
J_1$, the corresponding map $J_0^> \to J_1^>$ is a left-reflexive
full embedding by Lemma~\ref{cyl.le}, with the adjoint map $J_1^>
\to J_0^>$ sending the whole $J_1^> \ssetminus J_0^>$ to $o \in
J_0^>$.
\end{exa}

\begin{exa}\label{V.exa}
Let $q:\{0,1\} = \ppt \copr \ppt \to \ppt$ be the codiagonal
map. Then we have $\Cyl(q) \cong \V^o$, where $\V$ is the partially
ordered set \eqref{V.eq}. Lemma~\ref{cyl.le} shows that the
embeddings $\{0\} \subset \V^o$, $\{1\} \subset \V^o$ are
reflexive. Indeed, the embeddings $[1] \cong \ppt^> \to \V^o \cong
\{0,1\}^>$ induced by the two embeddings $\ppt \to \{0,1\}$ are
left-reflexive by Lemma~\ref{cyl.le}, and the embedding $s:\ppt =
[0] \to \ppt^> = [1]$ is right-reflexive by
Example~\ref{n.N.refl.exa}.
\end{exa}

In general, a full embedding $J \to J'$ is reflexive iff there
exists a finite sequences of subsets
\begin{equation}\label{j.filt}
J = J'_0 \subset J'_1 \subset \dots \subset J'_m = J'
\end{equation}
for some $m \geq 0$ equipped with maps $g_l:J'_l \to J'_{l-1}
\subset J'_l$, $1 \leq l \leq m$ such that $g_l = \id$ on $J'_{l-1}
\subset J'_l$, and either $g_l \leq \id$ or $g_l \geq \id$. Here is
a typical example of such a situation.

\begin{exa}\label{Z.exa}
Let $\ZZ_\infty$ be the partially ordered set of all integers $l
\geq 0$, with order relations $2l \leq 2l+1 \geq 2(l+1)$, and no
other non-trivial order relations. For any $m \geq 0$, let $\ZZ_m
\subset \ZZ_\infty$ be the subset $\{0,\dots,m\} \subset \ZZ_\infty
\cong \N$. Then the embedding $\ppt = \ZZ_0 \subset \ZZ_m$ is
reflexive, with the sequence \eqref{j.filt} consisting of subsets
$\ZZ_l \subset \ZZ_m$, $0 \leq l \leq m$, and with the map
$g_l:\ZZ_l \to \ZZ_{l-1}$ sending $l$ to $l-1$. For $m=2$, we have
$\ZZ_2 \cong \V^o \cong \{0,1\}^>$, where $\V$ is as in
\eqref{V.eq}, while the embedding is onto $\{0\} \subset \V^o$, thus
coincides with the reflexive embedding of Example~\ref{V.exa}. We
also note that we have a natural map
\begin{equation}\label{zeta.eq}
\zeta:\ZZ_\infty \to \N
\end{equation}
sending $2l,2l-1$ to $l$, and for any $m \geq 1$, $l=0,1$,
\eqref{zeta.eq} restricts to a map
\begin{equation}\label{z.m.eq}
  \zeta_{2m-l}:\ZZ_{2m-l} \to [m] = \N / m
\end{equation}
that provides an isomorphism $\ZZ_{2m} \cong \ZZ /_\zeta m$ if
$l=0$. The right comma-fiber $m \setminus_\zeta \ZZ_\infty$ is
naturally identified with $\ZZ_\infty$, and we have a standard
pushout square
\begin{equation}\label{zz.m.sq}
\begin{CD}
  \ppt @>{\eps(2m)}>> \ZZ_{2m}\\
  @V{\eps(0)}VV @VVV\\
  \ZZ_\infty @>{p^n}>> \ZZ_\infty,
\end{CD}
\end{equation}
where $p:\ZZ_\infty \to \ZZ_\infty$ is the shift map $l \mapsto
l+2$. If $q:\N \to \N$ is the shift map $l \mapsto l+1$ of
Example~\ref{sN.cl.exa}, then $\zeta \circ p = q \circ \zeta$. One
can use \eqref{zz.m.sq} to obtain $\ZZ_\infty$ by gluing copies of
$\ZZ_2 \cong \V^o$, and in fact, one can do the whole construction
at once: if one considers the map
\begin{equation}\label{x.z.eq}
x:\ZZ_\infty \to \V = \{0,1\}^<, \quad 2n \mapsto o, \ 4n+1 \mapsto
0,\ 4n+3 \mapsto 1,\ n \in \N,
\end{equation}
then the alternative form \eqref{pos.2.st.sq} of the corresponding
standard pushout square is a cartesian cocartesian square
\begin{equation}\label{Z.V.N.sq}
\begin{CD}
\overline{\N} \times \{0,1\} @>{a}>> \overline{\N} \times \V^o\\
@VVV @VVV\\
\overline{\N} @>{a'}>> \ZZ_\infty,
\end{CD}
\end{equation}
where $\overline{\N}$ is $\N$ with discrete order, and the maps are
$a(n \times l) = (n+l) \times l$, $n \geq 0$, $l=0,1$, and $a'(n) =
2(n+1)$.
\end{exa}

\begin{remark}
The set $\ZZ_\infty$ can be visualized as an infinite zigzag; then
$\ZZ_m \subset \ZZ_\infty$ is the zigzag of length $m$.
\end{remark}

\begin{lemma}\label{po.copr.le}
If the left-closed full embedding $J_{01} \to J_0$ in a standard
pushout square \eqref{pos.st.sq} is reflexive, then so is the map
$J_1 \to J$.
\end{lemma}

\proof{} Consider the filtration \eqref{j.filt} for $J_{01} \subset
J_0$, with the corresponding maps $g_l$, $1 \leq l \leq m$, and
define the filtration and the maps for $J_1 \subset J$ by taking
pushouts with respect to the left-closed embedding $J_{01} \subset
J$.
\endproof

\begin{lemma}\label{Z.m.le}
For any $n \geq 0$ and cofibration $J \to [n]$, with the induced
cofibration $\zeta_{2n}^*J \to \ZZ_{2n}$, the embedding $J_n \cong
(\zeta_{2n}^*J)_{2n} \to \zeta_{2n}^*J$ onto the fiber over $2n \in
\ZZ_{2n}$ is reflexive.
\end{lemma}

\proof{} Note that \eqref{zz.m.sq} for $m=1$ induces a standard
pushout square $\ZZ_{2n} \cong \ZZ_2 \copr_\ppt \ZZ_{2(n-1)}$ that
gives a corresponding standard pushout square decomposition
\begin{equation}\label{z.m.J.eq}
  \zeta_{2n}^*J \cong \zeta_2^*s^*J \copr_{J_1} \zeta_{2(n-1)}^*t^*J,
\end{equation}
where $s:[1] \to [n]$, $t:[n-1] \to [n]$ are the standard left
resp.\ right-closed embeddings. Therefore by induction on $n$ and
Lemma~\ref{po.copr.le}, it suffices to consider the case $n=1$. Then
we have $\ZZ_2 = \V^o = \{0,1\}^>$, the map $\zeta_2:\V^o \to [1]$
sends $0$ to $0$ and $1,o$ to $1$, and $\ppt = \{1\} \subset \{1,o\}
\cong [1] \subset \V^o$ gives a filtration \eqref{j.filt} for the
reflexive embedding $\ppt \to \V^o$ onto $1 \in \V^o$. This induces
a filtration $J_1 \subset J_1 \times [1] \subset \zeta_2^*J$ on
$\zeta_2^*J$, the embedding $J_1 \times [1] \subset \zeta_2^*J$ is
left-reflexive by Lemma~\ref{adm.cof.le}, and the embedding $J_1
\subset J_1 \times [1]$ is $\id \times s$, where the map
$s:\ppt = [0] \to [1]$ is right-reflexive by
Example~\ref{n.refl.exa}.
\endproof

\begin{lemma}\label{filt.le}
For any reflexive full embedding $f:J \to J'$, the full embedding
$s:J \to \Cyl(f)$ is reflexive, and for any map $\chi:J \to \N$, so
is the full embedding $s:J \to \Cyl(\chi \times f)$.
\end{lemma}

\proof{} For the first claim, since the embedding $s:J \to J \times
[1] \cong \Cyl(\id)$ is right-reflexive, it suffices to show that the
embedding $\Cyl(\id) \to \Cyl(f)$ is reflexive. To see this, take
the filtration \eqref{j.filt} for $f$, let $f_l:J \to J_l$ be the
embeddings, and filter $\Cyl(f)$ by $\Cyl(f_l) \subset \Cyl(f)$, $1
\leq l \leq m$, with $g_l:\Cyl(f_l) \to \Cyl(f_{l-1})$ the same as
for $J_l$ over $1 \in [1]$, and $g_l=\id$ over $0 \in [1]$. For the
second claim, note that $\id \times f:\N \times J \to \N \times J'$
is also a reflexive full embedding, and then $\chi \times \id:J \to
\N \times J$ is a reflexive full embedding by Example~\ref{J.N.refl.exa}.
\endproof

\subsection{Barycentric subdivision.}\label{pos.bary.subs}

For any set $S$, with the partially ordered set $P(S)$ of
Example~\ref{U.S.exa}, we let $V(S) \subset P(S)$ be the subset
spanned by finite subsets, and we let $W(S) \subset V(S)$ be the
subset spanned by finite non-empty subsets. Like $P(S)$, both $V(S)$
and $W(S)$ are directed, but unlike $P(S)$, they are not
contravariantly functorial with respect to $S$ (the preimage of a
finite set is not necessarily finite). However, the image of a
finite and/or non-empty set is finite and/or non-empty, so that all
the sets $W(S) \subset V(S) \subset P(S)$ are covariantly
functorial, with the functorialiy given by the adjoint maps $f_!$ of
Example~\ref{U.adj.exa}.

\begin{defn}\label{pos.bary.def}
The {\em barycentric subdivision} $BJ$ of a partially ordered set
$J$ is the subset $BJ \subset W(J)$ spanned by finite non-empty
subsets $S \subset J$ such that the induced order on $S$ is total.
\end{defn}

Since a surjective image of a finite totally ordered set is totally
ordered, $B(J)$ is functorial with respect to $J$, and since the
opposite to a total order is total, we have $B(J) \cong B(J^o)$. For
any injective map $J_0 \to J_1$, the corresponding map $BJ_0 \to
BJ_1$ is a left-closed full embedding. We also have functorial maps
\begin{equation}\label{B.xi}
\xi:BJ \to J, \qquad \xi_\perp:(BJ)^o \to J
\end{equation}
sending a totally ordered subset $S \subset I$ to its maximal
element $\max(S)$ resp.\ minimal element $\min(S)$. More generally,
denote by $B_\idot(J)$ the partially ordered set of pairs $\langle
S,s \rangle$, $S \in B(J)$, $s \in S$, with the order $\langle S,s
\rangle \leq \langle S',s' \rangle$ iff $S \subset S'$ and $s \leq
s'$. Then the forgetful map $\nu_\idot:B_\idot(J) \to B(J)$,
$\langle S,s \rangle \mapsto S$ is a cofibration corresponding to
the tautological functor $\phi:B(J) \to \Pos$ sending $S \in B(J)$
to $S$ considered as an object in $\Pos$. The cofibration
$\nu_\idot$ has a left resp.\ right-adjoint $\nu_\dg,\nu^\dg:B(J)
\to B_\idot(J)$ sending an element $S \in B(J)$ to $\langle
S,\min(S) \rangle$ resp.\ $\langle S,\max(S) \rangle$, and $\nu_\dg$
resp.\ $\nu^\dg$ has a left resp.\ right-adjoint
$\nu_>,\nu_<:B_\idot(J) \to B(J)$ sending $\langle S,s \rangle$ to
$S / s$ resp.\ $s \setminus S$, as in \eqref{nu.io.pos}. We then
also have the forgetful map
\begin{equation}\label{B.xi.dot}
\xi_\idot:B_\idot(J) \to J, \qquad \langle S,s \rangle \mapsto s,
\end{equation}
and it is related to \eqref{B.xi} by $\xi_\idot \cong \xi \circ
\nu_>$ and $\xi = \xi_\idot \circ \nu^\dg$.

\begin{exa}\label{BV.exa}
For any $n \geq 0$, $B([n])=W([n]) \cong [1]^n \ssetminus \{0^n\}$
consists of all non-empty subsets in $[n]$. The map $\xi:B([n]) \to
[n]$ is right-reflexive, with the adjoint $\xi_\dg:[n] \to B([n])$
sending $l \in [n]$ to $[n]/l = s([l]) \subset [n]$. For $n=1$, we
have $B([1]) \cong \V^o$, where $\V$ is the set of
Example~\ref{V.exa}, and the map $\xi:\V^o \cong B([1]) \to [1]$ of
\eqref{B.xi} is the map $\zeta_2:\V^o \cong \ZZ_2 \to [1]$ of
\eqref{z.m.eq}. For any $n \geq 0$, the map $\zeta_{2n}:\ZZ_{2n} \to
      [n]$ factors as
\begin{equation}\label{zeta.n.xi}
\begin{CD}
  \ZZ_{2n} @>{\beta_n}>> B([n]) @>{\xi}>> [n],
\end{CD}
\end{equation}
where $\beta_n:\ZZ_{2n} \to B([n])$ is a left-closed embedding that
identifies $\ZZ_{2n}$ with the subset in $B([n])$ formed by subsets
$\{l\} \subset [n]$, $0 \leq l \leq n$ and $\{l,l+1\} \subset [n]$,
$0 \leq l < n$. The map $\zeta$ of \eqref{zeta.eq} factors as
\begin{equation}\label{zeta.xi}
\begin{CD}
  \ZZ_\infty @>{\beta}>> B(\N) @>{\xi}>> \N,
\end{CD}
\end{equation}
where again, $\beta$ identifies $\ZZ_\infty$ with the left-closed
subset in $B(\N)$ formed by subsets $\{l\}$ and $\{l,l+1\}$, $l \geq
0$.
\end{exa}

\begin{exa}\label{W.exa}
The barycentric subdivision $B(\V)$ of the partially ordered set
$\V$ of \eqref{V.eq} is identified with the opposite $\W^o$ of the
partially ordered set
\begin{equation}\label{W.eq}
\W = ([1] \times \{0,1\}) \copr_{\{0\} \times \{0,1\}} \V^o,
\end{equation}
and under this identification, the map $\xi_\perp$ of
\eqref{B.xi} sends $1 \times l \in [1] \times \{0,1\}$ to $l \in
\V$, $l= 0,1$, and the whole $\V^o$ to $o \in \V$.
\end{exa}

\begin{exa}\label{B.S.exa}
More generally, for any set $S$ with discrete order, let $S^\hush$ be
the partially ordered set
\begin{equation}\label{B.S}
S^\hush = (S \times [1]) \copr_{\{0\} \times S} S^>.
\end{equation}
Then we have $S^\hush \cong B(S^<)^o$, and the projection
\begin{equation}\label{xi.B}
\xi_\perp:S^\hush \cong B(S^<)^o \to S^<
\end{equation}
of \eqref{B.xi} sends $S \times 1 \subset S \times [1]$ to $S
\subset S^<$, and the rest to $o \in S^<$. The projection $\xi^o:S^\hush
\to S^{<o} \cong S^>$ is equal to the identity on $S^>$ in
\eqref{B.S}, and sends $S \times [1]$ in \eqref{B.S} to $S \subset
S^>$ via the projection $S \times [1] \to S$ onto the first
factor. The product
\begin{equation}\label{xi.BS}
\xi^o \times \xi_\perp:S^\hush \to S^> \times S^<
\end{equation}
is a left-closed embedding identifying $S^\hush$ with $\delta(S)
\cup (S^> \times \{o\}) \subset S^> \times S^<$, where $\delta:S \to
S \times S \subset S^> \times S^<$ is the diagonal
map. Alternatively, $S^\hush \subset S^> \times S^< \cong S^> \times
S^{>o}$ is the twisted arrow set $\Tw(S^>)$.
\end{exa}

\begin{lemma}\label{B.cl.le}
For any map $f:J' \to J$ in $\Pos$, the induced map $B(f):B(J') \to
B(J)$ is a fibration.
\end{lemma}

\proof{} By definition, for any element in $B(J)$ represented by a
totally ordered subset $S \subset J$, the fiber $B(J')_S$ consists
of non-empty totally ordered subsets $S' \subset J'$ such that
$f(S') = S$. Then for any order relation $S_0 \subset S_1$ in
$B(J)$, the transition functor $B(J')_{S_1} \to B(J')_{S_0}$ sends
$S' \in B(J')_{S_1}$ to $S' \cap f^{-1}(S_0) \subset S'$.
\endproof

\begin{exa}\label{B.cl.exa}
For any $J \in \Pos$ equipped with a map $f:J \to [1]$ with fibers
$J_0$, $J_1$, we have $B([1]) = \V^o = \{0,1\}^>$, and then the
fibration $B(f)$ of Lemma~\ref{B.cl.le} has fibers $B(J)_0 \cong
B(J_0)$, $B(J)_1 \cong B(J_1)$, while the fiber $B(J)_o$ consist of
totally ordered subsets $S \subset J$ such that the induced map $f:S
\to [1]$ is surjective. The transition functor $B(J)_o \to B(J)_l$,
$l=0,1$ sends such a subset $S$ to the intersection $S \cap J_l$.
\end{exa}

\begin{exa}\label{B.X.exa}
Since the map \eqref{B.xi} is functorial, we have a commutative
square
\begin{equation}\label{B.X.sq}
\begin{CD}
  BJ' @>{\xi}>> J'\\
  @V{B(f)}VV @VV{f}V\\
  BJ @>{\xi}>> J
\end{CD}
\end{equation}
for any map $f:J' \to J$ in $\Pos$. If $f$ itself is a discrete
fibration, then the square \eqref{B.X.sq} is cartesian, and the
fibration $B(f)$ of Lemma~\ref{B.cl.le} is discrete.
\end{exa}

\begin{lemma}\label{B.corr}
The barycentric subdivision functor sends standard push\-out squares
\eqref{pos.st.sq} to standard pushout squares.
\end{lemma}

\proof{} By Example~\ref{closed.fib.exa}, any left-closed embedding
$J' \to J$ is a discrete fibration, so the square \eqref{B.X.sq} of
Example~\ref{B.X.exa} is cartesian -- that is, we have a natural
identification
\begin{equation}\label{B.pre.eq}
J' \times_J B(J) \cong B(J'),
\end{equation}
where $B(J) \to J$ is the map $\xi$ of \eqref{B.xi}. Now for any
square \eqref{pos.st.sq}, compose $\xi$ with
the joint characteristic function $\chi:J \to \V$ to obtain a
projection $B(J) \to \V$, and compute its left comma-fibers by
\eqref{B.pre.eq}.
\endproof

\begin{lemma}\label{B.le}
Assume given $J \in \Pos$, and let $\eps=\eps(o):\ppt \to J^>$ be
the embedding \eqref{eps.c} onto the new maximal element $o \in
J^>$. Then $B(\eps):\ppt \cong B(\ppt) \to B(J^>)$ is reflexive.
\end{lemma}

\proof{} Let $\chi:J^> \to [1]$ be the characteristic functor of the
left-closed embedding $J \subset J^>$, with the
fibration $B(\chi):B(J^>) \to \V^o \cong B([1])$ of
Example~\ref{B.cl.exa}, and note that the composition
\begin{equation}\label{B.le.eq}
\begin{CD}
  B(J^>) @>{B(\chi)}>> B([1]) @>{\xi}>> [1]
\end{CD}
\end{equation}
is a cofibration, with fibers $B(J^>)_0 \cong B(J)$, $B(J^>)_1
\cong B(J)^<$ and the transition functor $f:B(J) \to B(J)^<$ given
by the standard right-closed embedding. Therefore $B(J^>) \cong
\Cyl(f)$, and $B(\eps)$ is the composition of the left-reflexive
embedding $t:B(J)^< \to \Cyl(f)$ and the right-reflexive embedding
$\eps(o):\ppt \to B(J)^<$.
\endproof

\begin{corr}\label{xi.B.corr}
For any partially ordered set $J$ and functor $E:J \to \E$ to some
category $\E$, the Kan extension $\xi_!\xi^*E$ exists and is
isomorphic to $E$.
\end{corr}

\proof{} By \eqref{kan.eq}, we need to check that
$\colim_{BJ/j}\xi^*E \cong E(j)$ for any element $j \in J$. Since
$BJ / j \cong B(J/j)$, we may replace $J$ with $J/j$ and assume that
$j$ is the largest element. In other words, we are in the situation
of Lemma~\ref{B.le}: we have a partially ordered set $J$ equipped
with a functor $E:J^< \to \E$, and we need to check that
$\colim_{B(J^<)}\xi^*E \cong E(o)$. But then as in the proof of
Lemma~\ref{B.le}, we have the left-reflexive embedding $B(J^<)^1
\subset B(J^<)$, so that $\colim_{B(J^<)}\xi^*E \cong
\colim_{B(J^<)^1}\xi^*E$ by \eqref{adm.eq}, and then $\xi$ sends
$B(J^<)^1$ into $\{o\} \subset J^<$, so that $\xi^*E|_{B(J^<)^1}$ is
constant with value $E(o)$. The right-reflexive embedding $\ppt
\cong B(J^<)_1 \subset B(J^<)^1$ creates the initial element in
$B(J^<)^1$, so that $\colim_{B(J^<)^1}\xi^*E \cong E(o)$.  \endproof

\subsection{Dimension.}\label{pos.dim.subs}

For most of our intended applications, the category $\Pos$ is too
big; we will need to restrict our attention to some full subcategory
$\I \subset \Pos$ that is still large enough for our purposes. The
precise meaning of ``large enough'' is as follows.

\begin{defn}\label{amp.def}
A full subcategory $\I \subset \Pos$ is {\em ample} if it is closed
under finite limits and standard pushouts \eqref{pos.st.sq},
contains all finite partially ordered sets, and contains any subset
$J' \subset J$ of some $J \in \I$. An ample full subcategory $\I
\subset \Pos$ is {\em very ample} if it contains $\Cyl(f)$ and
$\Cyl^o(f)$ for any morphism $f:J' \to J$ in $\I$. The {\em
  extension} $\I^+$ of a full subcategory $\I \subset \Pos$ is the
full subcategory spanned by $J \in \Pos$ that admit a map $\chi:J
\to \N$ with $J/n \in \I$, $n \in \N$.
\end{defn}

\begin{lemma}\label{amp.le}
An ample full subcategory $\I \subset \Pos$ contains $J /_f J'$ and
$J' \setminus_f J$ for any map $f:J \to J'$ in $\I$, and
$\Pos(J',J)$ for any $J \in \I$ and finite $J' \in \Pos$. The
minimal saturated class of maps in $\I$ that contains the
projections $[1] \times J \to J$ for any $J \in \I$ also contains
all reflexive full embeddings $J' \to J$ in $\I$. The extension
$\I^+$ of an ample full subcategory $\I \subset \Pos$ is ample. The
extension $\I^+$ of a very ample full subcategory $\I \subset \Pos$
contains the sets $\ZZ_\infty$ and $\N$ of \eqref{zeta.eq}.
\end{lemma}

\proof{} For the first claim, we have $\Pos(J',J) \subset J^n$ for
$n=|J'|$, so $\Pos(J',J) \in \I$, and then $\Ar(J) \in \I$, and since
$\I$ is closed under fibered products, so are the comma-sets
$J/_fJ'$ and $J'\setminus_fJ$. For the second claim, the argument is
exactly as in Lemma~\ref{refl.sat.le} (where $\Cyl(f),\Cyl^o(f)
\subset J \times [1]$ must be in $\I$ for any full embedding $f:J'
\to J$). For the third claim, giving a map $J \to \N$ for some $J
\in \Pos$ is equivalent to representing $J$ as the colimit
\eqref{N.coli.eq} of an exhaustive filtration by left-closed subsets
$J(n) = J/n \subset J$; say that such a filtration is {\em good} if
$J(n) \in \I$, $n \in \N$. Then for any standard pushout $J = J_0
\cup J_1$, with good filtrations $J_l(-)$, $l = 0,1$, $J(n) = J_0(n)
\cup J_1(n) \subset J$, $n \in \N$ is a good filtration on $J$, so
that $\I^+$ is closed under standard pushouts, and conversely, for
any $J_0,J_1,J \in \I^+$ with good filtrations, and maps $f_l:J_l
\to J$, $l=0,1$, the filtrations by $J_l(n)' = J_l(n) \cap
f_l^{-1}(J(n)) \subset J_l$, $l=0,1$, $n \in \N$ are also good, and
then $J_0(n)' \times_{J(n)} J_1(n)' \subset J_0 \times_J J_1$ is a
good filtration, so that $J_0 \times_J J_1 \in \I^+$. For the last
claim, note that $\N/n$ and $\ZZ_\infty /_\zeta n$ are finite for
any $n \in \N$.
\endproof

The category $\Posff \subset \Pos$ of all finite partially ordered
sets is very ample in the sense of Definition~\ref{amp.def}, and by
Definition~\ref{amp.def}, it is the smallest possible very ample
subcategory in $\Pos$. Other useful very ample categories are
defined in terms of the following notion of ``dimension''.

\begin{defn}\label{dim.def}
A partially ordered set $J$ has {\em chain dimension $\leq n$} for
some integer $n \geq 0$ if for any injective map $f:[m] \to I$, $m
\geq 0$, we have $m \leq n$. A partially ordered set $J$ has {\em
  finite chain dimension} if it has chain dimension $\leq n$ for
some $n \geq 0$. In this case, $\dim J$ is the smallest such integer
$n$.
\end{defn}

In keeping with our general categorical usage, a non-degenerate
chain \eqref{c.idot} in $J$ is the same things as an injective map
$[n] \to J$, and then $\dim J \leq n$ iff all chains of length $>n$
are degenerate. Equivalently, we have $\dim J \leq n$ if there
exists a conservative map $J \to [n]$, and the converse is also true
--- for example, one can take the height map $\hht$ defined by
$\hht(j) = \dim J/j$. In fact, for any partially ordered set $J$, we
can define the {\em $n$-th skeleton} $\sk_nJ \subset J$ as the
left-closed subset of elements $j \in J$ such that $\dim J/j \leq
n$. Then $\dim \sk_n J = n$, and $\dim J = n$ iff $J = \sk_n J$. In
this case, we have the exhaustive {\em skeleton filtration}
\begin{equation}\label{sk.J}
\sk_0J \subset \sk_1 J \subset \dots \subset \sk_nJ = J,
\end{equation}
and this is the filtration \eqref{cl.J.eq} corresponding to the
height map $\hht$. Since we will be working mostly with partially
ordered sets of finite chain dimension, we denote the category they
form simply by $\Posf \subset \Pos$.

\begin{defn}\label{pos.fin.def}
A map $f:J' \to J$ between partially ordered sets is {\em
  left-finite} resp.\ {\em left-bounded} if its left comma-fibers
$J'/j$, $j \in J$ are finite resp.\ have finite chain dimension. A
partially ordered set $J$ is left-finite or left-bounded if so is
the identity map $\id:J \to J$.
\end{defn}

\begin{exa}\label{B.dim.exa}
The barycentric subdivision $BJ$ of any partially ordered set $J$ is
left-finite, and finite-dimensional iff so is $J$.
\end{exa}

\begin{exa}\label{f.b.exa}
For any cofibration $J \to J'$ with left-finite $J'$ and
left-bounded fibers $J_j$, $j \in J'$, the set $J$ is
left-bounded. Note that just asking for $J'$ to be left-bounded is
not enough.
\end{exa}

We denote the full subcategory of left-bounded partially ordered
sets by $\Posf^+ \subset \Pos$; this the the extensions of $\Posf$
in the sense of Definition~\ref{amp.def}, thus ample by
Lemma~\ref{amp.le}. We let $\Posf^- = \iota(\Posf^+) \subset \Pos$,
and we call sets $J \in \Posf^-$ {\em right-bounded}. We denote by
$\Posf^\pm \subset \Posf^+$ the full subcategory of left-finite
sets; it is also ample. For any left-bounded $J \in \Posf^+$, the
height function is still a well-defined conservative function
$\hht:J \to \N$, with $\sk_nJ = J /_{\hht} n$ of dimension $n$, and
we have the exhaustive filtration \eqref{sk.J} but it is now
infinite. A left-closed embedding $J \to J'$ in $\Posf^+$ preserves
the heights. Dually, for a right-bounded $J \in \Posf^-$, we have
the co-height function $\hht^\iota:J \to \N^o$ opposite to $\hht$
for $J^o$, and again, $\dim(n \setminus_{\hht^\iota} J) = n$ for any
$n \in \N^o$, and right-closed embeddings preserve the co-heights.

\begin{exa}
While the extension $\I^+$ of a very ample full subcategory $\I
\subset \Pos$ is ample by Lemma~\ref{amp.le}, it need not be very
ample -- indeed, the cylinder $\Cyl(f)$ of a map $f:J \to J'$ in
$\I^+$ is in $\I^+$ if $f$ is an embedding and/or $J \in \I$, but in
general, this need not be the case (for example, we have $\N^> \in
\I^+$ only if $\N \in \I$, and the latter is not true for
$\I=\Posff,\Posf_\kappa,\Posf$). Thus \eqref{cyl.deco} is not a
decomposition in $\I^+$. As an alternative, one can choose a map
$\chi:J \to \N$ with $J/n \in \I$, $n \in \N$, and consider the map
$\chi \times f:J \to \N \times J'$. Then $\Cyl(\chi \times f) \in
\I^+$, and $f$ decomposes as
\begin{equation}\label{cyl.N.deco}
\begin{CD}
  J @>{s}>> \Cyl(\chi \times f) @>>> J',
\end{CD}
\end{equation}
where the second map is the composition of $t_\dg:\Cyl(\chi \times
f) \to \N \times J'$ and the projection $\N \times J' \to J'$.
\end{exa}

For any $J \in \Posf$, we have $\dim J = \dim J^o = \dim BJ$, so
that the subcategory $\Posf \subset \Pos$ is closed both under the
barycentric subdivision functor $B$ and under the involution
$\iota:\Pos \to \Pos$, $J \mapsto J^o$. For any standard pushout
square \eqref{pos.st.sq}, we have $\dim J_{01} \leq \dim J_0,\dim
J_1$, and $\dim J = \max(\dim J_0,\dim J_1)$.

If $\dim J = 0$ then $J$ is discrete. There are various results for
partially ordered sets of small chain dimension; for example, we
have the following generalizations of Example~\ref{B.S.exa} and
Lemma~\ref{car.epi.le}.

\begin{lemma}
For any $J \in \Pos$ of dimension $\dim J = 1$, the product $\xi^o
\times \xi_\perp:(BJ)^o \to J^o \times J$ of the maps $\xi^o$ and
$\xi_\perp$ identifies $(BJ)^o$ with the twisted arrow set $\Tw(J^o)
\subset J^o \times J$.
\end{lemma}

\proof{} Since $\dim J = 1$, a totally ordered subset $S \subset J$
consists of its maximal and minimal element (that can
coincide). Therefore $\xi^o \times \xi_\perp$ is an embedding. Its
image consists of pairs $j \times j' \in J^o \times J$ such that $j
\geq j'$, and this is $\Tw(J^o)$ on the nose.
\endproof

\begin{lemma}\label{pos.epi.le}
Assume given a partially ordered set $J$ equipped with fibrations
$\C,\C' \to J$, and a functor $\gamma:\C \to \C'$ over $J$ such that
$\gamma(j):\C_j \to \C'_j$ is an equivalence resp.\ epivalence
resp.\ essentially surjective when $\hht j = 0$ resp.\ $1$
resp.\ $2$. Then $\gamma:\Sec^{\Tot}(\sk_n J,\C) \to
\Sec^{\Tot}(\sk_n J,\C')$ is an equivalence resp.\ epivalence
resp.\ essentially surjective for $n = 0$ resp.\ $1$ resp.\ $2$.
\end{lemma}

\proof{} The statement for $n=0$ is obvious. Lemma~\ref{car.epi.le}
is essentially the statement for $n=1,2$ and $J=\V=\{0,1\}^<$,
$J=\V^<=(\{0,1\}^<)^<$, and the statement is just as obvious for any
set $S$ and $S^<$, $(S^<)^<$. In the general case, for $n=1$,
replace $J$ with $\sk_1 J$, let $S \subset J$ be the discrete
right-closed subset of its maximal elements, and extend the identity
map $\id:S \to S$ to the universal map $p:J \to S^<$ of
\eqref{p.J.sq}. Then by \eqref{2kan.eq}, $p_*(\gamma):p_*\C \to
p_*\C'$ also satisfies the conditions of the lemma, so we are
reduced to the obvious case $J=S^<$. For $n=2$, again replace $J$
with $\sk_2 J$, with the subset $S \subset J$ of maximal elements,
and let $q:J \to (S^<)^<$ be the composition of the universal map
$p:J \to S^<$ and the embedding $S^< \to (S^<)^<$. Then the claim
for $n=0,1$ and \eqref{2kan.eq} shows that $q_*(\gamma)$ satisfies
the conditions of the lemma, and we are again reduced to the obvious
case $J=(S^<)^<$.
\endproof

A partially ordered set $J$ is {\em Noetherian} --- or equivalently,
satisfies the {\em descending chain condition} --- iff there are no
conservative maps $\N^o \to J$. Any $J \in \Posf^+$ is Noetherian,
and so is any well-ordered partially ordered set $Q$; consequently,
any $J \in \Pos$ that admits a conservative map $J \to Q$ to a
well-ordered $Q$ is Noetherian as well. Any full subset $J' \subset
J$ in a Noetherian $J$ is Noetherian. Any Noetherian set $J$
contains a minimal element $j \in J$.

\subsection{Size and direction.}\label{pos.size.subs}

For any partially ordered set $J$, we denote by $|J|$ the
cardinality of the underlying set, and we note that if $J$ is
infinite, it coincides with the cardinality of the corresponding
small category. For any regular cardinal $\kappa$, we say that $J
\in \Pos$ is {\em $\kappa$-bounded} if $|J| < \kappa$, and we denote
by $\Pos_\kappa \subset \Pos$ the full subcategory spanned by all
$\kappa$-bounded sets. We say that $J \in \Pos$ is {\em
  left-$\kappa$-bounded} if $|J/j| < \kappa$ for any $j \in J$.  We
let $\Posf_\kappa = \Posf \cap \Pos_\kappa \subset \Posf$, and we
note that this category is very ample in the sense of
Definition~\ref{amp.def}, with ample extension $\Posf^+_\kappa =
\Posf^+ \cap \Pos_\kappa$.

\begin{defn}\label{pos.filt.def}
For any regular cardinal $\kappa$, a partially ordered set $J$ is
{\em $\kappa$-directed} iff for any $\kappa$-bounded subset $S
\subset J$, $|S| < \kappa$, there exists $j \in J$ such that $S
\subset J/j$.
\end{defn}

\begin{exa}
For any set $S$, the subset $P_\kappa(S) \subset P(S)$ of
$\kappa$-bounded subsets $S' \subset S$, $|S'| < \kappa$ is
$\kappa$-directed.
\end{exa}

\begin{exa}\label{ord.exa}
For any infinite cardinal $\kappa$, there exists a
left-$\kappa$-bounded well-ordered partially ordered set $Q_\kappa$
of cardinality $|Q_\kappa|=\kappa$ -- take any well-ordered $S$ with
$|S|=\kappa$, and if it is not left-$\kappa$-bounded, replace it
with $S/'s$, where $s \in S$ is the minimal element such that
$|S/s|=\kappa$. If $\kappa$ is regular, then any such $Q_\kappa$ is
$\kappa$-directed. Indeed, for any $\kappa$-bounded $S \subset
Q_\kappa$, the left-closed subset $\Lambda(S) \subset Q_\kappa$ of
\eqref{La.eq} is then also $\kappa$-bounded, so that $Q_\kappa
\ssetminus \Lambda(S)$ is not empty. If $\kappa$ is the countable
cardinal, then $Q_\kappa=\N$.
\end{exa}

If $\kappa$ is the countable cardinal, then $\Posf_\kappa = \Posff$,
and an easy induction on $|S|$ shows that $J \in \Pos$ is
$\kappa$-directed iff it is directed. If $\kappa' \geq \kappa$, then
a $\kappa'$-directed $J \in \Pos$ is trivially $\kappa$-directed, so
that any $\kappa$-directed $J$ is directed. A directed partially
ordered set is connected. A full subset $J' \subset J$ in a
partially ordered set $J$ is {\em cofinal} if the embedding $J' \to
J$ is $0$-cofinal in the sense of Example~\ref{cofin.exa} --- that
is, if $j \setminus J'$ is connected for any $j \in J$ --- and if
$J$ is directed, it suffices to require that $j \setminus J'$ is not
empty.

\begin{defn}\label{sharp.def}
For any regular cardinals $\kappa < \mu$, $\mu$ is {\em sharply
  greater that} $\kappa$ if for any $\mu$-bounded set $S$, the
partially ordered set $P_\kappa(S)$ contains a cofinal $\mu$-bounded
subset $U$. We write $\kappa \trl \mu$ if either $\kappa=\mu$, or
$\mu > \kappa$ is sharply greater than $\kappa$,
\end{defn}

\begin{exa}\label{suc.exa}
For any regular cardinal $\kappa$ with successor cardinal
$\kappa^+$, we have $\kappa \trl \kappa^+$. Indeed, for any $S$ with
$|S|=\kappa$, we can equip $S$ with a minimal order $Q_\kappa$ of
Example~\ref{ord.exa}, and then the embedding $Q_\kappa \to
P_\kappa(S)$, $q \mapsto Q_\kappa/q$ induced by \eqref{yo.pos.eq} is
cofinal.
\end{exa}

\begin{exa}\label{P.ka.exa}
For any regular cardinal $\kappa$, we have $\kappa \trl
P(\kappa)^+$. Indeed, for any set $S$ with $|S|=\kappa$ and any set
$S'$, a map $S' \to P(S)$ is the same thing as a subset in $S'
\times S$, and $|S' \times S| = |S|$ are soon as $S'$ is non-empty
and $\kappa$-bounded, so that $|P_\kappa(P(S))| = |P(S)|$. More
generally, for any $\kappa' \leq \kappa$, we have
$|P_{\kappa'}(P(S))| \leq |P_\kappa(P(S))| < P(\kappa)^+$, so that
$\kappa' \trl P(\kappa)^+$.
\end{exa}

\begin{exa}\label{sharp.non.exa}
For any infinite cardinal $\kappa = \kappa_0$, define by induction
$\kappa_n = \kappa_{n-1}^+$, and let $\kappa_\infty =
\cup_n\kappa_n$. Then $\kappa_\infty$ is not regular. Its successor
$\kappa_\infty^+$ is regular; however, if $\kappa$ is regular but
uncountable, it is not true that $\kappa \trl \kappa_\infty^+$.
\end{exa}

\begin{lemma}\label{sharp.le}
Assume given regular cardinals $\kappa < \mu$. Then $\kappa \trl
\mu$ if and only if for any $\mu$-bounded set $S$, there exists a
left-closed embedding $S \subset J$ into a $\mu$-bounded
left-$\kappa$-bounded $\kappa$-directed partially ordered
set $J$. Moreover, one can choose $J$ to be Noetherian.
\end{lemma}

\proof{} The ``if'' part is clear: for any left-closed embedding
$\eps:S \to J$ with $\kappa$-directed left-$\kappa$-bounded $J$, the
map $\eps^* \circ \Y:J \to P_\kappa(S)$ has cofinal image.

For the ``only if'' part, assume that $\kappa \trl \mu$. Say that a
map $f:J' \to J$ of partially ordered sets is {\em weakly
  $\kappa$-directed} if if for any commutative diagram
$$
\begin{CD}
  S @>{e}>> S^>\\
  @V{a}VV @VV{g}V\\
  J @>{f}>> J'
\end{CD}
$$
with discrete $\kappa$-bounded $S$ and conservative $g$, there
exists a map $b:S^> \to J$ such that $b \circ e = a$ and $f \circ b
= g$. Note that if $J=\ppt$, then any map is trivially weakly
$\kappa$-directed (there are no conservative maps $S^> \to \ppt$
with non-empty $J$). Conversely, if $J$ is has no maximal elements,
then for any $S \subset J/j$, we can choose $j' > j$ and construct a
conservative cone $S^> \to J$ of the embedding $S \subset J$, with
vertex $j'$. Therefore if such a $J$ is $\kappa$-directed and $f:J'
\to J$ is weakly $\kappa$-directed, then $J'$ is $\kappa$-directed.

Now fix a set $S_\mu$ of cardinality $|S_\mu|=\mu$, let $Q_\kappa$
be a well-ordered set of Example~\ref{ord.exa} of cardinality
$|Q_\kappa|=\kappa$, and let $P$ be the set of pairs $\langle Q,J
\rangle$ of a left-closed subset $Q \subset Q_\kappa$, and a
$\mu$-bounded subset $J \subset S_\mu \times Q$ equipped with a
partial order such that $J$ is left-$\kappa$-bounded, and the
projection $J \to Q$ is a conservative order-preserving weakly
$\kappa$-directed map. In particular, for any $\mu$-bounded subset
$S \subset S_\mu$, $P$ contains the pair $\langle \{o\},S\rangle$,
where $o \in Q$ is the initial element, and $S$ carries the discrete
order. Moreover, order $P$ by setting $\langle Q',J' \rangle \leq
\langle Q'',J''\rangle$ iff $Q' \subset Q''$, and $J' = (S \times
Q') \cap J'' \subset S \times Q''$, with the order induced from
$J''$. Then since $\mu$ is regular, $P$ with this order has upper
bounds of all ascending chains given by the union of the chain, so
by the Zorn Lemma, it has a maximal element $\langle J,Q \rangle
\geq \langle \{o\},S \rangle$. If $Q$ is not the whole $Q_\kappa$,
then $Q = Q_\kappa/'q$ for some $q \in Q_\kappa$. Then consider $J$
with discrete order, take the cofinal $\mu$-bounded subset $U
\subset P_\kappa(J)$ provided by Definition~\ref{sharp.def}, forget
the order on $U$, choose an embedding $U \subset S_\mu$, and let
$\lambda:U \to L(J)$ be the map sending some $\kappa$-bounded $X
\subset J$ lying in $U \subset P_\kappa(J)$ to the union $\lambda(X)
= \cup_{x \in X}J/x \subset J$. Then $Q'=Q_\kappa/q \cong Q^>$, and
the union $J' = J \copr U \subset (S \times Q) \copr S \cong S
\times Q'$ with the glued order $J' = J \copr_\lambda U$ of
\eqref{pos.glue.eq} gives an element $\langle Q',J' \rangle$ in $P$
strictly bigger than $\langle J,Q\rangle$. This contradicts
maximality.
\endproof

\begin{corr}\label{sharp.corr}
Assume given regular cardinals $\kappa \trl \mu$. Then any
left-$\kappa$-bounded $\mu$-bound\-ed Noetherian partially ordered
set $J \in \Pos_\mu$ admits a left-closed embedding $J \subset J'$
into a left-$\kappa$-bounded $\mu$-bounded $\kappa$-directed
Noetherian $J' \in \Pos_\mu$.
\end{corr}

\proof{} Take a $\mu$-bounded cofinal subset $U \subset
P_\kappa(J)$, let $S$ be $U$ with discrete order, and take an
embedding $\eps:S \to J_\mu$ provided by Lemma~\ref{sharp.le}. Let
$\lambda = \Lambda \circ \ups \circ \eps^* \circ \Y:J_\mu \to
L_\kappa(J)$, where $\ups:P_\kappa(S) \to P_\kappa(J)$ sends a
subset $S' \subset S \subset P_\kappa(J)$ to the union of its
elements, $\kappa$-bounded if so was $S'$, and $\Lambda$ is the map
\eqref{La.eq} that sends $P_\kappa(J)$ into $L_\kappa(J)$ since $J$
is left-$\kappa$-bounded. Then $J' = J \copr_\lambda J_\mu$ does the
job.
\endproof

\subsection{Perfect colimits.}\label{pos.perf.subs}

In principle, the category $\Pos$ has all colimits (see
e.g. Remark~\ref{pos.col.rem} below) but they often
behave badly. Here is a minimal sanity check, and an additional
requirement that turns out to be very useful.

\begin{defn}\label{pos.perf.def}
A colimit $J' \cong \colim_J\phi(j)$ of a functor $\phi:J \to \Pos$
for some $J \in \Pos$ is {\em perfect} if it is preserved by the
forgetful functor $\Pos \to \Sets$ and the embedding $\Pos \subset
\Cat$. A perfect colimit $J'$ is {\em universal with respect to a
  map} $J' \to J''$ to some $J'' \in \Pos$ if for any map $J'_0 \to
J''$, we have a perfect colimit
\begin{equation}\label{pos.perf.eq}
  J'_0 \times_{J''} J' \cong \colim_j(\phi(j) \times_{J''} J'_0),
\end{equation}
and {\em universal} if it is universal with respect to the identity
map $\id:J' \to J'$.
\end{defn}

\begin{lemma}\label{pos.uni.le}
A perfect colimit $J' \cong \colim_J\phi(j)$ is universal iff it
satisfies the condition of Definition~\ref{pos.perf.def} for $J_0
\cong [1]$, and it is always universal with respect to any map $J'
\to J''$ with $\dim J'' \leq 1$. Moreover, for any perfect
colimit, \eqref{pos.perf.eq} holds if $J_0 \to J''$ is a fibration
or a cofibration.
\end{lemma}

\proof{} It is easy to see that all colimits of discrete sets are
universal, so we always have an isomorphism \eqref{pos.perf.eq} of
the underlying discrete sets, and to check that \eqref{pos.perf.eq}
is a colimit in $\Pos$, we just need to check that
\begin{itemize}
  \item for any cone of the functor $\phi(j) \times_{J''}
J'_0$, with some vertex $J_1$, the corresponding map $a:J'_0
\times_{J''} J'_0 \to J_1$ is order-preserving.
\end{itemize}
This can be checked after restricting to all subsets $[1] \subset
J'_0$, thus the first claim. Next, assume that $J'' \cong J'$ and
$J' \to J''$ is the identity map. Then if $J_0 \cong \ppt$, there is
nothing to check, and if $J' \cong \ppt$, then at least one of the
sets $J'(j)$ is non-empty, so that $J'(j) \times J_0$ contains a
copy of $J_0$, and \thetag{$\bullet$} is clear. By the combination
of these two observations, we see that \thetag{$\bullet$} holds when
$J_0 \to J' = J''$ factors through a one-element subset. In
particular, it always holds for the fibers of the map $J_0 \to
J'$. But if this map is a cofibration or a fibration, it is proper in
the sense of Definition~\ref{fl.def} by Lemma~\ref{fun.le} or its
dual, and the relative functor category $J_1' = \Fun(J_0|J',J_1)$ is
then a partially ordered set. The set-theoretic map $a$ of
\thetag{$\bullet$} gives rise to a set-theoretic map $a':J' \to
J_1'$, and to see that $a$ is order-preserving, it suffices to check
the same for $a'$; this holds by the universal property of the
colimit $J' = \colim_J\phi(j)$. Now to finish the proof, let $J''$
be arbitrary. Then if $J_0 \to J''$ is a fibration or a cofibration,
the same holds for $J' \times_{J''} J_0 \to J'$, and any map $[1]
\to J''$ of dimension $\leq 1$ is a composition of a fibration and a
cofibration.
\endproof

\begin{corr}\label{univ.corr}
For any perfect colimit $J' = \colim_J\phi(j)$ in $\Pos$, $\Ar(J') =
J' \setminus J' \cong \colim_J(J' \setminus \phi(j))$ is a
perfect colimit universal with respect to the projection
$\sigma:\Ar(J') \to J'$.
\end{corr}

\proof{} For any $J_0 \to J'$, $\tau:\Ar(J') \times_{J'} J_0 \cong
J_0 / J' \to J'$ is a cofibration.
\endproof

\begin{exa}\label{pos.perf.exa}
Since the embedding $\Pos \subset \Cat$ is fully faithful, a cone
$\phi_>:J^> \to \Pos$ of some functor $\phi:J \to \Pos$, $J \in
\Pos$ that becomes universal in $\Cat$ is also universal in
$\Pos$. Moreover, since the functor $\Cat \to \Sets$ sending a small
category to its set of objects has a right-adjoint $S \mapsto e(S)$,
such a cone also stays universal after applying the forgetful
functor $\Pos \to \Sets$, so it gives a perfect colimit in the sense
of Definition~\ref{pos.perf.def}. Conversely, if a functor $\phi:J
\to \Pos$ has a universal cone $\phi_>:J^> \to \Cat$, then its
composition $\Red \circ \phi_>:\phi_> \to \Pos$ with the reduction
functor \eqref{red.eq} is universal by adjunction; however, if
$\phi_>(o)$ is not in $\Pos$, the corresponding colimit is not
perfect. For example, consider the Kronecker category $\II$ of
Example~\ref{kron.exa}. Then we have two embeddings
$\eps(s),\eps(t):[1] \to \II$ of \eqref{eps.f}, and the square
\begin{equation}\label{II.red.sq}
\begin{CD}
\{0,1\} @>>> [1]\\
@VVV @VV{\eps(t)}V\\
[1] @>{\eps(s)}>> \II
\end{CD}
\end{equation}
is a cocartesian square in $\Cat$. However, since $\Red(\II) = [1]$,
the corresponding cocartesian square in $\Pos$ is not perfect (it is
not preserved by the embedding $\Pos \subset \Cat$).
\end{exa}

\begin{exa}
Let $\mu_l:[1] \to [3]$, $l=0,1$ be the map $\mu_l(n) = 2n+l$, and
let $e:[1],[3] \to [0] = \ppt$ be the tautological projections. Then
the commutative square
\begin{equation}\label{compl.sq}
\begin{CD}
[1] \copr [1] @>{\mu_0 \copr \mu_1}>> [3]\\
@V{e \copr e}VV @VV{e}V\\
[0] \copr [0] @>>> [0]
\end{CD}
\end{equation}
is cocartesian in $\Pos$ but not perfect. In fact, \eqref{compl.sq}
is a cocartesian square of categories in the sense of
Subsection~\ref{dia.subs} but not a cocartesian square in $\Cat$ ---
to turn it into one, one has to replace $[0]$ in the bottom right
corner with the equivalent category $e(\{0,1\})$.
\end{exa}

\begin{exa}\label{cone.sq.exa}
Colimits of arbitrary functors $\phi:J \to \Pos$ actually reduce to
cocartesian squares. Namely, let $J_\idot \to J$ be the cofibration
corresponding to $\phi$. Then a cone $\phi_>$ for $\phi$ with some
vertex $J'$ is the same thing as a map $J_\idot \to J \times J'$
cocartesian over $J$. Thus if we let $J_\idot^\sharp \subset J_\idot$
be the dense subset corresponding to maps cocartesian over $J$, then
by Example~\ref{IX.exa}, the cone is universal if and only if the
square
\begin{equation}\label{cone.sq}
\begin{CD}
  J_\idot^\sharp @>>> \pi_0(J_\idot^\sharp)\\
  @VVV @VVV\\
  J_\idot @>>> J'
\end{CD}
\end{equation}
is cocartesian. Moreover, the colimit $J'=\colim_J\phi$ is perfect
iff so is the cocartesian square \eqref{cone.sq}.
\end{exa}

A standard colimit \eqref{lc.coli.eq} is always perfect, and it is
also universal (for any map $J'_0 \to J'$, the corresponding colimit
is the standard colimit for the projection $J'_0 \to J' \to
J$). Alternative forms \eqref{pos.2.st.sq} of standard pushout
squares \eqref{pos.st.sq} are also perfect and universal. Let us
list some other examples of perfect colimits, some universal, some
not.

\begin{exa}\label{pos.triv.exa}
We start with two trivial examples. Firstly, $\Pos$ obviously has
all filtered colimits, so any colimit over a directed $J$ exists and
is perfect. Moreover, since $[1]$ is a compact object in $\Pos$,
such a colimit is universal by Lemma~\ref{pos.uni.le}. We also
observe that filtered colimits in $\Pos$ send pointwise left or
right-closed embeddings to left resp.\ right-closed
embeddings. Secondly, say that a map $g:J \to J'$ is {\em split} if
it is a full embedding, and $J' = J \copr (J' \setminus J)$. Then
for any split map $g:J \to J'$ and any map $f:J \to J_0$, the
coproduct $J_0' = J_0 \copr_J J' \cong J_0 \copr (J' \setminus J)$
exists, and this provides a perfect universal cocartesian square in
$\Pos$.
\end{exa}

\begin{exa}\label{J.S.sq.exa}
For any set $S$, the partially ordered set \eqref{B.S} of
Example~\ref{B.S.exa} fits into a cartesian square
\begin{equation}\label{S.S.sq}
\begin{CD}
S^> @>>> \ppt\\
@VVV @VV{\eps(o)}V\\
S^\hush @>{\xi_\perp}>> S^<,
\end{CD}
\end{equation}
where $\eps(o):\ppt \to S^<$ is the embedding onto the smallest
element $o \in S^<$. The square \eqref{S.S.sq} is cocartesian,
perfect, and since $\dim S^< \leq 1$, it is universal by
Lemma~\ref{pos.uni.le}. Moreover, by virtue of \eqref{B.S}, the
barycentric subdivision functor $B$ sends \eqref{S.S.sq} to a
cocartesian square.
\end{exa}

\begin{exa}\label{J.n.exa}
If $I$ in Lemma~\ref{cocart.le} is a partially ordered set, then the
square \eqref{cocart.sq} is a cartesian cocartesian square in
$\Pos$, and by Example~\ref{B.X.exa}, the barycentric subdivision
functor $B$ sends it to a cartesian cocartesian square (of the same
form). Since $BJ \cong BJ^o$, the same holds for the square formed
by opposite partially ordered sets. In particular, for any $n \geq 0$
and discrete cofibration $J \to [n]$, we have a cartesian
cocartesian square
\begin{equation}\label{J.n.sq}
\begin{CD}
  J_0 \times [n-1] @>{\id \times t}>> J_0 \times [n]\\
  @V{l \circ (f_0 \times \id)}VV @VV{l}V\\
  t^*J @>{t}>> J,
\end{CD}
\end{equation}
where $t:[n-1] \to [n]$ is the standard right-closed embedding of
Example~\ref{st.cl.exa}, $l$ are the functors of
Example~\ref{cyl.fib.exa}, and $f_0:J_0 \to J_1$ is the transition
map of the discrete cofibration $J \to [n]$. The square
\eqref{J.n.sq} is perfect, and barycentric subdivision functor $B$
sends it to a cartesian cocartesian square.
\end{exa}
  
\begin{exa}\label{pos.sq.exa}
For any standard pushout square \eqref{pos.st.sq}, the square
\begin{equation}\label{pos.st.gt.sq}
\begin{CD}
J^>_{01} @>{m_0}>> J^>_0\\
@V{m_1}VV @VVV\\
J^>_1 @>>> J^>
\end{CD}
\end{equation}
is also a cartesian cocartesian square in $\Pos$, it is perfect, and
the barycentric subdivision functor $B$ sends it to a standard
pushout square. We note that the embeddings $m_0$, $m_1$ in
\eqref{pos.st.gt.sq} are left-reflexive by
Example~\ref{cl.refl.exa}. The characteristic map $\chi:J \to \V$
extends to a map $\chi^>:J^> \to \V^> \cong [1]^2$. In the universal
case $J=\V$, \eqref{pos.st.gt.sq} becomes the square
\begin{equation}\label{pos.sq}
\begin{CD}
[1] @>{m}>> [2]\\
@V{m}VV @VV{a}V\\
[2] @>{b}>> [1]^2,
\end{CD}
\end{equation}
where the map $m:[1] \to [2]$ is as in Example~\ref{cyl.exa}, and
the map $a$ resp.\ $b$ sends $0$ to $0 \times 0$, $2$ to $1 \times
1$, and $1$ to $0 \times 1$ resp.\ $1 \times 0$. A general square
\eqref{pos.st.gt.sq} is induced from \eqref{pos.sq} by pullback with
respect to $\chi^>:J^> \to [1]^2$.
\end{exa}

\begin{exa}\label{seg.exa}
For any integers $n \geq l \geq 0$, let $s:[0] \to [n-l]$
resp.\ $t:[0] \to [l]$ be the embeddings onto the initial
resp.\ terminal element, as in Example~\ref{n.refl.exa}. Then the
square
\begin{equation}\label{seg.del.sq}
\begin{CD}
[0] @>{s}>> [n-l]\\
@V{t}VV @VV{t}V\\
[l] @>{s}>> [n]
\end{CD}
\end{equation}
in $\Pos$ is cocartesian and perfect. It is {\em not} universal as
soon as $n > l > 0$ (for a counterexample, take the map $[1] \to
[n]$ sending $0$ to $0$ and $1$ to $n$).
\end{exa}

\begin{exa}\label{cyl.pos.exa}
For any map $f:J' \to J$ in $\Pos$, the cylinder $\Cyl(f)$ and the
dual cylinder $\Cyl^o(f)$ fit into cocartesian squares
\eqref{cyl.o.sq}, and these are perfect cocartesian squares in
$\Pos$ by Example~\ref{pos.perf.exa}.
\end{exa}

\begin{exa}\label{Z.N.exa}
Let $\overline{\N}$ be $\N$ with discrete order. Then the map
$\zeta$ of \eqref{zeta.eq} fits into a cartesian square
\begin{equation}\label{Z.N.sq}
\begin{CD}
  \ZZ_\infty @>{\zeta}>> \N\\
  @A{q'}AA @AA{q}A\\
  \overline{\N} \times [1] @>>> \overline{\N},
\end{CD}
\end{equation}
where $q:\overline{\N} \to \N$ is the embedding $l \mapsto l+1$, and
the square \eqref{Z.N.sq} is also cocartesian and perfect. For any
$m \geq 1$, we can also consider the map \eqref{z.m.eq} and its
restriction to $\Z_{2m-1}$, and they also fit into cartesian
cocartesian perfect squares induced by \eqref{Z.N.sq}. The simplest
non-trivial one is
\begin{equation}\label{Z3.sq}
\begin{CD}
  \ZZ_3 @>>> [2]\\
  @AAA @AA{a}A\\
  [1] @>>> \ppt,
\end{CD}
\end{equation}
where $a:\ppt \to [2]$ is the embedding onto $1 \in [2]$.
\end{exa}

\begin{exa}\label{Z.NN.exa}
The map $q'$ in \eqref{Z.N.sq} actually fits into a cartesian square
\begin{equation}\label{Z.NN.sq}
\begin{CD}
  \overline{\N} \times [1] @>>> \ZZ_\infty\\
  @A{p'}AA @AA{q'}A\\
  \overline{\N} @>{p}>> \overline{\N} \times [1],
\end{CD}
\end{equation}
where $p$ and $p'$ are given by $p(2n+1-l)=n \times l$, $p'(2n+l-1)
= n \times l$, $n \in \N$, $l=0,1$ (note that in particular, $p$ is
dense). The square \eqref{Z.NN.sq} is cocartesian and perfect, and
$B$ sends it to a standard pushout square.
\end{exa}

\begin{exa}\label{BJ.coli.exa}
For any $J \in \Pos$, let $\overline{J} \subset J$ be $J$ with
discrete order, and let $\overline{B}_\idot(J) = B_\idot(J)^\sharp
\subset B(J)$ be the dense subset corresponding to maps cocartesian
with respect to the forgetful cofibration $\nu_\idot:B_\idot(J) \to
B(J)$. Then \eqref{B.xi.dot} induces a cartesian square
\begin{equation}\label{BJ.J.sq}
\begin{CD}
\overline{B}_\idot(J) @>{\xi_\idot}>> \overline{J}\\
@VVV @VVV\\
B_\idot(J) @>{\xi_\idot}>> J,
\end{CD}
\end{equation}
and this square is cocartesian and perfect. Therefore by
Example~\ref{cone.sq.exa}, we have
\begin{equation}\label{BJ.coli.eq}
J \cong \colim_{BJ}\phi,
\end{equation}
where $\phi:B(J) \to \Pos$, $S \mapsto S$ is the tautological
functor corresponding to the cofibration $\nu_\idot$, and the
colimit \eqref{BJ.coli.eq} is also perfect. We note that it is also
preserved by the barycentric subdivision functor $B$ --- indeed, the
cofibration corresponding to $B \circ \phi:BJ \to \Pos$ is
$\tau:\Ar(BJ) \to BJ$, so that $B$ sends \eqref{BJ.coli.eq} to the
standard colimit \eqref{lc.coli.eq} for $\id:BJ \to BJ$.
\end{exa}

\subsection{Cylindrical colimits.}\label{pos.cyl.subs}

Example~\ref{pos.sq.exa}, Example~\ref{seg.exa}, and
Example~\ref{cyl.pos.exa} for an injective $f$ follow the same
general pattern. We have a partially ordered set $J$ and two full
subsets $J_0,J_1 \subset J$ such that $J=J_0 \cup J_1$, with the
intersection $J_{01} = J_0 \cap J_1$. Example~\ref{Z.NN.exa} is
similar, but the subsets are not full. If the subsets are full, the
setup can conveniently encoded by considering the characteristic map
$\phi:J \to \V$ such that $J_l=\phi^{-1}(\V/l)$, $l=0,1$, and we
then have a functor $\V \to \Pos$ equipped with an cone $\V^> \cong
[1]^2 \to \Pos$ that sends the terminal element to $J$. If $\phi$ is
order-preserving, then $J_0,J_1 \subset J$ are left-closed, the cone
is universal, and we obtain a standard pushout square in
$\Pos$. However, even if $\phi$ is not order-preserving, it might
still happen that the cone is universal. This is exactly what
happens in Examples~\ref{pos.sq.exa}-\ref{cyl.pos.exa}.

More generally, one can replace $\V$ with some partially ordered set
$J'$, consider a map $\phi:J \to J'$, and define a functor
$F(\phi):J' \to \Pos$ sending $j \in J'$ to $\phi^{-1}(J'/j) \subset
J$ with the induced order. We again have a cone $F(\phi)_>$ with
vertex $J$, and we can ask whether it is universal. By
\eqref{lc.coli.eq}, this is automatic if $\phi$ is
order-preserving. However, even if $\phi$ is not order-preserving,
we still might have an isomorphism \eqref{lc.coli.eq}. Here is a
trivial but useful example.

\begin{exa}\label{N.coli.exa}
Take $J'=\N$ --- that is, assume given a partially ordered set $J$
equipped with an exhaustive filtration by full subsets $J(n) \subset
J$, $n \geq 0$. Then we have
\begin{equation}\label{N.coli.eq}
  J \cong \colim_nJ(n).
\end{equation}
Indeed, to check the universal property, assume given $J' \in \Pos$
and a compatible collection of order-preserving maps $J(n) \to J$,
$n \geq 0$. Then as in Lemma~\ref{pos.uni.le}, since the filtration
is exhaustive, these define a unique map $J \to J'$, and we just
need to check that it preserves order. But as in
Example~\ref{pos.triv.exa}, $[1] \in \Pos$ is compact, so a map $[1]
\to J$ factors through $J(n)$ for a sufficiently large integer $n$.
\end{exa}

For a less trivial example, consider Example~\ref{cyl.pos.exa} under
the additional assumption that $f:J' \to J$ is a left-closed
embedding corresponding to some map $\chi:J \to [1]$, so that
$\Cyl(f) \cong J/[1]$, and note that one can identify $\V$ with the
twisted arrow category $\Tw([1])$. Now take an arbitrary map $\chi:J
\to J'$ in $\Pos$, with the comma-set $J / J' = J/_{\chi}J'$,
identify $\Ar(J') \cong \Tw(J')$ as sets (without preserving the
order), and define a map $\phi:J/J' \to \Tw(J')$ as the composition
$$
\begin{CD}
J/J' @>{\chi/\id}>> J'/J' \cong \Ar(J') @>{\sim}>> \Tw(J').
\end{CD}
$$
Consider the corresponding functor $F(\phi):\Tw(J') \to \Pos$ with
its cone $F(\phi)_>$. Explicitly, for any $f \in \Tw(J')$
represented by a pair $j \leq j' \in J'$, $F(\phi)$ sends $f$ to the
subset
\begin{equation}\label{J.cyl.eq}
(J/J')_f = \{j_0 \leq j_0' \in J / J'\mid \chi(j_0) \leq j, j_0' \geq
j'\} \subset J/J',
\end{equation}
with the induced order, and one can identify $(J/J')_f \cong (J/j)
\times (j' \setminus J')$.

\begin{lemma}\label{pos.phi.le}
The cone $F(\phi)_>:\Tw(J')^> \to \Pos$ is universal, so
that $J/J' \cong \colim_{\Tw(J')}(J/J')_f$ in $\Pos$.
\end{lemma}

\proof{} If we equip $J/J'$ with the discrete order, then $\phi$
becomes order-preserving, so that the cone $F(\phi)_>$ is universal
by \eqref{lc.coli.eq}. Then as in Lemma~\ref{pos.uni.le}, for any
cone $F(\phi)'_>:\Tw(J')^> \to \Pos$ of $F(\phi)$, we have a unique
map $a:F(\phi)_> \to F(\phi)'_>$, and it suffices to check that $a$
preserves the order. In other words, for any $\langle
j_0,j_0'\rangle,\langle j_1,j_1' \rangle \in J/J'$ such that $j_0
\leq j_1$ and $j_0' \leq j_1'$, we have to check that $a(\langle
j_0,j_0' \rangle) \leq a(\langle j_1,j_1' \rangle)$. But we have $\langle
j_0,j_0' \rangle \leq \langle j_0,j_1' \rangle \leq \langle j_1,j_1'
\rangle$, so we may further assume that either $j_0=j_1$ or $j_0' =
j_1'$. In either case, both elements lie in a single subset
\eqref{J.cyl.eq}.
\endproof

\begin{remark}\label{K.rem}
While the colimits in Example~\ref{J.S.sq.exa},
Example~\ref{J.n.exa}, Example~\ref{pos.sq.exa},
Example~\ref{Z.NN.exa} and Example~\ref{N.coli.exa} are preserved by
the ba\-ry\-centric subdivision functor $B$, this is not the case in
the other examples. In Lemma~\ref{pos.phi.le}, since $B$ sends full
embeddings to left-closed full embeddings, the colimit
$\colim_{\Tw(J')}B(J/J')_f$ exists, and the comparison map
$\colim_{\Tw(J')}B(J/J')_f \to B(J/J')$ is injective, but it need
not be surjective: its image consists of subsets $[n] \subset J/J'$
that lie entirely in a single subset \eqref{J.cyl.eq}. Explicitly,
for any injective map $i:[n] \to J/J'$, define a left-closed subset
$P(i) \subset [n]$ by
\begin{equation}\label{K.ka}
P(i) = \{l \in [n]\mid \tau(i(l)) < \chi(\sigma(i(n)))\},
\end{equation}
where $\sigma:J/J' \to J$, $\tau:J/J' \to J'$ are the projections;
then the comple\-ment $B(J/J') \ssetminus
\colim_{\Tw(J')}B(J/J')_f$ consists of $i$ with non-empty $P(i)$.
\end{remark}

\begin{remark}\label{S.rem}
In the situation of Example~\ref{cyl.pos.exa}, we can consider the
map
\begin{equation}\label{S.B.eq}
  B(J' \times [1]) \copr_{B(J')} B(J) \to B(\Cyl(f))
\end{equation}
induced by the square \eqref{cyl.o.sq}. If $f$ is injective, then so
is \eqref{S.B.eq} --- if $f$ is left-closed, this is in fact covered
by Lemma~\ref{pos.phi.le} --- but in the general case,
\eqref{S.B.eq} is usually not injective. One exception to this is
the case when $J' = S$ is discrete. In this case, $f:S \to J$ splits
into a coproduct of maps $f_s:\ppt \to J$, and then $\Cyl(f)$ can be
expreseed as a co-standard colimit \eqref{lc.o.coli.eq}: we have
\begin{equation}\label{cyl.f.s}
\Cyl(f) \cong \colim_{S^<}\Cyl(f_\idot),
\end{equation}
where $f_o$ is the embedding $\emptyset \to J$. Since $B$ sends
co-standard colimits to standard ones, \eqref{S.B.eq} is then a
standard colimit over $S^<$ of the corresponding maps for
$f_s$. Since $f_s:\ppt \to J$ is injective, \eqref{S.B.eq} for $f_s$
is also injective, and then \eqref{S.B.eq} for $f$ is injective as
well. Moreover, the complement $\overline{B}(\Cyl(f))$ to its image
admits exactly the same description as in Remark~\ref{K.rem}: we
have
\begin{equation}\label{b.f.s}
\overline{B}(\Cyl(f)) \cong \coprod_{s \in S}\overline{B}(\Cyl(f_s))
\cong \coprod_{s \in S}\overline{B}(\Cyl(\overline{f}_s)),
\end{equation}
where $\overline{f}_s$, $s \in S$ is the left-closed embedding $\ppt
\to f_s(\ppt) \setminus J$.
\end{remark}

\subsection{Anodyne maps.}\label{ano.subs}

By Lemma~\ref{po.copr.le}, the class of reflexive full embeddings in
$\Posf$ is closed under pushouts with respect to standard pushout
squares \eqref{pos.st.sq} but it is not saturated, and its
saturation is no longer stable under standard pushouts. It will be
useful to distinguish a larger class of maps that is closed under
both operations. Here is the definition.

\begin{defn}\label{pos.ano.def}
The class of {\em anodyne maps} is the minimal saturated class of
maps in $\Posf$ such that
\begin{enumerate}
\item for any $J \in \Posf$, the projection $J \times [1] \to J$ is
  anodyne, and
\item for $J,J_0,J_1 \in \Posf$ equipped with maps $\chi_l:J_l \to
  J$, $l = 0,1$, and any map $f:J_0 \to J_1$ such that $f$ is co-lax
  over $J$ and the induced map $f/j:J_0 / j \to J_1 / j$ of
  \eqref{lax.pos.eq} is anodyne for any $j \in J$, the map $f$ is
  itself anodyne.
\end{enumerate}
A partially ordered set $J \in \Posf$ is {\em anodyne} if so is the
map $J \to \ppt$.
\end{defn}

The class of anodyne maps is manifestly closed with respect to
standard pushouts --- that is, for any map $J \to \V$ with anodyne
embedding $J_o \subset J/1$, the embedding $J/0 \subset J$ is
anodyne (apply Definition~\ref{pos.ano.def}~\thetag{ii} to the maps
$J/0 \to J \to \V$). On the other hand, by Lemma~\ref{amp.le},
Definition~\ref{pos.ano.def}~\thetag{i} insures that all reflexive
full embedding are anodyne. This has the following corollaries.

\begin{corr}\label{ano.corr}
The class of anodyne maps is the minimal saturated class of maps in
$\Posf$ that satisfies Definition~\ref{pos.ano.def}~\thetag{i} and
\begin{enumerate}
\renewcommand{\labelenumi}{(\roman{enumi}')}
\setcounter{enumi}{1}
\item for any cofibrations $J_0,J_1 \to J$ in $\Posf$, a map $f:J_0
  \to J_1$ cocartesian over $J$ is anodyne if so are all the fibers
  $f_i:J_{0,j} \to J_{1,j}$, $j \in J$.
\end{enumerate}
\end{corr}

\proof{} For any map $g:J' \to J$ in $\Posf$, we have the
factorization \eqref{lc.pos}, and since $\Posf$ is ample, $J' /_g J
\in \Posf$ by Lemma~\ref{amp.le}. Moreover, \eqref{lc.pos} is
functorial with respect to co-lax maps over $J$, and
Definition~\ref{pos.ano.def}~\thetag{ii} for a map $f:J_0 \to J_1$
is identically the same as \thetag{ii'} for $f / \id:J_0 / J \to J_1
/ J$. But $\eta$ in \eqref{lc.pos} is full and right-reflexive, so
that by Lemma~\ref{amp.le}, as soon as
Definition~\ref{pos.ano.def}~\thetag{i} is satisfied, \thetag{ii'}
implies \thetag{ii}. Conversely, if $J_0,J_1 \to J$ are
cofibrations, and $f:J_0 \to J_1$ is cocartesian over $J$, then
$J_{0,j} \subset J_0/j$ is full and right-admissible for any $j \in
J$, and similarly for $J_1$, so that \thetag{ii} implies
\thetag{ii'}.
\endproof

\begin{corr}\label{V.ano.corr}
Assume given two cofibrations $J,J' \to \V$ in $\Posf$, and a map
$f:J \to J'$ cocartesian over $\V$. Assume that $f$ is an
isomorphism over $o,1 \in \V$, and the embedding $J_0 \subset J$ is
anodyne. Then the same holds for the embedding $J'_0 \subset J'$.
\end{corr}

\proof{} Since $J \to \V$ is a cofibration, the embedding $J_0
\subset J/0$ is left-reflexive, so that $J_0 \subset J$ is anodyne
iff so is $J/0 \subset J$, and similarly for $J'$.

Let $\wV = \Cyl^o(\eps(0))$ be the dual cylinder of the embedding
$\eps(0):\ppt \to \V$ onto $0 \in \V$ (explicitly, $\wV$ is obtained
by adding one element $0'$ to $\V = \{o,0,1\}$, with order relation
$0 \leq 0'$). We have the embeddings $s:\V \to \wV$ and $e:\wV \to
\V \times [1] \cong \Cyl^o(\id)$ such that $e \circ s:\V \to \V
\times [1]$ is the embedding onto $\V \times \{0\}$. Composing $e$
with the projection $\V \times [1] \to \V$ gives a projection $p:\wV
\to \V$ right-adjoint to $s$, $p$ has a right-adjojnt embedding
$t:\V \to \wV$ onto $\{0',1,o\} \subset \wV$, and $t$ has a
right-adjoint $q:\wV \to \V$ sending $0'$ to $0$, $1$ to $1$, and
$o,0$ to $o$. Now consider the cylinder $\Cyl(f)$ of the map $f$,
with its cofibration $\Cyl(f) \to \V \times [1]$, and let $\wt{J} =
e^*\Cyl(f)$. Then we have a cofibration $\wt{J} \to \wV$ such that
$s^*\wt{J} \cong J$, $t^*\wt{J} \cong J'$, and the embedding $t^J:J'
\to \wt{J}$ is left-reflexive. Moreover, compose the cofibration
$\wt{J} \to \wV$ with the map $q:\wV \to \V$ to obtain a map $\wt{J}
\to \V$, thus representing $\wt{J}$ as a standard pushout. Then the
embedding $J'/0 \to \wt{J}/0$ is also left-reflexive, so it suffices
to show that $\wt{J}/0 \subset \wt{J}$ is anodyne. But $\wt{J}_o
\cong J/0$, $\wt{J}/1 \cong J$, and anodyne map are stable with
respect to standard pushout squares.
\endproof

\begin{corr}\label{B.ano.corr}
For any $J \in \Posf$, the map $\xi:BJ \to J$ of \eqref{B.xi} is
anodyne, and a map $f:J_0 \to J_1$ in $\Posf$ is anodyne iff so is
$f^o:J_0^o \to J_1^o$ iff so is $B(f):BJ_0 \to BJ_1$.
\end{corr}

\proof{} By Definition~\ref{pos.ano.def}~\thetag{ii}, to prove the
first claim, it suffices to check that the induced map $BJ/j \to
J/j$ is anodyne for any $j \in J$. But we obviously have $BJ/j \cong
B(J/j)$, the embedding $B(\eps(o)):\ppt \to B(J/j)$ onto the maximal
element $o \in J/j$ is reflexive by Lemma~\ref{B.le}, hence anodyne
by Lemma~\ref{amp.le}, and so is its composition $\eps(o) = \xi
\circ B(\eps(o)):\ppt \to J/j$ with the map $\xi$. The second claim
now immediately follows from the functoriality $B$ and the maps
\eqref{B.xi}, and the identifications $B(J) \cong B(J^o)$.
\endproof

\begin{lemma}\label{st.ano.le}
A saturated closed class $W$ of maps in $\Posf$ that satisfies
Definition~\ref{pos.ano.def}~\thetag{i} and is stable under standard
pushouts satifies Definition~\ref{pos.ano.def}~\thetag{ii} for
$J=\V$.
\end{lemma}

\proof{} By Lemma~\ref{amp.le}, $W$ contains all reflexive full
embeddings.  Any map $f:J_0 \to J_1$ has the factorization
\eqref{cyl.deco}, where $s:J_0 \to \Cyl(f)$ is left-closed, and
$t_\dg:\Cyl(f) \to J_1$ is inverse to a reflexive left-closed
embedding, thus in $W$. Therefore $f \in W$ if and only if $s \in
W$. In the setting of Definition~\ref{pos.ano.def}~\thetag{ii}, for
any $J$, we also have $\Cyl(f)/j \cong \Cyl(f/j)$, $j \in J$, so in
checking the claim, we may assume that $f$ is a left-closed full
embedding. Then if $J=\V$, the full embedding $J_0 \to J_1$ factors
as
$$
\begin{CD}
J_0 @>>> J_0 \copr_{J_0/0} J_1/0 @>>> J_1,
\end{CD}
$$
the first map is a standard pushout of a map in $W$, and the
second map fits into a standard pushout square
$$
\begin{CD}
J_0/1 \copr_{J_0/o} J_1/o @>>> J_1/1\\
@VVV @VVV\\
J_0 \copr_{J_0/0} J_1/0 @>>> J_1,
\end{CD}
$$
so it suffices to check that the top arrow is in $W$. But the
embedding $J_0/1 \to J_0/1 \copr_{J_0/o} J_1/o$ is the standard
pushout of a map in $W$, so we are done by the two-out-of-three
property.
\endproof

\begin{exa}\label{J.S.exa}
For any set $S$ and any $J \in \Posf$ equipped with a map $J \to
S^<$, the map
\begin{equation}\label{J.S}
J^\hush = S^\hush \times_{S^<} J \to J
\end{equation}
induced by \eqref{xi.B} is anodyne. Indeed, by \eqref{B.S},
we have a decomposition
\begin{equation}\label{B.S.J}
  J^\hush = S^\hush \times_{S^<} J \cong (J_o \times S^>) \copr_{J_o \times S}
  \coprod_{s \in S}J/s
\end{equation}
that induced a decomposition
\begin{equation}\label{J.pl}
J^\hush/s \cong (J_o \times S^>) \copr_{J_o} J/s
\end{equation}
for any $s \in S \subset S^<$, the projection $J^\hush/s \to J/s$ has a
section induced by the embedding $J_o \to J_o \times S^>$ onto $J_o
\times \{s\}$, and this section is a reflexive full embedding, hence
anodyne by Lemma~\ref{amp.le}, so that its one-sided inverse
\eqref{J.pl} is also anodyne. Moreover, we have $J^\hush/o \cong J^\hush_o
\cong S^> \times J_o \cong S^> \times J/o$, so that the projection
$J^\hush/o \to J/o$ also has a reflexive one-sided inverse, and then
\eqref{J.S} is anodyne by Definition~\ref{pos.ano.def}~\thetag{ii}.
\end{exa}

\begin{prop}\label{ano.prop}
The class of anodyne maps is the minimal saturated class of maps in
$\Posf$ that satisfies Definition~\ref{pos.ano.def}~\thetag{i} and
\thetag{ii} for $J = S^<$, $S \in \Sets$. Moreover, it is also the
minimal saturated class of maps in $\Posf$ that contains the maps
\eqref{J.S} for any $S \in \Sets$ and $J \to S^<$, and satisfies
Definition~\ref{pos.ano.def}~\thetag{ii} for $J = \V$ and for
discrete $J$. Furthermore, instead of
Definition~\ref{pos.ano.def}~\thetag{ii} for $J = \V$, it suffices
to require that the class of anodyne maps is closed with respect to
standard pushouts \eqref{pos.st.sq}.
\end{prop}

\proof{} For the first claim, assume
Definition~\ref{pos.ano.def}~\thetag{ii} for $S^<$, and let us prove
it for an arbirary $J$. Note that the claim for $J$ trivially
implies the claim for any full subset $J' \subset J$, so in
particular, we know it for a discrete $J$. By induction on
$n=\dim(J)$, we may assume it proven for the $(n-1)$-st skeleton
$\sk_{n-1}J$ of \eqref{sk.J}. Take $S = J \ssetminus \sk_{n-1}J$, and
as in \eqref{p.J.sq}, extend the right-closed embedding $S \to J$ to
a map $p:J \to S^<$ sending any $j \in J \ssetminus \sk_{n-1}J$ to
itself and the rest to $o \in S^<$. Then it suffices to check that
the map $f:J_0 \to J_1$ satifies the assumption of
Definition~\ref{pos.ano.def}~\thetag{ii} with respect to the map $p
\circ g:J_1 \to S^<$ --- that is, we need to show that the induced
map $J_0/s \to J_1/s$ is anodyne for any $s \in S^<$. If $s = o$,
this is the inductive assumption, and for any $s \in S \subset S^<$,
$J_0/o \subset J_0/s$ is an elementary extension that defines the
standard coproduct
\begin{equation}\label{sk.s}
J_0/s \cong \sk_{n-1}J \copr_{\sk_{n-1}J/j} J/j
\end{equation}
of Example~\ref{elem.exa}, where $j \in J \ssetminus \sk_{n-1}J = S$
corresponds to $s$. The same holds for $J_1$, so it suffices to
consider the characteristic maps $J_0/s,J_1/s \to \V$ and invoke
Definition~\ref{pos.ano.def}~\thetag{ii} for $J=\V$.

For the second claim, note that if $S=\ppt$, so that $S^< \cong
[1]$, and a map $J \to [1]$ factors through $0 \subset [1]$, then
$J^\hush = J \times [1]$, and the map \eqref{J.S} is the projection $J
\times [1] \to J$. Thus the minimal class in question satisfies
Definition~\ref{pos.ano.def}~\thetag{i}. Moreover, for any discrete
$S$ and $J \to S^<$, the decompositions \eqref{J.pl} together yield
a decomposition
\begin{equation}\label{J.pl.S}
J^\hush = (J_o \times S^>) \copr_{J_o \times S} \coprod_{s \in S}J/s,
\end{equation}
and the coproduct is again a standard pushout of \eqref{pos.st.sq}.
Thus to check that for any map $f:J_0 \to J_1$ over $S^<$ with
anodyne $J_0/s \to J_1/s$, $s \in S^<$, the map $f^\hush:J_0^\hush \to
J_1^\hush$ is anodyne, it suffices to use
Definition~\ref{pos.ano.def}~\thetag{ii} for $J=S$ and $J=\V$. Now
the first claim of the Proposition shows that all anodyne maps are
in our minimal class, and the converse immediately follows from
Example~\ref{J.S.exa}. Finally, the third claim is
Lemma~\ref{st.ano.le}.
\endproof

\subsection{Left-bounded sets.}

Let us now consider the category $\Posf^+$ of partially ordered sets
that are only left-bounded, not necessarily finite-dimensional.

\begin{defn}\label{pl.ano.def}
The class of {\em $+$-anodyne maps} is the minimal saturated class
of maps in $\Posf^+$ that satisfies the conditions \thetag{i},
\thetag{ii} of Definition~\ref{pos.ano.def}, with $\Posf$ replaced
by $\Posf^+$.
\end{defn}

Since $\Posf^+$ is ample by Lemma~\ref{amp.le}, all reflexive maps
are $+$-anodyne, and all anodyne maps in $\Posf \subset \Posf^+$ are
also trivially $+$-anodyne. Corollary~\ref{ano.corr} and
Lemma~\ref{st.ano.le} also hold for $\Posf^+$ with the same proof.

\begin{exa}\label{N.ano.exa}
Assume given a partially ordered set $J \in \Posf^+$ equip\-ped with
an exhaustive filtration $J(n) \subset J$ as in
Example~\ref{N.coli.exa}, and let $J_\idot \subset J \times \N$ be
the incidence subset -- namely, the subset of pairs $j \times n$
such that $j \in J(n)$. Then the projection $\tau:J_\idot \to \N$
is a cofibration with fibers $J(n)$ and transition functors given by
inclusions. On the other hand, for any $j \in J$, the fiber
$J_{\idot,j}$ of the projection $\sigma:J_\idot \to J$ is a right-closed
subset in $\N$, so it has the smallest element $n(j)$ and is
isomorphic to $n(j) \setminus \N$. Moreover, $\sigma$ is a cofibration
whose transition functor $J_{\idot,j} \to J_{\idot,j'}$ for any $j
\leq j'$ sends $n \in n(j) \setminus \N$ to $\max\{n,n(j')\}$. Since
as abstract partially ordered sets, $n \setminus \N \cong \N$ for
any $n$, and the embedding $\eps(0):\ppt \to \N$ is right-reflexive,
$\sigma:J_\idot \to J$ is $+$-anodyne, and $\kappa^+$-anodyne if $|J|
< \kappa$.
\end{exa}

\begin{lemma}\label{N.coli.le}
Assume given $J_0,J_1 \in \Posf^+$ equipped with exhaustive
filtrations as in Example~\ref{N.coli.exa}, and a map $f:J_0 \to
J_1$ that preserves the filtrations and induces a $+$-anodyne map
$f(n):J_0(n) \to J_1(n)$ for any $n \geq 0$. Then $f$ itself is
$+$-anodyne.
\end{lemma}

\proof{} By Example~\ref{N.ano.exa}, it suffices to show that the
map $f_\idot:J_{\idot,0} \to J_{\idot,1}$ between the corresponding
incidence sets is $+$-anodyne. But this map is cocartesian over $\N$
with fibers $f(n)$, so it is $+$-anodyne by
Definition~\ref{pos.ano.def}~\thetag{ii} in the equivalent form
\thetag{ii'} of Corollary~\ref{ano.corr}.
\endproof

\begin{remark}
Note that both $J_0$ and $J_1$ in Lemma~\ref{N.coli.le} may lie in
$\Posf$ -- e.g. one can take $\ZZ_\infty$ of Example~\ref{Z.exa}
with filtration $\ZZ_m \subset \ZZ_\infty$ -- but we do not claim
that the map $f$ is then anodyne.
\end{remark}

\begin{exa}\label{J.pl.exa}
The map \eqref{zeta.eq} of Example~\ref{Z.exa} is $+$-anodyne (apply
Definition~\ref{pos.ano.def}~\thetag{ii} over $J=\N$). Moreover, for
any cofibration $J \to \N$, the map $\zeta^*J \to J$ is $+$-anodyne
(again apply Definition~\ref{pos.ano.def}~\thetag{ii} over $\N$, and
use Lemma~\ref{Z.m.le}). Note that as in \eqref{z.m.J.eq},
\eqref{zz.m.sq} induces a standard coproduct decomposition
\begin{equation}\label{J.pl.dec}
\zeta^*J \cong (\zeta_2^*J) \copr_{J_1} \zeta^*q^*J,
\end{equation}
where $q:\N \to \N$ is the shift map, and $\zeta_2^*J \cong \ZZ_2
\times_{[1]} J/1 \cong B([1]) \times_{[1]} J/1$ is isomorphic to the
cylinder of the embedding $J_1 \to J/1$. In particular, for any $J
\in \Posf^+$ equipped with a map $\chi:J \to \N$, we have a map
\begin{equation}\label{J.pl.eq}
J^+ = \zeta^*(J/_\chi\N) \cong J \times_{\N} (\N \setminus \ZZ_\infty)
\to J,
\end{equation}
and since $\sigma:J / \N \to J$ is reflexive, the map
\eqref{J.pl.eq} is $+$-anodyne. The term $\zeta^*q^*(J/_\chi\N)$ in
the decomposition \eqref{J.pl.dec} for $J^+$ is of the same form: we
have $q^*(J/_\chi\N) \cong J/_{q_\dg \circ \chi} \N$, where
$q_\dg:\N \to \N$ sends $l \in \N$ to $\max\{0,l-1\}$.
\end{exa}

\begin{prop}\label{pl.ano.prop}
The class of $+$-anodyne maps is the minimal saturated class of maps
in $\Pos^+$ that is closed under coproducts and standard pushouts,
and contains all anodyne maps, the projection $J \times [1] \to J$
for any $J \in \Posf^+$, and the map \eqref{J.pl.eq} for any
left-bounded map $J \to \N$.
\end{prop}

\proof{} Let $W$ be the minimal class given in the Proposition. Then
it satisfies Definition~\ref{pl.ano.def}~\thetag{i}, and we need to
show that it also satisfies \thetag{ii}. For $J=\V$, this is
Lemma~\ref{st.ano.le}. Next, take $J=\N$ and assume that the maps
$\chi_l:J_l \to \N$, $l=0,1$ are left-bounded, and the comma-fibers
$f/n$, $n \in \Nn$ of the map $f:J_0 \to J_1$ are anodyne. Then by
\eqref{J.pl.eq}, to show that $f \in W$, we may instead consider the
map $f^+:J_0^+ \to J_1^+$ over $\ZZ_\infty$. Moreover, we can
compose the maps $\chi_l^+:J^+_l \to \ZZ_\infty$, $l=0,1$ with the
map $x$ of \eqref{x.z.eq}. Its comma-fibers are countable coproducts
of finite partially ordered sets isomorphic to $\V^o$ or $\ppt$, and
then comma-fibers of the map $f^+$ with respect to $x \circ
\chi_l^+:J_l^+ \to \V$, $l=0,1$ are countable coproducts of anodyne
maps. Therefore $f^+$, hence also $f$ are indeed in $W$.

Now, for any $J \in \Posf^+$ equipped with any map $\chi:J \to \N$,
the height function $\hht:J \to \N$ is left-bounded, and for any $n
\in \N$, the induced map $J /_{\hht} n \to \N$ is left-bounded, so
that the map $(J /_{\hht} n)^+ \to J /_{\hht} n$ of \eqref{J.pl.eq}
is in $W$. Moreover, $\hht$ induces a map $J^+ \to J \to \N$, and we
have $J^+/n \cong (J/_{\hht} n)^+$. Therefore the map
\eqref{J.pl.eq} for $J$ and $\chi$ is also in $W$. Analogously, for
any $J \in \Posf^+$, $S \in \Sets$, and a map $J \to S^<$, the
height function $\hht:J \to \N$ is left-bounded, and we have
$(J/_{\hht}n)^\hush \cong J^\hush/n$ for any $n \in \N$, so that the
map \eqref{J.S} is in $W$. Then exactly the same inductive argument
as in the proof of Proposition~\ref{ano.prop} shows that $W$
satisfies Definition~\ref{pl.ano.def}~\thetag{ii} for any $J_0,J_1
\in \Posf^+$ and $J \in \Posf$. Finally, for an arbitrary $J \in
\Posf^+$, composing the maps $\chi_l$, $l=0,1$ with the height
function $\hht:J \to \N$ reduces us to the case $J = \N$, then
\eqref{J.pl.eq} futher reduces us to $J = \ZZ_\infty$, and $\dim
\ZZ_\infty =1$ is finite.
\endproof

\begin{remark}
Definition~\ref{pos.ano.def}~\thetag{ii} for discrete $J$ is
equivalent to saying that the class of anodyne maps is closed under
those coproducts that exist in $\Posf$. Not all do, since the
dimension of the coproduct of an infinite number of
finite-dimensional sets may not be finite anymore. However, any
coproduct of left-bounded sets is left-bounded, so $\Posf^+$ has
arbitrary coproducts, and we can use simpler language in
Proposition~\ref{pl.ano.prop}.
\end{remark}

\begin{defn}\label{pos.S.ano.def}
For any regular cardinal $\kappa$, the class of {\em
  $\kappa$-anodyne maps} is the minimal saturated class of maps in
$\Posf_\kappa$ that satisfies the conditions \thetag{i}, \thetag{ii}
of Definition~\ref{pos.ano.def}, with $\Posf$ replaced by
$\Posf_\kappa$, and the class of $\kappa^+$-anodyne maps is the
minimal saturated class of maps in $\Posf_\kappa^+$ satisfying
\thetag{i}, \thetag{ii} with $\Posf$ replaced by $\Posf_\kappa^+$. A
partially ordered set $J \in \Posf_\kappa$ is {\em
  $\kappa$-anodyne}, and a partially ordered set $J \in
\Posf_\kappa^+$ is {\em $\kappa^+$-anodyne}, if so is the map $J \to
\ppt$.
\end{defn}

We note that all the results of the last two subsections about
anodyne and $+$-anodyne maps hold for $\kappa$-anodyne maps with the
same proofs, provided one replaces $\Posf$, $\Posf^+$ with
$\Posf_\kappa$, $\Posf^+_\kappa$ and requires that $|S| < \kappa$ in
\eqref{J.S}.

\section{Augmentations and biorders.}\label{bi.sec}

\subsection{Augmentations.}\label{pos.aug.subs}

The category $\Pos$ comes equipped with a full embedding $\phi:\Pos
\to \Cat$, so that the cofibration $\Cat_\idot \to \Cat$ of
Example~\ref{cat.exa} induces the cofibration \eqref{nu.pos}, and
for any category $I$, \eqref{C.bb} defines the category $\Pos \bb I$
of $I$-augmented partially ordered sets. Explicitly, objects in
$\Pos \bb I$ are pairs $\langle J,\alpha_J \rangle$, $J \in \Pos$,
$\alpha_J:J \to I$ a functor, with morphisms from $\langle
J,\alpha_J \rangle$ to $\langle J',\alpha_{J'}\rangle$ given by
pairs $\langle f,\alpha_f \rangle$ of an order-preserving map $f:J
\to J'$ and a map $\alpha_f:\alpha_J \to \alpha_{J'} \circ f$. We
will drop $\alpha_J$ or $\alpha_f$ from notation when it is clear
from the context. For any closed class $v$ of maps in $I$, a
morphism is {\em $v$-strict} if $\alpha_f$ is pointwise in $v$. A
morphism is {\em strict} if it is $\Iso$-strict, or equivalently, if
$\alpha_f$ an isomorphism (or in terms of Subsection~\ref{fun.subs},
if $\langle f,\alpha_f \rangle$ is a functor over $I$). Augmented
partially ordered set and strict maps form a dense subcategory $\Pos
\bbi I \subset \Pos \bb I$, and we denote by
\begin{equation}\label{pos.I.eq}
  \lambda:\Pos \bbi I \to \Pos \bb I
\end{equation}
the embedding functor. We also have the forgetful functors
\begin{equation}\label{pos.I.pi}
\pi:\Pos \bb I \to \Pos, \ \bpi = \pi \circ \lambda:\Pos \bbi I
\to \Pos
\end{equation}
induced by \eqref{fun.bb}, both are fibrations by
Example~\ref{bb.exa}, and if $I$ is rigid, then the fibration $\bpi$
is discrete. Since for any $I$-augmented partially ordered set
$\langle J,\alpha_J \rangle \in \Pos \bb I$ and any object $i \in
I$, the comma-fiber $i \setminus_\alpha J$ is a partially ordered
set, the extended Yoneda embedding \eqref{yo.cat.eq} of
Example~\ref{yo.cat.exa} induces a functor
\begin{equation}\label{E.eq}
\sE_I:\Pos \bb I \to I^o\Pos,
\end{equation}
and this functor is also a fully faithful embedding.

\begin{remark}\label{pos.aug.rem}
Since partially ordered sets are rigid, a pseudofunctor $I^o \to
\Pos \subset \Cat$ is the same thing as an honest functor, so that
if $I$ is essentially small, the Grothendieck construction defines a
full embedding $I^o\Pos \subset \Cat \mm I$ into the category
\eqref{fib.I}. Its essential image consists of fibrations whose
strict fibers are partially ordered sets. For any functor $X:I^o \to
\Pos$, we denote by
\begin{equation}\label{X.ora}
  \overrightarrow{I}X \to I
\end{equation}
the corresponding fibration with fibers $X(i)$ (if $X(i)$ is
discrete for any $i \in I$, then \eqref{X.ora} is the discrete
fibration $IX \to I$ of Example~\ref{IX.exa}). If $I$ itself is a
partially ordered set, then \eqref{X.ora} is a partially ordered set
for any $X$, we have $I^o\Pos \cong \Pos \mm I$, and the functors
\begin{equation}\label{pos.tri.dia}
\Pos \mm I \overset{a}{\longrightarrow} \Pos / I
\overset{b}{\longrightarrow} \Pos \bb I
\overset{\Y}{\longrightarrow} \Pos \mm I
\end{equation}
induced by \eqref{cat.tri.dia}, together with the maps
\eqref{coind.adj}, fall within the scope of Lemma~\ref{3.adj.le}, so
that $a$ is an honest right-adjoint to $\Y \circ b$, and $b$ is an
honest left-adjoint to $a \circ \Y$, while $\Y$ is the Yoneda
embedding \eqref{E.eq}. As another consequence of rigidity, one can
describe the category $\Pos \bb \Pos$ in terms of the Grothendieck
construction: it is equivalent to the subcategory $\Arc(\Pos)
\subset \Ar(\Pos)$ spanned by cofibrations $f:J' \to J$, with such
maps from $J_0' \to J_0$ to $J_1' \to J$ that the induced map $J_0'
\to J_1' \times_{J_1} J_0$ is a cocartesian functor over $J_0$. The
equivalence sends $\langle J,\alpha_J \rangle \in \Pos \bb \Pos$ to
the cofibration $\alpha_J^*\Pos_\idot \to J$.
\end{remark}

For any category $I$ and regular cardinal $\kappa$, we have the
full subcategories
\begin{equation}\label{subc.eq}
\Posff \bb I \subset \Posf_\kappa \bb I \subset \Posf\bb
I,\Posf^+_\kappa \bb I \subset \Posf^+ \bb I \subset \Pos \bb I
\end{equation}
spanned by pairs $\langle J,\alpha_J \rangle$ with $J$ in $\Posff
\subset \Posf_\kappa \subset \Posf,\Posf^+_\kappa \subset
\Posf^+$. We note that if a partially ordered set $J$ is
finite-dimensional, then the same is true for the comma-fiber $i
\setminus J$ for any $i \in I$, so that the functor \eqref{E.eq}
sends $\Posf \bb I$ into $I^o\Posf$. If $I$ is $\Hom$-bounded by
some regular cardinal $\kappa$, \eqref{E.eq} also sends
$\Posf_\kappa \bb I$ into $I^o\Posf_\kappa$.

\begin{defn}\label{good.def}
A {\em good filtration} on a category $I$ is a collection of full
subcategories $I_{\leq n} \subset I$, $n \geq 0$ such that $I_{\leq
  n} \subset I_{\leq n+1}$ and $I = \bigcup I_{\leq n}$.
\end{defn}

For any category $I$ equipped with a good filtration in the sense of
Definition~\ref{good.def}, we define a full subcategory $\Posf
\bb^\bB I \subset \Posf \bb I$ by
\begin{equation}\label{good.eq}
\Posf \bb^\bB I = \bigcup_{n \geq 0} \Posf \bb I_{\leq n} \subset
\Posf \bb I,
\end{equation}
and for any regular cardinal $\kappa$, we let $\Posf_\kappa
\bb^\bB I = \Posf_\kappa \bb I \cap \Posf \bb^\bB
I$. Explicitly, these subcategories consist of augmented partially
ordered sets $\langle J,\alpha \rangle$ such that $\alpha:J \to I$
factors through $I_{\leq n}$ for a large enough $n$; we will call
such sets {\em restricted}. If $I = I_{\leq n}$ for some $n$, this
condition is of course empty. We also say that $\langle J,\alpha
\rangle \in \Posf^+ \bb I$ is {\em locally restricted} if $J/j$ is
restricted for any $j \in J$, and we let
\begin{equation}\label{good.pl.eq}
  \Posf^+ \bb^\bB I \subset \Posf^+ \bb I
\end{equation}
be the full subcategory spanned by locally restricted augmented
partially ordered sets, with $\Posf^+_\kappa \bb^\bB I = \Posf^+
\bb^\bB I \cap \Posf^+_\kappa \bb I$.

\begin{defn}\label{aug.def}
A map $f = \langle f,\alpha_f \rangle:\langle J,\alpha_J \rangle \to
\langle J', \alpha_{J'} \rangle$ in $\Pos \bb I$ is a {\em full
  embedding} if it is strict, and $f$ is a full embedding. A full
embedding $f$ is {\em left} or {\em right} or {\em locally closed}
if $f$ is left or right or locally closed, {\em left-reflexive} if
$f$ is left-reflexive, {\em right-reflexive} if $f$ is
right-reflexive over $I$ (that is, the base change map $\alpha_J
\circ f^\dg \to \alpha_{J'} \circ f \circ f^\dg \to \alpha_{J'}$ for
the right-adjoint map $f^\dg:J' \to J$ is an isomorphism), and {\em
  reflexive} if it is a finite composition of right-reflexive and
left-reflexive maps.
\end{defn}

For any morphism $\langle f,\alpha_f \rangle:\langle J_0,\alpha_0
\rangle \to \langle J_1,\alpha_1 \rangle$ in $\Pos \bb I$, the
cylinder $\Cyl(f)$ has a natural augmentation given by $\alpha_l$ on
$J_l \subset \Cyl(f)$, $l=0,1$, and if $\langle J_l,\alpha_l
\rangle$ are restricted or locally restricted, then so is
$\Cyl(f)$. The embeddings $s:J_0 \to \Cyl(f)$, $t:J_1 \to \Cyl(f)$
are augmented full embeddings in the sense of
Definition~\ref{aug.def}, and the adjoint map $t_\dg:\Cyl(f) \to
J_1$ is augmented by adjunction, so that \eqref{cyl.deco} is a
decomposition in the category $\Pos \bb I$. If $f$ is strict, then
the dual cylinder $\Cyl^o(f)$ also carries a natural augmentation,
we have augmented full embeddings $s:J_1 \to \Cyl^o(f)$, $t:J_0 \to
\Cyl^o(f)$, and \eqref{cyl.o.deco} is a decomposition in the
category $\Pos \bbi I \subset \Pos \bb I$. The cylinder and dual
cylinder construction preserve each of the full subcategories
\eqref{subc.eq}. For any $I$-augmented right-reflexive full embedding
$f:J \to J'$, we also have the augmented version of the retraction
\eqref{refl.dia}, where we augment the product $J' \times [1]$ via
the projection $J' \times [1] \to J'$, and similarly for
left-reflexive full embeddings and dual cylinders. Therefore
Lemma~\ref{refl.sat.le} also holds for the categories $\Posf \bb I$,
$\Posf_\kappa \bb I$, $\Posf^+ \bb I$, $\Posf^+_\kappa \bb I$ and
their restricted versions \eqref{good.eq}, \eqref{good.pl.eq}, with
the same proof. So does Lemma~\ref{cyl.le}.

More generally, if we are given an $I$-augmented partially ordered
set $\langle J,\alpha_J \rangle \in \Pos \bb I$, a map $f:J'_0 \to
J_0$ in $\Pos$, and a map $\chi:J \to J_0$, then the fibered product
$J' = J \times_{J_0} J'_0$ comes equipped with a natural map $f':J'
\to J'_0$, and we will treat $J'$ as an object in $\Pos \bb I$ with
the structure map $\alpha_{J'} = \alpha_J \circ f'$, restricted or
locally restricted if so were $\langle J,\alpha_J \rangle$. In
particular, for any $j \in J_0$, we have the fiber $J_j = J
\times_{J_0} j$ of the map $\chi:J \to J_0$ and its left and right
comma-fibers $J/j = J \times_{J_0} J_0/j$, $j \setminus J = J
\times_{J_0} j \setminus J_0$. With these conventions, just as in
the absolute case $I = \ppt$, any left-closed full embedding
$\langle J',\alpha_{J'} \rangle \subset \langle J,\alpha_J \rangle$
in $\Pos \bb I$ is of the form $\langle J_0,\alpha_J|_{J_0} \rangle
\to \langle J,\alpha_J \rangle$ for a unique order-preserving
characteristic map $\chi:J \to [1]$. Moreover, for any standard
pushout square \eqref{pos.st.sq} with $I$-augmented $J_{01}$, $J_0$,
$J_1$ and $I$-augmented left-closed embeddings $J_{01} \subset J_0$,
$J_{01} \subset J_1$, the coproduct $J$ inherits the augmentation,
so that standard coproducts also exist in $\Posf \bb I$ and other
categories \eqref{subc.eq}, and correspond bijectively to
characteristic maps $\chi:J \to \V$. Lemma~\ref{po.copr.le} then
holds for augmented partially ordered sets with the same
proof. Moreover, if $I$ is equipped with a good filtration, then all
these statements also hold for restricted and locally restricted
$I$-augmented partially ordered sets.

\begin{defn}\label{bary.I.def}
For any $I$-augmented partially ordered set $\langle J,\alpha
\rangle$, its {\em bary\-centric subdivision} $B_IJ$ and {\em dual
  barycentric subdivision} $B^o_IJ$ are given by
\begin{equation}\label{bary.I.eq}
B_I(J) = \langle BJ,\alpha_J \circ \xi \rangle, \qquad B^o_I(J) = \langle
(BJ)^o,\alpha_J \circ \xi_\perp\rangle,
\end{equation}
where $\xi$ and $\xi_\perp$ are the maps \eqref{B.xi}.
\end{defn}

As in the non-augmented case, barycentric subdivision preserves
chain dimension, sends restricted sets to restricted sets, and due
to the way Definition~\ref{bary.I.def} is formulated, the maps
\eqref{B.xi} define strict augmented maps $\xi:B_IJ \to J$,
$\xi_\perp:B^o_IJ \to J$. We cannot say that $B_IJ \cong
B_IJ^o$ since in fact, $J^o$ is not well-defined (the opposite to an
$I$-augmented set has no natural augmentation). However, we still
have a version of Lemma~\ref{B.le}.

\begin{lemma}\label{B.I.le}
Assume given a partially ordered set $J$ with an augmentation
functor $\alpha:J^> \to I$, let $\eps(o):\ppt \to J^>$ be the
embedding \eqref{eps.c} onto the new maximal element $o \in J^>$,
and let $\ppt_\alpha = \langle \ppt,\alpha \circ \eps(o) \rangle \in
\Posf \bb I$. Then $\eps(o):\ppt_\alpha \to J^>$ is a left-reflexive
full embedding, and $B(\eps(o))$, $B(\eps(o))^o$ are augmented
reflexive full embeddings $\ppt_\alpha \to B_I(J^>)$, $\ppt_\alpha
\to B^o_I(J^>)$.
\end{lemma}

\proof{} For $B_I(J^>)$, use the same decomposition $\ppt \subset
B(J)^< \subset B(J^>)$ as in the proof of Lemma~\ref{B.le}, and note
that the first embedding is right-reflexive over $I$ with respect to
the augmentation $\alpha \circ \xi:B(J^>) \to I$. For $B^o_I(J^>)$, use
the opposite decomposition, and observe that the second embedding is
right-reflexive over $I$ with respect to $\alpha \circ \xi_\perp$.
\endproof

\begin{defn}\label{ano.aug.def}
Assume that a category $I$ is equipped with a good filtration in the
sense of Definition~\ref{good.def}. The class of {\em anodyne} maps
is the minimal saturated class of maps in $\Posf \bb^\bB I$
containing all the projections $\langle J,\alpha_J \rangle \times
[1] \to \langle J,\alpha_J \rangle$ and such that
\begin{enumerate}
\setcounter{enumi}{1}
\item for any map $\langle f,\alpha \rangle:\langle J_0,\alpha_0
  \rangle \to \langle J_1,\alpha_1 \rangle$ in $\Posf \bb^\bB I$,
  any $J \in \Posf$, and any maps $\chi_l:J_l \to J$, $l=0,1$ such
  that $f$ is co-lax over $J$ and $f/j:J_0/j \to J_1/j$ is anodyne
  for any $j \in J$, the same holds for $f$.
\end{enumerate}
The class of {\em $+$-anodyne} maps is the minimal saturated class
of maps in $\Posf^+ \bb^\bB I$ satisfying the above conditions
with $\Posf$ replaced by $\Posf^+$, and for any regular cardinal
$\kappa$, the classes of {\em $\kappa$-anodyne} resp.\ {\em
  $\kappa^+$-anodyne} maps are the minimal saturated classes of maps
in $\Posf_\kappa \bb^\bB I$ resp.\ $\Posf_\kappa^+ \bb^\bB I$
satisfying the above conditions with $\Posf$ replaced by
$\Posf_\kappa$ resp.\ $\Posf_\kappa^+$.
\end{defn}

With this definition, most the results of Subsection~\ref{ano.subs}
extend to the augmented setting. In particular, since \eqref{lc.pos}
has a obvious augmented version,
Definition~\ref{ano.aug.def}~\thetag{ii} can be replaced with the
appropriate version of \thetag{ii'} of Corollary~\ref{ano.corr}, and
Lemma~\ref{st.ano.le} also holds with the same proof. Analogously,
Lemma~\ref{B.I.le} implies that the map $\chi:B_I(J) \to J$ is
anodyne, and $\kappa$-anodyne if $|J| < \kappa$, so that a map $f$
is anodyne resp.\ $\kappa$-anodyne iff so is $B_I(f)$. The analogous
statement for $B^o_I(J)$ is also true but requires a proof.

\begin{lemma}\label{B.I.ano.le}
For any $J \in \Posf \bb^\bB I$, the map $\xi_\perp:B^o_I(J) \to J$ is
anodyne, and $\kappa$-anodyne if $|J| < \kappa$.
\end{lemma}

\proof{} Consider the barycentric subdivision $B_I(B^o_I(J))$, and
note that we have two natural maps $B_I(B^o_I(J)) \to J$, namely,
the compositions
\begin{equation}\label{B.I.xi.perp}
\begin{CD}
  B_I(B^o_I(J)) @>{\xi}>> B^o_I(J) @>{\xi_\perp}>> J,\\
  B_I(B^o_I(J)) @>{\xi_\perp}>> B^o_I(J)^o @>{\xi}>> J
\end{CD}
\end{equation}
Explicitly, an element in $B_I(B^o_I(J))$ is a collection $S_0
\subset \dots \subset S_n \subset J$ of totally ordered subsets; the
maps \eqref{B.I.xi.perp} send such a collection to the smallest
element in $S_0$ resp.\ the largest element in $S_n$. The map
$\xi_\perp \circ \xi$ is naturally $I$-augmented. The map $\xi \circ
\xi_\perp$ is not; however, $\xi \circ \xi_\perp \geq
\xi_\perp \circ \xi$, so that $\xi_\perp \circ \xi$ is co-lax over
$J$ with respect to $\xi \circ \xi_\perp:B_I(B^o_I(J)) \to
J$. Moreover, for any $j \in J$, we have $B_I(B^o_I(J)) /_{\xi \circ
  \xi_\perp} j \cong B_I(B^o_I(J/j))$, so by
Definition~\ref{ano.aug.def}, Lemma~\ref{B.I.le} immediately implies
that $\xi \circ \xi_\perp$ is anodyne. Since so is $\xi$, we are
done by the two-out-of-three property.
\endproof

Another useful result concerns a generalization of
Corollary~\ref{ano.corr} to fibrations. In the absolute case, it
would have been trivial by passing to the opposite partially ordered
sets, but in the augmented case, we need a separate proof.

\begin{lemma}\label{B.I.fib.le}
Assume given a map $f:J' \to J$ in $\Posf \bb^\bB I$, and a map
$\chi:J \to J_0$ for some $J_0 \in \Posf$ such that
  \begin{enumerate}
    \item $f:j \setminus J' \to j \setminus J$ is anodyne for any $j
      \in J_0$, or
    \item $\chi$ and $f \circ \chi$ are fibrations whose transition
      functors are strict $I$-aug\-ment\-ed maps, $f$ is cartesian over
      $J_0$, and $f:J'_j \to J_j$ is anodyne for any $j \in J_0$.
  \end{enumerate}
Then $f$ is anodyne, and $\kappa$-anodyne if $|J|,|J'| < \kappa$.
\end{lemma}

\proof{} In the case \thetag{i}, it suffices to note that for any $j
\in J_0$, we have $B_I(j \setminus_\chi J) \cong B_I(J) /_{\chi^o
  \circ \xi_\perp} j$, and similarly for $J'$, so that $B_I(f)$,
hence $f$ is anodyne. In the case \thetag{ii}, the conditions on the
fibration $\chi:J \to J_0$ insure that the embedding $J_j \to j
\setminus J$ is right-reflexive over $I$ for any $j \in J_0$, and
similarly for $\chi \circ f$.
\endproof

Finally, for any $I$-augmented $\langle J,\alpha_J \rangle$ and map
$J \to S^<$, $S \in \Sets$, \eqref{J.S} defines a $I$-augmented
partially ordered set $J^\hush$, and for any map $J \to \N$, so does
\eqref{J.pl.eq}. We then have the following version of
Proposition~\ref{ano.prop} and Proposition~\ref{pl.ano.prop}.

\begin{prop}\label{ano.aug.prop}
The class of anodyne maps is the minimal saturated class of maps in
$\Posf \bb^\bB I$ that contains the maps \eqref{J.S} for any $S
\in \Sets$ and $J \to S^<$, $\langle J,\alpha_J \rangle \in \Posf
\bb I$, satisfies Definition~\ref{pos.ano.def}~\thetag{ii} for
discrete $J$, and is closed with respect to standard pushouts
\eqref{pos.st.sq}. The class of $+$-anodyne maps is the minimal
saturated class of maps in $\Posf^+ \bb^\bB I$ that is closed
under coproducts and standard pushouts, and contains all anodyne
maps, the projection $J \times [1] \to J$ for any $\langle
J,\alpha_J \rangle \in \Posf^+ \bb^\bB I$, and the map
\eqref{J.pl.eq} for any left-bounded map $J \to \N$. The class of
$\kappa$-anodyne resp.\ $\kappa^+$-anodyne maps for some regular
cardinal $\kappa$ is the minimal saturated class of maps in
$\Posf_\kappa \bb^\bB I$ resp.\ $\Posf^+_\kappa \bb^\bB I$ with
the same properties.
\end{prop}

\proof{} Same as Proposition~\ref{ano.prop} and
Proposition~\ref{pl.ano.prop}. \endproof

\subsection{Biordered sets.}\label{rel.subs}

By Remark~\ref{pos.aug.rem}, $\Pos$-augmented partially ordered sets
can be described by cofibrations $J' \to J$, and in such a
situation, the trivial factorization system $\langle \Tot,\Id
\rangle$ on $J$ induces a factorization system of
Lemma~\ref{facto.fib.le}~\thetag{i} on $J'$. It is useful to
axiomatize the situation.

\begin{defn}\label{rel.def}
A {\em biordered set} is a partially ordered set $J$ equipped with
two additional orders $\leq^l$, $\leq^r$ such that $j \leq j'$ iff
there exists $j'' \in J$ such that $j \leq^l j'' \leq^r j'$, and
such a $j''$ is unique whenever it exists.
\end{defn}

Definition~\ref{rel.def} precisely means that $\langle \leq^l,\leq^r
\rangle$ is a factorization system on $J$. In particular, each of
the orders $\leq^l$, $\leq^r$ completely defines the other one. We
will denote by $J^l,J^r \subset J$ the dense subsets defined by
$\leq^l$, $\leq^r$, and we will sometimes define a structure of a
biordered set on some $J \in \Pos$ by specifying $J^l$ or $J^r$.

\begin{exa}\label{2.triv.rel.exa}
The set $[2] = \{0,1,2\}$ equipped with the orders $0
\leq^l 1$, $1 \leq^r 2$ is a biordered set in the sense of
Definition~\ref{rel.def}, denoted by $[2]^\dm$.
\end{exa}

\begin{defn}\label{rel.str.def}
A {\em morphism} between two biordered sets $J_0$, $J_1$ is a map
$f:J_0 \to J_1$ that preserves the orders $\leq$ and $\leq^l$. A
morphism is {\em strict} if it also preserves the order $\leq^r$.
\end{defn}

A biordered set $\langle J,\leq^l,\leq^r \rangle$ is {\em finite} or
{\em has finite chain dimension} or is {\em left} or {\em
  right-bounded} if so is $J \in \Pos$ (and then also $J^l$ and
$J^r$). We denote the category of biordered sets by $\Rel$, and we
let $\bRel \subset \Rel$ the dense subcategory defined by the closed
class of strict maps, with the dense embedding
\begin{equation}\label{rho.eq}
  \rho:\bRel \to \Rel.
\end{equation}
Every map $f:J_0 \to J_1$ of biordered sets defines a map $f^l:J_0^l
\to J_1^l$, and also $f^r:J_0^r \to J_1^r$ if it is strict. The
forgetful functor $U:\Rel \to \Pos$ has a left and right-adjoint
\begin{equation}\label{rel.pos}
L,R:\Pos \to \Rel
\end{equation}
defined by $L(J)=R(J)=L(J)^r=R(J)^l=J$, with discrete $L(J)^l$,
$R(J)^r$, and $R$ itself has a right-adjoint $U^l:\Rel \to \Pos$
sending $J$ to $J^l$. The functors \eqref{rel.pos} factor through
$\bRel \subset \Rel$, and when the functor $L$ is considered as a
functor to $\bRel$, it has a right-adjoint functor $U^r:\bRel \to
\Pos$ sending $J$ to $J^r$.

For any biordered set $J$, the opposite set $J^o$ is naturally
biordered, with $(J^o)^l = (J^r)^o$ and $(J^o)^r = (J^l)^o$, and
this is functorial with respect to strict maps. Since by definition,
the unique factorizations of Definition~\ref{rel.def} are preserved
by strict maps, the category $\bRel$ has all finite limits. Here is
one application of this.

\begin{exa}\label{rel.perf.exa}
Let $J' = \colim_J\phi(j)$ be a perfect colimit in $\Pos$ in the
sense of Definition~\ref{pos.perf.def}, and assume that it is
universal with respect to some map $J' \to J''$. Then any map $J_0
\to L(J'')$ in $\Rel$ is tautologically strict, and we have $L(J')
\times_{L(J'')} J_0 \cong \colim_J (L(\phi(j)) \times_{L(J'')}
J_0)$, both in $\Rel$ and in $\bRel$. Indeed, to check the universal
property, apply Definition~\ref{pos.perf.def} to the maps
$U(J_0),J_0^l,J_0^r \to J''$.
\end{exa}

We define the {\em unfolding} $\I^\dm \subset \Rel$ of a full
subcategory $\I \subset \Pos$ as its preimage $\I^\dm = U^{-1}(\I)$
under the forgetful functor $U$, and we denote $\overline{\I}^\dm =
\I^\dm \cap \bRel$, with the dense embedding $\rho:\overline{\I}^\dm
\to \I^\dm$ of \eqref{rho.eq}. In particular, we have the
ample subcategories $\Posf,\Posf^+,\Posf^- \subset \Pos$ of
Subsection~\ref{pos.dim.subs}, we let $\Relf = \Posf^\dm$, $\Relf^+
= \Posf^{+\dm}$, $\Relf^- = \Posf^{-\dm}$, and for any regular
cardinal $\kappa$, we let $\Relf_\kappa = \Posf_\kappa^\dm$,
$\Relf_\kappa^+ = \Posf_\kappa^{+\dm}$, with the dense subcategories
$\bRelf$, $\bRelf^+$, $\bRelf_\kappa$, $\bRelf^+_\kappa$.

\begin{exa}\label{Z.dm.exa}
Let $\ZZ_\infty$ be as in Example~\ref{Z.exa}. Then setting $2n
\leq^r 2n+1 \geq^l 2n+2$, $n \geq 0$ turns $\ZZ_\infty$ into a
biordered set $\ZZ^\dm_\infty$. For any $m \geq 0$, $\ZZ_m \subset
\ZZ_\infty$ acquires an induced biorder $\ZZ^\dm_m$, and we have
strict biordered maps
\begin{equation}\label{z.dm.eq}
\zeta^\dm:\ZZ_\infty^\dm \to L(\N), \qquad \zeta^\dm_{2m-k}:\ZZ_{2m-k}^\dm
\to L([m]), \quad m \geq 1, k=0,1
\end{equation}
induced by the maps \eqref{zeta.eq} and \eqref{z.m.eq}.
\end{exa}

\begin{exa}\label{bary.rel.exa}
For any partially ordered set $J$, the barycentric subdivision $BJ$
carries a biorder defined by $j \leq^l j'$ iff $j \leq j'$ and
$\xi(j) = \xi(j')$, where $\xi:BJ \to J$ is the map \eqref{B.xi}. If
we interpret elements in $BJ$ as subsets $S \subset J$, then $S
\leq^r S'$ iff $S \subset S'$ is left-closed, and for any $S \subset
S' \subset J$, the middle term $S'' \subset S'$ of the decomposition
$S \subset S'' \subset S'$ of Definition~\ref{rel.def} is the subset
$S'' \subset S'$ of all elements $s \in S'$ with non-empty $s
\setminus S$. If we let $B^\dm(J)$ be $BJ$ with this biorder, then
$\xi$ becomes a strict biordered map $\xi:B^\dm(J) \to L(J)$, and
any map $J \to J'$ of partially ordered sets induces a strict
biordered map $B^\dm(J) \to B^\dm(J')$, so that $\xi$ is a map of
functors $\xi:B^\dm \to L$. Alternatively, one can also consider the
biordered set $B_\dm(J) = B^\dm(J^o)$; this gives a different
biorder on the same partially ordered set $B(J) \cong B(J^o)$, and
we have a functorial map $\xi_\perp:B_\dm \to L \circ \iota$.
\end{exa}

\begin{exa}\label{cofib.rel.exa}
Assume given a cofibration $\pi:J'\to J$ in $\Pos$, and define an
order $\leq^r$ on $J'$ by setting $j \leq^r j'$ iff $j \leq j'$ and
$\pi(j) = \pi(j')$. Then $J'$ is a biordered set, with the order
$\leq^l$ corresponding to maps cocartesian over $J$. The map $\pi$
becomes a strict biordered map $\pi:J' \to R(J)$, and sending
$\pi:J' \to J$ to this map provides a fully faithful embedding
\begin{equation}\label{cofib.rel.eq}
  \Pos \mc J \to \Rel / R(J),
\end{equation}
where $R:\Pos \to \Rel$ is the functor \eqref{rel.pos}, and $\mc$
has the same meaning as in \eqref{fib.I}.
\end{exa}

\begin{defn}\label{rel.refl.def}
A map $f:J_0 \to J_1$ in $\Rel$ is {\em dense} if $U(f)$ is an
isomorphism while $f^l$ is dense, and a {\em full embedding} if both
$f$ and $f^l$ are full embeddings in $\Pos$. A full embedding $f$ is
{\em left} or {\em right} or {\em locally closed} if it is strict,
and $f$ (hence also $f^l$, $f^r$) is left or right or
locally-closed, and {\em left} resp.\ {\em right-reflexive} if both
$f$ and $f^l$ are left resp.\ right-reflexive, with the same adjoint
map $f_\dg=f^l_\dg$ resp.\ $f^\dg = f^{l\dg}$. A full embedding is
{\em reflexive} if it is a composition of left-reflexive and
right-reflexive full embeddings.
\end{defn}

\begin{exa}\label{m.rel.exa}
Let $[2]^\dm$ be as in Example~\ref{2.triv.rel.exa}, and consider
the embedding $m:[1] \to [2]$ onto $\{0,2\} \subset [2]$ of
Example~\ref{cyl.exa} as a biordered map $L([1]) \to [2]$. Then $m$
is a full embedding that is not strict.
\end{exa}

\begin{lemma}\label{bio.le}
Assume given a biordered set $J$ with a largest element $j$, and let
$J^R \subset J$ be the set of elements $j' \in J$ such that $j'
\leq^r j$. Then $J^R \subset J$ is right-closed and left-reflexive,
and the induced order $\leq^l$ on $J^R$ is discrete.
\end{lemma}

\proof{} The adjoint map is provided by
Example~\ref{facto.exa}.
\endproof

\begin{lemma}\label{le-re.le}
Assume given a biordered set $J$, and a left-admissible full subset
$J' \subset U(J)$, with the embedding map $\gamma:J' \to J$ and the
adjoint map $\gamma_\dg:J \to J' \subset J$. Moreover, assume that
$j \leq^l \gamma_\dg(j)$ for any $j \in J$. Then $J'$ with the
induced orders $\leq^l$, $\leq^r$ is a biordered set, $\gamma$ is a
strict full embedding, and $\gamma_\dg$ is a biordered map, so that
$\gamma$ is left-reflexive in the sense of
Definition~\ref{rel.refl.def}.
\end{lemma}

\proof{} Assume given $j \leq j'$ in $J' \subset J$, and consider the
unique factorization $j \leq^l j'' \leq^r j'$ in $J$. Then by
assumption, $j \leq^l j'' \leq^l \gamma_\dg(j'') \leq
\gamma_\dg(j') = j'$, and then by uniqueness, $\gamma_\dg(j'')
\leq^l j''$. Therefore $\gamma_\dg(j'')=j''$, so that $j'' \in J'$,
$J'$ is indeed a biordered set, and $\gamma$ is a strict full
embedding.  Now for any $j \leq^l j'$ in $J$, consider the unique
factorization $\gamma_\dg(j) \leq^l j'' \leq^r \gamma_\dg(j')$. Then
$j \leq^l j' \leq^l \gamma_\dg(j')$ by assumption, and $j \leq^l
\gamma_\dg(j) \leq^l j''$, so that $j'' = \gamma_\dg(j')$ by
uniqueness of the factorization $j \leq^l j'' \leq^r
\gamma_\dg(j')$. Therefore $\gamma_\dg$ is indeed a biordered map.
\endproof

For any biordered set $J$, a locally closed full subset $J' \subset
U(J)$ with the orders $\leq^l$, $\leq^r$ induced from $J$ is also a
biordered set, and the embedding $J' \to J$ is a strict locally
closed full embedding. In particular, for any $n \geq 0$ and
biordered set $J$, the $n$-th skeleton $\sk_n J \subset J$ carries a
natural biorder, and for any map $f:J \to J'$ of biordered sets, and
any element $j \in J'$, the comma-sets $J / j$, $j \setminus J$ and
the fiber $J_j = (J / j) \cap (j \setminus J)$ are biordered in a
natural way. Then as in $\Pos$, for any left-closed embedding $J_0
\to J$ in $\Rel$, $J_0 = \chi^{-1}(0)$ for a unique characteristic
map $\chi:U(J) \to [1]$ -- or equivalently, a unique map $\chi:J \to
R([1])$ -- and dually for right-closed embeddings. Moreover, $\Rel$
has standard pushouts \eqref{pos.st.sq} of left-closed embeddings,
these correspond bijectively to maps $U(J) \to \V$,
Lemma~\ref{po.copr.le} holds with the same proof, and morphisms
between standard pushout squares correspond to maps $f:J' \to J$
such that $U(f)$ is co-lax over $\V$. Since right-closed embeddings
are strict, $\Rel$ also has co-standard pushout squares
corresponding to maps $\chi:U(J) \to V^o$.

For any map $f:J_0 \to J_1$ in $\Rel$, the cylinder $\Cyl(f)$
becomes a biordered set by setting $\Cyl(f)^l = \Cyl(f^l) \subset
\Cyl(f)$. With this definition, \eqref{cyl.deco} induces a
decomposition
\begin{equation}\label{bicyl.deco}
\begin{CD}
  J_0 @>{s}>> \Cyl(f) @>{t_\dg}>> J_1
\end{CD}
\end{equation}
of $f$, where $s:J_0 \to \Cyl(f)$ is a left-closed full embedding,
and $t_\dg$ is left-adjoint to the right-closed left-reflexive full
embedding $t:J_1 \to \Cyl(f)$. We also have retractions
\eqref{refl.dia}, where we ignore dual cylinders and interpret $[1]$
as $R([1])$, so that $J \times R([1])$ is the cylinder of the
identity map $\id:J \to J$. With these conventions,
Lemma~\ref{refl.sat.le} and Lemma~\ref{cyl.le} hold for biordered
sets with the same proofs.

\subsection{From augmentations to biorders.}

A natural source of biordered sets is provided by augmented
partially ordered sets of Subsection~\ref{pos.aug.subs} via
Example~\ref{cofib.rel.exa} and Remark~\ref{pos.aug.rem}. Namely,
the biorder of Example~\ref{cofib.rel.exa} is obviously functorial
with respect to cocartesian maps, so we obtain a functor $\bT:\Pos
\bb \Pos \to \bRel$ and its composition $T = \rho \circ \bT$ with
the dense embedding \eqref{rho.eq}. Then for any category $I$
equipped with a functor $F:I \to \Pos$, we can define functors
\begin{equation}\label{T.F}
\begin{aligned}
  \bT_F&:\Pos \bb I \to \bRel,\\
  \T_F &\cong \rho \circ \bT_F:\Pos \bb I \to \Rel
\end{aligned}
\end{equation}
by composing $\bT$ and $\T$ with the functor $\Pos \bb I \to \Pos
\bb \Pos$ induced by $F$. Explicitly, for any $\langle J,\alpha
\rangle \in \Pos \bb I$, the biordered set $\bT_F(\langle J,\alpha
\rangle)$ is $J_\idot = (F \circ \alpha)^*\Pos_\idot$, with the
biorder of Example~\ref{cofib.rel.exa} corresponding to the
cofibration $J_\idot \to J$. If $F$ takes values in $\Posf \subset
\Pos$, then $I$ acquires a good filtration in the sense of
Definition~\ref{good.def} given by $I_{\leq n} = \{i \in I| \dim
F(i) \leq n \}$, and we have the categories \eqref{good.eq},
\eqref{good.pl.eq} of restricted and locally restricted
$I$-augmented partially ordered sets.

\begin{lemma}\label{T.le}
The functor $\T_F$ of \eqref{T.F} sends left resp.\ right-closed
embeddings to left resp.\ right-closed embeddings, reflexive full
embeddings to reflexive embeddings, cylinder to cylinders, and
standard resp.\ co-standard pushout square to standard
resp.\ co-standard pushout squares. Moreover, if $F$ takes values in
$\Posf_\kappa \subset \Pos$ for some regular $\kappa$, then $\T_F$
sends $\Posf \bb^\bB I$, $\Posf_\kappa \bb^\bB I$, $\Posf^+
\bb^\bB I$, $\Posf^+_\kappa \bb^\bB I$ into $\Relf$,
$\Relf_\kappa$, $\Relf^+$, $\Relf^+_\kappa$.
\end{lemma}

\proof{} Clear. \endproof

\begin{exa}\label{bary.bis.rel.exa}
For any $J \in \Pos$, let $B^\tr(J) \in \Pos \bb \Pos$ be the
barycentric subdivision $BJ$ augmented by the functor $\phi:BJ \to
\Pos$, $S \mapsto S$ corresponding to the cofibration
\eqref{B.xi.dot}. Then $B^\dm_\idot(J) = \T(B^\tr(J))$ is the
biordered set of pairs $\langle S,s \rangle$, $S \in BJ$, $s \in S$,
with the order given by setting $\langle S,s \rangle \leq \langle
S',s' \rangle$ iff $S \subset S'$ and $s \leq s'$, and $\leq^r$
resp.\ $\leq^l$ corresponding to the cases when $S = S'$
resp.\ $s=s'$. The required factorization is given by $\langle S,s
\rangle \leq^l \langle S',s \rangle \leq^r \langle S',s' \rangle$,
and the underlying partially ordered set is $B_\idot(J)$. Note that the
biordered set $B^\dm(J)$ of Example~\ref{bary.rel.exa} is the subset
$B^\dm(J) \subset B^\dm_\idot(J)$ spanned by $\langle S,s \rangle$
such that $s \in S$ is the largest element, and the embedding
$B^\dm(J) \to B^\dm_\idot(J)$ is full and right-reflexive, with the
strict adjoint map sending $\langle S,s \rangle$ to $S /s$.
\end{exa}
 
Under an additional assumption, the correspondence between
$I$-aug\-ment\-ed and biordered sets given by \eqref{T.F} can be
reversed. Namely, note that while the category $\Rel$ does not have
all finite limits, it has finite products, and these are also
products in $\bRel$.

\begin{lemma}
For any biordered sets $J$, $J'$, the set $\Rel(J,J')$ of biordered
maps $J \to J'$ and the subset $\bRel(J,J') \subset \Rel(J,J')$ of
strict maps are biordered in the sense of Definition~\ref{rel.def}
with respect to pointwise orders $\leq^l$, $\leq^r$, so that the
categories $\Rel$ and $\bRel$ are cartesian-closed in the sense of
Definition~\ref{cart.cl.def}.
\end{lemma}

\proof{} Let $f,f':J \to J'$ be two biordered maps such that $f \leq
f'$ pointwise. Since factorizations \eqref{facto.eq} in a
factorization system are functorial, we have a unique factorization
$f \leq f'' \leq f'$ for some map $f'':J \to J'$ preserving the
order $\leq$, and we need to check that it preserves $\leq^l$ and
$\leq^r$ if $f$, $f'$ are strict. To do this, assume given some $j
\leq j'$ in $J$, and consider the canonical commutative diagram
\begin{equation}\label{facto.map}
\begin{CD}
f(j) @>{l}>> f''(j) @>{r}>> f'(j)\\
@V{l}VV @V{l}VV @VV{l}V\\
f(j)'' @>{l}>> f''(j)'' @>{r}>> f'(j)''\\
@V{r}VV @V{r}VV @VV{r}V\\
f(j') @>{l}>> f''(j') @>{r}>> f'(j')
\end{CD}
\end{equation}
provided by the biorder in $J'$, where $l$ resp.\ $r$ in
\eqref{facto.map} indicate the orders $\leq^l$ resp.\ $\leq^r$. Then
$f(j) \leq^l f''(j)'' \leq^r f'(j')$ is the canonical factorization
of $f(j) \leq f'(j')$. If $j \leq^l j'$, so that $f(j) \leq^l
f(j')$, $f'(j) \leq^l f'(j')$, then $f(j)'' = f(j')$, $f'(j)'' =
f'(j')$, and we conclude that $f''(j)'' = f''(j')$, so that $f''(j)
\leq^l f''(j')$. If $j \leq^r j'$ and $f$, $f'$ are strict, then the
same argument shows that $f''(j)=f''(j)''$ and $f''(j) \leq^r
f''(j')$.
\endproof

\begin{lemma}\label{bT.le}
Assume that the functor $F:I \to \Pos$ is fully faithful, and $F(i)$
is connected for any $i \in I$. Then the functor $\bT_F$ of
\eqref{T.F} is also fully faithful. Moreover, the functor $\sE_I$ of
\eqref{E.eq} extends to the category $\bRel \supset \Pos \bb I$, in
that there exists a functor
\begin{equation}\label{E.F}
\bE_F:\bRel \to I^o\Pos
\end{equation}
equipped with an isomorphism $\bE_F \circ \bT_F \cong \sE_I$.
\end{lemma}

\proof{} Denote by $\bRel_F \subset \bRel$ the essential image of
the functor $\bT_F$. Then by definition, for any $U \cong
\bT_F(\langle J,\alpha \rangle) \in \bRel_F$, we have
\begin{equation}\label{U.co}
U^r = \coprod_{j \in J} F(\alpha(j)).
\end{equation}
By virtue of our assumptions on $F$, this means that $J \cong
\pi_0(U^r)$ as a set, with the natural reduction map $\pi:U \to J$,
and the order on $J$ is given by $j \leq j'$ iff there exists $u
\leq u' \in U$ such that $\pi(u) = j$, $\pi(u')=j'$ (recall that by
our convention, a connected category is non-empty, so that
$\pi_0(F(i)) \cong \ppt$ for any $i \in I$). Since $U \in \bRel_F$,
$\pi$ with this order becomes a cofibration, so that the
Grothendieck construction provides a functor $\alpha:J \to \Pos$
into the essential image of $F:I \to \Pos$. Sending $U$ to $\langle
J,\alpha \rangle$ then gives a functor from $\bRel_F$ to $\Pos \bb
I$ inverse to $\bT_F$. To define \eqref{E.F}, it suffices to set
\begin{equation}\label{bE.F.eq}
  \bE_F(J)(i) = \bRel(L(F(i)),J)^l, \qquad i \in I,
\end{equation}
and observe that as a set, $\bRel(L(F(i)),J) \cong \Pos(F(i),J^r)$
by adjunction, and since $F(i)$ is connected, $\Pos(F(i),-)$
commutes with coproducts. Then since $F$ is fully faithful,
\eqref{U.co} immediately provides an isomorphism $\bE_F \circ \bT_F
\cong \sE_I$.
\endproof

\subsection{Bicofibrations.}

As in the usual case, the cylinder is a simplest example of a
cofibration of biordered sets; here is the appropriate general
notion.

\begin{defn}\label{bicof.defn}
A map $f:J' \to J$ in $\Rel$ is a {\em bicofibration}
resp.\ a {\em biordered fibration} if $U(f)$ is a cofibration
resp.\ a fibration, and $f$ is strict.
\end{defn}

\begin{exa}\label{cyl.rel.exa}
For any biordered map $f:J_0 \to J_1$, the projection $\Cyl(f) \to
R([1])$ is a bicofibration.
\end{exa}

\begin{exa}\label{cofib.rel.bis.exa}
In the situation of Example~\ref{cofib.rel.exa}, the map $J'
\to R(J)$ is a bicofibration.
\end{exa}

\begin{exa}\label{R.groth.exa}
More generally, assume given $J \in \Pos$ equipped a functor $F:J
\to \Rel$, and let $J_\idot = (U \circ F)^*\Pos_\idot \to J$ be the
cofibration with fibers $J_{\idot,j} = U(F(j))$, $j \in J$
corresponding to the composition functor $U \circ F:J \to \Pos$ by
the Grothendieck construction. Then setting $J_\idot^l = (U^l \circ
F)^*\Pos_\idot$ turns $J_\idot$ into a biordered set, and the
forgetful map $J_\idot \to J$ becomes a bicofibration $J_\idot \to
R(J)$ with fibers $J_{\idot,j} = F(j)$. If $J=[1]$, this reduces to
Example~\ref{cyl.rel.exa}, and if the functor $F$ factors through
$L:\Pos \to \Rel$, this is Example~\ref{cofib.rel.bis.exa}. Dually,
if $F$ factors through $\bRel \subset \Rel$, we have the fibration
$J^\hdot = (\iota \circ U \circ F)^*\Pos^{\hdot o} \to J^o$, and
setting $(J^\idot)^r = (\iota \circ U^r \circ F)^*\Pos^{\hdot o}$
turns $J^\hdot$ into a biordered set; the projection $J^\hdot \to
L(J^o)$ then becomes a biordered fibration.
\end{exa}

\begin{exa}\label{facto.fib.exa}
Assume given a biordered set $J$ equipped with a cofibration $f:J'
\to U(J)$. Then the factorization system $\langle f^\sharp(l),f^*(r)
\rangle$ of Lemma~\ref{facto.fib.le}~\thetag{i} turns $J'$ into a
biordered set, and $f$ becomes a bicofibration in the sense of
Definition~\ref{bicof.defn}. Its fibers $J'_j$ have discrete order
$\leq^l$ (that is, $J'_j \cong L(U(J'_j))$, $j \in J$).
\end{exa}

\begin{exa}\label{bifib.lc.exa}
Dualizing Example~\ref{facto.fib.exa}, for any biordered set $J$ and
fibration $f:J' \to U(J)$, $J'$ carries a natural biordered
structure such that the fibers $J'_j$ have discrete order
$\leq^r$. In particular, for any partially ordered $J$ and any map
$g:J' \to R(J)$ in $\Rel$, the comma-set $U(J') /_g J$ comes
equipped with a fibration $\sigma:U(J') /_g J \to U(J')$, thus
acquires a biorder such that $(J' /_g R(J))^l \cong
{J'}^l /_{g^l} J$. Denote $J' /_g R(J)$ with this biorder by $J'
/^\dm_g R(J)$. Then \eqref{lc.pos} becomes a decomposition
\begin{equation}\label{lc.bipos}
\begin{CD}
  J' @>{\eta}>> J'/^\dm_g R(J) @>{\tau}>> R(J)
\end{CD}
\end{equation}
in $\Rel$, where $\eta$ is a right-reflexive full embedding, with
adjoint map $\sigma$, and $\tau$ is a bicofibration with fibers
$(J'/^\dm_g R(J))_j \cong J' /_g j$, $j \in J$.
\end{exa}

\begin{exa}\label{bary.rel.rel.exa}
As another application of the dual version of
Example~\ref{facto.fib.exa}, assume given a map $f:J' \to J$ of
partially ordered sets, and consider the fibration $B(f):B(J') \to
B(J)$ of Lemma~\ref{B.cl.le}. Then if we equip $B(J)$ with the
biorder $B^\dm(J)$ of Example~\ref{bary.rel.exa}, $B(J')$ acquires a
biorder such that $j \leq^l j'$ iff $j \leq j'$ and $\xi(f(j)) =
\xi(f(j'))$. We denote $B(J')$ with this biorder by $B^\dm(J'|J)$,
and we note that we have the tautological dense map $B^\dm(J') \to
B^\dm(J'|J)$ (in the two extreme cases, $B^\dm(J|J) \cong B^\dm(J)$
and $B^\dm(J|\ppt) \cong R(B(J))$). If $f$ itself is a fibration,
then $J'$ acquires a biorder $J'_\dm$ such that $j \leq^l j'$ iff $j
\leq j'$ and $f(j)=f(j')$, and we have a strict biordered map
\begin{equation}\label{xi.rel}
  \xi:B^\dm(J'|J) \to J'_\dm
\end{equation}
induced by \eqref{B.xi}.
\end{exa}

Since $\bRel$ has finite limits, for any bicofibration $J' \to J$
and strict map $f:J_0 \to J$, we have the induced bicofibration
$f^*J' = J' \times_J J_0 \to J_0$.

\begin{lemma}\label{l.rel.le}
For any bicofibration $J' \to J$ and strict left-reflexive full
embedding $f:J_0 \to J$, the map $f^*J' \to J'$ is a strict
left-reflexive full embedding.
\end{lemma}

\proof{} Immediately follows from Lemma~\ref{adm.cof.le}.
\endproof

\begin{lemma}\label{bifib.le}
For any map $f:J' \to J$ of biordered sets such that $U(f)$ is a
cofibration, the transition functors $J'_j \to J'_{j'}$ for all $j
\leq j'$ in $J$ are biordered maps. If we let $J''$ be $U(J')$ with
the biorder of Example~\ref{cofib.rel.exa}, then $f$ is a
bicofibration if and only if the identify map $\id:J'' \to J'$ is a
biordered map. In this case, $f^l$ is also a cofibration, and for
any $j \leq^l j'$, its transition functor ${J'}^l_j \to {J'}^l_{j'}$
is induced by the transition functor for $f$.
\end{lemma}

\proof{} For the first claim, take some $j \leq j'$ in $J$ and some
$j_0 \leq^l j_1$ in $J_j$, let $j_0 \leq j_0'$, $j_1 \leq j_1'$ be
the cocartesian lifting, and let $j_0' \leq^l j'' \leq^r j_1'$ be
the unique factorization of $j_0' \leq j_1'$. Then $j_0 \leq j_1,j''
\leq j_1'$ is a commutative square in $J'$, and since $j_0 \leq^l
j_1$, $j'' \leq^r j_1'$, this is a commutative square
\eqref{facto.sq}, so that we have $j_1 \leq j''$. But since $j_1
\leq j_1'$ is cocartesian, this means that $j_1' \leq j''$, so that
$j_1' = j''$ and $j_0' \leq^l j_1'$.

For the second claim, note that $f$ is strict if and only if for any
$j \leq^l j'$ in $J$, and any $j_0 \in J'_j$ with the cocartesian
lifting $j_0 \leq j_0'$, we actually have $j_0 \leq^l j_0$, and this
exactly means that $J'' \to J'$ is a biordered map. The third claim
is then clear.
\endproof

Lemma~\ref{bifib.le} easily implies that all bicofibrations $J' \to
R(J)$ are of the form given in Example~\ref{R.groth.exa}. The
situation is more complicated when the base $J$ of a bicofibration
$J' \to J$ has non-discrete order $J^r$. The simplest example is
$J=L([1])$, where $L$ is the functor \eqref{rel.pos}. In this case,
for any bicofibration $J \to L([1])$, the transition functor $g:J_0
\to J_1$ is a biordered map by Lemma~\ref{bifib.le}, but not all
biordered maps arise in this way.

\begin{defn}\label{bicyl.def}
A biordered map $f:J_0 \to J_1$ is {\em cylindrical} if $\Cyl(U(f))$
with $\Cyl(U(f))^l = J_0^l \copr J_1^l$ is a biordered set.
\end{defn}

For any cylindrical map $f$, we denote $\Cyl(U(f))$ with the biorder
of Definition~\ref{bicyl.def} by $\Cyl^\flat(f)$. We have the
cartesian square
\begin{equation}\label{cyl.fl.sq}
\begin{CD}
  \Cyl^\flat(f) @>>> \Cyl(f)\\
  @VVV @VVV\\
  L([1]) @>{\id}>> R([1]),
\end{CD}
\end{equation}
both vertical arrows are bicofibrations, and every bicofibration $J
\to L([1])$ arises in this way. Both horizontal arrows in
\eqref{cyl.fl.sq} are dense, but the map $\Cyl^\flat(f) \to \Cyl(f)$
is not an isomorphism unless $J_0$ is empty. Note that the
tautological dense map $\id:L([1]) \to R([1])$ is not strict, so
that \eqref{cyl.fl.sq} is a cartesian square in $\Rel$ and not in
$\bRel$ --- if $f$ is not cylindrical, the corresponding fibered
product does not exist.

\begin{exa}
For any biordered $J$, the projection $J \to \ppt$ is cylindrical.
\end{exa}

\begin{exa}\label{cl.rel.exa}
A right-closed full embedding $f:J_0 \to J_1$ is cylindrical, with
$\Cyl^\flat(f)$ given by the co-standard pushout square
\begin{equation}\label{cl.rel.sq}
\begin{CD}
  J_0 @>{f}>> J_1\\
  @V{t}VV @VVV\\
  J_0 \times L([1]) @>>> \Cyl^\flat(f).
\end{CD}
\end{equation}
\end{exa}

\begin{exa}\label{cl.hom.exa}
For any cylindrical map $g:J_0 \to J_1$ with $J=\Cyl^\flat(g)$, and any
connected $I \in \Pos$, the  induced map
$$
g_I:\Rel(R(I),J_0) \to \Rel(R(I),J_1)
$$
is cylindrical, with $\Cyl^\flat(g_I) \cong \Rel(R(I),J)$ (indeed,
since $I$ is connected, a map $I \to J^l \cong J_0^l \copr J_1^l$
must factor through $J_0^l$ or $J_1^l$).
\end{exa}

\begin{exa}\label{Z.cyl.exa}
Consider the map $\zeta^\dm$ of \eqref{z.dm.eq}, and compose it with the
bicofibration $L(s_\dg):L(\N) \to L([1])$, where $s_\dg:\N \to [1]$
is the cofibration of Example~\ref{sN.cl.exa} adjoint to the
left-closed embedding $s:[1] \to \N$. Then the composition $L(s_\dg)
\circ \zeta:\ZZ^\dm_\infty \to L([1])$ is a bicofibration, with
fibers $\ppt = \ZZ^\dm_0$ and $1 \setminus \ZZ^\dm_\infty$, and
$\zeta^\dm$ is cocartesian over $L([1])$. The corresponding cylindrical
map $\ppt \to 1 \setminus \ZZ^\dm_\infty$ is the right-closed
embedding onto $1 \in 1 \setminus \ZZ^\dm_\infty$. For any $m \geq
1$, the composition $L(s_\dg) \circ \zeta^\dm_m:\ZZ^\dm_m \to L([1])$ is
also a cofibration, and $\zeta^\dm_m$ is cocartesian over $L([1])$.

More generally, assume given a bicofibration $\chi:J \to L(\N)$,
with the induced bicofibration $\zeta^{\dm*}J \to \ZZ^\dm_\infty$. Then
the composition
\begin{equation}\label{s.dg.L.eq}
\begin{CD}
  \zeta^{\dm*}J @>{\zeta^\dm}>> J @>{\chi}>> L(\N) @>{L(s_\dg)}>> L([1])
\end{CD}
\end{equation}
is again a bicofibration, and $\zeta^\dm:\zeta^{\dm*}J \to J$ is
cocartesian over $L([1])$. For any integers $m \geq 1$, $k=0,1$, and
bicofibration $\chi:J \to L([m])$, the map $L(s_\dg) \circ \chi
\circ \zeta^\dm_{2m-k}:\zeta_{2m-k}^{\dm*}J \to L([1])$ is also a
bicofibration, and the map $\zeta^\dm_{2m-k}:\zeta_{2m-k}^{\dm*}J
\to J$ is cocartesian over $L([1])$.
\end{exa}

\begin{exa}\label{refl.rel.exa}
Any right-reflexive map $f:J_0 \to J_1$ with the adjoint map $f_\dg$
is cylindrical, with $\Cyl^\flat(f) \cong \Cyl(f^o_\dg)^o$, where
$f_\dg$ is automatically strict by adjunction, so that $f^o_\dg$ is
well-defined.
\end{exa}

\begin{exa}\label{2.rel.exa}
In particular, in the situation of Example~\ref{refl.rel.exa}, take
$J_0=\ppt$ and $J_1 = R([1])$, with the embedding $f=s:\ppt \to
R([1])$ onto $0 \in [1]$. This is right-reflexive, thus cylindrical,
and $\Cyl^\flat(s)$ is $[2]$ with the biorder $0 \leq^r 1,2$, $1
\leq^l 2$, so that $[2]^l \cong \ppt \copr [1]$ and $[2]^r \cong
\V$. Note that this is different from the biorder of
Example~\ref{2.triv.rel.exa}, so we denote $[2]$ with this biorder
by $[2]'$. The three embeddings $s,t,m:[1] \to [2]$ of
Example~\ref{n.refl.exa} and Example~\ref{cyl.exa} define strict
biordered full embeddings $t:R([1]) \to [2]'$ and $s,m:L([1]) \to
[2]'$, and $m$ is left-reflexive, with the adjoint strict biordered
map $[2]' \to L([1])$ given by $s_\dg$ of Example~\ref{n.refl.exa}.
\end{exa}

\begin{lemma}\label{bicyl.le}
For any cylindrical map $f:J_0 \to J_1$, the corresponding map
$s:J_0 \to \Cyl(f)$ is cylindrical.
\end{lemma}

\proof{} If $J_0=J_1=\ppt$, then $\Cyl(f) \cong R([1])$, and this is
Example~\ref{2.rel.exa}. In the general case, the complementary
order $\Cyl^\flat(s)^r$ is given by the standard pushout square
$$
\begin{CD}
  J_0^r @>{s}>> \Cyl^\flat(f)^r\\
  @V{\id \times s}VV @VVV\\
  J_0^r \times [1] @>>> \Cyl^\flat(s)^r,
\end{CD}
$$
and the projection $\Cyl(f) \to R([1])$ induces a projection
$p:\Cyl^\flat(s) \to [2]'$ that is compatible with both $\leq^l$ and
$\leq^r$ (this is the same projection that appears as the middle
vertical map in \eqref{m.t.dia}). Moreover, for all three strict
full embeddings $s,m:L([1]) \to [2]'$, $t:R([1]) \to [2]'$ of
Example~\ref{2.rel.exa}, $s^*\Cyl^\flat(s) \cong J_0 \times L([1])$,
$m^*\Cyl^\flat(s) \cong \Cyl^\flat(f)$ and $t^*\Cyl^\flat(s) \cong
\Cyl(f)$ are biordered, and the embeddings
\begin{equation}\label{s.t.m.rel.eq}
  s^*\Cyl^\flat(f),m^*\Cyl^\flat(f),t^*\Cyl^\flat(f) \to
  \Cyl^\flat(f)
\end{equation}
are full embeddings with respects to both orders $\leq^l$,
$\leq^r$. Then to construct the unique factorizations of
Definition~\ref{rel.def} for $\Cyl^\flat(s)$, it suffices to observe
that any pair $j \leq j'$ of elements in $\Cyl^\flat(s)$ lies inside
one of the subsets \eqref{s.t.m.rel.eq}, and since all
factorizations in $[2]'$ are trivial -- $p(j) \leq p(j')$ means that
either $p(j) \leq^l p(j')$ or $p(j) \leq^r p(j')$ -- and the
projection $p$ preserves both $\leq^l$ and $\leq^r$, all possible
factorizations of $j \leq j'$ must also lie inside the same
subset. Since all three subsets are biordered, such a factorization
exists and is unique.
\endproof

\begin{exa}\label{J.fl.exa}
For any bicofibration $J \to L([1])$, with biordered fibers $J_0,J_1
\subset J$ and cylindrical transition map $f:J_0 \to J_1$, let
$J^\flat = \Cyl(s)$ be the biordered set of
Lemma~\ref{bicyl.le}. Then the projection $p:J^\flat \to [2]' =
L([1])^\flat$ induced by $\Cyl(f) \to R([1])$ is a bicofibration,
with fibers $J^\flat_0=J^\flat_1 = J_0$, $J^\flat_2 = J_1$ and
transition functors $\id:J_0 \to J_0$, $f:J_0 \to J_1$. The strict
full embedding $m:L([1]) \to [2]'$ induces a strict full embedding
$J \cong m^*J^\flat \to J^\flat$, and this embedding is
left-reflexive by Lemma~\ref{l.rel.le}.
\end{exa}

\subsection{Bianodyne maps.}

For partially ordered sets, we have two ways to define anodyne maps:
Definition~\ref{pos.ano.def} and its equivalent form given in
Corollary~\ref{ano.corr}. For biordered sets, the two approaches are
not equivalent any more, and we use the second one.

\begin{defn}\label{biano.def}
The class of {\em bianodyne maps} is the smallest saturated class of
maps in $\Relf$ such that
\begin{enumerate}
\item the projection $R([1]) \times J \to J$ is bianodyne for any $J
  \in \Relf$, and
\item for any bicofibrations $J',J'' \to J$ in $\Relf$, a map
  $f:J' \to J''$ cocartesian over $J$ and such that the fibers
  $f:J'_j \to J''_j$ are bianodyne for any $j \in J$ is itself
  bianodyne.
\end{enumerate}
The class of {\em $+$-bianodyne maps} is the smallest saturated
class of maps in $\Relf^+$ satisfying \thetag{i}, \thetag{ii} with
$\Relf$ replaced by $\Relf^+$, and for any regular cardinal
$\kappa$, the classes of {\em $\kappa$-bianodyne} resp.\ {\em
  $\kappa^+$-bianodyne} maps are the smallest saturated classes
of maps in $\Relf_\kappa$ resp.\ $\Relf_\kappa^+$ satisfying
\thetag{i},  \thetag{ii} with $\Relf$ replaced by $\Relf_\kappa$
resp.\ $\Relf_\kappa^+$.
\end{defn}

\begin{exa}\label{refl.rel.ano.exa}
Lemma~\ref{refl.sat.le} holds for $\Relf$ and $\Relf_\kappa$ with
the same proof, so that all reflexive maps are bianodyne, and all
reflexive maps in $\Relf_\kappa$ are $\kappa$-bianodyne. All
bianodyne resp.\ $\kappa$-bianodyne maps are tautologically
$+$-bianodyne resp.\ $\kappa^+$-bianodyne.
\end{exa}

\begin{lemma}\label{rel.R.le}
Assume given a map $f:J_0 \to J_1$ in $\Relf^+$, some partially
ordered set $J \in \Posf^+$, and maps $\chi_l:U(J_l) \to J$, $l=0,1$
such that $U(f)$ is co-lax over $J$ and $f/j:J_0/j \to J_1/j$ is
$+$-bianodyne for any $j \in J$. Then $f$ is $+$-bianodyne,
$\kappa^+$-bianodyne if $|J_0|,|J_1|,|J| < \kappa$ and $f/j$, $j \in
J$ is $\kappa^+$-bianodyne, bianodyne if $\dim J_0,\dim J_1,\dim J <
\infty$ and $f/j$ are bianodyne, and $\kappa$-bianodyne if all of
these conditions hold.
\end{lemma}

\proof{} Immediately follows from Example~\ref{refl.rel.ano.exa} and
the decomposition \eqref{lc.bipos} of Example~\ref{bifib.lc.exa}.
\endproof

\begin{corr}\label{T.corr}
For any category $I$ and functor $F:I \to \Posf_\kappa$, with some
regular cardinal $\kappa$, the functor $T_F$ of \eqref{T.F} sends
anodyne resp.\ $+$-anodyne resp.\ $\kappa$-anodyne
resp.\ $\kappa^+$-anodyne maps to bianodyne resp.\ $+$-bianodyne
resp.\ $\kappa$-bianodyne resp.\ $\kappa^+$-bianodyne maps.
\end{corr}

\proof{} For any $J' \in \Pos \bb I$ equipped with a map $J' \to J$
to some $J \in \Pos$, $J'' = T_F(J')$ comes equipped with a map
$J'' \to R(J') \to R(J)$, and $J''/j \cong T_F(J'/j)$ for any $j \in
J$.
\endproof

\begin{corr}\label{B.biano.corr}
For any $J \in \Posf_\kappa$, the biordered map $\xi:B^\dm(J) \to
L(J)$ of Example~\ref{bary.rel.exa} is $\kappa$-bianodyne.
\end{corr}

\proof{} Lemma~\ref{rel.R.le} applied to the maps $\chi_0=\xi:B(J) =
U(B^\dm(J)) \to J$, $\chi_1=\id:J \to J$ immediately reduces us to
the case when $J$ has a largest element $j$, and by induction on
$\dim J$, we may assume that the claim is already proved for $J_0 =
J \ssetminus \{j\}$. Changing notation, we assume the claim known for
$J$, and we need to prove it for $J^>$. To do this, consider the
characteristic map $\chi:J^> \to [1]$ of $J \subset J^>$, and note
that \eqref{B.le.eq} provides a bicofibration
$$
\begin{CD}
B^\dm(J^>) @>{B^\dm(\chi)}>> B^\dm([1]) \cong \V^\dm @>{\xi}>>
L([1])
\end{CD}
$$
with fibers $B^\dm(J^>)_0 \cong B(J)$ and $B^\dm(J^>)_1 \cong
R(B(J)^<)$. Then the map $\xi:B^\dm(J^>) \to L(J^>)$ is cocartesian over
$L([1])$, its fiber $\xi_0$ is the map $\xi:B(J) \to J$, and its
fiber $\xi_1$ is the projection $R(B(J)^<) \to \ppt$ that has a
left-reflexive one-sided inverse $B(\eps(o)):\ppt \to R(B(J)^<)$, so
we done by Example~\ref{refl.rel.ano.exa} and
Definition~\ref{biano.def}~\thetag{ii} over $L([1])$.
\endproof

Corollary~\ref{B.biano.corr} is necessarily limited compared to
Corollary~\ref{B.ano.corr} since Example~\ref{bary.rel.exa} only
works for biordered sets of the form $L(J)$. We also have a more
general statement for the relative barycentric subdivisions of
Example~\ref{bary.rel.rel.exa}.

\begin{lemma}\label{bary.rel.rel.le}
For any fibration $J' \to J$ in $\Posf_\kappa$, the strict biordered
map \eqref{xi.rel} is $\kappa$-bianodyne.
\end{lemma}

\proof{} The same induction as in Corollary~\ref{B.biano.corr}
reduced us to the case when $J'$ and $J$ have largest elements $o'
\in J'$, $o \in J$, and $f(o')=o$. Then the embedding $\eps(o):\ppt
\to J'_o$ extends to a cartesian section $f_\dg:J \to J'$ of the
fibration $f$ that is right-closed and right-adjoint to $f$, and if
we let $J'_0 = (J' \ssetminus f_\dg(J))^>$, $J'_1 = f^{-1}(J
\ssetminus \{o\})^>$, and $J'_{01} = J'_0 \cap J'_1$, then we have a
cartesian cocartesian square
\begin{equation}\label{fib.pos.o.sq}
\begin{CD}
J'_{01} @>>> J'_1\\
@VVV @VVV\\
J'_0 @>>> J'
\end{CD}
\end{equation}
in $\Posf_\kappa$ that is obtained from \eqref{pos.sq} by pullback
with respect to the product $\chi_0 \times \chi_1:J' \to [1]^2$ of
the characteristic maps $\chi_0,\chi_1:J' \to [1]$ of the embeddings
$J' \ssetminus f_\dg(J) \subset J'$ and $f^{-1}(J \ssetminus \{o\})
\subset J'$. We let $J'_{l\dm} \subset J'_\dm$, $l=0,1,01$ be $J'_l$
with the induced biorder, and we note that since $f:J'_\dm \to L(J)$
is right-reflexive, it suffices to shows that $f \circ
\xi:B^\dm(J'|J) \to L(J)$ is $\kappa$-bianodyne. Indeed, recall that
the barycentric subdivision functor $B$ sends \eqref{fib.pos.o.sq}
to a standard pushout square, and let $B^\dm(J'|J)_l$, $l=0,1,01$ be
$B(J'_l) \subset B(J')$ with the biorder induced from
$B^\dm(J'|J)$. Then all the partially ordered sets in
\eqref{fib.pos.o.sq} have largest elements, we have bicofibrations
$\xi \circ B(\chi_1):B^\dm(J'|J)_{01},B^\dm(J'|J)_1 \to L([1])$,
$\xi \circ B(\chi_0):B^\dm(J'|J)_0 \to R([1])$, and by the same
argument as in Corollary~\ref{B.biano.corr}, induction on dimension
shows that \eqref{xi.rel} induces $\kappa$-bianodyne maps
$\xi:B^\dm(J'|J)_l \to J'_{l\dm}$, $l=0,1,01$. Since $f:J'_{1\dm}
\to L(J)$ is also right-reflexive, with the same adjoint, $f \circ
\xi:B^\dm(J'|J)_1 \to L(J)$ is then $\kappa$-bianodyne, and the
since the embedding $J'_{01\dm} \to J'_{0\dm}$ in
\eqref{fib.pos.o.sq} is left-reflexive by Example~\ref{cl.refl.exa},
the embedding $B^\dm(J'|J)_{01\dm} \to B^\dm(J'|J)_0$ is also
$\kappa$-bianodyne. The embedding $B^\dm(J'|J)_1 \to
B^\dm(J'|J)_\dm$ is then a standard pushout of a $\kappa$-bianodyne
map, thus $\kappa$-bianodyne, and $f \circ \xi:B^\dm(J'|J) \to L(J)$
is $\kappa$-bianodyne by the two-out-of-three property.
\endproof

\begin{corr}\label{bary.rel.rel.corr}
For any map $f:J' \to J$ in $\Pos$ with $|J'|,|J| < \kappa$ for some
regular cardinal $\kappa$, the map $B(e):B^\dm([1] \times J'|J) \to
B^\dm(J'|J)$ induced by the projection $e:[1] \times J' \to J'$ is
$\kappa$-bianodyne. If $f$ is right-reflexive, with a fully faithful
right-adjoint map $g:J \to J'$, then the map $B^\dm(J'|J) \to
B^\dm(J|J) \cong B^\dm(J)$ is also $\kappa$-bianodyne.
\end{corr}

\proof{} For the first claim, by Lemma~\ref{rel.R.le}, it suffices
to check that the left comma-fiber $B(e):B^\dm([1] \times J'|J)/S \to
B^\dm(J'|J)/S$ is $\kappa$-bianodyne for any $S \in B(J')$, $S \cong
[n] \subset J'$. However, this comma-fiber is the same map for $J' =
S$ and $J = f(S) \subset J$. Therefore we may assume right away that
$f$ is a surjective map $f:[n] \to [m]$ for some $n \geq m \geq
0$. But then $f$ is a fibration by Example~\ref{n.m.refl.exa}, and
so is $f \circ p:[1] \times [n] \to [m]$, so by
Lemma~\ref{bary.rel.rel.le}, it suffices to show that $p:([1] \times
[n])_\dm \to [n]_\dm$ is $\kappa$-bianodyne. This is obvious since
$([1] \times [n])_\dm \cong R([1]) \times [n]_\dm$. For the second
claim, note that the retraction \eqref{refl.dia} for the adjoint map
$g:J \to J'$ is actually a retraction in $\Posf_\kappa/J$, and that
$B^\dm(-|J)$, just as any other functor, preserves retractions, and
then conclude by exactly the same argument as in
Lemma~\ref{refl.sat.le}.
\endproof

While there is no barycentric subdivision for general biordered
sets, there are other useful subdivision functors. We will construct
one such much later in Section~\ref{seg.sp.sec}; for now, we limit
ourselves to the following simple general criterion.

\begin{lemma}\label{A.le}
Assume given a regular cardinal $\kappa$, an endofunctor $A$ of the
category $\bRelf_\kappa$ that preserves coproducts, and morphisms
$\alpha:A \to \id$, $\chi:U \circ A \to U$ such that $\chi \geq
U(\alpha)$. Moreover, assume the following.
\begin{enumerate}
\item For any left-closed embedding $J_0 \to J$ in $\Relf_\kappa$,
  the corresponding map $A(J_0) \to \chi^{-1}(J_0) \subset A(J)$ is
  $\kappa$-bianodyne.
\item For any projection $e:J \times R([1]) \to J$ in
  $\Relf_\kappa$, $A(e)$ is $\kappa$-bianodyne.
\item For any $J \in \Posf_\kappa$, the map $\alpha:A(L(J)) \to
  L(J)$ is $\kappa$-bianodyne.
\end{enumerate}
Then for any $J \in \Relf_\kappa$, the map $\alpha:A(J) \to J$ is
$\kappa$-bianodyne.
\end{lemma}

\proof{} As in Corollary~\ref{B.biano.corr}, Lemma~\ref{rel.R.le}
over $U(J)$ together with \thetag{i} reduce us to the
case when $J$ has a largest element. Now by \thetag{ii} and the
biordered version of Lemma~\ref{refl.sat.le}, $A(f)$ is
$\kappa$-bianodyne for any reflexive $f$, and Lemma~\ref{bio.le}
reduces us to the case $J=L(J_0^>)$. This is \thetag{iii}.
\endproof

\subsection{Examples.}

By the same argument as in Proposition~\ref{ano.prop},
Definition~\ref{biano.def}~\thetag{ii} for $J=R(\V)$ together with
Lemma~\ref{rel.R.le} immediately implies that the classes of
bianodyne, $+$-bianodyne, $\kappa$-bianodyne and
$\kappa^+$-bianodyne maps are closed with respect to standard
pushouts. Conversely, Lemma~\ref{st.ano.le} holds for biordered sets
with the same proof, and shows that any saturated class $W$ of maps
in $\Relf$, $\Relf^+$ or $\Relf_\kappa$, $\Relf^+_\kappa$ that is
closed under standard pushouts, and satisfies
Definition~\ref{biano.def}~\thetag{i}, also satisfies
Lemma~\ref{rel.R.le} for $J = \V$. We can also generalize
Example~\ref{J.S.exa} and Example~\ref{J.pl.exa}.

\begin{exa}\label{J.S.bi.exa}
Assume given a biordered set $J \in \Relf^+$, a set $S$, and a map
$U(J) \to S^<$. Equip the set $U(J)^\hush$ of \eqref{J.S} with a
biorder by upgrading \eqref{B.S.J} to
\begin{equation}\label{B.S.J.bi}
J^\hush = (J_o \times R(S^>)) \copr_{J_o \times R(S)} \coprod_{s \in
  S}J/s,
\end{equation}
where $J_o$ and $J/s \subset J$, $s \in S$ are equipped with the
induced biorders (equivalently, $J_o \times R(S^>)$ is the cylinder
of the projection $J_o \times S \to J_o$). Then we have $J^{l\hush}
\cong J^{\hush l}$, \eqref{J.S} becomes a strict biordered map
\begin{equation}\label{J.S.bi}
  J^\hush \to J,
\end{equation}
and by the same argument as in Example~\ref{J.S.exa}, it is
$+$-bianodyne, and bianodyne if $\dim J < \infty$. If $|J|,|S| <
\kappa$ for some regular cardinal $\kappa$, then \eqref{J.S.bi} is
$\kappa^+$-bianodyne, and $\kappa$-bianodyne if $\dim J < \infty$.
\end{exa}

\begin{exa}\label{J.pl.bi.exa}
For any biordered set $J \in \Relf^+$ and map $U(J) \to \N$, with
the corresponding map $J \to R(\N)$, we have the bicofibration $J
/^\dm R(\N)$ of Example~\ref{bifib.lc.exa}, and \eqref{J.pl.eq}
becomes a strict biordered map
\begin{equation}\label{J.pl.bi.eq}
J^+ = R(\zeta)^*(J /^\dm R(\N)) \cong J \times_{R(\N)} R(\N
\setminus_\zeta \ZZ_\infty) \to J.
\end{equation}
Then we have $J^{+l} \cong J^{l+}$, and Lemma~\ref{Z.m.le} works for
biordered sets with the same proof and shows that \eqref{J.pl.bi.eq}
is $+$-bianodyne, and $\kappa^+$-bianodyne if $|J| < \kappa$. Note
for any bicofibration $J \to R(\N)$, \eqref{J.pl.dec} induces a
standard coproduct decomposition
\begin{equation}\label{J.pl.dec.bi}
R(\zeta)^*J \cong (R(\zeta)^*J) \copr_{J_1} R(\zeta)^*R(q)^*J,
\end{equation}
where $R(\zeta_2)^*J \cong R(\V^o) \times_{R([1])} J/1$ can again be
identified with the cylinder of the embedding $J_1 \to J/1$, and for
any $\chi:J \to R(\N)$, we have $R(q)^*(J /_\chi^\dm R(\N)) \cong J
/^\dm_{R(q_\dg) \circ \chi} R(\N)$.
\end{exa}

\begin{lemma}\label{hom.bi.le}
Assume given a partially ordered set $I$ that has a largest element
$o \in I$. Then $\bRel(L(I),-)$ sends standard pushout squares to
standard pushout squares, and maps \eqref{J.S.bi}
resp.\ \eqref{J.pl.bi.eq} to maps \eqref{J.S.bi}
resp.\ \eqref{J.pl.bi.eq}.
\end{lemma}

\proof{} We have $\bRel(L(I),J) \cong \Pos(I,J^r)$, so the first
claim immediately follows from Lemma~\ref{hom.I.le}. Moreover,
$\bRel(L(I),-)$ tautologically sends products to products, and since
any strict map $f:L(I) \to R(J)$ for any $J \in \Pos$ must be
constant, it commutes with cylinders. Then \eqref{B.S.J.bi} provides
an identification
$$
\bRel(L(I),J^\hush) \cong \bRel(L(I),J)^\hush
$$
for any $\phi:J \to R(S^<)$, where the right-hand side is taken with
respect to the map $\phi \circ \ev_o:\bRel(L(I),J) \to R(S^<)$. For
\eqref{J.pl.bi.eq}, note that any map $f:L(I) \to J$ factors through
$J/n$ for some $n$ -- for example, one can take the image of $f(o)$
under the projection $J \to \N$ -- and then use the same argument
for the decomposition \eqref{J.pl.dec.bi}, and induction on $n$.
\endproof

\begin{lemma}\label{Z.cyl.le}
For any integers $m \geq 1$, $k=0,1$, and a bicofibration $\chi:J
\to L([m])$ in $\Relf$ resp.\ $\Relf^+$, the map
$\zeta_{2m-k}^{\dm*}J \to J$ of Example~\ref{Z.cyl.exa} is bianodyne
resp.\ $+$-bianodyne, and if $k=1$, then it in fact suffices to
require that $L(s_\dg) \circ \chi:J \to L([m-1])$ is a
bicofibration. For any bicofibration $J \to L(\N)$, the map
$\zeta^*J \to J$ is $+$-bianodyne.
\end{lemma}

\proof{} For the first claim, note that it trivially also holds for
$m=k=0$ -- indeed, $\zeta^\dm_0$ is an isomorphism -- and the case
$m=k=1$ is also trivial, since $\zeta^\dm_1$ is an isomorphism as
well. Assume that $2m-k \geq 2$. Then we have the bicofibration
$L(s_\dg) \circ \chi:J \to L([1])$ and the bicofibration
$\zeta_{2m-k}^{\dm*}J \to L([1])$ of Example~\ref{Z.cyl.exa}, and $f
= \zeta^\dm_{2m-k}:\zeta_{2m-k}^{\dm*}J \to J$ is cocartesian over
$L([1])$, so by definition, it suffices to check that its fibers
$f_0$ and $f_1$ are bianodyne resp.\ $+$-bianodyne. Since $f_0 =
\zeta^\dm_0$ is an isomorphism, the issue is $f_1$. We have a
standard pushout square decomposition $\ZZ^\dm_{2m-k} \cong
\ZZ^\dm_2 \copr_{\ppt} \ZZ^\dm_{2(m-1)-k}$ induced by
\eqref{zz.m.sq}, and then the fiber $(\zeta^{\dm*}_{2m-1}J)_1$
decomposes as
\begin{equation}\label{2.m.J}
(\zeta^{\dm*}_{2m-k}J)_1 \cong (R([1]) \times J_1) \copr_{J_1}
\zeta_{2(m-1)-k}^{\dm*}L(t)^*J,
\end{equation}
where $t:[m-1] \to [m]$ is the standard embedding of
Example~\ref{st.cl.exa}. In terms of \eqref{2.m.J}, the map
$f_1:(\zeta^{\dm*}_{2m-k}J)_1 \to L(t)^*J$ then factors as
$$
\begin{CD}
  (\zeta^{\dm*}_{2m-k}J)_1 @>>> \zeta_{2(m-1)-k}^{\dm*}L(t)^*J @>>>
  L(t)^*J,
\end{CD}
$$
where the second map is bianodyne resp.\ $+$-bianodyne by induction,
and the first map has a one-sided inverse that is a standard pushout
of a reflexive left-closed full embedding $(s \times \id):J_1 \to
R([1]) \times J_1$. This proves the first claim; for a
bicofibration $J \to L(\N)$, apply Lemma~\ref{rel.R.le} over $J=\N$.
\endproof

\begin{exa}\label{beta.rel.exa}
In the trivial case $J = L([m])$ and $\chi:L([m]) \to L([m])$ equal
to $\id$, Lemma~\ref{Z.cyl.le} shows that the map $\zeta^\dm_{2m}$
of \eqref{z.dm.eq} is bianodyne. As an immediate corollary, note
that for any $n \geq 0$, $\beta_n$ of \eqref{zeta.n.xi} induces a
strict biordered map
\begin{equation}\label{beta.n.dm}
  \beta_n^\dm:\ZZ_{2n}^\dm \to B^\dm([n]),
\end{equation}
where $B^\dm([n])$ is equipped with the biorder of
Example~\ref{bary.rel.exa}, and then Lem\-ma~\ref{Z.cyl.le},
Corollary~\ref{B.biano.corr} and the two-out-of-three property show
that \eqref{beta.n.dm} is bianodyne. As a further corollary, for any
$n \geq l \geq 0$, the cocartesian square \eqref{seg.del.sq} induces
a strict biordered map
\begin{equation}\label{b.dm.n.l}
  b^{l\dm}_n:B^\dm([l]) \copr_\ppt B^\dm([n-l]) \to B^\dm([n]),
\end{equation}
and this map is also bianodyne: the maps \eqref{beta.n.dm} for $n$,
$l$ and $n-l$ are bianodyne, and $\ZZ^\dm_{2n} \cong \ZZ^\dm_{2l}
\copr_\ppt \ZZ^\dm_{2(n-l)}$ by \eqref{zz.m.sq}.
\end{exa}

Another useful example of bianodyne map is obtained as follows. By
Example~\ref{R.groth.exa}, any functor $F:J \to \bRelf$ from some
partially ordered set $J \in \Posf$ defines a bicofibration
$\pi_\idot:J_\idot \to R(J)$ and a biordered fibration
$\pi^\hdot:J^\hdot \to L(J^o)$ in $\Relf$. We then have the
bicofibration $J_\idot /^\dm R(J) \to R(J)$ of
Example~\ref{bifib.lc.exa}, with the embedding $\eta:J_\idot \to
J_\idot /^\dm R(J)$ over $R(J)$, and since $J_\idot \to R(J)$ itself
is a bicofibration, $\eta$ admits a left-adjoint biordered map
$\eta_\dg:J_\idot /^\dm R(J) \to J_\idot$ cocartesian over
$R(J)$. Being a cofibration, $J_\idot /^\dm R(J) \to R(J)$ also
corresponds to some functor $F_\idot:J \to \bRelf$, and $\eta_\dg$
defines a map of functors $F_\idot \to F$. This in turns provides a
map
\begin{equation}\label{F.J.eq}
L(J^o) \setminus^\iota  J^\hdot \to J^\hdot,
\end{equation}
cartesian over $L(J)$, where $L(J^o) \setminus^\iota J^\hdot \to
L(J^o)$ denotes the biordered fibration corresponding to
$F_\idot$. We note that the construction is functorial: if we have
another functor $F':J \to \bRelf$, then maps $f:F \to F'$ correspond
bijectively to maps $L(J^o) \setminus^\iota f:L(J^o) \setminus^\iota
J^\hdot \to L(J^o) \setminus^\iota {J^\hdot}'$ cartesian over
$L(J^o)$, and for two maps $f,f':F \to F'$, we have
\begin{equation}\label{f.io}
L(J^o) \setminus^\iota f \leq^l L(J^o) \setminus^\iota f'
\end{equation}
as soon as $f \leq^l f'$.

\begin{lemma}\label{F.J.le}
Let $C$ be a closed class of maps in $\Relf$ consisting of strict
full embeddings. Then for any $J \in \Posf$ and functor $F:J \to
\Relf_C$, the corresponding map \eqref{F.J.eq} is bianodyne.
\end{lemma}

\proof{} By Lemma~\ref{rel.R.le}, it suffices to check that $(L(J^o)
\setminus^\iota J^\hdot)/j \to J^\hdot/j$ is bianodyne for any $j
\in J^\hdot$. The induced map $J^\hdot/j \to L(J^o)$ is a biordered
fibration that factors through $L(J^o/\pi^\hdot(j))$, and since the
transition functors of $\pi^\hdot$ are full embeddings, we have
$(L(J^o) \setminus^\iota J^\hdot)/j \cong L(J^o/\pi^\hdot(j))
\setminus^\iota (J^\hdot/j)$. Therefore we may replace $J^\hdot$
with $J^\hdot/j$ and assume right away that $J^\hdot$ has the
largest element $o \in J^\hdot$, and $\pi^\hdot(o)$ is the largest
element in $J^o$. We then have the right-closed left-reflexive full
subset $J^{\hdot R} \subset J^\hdot$ of Lemma~\ref{bio.le}, and since
the transition functors of the biordered fibration $\pi^\hdot$ are
strict, the induced map $J^{\hdot R} \to L(J^o)$ is also a biordered
fibration, the embedding $J^{\hdot R} \to J^\hdot$ is cartesian over
$L(J^o)$, and so is the adjoint map $J^\hdot \to J^{\hdot R}$ of
Lemma~\ref{bio.le}. Then by \eqref{f.io}, the map $L(J^o)
\setminus^\iota J^{\hdot R} \to L(J^o) \setminus^\iota J^\hdot$ is
also left-reflexive, and we are further reduced to the case when
$J^\hdot = J^{\hdot R}$ --- that is, the biorder $\leq^l$ on $J^\hdot$
is discrete. Moreover, $o \in J^\hdot$ defines a cartesian section
$o:L(J^o) \to J^\hdot$ of the biordered fibration $\pi^\hdot$, and
it is a right-closed embedding, so we have the complement
$\overline{J}^\hdot = J^\hdot \ssetminus o(L(J^o))$ and its
characteristic map $\chi:J^\hdot \to L([1])$. Then again, since the
transition functors of $\pi^\hdot$ are full embeddings, the map
$\chi$ inverts cartesian maps, and since it is a bicofibration over
each $j \in J$, it is a bicofibration, with fibers
$\overline{J}^\hdot$, $L(J^o)$ and transition map
$\bpi^\idot:\overline{J}^\hdot \to L(J^o)$ obtained by restricting
$\pi^\hdot$ to $\overline{J}^\idot \subset J^\hdot$ (note that
$\bpi^\hdot$ is trivially cylindrical since the biorder $\leq^l$ on
$J^\hdot$ is discrete). Moreover, the composition map $L(J^o)
\setminus^\iota J^\idot \to J^\hdot \to L([1])$ is also a
bicofibration, with fibers $L(J^o) \setminus^\iota
\overline{J}^\hdot$, $L(J^o) \setminus^\iota L(J^o)$ and transition
map $L(J^o) \setminus^\iota \bpi^\hdot$, and since $\dim
\overline{J}^\idot < \dim J^\idot$, induction on dimension and
Definition~\ref{biano.def}~\thetag{ii} over $L([1])$ further reduce
us to the case $J^\hdot = L(J^o)$, $\pi^\idot = \id$. In this case,
we have $(L(J^o) \setminus^\iota L(J^o))_j \cong R((j \setminus J /
o)^o)$ for any $j \in J$, and since all these sets have the smallest
element, the map \eqref{F.J.eq} admits a fully faithful
right-reflexive left-adjoint.
\endproof

\begin{exa}\label{F.J.exa}
The simplest example for Lemma~\ref{F.J.le} is $J=[1]$, with $F$
constant with value $\ppt$. In this case, $J^\hdot = L([1])$, and
$L(J^o) \setminus^\iota J^\hdot$ is $\V \cong B([1])^o$ with the
biorder opposite to $B^\dm([1])$ (that is, $0 \geq^r o \leq^l 1$).
\end{exa}

\subsection{Reconstruction.}

We now want to prove a reconstruction result saying that the classes
of bianodyne and $+$-bianodyne maps are generated by certain
standard elementary maps, just as in Proposition~\ref{ano.prop} and
Proposition~\ref{pl.ano.prop}. However, we need one more class of
elementary bianodyne maps.

\begin{exa}\label{Z.rel.exa}
Consider the embedding $\eps(0):\ppt \to \V = \{0,1\}^<$ onto $0 \in
\V$, and treat it as a right-closed biordered full embedding
$\eps:\ppt \to R(\V)$. Then it is cylindrical by
Example~\ref{cl.rel.exa}. As a partially ordered set,
$\Cyl(\eps(0))$ is the zigzag partially ordered set $\ZZ_3$ of
Example~\ref{Z.exa} -- that is, it is the set with elements $0$,
$1$, $2$, $3$ and order $0 \leq 1 \geq 2 \leq 3$ -- but the biorder
is different from that of Example~\ref{Z.dm.exa}, so we denote
$\Cyl^\flat(\eps) \in \Rel$ by $\ZZ_3'$ (explicitly, the two orders
are $0 \leq^r 1$ and $1 \geq^l 2 \leq^l 3$). We have the
characteristic map $\chi:\V \to [1]$ of the left-closed subset
$\{o,0\} \subset \V$ -- that is, $\chi(o)=\chi(0)=0$, $\chi(1)=1$ --
and the composition $R(\chi) \circ \eps:\ppt \to R([1])$ is the
embedding $s:\ppt \to R([1])$, so that $R(\chi)$ induces a strict
biordered map
\begin{equation}\label{Z.rel.eq}
  \ZZ_3' \cong \Cyl^\flat(\eps) \to [2]' \cong \Cyl^\flat(s),
  \qquad 0 \mapsto 0, \quad 1,2 \mapsto 1, \quad 3 \mapsto 2
\end{equation}
cocartesian over $L([1])$ (on the level of underlying partially
ordered sets, this is the map $\ZZ_3 \to [2]$ induced by
\eqref{zeta.eq}). Moreover, the map $\chi$ is left-reflexive, with
the adjoint $[1] \to \V$ given by the embedding onto $\{o,1\}
\subset \V$, so that by Example~\ref{refl.rel.ano.exa} and
Definition~\ref{biano.def}~\thetag{ii} for $J=L([1])$, the map
\eqref{Z.rel.eq} is bianodyne, and $\kappa$-bianodyne for any
regular $\kappa$. Note that it is not reflexive (in fact, it does
not even have a section). We also note that while $\ZZ_3'$ can be
obtained as the co-standard pushout square \eqref{cl.rel.sq} of
Example~\ref{cl.rel.exa}, it can also be represented as a standard
pushout square \eqref{zz.m.sq}, by considering left-closed subsets
$\{0,1,2\},\{2,3\} \subset \ZZ_3'$. This gives an identification
\begin{equation}\label{Z.V.rel}
  \ZZ_3' \cong \V^\dm \copr_{\ppt} R([1]),
\end{equation}
where $\V^\dm$ is $\V^o = \{0,1\}^>$ with the biorder $0 \leq^r o
\geq^l 1$ -- or equivalently, $V^\dm = B^\dm([1])$ with the biorder
of Example~\ref{bary.rel.exa}, as in Example~\ref{F.J.exa} -- and
the embeddings $\ppt \to \V^\dm,R([1])$ are onto $1 \in \V^\dm$
resp.\ $0 \in R([1])$.
\end{exa}

\begin{exa}\label{J.fl.hs.exa}
For any bicofibration $J \to L([1])$ in $\Relf^+$ with transition
functor $g:J_0 \to J_1$, let $J^\flat \to [2]'$ be the
bicofibration of Example~\ref{J.fl.exa}, and let $J^\hash =
\ZZ_3' \times_{[2]'} J^\flat$, where the product is taken with
respect to the strict map \eqref{Z.rel.eq}. Then we have a natural
strict map
\begin{equation}\label{J.fl.hs.eq}
  J^\hash \to J^\flat
\end{equation}
cocartesian over $L([1])$, and as in Example~\ref{Z.rel.exa}, it is
$+$-bianodyne by Definition~\ref{biano.def}~\thetag{ii} and
Lemma~\ref{l.rel.le}, and bianodyne if $\dim J < \infty$. If $|J| <
\kappa$ for some regular $\kappa$, then \eqref{J.fl.hs.eq} is
$\kappa^+$-bianodyne, and $\kappa$-bianodyne if $\dim J <
\infty$. In any case, \eqref{Z.V.rel} induces an identification
\begin{equation}\label{J.hs.eq}
J^\hash \cong (J_0 \times \V^\dm) \copr_{J_0} \Cyl(g),
\end{equation}
where the standard pushout in the right-hand side is taken with
respect to the left-closed embedding onto $J_0 \times \{1\} \subset
J_0 \times \V^\dm$ and the left-closed embedding $s:J_0 \to \Cyl(g)$.
\end{exa}

\begin{prop}\label{biano.prop}
The class of bianodyne maps in $\Relf$ is the minimal saturated
class that is closed under coproducts and standard pushouts, and
contains the map \eqref{J.S.bi} for any $J \to R(S^<)$ of
Example~\ref{J.S.bi.exa}, with the biorder \eqref{B.S.J.bi}, and the
map \eqref{J.fl.hs.eq} for any bicofibration $J \to L([1])$ of
Example~\ref{J.fl.hs.exa}. The class of $+$-bianodyne maps is the
minimal saturated class of maps in $\Relf^+$ that is closed under
coproducts and standard pushouts, and contains bianodyne maps, the
projection $R([1]) \times J \to J$ for any biordered set $J \in
\Relf^+$, and the map \eqref{J.pl.bi.eq} for any left-bounded map $J
\to R(\N)$. For any regular cardinal $\kappa$, the classes of
$\kappa$-bianodyne maps resp.\ $\kappa^+$-bianodyne maps are the
minimal saturated closed classes of maps in $\Relf_\kappa$
resp.\ $\Relf^+_\kappa$ with the same properties.
\end{prop}

\proof{} Let $W$ be the minimal saturated closed class of maps in
$\Relf$ described in the Proposition. We observe right away that for
any finite biordered set $J$, the functor $- \times J:\Relf \to
\Relf$ preserves coproducts, standard pushouts, maps \eqref{J.S.bi} and
maps \eqref{J.fl.hs.eq}, so that it also preserves the class $W$.
All maps in $W$ are bianodyne by Example~\ref{J.S.bi.exa} and
Example~\ref{J.fl.hs.exa}, and we need to prove that conversely, all
bianodyne maps are in $W$. As in Proposition~\ref{ano.prop}, the
projection $R([1]) \times J \to J$ is of the form \eqref{J.S.bi} for
any $J \in \Relf$, with $S = \ppt$, so by
Example~\ref{refl.rel.ano.exa}, $W$ contains all reflexive maps, and
it suffices to check that $W$ satisfies
Definition~\ref{biano.def}~\thetag{ii}.

To do this, use the same induction on $\dim J$ as in the proof of
Proposition~\ref{ano.prop}. If $\dim J = 0$, the claim just means
that $W$ is closed under coproducts. For $J=R(\V)$, the claim
follows from the fact that $W$ is closed under standard pushouts,
and \eqref{J.S.bi} coupled with \eqref{B.S.J.bi} then proves it for
$J=R(S^<)$ for any discrete set $S$. For a general $J$ of some
dimension $n=\dim J$, assume the claim proved for all biordered sets
of dimension $\leq n-1$, take $S = J \ssetminus \sk_{n-1}J$, and
consider the map $p:J \to R(S^<)$ sending any $s \in S \subset J$ to
itself and $\sk_{n-1}J$ to $o \in R(S^<)$. Then the claim for
$R(S^<)$ reduces us to considering the comma-sets $J / s$, $s \in
S$, and the standard pushout square \eqref{sk.s} further reduces us
to the case when $J$ has a largest element $j \in J$. In this
case, Lemma~\ref{bio.le} and Lemma~\ref{l.rel.le} further reduces us
to the case when $J = L(J_0^>)$ for some $J_0 \in \Posf$ of
dimension $\dim J_0 \leq n-1$, and then the characteristic map
$\chi:J_0^> \to [1]$ of $J_0 \subset J_0^>$ is a cofibration with
fibers of dimension $\leq n-1$, so we can compose the bicofibrations
$J',J'' \to J$ with the bicofibration $L(\chi):J \to L([1])$ and
reduce to the case $J=L([1])$. Thus we assume given bicofibrations
$J',J'' \to L([1])$ with fibers $J'_l,J''_l$, $l=0,1$ and transition
functors $g':J_0' \to J_1'$, $g'':J_0'' \to J_1''$, and a map $f:J'
\to J''$ cocartesian over $L([1])$ such that $f_l:J'_l \to J''_l$,
$l=0,1$ is in $W$, and we need to show that $f$ is in $W$.

Indeed, the maps ${J'}^\flat \to J'$, ${J''}^\flat \to J''$ are
reflexive by Example~\ref{J.fl.exa}, so it suffices to prove that
$f^\flat:{J'}^\flat \to {J''}^\flat$ is in $W$. Since the maps
\eqref{J.fl.hs.eq} are in $W$ by assumption, it further suffices to
check that $f^\hash:{J'}^\hash \to {J''}^\hash$ is in $W$. Then by
\eqref{J.hs.eq}, it suffices to check that $f_0 \times \id:J_0'
\times \V^\dm \to J_1' \times \V^\dm$ and the map $\wt{f}:\Cyl(g') \to
\Cyl(g'')$ induced by $f$ are in $W$. But on one hand, $\V^\dm$ is
finite, so that $W$ is closed under $- \times \V^\dm$, and $f_0$ is in
$W$, and on the other hand, the embeddings $t:J_1' \to \Cyl(g')$,
$t:J_1'' \to \Cyl(g'')$ are reflexive, thus in $W$, and $f_1$ is
also in $W$.

This finishes the proof in the bianodyne case. For the $+$-bianodyne
case, add the same additional argument as in
Proposition~\ref{pl.ano.prop}, and the proofs in the
$\kappa$-anodyne and $\kappa^+$-bianodynes cases are then exactly
the same (where we note that $- \times J$ for a finite $J$ sends
$\Relf_\kappa$ resp.\ $\Relf^+_\kappa$ into itself for any regular
cardinal $\kappa$).
\endproof

\subsection{Left-finite sets.}\label{lf.subs}

Let us now consider the ample full subcategory $\Posf^\pm \subset
\Pos$ of left-finite partially ordered sets, and denote by
$\Relf^\pm$ its unfolding -- that is, the category of left-finite
biordered sets. We then have the following version of
Definition~\ref{biano.def}.

\begin{defn}\label{pm.biano.def}
The class of {\em $\pm$-bianodyne maps} is the smallest saturated
closed class of maps in $\Relf^\pm$ that satisfies
Definition~\ref{biano.def}~\thetag{i} and \thetag{ii} for a
left-finite map $f$. A map in $\Relf^\pm \cap \Relf$ is {\em
  strongly $\pm$-bianodyne} if it is bianodyne and $\pm$-bianodyne.
\end{defn}

Just as in the situation of Definition~\ref{biano.def}, all
reflexive maps in $\Relf^\pm$ are $\pm$-bianodyne by
Lemma~\ref{refl.sat.le}, and Lemma~\ref{rel.R.le} holds for
$\pm$-biano\-dyne maps with the same proof, provided that $J$ and $f$
are left-finite. All $\pm$-bianodyne maps are automatically
$+$-bianodyne, and the converse is sometimes also true: for example,
in the situation of Lemma~\ref{Z.cyl.le}, the maps $\zeta_m^{\dm*}J
\to J$, $\zeta^{\dm*}J \to J$ are left-finite as soon as so is $J$,
and then both are $\pm$-bianodyne by the same argument. The same
goes for all the bianodyne maps of Example~\ref{beta.rel.exa}, the
maps \eqref{J.pl.bi.eq} of Example~\ref{J.pl.bi.exa}, and the maps
\eqref{J.fl.hs.eq} of Example~\ref{J.fl.hs.exa}. However, in the
situation of Example~\ref{J.S.bi.exa}, even if $J$ is left-finite,
$J^\hush$ is left-finite only if $S$ is finite, and for infinite
$S$, \eqref{J.S.bi} is not even a map in $\Relf^\pm$. This can be
corrected as follows. For any set $S$, define a partially ordered
set $S^\hhush$ by the standard coproduct
\begin{equation}\label{S.hh}
S^\hhush = (S \times [1]) \copr_S B(S^>),
\end{equation}
with the left-closed embedding $S \cong B(S) \to B(S^>)$ induced by
the embedding $S \to S^>$, and note that $\xi:B(S^>) \to S^>$
induces a map
\begin{equation}\label{S.hh.h}
\begin{CD}
S^\hhush @>{\id \copr \xi}>> S^\hush \cong (S \times [1]) \copr_S
S^> @>>> S^<
\end{CD}
\end{equation}
sending the second component in \eqref{S.hh} to $o \in S^<$. Then
for any biordered set $J$ equipped with a map $J \to R(S^<)$, modify
\eqref{B.S.J.bi} by setting
\begin{equation}\label{B.S.J.bis.bi}
J^\hhush = (J_o \times R(B(S^>))) \copr_{J_o
  \times R(S)} \coprod_{s \in S}J/s,
\end{equation}
so that $U(J^\hhush) \cong S^\hhush \times_{S^<} U(J)$, and
\eqref{S.hh} induces a strict biordered map
\begin{equation}\label{J.S.hh}
  J^\hhush \to J.
\end{equation}
Note that for any $S$, $B(S^>)$, hence also $S^\hhush$ are
left-finite by Example~\ref{B.dim.exa}, and if so is $J$, then
$J^\hhush$ is also left-finite and \eqref{J.S.hh} is a map in
$\Relf^\pm$.

\begin{lemma}\label{J.hh.le}
For any set $S$ and any $J \in \Relf^\pm$ equipped with a map $J \to
R(S^<)$, the map \eqref{J.S.hh} is $\pm$-bianodyne in the sense of
Definition~\ref{pm.biano.def}.
\end{lemma}

\proof{} Consider the arrow set $\Ar(S^<) \cong (S \times [1])^<$ of
Example~\ref{ar.S.exa}, let $\eps:S^< \to (S \times [1])^<$ be the
left-closed embedding induced by the embedding $S \to S \times [1]$
onto $S \times \{0\} \subset S \times [1]$, and note that the map
\eqref{S.hh} factors as
\begin{equation}\label{S.hh.dia}
\begin{CD}
  S^\hhush @>>> (S \times [1])^< \cong \Ar(S^<) @>{\tau}>> S^<,
\end{CD}
\end{equation}
where the first map is the embedding $S \times [1] \to (S \times
[1])^<$ resp.\ the composition map $\eps \circ \xi^\perp:B(S^>) \to
S^< \to (S \times [1])^<$ on the first resp.\ second component of
\eqref{S.hh}. Then for any $J \in \Relf^\pm$ equipped with a map $J
\to R(S^<)$, \eqref{S.hh.dia} induces a decomposition
\begin{equation}\label{J.hh.dia}
\begin{CD}
  J^\hhush @>>> J /^\dm R(S^<) @>{\tau}>> J
\end{CD}
\end{equation}
of the map \eqref{J.S.hh}, where the map $\tau$ is reflexive, thus
$\pm$-bianodyne by Lemma~\ref{refl.sat.le}, and the map $J^\hhush
\to J /^\dm R(S^<)$ is left-finite and has reflexive left
comma-fibers.
\endproof

\begin{prop}\label{pm.biano.prop}
The class of $\pm$-bianodyne maps in $\Relf^\pm$ is the minimal
saturated class that is closed under coproducts and standard
push\-outs, and contains the projection $R([1]) \times J \to J$ for
any $J \in \Relf^\pm$, the map \eqref{J.S.hh} for any $J \in
\Relf^\pm$ equipped with a map $J \to R(S^<)$, $S \in \Sets$, the
map \eqref{J.fl.hs.eq} for any bicofibration $J \to L([1])$ in
$\Relf^\pm$, and the map \eqref{J.pl.bi.eq} for any left-bounded map
$J \to R(\N)$, $J \in \Relf^\pm$. For any regular cardinal $\kappa$,
the class of $\kappa^\pm$-bianodyne maps is the minimal saturated
closed class of maps in $\Relf^\pm_\kappa$ with the same
properties.
\end{prop}

\proof{} Same as Proposition~\ref{biano.prop}. \endproof

More generally, by Example~\ref{B.dim.exa}, the barycentric
subdivision functor $B$ sends the whole $\Pos$ into $\Posf^\pm
\subset \Pos$, so that the two functors $B^\dm,B_\dm:\Pos \to \Rel$
of Example~\ref{bary.rel.exa} take values in $\Relf^\pm$. It is easy
to adapt Corollary~\ref{B.biano.corr} to the setting of
Definition~\ref{pm.biano.def} but we will need a more general
statement. Assume given a cartesian square
\begin{equation}\label{pm.sq}
\begin{CD}
J' @>>> J'_0\\
@VVV @VV{\gamma}V\\
J @>{\phi}>> J_0
\end{CD}
\end{equation}
in $\Pos$, and
consider the induced map
\begin{equation}\label{pm.eq}
B^\dm(J') \to B^\dm(J) \times_{L(J_0)} L(J_0'),
\end{equation}
where the product is taken with respect to $\phi \circ \xi:B^\dm(J)
\to L(J) \to L(J_0)$.

\begin{lemma}\label{pm.B.le}
For any cartesian square \eqref{pm.sq} with $J_0,J_0' \in \Posf^+$,
the corresponding map \eqref{pm.eq} is $+$-bianodyne. If $J_0,J_0'
\in \Posf^\pm$, then it is a left-finite $\pm$-bianodyne map in
$\Relf^\pm$, and it is strongly $\pm$-bianodyne if $J,J',J_0,J_0'
\in \Posf$.
\end{lemma}

\proof{} By definition, an element in the target of the map
\eqref{pm.eq} is given by a pair $\langle S,j \rangle$, $S \subset
J_0$, $j \in J_0'$ such that $S$ is finite non-empty totally
ordered, with some largest element $s \in S$, and $\gamma(j) =
\phi(s)$. Then we have
\begin{equation}\label{B.S.j}
(B(J) \times_{J_0} J_0')/\langle S,j \rangle \cong B(S) \times_S S',
\quad B(J')/\langle S,j \rangle \cong B(S'),
\end{equation}
where we denote $S' = S \times_{J_0} (J_0'/j)$, and since $S'$ is
then finite, the target of the map \eqref{pm.eq} and the map itself
are left-finite. Then by Lemma~\ref{rel.R.le}, it suffices to check
that the map \eqref{pm.eq} has $\pm$-bianodyne left comma-fibers. By
\eqref{B.S.j}, this amounts to proving the claim after replacing
$J$, $J_0$, $J_0'$ with $S$, $J_0/\gamma(j)$, $J_0'/j$. We can
therefore assume right away that \eqref{pm.sq} is a square in
$\Posff$. But then $\xi:B^\dm(J') \to L(J')$ factors as
$$
\begin{CD}
B^\dm(J') @>>> B(J) \times_{L(J_0)} L(J_0') @>{\xi \times \id}>>
L(J) \times_{L(J_0)} L(J_0') \cong L(J'),
\end{CD}
$$
where the first map is exactly the map \eqref{pm.eq}, so it suffices to
check that $\xi:B^\dm(J) \to L(J)$ is $\pm$-bianodyne for any finite
$J \in \Posff$. The argument for this is the same as in
Corollary~\ref{B.biano.corr}.
\endproof

\begin{exa}\label{B.i.exa}
Take $J_0=\ppt$ and $J_0' = [1]$. Then Lemma~\ref{pm.B.le} shows
that for any $J \in \Pos$, the natural map $B^\dm(J \times [1]) \to
B^\dm(J) \times L([1])$ is left-finite, $\pm$-bianodyne, and
strongly $\pm$-bianodyne if $J \in \Posf$. If we identify $[1] \cong
[1]^o$, then the same argument applied to $J^o$ provides a
left-finite $\pm$-bianodyne map $B_\dm([1] \times J) \to B_\dm(J)
\times L([1])$ that is also strongly $\pm$-bianodyne as soon as $J
\in \Posf$.
\end{exa}

In the situation of Example~\ref{B.i.exa}, we also have the relative
barycentric subdivision $B^\dm(J \times [1]|J)$ of
Example~\ref{bary.rel.rel.exa}, and we have a commutative square
\begin{equation}\label{B.l.r.sq}
\begin{CD}
B^\dm(J \times [1]) @>>> B^\dm(J) \times L([1])\\
@VVV @VVV\\
B^\dm(J \times [1]|J) @>>> B^\dm(J) \times R([1]),
\end{CD}
\end{equation}
where the vertical maps are dense, and the top map is
$\pm$-bianodyne, and strongly $\pm$-bianodyne if $J \in \Posf$. The
projection $B^\dm(J) \times R([1]) \to B^\dm(J)$ is trivially
$\pm$-bianodyne and bianodyne for $J \in \Posf$, so the bottom map
in \eqref{B.l.r.sq} is bianodyne for $J \in \Posf$ by
Corollary~\ref{bary.rel.rel.corr} and the two-out-of-three property,
and the same argument shows that it is always $\pm$-bianodyne. As
another application of Corollary~\ref{bary.rel.rel.corr}, note that
for any $J \in \Pos$, we have a commutative square
\begin{equation}\label{BB.sq}
\begin{CD}
B^\dm(B(J)|J) @>{\xi}>> B^\dm(J)\\
@V{B(\xi)}VV @VV{\xi}V\\
B^\dm(J) @>{\xi}>> L(J)
\end{CD}
\end{equation}
in $\Rel$, where $B^\dm(B(J)|J)$ is taken with respect to $\xi:B(J)
\to J$, and everything in \eqref{BB.sq} except $J$ is actually in
$\Relf^\pm$.

\begin{lemma}\label{B.B.le}
The maps $\xi,B(\xi):B^\dm(B(J)|J) \to B^\dm(J)$ in \eqref{BB.sq}
are $\pm$-bianodyne for any $J \in \Pos$.
\end{lemma}

\proof{} Since $B(\xi) \leq \xi$, both maps are co-lax over $B(J)$
with respect to $\xi:B^\dm(B(J)|J) \to B^\dm(J)$ and $\id:B^\dm(J)
\to B^\dm(J)$, so by Lemma~\ref{rel.R.le}, it suffices to prove that
they have $\pm$-bianodyne left comma-fibers over any $S \in
B(J)$. As in the proof of Lemma~\ref{pm.B.le}, these comma-fibers
over some $S \cong [n] \subset J$ are the same maps for $J = S$, so
we may assume right away that $J = [n]$, for some $n \geq 0$. Then
the whole square \eqref{BB.sq} is in $\Relf^\pm$, and since
$\xi:B^\dm([n]) \to L([n])$ is strongly $\pm$-bianodyne, it suffices
to prove that so is $\xi:B^\dm(B([n])|[n]) \to B^\dm([n])$. But
since $\xi:B([n]) \to [n]$ is right-reflexive by
Example~\ref{BV.exa}, we are done by
Corollary~\ref{bary.rel.rel.corr}.
\endproof

We will also need a corollary of Lemma~\ref{pm.B.le} describing what
the barycentric subdivision functors $B^\dm$, $B_\dm$ do to a map
\eqref{J.S}. Namely, modify \eqref{B.S.J.bis.bi} by setting
\begin{equation}\label{B.S.J.iii.bi}
\begin{aligned}
J^\fhush &= (J_o \times B^\dm(S^>)) \copr_{J_o
  \times R(S)} \coprod_{s \in S}J/s,\\
J^\hushf &= (J_o \times B_\dm(S^>)) \copr_{J_o
  \times R(S)} \coprod_{s \in S}J/s
\end{aligned}
\end{equation}
for any biordered set $J$ equipped with a map $J \to R(S^<)$, $S \in
\Sets$. Note that for any $J \in \Pos$ and $S \in \Sets$, we have
natural biordered maps
\begin{equation}\label{b.J.S}
  B^\dm(J \times S^>) \to B^\dm(J) \times B^\dm(S^>), \
    B_\dm(J \times S^>) \to B_\dm(J) \times B_\dm(S^>),
\end{equation}
and since $B^\dm$, $B_\dm$ send coproducts to coproducts and standard
pushout squares to standard pushout squares, \eqref{B.S.J} induces
standard coproduct decompositions
\begin{equation}\label{bb.J.S}
\begin{aligned}
B^\dm(J^\hush) &\cong B^\dm(J_o \times S^>) \copr_{B^\dm(J_o) \times
  R(S)} \coprod_{s \in S} B^\dm(J / s),\\
B_\dm(J^\hush) &\cong B_\dm(J_o \times S^>) \copr_{B_\dm(J_o) \times
  R(S)} \coprod_{s \in S} B_\dm(J / s).
\end{aligned}
\end{equation}
Then the results of applying $B^\dm$ and $B_\dm$ to \eqref{J.S}
factor as
\begin{equation}\label{b.j.S}
\begin{CD}
B^\dm(J^\hush) @>{a}>> B^\dm(J)^\fhush @>{b}>> B^\dm(J),\\
B_\dm(J^\hush) @>{a^\dg}>> B_\dm(J)^\hushf @>{b^\dg}>> B_\dm(J),
\end{CD}
\end{equation}
where $b$ and $b^\dg$ are \eqref{J.S.hh} with a different biorder,
and $a$, $a^\dg$ are standard pushouts of \eqref{b.J.S} for $J_o$
and $S$.

\begin{corr}\label{pm.B.corr}
For any $J \in \Pos$ resp.\ $\Posf$, the maps \eqref{b.J.S} are
$\pm$-bi\-ano\-dyne resp.\ strongly $\pm$-bianodyne, and the same
holds for the maps $a$, $a^\dg$ in \eqref{b.j.S}.
\end{corr}

\proof{} By Lemma~\ref{rel.R.le}, the classes of $\pm$-bianodyne and
strongly $\pm$-biano\-dyne maps are stable under standard pushouts, so
the second claim follows from the first one. For the first claim,
consider \eqref{b.J.S} as maps over $R(S^<)$ via the projection
\begin{equation}\label{pm.B.dia}
\begin{CD}
  B^\dm(J) \times B^\dm(S^>) @>>> B^\dm(S^>) @>{\xi^\perp}>> R(S^<),\\
  B_\dm(J) \times B_\dm(S^>) @>>> B_\dm(S^>) @>{\xi^\perp}>> R(S^<),
\end{CD}
\end{equation}
and apply again Lemma~\ref{rel.R.le}. Then it suffices to check that
all left comma-fibers of \eqref{b.J.S} are $\pm$-bianodyne
resp.\ strongly $\pm$-bianodyne. Over the largest element $o \in
S^<$, both maps are isomorphisms, and over some element $s \in S
\subset S^<$, the comma-fibers in question are the maps \eqref{b.J.S}
for $S = \ppt$. Thus we may assume right away that $S = \ppt$. Then
$S^> \cong [1]$, and the compositions
\begin{equation}\label{pm.B.1}
\begin{CD}
B^\dm(J \times [1]) @>>> B^\dm(J) \times B^\dm([1]) @>{\id \times
  \xi}>> B^\dm(J) \times L([1])\\
B_\dm(J \times [1]) @>>> B_\dm(J) \times B_\dm([1]) @>{\id \times
  \xi_\perp}>> B_\dm(J) \times L([1])
\end{CD}
\end{equation}
are the $\pm$-bianodyne resp.\ strongly $\pm$-bianodyne maps of
Example~\ref{B.i.exa}, while $\id \times \xi$ and $\id \times
\xi_\perp$ are reflexive.
\endproof

\subsection{Augmented biorders.}\label{aug.rel.subs}

Let us now combine augmented partially ordered sets of
Subsection~\ref{pos.aug.subs} and biordered sets of
Subsection~\ref{rel.subs}. Namely, consider the category $\Rel$
equipped with the composition $\Rel \to \Cat$ of the forgetful
functor $U:\Rel \to \Pos$ and the tautological embedding $\Pos
\subset \Cat$, and use \eqref{C.bb} to define the category $\Rel \bb
I$ of {\em $I$-augmented biordered sets}. Explicitly, object in
$\Rel \bb I$ are pairs $\langle J,\alpha \rangle$ of a biordered set
$J$ and a functor $\alpha:J \to I$, and maps from $\langle J,\alpha
\rangle$ to $\langle J',\alpha' \rangle$ are pairs $\langle
f,\alpha_f \rangle$ of a biordered map $f:J \to J'$ and a morphism
$\alpha_f:\alpha \to \alpha' \circ U(f)$. Such a map is {\em
  $I$-strict} if $\alpha_f$ is invertible, and the class of strict
maps defines a dense subcategory $\Rel \bbi I \subset \Rel \bb I$. We
can also consider maps that are strict in the biordered sense, as in
Definition~\ref{rel.str.def}, so that altogether, we obtain a
commutative diagram
\begin{equation}\label{rel.I.sq}
\begin{CD}
  \bRel / I @>{\rho}>> \Rel / I\\
  @V{\lambda}VV @VV{\lambda}V\\
  \bRel \bb I @>{\rho}>> \Rel \bb I,
\end{CD}
\end{equation}
where $\rho$ and $\lambda$ are dense embeddings. As in
\eqref{pos.I.pi}, we have the forgetful functors
\begin{equation}\label{rel.I.pi}
\pi:\Rel \bb I \to \Rel, \ \bpi = \pi \circ \lambda:\Rel / I
\to \Rel
\end{equation}
induced by \eqref{fun.bb}, both are fibrations by
Example~\ref{bb.exa}, and if $I$ is rigid, then the fibration $\bpi$
is discrete. For any $I$-augmented biordered set $\langle J,\alpha
\rangle$ and biordered set $J'$, the product $J \times J'$ is
naturally augmented via the projection $J \times J' \to J$. The
forgetful functor $U:\Rel \to \Pos$ tautologically extends to a
functor $U:\Rel \bb I \to \Pos \bb I$, so do its adjoints $L$, $R$
of \eqref{rel.pos} and the right-adjoint $U^l$ to $R$. As in
\eqref{subc.eq}, for any regular cardinal $\kappa$, we have
subcategories
\begin{equation}\label{rel.subc.eq}
\Relf_\kappa \bb I \subset \Relf \bb
I,\Relf^+_\kappa \subset \Relf^+ \bb I \subset \Rel \bb I
\end{equation}
spanned by pairs $\langle J,\alpha \rangle$ with $J$ in
$\Relf_\kappa \subset \Relf,\Relf^+_\kappa \subset \Relf^+$. If $I$
is equipped with a good filtration in the sense of
Definition~\ref{good.def}, we also have the subcategories $\Relf
\bb^\bB I \subset \Relf \bb^\bB I$, $\Relf^+ \bb^\bB I \subset
\Relf^+ \bb^\bB I$ of restricted resp.\ locally restricted
$I$-augmented biordered sets, and their $\kappa$-bounded
versions.

\begin{defn}\label{aug.rel.def}
A map $f = \langle f,\alpha_f \rangle:\langle J,\alpha \rangle \to
\langle J', \alpha' \rangle$ in $\Rel \bb I$ is a {\em full
  embedding} if it is $I$-strict, and $f$ is a full embedding in the
sense of Definition~\ref{rel.refl.def}. A full embedding $f$ is {\em
  left} or {\em right} or {\em locally closed} if $f$ is left or
right or locally closed in the sense of
Definition~\ref{rel.refl.def}, {\em left} resp.\ {\em
  right-reflexive} if both $U(f)$ and $U^l(f)$ are left
resp.\ right-reflexive in the sense of Definition~\ref{aug.def},
with the same adjoint, and {\em reflexive} if it is a finite
composition of right-reflexive and left-reflexive maps.
\end{defn}

As in Definition~\ref{aug.def}, for any morphism $\langle f,\alpha_f
\rangle:\langle J,\alpha \rangle \to \langle J',\alpha'\rangle$ in
$\Rel \bb I$, the biordered cylinder $\Cyl(f)$ has a natural
augmentation, and we have the augmented versions of the
decomposition \eqref{bicyl.deco} and retractions \eqref{refl.dia}.
Then the biordered versions of Lemma~\ref{refl.sat.le} and
Lemma~\ref{cyl.le} hold in the augmented case with the same proofs,
as do Lemma~\ref{bio.le} and Lemma~\ref{l.rel.le}. We can also
combine Definition~\ref{ano.aug.def} and Definition~\ref{biano.def}
in the following way.

\begin{defn}\label{biano.I.def}
Assume that the category $I$ is equipped with a good filtration.
The class of {\em bianodyne maps} is the smallest saturated class of
maps in $\Relf \bb^\bB I$ such that
\begin{enumerate}
\item for any $\langle J,\alpha \rangle \in \Relf \bb^\bB I$, the
  projection $R([1]) \times \langle J,\alpha \rangle \to \langle
  J,\alpha \rangle$ is bianodyne, and
\item for any $\langle J',\alpha' \rangle,\langle J'',\alpha''
  \rangle \in \Relf \bb^\bB I$ and bicofibrations $J',J'' \to J$
  in $\Relf$, a map $\langle f,\alpha_f \rangle:\langle J',\alpha'
  \rangle \to \langle J'',\alpha'' \rangle$ such that $f$ is
  cocartesian over $J$ and the fibers $\langle f,\alpha_f
  \rangle:\langle J'_j,\alpha' \rangle \to \langle J''_j,\alpha''
  \rangle$ are bianodyne for any $j \in J$ is itself bianodyne.
\end{enumerate}
The class of {\em $+$-bianodyne maps} is the smallest saturated
class of maps in $\Relf^+ \bb^\bB I$ satisfying \thetag{i},
\thetag{ii} with $\Relf$ replaced by $\Relf^+$, and for any regular
cardinal $\kappa$, the classes of {\em $\kappa$-bianodyne}
resp.\ {\em $\kappa^+$-bianodyne} maps are the smallest saturated
classes of maps in $\Relf_\kappa \bb^\bB I$ resp.\ $\Relf_\kappa^+
\bb^\bB I$ satisfying \thetag{i}, \thetag{ii} with $\Relf$
replaced by $\Relf_\kappa$ resp.\ $\Relf_\kappa^+$.
\end{defn}

If $J'$ in \eqref{lc.bipos} if $I$-augmented, then the comma-set $J'
/^\dm R(J)$ is also augmented via the projection $\sigma:J' /^\dm
R(J) \to J'$, and then Lemma~\ref{rel.R.le} holds in the augmented
case with the same proof. In paritucular, it again implies that the
classed of bianodyne maps of Definition~\ref{biano.I.def} are stable
under standard pushouts. Moreover, the standard bianodyne and
$+$-bianodyne maps \eqref{J.S.bi} and \eqref{J.pl.bi.eq} of
Proposition~\ref{biano.prop} are naturally $I$-augmented as soon as
so is $J$, and then they are bianodyne resp.\ $+$-bianodyne in the
sense of Definition~\ref{biano.I.def}. In the situation of
Example~\ref{J.fl.hs.exa}, if $J$ is augmented, we can augment
$J^\flat$ via the projection $t_\dg:J^\flat = \Cyl(s) \to J$, and
then \eqref{J.fl.hs.eq} with the induced augmentation on $J^\hash$
is also a bianodyne $I$-augmented map, while the embedding $J \to
J^\flat$ of Example~\ref{J.fl.hs.exa} is a left-reflexive
$I$-augmented biordered map.

\begin{prop}\label{biano.I.prop}
For any category $I$ equipped with a good filtration, the class of
bianodyne maps in $\Relf \bb^\bB I$ is the minimal saturated class
that is closed under coproducts and standard pushouts, and contains
the map \eqref{J.S.bi} for any $J \to R(S^<)$ of
Example~\ref{J.S.bi.exa}, with the biorder \eqref{B.S.J.bi}, and the
map \eqref{J.fl.hs.eq} for any bicofibration $J \to L([1])$ of
Example~\ref{J.fl.hs.exa}. The class of $+$-bianodyne maps is the
minimal saturated class of maps in $\Relf^+ \bb^\bB I$ that is
closed under coproducts and standard pushouts, and contains
bianodyne maps, the projection $R([1]) \times J \to J$ for any $J
\in \Relf^+ \bb^\bB I$, and the map \eqref{J.pl.bi.eq} for any
left-bounded map $J \to R(\N)$. For any regular cardinal $\kappa$,
the classes of $\kappa$-bianodyne maps resp.\ $\kappa^+$-bianodyne
maps are the minimal saturated closed classes of maps in
$\Relf_\kappa \bb^\bB I$ resp.\ $\Relf^+_\kappa \bb^\bB I$ with
the same properties.
\end{prop}

\proof{} Same as Proposition~\ref{biano.prop}.
\endproof

\begin{corr}\label{biano.I.corr}
For any object $i \in I$ in a category $I$ equipped with a good
filtration, the functor $\eps(i)_\trr:\Relf \to \Relf \bb^\bB I$
induced by \eqref{cat.phi} sends bianodyne maps to bianodyne maps.
For any regular cardinal $\kappa > |I|$,
$\eps(i)_\trr:\Relf_\kappa \to \Relf_\kappa \bb^\bB I$ sends
$\kappa$-bianodyne maps to $\kappa$-bianodyne maps.
\end{corr}

\proof{} Clear.
\endproof

\begin{exa}\label{bary.rel.I.exa}
For any $I$-augmented partially ordered set $\langle J,\alpha
\rangle$, the barycentric subdivision $B^\dm(J)$ of
Example~\ref{bary.rel.exa} is naturally augmented via $\alpha \circ
\xi:B(J) \to I$, and if we let $B^\dm_I(J) = \langle
B^\dm(J),\alpha \circ \xi\rangle$, then $\xi:B^\dm_I(J) \to L(J)$ is
a map in $\Rel \bb I$. If $I$ has a good filtration and $J \in \Relf
\bb^\bB I$, then $\xi$ is bianodyne, by exactly the same argument
as in Corollary~\ref{B.biano.corr}. The other biorder $B_\dm(J)$ on
$B(J)$ considered in Example~\ref{bary.rel.exa} defines another
functorial $I$-augmented biordered set $B_{I\dm}(J)$; augmentation
is still the same.
\end{exa}

\begin{lemma}\label{A.I.le}
Assume that the category $I$ is equipped with a good filtration, and
assume given a regular cardinal $\kappa$, an endofunctor $A$ of the
category $\bRelf_\kappa \bb^\bB I$ that preserves coproducts, and
morphisms $\alpha:A \to \id$, $\chi:\bpi \circ U \circ A \to \bpi
\circ U$ such that $\chi \geq \bpi(U(\alpha))$. Moreover, assume the
following.
\begin{enumerate}
\item For any left-closed embedding $J_0 \to J$ in $\Relf_\kappa
  \bb^\bB I$, the corresponding map $A(J_0) \to \chi^{-1}(J_0)
  \subset A(J)$ is $\kappa$-bianodyne.
\item For any projection $e:J \times R([1]) \to J$ in $\Relf_\kappa
  \bb^\bB I$, $A(e)$ is $\kappa$-bianodyne.
\item For any $J \in \Posf_\kappa \bb^\bB I$, the map
  $\alpha:A(L(J)) \to L(J)$ is $\kappa$-bianodyne.
\end{enumerate}
Then for any $J \in \Relf_\kappa \bb^\bB I$, the map $\alpha:A(J)
\to J$ is $\kappa$-bianodyne.
\end{lemma}

\proof{} Same as Lemma~\ref{A.le}.
\endproof

Finally, assume given two categories $I$, $I'$, and a functor $F:I'
\to \Pos$. Then \eqref{T.F} extends to functors
\begin{equation}\label{T.I.F}
\begin{aligned}
  \bT_{I,F}&:\Pos \bb (I \times I') \to \bRel \bb I,\\
  \T_{I,F} &\cong \rho \circ \bT_{I,F}:\Pos \bb (I \times I') \to
  \Rel \bb I.
\end{aligned}
\end{equation}
Moreover, assume that $I$ is equipped with a good filtration, and
note that if $F$ takes values in $\Posf \subset \Pos$, then $I'$
also acquires a good filtration.

\begin{lemma}\label{T.I.le}
The functor $\bT_{I,F}$ of \eqref{T.I.F} sends left
resp.\ right-closed embeddings to left resp.\ right-closed
embeddings, reflexive full embeddings to reflexive embeddings,
cylinder to cylinders, and standard resp.\ co-standard pushout
square to standard resp.\ co-standard pushout squares. Moreover, if
$I$ has a good filtration, and $F$ takes values in $\Posf_\kappa
\subset \Pos$ for some regular $\kappa$, then $\T_{I,F}$ sends
$\Posf \bb^\bB (I \times I')$, $\Posf_\kappa \bb^\bB (I \times
I')$, $\Posf^+ \bb^\bB (I \times I')$, $\Posf^+_\kappa \bb^\bB
(I \times I')$ into $\Relf \bb^\bB I$, $\Relf_\kappa \bb^\bB I$,
$\Relf^+ \bb^\bB I$, $\Relf^+_\kappa \bb^\bB I$, and it sends
anodyne, $\kappa$-anodyne, $+$-anodyne, $\kappa^+$-anodyne maps to
bianodyne, $\kappa$-bianodyne, $+$-bianodyne, $\kappa^+$-bianodyne
maps.
\end{lemma}

\proof{} Clear (for the last claim, use Lemma~\ref{rel.R.le} in the
same way as in Corollary~\ref{T.corr}).
\endproof

\begin{exa}\label{bary.bis.I.rel.exa}
As in Example~\ref{bary.bis.rel.exa}, for any $\langle J,\alpha
\rangle \in \Pos \bb I$, let $B^\tr_I(J) \in \Pos \bb (I \times
\Pos)$ be the barycentric subdivision $B_I(J)$ augmented by the
functor $B_I(J) \to I \times \Pos$ sending a totally ordered subset
$S \subset J$ with maximal element $s \in S$ to $\alpha(s) \times
S$. Then if we let $\T_I = \T_{I,\Id}$, $B^\dm_{I\idot}(J) =
\T(B^\tr_I(J))$ is the biordered set $B^\dm_\idot(J)$ of
Example~\ref{bary.bis.rel.exa} with the augmentation induced by the
projection $B_\idot(J) \to B(J)$. The augmented biordered set
$B^\dm_I(J)$ of Example~\ref{bary.rel.I.exa} is $B^\dm(J) \subset
B^\dm_\idot(J)$ with the induced augmentation, and the embedding
$B^\dm_I(J) \to B^\dm_{I\idot}(J)$ is right-reflexive over $I$.
\end{exa}

\chapter{Simplices.}\label{simp.ch}

This chapter is devoted to various simplicial notions that we will
need, and it is rather short since we will not need that much. We
start with standard material about the category $\Delta$, simplicial
sets and nerves of small categories. Then we discuss the Segal
condition and its interpretation in terms of horn embeddings. We
actually need several equivalent refomulations of the condition
given in Lemma~\ref{bva.le}, and also some facts about nerves of
small groupoids. All this is in Section~\ref{del.sec}, and nothing
is new.

Section~\ref{seg.sec} discusses Grothendieck fibrations over
$\Delta$ and Segal condition for those. Examples include simplicial
replacements of \cite{bou3} and various generalizations. We also
introduce the class of ``special maps'' in $\Delta$ --- namely, maps
between ordinals that send the last element to the last element ---
and prove a version of the standard result relating Grothendieck
fibrations over a category and over its simplicial replacement,
Proposition~\ref{spec.prop}. Again, nothing is new. In particular,
morally, if not technically, Proposition~\ref{spec.prop} goes back
at least to \cite{DKHS}.

Another very useful observation that we borrow from \cite{DKHS} is
that simplicial replacements of small categories are Reedy
categories. We discuss those in Section~\ref{reedy.sec}, and here
there is one result that is technically new. Namely, following an
idea from \cite{GZ}, we distinguish a special class of Reedy
categories that we call ``cellular''. These include categories of
simplices of simplicial sets, and share many combinatorial
properties of simplicial sets including a well-defined skeleton
filtration.

We then prove a general result, Proposition~\ref{ano.I.prop}, saying
that a cellular Reedy category admits an ``anodyne resolution'' ---
roughly speaking, a functor from a left-bounded partially ordered
set that has anodyne comma-fibers. For categories of simplices, the
construction is a version of Thomason's double barycentric
subdivision \cite{tho}, \cite{FL}, but it gives slightly more: a
partially ordered set rather than a general small
category. Eventually, for homotopical applications, this allows one
to replace arbitrary small categories with partially ordered sets.

\section{Simplices and nerves.}\label{del.sec}

\subsection{Ordinals and simplices.}\label{del.0.subs}

As usual, we denote by $\Delta \subset \Pos$ the full subcategory
spanned by finite ordinals $[n]=\{0,\dots,n\}$, $n \geq 0$, with the
usual order, and for any integers $n \geq l \geq 0$, we denote by
$s,t:[l] \to [n]$ the embeddings onto the initial resp.\ the
terminal segment of the ordinal $[n]$, as in
Example~\ref{n.refl.exa}. For any $n \geq l \geq 0$, we have the
cocartesian square \eqref{seg.del.sq} in $\Pos$ that is also a
cocartesian square in $\Delta \subset \Pos$. For any ordinal $[n]
\in \Delta$, we have a unique isomorphism $[n] \cong [n]^o$, so that
the involution $\iota:\Pos \to \Pos$ sends $\Delta$ into itself and
induces an involution
\begin{equation}\label{iota.eq}
\iota:\Delta \to \Delta, \qquad [n] \mapsto [n]^o.
\end{equation}
We also have a factorization system $\langle \dd,\ff \rangle$ on
$\Delta$ consisting of surjective resp.\ injective maps
(traditionally, $\ff$ is the class of ``face maps'', and $\dd$ is the
class of ``degeneracies''). For any $n \geq 0$, we denote by
\begin{equation}\label{D.n.eq}
  \D_n \subset \Delta
\end{equation}
the full subcategory spanned by $[l] \in \Delta$, $0 \leq l \leq n$.

For some purposes, it is useful to consider also the empty ordinal;
the corresponding full subcategory in $\Pos$ is naturally identified
with $\Delta^<$, with the empty ordinal corresponding to the initial
object $o \in \Delta^<$. Since a subset of an ordinal is a (possibly
empty) ordinal, the category $\Delta^<$ has pullbacks with respect
to injective maps, so that as in Example~\ref{comma.exa}, the
cofibration $\tau:\Ar(\Delta^<) \to \Delta^<$ is a bifibration over
$\Delta_\ff^< \subset \Delta^<$. Moreover, if we let $\eps =
\eps([0]):\ppt \to \Delta_\ff \subset \Delta^<_\ff$ be the embedding
\eqref{eps.c} onto $[0]$, then the functor
\begin{equation}\label{ar.i.dia}
\Ar(\Delta^<)_{\tau^*(i)} \to \eps_*\eps^*\Ar(\Delta^<)_{\tau^*(\ff)}
\cong \eps_*\Delta^<
\end{equation}
of \eqref{2.adj} is an equivalence over $\Delta^<_\ff$. Thus for for
any $[n] \in \Delta \subset \Delta^<$, we have a natural equivalence
\begin{equation}\label{del.n.pr}
(\Delta/[n])^< \cong \Delta^</[n] \cong (\eps_*\Delta^<)_{[n]} \cong
  \Delta^{<(n+1)},
\end{equation}
with the copies of $\Delta^<$ in the right-hand side numbered by $l
\in [n]$. Explicitly, the projection onto the $l$-th component in
\eqref{del.n.pr} is the transition functor $\eps(l)^*$ for the
embedding $\eps(l):[0] \to [n]$ onto $l \in [n]$, and it sends an
object $f:[m] \to [n]$ in $\Delta^</[n]$ to the ordinal $f^{-1}(l)
\subset [m]$.

Furthermore, \eqref{2.adj} for the constant fibration $\Delta^<_\ff
\times \Delta^< \to \Delta^<$ provides a functor $\Delta^<_\ff \times
\Delta^< \to \eps_*\Delta^<$, and we can consider the composition
\begin{equation}\label{edge.dia}
\begin{CD}
\Delta^<_\ff \times \Delta^< @>>> \eps_*\Delta^< \cong
\Ar(\Delta^<)_{\tau^*{\ff}} @>{\sigma}>> \Delta^<,
\end{CD}
\end{equation}
where the equivalence in the middle is inverse to \eqref{ar.i.dia},
and $\sigma$ is as in Example~\ref{ar.exa}. This sends $\Delta_\ff
\times \Delta$ into $\Delta$, thus induces a functor
\begin{equation}\label{edge.eq}
\edge:\Delta_\ff \times \Delta \to \Delta
\end{equation}
known as the {\em total edgewise subdivision functor}. Explicitly,
$\edge$ sends $[m] \times [n]$ to the product $[m] * [n]$ equipped
with the lexicographic order ($i \times j \leq i' \times j'$ if $i <
i'$ or $i=i'$ and $j \leq j'$). For any fixed $[m] \in \Delta$, we
have the {\em $m$-th edgewise subdivision functor}
\begin{equation}\label{edge.m.eq}
\edge_m:\Delta \to \Delta, \qquad [n] \mapsto [m] * [n]
\end{equation}
obtained by restricting \eqref{edge.eq} to $[m] \times \Delta
\subset \Delta_\ff \times \Delta$.

In keeping with both our convention and general usage, we denote by
$\Delta^o\Sets = \Fun(\Delta^o,\Sets)$ the category of simplicial
sets. For any non-negative integer $n$, the {\em elementary simplex}
$\Delta_n = \Y([n]) \in \Delta^o\Sets$ is the simplicial set
represented by $[n] \in \Delta$. By Example~\ref{cc.sets.exa}, the
category $\Delta^o\Sets$ is cartesian-closed, so that for any
simplicial sets $X$, $Y$, we have the mapping set $X^Y$ given by
\eqref{cc.XY}. All the sets $X^{\Delta_n}$, $n \geq 0$ together form
the functor $\delta^o_*X:\Delta^o \times \Delta^o \to \Sets$ of
\eqref{cc.Yi}. Alternatively, one can consider the full embedding
$\phi:\Delta \to \Pos$, and take the right Kan extension
$\phi^o_*X$; then the diagonal embedding $\delta:\Pos \to \Pos
\times \Pos$ has a right-adjoint cartesian product functor
$\mu=\delta_\dg:\Pos \times \Pos \to \Pos$, and we have
\begin{equation}\label{d.m.eq}
\delta^o_*X \cong (\phi^o \times \phi^o)^*(\phi^o \times
\phi^o)_*\delta^o_*X \cong (\phi^o \times \phi^o)^*\delta^o_*\phi^o_*X
\cong \theta^*\phi^o_*X,
\end{equation}
where we denote $\theta = \mu^o \circ (\phi^o \times
\phi^o):\Delta^o \times \Delta^o \to \Pos^o$. Explicitly, $\theta$
sends an object $[n] \times [m] \in \Delta^o \times \Delta^o$ to the
product $[n] \times [m] \in \Pos^o$. For any $n,m \geq 0$, the
tautogical map $[n] \times [m] \to [n] * [m]$ is order-preserving
($i \leq i'$ and $j \leq j'$ implies that $i \times j \leq i' \times
j'$ lexicographically), so that if we restrict $\theta$ to a functor
$\Delta_\ff^o \times \Delta^o \to \Pos^o$, then the total edgewise
subdivision functor \eqref{edge.eq} fits into a map $\phi^o \circ
\edge^o \to \mu$, and by \eqref{d.m.eq}, this induces a functorial
map
\begin{equation}\label{edge.map}
\delta^o_*X|_{\Delta_\ff^o \times \Delta^o} \to \edge^{o*}X.
\end{equation}
For any simplicial set $X \in \Delta^o\Sets$, its {\em category of
  simplices} $\Delta X$ is its category of elements of
Example~\ref{IX.exa} --- that is, the category of pairs $\langle
[n],x \rangle$, $[n] \in \Delta$, $x \in X([n])$. The forgetful
functor $\Delta X \to \Delta$ is a discrete fibration of
Example~\ref{IX.exa}, with fibers $(\Delta X)_{[n]} \cong
X([n])$. If $X = \Delta_n$ for some $n \geq 0$, then the category of
simplices $\Delta\Delta_n$ is naturally identified with the
comma-fiber $\Delta/[n]$ of the identity functor $\id:\Delta \to
\Delta$; in particular, it has the terminal object corresponding to
$\id:[n] \to [n]$. The {\em augmented category of simplices} is
given by $\Delta^<X = (\Delta X)^<$; for any $n \geq 0$, we have
$\Delta^<\Delta_n \cong \Delta^</[n]$, $[n] \in \Delta$.

A map $f:X \to X'$ between simplicial sets gives rise to a functor
$\Delta(f):\Delta X \to \Delta X'$ over $\Delta$, automatically
cartesian, and if $f$ is injective, $\Delta(f)$ is left-closed by
Example~\ref{cl.X.exa}. If we have another simplicial set $Y$ and a
map $X \to Y$, with $Y' = Y \copr_X X'$, then we have a cartesian
cocartesian square of categories provided by
Lemma~\ref{cocart.le}. The functor $\Delta^< X \to \Delta^< X'$ is
also left-closed. Moreover, if $f$ is an injective map $\Delta_m \to
\Delta_n$ corresponding to an injective map $[m] \to [n]$, then the
corresponding left-closed full embedding $\Delta^<\Delta_m \cong
\Delta^</[m] \to \Delta^<\Delta_n \cong \Delta^</[n]$ is
right-admissible, and then so is the embedding $\Delta^<(X \times
\Delta_m) \to \Delta^<(X \times \Delta_n)$ for any simplicial set
$X$. In particular, we have right-admissible left-closed full
embeddings
\begin{equation}\label{st.X.del}
s,t:\Delta^< X \to \Delta^<(X \times \Delta_1)
\end{equation}
induced by the maps $s,t:[0] \to [1]$.

A simplex $x \in X([n])$ of a simplicial set $X \in \Delta^o\Sets$
is {\em non-degenerate} iff it is not in the image of the map
$X([m]) \to X([n])$ for some non-identical surjective map $[n] \to
[m]$. We denote by $\overline{\Delta}X \subset \Delta X$ the full
subcategory spanned by non-degenerate simplices. A simplicial set $X
\in \Delta^o\Sets$ is {\em finite} if $\overline{\Delta}X$ is
finite.

For any map $\langle [n'],x' \rangle \to \langle [n],x \rangle$
between non-degenerate simplices, the underlying map $f:[n'] \to
[n]$ is automatically injective, so that the forgetful functor
$\overline{\Delta}X \to \Delta$ factors through the dense
subcategory $\Delta_\ff \subset \Delta$. For any simplex $x \in
X([n])$, there exists a unique pair of a non-degenerate simplex $x'
\in X([n'])$ and a surjective map $f:[n] \to [n']$ such that
$x=X(f)(x')$; the simplex $x'$ is the {\em normalization} of the
simplex $x$.

In general, normalization is not functorial. However, consider the
reduction $\Red(\Delta X)$ of the category $\Delta X$ in the sense
of \eqref{red.eq}. Then the functor
\begin{equation}\label{red.X.eq}
\overline{\Delta}X \to \Red(\Delta X)
\end{equation}
induced by \eqref{red.eq} is essentially bijective and full. This
identifies $\Red(\Delta X)$ with the set of non-degenerate simplices
in $X$ ordered by inclusion. Subsets $X' \subset X$ are then in
one-to-one correspondence with left-closed subsets in the partially
ordered set $\Red(\Delta X)$: if we let $P(X)$ be the partially
ordered set of subsets $X' \subset X$ ordered by inclusion, then we
have an isomorphism
\begin{equation}\label{P.X.J}
P(X) \cong L(\Red(\Delta X)), \qquad X' \mapsto \Red(\Delta X'),
\end{equation}
with the inverse isomorphism sending $J \subset \Red(\Delta X)$ to
the subset $X'$ corresponding to the left-closed subcategory $\Delta
X' = R^{-1}(J) \subset \Delta X$.

A simplicial set $X$ is {\em weakly regular} if for any
non-degenerate simplex $x \in X([n])$ and injective map $f:[n'] \to
[n]$, $X(f)(x) \in X([n'])$ is non-degenerate, and {\em regular} if
the Yoneda map $\Delta_n \to X$ corresponding to $x$ is injective
(in words, $X$ is weakly regular if any face of a non-degenerate
simplex is non-degenerate, and regular if all these faces are
distinct). For a weakly regular $X$, the projection
$\overline{\Delta}X \to \Delta_\ff$ is a discrete fibration
corresponding to some $\overline{X}:\Delta_\ff^o \to \Sets$, we have
\begin{equation}\label{ff.X}
  X \cong e^o_!\overline{X},
\end{equation}
where $e:\Delta_\ff \to \Delta$ is the embedding, and the full
embedding $\overline{\Delta}X \subset \Delta X$ is
left-reflexive. If $X$ is regular, then $\overline{\Delta}X$ is a
partially ordered set, the functor \eqref{red.X.eq} is an
isomorphism, and the reduction functor
\begin{equation}\label{red.del.eq}
R:\Delta X \to \overline{\Delta} X \cong \Red(\Delta X)
\end{equation}
induced by \eqref{red.eq} is a fibration left-adjoint to the
embedding $\overline{\Delta}X \subset \Delta X$.

\begin{exa}\label{Bn.exa}
For any $n \geq 0$, the elementary simplex $\Delta_n$ is regular. We
have $\overline{\Delta}\Delta_n \cong \Delta_\ff / [n]$, and sending
$i:[m] \to [n]$ to $i([m]) \subset [n]$ identifies $\Delta_\ff /
[n]$ with the barycentric subdivision $B([n])$ of
Definition~\ref{pos.bary.def}. The adjoint functor
\eqref{red.del.eq} is induced by the adjoint functor $\Ar(\Delta)
\to \Ar^\ff(\Delta)$ of Example~\ref{facto.exa}.
\end{exa}

We say that a subset $X' \subset X$ in a simplicial set $X$ is {\em
  elementary} if so is the corresponding left-closed subset
$\Red(\Delta X') \subset \Red(\Delta X)$ of \eqref{P.X.J} (that is,
$\Red(\Delta X) \ssetminus \Red(\Delta X')$ consists of a single
non-degenerate simplex $\langle [n],x \rangle$). If $X$ is finite,
then any increasing filtration $X_0 \subset \dots \subset X_n=X$ can
be refined to a filtration where all the subsets $X_i \subset
X_{i+1}$ are elementary. If $X = \Delta_n$, then $\Red(\Delta X)
\cong B([n])$ has the largest element $o \in B([n])$, and the
elementary subset $B([n]) \ssetminus o \subset B([n])$ corresponds
to the {\em simplicial sphere}
\begin{equation}\label{sph.eq}
\sigma_n:\SS_{n-1} \subset \Delta_n.
\end{equation}
Explicitly, it is the subset spanned by maps $f:[m] \to [n]$ that
are not surjective (that is, we remove the interior of the simplex,
and for $n=0$, we take $\SS_{-1}$ to be the empty set). In general,
an elementary extension $X' \subset X$ gives rise to a cocartesian
square
\begin{equation}\label{elem.X.sq}
\begin{CD}
\SS_{n-1} @>>> \Delta_n\\
@VVV @VVV\\
X' @>>> X,
\end{CD}
\end{equation}
where $\Delta_n \to X$ is the Yoneda map corresponding to the
unique non-dege\-ne\-rate simplex $\langle [n],x \rangle \in \Red(\Delta
X) \ssetminus \Red(\Delta X')$. Horizontal maps in \eqref{elem.X.sq}
are injective, and if $X$ is regular, so are the vertical maps.

\subsection{Nerves.}

The {\em nerve} $NI$ of a small category $I$ is the simplicial set
$NI:\Delta^o \to \Sets$ sending $[n] \in \Delta$ to the set of
functors $i_\idot:[n] \to I$ (or equivalently, chains \eqref{c.idot}
of length $n$ in $I$). In particular, $NI([0])$ is the set of
objects in $I$, and for any $i,i' \in I$, the fiber of the map
\begin{equation}\label{st.eq}
NI(s) \times NI(t):NI([1]) \to NI([0]) \times NI([0])
\end{equation}
over $i \times i'$ is the set of maps $I(i,i')$. The set $NI([2])$
is the set of composable pairs of maps in $I$, with the embeddings
$s,t:[1] \to [2]$ inducing an isomorphism
\begin{equation}\label{st2.eq}
NI(s) \times NI(t):NI([2]) \to NI([1]) \times_{NI([0])} NI([1]),
\end{equation}
where the fibered product in the right-hand side is taken with
respect to the projections $NI(t),NI(s):NI([1]) \to NI([0])$. The
map $m:[1] \to [2]$ of Example~\ref{cyl.exa} then induces the
composition map
\begin{equation}\label{m.I.eq}
NI(m):NI([2]) \to NI([1]),
\end{equation}
while the unique projection $e:[1] \to [0]$ induces the map
\begin{equation}\label{e.I.eq}
NI(e):NI([0]) \to NI([1])
\end{equation}
sending an object $i \in I$ to the identity arrow $\id:i \to i$.

For any small category $I$, the nerve $NI^o$ of the opposite
category $I^o$ is given by $NI^o \cong \iota^*NI$, where
$\iota:\Delta \to \Delta$ is the involution \eqref{iota.eq}. The
nerve of the category $[n]$, $[n] \in \Delta$ is the elementary
simplex $\Delta_n$. The nerve is functorial, and $N:\Cat \to
\Delta^o\Sets$ is fully faithful and commutes with arbitrary limits
(in particular, with products and fibered products). Moreover, for
any small categories $I_0$, $I_1$, we have $N(\Fun(I_0,I_1)) \cong
NI_1^{NI_0}$. For any small category $I$, we can promote the nerve
functor $N$ to a fully faithful {\em relative nerve} functor
\begin{equation}\label{N.I.eq}
N_I:\Cat \bb I \to I^o\Delta^o\Sets, \quad N_I(I')(i \times [n]) =
N(i \setminus I')([n])
\end{equation}
by applying the extended Yoneda embedding \eqref{yo.cat.eq} and then
taking the usual nerve pointwise.
  
\begin{remark}\label{N.ff.rem}
Since the nerve functor $N$ is fully faithful, a small category $I$
can be recovered from its nerve $NI$. However, one does not even
need the full simplicial set $NI \in \Delta^o\Sets$: if we consider
the full subcategory $\D_2 \subset \Delta$ of \eqref{D.n.eq}, and
let $\eps:\D_2 \to \Delta$ be the embedding, then $\eps^* \circ
N:\Cat \to \D_2^o\Sets$ is also fully faithful.
\end{remark}

\begin{exa}\label{BG.exa}
For any group $G$, the nerve $N(\ppt_G)$ of the groupoid $\ppt_G$ is
the stadard classifying simplicial set $BG:\Delta^o \to \Sets$ of
$G$, with terms $BG([n])=G^n$ and the universal $G$-torsor $EG \to
BG$. The simplicial set $EG$ is the nerve of the category $e(G)$
(where $G$ is considered as a set).
\end{exa}

The nerve $NI$ of a small category $I$ is weakly regular iff for any
two maps $f:i \to i'$, $f':i' \to i$ in $I$ such that $f' \circ f =
\id$, we have $i=i'$ and $f=f'=\id$. For example, this is satisfied
if $I = P$ is the category of Example~\ref{P.exa}. As this example
shows, $NI$ can be weakly regular but not regular. To insure that
$NI$ is regular, it suffices to require that for any two maps $f:i
\to i'$, $f':i' \to i$, we have $i=i'$ and $f=f'$. We will call such
a category $I$ {\em ordered}. Ordered categories do exists; in
particular, a partially ordered set $J$ is always ordered. In this
case, we have a natural identification
\begin{equation}\label{BJ.N}
\overline{\Delta}NJ \cong \Red(\Delta NJ) \cong BJ, \qquad
\langle [n],i_\idot \rangle \mapsto i_\idot([n]) \subset J,
\end{equation}
where $BJ$ is the barycentric subdivision of
Definition~\ref{pos.bary.def}. For $J = [n]$, this is
Example~\ref{Bn.exa}.

\begin{remark}\label{BJ.N.rem}
One consequence of \eqref{BJ.N} is that $BJ$ comes equipped with a
discrete fibration $BJ \to \Delta_\ff$ (if we compose it with the
embedding $\Delta_\ff \to \Pos$, we obtain the $\Pos$-augmentation
$B^\tr(J)$ that appeared in
Example~\ref{bary.bis.rel.exa}). Therefore the standard pushout
squares in $\Pos$, as well as those perfects colimits of
Subsection~\ref{pos.perf.subs} and Subsection~\ref{pos.cyl.subs}
that are preserved by $B$ --- namely, those of
Example~\ref{J.S.sq.exa}, Example~\ref{J.n.exa},
Example~\ref{pos.sq.exa}, Example~\ref{Z.NN.exa},
Example~\ref{BJ.coli.exa} and Example~\ref{N.coli.exa} --- produce
colimits in $\Delta_\ff^o\Sets$, and then by \eqref{BJ.N} and
\eqref{ff.X}, they are also preserved by the nerve functor $N:\Pos
\to \Delta^o\Sets$.
\end{remark}

\begin{lemma}\label{BB.le}
If a small ordered category $I$ has an initial object $i \in I$,
then the embedding $\ppt \to \overline{\Delta}NI$ induced by
$\eps(i):\ppt \to I$ is a reflexive map of partially ordered sets.
\end{lemma}

\proof{} Since $I$ is ordered and $i \in I$ is initial,
$\eps(i):\ppt \to I$ is a left-closed full embedding, so that we
have the characteristic functor $\chi:I \to [1]$ sending $i$ to $0
\in [1]$ and everything else to $1 \in [1]$. Then if we denote $BI =
\overline{\Delta}NI$, Example~\ref{B.cl.exa} and Lemma~\ref{B.le} apply
to $I$ with the same proofs.
\endproof
  
For any functor $\gamma:I \to I'$ between small categories, with the
corresponding map $N(\gamma):NI \to NI'$ between their nerves, the
maps \eqref{st.eq} fit into a commutative square
\begin{equation}\label{ngamma.sq}
\begin{CD}
NI([1]) @>{N(\gamma)([1])}>> NI'([1])\\
@V{NI(s) \times NI(t)}VV @VV{NI'(s) \times NI'(t)}V\\
NI([0]) \times NI([0]) @>{N(\gamma)([0]) \times N(\gamma)([0])}>>
NI'([0]) \times NI'([0]),
\end{CD}
\end{equation}
and $\gamma$ is fully faithful iff the square \eqref{ngamma.sq} is
cartesian. If $N(\gamma)$ is essentially surjective, then $\gamma$
is essentially surjective; if $\gamma$ is fully faithful, then
$N(\gamma)$ is surjective iff so is $N(\gamma)([0])$, and in this
case, $\gamma$ is an equivalence. If $\gamma$ is an isofibration in
the sense of Subsection~\ref{fib.subs}, then the converse it true:
$\gamma$ is essentially surjective iff it it surjective on the nose,
and this means that $N(\gamma)([0])$ is surjective. Moreover, say
that a map $f:X \to X'$ between simplicial sets is {\em dense} if
$f([0]):X([0]) \to X'([0])$ is an isomorphism; then $\gamma:I \to
I'$ is dense iff so is its nerve $N(\gamma)$.

For any map $F:X \to Y$ between simplicial sets, we can define the
cylinder $\Cyl(F)$ and the dual cylinder $\Cyl^o(F)$ by cocartesian
squares
\begin{equation}\label{cyl.X.sq}
\begin{CD}
X @>{\id \times t}>> X \times \Delta_1\\
@V{F}VV @VVV\\
Y @>>> \Cyl(F),
\end{CD}
\qquad\qquad
\begin{CD}
X @>{\id \times s}>> X \times \Delta_1\\
@V{F}VV @VVV\\
Y @>>> \Cyl^o(F),
\end{CD}
\end{equation}
and we can define the left and right comma-sets $X /_F Y$, $X
\setminus_F Y$ by cartesian squares
$$
\begin{CD}
X \setminus_F Y @>>> X\\
@VVV @VV{F}V\\
Y^{\Delta_1} @>{s^*}>> Y,
\end{CD}
\qquad\qquad
\begin{CD}
X /_F Y @>>> X\\
@VVV @VV{F}V\\
X^{\Delta_1} @>{t^*}>> Y,
\end{CD}
$$
where $s,t:\ppt \to \Delta_1$ are induced by $s,t:[0] \to [1]$. For
any simplicial set $X$, we also define
\begin{equation}\label{X.gt}
X^> = \Cyl(\pi), \qquad X^< = \Cyl^o(\pi),
\end{equation}
where $\pi:X \to \ppt$ is the tautological projection, and for any
two sets $X_0$, $X_1$ equipped with maps to $X$, we define
\begin{equation}\label{X.2.times}
X_0 \times^2_X X_1 = X_0 \times_X X^{\Delta_1} \times_X X_1.
\end{equation}
Then for any small category $I$, we have $(NI)^> \cong NI^>$,
$(NI)^< \cong NI^<$, and we have $N(I /_\gamma I') \cong NI
/_{N(\gamma)} NI'$, $N(I \setminus_\gamma I') \cong NI
\setminus_{N(\gamma)} NI'$ for any functor $\gamma:I \to I'$ between
small $I$, $I'$. If we have two functors $I_0 \to I$, $I_1 \to I$
and $I$ is a groupoid, we also have $N(I_0 \times_I I_1) \cong NI_0
\times^2_{NI} NI_1$. For any pair of functors $\lambda:I \to I'$,
$\rho:I' \to I$, an adjunction between $\rho$ and $\lambda$ is given
by an isomorphism
$$
NI /_{N(\lambda)} NI' \cong NI' \setminus_{N(\rho)} NI
$$
of simplicial sets that commutes with projections to $NI \times
NI'$.

\subsection{Segal condition.}\label{seg.X.subs}

It would be useful to characterize the image of the full embedding
$\Cat \subset \Delta^o\Sets$ given by the nerve functor, and it is
indeed possible to do so. An abstract simplicial set $X:\Delta^o \to
\Sets$ is the nerve of a small category if and only if it satisfies
the {\em Segal condition}: for any $n \geq l \geq 0$, the diagram of
sets
\begin{equation}\label{segal.sq}
\begin{CD}
X([n]) @>{X(s)}>> X([l])\\
@V{X(t)}VV @VV{X(t)}V\\
X([n-l]) @>{X(s)}>> X([0])
\end{CD}
\end{equation}
obtained by applying $X$ to the square \eqref{seg.del.sq} is
cartesian. Equivalently, a square \eqref{seg.del.sq} induces a
commutative square of elementary simplices, and this gives rise to
an injective map
\begin{equation}\label{b.ln.eq}
b^l_n:\Delta_l \copr_{\Delta_0} \Delta_{n-l} \to \Delta_n,
\end{equation}
where the coproduct is taken with respect to the embeddings
$s:\Delta_0 \to \Delta_{n-l}$, $t:\Delta_0 \to \Delta_l$. Then $X$
satisfies the Segal condition if and only if the functor
$\Hom(-,X):\Delta^o\Sets \to \Sets$ inverts the maps \eqref{b.ln.eq}
for all $n \geq l \geq 0$. In fact, it suffices to consider $n > l >
0$, since \eqref{b.ln.eq} is invertible for $l=0,n$.

It is also useful to have some additional equivalent reformulations
of the Segal condition. Firstly, by \eqref{P.X.J} and
Example~\ref{Bn.exa}, subsets in $\Delta_n$ correspond bijectively
to left-closed subsets in the barycentric subdivision $B([n])$. For
any $l$, $0 \leq l \leq n$, the subset $B([n])_l=\{S | S \cup \{l\}
\neq [n]\} \subset B([n])$ is left-closed, and the corresponding
subset $\V^l_n \subset \Delta_n$ is known as the {\em $l$-th horn},
with the embedding map
\begin{equation}\label{horn.eq}
v^l_n:\V^l_n \to \Delta_n.
\end{equation}
Note that the complement $\overline{B}([n])_l \subset B([n])$ to
$B([n])_l$ is right-closed and consists of the two subsets $[n],[n]
\ssetminus \{l\} \subset [n]$; it is therefore abstractly isomorphic
to $[1]$. If $l=0$ or $l=n$, we have
\begin{equation}\label{V.S.eq}
V^0_n \cong \SS_{n-1}^<, \quad V^n_n \cong \SS_{n-1}^>, \qquad
v^0_n = \sigma_n^<, \quad v^n_n = \sigma_n^>,
\end{equation}
where $\sigma_n$ is the sphere embedding \eqref{sph.eq}, and
$\SS_{n-1}^>$, $\SS_{n-1}^<$ are as in \eqref{X.gt}.

Secondly, for any map $\chi:J \to J'$ between partially ordered
sets, the comma-set $J/J'$ is represented as the colimit of
Lemma~\ref{pos.phi.le}, and this induces a map
\begin{equation}\label{a.chi.eq}
a_\chi:\colim_{\Tw(J')}N((J/J')_f) \to N(J/J')
\end{equation}
of the corresponding nerves. The map is injective but not
necessarily surjective. In particular, if $J'=[1]$, $J=[n]$ and
$\chi:[n] \to [1]$ is the characteristic function of a left-closed
embedding $s:[l] \to [n]$, \eqref{a.chi.eq} becomes the embedding
\begin{equation}\label{a.ln.eq}
a^l_n:\Cyl(N(s)) \to N(\Cyl(s))
\end{equation}
whose source is the simplicial cylinder \eqref{cyl.X.sq}. With this
notation, we have the following observation.

\begin{lemma}\label{bva.le}
Assume given a simplicial set $X$, and let $W(X)$ be the class of
maps in $\Delta^o\Sets$ that are inverted by $\Hom(-,X)$. Then the
following are equivalent: for any $n > l > 0$, \thetag{i} $b^l_n \in
W(X)$, \thetag{ii} $v^l_n \in W(X)$, \thetag{iii} $a^l_n \in
W(X)$. Moreover, these conditions are equivalent to
\begin{enumerate}
\setcounter{enumi}{3}
\item for any map $\chi:J \to J'$ between partially ordered sets of
  finite chain dimension, the map $a_\chi$ of \eqref{a.chi.eq} is
  in $W(X)$.
\end{enumerate}
\end{lemma}

\proof{} Note that by definition, the class $W(X)$ is saturated and
stable under pushouts and coproducts. To see that \thetag{ii}
implies \thetag{i}, take some integers $n > l > 0$. Then the source
of the map \eqref{b.ln.eq} corresponds to the left-closed subset
$B([n])^l \subset B([n])$ consising of subsets $S \subset B([n])$
that lie entirely in $s([l]) = \{0,\dots,l\} \subset [n]$ or in
$t([n-l]) = \{l,\dots,n\} \subset [n]$. For any $d \geq k \geq 0$,
let $A([n])^k_d \subset B([n])$ be the subset of $S \subset [n]$
such that $|S \cap \{0,\dots,l-1\}| \leq k$ and $|S \cap ([n]
\ssetminus \{l\})| \leq d$. Moreover, order pairs $\langle d,k
\rangle$ lexicographically --- that is, $\langle d,k \rangle \geq
\langle d',k' \rangle$ iff $d > d'$ or $d=d'$ and $k \geq k'$ ---
and let
\begin{equation}\label{B.km}
B([n])^k_d = B([n])^l \cup \bigcup_{\langle d',k' \rangle \leq
  \langle d,k \rangle} A([n])^{k'}_{d'}
\end{equation}
for any $d \geq k \geq 0$. Then $A([n])^0_0$ is empty, so that
$B([n])^0_0 = B([n])^l$, and $A([n])^l_n = B([n])$, so that we
obtain an increasing filtration
\begin{equation}\label{B.l.filt}
B([n])^l = B([n])^0_0 \subset \dots \subset B([n])^l_n =
B([n])
\end{equation}
where all the embeddings are moreover left-closed. For any $d \geq k
\geq 0$, the complement $\overline{B}([n])^k_d \subset B([n])$ to
the previous term in the filtration \eqref{B.l.filt} is right-closed
and consists of subsets $S \subset [n]$ whose intersection with
$\{0,\dots,l-1\}$ resp.\ $\{l+1,\dots,n\}$ contains exactly $k$
resp.\ $d-k$ elements. Therefore we have an isomorphism
\begin{equation}\label{B.S.km}
\overline{B}([n])^k_d \cong S(d,k) \times [1],
\end{equation}
where $S(d,k)$ is a discrete set whose elements are pairs of a
subset in $\{0,\dots,l-1\}$ of cardinality $k$ and a subset in
$\{l+1,\dots,n\}$ of cardinality $d-k$, and $S \in
\overline{B}([n])^k_d$ lies in $S(d,k) \times \{1\}$ or $S(d,k)
\times \{0\}$ depending on whether it contains $l$ or not. Moreover,
any $S \in \overline{B}([n])^k_d$ that does contain $l$ corresponds
to a map $f_S:[d] \to [n]$ such that $B(f_S)$ sends the horn subset
$B([d])_k \subset B([d])$ into the previous terms of the filtration
\eqref{B.l.filt}, so that if we consider the corresponding
filtration on $\Delta_n$, each embedding is a pushout of a finite
coproduct of maps $v^k_d$ of \eqref{horn.eq}. Therefore they are all
in $W(X)$, and so is their composition \eqref{b.ln.eq}.

Moreover, one observes that the last embedding in \eqref{B.l.filt}
is actually the horn embedding $B([n])_l \subset B([n])$ on the
nose, and all the previous ones are pushouts of horn emebdding
$v^k_d$ with $\langle d,k \rangle < \langle n,l \rangle$. Since
$W(X)$ is saturated, this means that \thetag{i} implies \thetag{ii}
by induction on $\langle n,l \rangle$.

Since $a^0_n = b^1_{n+1}$, \thetag{iii} implies \thetag{i} for
$l=1$, and then for any $n \geq l > 1$, we have a commutative square
$$
\begin{CD}
\Delta_1 \copr_{\Delta_0} \Delta_l \copr_{\Delta_0} \Delta_{n-l}
@>{\id \copr v^l_n}>> \Delta_1 \copr_{\Delta_0} \Delta_n\\
@V{v^1_l \copr \id}VV @VV{v^1_{n+1}}V\\
\Delta_{1+l} \copr_{\Delta_0} \Delta_{n-l}
@>{v^l_{n+1}}>> \Delta_{n+1},
\end{CD}
$$
so that \thetag{i} for $l=1$ implies \thetag{i} for any $l$ by
induction. Moreover, \thetag{iv} trivially implies \thetag{iii}, so
to finish the proof, it remains to prove that \thetag{ii} implies
\thetag{iv}.

To do this, take a map $\chi:J \to J'$ in $\Posf$, consider the
barycentric subdivision $B(J/J')$, and let $B(J/J')' =
\colim_{\Tw(J')}B(F(\phi)) \subset B(J/J')$. For any element in the
complement $B(J/J') \ssetminus B(J/J')'$ represented by an injective
map $i:[n] \to J/J'$, the left-closed subset $P(i) \subset [n]$ of
\eqref{K.ka} is non-empty by Remark~\ref{K.rem}, so that $P(i) =
s([p_i])$ for some $p_i \geq 0$. The composition $\sigma \circ i:[n]
\to J$ of the map $i$ with $\sigma:J/J' \to J$ need not be
injective, so the left-closed subset $K(i) = s([k_i]) \subset [n]$
given by $K(i) = \{l \in [n]|\sigma(i(l)) \leq \sigma(i(p_i))\}$ is
possibly larger. However, by \eqref{K.ka}, we still have
$\chi(\sigma(i(p_i))) \leq \tau(i(p_i)) < \chi(\sigma(i(n)))$, thus
$\sigma(i(p_i)) < \sigma(i(n))$, and since $\sigma(i(p_i)) =
\sigma(i(k_i))$, the largest element $n \in [n]$ is not in $K(i)$,
and $p_i \leq k_i < n$. Let $l_i=1$ if $\tau(i(k_i)) <
\tau(i(k_i+1))$ and $0$ otherwise, and note that since
$\tau(i(k_i+1)) \geq \tau(i(p_i+1)) \geq \chi(i(n)) > \tau(i(p_i))$,
the latter is only possible when $k_i > p_i$, so that in any case,
$k_i + l_i > 0$. Now for any $d \geq k \geq 0$, let $B(J/J')^k_d
\subset B(J/J')$ be the left-closed subset obtained by taking
$B(J/J')'$ and adding all elements $i$ such that $k_i +l_i \leq k$ and $n
+ l_i \leq d$. Then again, if we order the pairs $\langle d,k
\rangle$ lexicographically, we obtain a filtration \eqref{B.l.filt}
on $B(J/J')$, finite since $J/J' \in \Posf$, and the right-closed
complement $\overline{B}(J/J')^k_d \subset B(J/J')^k_d$ consists of
maps $i$ such that $k_i+l_i=k$ and $n+l_i=d$. Such maps can only
exist if $d > k > 0$, and the only possible order relation $i' \leq
i$ between two of them occurs when $l_{i'}=1$, $l_i=0$, and
$i'([d-1]) \subset B(J/J')$ is obtained by removing the element
$i(k) \in i([d])$. This can be done for any $i$ with $l_i=0$, and
conversely, the removed element $i(k) \in J/J'$ is represented by
the arrow $\sigma(i(k-1)) \leq \tau(i(k+1))$, so it only depends on
$i'$, and $i$ can be uniquely reconstructed from $i'$. We conclude
that again, $\overline{B}(J/J')^k_d$ is of the form \eqref{B.S.km}
for some discrete sets $S(d,k)$, and when we consider the filtration
on $N(J/J')$ induced by \eqref{B.l.filt}, each successive map is
again a pushout of a coproduct of horn embeddings \eqref{horn.eq}
(specifically, $N(J/J')^k_d$ is obtained as a pushout of several
copies of $v^k_d$, $d > k > 0$).  \endproof

\subsection{Nerves of groupoids.}

One can complement Lemma~\ref{bva.le} by a similar characterization
of nerves of small groupoids. Denote by $\Gamma$ the category of
finite non-empty sets $\{0,\dots,n\}$, and denote by $\gamma:\Delta
\to \Gamma$ the tautological functor sending $[n]$ to itself without
any order.

\begin{lemma}\label{gamma.le}
Assume given a small category $I$ with the nerve $NI$, and let
$W(I)$ be the class of injective maps $f:X_0 \to X_1$ of simplicial
sets such that $\Hom(-,NI)$ inverts the corresponding map
$f^>:X_0^> \to X_1^>$ of the simplicial sets \eqref{X.gt}. Then the
following conditions are equivalent.
\begin{enumerate}
\item $I$ is a groupoid.
\item $W(I)$ contains the sphere embedding $\sigma_1:\SS_0 \to \Delta_1$.
\item $W(I)$ contains all dense injective maps $X' \to X$ between
  finite simplicial sets.
\item There exists a functor $X:\Gamma^o \to \Sets$ such that
  $\gamma^*(X) \cong NI$.
\end{enumerate}
\end{lemma}

\proof{} By \eqref{V.S.eq}, \thetag{ii} amounts to saying that
$\Hom(-,NI)$ inverts $v^2_2$, or explicitly, that for any maps
$f_0:i_0 \to i$, $f_1:i_1 \to i$ in $I$, there exists a unique map
$f:i_0 \to i_1$ such that $f_0 = f_1 \circ f$. This is equivalent to
\thetag{i}.

To see that \thetag{ii} implies \thetag{iii}, note the functor $X
\mapsto X^>$ preserves colimits, so that $W(I)$ is saturated and
closed with respect to pushouts. Then by induction on the
cardinality $|\Red(\Delta X)|$, to check \thetag{iii}, it suffices
to consider elemenary dense embeddings $X' \to X$, and by
\eqref{elem.X.sq}, it further suffices to consider sphere embeddings
\eqref{sph.eq}. Then again by induction, it suffices to consider the
embeddings $\Delta_n([0]) \to \Delta_n$, and by the Segal condition
for $N(I \setminus I_{\Id})$, this reduces to the case $n=1$, that
is, \thetag{ii}.

Now, \thetag{iii} trivially implies \thetag{ii}, hence also
\thetag{i}, and to see that these equivalent conditions imply
\thetag{iv}, note that by the version of \eqref{kan.eq} for right
Kan extensions, we have $\gamma^*\gamma_*NI([n]) \cong
\Hom(Ne(\gamma([n]),NI)$ for any $[n] \in \Delta$, while the
adjunction map $NI \to \gamma^*\gamma_*NI$ evaluated at $[n]$ is
induced by the natural map $a_n:[n] \to e(\gamma([n]))$. If $I$ is a
groupoid, then $\Fun(-,I)$ inverts all the maps $a_n$, so that $NI
\to \gamma^*\gamma_*NI$ is an isomorphism.

Conversely, to see that \thetag{iv} implies \thetag{iii}, it
suffices to check that $f^>$ for any dense embedding $f:X' \to X$ of
finite simplicial sets is inverted by $\gamma_!$. By \eqref{kan.eq},
this amounts to checking that it is inverted by $\Hom(-,Ne(S))$ for
any $S \in \Gamma$, and $e(S)$ is a groupoid, thus satisfies
\thetag{iii}.
\endproof

\begin{corr}\label{edge.corr}
For any small groupoid $I$ with the nerve $X = NI$, the edgewise
subdivision map \eqref{edge.map} is an isomorphism.
\end{corr}

\proof{} For any $[m],[n] \in \Delta$, the morphism $[m] \times [n]
\to [m] * [n]$ is a dense map of finite partially ordered sets that
sends the largest element to the largest element; therefore by
Lemma~\ref{gamma.le}~\thetag{iii}, $\Hom(-,X)$ inverts the
corresponding map $N([m] \times [n]) \to N([m] * [n])$.
\endproof

\begin{lemma}\label{eck.le}
Assume given a small category $I$ with the nerve $NI$, and assume
that $NI$ admits a structure of a simplicial group. Then $I$ is a
groupoid.
\end{lemma}

\proof{} This is a version of the standard Eckmann-Hilton argument
claiming that a group object in the category of groups is an abelian
group, and both group structures coincide. The assumption means that
$I$ is a group object in the category of small categories. In
particular, its set of objects is a group, with some unit element $e
\in I$, and $I(e,e)$ is a monoid object in the category of groups,
with the group unit $e' \in I(e,e)$. Then
$$
e' = e' \cdot e' = (e' \circ \id) \cdot (\id \circ e') = (e' \cdot
\id) \circ (\id \cdot e') = \id \circ \id = \id
$$
and $f \cdot g = (f \circ \id) \cdot (\id \circ g) = (f \circ e')
\cdot (e' \circ g) = (f \cdot e') \circ (e' \cdot g) = f \circ g$
for any $f,g \in I(e,e)$, so that the monoid $I(e,e)$ is a
group. For any $i \in I$, group multiplication with $\id_i:i \to i$
identifies the monoids $I(e,e)$ and $I(i,i)$, and then for any map
$f:i \to i'$, we have a map $f' = \id_{i'} \cdot f^{-1} \cdot
\id_i:i' \to i$ with invertible compositions $f' \circ f \in
I(i,i)$, $f' \circ f \in I(i',i')$.  \endproof

\section{Segal categories.}\label{seg.sec}

\subsection{Simplicial replacements.}\label{del.subs}

To make sense of nerves of categories that are not necessarily
small, it is convenient to use the augmented categories of
Subsection~\ref{fun.subs}.

\begin{defn}\label{simpl.repl.def}
Assume given a category $I$ and a closed class $v$ of morphisms in
$I$. The {\em simplicial expansion} $\Delta^v I$ of the category $I$
with respect to the class $v$ is given by
$$
\Delta^v I = (\Delta \bb^\phi I)_{\phi^*\Cat(v)},
$$
where $\phi:\Delta \to \Cat$ is the tautological embedding, and
$\bb^\phi$ is as in \eqref{C.bb}. The {\em simplicial replacement}
$\Delta I = \Delta^{\Id} I$ of a category $I$ is its simplicial
expansion with respect to the class $\Id$ of all identity maps.
\end{defn}

The notion of a simplicial replacement goes back to \cite{bou3}.
Explicitly, for any category $I$, an object $\langle [n],i_\idot
\rangle$ in $\Delta I$ is given by a chain \eqref{c.idot} of length
$n$ in $I$ --- that is, a diagram
\begin{equation}\label{i.dot}
\begin{CD}
i_0 @>>> \dots @>>> i_n.
\end{CD}
\end{equation}
Objects in a simplicial expansion $\Delta^v I$ with respect to some
class $v$ are the same, but there are more morphisms between
them. The forgetful functor $\pi:\Delta^v I \to \Delta$, $\langle
[n],i_\idot \rangle \mapsto [n]$ is a fibration; if $v=\Id$, its
fibers are discrete. The transpose cofibration
\begin{equation}\label{trans.eq}
\Delta^v_\perp I = (\Delta^v I)_\perp \to \Delta^o
\end{equation}
is called the {\em cosimplicial expansion} of the category $I$ with
respect to the class $v$. Note that we have a tautological equivalence
\begin{equation}\label{d.io}
\iota^*(\Delta^v_\perp I)^o \cong \Delta^{v^o} I^o,
\end{equation}
where $\iota$ is the involution \eqref{iota.eq}, and $v^o$ is
opposite to $v$. If $v = \id$, this reads as $\Delta_\perp I \cong
(\Delta I)^o \cong \iota^*\Delta I^o$. If $I$ is small, then $\Delta
I = \Delta NI$ and $\pi$ corresponds to $NI \in \Delta^o\Sets$ by
the Grothendieck construction --- so that in particular, $\Delta [n]
\cong \Delta/[n]$ for any $[n] \in \Delta$ --- but $\Delta I$ is
perfectly well defined for an arbitrary $I$. If $v = \Tot$ is the
class of all maps in $I$, then by \eqref{cat.wr}, we have
\begin{equation}\label{del.wr}
\Delta^{\Tot} I \cong \Fun(\Delta_\idot | \Delta,I),
\end{equation}
where as in Example~\ref{cat.exa} and \eqref{nu.pos},
\begin{equation}\label{nu.eq}
\nu_\idot:\Delta_\idot = \phi^*\Cat_\idot \to \Delta
\end{equation}
is the tautological cofibration whose fiber over some $[n] \in
\Delta$ is $[n]$ itself (alternatively, $\Delta_\idot$ is the
category of pairs $\langle [n],l \rangle$, $[n] \in \Delta$, $l \in
[n]$, with maps $\langle [n],l \rangle \to \langle [n'],l' \rangle$
given by maps $f:[n] \to [n']$ such that $f(l) \leq l'$). By
\eqref{cat.I.sq}, the expansion $\Delta^v I$ for some other closed
class $v$ fits into a cartesian square
\begin{equation}\label{del.I.sq}
\begin{CD}
\Delta^v I @>>> \Delta^{\Tot} I\\
@VVV @VVV\\
\eps_*I_v @>>> \eps_*I,
\end{CD}
\end{equation}
where $\eps = \eps([0]):\ppt \to \Delta$ is the embedding
\eqref{eps.c} onto $[0] \in \Delta$.

The projection \eqref{nu.eq} has a fully faithful right-adjoint
$\nu^\dg:\Delta \to \Delta^\hdot$ sending $[n]$ to $\langle [n],n
\rangle$, and this in turn has a right-adjoint $\nu_<:\Delta^\hdot
\to \Delta$ sending $\langle [n],l \rangle$ to $[l]$ (this is the
restriction of \eqref{nu.io.pos} to $\Delta \subset \Pos$). Then the
relative evaluation functor \eqref{ev.I.rel.eq} for \eqref{del.wr}
induces a functor
\begin{equation}\label{xi.eq}
\xi = \ev \circ (\nu^\dg \times \id):\Delta^v I \to I,
\end{equation}
and conversely, we have $\ev \cong \xi \circ (\nu_< \times
\id)$, so that the functor $\ev$ can be recovered from the functor
$\xi$. Explicitly, $\xi$ sends a diagram \eqref{i.dot} to $i_n \in
I$. The comma-fibers of the functor \eqref{xi.eq} are given by
\begin{equation}\label{xi.comma}
(\Delta I) /_\xi i \cong \Delta (I/i), \qquad i \in I,
\end{equation}
so that if $I$ is locally small, \eqref{xi.eq} is bounded. We also
have its transpose version
\begin{equation}\label{xi.perp.eq}
\xi_\perp:\Delta^v_\perp I \to I
\end{equation}
sending \eqref{i.dot} to $i_0$.  If $v=\Tot$ is the class of all maps
in $I$, then $I \cong (\Delta^{\Tot}I)_{[0]}$, and the functor
\eqref{coev.eq} of Example~\ref{coev.exa} provides a section
\begin{equation}\label{del.coev}
r:\Delta \times I \to \Delta^{\Tot} I
\end{equation}
of the projection $\pi \times \xi:\Delta^{\Tot} I \to \Delta \times
I$. If $I=J$ is a partially ordered set and $v=\Id$, then we have
the identification \eqref{BJ.N}, and the maps \eqref{B.xi} are then
identified with the reductions of the functors \eqref{xi.eq} and
\eqref{xi.perp.eq}.

It is also useful to consider the case $I=\Delta$. Then by
\eqref{del.wr}, $\nu_<$ induces a functor
\begin{equation}\label{al.0.eq}
\delta:\Delta \to \Fun(\Delta^\hdot | \Delta,\Delta) \cong
\Delta^{\Tot}\Delta,
\end{equation}
and since $\nu^\dg$ is fully faithful, we have the isomorphisms
$\nu_\idot \circ \nu^\dg \cong \nu_< \circ \nu^\dg \cong \id$
that induce isomorphisms $\pi \circ \delta \cong \xi \circ \delta
\cong \id$, so that $\delta$ is a common section of the structural
fibration $\pi:\Delta^{\Tot}\Delta \to \Delta$ and the functor
\eqref{xi.eq}.

\begin{exa}\label{del.1.a.exa}
For an application of simplicial expansions with respect to a
non-trivial closed class of maps, assume given a functor $\gamma:I_1
\to I_0$, with the dual cylinder $I = \Cyl^o(\gamma)$ and the
fibration $\chi:I \to [1]$. Then we have a canonical cleavage of the
fibration $\chi$, that is, a choice of a cartesian lifting
$\gamma(i) \to i$ of the unique map $0 \to 1$ in $[1]$ for any $i
\in I_0 \subset I$. Adding the identity maps to these unique
cartesian liftings defines a closed class $c$ of maps in $I$, and we
have the projection
\begin{equation}\label{cyl.del}
\Delta(\chi):\Delta^c I \to \Delta^{\Tot}[1],
\end{equation}
where as usual, $\Tot$ is the class of all maps in $[1]$. One then
observes that \eqref{cyl.del} is a fibration. Conversely, since
$\Delta [1] \to \Delta$ is discrete, the projection
$\Delta(\chi):\Delta I \to \Delta [1]$ is a fibration for any
category $I$ and functor $\chi:I \to [1]$, and $\chi$ is a fibration
if and only if $\Delta I \to \Delta [1]$ extends to a fibration $\E
\to \Delta^{\Tot} [1]$. In this case we of course have $I \cong
\Cyl^o(\gamma)$, where $\gamma:I_1 \to I_0$ is the transition
functor of the fibration. To recover $\gamma$, one can use the
embedding \eqref{del.coev}; then $r^*\E \to \Delta \times [1]
\to [1]$ is a fibration with fibers $\Delta I_1$, $\Delta I_0$, and
its transition functor is $\Delta(\gamma)$.
\end{exa}

In general, it is not so easy to characterize fibrations over
arbitrary bases $I$ in terms of simplicial replacements, but there
is an easy way to do it for families of groupoids. Say that a family
of groupoids $\pi:\C \to I$ is {\em classical} if it is a classical
fibration in the sense of Subsection~\ref{fib.subs} (recall that it
suffices to assume that $\pi$ is an isofibration to insure this
property, and that for any family of groupoids $\C \to I$, the
equivalent family $I \times_I \C \to I$ is already classical).

\begin{lemma}\label{grp.v.le}
A functor $\gamma:\C \to I$ between two categories $\C$, $I$ is a
classical family of groupoids if and only if the corresponding
functor \eqref{pi.co.eq} is fully faithful, and its simplicial
replacement $\Delta(\pi \setminus \id)$ is essentially
surjective.
\end{lemma}

\proof{} The ``if'' part is clear: if $\Delta(\pi \setminus \id)$
is essentially surjective, $\pi \setminus \id$ itself is essentially
surjective, and then we are done by
Lemma~\ref{grp.c.le}. Conversely, if $\pi$ is a family of groupoids,
then $\pi \setminus \id$ is fully faithful by Lemma~\ref{grp.c.le},
and by the very definition of a classical fibration, $\Delta(\pi
\setminus \id)$ is essentially surjective over $[0] \in
\Delta$. Then by the Segal property of simplicial replacements, it
is essentially surjective over the whole $\Delta$.
\endproof

We note that in the situation of Lemma~\ref{grp.v.le}, the
simplicial replacements $\Delta(\C \setminus_{\id} \C)$, $\Delta(I
\setminus_\pi \C)$, $\Delta(\pi \setminus \id)$ are easily
constructed directly from the simplicial replacements $\Delta \C$,
$\Delta I$, $\Delta(\pi):\Delta \C \to \Delta I$, by the same
procedures as in Subsection~\ref{del.0.subs}. Moreover, the
condition that \eqref{pi.co.eq} is fully faithful can again be
checked purely in terms of its simplicial replacement. Thus
Lemma~\ref{grp.v.le} indeed gives an explicit description of
classical families of groupoids in terms of their simplicial
replacements. An analogous simplification occurs for essentially
surjective functors. In general, a functor $\gamma:\C \to \C'$ can
be essentially surjective even if its simplicial replacement
$\Delta(\gamma)$ is not; however if $\C$ and $\C'$ are groupoids,
then $\gamma$ is essentially surjective if and only if so so is the
projection $\Delta(\sigma):\Delta (\C' \setminus_\gamma \C) \to
\Delta \C'$. The same holds if $\C$, $\C'$ are families of groupoids
over some $I$, and $\gamma$ is a functor over $I$. In the latter
case, one can also replace $\C' \setminus_\gamma \C$ with its
relative version $\C' \setminus_\gamma^I \C$ of \eqref{rel.lc.sq}.

Alternatively, one can also characterize simplicial replacements of
classical families of groupoids in terms of horn embeddings
\eqref{horn.eq}. Namely, say that a map $f:X \to X'$ of simplicial
sets is {\em liftable} with respect to a functor $\pi:\C \to I$ if
for any commutative diagram
\begin{equation}\label{grp.lft.sq}
\begin{CD}
\Delta X @>{g}>> \Delta \C\\
@V{\Delta(f)}VV @VV{\Delta(\pi)}V\\
\Delta X' @>{g'}>> \Delta I
\end{CD}
\end{equation}
of categories and functors over $\Delta$, there exists a functor
$q:\Delta X' \to \Delta \C$ over $\Delta$ such that $g \cong q \circ
\Delta(f)$ and $\Delta(\pi) \circ q \cong g'$, and say that a
liftable $f$ is {\em uniquely liftable} if such a functor $q$ is
unique.

\begin{lemma}\label{horn.grp.le}
A functor $\pi:\C \to I$ is a classical family of groupoids if and
only if the horn embedding $v^n_n$ is liftable with respect to $\pi$
for $n=1,2$ and uniquely liftable for $n=2$. Moreover, in this case,
for any dense embedding $f:X_0 \to X_1$ of simplicial sets, the maps
$f_\dg = (\id \times f) \copr \id:\Cyl(f) \to \Delta_1 \times X_1$
and $f^>:X_0^> \to X_1^>$ are uniquely liftable with respect to
$\pi$.
\end{lemma}

\proof{} For the first claim, note that liftability for $v^1_1$
means that $\pi \setminus \id$ is essentially surjective, so that
for any $c \in \C$ and $f:i \to \pi(c)$, there exists a lifting
$f':c' \to c$, $\pi(f')=f$, and then unique liftability for $v^2_2$
means that every map in $\C$ is cartesian with respect to $\pi$. For
the second claim, note that the class of maps $f:X \to X'$ such that
$f^>$ and $f^\dg$ are uniquely liftable with respect to $\pi$ is saturated and
closed under pushouts, and use the same induction as in
Lemma~\ref{gamma.le} to reduce everything to the case $X_0=\SS_0$, $X_1
= \Delta_1$, $f$ the sphere embedding \eqref{sph.eq}. Then
$f^>=v^2_2$, and $f_\dg$ is a composition of pushouts of $v^2_2$
and $v^1_2$.
\endproof

\subsection{Special maps.}\label{spec.subs}

Another way to produce the cofibration \eqref{nu.eq} is to use a
factorization system on $\Delta$, as in
Lemma~\ref{facto.fib.le}~\thetag{ii}. Namely, say that a map $f$ in
$\Delta$ is {\em special} if it is left-reflexive in the sense of
Definition~\ref{po.refl.def}, and denote the class of special maps
by $+$, with the corresponding dense subcategory $\Delta_+ \subset
\Delta$. Note that if $f$ is left-reflexive, then the adjoint map
$f_\dg$ is tautologically right-reflexive, so that we have a natural
equivalence
\begin{equation}\label{jo.eq}
\Delta_+ \cong \Delta_+^o, \qquad [n] \mapsto [n]^o, f \mapsto f_\dg^o
\end{equation}
(this is a part of the so-called {\em Joyal duality}). Explicitly,
by Example~\ref{n.refl.exa}, $f:[n] \to [n']$ is special if and only
if $f(n) = n'$. Then the class $+$ fits into a factorization system
$\langle +,s \rangle$ on $\Delta$, with $s$ consisting of the
left-closed embeddings $s:[n] \to [n']$, $n \leq n'$ of
Example~\ref{st.cl.exa}. The corresponding full subcategory
$\Ar^s(\Delta) \subset \Ar(\Delta)$ is then naturally identified
with $\Delta_\idot$, and the functor $\nu^\dg \cong \eta:\Delta \to
\Ar^s(\Delta) \cong \Delta_\idot$ and its adjoints $\nu_< \cong
\sigma,\nu_\idot \cong \tau:\Ar^s(\Delta) \to \Delta$ correspond to
the functors $\sigma$, $\tau$, $\eta$ of Example~\ref{adj.exa}.

Among other things, this interpretation of the cofibration
$\nu_\idot$ shows that the functor $\delta$ of \eqref{al.0.eq} is
actually a functor
\begin{equation}\label{al.eq}
\delta:\Delta \to \Delta^+\Delta \subset \Delta^{\Tot}\Delta
\end{equation}
from $\Delta$ to its simplicial expansion $\Delta^+\Delta$ with
respect to the class of special maps. Another observation is that
for any category $I$, Example~\ref{ind.exa} provides an
identification
\begin{equation}\label{I.wr}
\Delta^{\Tot}I \cong \Ind^+(\Delta \times I),
\end{equation}
where $\Ind^+$ is the induction of Subsection~\ref{ind.subs} with
respect to the class of special maps, while the functor $\alpha$ of
Corollary~\ref{ind.corr} is given by
\begin{equation}\label{al.xi.eq}
\alpha = \pi \times \xi:\Ind^+(\Delta \times I) \cong \Delta^{\Tot}I
\to \Delta \times I,
\end{equation}
where $\pi:\Delta^{\Tot}I \to \Delta$ is the structural fibration,
and $\xi$ is the functor \eqref{xi.eq}. Since $\Delta \times I \to
\Delta$ is already a fibration, $\alpha$ admits a fully faithful
right-adjoint functor
\begin{equation}\label{be.eq}
\beta:\Delta \times I \to \Delta^{\Tot}I
\end{equation}
of Example~\ref{ab.exa} that coincides with $r$ of
\eqref{del.coev}. Explicitly, $\beta$ sends some $[n] \times i \in
\Delta \times I$ to the diagram \eqref{i.dot} with $i_0 = \dots =
i_n = i$, and all the maps equal to $\id_i$. Using the functor
\eqref{be.eq}, one can construct the following relative version of
simplicial expansions.

\begin{defn}\label{rel.simp.def}
The {\em relative simplicial expansion} $\Delta^v(\C | I)$ of a
fibration $\pi:\C \to I$ with respect to a closed class $v$ of maps
in $\C$ is given by the cartesian square
\begin{equation}\label{d.c.i}
\begin{CD}
\Delta^v(\C | I) @>>> \Delta^v \C\\
@VVV @VV{\Delta(\pi)}V\\
\Delta \times I @>{\beta}>> \Delta^{\Tot}I,
\end{CD}
\end{equation}
where $\beta$ is the functor \eqref{be.eq}.
\end{defn}

Explicitly, $\Delta^v(\C | I) \subset \Delta^v(\C)$ is the full
subcategory spanned by functors $c_\idot:[n] \to \C$ that factor
through $\C_{\pi^*\Id}$ (or in other words, functors such that $\pi
\circ c_\idot$ factors through the projection $[n] \to [0]$). As
soon as $v$ contains all maps in $\C$ cartesian over $I$, the
restriction $\C_v \to I$ of the fibration $\pi$ is also a fibration,
with the same transition functors $f^*$, and $\Delta^v(\C | I) \to I$ is
then a fibration with fibers $\Delta^v(\C_i)$, $i \in I$, and
transition functors $\Delta(f^*)$.

\medskip

Now assume given a category $I$ with a dense subcategory $I_v
\subset I$ defined by some class $v$. Say that a map $f$ in the
simplicial expansion $\Delta^vI$ is {\em special} if it is cartesian
with respect to the projection $\pi:\Delta^vI \to \Delta$, and
$\pi(f)$ is special in $\Delta$. Say that a fibration $\C \to
\Delta^vI$ is {\em special} if it is constant along all special
maps. Note that by definition, the functor \eqref{xi.eq} inverts all
special maps in $\Delta^vI$. If $v = \Id$, so that $\Delta^vI=\Delta
I$, then the class $\pi^\flat(+) = \pi^*(+)$ of special maps fits
into a factorization system $\langle \pi^*(+),\pi^\flat(s) \rangle$
of Lemma~\ref{facto.fib.le}~\thetag{i} on $\Delta I$ induced by the
factorization system $\langle +,s \rangle$ on $\Delta$.

\begin{remark}
If $J$ is a partially ordered set, then the factorization system
$\langle \pi^*(+),\pi^*(s) \rangle$ on $\Delta J$ restricts to the
biorder of Example~\ref{bary.rel.exa} on
$BJ=\overline{\Delta}J\subset \Delta J$.
\end{remark}

\begin{prop}\label{spec.prop}
Any special fibration $\C$ over $\Delta^v I$ is bounded with respect
to the functor $\xi$ of \eqref{xi.eq}, and the functor $\xi^*\xi_*\C
\to \C$ of \eqref{2.adj} is an equivalence. Conversely, for any
fibration $\C$ over $I$, the pullback $\xi^*\C$ is special, thus
bounded with respect to $\xi$, and $\C \to \xi_*\xi^*\C$ is an
equivalence.
\end{prop}

In order to prove this, we need some combinatorial
preliminaries. Denote by $\rho:\Delta_+ \to \Delta$ the embedding,
and consider the fibration $\rho^*\Delta^vI \to \Delta_+$. Since
$[0] \in \Delta_+$ is initial, the embedding $I_v \cong
(\Delta^vI)_{[0]} \to [0] \setminus \rho^*\Delta^vI \cong
\rho^*\Delta^vI$ admits right-adjoint functor $\rho^*\Delta^vI \to
I_v$. Let $\Delta_+^vI = \rho^*\Delta^vI \times_{I_v} I_{\Id}$, with
the embedding $\rho_I:\Delta_+^v I \to \rho^*\Delta^vI \to \Delta^v
I$. Then we also have the projection $\Delta^v_+I \to I_{\Id}$, or
equivalently, a decomposition
\begin{equation}\label{del.eq}
\Delta_+^v I \cong \coprod_{i \in I} (\Delta_+^v I)_i,
\end{equation}
where $(\Delta_+^v I)_i$ is the category of objects $\langle
[n],i_\idot \rangle \in \Delta^v I$ equipped with a special map
$\langle [0],i \rangle \to \langle [n],i_\idot \rangle$ (up to an
equivalence, $(\Delta^v_+ I)_i$ is the subcategory in $\Delta^v_+I$
spanned by diagrams \eqref{i.dot} with $i_n = i$ and maps between
them whose $i_n$-component is the identity map). The embedding
$\rho_I$ restricts to an embedding $\rho_i:(\Delta^v_+ I)_i \to
\Delta^v I$. Moreover, the functor $\rho:\Delta_+ \to \Delta$ has a
left-adjoint $\lambda:\Delta \to \Delta_+$, $[n] \mapsto [n]^>$, and
for any $i \in I$, we can consider the category $(\Delta^v_\dg I)_i =
\lambda^*(\Delta^v_+I)_i$ with the embedding $\lambda_i:(\Delta^v_\dg
I)_i \to (\Delta^v_+ I)_i$. The adjunction map $a:\id \to \rho \circ
\lambda$ then induces a functor $\nu_i = a^*:(\Delta^v_\dg I)_i \to
\Delta^v I$; explicitly, $\nu$ is the forgetful functor that takes a
diagram \eqref{i.dot} and removes the last term $i_n$. The functor
$\lambda_i$ has a right-adjoint $\overline{\rho}_i:(\Delta^v_+ I)_i
\to (\Delta^v_\dg I)_i$, and we have $\nu_i \circ \overline{\rho}_i
\cong \rho_i$.  On the other hand, the adjunction map $a:\id \to
\lambda \circ \rho$ lifts to a natural map $\alpha_i:\nu_i \to
\rho_i \circ \lambda_i$.

\begin{lemma}\label{del.le}
For any $i \in I$ and any special fibration $\C$ over $\Delta I$,
the fibrations $\nu_i^*\C$ resp.\ $\rho_i^*\C$ are bounded over
$(\Delta^v_\dg I)_i$ resp.\ $(\Delta^v_+I)_i$, and the functors
$$
\begin{aligned}
\alpha_i^* \circ \lambda_i^*:&\Sec^{\Tot}((\Delta^v_+I)_i,\rho_i^*\C) \to
\Sec^{\Tot}((\Delta^v_\dg I)_i,\nu_i^*\C),\\
\overline{\rho}_i^*:&\Sec^{\Tot}((\Delta^v_\dg I)_i,\nu_i^*\C) \to
\Sec^{\Tot}((\Delta^v_+I)_i,\overline{\rho}_i^*\nu_i^*\C) \cong
\Sec^{\Tot}((\Delta^v_+I)_i,\rho_i^*\C)
\end{aligned}
$$
are mutually inverse equivalences of categories.
\end{lemma}

\proof{} The functor $\rho_i$ sends all maps to special maps, so
that for any special fibration $\C \to \Delta^vI$, the fibration
$\rho_i^*\C \to (\Delta^v_+I)_i$ is constant, and the functor
\eqref{coev.eq} provides an equivalence $\rho_i^*\C \cong
(\Delta^v_+I)_i \times \C_{\langle [0],i \rangle}$. Therefore
$\rho_i^*\C$ is bounded, and we have
\begin{equation}\label{del.i.1}
\Sec^{\Tot}((\Delta^v_+ I)_i,\rho_i^*\C) \cong \C_{\langle [0],i \rangle},
\end{equation}
while $\overline{\rho}_i^*\alpha_i^*\lambda_i^* \cong \id$ on
$\Sec^{\Tot}((\Delta^v_+I)_i,\rho_i^*\C)$. Conversely, for any
given section $\sigma:(\Delta^v_\dg I)_i \to \nu_i^*\C$, the map $\alpha_i$
induces a natural map $\sigma \to
\alpha_i^*\lambda_i^*\overline{\rho}_i^*\sigma$, and if $\sigma$ is
cartesian, this map is an isomorphism. Therefore $\nu_i^*\C$ is
bounded, and we have $\alpha_i^*\lambda_i^*\overline{\rho}_i^* \cong
\id$ on $\Sec^{\Tot}((\Delta^v_\dg I)_i,\nu_i^*\C)$.
\endproof

\proof[Proof of Proposition~\ref{spec.prop}.]  For any $i \in I$,
the left comma-fiber $\Delta^vI/i$ of the functor \eqref{xi.eq} is
identified with the category $(\Delta^v_\dg I)_i$ of
Lemma~\ref{del.le}, and $\nu_i$ corresponds to the forgetful functor
$\sigma(i)$ of \eqref{p.i}. Therefore any special fibration $\C$
over $\Delta^v I$ is bounded with respect to $\xi$. Then the functor
$\xi^*\xi_*\C \to \C$ of \eqref{2.adj} is an equivalence by
Lemma~\ref{del.le} and \eqref{del.i.1}. Conversely, the functor
$\xi$ inverts special maps, so that for any fibration $\C$
over $I$, $\xi^*\C$ is special. We then have the functor $\C \to
\xi_*\xi^*\C$ of \eqref{2.adj}, and this is again an equivalence by
Lemma~\ref{del.le} and \eqref{del.i.1}.
\endproof

\subsection{Segal categories and Segal functors.}\label{seg.subs}

We now observe that the Segal condition of
Subsection~\ref{seg.X.subs} can be applied to fibrations over
$\Delta$ as well as to simplicial sets. The formal definition is as
follows.

\begin{defn}\label{seg.cat.def}
A {\em Segal category} is a fibration $\C \to \Delta$ that is
cartesian along the square \eqref{seg.del.sq} for any $n > l > 0$. A
{\em Segal functor} $\gamma:\C' \to \C$ between two fibrations
$\C',\C \to \Delta$ is a functor cartesian over $\Delta$. A {\em
  Segal fibration} $\C'$ over a Segal category $\C$ is a fibration
$\C' \to \C$ such that the composition $\C' \to \C \to \Delta$ turns
$\C'$ into a Segal category.
\end{defn}

\begin{exa}\label{I.seg.exa}
For any category $I$ and class of morphisms $v$ in $I$, the
simplicial expansion $\Delta^v I \to \Delta$ is a Segal
category. The product $I \times \Delta \to \Delta$ is a Segal
category as well.
\end{exa}

\begin{exa}\label{seg.pb.exa}
For any Segal functor $\gamma:\C_0 \to \C_1$ and Segal fibration
$\C_1' \to \C_1$, the pullback $\C_0' = \gamma^*\C_1 \to \C_0$ is
obviously a Segal fibration.
\end{exa}

\begin{exa}\label{opp.exa}
Since the involution $\iota:\Delta \to \Delta$, $[n] \mapsto [n]^o$
sends a square \eqref{seg.del.sq} to a square of the same form,
$\iota^*\C$ is a Segal category for any Segal category $\C$. We call
it the {\em $1$-opposite Segal category} to $\C$ and denote by
$\C^\iota$. The fibration $\C_\perp^o \to \Delta$ is also a Segal
category; we call it the {\em $2$-opposite Segal category} to $\C$
and denote by $\C^\tau$. The {\em $(1,2)$-opposite Segal category}
$\C^\sigma$ is $\C^\sigma = (\C^\tau)^\iota = \iota^*\C_\perp^o$. In
the situation of Example~\ref{I.seg.exa}, we have $(\Delta I)^\iota
\cong (\Delta I)^\sigma \cong \Delta I^o$, $(\Delta I)^\tau \cong
\Delta I$ and $(I \times \Delta)^\iota \cong I \times \Delta$, $(I
\times \Delta)^\tau \cong (I \times \Delta)^\sigma \cong I^o \times
\Delta$.
\end{exa}

\begin{exa}\label{2.cat.exa}
A Segal category with discrete $\C_{[0]}$ gives a convenient
packaging of the usual notion of a $2$-category (a.k.a.\ ``weak
$2$-category'', a.k.a.\ ``bicategory''). Explicitly, $\C_{[0]}$
consists of all the objects of the $2$-category, $\C_{[1]}$ is the
disjoint union of all of its categories of $1$-morphisms,
$s^*,t^*:\C_{[1]} \to \C_{[0]}$ send a $1$-morphism to its source
resp.\ target, the identity morphisms are provided by the functor
\begin{equation}\label{2.cat.e}
e^*:\C_{[0]} \to \C_{[1]},
\end{equation}
where as in \eqref{e.I.eq}, $e:[1] \to [0]$ is the unique
projection, and the composition is given by
\begin{equation}\label{2.cat.m}
m^*:\C_{[1]} \times_{\C_{[0]}} \C_{[1]} \cong \C_{[2]} \to \C_{[1]},
\end{equation}
where as in \eqref{m.I.eq}, $m:[1] \to [2]$ sends $0$ to $0$ and $1$
to $2$. All the higher associativity constraints are encoded in the
maps \eqref{fu.fib} for the fibration $\C \to \Delta$; for more
details on this, see e.g.\ \cite{ka.adj}.
\end{exa}

In practice, requiring that $\C_{[0]}$ is discrete is on one hand,
rather incovenient, and on the other hand, not really needed for
anything: one can treat a general Segal category as a basic
$2$-categorical object, and introduce appropriate versions of the
basic categorical constructions and notions. For example, every
cartesian functor $\gamma:\C \to \C'$ between fibrations $\C,\C' \to
\Delta$ induces a commutative square of categories
\begin{equation}\label{fib.fu.sq}
\begin{CD}
\C_{[1]} @>{\gamma_{[1]}}>> \C'_{[1]}\\
@V{s^* \times t^*}VV @VV{s^* \times t^*}V\\
\C_{[0]} \times \C_{[0]} @>{\gamma_{[0]} \times \gamma_{[0]}}>>
\C'_{[0]} \times \C'_{[0]},
\end{CD}
\end{equation}
a $2$-categorical version of the square \eqref{ngamma.sq}. If
$\gamma$ is a Segal functor between Segal categories, then
informally, the fibers of the vertical projections are the
categories of morphisms in $\C$, resp.\ $\C'$.

\begin{defn}\label{fib.fu.def}
A cartesian functor $\gamma:\C \to \C'$ between two fibrations
$\C,\C' \to \Delta$ is {\em $2$-fully faithful} resp.\ {\em
  $2$-full} resp.\ {\em $1$-full} if the square \eqref{fib.fu.sq} is
cartesian resp.\ semicartesian resp.\ weakly semicartesian. A full
subcategory $\C \subset \C'$ in a Segal category $\C$ is {\em
  $2$-full} if the embedding functor $\C \to \C'$ is $2$-fully
faithful.
\end{defn}

\begin{lemma}\label{2.ff.le}
Let $\eps = \eps([0]):\ppt \to \Delta$ be the embedding
\eqref{eps.c} onto $[0] \in \Delta$. Then a Segal functor $\gamma:\C
\to \C'$ is $2$-fully faithful if and only if the square
\begin{equation}\label{eps.sq}
\begin{CD}
\C @>{\gamma}>> \C'\\
@V{a}VV @VV{a'}V\\
\eps_*\eps^*\C @>{\eps_*\eps^*(\gamma)}>> \eps_*\eps^*\C'
\end{CD}
\end{equation}
is cartesian.
\end{lemma}

\proof{} By the Segal condition, it suffices to check that
\eqref{eps.sq} is cartesian over $[0]$ and $[1]$. Over $[0]$, the
adjuniction functors $a$, $a'$ are equivalences, so \eqref{eps.sq}
is tautologically cartesian, and over $[1]$, \eqref{eps.sq} is
exactly the square \eqref{fib.fu.sq}.
\endproof

\begin{exa}\label{v.01.exa}
For any Segal category $\C$ and closed classes of maps $v_0 \subset
v_1$ in $\C_{[0]}$, the embedding $\eps_*\C_{[0],v_0} \to
\eps_*\C_{[0],v_1}$ is trivially $2$-fully faithful in the sense of
Definition~\ref{fib.fu.def}, and then so is the Segal functor
\eqref{v.01.eq}. Both are also tautologically essentially
surjective.
\end{exa}

The classes of $1$-full, $2$-full and $2$-fully faithful Segal
functors are obviously closed under compositions and preserved by
pullbacks with respect to Segal functors. By the Segal condition, a
$2$-fully faithful Segal functor $\gamma:\C \to \C'$ is an
equivalence iff it is an equivalence over $[0] \in \Delta$, and if
$\gamma$ is only $2$-full, then it is essentially surjective
resp.\ an epivalence iff it is essentially surjective resp.\ an
epivalence over $[0]$. We also have the following corollary of
Lemma~\ref{car.epi.le}.

\begin{lemma}\label{2.epi.le}
Assume given fibrations $\C,\C' \to [1]^2 \times \Delta$ that are
cartesian over $[1]^2$, and a functor $\phi:\C' \to \C$ cartesian
over $[1]^2 \times \Delta$ that is $1$-full over $0 \times 1$, $1
\times 0$ and $2$-full over $0 \times 0$. Then $\phi$ is
$1$-full. Moreover, if $\phi$ is $2$-full over $0 \times 1$, $1
\times 0$ and $2$-fully faithful over $0 \times 0$, then $\phi$ is
$2$-full.
\end{lemma}

\proof{} Define a category $\C'(\phi)$ by the cartesian product
$$
\begin{CD}
\C'(\phi) @>>> \C_{[1]}\\
@VVV @VV{s^* \times t^*}V\\
\C'_{[0]} \times \C'_{[0]} @>{\phi_{[0]} \times \phi_{[0]}}>>
\C_{[0]} \times \C_{[0]}.
\end{CD}
$$
Then $\C'(\phi)$ is cartesian over $[1]^2$, $\phi_{[1]}:\C'_{[1]}
\to \C_{[1]}$ factors through a functor $\phi':\C'_{[1]} \to
\C'(\phi)$, and to prove the claim, it suffices to apply
Lemma~\ref{car.epi.le} to the functor $\phi'$.
\endproof

\begin{defn}\label{exa.seg.def}
A {\em $2$-family of groupoids} over a category $I$ is a Segal
fibration $\C \to \Delta \times I$ whose fibers are groupoids.
\end{defn}

For any $2$-family of groupoids $\C$ over $I$, let $\eps =
\eps_{[0]} \times \id:I \to \Delta \times I$ be the embedding onto
$\{[0]\} \times I \subset \Delta \times I$. Then \eqref{2.adj}
provides a functor
\begin{equation}\label{exa.seg.adj}
\C \to \eps_*\eps^*\C.
\end{equation}
We say that $\C$ is {\em reduced} if \eqref{exa.seg.adj} is a
discrete fibration, and define the {\em reduction}
\begin{equation}\label{seg.red}
\C^{red} = \pi_0(\C | \eps_*\eps^*\C)
\end{equation}
of $\C$ as the middle term of the factorization \eqref{grp.eq} of
the functor \eqref{exa.seg.adj}. By Lemma~\ref{sq.le}, $\C^{red}$ is
also a $2$-family of groupoids over $I$, and it is tautologically
reduced.

\begin{defn}\label{exa.bnd.def}
A $2$-family of groupoids $\C$ over a category $I$ is {\em bounded}
if the fibration \eqref{exa.seg.adj} has essentially small fibers.
\end{defn}

\begin{exa}\label{iso.seg.exa}
For any category $\C$, its simplicial expansion $\Delta^{\Iso}\C$
with respect to the class $\Iso$ of all isomorphisms, equipped with
the forgetful functor $\Delta^\Iso \C \to \Delta$, is a reduced
bounded $2$-family of groupoids over the point. We have
$\Delta^{\Iso}\C \cong \C \times \Delta$ if and only if $\C$ is a
groupoid. More generally, for any fibration $\pi:\C \to I$, let $\flat =
\pi^\flat(\Tot)$ be the class of maps in $\C$ cartesian over $I$. Then
the relative simplicial expansion $\Delta^\flat(\C|I)$ is a reduced
$2$-family of groupoids over $I$, and it is equivalent to $\C \times
\Delta$ if and only if $\C \to I$ is a family of groupoids.
\end{exa}

\begin{exa}\label{2.cat.fam.exa}
A $2$-category $\C$ in the sense of Example~\ref{2.cat.exa} is a
$2$-family of groupoids over a point if and only if all its
$2$-morphisms are invertible (nowadays, this is sometimes called a
{\em $(1,2)$-category}). For any $2$-category $\C$, we can consider
the dense subcategory $\C_{c} \subset \C$ defined by the class $c$
of maps cartesian over $\Delta$; this is the $(1,2)$-category
obtained by forgetting all non-invertible $2$-morphisms. The
reduction $\C^{red} \to \Delta$ of a $(1,2)$-category $\C$ is a
discrete fibration. If $\C$ is bounded, then its groupoids of
morphisms are essentially small, and $\C^{red} \cong \Delta \bC$ is
the simplicial replacement of the {\em truncation} $\bC$ of the
$2$-category $\C$ --- namely, the category with the same objects,
and with maps given by isomorphism classes of $1$-morphisms in $\C$.
\end{exa}

\subsection{Special functors.}

By Example~\ref{I.seg.exa}, every category $I$ produces two Segal
categories: the product $I \times \Delta$, and the simplicial
expansion $\Delta I$. We also have a natural functor
\begin{equation}\label{xi.seg}
\alpha = \xi \times \pi:\Delta I \to I \times \Delta
\end{equation}
induced by \eqref{al.xi.eq}. However, \eqref{xi.seg} is not a Segal
functor: it is a functor over $\Delta$ that is only cartesian over
$\Delta_+ \subset \Delta$ (that is, over special maps in $\Delta$ in
the sense of Subsection~\ref{spec.subs}). It will prove useful to
axiomatize the situation.

\begin{defn}\label{seg.spec.def}
A map in a fibration $\C \to \Delta$ is {\em special} if it is a
cartesian lifting of a special map in $\Delta$ A {\em special
  functor} between fibrations $\C,\C' \to \Delta$ is a functor
$\gamma:\C \to \C'$ over $\Delta$ that is cartesian over $\Delta_+
\subset \Delta$ (that is, sends special maps to special maps).
\end{defn}

\begin{defn}\label{simpl.rel.def}
The {\em $2$-simplicial expansion} $\Delta^{\Tot}(\C\|\Delta)$ of an
arbitrary fibration $\C \to \Delta$ is given by
$$
\Delta^{\Tot}(\C\|\Delta) = \Ind^+(\C),
$$
where $\Ind^+$ is the induction of Subsection~\ref{ind.subs} with
respect to the class of special maps. The {\em universal special
  functor} is the functor
\begin{equation}\label{2.del.al}
  \alpha:\Delta^\Tot(\C\|\Delta) \to \C
\end{equation}
induced by \eqref{ind.facto}.
\end{defn}

\begin{exa}
By \eqref{I.wr}, for any category $I$, we have a natural
identification $\Delta^{\Tot}(I \times \Delta\|\Delta) =
\Delta^{\Tot} I$. The universal special functor \eqref{2.del.al} is
$\alpha = \pi \times \xi$, where $\pi:\Delta^\Tot I \to \Delta$ is
the structural fibration, and $\xi$ is the functor \eqref{xi.eq}.
\end{exa}

One can also describe $2$-simplicial expansion in a way similar to
Definition~\ref{rel.simp.def} but with \eqref{be.eq} replaced with
\eqref{al.0.eq}. Namely, for any fibration $\C \to \Delta$, we have
a cartesian square
\begin{equation}\label{d.c.d}
\begin{CD}
\Delta^{\Tot}(\C \| \Delta) @>>> \Delta^{\Tot} \C\\
@VVV @VV{\Delta(\pi)}V\\
\Delta @>{\delta}>> \Delta^{\Tot}\Delta,
\end{CD}
\end{equation}
where $\delta$ is the functor \eqref{al.0.eq}. One can then combine
the two constructions: if we have a category $I$ and a fibration
$\pi:\C \to \Delta \times I$, we define the {\em relative
  $2$-simplicial expansion} $\Delta^{\Tot}(\C \| \Delta | I)$ by the
cartesian square
\begin{equation}\label{d.c.d.i}
\begin{CD}
\Delta^{\Tot}(\C \| \Delta | I) @>>> \Delta^{\Tot}(\C \| \Delta)\\
@VVV @VV{\Delta(\pi)}V\\
\Delta \times I @>{\beta}>> \Delta^{\Tot}I = \Delta^{\Tot}(I \times
\Delta \| \Delta),
\end{CD}
\end{equation}
or equivalently, by the cartesian square \eqref{d.c.d} with the
category $\Delta^{\Tot}\Delta$ replaced by $\Delta^{\Tot}(\Delta
\times I) \cong \Delta^{\Tot}\Delta \times_{\Delta} \Delta^{\Tot} I$
and $\delta$ replaced with $(\delta \circ \tau) \times \beta$, where
$\tau:\Delta \times I \to \Delta$ is the projection. Then
$\Delta^{\Tot}(\C \| \Delta | I) \to \Delta \times I \to I$ is a
fibration with fibers $\Delta^{\Tot}(\C_i \| \Delta)$, $i \in I$.

The $2$-simplicial expansion is functorial, in that a functor
$\pi:\C' \to \C$ cartesian over $\Delta_+$ induces a functor
$\Delta^{\Tot}(\pi):\Delta^{\Tot}(\C'\|\Delta) \to
\Delta^{\Tot}(\C\|\Delta)$ cartesian over $\Delta$. For any two
fibrations $\C,\C' \to \Delta$ with $\C$ essentially small,
\eqref{ind.2adj} provides an equivalence of categories
\begin{equation}\label{2exp.eq}
\Fun^+_\Delta(\C,\C') \cong
\Fun^{\Tot}_\Delta(\C,\Delta^{\Tot}(\C'\|\Delta)).
\end{equation}
Even if $\C$ is not essentially small, special functors from $\C$ to
$\C'$ are still naturally identified with functors from $\C$ to
$\Delta^{\Tot}(\C'\|\Delta)$ cartesian over $\Delta$. In particular,
this applies to Segal categories $\C$, $\C'$, and it is useful to
know when $\Delta^{\Tot}(\C'\|\Delta)$ is also a Segal
category. Here is a criterion for this.

\begin{defn}\label{seg.reg.def}
A Segal category $\C \to \Delta$ is {\em coregular} resp.\ {\em
  regular} if for any $[n] \in \Delta$, with the map $t:[0] \to
[n]$, the pullback functor $t^*:\C_{[n]} \to \C_{[0]}$ resp.\ the
opposite functor $t^{*o}:\C_{[n]}^o \to \C_{[0]}^o$ is a family of
groupoids.
\end{defn}

\begin{lemma}\label{rel.del.le}
The $2$-simplicial expansion $\Delta^{\Tot}(\C\|\Delta)$ of a
regular Segal category $\C$ is a Segal category.
\end{lemma}

\proof{} By induction, it suffices to check that
$\Delta^{\Tot}(\C\|\Delta)$ is cartesian along the squares
\eqref{seg.del.sq} with $l = n-1$. By Example~\ref{ind.exa}, for any
$[n] \in \Delta$, we have $\Delta^{\Tot}(\C\|\Delta)_{[n]} \cong
\Sec([n],\phi^*\C)$, where $\phi:[n] \cong \Delta^\hdot_{[n]} \to
\Delta$ sends $l \in [n]$ to $[l]$ and a map $l \leq l'$ to $s:[l]
\to [l']$. In particular, we have $\Delta^{\Tot}(\C\|\Delta)_{[0]}
\cong \C_{[0]}$, we have an equivalence
\begin{equation}\label{tot.del.1}
\Delta^{\Tot}(\C\|\Delta)_{[1]} \cong \C_{[0]} \setminus_{s^*}
\C_{[1]},
\end{equation}
and for any $n \geq 2$, we have a cartesian square
\begin{equation}\label{cor.C.sq}
\begin{CD}
\Delta^{\Tot}(\C\|\Delta)_{[n]} @>{s^*}>>
\Delta^{\Tot}(\C\|\Delta)_{[n-1]}\\
@VVV @VVV\\
\C_{[n-1]} \setminus_{s^*} \C_{[n]} @>{\tau}>> \C_{[n-1]}.
\end{CD}
\end{equation}
Therefore it suffices to prove that the square
$$
\begin{CD}
\C_{[n-1]} \setminus_{s*} \C_{[n]} @>{t^*}>> \C_{[0]}
\setminus_{s^*} \C_{[1]}\\
@V{\tau}VV @VV{\tau}V\\
\C_{[n-1]} @>{t^*}>> \C_{[0]}
\end{CD}
$$
is cartesian. By the Segal property of $\C$, checking this amounts to
checking that the functor $\C_{[n-1]}/_{\id}\C_{[n-1]} \to
\C_{[n-1]}/_{t^*}\C_{[0]}$ is an equivalence, and since $\C$ is
regular, this follows from Lemma~\ref{grp.c.le}.
\endproof

\begin{exa}\label{seg.reg.exa}
A constant Segal category $\C \times \Delta \to \Delta$ is both
regular and coregular, and so is a Segal category $\C \to \Delta$
that is itself a family of groupoids.
\end{exa}

\begin{exa}\label{spec.reg.exa}
For any fibration $\C \to I$, with the corresponding special
fibration $\xi^*\C \to \Delta I$, the category $\xi^*\C$ is both
regular and coregular.
\end{exa}

\begin{exa}\label{exa.reg.exa}
If a Segal category $\C$ is regular, then the $2$-opposite Segal
category $\C^\tau$ is coregular, and vice versa. In particular, for
any category $I$ and $2$-family of groupoids $\C \to I \times
\Delta$ in the sense of Definition~\ref{exa.seg.def}, $\C$ is
obviously coregular -- $t^*$ is a functor between families of
groupoids over $I$, thus itself a family of groupoids -- and then
$\C^\tau$ is regular. Since $\C^\iota$ is also a $2$-family of
groupoids over $I$, the $(1,2)$-opposite Segal category $\C^\sigma$
is regular as well.
\end{exa}

Just as in the absolute case of Definition~\ref{simpl.repl.def}, one
can refine Definition~\ref{simpl.rel.def} by considering a closed
class $v$ of maps in $\C_{[0]}$, and defining the $2$-simplicial
expansion $\Delta^v(\C\|\Delta)$ by the cartesian square
\begin{equation}\label{d.c.v}
\begin{CD}
\Delta^v(\C\|\Delta) @>>> \Delta^{\Tot}(\C\|\Delta)\\
@VVV @VVV\\
\eps_*\C_{[0],v} @>>> \eps_*\C_{[0]},
\end{CD}
\end{equation}
where as in \eqref{del.I.sq}, $\eps = \eps([0]):\ppt \to \Delta$ is
the embedding onto $[0] \in \Delta$. Then in particular, have
$\Delta^v(I \times \Delta\|\Delta) \cong \Delta^vI$ for any category
$I$. If we have two closed classes $v_0 \subset v_1$, then the
embedding $\C_{[0],v_0} \to \C_{[0],v_1}$ induces a Segal functor
\begin{equation}\label{v.01.eq}
\Delta^{v_0}(\C\|\Delta) \to \Delta^{v_1}(\C\|\Delta).
\end{equation}
As in Definition~\ref{simpl.repl.def}, we write
$\Delta(\C\|\Delta)=\Delta^{\Id}(\C\|\Delta)$ to simplify
notation. If we have a fibration $\pi:\C \to \Delta \times I$ and a
closed class of maps in $\C_{[0]}$ that contains all maps cartesian
over $I$, we define the relative $2$-simplicial expansion
$\Delta^v(\C \| \Delta | I)$ by replacing
$\Delta^{\Tot}(\C\|\Delta)$ in \eqref{d.c.v} with
$\Delta^{\Tot}(\C\|\Delta|I)$.

\subsection{Twisted $2$-simplicial expansion.}\label{tw.subs}

To illustrate Definition~\ref{simpl.rel.def}, let us write down
explicitly the Grothendieck construction in terms of nerves and
$2$-simplicial expansions. To do this, for any fibration $\C \to
\Delta$, define the {\em twisted $2$-simplicial expansion}
$\Delta^v_\sigma(\C\|\Delta)$ with respect to a closed class $v$ of
maps in $\C_{[0]}$ by
\begin{equation}\label{tw.rel.eq}
\Delta^v_\sigma(\C\|\Delta) = \Delta^{v^o}(\C^\sigma\|\Delta)^\sigma,
\end{equation}
where $v^o$ is the class of maps in $\C^\sigma_{[0]} \cong
\C^o_{[0]}$ opposite to $v$. If $\C^\sigma$ is a regular Segal
category, then \eqref{tw.rel.eq} is a Segal category by
Lemma~\ref{rel.del.le}.  Now assume given a category $I$ and a
discrete Segal fibration $\C \to I \times \Delta$. Then for any $i
\in I$, the restriction $\C_i = \C|_{\{i\} \times \Delta}$ is a
discrete Segal category, so that $\C_i \cong \Delta I'_i$ for some
category $I'_i$. Moreover, for any map $f:i \to i'$ in $I$, the
transition functors of the fibration $\C$ induce a functor
$f^*:I'_{i'} \to I'_i$, and this turns the collection $I'_i$ into a
fibration $I' \to I$ (in fact, $I' \to I$ is a classical fibration
equipped with a strict cleavage, as in Remark~\ref{cleav.rem}). Then
$\C^\sigma$ is regular by Example~\ref{exa.reg.exa}, and we have a
natural identification
\begin{equation}\label{str.del.eq}
\Delta I' \cong \Delta_\sigma(\C\|\Delta),
\end{equation}
where the right-hand side is the twisted $2$-simplicial expansion
\eqref{tw.rel.eq} with respect to $v = \Id$. For example, if $\C = I
\times \Delta$, so that $I'=I$, then indeed, $\Delta_\sigma(I \times
\Delta\|\Delta) = \Delta(I^o \times \Delta \| \Delta)^\sigma \cong
(\Delta I^o)^\sigma \cong \Delta I$.

More generally, for any category $I$ and Segal fibration $\pi:\C \to
I \times \Delta$, we can let $\flat = \pi^*\Id$ be the class of maps
$f$ in $\C_{[0]}$ such that $\pi(f)$ is an identity map in $I$. Then
we have the twisted $2$-simplicial expansion $\Delta^\flat_\sigma(\C
\| \Delta)$ and a fibration
$\Delta(\pi):\Delta^\flat_\sigma(\C\|\Delta) \to \Delta I$. If we
consider the class $\Tot$ of all maps, then we also have the
fibration $\Delta^\Tot(\pi):\Delta^\Tot_\sigma(\C\|\Delta) \to
\Delta^\Tot I$. If $\C^\sigma$ is regular, then $\Delta(\pi)$ and
$\Delta^\Tot(\pi)$ are Segal fibrations.

When $\pi$ is a family of groupoids, one can describe the twisted
simplicial expansion \eqref{tw.rel.eq} without passing to the
$(1,2)$-opposite Segal categories. To do this, consider the full
embedding $\beta$ of \eqref{be.eq}, and note that since any
section of a family of groupoids is automatically cartesian, we have
a natural identification
\begin{equation}\label{tw.beta.eq}
\beta^*\Delta^\Tot_\sigma(\C\|\Delta) \cong \C
\end{equation}
for any family of groupoids $\C \to I \times \Delta$.

\begin{lemma}\label{beta.le}
For any family of groupoids $\C \to I \times \Delta$, the functor
\begin{equation}\label{tw.rel.be.eq}
  \Delta^\Tot_\sigma(\C\|\Delta) \to \beta_*\C
\end{equation}
adjoint to \eqref{tw.beta.eq} is an equivalence of categories.
\end{lemma}

\proof{} For any $[n] \in \Delta$ and fibration $\C \to \Delta$,
\eqref{cor.C.sq} induces a cartesian square
\begin{equation}\label{cor.tw.sq}
\begin{CD}
\Delta^{\Tot}_\tau(\C\|\Delta)_{[n]} @>{t^*}>>
\Delta^{\Tot}_\tau(\C\|\Delta)_{[n-1]}\\
@VVV @VVV\\
\C^\tau_{[n-1]} \setminus_{t^*} \C_{[n]} @>{\tau}>> \C^\tau_{[n-1]}.
\end{CD}
\end{equation}
For the fibration $I \times \Delta \to \Delta$, we have
$\Delta^\Tot_\sigma(I \times \Delta \|\Delta) \cong \Delta^\Tot I$,
and cartesian sections of a corresponding square \eqref{cor.tw.sq}
are commutative squares
\begin{equation}\label{cor.I.sq}
\begin{CD}
\langle [n-1],i_0 \rangle @>{l \circ f_0}>> \langle [n-1],t^*i_\idot
\rangle\\
@V{t}VV @VV{t}V\\
\langle [n],i_0 \rangle @>>> \langle [n],i_\idot \rangle
\end{CD}
\end{equation}
in $\Delta^\Tot I$, where $\langle [n],i_\idot \rangle \in
\Delta^\Tot I$ is an arbitrary object represented by a diagram
\eqref{i.dot}, $f_0:i_0 \to i_1$ is the map in the diagram, and
$l:\langle [n],i_0 \rangle \to \langle [n],i_\idot \rangle$,
$l:\langle [n-1],i_1 \rangle \to \langle [n-1],t^*i_\idot \rangle$
are the tautological maps. For our given fibration $\C \to I \times
\Delta$, we also have the cartesian squares \eqref{cor.tw.sq}, and
then by Lemma~\ref{car.rel.le}, the fibration
$\Delta^\Tot_\sigma(\C\|\Delta) \to \Delta^\Tot I$ is cartesian along
each square \eqref{cor.I.sq}. Since \eqref{tw.rel.be.eq} is an
equivalence over the essential image of the full embedding $\beta$,
thus over both the source and the target of the left vertical arrow
in \eqref{cor.I.sq}, it suffices to check that $\beta_*\C$ is also
cartesian over \eqref{cor.I.sq} and apply induction on $n$. By
\eqref{2kan.eq}, this amounts to checking that the square
\begin{equation}\label{cor.co.sq}
\begin{CD}
(I \times \Delta) /_\beta \langle [n-1],i_0 \rangle @>>> (I \times
  \Delta) /_\beta \langle [n-1],t^*i_\idot
\rangle\\
@VVV @VVV\\
(I \times \Delta) /_\beta \langle [n],i_0 \rangle @>>> (I \times
\Delta) /_\beta \langle [n],i_\idot \rangle
\end{CD}
\end{equation}
induced by \eqref{cor.I.sq} is a cocartesian square of
categories. However, for any $\langle [n],i_\idot \rangle$, the
discrete fibration $(I \times \Delta) /_\beta \langle [n],i_\idot
\rangle \to I \times \Delta$ corresponds to the relative nerve
$N_I(\langle [n],i_\idot \rangle):I^o \times \Delta^o \to \Sets$ of
\eqref{N.I.eq}, so by Lemma~\ref{cocart.le}, it suffices to check
that $N_I$ sends a square \eqref{cor.I.sq} to a cocartesian
square. This is Example~\ref{J.n.sq} and Remark~\ref{BJ.N.rem}.
\endproof

\section{Reedy categories.}\label{reedy.sec}

\subsection{Cellular Reedy categories.}\label{cell.subs}

Recall that by Definition~\ref{good.def}, a {\em good filtration} on
a category $I$ is a collection of full subcategories $I_{\leq n}
\subset I$, $n \geq 0$ such that $I_{\leq n} \subset I_{\leq n+1}$
and $I = \bigcup I_{\leq n}$. For any small category $I$ equipped
with a good filtration $I_{\leq n}$, denote by $I_n = I_{\leq n}
\ssetminus I_{\leq n-1}$ the full subcategory spanned by objects in
$I_{\leq n}$ but not in $I_{\leq n-1}$, and let $\deg(i)$, $i \in I$
be the unique integer such that $i \in I_{\deg(i)}$. Note that the
degree function defines the filtration: $I_{\leq n} \subset I$ is
the full subcategory spanned by $i \in I$, $\deg(i) \leq n$.

\begin{defn}\label{dir.def}
A {\em directed category} is a small category $I$ equipped with a
good filtration $I_{\leq n}$ such that for any map $f:i \to i'$,
$\deg(i') \geq \deg(i)$, and the inequality is strict unless $i=i'$
and $f=\id$. A directed category $I$ is {\em locally finite} resp.\ {\em
  locally $\kappa$-bounded} for an infinite cardinal $\kappa$ iff
for any $i \in I$, $I/i$ is finite resp.\ $|I/i| < \kappa$.
\end{defn}

\begin{remark}\label{BB.rem}
A directed category $I$ is obviously ordered in the sense of
Lemma~\ref{BB.le}.
\end{remark}

\begin{defn}\label{reedy.def}
A {\em Reedy category} is a small category $I$ equipped with a good
filtration $I_{\leq n}$ and a factorization system $(M,L)$ such that
both $I_L$ and $I_M^o$ with the induced good filtrations are
directed. A Reedy category $I$ is {\em Reedy-finite} if so are the
directed categories $I_L$, $I_M^o$, and {\em Reedy-$\kappa$-bounded}
for some infinite cardinal $\kappa$ if $I_L$ is locally
$\kappa$-bounded and $I_M^o$ is locally $\kappa'$-bounded for some
$\kappa' \ll \kappa$ (where $\ll$ is as in
Example~\ref{sets.SS.exa}). A {\em Reedy functor} between two Reedy
categories $\langle I,M,L \rangle$, $\langle I',M',L' \rangle$ is a
functor $\gamma:I' \to I$ sending maps in $L'$ resp.\ $M'$ to maps
in $L$ resp.\ $M$.
\end{defn}

\begin{exa}\label{del.reedy.exa}
The category $\Delta$ is a locally finite Reedy category, with
$\deg([n])=n$, and $L = \ff$ resp.\ $M = \dd$ consisting of injective
resp.\ surjective maps.
\end{exa}

\begin{exa}\label{reedy.exa}
The product $I \times I'$ of two Reedy categories is a Reedy
category, with degree function $\deg(i \times i')=\deg(i)+\deg(i')$
and factorization system $(M \times M',L \times L')$, locally finite
if so are $I$, $I'$. The opposite $I^o$ of a Reedy category is a
Reedy category, with the same degree function and factorization
system $(L^o,M^o)$, locally finite if so is $I$. A partially ordered
set $I$ left-bounded in the sense of Definition~\ref{pos.fin.def} is
a Reedy category, with $\deg(i) = \dim(I/i)$, $I_L=I$, and discrete
$I_M$; $\langle I,L,M \rangle$ is locally finite iff $I$ is
left-finite.
\end{exa}

Definition~\ref{dir.def} implies that the subcategories $I_n \subset
I$, $n \geq 0$ in a directed category $I$ are discrete, so that for
any category $\C$, the functor category $\Fun(I_n,\C)$ is a product
of copies of $\C$. A Reedy category $I$ is rigid in the sense of
Definition~\ref{iso.def} but the subcategories $I_n$ are no longer
discrete --- they can have non-trivial non-invertible maps. However,
for any $n \geq 0$, we have a discrete dense subcategory
$\overline{I}_n = I_{L,n} = (I_{M,n}^o)^o \subset I_n$. It turns out
that the functor category $\C^I$ can be described inductively in
terms of the categories $\Fun(\overline{I}_n,\C)$, $n \geq
0$. Namely, assume that the target category $\C$ is complete and
cocomplete. Then for any $n \geq 0$, the restriction functor
$R:\Fun(I_{\leq n+1},\C) \to \Fun(I_{\leq n},\C)$ admits fully
faithful left and right-adjoint functors $R_\dg,R^\dg:\Fun(I_{\leq
  n},\C) \to \Fun(I_{\leq n+1},\C)$. Explicitly, for any $i \in I$, one
defines the {\em latching category} $L(i) \subset I_L/i \subset I/i$
as the full subcategory spanned by non-identical maps $i' \to i$
(that is, $L(i) \subset I_L(i)$ contains all objects except for the
terminal object $\id:i \to i$). Analogously, one defines the {\em
  matching category} $M(i)$ as $M(i) = (i \ssetminus I_M) \ssetminus
\{\id\}$. Then one has
\begin{equation}\label{LM.eq}
R_\dg(X)(i) = \colim_{L(i)}X|_{L(i)}, \ R^\dg(X)(i) =
\lim_{M(i)}X|_{M(i)}, \ i \in I_{n+1}
\end{equation}
for any functor $X:I_{\leq n} \to \C$. One then denotes
$L(X,i)=R_\dg(X)(i)$, $M(X,i)=R^\dg(X)(i)$, and calls them the {\em
  latching} and {\em matching object} of the functor $X$ with
respect to $i$. By adjunction, we have a functorial map $R_\dg \to
R^\dg$, thus a collection of maps $L(X,i) \to M(X,i)$, $i \in I_{n+1}$,
and one shows that extending $X$ to $I_{\leq n+1} \supset
I_{\leq n}$ is equivalent to providing objects $X(i)$, $i \in
I_{n+1}$, and factorizations
\begin{equation}\label{lm.eq}
\begin{CD}
L(X,i) @>{l}>> X(i) @>{m}>> M(X,i)
\end{CD}
\end{equation}
of the adjunction maps $L(X,i) \to M(X,i)$. The maps $l$, $m$ are
known as the {\em latching} and {\em matching} map. If $I$ is
Reedy-finite in the sense of Definition~\ref{reedy.def}, then it
suffices to require that $\C$ is finitely complete and finitely
cocomplete, and if $I$ is Reedy-$\kappa$-bounded, then it suffices
for $\C$ to be $\kappa$-cocomplete and $\kappa'$-complete for all
ordinals $\kappa' \ll \kappa$.

\begin{exa}\label{reedy.L.exa}
For any discrete fibration $\pi:I' \to I$ over a Reedy category $I$,
the category $I'$ with the classes $L'=\pi^*L$, $M'=\pi^*M$ and
degree function $\pi^*\deg$ is a Reedy category, $\pi$ is a Reedy
functor, and the induced functor $\pi:L(i) \to L(\pi(i))$ is an
equivalence for any $i \in I'$.
\end{exa}

For $I=\Delta$ and $I=\Delta^o$, this inductive description was
obtained by Reedy in \cite{reedy}. Definition~\ref{reedy.def} and
its various versions that appear in the literature is an abstraction
of the relevant properties of the categories $\Delta$ and
$\Delta^o$. We will also need to axiomatize some propeties such as
those listed in Subsection~\ref{del.0.subs} that are shared by
$\Delta$ and the categories of simplices $\Delta X$, but not
$\Delta^o$. To this end, we introduce the following definition that
is essentially borrowed from \cite[II.3.1]{GZ}.

\begin{defn}\label{ind.reedy.def}
A Reedy category $I=\langle I,\deg,M,L \rangle$ is {\em cellular}
if
\begin{enumerate}
\item all maps in $M$ resp.\ $L$ are surjective resp.\ injective, and
\item for any map $m:i \to i'$ in $M$, there exists a map $l:i' \to
  i$ in $L$ such that $m \circ l = \id$, and $m \in I_M(i,i')$ is
  uniquely determined by the subset $\Sec(m) = \{l|m \circ l = \id\}
  \subset I_L(i',i)$.
\end{enumerate}
\end{defn}

The category $\Delta$ with the standard Reedy structure is obviously
a cellular Reedy category. The same is true for the category of
simplices $\Delta X$ of a simplicial set $X \in \Delta^o\Sets$, by
virtue of the following result.

\begin{lemma}\label{ind.reedy.le}
Every left-bounded partially ordered set $J$ is a cellular Reedy
category. The product $I \times I'$ of two cellular Reedy categories
is a cellular Reedy category. For any cellular Reedy category $I$
and small discrete fibration $\pi:I' \to I$, $I'$ is a cellular
Reedy category.
\end{lemma}

\proof{} The first two claims are obvious (the Reedy structures in
both cases are those of Example~\ref{reedy.exa}). For the last
claim, define the classes $L'$, $M'$ by $f \in L',M'$ iff $\pi(f)
\in L,M$, and let $\deg(i') = \deg(\pi(i'))$, $i' \in I'$, as in
Lemma~\ref{facto.fib.le}~\thetag{i}. This turns $I'$ into a Ready
category, and we need to check that it is cellular. As noted in
Example~\ref{IX.exa}, a discrete fibration is always
faithful, so that Definition~\ref{ind.reedy.def}~\thetag{i} for $I$
immediately implies the corresponding statement for $I'$. Since
$\pi$ is a discrete fibraton, we have $\Sec(m) \cong \Sec(\pi(m))$
for any map $m:i \to i'$ in $I'_{M'}$, and since $\pi$ is
faithful, Definition~\ref{ind.reedy.def}~\thetag{ii} for $I$ implies
the corresponding statement for $I'$.
\endproof

Note that by Definition~\ref{ind.reedy.def}~\thetag{i}, for any
object $i \in I$ in a cellular Reedy category $I$, the
comma-category $i \setminus I_M$ is a partially ordered set. Here is
the main property of cellular Reedy categories that motivates the
definition.

\begin{lemma}\label{ez.le}
  For any object $i \in I$ in a cellular Reedy category $I$, the
  partially ordered set $i \setminus I_M$ has the largest element.
\end{lemma}

\proof{} We use the argument of \cite[II.3.1]{GZ} (the
``Eilenberg-Zilber Lem\-ma''). The degree function $\deg:i \setminus
I_M \to [n]$ is obviously conservative as a functor, so $\dim (i
\setminus I_M) \leq n$, the finite-dimensional partially ordered set
$i \setminus I_M$ has maximal elements, and it suffices to show that
they all coincide. If we have two such elements $m_0:i \to i_0$,
$m_1:i \to i_1$, then by Definition~\ref{ind.reedy.def}~\thetag{ii}
we can choose $l_0:i_0 \to i$ such that $m_0 \circ l_0 = \id$, and
then by maximality, $m_1 \circ l_0:i_0 \to i_1$ is in $L$. Dually,
for any $l_1 \in \Sec(m_1)$, the map $m_0 \circ l_1:i_1 \to i_0$ is
also in $L$.  But $I_L$ is ordered by the Reedy property, so that
$i_0=i_1$ and $m_1 \circ l_0 = \id$, $m_0 \circ l_1=\id$, for any
choices of $l_0$ and $l_1$. Therefore $\Sec(m_0)=\Sec(m_1)$, so that
$m_0=m_1$ by Definition~\ref{ind.reedy.def}~\thetag{ii}.
\endproof

\begin{remark}
For any cellular Reedy category $I$ and discrete fibration $\pi:I'
\to I$, the induced projection $\pi:I'_{M'} \to I_M$ is a
discrete fibration, and for any object $i \in I'$, the same holds
for the projection $\pi_i:i \setminus I'_{M'} \to \pi(i) \setminus
I_M$. Then Lemma~\ref{ind.reedy.le} and Lemma~\ref{ez.le} imply that
$i \setminus I'_{M'}$ has a largest element $m$, and by
Example~\ref{comma.discr.exa}, this means that $\pi_i$ identifies $i
\setminus I'_{M'}$ with the left-closed subset $(\pi(i) \setminus
I_M)/m \subset \pi(i) \setminus I_M$. This shows that if $\pi(i)
\setminus I_M$ satisfies some additional properties, such as being
finite and/or having finite coproducts, these are inherited by $i
\setminus I'_{M'}$. We note that both properties hold when
$I=\Delta$ (the easiest way to see it is to use the Joyal duality
\eqref{jo.eq} to identify $([n] \setminus \Delta_M)^o$ with the
partially ordered set of subsets $S \subset [n]$ that contain $n \in
[n]$).
\end{remark}

\begin{lemma}\label{I.corr}
Assume given a cellular Reedy category $I$, a Reedy functor $\chi:I
\to J$ to a left-bounded partially ordered set $J$ with the Reedy
structure of Example~\ref{reedy.exa}, and a subset $J' \subset
J$. Then the preimage $I' = \chi^{-1}(J') = I \times_J J'$ inherits
a cellular Reedy structure such that the embedding $\eps:I' \to I$
is a Reedy functor.
\end{lemma}

\proof{} Restrict the degree function on $I$ to $I' \subset I$, and
say that a map $f$ in $I'$ is in $L$ resp.\ $M$ iff $\eps(f)$ is in
$L$ resp.\ $M$. Since $\chi$ is a Reedy functor, $\chi(m) = \id$ for
any $m \in M$, so that for any map $f$ in $I'$, the canonical
factorization $\eps(f) = l \circ m$ lies entirely in $I'$. Thus
$(M,L)$ is a factorization system on $I'$. Moreover, for any $i \in
I'$, we have $M(i) = M(\eps(i))$ and $L(i) = L(\eps(i)) \cap
I'$. Therefore $I'$ is a Reedy category satisfying
Definition~\ref{ind.reedy.def}~\thetag{i}. For \thetag{ii}, observe
that $\chi(l)=\id$ for any $l \in \Sec(m)$, so that for any $m:i \to
i'$ in $I'_M \subset I_M$, the set $\Sec(m)$ computed in $I'$ and in
the ambient category $I$ is the same.
\endproof

\subsection{Skeleta.}

Assume given a cellular Reedy category $I$. Say that an object $i
\in I$ is {\em non-degenerate} if the matching category $M(i)$ is
empty, and denote by $\M(I) \subset I$ the full subcategory spanned
by all non-degenerate objects. Since $(M,L)$ is a factorization
system on $I$, any map $f:i \to i'$ in $I$ with non-degenerate $i$
is automatically in the class $L$, so that in particular, we have
$\M(I) \subset I_L$, and $\M(I)$ inherits a Reedy category
structure, with discrete $\M(I)_M \subset \M(I)$.

\begin{lemma}\label{I.M.le}
The embedding functor $\M(I)_M \subset I_M$ has a left-adjoint
$\lambda:I_M \to \M(I)_M$.
\end{lemma}

\proof{} By Lemma~\ref{ez.le}, for any $i \in I$, the category $i
\setminus I_M$ has the terminal object; $\lambda(i)$ is this object,
with the structure map $i \to \lambda(i)$ providing the adjunction
map.
\endproof

\begin{defn}\label{reg.I.def}
  A cellular Reedy category $I$ is {\em regular} if for any map $l:i
  \to i'$ in $L$ with non-degenerate $i'$, $i$ is also
  non-degenerate.
\end{defn}

\begin{exa}
  A cellular Reedy category $I$ that is directed -- that is, $I =
  I_L$ -- is trivially regular in the sense of
  Definition~\ref{reg.I.def}.
\end{exa}

\begin{exa}
The simplex category $\Delta X$ of a simplicial set $X$ is regular
in the sense of Definition~\ref{reg.I.def} iff $X$ is weakly regular
in the sense of Subsection~\ref{del.0.subs}.
\end{exa}

\begin{lemma}\label{reg.I.le}
Let $\langle I',L',M' \rangle$ and $\langle I,M,L \rangle$ be
cellular Reedy categories, and assume given a functor $\phi:I' \to
I^o$. Then the twisted comma-category $I \setminus_\phi^o I'$ of
Example~\ref{tw.comma.exa} is a cellular Reedy category. Moreover,
assume that $I'$ is regular, and $L' \subset \phi^*(L^o)$, $M'
\subset \phi^*(\Iso)$. Then $I \setminus_\phi^o I'$ is regular, and
we have
\begin{equation}\label{reg.M.eq}
  \M(I \setminus_\phi^o I') \cong \M(I') \setminus_\phi^o I_L,
\end{equation}
where since $L' \subset \phi^*(L^o)$, the functor $\phi:\M(I') \to
I^o$ factors through $I_L^o$.
\end{lemma}

\proof{} For the first claim, note that by
Example~\ref{tw.comma.exa}, we have a small discrete fibration $I
\setminus_\phi^o I' \to I \times I'$, and apply
Lemma~\ref{ind.reedy.le}. For the second claim, note that
explicitly, objects in $I \setminus_\phi^o I'$ are triples $\langle
i,i',\alpha \rangle$, $i \in I$, $i' \in I'$, $\alpha:i \to
\phi(i')$, and morphisms from $\langle i_1,i'_1,\alpha_1 \rangle$ to
$\langle i_0,i'_0,\alpha_0 \rangle$ are pairs $\langle f,f' \rangle$
of morphisms $f:i_1 \to i_0$, $f':i_1' \to i_0'$ such that the
diagram
\begin{equation}\label{tw.I.dia}
\begin{CD}
  i_0 @>{\alpha_0}>> \phi(i_0')\\
  @A{f}AA @VV{\phi(f')^o}V\\
  i_1 @>{\alpha_1}>> \phi(i_1')
\end{CD}
\end{equation}
is commutative. Then since $M \subset \phi^*(\Iso)$, the fibration
$I \setminus_\phi^o I' \to I \times I'$ is constant along maps of
the form $\id \times m$, $m \in M'$, so that if $\langle
i_0,i'_0,\alpha_0 \rangle$ is non-degenerate, then $i'_0$ is
non-degenerate in $I'$. Conversely, by the unique
factorization property, $\langle i_0,i_0',\alpha_0 \rangle$ with
non-degenerate $i_0'$ is non-degenerate if and only if $\alpha_0 \in
L$. Then we have some map $\langle f,f' \rangle:\langle
i_1,i_1',\alpha_1 \rangle \to \langle i_0,i_0',\alpha_0 \rangle$ in
the class $L$, then $f \in L$ and $f' \in L'$, and since $L' \in
\phi^*(L^o)$, we also have $\phi(f')^o \in L$, so that $\alpha_1$ in
\eqref{tw.I.dia} is also in $L$, and $\langle i_1,i_1',\alpha_1
\rangle$ is non-degenerate.
\endproof

\begin{corr}\label{reg.I.corr}
For any cellular Reedy category $I$ and object $i \in I$, the
category $I/i$ is regular.
\end{corr}

\proof{} Take $I' = \ppt$, and let $\phi = \eps(i):\ppt \to I^o$, as
in \eqref{eps.c}.
\endproof

Now consider a functor $X \in I^o\Sets$, with the category of
elements $IX$ of Example~\ref{IX.exa}. For any $n \geq 0$, denote by
$\eps_n:I_{\leq n} \to I$ the embedding functors and define the {\em
  $n$-th skeleton} $\sk_n X$ by
\begin{equation}\label{sk.eq}
\sk_nX = \eps_{n!}\eps_n^*X.
\end{equation}
By adjunction, we have a natural map $\sk_nX \to X$.

\begin{lemma}\label{sk.le}
For any $n \geq 0$ and $X:I^o \to \Sets$, the map $\sk_nX \to X$ is
injective.
\end{lemma}

\proof{} Let $\pi:IX \to I$ be the discrete fibration corresponding
to $X$, and let $x:IX^o \to \Sets$ be the constant functor with
value $\ppt$. Then by \eqref{kan.eq} and \eqref{IX.sets}, we have $X
\cong \pi^o_!x$, $\eps_n^*X \cong \pi^o_!\eps_n^*x$, and $\pi^o_!$
sends injective maps to injective maps. Therefore we may assume
right away that $IX = I$ and $X=x$ is the constant functor with
value $\ppt$. Then for any $i \in I$, we have
\begin{equation}\label{sk.M}
\sk_nX(i) = \colim_{M(i)^o \cap I^o_{\leq n}}X
\end{equation}
by \eqref{LM.eq}, and $M(i) \cap I_{\leq n}$ is a finite partially
ordered set that is either empty, or has a largest element.
\endproof

Lemma~\ref{sk.le} allows one to describe functors $X \in I^o\Sets$
inductively in terms of their skeleta. Namely, extend \eqref{sk.eq}
to $n=-1$ by setting $\sk_{-1}X = \emptyset$ (that is, the constant
functor with value $\emptyset \in \Sets$); note that
Lemma~\ref{sk.le} tautologically also holds for $n=-1$. Then for any
object $i \in I$, denote by $\Delta_i = \Y(i) \in I^o\Sets$ the
functor it represents -- that is, $\Delta_i(i')=I(i',i)$ -- and
define the {\em $i$-sphere} as $\SS_i = \sk_{\deg(i)-1}\Delta_i$ (if
$I = \Delta$ and $i = [n]$ for some $n \geq 0$, then the
$[n]$-sphere is just the usual simplicial sphere $\SS_{n-1}$, or
$\emptyset$ if $n=0$). We then have a natural map
\begin{equation}\label{sph.I.eq}
\sigma_i:\SS_i \to \Delta_i
\end{equation}
that is injective by Lemma~\ref{sk.le}, a generalization of the
sphere embedding \eqref{sph.eq}. Explicitly, by \eqref{sk.M}, we
have
\begin{equation}\label{sph.I.M}
  \SS_i = \colim_{i' \in M(i)^o} \Delta_{i'},
\end{equation}
and \eqref{sph.I.eq} is induced by the embedding $M(i) \subset I_M
\setminus i$. Moreover, for any $X \in I^o\Sets$ and $n \geq 0$,
\eqref{LM.eq} immediately shows that we have a cocartesian square
\begin{equation}\label{sk.I.sq}
\begin{CD}
\sk_{n-1}X @>>> \sk_nX\\
@AAA @AAA\\
\coprod_{i \in \M(IX)_n}\SS_{\pi(i)} @>{\copr \sigma_i}>>
\coprod_{i \in \M(IX)_n}\Delta_{\pi(i)},
\end{CD}
\end{equation}
where $\M(IX)_n \subset \M(IX)$ is the subset of objects of degree
$n$.

For a categorical interpretation of skeleta, denote by $\sk_nI
\subset I$, $n \geq -1$ the full subcategory spanned by objects $i
\in I$ with $\deg(\lambda(i)) \leq n$ (if $I=J$ is a left-bounded
partially ordered set with the cellular Reedy structure of
Lemma~\ref{ind.reedy.le}, then $\sk_nJ$ is the same as in
\eqref{sk.J}). Then for any $X \in I^o\Sets$, we have $I\sk_nX \cong
\sk_nIX$. In particular, the embedding $\sk_nI \to I$ is a discrete
fibration, so that $\sk_nI \subset I$ is left-closed in the sense of
Definition~\ref{cl.def}. Therefore we have the characteristic
functor
\begin{equation}\label{sk.chi.eq}
\chi_n:I \to [1]
\end{equation}
sending $\sk_nI \subset I$ to $0 \in [1]$ and everything else to $1
\in [1]$. Moreover, $\chi_n(m)=\id$ for any map $m:i \to i'$ in $M$,
so that if we equip $[1]$ with the Reedy structure of
Example~\ref{reedy.exa}, then $\chi_n$ is a Reedy functor. Just as
the squares \eqref{sk.I.sq}, the functors \eqref{sk.chi.eq} can be
used to set up an induction on skeleta and study $I$ by this
induction. This works especially well when $I$ is finite-dimensional
in the following sense.

\begin{defn}\label{I.dim.def}
A cellular Reedy category $I$ is {\em finite-dimensional} if we have
$I=\sk_nI$ for some $n \geq 0$, and the {\em dimension} $\dim I$ is
the smallest such integer $n$. For any cellular Reedy category $I$,
a functor $X \in I^o\Sets$ is {\em finite-dimensional} if so is
$IX$, and $\dim X = \dim IX$.
\end{defn}

\begin{exa}
A left-bounded partially ordered set $I$ with the cellular Reedy
structure of Example~\ref{reedy.exa} is finite-dimensional iff it
has finite chain dimension in the sense of
Definition~\ref{dim.def}, and $\dim I$ of Definition~\ref{dim.def}
coincides with $\dim I$ of Definition~\ref{I.dim.def}
\end{exa}

For any cellular Reedy category $I$, we denote by $I^o_f\Sets
\subset I^o\Sets$ the full subcategory spanned by finite-dimensional
functors $X$. If all functors are finite-dimensional then so is $I$,
but the converse is not necessarily true. For example, $\Delta$ is
finite-dimensional with $\dim \Delta = 0$, but $X \in \Delta^o\Sets$
is finite-dimensional iff it is finite-dimensional in the usual
sense (that is, the dimensions of non-degenerate simplices are
bounded by some $n \geq 0$). Even for $X \in I^o_f\Sets$, $\dim X$
is often strictly greater than $\dim I$.

\subsection{Anodyne resolutions.}\label{ano.res.subs}

As an application of the induction on ske\-leta, let us show that
cellular Reedy categories admit good ``resolutions'' in the
following sense.

\begin{defn}\label{ano.reso.def}
For any regular cardinal $\kappa$, a functor $J \to I$ is a {\em
  $\kappa$-anodyne resolution} of a small category $I$ if $J$ is a
partially ordered set, and for any $i \in I$, the right comma-fiber
$i \setminus J$ is $\kappa$-anodyne in the sense of
Definition~\ref{pos.S.ano.def}.
\end{defn}

\begin{exa}\label{xi.ano.exa}
For any regular cardinal $\kappa$ and $J \in \Posf_\kappa$, the
functor $\xi_\perp:B(J)^o \to J$ of \eqref{B.xi} is a
$\kappa$-anodyne resolution in the sense of
Definition~\ref{ano.reso.def}. Indeed, $j \setminus B(J)^o \cong B(j
\setminus J)^o$ for any $j \in J$, we have the embedding $\ppt \to j
\setminus J$ onto the smallest element $j$, and the induced
embedding $\ppt = B(\ppt)^o \to B(j \setminus J)^o$ is reflexive by
Lemma~\ref{B.le}.
\end{exa}

Fix a Reedy category $\langle I,M,L,\deg\rangle$ that is cellular in
the sense of Definition~\ref{ind.reedy.def}, and let $\Psi I =
\Delta^M_\perp I$ be the cosimplicial expansion of
\eqref{trans.eq}. Recall that objects in $\Psi I$ are diagrams
\eqref{i.dot}, and let $\Phi I \subset \Psi I$ be the full
subcategory spanned by diagrams such that $i_n \in I$ is
non-degenerate, and the maps $f_{d,d'}:i_d \to i_{d'}$, $0 \leq d
\leq d' \leq n$ are in $L$. Then \eqref{trans.eq} induces a
projection
\begin{equation}\label{pi.I.eq}
\Phi I \to \Delta^o,
\end{equation}
and on the other hand, composing \eqref{xi.perp.eq} with the
embedding $\rho:\Phi I \to \Psi I$ gives a functor
\begin{equation}\label{q.I.eq}
\xi^\Phi_\perp:\Phi I \to I.
\end{equation}
Any Reedy functor $\gamma:I' \to I$ between cellular Reedy
categories induces a functor $\Psi(\gamma):\Psi(I') \to \Psi(I)$
cocartesian over $\Delta^o$ that commutes with the projections
\eqref{xi.perp.eq}.

\begin{lemma}\label{Phi.le}
For any cellular Reedy category $I$, with the corresponding
categories $\Phi I \subset \Psi I$, the embedding $\rho:\Phi I
\subset \Psi I$ is left-reflexive over $\Delta^o$, with some adjoint
functor $\lambda:\Psi I \to \Phi I$, and the functor \eqref{pi.I.eq}
is a discrete cofibration. For any Reedy functor $\gamma:I' \to I$
between cellular Reedy categories, there exists a unique functor
$\Phi(\gamma):\Phi(I') \to \Phi(I)$ cocartesian over $\Delta^o$ such
that $\Phi(\gamma) \circ \lambda \cong \Psi(\gamma) \circ \lambda$,
and the isomorphism $\xi_\perp \circ \Psi(\gamma) \cong \gamma \circ
\xi_\perp$ induces a natural map
\begin{equation}\label{xi.Phi.map}
\gamma \circ \xi^\Phi_\perp \to \xi^\Phi_\perp \circ \Phi(\gamma)
\end{equation}
that is pointwise in the class $M$.
\end{lemma}

\proof{} For any $[n] \in \Delta$, the fiber $\Psi_{[n]}$ of the
cofibration $\pi:\Psi I \to \Delta^o$ is the category of diagrams
\eqref{i.dot} and maps between them that are pointwise in $M$. If
$[n]=0$, then $\Psi I_{[0]} = I_M$, and $\Phi I_{[0]} = \M(I)_M
\subset \Psi I_{[0]}$ is a discrete full left-reflexive subcategory
by Lemma~\ref{I.M.le}, with the functor $\lambda_{[0]}:\Psi I_{[0]}
\to \Phi I_{[0]}$ left-adjoint to the $\rho_{[0]}$, and the
adjunction map $a:\id \to \lambda$. For $n \geq 1$, note that the
maps $a \circ f_{d,n}:i_d \to \lambda(i_n)$, $0 \leq d \leq n$ admit
unique factorizations
\begin{equation}\label{l.m.facto}
\begin{CD}
i_d @>{m}>> i'_d @>{l}>> \lambda(i_n)
\end{CD}
\end{equation}
with $m \in M$, $l \in L$, and the middle terms $i'_d$ of these
factorizations form a diagram $i'_\idot$ that lies in $\Phi I_{[n]}
\subset \Psi I_{[n]}$. Thus the full subcategory $\Phi I_{[n]}
\subset \Psi I_{[n]}$ is also left-reflexive, with the adjoint
$\lambda_{[n]}:\Psi I_{[n]} \to \Phi I_{[n]}$ sending $i_\idot$ to
$i'_\idot$. Moreover, by the uniqueness of the factorizations
\eqref{l.m.facto}, the category $\Phi I_{[n]}$ is discrete. The
existence of the adjoint functors $\lambda_{[n]}$ then
tautologically implies that $\pi:\Phi I \to \Delta^o$ is a
precofibration, with the transition functor $f^\Phi_!$ for a map
$f:[n'] \to [n]$ given by $\lambda_{[n]} \circ f^\Psi_!$, where
$f^\Psi_!$ is the transition functor for the cofibration $\pi:\Psi I
\to \Delta^o$. Since the fibers of this precofibration are discrete,
it is a cofibration, and the functors $\lambda_{[n]}$ together
define a functor $\lambda:\Psi I \to \Phi I$ adjoint to the
embedding $\Phi I \to \Psi I$. Finally, if we have another cellular
Reedy category $I'$ and a Reedy functor $\gamma:I' \to I$, then for
any $i \in \Psi I'$, the functor $\lambda \circ \Psi(\gamma)$
inverts the adjunction map $i \to \lambda(i)$. Thus it indeed
factors as $\lambda \circ \Psi(\gamma) = \Phi(\gamma) \circ \lambda$
for a unique $\Phi(\gamma):\Phi I' \to \Phi I$ cocartesian over
$\Delta^o$, and by adjunction, we have the base change map
$\Psi(\gamma) \circ \rho \to \rho \circ \Phi(\gamma)$. To obtain
\eqref{xi.Phi.map}, apply $\xi^\Phi_\perp$ to this map.  \endproof

\begin{lemma}\label{Phi.red.le}
For a cellular Reedy category $I$,
$\xi^\Phi_\perp$ of \eqref{q.I.eq} factors as
$$
\begin{CD}
\Phi I @>{R}>> \Red(\Phi I) @>{q_I}>> I,
\end{CD}
$$
where $R$ is the reduction functor \eqref{red.eq}.
\end{lemma}

\proof{} By definition, a map $\langle [n'],i'_\idot \rangle \to
\langle [n],i_\idot \rangle$ in the category $\Phi I$ is given by a
map $g:[n] \to [n']$ and a collection of maps $g_d:i'_{g(d)} \to
i_d$ in $M$, $0 \leq d \leq n$, such that for any $0 \leq d \leq d'
\leq n$, we have $g_{d'} \circ f'_{g(d),g(d')} = f_{d,d'} \circ
g_d$. The functor $\xi^\Phi_\perp$ of \eqref{q.I.eq} sends $\langle
g,g_\idot \rangle$ to $g_0 \circ f'_{0,g(0)}:i'_0 \to i_0$. By
definition, $i_n$ is non-degenerate, so by Lemma~\ref{I.M.le}, for
any $i \in I$, there exists at most one map $m:i \to i_n$ in $M$,
and if it exists, then $i_n = \lambda(i)$, $m$ is the adjunction
map; moreover, in this case, for any other map $f = l' \circ m':i
\to i_n = \lambda(i)$, there must exist $m'' \in M$ such that $m''
\circ m' = m$, and then $m''=l'=\id$ for degree reasons, so that
$f=m$. So, for any map $\langle g',g'_\idot \rangle:\langle
      [n'],i'_\idot \rangle \to \langle [n],i_\idot \rangle$, say
      with $g'(n) \leq g(n)$, we must have $g'_n = g_n \circ
      f'_{g'(n),g(n)}$, so that
$$
f_{0,n} \circ \xi^\Phi_\perp(g) = g_n \circ f'_{0,g(n)} = g'_n
\circ f'_{0,g'(n)} = f_{0,n} \circ \xi^\Phi_\perp(g').
$$
Since $f_{0,n} \in L$ is injective, this implies that
$\xi^\Phi_\perp(g)=\xi^\Phi_\perp(g')$.
\endproof

If we denote by $Q(I) = \Red(\Phi I)$ the reduction of the category
$\Phi I$, then Lemma~\ref{Phi.red.le} provides a functor
\begin{equation}\label{Q.I.eq}
q_I:Q(I) \to I.
\end{equation}
By the universal property of reductions, for any Reedy functor
$\gamma:I' \to I$, the functor $\Phi(\gamma)$ of Lemma~\ref{Phi.le}
induces a functor $Q(\gamma):Q(I') \to Q(I)$, and the map
\eqref{xi.Phi.map} then reduces to a map
\begin{equation}\label{q.I.map}
\gamma \circ q_I \to q_{I'} \circ Q(\gamma)
\end{equation}
that is pointwise in the class $M$.  By Lemma~\ref{Phi.le}, $\Phi I
\cong (\Delta X(I))^o$ is opposite to the category of simplices of a
simplicial set $X(I) \in \Delta^o\Sets$, and $Q(I)$ is then opposite
to the set of non-degenerate simplices of $X(I)$ ordered by
inclusion.

\begin{exa}\label{q.I.exa}
If $I \in \Posf_\kappa$ is a partially ordered set, then $Q(I)$ is
the opposite $B(I)^o$ to the barycentric subdivion $B(I)$ of
Definition~\ref{pos.bary.def}, and $q_I$ is the $\kappa$-anodyne
resolution $\xi_\perp$ of Example~\ref{xi.ano.exa}.
\end{exa}

\begin{exa}\label{reg.I.exa}
More generally, if a cellular Reedy category $I$, $|I| < \kappa$ is
regular in the sense of Definition~\ref{reg.I.def}, then $\Phi I
\cong (\Delta \M(I))^o$, and since $\M(I)$ is directed, it is
ordered by Remark~\ref{BB.rem}. Therefore the simplicial set $X(I)$
is regular, we have
\begin{equation}\label{Q.reg}
Q(I) \cong (\overline{\Delta}\M(I))^o,
\end{equation}
and for any $i \in I$, we have $i \setminus Q(I) \cong
\overline{\Delta}(i \setminus \M(I))^o$. If $\dim I < \infty$, then
$\dim \overline{\Delta}(i \setminus \M(I))^o < \infty$ for any $i$,
and since the category $i \setminus \M(I)$ has an initial object $i
\to \lambda(i)$, $\overline{\Delta}(i \setminus \M(I))^o$ is
$\kappa$-anodyne by Lemma~\ref{BB.le}. Therefore $q_I$ is again a
$\kappa$-anodyne resolution.
\end{exa}

For a general cellular Reedy category $I$, the partially ordered set
$Q(I)$ still shares some of the properties of barycentric
subdivision.

\begin{lemma}\label{reedy.copr.le}
For any Reedy functor $\chi:I \to [1]$, the induced functor
$Q(\chi):Q(I) \to Q([1]) \cong B([1])^o \cong \V$ is a cofibration,
and the functors $Q(I_0) \to Q(I)$ and $Q(I_1) \to Q(I)$ induced by
the embeddings $I_0,I_1 \to I$ of Lemma~\ref{I.corr} are
right-closed embeddings that induce isomorphisms $Q(I_0) \cong
Q(I)_0$, $Q(I_1) \cong Q(I)_1$ with the fibers of the cofibration
$Q(\chi)$. Moreover, for any Reedy functor $\chi:I \to \V$, we have
$Q(I) \cong Q(I/0) \copr_{Q(I_o)} Q(I/1)$.
\end{lemma}

\proof{} For the first claim, let $p:\Delta [1] \to B([1]) \cong
\V^o$ be the fibration \eqref{red.del.eq}. Then $\Phi(\chi)$ is a
cocartesian functor between categories discretely cofibered over
$\Delta^o$, thus a cofibration, and $Q(\chi):Q(I) \to \V$ is the
reduction of the composition $p \circ \Phi(\chi)$. Therefore it is
also a cofibration. The rest of the argument is exactly the same as
in Example~\ref{B.cl.exa}. Moreover, if we compose the cofibration
$Q(\chi)$ with the projection $\xi_\perp:\V \cong B([1])^o \to [1]$,
then the fiber $Q(I)_0$ of the resulting map $Q(I) \to [1]$ is still
identified with $Q(I_0)$. Then for the second claim, as in the proof
of Lemma~\ref{B.corr}, compose the functor $Q(\chi):Q(I) \to Q(\V)
\cong B(\V)^o$ with the projection $\xi_\perp:B(\V)^o \to \V$, and
use this identification to see that $o \setminus Q(I) \cong Q(I_o)$,
$0 \setminus Q(I) \cong Q(I/0)$ and $1 \setminus Q(I) \cong Q(I/1)$.
\endproof

\begin{lemma}\label{Q.pushout.le}
Assume given a Reedy functor $\chi:I \to [1]$ and a discrete
fibration $\pi:I' \to I$ that is an isomorphism over $1 \in
[1]$. Then the induced map $Q(\pi):Q(I') \to Q(I)$ is an isomorphism
over $o,1 \in \V \cong Q([1])$.
\end{lemma}

\proof{} Since $I' \to I$ is a discrete fibration, any diagram
\eqref{i.dot} in $I_L$ that has non-empty intersection with $I_1
\subset I$ lifts uniquely to $I'_L$.
\endproof

\begin{lemma}\label{reedy.N.le}
Assume given a Reedy category $I$ equipped with a Reedy functor $I
\to \N$. Then $Q(I) \cong \colim_nQ(I/n)$.
\end{lemma}

\proof{} For any $n \geq 0$, the map $Q(I/n) \to Q(I)$ is a fully
faithful embedding by Lemma~\ref{reedy.copr.le}, so by
Example~\ref{N.coli.exa}, it suffices to check that the filtration
$Q(I/n) \subset Q(I)$, $n \geq 0$ is exhaustive. But any element in
$Q(I)$ is represented by a functor $[m] \to I$, and any such functor
factors through $I/n \subset I$ for a sufficiently large $n$.
\endproof

\begin{prop}\label{ano.I.prop}
Assume that the cellular Reedy category $I$ is
finite\--di\-men\-si\-onal in the sense of
Definition~\ref{I.dim.def}, and $|I| < \kappa$ for a regular
cardinal $\kappa$. Then $Q(I)$ has finite chain dimension $\dim Q(I)
\leq \dim I$, and the functor \eqref{Q.I.eq} is a $\kappa$-anodyne
resolution in the sense of Definition~\ref{ano.reso.def}.
\end{prop}

\proof{} The fact that $\dim Q(I) \leq \dim I$ is obvious from the
definition. To show that \eqref{Q.I.eq} is a $\kappa$-anodyne
resolution, note that if $I$ is regular, this is
Example~\ref{reg.I.exa}. In the general case, let $n = \dim I$, and
assume by induction that the claim is proved for categories of
dimension $\leq n-1$. Take an object $i \in I$. If
$\deg(\lambda(i))=n$, then $i \setminus Q(I)$ consists of a single
object $i \to \lambda(i)$, so it is anodyne. Otherwise we have $i
\in \sk_{n-1}I$. Consider the Reedy functor $\chi_{n-1}$ of
\eqref{sk.chi.eq} and the corresponding cofibration
$Q(\chi_{n-1})$. Compose it with the cofibration $i \setminus Q(I)
\to Q(I)$ of Example~\ref{comma.exa} to obtain a cofibration $i
\setminus Q(I) \to \V$. Moreover, let
$$
I' = \coprod_{i \in \M(I)_n} I/i,
$$
with the natural discrete fibration $\pi:I' \to I$, and consider the
induced functor
$$
\pi^Q:i \setminus Q(I') = \coprod_{i' \in I'_i}i' \setminus Q(I')
\to i \setminus Q(I).
$$
Then $\pi$ is a functor over $[1]$ with respect to the projection
$\chi_{n-1}:I \to [1]$ and an isomorphism over $1$, so by
Lemma~\ref{Q.pushout.le}, $\pi^Q$ is an isomorphism over $o,1 \in
\V$. Moreover, since $I'$ is regular by Corollary~\ref{reg.I.corr},
the embedding $i' \setminus Q(I'_{\leq n-1}) \to i' \setminus Q(I')$
is a map between $\kappa$-anodyne partially ordered sets for any $i'
\in I'_i$, thus $\kappa$-anodyne. Therefore the coproduct of these
embeddings is also $\kappa$-anodyne, and then by
Corollary~\ref{V.ano.corr}, the embedding $i \setminus Q(I_{\leq
  n-1}) \to i \setminus Q(I)$ is $\kappa$-anodyne, so we are done by
induction.
\endproof

\subsection{Barycentric subdivision.}

Let $\Delta X$ be the category of simplices of a simplicial set
$X:\Delta^o \to \Sets$. By Lemma~\ref{ind.reedy.le}, it is a
cellular Reedy category in the sense of
Definition~\ref{ind.reedy.def}. The set $\M(\Delta X)_M$ is the set
of non-degenerate simplices $\langle [n],x \rangle \in \Delta X$,
with $\M(\Delta I)_n$ consisting of simplices of degree $n$. The
functor $\lambda$ of Lemma~\ref{I.M.le} sends a simplex $\langle
[n],x \rangle$ to its normalization. The category $\Delta X$ is
regular in the sense of Lemma~\ref{reg.I.le} iff the simplicial set
$X$ is weakly regular in the sense of
Subsection~\ref{del.0.subs}. The $n$-th skeleton $\sk_n X$ of
\eqref{sk.eq} is the usual skeleton of a simplicial set, $\SS_{[n]}$
of \eqref{sph.I.eq} is the $(n-1)$-dimensional simplicial sphere
$\SS_{n-1}$ of \eqref{sph.eq}, and the gluing squares
\eqref{sk.I.sq} are obtained by iterating the squares
\eqref{elem.X.sq}. The cellular Reedy category $\Delta X$ has finite
dimension iff so does the set $X$, and the dimensions are the
same. Moreover, it turns out that similarly to
Example~\ref{q.I.exa}, the canonical anodyne resolution
\eqref{Q.I.eq} of Proposition~\ref{ano.I.prop} can be expressed in
terms of a barycentric subdivision.

\begin{defn}\label{bary.X.def}
  The {\em barycentric subdivision functor} $B$ from $\Delta^o\Sets$
  to itself is the left Kan extension of the functor $\Delta \to
  \Delta^o\Sets$, $[n] \mapsto N(B([n]))$ with respect to the Yoneda
  embedding $\Y:\Delta \to \Delta^o\Sets$.
\end{defn}

Definition~\ref{bary.X.def} is the standard definition of the
barycentric subdivision operation on simplicial sets. If $X = N(J)$
is the nerve of a partially ordered set $J$, then $B(N(J)) \cong
N(B(J))$, so that Definition~\ref{pos.bary.def} is consistent with
Definition~\ref{bary.X.def}.

\begin{lemma}\label{bary.X.le}
For any simplicial set $X \in \Delta^o\Sets$, we have an equivalence
$\Phi\Delta X \cong \Delta B(X)^o$, and this equivalence is
functorial in $X$.
\end{lemma}

\proof{} By Lemma~\ref{Phi.le}, for any $X \in \Delta^o\Sets$, we
have $\Phi\Delta X \cong \Delta B'(X)^o$ for some simplicial set
$B'(X)$, and this correspondence is functorial in $X$, so that we
have a functor $B':\Delta^o\Sets \to \Delta^o\Sets$. If $X$ is
regular, we have $\Phi \Delta X \cong (\Delta\overline{\Delta}X)^o$
and $B'(X) \cong N(\overline{\Delta}X)$. In particular, $\Delta_n
\subset \Delta^o\Sets$ is regular for any $n$, and
$\overline{\Delta}\Delta_n = \overline{\Delta}N([n]) = N(B([n]))$,
so that we have an isomorphism $B \cong B'$ on the image of
$\Y:\Delta \to \Delta^o\Sets$. By the universal property of the Kan
extension, it extends to a map of functors $B \to B'$, and as in
Example~\ref{yo.sets.exa}, to prove that it is an isomorphism, it
suffices to prove that $B'$ commutes with colimits. Moreover, $B'$
obviously commutes with filtered colimits, so by induction on
skeleta, it suffices to prove that it sends a square \eqref{sk.I.sq}
to a cocartesian square.

Indeed, by definition, for any $[m] \in \Delta$, we have
\begin{equation}\label{B.pr}
B'X([m]) = \coprod_{[n_\idot]}\overline{X}(\xi([n_\idot])),
\end{equation}
where the coproduct is over all functors $[n_\idot]:[m] \to
\Delta_\ff = \Delta_L \subset \Delta$, $\xi:\Delta^L\Delta \to
\Delta$ is the projection \eqref{xi.eq}, so that $\xi([n_\idot]) =
[n_m]$, and $\overline{X}([n]) \subset X([n])$, $[n] \in \Delta$ is
the subset of non-degenerate simplices. While the correspondence $X
\mapsto \overline{X}([n])$ is not functorial in $X$, it is
functorial with respect to injective maps $X \to X'$. Moreover, for
any injective map $X \to X'$ and any map $f:X \to Y$, we have an
injective map $Y \to Y' = X' \copr_X Y$, and for any $[n]$, $f$
induces an isomorphism
$$
\overline{X}'([n]) \ssetminus \overline{X}([n]) \cong
\overline{Y}'([n]) \ssetminus \overline{Y}([n]).
$$
Then \eqref{B.pr} shows that the natural map $B'X' \copr_{B'X} B'Y
\to B'Y'$ is an isomorphism, and this finishes the proof.
\endproof

As a corollary of Lemma~\ref{bary.X.le}, we see that the partially
ordered set $Q(\Delta X)$ is naturally opposite to the set of
non-degenerate simplices in the barycentric subdivision $B(X)$ of
$X$, with the order given by inclusion. If $X$ is regular, then
$B(X)$ is also regular, and we then have
\begin{equation}\label{Q.del.X}
Q(\Delta X) \cong (\overline{\Delta}BX)^o \cong B(\overline{\Delta}X)^o,
\end{equation}
functorially with respect to maps of regular simplicial sets.

\chapter{Localization.}\label{loc.ch}

The subject of this chapter is localization, in the original naive
sense of \cite{GZ}. We start with a bunch of standard and easy
observations, and some equally standard examples. Our notation for
the localization of a category $I$ in a class $W$ of morphisms is
$h^W(I)$. One application of the general localization machinery is a
left Kan extension operation $\gamma_!$ for families of groupoids, a
counterpart of the right Kan extension $\gamma_*$ constructed in
Section~\ref{fun.fib.sec}. All this material is in
Section~\ref{loc.sec}, and Section~\ref{seg.loc.sec} is devoted to
categories fibered over $\Delta$. Here the main result is
Lemma~\ref{del.prop} saying that, roughly speaking, a small category
is the localization of its simplicial explansion in the class of
special maps (just as Proposition~\ref{spec.prop}, this is again a
variation on the theme of \cite{DKHS}).

Another result that we will really need is
Proposition~\ref{exa.seg.prop} that characterizes families of
groupoids over $\Delta$ that appear as simplicial replacements; this
uses a toy version of a beautiful idea of \cite{rzk}. We complement
this material with the standard results on the localizations of the
categories of simplices of fibrant simplicial sets (this actually
goes back all the way to \cite{GZ}).

Section~\ref{mod.sec} introduces Quillen's model categories and
collects all we need from that theory (as in Chapter~\ref{simp.ch},
this is not that much). We omit the actual definition of a model
category and the proofs of the two main results, Quillen Adjunction
and the existence of the Reedy model structures. The main result
that we do prove is Proposition~\ref{mod.prop} saying that, roughly
speaking, a family of groupoids over a model category that is
constant along weak equivalences is completely determined by its
restriction to the subcategory of cofibrant objects and
cofibrations.

An application of this is the definition of a representable family
of group\-oids over a model category, and a construction of such
families; this is Definition~\ref{repr.def} and
Proposition~\ref{repr.prop}. It is exactly these representable
families that go into our version of Brown Representability proved
in Chapter~\ref{semi.ch}.

\section{Abstract localization.}\label{loc.sec}

\subsection{Relative categories.}\label{loc.subs}

It is convenient to start the discussion of localization with the
following notion introduced in \cite{BK}.

\begin{defn}
A {\em relative category} is a pair $\langle I,W \rangle$ of a
category $I$ and a class $W$ of morphisms in $I$ that contains all
isomorphisms.
\end{defn}

Letter $W$ here is traditional and stands for ``weak equivalence'',
but as far as the definition is concerned, it can be any class of
maps whatsoever, provided it contains all isomorphisms. The biggest
possible class is the class $\Tot$ of all morphisms in $I$, and
the smallest possible class is the class $\Iso$ of all
isomorphisms. As in the case of dense subcategories, for any functor
$\gamma:I' \to I$ to a relative category $\langle I,W \rangle$, we
denote by $\gamma^*W$ the class of morphisms $f$ in $I'$ such that
$\gamma(f) \in W$. If $I' \subset I$ is a subcategory and $\gamma$
is the embedding functor, we will shorten $\gamma^*W$ to $W$.  A
{\em morphism} from a relative category $\langle I,W \rangle$ to a
relative category $\langle I',W' \rangle$ is a functor $\gamma:I \to
I'$ such that $W \subset \gamma^*W'$.

\begin{exa}
A morphism $\langle I,W \rangle \to \langle I',\Iso \rangle$ is the
same thing as a functor $\gamma:I \to I'$ that inverts all maps $f
\in W$ (recall that if $I$ is essentially small, these form the
category $\Fun^W(I,I')$ of Example~\ref{fun.v.i.exa}).
\end{exa}

\begin{defn}
The {\em localization} of a relative category $\langle I,W \rangle$
is a category $h^W(I)$ equipped with a morphism $h:\langle I,W
\rangle \to \langle h^W(I),\Iso \rangle$ such that for any other
category $I'$ and a morphism $\gamma:\langle I,W \rangle \to \langle
I',\Iso \rangle$, there exists a functor $\gamma':h^W(I) \to I'$ and
an isomorphism $\alpha:\gamma' \circ h \cong \gamma$, and the pair
$\langle \gamma',\alpha \rangle$ is unique up to a unique
isomorphism.
\end{defn}

Alternatively, localization can be described by means of cocartesian
squares of categories of Subsection~\ref{dia.subs}: for any relative
category $\langle I,W \rangle$, the functors \eqref{eps.f} for all
$w \in W$ define a functor $W \times [1] \to I$, and the
localization, if it exists, fits into a cocartesian square
\begin{equation}\label{loc.sq}
\begin{CD}
W \times [1] @>>> I\\
@VVV @VV{h}V\\
W @>>> h^W(I),
\end{CD}
\end{equation}
where $W$ in the bottom left corner is understood as a discrete
category. In particular, the localization $h^W(I)$ is unique, if it
exists. If it does, we will say that the relative category $\langle
I,W \rangle$ is {\em localizable}, or that $I$ is {\em localizable
  with respect to $W$}. For any category $I$, we will say that $I$
is localizable if so is the relative category $\langle I,{\Tot}
\rangle$, and we will say that $I$ is {\em $1$-connected} if it is
localizable and $h^{\Tot}(I) \cong \ppt$. Being localizable is a
set-theoretic condition: one can always define $h^W(I)$ as the
category with the same objects as $I$ and morphisms given by formal
expressions
\begin{equation}\label{zig}
w_0^{-1} \circ f_1 \circ w_1^{-1} \circ \dots \circ w_{n-1}^{-1}
\circ f_n \circ w_n^{-1},
\end{equation}
with $w_\idot \in W$ and $f_\idot$ morphisms in $I$, modulo obvious
cancellation rules, and the only thing to check is that these
expressions form a set and not a proper class. This is automatically
the case if $I$ is small or at least essentially small. Moreover, if
$\|I\| \leq \kappa$ for a regular cardinal $\kappa$, then
$\|h^W(I)\| \leq \kappa$ for any $W$.

\begin{remark}\label{loc.cart.rem}
While localization is only a particular example of a cocartesian
square of categories, arbitrary such squares can in fact be
constructed by localization. Namely, it is immediate from the
definition that a cofibration $\pi:\C \to \V$ extends to a
cofibration $\C' \to [1]^2 = \V^>$ corresponding to a cocartesian
square if and only if $\C$ is localizable with respect to the class
$c=\pi^\sharp\Tot$ of maps cocartesian over $\V$, and in this case,
the fiber $\C'_o$ over the maximal element $o \in \V^>$ is the
localization $h^c(\C)$.
\end{remark}

\begin{remark}\label{pos.col.rem}
In particular, Remark~\ref{loc.cart.rem} shows that the category
$\Pos$ has all cocartesian squares (hence all colimits, since it has
all coproducts). Indeed, all partially ordered sets are by
definition small, thus localizable, so any functor $\V \to \Pos$
extends to a cocartesian square of small categories, and its
reduction is then a cocartesian square in $\Pos$, as in
Example~\ref{pos.perf.exa}.
\end{remark}

The most fundamental example of a localizable large category is the
category $\Cat$ of all small categories and functors between them,
with $W$ being the class of all equivalences.

\begin{lemma}\label{cat.loc.le}
The relative category $\langle \Cat,W \rangle$ is localizable, and
its localization $h^W(\Cat)$ is naturally equivalent to the category
$\Cat^0$ of small categories, with morphisms given by isomorphism
classes of functors.
\end{lemma}

\begin{remark}
If one promotes $\Cat$ to a $2$-category $\cCat$, then $\Cat^0 =
\overline{\cCat}$ is its truncation in the sense of
Example~\ref{2.cat.fam.exa}.
\end{remark}

\proof{} We have the obvious functor $h:\Cat \to \Cat^0$, and by
definition, it inverts equivalences, so that we have a morphism
$h:\langle \Cat,W \rangle \to \langle \Cat^0,\Iso \rangle$. If we
have two functors $\gamma_0,\gamma_1:I \to I'$ between small
categories and an isomorphism $\gamma_0 \cong \gamma_1$, then
these data define a functor $\gamma:e(\{0,1\}) \times I \to I'$,
where as in Subsection~\ref{pos.subs}, $e(\{0,1\})$ is the category with
two objects $0$, $1$ and exactly one map between any two
objects. The two embeddings $\eps(0),\eps(1):\ppt \to e(\{0,1\})$
have a common one-sided inverse given by the projection $e(\{0,1\})
\to \ppt$, and then $\eps(0) \times \id, \eps(1) \times \id:I \to
e(\{0,1\}) \times I$ are both inverse to the projection $e(\{0,1\})
\times I \to I$. Therefore for any morphism $f:\langle \Cat,W
\rangle \to \langle \C,\Iso \rangle$ to some $\C$, we have
$f(\eps(0) \times \id) = f(\eps(1) \times \id)$, and since $\gamma_l
= \gamma \circ (\eps(0) \times \id)$, $l = 0,1$, this implies
$f(\gamma_0)=f(\gamma_1)$, so that $f$ uniquely factors through $h$.
\endproof

For any regular cardinal $\kappa$, we have the similar statement for
the full subcategory $\Cat_\kappa \subset \Cat$ spanned by
$\kappa$-bounded small categories $I$ -- namely, $\langle
\Cat_\kappa,W \rangle$ is localizable, and $h^W(\Cat_\kappa) \cong
\Cat^0_\kappa \subset \Cat^0$ is the full subcategory spanned by $I
\in \Cat_\kappa$. More generally, for any category $I$, we have the
categories $\Cat \bbi I \subset \Cat \bb I$ of small categories and
functors resp.\ lax functor over $I$, as in
Subsection~\ref{fun.subs}. If we again let $W$ be the class of
equivalences, then both $\langle \Cat \bbi I,W \rangle$ and $\langle \Cat
\bb I,W \rangle$ are localizable, and the localizations are
$h^W(\Cat \bbi I) \cong (\Cat \bbi I)^0$, $h^W(\Cat \bb I) \cong (\Cat \bb
I)^0$ where $(\Cat \bbi I)^0$ resp.\ $(\Cat \bb I)^0$ is the category of
small categories $I' \in \Cat$ equipped with a functor $\pi:I' \to
I$, and with morphisms from $\pi_0:I'_0 \to I$ to $\pi_1:I_1' \to I$
given by functors resp.\ lax functors $\gamma:I'_0 \to I'_1$ over
$I$ modulo isomorphisms over $I$. If $I$ is small, thus corresponds
to an object $I \in \Cat^0$, then we have the comma-fiber
$\Cat^0/[I]$, and the tautological functor
\begin{equation}\label{cat.I.ev}
  (\Cat \bbi I)^0 \to \Cat^0/[I]
\end{equation}
is an epivalence.

\begin{exa}
The epivalence \eqref{cat.I.ev} is usually {\em not} an equivalence;
the reason for this is, maps in $\Cat^0/[I]$ are functors over $I$
modulo all isomorphisms, while maps in $(\Cat \bbi I)^0$ are functors
over $I$ modulo isomorphisms over $I$. For example, assume given a
group extension
\begin{equation}\label{gr.ex.dia}
\begin{CD}
  1 @>>> H @>>> \wt{G} @>{p}>> G @>>> 1
\end{CD}
\end{equation}
with abelian $H$, let $I=\ppt_G$, and let $\C = \ppt_H$,
$\C'=\ppt_{\wt{G}}$. Then $p$ in \eqref{gr.ex.dia} induces a functor
$\C' \to I$ making $\C'$ a category over $I$, and if we treat $\C$
as a category over $I$ via the tautological functor $\C \to \ppt \to
I$, then $(\Cat \bbi I)^0(\C,\C') = \Hom(H,H)$ is the set of group maps
$H \to H$, while $(\Cat^0/[I])(\C,\C')$ is the same set but modulo
conjugation by elements $g \in G$. If the adjoint $G$-action on $H$
is non-trivial, the latter set is strictly smaller.
\end{exa}

For any regular cardinal $\kappa$, the full subcategories $(\Cat
\bbi I)_\kappa \subset \Cat \bbi I$, $(\Cat \bb I)_\kappa \subset
\Cat \bb I$ spanned by $I' \in \Cat_\kappa$ are also localizable,
and the localizations are the full subcategories $h^W(\Cat_\kappa
\bbi I) \cong (\Cat_\kappa \bbi I)^0 \subset (\Cat/I)^0$,
$h^W(\Cat_\kappa \bb I) \cong (\Cat_\kappa \bb I)^0 \subset (\Cat
\bb I)^0$ spanned by $I' \in \Cat_\kappa$. Moreover, if $I$ is
small, then the subcategory $\Cat \mm I \subset \Cat \bbi I$ of
\eqref{fib.I} is also localizable with respect to the class $W$ of
equivalences over $I$, the localization $h^W(\Cat \mm I) \cong (\Cat
\mm I)^0 \subset (\Cat \bbi I)^0$ is the subcategory of fibrations
and isomorphism classes of cartesian functors, and for any regular
cardinal $\kappa$ with $|I| < \kappa$, the same holds for $(\Cat \mm
I)_\kappa = (\Cat \mm I) \cap (\Cat \bbi I)_\kappa$. The embeddings
$(\Cat \mm I)^0 \subset (\Cat \bbi I)^0$, $(\Cat \bbi I)^0 \subset
(\Cat \bb I)^0$ are left resp.\ right-reflexive, with the adjoint
functors in both cases sending $\pi:\C \to I$ to the comma-category
$\sigma:I \setminus_\pi \C \to I$. The extended Yoneda embedding
\eqref{yo.cat.eq} descends to a functor $(\C \bb I)^0 \to I^o\Cat^0$
that factors through a functor
\begin{equation}\label{cat.ev}
  (\Cat \mm I)^0 \to I^o\Cat^0
\end{equation}
provided by the Grothendieck construction (explicitly,
\eqref{cat.ev} corresponds to the functor $I^o \times (\Cat \mm I)^0
\to \Cat^0$, $\langle i,\C \rangle \mapsto \C_i$). The functor
\eqref{cat.ev} is obviously conservative, but usually not an
equivalence.

\begin{lemma}\label{cat.epi.le}
If $I = J$ is a partially ordered set with $\dim J \leq 1$, then
\eqref{cat.ev} is an epivalence.
\end{lemma}

\proof{} Since \eqref{cat.ev} is conservative, we just have to check
that it is essentially surjective and full. To lift a functor
$\C:J^o \to \Cat^0$ to an object in $(\Cat \mm I)^o$, just fix a
small category $\C(j) \in \Cat$ for any $j \in J$, and a functor
$\C(j) \to \C(j')$ representing the given map $\C(j) \to \C(j')$ in
$\Cat^0$, for any $j' < j$. Since $\dim J \leq 1$, this already
defines a functor $J^o \to \Cat$ that descends to an object in
$(\Cat \mm I)^0$. Finally, if we have two functors $\C,\C':J^o \to
\Cat$ and a morphism $g:\C \to \C'$ in $J^o\Cat^0$, then choosing a
representing functor $g(j):\C(j) \to \C'(j)$ for each $j \in J$, we
obtain a diagram
\begin{equation}\label{cat.sq.sq}
\begin{CD}
  \C(j) @>>> \C(j')\\
  @V{g(j)}VV @VV{g(j')}V\\
  \C'(j) @>>> \C'(j')
\end{CD}
\end{equation}
for any $j' < j$ in $J$ that commutes up to an isomorphism. To
lift $g$ to a map in $(\Cat \mm I)^0$, choose these isomorphisms.
\endproof

\begin{exa}
If $J=[1]$, then what \eqref{cat.ev} provides is a functor $(\Cat
\mm [1])^o \to \Ar(\Cat^0)$. It is a functor over $\Cat^0$ with
respect to projections $s^*:(\Cat \mm [1])^0 \to \Cat^0$,
$\tau:\Ar(\Cat^0) \to \Cat^0$, and its fiber over some $I \in
\Cat^0$ is the epivalence \eqref{cat.I.ev}.
\end{exa}

\begin{lemma}\label{adj.loc.le}
Assume given relative categories $\langle I_0,W_0 \rangle$ and
$\langle I_1,W_1 \rangle$, morphisms $\gamma_0:\langle I_0,W_0
\rangle \to \langle I_1,W_1 \rangle$ and $\gamma_1:\langle I_1,W_1
\rangle \to \langle I_0,W_0 \rangle$, and maps $a_0:\id \to \gamma_1
\circ \gamma_0$, $a_1:\gamma_0 \circ \gamma_1 \to \id$ such that
$a_0(i) \in W_0$ for any $i \in I_0$ and $a_1(i) \in W_1$ for any $i
\in I_1$.  Then $\langle I_0,W_0 \rangle$ is localizable if and only
if so is $\langle I_1,W_1 \rangle$, and in this case, $\gamma_0$ and
$\gamma_1$ induce mutually inverse equivalences between
$h^{W_0}(I_0)$ and $h^{W_1}(I_1)$.
\end{lemma}

\proof{} Clear. \endproof

\begin{lemma}\label{rel.loc.le}
Assume given a relative category $\langle I',W \rangle$ equipped
with a fibration $\pi:I' \to I$ such that $W \subset \pi^*(\Iso)$, the
fiber $\langle I'_i,W \rangle$ is localizable for any $i \in I$, and
the transition functor $f^*:\langle I'_{i'},W \rangle \to \langle
I'_i,W \rangle$ is a morphism for any map $f:i \to i'$ in $I$. Then
$\langle I',W \rangle$ is localizable, and the functor
$h(\pi):h^W(I') \to I = h^{\Iso}(I)$ is a fibration with fibers
$h^W(I')_i \cong h^W(I'_i)$. Moreover, if we have a functor
$\gamma:I_0 \to I$ from some category $I$, with the induced
fibration $I'_0 = \gamma^*I' \to I_0$, then $\langle I'_0,\gamma^*W
\rangle$ is localizable, and $h^{\gamma^*W}(I'_0) \cong
\gamma^*h^W(I')$.
\end{lemma}

\proof{} Clear. \endproof

As an application of Lemma~\ref{rel.loc.le}, assume given a functor
$\gamma:I' \to I$ between two categories $I$, $I'$, and consider the
right comma-category $I \setminus I'$ with its projections $\sigma:I
\setminus I' \to I$, $\tau:I \setminus I' \to I'$ and right
comma-fibers $i \setminus I'$. Say that a family of groupoids $\C'$
over $I'$ is {\em localizable with respect to $\gamma$} if for any
$i \in I$, the pullback $\tau(i)^*\C'$ with respect to the functor
\eqref{p.i} is a localizable category. Then by
Lemma~\ref{rel.loc.le}, the relative category $\langle
\tau^*\C',\sigma^*(\Iso) \rangle$ is localizable, and the natural
projection $\gamma_!\C' = h^{\sigma^*(\Iso)}(\tau^*\C') \to I$ is a
family of groupoids. Moreover, by the universal property of
localization, the operation $\gamma_!$ on families of groupoids is
left-adjoint to the pullback operation $\gamma^*$ in the same sense
as the pushforward operation $\gamma_*$ of
Subsection~\ref{2kan.subs} is right-adjoint to $\gamma^*$. In
particular, for any families of groupoids $\C \to I$, $\C' \to I$
such that $\C'$ and $\gamma^*\C$ are localizable with respect to
$\gamma$, we have natural functors
\begin{equation}\label{2.adj.bis}
\gamma_!\gamma^*\C \to \C, \qquad \C' \to \gamma^*\gamma_!\C',
\end{equation}
an analog of \eqref{2.adj}, and functors $\C' \to \gamma^*\C$ over
$I'$ are in a natural bijection with functors $\gamma_!\C' \to \C$
over $I$ (recall that since we are dealing with families of
groupoids, all functors are automatically cartesian). Note that if
$\gamma:I' \to I$ is a fibration, then by Lemma~\ref{adj.loc.le}, a
family of groupoids $\C'$ over $I'$ is localizable with respect to
$\gamma$ if and only if so is its restriction $\C'_i$ to any fiber
$I'_i \subset I$ of the fibration $\gamma$. In particular, if the
fibration $\gamma$ is itself a family of groupoids, then any family
of groupoids $\pi':\C' \to I'$ is trivially localizable with respect
to $\gamma$, and we have
\begin{equation}\label{pd.disc}
\gamma_!\C' \cong \C'
\end{equation}
considered as a family of groupoids over $I$ via the fibration $\pi
= \gamma \circ \pi':\C' \to I$. For a general fibration $\gamma$, if
we have a functor $\phi:I_0 \to I$ from some category $I_0$, with
the induced fibration $\gamma_0:I'_0 = \phi^*I' \to I_0$, then
Lemma~\ref{rel.loc.le} provides a natural equivalence
\begin{equation}\label{bc.fib}
  \gamma_{0!}\phi^*\C' \cong \phi^*\gamma_!\C'
\end{equation}
for any family of groupoids $\C' \to I'$ localizable with respect to
$\gamma$. If $I_0 = \ppt$, this reduces to equivalences
\begin{equation}\label{pd.fib}
(\gamma_!\C')_i \cong h^{\Tot}(C'_i), \qquad i \in I
\end{equation}
analogous to \eqref{pd.cof}

\begin{lemma}\label{fib.loc.le}
Assume given a functor $\gamma:I' \to I$ with $1$-connected
comma-fibers $i \setminus_\gamma I'$, $i \in I$. Then $I'$ is
localizable if and only if $I$ is localizable, and in this case,
$h(\gamma):h^{\Tot}(I') \to h^{\Tot}(I)$ is an equivalence.
\end{lemma}

\proof{} By Lemma~\ref{adj.loc.le}, we may replace $I'$ with the
comma-category $I' \setminus I$ and $\gamma$ with the fibration
$\sigma:I \setminus I' \to I$; in other words, we may assume that
$\gamma$ is a fibration with $1$-connected fibers. Then by
Lemma~\ref{rel.loc.le}, $\langle I',\gamma^*(\Iso) \rangle$ is
localizable, and $h^{\gamma^*(\Iso)}(I') \cong I$. If $I$ is
localizable, this proves the claim. If $I'$ is localizable, then we
have a natural functor $h':I \cong h^{\gamma^*(\Iso)}(I') \to
h^{\Tot}(I')$, and $h' \circ \gamma$ is isomorphic to the
localization functor $h:I' \to h^{\Tot}(I')$. Then the universal
property for $h$ implies the same property for $h'$, and this
finishes the proof.
\endproof

\subsection{Weak equivalences.}

Extending Lemma~\ref{cat.loc.le}, one can turn the category of small
relative categories into a relative category by designating a class
of morphisms as weak equivalences. There are several natural ways to
do it. The strongest one is the notion of a weak equivalence
considered in \cite{BK} (a morphism is a weak equivalence iff it
induces a weak equivalence of Dwyer-Kan localizations of
\cite{dw-ka}). We will need the following weaker notions.

\begin{defn}\label{equi.0.def}
A morphism $\gamma:\langle I',W' \rangle \to \langle I,W \rangle$ is
a {\em $0$-equivalence} if for any two morphisms $F,F':\langle I,W
\rangle \to \langle \C,\Iso \rangle$ to some category $\C$, any map
$\alpha':\gamma^*F \to \gamma^*F'$ is of the form $\gamma^*\alpha$
for a unique map $\alpha:F \to F'$. A morphism $\gamma$ is a {\em
  $1$-equivalence} if it is a $0$-equivalence, and moreover, any
morphism $F':\langle I',W' \rangle \to \langle \C,\Iso \rangle$ to
some category $\C$ is isomorphic to $\gamma^*F$ for some morphism
$F:\langle I,W \rangle \to \langle \C,\Iso \rangle$.
\end{defn}

\begin{defn}
A {\em family of groupoids} over a relative category $\langle I,W
\rangle$ is a family of groupoids $\C \to I$ constant along any map
$w \in W$.
\end{defn}

\begin{defn}\label{equi.def}
A morphism $\gamma:\langle I',W' \rangle \to \langle I,W \rangle$
between relative categories is a {\em $2$-equivalence} if
\begin{enumerate}
\item for any family of groupoids $\C'$ over $\langle I',W'
  \rangle$, there exists a family of groupoids $\C$ over $\langle
  I,W \rangle$ and an equivalence $\gamma^*\C \cong \C'$ over $I'$,
\item for any two families of groupoids $\C$, $\C'$ over $\langle
  I,W \rangle$ and a functor $F':\gamma^*\C \to \gamma^*\C'$ over
  $I'$, there exists a functor $F:\C \to \C'$ and an isomorphism $F'
  \cong \gamma^*F$, and
\item for any families of groupoids $\C$, $\C'$ over $\langle
  I,W \rangle$, functors $F,F':\C \to \C'$ over $I$, and a
  morphism $\alpha':\gamma^*F \to \gamma^*F'$ over $I'$, there
  exists a unique morphism $\alpha:F \to F'$ over $I$ such that
  $\alpha' = \gamma^*\alpha$.
\end{enumerate}
\end{defn}

If $I$ and $I'$ are essentially small, so that functor categories
are well-defined, then for any $\C$, a morphism $\gamma:\langle
I',W' \rangle \to \langle I,W \rangle$ induces a restriction functor
\begin{equation}\label{ga.W.eq}
\gamma^*:\Fun^W(I,\C) \to \Fun^{W'}(I',\C),
\end{equation}
where as Example~\ref{fun.v.i.exa}, $\Fun^W(-,-)$ stands for the
category of functors that invert maps in $W$, and similarly for
$W'$. Then $\gamma$ is a $0$-equivalence resp.\ a $1$-equivalence
iff for any $\C$, \eqref{ga.W.eq} is fully faithful resp.\ an
equivalence. Informally, the same holds for $2$-equivalences:
$\gamma$ is a $2$-equivalences iff it induces an equivalence between
$2$-categories of families of groupoids over $\langle I,W \rangle$
and $\langle I',W' \rangle$.

\begin{lemma}\label{loc.adj.le}
In the assumptions of Lemma~\ref{adj.loc.le}, the functors
$\gamma_0$ and $\gamma_1$ are $2$-equivalences.
\end{lemma}

\proof{} For any family of groupoids $\C_0$ over $\langle I_0,W_0
\rangle$, $a_0^*$ is an equivalence between
$\gamma_0^*\gamma_1^*\C_0$ and $\C_0$, and for any family of
groupoids $\C_1$ over $\langle I_1,W_1 \rangle$, $a_1^*$ is an
equivalence between $\C_1$ and $\gamma_1^*\gamma_0^*\C$.
\endproof

\begin{lemma}\label{2.1.eq.le}
A $2$-equivalence $\gamma:\langle I,W \rangle \to \langle I',W'
\rangle$ with essentially surjective $\gamma$ is a $1$-equivalence.
\end{lemma}

\proof{} Assume given a morphism $F:\langle I,W ,\rangle \to \langle
\C,\Iso \rangle$ for some category $\C$, and consider the Yoneda
pairing
\begin{equation}\label{Y.F.eq}
\Y(F):\C^o \times I \to \Sets, \qquad c \times i \mapsto \C(c,F(i))
\end{equation}
and the associated discrete fibration $\wt{\C} \to \C^o \times
I$. For any $c \in c$, its restriction $\wt{\C}_c \to I$ to $\{c\}
\times I \subset \C^o \times I$ is a family of discrete groupoids
over $\langle I,W \rangle$, and since $\gamma$ is a $2$-equivalence,
it comes from a unique family $\wt{\C}'_c$ over $\langle I',W'
\rangle$. Since the functor $\wt{\C}'_c \to \pi_0(\wt{\C}'_c)$
becomes an equivalence after restricting to $I$, it was an
equivalence to beging with, so that $\wt{\C}'_c$ is discrete, and by
virtue of uniqueness, the families $\wt{\C}'_c$ together define a
discrete family of groupoids $\wt{\C}' \to \C^o \times I'$, thus a
functor $\Y(F)':\C^o \times I' \to \Sets$ equipped with an
isomorphism $\Y(F) \cong (\id \times \gamma)^*\Y(F)'$, unique up to
a unique isomorphism. Then for any $i \in I$, $\Y(F)'$ restricts to
a representable functor $\C^o \to \Sets$ on $\C^o \times
\{\gamma(i)\}$, and since $\gamma$ is essentially surjective, any
$i' \in I'$ is of this form. By Yoneda Lemma, this means that
$\Y(F)' \cong \Y(F')$ for a functor $F':I' \to \C$ equipped with an
isomorphism $F \cong \gamma^*F'$, unique up to a unique
isomorphism. By construction, $F'$ is a morphism $\langle I',W'
\rangle \to \langle \C,\Iso \rangle$.
\endproof

\begin{remark}
Essential surjectivity in Lemma~\ref{2.1.eq.le} is needed because
the target category $\E$ in Definition~\ref{equi.0.def} need not be
Karoubi-closed. For instance, consider the category $P^=$ of
Example~\ref{kar.lim.exa} with its subcategory $P \subset P^=$; then
the embedding $\langle P,\Iso \rangle \to \langle P^=,\Iso \rangle$
is a $2$-equivalence but not a $1$-equivalence.
\end{remark}

\begin{lemma}\label{sat.le}
Let $\langle I,W \rangle$ be a relative category, and let $W' =
s(I,W)$ be the saturation of the class $W$ in the sense of
Definition~\ref{sat.def}. Then the natural morphism $\id:\langle I,W
\rangle \to \langle I,W \rangle$ is a $2$-equivalence.
\end{lemma}

\proof{} For any family of groupoids $\C \to I$ over $\langle I,W
\rangle$, the class $W(\C)$ of Lemma~\ref{sat.C.le} is saturated,
thus contains $W'$.
\endproof

\begin{lemma}\label{C.loc.le}
Assume given a $2$-equivalence $\gamma:\langle I',W' \rangle \to
\langle I,W \rangle$ and a family of groupoids $\pi:\C \to I$ over
$\langle I,W \rangle$, and let $\C'=\gamma^*\C$, with its fibration
$\pi':\C' \to I'$. Then $\gamma:\langle \C',{\pi'}^*W' \rangle \to
\langle \C,\pi^*W \rangle$ is a $2$-equivalence.
\end{lemma}

\proof{} A family of groupoids $\C'' \to \C$ over $\C$ is the same
thing as a family $\C'' \to I$ equipped with a functor $\C'' \to \C$
over $I$, and it is constant along $W$ iff it is constant along
$\pi^*W$. The same applies to $\C'$, $I'$ and $W'$.
\endproof

More-or-less by definition, a morphism $\gamma:\langle I',W' \rangle
\to \langle I,W \rangle$ between localizable categories is a
$1$-equivalence if and only if $h^{W'}(\gamma):h(I') \to h^W(I)$ is
an equivalence. For any category $I$, the tautological morphism
$\langle I,{\Tot} \rangle \to \ppt$ is a $0$-equivalence if and
only if $I$ is connected, and it is a $1$-equivalence if and only if
$I$ is $1$-connected.

\begin{defn}
A category $I$ is {\em $2$-connected} if the tautological morphism
$\langle I,{\Tot} \rangle \to \ppt$ is a $2$-equivalence. A functor
$\gamma:I' \to I$ is {\em $l$-connected}, $l=0,1,2$, if the morphism
$\gamma:\langle I',\gamma^*\Iso \rangle \to \langle I,\Iso \rangle$
is an $l$-equivalence.
\end{defn}

Lemma~\ref{2.1.eq.le} immediately implies that a
$2$-connected category is $1$-con\-nec\-ted, and an essentially
surjective $2$-connected functor is $1$-connected.

\begin{lemma}\label{triv.con.le}
Assume given a functor $\gamma:I \to I'$.
\begin{enumerate}
\item If a morphism $\gamma:\langle I,W \rangle \to \langle I,\Iso
  \rangle$ for some class $W \subset \gamma^*(\Iso)$ is an
  $l$-equivalence, $l=0,1,2$, then $\gamma$ is $l$-connected.
\item If $\gamma$ is fully faithful and left or right-reflexive,
  then the adjoint functor $\gamma_\dg:I' \to I$ is $l$-connected,
  $l = 0,1,2$.
\end{enumerate}
\end{lemma}

\proof{} \thetag{i} is obvious. In \thetag{ii}, since $\gamma$ is
fully faithful, $\gamma_\dg \circ \gamma \cong \id$, and as in
Example~\ref{kan.facto.exa}, the adjunction map between $\id$ and
$\gamma \circ \gamma_\dg$ is pointwise in $\gamma_\dg^*(\Iso)$, so
we are done by Lemma~\ref{loc.adj.le} and Lemma~\ref{adj.loc.le}.
\endproof

\begin{corr}\label{refl.con.corr}
Assume given a reflexive full embedding $f:J \to J'$ between
partially ordered sets. Then if $J$ is $l$-connected, $l=0,1,2$, so
is $J'$.
\end{corr}

\proof{} By induction and \eqref{j.filt}, it suffices to consider
the case when $f$ is left or right-reflexive; in this case, the
adjoint map $g:J' \to J$ is $l$-connected by
Lemma~\ref{triv.con.le}~\thetag{ii}.
\endproof

\begin{lemma}\label{2.con.le}
A category $I$ is $2$-connected if and only if either of the
following holds.
\begin{enumerate}
\item Any family of groupoids $\C$ over $\langle I,{\Tot} \rangle$ is
  bounded, and for any $i \in I$, the evaluation functor
  $\Sec^{\Tot}(I,\C) \to \C_i$ is an equivalence.
\item Any family of groupoids $\C$ over $\langle I,{\Tot} \rangle$
  is localizable, and for any $i \in I$, the functor $h_i:\C_i \to
  h^{\Tot}(\C)$ is an equivalence.
\end{enumerate}
\end{lemma}

\proof{} For \thetag{i}, note that for any groupoid $\C$, a
cartesian section $\sigma$ of the trivial fibration $\C \times I$ is
of the form $\sigma = \phi \times \id$ for some functor $\phi:I
\to \C$. If $I$ is $2$-connected, thus $1$-connected by
Lemma~\ref{2.1.eq.le}, then such a functor is constant, so that $\C
\times I \to I$ is bounded and $\Sec^{\Tot}(I,I \times \C) \cong
\C$. Conversely, \thetag{i} implies that for any family of groupoids
$\C$ over $\langle I,{\Tot} \rangle$ the total evaluation functor
$\Sec^{\Tot}(I,\C) \times I \to \C$ is an equivalence, so that $\C
\to I$ is trivial.

For \thetag{ii}, again take any groupoid $\C$, and note that if $I$
is $2$-connected, thus $2$-connected, then the projection $\C \times
I \to \C$ satisfies the assumptions of Lemma~\ref{rel.loc.le}. Since
$\C$ is a groupoid, this means that $\C \times I$ is localizable and
$h^{\Tot}(\C \times I) \cong \C$. Conversely, for any family of
groupoids $\C$ over $\langle I,{\Tot} \rangle$, \thetag{ii} implies
that the natural functor $\C \to h^{\Tot}(\C) \times I$ is an
equivalence, so that $\C \to I$ is trivial.
\endproof

\begin{lemma}\label{2.con.bis.le}
\begin{enumerate}
\item For $l=0,1,2$, a composition of $l$-connected functors is
  $l$-connected, and for any $l$-connected functor $\gamma:I' \to I$
  and relative category $\langle I,W \rangle$, the morphism $\langle
  I',\gamma^*W \rangle \to \langle I,W \rangle$ is an
  $l$-equivalence.
\item Assume given two functors $\pi':I' \to I$, $\pi'':I'' \to I$,
  and a functor $\gamma:I' \to I''$ over $I$. Moreover, assume that
  both $\pi'$ and $\pi''$ are either fibrations or cofibrations, and
  $\gamma$ is cartesian resp.\ cocartesian over $I$. Then if all the
  fibers $\gamma_i:I'_i \to I''_i$, $i \in I$ of the functor
  $\gamma$ are $l$-connected, $l=0,1,2$, $\gamma$ itself is
  $l$-connected.
\end{enumerate}
\end{lemma}

\begin{remark}\label{con.fib.rem}
Lemma~\ref{2.con.bis.le}~\thetag{ii} implies that any fibration or
cofibration $\gamma:I' \to I$ whose fibers are $l$-connected is
$l$-connected (take $I''=I$).
\end{remark}

\proof{} The first claim is clear. For the second, note that if
$\pi'$, $\pi''$ are cofibrations and $\gamma$ is cocartesian, the
embedding $I'_i/i'' \to I'/i$ is left-reflexive for any $i \in I$
and $i'' \in I''_i$. Therefore if $l=2$, a family of groupoids $\C$
over $I'$ is bounded with respect to $\gamma$ iff $\C|_{I'_i}$ is
bounded with respect to $\gamma_i$ for any $i \in I$, and in this
case, we have $\gamma_*\C|_{I''_i} \cong \gamma_{i*}\C|_{I'_i}$, $i
\in I$. Then if all the functors $\gamma_i$ are $2$-connected, a
family $\C'$ over $\langle I',\gamma^*(\Iso) \rangle$ is bounded with
respect to $\gamma$, and the functors \eqref{2.adj} are equivalences
over all the fibers $I''_i$, hence everywhere. This means that
$\gamma$ is $2$-connected. For $l=1$, the same argument shows that
for any morphism $E:\langle I',\gamma^*(\Iso) \rangle \to \langle
\E,\Iso \rangle$ to some category $\E$, the left Kan extension
$\gamma_!E$ exists by virtue of \eqref{kan.eq}, and the adjunction
map $E \to \gamma^*\gamma_!E$ is an isomorphism, so that $E$ factors
uniquely through $\gamma$, and $\gamma$ is therefore
$1$-connected. If $l=0$, the same holds for morphisms of the form
$\gamma^*E:I' \to \E$ for some $E:I'' \to \E$, and then
$\gamma_!\gamma^*E \to E$ is an isomorphism, so that $\gamma$ is
$0$-connected. Finally, if $\pi'$, $\pi''$ are fibrations,
interchange $\gamma_*$ and $\gamma_!$, replace \eqref{2.adj} with
\eqref{2.adj.bis}, and use the dual version of \eqref{kan.eq}.
\endproof

\subsection{Elementary examples.}\label{V.I.subs}

Lemma~\ref{2.con.bis.le} allows one to construct a bunch of simple
but non-trivial examples of $2$-equivalences and $2$-connected
functors. Here is the first one.

\begin{lemma}\label{Z.le}
The functor $\ZZ_\infty \to \N$ of \eqref{zeta.eq} is $2$-connected.
\end{lemma}

\proof{} The left comma-fibers $\ZZ_\infty/m = \ZZ_{2m}$ of the
functor \eqref{zeta.eq} are finite for any $m \in \N$, so for any
family $\C'$ over $\ZZ_\infty$, $\zeta_*\C'$ is well-defined, and we
have the functors \eqref{2.adj}. It suffices to prove that they are
equivalences for any $\C$ over $\N$ and $\C'$ over $\langle
\ZZ_\infty,\zeta^*\Iso \rangle$. This amounts to proving that
$\zeta_{2m}:\ZZ_{2m} \to [m] \subset \N$ is $2$-connected for any $m
\geq 0$. We have the left-reflexive embedding $\langle
\ZZ_{2m-1},\zeta^*(\Iso) \rangle \to \langle \ZZ_{2m},\zeta^*\Iso
\rangle$, and by Lemma~\ref{loc.adj.le}, this is a $2$-equivalence,
so it suffices to prove that $\zeta_{2m-1}:\ZZ_{2m-1} \to [m]$ is
$2$-connected. If $m=0$, $\ZZ_0 \cong [0]$ is the point and
$\zeta_0$ is an isomorphism. If $m \geq 1$, let $\chi = t_\dg:[m]
\to [1]$ be the characteristic function of the initial segment
$[m-1] \subset [m]$ of the ordinal $[m]$, and note that both $\chi$
and $\chi \circ \zeta_{2m-1}$ are fibrations, and $\zeta_{2m-1}$ is
cartesian over $[1]$ with fibers $\zeta_{2(m-1)}$, $\id$. It remains
to apply induction on $m$ and Lemma~\ref{2.con.bis.le}~\thetag{ii}.
\endproof

\begin{corr}\label{Z.corr}
The category $\ZZ_\infty$ is $2$-connected.
\end{corr}

\proof{} By Lemma~\ref{Z.le} and
Lemma~\ref{2.con.bis.le}~\thetag{i}, it suffices to prove that $\N$
is $2$-connected; since $\N$ has the smallest element $0 \in \N$,
this follows from Corollary~\ref{refl.con.corr}.
\endproof

\begin{remark}
By Lemma~\ref{2.1.eq.le}, Lemma~\ref{Z.le} implies that the map
$\zeta:\ZZ_\infty \to \N$ is a $1$-equivalence, so that $\N \cong
h^{\zeta^*(\Iso)}(\ZZ_\infty)$. The corresponding cocartesian square
\eqref{loc.sq} is then the cocartesian square \eqref{Z.N.sq}.
\end{remark}

For the second example, consider the connected groupoid $\ppt_{\Z}$,
with it unique object $o$ and maps $o \to o$ given by integers $n
\in \Z$, and let $\ppt^+_{\Z} \subset \ppt_{\Z}$ be the dense
subcategory spanned by non-negative integers.

\begin{lemma}\label{Z.pl.le}
The embedding $\gamma:\ppt^+_{\Z} \to \ppt_{\Z}$ is $2$-connected.
\end{lemma}

\proof{} The fiber of $\gamma$ over the unique object $o
\in \ppt_{\Z}$ is $\Z$ considered as a partially ordered set with
the standard order; it is $2$-connected by
Corollary~\ref{refl.con.corr} and Example~\ref{Z.refl.exa}, so we
are done by Remark~\ref{con.fib.rem}.
\endproof

To combine Corollary~\ref{Z.corr} and Lemma~\ref{Z.pl.le}, let $\II$
be the Kronecker category $\II$ of Example~\ref{kron.exa}, and
consider the functor $\II \to \ppt_{\Z}$ sending $s$ to $0$ and $t$
to $1$.

\begin{lemma}\label{II.le}
The functor $\II \to \ppt_{\Z}$ is $2$-connected.
\end{lemma}

\proof{} By its very definition, the functor factors through the
embedding $\gamma:\ppt^+_{\Z} \to \ppt_{\Z}$ of Lemma~\ref{Z.pl.le},
so it suffices to prove that the corresponding functor $\II \to
\ppt^+_{\Z}$ is $2$-connected. But its fiber $\II_o$ over the unique
object $o \in \ppt_{\Z}$ is $\ZZ_\infty$, so we are done by
Corollary~\ref{Z.corr} and Remark~\ref{con.fib.rem}.
\endproof

Explicitly, the equivalence $\II_o \cong \ZZ_\infty$ in the proof of
Lemma~\ref{II.le} sends $\langle l,m \rangle \in \II_o$, $l = 0,1 \in
\II$, $m:o \to o$ a non-negative integer, to $2m+l \in
\ZZ_\infty$. The functor
\begin{equation}\label{ZZ.II}
\ZZ_\infty \cong \II_o \to \II
\end{equation}
is actually a discrete fibration; both its fibers are identified
with the set $\Z_+$ of non-negative integers $n \geq 0$, and the
transition functors are given by $s^*(n)=n+1$, $t^*(n)=n$.

\begin{remark}\label{qui.rem}
The name ``Kronecker category'' comes from the more standard
``Kronecker quiver''. In fact, a general quiver is by definition a
functor $X:\II^o \to \Sets$; it has vertices $x \in X(0)$ and edges $x
\in X(1)$, and the maps $s,t:X(1) \to X(0)$ send an edge to its
source resp.\ target.
\end{remark}

The Kronecker category $\II$ is the smallest finite category with
infinite total localization $\ppt_{\Z} \cong h^{\Tot}(\II)$. For our
applications, what is inconvenient about $\II$ is that it is not a
partially ordered set. To achieve the same result with a finite
partially ordered set, compose the functor $\II \to \ppt_{\Z}$ with
the projection $\ppt_{\Z} \to \ppt_{\Z/2\Z}$ induced by the quotient
map $\Z \to \Z/2\Z$, and let $\VV$ be its fiber over $o \in
\ppt_{\Z/2\Z}$. We then have a commutative diagram
\begin{equation}\label{VV.dia}
\begin{CD}
\VV @>>> \ppt^+_{\Z} @>>> \ppt_{\Z},\\
@V{v}VV @VVV @VVV\\
\II @>>> \ppt^+_{\Z} @>>> \ppt_{\Z},
\end{CD}
\end{equation}
where the rightmost vertical arrow is induced by the map $\Z \to
\Z$, $n \mapsto 2n$, and both squares in \eqref{VV.dia} are
cartesian. The horizontal arrows at the top thus have the same
fibers as the arrows at the bottom, so they are again
$2$-connected. The projection $v:\VV \to \II$ is a discrete fibration
with fibers $\Z/2\Z = \{0,1\}$ and transition functors $s^*(n) = n+1
\mod 2$, $t^*(n) = n \mod 2$, and altogether, \eqref{ZZ.II} factors
into the composition
\begin{equation}\label{ZZ.VV}
\begin{CD}
\ZZ_\infty \cong \VV_o \cong \II_o @>{z}>> \VV @>{v}>> \II
\end{CD}
\end{equation}
of discrete fibrations with fibers $\Z_+$ resp.\ $\Z/2\Z$. However,
unlike $\II$, the category $\VV$ is a partially ordered set. It has
four elements, and can be alternatively described as a standard
coproduct \eqref{pos.st.sq}: we have
\begin{equation}\label{VV.eq}
\VV \cong \V^o \copr_{\{0,1\}} \V^o,
\end{equation}
where $\V$ is the partially ordered set \eqref{V.eq}, and the
coproduct is taken with respect to the left-closed embedding
$\{0,1\} \to \V^o \cong \{0,1\}^>$. On the first component $\V^o
\subset \VV$ of this decomposition, $v$ sends the map $0 \to o$ to
$s$ and $1 \to o$ to $t$, and on the second component, it is the
other way around.

\begin{remark}\label{V.I.Z.rem}
The characteristic map $\VV \to \V$ of the standard pushout square
\eqref{VV.eq} is a cofibration, with fibers $\VV_o \cong \{0,1\}$,
$\VV_0 \cong \VV_1 \cong \ppt$, and by Remark~\ref{loc.cart.rem},
the localization $h^{\Tot}(\VV) \cong \ppt_{\Z}$ of \eqref{VV.dia}
produces a cocartesian square of categories. In fact, we have a
diagram
\begin{equation}\label{V.I.Z.sq}
\begin{CD}
\{0,1\} @>>> [1] @>>> \ppt\\
@VVV @VVV @VVV\\
\ppt @>>> \ppt^+_{\Z} @>>> \ppt_{\Z},
\end{CD}
\end{equation}
where $\{0,1\} \to [1]$ is the standard embedding, and both squares
as well as the outer rectangle are cocartesian. If one wants to
include $\II$ in the picture, one has to replace $\ppt$ in the
bottom left corner of \eqref{V.I.Z.sq} with $[1]$; we then have $\II
\cong [1] \copr_{\{0,1\}} [1]$, as in \eqref{II.red.sq}.
\end{remark}

\begin{remark}\label{ine.V.rem}
Lemma~\ref{II.le} immediately implies that for any group\-oid $\C$,
with the inertia groupoid $\I(\C)=\ppt_\Z^o\C$ of
Lemma~\ref{ine.le}, we have a natural equivalence $\I(\C) \cong
\II^o\C$ (and this partially explains our notation for the
Kronecker category $\II$). Alternatively, one can use the
$2$-connected functors in the top row of \eqref{VV.dia} and obtain
an equivalence $\I(\C) \cong \VV^o\C$. The cocartesian square of
Remark~\ref{V.I.Z.rem} then provides a cartesian square
\begin{equation}\label{ine.cart.sq}
\begin{CD}
  \I(\C) @>>> \C\\
  @VVV @VVV\\
  \C @>>> \C \times \C,
\end{CD}
\end{equation}
where both functors $\C \to \C \times \C$ are diagonal embeddings.
\end{remark}

Let us now introduce a whole series of examples of $2$-connected
categories. Namely, recall from Subsection~\ref{lim.subs} that a
category $I$ is {\em filtered} if any functor $J \to I$ from a
finite category $J$ admits a cone, and it in fact suffices to
require this for $J=\emptyset,\{0,1\},\II$. Alternatively, it
suffices to require that cones exist when $J$ is a finite partially
ordered set (to find a cone for a functor $\II \to I$ from the
Kronecker category, compose it with the barycentric subdivision
functor $\xi:\VV \cong B(\II) \to \II$). It is not diffucult to show
that a filtered category is $2$-connected but we will need slightly
more.

\begin{defn}\label{ho.filt.def}
A category $I$ is {\em homotopy filtered} if for any finite artially
ordered set $J$ and functor $\gamma:J \to I$, there exists a functor
$\gamma':J^>\to I$ and a morphism $\gamma'|_J \to \gamma$.
\end{defn}

\begin{prop}\label{filt.prop}
Any homotopy filtered category $I$ is $2$-connected.
\end{prop}

\proof{} For any finite partially ordered set $J$, denote $J' = J^>
\times [1] \ssetminus \{1\}$, where $1$ stands for the terminal
element. We have a natural projection $J' \to J^>$, and this is a
fibration whose fiber $J'_j$ is $[0]$ for $j=1$ and $[1]$
otherwise. Since $J^>$ has a terminal element, it is $2$-connected
by Lemma~\ref{triv.con.le}~\thetag{ii}, and then $J'$ is also
$2$-connected by Lemma~\ref{2.con.bis.le}~\thetag{i}.

Now assume given a homotopy filtered category $I$ and a family of
group\-oids $\C$ over $\langle I,{\Tot} \rangle$. Take two objects
$i,i' \in I$ and objects $c \in \C_i$, $c' \in \C_{i'}$. Any zigzag
diagram $f$ of the form \eqref{zig} representing a map from $c$ to
$c'$ in the localization of $\langle \C,{\Tot} \rangle$ is finite,
and so is its image in $I$. Thus if we have two such zigzag diagrams
$f_0$, $f_1$, there exists a finite partially ordered set $J$ with
two elements $j,j' \in J$, a functor $\gamma:J \to I$ such that
$\gamma(j)=i$ and $\gamma(j')=i'$, and two diagrams $g_0$, $g_1$
representing maps from $c \in (\gamma^*\C)_j$ to $c' \in
(\gamma^*\C)_{j'}$ in $h^{\Tot}(\gamma^*\C)$ such that
$\gamma(g_0)=f_0$ and $\gamma(g_1)=f_1$. Since $I$ is homotopy
filtered, $\gamma$ factors as $\gamma = \gamma' \circ \eps$, where
$\eps:J \to J'$ is the embedding onto $J \times 1 \subset J'$, and
$\gamma'$ is a functor from $J'$ to $I$. However, $J'$ is
$2$-connected, so that $\gamma^{'*}\C \cong \C_i \times J'$ a
constant family. Thus $\eps(g_0)=\eps(g_1) \cdot g$ for some $g:c
\to c$ in $\C_i$, and then $f_0 = f_1 \cdot g$. Therefore $\langle
\C,{\Tot} \rangle$ is localizable, and the natural functor
$h_i:\C_i \to h^{\Tot}(\C)$ is full.

Moreover, if we have two maps $g,g':c \to c$ in $\C_i$ such that
$h_i(g) = h_i(g')$, then they are related by a finite chain of
cancellation rules involving a finite number of diagrams \eqref{zig}
in $\C$. Then again, we can choose a finite partially ordered set
$J$ with an element $j \in J$ and a functor $\gamma:J \to I$ such
that $\gamma(j)=i$, and $g$ becomes equal to $g'$ after projecting
to $h^{\Tot}(\gamma^*\C)$. Extending $\gamma$ to $\gamma':J' \to
I$, we see that $g$ and $g'$ are all the more equal after projection
to $h^{\Tot}(\gamma^{'*}\C) \cong \C_i$, so that $g=g'$. Therefore
the functor $h_i$ is also faithful.

Finally, for any $i' \in I$ and $c' \in \C_{i'}$, we can take
$J=\{j,j'\}$ with trivial order, and let $\gamma(j)=i$,
$\gamma(j')=i'$. Then extending $\gamma$ to $J'$, we obtain a
diagram $f$ in $I$ representing a map from $i$ to $i'$ in
$h(I,\Tot)$, and since $\C$ is a family over $\langle I,{\Tot}
\rangle$, it is actually a bifibration, so that $f$ admits a lifting
to a map $\wt{f}:c \to c'$ in $h^{\Tot}(\C)$ for some $c \in
\C_i$. Therefore $h_i$ is essentially surjective, and we are done by
Lemma~\ref{2.con.le}~\thetag{ii}.
\endproof

In particular, a filtered category is trivially homotopy filtered,
so Proposition~\ref{filt.prop} applies. Note that $I$ is {\em not}
required to be small. As an application of this, and motivated by
Remark~\ref{loc.cart.rem}, let us introduce the following version of
filtered colimits for large categories and fully faithful functors.

\begin{lemma}\label{filt.ff.le}
Assume given a cofibration $\C \to I$ such that $I$ is filtered, and
all transition functors $f_!:\C_i \to \C_{i'}$ are fully faithful.
Then $\C$ is localizable with respect to the class $\sharp$ of maps
cocartesina over $I$, and for any $i \in I$, the localization
functor $h:\C_i \to h^\sharp(\C)$ is fully faithful.
\end{lemma}

\proof{} For any $i,i' \in I$, the comma-fiber $i \times i'
\setminus I = i \times i' \setminus_{\delta} I$ of the diagonal
embedding $\delta:I \to I \times I$ is the category of objects $i''
\in I$ equipped with maps $f:i \to I''$, $f':i' \to i''$, and it is
obviously filtered, For any objects $c \in \C_i$, $c' \in \C_{i'}$,
we have a functor $\C(f_!(c),f'_!(c')):i \times i' \setminus I \to
\Sets$, thus a discrete cofibration over $i \times i' \setminus I$,
and since all the transition functors are fully faithful, the
cofibration is constant along all maps. Then by
Proposition~\ref{filt.prop}, it is trivial. Therefore setting
\begin{equation}\label{filt.ff.eq}
h^\sharp(\C)(c,c') = \colim_{(i \times i') \setminus
  I}\C(f_!(c),f_!(c'))
\end{equation}
gives a well-defined set, and the map $\C(f_!(c),f'_!(c)) \to
h^\sharp(\C)(c,c')$ is an isomorphism for any $\langle i'',f,f'
\rangle \in i \times i' \setminus I$. It is then an easy exercise
to define compositions for the sets \eqref{filt.ff.eq} turning
$h^\sharp(\C)$ into a category, and to show that the obvious
functor $h:\C \to h^\sharp(\C)$ is indeed a localization.
\endproof

\begin{defn}\label{filt.ff.def}
Under the assumptions of Lemma~\ref{filt.ff.le}, the {\em
  $2$-co\-limit} $\colim^2_I \C$ of the cofibration $\C$ is
$\colim^2_I\C = h^\sharp(\C)$.
\end{defn}

\begin{exa}\label{filt.ff.exa}
For any category $\C$, full small subcategories $\C' \subset \C$ and
inclusions between them form a category $I(\C)$ (were it small, it
would be a partially ordered set). Then if we let $\wt{\C} \subset \C
\times I(\C)$ be the full subcategory spanned by pairs $\langle
c,\C' \rangle$ such that $c \in \C'$, the projection $\wt{\C} \to
I(\C)$ is a cofibration, and $\C \cong \colim_{I(\C)}
\wt{\C}$. Alternatively, one can replace $I(\C)$ with any
subcategory $I'$ that is $0$-cofinal in the sense of
Example~\ref{cofin.exa}.
\end{exa}

\section{Localization over $\Delta$.}\label{seg.loc.sec}

\subsection{Simplicial expansions.}

Another series of examples of $2$-connected categories and functors
is provided by the simplicial enhancements and simplicial
replacements of Definition~\ref{simpl.repl.def}. Recall that for any
category $I$ and closed class $v$ of maps in $I$, we denote by $+$
the class of special maps in $\Delta^vI$.

\begin{lemma}\label{del.prop}
For any category $I$ and class of morphisms $v$ in $I$, the functor
\eqref{xi.eq} induces a $2$-equivalence $\xi:\langle \Delta^vI,+
\rangle \to \langle I,\Iso \rangle$, and in particular, $\xi$ is
$2$-connected. Moreover, the functor
$$
\gamma=\xi^o:(\Delta I^o)^o \to I
$$
is also $2$-connected.
\end{lemma}

\proof{} The first claim immediately follows from
Proposition~\ref{spec.prop}. For the second claim, it suffices to
show that for any family of groupoids $\C/I$ and any family of
groupoids $\C'/(\Delta I^o)^o$ constant along $\gamma^*(\Iso)$, the
functors \eqref{2.adj.bis} are equivalences. This amounts to
checking that for any $i \in I$, the comma-fiber $i\setminus I'$ is
$2$-connected. But as in the proof of Proposition~\ref{spec.prop},
we have a diagram
$$
\begin{CD}
\ppt @>>> (\Delta_+ I^o)_i^o @<<< i \setminus I',
\end{CD}
$$
and both functors in the diagram are right-reflexive. Since $\ppt$
is $2$-connected, the other two categories are then $2$-connected by
Lemma~\ref{loc.adj.le}.
\endproof

It turns out that Lemma~\ref{del.prop} has quite a lot of
applications and generalizations. For a first generalization, as in
Subsection~\ref{del.0.subs}, consider the factorization system
$\langle \dd,\ff \rangle$ on $\Delta$ consisting of surjective
resp.\ injective maps. For any category $I$ with a closed class of
maps $v$, define the {\em restricted simplicial expansion}
$\Delta^v_\ff I$ by $\Delta^v_\ff I = \Delta_\ff \times_\Delta
\Delta^v I$, and by abuse of notation, let $+$ be the class of maps
in $\Delta^v_\ff I$ that are special as maps in $\Delta^vI$.

\begin{corr}\label{del.corr}
For any category $I$ with a closed class of maps $v$, the morphism
$\xi:\langle \Delta_\ff^vI,+ \rangle \to \langle I,\Iso \rangle$
induced by \eqref{xi.eq} is a $2$-equivalence.
\end{corr}

\proof{} Let $I^p$ be the category with the same objects as $I$, and
whose morphisms are those of $I$ plus an extra morphism $p_i:i \to
i$ for any $i \in I$, subject to relations $p_{i'} \circ f = f \circ
p_i = f$ for any $f:i \to i'$ (in particular,
$p_i^2=p_i$). Informally, $I^p$ is obtained by formally adding a new
identity map to any object in $I$. Let $\nu:I^p \to I$ be the
forgetful functor sending $p_i$ to $\id_i$, obviously $2$-connected,
and note that $v$ also defines a closed class of maps in $I^p$ (not
containing any of the $p_i$). Moreover, the simplicial expansion
$\Delta^vI^p$ can be described as follows. Consider the fibration
$\sigma:\Ar^\dd(\Delta) \to \Delta$ of (the dual version of)
Lemma~\ref{facto.fib.le}~\thetag{ii}, and let $\Delta^\dd =
\Ar^\dd(\Delta)_{\sigma^\dg\Tot}$ be the underlying discrete
fibration. Since $\sigma^\flat\Tot = \tau^*\ff$, $\Delta^\dd$ is the
category of surjective arrows $p:[n] \to [m]$ in $\Delta$, with
morphisms given by diagrams \eqref{ar.sq} with injective arrow on
the right. Then we have $\Delta^vI^p \cong \Delta^\dd
\times_{\Delta_\ff} \Delta_\ff^vI$, with the fibration
$\sigma:\Delta^vI^p \to \Delta$ being the structural fibration.

Now, as in Example~\ref{adj.exa}, the projection $\tau:\Delta^\dd \to
\Delta_\ff$ has the fully faithful right-adjoint functor
$\eta:\Delta_\ff \to \Delta^\dd$, and this induces a fully faithful
right-adjoint $\eta:\Delta_\ff^vI \to \Delta^vI^p$ to the projection
$\tau:\Delta^v I^p \to \Delta_\ff^vI$. Then $\xi:\langle \Delta_\ff^vI,+
\rangle \to \langle I,\Iso \rangle$ factors as
\begin{equation}\label{v.p.dia}
\begin{CD}
\langle \Delta_\ff^vI,+ \rangle @>{\eta}>> \langle \Delta^vI^p,\tau^*(+)
\rangle @>{\xi}>> \langle I^p,\nu^*(\Iso) \rangle @>{\nu}>> \langle
I,\Iso \rangle.
\end{CD}
\end{equation}
Since $\nu$ is $2$-connected, the last morphisms in \eqref{v.p.dia}
is a $2$-equivalence, and the first morphism is a $2$-equivalence by
Lemma~\ref{loc.adj.le}, so it suffices to prove that the middle
morphism is also a $2$-equivalence. But since a surjective map in
$\Delta$ is special, we have $+ \subset \tau^*(+)$, so that by
Lemma~\ref{del.prop}, it suffices to check that a fibration $\C \to
I^p$ is constant along the maps $p_i$ if $\xi^*\C$ is constant along
all the maps in $\tau^*(+)$. This is clear: the arrow $p_i:i \to i$
defines an object $p'_i \in (\Delta I^p)_{[1]}$, and if we consider
the cartesian lifting $s_i:i \to p'_i$ of the initial embedding
$s:[0] \to [1]$, then $\tau(s_i) = \id$ is special, and $\xi(s_i) =
p_i$.
\endproof

For a first application, note that the notion of a $2$-connected
functor is not stable under base change: the pullback $\phi^*I' \to
I''$ of a $2$-connected functor $\gamma:I' \to I$ with respect to a
functor $\phi:I'' \to I$ need not be $2$-connected. To insure that a
functor $\gamma$ is $2$-connected together with all its pullbacks,
one needs to impose further conditions. For example, by
Remark~\ref{con.fib.rem}, it suffices to require that $\gamma$ is a
fibration or a cofibration with $2$-connected fibers. We will need
the following slight weakening of the latter.

\begin{defn}\label{loc.con.def}
A functor $\gamma:I' \to I$ is a {\em locally $2$-connected} if
\begin{enumerate}
\item for any object $i \in I$, the fiber $I'_i$ is $2$-connected,
  and
\item for any morphism $f:i \to i'$ in $I$ and object $\wt{i} \in
  I'_i$, the category $I'(f,\wt{i})$ of morphisms $\wt{f}:\wt{i} \to
  \wt{i}'$ such that $\gamma(\wt{f})=f$ is $2$-connected.
\end{enumerate}
\end{defn}

We note that if $\gamma:I' \to I$ is a cofibration, then the
category $I'(f,\wt{i})$ has an initial object given by a cocartesian
lifting of the map $f$, so that
Definition~\ref{loc.con.def}~\thetag{ii} is trivially satisfied by
Lemma~\ref{triv.con.le}~\thetag{ii}, and $\gamma$ is locally
$2$-connected iff it has $2$-connected fibers. We also note that
Definition~\ref{loc.con.def} is trivially stable under pullbacks.

\begin{prop}\label{reso.prop}
A functor $\gamma:I' \to I$ locally $2$-connected in the sense of
Definition~\ref{loc.con.def} is $2$-connected.
\end{prop}

\proof{} Let $v = \gamma^*\Id$ be the class of maps $f$ in $I'$ such
that $\gamma(v)$ is an identity map in $I$, and consider the
simplicial enhancement $\Delta^vI'$ and the simplicial replacement
$\Delta I$. Then $\gamma$ induces a functor $\Delta(\gamma):\Delta^v
I' \to \Delta I$, and by Lemma~\ref{del.prop}, it suffices to prove
that this functor is $2$-connected. It is obviously a fibration, so
that by Lemma~\ref{2.con.bis.le}~\thetag{ii}, it suffices to check
that the fiber $(\Delta^vI')_{\langle [n],i_\idot\rangle}$ is
$2$-connected for any $\langle [n],i_\idot \rangle \in \Delta I$.

If $n=0$, so that $i_\idot = i$ is a single object in $I$, then
$(\Delta^vI')_{\langle [0],i_\idot \rangle}$ is equivalent to
$I'_i$, so we are done by
Definition~\ref{loc.con.def}~\thetag{i}. If $n$ is strictly
positive, then we have the natural map $s:[n-1] \to [n]$ and the
corresponding transition functor
\begin{equation}\label{s.del}
s^*:(\Delta^vI')_{\langle [n],i_\idot \rangle} \to
(\Delta^vI')_{\langle [n-1],s^*i_\idot \rangle},
\end{equation}
and we may assume by induction that the right-hand side is
$2$-connected. Then to finish the proof, it suffices to observe that
\eqref{s.del} is a fibration, and its fiber over an object $\langle
[n'],i'_\idot \rangle \in (\Delta^vI')_{\langle
  [n-1],s^*i_\idot \rangle}$ is equivalent to $I'(f,i'_{n'})$, where
$f:i_{n-1} \to i_n$ is the map in the diagram \eqref{i.dot}. Thus we
are done by Definition~\ref{loc.con.def}~\thetag{ii} and
Remark~\ref{con.fib.rem}.
\endproof

\subsection{Regular Segal categories.}

One can also generalize Lemma~\ref{del.prop} to Segal categories in
the sense of Section~\ref{seg.sec} that are not simplicial
expansions. To do this, assume given a Segal category $\C$ that is
regular in the sense of Definition~\ref{seg.reg.def}, and let
$\pi:\C \to \Delta$ be the structural fibration.

\begin{defn}
The {\em skeleton} $\sk(\C)$ of a regular Segal category $\C$ is its
$2$-simplicial expansion $\sk(\C) = \Delta(\C\|\Delta)$ of
Definition~\ref{simpl.rel.def} with respect to the class of identity
maps in $\C_{[0]}$.
\end{defn}

By Lemma~\ref{rel.del.le}, $\sk(\C)$ is a $2$-category in the sense
of Example~\ref{2.cat.exa}. To see it explicitly, note that by
regularity of $\C$, for any objects $c_0 \in \C_{[0]}$, $c_1,c_1'
\in \C_{[1]}$ and morphisms $s:c_0 \to c_1$, $t:c_0 \to c_1'$ such
that $\pi(s) = s$, $\pi(t)=t$, and $t$ is cartesian over $\Delta$ we
have a cocartesian square in $\C$ of the form
$$
\begin{CD}
c_0 @>{s}>> c_1'\\
@V{t}VV @VVV\\
c_1 @>>> c_2
\end{CD}
$$
with some $c_2 \in \C_{[2]}$. Then the objects of $\sk(\C)$ are the
objects of $\C_{[0]}$, and for any $c,c' \in \C_{[0]}$, the category
of morphisms $\sk(\C)(c,c')$ is the category of diagrams
\begin{equation}\label{del.dom}
\begin{CD}
c @>{s}>> c_1 @<{t}<< c'
\end{CD}
\end{equation}
in $\C$ with $c_1 \in \C_{[1]}$, $\pi(s)=s$, $\pi(t)=t$, and $t$
cartesian over $\Delta$. The identity map $c_0 \to c_0$ corresponds
to the diagram \eqref{del.dom} with $c_1 = e^*c$, where $e:[1] \to
[0]$ is the projection, and the composition of morphisms $c_0 \to
c_0'$, $c_0' \to c_0''$ represented by diagrams \eqref{del.dom} with
middle terms $c_1$, $c_1'$ is the diagram
$$
\begin{CD}
c_0 @>>> m^*(c_1 \copr_{c_0'} c_1') @<<< c_0'',
\end{CD}
$$
where as in \eqref{2.cat.m}, $m:[1] \to [2]$ sends $0$ to $0$ and
$1$ to $2$. Moreover, regularity of $\C$ immediately implies that
all maps between diagrams \eqref{del.dom} are invertible, so that
$\sk(\C)$ is in fact a $(1,2)$-category (or equivalently, a
$2$-family of groupoids over a point).

\begin{lemma}\label{seg.loc.le}
Assume given a Segal category $\C$ that is regular in the sense of
Definition~\ref{seg.reg.def}, and assume that its skeleton $\sk(\C)$
is bounded in the sense of Definition~\ref{exa.bnd.def}.  Denote by
$+$ the class of maps in $\C$ special in the sense of
Definition~\ref{seg.spec.def}. Then $\langle \C,+ \rangle$ is
localizable, and $h^+(\C)$ is naturally equivalent to the truncation
$\overline{\sk(\C)}$ of the $(1,2)$-category $\sk(\C)$ in the sense
of Example~\ref{2.cat.fam.exa}. Moreover, assume given a Segal
category $\C'$ and a Segal functor $\pi:\C' \to \C$ that is an
epivalence, and an equivalence over $[0] \in \Delta$. Then $\langle
\C',+ \rangle$ is also localizable, and $h^+(\pi):h^+(\C') \to
h^+(\C)$ is an equivalence.
\end{lemma}

\proof{} For any $[n] \in \Delta$, we have a unique special map
$t:[0] \to [n]$, so that for any $c \in \C'$, we have a canonical
special map $t:c_0 = t^*c \to c$, and then for any map $f:c \to c'$
in $\C'$, we have, up to an isomorphism, a unique diagram
\eqref{del.dom} in $\C'$ equipped with a map $g:c_1 \to c'$ such
that $f \circ t = g \circ s$ and $g \circ t = t$. Moreover, since
$\pi$ is an equivalence over $[0]$ and an epivalence over $[1]$,
isomorphism classes of diagrams \eqref{del.dom} in $\C$ and $\C'$
are the same. Then sending $c \in \C'$ to $c_0 = t^*c$ and $f:c \to
c'$ to the class of this unique diagram defines a functor $h:\C' \to
\overline{\sk(\C)}$ that inverts special maps. Moreover, any functor
$\gamma:\C' \to \E$ to some $\E$ that inverts special maps factors
uniquely through $h$, with a map represented by a diagram
\eqref{del.dom} going to $\gamma(t)^{-1} \circ \gamma(s)$.  Thus
$\overline{\sk(\C)}$ indeed gives the localization $h^+(\C') \cong
h^+(\C)$.
\endproof

\begin{exa}\label{del.loc.exa}
For any category $I$, the product $I \times \Delta$ is a regular
Segal category, we have $\sk(I \times \Delta) \cong \Delta I$ and
$\overline{\sk(I \times \Delta)} \cong I$, and the localization
functor $\Delta I \to h^+(\Delta I) \cong I$ is the functor
\eqref{xi.eq}.
\end{exa}

In particular, consider a bounded Segal category $\C \to \Delta$
that is a family of groupoids (or equivalently, a $2$-family of
groupoids over a point in the sense of
Definition~\ref{exa.seg.def}). Then it is regular by
Example~\ref{seg.reg.exa}, thus lies within the scope of
Lemma~\ref{seg.loc.le}. Thus $\C$ is localizable with respect to the
class of special maps, with the localization $h^+(\C)$ described in
Lemma~\ref{seg.loc.le}. Note that in particular, the $1$-opposite
Segal category is also bounded, and we have
\begin{equation}\label{opp.loc}
h^+(\C^\iota) \cong h^+(\C)^o.
\end{equation}
Moreover, since $\C$ is a family of groupoids, every map in $\C$ is
cartesian over $\Delta$, so that the universal special functor
\eqref{2.del.al} is a Segal functor and an equivalence. Since
$\C_{[0]}$ is a groupoid, the localization functor $h:\C \to
h^+(\C)$ sends $\C_{[0]} \subset \C$ into the isomorphism groupoid
$h^+(\C)_{\Iso} \subset h^+(\C)$, and we have a commutative diagram
\begin{equation}\label{ind.h.sq}
\begin{CD}
\sk(\C) @>>> \pi_0(\sk(\C)) = \Delta h^+(\C)\\
@VVV @VVV\\
\C \cong \Delta^{\Tot}(\C\|\Delta) @>{\Ind(h \times \id)}>>
\Delta^{\Iso}(h^+(\C) \times \Delta\|\Delta) = \Delta^{\Iso}h^+(\C).
\end{CD}
\end{equation}

\begin{lemma}\label{exa.loc.le}
If a Segal functor $\gamma:\C \to \C'$ between bounded $2$-families
of groupoids $\C$, $\C'$ over a point is $1$-full resp.\ $2$-full in
the sense of Definition~\ref{fib.fu.def}, then the localization
$h^+(\gamma):h^+(\C) \to h^+(\C')$ is full resp.\ fully faithful,
and if $\gamma$ is essentially surjective, then so is $h^+(\gamma)$.
\end{lemma}

\proof{} By \eqref{ind.h.sq}, the special functor $\sk(\C) \to \C$
is a Segal functor of the form \eqref{v.01.eq}, thus $2$-fully
faithful and essentially surjective by Example~\ref{v.01.exa}. The
same applies to $\C'$, so that $\gamma$ is $1$-full or $2$-full if
and only if so is $\sk(\gamma)$. Then the essential surjectivity
claim immediately follows from \eqref{ind.h.sq}, and to deduce the
rest, it suffices to observe that if a functor $\G \to \G'$ between
two essentially small groupoids is essentially surjective resp.\ an
epivalence, then the induced map $\pi_0(\G) \to \pi_0(\G')$ is
surjective resp.\ an isomorphism.
\endproof

Now let $\C$ be a bounded $2$-family of groupoids over some category
$I$. Then by Example~\ref{exa.reg.exa}, $\C$ is co-regular but not
necessarily regular, so it lies outside of scope of
Lemma~\ref{seg.loc.le}. However, the $(1,2)$-opposite Segal category
$\C^\sigma$ is regular, so that for any closed class $v$ of maps in
$\C_{[0]}$, the twisted $2$-simplicial expansion
$\Delta^v_\sigma(\C\|\Delta)$ of \eqref{tw.rel.eq} is a Segal
category. We have a diagram of Segal functors
\begin{equation}\label{seg.loc.dia}
\begin{CD}
\C @>{\beta}>> \Delta^\Tot_\sigma(\C\|\Delta) @<{\lambda}<<
\Delta^v_\sigma(\C\|\Delta),
\end{CD}
\end{equation}
where $\beta$ corresponds to \eqref{tw.beta.eq}, and $\lambda$ is
the dense embedding \eqref{v.01.eq}. If $v \subset \Iso$ consists of
isomorphisms, $\Delta^v_\sigma(\C\|\Delta)$ is a bounded $2$-family
of groupoids over the point, thus regular by
Example~\ref{seg.reg.exa}.

\begin{lemma}\label{seg.tw.le}
Any bounded $2$-family of groupoids $\C$ over a category $I$ is
localizable with respect to the class of special maps, and the
projection $h^+(\C) \to I$ is a fibration. Moreover, the
localization functor $h:\C \to h^+(\C)$ factors as $h \cong h' \circ
\beta$, where $\beta$ is as in \eqref{seg.loc.dia} and
$h':\Delta^\Tot_\sigma(\C\|\Delta) \to h^+(\C)$ is a special
functor, and if $v \subset \Iso$ is a closed class $v \subset \Iso$
of isomorphisms in $\C_{[0]}$, then the special functor $h' \circ
\lambda:\Delta^v_\sigma(\C\|\Delta) \to h^+(\C)$ induces an
equivalence $h^+(\Delta^v_\sigma(\C\|\Delta)) \cong h^+(\C)$.
\end{lemma}

\proof{} Since for any $i \in I$, $\C|_{i \times \Delta}$ is regular
by Example~\ref{seg.reg.exa}, the first claim directly follows from
Lemma~\ref{rel.loc.le}. Moreover, the $(1,2)$-opposite Segal
category $\C^\sigma$ is localizable with respect to special maps by
Lemma~\ref{seg.loc.le}, and the equivalence \eqref{opp.loc} extends
to an equivalence $h^+(\C^\sigma) \cong h^+(\C)^o$. Then the
localization functor $h:\C^\sigma \to h^+(\C)^o$ is special, and we
can define $h'$ as the composition
$$
\begin{CD}
\Delta^\Tot_\sigma(\C\|\Delta) @>{\Delta(h)}>>
\Delta^\Tot_\sigma(h^+(\C) \times \Delta\|\Delta) \cong
\Delta^\Tot(h^+(\C)) @>{\xi}>> \C,
\end{CD}
$$
where $\xi$ is the special functor \eqref{xi.eq}. Then $h' \circ
\lambda$ is special, thus factors through the localization
$h^+(\Delta^v_\sigma(\C\|\Delta)$, and Lemma~\ref{seg.loc.le}
together with \eqref{tot.del.1} immediately show that the induced
functor $h^+(\Delta^v_\sigma(\C\|\Delta)) \to h^+(\C)$ is
essentially bijective and fully faithful.
\endproof

\subsection{Complete $2$-families of groupoids.}

A particular example of a $2$-family of groupoids $\C$ over some
category $I$ is given by $2$-families $\Delta^\flat(\C|I)$ of
Example~\ref{iso.seg.exa}. It turns out that one can use
localization to characterize precisely the $2$-families that arise
in this way.

\begin{defn}\label{exa.refl.def}
A family of groupoids $\C \to I \times \Delta$ is {\em reflexive} if
it is fully faithful along the map $e:[n] \to [0]$ for any $[n] \in
\Delta$.
\end{defn}

\begin{defn}\label{exa.compl.def}
A family of groupoids $\C \to I \times \Delta$ is {\em complete} if
the commutative square
\begin{equation}\label{exa.compl.sq}
\begin{CD}
\C_{[3]} @>{\mu_0 \times \mu_1}>> \C_{[1]} \times
\C_{[1]}\\
@A{e^*}AA @AA{e^* \times e^*}A\\
\C_{[0]} @>>> \C_{[0]} \times \C_{[0]}
\end{CD}
\end{equation}
induced by \eqref{compl.sq} is cartesian.
\end{defn}

\begin{prop}\label{exa.seg.prop}
\begin{enumerate}
\item For any fibration $\pi:\C \to I$, the $2$-family of groupoids
  $\Delta^\flat(\C|I)$ of Example~\ref{iso.seg.exa} is reflexive in the
  sense of Definition~\ref{exa.refl.def}, bounded in the sense of
  Definition~\ref{exa.bnd.def}, and complete in the sense of
  Definition~\ref{exa.compl.def}.
\item Conversely, assume given a complete $2$-family of groupoids
  $\C$ over a category $I$. Then $\C$ is reflexive. Moreover, if
  $\C$ is bounded, with the localization $\bC=h^+(\C)$ of
  Lemma~\ref{seg.tw.le} and the fibration $\bC \to I$, then we have
  $\C^{red} \cong \Delta^\flat(\bC|I)$, where $\C^{red}$ is the
  reduction \eqref{seg.red}.
\end{enumerate}
\end{prop}

\proof{} For \thetag{i}, reflexivity is obvious, and so is
boundedness: by the Segal condition, it suffices to check it over
$[1] \in \Delta$, and then for any $i \in I$, the fiber of
\eqref{exa.seg.adj} over $c \times c' \in \C_{[0] \times i} \times
\C_{[0] \times i} \cong \eps_*\eps^*\C_{[1] \times i}$ is the
$\Hom$-set $\C_i(c,c')$. To see that $\Delta^\flat(\C|I)$ is
complete, note that the fibered product in \eqref{exa.compl.sq} is
the groupoid of functors $c_\idot:[3] \to \C$ such that $\pi \circ
c_\idot:[3] \to I$ is isomorphic to the constant functor with value
$i$, and the maps $c_0 \to c_2$, $c_1 \to c_3$ are invertible. This
implies that $c_1 \to c_2$ is invertible, with right and
left-inverses given by $c_0 \to c_1$, $c_2 \to c_3$, and then the
whole $c_\idot$ is isomorphic to a constant functor. Since
$\Delta^\flat(\C/I)$ is reflexive, $e^* \times e^*$ in
\eqref{exa.compl.sq} is fully faithful, so this finishes the proof.

For \thetag{ii}, assume first that $I=\ppt$, so that $\C \to \Delta$
is a $2$-family over the point. Let $\C^e \subset \C$ be
the full subcategory spanned by essential images of the functors
$e^*:\C_{[0]} \to \C_{[n]}$, $[n] \in \Delta$. Then $\C^e$ is also a
complete $2$-family, and we have $\pi_0(\C^e|\Delta) \cong
\Delta \times \pi_0(\C_{[0]})$. Thus any object $c \in \C_{[0]}$
defines a section $c:\Delta \to \C^e$, and to show that $\C$ is
reflexive, it suffices to prove that $\pi_1(\C^e,c)$ is a constant
simplicial group. But by the Segal condition, it is actually the
nerve of the small category, this category must be a groupoid by
Lemma~\ref{eck.le}, and then since $\C^e$ is complete, the groupoid
must be discrete. Since it is connected, it is trivial.

Now assume further that $\C$ is bounded. Then \eqref{ind.h.sq}
provides a functor $\phi:\C \to \Delta^{\Iso}(\bC)$ whose target is
a reduced $2$-family of groupoids over $I = \ppt$. To prove that
$\phi$ induces an equivalence $\C^{red} \cong \Delta^\Iso\bC$, it
suffices to check that it is an equivalence over $[0]$ and an
epivalence over $[1]$. Over $[0]$, note that since $\C$ is complete,
a diagram \eqref{del.dom} represents an invertible map in $\bC$ only
if $c_1 \cong e^*c_0$. Therefore to compute $(\Delta^\Iso\bC)_{[0]}
\cong \bC_\Iso$, we may replace $\C$ with its full subcategory $\C^e
\subset \C$, and then since $\C$ is reflexive, $\C^e \cong \Delta
\times \C_{[0]}$ and the claim is obvious. Then over $[1]$, the
fiber $\C_{c_0 \times c'_0}$ of the projection \eqref{exa.seg.adj}
over $c_0 \times c_0' \in \C_{[0]} \times \C_{[0]}$ is the groupoid
of diagrams \eqref{del.dom}, and the corresponding fiber for
$\Delta^\Iso\bC$ is the set of isomorphism classes of its objects.

Finally, in the case of an arbitrary $I$, reflexivity is in any case
defined separately over each $i \in I$, and the localization functor
$h:\C \to \bC$ of Lemma~\ref{seg.tw.le} is cartesian over $I$, so
that we still have a functor
\begin{equation}\label{del.cmp.eq}
\phi:\C \cong \Delta^\flat(\C\|\Delta|I) \to \Delta^\flat(\bC|I)
\end{equation}
between $2$-families with reduced target. Moreover, it factors
through a functor $\C^{red} \to \Delta^\flat(\bC|I)$ that is an
equivalence over any $i \in I$, so it is itself an equivalence.
\endproof

\begin{remark}\label{ind.h.rem}
To see the functor \eqref{del.cmp.eq} more explicitly, one can use
the interpretation of $2$-simplicial expansions in terms of
cartesian squares \eqref{d.c.d}. Namely, the localization functor
$h$ is a morphism of relative categories $\langle \C,+ \rangle \to
\langle \bC,\flat \rangle$, thus induces a functor
$$
\Delta(h):\Delta^+(\C|I) \to \Delta^\flat(\bC|I).
$$
Then $\delta:\Delta \to \Delta^+\Delta$ of \eqref{al.eq}
lifts to a functor $\delta:\C \to \Delta^+(\C|I)$, and the
composition
$$
\begin{CD}
\C @>{\delta}>> \Delta^+(\C | I) @>{\Delta(h)}>>
\Delta^\flat(\bC|I),
\end{CD}
$$
is exactly the functor \eqref{del.cmp.eq}.
\end{remark}

\begin{remark}\label{fl.seg.rem}
As an application of Proposition~\ref{exa.seg.prop}, one can
characterize functors $\gamma:I' \to I$ that are proper in the sense
of Definition~\ref{fl.def}: an essentially small functor $\gamma$ is
proper iff for any $n > l > 0$ and functor $[n] \to I$, the
commutative square of categories
\begin{equation}\label{seg.rel.pr}
\begin{CD}
[0] \times_I I' @>>> [n-l] \times_I I'\\
@VVV @VVV\\
[l] \times_I I' @>>> [n] \times_I I'
\end{CD}
\end{equation}
induced by \eqref{seg.del.sq} is cocartesian. To construct the
relative functor category $\Fun(I'|I,\C)$ for some target category
$\C$, one first constructs the fibration $\Delta^\Tot\Fun(I'|I,\C)
\to \Delta^\Tot I$ with fibers $\Fun([n] \times_I I',\C)$, and then
the fact that \eqref{seg.rel.pr} is cocartesian insures that this is
a Segal fibration; to recover $\Fun(I'|I,\C)$, one then takes the
localization with respect to special maps. We do not give any
details since we will not need this.
\end{remark}

We will need a slight generalization of
Proposition~\ref{exa.seg.prop}. Namely, note that
\eqref{exa.seg.adj} and \eqref{seg.red} actually make sense for an
arbitrary family of groupoids $\C \to I \times \Delta$.

\begin{defn}\label{w.2.def}
A {\em weak $2$-family of groupoids} over a category $I$ is a family
of groupoids $\C \to I \times \Delta$, with reduction
\eqref{seg.red}, such that for any $n \geq 1$, the functor
\begin{equation}\label{ws.0.eq}
  \C_{[n]} \to \C^{red}_{[1]} \times_{\C_{[0]}} \C_{[n-1]}
\end{equation}
induced by \eqref{seg.del.sq} with $l=1$ is an epivalence. A weak
$2$-family of groupoids $\C$ over $I$ is {\em bounded} if $s^*
\times t^*:\C_{[1]} \to \C_{[0]} \times C_{[0]}$ is small.
\end{defn}

\begin{lemma}\label{w.2.le}
For any weak $2$-family of groupoids $\C$ over a category $I$, the
reduction $\C^{red}$ is a reduced $2$-family of groupoids, bounded
resp.\ reflexive if so is $\C$, and if $\C$ is reflexive, then
$\C^{red}$ is complete if and only if so is $\C$. If this happens,
then $\C$ is localizable with respect to the class of special maps,
and $h^+(\C) \to h^+(\C^{red})$ is an equivalence. Moreover,
$h^+(-)$ sends semicartesian squares of bounded reflexive complete
weak $2$-families over $[1]^2$ to semicartesian squares.
\end{lemma}

\proof{} By definition, we have an epivalence $\C_{[n]} \to
\C^{red}_{[n]}$ for any $[n] \in \Delta$, it is an equivalence for
$n=1$, and then the epivalence \eqref{ws.0.eq} gives an epivalence
$\C_{[n]} \to \C_{[1]}^{red} \times_{\C_{[0]}} \C^{red}_{[n-1]}$
that must factor through an equivalence $\C^{red}_{[n]} \to
\C_{[1]}^{red} \times_{\C_{[0]}} \C^{red}_{[n-1]}$, so that
$\C^{red}$ is a Segal category by induction on $n$. Since $\C_{[1]}
\to \C_{[1]}^{red}$ is an epivalence, $\C$ is bounded if and only if
so is $\C^{red}$. If $\C$ is reflexive, then $e^*:\C^{red}_{[0]}
\cong \C_{[0]} \to \C_{[1]} \to \C^{red}_{[1]}$ is full, and since
it has a section, it is also faithful, so that $\C^{red}$ is
reflexive over $[1]$, and then over any $[n]$ by the Segal
condition. If $\C$, hence also $\C^{red}$ is reflexive, then $e^*
\times e^*$ in \eqref{exa.compl.sq} is fully faithful, so the square
is cartesian if and only if it is weakly semicartesian, and by
Lemma~\ref{car.epi.le}, this holds for $\C$ if and only if it holds
for $\C^{red}$. If it does hold, then $h^+(\C^{red})$ exists by
Proposition~\ref{exa.seg.prop}, and $h^+(\C) \cong h^+(\C^{red})$ by
Lemma~\ref{seg.loc.le}.

Finally, for the last claim, assume given a semicartesian square of
weak $2$-families over $[1]^2$, and as in \eqref{cat.dia}, denote by
$\C,\C_0,\C_1,\C_{01} \to \Delta$ its fibers over $0 \times 0$, $1
\times 0$, $0 \times 1$, $1 \times 1$. Then we have a commutative
diagram
$$
\begin{CD}
  \C_{01} @>>> \C_{01}^{red}\\
  @VVV @VVV\\
  \C_0 \times_\C \C_1 @>>> \C^{red}_0 \times_{\C^{red}} \C^{red}_1
\end{CD}
$$
of families of groupoids over $\Delta$, the top and the left arrows
are $2$-full, and the bottom arrow is $1$-full by
Lemma~\ref{2.epi.le}, so that the right arrow is also $1$-full. But
since the simplicial expansion $\Delta^\flat(-)$ commutes with
cartesian products, we have $\C^{red}_0 \times_{\C^{red}} \C^{red}_1
\cong \Delta^\flat(h^+(\C_0) \times_{h^+(\C)} h^+(\C_1))$, so that
$h^+(\C_0) \times_{h^+(\C)} h^+(\C_1) \cong h^+(\C^{red}_0
\times_{\C^{red}} \C^{red}_1)$ by Proposition~\ref{exa.seg.prop},
and the claim immediately follows from Lemma~\ref{exa.loc.le}.
\endproof

\subsection{Fibrant simplicial sets.}\label{fibr.simp.subs}

Finally, let us drop the Segal condition altogether and consider
arbitrary fibrations $\C \to \Delta$. We will need some results for
the case when the fibration is discrete. In fact, we restrict our
attention to small discrete fibrations, or equivalently, to
categories of simplices $\Delta X$ of simplicial sets $X \in
\Delta^o\Sets$. Recall that for any integers $n \geq l \geq 0$, we
have the {\em $l$-th horn} $\V^l_n \in \Delta^o\Sets$, with the horn
embedding \eqref{horn.eq}, and we also have the simplicial sphere
$\SS_{n-1}$, with the sphere embedding \eqref{sph.eq}. In terms of
\eqref{sph.I.eq}, we have $\SS_{n-1} = \SS_{[n]}$ and
$\sigma_n=\sigma_{[n]}$, $[n] \in \Delta$. A map $f:X' \to X$ of
simplicial sets is a {\em Kan fibration} if it has the right lifting
proprety with respect to the horn embeddings \eqref{horn.eq} for all
$n \geq 1$ --- that is, for any $0 \leq l \leq n$ and commutative
diagram
$$
\begin{CD}
\V^l_n @>{a}'>> X'\\
@V{v^l_n}VV @VV{f}V\\
\Delta_n @>{a}>> X,
\end{CD}
$$
there exists a map $b:\Delta_n \to X'$ such that $a=f \circ b$ and
$a'=b \circ v^l_n$. A map $f$ is a {\em trivial Kan fibration} if it has
the same lifting property with respect to the sphere embeddings
\eqref{sph.eq} for all $n \geq 0$. A simplicial set $X$ is {\em
  fibrant} resp.\ {\em trivially fibrant} if the tautological map $X
\to \ppt$ is a Kan fibration resp.\ a trivial Kan fibration.

\begin{lemma}\label{kan.le}
For any fibrant simplicial set $X$ with the category of simplices
$\Delta X$, the groupoid $h^{\Tot}(\Delta X)$ is equivalent to the
following groupoid:
\begin{enumerate}
\item objects are elements $x \in X([0])$,
\item morphisms from $x$ to $x'$ are equivalence classes of elements
  $y \in X([1])$ such that $s^*y=x$ and $t^*y=x'$, where
\item two elements $y,y' \in X([1])$ lie in the same equivalence
  class iff there exists $z \in X([2])$ that restricts to $y$
  resp.\ $y'$ on $\{0,1\} \subset \{0,1,2\}$ resp.\ $\{0,2\} \subset
  \{0,1,2\}$, and to a degenerate simplex on $\{1,2\} \subset
  \{0,1,2\}$.
\end{enumerate}
For any trivial Kan fibration $f:X' \to X$ between fibrant
simplicial sets, the functor $h(f):h^{\Tot}(\Delta X') \to
h^{\Tot}(\Delta X)$ is an equivalence.
\end{lemma}

\proof{} Since for any $[n] \in \Delta$, there exists a non-trivial
map $[0] \to [n]$, all objects in $h^{\Tot}(\Delta X)$ can indeed be
represented by objects $\langle [0],x \rangle \in (\Delta X)_{[0]}$,
$x \in X([0])$, and moreover, any map between two such objects can
be represented by a diagram \eqref{zig} such that the sources of
$f_i$ and $w_i$ are in $(\Delta X)_{[0]}$. Since for any two maps
$a,b:[0] \to [n]$ there exists a map $c:[1] \to [n]$ such that $a =c
\circ s$ and $b = c \circ t$, one can further assume that the
targets of $f_i$ and $w_i$ are in $(\Delta X)_{[1]}$. Then the
lifting property for the horn $\V^1_2$ shows that any zigzag of
length $4$ reduces to one of length $2$, so that by induction, every
map is of the form $t_y^{-1} \circ s_y$ for some liftings
$t_y:\langle [0],x' \rangle \to \langle [1], y \rangle$,
$s_y:\langle [0],x \rangle \to \langle [1],y \rangle$ of the maps
$s,t:[0] \to [1]$. Moreover, sending $\langle [n],x \rangle$ to
$s^*x$ provides a full essentially surjective functor from
$h^{\Tot}(\Delta X)$ to the groupoid defined by \thetag{i},
\thetag{ii}, \thetag{iii}, and \thetag{iii} is indeed an equivalence
relation, while sending $x$ to $\langle [0],x \rangle$ and $y$ to
$t_y^{-1} \circ s_y$ gives a full essentially surjective inverse.
Finally, if $f:X' \to X$ is a trivial Kan fibration, then
\thetag{i}, \thetag{iii}, \thetag{ii} immediately imply that $h(f)$
is essentially surjective, faithful and full.
\endproof

\begin{corr}\label{semi.corr}
Assume given two maps $f_0:X_0 \to X$, $f_1:X_1 \to X$ of fibrant
simplicial sets, and assume that $f_0$ is a Kan fibration. Then the
natural functor
\begin{equation}\label{semi.eq}
h^{\Tot}(\Delta (X_0 \times_X X_1)) \to h^{\Tot}(\Delta X_0)
\times_{h^{\Tot}(\Delta X)} h^{\Tot}(\Delta X_1)
\end{equation}
is an epivalence.
\end{corr}

\proof{} To see that the functor \eqref{semi.eq} is essentially
surjective, note that by Lemma~\ref{kan.le}~\thetag{i},\thetag{ii},
an object in its target can be represented by a triple $\langle
x_0,x_1,y \rangle$ of elements $x_0 \in X_0([0])$, $x_1 \in
X_1([0])$, $y \in X([1])$ such that $f_0(x_0) = X(s)(y)$ and
$f_1(x_1)=X(t)(y)$. However, the embedding $s:[0] \to [1]$ is a horn
embedding, and since $f_0$ is a Kan fibration, we can modify $x_0$
so that $f_0(x_0)=X(t)(y)$, and we obtain an isomorphic object. But
then $f_0(x_1)=f_1(x_1)$, so that $x_0 \times x_1 \in (X_0 \times_X
X_1)([0])$ defines an object in the source of \eqref{semi.eq}.

Analogously, to see that \eqref{semi.eq} is full, note that by
Lemma~\ref{kan.le}~\thetag{ii} and \thetag{iii}, morphisms in the target
are represented triples $\langle y_0,y_1,z \rangle$, $y_l \in X_l([1])$,
$l=0,1$, $z \in X([2])$ restricts to $f_0(y_0)$ resp.\ $f_1(y_1)$ on
$\{0,1\}$ resp.\ $\{ 0,2 \}$, and use the liting property for the
horn embedding $v^0_2:V^0_2 \to \Delta_2$ to choose an equivalent
triple with $f_0(y_0)=f_1(y_1)$.
\endproof

\section{Model structures.}\label{mod.sec}

\subsection{Model categories.}

In this paper, a {\em model category} is what was called a closed
model category in \cite{qui.ho}, with no extra assumptions. By
definition, it is given by a relative category $\langle \C,W
\rangle$ equipped with two extra classes of maps $C$, $F$ satisfying
the usual axioms. We assume known standard terminology of
\cite{qui.ho} such as the notions of fibrant and cofibrant objects
(we will avoid using ``fibrations'' and ``cofibrations'' so as not
to create a clash of terminology with \cite{SGA}).

We will also assume known standard results about model
categories. In particular, the opposite $\C^o$ of a model category
$\C$ is a model category, with $C=F^o$, $F=C^o$, $W=W^o$. For any
object $c \in \C$ in a model category $\C$, the left comma-fiber
$\C/c$ is also a model category, with the classes $C$, $W$, $F$
induced from $\C$ via the forgetful functor $\sigma(c):\C \to
\C/c$. By duality, the same is true for the right comma-fiber $c
\backslash \C$.

It is a well-known property of the formalism of model categories ---
in fact, it is the main reason for its existence --- that for
any model category $\langle \C,C,W,F \rangle$, the underlying
relative category $\langle \C,W \rangle$ is localizable, and in a
very controlled way. Any map from $c' \in \C$ to $c \in \C$ in
$h^W(\C)$ can be represented by diagram \eqref{zig} of length at most
$3$, and one can take as $w_0:\wt{c}' \to c'$, $w_1:c \to \wt{c}$
fixed cofibrant resp.\ fibrant replacements of $c'$ resp.\ $c$. If
$c'$ itself is cofibrant, one can take $w_0=\id$, and if $c$ is
fibrant, one can take $w_1=\id$. If $c'$ is cofibrant and $c$ is
fibrant, then any map $c' \to c$ in $h^W(\C)$ can be represented by
a single map $f:c' \to c$. If one chooses a factorization
\begin{equation}\label{c.dec}
\begin{CD}
  c' \copr c' @>{i}>> \wt{c} @>{w}>> c'
\end{CD}
\end{equation}
of the co-diagonal map $c' \copr c' \to c'$ such that $w \in F \cap
W$ and $i \in C$, then two maps $f,f':c' \to c$ represent the same
map in $h^W(\C)$ iff $f \copr f':c' \copr c' \to c$ factors as $f
\copr f' = \wt{f} \circ i$ for some $\wt{f}:\wt{c} \to c$ (the maps
$f$, $f'$ are then called {\em homotopic}, and $\wt{f}$ is the
``homotopy'' between them). The localization functor $f:\C \to
h^W(\C)$ satisfies
\begin{equation}\label{h.W.C}
  h^*(\Iso) = W,
\end{equation}
that is, $h(f)$ is invertible iff $f \in W$. If we let
$\Cof(\C),\Fib(\C) \subset \C$ be the full subcategories spanned by
cofibrant resp.\ fibrant objects, then both $\langle \Cof(\C),W
\rangle$ and $\langle \Fib(\C),W \rangle$ are also localizable, and
the natural embeddings induce equivalences
\begin{equation}\label{mod.loc.eq}
  h^W(\Cof(\C)) \cong h^W(\C) \cong h^W(\Fib(\C))
\end{equation}
between the localizations. The main application of
\eqref{mod.loc.eq} is Quillen's Adjunction Theorem; we recall one of
its forms.

\begin{theorem}\label{qui.adj.thm}
Assume given model categories $\C$, $\C'$, and a pair of a functor
$\gamma:\C \to \C'$ and its left-adjoint $\gamma_\dg:\C' \to
\C$. Moreover, assume that $\gamma$ sends maps in $F$ resp.\ $F \cap
W$ to maps in $F$ resp.\ $F \cap W$, or equivalently, that
$\gamma^\dg$ sends maps in $C$ resp.\ $C \cap W$ to maps in $C$
resp.\ $C \cap W$. Then $\gamma|_{\Fib(\C)}$ and
$\gamma^\dg|_{\Cof(\C')}$ send maps in $W$ to maps in $W$, and the
induced pair of functors
$$
\begin{aligned}
  &R^\hdot\gamma:h^W(\C) \cong h^W(\Cof(\C)) \to h^W(\C'),\\
  &L^\hdot\gamma^\dg:h^W(\C') \cong h^W(\Fib(\C')) \to h^W(\C)
\end{aligned}
$$
is adjoint. Moreover, if $\gamma$ resp.\ $\gamma_\dg$ is fully
faithful, then so is the corresponding functor $R^\hdot\gamma$
resp.\ $L^\hdot\gamma_\dg$.\endproof
\end{theorem}

\begin{prop}\label{glue.prop}
In the situation of Theorem~\ref{qui.adj.thm}, consider the category
$\C' \setminus_\gamma \C$, and say that a map $\langle f,f' \rangle:
\langle c_0,c_0',\alpha_0 \rangle \to \langle c_1,c_1',\alpha_1
\rangle$ is in $F$ if so are $f$ and $f' \times \alpha_0:c_0' \to
c_1' \times_{\gamma(c_1)} \gamma(c_0)$, and in $C$ resp.\ $W$ iff so
are $f'$ and $f$. Then this is a model structure on $\C'
\setminus_\gamma \C$, and the functor
\begin{equation}\label{comma.mod.eq}
h:h^W(\C' \setminus_\gamma \C) \cong h^W(\Fib(\C' \setminus_\gamma
\C)) \to h^W(\C') \setminus_{R^\hdot\gamma} h^W(\C)
\end{equation}
induced by $h \setminus h:\Fib(\C' \setminus_\gamma \C) \subset
\Fib(\C') \setminus_\gamma \Fib(\C) \to h^W(\C')
\setminus_{R^\hdot\gamma} h^W(\C)$ is an epivalence.
\end{prop}

\proof{} The fact that we have a model structure on $\C'
\setminus_\gamma \C$ is very standard (see e.g.\ \cite[Proposition
  3.8]{ka.glue} for a more general statement). To see that
\eqref{comma.mod.eq} is essentially surjective, note that by
assumption, $\gamma$ sends fibrant objects to fibrant objects, so
that any object in its target can be represented by a triple
$\langle c,c',\alpha \rangle$ with fibrant $c$ and cofibrant $c'$
(and some $\alpha:c' \to \gamma(c)$). Since \eqref{comma.mod.eq} is
obviously conservative by virtue of \eqref{h.W.C}, it remains to
check that it is full. To do this, assume given a map $\langle
f_h,f'_h \rangle:c^h_0 \to c^h_1$ in $h^W(\C')
\setminus_{R^\hdot\gamma} h^W(\C)$. We need to check that for some
objects $\langle c_0,c_0',\alpha_0 \rangle$ and $\langle
c_1,c_1',\alpha_1 \rangle$ in $\C' \setminus_\gamma \C$ with
$h(\langle c_l,c'_l,\alpha_l \rangle) \cong c^h_l$, $l=0,1$, the map
$\langle f_h,f'_h \rangle$ lifts to a map $\langle f,f'
\rangle:\langle c_0,c_0',\alpha_0 \rangle \to \langle
c_1,c_1',\alpha_1 \rangle$. To do this, first choose some
representatives $\langle c_l,c'_l,\alpha_l \rangle$ for $c^h_l$,
$l=0,1$, and then take replacements so that both are fibrant and
cofibrant. Then $c_0$, $c_0'$ are cofibrant, $c_1$ is fibrant, and
$\alpha_1$ is in $F$, so that $c_1'$ is also fibrant. Therefore we
can represent maps $f_h$, $f'_h$ by maps $f:c_0 \to c_1$, $f':c_0
\to c'_1$ such that the diagram
\begin{equation}\label{comma.mod.sq}
\begin{CD}
  c_0' @>{f'}>> c_1'\\
  @V{\alpha_0}VV @VV{\alpha_1}V\\
  \gamma(c_0) @>{\gamma(f)}>> \gamma(c_1)
\end{CD}
\end{equation}
becomes commutative after applying the localization functor $h$ for
$\C'$ --- that is, we have $h(\gamma(f) \circ \alpha_0) = h(\alpha_1
\circ f')$. Then if we fix a decomposition \eqref{c.dec} for $c_0'$,
with some middle term $\wt{c}$, the two composition maps in
\eqref{comma.mod.sq} are related by a homotopy $\wt{g}:\wt{c} \to
\gamma(c_1)$. We can now replace $c'_1$ with $\wt{c} \copr_{c'_0}
c'_0$, with coproduct taken with respect to one of the two maps
$c'_0 \to \wt{c}$, replace $f'$ with the other map $c'_0 \to
\wt{c}$, and replace $\alpha_1$ with $\wt{g} \copr \alpha_1$. Then
\eqref{comma.mod.eq} becomes commutative in $\C' \setminus_\gamma
\C$ and gives the desired lifting.
\endproof

\begin{remark}
The generality in \cite{ka.glue} is as follows: $\gamma:\C \to \C'$
sends maps in $F$ resp.\ $F \cap W$ to maps in $F$ resp.\ $F \cap
W$, but is not required to have an adjoint, nor even to preserve
finite limits. In this setting, it is not even immediately obvious
--- but true --- that $\C' \setminus_\gamma \C$ is finitely
cocomplete.
\end{remark}

It is well-known that any two of the classes $C$, $W$, $F$ in a
model structure uniquely define the third one, so for example, to
describe the model structure of Proposition~\ref{glue.prop}, it
suffices to only describe $C$ and $W$. Here is another well-known
example of such a decription.

\begin{defn}\label{proj.mod.def}
For any model category $\C$ and essentially small category $I$, the
{\em projective model structure} on $\Fun(I,\C)$ is the one where a map
$f$ lies in $W$ resp.\ $F$ iff $f(i)$ is in $W$ resp.\ $F$ for any
$i \in I$.
\end{defn}

The projective model structure is unique, if it exists. Existence
requires imposing conditions either on $\C$ or on $I$. In
particular, it always exists if $I$ is a finite partially ordered
set: this is an easy example of the following general existence
theorem for Reedy categories.

\begin{theorem}\label{reedy.thm}
Assume given a Reedy category $I$ in the sense of
Definition~\ref{reedy.def}, and a model category $\langle \C,C,W,F
\rangle$. Assume that either $\C$ is complete and cocomplete, or $I$
is locally finite, or $I$ is Reedy-$\kappa$-bounded for some regular
cardinal $\kappa$, and $\C$ is $\kappa$-cocomplete and
$\kappa'$-complete for any $\kappa' \ll \kappa$. Say that a morphism
$f:X \to X'$ between two functors $X,X':I \to \C$ is in $W$ iff
$f(i):X(i) \to X'(i)$ is in $W$ for any $i \in I$, say that $f$ is
in $C$ iff for any $i \in I$ with the decompositions \eqref{lm.eq},
the map $L(X',i) \copr_{L(X,i)} X(i) \to X'(i)$ is in $C$, and say
that $f$ is in $F$ if for any $i \in I$, the map $X(i) \to X'(i)
\times_{M(X',i)} M(X,i)$ is in $F$. Then $\langle \Fun(I,\C),C,W,F
\rangle$ is a model category.\endproof
\end{theorem}

For a proof of this, see e.g.\ \cite{hovey} (keeping in mind that
although the author allows himself to redefine the notion of a model
category, the proof works {\em verbatim} in the original setting of
\cite{qui.ho}). The main idea is an iterated application of the
gluing construction given in Proposition~\ref{glue.prop}.

\begin{exa}\label{reedy.proj.exa}
If the Reedy category $I$ in Theorem~\ref{reedy.thm} is directed ---
that is, $I=I_L$ and $M=\Id$ --- then the Reedy model structure is
the projective model structure of Definition~\ref{proj.mod.def}
(and in particular, the latter exists). For example, this always
happens when $I=J$ is a left-bounded partially ordered set, with the
Reedy structure of Example~\ref{reedy.exa}.
\end{exa}

\begin{exa}\label{reedy.dim1.exa}
Assume given a finite partially ordered set $J$ of dimension $\dim J
= 1$. Let $J_l \subset J$ be the subset of elements of height $l$,
$l=0,1$, and let $J' \subset J_0 \times J_1$ be the set of
non-trivial order relations $j < j'$, with the projections $\pi_l:J'
\to J_l$, $l=0,1$. Then we have a Quillen-adjoint pair of functors
$\gamma = \pi_{1*} \circ \pi_0^*:J_0^o\C \to J_1^o\C$, $\gamma_\dg =
\pi_{0!} \circ \pi_1^*:J_1^o\C \to J_0^o\C$, and $J^o\C \cong
J_1^o\C \setminus_\gamma J_0^o\C \cong J_1^o\C /_{\gamma^\dg}
J_0^o\C$. The model structure of Proposition~\ref{glue.prop} is the
Reedy model structure on $J^o\C$ corresponding to the Reedy
structure on $J^o$ with $(J^o)_M=J^o$.
\end{exa}

\begin{lemma}\label{J.proj.le}
In the situation of Theorem~\ref{reedy.thm}, for any $i \in I$ and
cofibrant $X:I \to \C$, $X(i) \in \C$ is cofibrant, and for any map
$f:X \to X'$ in $C$, the map $f(i):X(i) \to X(i')$ is in
$C$. Moreover, if $I=J$ is a finite partially ordered set, with the
Reedy structure such that $J_L=J$, then for any cofibrant $X:J \to
\C$ and $j' \leq j$ in $J$, the map $X(j') \to X(j)$ is in $C$.
\end{lemma}

\proof{} For the first claim, note that the definition of the class
$C$ given in Theorem~\ref{reedy.thm} only depends on $I_L$, so we
may assume that $I=I_L$ is directed, as in
Example~\ref{reedy.proj.exa}. Then the evaluation $\ev_i:\Fun(I,\C)
\to \C$ has a right-adjoint given by the Kan extension $\eps(i)_*$,
and by \eqref{kan.eq}, it sends maps in $F$ resp.\ $F \cap W$ to
maps in $F$ resp.\ $F \cap W$, since both classes are closed under
products. Thus $X(i)$ is cofibrant resp.\ $f(i)$ is in $C$ by
adjunction.

For the second claim, let $\eps:[1] \to J$ be the embedding sending
$0$ to $j'$ and $1$ to $j$, and note that again by \eqref{kan.eq},
the right Kan extension $\eps_*:\Ar(\C) \to \Fun(J,\C)$ is given by
$$
\eps_*Y(j'')= \begin{cases}
Y(0), &\quad j'' \in J/j',\\
Y(1), &\quad j'' \in (J/j) \ssetminus (J/j'),\\
1, &\quad j'' \in J \ssetminus (J/j),
\end{cases}
$$
for any $Y:[1] \to \C$.  This functor again sends maps in $F$
resp.\ $W$ to maps in $F$ resp.\ $W$, so by adjunction, $\eps^*(X)$
is cofibrant in $\Ar(\C)$, and the claim reduces to the case
$J=[1]$. This is the situation of Proposition~\ref{glue.prop}.
\endproof

\begin{corr}\label{mod.corr}
For any model category $\C$ and object $c$, let $\C^W/c \subset
\C/c$ be the subcategory of cofibrant objects $c'$ equipped with
maps $f:c' \to c$ in the class $W$, and with maps between them that
are in the class $C$. Then $\C^W/c$ is homotopy filtered in the
sense of Definition~\ref{ho.filt.def}.
\end{corr}

\proof{} Assume given a functor $\gamma:J \to \C^W/c$ from a finite
partially ordered set $J$. Let $\gamma':J \to \C/c$ be a cofibrant
replacement of $\gamma$ in $(\C/c)^J$ with the projective model
structure, so that $\gamma'$ is cofibrant in $(\C/c)^J$ and equipped
with a map $f:\gamma' \to \gamma$ such that $f(j) \in W$ for any $j
\in J$. Since $W$ is closed under compositions,
Lemma~\ref{J.proj.le} implies that $\gamma'(j)$ lies in $\C^W/c
\subset \C/c$ for any $j \in J$, and moreover, maps $\gamma'(j) \to
\gamma(j')$, $j \leq j'$ are in the class $C$. Therefore $\gamma'$
factors through $\C^W/c$. Since $\C/c$ has finite coproducts,
$\gamma'$ extends to a functor $\gamma'':J^> \to \C/c$ with
$\gamma''(1) = \colim_J\gamma'$, and since $\gamma'$ is cofibrant in
$(\C/c)^J$, the extension $\gamma''$ is cofibrant in $(\C/c)^{J^>}$
by adjunction. In particular, $\gamma''(1)$ is cofibrant. Now
factorize the map $g:\gamma''(1) \to c$ as $f = c \circ w$,
$c:\gamma''(1) \to c'$ in $C$, $w:c' \to c$ in $W$, and let
$\wt{\gamma}:J^> \to \C/c$ be given by $\wt{\gamma}|_J = \gamma'$,
$\wt{\gamma}(1) = c'$. The functor $\wt{\gamma}$ then factors
through $\C^W/c \subset \C/c$ and provides the desired extension of
$\gamma'$.
\endproof

Dually, for any essentially small category $I$ and model category
$\C$, the {\em injective} model structure on $I^o\C = \Fun(I^o,\C)
\cong \Fun(I,\C^o)^o$ is dual to the projective model structure on
$\Fun(I,\C^o)$; explicitly, a map is in $W$ resp.\ $C$ with respect
to the injective model structure iff it is pointwise in $W$
resp.\ $C$. By Theorem~\ref{reedy.thm}, the injective model
structure always exists if $I=J$ is a left-finite partially ordered
set (where we treat $J^o$ as a Reedy category with $(J^o)_M = J^o$
and discrete $(J^o)_L$). A functor $X:J^o \to \C$ is fibrant with
respect to the injective model structure iff the map $X(j) \to
M(X,j)$ is in $F$ for any $j \in J$. By the universal property of
localizations, the evaluation functor \eqref{ev.eq} for $J^o\C$ fits
into a commutative square
\begin{equation}\label{ev.loc.sq}
\begin{CD}
  J^o \times J^o\C @>{\ev}>> \C\\
  @VVV @VVV\\
  J^o \times h^W(J^o\C) @>>> h^W(\C),
\end{CD}
\end{equation}
and the bottom arrow gives rise to a comparison functor
\begin{equation}\label{mod.ev}
  h^W(J^o\C) \to J^oh^W(\C).
\end{equation}
The functor \eqref{mod.ev} is conservative by virtue of
\eqref{h.W.C} but usually not an equivalence. If $J$ is finite and
$\dim J=1$, it is an epivalence by Example~\ref{reedy.dim1.exa} and
Proposition~\ref{glue.prop}.

\subsection{Cofibrant objects.}\label{repr.subs}

Assume given a model category $\langle I,C,W,F \rangle$, and use
\eqref{mod.loc.eq} to idenfity $h^W(I) \cong h^W(\Cof(I))$, so that
any object in $h^W(I)$ can be represented by a cofibrant $i \in
\Cof(I)$, and any morphism from $i' \in \Cof(I)$ to $i \in \Cof(I)$
in $h^W(I)$ can be represented by a diagram
\begin{equation}\label{codom}
\begin{CD}
i' @>{f}>> \wt{i} @<{w}<< i
\end{CD}
\end{equation}
with $w \in C \cap W$. Since $i$ is cofibrant, $\wt{i}$ is also
cofibrant. What is slightly less well-known is that one can also
arrange for $f$ to be in $C$. Namely, let $\Cof(I)_C = \Cof(I) \cap
I_C \subset \Cof(I)$ be the dense subcategory defined by the closed
class $C$, and by abuse of notation, let $C \cap W$ be the class of
maps in $\Cof(I)_C$ that are also in $W$. Denote by
\begin{equation}\label{eps.I}
\eps:\langle \Cof(I)_C,C \cap W \rangle \to \langle I,W \rangle
\end{equation}
the natural embedding (it is tautologically a morphism of relative
categories).

\begin{lemma}\label{mod.le}
The relative category $\langle \Cof(I)_C,C \cap W \rangle$ is
localizable, and the morphism \eqref{eps.I} is a $1$-equivalence in
the sense of Definition~\ref{equi.0.def}.
\end{lemma}

\proof{} For any two cofibrant objects $i,i' \in I$, denote by
$QI(i',i)$ the category of diagrams \eqref{codom}, and let
$QI_C(i',i) \subset QI(i',i)$ be the subcategory spanned by diagrams
such that $f \in C$, and morphisms between them that are pointwise
in $C$. It is well-known that two diagrams $d_0,d_1 \in QI(i',i)$
represent the same morphism in $h^W(I)$ if and only if they lie in
the same connected component of $QI(i',i)$, and even stronger, if
and only if there exists a third diagram $d_{01} \in QI(i_0,i_1)$
and maps $d_0 \to d_{01}$, $d_1 \to d_{01}$. On the other hand, for
the relative category $\langle \Cof(I)_C,C \cap W \rangle$, we note
that $\Cof(I) \subset I$ is closed under pushouts with respect to
maps in $C$, so that every zigzag diagram \eqref{zig} in $\langle
\Cof(I)_C,C \cap W \rangle$ is equivalent to a diagram $d \in
QI_C(i',i)$ by the cancellation rules. Moreover, every map in the
category $QI_C(i',i)$ is automatically pointwise in $C \cap W$, so
that by the same cancellation rules, $d_0,d_1 \in qI_C(i',i)$ are
equivalent if they are connected by a chain of maps. Thus to prove
that $\langle \Cof(I)_C,C \cap W \rangle$ is localizable and
$h(\eps)$ is fully faithful, it suffices to prove that for any $i,i'
\in \Cof(I)$, the functor
\begin{equation}\label{eps.i01}
\eps_{i',i}:QI_C(i',i) \to QI(i',i)
\end{equation}
induced by \eqref{eps.I} gives a bijection between connected
components of the two categories.

To prove that \eqref{eps.i01} gives a surjection, note that every
diagram $d$ in $QI(i',i)$ also defines an object $\wt{d}$ in the
model category $(i' \copr i)\backslash I$, and for any cofibrant
replacement $\wt{d}' \to \wt{d}$ of the object $\wt{d}$, the object
$\wt{d}'$ corresponds to a diagram $d \in QI_C(i',i)$.

To prove that \eqref{eps.i01} gives an injection, assume given two
diagrams $d_0,d_1 \in QI_C(i',i)$, and assume that $\eps(d_0)$ and
$\eps(d_1)$ represent the same map in $h^W(I)$. Then there exists a
diagram $d \in QI(i',i)$ and morphisms $w_0:\eps(d_0) \to d$,
$w_1:\eps(d_1) \to d$, and if we let $\wt{d} \in (i' \copr
i)\backslash I$ be the object corresponding to $d$, then the maps
$w_0$ and $w_1$ define two objects in the category $((i' \copr
i)\backslash I)^W/\wt{d}$ of Corollary~\ref{mod.corr}. But by
Proposition~\ref{filt.prop}, this category is $2$-connected, thus
connected by Lemma~\ref{2.1.eq.le}.

To finish the proof, it remains to observe that every object in $I$
has a cofibrant replacement, so that $h(\eps)$ is essentially
surjective.
\endproof

Using Corollary~\ref{mod.corr} in the proof of Lemma~\ref{mod.le}
looks like a huge overkill (and it provides much more than needed
for the argument). What we really need it for is the following
refinement of Lemma~\ref{mod.le}.

\begin{prop}\label{mod.prop}
For any model category $\langle I,C,W,F \rangle$, the morphism
$\eps$ of \eqref{eps.I} is a $2$-equivalence in the sense of
Definition~\ref{equi.def}.
\end{prop}

\proof{} Consider the comma-category $\Cof(I)_C /_\eps I$, and let
$\overline{I} \subset \Cof(I)_C /_\eps I$ be the full subcategory
spanned by triples $\langle i_0,i_1,\alpha \rangle$ such that
$\alpha:\eps(i_0) \to i_1$ is in the class $W$. We have the
projections $\tau:\overline{I} \to I$, $\sigma:\overline{I} \to
\Cof(I)_C$, and the embedding $\eta:\Cof(I)_C \to \Cof(I)_C /_\eps
I$ factors through $\overline{I} \subset \Cof(I)_C /_\eps I$, thus
gives a fully faithful embedding $\overline{\eta}:\Cof(I)_C \to
\overline{I}$ left-adjoint to the projection $\sigma$.

For any map $f$ in $\overline{I}$, we have $\tau(f) \in W$ if and
only if $\sigma(f) \in C \cap W$; denote by $\overline{W} = \tau^*W
= \sigma^*(C \cap W)$ the class of maps satisfying either of these
two equivalent conditions. Then $\overline{\eta}$ and $\sigma$ are
morphisms between $\langle \Cof(I)_C,C \cap W \rangle$ and $\langle
\overline{I},\overline{W} \rangle$, and $\tau$ is a morphism from
$\langle \overline{I},\overline{W} \rangle$ to $\langle I,W
\rangle$. Moreover, $\sigma \circ \overline{\eta} \cong \id$, and
the adjunction map $\overline{\eta} \circ \sigma \to \id$ is in
$\overline{W}$ by the definition of the category
$\overline{I}$. Therefore $\eta$ is a $2$-equivalence by
Lemma~\ref{loc.adj.le}, and it suffices to prove that $\tau$ is also
a $2$-equivalence.

But by Lemma~\ref{2.con.bis.le}~\thetag{i}, it then suffices to
prove that $\tau:\overline{I} \to I$ is $2$-connected, and by
Proposition~\ref{reso.prop}, it suffices to observe that it is locally
$2$-connected in the sense of Definition~\ref{loc.con.def}. Indeed,
for any object $i \in I$, the category $\overline{I}_i$ of
Definition~\ref{loc.con.def}~\thetag{i} is precisely the category
$I^W/i$ of Corollary~\ref{mod.corr}, and it is $2$-connected by
Proposition~\ref{filt.prop}. On the other hand, for any morphism
$f:i \to i'$ in $I$ and any object $\wt{i} \in \overline{I}_i$
represented by a morphism $w:c \to i$ in $W$ with cofibrant $c \in
I_C$, the category $\overline{I}(f,\wt{i})$ of
Definition~\ref{loc.con.def}~\thetag{ii} is the category of
factorizations $g = w' \circ f'$ of $g=f \circ w:c \to i'$ with
$f':c \to c'$ in the class $C$ and $w':c' \to i'$ in the class $W$,
with morphisms from $\langle w'_0,f'_0,c'_0 \rangle$ to $\langle
w'_1,f'_1,c'_1 \rangle$ given by maps $a:c'_0 \to c'_1$ in the class
$C$ such that $a \circ f'_0 = f'_1$ and $w'_1 \circ a = w'_0$. Then
$g:c \to i'$ defines an object $\wt{g}$ in the model category $c
\backslash I$, and $\overline{I}(f,\wt{i})$ is exactly the category
$(c \backslash I)^W/\wt{g}$ of Corollary~\ref{mod.corr}, again
$2$-connected by Proposition~\ref{filt.prop}.
\endproof

Proposition~\ref{mod.prop} shows that for any model category
$\langle I,C,W,F \rangle$, one loses no information when restricting
a family of groupoids over $\langle I,W \rangle$ to $\langle
\Cof(I)_C,C \cap W \rangle$. We will now use this to construct a
version of the Yoneda embedding for such families. Take two objects
$i,i' \in \Cof(I)$, and consider the category $QI_C(i',i)$ of
diagrams \eqref{codom} in $I_C$ of Lemma~\ref{mod.le}.

\begin{lemma}\label{i.loc.le}
The category $QI_C(i',i)$ is localizable.
\end{lemma}

\proof{} If for any model category $\langle I,C,W,F \rangle$ we
consider the dense subcategory $\Cof(I)_{C \cap W} \subset
\Cof(I)_C$, then the same argument as in Lemma~\ref{mod.le} shows
that $\Cof(I)_{C \cap W}$ is localizable. To prove the claim, it
remains to observe that $QI_C(i',i)$ is a union of connected
components of the category $\Cof((i' \copr i)\backslash I))_{C \cap
  W}$.
\endproof

Now fix one object $i \in I$, and let $QI_C(i)$ be the category of
diagrams \eqref{codom} in $\Cof(I)_C$ with fixed $i$ and arbitrary
$i'$.  Sending a diagram to $i'$ then provides a forgetful functor
$\phi:QI_C(i) \to \Cof(I)_C$, and this functor is obviously a
fibration with fibers $QI_C(i)_{i'} \cong QI_C(i',i)$. By
Lemma~\ref{i.loc.le} and \eqref{pd.fib}, the trivial family $QI_C(i)
\to QI_C(i)$ is localizable with respect to $\phi$. For any morphism
$w:i' \to i$ in $C \cap W$, composition with $w$ provides a functor
$w^*:QI_C(i) \to QI_C(i')$ cartesian over $I_C$.

\begin{defn}\label{repr.def}
The {\em representable family of groupoids} $\Hh(i)$ is the family
$\Hh(i) = \phi_!QI_C(i)$. The {\em universal object} $o \in
\Hh(i)_i$ is the object $o \in QI_C(i)$ represented by the diagram
\eqref{codom} with $i'=i=\wt{i}$ and $f=w=\id$.
\end{defn}

\begin{prop}\label{repr.prop}
Assume given a cofibrant object $i \in \Cof(I) \subset I$ in a model
category $\langle I,C,W,F \rangle$, and consider the corresponding
representable family $\Hh(i) \to \Cof(I)_C$ of
Definition~\ref{repr.def}.
\begin{enumerate}
\item $\Hh(i)$ is a family of groupoids over $\langle \Cof(I)_C,C
  \cap W \rangle$.
\item For any map $w:i' \to i$ in $C \cap W$, the functor
  $w^*:QI_C(i) \to QI_C(i')$ induces an equivalence $\Hh(i) \cong
  \Hh(i')$ over $I_C$.
\item For any family of groupoids $\C$ over $\langle \Cof(I)_C,C
  \cap W \rangle$ and object $c \in \C_i$, there exists a pair of a
  functor $\Hh(c):\Hh(i) \to \C$ over $I_C$ and an isomorphism
  $\Hh(c)(o) \cong c$, and such a pair is unique up to a unique
  isomorphism.
\end{enumerate}
\end{prop}

\proof{} For \thetag{i}, note that for any map $w:i_0 \to i_1$ in
$I_C$ in the class $C \cap W$, the corresponding transition functor
$w^*:QI_C(i_1,i) \to QI_C(i_0,i)$ of the fibration $QI_C(i) \to
\Cof(I)_C$ has a left-adjoint functor $w_!$ sending a diagram
$$
\begin{CD}
i_0 @>{f_0}>> \wt{i}_0 @<{w_0}<< i
\end{CD}
$$
to the diagram
$$
\begin{CD}
i_1 @>{f_1}>> \wt{i}_0 \copr_{i_0} i_1 @<{w' \circ w_0}<< i,
\end{CD}
$$
where $f_1:i_1 \to i_0' \copr_{i_0} i_1$, $w':i_0' \to i_0'
\copr_{i_0} i_1$ are the natural maps. Therefore the transition
functor $h(w)^*:h^{\Tot}(QI_C(i_1,i)) \to h^{\Tot}(QI_C(i_1,i))$ of the
fibration $\Hh(i) \to I_C$ is an equivalence by
Lemma~\ref{loc.adj.le}.

For \thetag{ii}, exactly the same argument shows that $w^*:QI_C(i)
\to QI_C(i')$ has an adjoint functor $w_!$, and then $w^*$ and $w_!$
descend to a pair of mutually inverse equivalences.

For \thetag{iii}, note that by the definition of localization, a
functor $QI_C(i) \to \C$ over $\Cof(I)_C$ uniquely factors through
$\Hh(i)$, so that to construct a functor $\Hh(c)$ from $\Hh(i)$ to
$\C$, it suffices to construct it on on $QI_C(i)$. This is
equivalent to constructing a section $\sigma$ of the family
$\phi^*\C \to QI_C(i)$. Let $\overline{Q}I_C(i) \subset QI_C(i)$ be
the full subcategory spanned by diagrams \eqref{codom} with
invertible $f$, and let $\overline{\phi}:\overline{Q}I_C(i) \to I_C$
be the restriction of the forgetful functor $\phi$. Then
$\overline{\phi}$ is actually a morphism from $\langle
\overline{Q}I_C(i),{\Tot} \rangle$ to $\langle \Cof(I)_C,C \cap W
\rangle$, so that $\overline{\phi}^*\C$ is a family over $\langle
\overline{Q}I_C(i),{\Tot} \rangle$. In particular,
$\overline{\phi}^*\C \to \overline{Q}I_C(i)$ is a bifibration, so
that $(\overline{\phi}^*\C)^o \to \overline{Q}I_C(i)$ is a
fibration. Since $o \in \overline{Q}I_C(i)$ is obviously the initial
object, \eqref{bnd.equi} shows that $(\overline{\phi}^*\C)^o$ is
bounded over $\overline{Q}I_C(i)^o$ and provides an equivalence
$$
\Sec^{\Tot}(\overline{Q}I_C(i)^o,(\overline{\phi}^*\C)^o) \cong
\C_i^o.
$$
Passing to the opposite categories and recalling that every section
of a family of groupoids is automatically cartesian, we deduce that
any object $c \in \C_i$ uniquely extends to a section
$\overline{\sigma}:\overline{Q}I_C(i) \to \overline{\phi}^*\C$. To
extend it further to a section $\sigma:QI_C(i) \to \phi^*\C$, it
remains to observe that the embedding $\overline{Q}I_C(i) \subset
QI_C(i)$ has a left-adjoint functor, and again apply
\eqref{bnd.equi}.
\endproof

For any cobirant object $i \in \Cof(I) \subset I$ of the model
category $I$, denote by $H(i):\Cof(I)^o_C \to \Sets$ the functor
given by
\begin{equation}\label{H.i.eq}
H(i)(i') = h^{C \cap W}(\Cof(I)_C)(h(i'),h(i)),
\end{equation}
that is, the pullback of the functor represented by $h(i) \in h^{C
  \cap W}(I_C)$. Then as we saw in Lemma~\ref{mod.le}, elements of
the set in the right-hand side of \eqref{H.i.eq} correspond
bijectively to connected components of the category $QI_C(i',i)$, so
that the fibers of the family of groupoids $\Hh(i)$ are essentially
small, and we have a natural isomorphism
\begin{equation}\label{H.HH}
H(i) \cong \pi_0(\Hh(i)).
\end{equation}
Proposition~\ref{repr.prop} is then a refinement of the usual Yoneda
property of representable functors.

\subsection{Explicit homotopies.}\label{exp.subs}

In practice, to study the representable functor $H(i):I_C^o \to
\Sets$, one chooses a fibrant replacement of the object $i$. This
does not change the functor $H(i)$ but reduces a diagram
\eqref{codom} to a single map. Namely, having assumed that $i$ is
fibrant, choose a factorization
\begin{equation}\label{cyl.eq}
\begin{CD}
i @>{w}>> \wt{i} @>{f}>> i \times i
\end{CD}
\end{equation}
of the diagonal map $i \to i \times i$, with $f \in F$ and $w \in C
\cap W$, and denote by $f_0,f_1:\wt{i} \to i$ the composition of the
map $f$ with the projections onto the two factors of $i \times
i$. Then just as for a decomposition \eqref{c.dec}, any map in
$H(i)(i')$ can be represented by a morphism $g:i' \to i$, and two
maps $g_0,g_1:i' \to i$ represent the same map if and only if there
exists a morphism $h:i' \to \wt{i}$ such that $g_l = f_l \circ h$,
$l=0,1$. We finish the section by giving a similar effective
description of the representable family $\Hh(i)$ of
Definition~\ref{repr.def}.

First, we need to recall one general categorical notion. Assume
given a cellular Reedy category $I$ in the sense of
Definition~\ref{ind.reedy.def}, and assume that $I$ is
Reedy-$\kappa$-bounded for some infinite cardinal $\kappa$ in the
sense of Definition~\ref{reedy.def}. Moreover, assume given a
$\kappa'$-complete category $\C$ for some $\kappa' \ll \kappa$, and
assume that $\C$ is equipped with a class of maps stable under
pullbacks, called {\em coverings}. For any $i \in I$, we have the
$i$-sphere embedding \eqref{sph.I.eq}. For any functor $X \in
I^o\C$, \eqref{hom.X} provides a well-defined object
$\Hom(\Delta_i,X) = X(i)$, and by \eqref{sph.I.M}, $\Hom(\SS_i,X)
\cong M(X,i)$ is also well-defined.

\begin{defn}
A map $f:Y' \to Y$ in $I^o\C$ is a {\em hypercovering} if for
any $i \in I$, the map
\begin{equation}\label{hyp.eq}
Y'(i) \to Y(i) \times_{\Hom(\SS_i,Y)} \Hom(\SS_i,Y')
\end{equation}
induced by $f$ and the sphere embedding \eqref{sph.I.eq} is a
covering.
\end{defn}

If $I=\Delta$ and $\C$ is the category $\Sets$ and coverings are
surjective maps, then $f$ is a hypercovering iff it is a trivial Kan
fibration, and the surjectivity of \eqref{hyp.eq} simply encodes the
lifting propety with respect to \eqref{sph.eq}. On the other hand,
for any Reedy-$\kappa$-bounded $I$, if $\C$ is a $\kappa$-cocomplete
$\kappa'$-complete model category and coverings are maps in $F$ or
in $F \cap W$, then $f$ is a hypercovering iff it is in $F$
resp.\ $F \cap W$ with respect to the Reedy model structure on
$I^o\C$. If $I=\Delta$, then it suffices to require that $\C$ is
finitely complete which it must be anyway simply by definition.
Moreover, in this case, $\C$ is also finitely cocomplete, so that
for any $c \in \C$ and finite simplicial set $X$, we have a
well-defined object $c \times X \in \Delta^o\C$, and for any $Y \in
\Delta^o\C$, we have a natural functorial identification
\begin{equation}\label{maps.du}
\Hom(c \times X,Y) \cong \Hom(c,\Hom(X,Y)) \cong \Hom(X,\Hom^\Delta(c,Y)),
\end{equation}
where $\Hom(X,Y) \in \C$ is given by \eqref{hom.X}, and
$\Hom^\Delta(c,Y) \in \Delta^o\Sets$ is given by
$\Hom^\Delta(c,Y)([n]) = \Hom(c,Y([n]))$, $[n] \in \Delta$. Then for
any map $w:c \to c'$ in the class $C \cap W$ and any $n \geq 1$, the
induced map $c \times \Delta_n \copr_{c \times \SS_{n-1}} c' \times
\SS_{n-1} \to c' \times \Delta_n$ is in $C \cap W$ with respect to
the Reedy structure, so that for any Reedy fibrant $Y \in
\Delta^o\C$, the map
$$
w^*:\Hom^\Delta(c',Y) \to \Hom^\Delta(c,Y)
$$
is a trivial Kan fibration of simplicial sets.

Now return to our original setting --- $I$ a model category, $i \in
I$ a fibrant object --- and choose a Reedy fibrant replacement
$i^\Delta \in \Delta^oI$ of the constant simplicial object $i \in I
\subset \Delta^oI$, with the corresponding map $e:i \to i^\Delta$ in
$C \cap W$. Replace $i$ with $i^\Delta([0])$, so that $i \to
i^\Delta$ is an isomorphism at $[0] \in \Delta$ (note that this does
not change the family $\Hh(i)$). Consider the product $\Cof(I)_C
\times \Delta$, with the projection $\theta:\Cof(I)_C \times \Delta
\to \Cof(I)_C$, and let $\Hh(i^\Delta) \to \Cof(I)_C \times \Delta$
be the family of discrete groupoids corresponding to the functor
$(\Cof(I)_C \times \Delta)^o \to \Sets$, $i' \times [n] \mapsto
\Hom(i',i^\Delta([n]))$. Explicitly, objects in $\Hh(i^\Delta)$ are
triples $\langle i',[n],f \rangle$, $i' \in \Cof(I)$, $[n] \in
\Delta$, $f$ a map from $i'$ to $i^\Delta([n])$. The composition
$\Hh(i^\Delta) \to \Cof(I)_C \times \Delta \to \Cof(I)_C$ is a
fibration, with fibers $\Hh(i^\Delta)_i \cong \Delta
\Hom^\Delta(i',i^\Delta)$, $i' \in \Cof(I)$.

\begin{prop}\label{Hh.prop}
There exists a natural equivalence $\theta_!\Hh(i^\Delta) \cong
\Hh(i)$ of families of groupoids over $\Cof(I)_C$.
\end{prop}

\proof{} Consider the category $QI_C(i)$ of
Subsection~\ref{repr.subs} with its forgetful functor $\phi:QI_C(i)
\to \Cof(I)_C$, let $QI(i)$ be the analogous category formed by
diagrams \eqref{codom} in the whole $I \supset \Cof(I)_C$, with the
forgetful functor $QI(i) \to I$, and consider the fibered product
$$
QI_C'(i) = QI(i) \times_I \Cof(I)_C,
$$
again with the natural forgetful functor $\phi':QI'_C \to
\Cof(I)_C$. Then $\phi'$ is a fibration with fibers $QI(i',i)$, $i'
\in \Cof(I)$, and we have a natural embedding $\eps_i:QI_C(i) \to
QI_C'(i)$, cartesian over $\Cof(I)_C$, whose fiber over some object $i' \in
\Cof(I)$ is the functor $\eps_{i',i}$ of \eqref{eps.i01}. For any
$d \in QI(i',i)$ represented by a diagram \eqref{codom}, the
comma-fiber $d \setminus QI_C(i',i)$ of the functor $\eps_{i',i}$ at
$d$ is the category $\C^W/c$ of Corollary~\ref{mod.corr}, for the
model category $\C = (i' \copr i)\setminus I$ and the object $c \in
\C$ corresponding to $d$. Therefore by Proposition~\ref{filt.prop},
Lemma~\ref{2.1.eq.le}, Lemma~\ref{fib.loc.le} and
Lemma~\ref{i.loc.le}, $QI(i',i)$ is localizable, and
$h(\eps_{i',i})$ is an equivalence. Thus $QI_C'(i)$ is
localizable with respect to $\phi'$, and the functor
$$
\eps_i:\Hh(i) = \phi_!QI_C(i) \to \phi'_!QI'_C(i)
$$
is an equivalence. On the other hand, sending $\langle i',[n],f \rangle \in
\Hh(i^\Delta)$ to the diagram
$$
\begin{CD}
i' @>{f}>> i^\Delta([n]) @<{w}<< i
\end{CD}
$$
provides a functor $\alpha:\Hh(i^\Delta) \to QI'_C(i)$ cartesian over
$\Cof(I)_C$, and it induces a functor
$$
\theta_!\Hh(i^\Delta) \to \phi'_!QI'_C \cong \Hh(i).
$$
To prove that this is an equivalence, it suffices to check that for
any object $i' \in \Cof(I)$, the functor
$$
\alpha_{i'}:\Delta \Hom^\Delta(i',i^\Delta) \cong \Hh(i^\Delta)_{i'}
\to QI(i',i) \cong QI'_c(i)_{i'}
$$
induces an equivalence of total localizations. By
Lemma~\ref{fib.loc.le}, it then suffices to check that for any $d
\in QI(i',i)$ represented by a diagram \eqref{codom}, the
comma-fiber $\Delta \Hom^\Delta(i',i^\Delta) \setminus_{\alpha_{i'}}
d$ is $1$-connected. But explicitly, this comma-fiber is the category of
pairs $\langle [n],g \rangle$ of an object $[n] \in \Delta$ and a
map $g:\wt{i} \to i^\Delta([n])$ such that $g \circ w = e:i \to
i^\Delta([n])$, and this is nothing but the category of simplices
$\Delta X$ of the simplicial set given by the fibered product
$$
\begin{CD}
X @>>> \Hom^\Delta(\wt{i},i^\Delta)\\
@VVV @VV{w^*}V\\
\ppt @>{e}>> \Hom^\Delta(i,i^\Delta).
\end{CD}
$$
But $w^*$ is a trivial Kan fibration, so that $X$ is trivially
fibrant, and we are done by Lemma~\ref{kan.le}.
\endproof

As a corollary of Proposition~\ref{Hh.prop}, for any $i' \in I_C$
and any choice of a Reedy-fibrant replacement $i^\Delta$, we have a
natural identification
\begin{equation}\label{H.del}
\Hh(i)_{i'} \cong h^{\Tot}(\Delta\Hom^\Delta(i',i^\Delta)).
\end{equation}
This helps in computations by virtue of the following.

\begin{lemma}\label{reedy.le}
For any morphism $f:i_0 \to i_1$ in $\Cof(I)_C \subset I$, fibrant
$i$ and Reedy-fibrant replacement $i^\Delta \in \Delta^oI$ of the
constant simplicial object $i$, the induced map
$f^*:\Hom^\Delta(i_1,i^\Delta) \to \Hom^\Delta(i_0,i^\Delta)$ is a
Kan fibration of simplicial sets. In particular,
$\Hom^\Delta(i',i^\Delta)$ is fibrant for any $i' \in \Cof(I)$.
\end{lemma}

\proof{} Since $i^\Delta$ is Reedy fibrant, the natural map
$i^\Delta([n]) \to \Hom(\SS_{n-1},i^\Delta)$ is in $F$ for any $n
\geq 0$, and then by induction, the map
\begin{equation}\label{e.e}
e^*:\Hom(X',i^\Delta) \to \Hom(X,i^\Delta)
\end{equation}
is in $F$ for any injective map $e:X \to X'$ of finite simplicial
sets. Let $E$ be the class of such maps $e$ for which \eqref{e.e} is
also in $W$. Then $E$ inherits the two-out-of-three property of the
class $W$, and since the class $F \cap W$ is stable under pullbacks,
$E$ is stable under pushouts. Moreover, since $i^\Delta$ is a
fibrant replacement of a constant simplicial object, $E$ contains
all the embeddings $t:\Delta_0 \subset \Delta_n$ corresponding to
right-closed embeddings $t:[0] \to [n]$. But as in
Lemma~\ref{bva.le}, for any $n \geq l \geq 0$, $s$ factors through
the source $\V^l_n$ of the horn embedding \eqref{horn.eq} via some
injective map $t':\Delta_0 \to \V^l_n$, and the map $t'$ is a
composition of pushouts of horn embeddings $v^{l'}_{n'}$ with $n' <
n$. Then by induction on dimension, $E$ contains all the horn
embeddings \eqref{horn.eq}, and $f^*$ is a Kan fibration by
\eqref{maps.du}.
\endproof

\begin{remark}
The isomorphism \eqref{H.HH} can be immediately observed in terms of
\eqref{H.del} --- in effect, a factorization \eqref{cyl.eq} is a
first step of a Reedy fibrant replacement $i^\Delta$, with $i =
i^\Delta([0])$ and $\wt{i}=i^\Delta([2])$, so we recover the usual
description of the functor $H(i)$. By Lemma~\ref{kan.le}, to compute
$\Hh(i)$, it actually suffices to consider $i^\Delta([0])$,
$i^\Delta([1])$ and $i^\Delta([2])$.
\end{remark}

\chapter{Semiexactness.}\label{semi.ch}

This chapter contains the first of our main results, a version of
Brown Representability for semiexact additive families of groupoids
over simplicial sets (or equivalently, over topological
spaces). This does seem to be new.

We start by preparing the ground. In Section~\ref{CW.sec}, we
introduce the notion of a ``CW-category'' --- roughly speaking, this
is a model category that has classes $C$ and $W$, but not $F$. It
also has pushouts of maps in the class $C$ (but is not required to
have all finite colimits). The notion is somewhat technical, even
more so than that of a model category, but it turns out that one can
prove something even in such a generality. The main examples are,
firstly, model categories, and secondly, various categories of
partially ordered sets --- in the latter case, $C$ and $W$ are
left-closed embeddings and reflexive maps of
Chapter~\ref{ord.ch}. Additional examples are finite-dimensional
simplicial sets, and chain complexes in an additive category.

One thing that survives in the setting of CW-categories is the Reedy
construction (this is Subsection~\ref{reedy.CW.subs}). Another is a
version of the small object argument, in
Subsection~\ref{ab.subs}. We actually develop the theory of
CW-categories further than strictly necessary, up to a non-trivial
result (Proposition~\ref{reg.prop}, a generalization of
Proposition~\ref{mod.prop}). This can be probably skipped, and the
same holds for an earlier Proposition~\ref{w.CW.prop}.

In Section~\ref{semi.sec}, we turn to the study of families of
groupoids on a CW-category. We introduce our two main conditions on
such families, additivity and semiexactness (the latter is a version
of the Mayer-Vietoris property for homology). The main surprise in
Section~\ref{semi.sec} is Subsection~\ref{fram.subs}, with
Lemma~\ref{ano.C.le} and Corollary~\ref{ano.C.corr}. It turns out
that under an additional assumption that happens to be satisfied in
examples of interest, if a semiexact additive family of groupoids
over a CW-category is constant along maps in the class $W$, it is
automatically constant along a much larger class of ``$W$-anodyne
maps'' (for $\Pos$, this corresponds to passing from reflexive maps
to anodyne maps of Subsection~\ref{ano.subs}). It is this phenomenon
that later allows us to reduce a homotopy invariance condition to a
much milder conditions such as e.g.\ excision of
Subsection~\ref{exci.subs}.

Then in the last part of Section~\ref{semi.sec}, we turn to
semiexact families over model categories, and we prove some further
results specific to this case. In particular,
Proposition~\ref{lift.prop} provides a simple criterion that allows
to recognize representable families of Subsection~\ref{repr.subs}.

Actual proof of representability is in
Section~\ref{repr.sec}. First, in Subsection~\ref{top.subs}, we give
a proof that follows very closely the standard proof found e.g.\ in
\cite{switz}. This uses topological spaces; the formal statement is
Theorem~\ref{coho.thm}. Then we move to simplicial sets, and provide
another proof, in Theorem~\ref{coho.bis.thm}, that actually gives
slightly more --- namely, we can control the size of the
representing object. Technically, both proofs depend on the
representability criterion of Proposition~\ref{lift.prop}; for
Theorem~\ref{coho.bis.thm}, the corresponding technical statement is
given separately in Subsection~\ref{repr.lft.subs}. As a free bonus,
we also obtain a representability theorem for families over the
category $I^o\Delta^o\Sets$, for any cellular Reedy category $I$.

In the last part of Section~\ref{repr.sec}, we restrict our
attention to finite-dimen\-sional simplicial sets. This is no longer a
model category, so this needs the technology of CW-categories of
Section~\ref{CW.sec}. However, we can still prove a version of
representability. As a corollary, we obtain
Proposition~\ref{X.ext.prop} saying that a family over
finite-dimensional simplicial sets extends to a family over all
simplicial sets, uniquely up to an equivalence unique up to a
non-unique isomorphism (we call this ``semicanonical
extension''). The same results hold for families over
$I^o\Delta^o\Sets$, for a cellular Reedy category $I$.

\section{CW-categories.}\label{CW.sec}

\subsection{Generalities.}\label{CW.subs}

By virtue of Proposition~\ref{mod.prop}, when studying families of
groupoids over a relative category $\langle I,W \rangle$ equipped
with a model structure, it suffices to consider cofibrant objects
and maps in the classes $C$ and $W$. We will also need to consider
more general situations when the ambient model category does not
exist, and we only have $C$, $W$, and not $F$. It is convenient to
axiomatize the situation. We start with following skeleton
definition, and we will add embellishments as we go along.

\begin{defn}\label{CW.def}
A {\em C-category} is a category $I$ equipped with a closed class of
maps $C$ such that $I$ has an initial object $0$, and for any maps
$f_0:i \to i_0$, $f_1:i \to i_1$ in $C$, there exists a cocartesian
square
\begin{equation}\label{st.CW.sq}
  \begin{CD}
    i @>{f_0}>> i_0\\
    @V{f_1}VV @VVV\\
    i_1 @>>> i_{01}
  \end{CD}
\end{equation}
in $I$, with both maps $i_0,i_1 \to i_{01}$ lying in $C$. A
{\em weak CW-category} is a C-category $I$ equipped with an
additional closed class of maps $W$ such that
\begin{enumerate}
\item for any square \eqref{st.CW.sq}, if $i \to i_0$ is in $W$,
  then so is $i_1 \to i_{01}$, and
\item for any $i_0 \in I$ such that $0 \to i_0$ is in $C$, any map
  $f:i_0 \to i_1$ in $I$ admits a decomposition
\begin{equation}\label{f.dec}
\begin{CD}
i_0 @>{c}>> i_2 @>{w}>> i_1
\end{CD}
\end{equation}
such that $c \in C$, and for any map $p:i_1 \to i_0$ in $C$ with $f
\circ p = \id$, we have $c \circ p \in W$.
\end{enumerate}
A {\em CW-category} $I$ is a weak CW-category such that
\eqref{f.dec} in \thetag{ii} can be chosen in such a way that $c
\circ p \in W$ for any $p:i_0' \to i_0$ in $C$ with $f \circ p \in
W$.
\end{defn}

\begin{exa}\label{CW.comma.exa}
For any C-category $\langle I,C \rangle$ and any discrete cofibration
$\pi:I' \to I$, $\langle I',\pi^*(C)\rangle$ is a C-category, and
the same holds if $\pi$ is a discrete fibration cartesian along each
square \eqref{st.CW.sq}. In both cases, if $I$ is a CW-category with
respect to some class $W$, then $\langle
I',\pi^*(C),\pi^*(W)\rangle$ is a CW-category, with the
decomposition \eqref{f.dec} induced from those in $I$. In
particular, for any object $i \in I$ in a CW-category $I$, the
comma-fibers $I/i$, $i \setminus I$ carry a natural CW-structure,
and the same holds for any two objects $i',i'' \in I$ and the
category $i' \setminus I /i$ of objects $i \in I$ equipped with maps
$i' \to i \to i''$. Note that $i' \setminus I / i$ decomposes into a
disjoint union according to the composition map $f:i' \to i''$, and
each of the components $(i' \setminus I / i)_f$ of this
decomposition is also a CW-category.
\end{exa}

\begin{exa}\label{CW.sets.exa}
The category $\Sets$ is trivially a $C$-category with respect to the
class $C$ of injective maps. This does not extend to a meaningful
CW-structure (not even to a weak one).
\end{exa}

\begin{exa}\label{mod.CW.exa}
A model category $\langle I, C,W,F \rangle$ is automatically a
CW-category in the sense of Definition~\ref{CW.def}, with the same
classes $C$, $W$, and \eqref{f.dec} provided by the model structure
factorization.
\end{exa}

\begin{exa}\label{CW.add.exa}
An additive category $\A$ is a $C$-category with respect to the
class $C$ of split injections. The category of complexes
$C_\idot(\A)$ is a C-category with respect to the class $C$ of
termwise-split injections, and it is a CW-category with respect to
the class $W$ of chain-homotopy equivalences; a decomposition
\eqref{f.dec} for a map $f:A_\idot \to B_\idot$ is given by
\begin{equation}\label{add.dec}
\begin{CD}
  A_\idot @>{a \oplus f}>> \Cone(\id_A) \oplus B_\idot @>{b}>>
  B_\idot,
\end{CD}
\end{equation}
where $a:A_\idot \to \Cone(\id_A)$ is the natural embedding, and $b$
is the projection.
\end{exa}

\begin{exa}\label{pos.CW.exa}
Let $I = \Pos$ be the category of partially ordered sets, with $C$
consisting of left-closed full embeddings, and $W$ consisting of
reflexive full embeddings of Definition~\ref{po.refl.def}. Then
$\langle \Pos,C,W \rangle$ is a CW-category in the sense of
Definition~\ref{CW.def}: \thetag{i} follows from
Lemma~\ref{po.copr.le}, and by Lemma~\ref{cyl.le} and
Lemma~\ref{filt.le}, \eqref{cyl.deco} gives a decomposition
\eqref{f.dec}. We call this the {\em standard} CW-structure on
$\Pos$. Dually, the {\em co-standard} CW-structure has the same
class $W$, and $C$ consists of right-closed embeddings. Any very ample
full subcategory $\I \subset \Pos$ in the sense of
Definition~\ref{amp.def} such as $\Posf$, $\Posf_\kappa$, $\Posff$
inherits both the standard and the co-standard CW-structures. Its
extension $\I^+$ also inherits both CW-structures, but for the
standard structure, we have to use \eqref{cyl.N.deco} instead of
\eqref{cyl.deco}.
\end{exa}

\begin{exa}\label{pos.I.CW.exa}
More generally, assume given a category $I$, and consider the
category $\Pos \bb I$ of $I$-augmented partially ordered sets, as in
Subsection~\ref{pos.aug.subs}. Then $\Pos \bb I$ with the classes
$C$ of left-closed embeddings and $W$ of reflexive full embeddings
of Definition~\ref{aug.def} is a CW-category, and $\I \bb I$
resp.\ $\I^+ \bb I$ for a very ample $\I \subset \Pos$ resp.\ its
extension $\I^+ \bb I$ inherits the CW-structure that we again call
{\em standard}. The {\em co-standard} CW-structure has the same $W$,
and $C$ consisting of right-closed embeddings. The decomposition
\eqref{f.dec} is given by \eqref{cyl.o.deco}, and since the latter
is only well-defined for strict maps, the co-standard CW-structure
exists only on $\Pos \bbi I \subset \Pos \bb I$.
\end{exa}

For any C-category $I$, we will call pushout squares
\eqref{st.CW.sq} {\em standard pushout squares} in $I$ (in the
situation of Example~\ref{pos.CW.exa}, these are standard pushout
squares \eqref{pos.st.sq} of partially ordered sets). In the
situation of Example~\ref{pos.I.CW.exa}, {\em co-standard pushout
  squares} are standard pushout squares in $\Pos \bbi I$ with respect
to the co-standard CW-structure. We will say that a functor
$\gamma:I \to I'$ between C-categories $\langle I,C \rangle$,
$\langle I',C' \rangle$ is a {\em C-functor} if $C \subset
\gamma^*C'$, and $\gamma$ sends standard pushout squares to standard
pushout squares. If $I$, $I'$ are weak CW-categories with classes
$W$, $W'$, then a $C$-functor $\gamma$ is a {\em weak CW-functor} if
$W \subset \gamma^*W'$. A weak CW-functor is a CW-functor if it
sends the initial object $0 \in I$ to the initial object $0 \in I'$.

\begin{exa}\label{pos.II.CW.exa}
For any functor $\gamma:I' \to I$ between some categories $I$, $I'$,
the induced functor $\Pos \bb I' \to \Pos \bb I$ is a CW-functor
with respect to the standard CW-structures of
Example~\ref{pos.I.CW.exa}.
\end{exa}

\begin{exa}
The involution $\iota:\Pos \to \Pos$, $J \mapsto J^o$ is a
CW-functor with respect to standard and co-standard CW-structures.
\end{exa}

\begin{exa}\label{rel.CW.exa}
Consider the category $\Rel$ of biordered sets of
Definition~\ref{rel.def}, with the classes $C$ resp.\ $W$ of
left-closed resp.\ reflexive full embeddings of
Definition~\ref{rel.refl.def}. Then $\Rel$ is a CW-category, and the
forgetful functor $U:\Rel \to \Pos$ is a CW-functor. As in
Example~\ref{pos.CW.exa} and Example~\ref{pos.I.CW.exa}, we call
this CW-structure on $\Rel$ {\em standard}; we note that the
standard CW-structure is inherited by the unfolding $\I^\dm =
U^{-1}(\I) \subset \Rel$ of any very ample full subcategory $\I
\subset \Pos$. The unfolding $\I^{+\dm}$ of the extension $\I^+$
then also inherits the standard CW-structure, with the same
modifications as in Example~\ref{pos.CW.exa}.
\end{exa}

As in the model category case, we will say that an object $i \in I$
in a C-category $I$ is {\em cofibrant} if $0 \to i$ is in $C$, and
we will denote by $\Cof(I) \subset I$ the full subcategory spanned
by cofibrant objects. It inherits a C-structure, and a CW-structure
if $I$ had one. In practice, this is the only part of $I$ that we
will need, but it is convenient to allow a more general definition.

\begin{exa}\label{CW.ar.exa}
Consider the arrow category $\Ar(I)$ of a C-category $I$. Say that
an arrow $i:[1] \to I$ is {\em cofibrant} if it factors through
$\Cof(I)_C$, and say that a map $f:i \to i'$ in $\Ar(I)$ between
cofibrant $i$, $i'$ is in the class $C$ if
\begin{enumerate}
\item $f(0):i(0) \to i'(0)$ is in $C$, so that $i(1) \copr_{i(0)}
  i'(0)$ exists by Definition~\ref{CW.def}~\thetag{i}, and
\item the induced map $f':i(1) \copr_{i(0)} i'(0) \to i'(1)$ is also in
  $C$.
\end{enumerate}
Then $\Ar(I)$ with this class $C$ is a $C$-category. Moreover, if
$I$ is a weak CW-category or a CW-category, say that $f$ is in $W$
if both $f(0)$ and $f'$ are in $C \cap W$. Then $\Ar(I)$ with these
classes $C$, $W$ is a weak CW-category resp.\ CW-category, the
subcategory $\Cof(\Ar(I)) \subset \Ar(I)$ of cofibrant arrows is a
CW-subcategory, and the evaluation functors
$\ev_0,\ev_1:\Cof(\Ar(I)) \to I$ are CW-functors. Alternatively, one
can define a smaller class of maps $C_0 \subset C$ by requiring that
$f(0) \in C$ and $f'$ is an isomorphism; then $\Ar(I)$ with the
classes $C_0$, $W$ is also a CW-category, and $\id:\langle I,C_0,W
\rangle \to \langle I,C,W \rangle$ is a CW-functor.
\end{exa}

\begin{exa}\label{CW.ar.non.exa}
Assume given a CW-functor $\gamma:I \to I'$ between CW-categories
$I$, $I'$, consider the comma-category $I /_\gamma I'$ with its
projections $\sigma:I /_\gamma I' \to I$, $\tau:I /_\gamma I' \to
I$, and say that a map $f$ in $I /_\gamma I'$ is in $C$ resp.\ $W$
if both $\sigma(f)$ and $\tau(f)$ are in $C$ resp.\ $W$. Then this
assignment does {\em not} turn $I /_\gamma I'$ into a CW-category,
nor even a weak one: while Definition~\ref{CW.def}~\thetag{i} is
satisfied, \thetag{ii} can fail. Already if $I'=I$ and
$\gamma = \id$, so that $I /_\gamma I' \cong \Ar(I)$, the classes
$C$ and $W$ defined above are different from the classes $C$ and $W$
of Example~\ref{CW.ar.exa}.
\end{exa}

\begin{lemma}\label{CW.W.le}
For any CW-category $\langle I,C,W \rangle$, the tautological
embedding $\langle \Cof(I),C \cap W \rangle \to \langle \Cof(I),W
\rangle$ is a $2$-equivalence.
\end{lemma}

\proof{} By Lemma~\ref{sat.le}, it suffices to show that $W$ lies
inside the saturation $W'=s(\Cof(I),C \cap W)$. To do this, one can
adapt a version of the cylinder construction quite standard in the
theory of model categories (where it is also known as ``Ken Brown's
Lemma'', \cite[Lemma 9.9]{DS}). Assume given a map $f:i \to i'$ in
$\Cof(I)$. Since both $i$ and $i'$ are cofibrant, the coproduct $i
\copr i'$ exists in $I$ by Definition~\ref{CW.def}, and we can take
a decomposition
\begin{equation}\label{f.i.dec}
\begin{CD}
i \copr i' @>{c \copr c'}>> \wt{i} @>{w}>> i'
\end{CD}
\end{equation}
of the map $f \copr \id:i \copr i' \to i'$ provided by
Definition~\ref{CW.def}~\thetag{ii}. Since $w \circ c' = \id$, $c'
\in C \cap W$, and then $w \in W'$. But if $f$ is in $W$, then $c$
is also in $C \cap W$, and then $f = w \circ c$ is in $W'$.
\endproof

\begin{defn}\label{abs.C.def}
A morphism $c:i \to i'$ in a $C$-category $I$ is {\em absolutely
  cofibrant} if any map $f:i \to i_1$ in $I$ fits into a cocartesian
square
\begin{equation}\label{abs.sq}
\begin{CD}
  i @>{c}>> i'\\
  @V{f}VV @VVV\\
  i_0 @>{c_0}>> i'_0
\end{CD}
\end{equation}
in $I$, and $c_0 \in C$.
\end{defn}

\begin{exa}\label{abs.exa}
Take $I = \Pos$, with the co-standard $C$-structure. Then any split
map $s:J \to J'$ in $I$ is absolutely cofibrant by
Example~\ref{pos.triv.exa}, and for any map $g:J_0 \to J_1$, the
right-closed embedding $t:J_1 \to \Cyl(f)$ is also absolutely
cofibrant, with \eqref{abs.sq} given by \eqref{cyl.o.sq} for the map
$f \circ g$.
\end{exa}

\subsection{The Reedy construction.}\label{reedy.CW.subs}

Example~\ref{CW.ar.exa} is the simplest case of the Reedy
construction of Theorem~\ref{reedy.thm}, and somewhat surprisingly,
the whole construction survives to some extent in the generality of
Definition~\ref{CW.def}. Namely, recall that any left-bounded
partially ordered set $J$ is a Reedy category by
Example~\ref{reedy.exa}, with latching categories $L(j) = J/'j$ also
being partially ordered sets.

\begin{defn}\label{reedy.CW.def}
A functor $X:J \to I$ from a Reedy category $J$ to a C\--ca\-te\-gory
$I$ is {\em cofibrant} if for any given object $j \in J$, the latching object
$L(X,j)=\colim_{L(j)}X$ exists in the category $I$, and the natural
map $L(X,j) \to X(j)$ is in the class $C$. A map $X \to X'$ between
two functors is in the class $C$ if both $X$ and $X'$ are cofibrant,
and for any $j \in J$, the map $L(X,j) \to L(X',j)$ is in $C$, and
the natural map
\begin{equation}\label{L.X}
L(X',j) \copr_{L(X,j)} X(j) \to X'(j)
\end{equation}
is in $C$. If $I$ is a weak CW-category, then
$f:X \to X'$ is in $W$ if both $X$ and $X'$ are cofibrant, and all
the maps \eqref{L.X} are in $C \cap W$.
\end{defn}

\begin{remark}
If $I$ is a model category, then $X:J \to I$ is cofibrant in the
sense of Definition~\ref{reedy.CW.def} iff it is cofibrant with
respect to the Reedy model structure of Theorem~\ref{reedy.thm}.
\end{remark}

\begin{remark}\label{reedy.CW.rem}
For any discrete fibration $\pi:J' \to J$, $J'$ inherits the Reedy
structure of Example~\ref{reedy.L.exa}, and since $\pi:L(j) \to
L(\pi(j))$ is an equivalence for any $j \in J'$, $\pi^*X:J' \to I$
is tautologically cofibrant for any cofibrant $X:J \to I$. In
particular, this applies when $J$ is directed and $J'=L(j)$ is the
latching category of an object $j \in J$, or when $J$ is a
left-bounded partially ordered set, and $J' \subset J$ is a
left-closed subset.
\end{remark}

\begin{lemma}\label{CW.J.le}
Assume given a functor $X:J \to I$ from a partially ordered set $J$
to a C-category $I$ that is cofibrant in the sense of
Definition~\ref{reedy.CW.def}.
\begin{enumerate}
\item If $J$ is finite, $\colim_JX$ exists and is cofibrant, and for
  any left-closed subset $J' \subset J$, the map $\colim_{J'}X \to
  \colim_JX$ is in the class $C$.
\item For any order-preserving left-finite map $\phi:J \to J_1$, the
  left Kan extension $\phi_!X:J_1 \to I$ exists and is cofibrant,
  and for any morphism $f:X \to X'$ in $C$, the morphism
  $\phi_!(f):\phi_!X \to \phi_! X'$ is in $C$.
\item The functor $X$ factors through through $\Cof(I)_C \subset I$,
  and if a map $f:X \to X'$ is in $C$, then $f(j):X(j) \to X'(j)$ is
  in $C$ for any $j \in J$.
\item For any left-closed subset $J_0 \subset J$ with the embedding
  map $\lambda:J_0 \to J$, the adjunction map $\lambda_!\lambda^*X
  \to X$ is in the class $C$.
\end{enumerate}
Moreover, if $I$ is a weak CW-category, then $\phi_!$ in \thetag{ii}
sends maps in $C \cap W$ to maps in $C \cap W$, and if $f$ in
\thetag{iii} is in $W$, then so is $f(j)$ for any $j \in J$.
\end{lemma}

\proof{} For \thetag{i}, note that for any left-closed $J' \subset
J$, the restriction $X|_{J'}$ is cofibrant by
Remark~\ref{reedy.CW.rem}. Then by induction on the cardinality
$|J|$, it suffices to consider $J' = J \ssetminus \{j\}$, $j \in J$ a
maximal element, and we may assume that $\colim_{J'}X$ exists and is
cofibrant. Then we have the standard pushout square of
Example~\ref{elem.exa}, and $\colim_JX$, if it exists, is given by
the cocartesian square
\begin{equation}\label{J.j.sq}
\begin{CD}
\colim_{L(j)}X @>{a}>> \colim_{J/j}X\\
@V{b}VV @VV{b'}V\\
\colim_{J'}X @>>> \colim_JX,
\end{CD}
\end{equation}
where $a$ is in $C$ by Definition~\ref{reedy.CW.def}, and $b$ is in
$C$ by induction. By Definition~\ref{CW.def}~\thetag{i}, the square
indeed exists, and $b'$ is in $C$.

For \thetag{ii}, the left Kan extension $\phi_!X$ is given by
\eqref{kan.eq} that reads as
\begin{equation}\label{pos.kan.eq}
\phi_!X(j) = \colim_{J/j}X, \qquad j \in J',
\end{equation}
and if $\phi$ is left-finite, all the comma-sets $J/j$ are finite,
so we may assume that $J$ is finite and $J_1=\ppt$, and then $\phi_!X
= \colim_JX$ exists by \thetag{i}. The condition of
Definition~\ref{reedy.CW.def} for $\phi_!X$ amounts to checking that
for any $j \in J'$, the map $\colim_{\phi^{-1}(L(j))}X \to
\colim_{J/j}X$ is in $C$, so that we may again assume that $J$ is
finite, and $\phi:J \to [1]$ is the characteristic function of a
left-closed subset. Then $\phi_!X$ is cofibrant by
\thetag{i}. Moreover, to check that $\phi_!(f)$ is in $C$ or in $C
\cap W$, we may again assume that $J$ is finite and $J_1=[1]$, and
then by the same induction as in \thetag{i}, the statement
immediately follows from Definition~\ref{CW.def}~\thetag{i} and the
cocartesian square \eqref{J.j.sq}.

For \thetag{iii}, the argument is essentially the same as
Lemma~\ref{J.proj.le}. Namely, if $J=[0]$ or $J=[1]$, the claim
immediately follows from Definition~\ref{CW.def}~\thetag{i}, so it
suffices to check that for any map $i:[1] \to J$, the
restriction functor $i^*$ sends cofibrant functors to cofibrant
functors and maps in $C$ resp.\ $C \cap W$ to maps in $C$ resp.\ $C
\cap W$. But $i$ corresponds to a pair $j',j \in J$, $j' \leq
j$, we can replace $J$ by $J/j$, thus assume that $j$ is the largest
element, and then the characteristic function $\chi = \chi_{J/j'}:J
\to [1]$ of $J/j' \subset j$ is adjoint to $i$, so that
$i^* \cong \chi_!$.

Finally, for \thetag{iv}, induction on $|J \ssetminus J_0|$ reduces
us to the case when $J \ssetminus J_0$ is discrete, and in this case,
the claim holds tautologically.
\endproof

\begin{defn}\label{thin.def}
A Reedy category $J$ is {\em thin} if it is directed (that is, $J =
J_L$) and $J/j$ is a finite partially ordered set for any $j \in J$.
\end{defn}

\begin{exa}
A left-bounded partially ordered set $J$ with the Reedy structure of
Example~\ref{reedy.exa} is thin if and only if it is left-finite.
\end{exa}

\begin{prop}\label{CW.prop}
Assume given a thin Reedy category $J$ and a C-cate\-go\-ry resp.\ a
CW-category $I$. Then the functor category $\Fun(J,I)$ with the
classes $C$, resp.\ $C$ and $W$ of Definition~\ref{reedy.CW.def} is
a C-category resp.\ a CW-category in the sense of
Definition~\ref{CW.def}. Moreover, for any morphism $\phi:J_0 \to
J_1$ between left-finite partially ordered sets, the left Kan
extension $\phi_!:\Cof(\Fun(J_0,I)) \to \Cof(\Fun(J_1,I))$ is a
C-functor resp.\ a CW-functor.
\end{prop}

\proof{} By definition, for any $j \in J$, the comma-category $J/j$
is a finite partially ordered set, with the discrete fibration
$\sigma(j):J/j \to J$ of \eqref{p.i}. Thus by
Remark~\ref{reedy.CW.rem}, a functor $X:J \to I$ is cofibrant iff so
is $\sigma(j)^*X$ for any $j \in J$, and a map $f$ between cofibrant
objects is in $C$ iff so is $\sigma(j)^*(f)$. In particular, maps in
$C$ are pointwise in $C$ by Lemma~\ref{CW.J.le}~\thetag{iii}, and
the existence of the coproducts in
Definition~\ref{CW.def}~\thetag{i} immediately follows from the
corresponding statement for $I$. Moreover, for any $j \in J$, let
$\chi(j):J/j \to [1]$ be the characteristic function of the latching
subset $L(j) \subset J/j$. Then by definition, a map $f:X \to X'$ is
in $C$ resp.\ $W$ iff so is $\chi(j)_!\sigma(j)^*(f)$ for any $j \in
J$, and $\sigma(j)^*$ obviously preserves pushouts and maps in $C$
resp.\ $W$, while $\chi(j)_!$ does so by
Lemma~\ref{CW.J.le}~\thetag{ii}. Thus checking that the pushouts of
maps in $C$ resp.\ $C \cap W$ are still in $C$ resp.\ $C \cap W$
reduces to the case $J=[1]$, where it again follows from
Definition~\ref{CW.def}~\thetag{i} for $I$. The functor $\phi_!$
preserves pushouts by adjunction, and $C$ and $C \cap W$ by
Lemma~\ref{CW.J.le}~\thetag{ii}. Finally, to construct a
decomposition
$$
\begin{CD}
X @>{c}>> X' @>{a}>> Y
\end{CD}
$$
of a map $g:X \to Y$ in $\Fun(J,I)$, $X \in \Cof(\Fun(J,I))$ with
the properties required in Definition~\ref{CW.def}~\thetag{ii}, note
that by induction on skeleta, we may assume this done on $\sk_{n-1}
J \subset J$, and we need to extend the decomposition to $\sk_n
J$. This can be done separately at each object $j$ in the
complement $\sk_n J \ssetminus \sk_{n-1} J$. Fix $j$, with the
projection $\bsigma(j):L(j) \to J/j \to J$, note that $L(X',j) =
\colim_{L(j)}\bsigma(j)^*X'$ exists by
Lemma~\ref{CW.J.le}~\thetag{i} and induction, and moreover, since
the functor $\bsigma(j)^*Y:L(j) \to I$ has a cone
$\sigma(j)^*Y:J/j = L(j)^> \to I$, the map $\bsigma(j)^*(a)$ gives
rise to a natural map $a(j):L(X',j) \to Y(j)$. Then one can define
$X'(j)$ and $a:X'(j) \to Y(j)$ by taking the decomposition
$$
\begin{CD}
L(X',j) \copr_{L(X,j)} X(j) @>{c}>> X'(j) @>{a}>> Y(j)
\end{CD}
$$
of the natural map $L(X',j) \copr_{L(X,j)} X(j) \to Y(j)$, where the
coproduct again exists by Lemma~\ref{CW.J.le}~\thetag{i}.
\endproof

\subsection{Abundant $C$-categories.}\label{ab.subs}

For any $C$-category $I$, finite set $S$ and functor $J:S^< \to I$
that factors through $I_C$, colimit $\colim_{S^<}J$ exists, and for
any subset $S_0 \subset S$, the map $\colim_{S_0^<}J \to
\colim_{S^<}J$ is in $C$ (for $S=\{0,1\}$ and $S_0 = \{0\}$, this is
Definition~\ref{CW.def}, and the general case immediately follows by
induction on $|S|$). It is useful to also consider $C$-categories
where this holds for an arbitrary set $S$.

\begin{defn}\label{ab.def}
A $C$-category $I$ is {\em abundant} if for any set $S$ and functor
$X:S^< \to I$ that factors through $I_C$, colimit $\colim_{S^<}X$
exists, and for any subset $S_0 \subset S$, the map $\colim_{S_0^<}X
\to \colim_{S^<}X$ is in $C$.
\end{defn}

\begin{exa}\label{ab.filt.exa}
If a $C$-category $I$ admits filtered colimits, and filtered
colimits preserve $C$, then $I$ is obviously abundant (and
conversely, if $I$ is abundant, then $\colim_{S^<}X$ is the filtered
colimit of $\colim_{S^<_0}X$ over all finite subsets $S_0 \subset
S$). In particular, this holds for the category $\Pos$ with either
the standard or co-standard $C$-structure of
Example~\ref{pos.CW.exa}. For the standard structure, colimits
$\colim_{S^<}J$ of Definition~\ref{ab.def} are standard colimits
\eqref{lc.coli.eq}. The subcategory $\Posf^+ \subset \Pos$ is also
abundant with respect to the standard structure, and dually,
$\Posf^- \subset \Pos$ is abundant with respect to the co-standard
one.
\end{exa}

\begin{lemma}\label{ab.CW.le}
For any abundant $C$-category $I$, Lemma~\ref{CW.J.le} holds for any
$J \in \Posf$ and any map $\phi:J \to J_0$ in $\Posf$, while
Proposition~\ref{CW.prop} holds for any directed cellular Reedy
category $J$.
\end{lemma}

\proof{} For Lemma~\ref{CW.J.le}~\thetag{i}, the claim immediately
follows from Definition~\ref{ab.def} if $J=S^<$ for some set $S$. In
the general case, assume by induction that everything is already
proved for any $J$ with $\dim J < l$ for some $l \geq 0$ (where
$\dim J = -1$ when $J$ is empty). Then for any $J$ with $\dim J =
l$, let $S = J \ssetminus \sk_{l-1}J$, and let $p:J \to S^<$ be the
standard extension \eqref{p.J.sq} of the right-closed embedding $S
\to J$. Then $p_!X$ exists by \eqref{pos.kan.eq} and induction, and
it is cofibrant by the same argument as in Lemma~\ref{CW.J.le}, so
the induction step reduces to the case $J=S^<$. The rest of the
proof is exactly the same as in Lemma~\ref{CW.J.le}. For
Proposition~\ref{CW.prop}, again, the same proof works.
\endproof

\begin{lemma}\label{pos.cof.le}
For any $J \in \Posf$, a functor $X:J \to \Posf^+$ is cofibrant with
respect to the standard $C$-structure of Example~\ref{pos.CW.exa}
iff $X \cong \sE_{J^o}({J'}^o)$ for some $J' \in \Posf^+ \bb J$,
where $\sE_{J^o}$ is the functor \eqref{E.eq}. Moreover, in this
case, we have $J' \cong \colim_JX$.
\end{lemma}

\proof{} Consider the cofibration $\pi:J_\idot \to J$ corresponding
to $X$ by the Grothendieck construction, and let $J' \subset
J_\idot$ be the full subset of all $j \in J_\idot$ that are not in
the image of the map $X(j') \to X(\pi(j))$ for any $j' \leq \pi(j)$
in $J$. Since $\pi$ is a cofibration, the embedding $J' \subset
J_\idot$ uniquely extends to a map $p:J' / J \to J_\idot$
cocartesian over $J$, and it suffices to check that it is an
isomorphism --- the identification $J' \cong \colim_JX$ is then
simply the standard colimit \eqref{lc.coli.eq}. Assume by induction
that everything is already proved when $\dim J < l$ for some $l \geq
0$, take some $J$ of dimension $\dim J = l$, and let $J_0 =
\sk_{l-1}J$, with the corresponding subset $J'_0 = \pi^{-1}(J_0)
\cap J' \subset J'$. Then by induction, $p$ is an isomorphism over
$J_0$, so it suffices to check that it is an isomorphism over any $j
\in J \ssetminus J_0$. But then $L(X,j) = \colim_{L(j)}X \cong J'_0
/ j$, so that since $X$ is cofibrant, $p:J' / j \to X(j)$ is a
left-closed embedding, and then by definition, $(J' \ssetminus J'_0)
/ j = X(j) \ssetminus L(X,j)$.
\endproof

\begin{lemma}\label{exc.CW.le}
Say that a $C$-category $I$ is {\em excisive} if all the maps $c$ in
$C$ are monomorphisms, and all the standard pushout squares
\eqref{st.CW.sq} are cartesian. Then for any compact object $i \in
I$ in an abundant excisive $C$-category $I$, and any cofibrant
functor $X:J \to I$ from some $J \in \Posf$, the functor $I(i,X):J
\to \Sets$ is cofibrant with respect to the trivial $C$-structure of
Example~\ref{CW.sets.exa}.
\end{lemma}

\proof{} Since all maps in $C$ are monomorphisms, $I(i,-):I \to
\Sets$ sends them to injective maps, and then by
Definition~\ref{reedy.CW.def}, it suffices to show that for any $J
\in \Posf$ and cofibrant functor $X:J \to I$, the natural map
$\colim_JI(i,X) \to I(i,\colim_JX)$ is injective. The same induction
on dimension as in Lemma~\ref{ab.CW.le} reduces this to the case $J
= S^<$, and then since $i$ is compact, $I(i,-)$ commutes with
filtered colimits, so by Example~\ref{ab.filt.exa}, we may assume
that $S$ is finite. But induction on $|S|$ further reduces us to the
case $S=\{0,1\}$, and this immediately follows from the definition.
\endproof

\begin{exa}
The category $\Pos$ with the standard or co-standard $C$-structure
of Example~\ref{pos.CW.exa} is excisive in the sense of
Lemma~\ref{exc.CW.le}, and so are the subcategories
$\Posf,\Posf^+,\Posf^- \subset \Pos$.
\end{exa}

One immediate consequence of Definition~\ref{ab.def} is that an
abundant $C$-category $I$ has arbitrary coproducts (this is
Lemma~\ref{ab.CW.le} for $J=S$). In particular, for any $i \in I$
and set $S$, we can define the product $S \times i$ as the colimit
of the constant functor $S \to I$ with value $i$. Moreover, for any
$i' \in I$, we can let $S = I(i,i')$, and we then have the
tautological evaluation map $S \times i \to i'$. This has the
following application. Assume given a set $S$ and a collection
$\{f_s\}$, $s \in S$ of maps $f_s:i_s \to i_s'$ in $I$ that are
absolutely cofibrant in the sense of Definition~\ref{abs.C.def}, and
define the {\em extension functor} $\Ex(\{f_s\}):I \to I$ by the
cocartesian square
\begin{equation}\label{ex.I.sq}
\begin{CD}
\coprod_{s \in S}I(i_s,i) \times i_s @>{\copr f_s}>> \coprod_{s \in
  S}I(i_s,i) \times i'_s\\
@VVV @VVV\\
i @>>> \Ex(\{f_s\})(i)
\end{CD}
\end{equation}
for any $i \in I$, where the vertical arrow on the left is the
coproduct of the evaluation maps. Note that the pushout
\eqref{ex.I.sq} indeed exists in $I$, since alternatively, it can be
expressed as
\begin{equation}\label{ex.I.eq}
\Ex(\{f_s\})(i) = \colim_{T^<}i_\idot,
\end{equation}
where $T = \copr_s I(i_s,i)$, and $i_\idot:T^< \to I$ sends $o$ to $i$ and
$t \in T$ corresponding to some map $i_s \to i$ to $i \copr_{i_s}
i_s'$. Since $\Ex(\{f_s\})$ is a functor, it can also be applied
pointwise to $I$-valued functors $X:J \to I$ from some $J \in
\Posf$, and the expression \eqref{ex.I.eq} is also functorial.

\begin{lemma}\label{ext.ab.le}
Assume that an abundant $C$-category $I$ is excisive in the sense of
Lemma~\ref{exc.CW.le}. Then for any $J \in \Posf$, cofibrant functor
$X:J \to I$ and collection $\{f_s\}$ of absolutely cofibrant maps
$f_s:i_s \to i'_s$ in $I$ such that $i_s$, $s \in S$ are compact,
the map $X \to X_1 = \Ex(\{f_s\})(X)$ is in $C$.
\end{lemma}

\proof{} By \eqref{ex.I.eq}, for any $j \in J$, the latching object
$L(X_1,j)$ is given by
$$
L(X_1,j) \cong \colim_{L(T,j)^<}L(X,j)_\idot
$$
where $L(T,j)$ is the latching set for the functor $T:J \to \Sets$
corresponding to $X$. But $T$ is cofibrant by
Lemma~\ref{exc.CW.le}, so the map $L(T,j) \to T$ is injective, and
then since $X$ is cofibrant, the latching map for $X_1$ is in $C$ by
Definition~\ref{ab.def}.
\endproof

\subsection{Relatively cofibrant functors.}

We note that the implication in Lemma~\ref{CW.J.le}~\thetag{iii}
only goes in one direction: a functor $X:J \to I$ may factor through
$\Cof(I)_C$ without being cofibrant. To measure the difference, it
is convenient to generalize Definition~\ref{reedy.CW.def} as
follows. Assume given a finite partially ordered set $J$ and a
subset $\J \subset L(J)$, with the incidence subset $\J_\idot
\subset J \times L(J)$ and its natural projections $\sigma(J)$,
$\tau(J)$ of \eqref{s.t.J}.

\begin{defn}\label{CW.J.def}
A functor $X:J \to I$ to a C-category $I$ is {\em $\J$-cofibrant} if
$\tau(J)_!\sigma(J)^*X:\J \to I$ exists and factors through
$\Cof(I)_C$.
\end{defn}

If $\J = L(J)$, then by Lemma~\ref{CW.J.le}, $X$ is $\J$-cofibrant
if and only if it is cofibrant in the sense of
Definition~\ref{reedy.CW.def}. On the other hand, if $\J = \Y(J)
\subset L(J)$ is the image of the Yoneda embedding
\eqref{yo.pos.eq}, then $X \cong \tau_!\sigma^*X$ is $\J$-cofibrant
iff it factors through $\Cof(I)_C$.

\begin{lemma}\label{CW.JJ.le}
Assume given a map $\phi:J \to J'$ of finite partially ordered sets,
and a subset $\J \subset L(J)$ such that $\J' = (\phi^*)^{-1}(\J)
\subset L(J')$ contains the Yoneda image $\Y(J') \subset L(J')$.
\begin{enumerate}
\item For any $\J$-cofibrant $X:J \to I$, $\phi_!X$ exists and is
  $\J'$-cofibrant.
\item If $\phi$ is a left-closed full embedding, then $X:J' \to I$
  is $\J'$-cofibrant iff $\phi^*X:J \to I$ is $\J$-cofibrant and for
  any $j \in J' \ssetminus J$, the latching object $L(X,j)$ exists
  and the map $L(X,j) \to X(j)$ is in $C$.
\end{enumerate}
\end{lemma}

\proof{} For \thetag{i}, by assumption, the map $\phi^* \circ \Y:J'
\to L(J)$ factors through a map $\Y':J' \to \J$, and we have a
commutative diagram
$$
\begin{CD}
J @<{\sigma}<< J /_\phi J' @>{\tau}>> J'\\
@| @VVV @VV{e'}V\\
J @<{\sigma(J)}<< \J_\idot @>{\tau(J)}>> \J
\end{CD}
$$
with cartesian square on the right, so that $\phi_!X \cong
\tau_!\sigma^*X \cong {\Y'}^*\tau(J)_!\sigma(J)^*X$ exists by
Lemma~\ref{bc.le}. Similarly, $\tau(J')_!\sigma(J')^*\phi_!X \cong
\tau(J)_!\sigma(J)^*X \circ \phi^*$, so it factors through
$\Cof(I)_C$.

For \thetag{ii}, the ``only if'' part is obvious: by assumtion we
have $J/j \in \J'$, and $\phi^*(L(j)) = \phi^*(J/j)$. For the ``if''
part, as in Lemma~\ref{CW.J.le}~\thetag{i}, we need to check that
$\colim_{J'}X$ exists and is cofibrant for any $J_1' \in \J'$, and
for any $J_0' \subset J_1'$, $J_0' \in \J'$, the map $\colim_{J_0'}X
\to \colim_{J_1'}X$ is in $C$. We may restrict to $J_1'$, so that
$J'_1 = J'$, and then for any maximal element $j \in J' \ssetminus
J$, $\phi^*(J' \ssetminus \{j\}) = \phi^*(J') \in \J$. By induction
on $|J'|$, it suffices to consider the case $J_0' = J' \ssetminus
\{j\}$, and assume that $\colim_{J_0'}X$ exists and is
cofibrant. The claim then follows from the cocartesian square
\eqref{J.j.sq}.
\endproof

By Lemma~\ref{P.adj.le}, an example of the situation of
Lemma~\ref{CW.JJ.le} occurs when a left-closed embedding $\phi:J \to
J'$ is right-reflexive, and $\J = \Y(J)$. In this case, we will say
that a functor $X:J' \to I$ is {\em cofibrant relatively to $J$} if
it is $\J'$-cofibrant. More generally, a left-closed full Reedy
subcategory $J_0 \subset J$ in a thin Reedy category $J$ is
automatically thin, and if $J_0 \subset J$ is right-admissible, with
the adjoint functor $\gamma:J \to J_0$, then for any $j \in J$,
$J_0/\gamma(j) \subset J/j$ is a left-closed right-admissible subset
in the finite partially ordered set $J/j$. Say that a functor $X:J
\to X$ is {\em cofibrant relatively to $J_0 \subset J$} if for any
$j \in J$, $\sigma(j)^*X$ is cofibrant relatively to $J'/\gamma(j)
\subset J/j$.

\begin{lemma}\label{thin.le}
Assume given a thin Reedy category $J$, a left-closed full Reedy
subcategory $J_0 \subset J$, and a functor $X:J \to I$ to a
CW-category $I$ whose restriction to $J_0$ is cofibrant. Then there
exists a cofibrant functor $X':J \to I$ and a map $a:X' \to X$ such
that $a$ is an isomorphism over $J_0$, and for any map $f:j_0 \to j$
in $J$, $j_0 \in J_0$ such that $X(f) \in W$, $X'(f)$ is in $C \cap
W$. Moreover, if we have a left-closed Reedy subcategory $J_{00}
\subset J_0$ that is right-admissible in $J$, hence also in $J_0$ by
Example~\ref{cl.adj.exa}, then the same holds for functors cofibrant
relatively to $J_{00}$.
\end{lemma}

\proof{} Since $J_0 \subset J$ is left-closed, we can build $X'$ by
induction on skeleta, starting with $X'|_{J_0}$. Then as in
Proposition~\ref{CW.prop}, we may assume by induction that $X'$ and
$a$ has been already constructed on a left-closed subcategory $J'
\subset J$, $J \supset J_0$ such that $J \ssetminus J'$ consists of a
single object $j \in J$, consider the map $a(j):L(X',j) \to X(j)$
that exists by Lemma~\ref{CW.J.le}~\thetag{i}, and construct $X'(j)$
and $a:X'(j) \to X(j)$ by taking its decomposition
$$
\begin{CD}
L(X',j) @>{c}>> X'(j) @>>> X(j)
\end{CD}
$$
of Definition~\ref{CW.def}~\thetag{ii}. Then for any $f:j_0 \to j$,
$j_0 \in J_0$, the natural map $c:X(j_0) \cong X'(j_0) \to L(X',j)$
is in $C$ by Lemma~\ref{CW.J.le}~\thetag{iv}, and $a(j) \circ c =
X(f)$, so that if $X(f)$ is in $W$, then $X'(f)$ must be in $C \cap
W$ by Defi\-nition~\ref{CW.def}~\thetag{ii}. This finishes the proof
in the absolute case; in the relative case, replace
Lemma~\ref{CW.J.le}~\thetag{i} with
Lemma~\ref{CW.JJ.le}~\thetag{ii}.
\endproof

\begin{exa}\label{pos.thin.exa}
Take $I = \Pos$, with the standard CW-structure of
Example~\ref{pos.CW.exa}, and assume that $J_0 = \sk_0J$ is
discrete, and there is no $J_{00}$. Then a functor $X:J \to \Pos$
defines a cofibration $J' \to J$ with fibers $J'_j = X(j)$, $j \in
J$, and if we use the decompositions \eqref{f.dec} given by
\eqref{cyl.deco}, the functor $X':J \to \Pos$ provided by
Lemma~\ref{thin.le} is $X'(j) = J'/j$. Note that since $J/j$ is a
partially ordered set for any $j \in J$ by
Definition~\ref{thin.def}, so is $J'/j$, and the functor $X'$ is
cofibrant by Lemma~\ref{pos.cof.le}.
\end{exa}

\begin{exa}\label{thin.exa}
Let $J = \II$ be the Kronecker category of
Subsection~\ref{V.I.subs}, let $J_0 = \ppt \subset \II$ be the full
subcategory spanned by $0 \in \II$, and let $X:\II \to I$ be a
constant functor with some value $i \in I$. Then $\II$ is a thin
Reedy category with respect to the degree function $\deg(l)=l$,
$l=0,1$, and as soon as $i \in I$ is cofibrant, we are within the
scope of Lemma~\ref{thin.le}. Then the latching object $L(X,1) \in
I$ is $i \copr i$, and to construct the cofibrant $X':\II \to I$,
one need to apply Definition~\ref{CW.def}~\thetag{ii} to the
codiagonal map $i \copr i \to i$. Explicitly, \eqref{f.dec} in this
case reads as
\begin{equation}\label{cod.dec}
\begin{CD}
i \copr i @>{c}>> i' @>{w}>> i,
\end{CD}
\end{equation}
and $X'$ sends $0$ to $i$, $1$ to $\wt{i}$, and $s$ resp.\ $t$ to
the composition of $c$ with the embedding $i \to i \copr i$ onto the
first resp.\ second summand, while $a:X \to X'$ is $\id$ at $0$ and
$w$ at $1$. We note that diagrams \eqref{cod.dec} are in one-to-one
correspondence with functors $\wt{i}:\D_1 \to I$, where $\D_1$ is as
in \eqref{D.n.eq}; the correspondence is given by $i = \wt{i}([0])$,
$i' = \wt{i}([1])$, $c = \wt{i}(s) \copr \wt{i}(t)$, and
$w=\wt{i}(e)$, where $e:[1] \to [0]$ is the unique projection. If we
equip $\D_1 \subset \Delta$ with the induced Reedy structure, its
latching category is exactly $\II$, and the functor $X':\II \to I$
is $l^*\wt{i}$, where $\wt{i}:\D_1 \to I$ corresponds to
\eqref{cod.dec}, and $l:\II \to \D_1$ is the embedding functor.
\end{exa}

\begin{exa}\label{thin.bis.exa}
In the setting of Example~\ref{thin.exa}, one can also consider the
restriction $v:\V^o \subset \VV \to \II$ of the discrete fibration
$v$ of \eqref{VV.dia} to the first component of \eqref{VV.eq}. Then
$v:\V^o \to \II$ is still a discrete fibration (in fact, one has
$\V^o \cong \II/1$, and $v$ is the fibration $\sigma(1):\II/1 \to
\II$ of \eqref{p.i}), $\V^o$ is a thin Reedy category with respect
to the induced Reedy structure, and $v^*X':\V^o \to I$ is also
cofibrant and obtained by applying Lemma~\ref{thin.le} to the
constant functor $v^*X:\V^o \to I$ with value $i$. More generally,
one take a map $f:i \to i'$ in $\Cof(I)$, and define a functor
$X:\V^o \to I$ by $X(0)=i$, $X(1)=X(o)=i'$, with the maps $f:X(0)
\to X(o)$, $\id:X(1) \to X(o)$. Then $L(X,o) \cong i \copr i'$, and
if one applies Lemma~\ref{thin.le}, then instead of \eqref{cod.dec}
one has to consider a decomposition \eqref{f.i.dec} used in
Lemma~\ref{CW.W.le}.
\end{exa}

\subsection{CW-augmentations.}

Since a finite partially ordered set $J$ is trivially left-finite,
Proposition~\ref{CW.prop} immediately provides a CW-structure on the
category $\Cof(\Fun(J,I))$ for any CW-category $I$. It turns out that one
can collect all these structures together into a single CW-structure
on the category $\Posff \bb I$. To do this, we modify
Example~\ref{pos.I.CW.exa} in the following way. Say that an
$I$-augmented finite partially ordered set $\langle J,\alpha_J
\rangle$ is {\em cofibrant} if $\alpha_J:J \to I$ is cofibrant in the
sense of Definition~\ref{reedy.CW.def}, and let $\Cof(\Posff \bb I)
\subset \Posff \bb I$ be the full subcategory spanned by cofibrant
objects.

\begin{exa}\label{pos.C.pos.exa}
If $I=\Posff$ with the CW-structure of Example~\ref{pos.CW.exa},
then by Lemma~\ref{pos.cof.le}, $\Cof(\Posff \bb \Posff)$ is the
category $\Arr(\Posff)$ whose objects are maps $f:J \to J'$ in
$\Posff$, and whose maps are given by lax commutative squares
\begin{equation}\label{arr.sq}
\begin{CD}
J_0 @>{f_0}>> J_0'\\
@V{g}VV @VV{g'}V\\
J_1 @>{f_1}>> J_1',
\end{CD}
\end{equation}
or explicitly, squares \eqref{arr.sq} such that $f_1 \circ g \leq g'
\circ f_0$. In particular, we have a dense embedding $\Ar(\Posff)
\subset \Arr(\Posff) \cong \Cof(\Posff \bb \Posff)$.
\end{exa}

\begin{defn}\label{w.CW.def}
A map $\langle f,\alpha_f \rangle:\langle J,\alpha_J \rangle \to
\langle \langle J',\alpha_{J'} \rangle$ between cofibrant
$I$-augmented finite partially ordered sets $\langle J,\alpha_J
\rangle$, $\langle J',\alpha_{J'} \rangle$ is {\em CW-left-closed}
resp.\ {\em CW-left-reflexive} if $f:J \to J'$ is a left-closed
resp.\ left-reflexive full embedding, and $\alpha_f:\alpha_J \to
f^*\alpha_{J'}$ is in the class $C$ resp.\ $W$ with respect to the
CW-structure of Proposition~\ref{CW.prop}. A map is {\em
  CW-right-reflexive} if it is a right-reflexive full embedding in
the sense of Definition~\ref{aug.def}, and {\em $CW$-reflexive} if
it is a finite composition of CW-right-reflexive and
CW-left-reflexive full embeddings.
\end{defn}

By definition, a CW-right-reflexive full embedding of augmented
partially ordered sets is strict, and by
Lemma~\ref{CW.J.le}~\thetag{iii}, a CW-left closed
resp.\ CW-left-reflexive full embedding is $C$-strict
resp.\ $W$-strict. Note that since the pullback $f^*X$ of a
cofibrant functor $X:J' \to I$ with respect to a map $f:J \to J'$ of
finite partially ordered sets is not necessarily cofibrant, the
projection
\begin{equation}\label{CW.bb}
\Cof(\Posff \bb I) \subset \Posff \bb I \to \Posff
\end{equation}
induced by the fibration of Example~\ref{bb.exa} is not a
fibration. However, by virtue of Lemma~\ref{CW.J.le}~\thetag{ii}, it
is a cofibration. The embedding \eqref{ups.eq} restricts to an
embedding $\ups:\Cof(I) \cong \Cof(\Posff \bb I)_{\ppt} \to
\Cof(\Posff \bb I)$, and it has a left-adjoint functor
$\chi:\Cof(\Posff \bb I) \to \Cof(I)$. By Lemma~\ref{bc.le} and
Lemma~\ref{CW.J.le}~\thetag{i}, $\chi$ can also be described as
\begin{equation}\label{chi.pos}
\chi = (\nu_\idot)_!\ev:\Cof(\Posff \bb I) \to I,
\end{equation}
where $\ev$ and $\nu_\idot$ are obtained by restricting the functors
of Example~\ref{cat.bb.exa} to $\Cof(\Posff_\idot \bb I) =
\Posff_\idot \times_{\Posff} \Cof(\Posff \bb I) \subset \Posff_\idot
\bb I$.

\begin{prop}\label{w.CW.prop}
For any weak CW-category $I$, the corresponding category
$\Cof(\Posff \bb I)$ with the classes $C$ resp.\ $W$ of
CW-left-closed resp.\ CW-reflexive full embeddings is a weak
CW-category. The functor \eqref{chi.pos} is a CW-functor, and for
any family of groupoids $\C$ over $\langle I,W \rangle$,
$(\nu_\idot)_*\ev^*\C$ is a family of groupoids over $\langle
\Cof(\Posff \bb I),W \rangle$.
\end{prop}

\proof{} For any morphism $\langle f,\alpha_f\rangle:\langle J,\alpha_J
\rangle \to \langle J',\alpha_{J'}\rangle$ in $\Cof(\Posff \bb I)$,
the left Kan extension $f_!\alpha_J:J' \to I$ exists and is
cofibrant by Lemma~\ref{CW.J.le}~\thetag{i}, and if $\langle
f,\alpha_f \rangle$ is CW-left-closed, then the map
$\alpha^\dg_f:f_!\alpha_J \to \alpha_{J'}$ adjoint to $\alpha_f$ is
in the class $C$ by Lemma~\ref{CW.J.le}~\thetag{ii},\thetag{iv}. To
construct the coproduct of $\langle f,\alpha_f \rangle$ with another
CW-left-closed map $\langle g,\alpha_g\rangle:\langle J,\alpha_J
\rangle \to \langle J'',\alpha_{J''}\rangle$, take the finite
partially ordered set $J' \copr_J J''$, and augment it by
\begin{equation}\label{tau.J.copr}
\alpha_{J' \copr_J J''} = \eps'_!\alpha_{J'} \copr_{\eps_!\alpha_J}
\eps''_!\alpha_{J''},
\end{equation}
where $\eps$, $\eps'$, $\eps''$ are the embeddings of $J$, $J'$,
$J''$ into $J' \copr_J J''$. Moreover, if $\langle f,\alpha_f \rangle$
is CW-reflexive, then $J'$ admits a filtration \eqref{j.filt}, with
all the maps $\langle f_l,\alpha_l \rangle:J'_{l-1} \to J'_l$ CW-left
or right-reflexive, and to prove that the pushout $\langle
f',\alpha_f'\rangle:J'' \to J' \copr_J J''$ of the map $\langle f,\alpha_f
\rangle$ is CW-reflexive, it suffices to prove the same for the
pushouts $\langle f'_l,\alpha'_l \rangle:J'_{l-1} \copr_J J'' \to J'_l
\copr_J J''$ of the maps $\langle f_l,\alpha_l \rangle$. But if
$\langle f_l,\alpha_l \rangle$ is CW-right-reflexive, the situation is
the same as in Example~\ref{pos.I.CW.exa}, and in the
CW-left-reflexive case, the map $\alpha_l'$ is a pushout
\eqref{tau.J.copr} of the map $\alpha_l$ with respect to a map in $C$,
thus lies in $W$ (where we recall that by
Definition~\ref{reedy.CW.def}, $\alpha_l \in W$ means that $\alpha_l \in
C \cap W$). Therefore $\Cof(\Posff \bb I)$ satisfies
Definition~\ref{CW.def}~\thetag{i}.

For Definition~\ref{CW.def}~\thetag{ii}, consider a map $\langle
f,\alpha_f \rangle:\langle J,\alpha_J\rangle \to \langle
J',\alpha_{J'}\rangle$ in $\Cof(\Posff \bb I)$, let $J'' = \Cyl(f)$,
with the decomposition \eqref{cyl.deco}, choose a decomposition
\begin{equation}\label{tau.dec}
\begin{CD}
f_!\alpha_J @>{c}>> \alpha'_{J'} @>{w}>> \alpha_{J'}
\end{CD}
\end{equation}
of the adjoint map $\alpha^\dg_f:f_!\alpha_J \to \alpha_{J'}$ provided by
Definition~\ref{CW.def}~\thetag{ii}, and let $\alpha_{J''} =
t_\dg^*\alpha'_J:J'' \to I$. Then \eqref{cyl.deco} lifts to a
decomposition
\begin{equation}\label{tau.J.dec}
\begin{CD}
\langle J,\alpha_J \rangle @>{\langle a,\alpha_a \rangle}>> \langle
J'',\alpha_{J''} \rangle @>{\langle b,\alpha_b \rangle}>> \langle
J',\alpha_{J'} \rangle
\end{CD}
\end{equation}
of the map $\langle f,\alpha_f \rangle$, where $a = s:J \to J''$, $b
= t_\dg:J'' \to J'$, the map $\alpha_a:\alpha_J \to a^*\alpha_{J''}
\cong f^*\alpha'_{J}$ is adjoint to the map $c$ in \eqref{tau.dec},
and the map $\alpha_b=t_\dg^*w$ is the pullback of the map $w$. Then
for any CW-left-closed morphism $\langle p,\alpha_p \rangle:\langle
J',\alpha_{J'} \rangle \to \langle J,\alpha_J \rangle$ such that
$\langle f,\alpha_f \rangle \circ \langle p,\alpha_p \rangle = \id$,
the map $\id:\alpha_{J'} \cong (f \circ p)_!\alpha_{J'} \to
\alpha_{J'}$ factors as
$$
\begin{CD}
(f \circ p)_!\alpha_{J'} @>{f_!(\alpha^\dg_p)}>> f_!\alpha_J
  @>{\alpha_f^\dg}>> \alpha_{J'}.
\end{CD}
$$
The cylinder $\Cyl(\id)$ of the identity map $\id:\langle
J',\alpha_{J'} \rangle \to \langle J',\alpha_{J'} \rangle$ is the
product $J' \times [1]$ equipped with the augmentation $\alpha_{J'
  \times [1]} = t_\dg^*\alpha_{J'}$, Lemma~\ref{cyl.le} provides a
left-reflexive full embedding $q:J' \times [1] \to J''$, and if we
let
$$
\alpha_q = t_\dg^*(c \circ f_!(\alpha^\dg_p)):\alpha_{J' \times [1]} \to
q^*\alpha_{J''} \cong t_\dg^*\alpha'_{J'},
$$
then the map $\langle q,\alpha_q \rangle$ is CW-left-reflexive by
virtue of Definition~\ref{CW.def}~\thetag{ii}, while the map
$\langle \sigma,\id \rangle:\langle J',\alpha_{J'} \rangle \to \langle
J' \times [1],\alpha_{J' \times [1]}\rangle$ is tautologically
CW-right-\-refle\-xive. Therefore $\langle a,\alpha_a \rangle \circ
\langle p,\alpha_p \rangle \cong \langle q,\alpha_q \rangle \circ
\langle \sigma,\id \rangle$ is CW-reflexive, and \eqref{tau.J.dec}
satisfies the condition of Definition~\ref{CW.def}~\thetag{ii}.

By Lemma~\ref{CW.J.le}~\thetag{ii},\thetag{iv}, the left Kan
extension \eqref{chi.pos} sends maps in $C$ to maps in $C$, it
trivially sends $0$ to $0$, and \eqref{tau.J.copr} immediately shows
that it sends pushout squares of Definition~\ref{CW.def}~\thetag{i}
to pushout squares. Thus to check that it is a CW-functor, it
suffices to show that for any CW-left-reflexive or
CW-right-reflexive map $\langle f,\alpha_f \rangle:\langle J,\alpha_J
\rangle \to \langle J',\alpha_{J'} \rangle$, the composition map
\begin{equation}\label{kkk}
\begin{CD}
\colim_J\alpha_J @>{\colim_J\alpha_f}>> \colim_Jf^*\alpha_{J'} @>>>
\colim_{J'}\alpha_{J'}
\end{CD}
\end{equation}
is in $W$. In the left-reflexive case, the first map in \eqref{kkk}
is in $W$ by Lemma~\ref{CW.J.le}~\thetag{i}, and the second map is
an isomorphism by \eqref{adm.eq}. In the right-reflexive case, the
whole map is also $\colim_{J'}\alpha^\dg_f$, and $\alpha^\dg_f$ is an
isomorphism by Definition~\ref{aug.def}. Finally, to see that
$(\nu_\idot)_*\ev^*\C$ is constant over $W$, use the same argument with
Lemma~\ref{bc.le} and \eqref{adm.eq} replaced by \eqref{2kan.eq} and
Lemma~\ref{2adm.le}.
\endproof

\begin{remark}\label{C.bb.rem}
The category $(\nu_\idot)_*\ev^*\C$ of Proposition~\ref{w.CW.prop}
can be also described by \eqref{nu.ev.bb}: it is the full
subcategory $\Cof(\Posff\bb\C) \subset \Posff\bb\C$ spanned by
$\C$-augmented partially ordered sets whose projection to $\Posff\bb
I$ is cofibrant.
\end{remark}

\begin{remark}\label{w.CW.rem}
Definition~\ref{w.CW.def} immediately implies that the cofibration
\eqref{CW.bb} is a CW functor, with fibers $\Cof(\Posff \bb I)_J
\cong \Cof(\Fun(J,I))$, and for any $J \in \Posff$, the embedding
$\ups_J:\Cof(\Fun(J,I)) \to \Cof(\Posff \bb I)$ of the fiber over
$J$ is a weak CW-functor (this includes the embedding $\ups =
\ups_{\ppt}:\Cof(I) \to \Cof(\Posf \bb I)$ of
\eqref{ups.eq}). However, these are only weak CW-functors, since $0
\in \Fun(J,I)$ goes to $\langle J,0 \rangle$ and not to $\langle
\emptyset,0 \rangle$.
\end{remark}

\begin{remark}\label{w.CW.thin.rem}
Alternatively, one can construct decompositions \eqref{f.dec} for
the CW-category $\Cof(\Posff \bb I)$ by using
Lemma~\ref{thin.le}. Namely, the augmented cylinder $\langle
J'',\alpha_{J''} \rangle$ of a map $\langle f,\alpha_f \rangle:\langle
J,\alpha_J\rangle \to \langle J',\alpha_{J'}\rangle$ between cofibrant
$I$-augmented finite partially ordered sets is not in general
cofibrant, but $s:J \to J''$ is a left-closed embedding, and
$s^*\alpha_{J''} \cong \alpha_J:J \to I$ is cofibrant. If we treat $J''$
as a Reedy category, by Example~\ref{reedy.exa}, then
Lemma~\ref{thin.le} provides a cofibrant functor $\alpha'_{J''}:J''
\to I$ and a map $\alpha'_{J''} \to \alpha_{J''}$ that is an isomorphism
over $J$. Replacing $\alpha_{J''}$ with $\alpha'_{J''}$ gives
\eqref{tau.J.dec}. The construction used in the proof of
Proposition~\ref{w.CW.prop} is essentially the same but shorter:
since we know that $t^*\alpha_{J''}:J' \to I$ is also cofibrant, it
suffices to take the decomposition \eqref{tau.dec} once instead of
doing the induction over $j \in J' = J'' \ssetminus J$ as in
Lemma~\ref{thin.le}.
\end{remark}

\begin{exa}\label{w.CW.exa}
To illustrate how \eqref{f.dec} and \eqref{cod.dec} work in the
CW-category $\Cof(\Posff \bb I)$, take a cofibrant object $i \in I$,
and consider its image $\ups(i) = \langle \ppt,i \rangle$ with
respect to the embedding \eqref{ups.eq}. Note that $\ups(i) \copr
\ups(i)$ is $\{0,1\} \cong \ppt \copr \ppt$ augmented by the
constant functor $i:\{0,1\} \to I$ with value $i$, and not $\eps(i
\copr i) \cong \langle \ppt, i \copr i \rangle$. To construct
\eqref{cod.dec} for $\ups(i)$, one first takes the cylinder $\V^o$
of the codiagonal map $\{0,1\} \cong \ppt \copr \ppt \to \ppt$, and
then chooses a decomposition \eqref{cod.dec} for $i$ in $I$, with
the corresponding functor $\wt{i}:\II \to I$. The resulting
decomposition is
\begin{equation}\label{w.cod.dec}
\begin{CD}
\ups(i) \copr \ups(i) \cong \langle \{0,1\},i \rangle @>{c}>> \langle
\V^o,v^*\wt{i} \rangle @>{w}>> \ups(i) = \langle \ppt,i \rangle,
\end{CD}
\end{equation}
where $v:\V^o \to \II$ is as in Example~\ref{thin.bis.exa}.
\end{exa}

\subsection{Regular CW-categories.}

One other result of Section~\ref{mod.sec} that survives, with some
modifications, in the generality of Definition~\ref{CW.def} is
Proposition~\ref{mod.prop}. Namely, assume given a CW-category $I$,
with the corresponding full subcategory $\Cof(I) \subset I$ of
cofibrant objects and classes $C$ and $W$.

\begin{defn}\label{reg.CW.def}
A CW-category $I$ is {\em regular} if the natural embedding
$\nu:\langle \Cof(I)_C,C \cap W\rangle \to \langle I,W \rangle$ is
a $2$-equivalence.
\end{defn}

If $I$ is a model category $I$ with the CW-structure of
Example~\ref{mod.CW.exa}, then $I$ is regular in the sense of
Definition~\ref{reg.CW.def} by Proposition~\ref{mod.prop}. One
cannot expect a general CW-category $I$ to be regular: at the very
least, nothing in Definition~\ref{reg.CW.def} forbids adding a new
non-cofibrant object $i$ to a CW-category $I$ without changing the
classes $C$ and $W$ (nor in fact having a CW-category whose only
cofibrant object is $0 \in I$). But as it happens, excluding trivial
possibilities of this sort does the job.

\begin{prop}\label{reg.prop}
Assume given a CW-category $\langle I,C,W \rangle$ such that
$\Cof(I)=I$. Then $I$ is regular in the sense of
Definition~\ref{reg.CW.def}.
\end{prop}

\proof{} For any category $I$, let $e:[1] \times I \to I$ be the
projection, with its sections $s,t:I \to [1] \times I$ given by
embeddings onto $\{0\} \times I$ resp.\ $\{1\} \times I$. By
Lemma~\ref{triv.con.le}~\thetag{ii}, $e:\langle [1] \times I,e^*\Iso
\rangle \to \langle I,\Iso \rangle$ is a $2$-equivalence, with
inverse $2$-equivalence given by either $s$ or $t$. The restricted
simplicial expanion $\Delta_\ff I$ of Corollary~\ref{del.corr} is a
thin Reedy category in the sense of Definition~\ref{thin.def}, and
we have a commutative diagram
\begin{equation}\label{2.dia.eq}
\begin{CD}
\langle \Delta_\ff I,+ \rangle @>{\Delta_\ff(s)}>> \langle \Delta_\ff([1]
\times I),\Delta_\ff(e)^*+ \rangle @<{\Delta_\ff(t)}<< \langle
\Delta_\ff I,+ \rangle\\
@| @VV{\xi'}V @|\\
\langle \Delta_\ff I,+ \rangle @>{s}>> \langle [1] \times
\Delta_\ff I,e^*(+) \rangle @<{t}<< \langle \Delta_\ff I,+ \rangle\\
@V{\xi}VV @VV{\xi}V @VV{\xi}V\\
\langle I,\Iso \rangle @>{s}>> \langle [1] \times I,e^*\Iso \rangle
@<<{t}< \langle I,\Iso \rangle,
\end{CD}
\end{equation}
where all the morphisms $\xi$ are $2$-equivalences by
Corollary~\ref{del.corr}, and so is the composition $\xi \circ
\xi'$, so that all the morphisms in \eqref{2.dia.eq} are
$2$-equivalences.

Now let $I$ be our CW-category, consider the thin Reedy category
$\Delta_\ff I$, and note that since $I = \Cof(I)$, the functor
$\xi:\Delta_\ff I \to I$ is cofibrant on the skeleton $\sk_0\Delta_\ff I
= (\Delta I)_{[0]} = I_{\Id} \subset \Delta_\ff I$. Thus if we take
$J=\Delta_\ff I$, $J_0=\sk_0\Delta_iI$, $X=\xi$, we can apply
Lemma~\ref{thin.le} and obtain a cofibrant functor $X':\Delta_\ff I
\to I$ and a morphism $a:X' \to X$. Since $X'$ is cofibrant, it
factors through $I_C$. Moreover, $X = \xi$ inverts all special maps,
so that for any $\langle [n],i_\idot \rangle \in \Delta_\ff I$, $X'$
sends the special map $t:\langle [0],i_n \rangle \to \langle
[n],i_\idot \rangle$ to a map in $C \cap W$. Let $W' = s(I,W)$, $(C
\cap W)'=s(I_C,C \cap W)$ be the saturations of the classes $W$ and
$C \cap W$ (note that $(C \cap W)' \subset C \cap W'$, so that $\nu$
extends to a morphism $\nu':\langle I_C,(C \cap W)' \rangle \to
\langle I,W' \rangle$, but the inclusion might be strict). Then
since $a \circ X'(t) = X(t) \circ a$ is an isomorphism, $a:X' \to X$
is pointwise in $W'$ and defines a morphism
\begin{equation}\label{st.X}
X'':\langle [1] \times \Delta_\ff I,e^*\xi^*W \rangle \to \langle
I,W' \rangle
\end{equation}
such that $X'' \circ s \cong X'$ and $X'' \circ t \cong X$. On the
other hand, for any map $f:\langle [n],i_\idot \rangle \to \langle
[n'],i'_\idot \rangle$ such that $X(f) \in W$, we have $X'(f) \circ
X'(t) \in C \cap W$ by Lemma~\ref{thin.le}, so that $X(f) \in (C
\cap W)'$. Therefore $X'$ factors as
\begin{equation}\label{XY.eq}
\begin{CD}
\langle \Delta_\ff I,\xi^*W) \rangle @>{Y}>> \langle I_C,(C \cap W)')
\rangle @>{\nu'}>> \langle I,W' \rangle
\end{CD}
\end{equation}
for some $Y$. Since $s$, $t$ and $\xi$ in \eqref{2.dia.eq} are
$2$-equivalences, the morphisms
$$
s,t:\langle \Delta_\ff I,\xi^*W \rangle \to \langle [1] \times
\Delta_\ff I,e^*\xi^*W\rangle, \ \xi:\langle \Delta_\ff I,\xi^*W
\rangle \to \langle I,W \rangle \to \langle I,W' \rangle
$$
are $2$-equivalence by Lemma~\ref{2.con.bis.le}~\thetag{i} and
Lemma~\ref{sat.le}, and since we have $\xi \cong X'' \circ t$, so is
$X''$ and therefore also $X' \cong X'' \circ t$. Thus $Y$ in
\eqref{XY.eq} is a one-sided inverse to $\nu'$.

To construct the inverse on the other side, consider the diagram
\eqref{2.dia.eq} for the category $I_C$, and note that $\Delta_i(s)$
and $\Delta_i(t)$ define left-closed embeddings
$$
\begin{CD}
\Delta_\ff I_C @>{b}>> \Delta_\ff I_C \copr \Delta_\ff I_C @>{\Delta_\ff(s)
  \copr \Delta_\ff(t)}>> \Delta_\ff([1] \times I_C),
\end{CD}
$$
where $b$ is the embedding onto the second component. Moreover, if
we augment all the categories by adding initial objects, then the
embeddings
\begin{equation}\label{b.bb.dia}
\begin{CD}
  \Delta^<_\ff I_C @>{b^<}>> (\Delta_\ff I_C \copr \Delta_\ff I_C)^<
  @>{\Delta^<_\ff(s) \copr \Delta^<_\ff(t)}>> \Delta^<_\ff([1] \times I_C)
\end{CD}
\end{equation}
are still left-closed embeddings of thin Reedy categories, and their
composition $\Delta^<(t)$ is the right-admissible embedding
\eqref{st.X.del}, with some right-adjoint $\Delta^<(t)^\dg$. On the
other hand, the morphism \eqref{st.X} fits into a commutative
diagram
$$
\begin{CD}
\langle [1] \times \Delta_ff I_C,e^*\xi^*(C \cap W) \rangle @>{Y'}>>
\langle I_C,(C \cap W)'\rangle\\
@V{\Delta\nu}VV @VV{\nu'}V\\
\langle [1] \times \Delta_\ff I,e^*\xi^*W \rangle @>{X''}>> \langle
I,W' \rangle
\end{CD}
$$
for some morphism $Y'$, and $s^*Y' \cong Y \circ \Delta\nu$ is a
cofibrant functor $\Delta_\ff I_C \to I$, while $t^*Y' \cong \nu \circ
\xi$, while not cofibrant, still factors through $I_C \subset I$. We
can then consider the restriction $Y'' \circ \xi'$, where $\xi'$ is
as in \eqref{2.dia.eq}, and extends it to a functor
$Z:\Delta_\ff^<([1] \times I_C) \to I_C$ by sending $o$ to $0 \in
I$. If we now take $J = \Delta_\ff^<([1] \times I_C)$, $J_0 =
(\Delta_\ff I_C \copr \Delta_\ff I_C)^<$, and $J_{00}=\Delta^<_\ff I_C$,
with the embeddings \eqref{b.bb.dia}, then we can apply
Lemma~\ref{thin.le} to the functor $Z$ and obtain a functor
$Z':\Delta^<_\ff([1] \times I_C) \to I$ cofibrant over
$\Delta^<_\ff I_C$, and a map $a:Z' \to Z$. Since $Z'$ is relatively
cofibrant, it factors through $I_C$. Moreover, $Z \cong {\xi'}^*Y'$
defines a morphism
\begin{equation}\label{Z.eq}
Z:\langle \Delta_\ff([1] \times I_C),\xi^*e^*(C \cap W)\rangle \to \langle
I_C,(C \cap W)'\rangle,
\end{equation}
and for any $j \in J \ssetminus J_0$, the adjunction map
$a(j):\Delta^<(t)\Delta^<(t)^\dg(j) \to j$ is inverted by $e \circ
\xi$. Then $a \circ Z'(f) \circ Z'(a(j)) = Z(f)$ for any map $f:j
\to j'$, so that by the same argument as before, if $e(\xi(f)) \in C
\cap W$, then $Z'(f)$ is in $(C \cap W)'$. Thus $Z'$ also defines a
morphism of the form \eqref{Z.eq}. Then $Z' \circ \Delta(t) \cong
\xi$ is a $2$-equivalence by Lemma~\ref{sat.le} and
\eqref{2.dia.eq}, and since both $\Delta(t)$ and $\Delta(s)$ are
also $2$-equivalences of \eqref{2.dia.eq}, $Z' \circ \Delta(s) \cong
Y \circ \Delta(\nu)$ is then also a $2$-equivalence. We thus have a
sequence of morphisms
$$
\langle \Delta_\ff I_C,\xi^*(C \cap W) \rangle
\overset{\Delta(\nu)}{\longrightarrow} \langle 
\Delta_\ff I,\xi^*W\rangle \overset{Y}{\longrightarrow} \langle
I_C,C \cap W' \rangle \overset{\nu'}{\longrightarrow} \langle I,W'
\rangle
$$
where both compositions are $2$-equivalences. Therefore all the
morphisms in the diagram are $2$-equivalences, and this includes
$\nu'$. To finish the proof, it remains to again invoke
Lemma~\ref{sat.le}.
\endproof

\begin{remark}
The heuristic idea behind the proof of Proposition~\ref{reg.prop} is
very simple, and it is perhaps useful to summarize it. We have a
commutative diagram
\begin{equation}\label{xi.c}
\begin{CD}
\langle \Delta_\ff I_C,\xi^*(C \cap W) \rangle @>{\xi}>> \langle I,C
\cap W \rangle\\
@V{\Delta(\nu)}VV @VV{\nu}V\\
\langle \Delta_\ff I,\xi^*W \rangle @>{\xi}>> \langle I,W \rangle
\end{CD}
\end{equation}
where the horizontal arrows are $2$-equivalences. To prove the
claim, it would suffice to construct a diagonal arrow $Y:\langle
\Delta_\ff I,\xi^*W \rangle \to \langle I,C \cap W \rangle$ such that
$\nu \circ Y \cong \xi$ and $Y \circ \Delta(\nu) \cong \xi$. This is
not possible --- $\xi$ in the bottom arrow in \eqref{xi.c} does not
factor through $I_C$. However, we can use Lemma~\ref{thin.le} to
construct a ``cofibrant replacement'' $X'$ of $X = \xi$, and a map
$a:X' \to X$. Being cofibrant, $X'$ does factor through $I_C$, and
since the map $a$ is pointwise in $W'$, the $2$-functors induced by
$\xi$ and by $X' = \nu \circ Y$ on relevant $2$-categories are the
same. The simplest way to make this rigourous is to note that both
$\xi$ and $\nu \circ Y$ factor through the product $[1] \times
\Delta_\ff I$ via the $2$-equivalences given by $s$ and
$t$. Unfortunately, $a$ is not pointwise in $C$, so one cannot apply
the same argument to identify the $2$-functors induced by
$\Delta(\nu) \circ Y$ and $\xi$ --- there is no map between them. To
correct this, we take the simplicial replacement $\Delta_i([1]
\times I_C)$ instead of the product $[1] \times \Delta_i I_C$, and
apply Lemma~\ref{thin.le} once more to obtain a cofibrant
replacement for the non-existent map.
\end{remark}

\begin{remark}
Proposition~\ref{reg.prop} can be thought of as a strengthening of
Lemma~\ref{CW.W.le} (and combining the two, we deduce that the
embedding $\langle \Cof(I)_C,C \cap W \rangle \to \langle \Cof(I),C
\cap W \rangle$ is also a $2$-equivalence for any CW-category
$I$). Our proof also elaborates on the same general idea. In fact,
the construction of the cofibrant replacement $X'$ provided by
Lemma~\ref{thin.le} works by induction on degree, and the first step
of this inductive construction is exactly the same as in
Example~\ref{thin.bis.exa}.
\end{remark}

\begin{exa}
Take $I=\Pos$, with the standard CW-structure of
Example~\ref{pos.CW.exa}. Then objects in $\Delta I$ correspond to
cofibrations $J \to [n]$, $[n] \in \Delta$, $J \in \Pos$. The
functor $X=\xi:\Delta I \to I$ sends such a cofibration to the fiber
$J_n \in \Pos$, while its ``cofibrant replacement'' $X':\Delta_\ff I
\to I$ sends it to the product $J \times_{[n]} B([n])$, where
$B([n])$ is as in Definition~\ref{pos.bary.def}.
\end{exa}

\subsection{CW-replacements.}

Unfortunately, Proposition~\ref{reg.prop} does not apply to weak
CW-categories (and in particular, to weak CW-categories of
Proposition~\ref{w.CW.prop}). One way to cure weakness is to enlarge
the class $W$ so that it becomes saturated; the prototype here is
passing from reflexive maps of Subsection~\ref{pos.exa.subs} to
anodyne maps of Subsection~\ref{ano.subs}. Here is the formal
definition.

\begin{defn}\label{ano.W.def}
The class $a(W)$ of {\em $W$-anodyne maps} in a CW-category $I$ is
the smallest class of maps in $I$ that contains $W$, is saturated in
the sense of Definition~\ref{sat.def}, and such that for any two
maps $f_0:i \to i_0$, $f_1:i \to i_1$ in $C$ such that $f_0 \in
a(W)$, the pushout map $i_1 \to i_0 \copr_i i_1$ is in $a(W)$.
\end{defn}

\begin{exa}\label{W.ano.exa}
Assume given an arrow $c:i_0 \to i_1$ in $\Cof(I)_C$ considered as a
cofibrant object $i \in \Ar(I)$, with the CW-structure of
Example~\ref{CW.ar.exa}, and take a decomposition \eqref{cod.dec} of
the map $i \copr i \to i$. Then $i'$ corresponds to an arrow $i_0'
\to i_1'$ in $C$, and either of the two maps $i_1 \to i_1'$ factors
as
\begin{equation}\label{CW.ano}
\begin{CD}
  i_1 @>{l}>> i_0' \copr_{i_0} i_1 @>{r}>> i_1'.
\end{CD}
\end{equation}
Both maps $i_0 \to i_0'$, $i_1 \to i_1'$ are in $C \cap W$ by
assumption, and then by Definition~\ref{CW.def}~\thetag{i}, so is
the map $l$ in \eqref{CW.ano}. However, the best one can say about
the map $r$ is that it is $W$-anodyne: it need not lie in $W$
(e.g. if $I$ is the category $\Pos$ with the CW-structure of
Example~\ref{pos.CW.exa}, it usually does not).
\end{exa}

\begin{lemma}\label{W.ano.W.le}
For any weak CW-category $\langle I,C,W \rangle$,
$\langle \Cof(I),C,a(W) \rangle$ is a CW-category.
\end{lemma}

\proof{} Definition~\ref{CW.def}~\thetag{i} is immediate. For
Definition~\ref{CW.def}~\thetag{ii}, note that since $a(W)$ is
saturated, a factorization \eqref{f.dec} satisfies the condition as
soon as $w \in a(W)$. To obtain such a factorization for any map
$f:i \to i'$, it suffices to use the cylinder construction of
Lemma~\ref{CW.W.le} and consider the factorization \eqref{f.i.dec}.
\endproof

In most examples, the class of $W$-anodyne maps is too large to
allow for an effective description but has better properties. In
particular, we have the following strengthening of
Lemma~\ref{CW.J.le}~\thetag{iii}.

\begin{lemma}\label{ano.W.le}
Assume given a left-finite map $\phi:J \to J'$ between partially
ordered sets, two functors $X_0,X_1:J \to I$ cofibrant in the sense
of Definition~\ref{reedy.CW.def}, and a map $f:X_0 \to X_1$ that is
pointwise in $a(W)$. Then so is $\phi_!(f):\phi_!X_0 \to \phi_!X_1$.
\end{lemma}

\proof{} As in the proof of Lemma~\ref{CW.J.le}~\thetag{i},
\eqref{pos.kan.eq} immediately shows that it is enough to consider
the case $J'=\ppt$, so that $J$ is left-finite. Then it is a thin
Reedy category in the sense of Definition~\ref{thin.def}, and
Proposition~\ref{CW.prop} provides a CW-structure on the functor
category $\Fun(J,I)$. Consider the decomposition
$$
\begin{CD}
X_0 \copr X_1 @>{c_0 \copr c_1}>> X @>{w}>> X_1
\end{CD}
$$
of the map $f \copr \id:X_0 \copr X_1 \to X_1$ provided by
Definition~\ref{CW.def}~\thetag{ii}. Then $c_1 \in W$, so that $w
\in a(W)$ since $w \circ c_1 = \id$, and then $c_0 \in a(W)$ since
$f = w \circ c_0 \in a(W)$. Therefore it suffices to prove the claim
for the maps $c_0$ and $c_1$, or in other words, we may assume that
$f$ is in the class $C$ of Definition~\ref{reedy.CW.def}. Finally,
let $S \subset J$ be the right-closed subset of maximal elements,
with discrete order, and as in the proof of Lemma~\ref{pos.epi.le},
extend the identity map $\id:S \to S$ to the universal map $p:J \to
S^<$ of \eqref{p.J.sq}. Then by induction on $\dim J$, we may assume
the claim proved for $\phi=p$, so that we are reduced to the case
$J=S^<$, $J'=\ppt$, $f \in C$. But then $f$ factors as
$$
\begin{CD}
X_0 @>{f_0}>> X_0' @>{f_1}>> X_1,
\end{CD}
$$
where $X_0'(o) = X_1(o)$ and $X'_0(s) = X_0(s) \copr_{X_0(o)}
X_1(o)$, $s \in S$, and it suffices to prove that $\colim_{S^<}f_0$
and $\colim_{S^<}f_1$ are $W$-anodyne. But $\colim_{S^<}f_0$ is a
standard pushout of the map $f(o):X_0(o) \to X_1(o)$, and
$\colim_{S^<}f_1$ is a composition of standard pushouts of maps
$f_1(s)$, $s \in S$.
\endproof

\begin{exa}\label{W.exc.exa}
Assume given a standard pushout square \eqref{st.CW.sq} of cofibrant
objects in $I$, treat it as a cofibrant object in $\Fun(\V,I)$, and again
take the decomposition \eqref{cod.dec}. Then the map $w:i' \to i$
induces a map
\begin{equation}\label{CW.exc}
  i_{01}' = i_0 \copr_i i' \copr_i i_1 \to i_{01} = i_0 \copr_i i_1,
\end{equation}
and this map is again $W$-anodyne: it factors as
$$
\begin{CD}
i_{01}' @>>> i_0' \copr_i i_1 @>>> i_0' \copr_{i'} i_1' @>>> i_0
\copr_i i_1 = i_{01},
\end{CD}
$$
the first two maps are standard pushouts of the anodyne maps of
Example~\ref{W.ano.exa}, and the last map is $W$-anodyne by
Lemma~\ref{ano.W.le}.
\end{exa}

As a corollary of Lemma~\ref{ano.W.le}, one can also prove a version
of Corollary~\ref{B.ano.corr}. The setup generalizes
Example~\ref{thin.bis.exa} that can reinterpreted in terms of
barycentric subdivision: if we identify $\V^o \cong B([1])$, with
the functor $\xi:\V^o \cong B([1]) \to [1]$ of \eqref{B.xi}, then we
have $X = \xi^*\eps(f)$, where $\eps(f):[1] \to \Cof(I)$ is the
embedding \eqref{eps.f} corresponding to the morphism $f$. More
generally, assume given a left-finite partially ordered set $J$ and
a functor $X:J \to \Cof(I)$. Then $BJ$ is left-finite, thus a thin
Reedy category, the $0$-skeleton $BJ_0 = \sk_0 BJ \subset BJ$ is a
left-closed subcategory, and since $J_0$ is discrete,
$\xi^*X|_{BJ_0}:BJ_0 \to \Cof(I)$ is cofibrant, so that
Lemma~\ref{thin.le} provides a cofibrant functor $X':BJ \to I$ and a
map $a:X' \to \xi^*X$ that restricts to an isomorphism over
$BJ_0$. By adjunction, $a$ induces a map
\begin{equation}\label{A.X.W}
\xi_!X' \to X,
\end{equation}
where the Kan extension on the left exists by
Lemma~\ref{CW.J.le}~\thetag{i}.

\begin{corr}\label{ano.W.corr}
In the assumptions above, the map \eqref{A.X.W} is pointwise
$W$-anodyne.
\end{corr}

\proof{} As in the proof of Corollary~\ref{xi.B.corr}, since $BJ/j
\cong B(J/j)$ for any $j \in J$, \eqref{kan.eq} shows that to prove
the claim at some $j \in J$, we may replace $J$ with $J/j$ and place
ourselves in the situation of Lemma~\ref{B.le}. Namely, we have a
finite partially ordered set $J$, a functor $X:J^< \to \Cof(I)$, and
its replacement $X':B(J^<) \to I$ provided by Lemma~\ref{thin.le}
that is equipped with the map $a:X' \to \xi^*X$, and we need to
check that the map $\colim(a):\colim_{B(J^<)}X' \to
\colim_{B(J^<)}\xi^*X \cong X(o)$ is $W$-anodyne, where the
isomorphism $\colim_{B(J^<)}\xi^*X \cong X(o)$ is
Corollary~\ref{xi.B.corr}. Then as in the proof of
Corollary~\ref{xi.B.corr}, we have a left-admissible full embedding
$i:B(J^<)^1 \to B(J^<)$, with the adjoint map $i^\dg:B(J^<) \to
B(J^<)^1$, and a right-admissible left-closed full embedding $p:\ppt
\cong B(J^<)_1 \to B(J^<)^1$, and it suffices to take the colimits
over $B(J^<)^1 \subset B(J^<)$. Over this subset, $i^*X' \cong
i^\dg_!X'$ is cofibrant by Lemma~\ref{CW.J.le}~\thetag{i}, and since
$\xi \circ i$ factors through $\{o\} \subset J^<$, $i^*\xi^*X$ is
constant. Thus $i^*\xi^*X \cong p_!p^*i^*\xi^*X \cong p_!X(o)$ is
cofibrant by Lemma~\ref{CW.J.le}~\thetag{iv}, and $i^*X'$ factors
through $I_{C \cap W}$ by the conditions of Lemma~\ref{thin.le}. But
the map $a$ is an isomorphism over $p(\ppt) \subset B(J^<)^1$, and
since $i^*\xi^*X$ is constant and $i^*X'$ factors through $I_{C \cap
  W}$, $a$ is pointwise $W$-anodyne. Therefore $\colim_{B(J^<)^1}a$
is $W$-anodyne by Lemma~\ref{ano.W.le}.
\endproof

Now assume given an arbitrary CW-category $I$, and consider the
CW-category $\Cof(\Posff \bb I)$ with the weak CW-structure of
Proposition~\ref{w.CW.prop}. To turn it into a CW-structure,
automatically regular, one can simply pass to the anodyne completion
of Lemma~\ref{W.ano.W.le}, but we can go for an even larger class
--- namely, the class $\chi^*a(W)$, where $\chi:\Cof(\Posff \bb I)
\to \Cof(I)$ is the functor \eqref{chi.pos} left-adjoint to the
embedding $\ups:\Cof(I) \to \Cof(\Posff \bb I)$. To make for a
stronger result, we can also replace $C$ with a smaller class ---
namely, we take the class of left-closed embeddings that are
$\Id$-strict in the sense of Subsection~\ref{pos.aug.subs}. The
formal definition is as follows.

\begin{defn}\label{s.C.def}
The {\em class $s(C)$} is the closed class spanned by morphisms
$\langle f,\alpha_f \rangle:\langle J,\alpha_J \rangle \to \langle
J',\alpha_{J'} \rangle$ in $\Cof(\Posff \bb I)$ such that $f$ is a
left-closed full embedding, and $\alpha_f$ is pointwise an identity
map.
\end{defn}

One can also rephrase Definition~\ref{s.C.def} in a way similar to
Definition~\ref{simpl.repl.def}. Namely, observe that since the
pullback $f^*$ with respect to a left-closed embedding $f:J \to J'$
sends cofibrant functors to cofibrant functors, the cofibration
\eqref{CW.bb} restricts to a bifibration $\Cof(\Posff \bb I)_C \to
\Posff_C$ with fibers $\Cof(J^I)_C$, $J \in \Posff$. In particular,
the fiber over $\ppt \in \Posff$ is $\Cof(I)$, and then the dense
subcategory $\Cof(\Posff \bb I)_{s(C)} \subset \Cof(\Posff \bb I)_C$
fits into a cartesian square
\begin{equation}\label{s.C.sq}
\begin{CD}
\Cof(\Posff \bb I)_{s(C)} @>>> \Cof(\Posff \bb I)_C\\
@VVV @VVV\\
\eps(\ppt)_*\Cof(I)_{\Id} @>>> \eps(\ppt)_*\Cof(I),
\end{CD}
\end{equation}
where $\eps(\ppt):\ppt \to \Posff$ is the embedding \eqref{eps.c}
onto $\ppt \in \Posff$. Because of this interpretation, we call
$\Cof(\Posff \bb I)_{s(C)}$ the {\em CW-replacement} of the
CW-category $I$. It is fibered over $\Posff$, and the fiber over
some $J \in \Posff$ is the discrete category
$\Cof(\Fun(J,I))_{\Id}$.

\begin{lemma}\label{CW.repl.le}
For any CW-category $I$, the category $\Cof(\Posff \bb I)$ with the
classes of maps $s(C)$, $\chi^*a(W)$ is a CW-category, and
the CW-functor $\chi:\langle \Cof(\Posff \bb I),\chi^*a(W) \rangle
\to \langle I,a(W) \rangle$ is a $2$-equivalence.
\end{lemma}

\proof{} As in Lemma~\ref{W.ano.W.le},
Definition~\ref{CW.def}~\thetag{i} is immediate. For
Definition~\ref{CW.def}~\thetag{ii}, use the same construction as in
Proposition~\ref{w.CW.prop}, but as in Lemma~\ref{W.ano.W.le},
choose a decomposition \eqref{tau.dec} induced from a decomposition
\eqref{f.i.dec}. Then $w$ is $W$-anodyne and pointwise $W$-anodyne,
so that it is in $\chi^*a(W)$ by Lemma~\ref{ano.W.le}, and then as
in Lemma~\ref{W.ano.W.le}, since the class $\chi^*a(W)$ is
saturated, the decomposition satisfies the condition of
Definition~\ref{CW.def}~\thetag{ii}. Finally, to see that $\chi$ is a
$2$-equivalence, note that the adjunction map $\chi \circ \ups \to
\id$ is an isomorphism, and the adjuntion map $\id \to \ups \circ
\chi$ is sent to an isomorphism by $\chi$, thus lies in
$\chi^*a(W)$. Thus we are done by Lemma~\ref{loc.adj.le}.
\endproof

Lemma~\ref{CW.repl.le} is somewhat analogous to
Proposition~\ref{spec.prop}, in that it reduces the study of
families of groupoids over $\langle I,W \rangle$ to families over
the CW-replacement $\Cof(\Posff \bb I)_{s(C)}$ of the category
$I$.

\subsection{Saturated CW-categories.}

Finally, let us introduce a class of CW-categories that is so close
to model categories that Proposition~\ref{mod.prop} works with the
same proof.

\begin{defn}\label{CW.sat.def}
A CW-category $\langle I,C,W \rangle$ is {\em saturated} if the
class $W$ is saturated, and for any map $i_0 \to i_1$ in $I$, there
exists a decomposition \eqref{f.dec} with $c \in C$, $w \in W$.
\end{defn}

\begin{exa}
In Example~\ref{CW.comma.exa}, $I'$ is saturated in the sense of
Definition~\ref{CW.sat.def} if so is $I$. In
Example~\ref{mod.CW.exa}, $I$ is saturated by definition. In
Example~\ref{CW.add.exa}, $C_\idot(\A)$ is saturated by virtue of
\eqref{cone.oct}.
\end{exa}

\begin{exa}\label{CW.sat.fun.exa}
For any saturated CW-category $I$ and thin Reedy category $J$, the
functor category $\Fun(J,I)$ with the C-structure of
Proposition~\ref{CW.prop} and $W$ consisting of maps that are
pointwise in $W$ is saturated (the same inductive procedure as in
the proof of Proposition~\ref{CW.prop} provides decompositions
\eqref{f.dec} with $w \in W$).
\end{exa}

\begin{exa}
For any weak CW-category $\langle I,C,W \rangle$, the CW-category
$\langle \Cof(I),C,a(W) \rangle$ of Lemma~\ref{W.ano.W.le} is
saturated.
\end{exa}

\begin{prop}\label{CW.sat.prop}
A CW-category $I$ saturated in the sense of
Definition~\ref{CW.sat.def} is regular in the sense of
Definition~\ref{reg.CW.def}.
\end{prop}

\proof{} As in the proof of Proposition~\ref{mod.prop}, consider the
subcategory $\overline{I} \subset \Cof(I)_C / I$ consisting of
triples $\langle i_0,i_1,\alpha\rangle$ such that $\alpha \in W$. Then it
again suffices to prove that $\tau:\overline{I} \to I$ is
$2$-connected, and to do this, it suffices to prove a version of
Corollary~\ref{mod.corr}: for any saturated CW-category $I$ and
object $i \in I$, $I^W/i = \overline{I}_i$ is homotopy filtered (and
then apply this result first to $I$, then to the appropriate
comma-fibers $i' \setminus I$). To prove this, assume given a finite
partially ordered set $J$ and a functor $\gamma:J \to I^W/i$, use
\eqref{f.dec} in $\Fun(J,I)$ to construct a cofibrant replacement
$\wt{\gamma} \to \gamma$, take the universal cone $\wt{\gamma}^>:J^>
\to I$, and then take a decomposition \eqref{f.dec} of
$\colim_J\wt{\gamma} = \wt{\gamma}(o) \to i$ to construct a cone
$J^> \to I^W/i$ of the same $\wt{\gamma}$ with values in $I^W/i$.
\endproof

A useful example of a saturated CW-category is the category
$C_\idot(\A)$ of chain complexes in an additive category $\A$, with
the CW-structure of Example~\ref{CW.add.exa}. This category is
usually not a model category (it need not even be finitely complete
and cocomplete). However, one can easily show that $\langle
C_\idot(\A),W \rangle$ is localizable, with the localization given
by the chain-homotopy category $\Ho(\A)$.

\begin{lemma}\label{add.loc.le}
The relative category $\langle C_\idot(\A),W \rangle$ for any
additive category $\A$ is localizable, the natural functor
$h:C_\idot(\A) \to \Ho(\A)$ factors through an equivalence
$h^W(C_\idot(\A)) \cong \Ho(\A)$, and we have $W = h^*(\Iso)$.
\end{lemma}

\proof[Sketch of a proof.] This is very standard, so we only sketch
a proof for the convenience of the reader. If a map $f:A_\idot \to
B_\idot$ is in $W$, then the contracting homotopy for $\Cone(f)$
provides, among other things, a map $g:B_\idot \to A_\idot$ and
chain homotopies between $\id_A$ and $g \circ f$, and $\id_B$ and $f
\circ g$. Therefore $h$ inverts maps in $W$, and then as in
Lemma~\ref{cat.loc.le}, it suffices to prove that any functor
$C_\idot(\A) \to \E$ that inverts maps in $W$ factors through
$h$. Again as in Lemma~\ref{cat.loc.le}, this immediately follows
from the facts that either of the two embeddings $A_\idot \to
\Cone(\delta_A)$ is in $W$.
\endproof

Lemma~\ref{add.loc.le} is of course completely analogous to
Lemma~\ref{cat.loc.le}, and we also have an analog of
Lemma~\ref{cat.epi.le} --- and more generally, of
Proposition~\ref{glue.prop}.

\begin{lemma}\label{add.glue.le}
For any additive functor $\gamma:\A \to \A'$ between additive
categories $\A$, $\A'$, the comma-category $\A' \setminus_\gamma \A$
is additive, and the functor
\begin{equation}\label{ho.A.epi}
\Ho(A' \setminus_\gamma A) \to \Ho(\A') \setminus_{\Ho(\gamma)} \Ho(\A)
\end{equation}
is an epivalence.
\end{lemma}

\proof{} The first claim is obvious, and it is also obvious that
\eqref{ho.A.epi} is conservative and essentially surjective, so it
suffices to check that it is full. Then as in the proof of
Proposition~\ref{glue.prop}, if we have two objects in the target of
\eqref{ho.A.epi} represented by $\langle A_\idot,A'_\idot,\alpha
\rangle$ and $\langle B_\idot,B'_\idot,\beta \rangle$ in
$C_\idot(\A' \setminus_\gamma \A)$, a map in
$\Ho(\A)/_{\Ho(\gamma)}\Ho(\A')$ between them gives rise to a
diagram
\begin{equation}\label{comma.add.sq}
\begin{CD}
  A'_\idot @>{f'}>> B'_\idot)\\
  @V{\alpha}VV @VV{\beta}V\\
  \gamma(A_\idot) @>{f'}>> \gamma(B_\idot)
\end{CD}
\end{equation}
that commutes up to a chain homotopy $\Cone(\delta_{A'}) \to
\gamma(B_\idot)$. We can now replace $B'_\idot$ with
the quotient $(\Cone(\delta_{A'}) \oplus B'_\idot)/A'_\idot$ to make
\eqref{comma.add.sq} commute on the nose; the quotient exists since
the embedding $A'_\idot \to \Cone(\delta_{A'})$ is termwise-split.
\endproof

More generally, assume given a thin Reedy category $J$, and consider
the functor category $\Fun(J,C_\idot(\A))$, with the saturated
CW-structure of Example~\ref{CW.sat.fun.exa}. Then we have
$C_\idot(\Fun(J,\A)) \cong \Fun(J,C_\idot(\A))$, and $\Fun(J,\A)$ is
also additive, but the two CW-structures are different: a complex
$A_\idot$ of functors $J \to \A$ can be contractible pointwise but
not contractible. Nevertheless, we have the following.

\begin{lemma}\label{A.cof.le}
For any cofibrant functor $A_\idot:J \to C_\idot(\A)$, the
representable functor $\Ho(\Fun(J,\A))(A_\idot,-)$ inverts maps in
the class $W$.
\end{lemma}

\proof{} For any two complexes $A_\idot,B_\idot \in \A$ in an
additive category $\A$, we have the complex of abelian groups
$\Hom^\hdot(A_\idot,B_\idot)$, maps of complexes $A_\idot \to
B_\idot$ correspond to cocycles in $\Hom^0(A_\idot,B_\idot)$, and
two maps are chain-homotopic iff their difference is a
coboundary. Thus it suffices to show that for any $A_\idot,B_\idot
\in C_\idot(\A^J)$ such that $A_\idot$ is cofibrant and $B_\idot(j)$
is contractible for any $j \in J$, the complex
$\Hom^\hdot(A_\idot,B_\idot)$ is acyclic. But an easy induction
shows that $A:J \to \A$ is cofibrant iff it is a sum of objects of
the form $\eps(j)_!A^j$, $A \in \A$, $j \in J$. Then any cofibrant
$A_\idot \in C_\idot(\A^J)$ has a termwise-split skeleton filtration
$\sk_n A_\idot$ that induces a convergent decreasing filtration on
$\Hom^\hdot(A_\idot,B_\idot)$, and its $n$-th associated graded
quotient is a products of complexes of the form
$$
\Hom^\hdot(\eps(j)_!A^j_\idot,B_\idot) \cong
\Hom^\hdot(A^j_\idot,B_\idot(j)), \qquad \deg j = n
$$
that are acyclic since $B_\idot(j)$ is contractible.
\endproof

\begin{corr}\label{add.A.corr}
For any thin Reedy category $J$ and additive category $\A$, the
relative category $\langle C_\idot(\Fun(J,\A)),W \rangle$ is
localizable, and if we denote by
$\Ho(\A,J)=h^W(C_\idot(\Fun(J,\A)))$ the localization, then the
embedding $\Cof(C_\idot(\Fun(J,\A))) \subset C_\idot(\Fun(J,\A))$
induces an equivalence
\begin{equation}\label{ho.A.J.eq}
\Ho(\Cof(\Fun(J,\A))) \cong \Ho(\A,J).
\end{equation}
Moreover, for any left-finite map $f:J' \to J$ between left-finite
partially ordered sets, the pullback functor $f^*:\Ho(\A,J) \to
\Ho(\A,J')$ admits a left-adjoint $f_!:\Ho(\A,J') \to \Ho(\A,J)$.
\end{corr}

\proof{} Lemma~\ref{A.cof.le} immediately implies that a cofibrant
complex $A_\idot$ in $\Fun(J,\A)$ is contractible iff $A_\idot(j)$
is contractible for any $i$ (to construct a contracting homotopy,
note that $\Cone(\id_A)$ is still cofibrant, and apply
Lemma~\ref{A.cof.le} to $\Cone(\id_A)$ and the map $A_\idot \to
0$). Thus a map $f:A_\idot \to A_\idot'$ between two cofibrant
complexes in $\Fun(J,\A)$ is in $W$ iff it is a chain-homotopy
equivalence. Then $\Ho(\Fun(J,\A)) \cong
h^W(\Cof(C_\idot(\Fun(J,\A))))$ by Lemma~\ref{add.loc.le}, and to
deduce \eqref{ho.A.J.eq}, it suffices to apply
Proposition~\ref{CW.sat.prop}. Moreover, Lemma~\ref{A.cof.le} then
shows that as soon as $A_\idot$ is cofibrant, we have
$\Ho(\A,J)(A_\idot,B_\idot) \cong \Ho(\Fun(J,\A))(A_\idot,B_\idot)$.
Now to construct $f_!$, it suffices to restrict our attention to
cofibrant complexes and maps in $C$, and apply Lemma~\ref{CW.J.le}.
\endproof

\begin{remark}
We use Proposition~\ref{CW.sat.prop} in Corollary~\ref{add.A.corr}
because we have it anyway; however, just as in Lemma~\ref{mod.le},
this is probably an overkill, and one can deduce
Corollary~\ref{add.A.corr} directly from Lemma~\ref{A.cof.le}.
\end{remark}

\section{Semiexact families.}\label{semi.sec}

\subsection{Semiexactness and additivity.}

We now turn to studying familes of groupoids over
CW-categories. Assume given a C-category $\langle I,C \rangle$ in
the sense of Definition~\ref{CW.def}, and a family of groupoids $\C$
over $I$.

\begin{defn}\label{semi.def}
The family $\C$ is {\em semiexact} resp.\ {\em exact} if for any two
maps $f_0:i \to i_0$, $f_1:i \to i_1$ in $I$ that lie in the class
$C$, the family $\C$ is semicartesian resp.\ cartesian along the
standard pushout square \eqref{st.CW.sq}, and cartesian if if $i=0$
is the initial object.
\end{defn}

We recall that explicitly, Definition~\ref{semi.def} means that for
any $f_0:i \to i_0$, $f_1:i \to i_1$ in $C$, the functor
\begin{equation}\label{exa.eq}
f_0^* \times f_1^*:\C_{i_0 \copr_i 1_1} \to \C_{i_0} \times_{\C_i}
\C_{i_1}
\end{equation}
is an epivalence, and an equivalence if $i=0$.

\begin{defn}\label{add.def}
The family $\C$ is {\em additive} if for any set $S$ and a family
of object $i_s \in I$ indexed by elements $s \in S$ that admits a
coproduct $\copr_s i_s \in I$, the natural functor
\begin{equation}\label{conti.eq}
C_{\copr_s i_s} \to \prod_s \C_{i_s}
\end{equation}
is an equivalence.
\end{defn}

\begin{remark}
If $I$ has all small coproducts, then Definition~\ref{add.def} is
equivalent to saying that $\C$ also has all small coproducts (and
the projection $\C \to I$ then automatically preserves them). For
finite coproducts, the claim is already included in
Definition~\ref{semi.def}.
\end{remark}

\begin{exa}\label{semi.II.exa}
For any family of groupoids $\C \to I$ over any category $I$, the
CW-functor $\Pos \bb \C \to \Pos \bb I$ of
Example~\ref{pos.II.CW.exa} is a family of groupoids, and it is
exact and additive with respect to the CW-structure of
Example~\ref{pos.I.CW.exa} on $\Pos \bb I$. The same holds for
$\Posf$ and $\Posff$.
\end{exa}

\begin{lemma}\label{exa.CW.le}
Assume that $\pi:\C \to I$ is an exact family of groupoids over a
CW-category $\langle I,C,W \rangle$ in the sense of
Definition~\ref{semi.def}, and as in Definition~\ref{lft.v.def}, let
$\pi^\flat(C)$ resp.\ $\pi^\flat(W)$ be the classes of cartesian lifting
of maps in $C$ resp.\ $W$. Then $\langle \C,\pi^\flat(C),\pi^\flat(W)
\rangle$ is a CW-category, $\pi$ is a CW-functor, and a family of
groupoids $\gamma:\C' \to \C$ is exact resp.\ semiexact if and only
if so is $\pi \circ \gamma:\C' \to I$.
\end{lemma}

\proof{} Clear.\endproof

Definition~\ref{semi.def} has one useful corollary that can be
proved right away. Fix a CW-category $I$ and a semiexact family of
groupoids $\C$ over $\langle I,W \rangle$. If we have a CW-category
$I'$ and a CW-functor $\gamma:I' \to I$, then $\gamma^*\C$ is a
semiexact family over $I'$. Consider the CW-category $\Cof(\Posff
\bb I)$ of Proposition~\ref{w.CW.prop}. As in \eqref{chi.pos},
denote
$$
\Cof(\Posff_\idot \bb I) = \Posff_\idot \times_{\Posff}
\Cof(\Posff_\idot \bb I) \subset \Posff_\idot \bb I,
$$
with the cofibration $\nu_\idot$ and the evaluation functor $\ev$
induced by those of Example~\ref{cat.bb.exa}. By \eqref{chi.pos},
the left Kan extension $\chi=(\nu_\idot)_!\ev$ exists, and then $\C$
induces two families of groupoids over $\Cof(\Posff \bb I)$, namely,
$\chi^*\C$ and $(\nu_\idot)_*\ev^*\C$, and a functor
\begin{equation}\label{coli.X.exa}
\chi^*\C \to (\nu_\idot)_*\ev^*\C
\end{equation}
over $\Cof(\Posff \bb I)$ provided by \eqref{coli.I}. Since $\chi$
is a CW-functor by Proposition~\ref{w.CW.prop}, the source
$\chi^*\C$ of the functor \eqref{coli.X.exa} is a semiexact family
over $\langle I,W \rangle$. Its target $(\nu_\idot)_*\ev^*\C$ need not
be semiexact, but by Proposition~\ref{w.CW.prop}, it is still a
family over $\langle I,W \rangle$. Explicitly, for any finite
partially ordered set $J$ equipped with a cofibrant functor
$i_\idot:J \to I$, the fiber of the functor \eqref{coli.X.exa} over
$\langle J,i_\idot \rangle \in \Cof(\Posff \bb I)$ is a functor
\begin{equation}\label{coli.exa}
\C_{\colim_J i_\idot} \to \Sec(J,i_\idot^*\C),
\end{equation}
and the functor \eqref{coli.exa} exists even if $J$ is not finite
(but $\colim_J i$ exists in $I$). If $J$ has a largest element, this
functor is tautologically an equivalence, and
Definition~\ref{semi.def} says that this is an epivalence for $J=\V$
of \eqref{V.eq} and any cofibrant $i_\idot:\V \to I$. For $J$ of
small dimension $\dim J$, we have the following general result.

\begin{lemma}\label{cube.le}
Assume given a finite partially ordered set $J$ and a cofibrant
functor $i_\idot:J \to I$. Moreover, assume that $\dim J \leq l$ for
some $l$. Then and for any semiexact family $\C$ over $I$, the
functor \eqref{coli.exa} is an equivalence resp.\ an epivalence
resp.\ essentially surjective if $l=0$ resp.\ $1$ resp.\ $2$.
\end{lemma}

\proof{} If $J=S$ or $J=S^<$ for a finite discrete set $S$, then the
statement follows from Definition~\ref{semi.def} by induction on the
cardinality $|S|$. If we have a left-reflexive subset $J' \subset
J$, then the embedding $J' \to J$ induces an equivalence between
both the sources and the targets of the functor \eqref{coli.exa},
and $\dim J' \leq \dim J$, so the claim for $J'$ implies the claim
for $J$. In particular, the claim holds if $J$ has a largest
element. In general, by induction on $|J|$, we may assume that we
have a finite set $S$ and a map $\phi:J \to J' = S^>$ such that the
claim holds for all the comma-sets $J / j$, $j \in J'$ (with any $l$
and $i_\idot$). But then $\phi_! i_\idot:J' \to I$ is cofibrant by
Lemma~\ref{CW.J.le}~\thetag{i}, and we have $\colim_J i_\idot \cong
\colim_{J'}\phi_!i_\idot$, $\Sec(J,i_\idot^*\C) \cong
\Sec(J',\phi_*i_\idot^*\C)$, so that \eqref{coli.exa} factors as
$$
\begin{CD}
\C_{\colim_{J'} \phi_!i_\idot} @>{a}>>
\Sec(J',\phi_!i_\idot^*\C) @>{b}>> \Sec(J',\phi_*i_\idot^*\C)
\end{CD}
$$
where we already know the claim for $a$, so it suffices to check it
for $b$. By inductive assumption and \eqref{kan.eq},
\eqref{2kan.eq}, the claim for $b$ follows from
Lemma~\ref{pos.epi.le}.
\endproof

\begin{corr}\label{semi.epi.corr}
Assume given a CW-category $I$ regular in the sense of
Definition~\ref{reg.CW.def}, and a functor $\gamma:\C' \to \C$ over
$I$ between two semiexact families of groupoids over $\langle
I,W\rangle$. Then if $\gamma$ is an epivalence, it is an
equivalence.
\end{corr}

\proof{} We have to prove that for any $i \in I$, the epivalence
$\gamma_i:\C'_i \to \C_i$ is an equivalence, and since $I$ is
regular, it suffices to consider cofibrant $i$. Fix such an $i \in
\Cof(I)$, construct a decomposition \eqref{cod.dec} with the
corresponding cofibrant functor $X':\II \to I$, as in
Example~\ref{thin.exa} and let $Y=v^*X':\VV \to I$, where $\VV$ is
as in \eqref{VV.eq} and $v$ is as in \eqref{VV.dia}. Then $Y$ is
cofibrant, so that $\langle \VV,Y \rangle$ is an object in
$\Cof(\Posff \bb I)$, and since $\dim\VV = 1$, Lemma~\ref{cube.le}
provides an epivalence
\begin{equation}\label{C.VV}
\C_{i'} \to \Sec(\VV,Y^*\C) \cong \C_i^\VV \cong \I(\C_i),
\end{equation}
where the last equivalence is as in Remark~\ref{ine.V.rem}, we
denote $i'=\colim_{\VV}Y$, and we observe that since $Y$
is equipped with a map $v^*a:Y=v^*X' \to v^*X$ in $W$ to the
constant functor $v^*X:\VV \to I$ with value $i$, $Y^*\C \to \VV$ is
a constant family with fiber $\C_i$. We also have the
counterpart of the epivalence \eqref{C.VV} for $\C'$, and since by
assumption, $\gamma_{i'}:\C'_{i'} \to \C_{i'}$ is an epivalence, so
is $\I(\gamma_i):\I(\C'_i) \to \I(\C_i)$. We are then done by
Lemma~\ref{ine.le}.
\endproof

\begin{exa}\label{efmid.exa}
The dimension estimate in Lemma~\ref{cube.le} is certainly not
sharp. For example, in the assumptions of Lemma~\ref{cube.le},
assume given a full left-reflexive subset $J' \subset J$. Then by
adjunction, $\colim_J i_\idot \cong \colim_{J'} i_\idot$ and
$\Sec(J,i_\idot^*\C) \cong \Sec(J',i_\idot^*\C)$, so that the
functors \eqref{coli.exa} for $J$ and $J'$ are the same. Thus the
conclusion of Lemma~\ref{cube.le} holds for $J$ as long as $\dim J'
\leq l$; the dimension of $J$ itself can be
arbitrary. Informally, one can say that such a $J$ has ``effective
dimension'' $\leq l$.
\end{exa}

\subsection{Framings.}\label{fram.subs}

Somewhat surprisingly, one can show that if a regular CW-category
$\langle I,C,W \rangle$ is good enough, an additive semiexact family
of groupoids $\C$ over $\langle I,W \rangle$ is automatically
constant along the class $a(W)$ of $W$-anodyne maps in $I$ in the
sense of Definition~\ref{ano.W.def}. This uses the same general idea
as in the proof of Corollary~\ref{semi.epi.corr}, and requires
putting an additional structure on $I$. Roughly speaking, we need to
choose a decomposition \eqref{cod.dec} for any object $i \in I$ in a
functorial way, or failing that, at least a contractible family of
such decompositions.

As in Example~\ref{thin.exa}, encode decompositions \eqref{cod.dec}
for objects $i \in I$ by functors $\wt{i}:\D_1 \to I$ such that
$l^*\wt{i}:\II \to I$ is cofibrant and factors through $I_{C \cap W}
\subset I$. Denote by $\Dec(I) \subset I^{\D_1}$ the full
subcategory spanned by all such $\wt{i}$, and let $C$ be the class
of maps $f$ in $\Dec(I)$ such that $l^*(f)$ is in $C$ with respect
to the CW-structure of Lemma~\ref{thin.le}. Denote by $\phi =
\ev_{[0]}:\Dec(I)_C \to \Cof(I)_C$ the evaluation functor sending a
diagram \eqref{cod.dec} to $i$, and let $\psi:\Dec(I)_C \to
\Cof(I)_C$ be the functor sending such a diagram $\wt{i}$ to
\begin{equation}\label{cod.psi.eq}
\colim_{\VV}v^*l^*\wt{t} \cong i' \copr_{i \copr i} i',
\end{equation}
where $v:\VV \to \II$ is as in \eqref{VV.dia}, and the colimit
exists since $v^*l^*\wt{i}$ is cofibrant.

\begin{defn}\label{fr.0.def}
A {\em framing} on a regular CW-category $I$ is a subcategory
$\wt{I} \subset \Dec(I)_C$ such that the induced functor
$\phi:\wt{I} \to \Cof(I)_C$ is locally $2$-connected in the sense of
Definition~\ref{loc.con.def}.
\end{defn}

If we are given a CW-category $I$ equipped with a framing $\wt{I}$
in the sense of Definition~\ref{fr.0.def}, then
Definition~\ref{loc.con.def} implies that any object $i \in \Cof(I)$
admits a lifting $\wt{i} \in \wt{I}_i$ to $\wt{I}$. Moreover, if we
have a map $f:i_0 \to i_1$ in $\Cof(I)_C$ and a lifting $\wt{i}_0
\in \wt{I}_{i_0}$ of its source $i_0$, then we can find a map
$\wt{f}:\wt{i}_0 \to \wt{i}_1$ such that $\phi(\wt{f}) = f$. If the
class $W$ is saturated, then it is easy to check that once $f$ is in
$C \cap W$, the induced maps $f':i_0' \to i_1'$,
$\psi(f):\psi(\wt{i}_0) \to \psi(\wt{i}_1)$ are also in $C \cap W$,
but in general, this need not be the case. Moreover, if we have a
standard pushout square
\begin{equation}\label{fr.sq}
\begin{CD}
i @>{f_0}>> i_0\\
@V{f_1}VV @VVV\\
1_1 @>>> i_{01} = i_0 \copr_i i_1
\end{CD}
\end{equation}
in $I$, and choose a lifting $\wt{i} \in \wt{I}_i$ and then
$\wt{f}_l:\wt{i} \to \wt{i}_l$, $l=0,1$, then the coproduct
$\wt{i}_{01} = \wt{i}_0 \copr_{\wt{i}} \wt{i}_1$ exists in
$I^{\D_1}$ since $I^{\II}$ is a CW-category, but nothing insures
that it lies in $\wt{I}$ (nor in fact in $\Dec(I) \subset
I^{\D_1}$). For our purposes, we only need framings that are free
from these two problems.

\begin{defn}\label{fr.def}
A framing $\wt{I}$ on a regular CW-category $I$ is {\em perfect} if
\begin{enumerate}
\item $\psi^*(C \cap W) \supset \phi^*(C \cap W)$, and
\item for any standard pushout square \eqref{fr.sq} in $I$, we can
  choose liftings $\wt{i} \in \wt{I}_i$, $\wt{f}_l:\wt{i} \to
  \wt{i}_l$, $l=0,1$ in such a way that $\wt{i}_{01} = \wt{i}_0
  \copr_{\wt{i}} \wt{i}_1$ is in $\wt{I} \subset I^{\D_1}$.
 \end{enumerate}
A regular CW-category $I$ is {\em strong} if it admits a perfect
framing.
\end{defn}

\begin{exa}\label{mod.cod.exa}
If $I$ is a model category with the CW-structure of
Example~\ref{mod.CW.exa}, consider $\Fun(\D_1,I)$ with the Reedy
model structure, and note that $\Dec(I)_C \subset
\Cof(\Fun(\D_1,I))_C$ is then the full subcategory spanned by
decompositions \eqref{cod.dec} with $w \in W$. Then $\wt{I} =
\Dec(I)$ is a perfect framing for $I$ (the fact that $\phi:\Dec(I)
\to \Cof(I)_C$ is locally $2$-connected easily follows from
Corollary~\ref{mod.corr}).
\end{exa}

\begin{exa}\label{mod.cod.nonexa}
More generally, one can show that for an arbitrary regular
$CW$-category $I$, the projection $\phi:\Dec(I)_C \to \Cof(I)_C$ is
locally $2$-connected, so that $\wt{I}=\Dec(I)$ is a framing for
$I$. If $W$ is saturated, then the framing is perfect, but in
general, thus need not be the case.
\end{exa}

\begin{exa}\label{adm.CW.exa}
Say that a CW-category $I$ is {\em distributive} if it has finite
products, the terminal object $1 \in I$ is cofibrant, $i \times -:I
\to I$ is a CW-functor for any cofibrant $i \in I$, and for any two
cofibrant $f,f' \in \Cof(\Ar(I))$, the product square $f \times
f':[1]^2 \to I$ is cocartesian in $I$. Then for any distributive
CW-category $I$, a choice of a decomposition \eqref{cod.dec} for the
object $1$ defines a framing for $I$. Namely, denote the
decomposition by
\begin{equation}\label{univ.cod}
\begin{CD}
1 \copr 1 @>{c}>> 1' @>>> 1,
\end{CD}
\end{equation}
with the corresponding functor $\wt{1}:\D_1 \to I$, and then let
$\wt{I}=\Cof(I)_C$, with the functor $- \times \wt{1}:\Cof(I)_C \to
\Fun(\D_1,I)$. Note that the evaluation functor $\phi:\wt{I} \to
\Cof(I)_C$ in this framing is simply the identity functor, while
$\psi$ sends $i$ to $i \times (1' \copr_{1 \copr 1} 1')$, so that
the framing is perfect. Here are two specific examples of such a
situation:
\begin{enumerate}
\item Take $I=\Pos$ with the standard CW-structure of
  Example~\ref{pos.CW.exa}. Then it is distributive, and the
  standard decomposition \eqref{univ.cod} is given by the cylinder
  construction (in particular, $1' = \V^o$). This is also inherited
  by very ample full subcategories $\I \subset \Pos$ and their
  extensions $\I^+$; all these CW-categories are distributive. The
  category $I=\Rel$ with the standard CW-structure of
  Example~\ref{rel.CW.exa} is distributive, with $1' =
  R(\V^o)$, and so is the unfolding $\I^\dm \subset \Rel$ of a very
  ample subcategory $\I \subset \Pos$, and its extension
  $\I^{+\dm}$.
\item Let $I = \Delta^o\Sets$ be the model category of simplicial
  sets, with the CW-structure of Example~\ref{mod.CW.exa}. Then it
  is also distributive, so to construct a perfect framing for $I$
  smaller that than of Example~\ref{mod.cod.exa}, it suffices to
  choose a decomposition \eqref{univ.cod}. For example, one can take
  $1' = \Delta_1$, the elementary $1$-simplex. Analogously, if
  $I=\Top$ is the model category of compactly generated topological
  spaces, then it is also distributive, and the unit interval $1' =
  [0,1]$ gives a decomposition \eqref{univ.cod} for $1 = \ppt$.
\end{enumerate}
\end{exa}

\begin{exa}\label{adm.CW.I.exa}
For any category $I$ and very ample full subcategory $\I \subset \Pos$,
the category $\I \bb I$ with its standard CW-structure of
Example~\ref{pos.I.CW.exa} is not distributive. However, we still
have the product functor $\I \times (\I \bb I) \to \I \bb I$ that is
distributive in the same sense as in Example~\ref{adm.CW.exa}, and
choosing a decomposition \eqref{univ.cod} for $\I$ provides a
perfect framing for $\I \bb I$. The same goes for the extension
$\I^+$.
\end{exa}

Assume given a regular CW-category $I$ with a framing $\wt{I}
\subset \Dec(I)$, and note that by \eqref{cod.psi.eq}, the functor
$\psi:\Dec(I)_C \to \Cof(I)_C$ factors as $\psi \cong \chi \circ
\wt{\psi}$, where $\chi:\Cof(\Posff \bb I)_C \to \Cof(I)_C$ is the
functor \eqref{chi.pos}, and $\wt{\psi}$ sends $\wt{i} \in \Dec(I)$
to $\langle \VV,v^*l^*\wt{i} \rangle \in \Cof(\Posff \bb I)$. Then
for any family of groupoids $\C$ over $I$, \eqref{coli.X.exa}
induces a functor
\begin{equation}\label{coli.psi.exa}
\wt{\psi}^*\chi^*\C \to \wt{\psi}^*\nu_*\ev^*\C,
\end{equation}
and if $\C$ is a family over $\langle I,W \rangle$, then as in the
proof of Corollary~\ref{semi.epi.corr}, the fiber
$(\wt{\psi}^*\nu_*\ev^*\C)_{\wt{i}}$ of its target over some $\wt{i}
\in \Dec(I)$ can be canonically identified with the inertia groupoid
$\I(\C_i)$, $i = \phi(\wt{i})$. Indeed, by definition, $\wt{i}:\D_1
\to I$ sends all maps in $\D_1$ to maps in the saturation $s(I,W)$
of the class $W$, and since $\D_1$ has a terminal object, we have a
canonical trivialization $\wt{i}^*\C \cong \C_i \times \D_1$ that
provides an identification
$$
(\wt{\psi}^*\nu_*\ev^*\C)_{\wt{i}} \cong \Sec(\VV,v^*l^*\wt{i}^*\C)
\cong \C_i^{\VV} \cong \I(\C_i),
$$
where the last identification is Remark~\ref{ine.V.rem}. Then
\eqref{coli.psi.exa} restricts to a functor
\begin{equation}\label{coli.psi.i.exa}
\C_{\psi(\wt{i})} \to \I(\C_i),
\end{equation}
and if $\C$ is semiexact, this is an epivalence by
Lemma~\ref{cube.le}.

\begin{lemma}\label{ano.C.le}
  Assume given a regular CW-category $I$ equipped with a perfect
  framing $\wt{I}$, a semiexact family of groupoids $\C$ over
  $\langle I,W \rangle$, and a standard pushout square \eqref{fr.sq}
  in $I$. Moreover, assume that both $\C$ and $\phi_*\psi^*\C$ are
  constant resp.\ fully faithful along $f_0$. Then $\C$ is also
  constant resp.\ fully faithful along the pushout $i_1 \to i_{01}$
  of the map $f_0$.
\end{lemma}

\proof{} Lift the square to a pushout square in $\wt{I}$ by
Definition~\ref{fr.def}~\thetag{ii}, and note that $\psi$ then sends
it to a standard pushout square in $I$. Since $\phi$ is
$2$-connected, $\psi^*\C \cong \phi^*\phi_*\psi^*\C$ is constant
resp.\ fully faithful along $\wt{f}_0$, so that $\C$ is constant
resp.\ fully faithful along $\psi(\wt{f}_0)$, and then
Definition~\ref{semi.def} implies that $\C_{\psi(\wt{i}_{01})} \to
\C_{\psi(\wt{i}_1)}$ is an epivalence resp.\ full. Together with the
epivalences \eqref{coli.psi.i.exa}, we have a commutative square
$$
\begin{CD}
  \C_{\psi(\wt{i}_{01})} @>>> \C_{\psi(\wt{i}_1)}\\
  @VVV @VVV\\
  \I(\C_{i_{01}}) @>>> \I(\C_{i_1})
\end{CD}
$$
where both vertical arrows are epivalences, and the top horizontal
arrow is an epivalence resp.\ full. Then the bottom horizontal arrow
is an epivalence resp.\ full by Lemma~\ref{epi.epi.le}, and we
are done by Lemma~\ref{ine.le}.
\endproof

\begin{corr}\label{ano.C.corr}
  Assume given a strong regular CW-category $\langle I,C,W \rangle$,
  and a semiexact family of groupoids $\C$ over $\langle I,W
  \rangle$. Then $\C$ is constant along the class $a(W) \supset W$
  of $W$-anodyne maps of Definition~\ref{ano.W.def}.
\end{corr}

\proof{} Choose a perfect framing $\wt{I}$ for $I$ in the sense of
Definition~\ref{fr.def}. Let $W_0(\C)$ be the class of maps $f$ in
$I$ such that $\C$ is constant along $f$, and then by induction,
define $W_n(\C)$, $n \geq 1$ as the class of maps $f$ in
$W_{n-1}(\C)$ such that $\psi(\wt{f}) \in W_{n-1}(\C)$ for any
$\wt{f}$ in $\wt{I}$ with $\phi(\wt{f}) = f$. Denote by $W(\C) =
\cap_n W_n(\C) \subset W_0(\C)$ the intersection of all these
classes. Then by Definition~\ref{fr.def}~\thetag{i}, we have $C \cap
W \subset W(\C)$, and by Lemma~\ref{ano.C.le}, standard pushouts of
maps in $W_n(\C)$ lie in $W_{n-1}(\C)$, so that $W(\C)$ is closed
under standard pushouts. Since all the classes $W_n(\C)$, hence also
$W(\C)$ are saturated by Lemma~\ref{sat.C.le}, we are done.
\endproof

\subsection{Mayer-Vietoris presentations.}\label{mv.subs}

It turns out that as an application of Corollary~\ref{ano.C.corr},
one can measure, to some extent, the failure of an additive semiexact
family of groupoids over a strong CW-category to be exact. Namely,
one can construct a natural presentation for an epivalence
\eqref{exa.eq} in the sense of Definition~\ref{pres.def}.

Assume given an arbitrary weak CW-category $I$, and an additive
semiexact family of groupoids $\C$ over $\langle I,W \rangle$. For
any map $q:i \to i'$ in $I$, denote by $\C_q \subset \C_i$ the
essential image of the transition functor $q^*$ (that is, the full
subcategory spanned by objects isomorphic to $q^*c$, $c \in
\C_{i'}$).

\begin{lemma}\label{retra.le}
Assume given two maps $s:i_0 \to i'_0$, $p:i_0 \to i_1$ in $C$, and
a map $q:i_1 \to i_0$ such that $q \circ p =\id$.  Let $i_1' = i_0'
\copr_{i_0} i_1$, with the natural maps $p' = \id \times p:i_0' \to
i_1'$, $s_1 = s \times \id:i_1 \to i_1'$, $q'=\id \times q:i_1' \to
i_0'$. Then for any semiexact family $\C$ over $\langle I,W
\rangle$, the functor \eqref{exa.eq} induces a functor
\begin{equation}\label{exa.rel}
\nu:\C_{q'} \to \C_{i_0'} \times_{\C_{i_0}} \C_q,
\end{equation}
and this functor is an equivalence of categories.
\end{lemma}

\proof{} Since $q' \circ s_1 = s \circ q$, $s_1^*$ sends the
essential image of $q'$ into the essential image of $q^*$, so that
\eqref{exa.eq} indeed induces a functor $\nu$ of
\eqref{exa.rel}. Since \eqref{exa.eq} is an epivalence, the functor
$\nu$ is full. Moreover, since $q \circ p = \id$, any object in its
target is isomorphic to $c \times q^*s^*c$ for some $c \in \C_{i'}$,
and since we have $c \times q^*s^*c \cong \nu(q^{'*}c)$, $\nu$ is
essentially surjective. Therefore it is an epivalence. Now choose a
decomposition \eqref{f.dec} of the map $q$, with $c:i_1 \to i_2$ in
$C$ and $c \circ p:i \to i_2$ in $W$, and note that $i_2' = i_2
\copr_{i_1} i_1' = i_2 \copr_{i_0} i_0'$ with the natural maps
$c':i_1' \to i_2'$, $w':i_2' \to i_0'$ gives a decomposition
\eqref{f.dec} of the map $q':i_1' \to i_0'$. Then
Definition~\ref{semi.def} applies to the coproduct $i_2' = i_1'
\copr_{i_1} i_2$, so that altogether, we obtain a commutative
diagram
$$
\begin{CD}
\C_{i_0'} \cong \C_{i_2'} @>{c^{'*}}>> \C_{q'} @>{p^{'*}}>> \C_{i_0'}\\
@V{s^*}VV @VV{s_1^*}V @VV{s^*}V\\
\C_{i_0} \cong \C_{i_2} @>{c^*}>> \C_q @>{p^*}>> \C_{i_0}
\end{CD}
$$
with semicartesian squares. Since $c^*:\C_{i_0} \to \C_q$ is
essentially surjective by definition, we are done by
Corollary~\ref{retra.corr}.
\endproof

Now assume that $I$ is a CW-category that is strong in the sense of
Definition~\ref{fr.def}, and fix a perfect framing $\wt{I}$ on
$I$. For any additive semiexact family of groupoids $\C$ over
$\langle I,W \rangle$, and any object $i \in I$, define the {\em
  loop groupoid} $\Omega(\C_i)$ by
the cartesian square
\begin{equation}\label{loop.sq}
\begin{CD}
\Omega(\C_i) @>>> (\phi_*\psi^*\C)_i\\
@VVV @VVV\\
\C_i @>>> \I(\C_i),
\end{CD}
\end{equation}
where the vertical arrow on the right is the epivalence
\eqref{coli.psi.i.exa}, and the bottom arrow is the tautological
fully faithful embedding. Therefore the top arrow is also fully
faithful. Moreover, for any decomposition \eqref{cod.dec} for $i$
corresponding to an object $\wt{i} \in \wt{I}$, we have
$(\phi_*\psi^*\C)_i \cong \C_{\psi(\wt{i})}$, and then $\Omega(\C_i)
\cong \C_q \subset \C_{\psi(\wt{i})}$, where $q:i'' = i' \copr_{i
  \copr i} i' \to i'$ is the codiagonal map.

\begin{lemma}\label{CW.pres.le}
For any standard pushout square \eqref{st.CW.sq}, the epivalence
\eqref{exa.eq} admits a presentation $\C'_{i_{01}}$ in the sense of
Definition~\ref{pres.def} such that
\begin{equation}\label{loop.pres.eq}
\C'_{i_{01}} \cong \C_{i_{01}} \times_{\C_i} \Omega(\C_i),
\end{equation}
where the projection $\C_{i_{01}} \to \C_i$ is induced by the
embedding $i \to i_{01}$.
\end{lemma}

\proof{} To simplify notation, let $\C_o=\C_i$ and $\C_l =
\C_{i_l}$, $l=0,1,01$, so that \eqref{exa.eq} becomes an epivalence
\begin{equation}\label{exa.C.eq}
  \C_{01} \to \C_0 \times_{\C_o} \C_1.
\end{equation}
Choose a decomposition \eqref{cod.dec} for $i$ corresponding to some
$\wt{i} \in \wt{I}$, and consider the associated anodyne map
$i'_{01} \to i_{01}$ of Example~\ref{W.exc.exa}. Then by
Corollary~\ref{ano.C.corr}, we also have $\C_{01} \cong
\C_{i'_{01}}$. The natural map $i_0 \copr i_1 \to i_{01}$ is in the
class $C$, so we can form the standard pushout $i''_{01} = i_{01}
\copr_{i_0 \copr i_1} i'_{01}$, and then \eqref{exa.eq} for this
standard pushout provides a presentation $\C'_{01} = \C_{i''_{01}}$
for the projection $\C_{01} \to \C_0 \times \C_1$. This projection
factors through \eqref{exa.C.eq} as
$$
\begin{CD}
  \C_{01} @>>> \C_0 \times_{\C_o} \C_1 @>>> \C_0 \times \C_1,
\end{CD}
$$
so by Example~\ref{pres.exa}, $\C'_{01}$ induces a presentation
$\delta^*\C'_{01}$ for \eqref{exa.C.eq}. Now observe that by
\eqref{ine.cart.sq}, we have
$$
(\C_0 \times_{\C_o} \C_1) \times_{\C_0 \times \C_1} (\C_0
\times_{\C_o} \C_1) \cong (\C_0 \times_{\C_o} \C_1) \times_{\C_o}
\I(\C_o),
$$
and the diagonal embedding $\delta$ in our application of
Example~\ref{pres.exa} is induced by the tautological embedding
$\C_o \to \I(\C_o)$. Therefore it is fully faithful, and we have
$\delta^*\C'_{01} \cong \C_{q'} \subset \C'_{01}$, where
$q':i''_{01} \to i'_{01}$ is the codiagonal map. But we also have
$i''_{01} \cong i'_{01} \copr_{i'} i''$, and $q'$ is
the pushout of the codiagonal map $q:\psi(\wt{i}) = i' \copr_{i \copr i}
i' \to i'$. Then \eqref{loop.pres.eq} follows from
Lemma~\ref{retra.le}.
\endproof

\begin{corr}\label{CW.pres.corr}
Let $\C$, $\C'$ be additive semiexact families of groupoids over a
strong CW-category $I$, and assume given a functor $\gamma:\C \to
\C'$ over $I$, and a standard pushout square \eqref{st.CW.sq} in $I$
such that $\gamma$ is an equivalence over $i$, $i_0$, $i_1$, and
$\phi_*\psi^*(\gamma)$ is an equivalence over $i$. Then $\gamma$ is
an epivalence over $i_{01}$.
\end{corr}

\proof{} Since $\gamma$ is an equivalence over $i$, $i_0$ and $i_1$,
$\gamma:\C_{i_{01}} \to \C'_{i_{01}}$ can be treated as a functor
over $\C_{i_0} \times_{\C_i} \C_{i_1} \cong \C'_{i_0} \times_{\C'_i}
\C'_{i_1}$. The presentation of Lemma~\ref{CW.pres.le} is functorial
with respect to $\C$, so that we actually obtain presentations for
the epivalences \eqref{exa.eq} for both $\C$ and $\C'$, and also a
presentation $\gamma'$ for $\gamma$. Moreover, since
$\phi_*\psi^*(\gamma)$ is also an equivalence over $i$, the functor
$\Omega(\gamma):\Omega(\C_i) \to \Omega(\C'_i)$ is an equivalence,
and then \eqref{loop.pres.eq} shows that the presentation $\gamma'$
is strict. We are then done by Lemma~\ref{pres.le}.
\endproof

\subsection{Families over model categories.}\label{mod.CW.subs}

Let us now step back to the situation of Section~\ref{mod.sec}: $I =
\langle I,C,W,F \rangle$ is a model category, with the CW-structure
of Example~\ref{mod.CW.exa}. Here we have an ample source of
examples provided by the following.

\begin{lemma}\label{repr.le}
For any object $i \in I$, the representable family $\Hh(i)$ of
Definition~\ref{repr.def} is additive and semiexact.
\end{lemma}

\proof{} By Proposition~\ref{Hh.prop}, we may assume that $i$ is
fibrant and cofibrant, and $\Hh(i)$ is given by \eqref{H.del} for
some Reedy fibrant $i^\Delta$. Then to prove that $\Hh(i)$ is
additive, it suffices to check that for any collection $X_s$, $s
\in S$ of fibrant simplicial sets with product $X = \prod_s X_s$,
the functor
$$
h^{\Tot}(\Delta X) \to \prod_s h^{\Tot}(\Delta X_s)
$$
is an equivalence. This immediately follows from
Lemma~\ref{kan.le}. To prove that $\Hh(i)$ is semiexact, again use
\eqref{H.del}, combined with Lemma~\ref{reedy.le} and
Corollary~\ref{semi.corr}.
\endproof

There are also two further properties of semiexact families that are
automatically satisfied in the model category settings. We start
with the following standard remark.

\begin{lemma}\label{Q.adj.le}
Assume given a regular cardinal $\kappa$, a $\kappa$-cocomplete
model category $\C$, and a Reedy category $J$ such that $|J| <
\kappa$ and the matching categories $M(j)$, $j \in J$ are finite
ordinals. Then $\colim_J:\Fun(J,I) \to I$ sends weak equivalences between
Reedy-cofibrant objects to weak equivalences.
\end{lemma}

\proof{} The assumption of the Lemma insures that the tautological
functor $I \to \Fun(J,I)$ sending $i \in I$ to the constant functor with
value $i$ sends maps in $F$ resp.\ $F \cap W$ to maps in $F$
resp.\ $F \cap W$ with respect to the Reedy model structure on
$\Fun(J,I)$. Then the colimit functor $\colim_J$ sends maps in $C$
resp.\ $C \cap W$ to maps in $C$ resp.\ $C \cap W$ --- that is, we
have a Quillen adjunction --- and by Ken Brown's
Lemma~\ref{CW.W.le}, $\colim_J$ indeed sends maps in $\Cof(\Fun(J,I))$ that
are in $W$ to maps in $W$.
\endproof

The first additional property of semiexact additive families over a
model category is the following one.

\begin{defn}\label{str.semi.def}
  A family of groupoids $\C$ over a $C$-category $I$ is {\em
    strongly semiexact} if \eqref{exa.eq} is an epivalence for any
  cocartesian square \eqref{st.CW.sq} in $\Cof(I)$ with $f_1 \in C$,
  and an equivalence if $i=0$.
\end{defn}

\begin{lemma}\label{semi.mod.le}
Assume given a family of groupoids $\C \to I$ over $\langle I,W
\rangle$ that is semiexact and additive in the sense of
Definition~\ref{semi.def} and Definition~\ref{add.def}. Then $\C$ is
strongly semiexact in the sense of Definition~\ref{str.semi.def}.
\end{lemma}

\proof{} By virtue of the factorization property of model
categories, it suffices to prove the claim for $f_0 \in C$ and for
$f_0 \in F \cap W$. In the first case, this is
Definition~\ref{semi.def}, so it suffices to show that if $f_0 \in
W$, then its pushout $f_0':i_1 \to i_0 \copr_i i_1$ with respect to
$f_1$ is also in $W$. This is a well-known fact from the theory of
model categories; let us reproduce the proof for the convenience of
the reader. Consider the partially ordered set $\V = \{0,1\}^<$ as a
Reedy category by taking $\deg 0 = 0$, $\deg o = 1$, $\deg 1 = 2$,
with the map $o \to 0$ resp.\ $o \to 1$ as the only non-trivial
matching resp.\ latching map. Then a functor $i_\idot:\V \to I$ is
cofibrant iff $i = i_\idot(o)$, $i_0 = i_\idot(0)$ are cofibrant and
the map $i \to i_1=i_\idot(1)$ is in $C$, and moreover, $\V$
satisfies the assumptions of Lemma~\ref{Q.adj.le}. It remains to let
$i'=i'_0=i$, $i'_1=i_1$, with the maps $\id:i' \to i'_0$, $f_1:i'_0
\to i_1$, consider the map $f_\idot:i'_\idot \to i_\idot$ given by
$\id$ at $o$ and $1$ and by $f_0$ at $0$, and note that $f'_0 =
\colim_{\V}f_\idot$.
\endproof

The second result concerns a version of the telescope
construction. Assume given a functor $i_\idot:\N \to I$ to some
category $I$, or explicitly, a collection of objects $i_n \in I$ and
maps $b_n:i_n \to i_{n+1}$, $n \geq 0$.

\begin{defn}\label{count.def}
  If $\colim_ni_n$ exists, then the map $b:i_0 \to \colim_ni_n$ is
  the {\em countable composition} of the maps $b_\idot$. A family
  $\C$ of groupoids over $I$ is {\em continuous} resp.\ {\em
    semicontinuous} along the countable composition $b$ if the
  functor \eqref{coli.exa} for $J=\N$ is an equivalence resp.\ an
  epivalence.
\end{defn}

To analyse continuity along countable compositions in CW-categories,
it is convenient to use the partially ordered set $\ZZ_\infty$ of
Example~\ref{Z.exa}, with the map \eqref{zeta.eq}. Note that $\dim
\ZZ_\infty = 1$, and $\ZZ_\infty$ is left-finite, thus a thin Reedy
category in the sense of Definition~\ref{thin.def}. The $0$-th
skeleton $\sk_0\ZZ_\infty$ is the subset $\{2l| l \geq 0\} \subset
\ZZ_\infty = \{l \geq 0\}$ of even integers, and it is abstractly
isomorphic to the set $\overline{\N}$ with discrete order. For any
CW-category $I$ and functor $Y:\N \to \Cof(I)_C$, let $X =
\zeta^*Y:\ZZ_\infty \to \Cof(I)_C$, and let $X':\ZZ_\infty \to
\Cof(I)_C$ and $a:X' \to X$ be its cofibrant replacement provided by
Lemma~\ref{thin.le}, with $J=\ZZ_\infty$ and
$J_0=\sk_0\ZZ_\infty$. Then $a$ is pointwise in the saturation
$s(W)$ of the class $W$, and for any family of groupoids $\C$ over
$\langle I,W \rangle$, we have functors
\begin{equation}\label{tel.C.eq}
\begin{CD}
\Sec(\N,Y^*\C) @>{\zeta^*}>> \Sec(\ZZ_\infty,X^*\C) @>{a^*}>>
\Sec(\ZZ_\infty,{X'}^*\C),
\end{CD}
\end{equation}
where $a^*$ is an equivalence. Moreover, by Lemma~\ref{Z.le},
$\zeta^*$ in \eqref{tel.C.eq} is an equivalence as well. Now define
a {\em telescope} $\Tel(Y)$ of the functor $Y$ by
\begin{equation}\label{tel.Y.eq}
  \Tel(Y) = \colim_{\ZZ_\infty}X' \in I.
\end{equation}
The telescope need not exist ($\ZZ_\infty$ is not finite, so
Lemma~\ref{CW.J.le} does not apply), and it need not be unique
(because $X'$ and $a:X' \to X$ are not unique). However, we have the
following general result.

\begin{lemma}\label{tel.I.le}
Assume that a CW-category $I$ has countable coproducts, and the
class $C$ is stable under countable coproducts. Then for any
functor $Y:\N \to \Cof(I)_C$ and cofibrant replacement $X$', $a:X'
\to X=\zeta^*Y$ of Lemma~\ref{thin.le}, the telescope
\eqref{tel.Y.eq} exists, and for any additive semiexact family of
groupoids $\C$ over $\langle I,W \rangle$, the functor
$$
\C_{\Tel(Y)} \to \Sec(\ZZ_\infty,{X'}^*\C) \cong \Sec(\N,Y^*\C)
$$
of \eqref{coli.exa} is an epivalence.
\end{lemma}

\proof{} Consider the discrete fibration $z:\ZZ_\infty \to \VV$ of
\eqref{ZZ.VV}. Each of its left comma-fibers $\ZZ_\infty /_z v$, $v
\in \VV$ is a countable coproduct of partially ordered sets
isomorphic to $\ppt$ or $\V^o$, and both $\ppt$ and $\V^o$ have a
largest element, thus are of effective dimension $0$ in the sense of
Example~\ref{efmid.exa}. Since $I$ has countable products,
\eqref{kan.eq} then shows that $z_!X':\VV \to I$ exists, and since
the family $\C$ is additive, we have an equivalence
$$
(z_!X')^*\C \cong z_*({X'}^*\C),
$$
so it remains to prove that \eqref{coli.exa} is an epivalence for
$J=\VV$ and the functor $z_!X':\VV \to I$. But by \eqref{kan.eq} and
Lemma~\ref{CW.J.le}, this functor is cofibrant and $\dim\VV=1$, so
we are done by Lemma~\ref{cube.le}.
\endproof

\begin{remark}\label{tel.I.rem}
Alternatively, instead of the map $z$, we could have proved
Lemma~\ref{tel.I.le} by using the map $x:\ZZ_\infty \to \V$ of
\eqref{x.z.eq} (that happens to be the composition of $z$ and the
cofibration $\VV \to \V$ of Remark~\ref{V.I.Z.rem}).
\end{remark}

Now assume that $I$ is a model category, and assume that it has
countable coproducts. Then is also has all countable colimits, thus
countable compositions, and $C$ is stable under countable
coproducts, so Lemma~\ref{tel.I.le} applies.

\begin{lemma}\label{tele.le}
Assume that a model category $I$ has countable coproducts. Then any
additive semiexact family of groupoids $\C$ over $\langle I,W
\rangle$ is semicontinuous along any countable compositions of maps
in $\Cof(I)_C$.
\end{lemma}

\proof{} Take a functor $Y:\N \to \Cof(I)_C$ and its telescope
$\Tel(Y)$ of \eqref{tel.Y.eq}. Then $a:X' \to X = \zeta^*Y$ induces
by adjunction a map $a_\dg:\zeta_!X' \to Y$, and by
Lemma~\ref{tel.I.le}, it suffices to check that the corresponding
morphism $\Tel(Y) \to \colim_{\N}Y$ is in $W$. But $Y$ is
Reedy-cofibrant, and by adjunction, so is $\zeta_!X'$. Therefore by
Lemma~\ref{Q.adj.le}, it suffices to check that $a_\dg:\zeta_!X' \to
Y$ is pointwise in $W$. By \eqref{kan.eq}, this amounts to checking
that the morphism $a_\dg:\Tel_m(Y) = \colim_{\ZZ_{2m}}X' \to X(2m)
\cong Y(m)$ is in $W$ for any integer $m \geq 0$. Then the
square \eqref{zz.m.sq} provides an isomorphism
\begin{equation}\label{tel.X.eq}
\Tel_m(Y) \cong X'(1) \copr_{Y(1)} \Tel_{m-1}(q^*Y),
\end{equation}
and by the same induction as in Lemma~\ref{Z.m.le}, we are reduced
to $m=1$, where $\Tel_1(Y) = X'(1)$, and $a_\dg:\Tel_1(Y) =
X'(1) \to Y(1) = X(2) = X(1)$ is the map $a$ on the nose.
\endproof

\subsection{Representability criterion.}\label{repr.cr.subs}

To complement Lemma~\ref{repr.le}, it will be useful to have a way
to check whether an object $i \in I$ in a category $I$ represents a
given family $\C \to I$.

\begin{defn}\label{lift.def}
Let $\C$ be a family of groupoids over a C-category $I$, and assume
given an object $c \in \C$ and a map $f:i \to i'$ in the class
$C$. Then $c$ is {\em liftable} with respect to $f$ if for any $c'
\in \C_{i'}$, any map $a:f^*c' \to c$ factors through the natural
map $\wt{f}:f^*c' \to c'$ by a map $b:c' \to c$.
\end{defn}

For any $c \in \C$, we denote by $C(c) \subset C$ the class of all
maps $f$ in $\Cof(I)_C$ such that $c$ is liftable with respect to
$f$. The class $C(c)$ is obviously closed under compositions and
contains all identity maps.

\begin{lemma}\label{lift.le}
Assume given a family of groupoids $\C$ over a C-category $I$, and
an object $c \in \C$.
\begin{enumerate}
\item If the family $\C \to I$ is strongly semiexact in the sense of
  Definition~\ref{str.semi.def}, then the class $C(c)$ is closed under
  pushouts.
\item If the family $\C \to I$ is additive in the sense of
  Definition~\ref{add.def}, then the class $C(c)$ is closed under
  coproducts.
\item If the family $\C \to I$ is semicontinuous in the sense of
  Definition~\ref{count.def} over a countable composition $b$ of
  maps $b_n \in C(c)$, then $b \in C(c)$.
\end{enumerate}
\end{lemma}

\proof{} Denote by $\pi:\C \to I$ the projection. Note that since
$\pi$ is a fibration, giving a map $b:c' \to c$ in
Definition~\ref{lift.def} is equivalent to giving the map $\pi(b):i'
\to \pi(c)$ in $I$ such that $\pi(a) = \pi(b) \circ f$, and
moreover, an isomorphism $\alpha:c' \cong \pi(b)^*c$. In all the
cases \thetag{i}, \thetag{ii}, \thetag{iii}, the map $\pi(b)$ is
immediately provided by the universal properties of colimits, and
the issue is the isomorphism $\alpha$. In \thetag{i}
resp.\ \thetag{ii}, its existence follows from \eqref{exa.eq}
resp.\ \eqref{conti.eq}. For \thetag{iii}, use
Definition~\ref{count.def} directly.
\endproof

\begin{prop}\label{lift.prop}
Assume given a model category $I$, a cofibrant object $i \in I$, an
additive semiexact family of groupoids $\C$ over $\langle I,W
\rangle$, and an object $c \in \C_i \subset \C$. Then \thetag{i}
$C(c)=C$ if and only if \thetag{ii} the object $i \in I$ is fibrant,
and the comparison functor $\Hh(c):\Hh(i) \to \C$ of
Pro\-po\-si\-tion~\ref{repr.prop}~\thetag{iii} is an equivalence
over $\Cof(I)_C$.
\end{prop}

In order to prove this, we need to modify slightly the description
of the family $\Hh(i)$ given in Proposition~\ref{Hh.prop}. Assume
given two objects $i',i \in I$, $i$ fibrant and $i'$ cofibrant. As
in Subsection~\ref{exp.subs}, choose a Reedy fibrant replacement
$i^\Delta:\Delta^o \to I$ of the constant simplicial object in $I$
with value $i$, and assume that $i \cong i^\Delta([0])$. Moreover,
the category of cosimplicial objects in $I$ also has a Reedy model
structure, so that we can choose a cofibrant replacement
$i'_\Delta:\Delta \to I$ of the constant cosimplicial object $i'$,
with $i'_\Delta([0]) \cong i'$. We then have a simplicial set
$\Hom^\Delta(i'_\Delta,i):\Delta^o \to \Sets$ given by $[n] \mapsto
\Hom(i'_\Delta([n]),i)$, and its category of simplices $\Delta
\Hom^\Delta(i'_\Delta,i)$.

\begin{lemma}\label{del.del.le}
In the assumptions above, there is a natural equivalence of
groupoids $h^{\Tot}(\Delta \Hom^\Delta(i',i^\Delta)) \cong
h^{\Tot}(\Delta \Hom^\Delta(i'_\Delta,i))$.
\end{lemma}

\proof{} Consider the bisimplicial set
$\Hom^\Delta(i'_\Delta,i^\Delta):\Delta^o \times \Delta^o \to
\Sets$ given by $[n'] \times [n] \mapsto
\Hom(i'_\Delta([n']),i^\Delta([n])$, and denote by
$\Delta^2\Hom^\Delta(i'_\Delta,i^\Delta) \to \Delta \times \Delta$
its category of elements, with objects $\langle [n']
\times [n],f \rangle$, $f$ a map from $i'_\Delta([n'])$ to
$i^\Delta([n])$. Denote by
$\tau',\tau:\Delta^2\Hom^\Delta(i'_\Delta,i^\Delta) \to \Delta$ the
fibrations sending $\langle [n'],[n],f \rangle$ to $[n']$
resp.\ $[n]$. Then by Proposition~\ref{Hh.prop}, the family of
groupoids
$$
\tau'_!\Delta^2\Hom^\Delta(i'_\Delta,i^\Delta) \to \Delta
$$
has fibers $\Hh(i)_{i'_\Delta([n])}$, $[n] \in \Delta$. Since
$i'_\Delta$ is a replacement of a constant functor, all the
transition functors of this family are therefore equivalences of
groupoids, and since $\Delta$ has the terminal object $[0]$, the
embedding
\begin{equation}\label{del.del2}
\Delta \Hom^\Delta(i',i^\Delta) \cong
\Delta^2\Hom^\Delta(i'_\Delta,i^\Delta)_{[0]} \subset
\Delta^2\Hom^\Delta(i'_\Delta,i^\Delta)
\end{equation}
induces an equivalence of total localizations. To finish the proof,
apply the same argument to the opposite model category $I^o$ and the
projection $\tau$, and deduce that the embedding
\begin{equation}\label{codel.del2}
\Delta \Hom^\Delta(i'_\Delta,i) \subset
\Delta^2\Hom^\Delta(i'_\Delta,i^\Delta)
\end{equation}
also induces an equivalence of total localizations.
\endproof

Combining Lemma~\ref{del.del.le} and \eqref{H.del}, we see that for
any $i$, $i^\Delta$, $i'$, $i'_\Delta$ as above, we have a
canonical equivalence
\begin{equation}\label{H.del.bis}
\Hh(i)_{i'} \cong h^{\Tot}(\Delta \Hom^\Delta(i'_\Delta,i)).
\end{equation}
By Lemma~\ref{reedy.le} applied to the opposite model category
$I^o$, the simplicial set $\Hom^\Delta(i'_\Delta,i)$ is fibrant, so
that the right-hand side of \eqref{H.del.bis} can be described by
Lemma~\ref{kan.le}. Explicitly, let $\wt{i} = i'_\Delta([1])$, and
consider the corresponding factorization
$$
\begin{CD}
i' \copr i' @>{p}>> \wt{i} @>{q}>> i
\end{CD}
$$
of the codiagonal map $i \copr i \to i$, with $p = p_0 \copr p_1$ in
$C$ and $q$ in $W$ ($p_0$ and $p_1$ are induced by the maps $s,t:[0]
\to [1]$, and $q$ is induced by the projection $[1] \to [0]$). Then
objects of the groupoid $\Hh(i)_{i'}$ are represented by maps $f:i'
\to i$, and morphisms from $f_0$ to $f_1$ are represented by maps
$g:\wt{i} \to i$ such that $f_l = f \circ p_l$, $l=0,1$.

\begin{lemma}\label{hh.exp.le}
Assume given a family of groupoids $\C$ over $\langle I,W \rangle$
and an object $c \in \C_i$. Then in terms of the identification
\eqref{H.del.bis}, the functor $\Hh(c)_{i'}:\Hh(i)_{i'} \to \C_{i'}$
of Proposition~\ref{repr.prop} sends an object represented by a map
$f:i' \to i$ to $f^*c \in \C_{i'}$, and a morphism represented by a
map $g:\wt{i} \to i$ to the composition
$$
\begin{CD}
(g \circ p_0)^* c @>{p_0^*(\alpha^{-1})}>> p_0^*q^*q_!g^*c \cong
  q_!g^*c \cong p_1^*q^*q_!g^*c @>{p_1^*(\alpha)}>> (g \circ p_1)^*c,
\end{CD}
$$
where $\alpha:g^*c \to q^*q_!g^*c$ is the adjunction map, and the
isomorphisms are induced by the identities $q \circ p_0 = q \circ
p_1 = \id$.
\end{lemma}

\proof{} For any $[n] \in \Delta$, denote the maps corresponding to
the projection $[n] \to [0]$ by $e:i = i^\Delta([0]) \to
i^\Delta([n])$ and $e':i'_\Delta([n]) \to i'_\Delta([0]) = i'$. Then
for any $[n],[n'] \in \Delta$ and $f:i'_\Delta([n']) \to
i^\Delta([n])$, we have a functor
$$
\begin{CD}
\C_i @>{e_!}>> \C_{i^\Delta([n])} @>{f^*}>> \C_{i'_\Delta([n'])}
@>{e^{'*}}>> \C_{i'}.
\end{CD}
$$
These functors are compatible with maps in
$\Delta^2\Hom^\Delta(i'_\Delta,i^\Delta)$, thus glue together to a
functor $\C_i \times \Delta^2\Hom^\Delta(i'_\Delta,i^\Delta) \to
\C_{i'}$.  Restricting this to $c \in \C_i$, we obtain a functor
$$
\Hh'(c):\Delta^2\Hom^\Delta(i'_\Delta,i^\Delta) \to \C_{i'}.
$$
When we compose $\Hh(c)'$ with the embedding \eqref{del.del2}, we
obtain the functor $\Hh(c)_{i'}$ of Proposition~\ref{Hh.prop}, and
since \eqref{del.del2} induces an equivalence of total
localizations, we have $\Hh'(c) \cong \Hh(c)_{i'}$. To finish the
proof, it remains to compose $\Hh'(c)$ with the embedding
\eqref{codel.del2} and use Lemma~\ref{kan.le}.
\endproof

\proof[Proof of Proposition~\ref{lift.prop}.] The fact that
\thetag{ii} implies \thetag{i} is obvious. Conversely, assume that
$C = C(c)$. Then for any morphism $w:i_0 \to i_1$ in $C \cap W$, the
transition functor $w^*:\C_{i_1} \to \C_{i_0}$ is an equivalence, so
that for any map $f:i_0 \to i$, we have $f^*c \cong w^*c'$ for some
$c' \in \C_{i_1}$. Then Definition~\ref{lift.def} immediately
implies that $f$ factors through $w$, so that $i$ is fibrant.

It remains to prove that $\Hh(c)_{i'}:\Hh(i)_{i'} \to \C_{i'}$ is an
equivalence for any $i' \in I$, and to do this, we use
Lemma~\ref{hh.exp.le}. To prove that $\Hh(c)_{i'}$ is essentially
surjective, we need to check that any object $c' \in \C_{i'}$ is of
the form $f^*c$ for some $f:i' \to i$; this follows from
Definition~\ref{lift.def} applied the map $0 \to i'$ (note that
\eqref{conti.eq} for $S = \emptyset$ means that $\C_0 \cong
\ppt$). To prove that $\Hh(c)_{i'}$ is full, assume given two maps
$f_0,f_1:i' \to i$ and a morphism $\phi:f_0^*c \to f_1^*c$, let $c'
= q^*f_0^*c \in \C_{\wt{i}}$, identify $p^*c' \cong p_0^*c' \copr
p_1^*c' \cong f_0^*c \copr f_0^*c$, and apply Definition~\ref{lift.def}
to the composition of the canonical map $f_0^*c \copr f_1^*c \to c$
and the map $\id \copr \phi$. This gives a map $g:\wt{i} \to i$, and
by Lemma~\ref{hh.exp.le}, the functor $\Hh(c)_{i'}$ sends the
morphism it represents to $\phi$.

Finally, to prove that $\Hh(c)_{i'}$ is faithful, let $J = L([2])$
be the latching category of the object $[2] \in \Delta$. Explicitly,
$J$ is the partially ordered set of non-empty proper subsets in $[2]
= \{0,1,2\}$, so it has chain dimension $1$. We denote the objects
in $J$ by $0$, $1$, $2$, $01$, $02$, $12$, and we let $i_\idot:J \to
I$ be the composition of the functor $i'_\Delta$ with the forgetful
functor $J \to \Delta$ that sends $[n] \to [2]$ to $[n]$. Then since
since $i'_\Delta$ is Reedy cofibrant, $i_\idot \in \Fun(J,I)$ is cofibrant
with respect to the projective model structure, and moreover, the
natural map $l:\colim i_\idot \to i'_\Delta([2])$ is in $C$ (it is
in fact one of the latching maps \eqref{lm.eq}). Composing it with
the map $e':i'_\Delta([2]) \to i'$ induced by the projection $[2]
\to [0]$, we obtain a map $a:\colim i_\idot \to i'$.

Now assume given two maps $f_0,f_1:i' \to i$, and two maps
$g_0,g_1:\wt{i} \to i$ such that $g_0 \circ p_0 = g_1 \circ p_0$,
$g_0 \circ p_1 = g_1 \circ p_1$, and $\Hh(c)_{i'}(g_0) =
\Hh(c)_{i'}(g_1)$. Then we have a map $i_\idot \to i$ from $i_\idot$
to the constant functor with value $i$ that is equal to $f_0$,
$f_0$, $f_1$, $f_0 \circ q$, $g_0$, $g_1$ on $i_0$, $i_1$, $i_2$,
$i_{01}$, $i_{02}$, $i_{12}$, and by adjunction, it induces a map
$b:\colim i_\idot \to i$. Moreover, by Lemma~\ref{hh.exp.le}, the
objects $b^*c$ and $a^*f_0^*c$ become isomorphic after applying the
comparison functor \eqref{coli.exa}, and then by
Lemma~\ref{cube.le}, they were isomorphic to begin with. Thus $b^*c
\cong l^*c'$ for $c' = (f_0 \circ e')^*c \in
\C_{i'_\Delta([2])}$. Applying Definition~\ref{lift.def} to the map
$l$, we conclude that $b$ factors through $l$, and then by
Lemma~\ref{kan.le}, $g_0$ and $g_1$ represent the same map already
in $\Hh(i)_{i'}$.
\endproof

\section{Representability theorems.}\label{repr.sec}

\subsection{Relative CW complexes.}\label{top.subs}

Consider the category $\Top$ of compactly generated topological
spaces, and equip it with Quillen's model structure of
\cite{qui.ho}: a map $f:X \to Y$ is in $W$ iff $f:\pi_n(X,x) \to
\pi_n(Y,f(x))$ is an isomorphism for any $x \in X$, $n \geq 0$, and
a map $f:X \to Y$ is in $F$ iff it has the right lifting property
with respect to maps $D^{n-1} \times \{0\} \subset D^n$, $n \geq 1$,
where $D^n = [0,1]^n$ is the $n$-th self-product of the unit
interval $[0,1]$. The model category $\Top$ is cofibrantly generated
in the sense of \cite{hovey}, with generating cofibrations given by
the embeddings $S^{n-1} \to D^n$, $n \geq 0$, where $S^{n-1}$ is the
$(n-1)$-dimensional sphere, embedded as the boundary of the
$n$-dimensional disc $D^n$, and $S^{-1}$ is taken to be the empty
set. As a corollary of this, for any essentially small category $I$,
the category $\Fun(I,\Top)$ admits the projective model
structure. Cofibrant objects in $\Top$ are exactly CW complexes, so
it is natural to introduce the following.

\begin{defn}
For any essentially small category $I$, an {\em $I$-CW complex} is a
cofibrant object $X \in \Fun(I,\Top)$.
\end{defn}

We will denote by $I\CW = \Cof(\Fun(I,\Top))_C$ the category of
$I$-CW complexes and maps between them that are in the class
$C$. For any $i \in I$, the evaluation functor $\ev_i:\Fun(I,\Top)
\to \Top$ admits an obvious left-adjoint functor $i_!:\Top \to
\Fun(I,\Top)$, $i_!(X)(i') = X \times I(i',i)$, and by adjunction,
$i_!$ sends cofibrant objects to cofibrant objects. For any $n \geq
0$, we will denote $D^n_i = i_!D^n \in I\CW$ and $S^n_i=i_!S^n \in
I\CW$. By adjunction, discs $D^n_i$ and spheres $S^n_i$ are compact
objects in $\Fun(I,\Top)$. Every object in $\Top$, hence also in
$\Fun(I,\Top)$ is fibrant, and for any $X \in \Fun(I,\Top)$, a
convenient fibrant replacement of the diagonal $X \to X \times X$ is
given by
\begin{equation}\label{P.eq}
\begin{CD}
X @>>> \Maps(D^1,X) @>>> X \times X,
\end{CD}
\end{equation}
where for any $i \in I$, $\Maps(D^1,X)(i) = \Maps(D^1,X(i))$ is the
space of maps from the unit interval $D^1$ to $X(i)$ with the
compact-open topology, and the projection to $X \times X$ is given
by evaluation at $0,1 \in D^1$. Therefore $h^W(\Fun(I,\Top))$
coincides with the naive homotopy category: morphisms from $Y$ to
$X$ are represented by maps $f:Y \to X$ in $\Fun(I,\Top)$, and two
maps $f_0$, $f_1$ represent the same morphism if and only if there
exists a map $g:Y \times D^1 \to X$ that restricts to $f_0$ on $Y
\times \{0\} \subset Y \times D^1$ and to $f_1$ on $Y \times \{1\}
\subset Y \times D^1$. For any $X \in I\CW$, we have the
representable family of groupoids $\Hh(X)$ over $I\CW$ of
Definition~\ref{repr.def}, and again, it can be explicitly described
by Proposition~\ref{Hh.prop} using \eqref{P.eq}: objects of the
fiber $\Hh(X)_Y$ for some $I$-CW complex $Y$ are represented by maps
$f:Y \to X$, and morphisms between objects represented by $f_0,f_1:Y
\to X$ are homotopy classes of maps $g:Y \times D^1 \to X$ that
restrict to $f_0$ resp.\ $f_1$ on $Y \times \{0\}$ resp.\ $Y \times
\{1\}$. If $I = \ppt$ is the point, then $\Hh(X)_Y$ is the
fundamental groupoid of the space $\Maps(Y,X)$.

\begin{theorem}\label{coho.thm}
Assume given an essentially small category $I$ and a small family
$\C$ of group\-oids over $\langle I\CW, W \rangle$ such that $\C$ is
semiexact in the sense of Definition~\ref{semi.def} and additive in
the sense of Definition~\ref{add.def}. Then there exists an object
$X \in I\CW$ and an equivalence $\Hh(X) \cong \C$, and $X$ is unique
up to a weak equivalence.
\end{theorem}

This result is an analog of the ``cohomological'' version of the
Brown Theorem \cite{brown}, and to prove it, we follow the same
plan, see e.g.\ \cite{switz}.

\begin{defn}\label{univ.def}
For any family $\C$ as in Theorem~\ref{coho.thm}, any $I$-CW complex
$X \in I\CW$, and any integer $n \geq 0$, an object $x \in \C_X$ is
{\em $n$-universal} if it is liftable in the sense of
Definition~\ref{lift.def} with respect to the natural embedding
$S^{m-1}_i \to D^m_i$ for any $i \in I$ and integer $m \leq n$. An
object $x$ is {\em universal} if it is $n$-universal for any $n \geq
0$.
\end{defn}

\begin{lemma}\label{univ.exi.le}
\begin{enumerate}
\item Assume given an integer $n \geq 0$, some $X \in I\CW$, and an
  object $x \in \C_X$, and assume that either $n=0$, or the object
  $x$ is $(n-1)$-universal. Then there there exists a map $e:X \to
  X'$ in $I\CW$ and an $n$-universal object $x' \in \C_{X'}$ such
  that $e^*(x') \cong x$.
\item Assume given $X \in I\CW$ and an object $x \in \C_X$. Then
  there exists a map $e:X \to X'$ in $I\CW$ and a universal $x' \in
  \C_{X'}$ such that $e^*x' \cong x$.
\end{enumerate}
\end{lemma}

\proof{} For \thetag{i}, assume given $X$, $x \in \C_X$ and $n$, and
for any $i \in I$, let $A_i$ be the set of triples $\langle f,u,\phi
\rangle$ of a map $f:S^{n-1}_i \to X$, an object $u \in \C_{D^n_i}$,
and an isomorphism $f^*x \cong u|_{S^{n-1}_i}$. Since $\C$ has small
fibers, $A_i$ is indeed a set, and we can form the disjoint unions
$$
Z = \coprod_{i \in I} S^{n-1}_i \times A_i \subset Y = \coprod_{i
  \in I} D^n_i \times A_i.
$$
We then have a natural map $f:Z \to X$, and moreover, since $\C$ is
additive, we have an object $u \in \C_Y$ and an isomorphism $f^*x
\cong u|_Z$. Since $\C$ is also semiexact, we can then set $X' = Y
\copr_Z X$, and choose an object $x' \in \C_{X'}$ that restricts to
$u$ resp.\ $x$ on $Y$ resp.\ $X$.

To check that $x'$ is $n$-universal, we need to prove that for any
$i \in I$, integer $m \leq n$ and object $u \in \C_{\D^m_i}$, any
map $f_{m-1}:S^{m-1}_i \to X'$ such that $f_{m-1}^*x' \cong
u|_{S^{m-1}_i}$ extends to a map $f_m:D^m_i \to X'$ such that $u
\cong f_m^*x'$. But this claim only depends on the homotopy class of
$f_{m-1}$, and by adjunction and standard cellular approximation
theorem, we may therefore assume that $f_{m-1}$ factors through $X
\subset X'$. Then the claim follows by assumption if $m < n$ and by
construction if $m = n$.

For \thetag{ii}, let $X_{-1}=X$, $x_{-1} = x$, use \thetag{i} and
induction to construct a collection of $I$-complexes $X_n$, objects
$x_n \in \C_{X_n}$, maps $X_{n-1} \to X_n$ in $I\CW$, and
isomorphisms $x_n |_{X_{n-1}} \cong x_{n-1}$ such that $x_n$ is
$n$-universal for any $n \geq 0$, let $X' = \colim X_n$, and use
Lemma~\ref{tele.le} to obtain $x \in \C_{X'}$ that restricts to
$x_n$ on each $X_n$. Then since the spheres are compact, for any $i
\in I$ and $m \geq 0$, any map $f:S^{m-1}_i \to X'$ factors through
$X_n$ for some $n$, and increasing $n$ if necessary, we may assume
that $x_n \cong x'|_{X_n}$ is $m$-universal. Therefore $x'$ is
liftable with respect to the embedding $S^{m-1}_i \to D^m_i$, thus
universal.
\endproof

\begin{lemma}\label{univ.uni.le}
Assume given two $I$-CW complexes $X,X' \in I\CW$, universal objects
$x \in \C_X$, $x' \in \C_{X'}$, and a map $e:X \to X'$ such that
$e^*(x') \cong x$. Then $e$ is a weak equivalence.
\end{lemma}

\proof{} By the definition of the model structure on $\Fun(I,\Top)$, it
suffices to prove that $e$ becomes a weak equivalence after
evaluation at any object $i \in I$, so that we may assume that $I =
\ppt$ and $X$, $X'$ are usual CW complexes. For any $n \geq 0$, the
liftability of $x$ with respect to the embedding $S^{n-1} \subset
D^n$ shows that any map $f:S^{n-1} \to X$ such that $e \circ f$
extends to $D^n$ itself extends to $D^n$; this shows that
$\pi_n(X,o) \to \pi_n(X',e(o))$ is injective for any choice of a
distinguished point $o \in X$. To prove that $\pi_0(X) \to
\pi_0(X')$ is surjective, take a map $f:\ppt \to X'$, and note that
by the universality of $x$, $f^*x' \cong g^*x$ for some map $g:\ppt
\to X$. Then by the universality of $x'$, the map $f \copr (e \circ
g):S^0 = \ppt \copr \ppt \to X'$ extends to $D^1$, so that $f$ is
homotopic to $e \circ g$. Finally, to prove that $\pi_n(X,o) \to
\pi_n(X',e(o))$ is surjective for some $n \geq 1$ and $o \in X$,
note that by Lemma~\ref{lift.le}~\thetag{i}, $x$ and $x'$ are also
liftable with respect to the embedding $S^n \copr S^n \to S^n \times
D^1$ and the embedding $\ppt \to S^n$ onto a distinguished point $o
\in S^n$. Then for any pointed map $f:\langle S^n,o \rangle \to
\langle X',e(o) \rangle$, there exists a pointed map $g:\langle
S^n,o \rangle \to \langle X,o \rangle$ such that $f^*x' \cong g^*x$,
and then the map $f \copr e \circ g:S^n \copr S^n \to X'$ extends to
$S^n \times D^1$, so that $f$ is homotopic to $e \circ g$.
\endproof

\proof[Proof of Theorem~\ref{coho.thm}.] Since $\C$ is additive,
$\C_{\emptyset}$ consists of a single object, and using
Lemma~\ref{univ.exi.le}, we can start from this object and construct
a universal object $x \in \C_X$ over some $X \in I\CW$. Therefore
universal objects exist. They are unique by Lemma~\ref{univ.uni.le},
so by Proposition~\ref{lift.prop}, it suffices to prove that for any
$X \in I\CW$, a universal object $x \in \C_X$ is liftable with
respect to any map $e:Y \to Y'$ in $I\CW$.

Indeed, assume given such a map, an object $y \in \C_{Y'}$, a map
$f:Y \to X'$ and an isomorphism $f^*x' \cong e^*y$. Then by
\eqref{exa.eq}, we can construct an object $y' \in \C_{Y' \copr_Y
  X}$ that restricts to $y$ resp.\ $x$ on $Y'$ resp.\ $X$, and then
again by Lemma~\ref{univ.exi.le}, we can construct a map $g:Y'
\copr_Y X \to X'$ in $I\CW$, a universal object $x' \in \C_{X'}$,
and an isomorphism $g^*x' \cong y'$. Then the map $h:X \to X'$ is in
$C$, and since $h^*x' \cong x$, it is also in $W$ by
Lemma~\ref{univ.uni.le}. Since $X$ is fibrant, we have a map $p:X'
\to X$ such that $p \circ h = \id$, and then $p^*$ is an
equivalence, so that $x \cong h^*x'$ implies $x' \cong p^*x$. This
exactly means that $x$ is liftable with respect to $e$.
\endproof

\subsection{Simplicial sets.}\label{simp.repr.subs}

Consider now the category $\Delta^o\Sets$ of simplicial sets, also
equipped with the standard model structure of \cite{qui.ho}: $C$ is
the class of all injective maps, $F$ is the class of all Kan
fibrations. This completely defines the model structure; however, it
is also known that $F \cap W$ consists of trivial Kan fibrations. It
immediately follows from the definitions that $C \cap W$ contains
all the horn embeddings \eqref{horn.eq}, and is stable under
coproducts, pushouts and countable compositions.

Apart from its mere existence --- that is rather non-trivial, but
can be used a black box --- we will need two general results on the
model structure on $\Delta^o\Sets$ concerning the size of things.
As in Example~\ref{sets.S.exa}, for any regular cardinal $\kappa$,
let $\Sets_\kappa \subset \Sets$ be the category of sets $S$ such
that $|S| < \kappa$, and let $\Delta^o\Sets_\kappa \subset
\Delta^o\Sets$ be the full subcategory spanned by simplicial sets
$X:\Delta^o \to \Sets$ that take values in $\Sets_\kappa \subset
\Sets$. Recall that as in Example~\ref{sets.SS.exa}, $\Sets_\kappa$
and $\Delta^o\Sets_\kappa$ are $\kappa$-cocomplete and
$\kappa'$-complete for any infinite cardinal $\kappa' \ll \kappa$.

\begin{lemma}\label{simp.S.le}
For any uncountable regular cardinal $\kappa$, the full subcategory
$\Delta^o\Sets_\kappa \subset \Delta^o\Sets$ is a model subcategory.
\end{lemma}

\proof{} Say that a map $f$ in $\Delta^o\Sets_\kappa$ lies in $C$,
$W$, $F$ iff it does so as a map in $\Delta^o\Sets$. To check that
this defines a model structure, we have to check that the
factorizations of maps in $\Delta^o\Sets_\kappa$ can be chosen so
that they also lie in $\Delta^o\Sets_\kappa$. To see that this is
the case, let us recall the standard construction of these
factorizations. Assume given such a map $f:X \to Y$, treat it as an
object in the arrow category $\Ar(\Delta^o\Sets_\kappa)$, and note
that for any $n \geq k \geq 0$, with the corresponding horn
embedding $v^k_n \in \Ar(\Delta^o\Sets_\kappa)$ of \eqref{horn.eq},
the set $S(f)^k_n$ of maps $v^k_n \to f$ in the arrow category
satisfies $|S(f)^k_n| < \kappa$. Since $\kappa$ is uncountable, the
same holds for the disjoint union $\copr_{n,k}S(f)^k_n$. Define
$\X_1 \in \Delta^o\Sets_\kappa$ by the pushout square
\begin{equation}\label{ex.sq}
\begin{CD}
\coprod_{n \geq k \geq 0}S(f)^k_n \times \V^k_n @>{\copr \id \times
  v^k_n}>> \coprod_{n \geq k \geq 0}S(f)(n,k) \times \Delta_n\\
@VVV @VVV\\
X @>{a_0}>> X_1,
\end{CD}
\end{equation}
with the natural map $f_1:X_1 \to Y$ such that $f=f_1 \circ a_0$,
and iterate the construction using $f_n:X_n \to Y$ at step $n$ to
obtain a system of simplicial sets $X_n \in \Delta^o\Sets_\kappa$,
$n \geq 1$ and maps $a_n:X_n \to X_{n+1}$, $f_{n+1}:X_n \to Y$ such
that $f_n=f_{n+1} \circ a_n$, and $a_n$ is a pushout of a coproduct
of horn embeddings. Then $X_\infty = \colim_nX_n$ also lies in
$\Delta^o\Sets_\kappa$, and we have a factorization
\begin{equation}\label{ab.facto}
\begin{CD}
X @>{a_\infty}>> X_\infty @>{f_\infty}>> Y
\end{CD}
\end{equation}
of the map $f$, where $a_\infty$ is the countable composition of the
maps $a_n$. By Quillen's small object argument, since $\V^k_n$ is a
compact object in $\Delta^o\Sets_\kappa$ for any $n \geq k \geq 0$,
$f_\infty$ has a lifting property with respect to all horn
embeddings, thus lies in $F$, and $a_\infty$ lies in $C \cap
W$. This provides one of the factorizations for $f$, and the other
one is obtained in the same way with horns replaced by spheres.
\endproof

\begin{remark}
The square \eqref{ex.sq} is the motivating example for the extension
construction of Subsection~\ref{ab.subs}, where it appears as
\eqref{ex.I.sq}.
\end{remark}

For our second general result, recall that a simplicial set $X$ is
{\em contractible} if the projection $X \to \ppt$ is a trivial Kan
fibration, and say that a map $f:X \to Y$ between simplicial sets is
{\em homotopy trivial} if the induced map $h(f):h(X) \to h(Y)$ in
the localization $h^W(\Delta^o\Sets)$ factors through $\ppt$. For
any simplicial set $Y$ and regular cardinal $\kappa$, let
$I_\kappa(Y)$ be the set of simplicial subsets $X \subset Y$ such
that $X \in \Delta^o\Sets_\kappa$, and $X$ is finite if $\kappa$ is
the countable cardinal.

\begin{lemma}\label{kappa.X.le}
Assume given a contractible simplicial set $Y$ and a subset $X \in
I_\kappa(Y)$. Then there exists $X' \in I_\kappa(Y)$ containing $X$
such that the embedding map $X \to X'$ is homotopy
trivial. Moreover, if $X=N(J)$ and $Y=N(I)$ are nerves of partially
ordered sets, one can also choose $X' = N(J')$, for a subset $J'
\subset I$.
\end{lemma}

\proof{} It is convenient to fix a point $x \in X([0]) \subset
Y([0])$. We have the geometric realization functor $R:\Delta^o\Sets
\to \Top$, and to insure that the embedding map $f:X \to X'$ is
homotopy trivial, it suffices to find a homotopy $h$ between $R(f)$
and the map $R(X) \to \ppt \cong \{x\} \subset R(X')$. Explicitly,
let $D^1$ be the unit interval, and for any topological space $Z$,
let $C(Z) = (Z \times D^1) \copr_{Z \times \{1\}} \ppt$ be the cone
over $Z$, with the vertex $o \in C(Z)$ given by the image of $Z
\times \{1\}$; then we need a map $h:C(R(X)) \to R(X')$ sending $o$
to $x$. The category $\Delta$ is a cellular Reedy category in the
sense of of Definition~\ref{ind.reedy.def}, so for any $n \geq 0$,
we have the $n$-th skeleton $\sk_n X \subset X$ of \eqref{sk.eq},
and the tautological map $\sk_{-1}X = \emptyset \to X$ is the
countable composition
\begin{equation}\label{sk.X}
X \cong \colim_n\sk_nX
\end{equation}
of the embeddings $\sk_nX \to \sk_{n+1}X$ of \eqref{sk.I.sq}. Assume
by induction on $n$ that we already have $X'_{n-1} \in I_\kappa(Y)$
containing $X$ and a homotopy $h_{n-1}:C(R(\sk_{n-1}(X))) \to
R(X'_{n-1})$. Then $R$ commutes with colimits, and by
\eqref{sk.I.sq}, $R(\sk_n X) \copr_{R(\sk_{n-1}(X))}
C(R(\sk_{n-1}(X)))$ is homotopy equivalent to a bouquet of spheres
$S^n$ indexed by the set $X_n \subset X([n])$ of non-degenerate
$n$-simplices in $X$. Thus effectively, we have a bunch of elements
$\alpha_s \in \pi_n(R(X'),x)$, $s \in X_n$ in the $n$-th homotopy
group of $R(X')$. Now let $I_n \in I_\kappa(Y)$ be the subset of $X'
\subset Y$ containing $X'_{n-1}$, and order it by inclusion. Then
$I_n$ is filtered, $Y \cong \colim_{X' \in I_n}X'$, $R$ preserves
colimits, and since spheres are compact, $\pi_n(-,x)$ commutes with
filtered colimits. Therefore for each $s \in X_n$, $\alpha_s$
becomes trivial in $\pi_n(X',x)$ for some $X'_s \in I_n$, and since
$\kappa$ is regular, the union $X'_n = \cup_sX'_s \subset Y$ is
still in $I_n$.  The map $R(\sk_n X) \copr_{R(\sk_{n-1}(X))}
C(R(\sk_{n-1}(X))) \to R(X'_n)$ is then homotopy trivial, so we can
extend $h_{n-1}$ to a continuous map $h_n:C(R(\sk_nX)) \to R(X'_n)$.

If $X$ is finite, $X=\sk_nX$ for some $n$, and we can take $X'=X'_n$
and finish the proof. If $\kappa$ is uncountable, take
$X'=\cup_nX'$. Finally, if $X=N(J)$, $Y=N(I)$, that at each step in
the construction, we can restrict our attention to subsets $X'
\subset Y$ of the form $N(J')$; then the subset $I'_n \subset I_n$
these form is still filtered, and $Y=\colim_{X' \in I'_n}X'$, so the
the same argument works.
\endproof

\begin{exa}\label{kappa.conn.exa}
Take a partially ordered set $I$ with contractible nerve
$N(I)$. Then Lemma~\ref{kappa.X.le} for the two-element discrete set
$J=\{0,1\}$ says that any two elements $i_0,i_1 \in I$ can be
connected by a finite zigzag of maps in $I$ --- that is, $I$ is
connected as a category.
\end{exa}

\begin{remark}
Instead of using geometric realizations in the proof of
Lemma~\ref{kappa.X.le}, one can notice that $X_\infty$ of
\eqref{ab.facto} is actually functorial in $X$, thus gives a
functorial fibrant replacement in $\Delta^o\Sets$ that commutes with
filtered colimits, and use homotopies between maps to these fibrant
replacements. In any case, as Example~\ref{kappa.conn.exa} shows,
any proof would involve some sort of inexplicit and non-constructive
procedure: the size of a zigzag in Example~\ref{kappa.conn.exa} is
impossible to control.
\end{remark}

The geometric realization functor $R:\Delta^o\Sets \to \Top$ that we
have used in Lemma~\ref{kappa.X.le} is left-adjoint to functor
$S:\Top \to \Delta^o\Sets$ sending a topological space to its
singular simplicial set, both functors send weak equivalences to
weak equivalences, and the adjunction maps $\id \to S \circ R$, $S
\circ R \to \id$ are pointwise weak equivalences. Therefore by
Lemma~\ref{loc.adj.le}, families of groupoids over $\langle
\Delta^o\Sets,W \rangle$ and $\langle \Top,W \rangle$ are naturally
identified, and Theorem~\ref{coho.thm} for $I=\ppt$ immediately
implies a corresponding statement for families over
$\Delta^o\Sets$. However, one can also give an alternative proof
that gives the following slightly more precise result.

\begin{defn}\label{repr.bnd.def}
For any uncountable regular cardinal $\kappa$, a family of
group\-oids $\C$ over $\Delta^o\Sets$ or $\Delta^o\Sets_\kappa$ is
{\em $\kappa$-bounded} iff $\|\C_X\| < \kappa$ for any finite $X \in
\Delta^o\Sets_\kappa \subset \Delta^o\Sets_\kappa$.
\end{defn}

\begin{theorem}\label{coho.bis.thm}
Assume given an uncountable regular cardinal $\kappa$ and a
$\kappa$-bounded family $\C$ of groupoids over $\langle
\Delta^o\Sets_\kappa, W \rangle$ such that $\C$ is semiexact in the
sense of Definition~\ref{semi.def} and additive in the sense of
Definition~\ref{add.def}. Then there exists a simplicial set $X \in
\Delta^o\Sets_\kappa$ and an equivalence $e:\Hh(X) \cong \C$, so
that $\C$ extends to a $\kappa$-bounded additive semiexact family of
groupoids over $\langle \Delta^o\Sets,W \rangle$. The simplicial set
$X$ is unique up to a unique isomorphism in
$h^W(\Delta^o\Sets)$. Moreover, for any two $\kappa$-bounded
additive semiexact families of groupoids $\C_0$, $\C_1$ over
$\langle \Delta^o\Sets,W \rangle$, a functor between their
restrictions to $\Delta^o\Sets_\kappa$ extends to a functor
$\gamma:\C_0 \to \C_1$ over the whole $\Delta^o\Sets$, and for any
two such functors $\gamma_0,\gamma_1:\C_0 \to \C_1$, an isomorphism
between their restrictions to $\Delta^o\Sets_\kappa$ uniquely
extends to an isomorphism $\gamma_0 \cong \gamma_1$.
\end{theorem}

\begin{remark}\label{sm.rem}
Any small family of groupoids $\C \to \Delta^o\Sets$ is obviously
$\kappa$-bounded for a large enough $\kappa$, thus falls within the
scope of Theorem~\ref{coho.bis.thm}.
\end{remark}

Just as for Theorem~\ref{coho.thm}, our proof of
Theorem~\ref{coho.bis.thm} depends on the representability criterion
of Proposition~\ref{lift.prop}, and it will be useful to do the main
construction in larger generality. Namely, assume given a general
cellular Reedy category $I$ in the sense of
Definition~\ref{ind.reedy.def}, and consider the category
$I^o\Sets$. It has all colimits, and we can turn it into a
$C$-category in the sense of Definition~\ref{CW.def} by letting $C$
be the class of injective maps. For any $X \in I^o\Sets$, we still
have the skeleton filtration, and the tautological map $\sk_{-1}X =
\emptyset \to X$ decomposes into the countable composition
\eqref{sk.X}. Let $I^o_f\Sets \subset I^o\Sets$ be the full
subcategory spanned by objects that are finite-dimensional in the
sense of Definition~\ref{I.dim.def} (that is, $X \cong \sk_nX$ for
some $n$), and note that $I^o_f\Sets$ inherits a C-category
structure. Generalizing Definition~\ref{repr.bnd.def}, say that $X
\in I^o\Sets$ is {\em finite} if $X \in I^o_f\Sets$ and $X(i)$ is
finite for any $i \in I$, and say that a family of groupoids $\C$
over $I^o\Sets$ is {\em $\kappa$-bounded} if $\|\C_X\| < \kappa$ for
any finite $X \in I^o\Sets$. With this terminology, here is our main
liftability result.

\begin{prop}\label{lft.I.prop}
  Assume given a cellular Reedy category $I$ and a small family of
  groupoids $\C$ over the $C$-category $I^o\Sets$ such that
\begin{enumerate}
\item over $I^o_f\Sets \subset I^o\Sets$, $\C$ is additive in the
  sense of Definition~\ref{add.def} and strongly semiexact in the
  sense of Definition~\ref{str.semi.def}, and
\item for any $X \in I^o\Sets$, $\C$ is semicontinuous over the
  countable composition \eqref{sk.X} in the sense of
  Definition~\ref{count.def}.
\end{enumerate}
Then there exists $X \in I^o\Sets$ and an object $c \in \C_X$ that
is liftable in the sense of of Definition~\ref{lift.def} with
respect to all the $i$-sphere embeddings \eqref{sph.I.eq}. Moreover,
if $\C$ is $\kappa$-bounded for some regular cardinal $\kappa$, and
$I$ is $\Hom$-finite, then one can choose $X$ in $I^o\Sets_\kappa$,
and it only depends on the restriction of $\C$ to $I^o\Sets_\kappa$.
\end{prop}

We postpone the proof of Proposition~\ref{lft.I.prop} until
Subsection~\ref{repr.lft.subs}. For now, we just note that while it
implies Theorem~\ref{coho.bis.thm} almost immediately, in fact it
gives more. Indeed, by Lemma~\ref{ind.reedy.le}, for any cellular
Reedy category $I$, the product $I \times \Delta$ is a cellular
Reedy category, thus lies within scope of
Proposition~\ref{lft.I.prop}. Let us define {\em $I$-simplicial
  sets} as functors $I^o \to \Delta^o\Sets$. Then the category
$I^o\Delta^o\Sets \cong \Fun(I^o,\Delta^o\Sets)$ they form carries a
Reedy model structure induced by the Kan-Quillen model structure on
$\Delta^o\Sets$. Moreover, if $I$ is Reedy-$\kappa$-bounded in the
sense of Definition~\ref{reedy.def} for some uncountable regular
cardinal $\kappa$, then category $I^o\Delta^o\Sets_\kappa$ also
carries a Reedy model structure by Lemma~\ref{simp.S.le}. In fact,
these Reedy model structures have a description very similar to the
Kan-Quillen model structure. Namely, for any functor $X:I^o \to
\Sets$ and simplicial set $Y:\Delta^o \to \Sets$, denote by $X
\boxtimes Y:I^o \times \Delta^o \to \Sets$ the $I$-simplicial set
given by $X \boxtimes Y(i \times [n]) = X(i) \times Y([n])$. For any
object $i \in I$ and integer $n \geq 0$, let $\Delta_{i,n} \in
I^o\Delta^o\Sets$ be the functor represented by $i \times [n] \in I
\times \Delta$, and for any $i \in I$, $n \geq k \geq 0$, let
\begin{equation}\label{s.v.eq}
\begin{aligned}
\SS_{i,n} &= (\Delta_i \boxtimes \SS_{n-1}) \copr_{\SS_i
  \boxtimes \SS_{n-1}} (\SS_i \boxtimes \Delta_n) \subset
\Delta_{i,n} = \Delta_i \boxtimes \Delta_n,\\
\V_{i,n}^k &= (\Delta_i \boxtimes \V^k_n) \copr_{\SS_i
  \boxtimes \V^k_n} (\SS_i \boxtimes \Delta_n) \subset
\Delta_{i,n} = \Delta_i \boxtimes \Delta_n,
\end{aligned}
\end{equation}
where $\SS_i$ and $\Delta_i$ are as in \eqref{sph.I.eq}, and
$\V^k_n$, $\SS_{n-1}$ are as in \eqref{horn.eq}, \eqref{sph.eq}. We
then have the natural injective maps
\begin{equation}\label{bi.sph.eq}
\sigma_{i,n}:\SS_{i,n} \hookrightarrow \Delta_{i,n}
\end{equation}
and
\begin{equation}\label{bi.horn.eq}
v^k_{i,n}:\V_{i,n}^k \hookrightarrow \Delta_{i,n}.
\end{equation}
Note that in terms of the cellular Reedy structure on $I \times
\Delta$ and the sphere embeddings \eqref{sph.I.eq}, we have
$\Delta_{i,n} = \Delta_{i \times [n]}$, $\SS_{i,n} = \SS_{i \times
  [n]}$, $\sigma_{i,n} = \sigma_{i \times [n]}$, and
\eqref{bi.sph.eq} is exactly \eqref{sph.I.eq}, while
\eqref{bi.horn.eq} is new.

\begin{lemma}\label{bisimp.le}
A map $f$ in the category $I^o\Delta^o\Sets$ is the class in $F$
resp.\ $F \cap W$ iff it it has the right lifting property with
respect to the horn embeddings \eqref{bi.horn.eq} resp.\ sphere
embeddings \eqref{bi.sph.eq}, and it is in $C$ iff it is injective.
The same is true for the full subcategory $I^o\Delta^o\Sets_\kappa
\subset I^o\Delta^o\Sets$ for any uncountable regular cardinal
$\kappa$ such that $I$ is Reedy-$\kappa$-bounded.
\end{lemma}

\proof{} For any cocomplete model category $\C$, cellular Reedy
category $I$, object $i \in I$ and functor $X:I^o \to \C$, we have
$\Hom(\Delta_i,X) \cong X(i)$ and, as we mentioned already in
Subsection~\ref{exp.subs}, \eqref{sph.I.M} immediately implies that
$\Hom(\SS_i,X) \cong M(X,i)$, so that a map $f:X \to Y$ in $I^o\C$
is in $F$ resp.\ $F \cap W$ with respect to the Reedy structure if
and only if the natural map
$$
X(i) \to Y(i) \times_{\Hom(\SS_i,Y)} \Hom(\SS_i,X)
$$
is in $F$ resp.\ $F \cap W$ for any $i \in I$. Taking
$\C=\Delta^o\Sets$ and plugging in \eqref{s.v.eq}, we obtain the
claim about $F \cap W$ and $F$. For the last claim, note that
exactly the same argument as in Lemma~\ref{simp.S.le} provides a
factorization \eqref{ab.facto} of a map $f$ such that $b$ is in $F
\cap W$, and $a$ is a countable composition of pushouts of
coproducts of the embeddings \eqref{bi.sph.eq}. Then $a$ is in $C$,
and if $f$ is in $C$, it is a retract of the map $a$. By
\eqref{sk.X} and \eqref{sk.I.sq}, the latter holds iff $f$ is
injective.
\endproof

\begin{corr}\label{sk.corr}
For any non-negative integer $n \geq 0$, the $n$-th skeleton functor
$\sk^I_n:I^o\Delta^o\Sets \to I^o\Delta^o\Sets$ taken with respect
to $I$ sends maps in $W$ to maps in $W$, and the same holds for
$I^o\Delta^o\Sets_\kappa \subset I^o\Delta^o\Sets$ for any
uncountable regular cardinal $\kappa$ such that $I$ is
Reedy-$\kappa$-bounded.
\end{corr}

\proof{} Since all objects in $I^o\Delta^o\Sets$ or
$I^o\Delta^o\Sets_\kappa$ are cofibrant, the same argument as in
Lemma~\ref{Q.adj.le} shows that it suffices to check that $\sk_n^I$
sends maps in $C \cap W$ to maps in $C \cap W$. By \eqref{sk.eq} and
adjunction, this amounts to checking that $\eps_n^*$ preserves $C
\cap W$ and $F$, and this is immediately obvious from the definition
of the Reedy model structure.
\endproof

\begin{exa}\label{holim.I.exa}
One immediately corollary of Lemma~\ref{bisimp.le} is that for any
functor $\gamma:I_0 \to I_1$ between cellular Reedy categories, the
pullback functor $\gamma^*:I_1^o\Delta^o\Sets \to
I_0^*\Delta^o\Sets$ sends maps in $C$ resp.\ $W$ to maps in $C$
resp.\ $W$. Thus it descends to a functor
$h^W(I^o_1\Delta^o\Sets) \to h^W(I^o_0\Delta^o\Sets)$, and this
functor has a right-adjoint provided by the Quillen Adjunction
Theorem~\ref{qui.adj.thm}. For any Reedy cellular $I_0=I$ and $I_1 =
\ppt$, this is the homotopy limit functor
\begin{equation}\label{holim.I.eq}
\holim_{I^o}:h^W(I^o\Delta^o\Sets) \to h^W(\Delta^o\Sets)
\end{equation}
right-adjoint to the tautological functor $h^W(\Delta^o\Sets) \to
h^W(I^o\Delta^o\Sets)$.
\end{exa}

\begin{theorem}\label{coho.bisi.thm}
Assume given an uncountable regular cardinal $\kappa$ and a
$\Hom$-finite Reedy-$\kappa$-bounded cellular Reedy category
$I$. Then for any $\kappa$-bounded family $\C$ of groupoids over
$\langle I^o\Delta^o\Sets_\kappa, W \rangle$ that is semiexact in
the sense of Definition~\ref{semi.def} and additive in the sense of
Definition~\ref{add.def}, there exists an $I$-simplicial set $X \in
I^o\Delta^o\Sets_\kappa$ equipped with an equivalence $\Hh(X) \cong \C$, so
that $\C$ extends to a $\kappa$-bounded additive semiexact family of
groupoids over $\langle I^o\Delta^o\Sets,W \rangle$.  The set $X$ is
unique up to a unique isomorphism in
$h^W(I^o\Delta^o\Sets)$. Moreover, for any two $\kappa$-bounded
additive semiexact families of groupoids $\C_0$, $\C_1$ over
$\langle I^o\Delta^o\Sets,W \rangle$, a functor between their
restrictions to $I^o\Delta^o\Sets_\kappa$ extends to a functor
$\gamma:\C_0 \to \C_1$ over $I^o\Delta^o\Sets$, and for any two such
functors $\gamma_0,\gamma_1:\C_0 \to \C_1$, an isomorphism between
their restrictions to $\Delta^o\Sets_\kappa$ uniquely extends to an
isomorphism $\gamma_0 \cong \gamma_1$.
\end{theorem}

\begin{remark}\label{sm.I.rem}
As in Remark~\ref{sm.rem}, for any cellular Reedy category $I$, any
additive semiexact small family of groupoids over $\langle
I^o\Delta^o\Sets, W \rangle$ falls within the scope of
Theorem~\ref{coho.bisi.thm} for a large enough $\kappa$.
\end{remark}

\proof[Proof of Theorem~\ref{coho.bis.thm} and
  Theorem~\ref{coho.bisi.thm}.] Note that by \eqref{sk.X} and
\eqref{sk.I.sq}, every injective map in $I^o\Sets$ for any cellular
Reedy category $I$ is a countable composition of pushouts of
$i$-sphere embeddings \eqref{sph.I.eq}. Then take $I = \Delta$ for
Theorem~\ref{coho.bis.thm}, or $I = I \times \Delta$ for
Theorem~\ref{coho.bisi.thm}, and combine
Proposition~\ref{lft.I.prop}, Proposition~\ref{lift.prop},
Lemma~\ref{lift.le}, Lemma~\ref{tele.le} and
Lemma~\ref{semi.mod.le}.
\endproof

\subsection{Liftability.}\label{repr.lft.subs}

Our proof of Proposition~\ref{lft.I.prop} works by induction on
skeleta of \eqref{sk.X}, and we start by a general liftablity result
needed for the induction step. Assume given two families of
groupoids $\pi:\C \to I$, $\pi':\C' \to I$ over a category $I$, and
a functor $\gamma:\C \to \C'$ over $I$. For any object $c' \in \C'$,
we then have a family of groupoids $\C /_\gamma c' \to I / i$, where
$i = \pi'(c') \in I$, and it fits into a cartesian square
\begin{equation}\label{C.I.sq}
\begin{CD}
  \C /_\gamma c' @>>> \sigma(i)^*\C\\
  @VVV @VV{\sigma(i)^*(\gamma)}V\\
  I / i @>{s(c')}>> \sigma(i)^*\C',
\end{CD}
\end{equation}
where $\sigma(i):I / i \to I$ is the discrete fibration \eqref{p.i},
and $s(c')$ is the tautological section sending $f:i' \to i$ to
$f^*c' \to c'$.

\begin{lemma}\label{lft.rel.le}
In the assumptions above, assume given an object $c \in \C$ such
that $\gamma(c) \in \C'$ is liftable with respect to a map $f$ in
$I$, and the tautological object $\langle c,\id \rangle \in \C
/_\gamma \gamma(c)$ is liftable with respect to any map $g$ in $I /
\pi(c)$ such that $\sigma(g)=f$. Then $c$ is liftable with respect
to $f$.
\end{lemma}

\proof{} Straightforward diagram chasing.
\endproof

Now as in Proposition~\ref{lft.I.prop}, assume given a cellular Reedy
category $I$, and for any $n \geq 0$, let $C_n = C \cap
\sk_n^*(\Iso)$ be the class of injective maps $f$ in $I^o\Sets$ such
that $\sk_n(f)$ is invertible. Then $\langle I^o\Sets,C_n \rangle$
is a $C$-category in the sense of Definition~\ref{CW.def}.

\begin{lemma}\label{lft.n.le}
Let $\C$ be a small family of groupoids over the relative category
$\langle I^o\Sets,C_n \rangle$ that is additive and strongly
semiexact with respect to the $C$-category structure on $I^o\Sets$
given by $C_{n-1}$, and assume that $\sk_{n-1}^*\C$ is
constant. Then there exists a functor $X \in I^o\Sets$ such that
$X(i)=\ppt$ whenever $\deg i < n$, and an object $c \in \C_X$
liftable with respect to all sphere embeddings
\eqref{sph.I.eq}. Moreover, if $I$ is $\Hom$-finite, and $\C$ is
$\kappa$-bounded for some regular cardinal $\kappa$, then $X$ can be
chosen in $I^o\Sets_\kappa \subset I^o\Sets_\kappa$, and it only
depends on the restriction of $\C$ to $I^o\Sets_\kappa$.
\end{lemma}

\proof{} Let $\Y:I \to I^o\Sets$ be the Yoneda embedding, and take
the decomposition \eqref{grp.eq} for the family of groupoids $\Y^*\C
\to I$. Then $\pi_0(\Y^*\C | I) \cong IX$ for some $X \in
I^o\Sets$. If $I$ is $\Hom$-finite, then $\Delta_i$, hence also
$\SS_i \subset \Delta_i$ is finite for any $i \in I$, and if $\C$ is
$\kappa$-bounded, then $|X(i)| < \kappa$. Since $\C$ is additive,
$\C_\emptyset = \ppt$, and since $\sk_{n-1}^*\C$ is constant, $\C_Y
\cong \ppt$ for any $Y$ of dimension $\dim Y < n$. In particular,
$X(i) = \ppt$ for any $i$ with $\deg i < n$. Moreover, the embedding
$\sk_nX \to X$ is in $C_n$, so $\C_{\sk_nX} \cong \C_X$, and all the
horizontal maps in \eqref{sk.I.sq} are in $C_{n-1}$, so $\C$ is by
assumption semiexact along this square. Then \eqref{sk.I.sq} and
\eqref{grp.eq} give epivalences
\begin{equation}\label{sk.pi}
\C_X \cong \C_{\sk_n X} \to \prod_{i \in \M(IX)_n}\Y^*\C_{\pi(i)}
\to \prod_{i \in \M(IX)_n}X(\pi(i)),
\end{equation}
where $\pi:IX \to I$ is the projection. By definition, an object $i
\in IX$ is of the form $\langle \pi(i),x(i) \rangle$, $\pi(i) \in
I$, $x(i) \in X(\pi(i))$, so the discrete groupoid in the right-hand
side of \eqref{sk.pi} has a tautological object $\prod x(i)$, and we
can choose an object $c \in \C_X$ that is sent to $\prod x(i)$ by
the epivalence \eqref{sk.pi}. We claim that $c$ is liftable with
respect to all embeddings \eqref{sph.I.eq}.

Indeed, assume given an embedding \eqref{sph.I.eq}, an object $c' \in
\C_{\Delta_i} \cong \Y^*\C_i$, and a map $a:\sigma_i^*c' \to c$. Let
$\pi_\C:\C \to I^o\Sets$ be the projection. Then by the construction
of $X$, there exists a unique map $f:\Delta_i \to X$ such that $f^*c
\cong c'$, we have $f \circ \sigma_i = \pi_\C(a)$ by virtue of this
uniqueness, so that $a:\sigma_i^*c' \to c$ factors as
\begin{equation}\label{lft.facto}
\begin{CD}
  \sigma_i^*c' @>{\alpha}>> \sigma_i^*f^*c \cong
  \pi_\C(a)^*c \to c
\end{CD}
\end{equation}
for some isomorphism $\alpha$ in $\C_{\SS_i}$, and we need to show
that $\alpha = \sigma_i^*(\beta)$ for some isomorphism $\beta:c'
\cong f^*c$ in $\C_{\Delta_i}$. But $\deg i \leq n$, then
$\C_{\SS_i} \cong \ppt$ is discrete, so any $\beta$ will do, and if
$\deg i > n$, then $\sigma_i$ is in $C_n$, so $\sigma_i^*$ is an
equivalence.
\endproof

\begin{remark}\label{post.can.rem}
While there is no uniqueness claim in Lemma~\ref{lft.n.le}, the
construction is pretty canonical: $X$ is uniquely defined by $\C$,
and $c \in \C_X$ is unique up to a non-unique isomorphism.
\end{remark}

Now for any family of groupoids $\C$ over $I^o\Sets$ and integer $n
\geq -1$, let $\C^{\leq n} = \sk_n^*\C$. The functorial maps
$\sk_{n-1} \to \sk_n$ of \eqref{sk.I.sq} induce functors $q:\C^{\leq
  n} \to \C^{n-1}$, so that effectively, we have a family of
groupoids $\C^\hdot \to \N^o \times I^o\Sets$ whose restriction to
$n \times I^o\Sets$ is $\C^{\leq n}$, and we have $\C \cong
p_*\C^\hdot$, where $p:\N^o \times I^o\Sets \to I^o\Sets$ is the
projection. If $\C$ is additive, then $\C_\emptyset = \ppt$, so that
$\C^{\leq -1}$ is the trivial family $I^o\Sets \to I^o\Sets$.

\begin{lemma}\label{sk.lift.le}
Assume given a small family of groupoids $\C \to I^o\Sets$ that
satisfies the assumption \thetag{i} of
Proposition~\ref{lft.I.prop}. For some $n \geq 0$, assume given a
functor $X:I^o \to \Sets$ and an object $c \in \C^{\leq n-1}_X$
liftable with respect to all the embeddings \eqref{sph.I.eq}. Then
there exists $X' \in I^o\Sets$, a map $r:X' \to X$, and an object
$c' \in \C^{\leq n}_{X'}$ equipped with an isomorphism $a:q(c')
\cong r^*c$ such that $c'$ is liftable with respect to all the
embeddings \eqref{sph.I.eq}, and $\sk_{n-1}(r):\sk_{n-1}X' \to
\sk_{n-1}X$ is an isomorphism. Moreover, if $\C$ is $\kappa$-bounded
for some regular cardinal $\kappa$, then $|r^{-1}(x)| < \kappa$ for
any $i \in I$ and element $x \in X(i)$.
\end{lemma}

\proof{} Identify $(IX)^o\Sets \cong I^o\Sets / X$ as in
\eqref{IX.sets}, and consider the family of groupoids $\C' =
\C^{\leq n} /_q c \to (IX)^o\Sets$. By Lemma~\ref{lft.rel.le}, to
lift $c$ to a liftable object in $\C^{\leq n}$, it suffices to
construct a liftable object in $\C'$, and by
Lemma~\ref{ind.reedy.le}, $IX$ is a cellular Reedy category, so it
suffices to check that $\C'$ satisfies the assumption of
Lemma~\ref{lft.n.le}. Indeed, by definition, $\C^{\leq n}$ and
$\C^{\leq n-1}$ are both constant along $C_n$, and $\sk_{n-1}^*(q)$
is an equivalence, so by \eqref{lft.facto} $\C'$ is constant along
$C_n$ and $\sk_{n-1}^*\C'$ is constant. Moreover, again by
definition, $\C^{\leq n-1}$ is constant along $C_{n-1}$, so it is
not only strongly semiexact but actually strongly exact with respect
to the $C$-category structure defined by $C_{n-1}$. Therefore
\eqref{lft.facto} also implies that $\C'$ is additive and strongly
semiexact.
\endproof

\proof[Proof of Proposition~\ref{lft.I.prop}.] By induction on $n$,
Lemma~\ref{sk.lift.le} provides a collection of functors $X_n:I^o
\to \Sets$, maps $r:X_n \to X_{n-1}$ and liftable objects $c_n \in
\C^{\leq n}_{X_n} \cong \C_{\sk_nX_n}$. Moreover, if we take $X =
\lim_nX$, then the natural map $\sk_nX \to \sk_nX_n$ is an
isomorphism, so $c_n$ is an object in $\C_{\sk_nX}$. Isomorphisms
$a$ of Lemma~\ref{sk.lift.le} then turn the collection of objects
$c_n$, $n \geq 0$, into an object in the target of the functor
\eqref{coli.exa} for the countable composition \eqref{sk.X}. Since
$\C$ is semicontinuous over this countable composition by
Proposition~\ref{lft.I.prop}~\thetag{ii}, there exists an object $c
\in \C_X$ equipped with isomorphisms $c|_{\sk_nX} \cong c_n$. To
check that $c$ is liftable with respect to an embedding
\eqref{sph.I.eq} for some $i \in I$ of degree $n = \deg i$, note
that a map from $\Delta_i$ or $\SS_i$ to $X$ automatically factors
through $\sk_nX$, and the liftability then follows from the
liftability of $c_n$. Finally, if $\C$ is $\kappa$-bounded, then
again by Lemma~\ref{sk.lift.le} and induction, $|\sk_nX(i)| =
|\sk_nX_n(i)| < \kappa$ for any $i \in I$ of degree $n = \deg i$,
and then $|X(i)| < \kappa$ for any $i \in I$.
\endproof

\begin{remark}\label{exa.rem}
For any cellular Reedy category $I$ and small family of groupoids
$\C \to I$, the family $\Y_*\C \to I^o\Sets$ trivially satifies the
assumptions of Proposition~\ref{lft.I.prop} -- in effect, we have
$\Y_*\C_X \cong \Sec(IX,\pi^*\C)$ by \eqref{2kan.eq}, so by
Lemma~\ref{cocart.le}, $\Y_*\C$ is not only strongly semiexact but
also strongly exact. However, even in this simple case the
construction is somewhat non-trivial. In particular, if $I = \Delta$
and $\C = \ppt_G \times \Delta$ is the constant family with fiber
$\ppt_G$ for some group $G$, then what our proof of
Proposition~\ref{lft.I.prop} produces is $X = BG$, with the liftable
object $c$ given by the functor $\xi \times \pi:\Delta X = \Delta
\ppt_G \to \ppt_G \times \Delta$, where $\xi$ is as in
\eqref{xi.eq}. Another illuminating example is to take a functor
$M:\Delta^o \to \Ab$ to abelian groups, and to consider the
corresponding family $\C \to \Delta$ with fibers $\C_{[n]} =
\ppt_{M([n])}$. Then $X$ is the underlying simplicial set of
$BM:\Delta^o \to \Ab$, where, if $M$ corresponds to a complex
$M_\idot$ by the Dold-Kan equivalence, $BM$ corresponds to its
homological shift $M_\idot[1]$.
\end{remark}

\begin{remark}\label{post.rem}
The proof of Theorem~\ref{coho.thm} given in
Subsection~\ref{top.subs} and the proof of
Theorem~\ref{coho.bis.thm} based on Proposition~\ref{lft.I.prop} are
related by a sort of Eckman-Hilton duality: in
Thereom~\ref{coho.thm}, we construct the representing object $X$ as
the colimit of elementary cofibrations obtained by attaching cells,
while in Theorem~\ref{coho.bis.thm}, $X$ appears as the limit of a
tower of Kan fibrations. The second construction is more canonical,
in that by Remark~\ref{post.can.rem}, it produces $X$ uniquely up to
a (non-unique) isomorphism, and any equivalence $\C \cong \C'$
between families is induced by a (non-canonical) isomorphism $X
\cong X'$ between the representing objects. In fact, it is not
difficult to see that while the maps $r:X_n \to X_{n-1}$ in the
tower are not surjective, it can be corrected by replacing
$X_n$ with $Y_n = r(X_{n+1}) \subset X_n$, and $Y_\idot$ is then the
minimal Postnikov tower of the simplicial set $X$. Namely, $Y_0 =
X_0$ is $\pi_0(X)$, the map $\pi_0(Y_n) \to Y_0$ is an isomorphism
for any $n$, and if for some fixed $o \in \pi_0(X) \cong
\pi_0(Y_n)$, we let $Y^o_n \subset Y_n$ be the corresponding
connected component, then $Y^o_1 \cong B\pi_1(Y,o)$, and $Y^o_n \to
Y^o_{n-1}$, $n \geq 2$ is a torsor with respect to $B^n\pi_n(Y,o)$.

For an arbitrary small category $I$, one cannot expect a similar
nice minimal Postnikov tower construction for objects in
$I^o\Delta^o\Sets$ -- among other things, this would provide a
canonical minimal projective resolution for any functor $M$ from
$I^o$ to abelian groups, and this is clearly not possible without
additional assumptions. The assumption we impose in
Theorem~\ref{coho.bisi.thm} is that $I$ is a cellular Reedy
category.
\end{remark}

\subsection{Finite dimension.}

Now again, assume given a cellular Reedy category $I$, and consider
the category $I^o_f\Delta^o_f\Sets = (I \times \Delta)^o_f\Sets$ of
$I$-simplicial sets that are finite-dimensional in the sense of
Definition~\ref{I.dim.def}. For any regular cardinal $\kappa$, let
$I^o_f\Delta^o_f\Sets_\kappa = I^o_f\Delta^o_f\Sets \cap
I^o\Delta^o\Sets_\kappa \subset I^o\Delta^o\Sets$, and denote by
\begin{equation}\label{i.i.k}
i:I^o_f\Delta^o_f\Sets \to I^o\Delta^o\Sets, \qquad
i_\kappa:I^o_f\Delta^o_f\Sets_\kappa \to I^o\Delta^o\Sets_\kappa
\end{equation}
the embedding functors.

Note that even if $I=\ppt$ is the point, the category
$I^o_f\Delta^o_f\Sets$ is {\em not} a model category: when a
finite-dimensional simplicial set $X$ has an infinite number of
non-trivial homotopy groups, its fibrant replacement cannot be
finite-dimensional. Nevertheless, $I^o_f\Delta^o_f\Sets$ has a
CW-structure in the sense of Definition~\ref{CW.def}. In fact, it
has several; the one that we will need is defined as follows.

\begin{defn}\label{f.CW.def}
The {\em class $C$} of maps in $I^o_f\Delta^o_f\Sets$ or
$I^o_f\Delta^o_f\Sets_\kappa$ consists of all injective maps. The
class $W$ of {\em anodyne maps} is the smallest saturated class of
maps in $I^o_f\Delta^o_f\Sets$ that is closed under coproducts,
contains the projections $\Delta_{i,n} \to \Delta_{i,0}$ for any $i
\in I$, $[n] \in \Delta$, and is such that $C \cap W$ is closed
under pushouts -- that is, for any $f:X \to X'$ in $C \cap W$ and
any map $X \to Y$, the map $f':Y \to Y' = Y \copr_X X'$ is in $C
\cap W$. The class of {\em $\kappa$-anodyne maps} is the smallest
class of maps in $I^o_f\Delta^o_f\Sets_\kappa$ with the same
properties.
\end{defn}

\begin{exa}
The class of all weak equivalences of $I$-simplicial sets obviously
satisfies the conditions of Definition~\ref{f.CW.def}, so that all
anodyne maps are weak equivalences, and the embedding functors
\eqref{i.i.k} are CW-functors. The same hold for $\kappa$-anodyne
maps for a large enough $\kappa$ (namely, $\kappa$ has to be
uncountable, and $I$ has to be Reedy-$\kappa$-bounded).
\end{exa}

By abuse of notation, for any $I$-simplicial set $X$ and simplicial
set $Z$, denote by $X \times Z = X \times (\ppt \boxtimes Z)$ the
$I$-simplicial set given by
$$
(X \times Z)(i \times [n]) = X(i \times [n]) \times Z([n]), \qquad i
\in I, [n] \in \Delta.
$$

\begin{lemma}\label{t.horn.le}
For any $i \in I$ and $n \geq l \geq 0$, the horn embedding
$v^l_{i,n}$ of \eqref{bi.horn.eq} is anodyne in the sense of
Definition~\ref{f.CW.def}, and $\kappa$-anodyne for any regular
cardinal $\kappa$ such that $|L(I,i)| < \kappa$. Moreover, for any
$m,n \geq 0$, the same is true for the embedding
\begin{equation}\label{t.horn.eq}
t^m_{i,n}:\SS_{i,m,n} = (\SS_{i,m} \times \Delta_n)
\copr_{\SS_{i,m}} \Delta_{i,m} \to \Delta_{i,m} \times \Delta_n,
\end{equation}
where the embedding $\SS_{i,m} \to \SS_{i,m} \times \Delta_n$ is
induced by $t:\ppt = \Delta_0 \to \Delta_m$.
\end{lemma}

\proof{} For any $[n] \in \Delta$, the projection $[n] \to [0]$ has
a section $t:[0] \to [n]$, and for any $i \in I$, this gives a
section $t_i:\Delta_{i,0} \to \Delta_{i,n}$ of the projection
$\Delta_{i,n} \to \Delta_{i,0}$. Since the class of $W$ of anodyne
maps is saturated, $t_{i,n}$ is in $W$, and it is
injective. Moreover, it factors through the source $\SS_{i,0,n}$ of
the embedding $t^0_{i,n}$ of \eqref{t.horn.eq}, and the induced map
$\Delta_{i,0} \to \SS_{i,0,n}$ is a composition of pushouts of
coproducts of maps $t^0_{i',n}$ with $i' \in L(I,i)$. Since $W$ is
anodyne, $t^0_{i,n} \in W$ by induction on $\deg i$. Now assume for
a moment that $I = \ppt$, and drop $i$ from notation. Then as in
Lemma~\ref{reedy.le}, $t^0_n = t_n:\Delta_0 \to \Delta_n$ factors
through the source $\V^l_n$ of the horn embedding \eqref{horn.eq},
and the induced map $\Delta_0 \to \V^l_n$ is a composition of
pushouts of horn embeddings $v^l_{n'}$ with $n' < n$, so by
induction on $n$, $v^l_n \in W$. Moreover, for any $m \geq 0$, the
map $t^m_n$ of \eqref{t.horn.eq} is a composition of pushouts of
horn embeddings, so it is also anodyne.

Returning to the general case, note that for any $Z \in
\Delta^o_f\Sets$, subsets $X' \subset \Delta_i \boxtimes Z$ that
contain $\SS_i \boxtimes Z$ are in one-to-one correspondence with
subsets $Z' \subset Z$, the embedding $t^m_{i,n}$ corresponds to the
embedding $t^m_n$, and if $Z' \subset Z'' \subset Z$ is a pushout of
coproducts of horn embeddings $v^l_n$ of \eqref{horn.eq}, then the
corresponding embedding $X' \subset X''$ is a pushout of coproducts
of horn embeddings $v^l_{i,n}$ of \eqref{bi.horn.eq}. Then exactly
the same induction as for $I=\ppt$ proves that $v^l_{i,n}$ and
$t^m_{i,n}$ are anodyne as required. In the $\kappa$-anodyne case,
repeat the same proof, and observe that since $|L(I,i)| < \kappa$,
all the coproducts that occur exist in $\Sets_\kappa$.
\endproof

\begin{corr}\label{Y.corr}
For any $Y \in I^o_f\Delta^o_f\Sets$ and any map $f:[n] \to [m]$ in
$\Delta$, the map $\id_Y \times f:Y \times \Delta_n \to Y \times
\Delta_m$ is anodyne, and $\kappa$-anodyne if $Y \in
I^o_f\Delta^o_f\Sets_\kappa$ and $I$ is Reedy-$\kappa$-bounded.
\end{corr}

\proof{} Since $W$ is saturated, it suffices to consider the maps
$t:[0] \to [n]$, and then by \eqref{sk.I.sq}, $\id_Y \times t$ is a
composition of pushouts of coproducts of the embeddings
\eqref{t.horn.eq}.
\endproof

\begin{corr}\label{f.CW.corr}
The category $I^o_f\Delta^o_f\Sets$ with the classes $C$ resp.\ $W$
of injective resp.\ anodyne maps of Definition~\ref{f.CW.def} is a
CW-category in the sense of Definition~\ref{CW.def}, and for any
regular cardinal $\kappa$ such that $I$ is Reedy-$\kappa$-bounded,
so is $I^o_f\Delta^o_f\Sets_\kappa$, with the classes $C$ resp.\ $W$
of injective resp.\ $\kappa$-anodyne maps.
\end{corr}

\proof{} Definition~\ref{CW.def}~\thetag{i} is a part of
Definition~\ref{f.CW.def}. For \thetag{ii}, note that as in
\eqref{cyl.deco}, any map $f:X \to Y$ in $I^o_f\Delta^o_f\Sets$
admits a decomposition
\begin{equation}\label{ff.dec}
\begin{CD}
  X @>{c}>> \Cyl(f) @>{w}>> Y,
\end{CD}
\end{equation}
where $\Cyl(f)$ is defined by the cocartesian square
\eqref{cyl.X.sq}. Then $w$ has an injective one-sided inverse $t:Y
\to \Cyl(f)$ that is a pushout of the embedding $t:X \to X \times
\Delta_1$, thus anodyne by Corollary~\ref{Y.corr}. Since $W$ is
saturated, $w \in W$, and then \eqref{ff.dec} works as \eqref{f.dec}
for the same reason as in Lemma~\ref{W.ano.W.le}.
\endproof

We will call the CW-structure of Corollary~\ref{f.CW.corr} the {\em
  standard CW-structure} on $I^o_f\Delta^o_f\Sets$ or
$I^o_f\Delta^o_f\Sets_\kappa$. We note that it is distributive in
the sense of Example~\ref{adm.CW.exa}, thus strong in the sense of
Definition~\ref{fr.def}. To fix a perfect framing, it suffices to
fix a decomposition \eqref{univ.cod}, and if $I=\ppt$, the simplest
one is that of Example~\ref{adm.CW.exa}~\thetag{ii}, that is,
\begin{equation}\label{D.univ}
\begin{CD}
  \ppt \copr \ppt @>{s \copr t}>> \Delta_1 @>>> \ppt,
\end{CD}
\end{equation}
where $\Delta_1$ is the standard $1$-simplex. We denote
$\overline{\SS}_1 = \Delta_1 \copr_{\ppt \copr \ppt} \Delta_1$
(equivalently, $\overline{\SS}_1$ is the $1$-sphere $\SS_1$ with
contracted edge $s(\Delta_1) \subset \SS_1 \subset \Delta_2$). If
$I$ is non-trivial, one has to pull back \eqref{D.univ} and
$\overline{\SS}_1$ to $I^o_f\Delta^o_f\Sets$ via the projection $I^o
\times \Delta^o \to \Delta^o$.

\begin{lemma}\label{f.CW.le}
Let $\I$ be $I^o_f\Delta^o_f\Sets$ or $I^o_f\Delta^o_f\Sets_\kappa$,
with the standard CW-struc\-ture, and let $\C$ be a semiexact
additive family of groupoids over $\langle I,W \rangle$. Then $\C$
is strongly semiexact in the sense of Definition~\ref{str.semi.def}
-- that is, it is semicartesian over any cocartesian square
\begin{equation}\label{f.CW.sq}
\begin{CD}
  X @>{f}>> Y\\
  @VVV @VVV\\
  X' @>{f'}>> Y'
\end{CD}
\end{equation}
in $\I$ with injective $f$, and cartesian if $X$ is empty. Moreover,
if $\C$ is constant along $f$ and $f \times \id_{\overline{\SS}_1}$,
then it is constant along $f'$. Finally, if $f$ admits a one-sided
inverse $e:Y \to X$, $e \circ f = \id$, with the induced inverse
$e':Y' \to X'$, then \eqref{f.CW.sq} yields an equivalence
$$
\C_{Y'}^e \cong \C_{X'} \times_{\C_X} \C^e_Y,
$$
where $\C^e_Y \subset \C_Y$ resp.\ $\C^e_{Y'} \subset \C_{Y'}$ are
the essentialy images of the functors $e^*$ resp.\ ${e'}^*$.
\end{lemma}

\proof{} For any square \eqref{f.CW.sq}, the embeddings $Y \to
\Cyl(f)$, $Y' \to \Cyl(f')$ are in $W$, and Corollary~\ref{Y.corr}
shows that the embedding $X' \copr_X \Cyl(f) \to \Cyl(f')$ is also
in $W$. Then $\C_{Y'} \cong \C_{X' \copr_X \Cyl(f)}$, $\C_{\Cyl(f)}
\cong \C_Y$, and it remains to apply Definition~\ref{semi.def},
Lemma~\ref{ano.C.le} and Lemma~\ref{retra.le} to $X' \copr_X
\Cyl(f)$.
\endproof

\begin{remark}
As an alternative, one can equip $I^o_f\Delta^o_f\Sets$ with a
smaller class $W'$ formed by compositions of pushouts of coproducts
of horn extensions \eqref{bi.horn.eq}. With slightly more care, one
can show that \eqref{ff.dec} would still works as \eqref{f.dec}, so
we would still obtain a CW-structure. However, Lemma~\ref{f.CW.le}
would also work and easily imply that a semiexact additive family of
groupoids $\C \to I^o_f\Delta^o_f\Sets$ constant along $W'$ is also
constant along anodyne maps of Definition~\ref{f.CW.def}. Thus the
end result would be the same.
\end{remark}

\subsection{Semicanonical extensions.}

Let us now turn back to the model categories $I^o\Delta^o\Sets$,
$I^o\Delta^o\Sets_\kappa$. It turns out that both also carry a
CW-category structure with a smaller class $W$ that is still large
enough to allow one to prove Theorem~\ref{coho.bisi.thm}. Namely,
let us introduce the following generalization of
Definition~\ref{f.CW.def}.

\begin{defn}\label{f.pl.CW.def}
The class $W^+$ of {\em $+$-anodyne maps} is the smallest saturated
class of maps in $I^o\Delta\Sets$ that is closed under coproducts,
contains the projections $\Delta_{i,n} \to \Delta_{i,0}$ for any $i
\in I$, $[n] \in \Delta$, and is such that $C \cap W^+$ is closed
under pushouts in the same sense as in Definition~\ref{f.CW.def},
and for any functors $X_\idot,Y_\idot:\N \to I^o\Delta^o\Sets$ that
factor through $(I^o\Delta^o\Sets)_C$, and any map $f_\idot:X_\idot
\to Y_\idot$ such that $f_n$ is $+$-anodyne for any $n$, the map
$f:X = \colim_nX_n \to Y = \colim_nY_n$ is $+$-anodyne. The class
$W^+$ of {\em $\kappa^+$-anodyne maps} is the smallest class of maps
in $I^o_f\Delta^o_f\Sets_\kappa$ with the same properties.
\end{defn}

\begin{exa}
All anodyne resp.\ $\kappa$-anodyne maps are obviously $+$-ano\-dyne
resp.\ $\kappa^+$-anodyne, so that the embeddings \eqref{i.i.k} send
$W$ into $W^+$. On the other hand, weak equivalences satisfy the
extra condition of Definition~\ref{f.pl.CW.def} by
Lemma~\ref{Q.adj.le}, so $W^+ \subset W$. The inclusion is probably
strict (we haven't checked).
\end{exa}

\begin{exa}\label{tele.X.exa}
Recall that for any functor $X_\idot:\N \to I^o_f\Delta^o_f\Sets$
that factors through $(I^o_f\Delta^o_f\Sets)_C$, with $X =
\colim_nX_n \in I^o\Delta^o\Sets$, a choice of a cofibrant
replacement of the functor $\zeta^*X_\idot:\ZZ_\infty \to
I^o_f\Delta^o_f\Sets$ provides a telescope $\Tel(X_\idot) \in
I^o\Delta^o\Sets$ of \eqref{tel.Y.eq} equipped with a map
\begin{equation}\label{tele.X}
  a:\Tel(X_\idot) \to X = \colim_{\N}X_\idot.
\end{equation}
Then exactly the same argument as in Lemma~\ref{tele.le} shows that
the map \eqref{tele.X} is $+$-anodyne in the sense of
Definition~\ref{f.pl.CW.def}, and it is $\kappa^+$-anodyne if
$X_\idot$ takes values in $I^o_f\Delta^o_f\Sets_\kappa$.
\end{exa}

\begin{lemma}\label{f.pl.CW.le}
The category $I^o\Delta^o\Sets$ with the classes of injective and
$+$-anodyne maps is a strong CW-category, and for any regular
cardinal $\kappa$, so is the category $I^o\Delta^o\Sets_\kappa$,
with the classes of injective and $\kappa^+$-anodyne maps.
\end{lemma}

\proof{} Exactly the same argument as in Corollary~\ref{f.CW.corr}
works; the only thing to check is that the embedding $t:X \to
\Delta^1 \times X$ is $+$-anodyne for any $X$, and
$\kappa^+$-anodyne for $X \in I^o\Delta^o\Sets_\kappa$. This
immediately follows by presenting $X = \colim_n\sk_nX$ as the
countable composition \eqref{sk.X}. A perfect framing in both cases
is provided by \eqref{D.univ}.
\endproof

\begin{defn}\label{X.ext.def}
A {\em semicanonical extension} of an additive semiexact family of
groupoids $\C$ over $\langle I^o_f\Delta^o_f\Sets_f,W \rangle$ is an
additive semiexact family $\C^+$ over $\langle I^o\Delta^o\Sets,W^+
\rangle$ equipped with an equivalence $i^*\C^+ \cong \C$. For any
regular cardinal $\kappa$, a {\em semicanonical extension} of an
additive semiexact family of groupoids $\C$ over $\langle
I^o_f\Delta^o_f\Sets_\kappa,W \rangle$ is an additive semiexact
family $\C^+$ over $\langle I^o\Delta^o\Sets_\kappa,W^+ \rangle$
equipped with an equivalence $i_\kappa^*\C^+ \to \C$. A
semicanonical extension $\C^+$ is {\em strong} if it is constant
along all maps that are weak equivalences in $I^o\Delta^o\Sets$.
\end{defn}

\begin{lemma}\label{can.le}
For any additive semiexact family of groupoids $\C$ over $\langle
I^o_f\Delta^o\Sets_f,W \rangle$, and any semicanonical extension
$\C^+$ in the sense of Definition~\ref{X.ext.def}, the functor $\C^+
\to i_*\C$ adjoint to the equivalence $i^*\C^+ \to \C$ is an
epivalence, and the same holds for semicanonical extensions of
families over $\langle I^o_f\Delta^o_f\Sets_\kappa,W \rangle$, for
any uncountable regular $\kappa$.
\end{lemma}

\proof{} For any $X \in I^o\Delta^o\Sets$, the family $\C^+$ is
semicontinous along the countable composition \eqref{sk.X} by
Lemma~\ref{tel.I.le}, and $X_\idot:\N \to I^o_f\Delta^o_f\Sets$
gives an embedding $\N \to I^o_f\Delta^o_f\Sets / X$ that is
left-admissible, with the adjoint functor sending $f:Y \to X$ to
$\dim f(Y) \subset X$. It remains to apply \eqref{2kan.eq} and
\eqref{2adm.eq}. The argument for $I^o\Delta^o\Sets_\kappa$ is
exactly the same.
\endproof

\begin{corr}\label{can.corr}
Assume given two additive semiexact small families of group\-oids
$\C$, $\C'$ over $\langle I^o_f\Delta^o_f\Sets,W \rangle$ equipped
with semicanonical extensions $\C^+$, ${\C'}^+$ in the sense of
Definition~\ref{X.ext.def}. Moreover, assume that ${\C'}^+$ is
strong. Then any functor $\gamma:\C' \to \C$ over
$I^o_f\Delta^o_f\Sets$ extends to a functor $\gamma^+:{\C'}^+ \to
\C^+$ over $I^o\Delta^o\Sets$, and a morphism $\gamma_0 \to
\gamma_1$ between two such functors with some extensions
$\gamma_0^+$, $\gamma_1^+$ extends to a morphism $\gamma_0^+ \to
\gamma_1^+$. Moreover, if $\gamma$ is an equivalence, then so is
$\gamma^+$. For any uncountable regular cardinal $\kappa$ such that
$I$ satisfies the assumptions of Theorem~\ref{coho.bisi.thm}, the
same is true for $\kappa$-bounded additive semiexact families of
groupoids over $\langle I^o_f\Delta^o_f\Sets_\kappa,W \rangle$.
\end{corr}

\proof{} By Theorem~\ref{coho.bisi.thm}, we have ${\C'}^+ \cong
\Hh(X)$ for some $I$-simplicial set $X$. The projection ${\C'}^+
\cong \Hh(X) \to i_*\C'$ corresponds to an object $c' \in
(i_*\C')_X$ that gives an object $c = i_*(\gamma)(c') \in (i_*\C)_X
\subset i_*\C$, and different functors $\gamma^+:{\C'}^+ \to \C^+$
that extend $\gamma$ correspond to objects $c^+ \in \C^+_c$, where
the fiber is taken with respect to the projection $\C^+ \to
i_*\C$. Since $\C^+ \to i_*\C$ is an epivalence by
Lemma~\ref{can.le}, its fibers are connected groupoids, so this
proves the first claim for $I^o_f\Delta^o_f\Sets$.

For the second claim, for any $X \in I^o\Delta^o\Sets$, we
have the $+$-anodyne map \eqref{tele.X} for some cofibrant
replacement $X':\ZZ_\infty \to I^o_f\Delta^o_f\Sets$ of
$\zeta^*X_\idot$, and both $\C^+$ and ${\C'}^+$ are constant over
\eqref{tele.X}. As in Lemma~\ref{tel.I.le}, consider the fibration
$z:\ZZ_\infty \to \VV$ of \eqref{ZZ.VV}, and compose it with the
cofibration $\VV \to \V$ of Remark~\ref{V.I.Z.rem} corresponding to
\eqref{VV.eq} to obtain a functor $z':\ZZ_\infty \to \V$. Then
$\Tel(X_\idot) \cong \colim_{\V}z'_!X'$, and as in
Lemma~\ref{tel.I.le}, this presents $\Tel(X_\idot)$ as a standard
pushout square of coproducts of finite-dimensional $I$-simplicial
sets. Moreover, $z'_!(X')(o) \times \overline{\SS}_1$ has the same
form. Since both $\C^+$ and ${\C'}^+$ are additive, and $\gamma^+$
is an equivalence over $I^o_f\Delta^o_f\Sets$, we can apply
Corollary~\ref{CW.pres.corr} and conclude that $\gamma^+$ is an
epivalence over $\Tel(X_\idot)$, hence also over $X$. Since $X$ was
arbitrary, $\gamma^+$ is an epivalence everythere, and then it is an
equivalence by Corollary~\ref{semi.epi.corr}.

This finishes the proof for $I^o\Delta^o\Sets$; for
$I^o\Delta^o\Sets_\kappa$, note that the object $X$ representing
${\C'}^+$ lies in $I^o\Delta^o\Sets_\kappa$, and use the
same argument.
\endproof

\begin{remark}
The map $\gamma_0^+ \to \gamma_1^+$ in Corollary~\ref{can.corr} is
{\em not} unique -- that is why the extension is semicanonical and
not canonical.
\end{remark}

\begin{prop}\label{X.ext.prop}
For any small additive semiexact family of groupoids $\C$ over
$\langle I^o_f\Delta^o_f\Sets,W \rangle$, a strong semicanonical
extension $\C^+$ in the sense of Definition~\ref{X.ext.def}
exists. The same is true for any uncountable regular cardinal
$\kappa$ such that $I$ satisfies the assumptions of
Theorem~\ref{coho.bisi.thm}, and $\kappa$-bounded additive semiexact
family $\C$ over $\langle I^o_f\Delta^o_f\Sets_\kappa,W \rangle$.
\end{prop}

In order to prove Proposition~\ref{X.ext.prop}, we use the
description of the representable family $\Hh(X)$, $X \in
I^o\Delta^o\Sets$ in terms of Reedy cofibrant replacements
$i_\Delta$ given in Subsection~\ref{repr.cr.subs}. However, while
Subsection~\ref{repr.cr.subs} dealt with arbitrary model categories,
we are now only interested in $I^o\Delta^o\Sets$, and here we can do
better: we can choose cofibrant replacements $Y_\Delta$ for all $Y
\in I^o\Delta^o\Sets$ in a functorial way. Namely, for any cellular
Reedy category $I$, one can consider the Yoneda embedding $\Y:I \to
I^o\Sets$ as an object in $\Fun(I,I^o\Sets)$, and then by
\eqref{sph.I.M}, for any $i \in I$, the latching object $L(\Y,i)$ is
exactly the $i$-sphere $\SS_i$, and the latching map $L(\Y,i) \to
\Y(i)$ is the sphere embedding \eqref{sph.I.eq}. In particular, it
is injective. For $I = \Delta$, this means that $\Y \in
\Fun(\Delta,\Delta^o\Sets)$ is cofibrant with respect to the Reedy
model structure. Moreover, since for any $I$, the class $C$ in
$I^o\Delta^o\Sets$ consists of injective maps by
Lemma~\ref{bisimp.le}, the same is true for $\Y \times Y \in
\Fun(\Delta, I^o\Delta^o\Sets)$, for any $Y \in
I^o\Delta^o\Sets$. Now consider the product $I^o\Delta^o\Sets \times
\Delta$, denote by $\theta:I^o\Delta^o\Sets \times \Delta \to
I^o\Delta^o\Sets$ the projection, and let
\begin{equation}\label{X.ups}
  \ups:I^o\Delta^o\Sets \times \Delta \to I^o\Delta^o\Sets
\end{equation}
be the functor given by $\ups(X \times [n]) = X \times
\Y([n])$. Then for any $Y \in I^o\Delta^o\Sets$, $Y_\Delta =
\ups|_{Y \times \Delta}:\Delta \to I^o\Delta^o\Sets$ is Reedy
cofibrant, and we have the following.

\begin{lemma}\label{I.cosimp.le}
  For any fibrant $I$-simplicial set $X \in I^o\Delta^o\Sets$, we
  have $\Hh(X) \cong \theta_!\ups^*h(X)$, where $\ups$ is the
  functor \eqref{X.ups}, and $h(X) \to I^o\Delta^o\Sets$ is the
  category of elements of $\Y(X):(I^o\Delta^o\Sets)^o \to \Sets$.
\end{lemma}

\proof{} Choose a Reedy fibrant replacement $X^\Delta \in
\Delta^oI^o\Delta^o\Sets$ of the constant simplicial $I$-simplicial
set $X$, starting with $X^\Delta([0]) = X$, use the description of
$\Hh(X)$ in terms of $X^\Delta$ given in Proposition~\ref{Hh.prop},
and then apply the same argument as in Lemma~\ref{del.del.le}.
\endproof

\begin{remark}\label{I.cosimp.rem}
By \eqref{2kan.eq}, the discrete fibration $h(X) \to
I^o\Delta^o\Sets$ in Lemma~\ref{I.cosimp.le} can be identified with
$\Y_*I\Delta X$, where $I\Delta X \to I \times \Delta$ is the
category of elements of $X$, and $\Y:I \times \Delta \to
I^o\Delta^o\Sets$ is the Yoneda embedding. Analogously,
$v^*h(X) \cong (\Y \times \pi)_*I\Delta X$, where $\pi:I \times
\Delta \to \Delta$ is the projection.
\end{remark}

\proof[Proof of Proposition~\ref{X.ext.prop}.] We only do the
absolute case; the proof in the $\kappa$-bounded case is exactly the
same. Let $\C' = i_*\C$, and note that it satisfies
Proposition~\ref{lft.I.prop}~\thetag{i} by
Lemma~\ref{f.CW.le}. Moreover, by the same argument as in
Lemma~\ref{can.le}, it is continuous along the skeleton filtration
\eqref{sk.X} for any $X$.  In particular, $\C'$ also satisfies
Proposition~\ref{lft.I.prop}~\thetag{ii}, and we can choose $X \in
I^o\Delta^o\Sets$ and a liftable object $c \in \C'_X$. Since horn
embeddings \eqref{bi.horn.eq} are anodyne by Lemma~\ref{t.horn.le},
and $c$ is liftable, $X$ is fibrant by Lemma~\ref{bisimp.le}, and we
have $\Hh(X) \cong \theta_!\ups^*h(X)$ by Lemma~\ref{I.cosimp.le},
while $c$ induces a functor $\lambda:h(X) \to \C'$ over
$I^o\Delta^o\Sets$. But for any $Y \in I^o_f\Delta^o_f\Sets$,
$Y_\Delta^*\C' \cong Y_\Delta^*\C$ is constant along all maps in
$\Delta$ by Corollary~\ref{Y.corr}. Since $\C'$ is continuous along
the skeleton filtrations, $Y_\Delta^*\C'$ is then constant for any
$Y \in I^o\Delta^o\Sets$. We conclude that $v^*\C' \cong
\theta^*\C'$, so that $\lambda$ induces a functor
$\ups^*(\lambda):\ups^*h(X) \to \theta^*\C'$ that
provides a functor $\lambda_\dg:\Hh(X) \cong \theta_!\eps^*h(X) \to
\C'$ over $I^o\Delta^o\Sets$ by adjunction. By the same argument as
in the proof of Proposition~\ref{lift.prop}, $\lambda_\dg$ is an
equivalence over $I^o_f\Delta^o_f\Sets$. Moreover, since for any $Y
\in I^o\Delta^o\Sets$, $\Hh(X)$ is semicontinuous along the skeleton
filtration $\sk_n Y$ by Lemma~\ref{repr.le} and Lemma~\ref{tele.le},
and $\C'$ is continuous along the same filtration, $\lambda_\dg$ is
an epivalence over the whole $I^o\Delta^o\Sets$. Then $\C^+=\Hh(X)$
with the functor $\lambda_\dg$ is our semicanonical extension.
\endproof

\begin{corr}\label{W.CW.corr}
Any small additive semiexact family of groupoids $\C$ over $\langle
I^o_f\Delta^o_f\Sets,W \rangle$ is constant along maps $w:X \to
X'$ that are weak equivalences of $I$-simplicial sets. For any
uncountable regular cardinal $\kappa$ such that $I$ satisfies the
assumptions of Theorem~\ref{coho.bisi.thm}, the same holds for any
$\kappa$-bounded additive semiexact family of groupoids $\C$ over
$\langle I^o_f\Delta^o_f\Sets_\kappa,W \rangle$. All semicanonical
extensions are strong and unique up to a unique equivalence.
\end{corr}

\proof{} A strong semicanonical extension exist by
Proposition~\ref{X.ext.prop}, and then it is unique up to a unique
equivalence by Corollary~\ref{can.corr}.
\endproof

\begin{defn}\label{st.I.def}
Let $\I$ be the category $I^o\Delta^o\Sets$, $I^o_f\Delta^o_f\Sets$,
or either of the categories $I^o\Delta^o\Sets_\kappa$,
$I^o_f\Delta^o_f\Sets_\kappa$ for a regular cardinal $\kappa$. Then
a family of groupouds $\C$ over $\I$ is {\em stably constant} along
a map $f$ in $\I$ if it is constant along $f \times \id_Z$ for any
finite simplicial set $Z$.
\end{defn}

\begin{corr}\label{W.C.corr}
Let $\C$ be an additive semiexact small family of group\-oids over
$\langle I^o_f\Delta^o_f\Sets,W\rangle$, with semicanonical
extension $\C^+$, and let $W(\C)$ be the class of maps $f$ in
$I^o\Delta^o\Sets$ such that $\C^+$ is stably constant along
$f$. Then for any object $J \in \Posf^+$, functors $X,X':J \to
I^o\Delta^o\Sets$ cofibrant with respect to the projective model
structure, and a map $f:X \to X'$ such that $f(j) \in W(\C)$ for any
$j \in J$, the map $\colim_J(f):\colim_J(X) \to \colim_J(X')$ is
also in $W(\C)$. Moreover, if $\C^+ \cong \Hh(X)$ for some fibrant
$X \in I^o\Delta^o\Sets$, then for any injective map $f:Y \to Y'$ in
$I^o\Delta^o\Sets$, \thetag{i} $f \in W(\C)$ iff \thetag{ii} for any
injective map $Z \to Z'$ of finite simplicial sets, $\Hom(-,X)$
sends the induced map $f(Z,Z'):(Z \times Y') \copr_{Z \times Y} (Z'
\times Y) \to Z' \times Y'$ to a surjection. Furthermore,
\thetag{ii} holds for all injective maps $Z \to Z'$ iff it holds all
the sphere embeddings \eqref{sph.eq}.
\end{corr}

\proof{} For the first claim, note that we may assume that $f$ is in
the class $C$ with respect to the projective model structure --
indeed, we can factor it as
$$
\begin{CD}
  X @>{c}>> X'' @>{f'}>> X',
\end{CD}
$$
with $c \in C$ and $f' \in F \cap W$, and then $\colim_J(f) \in W$
by Lemma~\ref{Q.adj.le}. If $J=S^<$ and $X(o)=X'(o)=\emptyset$, then
the claim simply means that $W(\C)$ is closed under coproducts, and
this is immediate from additivity. If $X$ and $X'$ are arbitrary but
$S = \{0,1\}$ so that $J=\V$, the claim follows from
Lemma~\ref{f.CW.le}. If $J=\N$, then by Example~\ref{tele.X.exa}, it
suffices to prove that $\Tel(f) \in W(\C)$; by the same argument as
in Lemma~\ref{tel.I.le} and Remark~\ref{tel.I.rem}, this reduces to
additivity and the case $J=\V$. Finally, for a general $J$, we can
first consider the height function $\hht:J \to \N$, and apply the
claim for $\hht_!(f):\hht_!X \to \hht_!X'$ to reduce to the case
when $J \in \Posf$, and then as in Lemma~\ref{ab.CW.le}, apply
induction on dimension to reduce to the case $J = S^<$. Then in this
case, by \eqref{sk.X}, it suffices to consider the case when both
$X$ and $X'$ factor through $I^o_f\Delta^o_f\Sets$, and the colimits
$\colim_J X$, $\colim_JX'$ also lie in $I^o_f\Delta^o\Sets$. But
then we have the cocartesian square
$$
\begin{CD}
  X(o) \times S @>>> X(o)\\
  @VVV @VVV\\
  \coprod_{s \in S}X(s) @>>> \colim_JX
\end{CD}
$$
in $I^o_f\Delta^o_f\Sets$, and similarly for $X'$, so we are
done by Lemma~\ref{f.CW.le}.

For the second claim, if \thetag{ii} holds for all sphere
embeddings, it holds for all injective maps by induction on the
elementary extensions \eqref{elem.X.sq}. By Lemma~\ref{I.cosimp.le},
we have $\Hh(X)_Y \cong h^\Tot(\Delta \Hom^\Delta(Y,X))$, where the
simplicial set $\Hom^\Delta(Y,X):\Delta^o \to \Sets$ sends $[n] \in
\Delta$ to $\Hom(Y \times \Delta_n,X)$, and similarly for $Y'$. Then
\thetag{ii} for sphere embeddings is equivalent to saying that the
pullback map $f^*:\Hom^\Delta(Y',X) \to \Hom^\Delta(Y,X)$ is a
trivial Kan fibration, and by Lemma~\ref{kan.le}, this means that
$\Hh(X)$ is constant along $f$. Moreover, \thetag{ii} for $f$
contains \thetag{ii} for $f \times \id_Z$, for any finite $Z \in
\Delta^o\Sets$, so the same then holds for any such $f \times
\id_Z$, so that \thetag{ii} implies \thetag{i}. In addition to this,
since $X$ is fibrant, $f^*$ is in any case a Kan fibration of
fibrant simplicial sets, and then if $\Hh(X)$ is constant along $f$,
so that $f^*$ induces an isomorphism on $\pi_0$, Lemma~\ref{kan.eq}
immediately implies that its component $f^*:\Hom(Y',X) \to
\Hom(Y,X)$ over $[0] \in \Delta$ is surjective. But then, since
$W(\C)$ is saturated and closed under coproducts and pushouts by
Lemma~\ref{f.CW.le}, the same induction on elementary extensions
shows that it contains $f(Z,Z')$ as soon as it contains $f$, so
\thetag{i} implies \thetag{ii}.
\endproof

\begin{remark}\label{W.C.rem}
In particular, Corollary~\ref{W.C.corr}~\thetag{ii} applies to the
embedding $\emptyset \to \ppt = \Delta_0$, and says that $\Hom(Y',X)
\to \Hom(Y,X)$ is surjective as soon as $f \in W(\C)$.
\end{remark}

\chapter{Between simplices and orders.}\label{nerve.ch}

As we have mentioned in the Introduction, this chapter is the
``technical heart'' of the text, something that has no conceptual
meaning but needs to be done. What happens is, to move from
Representability Theorems of Chapter~\ref{semi.ch} to the enhanced
category theory of Chapter~\ref{enh.ch}, we need to do two
things. Firstly, we need to replace families over simplicial sets
with families over partially ordered sets. Secondly, we need to
replace the homotopy invariance condition of
Theorem~\ref{coho.bis.thm} with something more manageable and
concrete.

We start with enhanced groupoids whose theory is easier and can
serve as a warm-up. This is Section~\ref{grp.sec}. Homotopy types
represent families over $\Delta^o\Sets$. On the other side, we have
the category $\Posf^+$ of left-bounded partially ordered
sets. Moreover, we can also restrict our attention to the
subcategory $\Delta^o_f\Sets$ of finite-dimensional simplicial sets;
these correspond to partially ordered sets of finite chain
dimension. To move from the category $\Posf$ of such sets to
$\Delta^o_f\Sets$, we use the nerve functor $\Nn:\Posf \to
\Delta^o_f\Sets$. In the other direction, we adapt the functorial
anodyne resolution construction of
Proposition~\ref{ano.I.prop}. This gives a functor $\QQ:\Posf \to
\Delta^o_f\Sets$.

We then show that $\Nn$ and $\QQ$ establish a bijection between
homotopy invariant semiexact additive families of groupoids over
$\Delta^o_f\Sets$ --- that is, the families representing homotopy
types --- and ``enhanced groupoids'', defined as additive semiexact
families of groupoids over $\Posf$ satisfying one simple additional
property (``excision'' of Definition~\ref{pos.def}). This works as a
replacement for homotopy invariance, and the reason it is so simple
is the general results of Subsection~\ref{fram.subs}.

Then we also consider semicanonical extensions of our families in
the sense of Proposition~\ref{X.ext.prop}, and we define the
corresponding notion on the enhanced groupoid side --- roughly
speaking, it corresponds to families over $\Posf^+$ satisfying an
additional axiom (``semicontinuity'' of
Definition~\ref{pos.pl.def}). This finishes the story for groupoids;
the end result is Proposition~\ref{pos.prop}. In the last part of
Section~\ref{grp.sec}, we give an $I$-augmented version of the
theory. Unfortunately, this needs a rather strong assumption on $I$
--- effectively, it has to be a partially ordered set.

What we really want for general enhanced categories is to consider
bisimplicial sets, that is, to take $I=\Delta$. This is done in
Section~\ref{seg.sp.sec}, and the result here is more restrictive
than in Section~\ref{grp.sec}. We cannot handle general families
over $\Delta^o\Delta^o\Sets$, and we only consider those that
satisfy an additional condition --- namely, a version of the Segal
condition (meaning that the family is represented by a Segal
space). On the other side, we use the categories $\Relf$, $\Relf^+$
of biordered sets of Section~\ref{bi.sec}; we introduce versions of
excision and semicontinuity for families over these categories, and
one additional axiom called ``the cylinder axiom''. The end result
is Proposition~\ref{seg.prop}.

Having done all that, we finish with two more things. Firstly ---
this is Section~\ref{css.sec} --- we show that under an additional
assmption, a ``Segal family'' of groupoids over $\Relf^+$ is
completely determined by its restriction to $\Posf^+ \subset
\Relf^+$, and we characterize families over $\Posf^+$ that appear in
this way. This takes up back to partially ordered sets that we
started with, and it requires the representing Segal space to be
complete in the sense of \cite{rzk}. Secondly, in
Section~\ref{aug.sec}, we study the $I$-augmented situation, with
$I$ again a partially ordered set. The main result here is the
rather non-trivial Proposition~\ref{tw.seg.prop}.

\section{Enhanced groupoids.}\label{grp.sec}

\subsection{Excision and continuity.}\label{exci.subs}

If one wants to find a model category representing all homotopy
types, then there are two standard options: one can use either
topological spaces, or simplicial sets. However, if one is satisfied
with a CW-category in the sense of Subsection~\ref{CW.subs}, then
there are other possibilities. For example, one can use the category
$\Pos$ of partially ordered sets, with the CW-structure of
Example~\ref{pos.CW.exa}.

Recall that as in Subsection~\ref{pos.dim.subs}, we have the notion
of an ample full subcategory $\I \subset \Pos$ given in
Definition~\ref{amp.def}, and the main example is the subcategory
$\Posf \subset \Pos$ of partially ordered sets of finite chain
dimension. For any infinite cardinal $\kappa$, $\Posf_\kappa \subset
\Posf$ is the ample full subcategory spanned by $J \in \Posf$ with
$|J| < \kappa$. Moreover, for any full subcategory $\I \subset
\Pos$, we have the extension $\I^+ \subset \Pos$, ample if so were
$\I$, and in particular, $\Posf^+ \subset \Pos$ is the category of
left-bounded partially ordered sets, and $\Posf^+_\kappa \subset
\Posf^+$ is the full subcategory of $J \in \Posf^+$ with $|J| <
\kappa$. We have the tautological full embedding $i:\Posf \to
\Posf^+$ that induces a full embedding $i:\Posf_\kappa \to
\Posf^+_\kappa$ for any $\kappa$.  Any set $S \in \Sets$ with
discrete order is an object in $\Pos$, we have the map \eqref{J.S}
for any $J \in \Pos$ equipped with a map $J \to S^<$, and if
$J,S^<,S^\hush \in \I$ for some ample $\I \subset \Pos$, then
$J^\hush$ and \eqref{J.S} are also in $\I$ (if $\I = \Posf_\kappa$,
then $S^<,S^\hush \in \I$ iff $|S| < \kappa$). For any $J \in \Pos$
equipped with a map $J \to \N$, we have the map \eqref{J.pl.eq}, and
if $J \in \I^+$ for an ample $\I$, then $J^+$, $J/\N$ and
\eqref{J.pl.eq} are all in $\I^+$.

\begin{defn}\label{pos.def}
Let $\I \subset \Pos$ be an ample full subcategory in the sense of
Definition~\ref{amp.def}, and let $\C$ be a family of groupoids over
$\I$.
\begin{enumerate}
\item The family $\C$ is {\em $[1]$-invariant} if for any $J \in
  \I$, $\C$ is constant along the projection $J \times [1] \to J$.
\item The family $\C$ is {\em semiexact} if for any map $\chi:J \to
  \V = \{0,1\}^<$ in $\I$, the functor
$$
\C_J \to \C_{J/0} \times_{\C_{J_o}} \C_{J/1}
$$
is an epivalence.
\item The family $\C$ is {\em additive} if for any $J \in \I$
  and a map $\chi:J \to S$ to a set $S \in \Sets \cap \I
  \subset \Pos$, the functor
$$
\C_J \to \prod_{s \in S}\C_{J_s}
$$
is an equivalence.
\item The family $\C$ {\em satisfies excision} if for any $J \in
  \I$, $S \in \Sets$ such that $S^<,S^\hush \in \I$, and a map $J
  \to S^<$, $\C$ is constant along the map $J^\hush \to J$ of
  \eqref{J.S}.
\end{enumerate}
\end{defn}

\begin{defn}\label{pos.pl.def}
A family of groupoids $\C$ over the extension $\I^+$ of an ample
subcategory $\I \subset \Pos$ is {\em semicontinuous} if for any $J
\in \I^+$ equipped with a map $J \to \N$, $\C$ is constant along the
map $J^+ \to J$ of \eqref{J.pl.eq}.
\end{defn}

\begin{defn}\label{grp.enh.def}
A {\em small enhanced groupoid} is a small family of group\-oids $\C$
over $\Posf$ that is additive, semiexact, and satisfies
excision. For any regular cardinal $\kappa$, a {\em $\kappa$-bounded
  enhanced groupoid} is an additive semiexact family
$\C/\Posf_\kappa$ satisfying excision such that $\|\C_J\| < \kappa$
for any finite partially ordered set $J$. For any small
resp.\ $\kappa$-bounded enhanced groupoid $\C$, its {\em semicanonical
  extension} is family of groupoids $\C^+$ over $\Posf^+$
resp.\ $\Posf^+_\kappa$ that is $[1]$-invariant, additive and
semiexact in the sense of Definition~\ref{pos.def}, semicontinuous
in the sense of Definition~\ref{pos.pl.def}, and is equipped with an
equivalence $\C \cong i^*\C^+$.
\end{defn}

Recall that for any small category $I$, we have the category
$(\Cat/I)^0$ whose objects are small categories over $I$, and whose
morphisms are functors over $I$ modulo isomorphisms over $I$. For
any regular cardinal $\kappa$, $\Posf_\kappa$ is essentially small,
and we denote by
\begin{equation}\label{grp.S.eq}
\Sets^h_\kappa \subset (\Cat/\Posf_\kappa)^0
\end{equation}
the full subcategory spanned by $\kappa$-bounded enhanced
groupoids. Any small enhanced groupoid gives a $\kappa$-bounded
enhanced groupoid after restriction to $\Posf_\kappa \subset \Posf$
for a sufficiently large $\kappa$. However, since $\Posf$ is not
small, it is not {\em a priori} clear that small enhanced groupoids
also form a category.

\begin{lemma}\label{fin.enh.grp.le}
If $\kappa$ is the smallest infinite cardinal, then any
$\kappa$-bounded enhanced groupoid $\C \to \Posf_\kappa = \Posff$ is
a discrete fibration $\Posff / S \to \Posff$ corresponding to some
finite set $S \in \Sets_\kappa$, so that $\Sets^h_\kappa \cong
\Sets_\kappa$.
\end{lemma}

\proof{} By assumption, $\|\C_J\|$ is finite for any $J \in \Posff$,
so that the groupoid $\C_J$ is rigid, thus discrete by
Example~\ref{fin.rig.exa}. Since epivalences between rigid
categories are equivalences, it is then cartesian along all standard
pushout squares in $\Posff$, and then if we let $S=\C_{\ppt}$, $\C
\to \eps(\ppt)_*S \cong \Posff / S$ is an equivalence over any $J
\in \Posff$, by $[1]$-invariance and induction on $|J|$.
\endproof

\begin{remark}
We use notation $\Sets^h$ in \eqref{grp.S.eq} to emphasize the fact
that in the enhanced context, it is enhanced groupoids that play the
role of enhanced sets (or ``spaces'', whatever that means). By
Lemma~\ref{fin.enh.grp.le}, for finite sets, there is no
enhancement.
\end{remark}

\begin{lemma}\label{ano.grp.le}
  A small enhanced groupoid $\C \to \Posf$ is constant along maps
  anodyne in the sense of Definition~\ref{pos.ano.def}, and a
  semicanonical extension $\C^+$ of $\C$ is constant along maps
  $+$-anodyne in the sense of Definition~\ref{pl.ano.def}. For any
  uncountable regular cardinal $\kappa$, a $\kappa$-bounded enhanced
  groupoid $\C \to \Posf_\kappa$ is constant along maps
  $\kappa$-anodyne in the sense of Definition~\ref{pos.S.ano.def},
  and any of its semicanonical extensions $\C^+$ is constant along
  $\kappa^+$-anodyne maps. Moreover, for any $J \in \Pos^+_\kappa$
  equipped with a map $J \to \N$, the semicanonical extension $\C^+$
  is semicontinuous in the sense of Definition~\ref{count.def} along
  the countable compositions of embeddings $J/n \to J/(n+1)$.
\end{lemma}

\proof{} Assume given a small enhanced groupoid $\C$, and let
$W(\C)$ be the class of maps $f:J' \to J$ in $\Posf$ such that $\C$
is constant along the map $f \times \id:J' \times J_0 \to J \times
J_0$ for any finite $J_0 \in \Posff \subset \Posf$. It suffices to
prove that $W(\C)$ contains all the anodyne maps. Note that for any
$J \to S^<$ of Definition~\ref{pos.def}~\thetag{iv} and any
$J_0 \in \Posff$, we have $(J \times J_0)^\hush \cong J^\hush \times
J_0$, so that $W(\C)$ contains the map $J^\hush \to J$. Then by the last
claim of Proposition~\ref{ano.prop}, it suffices to check that
$W(\C)$ is stable under standard pushouts \eqref{pos.st.sq}. Indeed,
equip $\Posf$ with the standard CW-structure of
Example~\ref{pos.CW.exa}. Then by
Definition~\ref{pos.def}~\thetag{ii}, $\C$ is semiexact with respect
to this CW-structure, and by the same argument as in the proof of
Proposition~\ref{ano.prop}, Definition~\ref{pos.def}~\thetag{iv}
implies Definition~\ref{pos.def}~\thetag{i}, so that $W(\C)$
contains all the projections $J \times [1] \to J$, $J \in
\Posf$. Thus by Lemma~\ref{amp.le}, $W(\C)$ also contains all
reflexive full embeddings. However, by
Example~\ref{adm.CW.exa}~\thetag{i}, the CW-structure on $\Posf$ is
distributive, thus strong in the sense of Definition~\ref{fr.def},
and $1'$ in the corresponding decomposition \eqref{univ.cod} is the
finite partially ordered set $\V^o \cong \{0,1\}^>$, so that $1'
\copr_{1 \copr 1} 1' \cong \VV$ is also finite. Then the functor
$\psi:\wt{\Posf} \cong \Posf_C \to \Posf$ for the corresponding
perfect framing on $\Posf$ is given by $J \mapsto J \times \VV$, and
for any standard pushout square \eqref{pos.st.sq} with $f_0:J \to
J_0$ in $W(\C)$, $\psi(f)$ is also in $W(\C)$. Then to show that the
pushout $J_1 \to J_{01}$ of the map $f_0$ is in $W(\C)$, it suffices
to invoke Lemma~\ref{ano.C.le}. This finishes the proof for small
enhanced groupoids, and the proof for $\kappa$-bounded enhanced
groupoids is exactly the same. For $+$-anodyne and
$\kappa^+$-anodyne maps, replace Proposition~\ref{ano.prop} with
Proposition~\ref{pl.ano.prop}. Finally, for semicontinuity, note
that by Example~\ref{pos.thin.exa} and \eqref{lc.coli.eq}, the set
$J^+$ is a telescope for the countable composition of embeddings
$J/n \to J/(n+1)$ in the sense of Lemma~\ref{tel.I.le}.
\endproof

\begin{corr}\label{X.exci.corr}
Any semicanonical extension $\C^+$ of a small or $\kappa$-bound\-ed
enhanced groupoid $\C$ satisfies excision.
\end{corr}

\proof{} The map \eqref{J.S} is $+$-anodyne, and $\kappa^+$-anodyne
if $|J| < \kappa$.
\endproof

\begin{corr}\label{ano.grp.corr}
  Let $\iota:\Posf \to \Posf$ be the involution sending $J \in
  \Posf$ to $J^o$. Then for any small enhanced groupoid $\C \to
  \Posf$, $\iota^*\C$ is a small enhanced groupoid, and for any
  $\kappa$-bounded enhanced groupoid $\C \to \Posf_\kappa$, so is
  $\iota^*\C$. Moreover, any semicanonical extension $\C^+$ of the
  enhanced groupoid $\C$ is semiexact with respect to co-standard
  pushout squares.
\end{corr}

\proof{} Definition~\ref{pos.def}~\thetag{iii} is manifestly
invariant under $\iota^*$. By virtue of Lemma~\ref{ano.grp.le}, one
can replace Definition~\ref{pos.def}~\thetag{iv} with requiring
$\C$ to be constant along all anodyne maps, and this is also
invariant under $\iota^*$ by Corollary~\ref{B.ano.corr}. Finally,
Definition~\ref{pos.def}~\thetag{ii} literally means that $\C$ is
semiexact over standard pushout squares \eqref{pos.st.sq}, and to
check this for $\iota^*\C$, it suffices to invoke again
Lemma~\ref{ano.grp.le} together with Corollary~\ref{B.ano.corr}, and
note that by Lemma~\ref{B.corr}, $B \circ \iota \cong B$ sends a
standard pushout square to a standard pushout square. For $\C^+$,
while we do not have the involution $\iota$ anymore, the barycentric
subdivision functor $B:\Posf^+ \to \Posf^+$ still sends both
standard and co-standard pushout squares to standard pushout
squares, and the map $\xi:BJ \to J$ is $+$-anodyne for any $J \in
\Posf^+$, and $\kappa^+$-anodyne if $|J| < \kappa$.
\endproof

To obtain examples of small enhanced groupoids in the sense of
Definition~\ref{grp.enh.def}, consider the nerve functor
\begin{equation}\label{N.Q}
\Nn:\Pos \to \Delta^o\Sets, \qquad J \mapsto N(J).
\end{equation}
If we restrict it to $\Posf^+ \subset \Pos$, then for any regular
cardinal $\kappa$, it sends $\Posf^+_\kappa \subset \Posf^+$ into
$\Delta^o\Sets_\kappa \subset \Delta^o\Sets$. Moreover, for any
partially ordered set $J$ of finite chain dimension, its nerve
$N(I)$ is a finite-dimensional simplicial set, so that \eqref{N.Q}
also sends $\Posf$ into $\Delta^o_f\Sets$, and $\Posf_\kappa$ into
$\Delta^o_f\Sets_\kappa$.

\begin{lemma}\label{N.exa.le}
The functor \eqref{N.Q} sends full embeddings to injective maps, and
for any map $\phi:J' \to J$ in $\Pos$, the map
\begin{equation}\label{lc.X.coli.eq}
  \colim_{j \in J}\Nn(J'/j) \to \Nn(J')
\end{equation}
induced by \eqref{lc.coli.eq} is an isomorphism.
\end{lemma}

\proof{} By definition, we have $\Nn(J)([n]) = \Pos([n],J)$ for any
$[n] \in \Delta$ and $J \in \Pos$. Then the first claim is clear,
and the second one immediately follows from Lemma~\ref{hom.I.le}.
\endproof

\begin{lemma}\label{NQ.le}
The functor \eqref{N.Q} commutes with coproducts, sends left-closed
embeddings to injective maps, and sends standard pushout
squares \eqref{pos.st.sq} to pushout squares. Moreover, for any
$+$-anodyne map in the sense of Definition~\ref{pl.ano.def}, $\Nn(f)$
is a weak equivalence, and the same holds for anodyne,
$\kappa$-anodyne and $\kappa^+$-anodyne maps.
\end{lemma}

\proof{} The first claim immediately follows from
Lemma~\ref{N.exa.le}. For the second, note that $\Nn(J \times [1])
\cong \Nn(J) \times \Nn([1]) \cong \Nn(J) \times \Delta_1$, so the
class $\Nn^*W$ of maps $f$ such that $\Nn(f)$ is a weak equivalence
trivially satisfies Definition~\ref{pos.ano.def}~\thetag{i}.
Moreover, since by Lemma~\ref{N.exa.le}, $\Nn$ preserves colimits
\eqref{lc.coli.eq}, it sends functors $J \to \Pos$, $J \in
\Posf^+$ cofibrant in the sense of Definition~\ref{reedy.CW.def} for
the standard CW-structure of Example~\ref{pos.CW.exa} to
Reedy-cofibrant functors $J \to \Delta^o\Sets$, and then by
Lemma~\ref{pos.cof.le}, the functor in the left-hand side of
\eqref{lc.X.coli.eq} is Reedy-cofibrant. Then
Definition~\ref{pos.ano.def}~\thetag{ii} for $\Nn^*W$ reduces to
Lemma~\ref{Q.adj.le}.
\endproof

For any uncountable regular cardinal $\kappa$ and simplicial set $X
\in \Delta^o\Sets_\kappa$, the representable family of groupoids
$\Hh_X \to \Delta^o\Sets_\kappa$ is additive and semiexact by
Lemma~\ref{repr.le}, and it is constant along weak equivalences, so
that $\Nn(f)^*\Hh_X$ is a $\kappa$-bounded enhanced groupoid by
Lemma~\ref{NQ.le}. This gives an ample supply of $\kappa$-bounded
enhanced groupoids. It is also functorial with respect to $X$, so we
obtain a functor
\begin{equation}\label{grp.ka.eq}
h^W(\Delta^o\Sets_\kappa) \to \Sets^h_\kappa, \qquad X \mapsto
\Nn^*\Hh_X,
\end{equation}
where $\Sets^h_\kappa$ is as in \eqref{grp.S.eq}.

\subsection{Comparison.}

We now want to show that the functor \eqref{grp.ka.eq} is an
equivalence for any uncountable regular $\kappa$, and moreover, all
these equivalences together provide an equivalence between
$h^W(\Delta^o\Sets)$ and the category of small enhanced groupoids
and isomorphism classes of functors between them (so that {\em a
  posteriori}, the latter category is well-defined). By virtue of
Theorem~\ref{coho.bis.thm}, homotopy classes of simplicial sets
correspond bijectively to additive semiexact families of groupoids
over $\Delta^o\Sets$, so what we need to do is to find a way to pass
from families over $\Posf$ to families over $\Delta^o\Sets$ that is
inverse to $\Nn^*$. To do this, we use Reedy categories and anodyne
resolutions of Section~\ref{reedy.sec}.

Consider the category $\Delta$ with the Reedy structure of
Example~\ref{del.reedy.exa}, so that in particular, $M$ is the class
of surjective maps, and let $\Delta^M\Delta$ be its simplicial
expansion in the sense of Definition~\ref{simpl.repl.def}. We then
have the fibration $\pi:\Delta^M\Delta \to \Delta$, and the functor
$\xi:\Delta^M\Delta \to \Delta$ of \eqref{xi.eq}. Explicitly,
objects in $\Delta^M\Delta$ are pairs $\langle [m],[n_\idot]
\rangle$ of an object $[m] \in \Delta$ and a functor $[n_\idot]:[m]
\to \Delta$ represented by a diagram
\begin{equation}\label{n.dot}
\begin{CD}
[n_0] @>{f_1}>> \dots @>{f_m}>> [n_m]
\end{CD}
\end{equation}
in $\Delta$. The functors $\pi$, $\xi$ send such a pair to $[m]$
resp.\ $[n_m]$. We also have a morphism
\begin{equation}\label{pi.xi}
b:\pi \to \xi
\end{equation}
sending $l \in [m] = \pi([n_\idot])$ to $\overline{f}_l(n_l) \in [n_m]
= \xi([n_\idot])$, where $\overline{f}_l:[n_l] \to [n_m]$ is the map
in \eqref{n.dot}. Since every surjective map is special in the sense
of Subsection~\ref{spec.subs}, \eqref{pi.xi} is indeed well-defined
on the whole $\Delta^M\Delta$.

For any simplicial set $X$, the category of simplices $\Delta X$
inherits the Reedy structure, and we have $\Delta^M\Delta X \cong
\Delta^M\Delta \times_\Delta \Delta X$, so that objects in
$\Delta^M\Delta X$ are triples $\langle [m],[n_\idot],x \rangle
\rangle$, $\langle [m],[n_\idot] \rangle \in \Delta^M\Delta$, $x \in
X([n_m])$. Sending such a triple to $\langle [m],b^*x \rangle$ then
gives a functor
\begin{equation}\label{beta.eq}
\beta:\Delta^M\Delta X \to \Delta X,
\end{equation}
cartesian over $\Delta$. If $X = NJ$ is the nerve of a partially
ordered set $J$, then $\Delta NJ = \Delta J$ is its simplicial
replacement, and we have
\begin{equation}\label{be.de.eq}
  \beta \cong \Delta(\xi):\Delta^M\Delta J \to \Delta J,
\end{equation}
where $\Delta(\xi)$ is well-defined on $\Delta^M\Delta J$ since the
functor $\xi:\Delta J \to J$ of \eqref{xi.eq} inverts all maps in
$M$.

We now recall that the Reedy category $\Delta X$ is cellular in the
sense of Definition~\ref{ind.reedy.def}. If $X$ is
finite-dimensional, then $\Delta X$ is also finite-dimensional, and
it has the anodyne resolution $Q(\Delta X)$ of
Proposition~\ref{ano.I.prop}. This is functorial in $X$ and has
finite chain dimension, so that we obtain a functor
\begin{equation}\label{Q.N}
\QQ:\Delta^o_f\Sets \to \Posf, \qquad X \mapsto Q(\Delta X).
\end{equation}
For any uncountable regular cardinal $\kappa$, $\QQ$ sends
$\Delta^o_f\Sets_\kappa$ into $\Posf_\kappa$, and it sends finite
simplicial sets to finite partially ordered sets.

If a simplicial set $X$ is regular, $\QQ(X)$ can be described by the
functorial isomorphism \eqref{Q.del.X}. In particular, this applies
to nerves of partially ordered sets, and provides a functorial
isomorphism
\begin{equation}\label{Q.B}
\QQ(\Nn(J)) \cong B(B(J))^o,
\end{equation}
where $B:\Posf \to \Posf$ is the barycentric subdivision functor
of Definition~\ref{pos.bary.def}. Composing the maps $\xi$ and
$\xi_\perp$ of \eqref{B.xi} gives a functorial map
\begin{equation}\label{QN.id}
a(J):\QQ(\Nn(J)) \cong B(B(J))^o \to J
\end{equation}
for any $J \in \Posf$. Equivalently, \eqref{QN.id} is the
composition
\begin{equation}\label{QN.bis.id}
\begin{CD}
\QQ(\Nn(J)) @>{q_{\Delta J}}>> \Delta \Nn(J) = \Delta J @>{\xi}>> J
\end{CD}
\end{equation}
of the functors \eqref{Q.I.eq} and \eqref{xi.eq}. Conversely, for
any $X \in \Delta^o_f\Sets$, the partially ordered set $\QQ(X)$
comes equipped with the functor $q_X = q_{\Delta X}$ of
\eqref{Q.I.eq}, and composing its simplicial replacement
$\Delta(q_X)$ with the functor \eqref{beta.eq}, we obtain a functor
$\beta \circ \Delta(q_X):\Delta \QQ(X) \to \Delta X$ cartesian over
$\Delta$, or equivalently, a map
\begin{equation}\label{NQ.id}
a(X):\Nn(\QQ(X)) \to X.
\end{equation}
If $X = NJ$ for some $J \in \Posf$, then \eqref{be.de.eq}
immediately shows that we have
\begin{equation}\label{a.n.j.eq}
  a(\Nn(J)) = \Nn(a(J)).
\end{equation}
For any map $f:X \to X'$ between simplicial sets, the map
\eqref{q.I.map} induces a possibly non-trivial map $\Delta(q_{X'})
\circ \Delta \QQ(f) \to \Delta^M(\Delta(f)) \circ
\Delta(q_X)$. However, $\beta$ inverts it, so that the map
\eqref{NQ.id} is functorial with respect to $X$.

\begin{lemma}\label{QN.le}
\begin{enumerate}
\item The functor $\QQ$ of \eqref{Q.N} commutes with coproducts,
  sends injective maps to right-closed embeddings, and for any
  cocartesian square
\begin{equation}\label{QN.sq}
\begin{CD}
  X @>{f}>> X'\\
  @V{g}VV @VV{g'}V\\
  Y @>{f'}>> Y'
\end{CD}
\end{equation}
with injective $f$, $g$, we have $\QQ(Y') \cong \QQ(Y)
\copr_{\QQ(X)} \QQ(X')$.
\item The map \eqref{QN.id} is anodyne for any $J \in \Posf$, and
  $\kappa$-anodyne for any uncountable regular cardinal $\kappa$
  such that $J \in \Posf_\kappa$.
\item For any uncountable regular cardinal $\kappa$, $\QQ$ sends
  anodyne resp.\ $\kappa$-anodyne maps in $\Delta^o_f\Sets_\kappa$
  in the sense of Definition~\ref{f.CW.def} to anodyne
  resp.\ $\kappa$-anodyne maps in $\Posf_\kappa$.
\item For any cocartesian square \eqref{QN.sq} with injective $f$,
  the map
\begin{equation}\label{NQ.sq.eq}
  \Nn(\QQ(Y)) \copr_{\Nn(\QQ(X))} \Nn(\QQ(X')) \to \Nn(\QQ(Y'))
\end{equation}
is a weak equivalence.
\end{enumerate}
\end{lemma}

\proof{} For \thetag{i}, note that an injective map $X \to X'$ gives
rise to a fully faithful fibration $\Delta X \to \Delta X'$, thus a
characteristic Reedy functor $\chi:\Delta X' \to [1]$, and if we
have a square \eqref{QN.sq} with injective $f$ and $g$, we have the
joint characteristic Reedy functor $\chi:\Delta(Y') \to \V$. The
claims then follow from Lemma~\ref{reedy.copr.le}. For \thetag{ii},
use Corollary~\ref{B.ano.corr}.

For \thetag{iii}, let $W$ be the class of maps $f$ in
$\Delta^o_f\Sets_\kappa$ such that $\QQ(f)$ is
$\kappa$-anodyne. Then $W$ is saturated and closed under
coproducts. Moreover, since the projection $[n] \to [0]$ is
reflexive, thus $\kappa$-anodyne, $W$ contains all the projections
$\Delta_n = \Nn([n]) \to \Delta_0 = \Nn([0])$ by \thetag{ii}. It
remains to show that $W \cap C$ is closed under pushouts. Indeed,
assume given a cocartesian square \eqref{QN.sq} in
$\Delta^o_f\Sets_\kappa$ with injective $\kappa$-anodyne $f$ and
arbitrary $g$. Then we still have the characteristic Reedy functor
$\chi:\Delta Y' \to [1]$ of the embedding $\Delta(f'):\Delta Y \to
\Delta Y'$, the discrete fibration $\Delta(g'):\Delta X' \to \Delta
Y'$ satisfies the assumptions of Lemma~\ref{Q.pushout.le}, and
Corollary~\ref{V.ano.corr} then proves that $f' \in W$. The proof in
the anodyne case is the same.

Finally, for \thetag{iv}, apply again Lemma~\ref{Q.pushout.le}: if
we let $J = \QQ(Y')$ and $J' = \QQ(X')$, then we have cofibrations
$J,J' \to \V$, and $\QQ(g')$ is a map $J' \to J$ cocartesian over
$\V$ that is an isomorphism over $o,1 \in \V$. Then as in
Corollary~\ref{V.ano.corr}, take the partially ordered set $\wt{J}$,
and note that we have anodyne embeddings $J \subset \wt{J}$, $J_0
\subset \wt{J}/0$, $J'_0 \subset \wt{J}_o$, and an isomorphism
$\wt{J}/1 \cong J'$. By Lemma~\ref{NQ.le}, the horizontal arrows in
the commutative square
$$
\begin{CD}
\Nn(J_0) \copr_{\Nn(J'_0)} \Nn(J') @>>>
\Nn(\wt{J}/0) \copr_{\Nn(\wt{J}_o)} \Nn(\wt{J}/1)\\
@VVV @|\\
\Nn(J) @>>> \Nn(\wt{J})\\
\end{CD}
$$
are weak equivalences, and the vertical arrow on the left is
\eqref{NQ.sq.eq}.
\endproof

\begin{corr}\label{NQ.corr}
For any $X \in \Delta^o_f\Sets$, \eqref{NQ.id} is a weak
equivalence.
\end{corr}

\proof{} Say that $X$ is good if \eqref{NQ.id} is a weak
equivalence. Then by Lemma~\ref{QN.le}~\thetag{ii} and
\eqref{a.n.j.eq}, $\Nn(J)$ is good for any $J \in \Posf$, so that
elementary simplices $\Delta_n$ are good. Moreover, by
Lemma~\ref{QN.le}~\thetag{iv} and Lemma~\ref{Q.adj.le} applied as in
Lemma~\ref{semi.mod.le}, for any cocartesian square \eqref{QN.sq}
with good $X$, $Y$, $X'$ and injective $f$, $Y'$ is also good. But
coproducts of good simplicial sets are obviously good, so by
induction on $\dim X$ and \eqref{sk.I.sq}, any $X \in
\Delta^o_f\Sets$ is good.
\endproof

Extending the functor \eqref{Q.N} to arbitrary simplicial sets is
slightly non-trivial: while by Lemma~\ref{reedy.N.le}, $\QQ$
commutes with countable compositions of injective maps, so that we
have $\QQ(X) = Q(\Delta X) \cong \colim_n\QQ(\sk_nX)$ for any $X \in
\Delta^o\Sets$, the maps $\QQ(\sk_nX) \to \QQ(\sk_{n+1}X)$ are not
left closed --- by Lemma~\ref{QN.le}~\thetag{i}, they are
right-closed --- and $\QQ(X)$ is typically not left-bounded (in
fact, it is right-bounded). To circumvent this, we use the incidence
subset construction of Example~\ref{N.ano.exa}. We define a functor
\begin{equation}\label{Q.pl}
\begin{aligned}
&\QQ_+:\Delta^o\Sets \to \Posf^+,\\
&\QQ_+(X) = \{j \times n|j \in \QQ(\sk_nX)\} \subset \QQ(X) \times
  \N
\end{aligned}
\end{equation}
that sends $\Delta^o\Sets_\kappa$ into $\Posf_\kappa^+$ for any
uncountable regular $\kappa$, and comes equipped with a functorial
map
\begin{equation}\label{Q.pl.Q}
  \QQ_+(X) \to \QQ(X) = \colim_n\QQ(\sk_nX).
\end{equation}
By Example~\ref{N.ano.exa}, \eqref{Q.pl.Q} is $+$-anodyne for $X \in
\Delta^o_f\Sets$, and $\kappa^+$-anodyne for $X \in
\Delta^o_f\Sets_\kappa$. Moreover, combining \eqref{Q.pl.Q} with
\eqref{QN.id} and \eqref{NQ.id} provides functorial maps
\begin{equation}\label{Q.N.pl}
\begin{aligned}
a_+(J)&:\QQ_+(\Nn(J)) \to J, \quad J \in \Posf^+,\\
a_+(X)&:\Nn(\QQ_+(X)) \to X, \quad X \in \Delta^o\Sets,
\end{aligned}
\end{equation}
where $a_+(X)$ is a weak equivalence by Corollary~\ref{NQ.corr} and
Lemma~\ref{Q.adj.le}, and $a_+(J)$ is $+$-anodyne by
Lemma~\ref{QN.le}~\thetag{ii} and Lemma~\ref{N.coli.le}, and
$\kappa^+$-anodyne if $J \in \Posf^+_\kappa$.

\begin{prop}\label{pos.prop}
\begin{enumerate}
\item Assume given a small enhanced groupoid $\C$ in the sense of
  Definition~\ref{grp.enh.def}. Then $\QQ^*\C$ is an additive
  semiexact family over $\Delta^o_f\Sets$ constant along anodyne
  maps of Definition~\ref{f.CW.def}, and the functor $\C \to
  \Nn^*\QQ^*\C$ induced by the map \eqref{QN.id} is an
  equivalence. Moreover, for any semicanonical extension $\C^+$ of $\C$,
  $\QQ_+^*\C^+$ is a semicanonical extension of $\QQ^*\C$ in the sense of
  Definition~\ref{X.ext.def}, and the functor $\C^+ \to
  \Nn^*\QQ_+^*\C^+$ induced by \eqref{Q.N.pl} is an equivalence.
\item Assume given a small family of groupoids $\C$ over
  $\Delta^o_f\Sets$ that is additive, semiexact, and constant along
  anodyne maps. Then $\Nn^*\C$ is a small enhanced groupoid, and the
  functor $\C \to \QQ^*\Nn^*\C$ induced by the map \eqref{NQ.id} is
  an equivalence. Moreover, for any semicanonical extension $\C^+$ of
  $\C$, $\Nn^*\C^+$ is a semicanonical extension of $\Nn^*\C$, and the
  functor $\C \to \QQ_+^*\Nn^*\C$ induced by \eqref{Q.N.pl} is an
  equivalence.
\end{enumerate}
Moreover, for any uncountable regular cardinal $\kappa$, the same
claims hold for $\kappa$-bounded enhanced groupoids and
$\kappa$-bounded families over $\Delta^o_f\Sets_\kappa$ constant
along $\kappa$-anodyne maps.
\end{prop}

\proof{} For \thetag{i}, $\QQ^*\C$ is additive by Lemma~\ref{QN.le},
it is constant along anodyne maps by Lemma~\ref{QN.le} and
Lemma~\ref{ano.grp.le}, and semiexact by Lemma~\ref{QN.le} and
Corollary~\ref{ano.grp.corr} (where the latter implies that $\C$ is
semiexact with respect to co-standard pushout squares in
$\Posf$). To see that $\C \cong \Nn^*\QQ^*\C$, combine
Lemma~\ref{QN.le}~\thetag{ii} and Lemma~\ref{NQ.le}, and to see that
$\QQ_+^*\C$ is a semicanonical extension, note that $\QQ_+$ sends
standard pushout squares to co-standard pushout squares by
Lemma~\ref{QN.le}~\thetag{i} and \eqref{Q.pl}, and by
Lemma~\ref{coli.exa}, it also sends all $+$-anodyne maps of
Definition~\ref{f.pl.CW.def} to $+$-anodyne maps.

For \thetag{ii}, note that by Corollary~\ref{W.CW.corr}, the family
$\C$ is constant along all weak equivalences in
$\Delta^o_f\Sets$. Then $\Nn^*\C$ is a small enhanced groupoid with
semicanonical extension $\Nn^*\C^+$ by Lemma~\ref{NQ.le}, and
\eqref{NQ.id}, \eqref{Q.N.pl} induce equivalences by
Corollary~\ref{NQ.corr} and Lemma~\ref{Q.adj.le}.

For any uncountable regular cardinal $\kappa$, the proofs in the
$\kappa$-bounded case are exactly the same.
\endproof

\begin{corr}\label{grp.X.ext.corr}
Any small resp.\ $\kappa$-bounded enhanced groupoid $\C$ admits a
semicanonical extension $\C^+$, a functor $\gamma:\C_0 \to \C_1$
between small resp.\ $\kappa$-bounded enhanced groupoids $\C_0$,
$\C_1$ with semicanonical extensions $\C_0^+$, $\C_1^+$ extends to a
functor $\gamma^+:\C_0^+ \to \C_1^+$, and a morphism $\gamma_0 \to
\gamma_1$ between two such functors with some extensions
$\gamma_0^+$, $\gamma_1^+$ extends to a (non-unique) morphism
$\gamma_0^+ \to \gamma_1^+$.
\end{corr}

\proof{} Proposition~\ref{pos.prop},
Proposition~\ref{X.ext.prop} and Lemma~\ref{can.le}.
\endproof

\begin{corr}\label{grp.corr}
  For any uncountable regular cardinal $\kappa$, the functor
  \eqref{grp.ka.eq} is an equivalence, and any $\kappa$-bounded
  enhanced group\-oid $\C \in \Sets^h_\kappa$ extends to a small
  enhanced groupoid $\C'/\Posf$, uniquely up to an equivalence
  unique up to an isomorphism. Moreover, any functor $\gamma:\C_0
  \to \C_1$ over $\Posf_\kappa$ between $\kappa$-bounded enhanced
  groupoids extends to a functor over $\Posf$, uniquely up to an
  isomorphism, so that small enhanced groupoids and isomorphism
  classes of functors over $\Posf$ between them form a well-defined
  category
$$
\Sets^h = \bigcup_\kappa \Sets^h_\kappa.
$$
We have $\Sets^h \cong h^W(\Delta^o\Sets)$, and for any regular
cardinal $\kappa$, this equivalence reduces to the equivalence
\eqref{grp.ka.eq}.
\end{corr}

\proof{} Immediately follows from Proposition~\ref{pos.prop},
Theorem~\ref{coho.bis.thm} and
Proposition~\ref{X.ext.prop}. \endproof

\subsection{Augmentations.}\label{aug.subs}

Recall that for any category $I$, we have the categories $\Posf \bb
I \subset \Posf^+ \bb I$ of $I$-augmented finite-dimensional
left-bounded partially ordered sets described in
Subsection~\ref{pos.aug.subs}, with the full subcategories
$\Posf_\kappa \bb I \subset \Posf^+_\kappa \bb I$ for any regular
cardinal $\kappa$. Moreover, assume that $I$ is equipped with a good
filtration in the sense of Definition~\ref{good.def}. Then we also
have the categories $\Posf \bb^\bB I$, $\Posf^+ \bb^\bB I$ of
\eqref{good.eq}, \eqref{good.pl.eq} formed by restricted and locally
restricted $I$-augmented partially ordered sets, and their versions
$\Posf_\kappa \bb^\bB I$, $\Posf^+_\kappa \bb^\bB I$.  We have
the full embedding $i:\Posf \bb^\bB I \to \Posf^+ \bb^\bB I$
that sends $\Posf_\kappa \bb^\bB I$ into $\Posf^+_\kappa \bb^\bB
I$ for any $\kappa$. By our convention, for any map $J'_0 \to J_0$
of partially ordered sets, an $I$-augmented partially ordered set
$\langle J,\alpha_J \rangle$, and a map $J \to J_0$, we augment the
fibered product $J' = J \times_{J_0} J_0'$ via the projection $J'
\to J$. With this convention, we have the augmented comma-sets $J /
j$, $j \in J_0$. Also, for any set $S$ and map $J \to S^<$, the set
$J^+ = S^+ \times_{S^>} J$ is augmented, and the map \eqref{J.S} is
a strict $I$-augmented map, while for any map $J \to \N$, so are the
set $J^+$ and the map $J^+ \to J$ of \eqref{J.pl.eq}.  With this
convention, Definition~\ref{pos.def} and Definition~\ref{pos.pl.def}
make sense for families of groupoids $\C$ over $\Posf \bb I$,
$\Posf^+ \bb I$ and $\Posf_\kappa \bb I$, $\Posf^+_\kappa \bb I$, so
we can introduce the following version of
Definition~\ref{grp.enh.def}.

\begin{defn}\label{grp.I.enh.def}
A {\em small $I$-augmented enhanced groupoid} is a small family of
groupoids $\C$ over $\Posf \bb^\bB I$ that is additive, semiexact,
and satisfies excision in the sense of Definition~\ref{pos.def}. For
a regular cardinal $\kappa$, a {\em $\kappa$-bounded $I$-augmented
  enhanced groupoid} is an additive semiexact family $\C$ over
$\Posf_\kappa \bb^\bB I$ satisfying excision such that $\|\C_J\| <
\kappa$ for any finite $I$-augmented partially ordered set $\langle
J,\alpha_J \rangle$. A {\em semicanonical extension} $\C^+$ of a small
$I$-augmented enhanced groupoid resp.\ a $\kappa$-bounded
$I$-augmented enhanced groupoid $\C$ is a family of groupoids $\C^+$
over $\Posf^+ \bb^\bB I$ resp.\ $\Posf^+_\kappa \bb^\bB I$ that
is $[1]$-invariant, additive, semiexact and semicontinous in the
sense of Definition~\ref{pos.pl.def}, and is equipped with an
equivalence $i^*\C^+ \cong \C$.
\end{defn}

\begin{lemma}\label{ano.I.le}
Any small $I$-augmented enhanced groupoid $\C$ is constant along
maps anodyne in the sense of Definition~\ref{ano.aug.def}, and its
restriction to $\Posf \bbi^\bB I = \Posf \bbi I \cap \Posf \bb^\bB I
\subset \Posf \bb^\bB I$ is semiexact with respect to the
co-standard CW-structure of Example~\ref{pos.I.CW.exa}. Any
semicanonical extension $\C^+$ of $\C$ is constant along $+$-anodyne
maps, and its restriction to $\Posf^+ \bbi^\bB I$ is semiexact with
respect to the co-standard CW-structure. For any uncountable regular
cardinal $\kappa$, the same holds for $\kappa$-bounded augmented
enhanced groupoids and $\kappa$-anodyne resp.\ $\kappa^+$-anodyne
maps.
\end{lemma}

\proof{} For the first claim, use the same argument as in
Lemma~\ref{ano.grp.le}, with Proposition~\ref{ano.prop} replaced by
Proposition~\ref{ano.aug.prop}, and with the perfect framing on
$\Posf \bb^\bB I$ resp.\ $\Posf_\kappa \bb^\bB I$ provided by
Example~\ref{adm.CW.I.exa}. For the second claim, as in
Corollary~\ref{ano.grp.corr}, note that the barycentric subdivision
functor $B_I$ of Definition~\ref{bary.I.def} sends co-standard
pushout squares to standard pushout squares, and the map $\xi:B_I(J)
\to J$ is anodyne resp.\ $+$-anodyne for any $J \in \Posf \bb I$
resp.\ $\Posf^+ \bb I$, and $\kappa$-anodyne
resp.\ $\kappa^+$-anodyne if $|J| < \kappa$.
\endproof

\begin{corr}\label{X.I.exci.corr}
A semicanonical extension $\C^+$ of a small or $\kappa$-bounded
$I$-augmented enhanced groupoid $\C$ satisfies excision.
\end{corr}

\proof{} Same as Corollary~\ref{X.exci.corr}.
\endproof

Now assume that $I$ is a cellular Reedy category in the sense of
Definition~\ref{ind.reedy.def}, with the good filtration given by
the degree filtration. Then for any $I$-simplicial set $X \in
I^o\Delta^o\Sets$, we have the corresponding discrete fibration
$I\Delta X \to I \times \Delta$, and $I\Delta X$ is a cellular Reedy
category by Lemma~\ref{ind.reedy.le}, finite-dimensional if so is
$X$. In this case, we have the anodyne resolution $Q(I\Delta X)$
equipped with the projection $q_{I\Delta X}$ of
\eqref{Q.I.eq}. Composing $q_{I\Delta X}$ with the projection
$I\Delta X \to I$ turns $Q(I\Delta X)$ into a restricted
finite-dimensional $I^o$-augmented partially ordered set. For any
map $f:X \to X'$ in $I^o_f\Delta^o_f\Sets$, we have the
corresponding map $Q(f):Q(I\Delta X) \to Q(I\Delta X')$, and
together with the maps \eqref{q.I.map}, this defines a map in $\Posf
\bb^\bB I$, so that we have a functor
\begin{equation}\label{Q.N.I}
\QQ_I:I^o_f\Delta^o_f\Sets \to \Posf \bb^\bB I, \qquad X \mapsto
Q(I\Delta X).
\end{equation}
As in \eqref{Q.pl}, we extend it to $I^o\Delta^o\Sets$ by defining a
functor
\begin{equation}\label{Q.I.pl}
\begin{aligned}
&\QQ_{I+}:I^o\Delta^o\Sets \to \Posf^+ \bb^\bB I,\\
&\QQ_{I+}(X) = \{j \times n|j \in \QQ_I(\sk_nX)\} \subset \QQ_I(X) \times
  \N,
\end{aligned}
\end{equation}
where $\sk_\idot X$ is the skeleton filtration for the cellular
Reedy category $I \times \Delta$, and we have a map
\begin{equation}\label{Q.I.pl.Q}
  \QQ_{I+}(X) \to \QQ_I(X) = \colim_n\QQ_I(\sk_nX)
\end{equation}
that is $+$-anodyne for $X \in I^o_f\Delta^o_f\Sets$. Conversely,
applying the functor \eqref{N.Q} pointwise, as in \eqref{N.I.eq}, we
obtain a functor
\begin{equation}\label{N.Q.I}
\begin{aligned}
\Nn_I:\Posf^+ \bb^\bB I &\to I^o\Delta^o\Sets,\\
\Nn_I(J)(i \times [n]) &= \Nn(\sE_I(J)(i))([n]) = \Nn(i \setminus
J)([n])
\end{aligned}
\end{equation}
for any $i \in I$ and $[n] \in \Delta$, where $\sE_I$ is the functor
\eqref{E.eq}. If $|I| < \kappa$ for some uncountable regular
cardinal $\kappa$, then \eqref{Q.N.I}, \eqref{Q.I.pl} and
\eqref{N.Q.I} restrict to functors between
$I^o_f\Delta^o_f\Sets_\kappa$, $I^o\Delta^o\Sets_\kappa$ and
$\Posf_\kappa \bb^\bB I$, $\Posf^+_\kappa \bb^\bB I$. Moreover,
for any $X \in I^o_f\Delta^o_f\Sets$, the projection $I\Delta X \to
I$ is a fibration, so that for any $i \in I$, the projection $q_{I
  \Delta X}$ induces a functor
\begin{equation}\label{q.i}
q_i:i \setminus \QQ_I(X) \to \Delta X(i),
\end{equation}
where $X(i) \in \Delta^o_f\Sets$ is given by $X(i)([n]) = X(i \times
[n])$. As in the absolute case, $q_i$ gives rise to a map $\Nn(i
\setminus \QQ_I(X)) \to X(i)$, and these maps are functorial in $i$,
thus glue together to a map
\begin{equation}\label{NQI.id}
a(X):\Nn_I(\QQ_I(X)) \to X.
\end{equation}
For the same reasons as in the absolute case, this map is functorial
with respect to $X$. Conversely, for any $I$-augmented partially
ordered set $\langle J,\alpha \rangle$, we have the identification
\begin{equation}\label{NJ.eq}
  I\Delta \Nn_I(J) \cong I \setminus^o_{\alpha^o \circ \xi_\perp}
  \Delta J,
\end{equation}
where the twisted comma-category of Example~\ref{tw.comma.exa} is
taken with respect to the composition
$$
\begin{CD}
(\Delta J)^o \cong \Delta_\perp J @>{\xi_\perp}>> J^o @>{\alpha^o}>>
  I^o
\end{CD}
$$
of the functor \eqref{xi.perp.eq} and the functor $\alpha^o:J^o \to
I^o$. Then the composition
$$
\begin{CD}
\QQ_I(\Nn_I(J)) @>{q_{I\Delta \Nn_I(J)}}>> I\Delta \Nn_I(J) \cong
I \setminus^o_{\alpha^o \circ \xi_\perp} \Delta J @>{\tau}>> \Delta
J @>{\xi}>> J
\end{CD}
$$
gives an $I$-augmented map
\begin{equation}\label{QNI.id}
a(J):\QQ_I(\Nn_I(J)) \to J
\end{equation}
functorial with respect to $\langle J,\alpha \rangle \in \Posf \bb
I$, and we have
\begin{equation}\label{a.n.j.I.eq}
  a(\Nn_I(J)) = \Nn_I(a(J)),
\end{equation}
an augmented version of \eqref{a.n.j.eq}. We also obtain
functorial maps
\begin{equation}\label{NQI.id.pl}
\begin{aligned}
a_+(X)&:\Nn_I(\QQ_{I+}(X)) \to X, \qquad X \in I^o\Delta^o\Sets\\
a_+(J)&:\QQ_{I+}(\Nn_I(J)) \to J, \qquad J \in \Posf^+ \bb I
\end{aligned}
\end{equation}
by composing \eqref{NQI.id} and \eqref{QNI.id} with
\eqref{Q.I.pl.Q}.

\begin{lemma}\label{QNI.le}
\begin{enumerate}
\item The functor \eqref{N.Q.I} sends restricted $I$-augmented
  partially ordered sets to finite-dimensional $I$-simplicial
  sets. Moreover, it commutes with coproducts, sends left-closed
  embeddings to injective maps, and standard pushout squares in
  $\Posf \bb I$ to pushout squares in $I^o_f\Delta^o_f\Sets$. For
  any uncountable regular cardinal $\kappa > |I|$ and any map $f$ in
  $\Posf_\kappa \bb^\bB I$ that is $\kappa$-anodyne in the sense of
  Definition~\ref{ano.aug.def}, $\Nn_I(f)$ is a weak equivalence,
  and the same holds for any map $f$ in $\Posf_\kappa^+ \bb^\bB I$
  that is $\kappa^+$-anodyne.
\item For any uncountable regular $\kappa > |I|$ and any $\langle J,\alpha
  \rangle \in \Posf_\kappa \bb^\bB I$ such that $\alpha:J \to I$
  factors through $I_L \subset I$, the map \eqref{QNI.id} is
  $\kappa$-anodyne, and for any $\langle J,\alpha \rangle \in
  \Posf_\kappa^+ \bb^\bB I$ satisfying the same assumption, the
  map $a_+(J)$ of \eqref{NQI.id.pl} is $\kappa^+$-anodyne.
\item The functor \eqref{Q.N.I} commutes with coproducts, sends
  injective maps to right-closed embeddings, and standard pushout
  squares in $I^o_f\Delta^o_f\Sets$ to co-standard pushout squares
  in $\Posf \bb^\bB I$. Moreover, for any uncountable regular
  $\kappa > |I|$, it sends $\kappa$-anodyne maps of
  Definition~\ref{f.CW.def} to $\kappa$-anodyne maps. The functor
  \eqref{Q.I.pl} commutes with coproducts, sends injective maps to
  right-closed embeddings, $\kappa^+$-ano\-dyne maps of
  Definition~\ref{f.pl.CW.def} to $\kappa^+$-anodyne maps, and
  standard pushout squares in $I^o\Delta^o\Sets$ to co-standard
  pushout squares in $\Posf^+ \bb^\bB I$.
\item For any $X \in I^o_f\Delta^o_f\Sets$, the map \eqref{NQI.id} is a
  weak equivalence, and for any $X \in I^o\Delta^o\Sets$, so is the
  map $a_+(X)$ of \eqref{NQI.id.pl}.
\end{enumerate}
\end{lemma}

\proof{} For \thetag{i}, note that for any $\langle J,\alpha \rangle
\in \Posf \bb^\flat I$, a simplex in $\Nn_I(J)(i)$ is explicitly
given by a map $f:[n] \to J$ and a map $g:i \to \alpha(f(0))$, and
if the simplex is non-degenerate, then $f$ is injective, so that $n
\leq \dim J$, and $g$ is in $L$, so that $\deg(i) \leq
\deg(\alpha(f(0)))$, and the latter is bounded by a constant since
$\langle J,\alpha \rangle$ is restricted. This proves the first
claim, and the rest of the argument is exactly the same as in
Lemma~\ref{NQ.le} (combined with Lemma~\ref{Q.adj.le} in the
$+$-anodyne case).

For \thetag{ii}, note that under our assumptions, the twisted
comma-category of \eqref{NJ.eq} is regular by Lemma~\ref{reg.I.le},
and then $\QQ_I(\Nn_I(J))$ is given by \eqref{reg.M.eq}. Explicitly,
this reads as
\begin{equation}\label{Q.I.N.I.B}
\QQ_I(\Nn_I(J)) \cong B^o_I(I \setminus^o_{\xi_\perp} B_I(J)),
\end{equation}
where the twisted comma-category $I \setminus^o_{\xi_\perp} B_I(J)$
is fibered over $B_I(J)$ with fibers $I_L/j$, $j \in J$ that are
partially ordered sets by Definition~\ref{ind.reedy.def}~\thetag{i},
so that $I \setminus^o_{\xi_\perp} B_I(J)$ is a partially ordered
set. Then the map \eqref{QNI.id} factors as
$$
\begin{CD}
B^o_I(I \setminus^o_{\xi_\perp} B_I(J)) @>{\xi_\perp}>> I
\setminus^o_{\xi_\perp} B_I(J) @>{\tau}>> B_I(J) @>{\xi}>> J,
\end{CD}
$$
the map $\xi$ is $\kappa$-anodyne, the map $\xi_\perp$ is
$\kappa$-anodyne by Lemma~\ref{B.I.ano.le}, and the fibration $\tau$
is $\kappa$-anodyne by Lemma~\ref{B.I.fib.le}~\thetag{ii}. This
finishes the proof for \eqref{QNI.id}; to deduce the statement for
\eqref{NQI.id.pl}, use Lemma~\ref{N.coli.le}.

For \thetag{iii}, for \eqref{Q.N.I}, use the same argument as in
Lemma~\ref{QN.le} to prove all the claims but the last, and reduce
the last one to the case when $f$ is an elementary simplex
projection $\Delta_{i,n} \to \Delta_{i,0}$ for some $i \in I$ and $n
\geq 0$. Then for any integer $m \geq 0$, we have $\Delta_{i,m} =
\Nn_I([m]_i)$, where $[m]_i \in \Posf \bb^\flat I$ is $[m]$
augmented via the constant functor $[m] \to I$ with value $i$, and
$[m]_i$ trivially satisfies the assumptions of \thetag{ii}.  Since
the projection $[n]_i \to [0]_i$ is obviously $\kappa$-anodyne for
any $\kappa$, we are done by \thetag{ii}. To deduce the claim for
\eqref{Q.I.pl}, use Lemma~\ref{reedy.N.le}.

Finally, for \thetag{iv}, the same argument as in the proof of
Corollary~\ref{NQ.corr} reduces us to the case $X = \Delta_{i,n} =
\Nn_I([n]_i)$, $i \in I$, $n \geq 0$, and then since $[n]_i$
satisfies the assumptions of \thetag{ii}, we are done by
\eqref{a.n.j.I.eq}, \thetag{ii}, \thetag{i}, and
Lemma~\ref{Q.adj.le} for $a_+(X)$.
\endproof

\begin{lemma}\label{a.Q.le}
For any $J \in \Posf \bb I$ and left-closed $J_0 \subset J$, the
embedding $\QQ_I(\Nn_I(J_0)) \to a(J)^{-1}(J_0) \subset
\QQ_I(\Nn_I(J))$ is left-reflexive in $\Posf \bb I$.
\end{lemma}

\proof{} Let $\chi:I\Delta\Nn_I(J) \to [1]$ be characteristic Reedy
functor of the left-closed embedding $I\Delta\Nn_I(J_0) \to
I\Delta\Nn_I(J)$, and let $Q(\chi):\QQ_I(\Nn_I(J)) \to \V$ be the
induced cofibration of Lemma~\ref{reedy.copr.le}. Then
$\QQ_I(\Nn_I(J_0)) \cong \QQ_I(\Nn_I(J))_0$, and $a(J)^{-1}(J_0)
\cong \QQ_I(\Nn_I(J))/0$.
\endproof

Lemma~\ref{QNI.le}~\thetag{i} is an $I$-augmented version of
Lemma~\ref{NQ.le}, and just as in the absolute case, it immediately
shows that for any cellular $I$ with $|I| < \kappa$ and
$I$-simplicial set $X \in I^o\Delta^o\Sets_\kappa$, $\Nn_I^*\Hh_X$
is a $\kappa$-bounded $I$-augmented enhanced groupoid, so that as in
\eqref{grp.ka.eq}, we have a functor
\begin{equation}\label{grp.ka.I.eq}
h^W(I^o\Delta^o\Sets_\kappa) \to \Sets^h_\kappa(I), \qquad X \mapsto
\Nn_I^*\Hh_X,
\end{equation}
where $\Sets^h_\kappa(I)$ stands for the category of $\kappa$-bounded
$I$-augmented enhanced groupoids and isomorphism classes of functors
over $\Posf_\kappa \bb^\bB I$. Ideally, we would like to show that
\eqref{grp.ka.I.eq} is an equivalence, just like in the absolute
case, but this is not possible because of the additional assumption
we have to impose in Lemma~\ref{QNI.le}~\thetag{ii}. This needs to
be taken care of by imposing an assumption on the cellular Reedy
category $I$.

\begin{prop}\label{pos.I.prop}
  Assume that a cellular Reedy category $I$ is directed (that is, $I
  = I_L$) and $\Hom$-finite. Then for any uncountable regular
  cardinal $\kappa$ such that $|I| < \kappa$, the functor
  \eqref{grp.ka.I.eq} is an equivalence of categories. Moreover, in
  this case, for any regular cardinal $\kappa' \geq \kappa$, any
  $\kappa'$-bounded $I$-augmented enhanced groupoid $\C_\kappa$
  extends to a small $I$-augmented enhanced groupoid $\C$, any functor
  $\gamma_{\kappa'}:\C_{\kappa'} \to \C'_{\kappa'}$ over
  $\Posf_{\kappa'} \bb^\bB I$ extends to a functor $\gamma:\C \to
  \C'$ over $\Posf \bb^\bB I$, and a morphism $\gamma_{\kappa'}
  \to \gamma'_{\kappa'}$ extends to a (non-unique) morphism $\gamma
  \to \gamma$. In particular, $\C$ is unique up to an equivalence
  unique up to an isomorphism, so that small $I$-augmented enhanced
  groupoids and isomorphisms classes of functors over $\Posf
  \bb^\bB I$ form a well-defined category
$$
\Sets^h(I) \cong \bigcup_{\kappa' \geq \kappa} \Sets^h_{\kappa'}(I),
$$
and the functors \eqref{grp.ka.I.eq} give an equivalence
$h^W(I^o\Delta^o\Sets) \cong \Sets^h(I)$. Finally, any small
$I$-augmented enhanced groupoid $\C$ admits a semicanonical extension
$\C^+$, any functor $\gamma:\C_0 \to \C_1$ between small
$I$-augmented enhanced groupoids extends to a functor
$\gamma^+:\C_0^+ \to \C_1^+$, and any morphism $\gamma_0 \to
\gamma_1$ between such functors extends to a (non-unique) morphism
$\gamma_0^+ \to \gamma_1^+$.
\end{prop}

\proof{} The assumptions on $I$ contain the
assumptions of Theorem~\ref{coho.bisi.thm}, and since $I=I_L$, every
$J \in \Posf_\kappa \bb I$ trivially satisfies the assumptions of
Lemma~\ref{QNI.le}~\thetag{iv} since $I=I_L$. As in
Proposition~\ref{pos.prop}, Corollary~\ref{grp.X.ext.corr} and
Corollary~\ref{grp.corr}, everything then directly follows from
Lemma~\ref{ano.I.le}, Lemma~\ref{QNI.le},
Theorem~\ref{coho.bisi.thm} and Proposition~\ref{X.ext.prop}.
\endproof

\begin{remark}
Under the assumptions of Proposition~\ref{pos.I.prop}, all
$I$-aug\-mented left-bounded partially ordered sets are automatically
locally restricted, so that $\Posf^+ \bb^\bB I = \Posf^+ \bb
I$. If $I$ is finite-dimensional in the sense of
Definition~\ref{I.dim.def}, then we also have $\Posf \bb^\bB I =
\Posf \bb I$.
\end{remark}

\section{Segal spaces.}\label{seg.sp.sec}

\subsection{Segal families.}\label{del.seg.subs}

While it might be possible to generalize
Proposition~\ref{pos.I.prop} to arbitrary cellular Reedy categories,
it is not easy. Therefore we now concentrate on one specific
category, and only consider families of groupoids satisfying a
certain condition. The category is $\Delta$ itself, and
$\Delta$-simplicial sets $X \in \Delta^o\Delta^o\Sets$ are usually
called {\em bisimplicial sets}. We keep this usage but we need to
distinguish between the two factors in the product $\Delta \times
\Delta$. We will call the first one the {\em Reedy factor}, since
$\Delta$ here plays the role of a cellular Reedy category $I$, and
we will call the second one the {\em simplicial factor}. The Reedy
factor will always come first. For any $n,m \geq 0$, we will denote
by $\Delta_{n,m} = \Delta_{[n],m} = \Delta_n \boxtimes \Delta_m$ the
corresponding elementary bisimplex, where per our convention, the
first integer $n$ corresponds to the Reedy factor, and the second
integer $m$ corresponds to the simplicial one. We also have two
functors
\begin{equation}\label{LR.simp}
  L = l^*,R = r^*:\Delta^o\Sets \to \Delta^o\Delta^o\Sets,
\end{equation}
where $l,r:\Delta^o \times \Delta^o \to \Delta^o$ are the
projections onto the Reedy resp.\ simplicial factor, and for any
$X,Y \in \Delta^o\Sets$, we denote $X \boxtimes Y = L(X) \times
R(Y)$. As in Example~\ref{holim.I.exa}, the functor $R$ sends maps in
the class $C$ resp.\ $W$ to maps in $C$ resp.\ $W$, but specifically
in the case $I=\Delta$, it also sends maps in $F$ to maps in $F$, so
the induced functor $R^h:h^W(\Delta^o\Sets) \to
h^W(\Delta^o\Delta^o\Sets)$ has both a left and a right-adjoint
\begin{equation}\label{colim.del.eq}
R^h_!=\hocolim_{\Delta^o},R^h_*=\holim_{\Delta^o}:h^W(\Delta^o\Delta^o\Sets)
\to h^W(\Delta^o\Sets)
\end{equation}
provided by the Quillen Adjunction Theorem~\ref{qui.adj.thm} (where
both $\holim_{\Delta^o}$ and $\hocolim_{\Delta^o}$ in
\eqref{colim.del.eq} are taken in the Reedy direction). Since $\Delta
\times \Delta$ is a cellular Reedy category, any bisimplicial set
$X$ carries the skeleton filtration $\sk_\idot X$ of \eqref{sk.X},
and we can also take skeletons only along the Reedy
resp.\ simplicial factor; we denote those by $\sk_\idot^L X$
resp.\ $\sk^R_\idot X$. For any bisimplicial set $X$, we let
$X^\iota = (\iota \times \id)^*X$, where the involution
$\iota:\Delta \to \Delta$, $[n] \mapsto [n]^0$ acts via the Reedy
factor.

Our condition on families is a version of the Segal condition of
Subsection~\ref{seg.X.subs}. Recall that for any $n \geq l \geq 0$,
the cocartesian square \eqref{seg.del.sq} induces a map
\eqref{b.ln.eq} of simplicial sets. If $l=0$ or $l=n$, this map is
an isomorphism, but if $n > l > 0$, it is not.

\begin{defn}\label{bi.seg.def}
A family of groupoids $\C$ over the category
$\Delta^o_f\Delta^o_f\Sets$ or its full subcategory
$\Delta^o_f\Delta^o_f\Sets_\kappa \subset \Delta^o_f\Delta^o_f\Sets$
for some regular cardinal $\kappa$ is a {\em Segal family} iff for
any integer $n > l > 0$, with the corresponding map $b_n^l$ of
\eqref{b.ln.eq}, the family $\C$ is stably constant along the map
$L(b_n^l)$ in the sense of Definition~\ref{st.I.def}.
\end{defn}

We recall that explicitly, ``stably constant'' in
Definition~\ref{bi.seg.def} means that $\C$ is constant along
$b_n^l \times \id_{R(Y)}$ for any finite simplicial set $Y \in
\Delta^o\Sets$.

\begin{exa}\label{seg.sp.exa}
For any model category $\C$, Reedy-fibrant simplicial object $X \in
\Delta^o\C$, and integers $n > l > 0$, the map
\begin{equation}\label{seg.F}
X([n]) \to X([l]) \times_{X([0])} X([n-l])
\end{equation}
induced by \eqref{seg.del.sq} is in $F$, and $X$ is a {\em homotopy
  Segal object} if all the maps \eqref{seg.F} are also in $W$. Say
that a bisimplicial set $X$ is a {\em Segal space} if all the maps
\eqref{seg.F} are trivial Kan fibrations. Then a fibrant $X$ is a
Segal space iff the corresponding represented family $\Hh_X$ is a
Segal family in the sense of Definition~\ref{bi.seg.def}. We note
that since the involution $\iota:\Delta \to \Delta$ sends squares
\eqref{seg.del.sq} to squares of the same form, a bisimplicial set
$X$ is a Segal space if and only if so is $X^\iota$.
\end{exa}

\begin{remark}
While we do allow Segal spaces $X \in \Delta^o\Delta^o\Sets$ not to
be fibrant, this notion in such a generality should be taken with a
grain of salt. Among other things, it is not invariant under weak
equivalences --- in fact, as we shall see, {\em any} bisimiplicial
set is weakly equivalent to a Segal space.
\end{remark}

\begin{exa}
If $\C \to \Delta^o_f\Delta^o_f\Sets$ is exact, so that as in
Remark~\ref{exa.rem}, $\C \cong \Y_*\Y^*\C$ is completely defined by
its restriction $Y^*\C \to \Delta \times \Delta$, and constant along
weak equivalences, so that $\Y^*\C$ is constant along the simplicial
factor, then $\C$ is a Segal family if and only if $\Y^*\C$ is a
Segal category in the sense of Definition~\ref{seg.cat.def}.
\end{exa}

We will need an equivalent reformulation of
Definition~\ref{bi.seg.def} in the spirit of
Lemma~\ref{bva.le}. Firstly, for any $m \geq 0$ and $n \geq l \geq
0$, we have the natural embedding
\begin{equation}\label{bi.horn.bis.eq}
v^l_{n,m}:\V_{n,m}^l \hookrightarrow \Delta_{n,m},
\end{equation}
where as in \eqref{s.v.eq}, we denote
\begin{equation}\label{s.v.bis.eq}
\V_{n,m}^l = (\V^l_n \boxtimes \Delta_m) \copr_{\V^l_n \boxtimes
  \SS_m} (\Delta_n \boxtimes \SS_m) \subset \Delta_{n,m} = \Delta_n
\boxtimes \Delta_m.
\end{equation}
The map \eqref{bi.horn.bis.eq} is essentially the horn embedding
\eqref{bi.horn.eq} for $I = \Delta$, but along the Reedy factor
rather than the simplicial one. Secondly, for any map $\chi:J \to
J'$ of partially ordered sets, we have the map
\eqref{a.chi.eq}. Moreover, if $J \in \Delta^o\Pos$ is a simplicial
partially ordered set, then we can take its nerve pointwise and
obtain functors
\begin{equation}\label{bi.nerve}
N^L,N^R:\Delta^o\Pos \to \Delta^o\Delta^o\Sets,
\end{equation}
where the nerve is taken along the Reedy resp.\ simplicial
factor. For any $J:\Delta^o \to \Pos$, we denote $J^\iota = \iota
\circ J$, so that $J^\iota([n]) = J([n])^o$, $[n] \in \Delta$, and
then we tautologically have $N^L(J)^\iota \cong N^L(J^\iota)$. For
any $J \in \Delta^o\Pos$, the functor \eqref{xi.eq} induces a
functor
\begin{equation}\label{J.ora}
  \Delta\Delta N^R(J) \to \overrightarrow{\Delta}J
\end{equation}
whose target is the category \eqref{X.ora} corresponding to
$J:\Delta^o \to \Pos$, and similarly for $N^L$. Now, if $J$ is
equipped with a map $\chi:J \to J'$ to a partially ordered set $J'$
understood as a constant simplicial partially ordered set, the map
\eqref{a.chi.eq} induces a map
\begin{equation}\label{a.chi.bis.eq}
a^\dm_\chi:\colim_{\Tw(J')}N^L((J/J')_f) \to N^L(J/J')
\end{equation}
of bisimplicial sets. For any regular cardinal $\kappa$, the functors
\eqref{bi.nerve} sends $\Delta^o\Pos_\kappa$ into
$\Delta^o\Delta^o\Sets_\kappa$ but controlling the dimension
needs more care.

\begin{defn}\label{glob.def}
A simplicial partially ordered set $J$ is {\em globally
  finite-dimensional} if its underlying simplicial set is
finite-dimensional, and we have $\dim J([m]) \leq d$ for all $[m]
\in \Delta$ and some fixed constant $d \geq 0$ independent of
$[m]$. Moreover, $J \in \Delta^o\Pos$ is {\em globally left-bounded}
if it admits a height function $\hht:J \to \N$ with globally
finite-dimensional comma-fibers $J /_{\hht} n$, $n \in \N$, where
$\N$ is considered as a constant simplicial partially ordered set,
and $J$ is {\em globally right-bounded} if $J^\iota$ is globally
left-bounded, so that $J$ admits a co-height function $\hht^\iota:J
\to \N^o$ with globally finite-dimensional comma-fibers $n
\setminus_{\hht^\iota} J$, $n \in \N^o$.
\end{defn}

In particular, if $J \in \Delta^o\Pos$ is globally
finite-dimensional, $J$ factors through $\Posf$, and we denote by
$\Delta^o_{fd}\Posf \subset \Delta^o_f\Posf$ the full subcategory
spanned by such $J$. If $J$ is globally left-bounded, it factors
through $\Delta^o\Posf^+$, and we let $\Delta^o_{d+}\Posf^+ \subset
\Delta^o\Posf^+$ be the full subcategory spanned by such $J$, and
dually, globally right-bounded $J$ factors through $\Posf^-$, and we
let $\Delta^o_{d-}\Posf^- \subset \Delta^o\Posf^-$ be the full
subcategory spanned by such $J$. Note that all inclusions is strict.

\begin{exa}
For any partially ordered set $J$, the simplicial expansion
$\Delta^\Tot J$ with respect to the class of all maps corresponds to
a simplicial partially ordered set $N^\Tot(J)$ by the Grothendieck
construction, and the underlying simplicial set of $N^\Tot(J)$ is
the nerve $N(J)$. If $\dim J < \infty$, then $N(J)$ is also
finite-dimensional, and $N^\Tot(J)([m]) = \Posf([m],J)$ is
finite-dimensional for any $[m]$, so that $N^\Tot(J) \in
\Delta^o_f\Posf$. However, $N^\Tot(J) \in \Delta^o_{fd}\Posf$ only
if $\dim J = 0$. Already for $J = [1]$, we have $N^\Tot(J)([m])
\cong [m+1]$, and $N^\Tot([1])$ is not even left or right-bounded:
as soon as a subset $J \subset N^\Tot([1])$ contains
$N^\Tot([1])([1]) \cong [2]$, it contains the whole $N^\Tot([1])$.
\end{exa}

\begin{lemma}\label{binerve.le}
For any $J \in \Delta^o_f\Pos$ and integer $d \geq 0$, both
$\sk^L_dN^L(J)$ and $\sk^R_dN^R(J)$ are in
$\Delta^o_f\Delta^o_f\Sets$, and if $J$ is globally finite
dimensional, then $N^L(J) = \sk^L_dN^L(J)$, $N^R(J) = \sk^R_dN^R(J)$
for some $d \geq 0$, so that the nerve functors \eqref{bi.nerve}
send $\Delta^o_{fd}\Posf$ into $\Delta^o_f\Delta^o_f\Sets$.
\end{lemma}

\proof{} The proofs for $N^L$ and $N^R$ are identical, so let us
only consider $N^L$. For any $d \geq 0$, let $i_d:\Delta \to \Delta
\times \Delta$ be the embedding onto $[d] \times \Delta$. Then
$i_d^*N^L(J) \cong J^{[d]}$ is a subset of $J^d$, with the embedding
$J^{[d]} \to J^d$ given by restriction to the discrete set
$\{0,\dots,d\} \subset [d]$. Therefore
$$
\dim \sk^L_dN^L(J) = \max\{p \leq d|\dim i_p^*N^L(J)\} \leq d \dim
U(J),
$$
where $U(J) \in \Delta^o_f\Sets$ is the underlying simplicial set of
$J$. This proves the first claim. For the second, note that if $\dim
J([m]) \leq d$ for all $[m] \in \Delta$ and some fixed $d \geq 0$,
then $N^L(J) \cong \sk^L_dN^L(J)$.
\endproof

\begin{prop}\label{bi.seg.prop}
Assume given an additive semiexact family of group\-oids $\C$ over
$\Delta^o_f\Delta^o_f\Sets$ or $\Delta^o_f\Delta^o_f\Sets_\kappa
\subset \Delta^o_f\Delta^o_f\Sets$ constant along anodyne
resp.\ $\kappa$-anodyne maps of Definition~\ref{f.CW.def}. Then the
following conditions are equivalent.
\begin{enumerate}
\item The family $\C$ is a Segal family in the sense of
  Definition~\ref{bi.seg.def}.
\item For any $n > l > 0$, $m \geq 0$, $\C$ is stably constant along the
  map \eqref{bi.horn.bis.eq}.
\item For any $J' \in \Posf$ resp.\ $\Posf_\kappa$ and $J \in
  \Delta^o_{fd}\Posf$ resp.\ $\Delta^o_{fd}\Posf_\kappa$ equipped
  with a map $\chi:J \to J'$, $\C$ is constant along the map
  $a^\dm_\chi$ of \eqref{a.chi.bis.eq}.
\item For any $J \in \Posf$ resp.\ $\Posf_\kappa$ and a map $\chi:J
  \to [1]$, $\C$ is stably constant along the map $L(a_\chi)$, where
  $a_\chi$ is the map \eqref{a.chi.eq}.
\end{enumerate}
\end{prop}

\proof{} Let $W(\C)$ be the class of injective maps $w$ in
$\Delta^o_f\Delta^o\Sets^f$
resp.\ $\Delta^o_f\Delta^o\Sets^f_\kappa$ such that $\C$ is stably
constant along $w$. Then $W(\C)$ is saturated, and closed under
coproducts and pushouts by Lemma~\ref{f.CW.le}. The same induction
as in Lemma~\ref{bva.le} then shows that if $\C$ is a Segal family,
it is stably constant along the horn extensions $L(v^l_n)$, $n > l >
0$ in the Reedy direction, and an additional induction on the
skeleta of an elementary simplex $\Delta_m$ shows that it is also
stably constant along $v^l_{n,m}$, so that \thetag{i} implies
\thetag{ii}. Next, note that any $Y \in \Delta^o_f\Sets$ with
discrete order is a globally finite-dimensional simplicial partially
ordered set, and if for some map $\chi:J \to J'$ in $\Posf$, we let
$p:Y \times J \to J'$ be the projection, then $L(a(\chi)) \times
\id_{R(Y)} = a^\dm(\chi \circ p)$, so that \thetag{iii} implies
\thetag{iv}. Since \thetag{iv} implies \thetag{i} by the same
argument as in Lemma~\ref{bva.le}, it remains to show that
\thetag{ii} implies \thetag{iii}.

To do this, it suffices to show that the map $a^\dm_\chi$ of
\eqref{a.chi.bis.eq} is a composition of pushouts of coproducts of
horn embeddings \eqref{bi.horn.bis.eq}. By induction on $m$, we may
assume that this is already proved for $\sk^R_{m-1}(a^\dm_\chi)$,
and we need to prove it for $\sk^R_m(a^\dm_\chi)$. Subsets in
$\sk_m^RN^L(J/J')$ that contain $\sk_{m-1}^RN^L(J/J')$ correspond
bijectively to left-closed subsets in $B(J([m])/J')$ that contain
the left-closed subset $B(J([m])/J')^\dg \subset B(J([m])/J')$ of
degenerate bisimplices. Explicitly, any bisimplex in $B(J([m])/J')$
is represented by an injective map $i:[n] \to J([m])/J'$, and it is
degenerate iff for some $m' < m$, the map $i$ factors through the
map $J([m'])/J' \to J([m])/J'$ corresponding to a surjective map
$[m] \to [m']$. Since $J'$ is constant in the simplicial direction,
one can equivalently require that the factorization exists for
$\sigma \circ i:[n] \to J([m])$.

Now, as in Lemma~\ref{bva.le}, the barycentric subdivision
$B(J([m])/J')$ has a filtration $B(J([m])/J')^k_d$, and we need to
show that if we enlarge its terms by replacing $B(J([m])/J')^k_d$
with the union $B(J([m])/J')^k_d \cup B(J([m])/J')^\dg$, then the
complements $\overline{B}(J([m])/J')^k_d$ still have the form
\eqref{B.S.km}. In other words, we have to show that for any term
$[1]$ in the right-hand side of the decomposition \eqref{B.S.km} for
$\overline{B}(J([m])/J')^k_d$, bisimplices corresponding to both
$0,1 \in [1]$ are degenerate or non-degenerate at the same time.
But these bisimplices are represented by injective maps $i_0:[d] \to
J([m])/J'$ and $i_1:[d+1] \to J([m])/J'$ described in detail in the
proof of Lemma~\ref{bva.le}, and in particular, the composition
$\sigma \circ i_1:[d+1] \to J([m])$ is equal to $\sigma \circ i_0
\circ p$ for some map $p:[d+1] \to [d]$. Therefore the factorization
condition for $\sigma \circ i_0$ is indeed equivalent to the same
condition for $\sigma \circ i_1$.
\endproof

\begin{corr}\label{bi.seg.corr}
Let $X \in \Delta^o_f\Sets_\kappa$ be a simplicial set, and let $\C$
be an additive semiexact Segal family of groupoids over
$\Delta_f^o\Delta^o_f\Sets_\kappa$ constant along $\kappa$-anodyne
maps of Definition~\ref{f.CW.def}. Let
$x:\Delta^o_f\Delta^o_f\Sets_\kappa \to
\Delta^o_f\Delta^o_f\Sets_\kappa$ be the functor $Y \mapsto L(X)
\times Y$. Then $x^*\C$ is also an additive semiexact Segal family
constant along $\kappa$-anodyne maps.
\end{corr}

\proof{} We need to check that $x^*\C$ is a Segal family, the rest
is automatic. As in the proof of Proposition~\ref{bi.seg.prop}, let
$W(\C)$ be the class of injective maps $f$ in
$\Delta^o\Delta^o_f\Sets_\kappa$ such that $\C$ is stably constant
along $f$. Then again, $W(\C)$ is saturated and closed under
coproducts and pushouts by Lemma~\ref{f.CW.le}, and by
Proposition~\ref{bi.seg.prop}~\thetag{iv}, it suffices to check that
for any $J \in \Posf_\kappa$ and $\chi:J \to [1]$, $W(\C)$ contains
$L(\id_X \times a_\chi)$. To simplify notation, denote the source and
target of the map $a_\chi$ by $Y$ and $Y'$, and more generally, let
$W$ be the class of injective maps $f:Z \to Z'$ in
$\Delta_f^o\Sets_\kappa$ such that the induced map $(Z' \times Y)
\copr_{Z \times Y} (Z \times Y') \to Z' \times Y'$ is in
$L^*W(\C)$. Then $W$ is also saturated and closed under coproducts and
pushouts, and we need to check that the tautological embedding
$\emptyset \to X$ is in $W$. By \eqref{sk.eq}, it then suffices to
show that $W$ contains all the sphere embeddings \eqref{sph.eq}, and
by induction on dimension, it further suffics to show that it
contains the embeddings $\emptyset \to \Delta_n$, $n \geq 0$. In
other words, we may assume that $X = \Delta_n$.  But then $X =
N([n])$, and $\id_X \times a_\chi = a_{\chi'}$, where $\chi'$ is the
projection $[n] \times J \to J \to [1]$.
\endproof

\subsection{Nerves and subdivisions.}\label{ner.div.subs}

Our version of Proposition~\ref{pos.prop} for Segal families uses
biordered sets of Section~\ref{bi.sec}, and we start by
constructing versions of the functors \eqref{N.Q}, \eqref{Q.N}
relating bisimplicial sets and biordered sets. By definition, the
category $\Delta$ comes equipped with the tautological embedding
$F:\Delta \to \Posf \subset \Pos$, $F([n]) = [n]$ that satisfies all
the assumptions of Lemma~\ref{bT.le}. Moreover, the good filtration
on $\Delta$ induced by $F$ coincides with the good filtration that
is the part of the Reedy structure: $\Delta_{\leq n}$ consists
exactly of $[m] \in \Delta$ with $\dim F([m]) = m \leq n$. Using
this embedding $F$ and the functors \eqref{rel.pos}, we can define
the nerve functors
\begin{equation}\label{bi.N.Q}
\Nn_\dm:\Rel \to \Delta^o\Delta^o\Sets, \qquad
\bN_\dm:\bRel \to \Delta^o\Delta^o\Sets
\end{equation}
by setting
$$
\begin{aligned}
\Nn_\dm(J)([n] \times [m]) &= \Rel(L([n]) \times R([m]),J),\\
\bN_\dm(J)([n] \times [m]) &= \bRel(L([n]) \times R([m]),J),
\end{aligned}
$$
and the embedding $\rho:\bRel \to \Rel$ induces a morphism
\begin{equation}\label{bN.N}
\bN_\dm \to \Nn_\dm \circ \rho.
\end{equation}
One can also construct the nerve functors \eqref{bi.N.Q} in two
steps, and there are two ways to do it: one can start with either
the Reedy or the simplicial direction in bisimplicial sets. Thus on
one hand, we have natural isomorphisms
\begin{equation}\label{N.E}
\Nn_\dm \cong N^R \circ \E_F, \qquad \bN_\dm \cong N^R \circ
\bE_F,
\end{equation}
where $N^R$ is one of the nerve functors \eqref{bi.nerve}, $\bE_F$
is the functor \eqref{E.F} associated to the embedding $F:\Delta \to
\Posf$, and $\E_F:\Rel \to \Delta^o\Pos$ is its version given
by
\begin{equation}\label{E.1.eq}
\E_F(J)([n]) = \Rel(L([n]),J)^l, \qquad J \in \Rel, \ [n] \in
\Delta.
\end{equation}
In this approach, the map \eqref{bN.N} is induced by a natural map
\begin{equation}\label{bE.E}
\bE_F \to \E_F \circ \rho
\end{equation}
corresponding to the embedding $\rho:\bRel \to \Rel$. On the
other hand, we also have
\begin{equation}\label{N.sE}
\Nn_\dm \cong N^L \circ \sE_F,
\end{equation}
where $\sE_F:\Rel \to \Delta^o\Pos$ and $\bsE_F:\bRel \to
\Delta^o\Pos$ are given by
\begin{equation}\label{E.2.eq}
\begin{aligned}
\sE_F(J)([m]) &= \Rel(R([m]),J),\\
\bsE_F(J)([m]) &= \Rel(R([m]),J)^r
\end{aligned}
\end{equation}
for any $J \in \Rel$, $[m] \in \Delta$, and $N^L$ is the other nerve
functor of \eqref{bi.nerve}. We note that for any $J \in \Rel$, we
have a natural identification
\begin{equation}\label{bsE.bE}
    \bsE_F(J^o) \cong \bE_F(J)^\iota,
\end{equation}
where $\bE_F(J)^\iota([n]) = \bE_F(J)([n])^o$, $[n] \in \Delta$, and
$J^o$ is $J$ with the opposite biorder.

\begin{lemma}\label{N.bi.exa.le}
The functors $\bN_\dm$ and $\Nn_\dm$ send coproducts to coproducts,
left-closed and right-closed full embeddings to injective maps, and
standard and co-standard pushout squares to standard pushout
squares.
\end{lemma}

\proof{} For any connected partially ordered ser $I$, and in
particular for $I = [n]$, $n \geq 0$, $\bRel(L(I),-)$ and
$\Rel(L(I),-)$ obviously preserve coproducts, and send left
resp.\ right-closed embeddings to left resp.\ right-closed
embeddings. If $I$ has a largest element, then both functors also
send standard pushout squares to standard pushout squares (for
$\bRel$, this is Lemma~\ref{hom.bi.le}, and for $\Rel(L(I),-) \cong
\Pos(I,-)$, the claim similarly immediately follows from
Lemma~\ref{hom.I.le}). If $I$ has a smallest element $o$, then both
functors also preserve co-standard pushout squares, by applying
Lemma~\ref{hom.I.le} to opposite partially ordered sets. Finally,
the correspondence $J \mapsto J^l$ also preserves standard and
co-standard pushout squares.  To finish the proof, use \eqref{N.E},
and apply Lemma~\ref{N.exa.le} pointwise to $N^R$.
\endproof

\begin{lemma}\label{bN.N.le}
For any biordered set $J$, the morphism $\bN_\dm(J) \to \Nn_\dm(J)$
of \eqref{bN.N} is a weak equivalence in $\Delta^o\Delta^o\Sets$.
\end{lemma}

\proof{} For any $[n] \in \Delta$, we have $[n] \cong I^>$, with $I
= [n-1]$ if $n \geq 1$ and $I=\emptyset$ if $n=0$. Then by
Lemma~\ref{fun.i.ar.le} and Lemma~\ref{le-re.le}, the embedding
$\bE_F(J)([n]) \to \E_F(J)([n])$ of \eqref{bE.E} is left-reflexive
for any $[n]$. Since by the same argument as in Lemma~\ref{NQ.le},
the nerve of a left-reflexive embedding is a weak equivalence, we
are done by \eqref{N.E}.
\endproof

\begin{lemma}\label{loc.N.le}
For any biordered ser $J \in \Rel$, with the adjunction map $a:J \to
R(U(J))$, the functor $R^h_!=\hocolim_{\Delta^o}$ of
\eqref{colim.del.eq} sends the map $\Nn_\dm(a):\Nn_\dm(J) \to
\Nn_\dm(R(U(J)))$ to a weak equivalence.
\end{lemma}

\proof{} For any map $f:X \to X'$ of bisimplicial sets, $R^h_!(f)$
is a weak equivalence if for any $[m] \in \Delta$, so is the
restriction $i_m^*(f)$ with respect to the embedding $i_m:\Delta^o
\to \Delta^o \times \Delta^o$ onto $\Delta^o \times [m]$. For our
map $s$, $i_0(a)$ is not only a weak equivalence but an isomorphism,
so it suffices to check that for any $m$ and $J \in \Rel$, the map
$e^*:i_0^*\Nn_\dm(J) \to i_m^*\Nn_\dm(J)$ induced by $e:[m] \to [0]$
is a weak equivalence. But by \eqref{N.sE}, $i_m^*\Nn_{\dm}(J)$ is
the nerve of the partially ordered set $\Rel(R([m]),J)$, and
$R(e):R([m]) \to R([0])$ is both left and right-reflexive, so we are
done by Example~\ref{fun.adj.exa} and the same argument as in
Lemma~\ref{bN.N.le}.
\endproof

Going in the other direction, we have functors
\begin{equation}\label{bi.Q.N}
\begin{aligned}
\bQQ_\dm &= \bT_F \circ \QQ_{\Delta}:\Delta^o\Delta^o\Sets \to
\bRel,\\
\QQ_\dm &= \rho \circ \bQQ_\dm \cong \T_F \circ
\QQ_{\Delta}:\Delta^o\Delta^o\Sets \to \Rel,
\end{aligned}
\end{equation}
where $\bT_F$ and $\T_F$ are the functors \eqref{T.F}, and
$\QQ_{\Delta}$ is the functor \eqref{Q.N.I} for $I = \Delta$.
Spelling out the definitions, we see that for any $X \in
\Delta^o\Delta^o\Sets$, we have a cartesian diagram
\begin{equation}\label{Q.rel.sq}
\begin{CD}
\QQ_\dm(X) @>{q_\idot}>> \Delta_\idot\Delta X\\
@V{p}VV @VVV\\
\QQ_{\Delta}(X) @>{q}>> \Delta\Delta X,
\end{CD}
\end{equation}
where $q$ is the augmentation functor, and $\Delta_\idot\Delta X =
\Delta_\idot \times_{\Delta} \Delta\Delta X$ is the cofibration over
the bisimplex category $\Delta\Delta X$ induced by the cofibration
$\nu_\idot$ of \eqref{nu.eq}. For any $J \in \Rel$ and $[n] \in
\Delta$, we have the evaluation map
\begin{equation}\label{q.del.ev}
[n] \times \E_F(J)([n]) \to U(J), \qquad i \times f \mapsto
f(i),
\end{equation}
where $i \times f \in [n] \times \E_F(J)([n])$ is a pair of an
element $i \in [n]$ and a map $f:L([n]) \to J$. This is functorial
with respect to $[n]$, thus defines a functor
\begin{equation}\label{E.ora}
\overrightarrow{\Delta}_\idot \E_F(J) = \Delta_\idot \times_\Delta
\overrightarrow{\Delta}\E_F(J) \to U(J)
\end{equation}
that inverts all maps cocartesian over
$\overrightarrow{\Delta}\E_F(J)$. Then combining \eqref{E.ora} and
\eqref{J.ora}, we obtain an evaluation functor
$$
\ev:\Delta_\idot\Delta \Nn_\dm(J) \cong \Delta_\idot \Delta
N^R(\E_F(J)) \to \overrightarrow{\Delta}^\hdot \E_F(J) \to U(J)
$$
that inverts all maps cocartesian over $\Delta\Delta
\Nn_\dm(J)$. Composing it with the functor $q_\idot$ of
\eqref{Q.rel.sq} gives a functor $\ev \circ
q_\idot:\QQ_\dm(\Nn_\dm(J)) \to U(J)$ that inverts all maps
cocartesian over $\QQ_\Delta(\Nn_\dm(J))$, and by
Example~\ref{cofib.rel.exa}, this means that $\ev \circ q_\idot$ is
actually a biordered map. It is functorial in $J$, so we obtain a
functorial map
\begin{equation}\label{bi.QN.id}
a_\dm:\QQ_\dm \circ \Nn_\dm \to \id.
\end{equation}
On the other hand, we have natural identifications
$$
\bN_\dm \circ \bQQ_\dm \cong \bN_\dm \circ \bT_F \circ \QQ_\Delta
\cong N^R \circ \bE_F \circ \bT_F \circ \QQ_\Delta \cong N^R \circ
\sE_\Delta \circ \QQ_\Delta \cong \Nn_\Delta \circ \QQ_\Delta,
$$
where $\Nn_\Delta$ is \eqref{N.Q.I} for $I=\Delta$, and the maps
\eqref{NQI.id} and \eqref{bN.N} induce a diagram
\begin{equation}\label{bi.NQ.id}
\begin{CD}
\Nn_\dm \circ \QQ_\dm @<<< \bN_\dm \circ \bQQ_\dm \cong \Nn_\Delta
\circ \QQ_\Delta @>>> \id.
\end{CD}
\end{equation}
For any uncountable regular cardinal $\kappa$, the functors
\eqref{bi.N.Q} and \eqref{bi.Q.N} restrict to functors between
$\Rel_\kappa$ resp.\ $\bRel_\kappa$ and
$\Delta^o\Delta^o\Sets_\kappa$. Moreover, the functor $\QQ_{\Delta}$
sends $\Delta^o_f\Delta^o_f\Sets$ into the category $\Posf \bb^\bB
\Delta$ of restricted finite-dimensional $\Delta$-augmented
partially ordered sets, and then \eqref{bi.Q.N} and Lemma~\ref{T.le}
immediately imply that $\bQQ_\dm$ resp.\ $\QQ_\dm$ sends
$\Delta^o_f\Delta^o_f\Sets$ into $\bRelf$ resp.\ $\Relf$. We can
also modify the two functors to send $\Delta^o\Delta^o\Sets$ into
$\Relf^+$ by setting
\begin{equation}\label{bi.Q.N.pl}
  \bQQ_{\dm+} = \bT_F \circ \QQ_{\Delta+}, \qquad \QQ_{\dm+} = \rho
  \circ \bQQ_{\dm+} \cong \T_F \circ \QQ_{\Delta+},
\end{equation}
where $\QQ_{\Delta+}$ is the functor \eqref{Q.I.pl} for $I =
\Delta$, and then \eqref{Q.I.pl.Q} induces functorial maps
\begin{equation}\label{bi.Q.N.pl.Q}
  \QQ_{\dm+}(X) \to \QQ_\dm(X), \qquad X \in \Delta^o\Delta^oX
\end{equation}
that are $\kappa^+$-bianodyne as soon as $X \in
\Delta^o_f\Delta^o_f\Sets_\kappa$ by Corollary~\ref{T.corr}.

\begin{lemma}\label{N.d.le}
For any $J \in \Relf$, the nerve $\bN_\dm(J)$ lies in
$\Delta^o_f\Delta^o_f\Sets$, and so does the Reedy skeleton
$\sk^L_d\Nn_\dm(J)$ for any $d \geq 0$.
\end{lemma}

\proof{} For the first claim, take $[n] \in \Delta$, write
$[n] = I^>$ as in Lemma~\ref{bN.N.le}, and note that
\eqref{fun.i.ar.eq} provides a cofibration $\bE_F(J)([n]) \to J$
whose fibers are discrete by Lemma~\ref{bio.le}. Therefore $\dim
\bE_F(J)([n]) \leq \dim J$, for any $[n]$, so that $\bE_F(J)$ lies
in $\Delta^o_{fd}\Posf$, and we are done by Lemma~\ref{binerve.le}
and \eqref{N.E}. For the second claim, use \eqref{N.sE}: again by
Lemma~\ref{binerve.le}, it suffices to check that $\sE_F(J)$ lies in
$\Delta^o_f\Posf$, and this is clear since its underlying simplicial
set is the nerve $N(J^l)$.
\endproof

By Lemma~\ref{bN.N.le}, the first map in \eqref{bi.NQ.id} is a
pointwise weak equivalence, and if we restrict our attention to
finite-dimensional bisimplicial sets, then the second map is a
pointwise weak equivalence by Lemma~\ref{QNI.le}~\thetag{iv}.
Analysing the map \eqref{bi.QN.id} is more diffucult: the nerve
$\Nn_\dm(J)$ of a biorder\-ed set $J \in \Relf$ is usually not a
regular bisimplicial set, so there is no analog of the isomorhism
\eqref{Q.B}. However, assume given $J \in \Posf$, and consider its
image $L(J) \in \Relf$. Then $\sE_F(L(J))$ is a constant simplicial
partially ordered set, so that $\Nn_\dm(L(J))$ is constant in the
simplicial direction, we have $\Delta\Delta\Nn_\dm(J) \cong \Delta J
\times \Delta$, and \eqref{Q.rel.sq} reduces to a cartesian square
\begin{equation}\label{Q.L.sq}
\begin{CD}
\QQ_\dm(\Nn_\dm(L(J))) @>{q^\hdot}>> B^\dm_\idot(J)\\
@VVV @VVV\\
B(B(J))^o @>{q}>> B(J),
\end{CD}
\end{equation}
where $B^\dm_\idot(J)$ is as in Example~\ref{bary.bis.rel.exa}, and we
identify $\QQ(\Nn(J)) \cong B(B(J))^o$, $B(J) \cong
\overline{\Delta}N(J)$. If we interpret $B^\tr(J)$ of
Example~\ref{bary.bis.rel.exa} as a $\Delta$-augmented partially
ordered set, then we have $\QQ_\dm(\Nn_\dm(L(J))) \cong
\T_F(B^o_\Delta(B^\tr(J)))$, while $B^\dm_\idot(J) \cong \T_F(B^\tr(J))$
and $q=\xi_\perp$, $q^\hdot = \T_F(\xi_\perp)$.

\begin{lemma}\label{A.QN.le}
For any uncountable regular cardinal $\kappa$ and $J \in
\Relf_\kappa$, the composition $\bQQ_\dm(\bN_\dm(J)) \to J$ of the
maps \eqref{bi.QN.id} and \eqref{bN.N} is $\kappa$-bianodyne, and
for any $J \in \Relf^+_\kappa$, the composition
$\bQQ_{\dm+}(\bN_\dm(J)) \to J$ of the maps \eqref{bi.QN.id},
\eqref{bN.N} and \eqref{bi.Q.N.pl.Q} is $\kappa^+$-bianodyne.
\end{lemma}

\proof{} Let $A = \bQQ_\dm \circ \bN_\dm:\Relf_\kappa \to
\Relf_\kappa$, and let $\alpha:A \to \Id$ be the composition of the
maps \eqref{bi.QN.id} and \eqref{bN.N}. Moreover, let $\chi = a
\circ p$, where $p$ is as in \eqref{Q.rel.sq}, and $a$ is the
map \eqref{QNI.id} for $I=\Delta$. Then we have $\chi \geq
U(\alpha)$, and $A$ preserves coproducts by Lemma~\ref{N.bi.exa.le},
Lemma~\ref{QNI.le}~\thetag{i} and Lemma~\ref{T.le}, so for the first
claim, it suffices to check that the triple $\langle A,\alpha,\chi
\rangle$ satisfies the conditions \thetag{i}-\thetag{iii} of
Lemma~\ref{A.le}. For \thetag{i}, note that $\chi^{-1}(J_0) \cong
\T_F(a^{-1}(J_0))$, and apply Lemma~\ref{T.le} and
Lemma~\ref{a.Q.le} for $I=\Delta$. For \thetag{ii}, note that
$\bN_\dm(J \times R([1])) \cong \bN_\dm(J) \times \bN_\dm(R([1]))
\cong \bN_\dm \times R(\Delta_1)$, and for any $X \in
\Delta^o_f\Delta^o_f \Sets_\kappa$, the projection $X \times
R(\Delta_1) \to X$ is $\kappa$-anodyne in the sense of
Definition~\ref{f.CW.def}, so we are done by
Lemma~\ref{QNI.le}~\thetag{iii} and Corollary~\ref{T.corr}.
Finally, for \thetag{iii}, note that by \eqref{Q.L.sq},
the map $\alpha:\bQQ_\dm(\bN_\dm(L(J))) \cong \T_F(B^o_\Delta(B^\tr(J))) \to
L(J)$ factors as
$$
\begin{CD}
\T_F(B^o_\Delta(B^\tr(J))) @>{\T_F(\xi_\perp)}>> B^\dm_\idot(J) @>>>
B^\dm(J) @>{\xi}>> L(J),
\end{CD}
$$
and $\T_F(\xi_\perp)$ is $\kappa$-bianodyne by Corollary~\ref{T.corr}
and Lemma~\ref{B.I.ano.le}, $\xi$ is $\kappa$-bianodyne by
Corollary~\ref{B.biano.corr}, and the map in the middle is the
strict adjoint map to the right-reflexive embedding $B^\dm(J) \to
B^\dm_\idot(J)$ of Example~\ref{bary.bis.rel.exa}. For the second claim,
since \eqref{bi.Q.N.pl.Q} is $+$-anodyne for
a finite-dimensional $X$, everything is already proved if $J
\in \Relf$, and then for any $J \in \Relf^+$, use
Lemma~\ref{rel.R.le} for the height map $\hht:J \to R(\N)$.
\endproof

\subsection{The cylinder axiom.}

We now turn to the study of families of groupoids over categories of
biordered sets. Assume given a full subcategory $\I \subset \Pos$
ample in the sense of Definition~\ref{amp.def}, and as in
Subsection~\ref{rel.subs}, let $\I^\dm = U^{-1}(\I) \subset \Rel$ be
its unfolding. In particular, $\I^\dm$ can be $\Relf$, $\Relf^+$ or
one of the categories $\Relf_\kappa$, $\Relf_\kappa^+$ for a regular
cardinal $\kappa$. As in Section~\ref{grp.sec}, we denote by
$i:\Relf \to \Relf^+$ the tautological embedding, together with the
induced embeddings $i:\Relf_\kappa \to \Relf^+_\kappa$. We now
consider a families of groupoids over $\I^\dm$. We note that
Definition~\ref{pos.def}~\thetag{ii},\thetag{iii} extend to such
families literally, and so does \thetag{i} if we replace $[1]$ with
$R([1])$. For \thetag{iv}, we need to find a biordered version of
the map \eqref{J.S}, and we use the map \eqref{J.S.bi} of
Example~\ref{J.S.bi.exa}. Definition~\ref{pos.pl.def} also makes
sense if we use \eqref{J.pl.bi.eq} instead of \eqref{J.pl.eq}.

\begin{defn}\label{rel.fam.def}
Let $\I^\dm \subset \Rel$ be the unfolding of an ample full
subcategory $\I \subset \Pos$, with extension $\I^{+\dm}$.  A family
of groupoids $\C$ over $\I^\dm$ is {\em $[1]$-invariant}, {\em
  additive}, {\em semiexact}, or {\em satisfies excision} if it is
additive, semiexact, or satisfies excision in the sense of
Definition~\ref{pos.def}, where for $[1]$-invariance, we replace
$[1]$ with $R([1])$, and for excision, we assume that $S^<,S^\hush
\in \I$, we take a map $J \to R(S^<)$ in $\I^\dm$, and we biorder
$J^\hush$ in \eqref{J.S} as in \eqref{B.S.J.bi}. A family $\C$ over
$\I^{+\dm}$ is {\em semicontinuous} if for any $J \in \I^{+\dm}$
equipped with a map $J \to R(\N)$, $\C$ is constant along the map
\eqref{J.pl.bi.eq}.
\end{defn}

Next, we want to have a version of Lemma~\ref{ano.grp.le} for
families over $\Relf$ and bianodyne maps, and as we saw in
Proposition~\ref{biano.prop}, one thing is missing: the bianodyne
maps \eqref{J.fl.hs.eq} of Example~\ref{J.fl.hs.exa}. This has to
required separately.

\begin{defn}\label{cont.def}
A family of groupoids $\C$ over $\I^\dm$ {\em satisfies the cylinder
  axiom} if for any $J \in \I$ equipped with a bicofibration $J \to
L([1])$, $\C$ is constant along the corresponding map
\eqref{J.fl.hs.eq}.
\end{defn}

\begin{defn}\label{rel.seg.def}
For any uncountable regular cardinal $\kappa$, a {\em
  $\kappa$-bounded Segal family} of groupoids over $\Relf_\kappa$ is
a family of groupoids $\C \to \Relf_\kappa$ that is additive,
semiexact and satisfies excision in the sense of
Definition~\ref{rel.fam.def}, satisfies the cylinder axiom of
Definition~\ref{cont.def}, and such that $\|\C_J\| < \kappa$ for any
finite $J \in \Relf$. A {\em Segal family} of groupoids over $\Relf$
is a family of groupoids $\C \to \Relf$ that is additive, semiexact,
satisfies excision in the sense of Definition~\ref{rel.fam.def}, and
satisfies the cylinder axiom of Definition~\ref{cont.def}. A {\em
  semicanonical extension} of a Segal family resp.\ $\kappa$-bounded
Segal family $\C$ is a family of groupoids $\C^+$ over $\Relf^+$
resp.\ $\Relf^+_\kappa$ equipped with an equivalence $i^*\C^+ \cong
\C$ that is $[1]$-invariant, additive, semiexact, and semicontinuous
in the sense of Definition~\ref{rel.fam.def}.
\end{defn}

As in the situation of Definition~\ref{grp.enh.def}, for any
uncountable regular cardinal $\kappa$, $\kappa$-bounded Segal
families over $\Relf_\kappa$ form a full subcategory in
$(\Cat/\Relf_\kappa)^0$ that we denote by $\Seg_\kappa$. Any small
Segal family $\C/\Relf$ restricts to a $\kappa$-bounded family over
$\Relf_\kappa \subset \Relf$ for a large enough $\kappa$.

\begin{lemma}\label{bi.del.le}
Let $\C$ be a Segal family of groupoids over $\Relf$, or a
$\kappa$-bounded Segal family of groupoids over $\Relf_\kappa$ in
the sense of Definition~\ref{rel.seg.def}. Then $\C$ is constant
along bianodyne resp.\ $\kappa$-bianodyne maps, and if $\C \to
\Relf$ is small, then $\T_F^*\C$ is a $\Delta$-augmented enhanced
groupoid in the sense of Definition~\ref{grp.I.enh.def}. Moreover,
any semicanonical extension $\C^+$ of such a family $\C$ is constant
along $+$-bianodyne resp.\ $\kappa^+$-bianodyne maps, and
$\T_F^*\C^+$ is a semicanonical extension of $\T_F^*\C$.
\end{lemma}

\proof{} For Segal families, the first claim immediately follows from
Proposition~\ref{biano.prop} by exactly the same argument as in
Lemma~\ref{ano.grp.le}. For the second claim, use Lemma~\ref{T.le}
and Corollary~\ref{T.corr}. For semicanonical extensions, the proof is
the same.
\endproof

\begin{corr}\label{X.bi.exci.corr}
A semicanonical extension $\C^+$ of a small or $\kappa$-bounded Segal
family in the sense of Definition~\ref{rel.seg.def} satisfies
excision and the cylinder axiom.
\end{corr}

\proof{} As in Corollary~\ref{X.exci.corr} and
Corollary~\ref{X.I.exci.corr}, the maps \eqref{J.S.bi} and
\eqref{J.fl.hs.eq} are $+$-bianodyne, and $\kappa^+$-bianodyne if
$|J| < \kappa$.
\endproof

\begin{corr}\label{cyl.2.corr}
Assume given a bicofibration $\chi:J \to L([2])$ in $\Relf$, and
consider the commutative square
\begin{equation}\label{s.1.t.sq}
\begin{CD}
  J_1 @>>> s^*J\\
  @VVV @VVV\\
  t^*J @>>> J
\end{CD}
\end{equation}
induced by \eqref{seg.del.sq} with $n=2$ and $l=1$. Then any Segal
family of groupoids over $\Relf$ is semicartesian along
\eqref{s.1.t.sq}.
\end{corr}

\proof{} Consider the map $\zeta^\dm_3:\ZZ^\dm_3 \to L([2])$ of
\eqref{z.dm.eq} and its restriction $\zeta^\dm_2:\ZZ^\dm_2 \to L([1])$
to $\ZZ^\dm_2 \subset \ZZ^\dm_3$. Then \eqref{zz.m.sq} induces a
standard coproduct decomposition
\begin{equation}\label{cyl.2.eq}
  \zeta_3^{\dm*}J \cong \zeta_2^{\dm*}s^*J \copr_{J_1} t^*J,
\end{equation}
and $\C$ is constant along the maps $\zeta_3^{\dm*}J \to J$,
$\zeta_2^{\dm*}s^*J \to s^*J$ by Lemma~\ref{Z.cyl.le} and
Lemma~\ref{bi.del.le}.
\endproof

We note that since all left-closed embeddings in $\Rel$ are strict,
and so are the maps of Example~\ref{J.S.bi.exa},
Example~\ref{J.pl.bi.exa} and Example~\ref{J.fl.hs.exa}, the
property of being a Segal family for some $\C$ only depends on its
restriction $\rho^*\C$ with respect to the embedding \eqref{rho.eq},
and it makes sense to speak about Segal families over $\bRelf$,
$\bRelf_\kappa$ and their semicanonical extensions to $\bRelf^+$,
$\bRelf^+_\kappa$. We then have the following result.

\begin{lemma}\label{NQ.bi.le}
Assume given an additive semiexact family of groupoids $\C$ over
$\Delta^o_f\Delta^o_f\Sets_\kappa$ for some uncountable regular
cardinal $\kappa$, and assume that $\C$ is $\kappa$-bounded,
constant along the class $W$ of $\kappa$-anodyne maps of
Definition~\ref{f.CW.def}, and is a Segal family in the sense of
Definition~\ref{rel.seg.def}. Then $\bN_\dm^*\C$ is a Segal family
over $\bRelf_\kappa$ in the sense of
Definition~\ref{rel.seg.def}. Moreover, for any semicanonical
extension $\C^+$ of the family $\C$ in the sense of
Definition~\ref{X.ext.def}, $\bN_\dm^*\C^+$ is a semicanonical
extension of $\bN_\dm^*\C$.
\end{lemma}

\proof{} By Lemma~\ref{N.bi.exa.le}, $\bN_\dm^*\C$ and
$\bN_\dm^*\C^+$ are additive and semiexact. Moreover, by
Corollary~\ref{W.CW.corr}, $\C$ and $\C^+$ are constant along weak
equivalences, and since $\bN_\dm(R([1])) \cong \Delta_1$, both
$\bN_\dm^*\C$ and $\bN_\dm^*\C^+$ are $[1]$-invariant. To see that
$\bN_\dm^*\C$ satisfies excision, note that for any $J \to R(S^<)$
in $\Relf_\kappa$, and any $[n] \in \Delta$, we have
$\bE_F(J^\hush)([n]) \cong \bE_F(J)([n])^\hush$ by
Lemma~\ref{hom.bi.le}, so by \eqref{N.E} and Lemma~\ref{NQ.le},
$\bN_\dm$ sends excision maps \eqref{J.S.bi} to weak equivalences in
$\Delta^o\Delta^o\Sets_\kappa$. Exactly the same argument shows that
$\bN_\dm^*\C^+$ is semicontinous, so it remains to check the
cylinder axiom for $\bN_\dm^*\C$. However, the nerve of a map
\eqref{J.fl.hs.eq} factors as
\begin{equation}\label{N.J.fl.dia}
\begin{CD}
\bN_\dm(J^\hash) \cong \bN_\dm(J_0 \times \V^\dm)
\copr_{\bN_\dm(J_0)} \bN_\dm(\Cyl(g_J))\\
@V{b}VV\\
\bN_\dm(J_0 \times [1]) \copr_{\bN_\dm(J_0)} \bN_\dm(\Cyl(g_J))\\
@V{a}VV\\
\bN_\dm(J^\flat),
\end{CD}
\end{equation}
where the first isomorphism is induced by \eqref{J.hs.eq}. Then $b$
in \eqref{N.J.fl.dia} is a pushout of the nerve of a reflexive map
$J_0 \times \V^\dm \to J_0 \times [1]$, thus a weak equivalence, and
by Example~\ref{cl.hom.exa} and \eqref{N.sE}, $a$ in
\eqref{N.J.fl.dia} is the map $a^\dm(\chi)$ of \eqref{a.chi.bis.eq}
for the characteristic map $\chi:\sE_F(\Cyl(g_J)) \to [1]$ of the
left-closed embedding $\sE_F(s):\sE_F(J_0) \to
\sE_F(\Cyl(g_J))$. Therefore $\C$ is constant along $a$ by
Proposition~\ref{bi.seg.prop}~\thetag{iii}.
\endproof

\subsection{Comparison.}

We can now prove a version of Proposition~\ref{pos.prop} comparing
Segal families of Definition~\ref{rel.seg.def} with those of
Definition~\ref{bi.seg.def}. One difficulty here is that the nerve
functor $\Nn_\dm$ of \eqref{bi.N.Q} does not send $\Relf$ into
$\Delta^o_f\Delta^o_f\Sets$. To circumvent this, look at the
product $\N \times \Rel$. For any category $\I$, the
projection $\pi:\N \times \I \to \I$ has a left-adjoint section $\I
\to \N \times \I$ onto $\{0\} \times \I$, so that $\C \cong
\pi_*\pi^*\C$ for any family of groupoids $\C$ over $\I$, and for
any functor $\gamma:\I_0 \to \I_1$, we have $\pi_*(\id \times
\gamma)^*\C \cong \gamma^*\pi_*\C$ for any family of groupoids $\C$
over $\N \times \I_1$. Moreover, for any $d \geq 0$, the shift
functor $q^d:\N \to \N$, $l \mapsto l + d$ is cofinal, and we have
\begin{equation}\label{shift.eq}
  \pi_*\C \cong \pi_*(q^d \times \id)^*\C
\end{equation}
for any family of groupoids $\C$ over $\N \times \I$. We then
consider the functors
\begin{equation}\label{bi.N.dot}
\begin{aligned}
&\bN_\idot:\N \times \Rel \to \Delta^o\Delta^o\Sets, \qquad d \times
J \mapsto \sk_d^L\bN_\dm(J),\\
&\Nn_\idot:\N \times \Rel \to \Delta^o\Delta^o\Sets, \qquad d \times
J \mapsto \sk_d^L\Nn_\dm(J).
\end{aligned}
\end{equation}
We have tautological maps
\begin{equation}\label{N.d.N}
  \bN_\idot \to \bN_\dm \circ \pi, \qquad \Nn_\idot \to \Nn_\dm
  \circ \pi
\end{equation}
defined by the skeleton inclusons, and by Lemma~\ref{N.d.le}, both
functors in \eqref{bi.N.dot} send $\N \times \Relf$ into
$\Delta^o_f\Delta^o_f\Sets$, while the map \eqref{N.d.N} for
$\bN_\idot(d \times J)$ is an isomorphism for any fixed $J \in
\Relf$ and sufficiently large $d$.

\begin{prop}\label{seg.prop}
\begin{enumerate}
\item For any additive semiexact family of group\-oids $\C$ over
  $\Delta^o_f\Delta^o_f\Sets$ that is small, constant along anodyne maps of
  Definition~\ref{f.CW.def} and is a Segal family in the sense of
  Definition~\ref{bi.seg.def}, the family of groupoids
  $\pi_*\Nn_\idot^*\C$ is a small Segal family over $\Relf$ in the
  sense of Definition~\ref{rel.seg.def}, and the functors
$$
\begin{CD}
\C \cong \pi_*\pi^*\C @>>> \pi_*(\id \times \bQQ_\dm)^*\bN_\idot^*\C
\cong \bQQ_\dm^*\pi_*\bN_\idot^* @<<< \QQ_\dm^*\pi_*\Nn_\idot^*\C
\end{CD}
$$
induced by the maps \eqref{bi.NQ.id} and \eqref{N.d.N} are
equivalences. Moreover, for any semicanonical extension $\C^+$ of the
family $\C$ in the sense of Definition~\ref{X.ext.def},
$\Nn_\dm^*\C^+$ with the functor $i^*\Nn_\dm^*\C^+ \cong
i^*\pi_*\pi^*\Nn_\dm^*\C^+ \to \pi_*\Nn_\idot^*\C$ induced by
\eqref{N.d.N} is a semicanonical extension of the family
$\pi_*\Nn_\idot^*\C$, and the functors
$\C^+,\QQ_{\dm+}^*\Nn_\dm^*\C^+ \to \bQQ_{\dm+}^*\bN_\dm^*\C^+$
induced by the maps \eqref{bi.Q.N.pl.Q} and \eqref{bi.NQ.id} are
equivalences.
\item For any small Segal family $\C$ over $\Relf$, the pullback
  $\QQ_\dm^*\C$ with respect to the functor \eqref{bi.Q.N} is an
  additive semiexact family over the category
  $\Delta^o\Delta^o\Sets$ constant along anodyne maps of
  Definition~\ref{f.CW.def}. Moreover, $\QQ_\dm^*\C$ is a Segal
  family in the sense of Definition~\ref{bi.seg.def}, and the
  functor $\C \cong \pi_*\pi^*\C \to \pi_*\Nn_\idot^*\QQ_\dm^*\C$
  induced by the composition of the maps \eqref{N.d.N} and
  \eqref{bi.QN.id} is an equivalence. For any semicanonical extension
  $\C^+$ of the Segal family $\C$, $\QQ_{\dm+}^*\C^+$ with the
  functor $i^*\QQ_{\dm+}^*\C^+ \to \QQ_\dm^*\C$ induced by the map
  \eqref{bi.Q.N.pl.Q} is a semicanonical extension of the Segal family
  $\QQ_\dm^*\C$, and the functor $\C^+ \to
  \Nn_\dm^*\QQ_{\dm+}^*\C^+$ induced by the maps
  \eqref{bi.Q.N.pl.Q}, \eqref{bi.QN.id} and \eqref{bN.N} is an
  equivalence.
\end{enumerate}
Moreover, for the any uncountable regular cardinal $\kappa$, the
same statements hold for $\kappa$-bounded Segal families over
$\Delta^o_f\Delta^o_f\Sets_\kappa$ and $\Relf_\kappa$.
\end{prop}

\proof{} It suffices to prove the $\kappa$-bounded version of both
statements; once everything is proved for an arbitrary $\kappa$, the
absolute statements follow. For \thetag{i}, the maps \eqref{bN.N}
and \eqref{N.d.N} induce functors
\begin{equation}\label{N.b.N.dia}
\begin{CD}
\rho^*\Nn_\idot^*\C @>>> \bN_\idot^*\C @<<< \pi^*\bN_\dm^*\C,
\end{CD}
\end{equation}
and since $\C$ is constant along weak equivalences by
Corollary~\ref{W.CW.corr}, the first functor is an equivalence by
Lemma~\ref{bN.N.le} and Corollary~\ref{sk.corr}. For any fixed $J
\in \Relf_\kappa$, the second functor becomes an equivalence over
$\N \times \{J\}$ once we restrict with respect to $q^d$ for a
sufficiently large $d$, so by \eqref{shift.eq}, \eqref{N.b.N.dia}
induces equivalences
\begin{equation}\label{N.b.N.equi}
\rho^*\pi_*\Nn_\idot^*\C \cong \pi_*\bN_\idot^*\C \cong \bN_\dm^*\C.
\end{equation}
Then $\pi_*\Nn_\idot^*\C$ is a Segal family by Lemma~\ref{NQ.bi.le},
and $\pi_*\Nn_\idot^*\C^+$ is its semicanonical extension. To finish the
proof, it remains to observe that the map $\bN_\dm \circ \bQQ_\dm
\cong \Nn_\Delta \circ \QQ_\Delta \to \id$ and its $+$-version are
weak equivalences by Lemma~\ref{QNI.le}~\thetag{iv}.

For \thetag{ii}, $\QQ_\dm^*\C$ is semiexact, additive and constant
along $\kappa$-anodyne maps of Definition~\ref{f.CW.def} by
Lemma~\ref{QNI.le}~\thetag{iii} and Lemma~\ref{bi.del.le}, and
$\QQ_{\dm+}^*\C^+$ is its semicanonical extension for the same reason,
while $\C \cong \pi_*\Nn_\dm^*\QQ_\dm^*\C$ and $\C^+ \cong
\Nn_\dm^*\QQ_{\dm+}^*\C^+$ by Lemma~\ref{A.QN.le} and
\eqref{N.b.N.equi}. To finish the proof, it remains to check that
$\QQ_\dm^*\C$ is a Segal family in the sense of
Definition~\ref{bi.seg.def}. By
Proposition~\ref{bi.seg.prop}~\thetag{iv}, we may assume given some
$J_l \in \Posf_\kappa$, a map $\chi:J_l \to [1]$, and a simplicial
set $Y \in \Delta^o_f\Sets_\kappa$, and we need to check that $\C$
is constant along the corresponding map $L(a_\chi) \times
\id_{R(Y)}$. Moreover, by Corollary~\ref{NQ.corr}, we may replace
$Y$ with $\Nn(\QQ(Y))$, so that $Y=\Nn(J_r)$ for some $J_r \in
\Posf_\kappa$. Then $L(\Nn(I)) \times R(Y) \cong \bN_\dm(L(I) \times
R(J_r))$ for any $I \in \Posf$, and by \eqref{a.chi.eq}, we need to
check that $\QQ_\dm^*\C$ is constant along the map
\begin{equation}\label{Q.a.chi}
\bN_\dm(J_0 \times L([1])) \copr_{\bN_\dm(J_0)} \bN_\dm(J_1) \to
\bN_\dm(J),
\end{equation}
where $J = L(J_l/[1]) \times R(J_r)$, with the bicofibration $J \to
L(J_l/[1]) \to L([1])$ induced by $\tau:J_l/[1] \to [1]$. Moreover,
since $J^\flat \to J$ is reflexive, we may replace $J$ resp.\ $J_1$
in \eqref{Q.a.chi} with $J^\flat$ resp.\ $J^\flat_1 = \Cyl(g_J)$,
where $g_J:J_0 \to J_1$ is the cylindrical map corresponding to the
bicofibration $J \to L([1])$. Then the claim immediately follows
from the cylinder axiom for $J$, Lemma~\ref{A.QN.le}, and the
decomposition \eqref{N.J.fl.dia}.
\endproof

\begin{corr}\label{seg.corr}
For any uncountable regular cardinal $\kappa$, any $\kappa$-bounded
Segal family $\C_\kappa \in \Seg_\kappa$ over $\Relf_\kappa$ extends
to a small Segal family $\C$ over $\Relf$. Moreover, for any two
$\kappa$-bounded Segal families $\C_\kappa,\C'_\kappa \in
\Seg_\kappa$ with extensions $\C$, $\C'$, any functor
$\gamma_\kappa:\C_\kappa \to \C_\kappa$ over $\Relf_\kappa$ extends
to a functor $\gamma:\C \to \C'$ over $\Relf$, and any morphism
$\gamma_\kappa \to \gamma'_\kappa$ between two such functors
$\gamma_\kappa$, $\gamma_\kappa'$ with extensions $\gamma$,
$\gamma'$ extends to a (non-unique) morphism $\gamma \to
\gamma'$. In particular, small Segal families and isomorphism
classes of functors over $\Relf$ between them form a well-defined
category
\begin{equation}\label{Seg.eq}
\Seg = \bigcup_\kappa \Seg_\kappa.
\end{equation}
We have an equivalence between the category $\Seg$ and the full
subcategory in $h^W(\Delta^o\Delta^o\Sets)$ spanned by Segal spaces,
and for any uncountable regular cardinal $\kappa$, this equivalence
induces an equivalence between $\Seg_\kappa$ and the full
subcategory of Segal spaces in
$h^W(\Delta^o\Delta^o\Sets_\kappa)$. Moreover, any small or
$\kappa$-bounded Segal family $\C$ admits a semicanonical extension
$\C^+$, for two families $\C_0$, $\C_1$ with semicanonical
extensions $\C_0^+$, $\C_1^+$, a functor $\gamma:\C_0 \to \C_1$
extends to a functor $\gamma^+:\C_0^+ \to \C_1^+$, and for two
functors $\gamma_0,\gamma_1:\C_0 \to \C_1$ with extensions
$\gamma_0^+$, $\gamma_1^+$, a morphism $\gamma_0 \to \gamma_1$
extends to a (non-unique) morphism $\gamma_0^+ \to \gamma_1^+$, so
that in particular, the semicanonical extension is unique up to a
unique equivalence.
\end{corr}

\proof{} Clear. \endproof

\section{Complete Segal spaces.}\label{css.sec}

\subsection{Reflexive families.}\label{refl.grp.subs}

As it happens, one can show that under an additional assumption, a
Segal family of groupoids $\C$ over $\Relf$ in the sense of
Definition~\ref{rel.seg.def} and its semicanonical extension $\C^+$
can be completely reconstructed from their restriction $L^*\C$,
$L^*\C^+$ to $\Posf$, $\Posf^+$ with respect to the functor $L$ of
\eqref{rel.pos}. Here is the assumption.

\begin{defn}\label{rel.refle.def}
Let $\I^\dm \subset \Rel$ be the unfolding of an ample full
subcategory $\I \subset \Pos$. A family of groupoids $\C$ over
$\I^\dm$ is {\em reflexive} if it is fully faithful along the
projection $e:J \times L([1]) \to J$ for any $J \in \I^\dm$.
\end{defn}

To state and prove the reconstruction result, we need some
notation. For any reflexive family of groupoids $\C \to \I^\dm$,
biordered set $J \in \I^\dm$, and a set $W \subset \I^\dm(L([1]),J)$ of
maps $w:L([1]) \to J$, let $\C_J(W) \subset \C_J$ be the full
subcategory spanned by objects $c \in \C_J$ such that $w^*c$ lies in
$e^*(\C_\ppt) \subset \C_{L([1])}$ for any $w \in W$. Moreover, for
any map $f:J' \to J$ in $\I^\dm$, let $W(f) \in \I^\dm(L([1]),J)$ be the set
of maps $w:L([1]) \to J'$ such that $f \circ w:L([1]) \to J$ factors
through $L([1]) \to R([1])$. Finally, note that for any map $f:J
\to J'$ in $\I^\dm$, the transition functor $f^*:\C_{J'} \to \C_J$
factors through a functor
\begin{equation}\label{ol.f.eq}
\C_{J'} \to \C_J(W(f)).
\end{equation}
Say that $\C$ is {\em adapted to} the map $f$ iff \eqref{ol.f.eq} is
an equivalence of categories, and denote by $E(\C)$ the class of
dense maps $f:J \to J'$ in $\I^\dm$ such that $\C$ is adapted to $\id
\times f:J_0 \times J \to J_0 \times J'$ for any $J_0 \in \I^\dm$.

\begin{lemma}\label{E.0.le}
Assume given a reflexive family of groupoids $\C$ over $\I^\dm$
constant along a morphism $f:J \to J'$, and subsets $W \subset
\I^\dm(L([1]),J)$, $W' \subset \I^\dm(L([1]),J')$ such that $f$
induces a surjective map $W \to W'$. Then the equivalence
$f^*:\C_{J'} \to \C_J$ induces an equivalence $\C_{J'}(W') \cong
\C_J(W)$.
\end{lemma}

\proof{} Since $f(W') \subset W$, the equivalence $f^*$ indeed sends
$\C_{J'}(W')$ into $\C_J(W)$, and it suffices to prove that the
induced fully faithful functor is essentially surjective. But since
$W \to W'$ is surjective, any $w' \in W'$ lifts to $w \in W$, and
then for any $c \in \C_{J'}$, we have ${w'}^*c \cong w^*f^*c$.
\endproof

\begin{lemma}\label{E.3.le}
For any reflexive Segal family $\C$ over $\I^\dm=\Relf$, the tautological
map $p:L([1]) \to R([1])$ is in $E(\C)$.
\end{lemma}

\proof{} By Lemma~\ref{bi.del.le}, $\C$ is constant along the
projection $J \times R([1]) \to J$ for any $J \in \Relf$. Since $\C$
is reflexive, the functor \eqref{ol.f.eq} for the morphism $f = f(J)
= \id_J \times p:J \times L([1]) \to J \times R([1])$ is always
fully faithful. Say that $J \in \Relf$ is {\em good} if this functor
is essentially surjective. Then we need to prove that any $J$ is
good. For any left-reflexive full embedding $g:J' \to J$, with the
adjoint map $g_\dg:J \to J$, $g_1 = g \times \id_{L([1])}$ is also
left-reflexive, with the adjoint map $g_{1\dg} = g_\dg \times
\id_{L([1])}$, $\C$ is constant along $g_\dg$ and $g_{1\dg}$ by
Lemma~\ref{bi.del.le}, and the subsets $W(f(J'))$, $W(f(J))$ satisfy
the assumptions of Lemma~\ref{E.0.le} with respect to $g_{1\dg}$, so
that $J$ is good iff so is $J'$. If $J$ admits a map $\chi:J \to \V$
with good comma-fibers $J/0$, $J/1$, $J_o$, then $J$ is good by
semiexactness of $\C$, and if $J$ admits a map $J \to R(S^<)$ for
some set $S$, then the map $(J \times L([1]))^\hush \cong J^\hush
\times L([1]) \to J \times L([1])$ of \eqref{J.S.bi} also satisfies
the assumptions of Lemma~\ref{E.0.le} with respect to
$W(f(J^\hash))$ and $W(f(J))$, so that $J$ is good iff so is
$J^\hush$. By additivity of $\C$, a coproduct of good biordered sets
is good, and then by excision and \eqref{B.S.J.bi}, any $J$ that
admits a map $J \to R(S^<)$, $S \in \Sets$ with good comma-fibers is
itself good. Now the same induction on dimension as in
Lemma~\ref{A.le} reduces us to the situation when we have a
cofibration $\chi:J \to [1]$ in $\Posf$, the fibers $L(J_0)$,
$L(J_1)$ are good, and we need to check that $L(J)$ is good. In this
case, we have the cofibration $\chi \times \id:J \times [1] \to
[1]^2$, the diagonal embedding $\delta:[1] \to [1]^2$ lifts to a
section $\delta:J \to J \times [1]$ of the projection $e:J \times
[1] \to J$, and it suffices to show that for any $c \in \C_{L(J
  \times [1])}(W(f(J)))$, there exists an isomorphism $f:c \to
e^*\delta^*c$ such that $\delta^*(f)=\id$.

To do this, recall that we have the two embeddings $a,b:[2] \to
[1]^2$ that fit into the cocartesian square \eqref{pos.sq} of
Example~\ref{pos.sq.exa}, and $\delta = a \circ m = b \circ m$,
where $m:[1] \to [2]$ sends $0$ to $0$ and $1$ to $2$. Then if we
denote $J^a = a^*(J \times [1])$, $J^b = b^*(J \times [1])$,
\eqref{pos.sq} induces a commutative diagram
\begin{equation}\label{ab.sq}
\begin{CD}
J^a @>{a}>> J \times [1])\\
@A{m_a}AA @AA{b}A\\
J @>{m_b}>> J^b
\end{CD}
\end{equation}
in $\Posf$ that becomes a standard pushout diagram in $\Relf$ after
we apply the barycentric subdivision functor $B^\dm$ of
Example~\ref{bary.rel.exa}. Therefore by
Corollary~\ref{B.biano.corr}, Lemma~\ref{bi.del.le} and
semiexactness of $\C$, the family $L^*\C$ is semicartesian along
\eqref{ab.sq}, so it suffices to prove that there exist isomorphisms
$a^*c \cong m_a^*\delta^*c$ in $\C_{L(J^a)}$ and $b^* \cong
m_b^*\delta^*c$ in $\C_{L(J^b)}$ that restrict to $\id$ on $J
\subset J^a,J^b$. This immediately follow from
Corollary~\ref{cyl.2.corr} applied to $L(J^a),L(J^b) \to L([2])$,
and the assumption that $L(J_0)$, $L(J_1)$ are good.
\endproof

\begin{lemma}\label{E.1.le}
Assume given a left-reflexive full embedding $g:J' \to J$ in
$\I^\dm$, and a $[1]$-invariant reflexive family of groupoids $\C$
over $\I^\dm$ such that the map $p:L([1]) \to R([1])$ is in
$E(\C)$. Then the natural map $f:L(U(J)) \to J$ is in $E(\C)$ if and
only if so is the map $f':L(U(J')) \to J'$.
\end{lemma}

\proof{} Since the family $\C$ is $[1]$-invariant, it is constant
along reflexive maps by Lemma~\ref{amp.le}. Moreover, $g \times
\id:J' \times J_0 \to J \times J_0$ is left-reflexive for any $J_0
\in \Relf$, so it suffices to prove that $g^*:\C_{L(U(J))} \to
\C_{L(U(J'))}$ induces an equivalence $\C_{L(U(J))}(W(f)) \cong
\C_{L(U(J'))}(W(f'))$. Let $g_\dg:J \to J'$ be the adjoint map, and
let $h:J \times [1] \to J$ be the map equal to $\id$ on $J \times
\{0\}$ and $g \circ g_\dg$ on $J \times [1]$. Then for any $c \in
\C_{L(U(J))}(W(f))$, $h^*c$ lies inside $\C_{L(U(J) \times
  [1])}(W(\id \times p))$. Since $\C$ is reflexive and $p \in
E(\C)$, we have the equivalence
\begin{equation}\label{p.eq.eq}
(\id \times p)^*:\C_{L(U(J))} \cong \C_{L(U(J)) \times R([1])} \to
\C_{L(U(J) \times [1])}(W(\id \times p)),
\end{equation}
and the functors $s^*,t^*:\C_{L(U(J) \times [1])}(W(\id \times p))
\to \C_{L(U(J))}$ induced by the embeddings $s,t:[0] \to [1]$ onto
$0,1 \in [1]$ are both one-sided inverses to
\eqref{p.eq.eq}. Therefore they are both equivalences, and they are
isomorphic. Thus $\id \cong g_\dg^* \circ g^*:\C_{L(U(J))}(W(f)) \to
\C_{L(U(J))}(W(f))$, and $g^*$, $g_\dg^*$ indeed induce mutually
inverse equivalences of categories between $\C_{L(U(J))}(W(f))$ and
$\C_{L(U(J'))}(W(f'))$.
\endproof

\begin{lemma}\label{E.2.le}
Assume that $\I$ is either a very ample full subcategory, or an
extension of such, and assume given a reflexive $[1]$-invariant
additive semiexact family of groupoids $\C$ over $\I^\dm$, a dense
map $f:J' \to J$, and a map $\chi:J \to R(\V)$ such that the
comma-fibers $f/j:J'/j \to J/j$, $j \in \V$ are in $E(\C)$. Then $f
\in E(\C)$.
\end{lemma}

\proof{} For any $J_0 \in \I^\dm$, $f \times \id:J' \times J_0 \to J
\times J_0$ also satisfies the assumptions of the Lemma, so it
suffices to prove that the functor \eqref{ol.f.eq} for $f$ is an
equivalence. It is an epivalence by the semiexactness of $\C$, so it
further suffices to check that it is fully faithful. More generally,
let $E'(\C) \supset E(\C)$ be the class of all maps $g$ in $\I^\dm$
such that $\C$ is fully faithful along $g \times \id_{J_0}$ for any
$J_0 \in \I^\dm$. Then $E'(\C)$ is left-saturated by
Lemma~\ref{sat.C.le}, so it has ``one-half'' of the two-out-of-three
property: if $g_0 \circ g_1,g_0 \in E'(\C)$, then $g_1 \in
E'(\C)$. Then since $\C$ is $[1]$-invariant, it is constant along
all reflexive maps by the biordered version of
Lemma~\ref{refl.sat.le}, and by the same argument as in
Lemma~\ref{ano.grp.le}, Lemma~\ref{ano.C.le} implies that $E'(\C)$
is closed under standard pushouts. Now exactly the same
decomposition as in Lemma~\ref{st.ano.le} shows that $f \in E'(\C)$
as soon as $f/j \in E'(\C)$ for all $j \in \V$.
\endproof

\begin{lemma}\label{E.le}
Assume given a reflexive Segal family $\C$ over $\Relf$. Then any
dense map $f:J_0 \to J$ in $\Relf$ lies in the class
$E(\C)$. Moreover, any semicanonical extension $\C^+$ of the family $\C$
is also reflexive, and $E(\C^+)$ also contains all dense maps.
\end{lemma}

\proof{} For the first claim, note that as in the proof of
Lemma~\ref{E.3.le}, for any set $S$ and map $J \to R(S^<)$,
Lemma~\ref{E.0.le} shows that $f:J_0 \to J$ is in $E(\C)$ iff so is
$f^\hush:J_0^\hush \to J^\hush$, and then Lemma~\ref{E.2.le} and
induction on dimension reduce the proof to the
case when $J$ has a largest element $j$. Then it suffices to prove
that the maps $L(U(J_0)) = L(U(J) \to J_0,J$ are in $E(\C)$, so that we
may assume right away that $J_0 = L(U(J))$. In this case, as in the
proof of Lemma~\ref{A.le}, take the reflexive subset $J^R \subset J$
of Lemma~\ref{bio.le}, and apply Lemma~\ref{E.3.le} and
Lemma~\ref{E.1.le}. For the second claim, it suffices to prove that
any dense map $f:J' \to J$ in $\Relf^+$ is in $E(\C^+)$. To do this,
consider the height map $\hht:J \to \N$, and note that again,
Lemma~\ref{E.0.le} together with the semicontinuity of $\C^+$ show
that $f$ lies in $E(\C^+)$ iff so does $f^+:J_0^+ \to J^+$. But by
additivity, $E(\C^+)$ is closed under coproducts, and then the same
additional argument as in Proposition~\ref{pl.ano.prop} reduces us to
Lemma~\ref{E.2.le}.
\endproof

\subsection{Weak excision and weak semicontinuity.}\label{unf.subs}

To describe reflexive Segal families over $\Relf$ and their
semicanonical extensions to $\Relf^+$ purely in terms of their
restrictions to $\Posf$, resp.\ $\Posf^+$, we start by repeating
Definition~\ref{rel.refle.def} literally.

\begin{defn}\label{pos.refle.def}
Let $\I \subset \Pos$ be an ample full subcategory.  A family of
groupoids $\C$ over $\I$ is {\em reflexive} if for any $J \in \I$,
$\C$ is fully faithful along the projection $e:J \times [1] \to J$.
\end{defn}

As in the biordered case, for any reflexive family of groupoids $\C$
over $\I$ and any $J \in \I$ equipped with a subset $W \subset
\I([1],J)$, we denote by $\C_J(W) \subset \C_J$ the full subcategory
spanned by objects $c$ such that $w^*c \in \C_{[1]}$ lies in the
image of the embedding $\C_\ppt \to \C_{[1]}$ for any $w \in W$. In
particular, if we let $\I^\dm$ be the unfolding of $\I$, then for
any set $J \in \I^\dm$, we can let $W(J) = \I([1],J^l) \subset
\I([1],U(J))$ and $\C^\dm_J = \C_{U(J)}(W(J)) \subset \C_{U(J)}$, and
this defines a family of groupoids $\C^\dm$ over $\I^\dm$. We call
$\C^\dm$ the {\em unfolding} of the family $\C$. We have $L^*\C^\dm
\cong \C$, where $L:\I \to \I^\dm$ is the functor \eqref{rel.pos}.

\begin{defn}\label{pos.nondege.def}
A reflexive family of groupoids $\C$ over an ample $\I$ with
unfolding $\C^\dm$ is {\em non-de\-ge\-ne\-rate} if for any $J \in
\I$, $\C^\dm$ is constant along the projection $L(J) \times R([1])
\to L(J)$
\end{defn}

\begin{lemma}\label{E.1.bis.le}
The unfolding $\C^\dm$ of a non-degenerate reflexive family of
groupoids $\C$ is $[1]$-invariant and constant along all reflexive
full embeddings.
\end{lemma}

\proof{} By the very definition of the unfolding $\C^\dm$,
Definition~\ref{pos.nondege.def} implies that is is $[1]$-invariant,
and then it is constant along reflexive full embeddings by the
biordered version of Lemma~\ref{refl.sat.le}.
\endproof

We also need appropriate versions of excision, semicontinuity and
the cylinder axiom. To state them, recall that the maps \eqref{xi.B}
resp.\ \eqref{zeta.eq} fit into perfect cartesian cocartesian
squares \eqref{S.S.sq} resp.\ \eqref{Z.N.sq} in $\Pos$. We can treat
\eqref{S.S.sq} as a square in $\I/S^<$, and \eqref{Z.N.sq} is a
square in $\I / \N$. Thus for any partially ordered set $J$ equipped
with a map $J \to S^<$, the map \eqref{J.S} fits into a cartesian
square
\begin{equation}\label{J.S.sq}
\begin{CD}
S^> \times J_o @>>> J_o\\
@VVV @VVV\\
J^\hush @>>> J,
\end{CD}
\end{equation}
where we identify $S^> \times J_o \cong S^> \times_{S^<} J^\hush$,
$J_o \cong \ppt \times_{S^<} J$, while for any $J$ equipped with a
map $J \to \N$, the map \eqref{J.pl.eq} fits into a diagram
\begin{equation}\label{J.pl.sq}
\begin{CD}
[1] \times \coprod_{n \geq 1} J/n @>>> \coprod_{n \geq 1} J/n\\
@VVV @VVV\\
J^+ @>>> J / \N @>{\sigma}>> J,
\end{CD}
\end{equation}
where $\sigma$ is reflexive, and the square is a cartesian square
induced by the cartesian square \eqref{Z.N.sq}. Since the
cocartesian square \eqref{S.S.sq} is universal, the square
\eqref{J.S.sq} is also cocartesian, and the square \eqref{J.pl.sq}
is cocartesian by Corollary~\ref{univ.corr}. If $J \in \I$ and
$S^<,S^\hush \in \I$, then \eqref{J.S} is also in $\I$, while if
$\ZZ_\infty,\N \in \I$, then \eqref{J.pl.sq} is in $\I$ as
well. Moreover, the map $\zeta:\ZZ_3 \to [2]$ induced by
\eqref{zeta.eq} fits into the cartesian square \eqref{Z3.sq} that
can be understood as a square in $\Posf/[2]$, and we can treat
\eqref{Z3.sq} as a square in $\Posf/[1]$ by composing the structure
maps with the projection $t_\dg:[2] \to [1]$ of
Example~\ref{n.refl.exa} (explicitly, $t_\dg(0) = t_\dg(1) = 0$,
$t_\dg(2) = 1$). Then for any partially ordered set $J$ equipped
with a map $J \to [1]$, \eqref{Z3.sq} induces a cartesian square
\begin{equation}\label{J.cyl.sq}
\begin{CD}
J_0 \times [1] @>>> J_0\\
@VVV @VVV\\
\ZZ_3 \times_{[1]} J @>>> J/[1] \cong t_\dg^*J.
\end{CD}
\end{equation}
Since by Lemma~\ref{pos.uni.le}, \eqref{Z3.sq} is universal with
respect to $t_\dg:[2] \to [1]$, the square \eqref{J.cyl.sq} is also
cocartesian, and lies in $\I$ if so does $J$ and $\ZZ_3$, $[2]$. If
the map $J \to [1]$ is a bicofibration $J^\dm \to L([1])$ for some
biorder $J^\dm$ on $J$, then the bottom arrow in \eqref{J.cyl.sq} is
the map underlying \eqref{J.fl.hs.eq}.

\begin{defn}\label{refl.exci.def}
\begin{enumerate}
\item A family of groupoids $\C$ over an ample subcategory $\I
  \subset \Pos$ {\em satisfies weak excision} if it is semicartesian
  along the square \eqref{J.S.sq} for any set $S$ such that
  $S^<,S^\hush \in \I$, and any $J \in \I$ equipped with a map $J
  \to S^<$.
\item A family of groupoids $\C$ over $\I$ {\em satisfies the
  cylinder axiom} if it is semicartesian along the square
  \eqref{J.cyl.sq} for any $J \in \I$ equipped with a map $J \to
        [1]$.
\item A family of groupoids $\C$ over $\I$ with $\ZZ_\infty,\N \in
  \I$ is {\em weakly semicontinuous} if it is semicartesian along
  the square \eqref{J.pl.sq} for any $J \in \I$ equipped with a map
  $J \to \N$.
\end{enumerate}
\end{defn}

\begin{defn}\label{red.seg.def}
A {\em restricted Segal family} is a family of groupoids $\C$ over
$\Posf$ that is reflexive in the sense of
Definition~\ref{pos.refle.def}, non-degenerate in the sense of
Definition~\ref{pos.nondege.def}, additive and semiexact in the
sense of Definition~\ref{pos.def}~\thetag{ii},\thetag{iii}, and
satisfies weak excision and the cylinder axiom in the sense of
Definition~\ref{refl.exci.def}. A {\em semicanonical extension} of a
restricted Segal family $\C$ is a family of groupoids $\C^+$ over
$\Posf^+$ that is additive and semiexact in the sense of
Definition~\ref{pos.def}~\thetag{ii},\thetag{iii}, reflexive in the
sense of Definition~\ref{pos.refle.def}, non-degenerate in the sense
of Definition~\ref{pos.nondege.def}, weakly semicontinuous in the
sense of Definition~\ref{refl.exci.def}, and equipped with an
equivalence $i^*\C^+ \cong \C$.
\end{defn}

\begin{exa}\label{I.triv.seg.exa}
Since all the squares \eqref{J.S.sq}, \eqref{J.cyl.sq},
\eqref{J.pl.sq} are cocartesian squares of categories, and so are
standard pushout squares, the fibration $\bpi:\Posf \bbi I \to \Posf$
of \eqref{pos.I.pi} for any category $I$ is a restricted Segal
family, and $\bpi:\Posf^+ \bbi I \to \Posf^+$ is its semicanonical
extension. Explicitly, for any $J \in \Posf^+$, the fiber $(\Posf^+
\bbi I)_J$ is the isomorphism groupoid $\Fun(J,I)_{\Iso}$, and for any
biorder $J^\dm$ on $J$, the fiber $(\Posf^+ \bbi I)^\dm_{J^\dm}$ of the
unfolding is spanned by functors $J \to I$ that send each order
relation $j \leq^l j'$ to an invertible map in $I$. If $I$ is rigid,
then $\Posf^+ \bbi I \to \Posf^+$ is discrete.
\end{exa}
  
\begin{prop}\label{restr.seg.prop}
The unfolding $\C^\dm$ of a restricted Segal family $\C$ over $\Posf$ in
the sense of Definition~\ref{red.seg.def} is a Segal family in the
sense of Definition~\ref{rel.seg.def} and reflexive in the sense of
Definition~\ref{rel.refle.def}. Conversely, the restriction $L^*\C$
of a reflexive Segal family $\C$ over $\Relf$ in the sense of
Definition~\ref{rel.seg.def} and Definition~\ref{rel.refle.def} is a
restricted Segal family in the sense of
Definition~\ref{red.seg.def}, and we have a natural equivalence $\C
\cong (L^*\C)^\dm$. Moreover, for any semicanonical extension $\C^+$ of
the family $\C$, the unfolding $\C^{+\dm}$ is a semicanonical extension
of the unfolding $\C^\dm$, and conversely, for any semicanonical
extension $\C^{\dm+}$ of the unfolding $\C^\dm$, $\C^+ =
L^*\C^{\dm+}$ is a semicanonical extension of $\C$, it is reflexive,
safisfies weak excision and the cylinder axiom of
Definition~\ref{refl.exci.def}, and we have $\C^{\dm+} \cong
\C^{+\dm}$.
\end{prop}

\proof{} For the first claim, note that $\C^\dm$ is tautologically
additive and semiexact, and it is constant along reflexive maps by
Lemma~\ref{E.1.bis.le}. In particular, it is constant along the map
$J \times R(S^>) \to J$ for any $J \in \Relf$, and then the weak
excision property for $\C$ implies that $\C^\dm$ is semiconstant
along a map $J^\hush \to J$ with $J^\hush$ as in
\eqref{B.S.J.bi}. Moreover, for any $J' \in \Relf$, we have $(J'
\times J)^\hush \cong J' \times J^\hush$, so if we let
$j,j_\hush:\Relf \to \Relf$ be the functors $J' \mapsto J' \times
J,J' \times J^\hush$, then the transition functor $j^*\C^\dm \to
j_\hush^*\C^\dm$ is an epivalence. Therefore it is an equivalence by
Corollary~\ref{semi.epi.corr}, so that $\C^\dm$ satisfies
excision. For the same reason, the cylinder axiom for $\C^\dm$ for a
bicofibration $J \to L([1])$ immediately follows from the cylinder
axiom for $\C$, for the underlying map $U(J) \to [1]$. Conversely,
for a reflexive Segal family $\C$ over $\Relf$, the restriction
$L^*\C$ is tautologically additive, semiexact and reflexive, and it
is non-degenerate by Lemma~\ref{E.3.le}. Weak excision for $L^*\C$
then immediately follows from excision for $\C$. Moreover, we have
$\C \cong (L^*\C)^\dm$ by Lemma~\ref{E.le}. As for the cylinder
axiom, note that for any map $\chi:J \to [1]$ in $\Posf$, the map
$L(J/[1]) \cong t_\dg^*L(J) \to L([2])$ satisfies the assumption of
Lemma~\ref{Z.cyl.le}, so that the map $\zeta_3^{\dm*}L(J/[1]) \to
L(J/[1])$ is bianodyne, and $\C^\dm$ is constant along this map by
Lemma~\ref{bi.del.le}. Again by Lemma~\ref{E.1.bis.le}, this is
equivalent to the cylinder axiom for $\C$.

For $\C^+$, the family $\C^{+\dm}$ is $[1]$-invariant since $\C^+$
is non-degenerate, and then again, by Lemma~\ref{E.1.bis.le} and
Corollary~\ref{semi.epi.corr}, weak semicontinuity for $\C^+$
immediately implies semicontinuity for $\C^{+\dm}$. Conversely,
$\C^{\dm+}$ satisfies excision and the cylinder axiom by
Corollary~\ref{X.bi.exci.corr}, and then $\C^+$ is reflexive by
Lemma~\ref{E.le} and non-degenerate since $\C^{\dm+}$ is
$[1]$-invariant. Then again, Lemma~\ref{E.1.bis.le} applies, so that
weak excision for $\C^+$ follows from excision for $\C^{\dm+}$, weak
semicontinuity resp.\ the cylinder axiom follow from
Lemma~\ref{bi.del.le} and Lemma~\ref{Z.cyl.le} applied to the map
$\zeta^{\dm*}L(J/\N) \to L(J/\N)$ resp.\ $\zeta_3^{\dm*}L(J) \to
L(J/[1])$, and we have $\C^{\dm+} \cong \C^{+\dm}$ by
Lemma~\ref{E.le}.
\endproof

\begin{corr}\label{seg.exa.corr}
A restricted Segal family $\C \to \Posf$ is cartesian along all
squares \eqref{J.S.sq}, \eqref{J.cyl.sq}, and its canonical
extension $\C^+$ is cartesian along all squares \eqref{J.pl.sq}.
\end{corr}

\proof{} Clear. \endproof

\begin{corr}\label{iota.corr}
Assume given a small restricted Segal family $\C \to \Posf$, and let
$\iota:\Posf \to \Posf$ be the involution $J \mapsto J^o$. Then
$\iota^*\C$ is also a restricted Segal family.
\end{corr}

\proof{} By virtue of Proposition~\ref{restr.seg.prop} and
Proposition~\ref{seg.prop}~\thetag{ii}, we have $\C \cong
L^*(\C^\dm) \cong L^*\pi_*\Nn_\idot^*\C^\Delta$ for some small
additive semiexact family of groupoids $\C^\Delta$ over $\langle
\Delta_f^o\Delta_f^o\Sets,W \rangle$ that is a Segal family in the
sense of Definition~\ref{bi.seg.def}. Equivalently, $\C \cong
\Nn^*L^*\C^\Delta$, where $L$ is the functor \eqref{LR.simp}. But if
we let $\iota_L:\Delta \times \Delta \to \Delta \times \Delta$ be
the involution $\iota$ acting on the Reedy factor, and denote
$\iota_\Delta = \iota_L^*:\Delta^o_f\Delta^o_f\Sets \to
\Delta^o_f\Delta^o_f\Sets$, then $L \circ \Nn \circ \iota \cong
\iota_\Delta \circ \Nn \circ \iota$, so that $\iota^*\C \cong
L^*\pi_*\Nn_\idot^*\iota_\Delta^*\C^\Delta$. It remains to observe
that Definition~\ref{bi.seg.def} is manifestly invariant under
$\iota_\Delta$, so that $\iota_\Delta^*\C^\Delta$ is also a small
Segal family over $\langle \Delta_f^o\Delta^o_f\Sets,W \rangle$, and
we are done by Proposition~\ref{seg.prop}~\thetag{i} and
Proposition~\ref{restr.seg.prop}.
\endproof

\subsection{Barycentric dualization.}\label{bary.dua.subs}

While our proof of Corollary~\ref{iota.corr} is short, it is
somewhat unsatisfactory. Using the full force of
Proposition~\ref{seg.prop} looks like an overkill, and comes with a
penalty: we have to assume that $\C$ has essentially small
fibers. Let us now extract the relevant parts of
Proposition~\ref{seg.prop} to give a self-contained proof of a
version of Corollary~\ref{iota.corr}; as a bonus, it works for
families with large fibers, and includes a statement about semicanonical
extensions.

The main technology used are the barycentric subdivision functors
$B^\dm$, $B_\dm$ of Example~\ref{bary.rel.exa}. We recall that
$B^\dm$ and $B_\dm$ send $\Posf$ into $\Relf$ and the whole $\Pos$
into the subcategory $\Relf^\pm \subset \Rel$ of left-finite
biordered sets, and for any $J$ in $\Posf^\pm$, the biordered map
$\xi:B^\dm(J) \to L(J)$ is $\pm$-bianodyne in the sense of
Definition~\ref{pm.biano.def}, and strongly $\pm$-bianodyne if $J
\in \Posf$. For any partially ordered set $J$ equipped with a map
$\chi:J \to [1]$, the comma-set $J/[1]$ fits into the
cocartesian square \eqref{cyl.o.sq}. We denote by
\begin{equation}\label{B.J.1}
B(J,\chi) = B(J_0 \times [1]) \copr_{B(J_0)} B(J)
\subset B(J/[1])
\end{equation}
the corresponding left-closed subset in its barycentric subdivision,
and let
\begin{equation}\label{B.J.1.dm}
  B^\dm(J,\chi) \subset B^\dm(J/[1]), \quad
  B_\dm(J,\chi) \subset B_\dm(J/[1])
\end{equation}
be the subset $B(J,\chi)$ with the biorders induced by $B^\dm$ and
$B_\dm$.

\begin{lemma}\label{bary.cyl.le}
For any $J$ in $\Posf$ equipped with a map $\chi:J \to [1]$, the
embeddings \eqref{B.J.1.dm} are strongly $\pm$-bianodyne.
\end{lemma}

\proof{} For any $n > l > 0$, denote by $B([n])_l \subset B([n])$ the
left-closed subset corresponding to the horn embedding
\eqref{horn.eq}, let $B^\dm([n])_l$ resp.\ $B_\dm([n])_l$ be
$B([n])_l$ with the biorder induced by $B^\dm([n])$
resp.\ $B_\dm([n])$, and let $v^{l\dm}_n:B^\dm([n])_l \to
B^\dm([n])$, $v^l_{n\dm}:B_\dm([n])_l \to B_\dm([n])$ be the
embeddings. Then we have the strongly $\pm$-bianodyne map
\eqref{b.dm.n.l} of Example~\ref{beta.rel.exa}, and the same
induction as in Lemma~\ref{bva.le} shows that $v^{l\dm}_n$ is also
strongly $\pm$-bianodyne. Moreover, we can identify $B_\dm([n]) =
B^\dm([n]^o) \cong B^\dm([n])$, and then \eqref{b.dm.n.l} for the
opposite ordinals is also strongly $\pm$-bianodyne, so that
$v^l_{n\dm}$ is strongly $\pm$-bianodyne as well.  Now, the same
argument as in the proof of Lemma~\ref{bva.le}~\thetag{iv} provides
a finite filtration $B(J/[1])^d_k$ by left-closed subsets on the
target of the embedding \eqref{B.J.1} that starts with its source
and ends with the whole target, and this provides filtrations on the
biordered embeddings \eqref{B.J.1.dm}. The successive right-closed
complements in the filtration are of the form \eqref{B.S.km}, and
the corresponding projections $\overline{B}^\dm(J/[1])^d_k \to S$
resp.\ $\overline{B}_\dm(J/[1])^d_k \to S$ extend to the characteristic
maps $p:B^\dm(J/[1])^d_k \to R(S^<)$ resp.\ $p:B_\dm(J/[1])^d_k \to
R(S^<)$ of \eqref{p.J.sq}. Then with respect to these maps, the left
comma-fibers of the embedding of the previous term in the filtration
are $\id$ over $o \in S^<$ and a standard pushout of the horn
embedding $v^{d\dm}_k$ resp.\ $v^d_{k\dm}$ over each $s \in S
\subset S^<$, so we are done by Lemma~\ref{rel.R.le}.
\endproof

\begin{corr}\label{bary.cyl.corr}
For any map $\chi:J \to \N$ in $\Posf^+$, define a left-closed
subset $B(J/\N,\chi) \subset B(J/\N)$ by the cartesian square
\begin{equation}\label{B.N.sq}
\begin{CD}
  B(J/\N,\chi) @>>> B(J/\N)\\
  @V{p}VV @VVV\\
  \ZZ_\infty @>{\beta}>> B(\N),
\end{CD}
\end{equation}
where $\beta$ is embedding \eqref{zeta.xi}, and let
$B_\dm(J/\N,\chi) \subset B_\dm(J/\N)$ be $B(J/\N,\chi)$ with the
induced biorder. Then the embedding $B_\dm(J/\N,\chi) \to
B_\dm(J/\N)$ is $\pm$-bianodyne.
\end{corr}

\proof{} The height map $\hht:J \to \N$ is left-bounded, it induces
a left-finite map $\hht \circ \sigma \circ \xi:B(J/\N) \to J/\N \to
J \to \N$, and for any $n \in \N$, we have $B_\dm(J/\N)/n \cong
B_\dm((J/_{\hht}n)/\N)$, so by Lemma~\ref{rel.R.le}, we may assume
right away that $J \in \Posf$. Then $\tau:J/\N \to \N$ is
left-bounded and induces a left-finite map $\tau \circ \xi:B(J/\N)
\to J/\N \to \N$, and with respect to this map, we have $B(J/\N)/n
\cong B((J/n)/[n])$ for any $n \in \N$. Therefore it suffices to
show that for any $J \in \Posf$ equipped with a map $\chi:J \to [n]$
for some $[n]$, if we define $B(J/[n],\chi) \subset B(J/[n])$ by
\eqref{B.N.sq} with $\beta$ replaced by $\beta_n$ of
\eqref{zeta.n.xi}, then the embedding $B_\dm(J/[n],\chi) \to
B_\dm(J/[n])$ is $\pm$-anodyne. For $n=1$, this is
Lemma~\ref{bary.cyl.le}, and the general case immediately follows by
induction, using the decomposition \eqref{J.pl.dec}.
\endproof

\begin{remark}
To visualize better the biordered set $B_\dm(J/\N,\chi)$ of
Corollary~\ref{bary.cyl.corr}, we note that the left comma-fibers of
the projection $p$ of \eqref{B.N.sq} over the elements $2n,2n\!+\!1
\in \ZZ_\infty$, $n \geq 0$ are given by
\begin{equation}\label{B.N.f}
B_\dm(J/\N,\chi)/2n \cong B_\dm(J/n), \ B_\dm(J/\N,\chi)/(2n\!+\!1)
\cong B_\dm(J / e_n([1])),
\end{equation}
where $e_n:[1] \to \N$ is the embedding onto $\{n,n+1\} \subset \N$.
\end{remark}

\begin{prop}\label{iota.prop}
For any restricted Segal family $\C$ over $\Posf$, with unfolding
$\C^\dm$, $B_\dm^*\C$ is a restricted Segal family over $\Posf$, and
the bianodyne map $\xi$ of Corollary~\ref{B.biano.corr} induces an
equivalence $B_\dm^*\C^\dm \cong \iota^*\C$. Moreover, for any
semicanonical extension $\C^+$ of a restricted Segal family $\C$,
$B_\dm^*\C^{+\dm}$ is a semicanonical extension of $\iota^*\C \cong
B_\dm^*\C^\dm$.
\end{prop}

\proof{} The fact that $\C^\dm$ is constant along the map $\xi$
immediately follows from Proposition~\ref{restr.seg.prop} and
Lemma~\ref{bi.del.le}, so we indeed have $\iota^*\C \cong
B_\dm^*\C^\dm$. Since $B_\dm$ preserves coproducts and standard
pushout square, the family $B_\dm^*\C$ is additive and semiexact. To
see that it is reflexive and non-degenerate, recall that for any $J
\in \Posf$, we have the functorial strongly $\pm$-bianodyne map
$B_\dm(J \times [1]) \to B_\dm(J) \times L([1])$ of
Example~\ref{B.i.exa}. To see that $B_\dm^*\C^\dm$ satisfies the
cylinder axiom, consider a square \eqref{J.cyl.sq} in $\Posf$, and
use Lemma~\ref{bary.cyl.le} and Lemma~\ref{pm.B.le} to reduce the
question to showing that $\C^\dm$ is cartesian along the square
$$
\begin{CD}
  B_\dm(J_0) \times L([1]) @>>> B_\dm(J_0)\\
  @VVV @VVV\\
   (B_\dm(J_0) \times L(\V)) \copr_{B_\dm(J_0)} B_\dm(J) @>>>
  (B_\dm(J_0) \times L([1])) \copr_{B_\dm(J_0)} B_\dm(J).
\end{CD}
$$
Then since $\C$ is reflexive and non-degenerate, $\C^\dm$ is fully
faithful along the top arrow, so as in
Proposition~\ref{restr.seg.prop}, it suffices to observe that
$\C^\dm$ is constant along the map
$$
(B_\dm(J_0) \times \V^{\dm o}) \copr_{B_\dm(J_0)} B_\dm(J) \to
  (B_\dm(J_0) \times L([1])) \copr_{B_\dm(J_0)} B_\dm(J)
$$
that happens to be a standard pushout of a reflexive
map. Analogously, to check that $B_\dm^*\C^\dm$ satisfies
weak excision, assume given a square \eqref{J.S.sq} in $\Posf$, and
note that by Corollary~\ref{pm.B.corr}, it suffices to check that
$\C^\dm$ is cartesian along the square
$$
\begin{CD}
  B_\dm(J_o) \times B_\dm(S^>) @>>> B_\dm(J_o)\\
  @VVV @VVV\\
  B_\dm(J)^\hushf @>>> B_\dm(J),
\end{CD}
$$
where $B^\dm(J^o)^\hushf$ is as in \eqref{B.S.J.bi}. As in
Proposition~\ref{restr.seg.prop}, this immediately follows from
Lemma~\ref{J.hh.le}. This proves that $B_\dm^*\C^\dm \cong
\iota^*\C$ is indeed a restricted Segal family. Finally, if $\C^+$
is a semicanonical extension of the family $\C$, then $B^*_\dm\C^{+\dm}$ is
again reflexive and non-degenerate by Example~\ref{B.i.exa}, so it
suffices to prove that it is weakly semicontinuous. We thus need to
show that for any $J \in \Posf^+$ equipped with a map $\chi:J \to
\N$, the unfolding $\C^{+\dm}$ is cartesian along the square
\begin{equation}\label{B.J.pl.sq}
\begin{CD}
\coprod_{n \geq 1} B_\dm([1] \times J/n) @>>> \coprod_{n \geq
  1}B_\dm(J/n)\\
@VVV @VVV\\
B_\dm(J^+) @>>> B_\dm(J/\N)
\end{CD}
\end{equation}
obtained by applying $B_\dm$ to the square \eqref{J.pl.sq}. The left
vertical arrow in \eqref{B.J.pl.sq} fits into a standard pushout
square
$$
\begin{CD}
\coprod_n B_\dm(J/n) @>>> \coprod_{n \geq 1} B_\dm([1] \times
J/n)\\
@VVV @VVV\\
\coprod_n B_\dm(J/ e_n([1])) @>>> B_\dm(J^+)
\end{CD}
$$
induced by \eqref{Z.NN.sq}, where $e_n$ is as in \eqref{B.N.f},
while in the top arrow in \eqref{B.J.pl.sq}, each term in the
coproduct factors as
$$
\begin{CD}
B_\dm([1] \times J/n) @>{a}>> B_\dm([1] \times J/n|J/n) @>{b}>>
B_\dm(J/n),
\end{CD}
$$
where $a$ is the natural dense map, and $b$ is the $\pm$-bianodyne
map of (the $\pm$-bianodyne version of)
Corollary~\ref{bary.rel.rel.corr}. Thus if we simply notation by
writing $B_\dm(J/n)'=B_\dm([1] \times J/n|J/n)$, and define a
biordered set $B_\dm(J^+)'$ by a standard pushout square
$$
\begin{CD}
\coprod_n B_\dm(J/n) @>>> \coprod_{n \geq 1} B_\dm(J/n)'\\
@VVV @VVV\\
\coprod_n B_\dm(J/e_n([1])) @>>> B_\dm(J^+)',
\end{CD}
$$
then \eqref{B.J.pl.sq} fits into a commutative diagram
\begin{equation}\label{B.big.dia}
\begin{CD}
\coprod_{n \geq 1} B_\dm([1] \times J/n) @>>> B_\dm(J^+)\\
@V{a}VV @VV{a'}V\\
\coprod_{n \geq 1} B_\dm(J/n)' @>>> B_\dm(J^+)'\\
@V{b}VV @VV{b'}V\\
\coprod_{n \geq 1} B_\dm(J/n) @>>> B_\dm(J/\N),
\end{CD}
\end{equation}
where $a$ and $a'$ are dense, the unfolding $\C^{+\dm}$ is cartesian
along the top square by the same argument as in Lemma~\ref{E.0.le},
and $b$ is a coproduct of $\pm$-bianodyne maps, thus
$\pm$-bianodyne. Thus to finish the proof, it suffices to check that
$b'$ is also $\pm$-bianodyne. But it factors as
$$
\begin{CD}
  B_\dm(J^+)' @>{c}>> B_\dm(J/\N,\chi) @>{d}>> B_\dm(J/\N),
\end{CD}
$$
the map $d$ is $\pm$-bianodyne by Corollary~\ref{bary.cyl.corr}, so
it further suffices to check that so is $c$. Now, if we treat it as
a map over $\ZZ_\infty$ via the projection $p$ of \eqref{B.N.sq},
then in terms of \eqref{B.N.f}, its left comma fiber $c/2n$ is the
$\pm$-bianodyne projection $B_\dm(J/n)' \to
B_\dm(J/n)$, while its left comma-fiber $c/(2n\!+\!1)$ is the
projection
$$
\begin{CD}
B_\dm(J/n)' \copr_{B_\dm(J/n)} B_\dm(J/e_n([1]))
  \copr_{B_\dm(J/(n\!+\!1))} B_\dm(J/(n\!\!+\!\!1))'\\ 
  @VVV\\
B_\dm(J/ e_n([1])),
\end{CD}
$$
so that it is inverse to a composition of two standard pushouts of
bianodyne maps. Thus we are done by by Lemma~\ref{rel.R.le}.
\endproof

\subsection{Complete Segal families.}

We now introduce one more condition on families of groupoids over an
ample full subcategory $\I \subset \Pos$. By
Definition~\ref{amp.def}, $\I$ contains $\Delta \subset \Posff
\subset \Pos$. As in \eqref{compl.sq}, denote by $\mu:\{0,1\} \times [1]
\to [3]$ the map given by $\mu(l \times l') = l + 2l'$,
$l,l'=0,1$, and for any $[n] \in \Delta$, let $e:[n] \to [0]$ be the
tautological projection. We then have the following version of
Definition~\ref{exa.compl.def}.

\begin{defn}\label{pos.compl.def}
A family of groupoids $\C$ over an ample $\I \subset \Pos$ is {\em
  complete} if it is cartesian over the square
\begin{equation}\label{pos.compl.sq}
\begin{CD}
\{0,1\} \times [1] \times J @>{m \times \id}>> [3] \times J\\
@V{\id \times e \times \id}VV @VV{e \times \id}V\\
\{0,1\} \times J @>>> J
\end{CD}
\end{equation}
for any partially ordered set $J \in \I$.
\end{defn}

\begin{exa}\label{I.triv.compl.exa}
For any category $I$, the family $\bpi:\Posf^+ / I \to
\Posf^+$ of Example~\ref{I.triv.seg.exa} is complete.
\end{exa}

We also need a version of Definition~\ref{pos.compl.def} for
bisimplicial sets. Recall that for any $[n] \in \Delta$, $\Nn([n])
\cong \Delta_n$ is an elementary $n$-simplex, and then the square
\eqref{pos.compl.sq} with $J = \ppt$ induces a cocartesian square
\begin{equation}\label{bi.compl.sq}
\begin{CD}
\{0,1\} \times \Delta_1 @>>> \Delta_3\\
@VVV @VVV\\
\{0,1\} @>>> \V_3
\end{CD}
\end{equation}
for some $\V_3 \in \Delta^o_f\Sets$. Geometrically, the simplicial
set $\V_3$ is a tetrahedron with two edges contracted to points.

\begin{defn}\label{bi.compl.def}
A family of groupoids $\C$ over $\Delta^o_f\Delta^o_f\Sets$ is {\em
  complete} if for any $n \geq 0$, $\C$ is stably constant along the
projection $L(\V_3 \times \Delta_n) \to L(\Delta_n)$ in the sense of
Definition~\ref{st.I.def}. A fibrant Segal space $X \in
\Delta^o\Delta^o\Sets$ in the sense of Example~\ref{seg.sp.exa} is
      {\em complete} if so is the family $\Hh(X)$.
\end{defn}

\begin{exa}\label{css.sp.exa}
If we let $\Hhom(-,-)$ be the internal $\Hom$ in the
cartesi\-an-closed category $\Delta^o\Delta^o\Sets$, then a fibrant
Segal space $Y \in \Delta^o\Delta^o\Sets$ is complete if and only if
$\Hhom(-,Y)$ sends the projection $L(\V_3) \to \ppt = L(\Delta_0)$
to a weak equivalence. Indeed, choose a section $v:\ppt \to \V_3$ of
the projection $\V_3 \to \ppt$; then since $v$ is injective and $Y$
is fibrant, $\Hhom(-,Y)$ always sends $v$ to a pointwise Kan
fibration, and by Corollary~\ref{W.C.corr}, this fibration is
pointwise-trivial iff $\Hh(Y)$ is constant along $L(v
\times \id):L(\Delta_n) \to L(\Delta_n \times \V_3)$ for any $n \geq
0$.
\end{exa}

\begin{exa}\label{compl.sq.exa}
For another interpreration of complete Segal spaces, note that the
square \eqref{bi.compl.sq} is also homotopy cocartesian in
$\Delta^o\Delta^o\Sets$. Let $\delta:\Delta \to \Delta \times
\Delta$ be the diagonal embedding, and for any $X:\Delta^o \to
\Delta^o\Sets$ and $[n],[m] \in \Delta$, denote $X_{[n],[m]} =
\delta^o_*X([n] \times [m]) \in \Delta^o\Sets$. Then a fibrant Segal
space $X$ is complete if and only if for any $[n] \in \Delta$, the
square
\begin{equation}\label{bi.compl.bis.sq}
\begin{CD}
X_{[n],[0]}  @>>> X_{[n],[0]} \times X_{[n],[0]}\\
@VVV @VVV\\
X_{[n],[3]} @>>> X_{[n],[1]} \times X_{[n],[1]}
\end{CD}
\end{equation}
induced by \eqref{bi.compl.bis.sq} is homotopy cartesian in
$\Delta^o\Sets$.
\end{exa}

\begin{exa}\label{css.ni.exa}
Let $NI$ be the nerve of a small category $I$, constant in the
simplicial direction. Then a map $\V_3 \to I$ is a diagram
\eqref{i.dot} of length $4$ such that the maps $i_0 \to i_2$, $i_1
\to i_3$ are invertible, and then all maps in the diagram must be
invertible. If $I$ is rigid, any map $\V_3 \to I$ must be therefore
constant. In this case, the fibrant Segal space $L(NI)$ is complete:
indeed, since the product of \eqref{bi.compl.sq} with any
$\Delta_n$, $n \geq 0$ is still cocartesian, a map $\V_3 \times
\Delta_n \to NI$ is then the same thing as a functor $[3] \times [n]
\to I$ whose restriction to any $[3] \times \{l\}$, $l \in [n]$ is
constant, and such a functor factors through the projection $[3]
\times [n] \to [n]$.
\end{exa}

\begin{lemma}\label{css.le}
Assume given a small additive semiexact Segal family of group\-oids
$\C \to \Delta_f^o\Delta^o_f\Sets$ constant along anodyne maps of
Definition~\ref{f.CW.def} that is complete in the sense of
Definition~\ref{bi.compl.def}. Then $\C$ is stably constant along
the projection $L(\V_3 \times X) \to L(X)$ for any simplicial set $X
\in \Delta^o_f\Sets$
\end{lemma}

\proof{} Same as Corollary~\ref{bi.seg.corr}.
\endproof

\begin{prop}\label{css.prop}
Assume given a small semiexact additive family of groupoids $\C \to
\Delta^o_f\Delta^o_f\Sets$ that is constant along weak equivalences
and satifies the Segal condition of
Definition~\ref{bi.seg.def}. Then $\C$ is complete in the sense of
Definition~\ref{bi.compl.def} if and only if $\Nn^*L^*\C$ is a
restricted Segal family in the sense of Definition~\ref{red.seg.def}
and complete in the sense of
Definition~\ref{bi.compl.def}. Moreover, if this happens, then
$\pi_*\Nn_\idot^*\C$ is reflexive in the sense of
Definition~\ref{rel.refle.def}.
\end{prop}

\proof{} Note that $L \circ \Nn \cong \Nn_\dm \circ L \cong \bN_\dm
\circ L$, and $\bN_\dm^*\C$ is a Segal family over $\Relf$ by
Lemma~\ref{NQ.bi.le}. Assume first that $\Nn^*L^*\C \cong
L^*\bN_\dm^*\C$ is a complete restricted Segal family. Then
$\Nn^*L^*\C \cong L^*\pi_*\Nn_\idot^*\C$ by
Proposition~\ref{seg.prop}, and $\pi_*\Nn_\idot^*\C \cong
(\Nn^*L^*\C)^\dm$ by Proposition~\ref{restr.seg.prop}, so that in
particular, $\pi_*\Nn_\idot^*\C$ is cartesian along the square
\eqref{pos.compl.sq} for any $J \in \Relf$ (where we replace $[0]$,
$\{0,1\}$, $\{0,1\} \times [1]$ and $[3]$ with their images under
the functor $L$). Moreover, let $v:\Delta^o_f\Delta^o_f\Sets \to
\Delta^o_f\Delta^o_f\Sets$ be the functor sending $X$ to $X \times
L(\V_3)$. Then by Corollary~\ref{bi.seg.corr}, $v^*\C$ is also an
additive semiexact Segal family constant along weak equvalences, so
that $\pi_*\Nn_\idot^*v^*\C$ is also a Segal family over
$\Relf$. Now note that by Lemma~\ref{f.CW.le}, $\C$ is strongly
semiexact over $\Delta^o_f\Delta^o_f\Sets$, and in particular, it is
semiexact over the square \eqref{bi.compl.sq} and its cartesian
product with any $Y \in \Delta^o_f\Delta^o_f\Sets$. Since
$\pi_*\Nn_\dm^*\C$ is cartesian along \eqref{pos.compl.sq}, this
implies that the functor $\pi_*\Nn_\dm^*\C \to \pi_*\Nn_\dm^*v^*\C$
induced by the tautological projection $v \to \id$ has a one-sided
inverse $\pi_*\Nn_\dm^*v^*\C \to \pi_*\Nn_\dm^*\C$ that is an
epivalence. Since the source and target of this epivalence are
semiexact families over $\Relf$ with respect to the standard
CW-structure, it is an equivalence by Corollary~\ref{semi.epi.corr},
and its inverse $\pi_*\Nn_\dm^*\C \to \pi_*\Nn_\dm^*v^*\C$ is then
also an equivalence. Therefore $\C \to v^*\C$ is an equivalence by
Proposition~\ref{seg.prop}~\thetag{i}, and $\C$ is complete.

Conversely, assume that $\C$ is complete. Then again, by strong
semiexactness, and Lemma~\ref{css.le}, $\Nn^*L^*\C$ is semicartesian
along \eqref{pos.compl.sq} for any $J \in \Posf$. However, since
$e:[3] \to [0]$ admits a one-sided inverse, $(\id \times e)^*$ in
\eqref{pos.compl.sq} is faithful, so that $\Nn_L^*\C$ is in fact
cartesian along \eqref{pos.compl.sq}. By
Proposition~\ref{restr.seg.prop}, to finish the proof, it suffices
to show that the Segal family $\pi_*\Nn_\idot^*\C$ is reflexive in
the sense of Definition~\ref{rel.refle.def}.

To do this, take any $J \in \Relf$, consider the functor
$j:\Delta \to \Relf$, $[n] \mapsto J \times L([n])$, and let $\C(J)
= j^*\pi_*\Nn_\idot^*\C$ be the induced family of groupoids over
$\Delta$. Since $\pi_*\Nn_\idot^*\C$ is a Segal family, $\C(J)$ is
semicartesian along the square \eqref{seg.del.sq} for any $n \geq 0$
and $l=1$ (let $s_\dg:[n] \to [2]$ be the cofibration of
Example~\ref{n.refl.exa}, and apply Corollary~\ref{cyl.2.corr} to
the corresponding bicofibration $J \times L([n]) \to L([n]) \to
L([2])$). Moreover, as in the proof of
Proposition~\ref{exa.seg.prop}, let $\C(J)^e \subset \C(J)$ be the
full subcategory spanned by essential images of the functors
$e^*:\C(J)_{[0]} \to \C(J)_{[n]}$. Then by Lemma~\ref{retra.le}
applied to the standard pushout square \eqref{cyl.2.eq}, $\C(J)^e$
is actually cartesian along \eqref{seg.del.sq} for $l=1$ and any $n
\geq 0$. By induction on $n$, $\C(J)^e$ is then cartesian along
\eqref{seg.del.sq} for any $n \geq l \geq 0$, so that it is a Segal
category in the sense of Definition~\ref{seg.cat.def}. Now the same
argument as for $\Nn^*L^*\C$ shows that $\C(J)$ is complete, and
$\C(J)^e \subset \C(J)$ is then also complete. Therefore by
Proposition~\ref{exa.seg.prop}~\thetag{ii}, it is reflexive, and
since $J$ was arbitrary, $\pi_*\Nn_\dm^*\C$ is also reflexive.
\endproof

As a corollary of Proposition~\ref{restr.seg.prop} and
Proposition~\ref{css.prop}, we note that if we consider the category
$\Seg$ of \eqref{Seg.eq}, and let $\CSS \subset \Seg$ be the full
subcategory spanned by unfoldings of small complete restricted Segal
families, then the equivalence of Corollary~\ref{seg.corr}
identifies $\CSS$ with the full subcategory in
$h^W(\Delta^o\Delta^o\Sets)$ formed by complete Segal spaces. In
fact, one can do much better: by a beautiful discovery of Ch. Rezk,
there is another model structure on $\Delta^o\Delta^o\Sets$ whose
homotopy category is $\CSS$ on the nose. We will not really need
this construction, but we will need one of its corollaries: the
subcategory $\CSS \subset h^W(\Delta^o\Delta^o\Sets)$ is
left-admissible. For completeness, let us reprove this. We will
actually prove slightly more (namely, a weak version of the
so-called ``Bousfield localization'', \cite{bou1}).

Assume given a cellular Reedy category $I$, and a set $T$ of arrows
$Y \to Y'$ in $I^o\Sets$ whose sources $Y$ are compact. Say that a
semicanonical extension $\C^+$ of a small additive semiexact family
of groupoids $\C \to I^o_f\Delta^o_f\Sets$ is a {\em $T$-family} if
it is stably constant along any $f \in T$ in the sense of
Definition~\ref{st.I.def}, say that a fibrant object $X \in
I^o\Delta^o\Sets$ is a {\em $T$-space} if $\Hh(X)$ is a $T$-family,
and say that a map $g$ in $I^o\Delta^o\Sets$ is a {\em
  $T$-equivalence} if for any $T$-space $X$, $\Hh(X)$ is stably
constant along $g$. By definition, the class $W(T)$ of injective
$T$-equivalences contains $T$, but it is much larger; in particular,
by Lemma~\ref{f.CW.le} and Corollary~\ref{W.C.corr}, it is saturated
and closed under coproducts, pushouts and countable compositions,
and we also have the following easy observation.

\begin{lemma}\label{T.pb.le}
For any cellular Reedy category $I'$ and Reedy-fibrant functor
$X:{I'}^o \to I^o\Delta^o\Sets$ such that $X(i)$ is a $T$-space for
any $i \in I'$, $\lim_{{I'}^o}X$ is a $T$-space.
\end{lemma}

\proof{} By Corollary~\ref{W.C.corr}~\thetag{ii}, a fibrant $X \in
I^o\Delta^o\Sets$ is a $T$-space iff for any $f:Y \to Y'$ in $T$,
the Kan fibration $f^*:\Hom^\Delta(Y',X) \to \Hom^\Delta(Y,X)$ is
trivial, and then as in Example~\ref{holim.I.exa}, by Quillen
Adjunction Theorem~\ref{qui.adj.thm}, the functor $\lim_{{I'}^o}$
sends weak equivalences between Reedy-fibrant objects to weak
equivalences.
\endproof

\begin{exa}\label{seg.eq.exa}
Let $I = \Delta$, and let $T$ be the class of horn extensions
$v^l_n$ of \eqref{horn.eq}, with $n > l > 0$. Then a $T$-family is a
Segal family in the sense of Definition~\ref{bi.seg.def}, and a
$T$-space is a fibrant Segal space of Example~\ref{seg.sp.exa}. A
{\em Segal equivalence} is a $T$-equivalence for this particular set
$T$.
\end{exa}

\begin{exa}\label{css.eq.exa}
Let $I = \Delta$, consider the set $\V_3$ of \eqref{bi.compl.sq},
choose a map $v:\ppt \to \V_3$, and let $T'$ be the set $T$ of
Example~\ref{seg.eq.exa} enlarged by adding the embeddings $\id
\times v:\Delta_n \to \Delta_n \times \V_3$ for all $n \geq 0$. Then
$T'$-spaces are complete Segal spaces of
Definition~\ref{bi.compl.def}.
\end{exa}

\begin{lemma}\label{bous.le}
Assume given a cellular Reedy category $I$, and a set $T$ of arrows
$f:Y \to Y'$ in $I^o\Sets$ whose source $Z$ is compact.
\begin{enumerate}
\item Any $T$-equivalence $f:X \to X'$ between $T$-spaces is a weak
  equivalence.
\item The full subcategory $h_T^W(I^o\Delta^o\Sets) \subset
  h^W(I^o\Delta^o\Sets)$ spanned by $T$-spaces is left-admissible.
\end{enumerate}
\end{lemma}

\proof{} For \thetag{i}, since $f$ is a $T$-equivalence, the family
$\Hh(X)$ is constant along $f$, and so is the represented functor
$h(X) = \pi_0(\Hh(X))$. This provides a map $g:X' \to X$ such that
$g \circ f = \id$ in $h^W(\Delta^o\Delta^o\Sets)$. The same argument
applied to $g$ then provides a map $f':X \to Y'$ such that $f' \circ
g = \id$, and then $f = f' \circ g \circ f = f'$.

For \thetag{ii}, by Lemma~\ref{p.le}, we have to construct a functor
\begin{equation}\label{p.h.eq}
p:h^W(I^o\Delta^o\Sets) \to h^W_T(I^o\Delta^o\Sets) \subset
h^W(I^o\Delta^o\Sets)
\end{equation}
and a functorial map $\alpha:X \to p(X)$ invertible when $X \in
h^W_T(I^o\Delta^o\Sets)$. Let $T_\idot$ be the set of maps
$f(\SS_{n-1},\Delta_n)$ of Corollary~\ref{W.C.corr}~\thetag{ii},
for all $f \in T$ and $n \geq 0$, and enlarge $T_\idot$ further to a
set $T'_\idot$ by adding all the horn embeddings $\V^k_{i,n}$ of
\eqref{s.v.eq}. Then since all injective maps in $I^o\Delta^o\Sets$
are trivially absolutely cofibrant in the sense of
Definition~\ref{abs.C.def} for the standard $C$-category structure,
we have the extension functor $\Ex(T'_\idot)$ of \eqref{ex.I.sq},
and as in Lemma~\ref{simp.S.le}, for any $X \in I^o\Delta^o\Sets$,
we can consider $X_1 = \Ex(T'_\idot)(X)$ equipped with a map $a_0:X
= X_0 \to X_1$, and then iterate the construction to obtain
functorial maps $a_n:X_n \to X_{n+1}$. Then we take $X_\infty=
\colim_nX_n$, with the functorial map $a_\infty:X \to X_\infty$, and
we note that by Corollary~\ref{W.C.corr} and \eqref{ex.I.eq}, all
the maps $a_n$ are $T$-equivalences, and then again by
Corollary~\ref{W.C.corr}, so is the map $a_\infty$. Moreover, the
sources of all maps $f \in T'_\idot$ are compact in
$I^o\Delta^o\Sets$, so the small object argument applies and shows
that $X_\infty$ has the lifting property with respect to all such
maps $f$. Then by Lemma~\ref{bisimp.le} and
Corollary~\ref{W.C.corr}~\thetag{ii}, $X_\infty$ is a $T$-space. If
$X$ is also a $T$-space, then $a_\infty$ is a weak equivalence by
\thetag{i}, so to finish the proof, it suffices to show that the
functor $X \mapsto X_\infty$ sends weak equivalences to weak
equivalences, thus descends to a functor \eqref{p.h.eq} (with
$\alpha$ then given by $a_\infty$). But since the class of
$T$-equivalences is saturated, for any weak equivalence $w:X \to
X'$, the corresponding map $w_\infty:X_\infty \to X'_\infty$ is a
$T$-equivalence, and then it must be a weak equivalence by
\thetag{i}.
\endproof

\section{Augmentations.}\label{aug.sec}

\subsection{Augmented Segal spaces.}

For any category $I$, define an {\em $I$-bisimp\-li\-cial set} as a
functor $I^o \times \Delta^o \times \Delta^o \to \Sets$, or
equivalently, a functor $I^o \to \Delta^o\Delta^o\Sets$. Assume that
$I$ is a directed $\Hom$-finite cellular Reedy category $I$, as in
Proposition~\ref{pos.I.prop}. Then the category
$I^o\Delta^o\Delta^o\Sets$ carries a Reedy model structure. As in
Subsection~\ref{del.seg.subs}, we distinguish the Reedy and the
simplicial factor $\Delta$ in $I \times \Delta \times \Delta$, and
we let $r:I^o \times \Delta^o \times \Delta^o \to \Delta^o$ be the
projection onto the simplicial factor, with the corresponding
functor $R = r^*:\Delta^o\Sets \to I^o\Delta^o\Delta^o\Sets$. Each
$I$-bisimplicial set $X$ has three skeleton filtrations
$\sk^I_\idot$, $\sk^L_\idot$, $\sk^R_\idot$ along $I$ and the Reedy
resp.\ simplicial copy of $\Delta$, and we let $\sk^{I,L}_\idot$
resp.\ $\sk^{I,R}_\idot$ be the skeleton filtration along the
product $I \times \Delta$ of $I$ and the Reedy resp.\ simplicial
copy of $\Delta$.  For any $I$-bisimplicial set $X$ and any object
$i \in I$, evaluation at $i$ gives a bisimplicial set $X(i) =
\eps(i)^{o*}X$, and since we also have the left Kan extension
functor $\eps(i)^o_!:\Delta^o\Delta^o\Sets \to
I^o\Delta^o\Delta^o\Sets$, a family of groupoids $\C \to
I^o\Delta^o\Delta^o\Sets$ restricts to a family
$\C(i)=(\eps(i)^o_!)^*\C \to \Delta^o\Delta^o\Sets$. We note that
explicitly, by \eqref{kan.eq}, we have
\begin{equation}\label{kan.X}
  \eps(i)^o_!X \cong \Delta_i \boxtimes X,
\end{equation}
where $\Delta_i = \Y(i):I^o \to \Sets$ is represented by $i$, as in
\eqref{sph.I.eq}. Since $I \times \Delta$ is also a cellular Reedy
category, the category of $I$-bisimplicial sets carries the Reedy
model structure, and then Proposition~\ref{Hh.prop} immediately
implies that for any fibrant $I$-bisimplicial set $X$ and object $i
\in I$, we have $\Hh(X)(i) \cong \Hh(X(i))$, where $X(i)$, being an
evaluation of a Reedy-fibrant object, is automatically
fibrant. Moreover, since $I$ is directed, $\eps(i)^o_!$ sends
$\Delta^o_f\Delta^o_f\Sets$ into $I^o_f\Delta^o_f\Delta^o_f\Sets$,
and for any uncountable regular $\kappa > |I|$, it also
sends $\Delta^o\Delta^o\Sets_\kappa$ into
$I^o\Delta^o\Delta^o\Sets_\kappa$. Therefore for any family of
groupoids $\C$ over $I^o_f\Delta^o_f\Delta^o_f\Sets$
resp.\ $I^o_f\Delta^o_f\Delta^o_f\Sets_\kappa$, we have a
well-defined family of groupoids $\C(i)$ over
$\Delta^o_f\Delta^o_f\Sets$
resp.\ $\Delta^o_f\Delta^o_f\Sets_\kappa$.

\begin{defn}\label{I.bi.seg.def}
An {\em $I$-Segal space} is an $I$-bisimplicial set $X$ such that
for any $i \in I$, $X(i)$ is a Segal space in the sense of
Example~\ref{seg.sp.exa}. A family of groupoids $\C$ over
$I^o_f\Delta^o_f\Delta^o_f\Sets$, or
$I^o_f\Delta^o_f\Delta^o_f\Sets_\kappa$ for some uncountable regular
cardinal $\kappa > |I|$, is a {\em Segal family} if for any $i \in
I$, $\C(i)$ is a Segal family in the sense of
Definition~\ref{bi.seg.def}. A fibrant $I$-Segal space $X$ is {\em
  complete} if so is $X(i)$ for any $I \in I$, and a Segal family
$\C$ is {\em complete} if for any $i \in I$, the Segal family
$\C(i)$ is complete in the sense of Definition~\ref{bi.compl.def}.
\end{defn}

\begin{remark}\label{I.seg.eq.rem}
As in Example~\ref{seg.eq.exa}, let $T$ be set of maps in
$I^o\Delta^o\Sets$ of the form $\id \times v^l_n:\Delta_i \times
\V^l_n \to \Delta_i \times \Delta_n$, for any $i \in I$ and $n > l >
0$. Then a fibrant $I$-Segal space is the same thing as a
$T$-space. As in Example~\ref{css.eq.exa}, to describe complete
$I$-Segal spaces, choose a map $v:\ppt \to \V_3$, and enlarge $T$ by
adding all the maps $\id \times v:\Delta_{i,n} \to \Delta_{i,n}
\times \V_3$.
\end{remark}

\begin{exa}\label{I.biseg.triv.exa}
Let $\C = (I \times \Delta)X \to I \times \Delta$ be a discrete
Segal fibration corresponding to an $I$-simplicial set $X: I^o
\times \Delta^o \to \Sets$. Then if we treat $X$ as an
$I$-bisimplicial set $X'$ constant in the simplicial direction, it
is a fibrant $I$-Segal space in the sense of
Definition~\ref{I.bi.seg.def}.
\end{exa}

\begin{lemma}\label{I.seg.le}
Assume given an additive semiexact family of groupoids $\C$ over
$I^o_f\Delta^o_f\Delta^o_f\Sets_\kappa$, for some uncountable
regular cardinal $\kappa > |I|$, and assume that $\C$ is constant
along $\kappa$-anodyne maps of Definition~\ref{f.CW.def}. Moreover,
assume given a map $f:X \to Y$ in
$I^o_f\Delta^o_f\Delta^o_f\Sets_\kappa$ such that for any morphism
$i \to i'$ in $I$, $\C(i)$ is stably constant along $f(i'):X(i') \to
Y(i')$ in the sense of Definition~\ref{st.I.def}. Then $\C$ is
constant along $f$.
\end{lemma}

\proof{} Replacing $Y$ with the cylinder $\Cyl(f)$, we may assume
right away that $f$ is injective. Let $W(\C)$ be the class of
injective maps $g$ in $I^o_f\Delta^o_f\Delta^o_f\Sets_\kappa$ such
that $\C$ stably constant along $g$ (that is, constant along $g
\times \id_{R(Z)}$ for any $Z \in \Delta^o\Sets_\kappa$). Then as in
Proposition~\ref{bi.seg.prop}, $W(\C)$ is closed under coproducts
and pushouts by Lemma~\ref{f.CW.le}, and by \eqref{kan.X}, it
contains the map $\id \times f(i'):\Delta_i \boxtimes X(i') \to
\Delta_i \boxtimes Y(i')$ for any map $i \to i'$ in $I$. But any such
map factors through a map
\begin{equation}\label{g.f.eq}
(\Delta_i \boxtimes X(i')) \copr_{\SS_i \times X(i')} (\SS_i
  \boxtimes Y(i')) \to \Delta_i \boxtimes Y(i'),
\end{equation}
where $\SS_i \subset \Delta_i$ is the sphere embedding
\eqref{sph.I.eq}, and then by induction on $\deg i$, $W(\C)$ also
contains all the maps \eqref{g.f.eq}. It remains to observe that by
induction on skeleta, the map $f$ is a finite composition of
pushouts of coproducts of maps \eqref{g.f.eq}.
\endproof

\begin{corr}\label{I.seg.corr}
Assume given an $I$-bisimplicial set $Y \in
I^o_f\Delta^o_f\Delta^o_f\Sets_\kappa$, with $\kappa > |I|$, and
consider the functor
\begin{equation}\label{F.Y}
F_Y:\Delta^o_f\Delta^o_f\Sets_\kappa \to
I^o_f\Delta^o_f\Delta^o_f\Sets_\kappa, \quad Z \mapsto Y \times Z.
\end{equation}
Then for any additive semiexact family of groupoids $\C$ over
$I^o_f\Delta^o_f\Delta^o_f\Sets_\kappa$ that is constant along
$\kappa$-anodyne maps and Segal in the sense of
Definition~\ref{I.bi.seg.def}, $F_Y^*\C$ is an additive semiexact
Segal family that is constant along $\kappa$-anodyne maps, and it is
complete if so is $\C$.
\end{corr}

\proof{} The functor \eqref{F.Y} commutes with products and standard
pushout squares, it sends $\kappa$-anodyne maps to $\kappa$-anodyne
maps, and for any simplicial set $Z \in \Delta^o_f\Sets_\kappa$ and
integers $n > l > 0$, with the corresponding map $b^l_n$ of
\eqref{b.ln.eq}, $\C$ is constant along $F_Y(L(b^l_n) \times
\id_{R(Z)})$ by Lemma~\ref{I.seg.le}. If $\C$ is complete, then for
any simplicial set $X \in \Delta^o_f\Sets_\kappa$, with the
projection $p_X:L(\V_3 \times X) \to L(X)$ of
Definition~\ref{bi.compl.def}, Lemma~\ref{I.seg.le} also shows that
$\C$ is constant along $F_Y(p_X \times \id_{R(Z)}$.
\endproof

In the situation of Example~\ref{I.biseg.triv.exa}, we also have the
twisted $2$-simplicial expansion $\Delta_\sigma(\C\|\Delta)$ of
\eqref{tw.rel.eq}, and it corresponds to a (non-augmented) Segal
space by \eqref{str.del.eq}. To generalize this to arbitrary
$I$-Segal spaces, one can use the alternative description of twisted
$2$-simplicial expansions given in Subsection~\ref{tw.subs}. Namely,
consider the total simplicial expansion $\Delta^\Tot I$, with the
fully faithful embedding $\beta:I \times \Delta \to \Delta^\Tot I$
of \eqref{be.eq}. Denote by $\lambda:\Delta I \to \Delta^\Tot I$
the natural dense embedding, and let $\bpi:\Delta I \to \Delta$ be
the structural cofibration. Then for any $I$-bisimplicial set $X \in
I^o\Delta^o\Delta^o\Sets$, we can define its {\em twisted simplicial
  expansion} by
\begin{equation}\label{tw.seg.eq}
  N_\sigma(X|I) = (\bpi^o \times \id)_!(\lambda^o \times
  \id)^*(\beta^o \times \id)_*X,
\end{equation}
and this actually makes sense for any small $I$. In the situation of
Example~\ref{I.biseg.triv.exa}, Lemma~\ref{beta.le} and
\eqref{str.del.eq} provide an isomorphism
\begin{equation}\label{str.seg.eq}
N_\sigma(X'|I) \cong NI',
\end{equation}
where we interpret $X:I^o \times \Delta^o \to \Sets$ as a functor
$I^o \to \Cat \subset \Delta^o\Sets$, and $I' \to I$ is the
fibration corresponding to this functor by the Grothendieck
construction.

\begin{defn}\label{t.fib.def}
An $I$-bisimiplicial set $X:I^o \times \Delta^o \to \Delta^o\Sets$
is {\em $t$-fibrant} if for any $i \in I$, $X(i \times [0])$ is a
fibrant simplicial set, and for any integer $n \geq 1$, the map
$X(\id \times t):X(i \times [n]) \to X(i \times [n-1])$ is a Kan
fibration.
\end{defn}

\begin{lemma}\label{tw.seg.le}
For any $t$-fibrant $I$-Segal space $X$, the twisted simlicial
expansion $Y = N_\sigma(X|I)$ is a $t$-fibrant Segal space, fibrant
and complete if so is $X$, and for any weak equivalence $w:X \to X'$
of $t$-fibrant $I$-Segal spaces, the induced map $Y \to
Y'=N_\tau(X'|I)$ is a weak equivalence.
\end{lemma}

\proof{} If we treat the bisimplicial set $Y$ as a functor $\Delta^o
\to \Delta^o\Sets$, then \eqref{hom.X} provides a functor
$\Hom(-,Y):\Delta^o\Sets \to \Delta^o\Sets$, and as in
Lemma~\ref{bisimp.le}, to show that $Y$ is a Segal space, we need to
check that $\Hom(-,Y)$ sends the injective maps \eqref{b.ln.eq} to
maps $F \cap W$, and to check that $Y$ is fibrant, we need to check
that is sends embeddings \eqref{sph.eq} to maps in $F$. Moreover, by
the same argument as in Lemma~\ref{bva.le}, it suffices to check the
former for $b^l_n$ with $l=1$. Altogether, for any $n \geq 0$,
$b^1_n$ factors through \eqref{sph.eq}, we have maps
\begin{equation}\label{I.seg.dia}
\Hom(\Delta_n,Y) \longrightarrow \Hom(\SS_{n-1},Y) \longrightarrow
\Hom(\Delta_1 \copr_{\Delta_0} \Delta_{n-1},Y),
\end{equation}
and we need to show that the composition is in $F \cap W$, and the
first map is in $F$ if $X$ is fibrant. Since $\pi^o:(\Delta I)^o \to
\Delta^o$ is a discrete cofibration, $Y([n]) = \Hom(\Delta_n,Y)$
splits into a disjoint union
\begin{equation}\label{Y.spl}
  Y([n]) = \coprod_{i_\idot}Y([n])(i_\idot)
\end{equation}
indexed by all the functors $i_\idot:[n] \to I$ given by diagrams
\eqref{i.dot}, and \eqref{Y.spl} induces the corresponding
splittings for all the other simplicial sets in
\eqref{I.seg.dia}. Then spelling out the definitions, we observe
that as in Lemma~\ref{beta.le}, we have functorial maps
$a(i_\idot):Y([n])(i_\idot) \to X(i_0 \times [n])$ that are
isomorphisms for $n=0$, and for any $n \geq 1$, we have a cartesian
square
\begin{equation}\label{Y.seg.sq}
\begin{CD}
  Y([n])(i_\idot) @>{a(i_\idot)}>> X(i_0 \times [n])\\
  @VVV @VV{X(i_0)(t^o)}V\\
  Y([n-1])(t^*i_\idot) @>{X((f_0 \times \id)^o) \circ
    a(t^*i_\idot)}>> X(i_0 \times [n-1]),
\end{CD}
\end{equation}
where $t:[n-1] \to [n]$ is the standard embedding, and $f_0:i_0 \to
i_1$ is the map in \eqref{i.dot}. Since $X:I^o \times \Delta^o \to
\Delta^o\Sets$ is assumed to be $t$-fibrant, both simplicial sets on
the right are fibrant, and the vertical arrow on the right is in
$F$. Moreover, by induction on $n$, $Y([n-1])(t^*i_\idot)$ is also
fibrant, and then so is $Y([n])(i_\idot)$ and the vertical arrow on
the left. Therefore $Y$ is $t$-fibrant. Moreover, by the same
argument as in Lemma~\ref{semi.mod.le} and induction on $n$, the map
$Y([n])(i_\idot) \to Y'([n])(i_\idot)$ induced by a weak equivalence
$w:X \to X'$ is itself a weak equivalence. Finally, the embedding
$t:\Delta_{n-1} \to \Delta_n$ factors through the map $b^1_n$, so
that the maps \eqref{I.seg.dia} then fit into a commutative diagram
\begin{equation}\label{Y.seg.1.sq}
\begin{CD}
  Y([n])(i_\idot) @>>> X(i_0 \times [n])\\
  @VVV @VVV\\
  \Hom(\SS_{n-1},Y)(i_\idot) @>>> \Hom(\SS_{n-1},X(i_0))\\
  @VVV @VVV\\
  \Hom(\Delta_1 \copr_{\Delta_0} \Delta_{n-1},Y)(i_\idot)
  @>>> \Hom(\Delta_1 \copr_{\Delta_0} \Delta_{n-1},X(i_0))\\
  @VVV @VVV\\
  Y([n-1])(t^*i_\idot) @>>> X(i_0 \times [n-1]),
\end{CD}
\end{equation}
with cartesian squares. Since $X(i_0)$ is assumed to be a Segal
space, the map $X(i_0 \times [n]) \to \Hom(\Delta_1 \copr_{\Delta_0}
\Delta_{n-1},X(i_0))$ is in $F \cap W$, and the same then holds for
the corresponding map on the left. This proves that $Y$ is a Segal
space. Moreover, if $X$ is fibrant, then the same argument as in in
Lemma~\ref{J.proj.le} shows that $X((f_0 \times \id)^o)$ in
\eqref{Y.seg.sq} is in $F$, and then by induction on $n$, so is
$a(i_\idot)$. In addition to this, the vertical arrows on the right
in \eqref{Y.seg.1.sq} are in $F$, and then so are the arrows on the
left. This proves that $Y$ is fibrant. Finally, for any bisimplicial
set $X$, let $X^v = \Hhom(L(\V_3),X)$, as in
Example~\ref{css.sp.exa}. Then this is functorial in $X$, thus can
be applied pointwise to $I$-bisimplicial sets. Moreover, for any
fibrant $I$-Segal space $X$, $X^v$ is a fibrant $I$-Segal space, and
by Example~\ref{css.sp.exa}, $X$ is complete iff the map $X \to X^v$
induced by the projection $\V_3 \to \ppt$ is a weak equivalence. But
since $I$ is rigid, Example~\ref{css.ni.exa} immediately implies
that for any $I$-bisimplicial set $X$, we have $N_\sigma(X|I)^v \cong
N_\sigma(X^v|I)$. Thus if a fibrant $I$-Segal space is complete, the
weak equivalence $X \to X^v$ induces a weak equivalence $N_\sigma(X|I)
\to N_\sigma(X^v|I) \cong N_\tau(X|I)^v$, and $N_\sigma(X|I)$ is
complete.
\endproof

\subsection{Nerves.}

Our version of Proposition~\ref{seg.prop} for $I$-Segal spaces uses
the categories $\Relf \bb^\bB I$, $\Relf^+ \bb^\bB I$ of
restricted $I$-augmented biordered sets of
Subsection~\ref{aug.rel.subs}, where the good filtration of
Definition~\ref{good.def} on $I$ is given by the Reedy category
structure. We note right away that since $I=I_L$ is assumed to be
directed, any $J \in \Relf^+ \bb I$ is automatically locally
restricted with respect to this filtration, so that $\Relf^+
\bb^\bB I = \Relf^+ \bb I$. As in Subsection~\ref{ner.div.subs},
we start by constructing nerve functors relating augmented biordered
sets and $I$-bisimplicial set. First, consider the category
$I^o\Delta^o\Pos$, and as in \eqref{bi.nerve}, define functors
\begin{equation}\label{bi.I.nerve}
N^L,N^R:I^o\Delta^o\Pos \to I^o\Delta^o\Delta^o\Sets
\end{equation}
by taking the nerve pointwise along the Reedy resp.\ simplicial
factor. Say that $J \in I^o\Delta^o\Pos$ is {\em globally
  finite-dimensional} if its underlying $I$-simp\-li\-cial set is
finite-dimensional, and $\dim J(i \times [m]) \leq d$ for some $d$
independent of $i \times [m] \in I \times \Delta$.

\begin{lemma}\label{binerve.I.le}
For any $J \in \Delta^o_f\Pos$ and $d \geq 0$, both
$\sk^{I,L}_dN^L(J)$ and $\sk^{I,R}_dN^R(J)$ are in
$I^o_f\Delta^o_f\Delta^o_f\Sets$, and if $J$ is globally finite
dimensional, then $N^L(J) = \sk^{I,L}_dN^L(J)$, $N^R(J) =
\sk^{I,R}_dN^R(J)$ for some $d \geq 0$.
\end{lemma}

\proof{} Since the Reedy category $I$ is assumed to be directed, all
its objects $i \in I$ are non-degenerate. Therefore for any $X \in
I^o_f\Delta^o_f\Sets$, we have $\dim X = \max_i\{\deg i + \dim
X(i)\}$, where $\max$ is taken over all $i \in I$ with non-empty
$X(i)$. Thus for any $d \geq 1$, we have $\dim X^d \leq d\dim
X$, and we can then use exactly the same argument as in
Lemma~\ref{binerve.le}.
\endproof

Now as in \eqref{bi.N.Q}, we can define nerve functors
\begin{equation}\label{bi.I.N.Q}
\begin{aligned}
\Nn_\dm(I):&\Rel\bb I \to I^o\Delta^o\Delta^o\Sets, \\
\bN_\dm(I):&\bRel\bb I \to I^o\Delta^o\Delta^o\Sets
\end{aligned}
\end{equation}
by setting
\begin{equation}\label{bi.I.N}
\begin{aligned}
\Nn_\dm(I)(J)(i \times [n] \times [m]) &= \Rel(L([n]) \times
R([m]),i \setminus J),\\
\bN_\dm(I)(J)(i \times [n] \times [m]) &= \bRel(L([n]) \times
R([m]),i \setminus J),
\end{aligned}
\end{equation}
where $i \setminus J$ is biordered via the discrete cofibration $i
\setminus J \to J$, and as in \eqref{bN.N}, the embedding
$\rho:\bRel \to \Rel$ induces a map
\begin{equation}\label{bN.N.I}
\bN_\dm(I) \to \Nn_\dm(I) \circ \rho.
\end{equation}
Moreover, as in \eqref{bi.N.dot}, we can extend the functors
\eqref{bi.I.N.Q} to functors
\begin{equation}\label{bi.N.I.dot}
\begin{aligned}
&\bN_\idot(I):\N \times \Rel \bb I \to I^o\Delta^o\Delta^o\Sets,\\
&\Nn_\idot(I):\N \times \Rel \bb I \to I^o\Delta^o\Delta^o\Sets,
\end{aligned}
\end{equation}
by taking the skeleton $\sk_d^{I,L}$, the skeleton inclusion define
tautological maps
\begin{equation}\label{N.d.N.I}
  \bN_\idot(I) \to \bN_\dm(I) \circ \pi, \qquad \Nn_\idot(I) \to
  \Nn_\dm(I) \circ \pi,
\end{equation}
and by Lemma~\ref{N.d.le} combined with Lemma~\ref{binerve.I.le},
both functors \eqref{bi.N.I.dot} send $\N \times \Relf \bb^\bB I$
into $I^o_f\Delta^o_f\Delta^o_f\Sets$, while the map \eqref{N.d.N.I}
for $\bN_\idot(d \times J)$ is an isomorphism for any fixed $J \in
\Relf \bb^\bB I$ and sufficiently large $d$.

Going in the other direction, we can upgrade \eqref{bi.Q.N} to
define functors
\begin{equation}\label{bi.Q.N.I}
\begin{aligned}
\bQQ_\dm(I) &= \bT_{I,F} \circ \QQ_{I \times
  \Delta}:I^o\Delta^o\Delta^o\Sets \to \bRel \bb I,\\
\QQ_\dm(I) &= \T_{I,F} \circ \QQ_{I \times
  \Delta}:I^o\Delta^o\Delta^o\Sets \to \Rel \bb I,
\end{aligned}
\end{equation}
where $\bT_{I,F}$ and $\T_{I,F}$ are the functors
\eqref{T.I.F}. Then functors \eqref{bi.Q.N.I} send
$I^o_f\Delta^o_f\Delta^o_f\Sets$ into $\bRelf \bb^\bB I$
resp.\ $\Relf \bb^\bB I$, and we can modify them to send
$I^o\Delta^o\Delta^o\Sets$ into $\Relf^+ \bb^\bB I$ by setting
\begin{equation}\label{bi.Q.N.I.pl}
  \bQQ_{\dm+}(I) = \bT_{I,F} \circ \QQ_{(I \times \Delta)+}, \quad
  \QQ_{\dm+}(I) = \T_{I,F} \circ \QQ_{(I \times \Delta)+}.
\end{equation}
Then \eqref{Q.I.pl.Q} induces functorial maps
\begin{equation}\label{bi.Q.N.I.pl.Q}
  \QQ_{\dm+}(I)(X) \to \QQ_\dm(I)(X), \qquad X \in I^o\Delta^o\Delta^oX
\end{equation}
that are $+$-bianodyne as soon as $X \in
I^o\Delta^o_f\Delta^o_f\Sets$ by Lemma~\ref{T.I.le}. The same
construction as in Subsection~\ref{ner.div.subs} then provides a
functorial map
\begin{equation}\label{bi.QN.I.id}
a_\dm(I):\QQ_\dm(I) \circ \Nn_\dm(I) \to \id,
\end{equation}
an $I$-augmented version of \eqref{bi.QN.id}, and the maps
\eqref{NQI.id} and \eqref{bN.N.I} induce a diagram
\begin{equation}\label{bi.NQ.I.id}
\begin{CD}
\Nn_\dm(I) \circ \QQ_\dm(I) @<<< \bN_\dm(I) \circ \bQQ_\dm(I) \cong
\Nn_{I \times \Delta} \circ \QQ_{I \times \Delta} @>>> \id.
\end{CD}
\end{equation}
For any uncountable regular cardinal $\kappa$, the functors
\eqref{bi.I.N.Q}, \eqref{bi.N.I.dot}, \eqref{bi.Q.N.I},
\eqref{bi.Q.N.I.pl} sends $\kappa$-bounded $I$-bisimplicial sets to
$\kappa$-bounded augmented biordered sets and vice versa, and the
map \eqref{bi.Q.N.I.pl.Q} is $\kappa^+$-bianodyne for $X \in
I^o_f\Delta^o_f\Delta^o_f\Sets_\kappa$.

\begin{lemma}\label{A.QN.I.le}
Let $\kappa$ be an uncountable regular cardinal. Then for any $J \in
\Relf_\kappa \bb^\bB I$, the composition $\bQQ_\dm(I)(\bN_\dm(I)(J))
\to J$ of the morphisms \eqref{bi.QN.I.id} and \eqref{bN.N.I} is
$\kappa$-bi\-anodyne, and for any $J \in \Relf^+_\kappa \bb^\bB I$,
the composition $\bQQ_{\dm+}(I)(\bN_\dm(I)(J)) \to J$ of the
morphisms \eqref{bi.QN.I.id}, \eqref{bN.N.I} and
\eqref{bi.Q.N.I.pl.Q} is $\kappa^+$-bianodyne.
\end{lemma}

\proof{} As in Lemma~\ref{A.QN.le}, consider the endofunctor $A =
\bQQ_\dm(I) \circ \bN_\dm(I)$ of the category $\Relf_\kappa \bb^\bB
I$, and let $\alpha:A \to \Id$ be the composition of the maps
\eqref{bi.QN.I.id} and \eqref{bN.N.I}. Moreover, let $\chi = a \circ
\rho$, where $p:\QQ_\dm \to Q_{I \times \Delta}$ is the
projection, and $a$ is the map \eqref{QNI.id} for $I \times
\Delta$. Then we have $\chi \geq U(\alpha)$, and $A$ preserves
coproducts, so for the first claim, it suffices to check that the
triple $\langle A,\alpha,\chi \rangle$ satisfies the conditions
\thetag{i}-\thetag{iii} of Lemma~\ref{A.le}. For \thetag{i} and
\thetag{ii}, the argument is the same as in Lemma~\ref{A.QN.le}. For
\thetag{iii}, note that as in \eqref{Q.L.sq}, for any $J \in \Posf
\bb^\flat I$, we have a cartesian square
\begin{equation}\label{Q.L.I.sq}
\begin{CD}
\QQ_\dm(I)(\Nn_\dm(I)(L(J))) @>{q^\hdot}>> B^\dm_{I\idot}(J)\\
@VVV @VVV\\
B^o_I(I \setminus^o_{\xi_\perp} B(J))^o @>{q}>> B_I(J),
\end{CD}
\end{equation}
where $B^\dm_{I\idot}(J)$ is as in Example~\ref{bary.bis.I.rel.exa}, and
the bottom line is \eqref{Q.I.N.I.B}. Then $q$ in \eqref{Q.L.I.sq}
is $\kappa$-anodyne by the same argument as in Lemma~\ref{QNI.le},
so $q^\hdot$ is $\kappa$-bianodyne by Lemma~\ref{T.I.le}, and then
the rest of the proof is exactly as in Lemma~\ref{A.QN.le}. The same
goes for the second claim.
\endproof

Now assume given an object $i \in I$ in the category $I$, and
consider the embedding $\eps(i)_\dg:\Relf \to \Relf \bb^\bB I$
induced by \eqref{cat.phi}. Then \eqref{bi.I.N} together with
\eqref{kan.X} immediately provides isomorphisms
\begin{equation}\label{N.eps.i}
  \eps(i)^o_! \circ \bN_\dm \cong \bN_\dm(I) \circ \eps(i)_\dg,\quad
  \eps(i)^o_! \circ \Nn_\dm \cong \Nn_\dm(I) \circ \eps(i)_\dg,
\end{equation}
and on the other hand, for any bisimplicial set $X$, the embedding
$\eps(i)$ extends to an embedding $\eps(i):\Delta\Delta X \to
I\Delta\Delta\eps(i)^o_!X$, and this is a Reedy functor, so we
obtain a map
\begin{equation}\label{Q.eps.i}
  Q(\eps(i)):\eps(i)_\dg(\QQ_\dm(X)) \to \QQ_\dm(I)(\eps(i)_!^oX)
\end{equation}
functorial with respect to $X$.

\begin{lemma}\label{Q.eps.i.le}
For any uncountable regular cardinal $\kappa$ and $J \in
\Relf_\kappa$, the map \eqref{Q.eps.i} for $X=\bN_\dm(J)$ is
$\kappa$-bianodyne.
\end{lemma}

\proof{} By \eqref{N.eps.i}, the target of the map \eqref{Q.eps.i}
is naturally identified with
$\QQ_\dm(I)(\bN_\dm(I)(\eps(i)_\dg(J)))$, so we are done by
Corollary~\ref{biano.I.corr}, Lemma~\ref{A.QN.le} and
Lemma~\ref{A.QN.I.le}.
\endproof

\subsection{Augmented Segal families.}

For the moment, drop the assumption that $I$ is a directed Reedy
category, and let it be any category equipped with a good
filtration, as in Subsection~\ref{aug.rel.subs}. Note that just as
in Subsection~\ref{aug.subs}, the maps \eqref{J.S.bi},
\eqref{J.pl.bi.eq} and \eqref{J.fl.hs.eq} are naturally
$I$-augmented as soon as so is $J$, and the same goes for standard
pushout squares $J \to R(\V)$, so that Definition~\ref{rel.fam.def}
and Definition~\ref{cont.def} make sense for $I$-augmented biordered
sets. We can then repeat {\em verbatim}
Definition~\ref{rel.seg.def}.

\begin{defn}\label{rel.seg.I.def}
For any uncountable regular cardinal $\kappa > |I|$, a {\em
  $\kappa$-bounded Segal family} of groupoids over $\Relf_\kappa
\bb^\bB I$ is a family of groupoids $\C \to \Relf_\kappa \bb^\bB
I$ that is additive, semiexact and satisfies excision in the sense
of Definition~\ref{rel.fam.def}, satisfies the cylinder axiom of
Definition~\ref{cont.def}, and such that $\|\C_J\| < \kappa$ for any
finite $J \in \Relf$. A {\em Segal family} of groupoids over $\Relf
\bb^\bB I$ is a family of groupoids $\C \to \Relf \bb^\bB I$
that is additive, semiexact, satisfies excision in the sense of
Definition~\ref{rel.fam.def}, and satisfies the cylinder axiom of
Definition~\ref{cont.def}. A {\em semicanonical extension} of a
Segal family resp.\ $\kappa$-bounded Segal family $\C$ is a family
of groupoids $\C^+$ over $\Relf^+ \bb^\bB I$ resp.\ $\Relf^+_\kappa \bb^\bB
I$ equipped with an equivalence $i^*\C^+ \cong \C$ that is
$[1]$-invariant, additive, semiexact, and semicontinuous in the
sense of Definition~\ref{rel.fam.def}.
\end{defn}

\begin{lemma}\label{bi.del.I.le}
Any Segal family of groupoids $\C$ over $\Relf \bb^\bB I$ is
constant along bianodyne maps, and for any $\kappa > |I|$, a
$\kappa$-bounded Segal family $\C$ is constant along
$\kappa$-bianodyne maps. A semicanonical extension $\C^+$ of such a
family $\C$ is constant along $+$-bianodyne
resp.\ $\kappa^+$-bianodyne maps.
\end{lemma}

\proof{} As in Lemma~\ref{bi.del.le}, everything immediately follows
from Proposition~\ref{biano.I.prop}, by the same argument as in
Lemma~\ref{ano.grp.le}.
\endproof

\begin{corr}\label{cyl.2.I.corr}
For any $I$-augmented $\langle J,\alpha \rangle \in
\Relf \bb^\bB I$ equipped with a bicofibration $J \to L([2])$, any
Segal family $\C$ over $\Relf \bb^\bB I$ is semicartesian along
the square \eqref{s.1.t.sq}.
\end{corr}

\proof{} Same as Corollary~\ref{cyl.2.corr}.
\endproof

Next, we need an augmented version of the material of
Section~\ref{css.sec}. Again, we start with repeating
Definition~\ref{rel.refle.def} and Definition~\ref{pos.refle.def}.

\begin{defn}\label{rel.refle.I.def}
For any ample subcategory $\I \subset \Pos$, with unfolding $\I^\dm
\subset \Rel$, a family of groupoids $\C$ over $\I^\dm \bb^\bB I$
resp.\ $\I \bb^\bB I$ is {\em reflexive} if it is fully faithful
along the projection $L([1]) \times J \to J$ resp.\ $[1] \times J
\to J$ for any $J \in \I^\dm \bb^\flat I$ resp.\ $J \in \I \bb^\flat
I$.
\end{defn}

However, in the augmented case this is not enough: we also need to
consider more general augmentations on the product $J \times
[1]$. Explicitly, for any $J \in \Pos$, an $I$-augmentation
$\alpha:J \times [1] \to I$ is given by two augmentations
$\alpha_0,\alpha_1 \to I$ and a morphism $f:\alpha_0 \to
\alpha_1$. If we denote $J_l = \langle J,\alpha_l \rangle$, $l=0,1$,
and let $g=\langle \id,f \rangle:J_0\to J_1$ be the corresponding
augmented map, then $\langle J \times [1],\alpha \rangle \cong
\Cyl(g)$, and we have the embedding $t:J_1 \to \Cyl(g)$. We also
have an augmented map $\wt{g}:J_0 \times [1] \to \Cyl(g)$ given by
$\id$ resp.\ $g$ over $0$ resp.\ $1$, and a commutative square
\begin{equation}\label{a.cyl.sq}
\begin{CD}
\Cyl(g) @<{t}<< J_1\\
@A{\wt{g}}AA @AA{g}A\\
J_0 \times [1] @<{\id \times t}<< J_0.
\end{CD}
\end{equation}
Moreover, we have the projection $e_0:J_0 \times [1] \to J_0$ and the
adjoint map $e=t_\dg:\Cyl(g) \to J_1$. Then for any ample subcategory
$\I \subset \Pos$ such that $\Cyl(g) \in \I \bb^\bB I$, and for
any family of groupoids $\C$ over $\I \bb^\bB I$, we can denote by
$\C^e_{\Cyl(g)} \subset \C_{\Cyl(g)}$ resp.\ $\C^e_{J_0 \times [1]}
\subset \C_{J_0 \times [1]}$ the essential image of the functor
$e^*$ resp.\ $e_0^*$, and \eqref{a.cyl.sq} induces a functor
\begin{equation}\label{a.cyl.eq}
  \C^e_{\Cyl(g)} \to \C^e_{J_0 \times [1]} \times_{\C_{J_0}} \C_{J_1}.
\end{equation}
Moreover, consider the projection $s^* \times t^*:\C_{J_0 \times [1]} \to
\C_{J_0}^2 = \C_{J_0} \times \C_{J_0}$ and its decomposition
$$
\begin{CD}
\C_{J_0 \times [1]} @>{a}>> \pi_0(\C_{J_0 \times [1]} | \C_{J_0}^2)
@>>> \C_{J_0}^2
\end{CD}
$$
of \eqref{grp.eq}; then \eqref{a.cyl.eq} also induces a commutative
square
\begin{equation}\label{a.red.sq}
\begin{CD}
\C_{\Cyl(g)} @>>> \C_{J_1}\\
@V{a \circ \wt{g}^*}VV @VV{g^*}V\\
\pi_0(\C_{J_0 \times [1]} | \C_{J_0}^2) @>{t^*}>> \C_{J_0}.
\end{CD}
\end{equation}
If $J$ is equipped with a biorder, then \eqref{a.cyl.eq} with $[1]$
replaced by $L([1])$ becomes a commutative square of augmented
biordered sets, and induces a functor \eqref{a.cyl.eq} and a
commutative square \eqref{a.red.sq} for any family of groupoids $\C$
over $\I^\dm \bb^\bB I$.

\begin{defn}\label{cohe.def}
A family of groupoids over $\I \bb^\bB I$ resp.\ $\I^\dm \bb^\bB
I$ is {\em coherent} if for any $\langle J \times [1],\alpha \rangle
= \Cyl(g) \in \I \bb^\bB I$ resp.\ $\langle J \times L([1]),\alpha
\rangle \in \I^\dm \bb^\bB I$, the corresponding functor
\eqref{a.cyl.eq} is an equivalence, and the square \eqref{a.red.sq}
is semicartesian.
\end{defn}

Comparing Definition~\ref{rel.refle.I.def} and
Definition~\ref{cohe.def}, we see that if a family of groupoids $\C$
over $\I \bb^\bB I$ is both coherent and reflexive, then it is
also fully faithful along the projection $e:\langle J \times
[1],\alpha \rangle \to \langle J,\alpha_1 \rangle$ for any $\langle
J \times [1],\alpha \rangle$ in $\I \bb^\bB I$, and similarly in
the biordered case. In particular, take $J = \ppt$; then for
any $\alpha:[1] \to I$, a coherent reflexive family of groupoids
$\C$ over $\I^\dm \bb^\bB I$ is fully faithful along the
projection $e:\langle [1],\alpha \rangle \to \langle
\ppt,\alpha_1\rangle$. Now for any $J \in \I^\dm \bb^\bB I$ and a
collection $W$ of maps $w:L([1]) \to J$, we have the full
subcategory $\C_J(W) \subset \C_J$ of $c \in \C_J$
such that $w^*c \in e^*\C_\ppt \subset \C_{L[1]}$ for all $w \in
W$. With this convention, and with Lemma~\ref{bi.del.le} and
Corollary~\ref{cyl.2.corr} replaced by Lemma~\ref{bi.del.I.le} and
Corollary~\ref{cyl.2.I.corr}, all the material in
Subsection~\ref{refl.grp.subs} carries over to the $I$-augmented
setting; in particular, if we use Definition~\ref{rel.refle.I.def}
instead of Definition~\ref{rel.refle.def}, then
Lemmas~\ref{E.0.le}--\ref{E.le} hold for families over $\Relf
\bb^\bB I$, $\Relf^+ \bb^\bB I$ and their $\kappa$-bounded versions
with the same proof.

Then for any family of groupoids $\C$ over $\I \bb^\bB I$, we
define the {\em relative unfolding} $(\C|I)^\dm \subset U^*\C$ as
the full subcategory spanned by $\C_J(W(J))$ for all $J \in \I^\dm
\bb^\bB I$, and as in Definition~\ref{pos.nondege.def}, we say
that the family $\C$ is {\em non-degenerate} if $(\C|I)^\dm$ is
constant along the map $J \times L([1]) \to J \times R([1])$ for any
$J \in \I \bb^\bB I$. Moreover, we observe that an augmentation on
some $J \in \Pos$ induces augmentations on all the squares
\eqref{J.S.sq}, \eqref{J.cyl.sq}, \eqref{J.pl.sq}, so that
Definition~\ref{refl.exci.def} also makes sense in the augmented
setting.

\begin{defn}\label{red.seg.I.def}
An {\em $I$-augmented restricted Segal family} is a family of
groupoids $\C$ over $\Posf \bb^\bB I$ that is additive and semiexact
in the sense of Definition~\ref{pos.def}~\thetag{ii},\thetag{iii},
non-degenerate, reflexive in the sense of
Definition~\ref{rel.refle.I.def}, coherent in the sense of
Definition~\ref{cohe.def}, and satisfies weak excision and the
cylinder axiom in the sense of Definition~\ref{refl.exci.def}. A
{\em semicanonical extension} of a restricted Segal family $\C$ is a
family of groupoids $\C^+$ over $\Posf^+ \bb^\bB I$ that is additive
and semiexact in the sense of
Definition~\ref{pos.def}~\thetag{ii},\thetag{iii}, reflexive in the
sense of Definition~\ref{rel.refle.I.def}, non-degene\-rate, coherent
in the sense of Definition~\ref{cohe.def}, weakly semicontinuous in
the sense of Definition~\ref{refl.exci.def}, and equipped with an
equivalence $i^*\C^+ \cong \C$.
\end{defn}

\begin{prop}\label{restr.seg.I.prop}
The unfolding $(\C|I)^\dm$ of a $I$-augmented restricted Segal
family $\C \to \Posf \bb^\bB I$ in the sense of
Definition~\ref{red.seg.I.def} is a Segal family in the sense of
Definition~\ref{rel.seg.I.def}, reflexive in the sense of
Definition~\ref{rel.refle.I.def} and coherent in the sense of
Definition~\ref{cohe.def}. Conversely, the restriction $L^*\C$ of a
coherent reflexive Segal family $\C \to \Relf \bb^\flat I$ in the
sense of Definition~\ref{rel.seg.I.def},
Definition~\ref{rel.refle.I.def} and Definition~\ref{cohe.def} is a
restricted Segal family in the sense of
Definition~\ref{red.seg.I.def}, and we have a natural equivalence
$\C \cong (L^*\C|I)^\dm$. Moreover, for any semicanonical extension
$\C^+$ of the family $\C$, the unfolding $(\C^+|I)^\dm$ is a
semicanonical extension of the unfolding $(\C|I)^\dm$, and conversely,
for any semicanonical extension $(\C|I)^{\dm+}$ of the unfolding
$(\C|I)^\dm$, $\C^+ = L^*(\C|I)^{\dm+}$ is a semicanonical extension of
$\C$, it is reflexive, coherent, safisfies weak excision and the
cylinder axiom of Definition~\ref{refl.exci.def}, and we have
$(\C|I)^{\dm+} \cong (\C^+|I)^\dm$.
\end{prop}

\proof{} Same as Proposition~\ref{restr.seg.prop}.
\endproof

\begin{lemma}\label{iota.I.le}
Let $B_{I\dm}:\Pos \bb I \to \Relf^{\pm} \bb I$ be the $I$-augmented
barycentric subdivision functor of
Example~\ref{bary.rel.I.exa}. Then for any $I$-augment\-ed retricted
Segal family of groupoids $\C \to \Posf \bb^\bB I$, with semicanonical
extension $\C^+$, the family $\C^\iota = B_{I\dm}^*(\C|I)^{\dm
  o}_\perp$ is an $I$-augmented restricted Segal family, and
$\C^{\iota+} = B_{I\dm}^*(\C^+|I)^{\dm o}_\perp$ is its semicanonical
extension.
\end{lemma}

\proof{} Everything except for the fact that $\C^\iota$ is coherent
follows by exactly the same argument as in
Proposition~\ref{iota.prop}. To see that $\C^\iota$ is coherent,
note that as in the proof of Proposition~\ref{iota.prop}, by
Lemma~\ref{pm.B.le} and Lemma~\ref{bary.cyl.le}, if we apply
$B_{I\dm}$ to a square \eqref{a.cyl.sq}, then the resulting square
in $\Relf^{\pm} \bb I$ maps to the square \eqref{a.cyl.eq} for
$B_{I\dm}(g)$ by an $I$-augmented
$\pm$-bianodyne map.
\endproof

\begin{defn}\label{pos.compl.I.def}
A family of groupoids $\C$ over $\I \bb^\bB I$ for an ample
subcategory $\I \subset \Pos$ is {\em complete} if it is cartesian
over \eqref{pos.compl.sq} for any $J \in \I \bb^\bB
I$. A family of groupoids $\C$ over $\I^\dm \bb^\bB I$ is {\em
  complete} if so is $L^*\C$.
\end{defn}

Now recall that we have the embedding $\lambda:\Pos \bbi I \to \Pos \bb
I$ of \eqref{pos.I.eq} and the fibration $\bpi:\Pos \bbi I \to \Pos$ of
\eqref{pos.I.pi}, and if $I$ is rigid, the fibration $\bpi$ is
discrete. Then $\lambda$ and $\bpi$ induce functors between $\Posf^+
\bbi I$, $\Posf^+ \bb I$ and $\Posf^+$, and any family of groupoids
$\C$ over $\Posf^+ \bb I$ induces a family of groupoids
$\bpi_!\lambda^*\C$ over $\Posf^+$.

\begin{lemma}\label{seg.I.le}
A coherent family of groupoids $\C^+$ over $\Posf^+ \bb I$ is a
semicanonical extension of a restricted Segal family if and only if so
is $\bpi_!\lambda^*\C^+$, and then $\C^+$ is complete if and only if
so is $\bpi_!\lambda^*\C^+$.
\end{lemma}

\proof{} Note that $\bpi:\Posf^+ \bbi I \to \Posf^+$ is a semicanonical
extension of a restricted Segal family by
Example~\ref{I.triv.seg.exa}. In particular, it is fully faithful
along the projection $[1] \times J \to J$ for any $J \in
\Posf^+$. Then $\C^+$ is obviously reflexive and/or non-degenerate
and/or additive iff so is $\bpi_!\lambda^*\C^+$, and the rest
immediately follows from Lemma~\ref{car.rel.le}.
\endproof

\begin{lemma}\label{seg.bc.le}
Assume that the category $I$ is small and rigid. Then for any
$I$-bisimplicial set $X$, with twisted simplicial expansion
$Y=N_\sigma(X|I)$ of \eqref{tw.seg.eq}, we have
$\bpi_!\lambda^*L^*\Nn_\dm(I)^*\Hh(X) \cong L^*\Nn_\dm^*\Hh(Y)$,
functorially with respect to $X$.
\end{lemma}

\proof{} By Proposition~\ref{Hh.prop}, it suffices to consider
$I$-bisimplicial sets $X$ that are constant in the simplicial
direction; in other words, for any $I$-simplicial set $X$, we have
to construct a functorial isomorphism
\begin{equation}\label{I.X.Y}
\bpi^o_!\lambda^{o*}\Nn_I^{o*}\Y(I)^o_*(X) \cong
\Nn^{o*}\Y^o_*(\bpi^o_!\lambda^{o*}\beta^o_*X),
\end{equation}
where $\Y(I):I \times \Delta \to I^o\Delta^o\Sets$, $\Y:\Delta \to
\Delta^o\Sets$ are the Yoneda embeddings \eqref{yo.yo.eq}, and
$\beta$ is the embedding \eqref{be.eq}, as in
\eqref{tw.seg.eq}. Denote by $F:\Delta \to \Posf^+$ and
$F(I):\Delta^\Tot I = F^*(\Posf^+ \bb I) \to \Posf^+ \bb I$ the
standard full embeddings. Then we have the full embedding $F(I)
\circ \beta:I \times \Delta \to \Posf^+ \bb I$, and $\Y(I) \cong
\Nn_I \circ F(I) \circ \beta$. Since the nerve functor
$\Nn_I:\Posf^+ \bb I \to I^o\Delta^o\Sets$ is fully faithful, we
have $\Nn_I^{o*} \circ \Y(I)^o_* \cong (F(I) \circ \beta)^o_* = F(I)^o_*
\circ \beta^o_*$, and analogously, $\Nn^{o*} \circ \Y^o_* \cong
F^o_*$. We then have the composition
\begin{equation}\label{I.X.Y.dia}
\begin{CD}
\bpi^o_!\lambda^{o*}F(I)^o_*\beta^o_* @>>>
\bpi^o_!F(I)^o_*\lambda^{o*}\beta^o_* @>>>
F^o_*\bpi^o_!\lambda^{o*}\beta^o_*
\end{CD}
\end{equation}
of the base change isomorphism \eqref{bc.adj.eq} and the base change
map \eqref{bc.2.eq}. But the discrete fibration $\bpi:\Delta I \to
\Delta$ corresponds to $N(I):\Delta^o \to \Sets$, and
since the nerve functor $N:\Cat \to \Delta^o\Sets$ is fully
faithful, the discrete fibration $\Posf^+ \bbi I \to \Posf^+$
corresponds to $F^o_*N(I):\Posf^+ \to \Sets$. Then by
Lemma~\ref{bc.2.le}, the second map in \eqref{I.X.Y.dia} is also an
isomorphism, and the composition provides a functorial isomorphism
\eqref{I.X.Y}.
\endproof

\subsection{Comparison over $I$.}

We can now prove an $I$-augmented version of
Proposition~\ref{seg.prop}, combined with
Proposition~\ref{css.prop}.

\begin{lemma}\label{cohe.le}
For any $\Hom$-finite directed cellular Reedy category $I$ and
$I$-Segal space $X$, the family of groupoids $\Nn_\dm(I)^*\Hh(X)$
is coherent in the sense of Definition~\ref{cohe.def}.
\end{lemma}

\proof{} Observe that the nerve functor $\Nn_\dm(I)$ sends a
commutative square \eqref{a.cyl.sq} into a cocartesian square of the
form \eqref{f.CW.sq} in $I^o\Delta^o\Delta^o\Sets$, and apply
Lemma~\ref{f.CW.le}.
\endproof

\begin{prop}\label{seg.I.prop}
Assume given a $\Hom$-finite directed celullar Reedy category $I$.
\begin{enumerate}
\item For any small additive semiexact family of groupoids $\C$ over
  the category $I^o_f\Delta^o_f\Delta^o_f\Sets$ that is constant
  along anodyne maps of Definition~\ref{f.CW.def} and is a complete
  Segal family in the sense of Definition~\ref{I.bi.seg.def}, the
  family of groupoids $\pi_*\Nn_\idot(I)^*\C$ is a small Segal
  family over $\Relf \bb^\bB I$ in the sense of
  Definition~\ref{rel.seg.I.def}, reflexive in the sense of
  Definition~\ref{rel.refle.I.def} and complete in the sense of
  Definition~\ref{pos.compl.I.def}, while the functors
\begin{equation}\label{fn.I}
\begin{CD}
\C  @>>> \bQQ_\dm(I)^*\pi_*\bN_\idot(I)^* @<<<
\QQ_\dm(I)^*\pi_*\Nn_\idot(I)^*\C
\end{CD}
\end{equation}
induced by the maps \eqref{bi.NQ.I.id} and \eqref{N.d.N.I} are
equivalences. For any semicanonical extension $\C^+$ of the
family $\C$ in the sense of Definition~\ref{X.ext.def},
$\Nn_\dm(I)^*\C^+$ with the functor $i^*\Nn_\dm(I)^*\C^+ \to
\pi_*\Nn_\idot(I)^*\C$ induced by \eqref{N.d.N.I} is a semicanonical
extension of the family $\pi_*\Nn_\idot(I)^*\C$, and the functors
$\C^+,\QQ_{\dm+}(I)^*\Nn_\dm(I)^*\C^+ \to
\bQQ_{\dm+}(I)^*\bN_\dm(I)^*\C^+$ induced by the maps
\eqref{bi.Q.N.I.pl.Q} and \eqref{bi.NQ.I.id} are equivalences.
\item For any small Segal family $\C$ over $\Relf \bb^\bB I$,
  reflexive in the sense of Definition~\ref{rel.refle.I.def},
  coherent in the sense of Definition~\ref{cohe.def} and complete in
  the sense of Definition~\ref{pos.compl.I.def}, the pullback
  $\QQ_\dm(I)^*\C$ with respect to the functor \eqref{bi.Q.N.I} is
  an additive semiexact family over the category
  $I^o_f\Delta^o\Delta^o\Sets$ that is constant along anodyne maps
  of Definition~\ref{f.CW.def}. Moreover, $\QQ_\dm(I)^*\C$ is a
  complete Segal family in the sense of
  Definition~\ref{I.bi.seg.def}, and the functor $\C \cong
  \pi_*\pi^*\C \to \pi_*\Nn_\idot(I)^*\QQ_\dm(I)^*\C$ induced by the
  composition of the maps \eqref{N.d.N.I} and \eqref{bi.QN.I.id} is
  an equivalence. For any semicanonical extension $\C^+$ of the
  Segal family $\C$, $\QQ_{\dm+}(I)^*\C^+$ with the functor
  $i^*\QQ_{\dm+}(I)^*\C^+ \to \QQ_\dm(I)^*\C$ induced by the map
  \eqref{bi.Q.N.I.pl.Q} is a semicanonical extension of the Segal
  family $\QQ_\dm^*\C$, and the functor $\C^+ \to
  \Nn_\dm(I)^*\QQ_{\dm+}(I)^*\C^+$ induced by the maps
  \eqref{bi.Q.N.I.pl.Q}, \eqref{bi.QN.I.id} and \eqref{bN.N.I} is an
  equivalence.
\end{enumerate}
Moreover, for the any uncountable regular cardinal $\kappa > |I|$,
the same statements hold for $\kappa$-bounded Segal families over
$I^o_f\Delta^o_f\Delta^o_f\Sets_\kappa$ and $\Relf_\kappa \bb^\bB
I$.
\end{prop}

\proof{} For \thetag{i}, by Proposition~\ref{X.ext.prop} and
Theorem~\ref{coho.bisi.thm}, we may assume that $\C$ comes equipped
with a semicanonical extension $\C^+ \cong \Hh(X)$ for some complete
$I$-Segal space $X$. In particular, it is constant along weak
equivalences in $I^o\Delta^o\Delta^o\Sets$, so that the functors
\eqref{fn.I} are equivalences by exactly the same argument as in
Proposition~\ref{seg.prop}~\thetag{i}, and moreover,
$\rho^*\pi_*\Nn_\idot(I)^*\C \cong \bN_\dm(I)^*\C$. Since
$\bN_\dm(I)$ commutes with coproducts and sends standard pushout
squares to standard pushout squares, $\bN_\dm(I)^*\C$ is additive
and semiexact. Moreover, for any $J \in \Relf \bb^\flat I$ and $i
\in I$, taking the comma-fiber $i \setminus J$ commutes with
constructing either of the maps \eqref{J.S.bi}, \eqref{J.fl.hs.eq},
and then by \eqref{bi.I.N}, $\bN_\dm(I)^*\C$ satisfies excision and
the cylinder axiom by Proposition~\ref{seg.prop}~\thetag{i} and
Lemma~\ref{I.seg.le}. It is coherent by Lemma~\ref{cohe.le} and
reflexive and complete by Corollary~\ref{I.seg.corr} and
Proposition~\ref{css.prop}. Moreover, forming a map
\eqref{J.pl.bi.eq} also commutes with taking the comma-fiber $i
\setminus J$, so that $\bN_\dm(I)_\dm$ sends such a map to a weak
equivalence, and therefore $\Nn_\dm(I)^*\C^+$ is semicontinuous, and
$\C^+ \cong \QQ_{\dm+}(I)^*\Nn_\dm(I)^*\C^+$.

For \thetag{ii}, by Lemma~\ref{QNI.le}~\thetag{iii},
Lemma~\ref{T.I.le} and Lemma~\ref{bi.del.I.le}, the family
$\QQ_\dm(I)^*\C$ is semiexact, additive and constant along
$\kappa$-anodyne maps, and then it is constant along weak
equivalences by Corollary~\ref{W.CW.corr}. Since for any $i \in I$,
$\eps(i)^o_!$ sends weak equivalences to weak equivalences by
\eqref{kan.X}, while by Lemma~\ref{QNI.le}~\thetag{iv}, every
bisimplicial set $X \in \Delta^o_f\Delta^o_f\Sets$ is weakly
equivalent to the nerve $\bN_\dm(J)$ of a biordered set $J \in
\Relf$, Lemma~\ref{Q.eps.i.le} and Lemma~\ref{bi.del.I.le}
immediately imply that $\C$ is constant along the map
\eqref{Q.eps.i} for any $X \in \Delta^o_f\Delta^o_f\Sets$ and $i \in
I$. Therefore $\QQ_\dm(I)^*\C(i) \cong \QQ_\dm^*\eps(i)_\dg^*\C$ is
a complete Segal family by Proposition~\ref{seg.prop} and
Proposition~\ref{css.prop}, and then as in
Proposition~\ref{seg.prop}~\thetag{ii}, the remaining claims
immediately follow from Lemma~\ref{A.QN.I.le}.
\endproof

\begin{remark}
Proposition~\ref{seg.I.prop}~\thetag{i} together with
Proposition~\ref{restr.seg.I.prop} imply that
$L^*\Nn_\dm(I)^*\C^+$ is a complete restricted Segal family over
$\Posf^+ \bb I$. Alternatively, this can be seen directly: the
twisted simplicial expansion $Y = N_\sigma(X|I)$ is a complete Segal
space by Lemma~\ref{tw.seg.le}, and then it suffices to apply
Lemma~\ref{seg.I.le} and Lemma~\ref{seg.bc.le}. Unfortunately, this
is not enough to conclude anything about $\Nn_\dm(I)^*\C^+$ --- in
order to apply Proposition~\ref{restr.seg.I.prop} in the other
direction, we need to know already that $\Nn_\dm(I)^*\C^+$ is a
Segal family. An obvious counterpat of Lemma~\ref{seg.bc.le} for
biordered sets is wrong, already for the trivial family with fiber
$\ppt$. Indeed, the unfolding $(\Posf^+ \bbi I)^\dm$ of the reflexive
family $\Posf^+ \bbi I \to \Posf^+$ is not the whole $\Relf^+ \bbi I$ ---
it is the full subcategory spanned by $I$-augmented biordered sets
$\langle J,\alpha \rangle$ such that $\alpha(j) \to \alpha(j')$ is
an isomorphism for any $j \leq^l j'$ in $J$. Since it is a Segal
family by Proposition~\ref{restr.seg.prop}, the whole $\Relf^+ \bbi I
\to \Relf^+$ cannot be Segal.
\end{remark}

We also have an $I$-augmented version of
Corollary~\ref{seg.corr}. For any uncountable regular $\kappa >
|I|$, denote by $\CSS_\kappa(I) \subset (\Cat \bbi (\Relf_\kappa
\bb^\bB I))^0$ the full subcategory spanned by reflexive complete
$\kappa$-bounded Segal families.

\begin{corr}\label{seg.I.corr}
In the assumptions of Proposition~\ref{seg.I.prop}, for any
uncountable regular cardinal $\kappa > |I|$, any Segal family
$\C_\kappa \in \CSS_\kappa(I)$ extends to a small Segal family $\C
\to \Relf \bb^\flat I$. Moreover, for any two such
$\C_\kappa,\C'_\kappa \in \CSS_\kappa(I)$ with extensions $\C$,
$\C'$, any functor $\gamma_\kappa:\C_\kappa \to \C'_\kappa$ over
$\Relf_\kappa \bb^\flat I$ extends to a functor $\gamma:\C \to \C'$
over $\Relf \bb^\flat I$, and any morphism $\gamma_\kappa \to
\gamma'_\kappa$ between two such functors $\gamma_\kappa$,
$\gamma_\kappa'$ with extensions $\gamma$, $\gamma'$ extends to a
(non-unique) morphism $\gamma \to \gamma'$. In particular, small
complete reflexive Segal families and isomorphism classes of
functors over $\Relf \bb^\flat I$ between them form a well-defined
category
$$
\CSS(I) = \bigcup_\kappa \CSS_\kappa(I).
$$
We have an equivalence between the category $\CSS(I)$ and the full
subcategory in $h^W(I^o\Delta^o\Delta^o\Sets)$ spanned by complete
Segal spaces, and for any uncountable regular cardinal $\kappa >
|I|$, this equivalence induces an equivalence between
$\CSS_\kappa(I)$ and the full subcategory of complete Segal spaces
in $h^W(I^o\Delta^o\Delta^o\Sets_\kappa)$. Moreover, any small or
$\kappa$-bounded reflexive Segal family $\C$ with complete $L^*\C$
admits a semicanonical extension $\C^+$, for two families $\C_0$,
$\C_1$ with semicanonical extensions $\C_0^+$, $\C_1^+$, any functor
$\gamma:\C_0 \to \C_1$ extends to a functor $\gamma^+:\C_0^+ \to
\C_1^+$, and for two functors $\gamma_0,\gamma_1:\C_0 \to \C_1$ with
extensions $\gamma_0^+$, $\gamma_1^+$, a morphism $\gamma_0 \to
\gamma_1$ extends to a (non-unique) morphism $\gamma_0^+ \to
\gamma_1^+$, so that in particular, the semicanonical extension is
unique up to a unique equivalence.
\end{corr}

\proof{} Clear. \endproof

\begin{remark}
By Remark~\ref{I.seg.eq.rem}, Lemma~\ref{bous.le}~\thetag{ii}
immediately implies that $\CSS(I) \subset h^W(I^o\Delta^o\Sets)$ is
left-admissible for any $\Hom$-finite directed cellular Reedy
category $I$.
\end{remark}

\subsection{Comparison for expansions.}

For any $\Hom$-finite directed cellular Reedy category $I$ and
$I$-bisimplicial set $X$, Proposition~\ref{seg.I.prop} provides an
$I$-augmented biordered set $J \in \Rel \bb I$ equipped with a weak
equivalence $w:\Nn_\dm(I)(J) \to X$, and $w$ induces a map
\begin{equation}\label{tw.w.eq}
  N_\sigma(w):N_\sigma(\Nn_\dm(I)(J)|I) \to N_\sigma(X|I).
\end{equation}
Moreover, \eqref{str.seg.eq} induces an isomorphism
\begin{equation}\label{tw.I.eq}
  N_\sigma(\Nn_\dm(I)(J)|I) \cong \Nn_\dm(I \setminus J),
\end{equation}
where since $I$ is cellular, the comma-category $I \setminus J$ is a
partially ordered set, cofibered over $J$ via $\tau:I \setminus J
\to J$ and we biorder it as in Example~\ref{facto.fib.exa}. It would
be nice to know that under appropriate assumptions --- for example,
if $X$ is a fibrant $I$-Segal space --- \eqref{tw.w.eq} is a weak
equivalence. Unfortunately, this is not at all obvious: while
$\Nn_\dm(J)$ is an $I$-Segal space, it is not fibrant, nor even
$t$-fibrant unless the biorder on $J$ is trivial, so
Lemma~\ref{tw.seg.le} does not apply. If $I$ itself is a partially
ordered set, then more generally, for any functor $J:I^o \to \Rel$
that factors through $\bRel \subset \Rel$, the composition $\Nn_\dm
\circ J:I^o \to \Delta^o\Delta^o\Sets$ is an $I$-bisimplicial
set. If we have another $I$-bisimplicial set $X$, a map $w:\Nn_\dm
\circ J \to X$ induces a map
\begin{equation}\label{N.w}
  N_\sigma(w):\Nn_\dm(J^\hdot) \to N_\sigma(X|I),
\end{equation}
where $J^\hdot \to L(I)$ is the biordered fibration of
Example~\ref{R.groth.exa} corresponding to $J$, and as in
\eqref{tw.I.eq}, we use \eqref{str.seg.eq} to identify the source of
the map. If $J = \sE_I(J')$, $i \mapsto i \setminus J'$ corresponds
to an $I$-augmented biordered set $J' \in \Rel \bb I$, then we have
$J^\hdot \cong I \setminus J'$, and \eqref{N.w} reduces to
\eqref{tw.w.eq} combined with \eqref{tw.I.eq}. We then have the
following general result.

\begin{prop}\label{tw.seg.prop}
Assume given a left-bounded partially ordered set $I$, a functor
$J:I^o \to \Relf^-_C$, where $\Relf^-$ is the category of
right-bounded biordered sets equiped with the co-standard
$C$-category structure, a fibrant $I$-Segal space $X$ and a weak
equivalence $w:\Nn_\dm \circ J \to X$. Then the corresponding map
\eqref{N.w} is a Segal equivalence in the sense of
Example~\ref{seg.eq.exa}.
\end{prop}

Our main technical tool in proving Proposition~\ref{tw.seg.prop} are
the isomorphisms \eqref{N.E} that express the biordered nerve
$\Nn_\dm$ as the nerve of a simplicial partially ordered set, so we
need some preliminaries on $\Delta^o\Pos$. Firstly, recall that for
any partially ordered set $J' \in \Pos$ and simplicial partially
ordered set $J \in \Delta^o\Pos$ equipped with a map $\chi:J \to
J'$, we have the corresponding map $a^\dm_\chi$ of
\eqref{a.chi.bis.eq}. Moreover, if we identify $(J/J')^\iota \cong
({J'}^o \setminus J^\iota)$, then we also have a dual map
\begin{equation}\label{b.chi.bis.eq}
\colim_{\Tw(J')}N^L((J' \setminus J)_f) \to N^L(J'
\setminus J),
\end{equation}
where $(J' \setminus J)_f \cong (j' \setminus J) \times (J' / j)$
for any $f \in \Tw(J')$ given by a pair $j \leq j'$.

\begin{lemma}\label{b.chi.le}
Assume that $J' \in \Posf$, and $J \in \Delta^o_{d+}\Posf^+$ is
globally left-bounded in the sense of
Definition~\ref{glob.def}. Then the corresponding map
\eqref{a.chi.bis.eq} is a Segal equivalence. Moreover, if $J \in
\Delta^o_{d-}\Posf^-$ is globally right-bounded, then so is the dual
map \eqref{b.chi.bis.eq}.
\end{lemma}

\proof{} For the first claim, note that by Corollary~\ref{W.C.corr},
the class of injective Segal equivalences is closed under colimits
over $\N$, so it suffices to prove the claim for the comma-fibers $J
/ n$, $n \geq 0$ of some map $J \to \N$. If we take the height
function $\hht:J \to \N$, then all these comma-fibers are globally
finite-dimensional, and we are done by
Proposition~\ref{bi.seg.prop}~\thetag{iii}. For the second claim,
recall that $N^L(J^\iota) \cong N^L(J)^\iota$ for any $J \in
\Delta^o\Pos$, and by Example~\ref{seg.sp.exa}, $(\iota \times
\id)^*:\Delta^o\Delta^o\Sets \to \Delta^o\Delta^o\Sets$ sends Segal
equivalences to Segal equivalences.
\endproof

Next, we consider the category $\Delta^o\Posf^-$ of simplicial
partially ordered sets that are pointwise right-bounded. The
co-standard $C$-structure on $\Posf^-$ of Example~\ref{pos.CW.exa}
induces a $C$-category structure on $\Delta^o\Posf^-$, with the
class $C$ defined pointwise, and as a $C$-category,
$\Delta^o\Posf^-$ is abundant in the sense of
Definition~\ref{ab.def} and excisive in the sense of
Lemma~\ref{exc.CW.le}. Say that a map $f$ in $\Delta^o\Pos$ is a
{\em Segal equivalence} if so is the induced map $N^L(f)$, and say
that a commutative square
$$
\begin{CD}
  J @>>> J_0\\
  @VVV @VVV\\
  J_1 @>>> J_{01}
\end{CD}
$$
in $\Delta^o\Pos$ is {\em Segal-cocartesian} if the induced map
\begin{equation}\label{N.L.J.eq}
N^L(J_0) \copr_{N^L(J)} N^L(J_1) \to N^L(J_{01})
\end{equation}
is an injective Segal equivalence.

\begin{lemma}\label{seg.pos.le}
For any finite simplicial set $X$ and subset $Y \subset X$, the
embedding $f:Y \to X$ is absolutely cofibrant in $\Delta^o\Posf^-
\supset \Delta^o\Sets$, and so is the embedding $\wt{f}:\Cyl(f) \to
X \times [1]$ equal to $\id$ resp.\ $f$ over $1$
resp.\ $0$. Moreover, in both cases, all the corresponding squares
\eqref{abs.sq} are Segal-cocartesian.
\end{lemma}

\proof{} The claim for $f$ is obvious --- since $X$ and $Y$ are
equipped with a discrete order, $f$ is pointwise-split in the sense
of Example~\ref{pos.triv.exa}, thus absolutely cofibrant by
Example~\ref{abs.exa}, and then $N^L$ sends any square
\eqref{abs.sq} for $f$ to a cocartesian square in
$\Delta^o\Delta^o\Sets$, so that the map \eqref{N.L.J.eq} is an
isomorphism. To check that $\wt{f}$ is absolutely cofibrant, it
suffices to check it pointwise. But for any $[n] \in \Delta$, we can
let $X([n])' = X([n]) \ssetminus Y([n])$, and then we have
$\Cyl(f)([n]) \cong (Y([n]) \times [1]) \copr X([n])'$, and
$\wt{f}([n])$ is the disjoint union of an identity map and the
embedding $\id \times t:X([n])' \to X([n])' \times [1]$; this is
also absolutely cofibrant by Example~\ref{abs.exa}. Moreover, for
each $[n] \in \Delta$, the square \eqref{abs.sq} corresponding to
some map $f:\Cyl(f) \to J$ evaluates to the square \eqref{cyl.o.sq}
for the induced map $f \circ t:X([n])' \to J([n])$. Therefore
\eqref{N.L.J.eq} is indeed injective by Remark~\ref{S.rem}, and
moreover, by \eqref{b.f.s}, exactly the same argument as in
Lemma~\ref{bva.le} and Proposition~\ref{bi.seg.prop} shows that
\eqref{N.L.J.eq} is a countable composition of pushouts of
coproducts of horn embeddings \eqref{bi.horn.bis.eq}. Since by
Lemma~\ref{f.CW.le} and Corollary~\ref{W.C.corr}, the class of
injective Segal equivalences is stable under coproducts, pushouts
and countable compositions, \eqref{N.L.J.eq} is indeed an injective
Segal equivalence.
\endproof

Now for any horn embedding $v^k_n$ of \eqref{horn.eq}, let $\W^k_n =
\Cyl(v^k_n) \in \Delta^o\Posf$ be its pointwise cylinder, where as
in Lemma~\ref{seg.pos.le}, we equip $\V^k_n$ and $\Delta_n$ with
discrete order, and let
\begin{equation}\label{pos.horn.eq}
  w^k_n:\W^k_n \to \Delta_n \times [1]
\end{equation}
be the embedding given by $\id$ resp.\ $v^k_n$ over $1 \in [1]$
resp.\ $0 \in [1]$.

\begin{lemma}\label{t.fib.le}
Say that $J \in \Delta^o\Pos$ is {\em $t$-fibrant} if the nerve
$N^L(J)$ is $t$-fibrant in the sense of
Definition~\ref{t.fib.def}. Then $J$ is $t$-fibrant if and only if
it has the lifting property with respect to all horn embddings
\eqref{horn.eq} and their cylinders \eqref{pos.horn.eq}.
\end{lemma}

\proof{} Being a nerve, the bisimplicial set $X=N^L(J)$ is in any
case a Segal space, so by induction on $n$, it is $t$-fibrant if and
only if $X([0])$ is fibrant, and $X(t):X([1]) \to X([0])$ is a Kan
fibration. The former means the lifting property
w.r.t. \eqref{horn.eq}, and since $N([1]) \cong \Delta_1$, the
latter means the lifting property w.r.t. \eqref{pos.horn.eq}.
\endproof

We can now construct $t$-fibrant replacement for simplicial
partially ordered sets. The procedure is the same as in
Lemma~\ref{simp.S.le} but we use both the horn embeddings
\eqref{horn.eq} and their cylinders \eqref{pos.horn.eq}. We apply
the extension construction of Subsection~\ref{ab.subs}. By
Lemma~\ref{seg.pos.le}, the maps $v^k_n$ and $w^k_n$ of
Lemma~\ref{t.fib.le} are absolutely cofibrant in the abundant
$C$-category $\Delta^o\Posf^-$. For any $J = J_0 \in
\Delta^o\Posf^-$, we can then define
$$
J_1 = \Ex(\{v^k_n,w^k_n\})(J),
$$
where $\Ex$ is the extension \eqref{ex.I.sq} with respect to the
collection $\{v^k_n,w^k_n\}$ of maps $v^k_n$, $w^k_n$ for all
integers $n \geq k \geq 0$. This comes equipped with a map $a_0:J_0
\to J_1$, and we can then iterate the construction to
obtain maps $a_n:J_n \to J_{n+1}$, $n \geq 0$, and let $J_\infty =
\colim_nJ_n$, together with the map
\begin{equation}\label{a.inf.eq}
  a_\infty:J \to J_\infty
\end{equation}
obtained as the countable composition of the maps $a_n$. Note that
since the maps $a_n$ are right-closed, we have $J_\infty \in
\Delta^o\Posf^-$. Moreover, $J_\infty$ is functorial in $J$, so that
for any category $I$ and functor $J:I \to \Delta^o\Posf^-$, we can
apply it pointwise and obtain a functor $J_\infty:I \to
\Delta^o\Posf^-$ and a map $a_\infty:J \to J_\infty$.

\begin{lemma}\label{t.repl.le}
For any $J \in \Delta^o\Posf^-$, $J_\infty \in \Delta^o\Posf^-$ is
$t$-fibrant, globally right-bounded if so is $J$, and $a_\infty$ is
a Segal equivalence. Moreover, for any $I \in \Posf$ and $J:I \to
\Delta^o\Posf^-$ cofibrant with respect to the co-standard
$C$-structure, $a_\infty:J \to J_\infty$ is in $C$.
\end{lemma}

\proof{} Since both $\V^k_n$ and $W^k_n$ are compact as objects in
$\Delta^o\Posf^-$, the fact that $J_\infty$ is $t$-fibrant
immediately follows from the small object argument and
Lemma~\ref{t.fib.le}. Moreover, the map $a_1$ is right-closed, and
if $\dim J([n]) \leq m$ for any $[n] \in \Delta$ and some $m \geq
0$, then $\dim J_1([n]) \leq m+1$, again for any $[n] \in
\Delta$. Then if $J$ is equipped with a co-height function $\hht^i:J
\to \N^o$ with globally finite-dimensional comma-fibers $n
\setminus_{\hht^\iota} J$, then $J_\infty$ is filtered by globally
finite-dimensional right-closed subsets $(n \setminus_{\hht^\iota}
J)_n \subset J_\infty$, $n \in \N^o$, and this defines a co-height
function $\hht^\iota:J_\infty \to \N^o$. For the second claim, since
injective Segal equivalences are closed under countable compositions
by Corollary~\ref{W.C.corr}, it suffices to prove that $a_1:J \to
J_1$ is a Segal equivalence. By \eqref{ex.I.eq}, $J_1$ can be expressed
as a co-standard colimit of some functor $J_\idot:T^< \to
\Delta^o\Posf^-$ sending $o$ to $J$ and $t \in T$ to some
right-closed embedding $a_t:J \to J_t$ obtained as a pushout of a
horn embedding $v^k_n$ or its cylinder $w^k_n$. Since $N^L$ sends
co-standard colimits to standard colimits, Corollary~\ref{W.C.corr}
again shows that it suffices to prove that each $a_t$ is a Segal
equivalence. Then by Lemma~\ref{seg.pos.le}, it actually suffices to
check this for the pushout of the corresponding nerves, and this is
obvious, since both $N^L(v^k_n)$ and $N^L(w^k_n)$ are not only Segal
equivalences but weak equivalences in
$\Delta^o\Delta^o\Sets$. Finally, since $\Delta^o\Posf^-$ is
excisive, the last claim immediately follows from
Lemma~\ref{ext.ab.le}.
\endproof

\proof[Proof of Proposition~\ref{tw.seg.prop}.] First of all, since
right-closed embeddings are strict, $J$ factors through $\bRelf^-$,
and then by Lemma~\ref{bN.N.le}, we may replace the nerves
$\Nn_\dm$, $\Nn_\dm(I)$ with their strict versions $\bN_\dm$,
$\bN_\dm(I)$. Moreover, since by Corollary~\ref{W.C.corr}, the class
of injective Segal equivalences is closed under countable
compositions, it suffices to prove the claim after restricting to
the skeleta $\sk_nI$, $n \geq 0$; thus we may assume right away that
$I$ is finite-dimensional. We have the map
\begin{equation}\label{F.J.I.eq}
p:L(I) \setminus^\iota J^\hdot \to J^\hdot
\end{equation}
of \eqref{F.J.eq}. Since $J$ factors through $\Relf^-_C$ and
right-closed embeddings preserve the co-heights, the co-height
functions $\hht^\iota:J(i) \to R(\N^o)$, $i \in I$ taken together
define a co-height function $\hht^\iota:J^\hdot \to R(\N^o)$, with
finite-dimensional comma-fibers $n \setminus J^\hdot$. Then if we
replace $J^\hdot$ in \eqref{F.J.I.eq} with $n \setminus J^\hdot$,
the map \eqref{F.J.I.eq} is fiberwise-reflexive by construction and
bianodyne by Lemma~\ref{F.J.le}. Thus $\bN_\dm(I)(p)$ is a weak
equivalence, and $\bN_\dm(p)$ is a Segal equivalence by
Lemma~\ref{bi.del.le} and Lemma~\ref{NQ.bi.le}. Since the functors
$\N \to \Relf^-$ given by $n \mapsto n \setminus J^\hdot$ and $n
\mapsto L(I) \setminus^\iota (n \setminus J)$ are cofibrant, and by
Lemma~\ref{Q.adj.le} and Corollary~\ref{W.C.corr}, both weak
equivalences and Segal equivalences between cofibrant functors from
$\N$ are closed under colimits, the same holds for the map
\eqref{F.J.I.eq} itself. Therefore we may replace $J^\hdot$ with
$L(I) \setminus^\iota J^\hdot$ and assume right away that $J:I^o \to
\Relf^-$ is cofibrant with respect to the co-standard $C$-category
structure.

Moreover, by Lemma~\ref{hom.bi.le}, the functor $\bsE_F$ of
\eqref{E.2.eq} is a $C$-functor with respect to the co-standard
$C$-category structures, so that the functor $\bsE_F(J):I^o \to
\Delta^o\Posf^-$ is then also cofibrant, and Lemma~\ref{t.repl.le}
provides a cofibrant functor $\bsE_F(J)_\infty:I^o \to
\Delta^o\Posf^-$ such that $\bsE_F(J)_\infty(i)$ is $t$-fibrant for
any $i \in I$, and a map $a_\infty:\bsE_F(J) \to \bsE_F(J)_\infty$
in the class $C$ such that $a_\infty(i)$ is a Segal equivalence for
any $i \in I$. By \eqref{bsE.bE} and Lemma~\ref{N.d.le},
$\bsE_F(J)(i) \in \Delta^o\Posf^-$ is globally right-bounded for any
$i \in I$, and then by Lemma~\ref{t.repl.le}, so is
$\bsE_F(J)_\infty(i)$.

Now, the map $N^L(a_\infty)$ is pointwise an injective Segal
equivalence, and then by Lemma~\ref{I.seg.le}, the family $\Hh(X)$
is constant along $N^L(a_\infty)$, so by Remark~\ref{W.C.rem} and
Corollary~\ref{W.C.corr}, we have $w = w' \circ N^L(a_\infty)$ for
some map $w':N^L(\bsE_F(J)_\infty) \to X$. Then $w'$ is a pointwise
Segal equivalence. Moreover, factor it as
\begin{equation}\label{w.facto}
\begin{CD}
  N^L(\bsE_F(J)_\infty) @>{c}>> X' @>{f}>> X,
\end{CD}
\end{equation}
with $c \in C \cap W$ and $f \in F$. Then $f$ is a Segal
equivalence, and since $X$ is fibrant, so is $X'$. However, for any
$t$-fibrant bisimplicial set $Y$ and integers $n > l > 0$, the
vertical maps in the corresponding square \eqref{segal.sq} are in
$F$, so that if we have a weak equivalence $Y \to Y'$, with
$t$-fibrant $Y$ and fibrant $Y'$, then $Y'$ is a Segal space as soon
as so is $Y$. We conclude that $X'$ in \eqref{w.facto} is pointwise
a fibrant Segal space, and then the Segal equivalence $f$ is a weak
equivalence by Lemma~\ref{bous.le}~\thetag{i}. Therefore the whole
map $w'$ is a weak equivalence between $t$-fibrant $I$-Segal spaces,
and then by Lemma~\ref{tw.seg.le}, it becomes a Segal equivalence
after taking the twisted simplicial expansion. Thus to finish the
proof, it remains to show the same for the pointwise Segal
equivalence $N^L(a_\infty)$.

But recall that $\bsE_F(J)$ and $\bsE_F(J)_\infty$ are
cofibrant. Then for any cofibrant functor $X:I^o \to
\Delta^o\Posf^-$, Lemma~\ref{pos.cof.le} applied to $X^\iota$
pointwise over $\Delta$ provides an isomorphism $X \cong \sE_I(X')$,
where $X' = \colim_IX$ is equipped with a natural augmentation
$\chi:X' \to I$. Thus if we denote $J' = \colim_{I^o}\sE_F(J)$,
$J'_\infty = \colim_{I^o}\sE_F(J)_\infty$, $a'_\infty =
\colim_{I^o}a_\infty$, we have a commutative diagram
\begin{equation}\label{J.inf.dia}
\begin{CD}
\colim_{\Tw(I)}N^L((I \setminus J')_f) @>{\colim N^L(a'_\infty)_f}>>
\colim_{\Tw(I)}N^L((I \setminus J'_\infty)_f)\\
@VVV @VVV\\
N^L(I \setminus J') @>{N^L(a'_\infty)}>> N^L(I \setminus J'_\infty),
\end{CD}
\end{equation}
where the vertical arrows on the left and on the right are the Segal
equivalences \eqref{b.chi.bis.eq} of Lemma~\ref{b.chi.le}, and the
bottom map is the twisted simplicial expansion of the map
$N^L(a_\infty)$. Then to finish the proof, it remains to show that
the top arrow in \eqref{J.inf.dia} is a Segal equivalence. Moreover,
by Corollary~\ref{W.C.corr}, it suffices to check it for
$N^L(a'_\infty)_f$ for each $f \in \Tw(I)$ represented by a pair $i
\leq i'$. Since $N^L(a'_\infty)_f = N^L(a_\infty(i')) \times
\id_{L(N(I / i))}$, we are done by Corollary~\ref{bi.seg.corr}.
\endproof

\chapter{Enhanced category theory.}\label{enh.ch}

In this chapter, we are finally able to get to our main task --- we
build enhanced category theory.

We start with the notion of a {\em reflexive family}
(Section~\ref{refl.sec}). This is a very simple technical gadget ---
roughly speaking, a reflexive family of categories over a big enough
full subcategory $\I \subset \Pos$ is a Grothendieck fibration $\C
\to \I$ such that for any left-reflexive full embedding $f:J \to J'$
in $\I$, with adjoint map $f_\dg:J' \to J$, the isomorphism $f^*
\circ f_\dg^* \cong \id$ induces an adjunction between the
transition functors $f^*:\C_{J'} \to \C_J$ and $f_\dg^*:\C_{J'} \to
\C_J$. A trivial example of such a family is the family
$\Unf(\E,\I)$ with fibers $\Unf(\E,\I)_J \cong J^o\E$, for some
category $\E$. Suprisingly, even with such a naive definition, one
can already prove something --- for example,
Proposition~\ref{trunc.prop} shows that $\Unf(-,\I)$ is
right-adjoint to the truncation operation $\C \mapsto \C_\ppt$. We
also define reflexive families augmented by a partially ordered set
$I$, and we prove that in good cases, the Kan extension operation
$f_*$ associated to a map $f:I \to I'$ sends reflexive families to
reflexive families (the precise statement is
Corollary~\ref{pf.corr}).

Section~\ref{refl.sec} ends with Subsection~\ref{expo.subs} that
contains one reconstruction result, Proposition~\ref{exp.prop},
saying that --- again in good cases --- a reflexive family is
completely determined by the underlying family of groupoids; this
establishes a link between enhanced categories and families of
groupoids studied in Chapter~\ref{nerve.ch}.

Enhanced categories themselves are defined in Section~\ref{enh.sec},
as reflexive families semicartesian over certain standard
commutative squares, and the remainder of Section~\ref{enh.sec}
collects definitions and results that do not need representability
theorems. Again somewhat suprisingly, this is more than one would
expect. In particular, two notions that are simply borrowed {\em
  verbatim} from the usual category theory are fully faithful
embeddings (Corollary~\ref{enh.ff.corr}) and adjoint pairs
(Definition~\ref{enh.refl.def}). We also construct basic examples of
enhanced categories such as enhanced localizations of a model
category (Lemma~\ref{enh.mod.le}, and enhanced localization of the
chain-homotopy category of an additive category
(Lemma~\ref{enh.ho.le}).

In Section~\ref{enh.repr.sec} we introduce our representability
result, Theorem~\ref{enh.thm}, and see what it gives us. The theorem
itself is an easy combination of Proposition~\ref{exp.prop} and
representability theorems of Chapter~\ref{nerve.ch} (that in turn
depend on Brown Representability of Chapter~\ref{semi.ch}). The main
corollary of Theorem~\ref{enh.thm} is Proposition~\ref{univ.prop}
saying that, roughly speaking, any small enhanced category is a
localization of a partially ordered set $I$, and in a rather
controlled way (in particular, the class $W$ of maps in $I$ that we
use for localization fits into a biorder on $I$ in the sense of
Section~\ref{bi.sec}). This fact itself has two main corollaries ---
firstly, Corollary~\ref{fun.h.corr} that constructs enhanced functor
categories, and secondly, Corollary~\ref{sq.h.corr} and
Lemma~\ref{sq.h.le} showing that enhanced categories admit
semicartesian products, and those products have a universal
property. This is our main tool for working with enhanced
categories, and it is heavily used throughout the remainder of this
chapter. As a less immediate corollary of
Proposition~\ref{univ.prop}, we also construct an enhancement for
the category of enhanced small categories
(Subsection~\ref{cath.subs}) and enhanced versions of the lax
functor categories $\Cat \bb I$, $I \bbd \Cat$ of
Subsection~\ref{fun.subs} (this is Subsection~\ref{lax.subs}).

Starting from Section~\ref{yo.sec}, we develop enhanced category
theory in proper sense, following the template set in
Chapter~\ref{cat.ch}. The slogan is, replace cartesian products with
semicartesian products of Subsection~\ref{univ.subs}, and then
everything works. Section~\ref{yo.sec} starts with enhanced versions
of cylinders and comma-categories; an important technical result is
Lemma~\ref{adj.eq.le} that gives a characterization of adjunction
stable under semicartesian products. We then build a theory of
enhanced fibrations and cofibrations, in parallel to
Section~\ref{fib.sec}, prove an enhanced version of the Grothendieck
construction (Proposition~\ref{groth.prop}), and use it to derive
the enhanced version of the Yoneda Lemma
(Proposition~\ref{yo.prop}). We also construct enhanced relative
functor categories, in parallel to Section~\ref{fun.fib.sec}; this
is Subsection~\ref{enh.rel.subs}.

Section~\ref{enh.lim.sec} deals with enhanced limits, colimits and
Kan extensions, and here again, the story is largely parallel to its
unenhanced version of Subsection~\ref{lim.subs}. Representability
and representing objects only enter through semicartesian squares,
with one important exception --- the fact that the category of small
enhanced categories is enhanced-cocomplete
(Proposition~\ref{ccat.coco.prop}). The two things that have no
direct parallel in Subsection~\ref{lim.subs} are an enhanced version
of nerves and simplicial replacements given in
Subsection~\ref{enh.repl.subs}, and an envelope construction of
Subsection~\ref{env.subs} that allows one to formally add enhanced
colimits of a specified kind via a universal construction
(Lemma~\ref{env.le}).

Finally, in Section~\ref{enh.large.sec}, we discuss what to do with
large enhanced categories. This is even more problematic than the
situation in the usual category theory, since in the enhanced
setting, even semicartesian products only exist for small enhanced
categories, and this is an undispensable tool for building a
workable theory. As a solution, we adapt the machinery of filtered
colimits and accessible categories of \cite{sga4} and \cite{GU}. A
standard reference for this material is a wondeful book \cite{AR},
and the enhanced theory that we build is again very similar to its
unenhanced prototype (for comparison, there is a short overview of
the unenhanced story in \cite{ka.acc} written in parallel to
Section~\ref{enh.large.sec}). We construct enhanced Karoubi
completions in Subsection~\ref{enh.kar.sus}, then develop the basic
machinery of filtered colimits and inductive completions, and then
define accessible enhanced categories and show how to adapt our main
technical tools such as semicartesian products to the accessbile
world. We end the text with a brief discussion of accessible
enhanced fibrations, and an accessible version of relative enhanced
functor categories of Subsection~\ref{enh.rel.subs}.

\section{Reflexive families.}\label{refl.sec}

\subsection{Definitions and basic properties.}

We now turn to families of categories that are not necessarily
families of groupoids. Assume given a full subcategory $\I \subset
\Pos$ that is ample in the sense of Definition~\ref{amp.def}. From
now on, a {\em family of categories} over $\I$ means a fibration $\C
\to \I$, and it is small if it is small as a functor (equivalently,
has essentially small fibers $\C_i$). Since $\I$ is ample, it
contains $\ppt \in \Pos$, and for any family of categories $\C \to
\I$, \eqref{2.adj} provides a functor
\begin{equation}\label{nondeg.eq}
  \C \to \eps_*\eps^*\C,
\end{equation}
cartesian over $\I$, where $\eps = \eps(\ppt):\ppt \to \I$ is the
embedding onto $\ppt \in \I$. Explicitly, for any $J \in \I$, and an
element $j \in J$, we denote by
\begin{equation}\label{ev.J.eq}
  \ev_j = \eps(j)^*:\C_J \to \C_\ppt
\end{equation}
the transition functor for the embedding $\eps(j):\ppt \to J$ onto
$j \in J$. Then \eqref{nondeg.eq} over some $J \in \I$ is the
product of the functors $\ev_j$ for all $j \in J$.

\begin{defn}\label{cat.nondege.def}
A family of categories $\C \to \I$ is {\em non-degenerate} if the
corresponding functor \eqref{nondeg.eq} is conservative.
\end{defn}

Another consequence of Definition~\ref{amp.def} is that $[1] \in
\I$, and $[1] \times J \in \I$ for any $J \in \I$. Denote by $s,t:J
\to J \times [1]$ the embeddings onto $J \times \{0\}$ resp.\ $J
\times \{1\}$, and let $e:J \times [1] \to J$ be the projection.

\begin{defn}\label{cat.refl.def}
A family of categories $\C \to \I$ is {\em reflexive} resp.\ {\em
  co-reflexive} if for any $J \in \I$, with the transition functors
$s^*,t^*:\C_{J \times [1]} \to \C_J$ and $e^*:\C_J \to \C_{J \times
  [1]}$, the isomorphism $s^* \circ e^* \cong \id$ resp.\ $t^* \circ
e^* \cong \id$ induced by the isomorphism $e \circ s \cong \id$
resp.\ $e \circ t \cong \id$ defines an adjunction between $e^*$ and
$s^*$ resp.\ $e^*$ and $t^*$.
\end{defn}

\begin{exa}
For any ample $\I \subset \Pos$, $\iota(\I) \subset \Pos$ is also
ample, and for any reflexive resp.\ co-reflexive family $\C$ over
$\I$, $\iota^*\C$ is tautologically co-reflexive resp.\ reflexive
over $\iota(\I)$.
\end{exa}

\begin{exa}\label{sub.refl.exa}
Assume given a reflexive resp.\ coreflexive family $\C \to \I$, and
a full subcategory $\C' \subset \C$ such that the induced functor
$\C' \to \I$ is a fibration, and the embedding $\C' \to \C$ is
cartesian over $\I$. Then $\C'$ is reflexive resp.\ co-reflexive.
\end{exa}

If $\C$ is a reflexive family of categories over
$\I$, then for any $J \in \I$ and $c \in \C_{J \times [1]}$, the
adjunction map $a(c):c \to e^*s^*c$ is functorial with respect to
$c$, and sending $c$ to $t^*a(c)$ provides a functor
\begin{equation}\label{nu.J.eq}
\nu_J:\C_{J \times [1]} \to \Ar(\C_J) = \Fun([1],\C_J)
\end{equation}
isomorphic to $t^*$ resp.\ $s^*$ after evaluation at $0 \in [1]$
resp.\ $1 \in [1]$. For any two maps $f_0,f_1:J \to J'$ in $\I$ such
that $f_0 \leq f_1$, \eqref{nu.J.eq} induces a natural map
$\nu(f_0,f_1):f_1^* \to f_0^*$ (factor $f_0 \copr f_1:J \copr J \to
J'$ through $J \times [1]$), we have $f^*(\nu(f_0,f_1)) = \nu(f_0
\circ f,f_1 \circ f)$ for any map $f:J'' \to J$, and we have
$\nu(f_0,f_1) = \id$ if $f_0=f_1$. For any map $f:J \to J'$ in $\I$
such that $\Cyl(f),\Cyl^o(f) \in \I$ --- in particular, for any full
embedding $f$ --- we then have functors
\begin{equation}\label{cyl.J.eq}
\nu_f:\C_{\Cyl(f)} \to \C_{J'} /_{f^*} \C_J, \qquad
\nu^\perp_f:\C_{\Cyl^o(f)} \to \C_J \setminus_{f^*} \C_{J'}
\end{equation}
given by $t^* \times s^*$ with the maps $\nu(s,t \circ f)$,
resp.\ $\nu(s \circ f,t)$. If $f$ is right-reflexive with some
adjoint $f^\dg$, so that $\Cyl(f) \cong \Cyl^o(f^\dg)$, then $\nu_f
\cong \nu^\perp_{f^\dg}$. For $f = \id:J \to J$, both functors
\eqref{cyl.J.eq} reduce to \eqref{nu.J.eq}.

\begin{exa}\label{unf.E.exa}
For any category $\E$, the category $\I \bb \E$ with the fibration
$\I \bb \E \to \I$ of \eqref{fun.bb} is a co-reflexive family over
$\I$ in the sense of Definition~\ref{cat.refl.def}. Dually, the
category
\begin{equation}\label{unf.eq}
  \Unf(\E,\I) = \I \bb^{\iota} \E \cong \iota^*(\I \bb \E)
\end{equation}
equipped with the composition $\I^o \bb \E \to \iota(\I) \to \I$ of
the fibration \eqref{fun.bb} and the equivalence $\iota:\iota(\I)
\to \I$ is a reflexive family, with fibers given by $\Unf(\E,\I)_J =
J^o\E = \Fun(J^o,\E)$. The families $\I \bb \E$ and $\Unf(\E,\I)$
are non-degenerate in the sense of
Definition~\ref{cat.nondege.def}. The functor \eqref{ev.J.eq} for
$\Unf(\E,\I)$ is the evaluation functor, \eqref{nu.J.eq} is the
equivalence $\Fun(J^o \times [1]^o,\E) \to \Fun([1],J^o\E)$, where
we identify $[1] \cong [1]^o$, and \eqref{cyl.J.eq} are the
equivalences \eqref{cyl.com}. We have
\begin{equation}\label{unf.perp}
  \Unf(\E^o,\I)^o_\perp \cong \I \bb \E^o,
\end{equation}
where $\Unf(\E^o,\I)^o_\perp \to \I$ is the transpose-opposite
fibration. In terms of relative functor categories of
Definition~\ref{fl.def}, we have
\begin{equation}\label{refl.fun}
\Unf(\E,\I) \cong \Fun(\I^\hdot | \I,\E), \qquad \I \bb \E \cong
\Fun(\I_\idot | \I,\E),
\end{equation}
where $\nu_\idot:\I_\idot \to \I$, $\nu^\hdot:\I^\hdot \cong
\I_\idot^{o\perp} \to \I$ are the cofibrations with fibers $J$,
$J^o$ induced by $\Pos_\idot,\Pos^\hdot \to \Pos$.
\end{exa}

\begin{lemma}\label{refl.J.le}
Assume given a reflexive family of categories $\C$ over an ample
full subcategory $\I \subset \Pos$. Then for any left
resp.\ right-reflexive full embedding $f:J' \to J$ in $\I$, with the
adjoint map $f_\dg:J \to J'$, the isomorphism $f^* \circ f_\dg^*
\cong \id$ resp.\ $\id \cong f^* \circ f_\dg^*$ defines an
adjunction between $f_\dg^*$ and $f^*$ resp.\ $f^*$ and $f_\dg^*$.
\end{lemma}

\proof{} The other adjunction map is $\nu(f \circ f_\dg,\id)$
resp.\ $\nu(\id,f \circ f_\dg)$.
\endproof

\begin{corr}\label{refl.refl.corr}
A family of categories $\C \to \I$ is reflexive if and only if for
any right-reflexive full embedding $f:J \to J'$, with adjoint
$f_\dg:J' \to J$, the isomorphism $f^* \circ f_\dg^* \cong \id$
defines an adjunction between $f^*$ and $f_\dg^*$.
\end{corr}

\proof{} For any $J \in \I$, $s:J \to J \times [1]$ is
right-reflexive, and $e = s_\dg$.
\endproof

\begin{exa}\label{corefl.exa}
For any $J \in \I$, the embedding $t:J \to J \times [1]$ is
left-reflexive, and Lemma~\ref{refl.J.le} shows that for any
reflexive family of categories $\C$, the isomorphism $\id \cong e
\circ t$ defines an adjunction between $t^*$ and $e^*$. In other
words, $e^*$ actually has both a left-adjoint $s^*$ and a
right-adjoint $t^*$, and \eqref{nu.J.eq} is the corresponding map
\eqref{adj.fu}. This shows that for any reflexive family $\C \to
\I$, the transpose-opposite fibration $\C^o_\perp \to \I$ is a
co-reflexive family. The converse also holds (apply the same
argument to $\iota^*\C$).
\end{exa}

\begin{defn}\label{cat.prp.def}
A reflexive family of categories $\C$ over $\I$ is {\em separated}
if the functor \eqref{nu.J.eq} is an epivalence for any $J \in \I$,
and {\em proper} if so is $\nu^\perp_f$ of \eqref{cyl.J.eq} for any
right-closed full embedding $f$ in $\I$. A co-reflexive family $\C
\to \I$ is {\em separated} or {\em proper} if so the reflexive
family $\C^o_\perp \to \I$.
\end{defn}

\begin{exa}
For any category $\E$, the reflexive family $\Unf(\E,\I)$ of
Example~\ref{unf.E.exa} is proper in the sense of
Definition~\ref{cat.prp.def}.
\end{exa}

\begin{lemma}\label{cyl.J.le}
A reflexive family of categories $\C$ over $\I$ is separated if and
only if for any left-reflexive full embedding $f:J \to J'$ in $\I$,
with adjoint map $g:J' \to J$, the functor $\nu^\perp_f \cong \nu_g$
of \eqref{cyl.J.eq} is an epivalence. A reflexive family of
categories $\C$ over $\I$ is proper if and only if $\nu^\perp_f$ is
an epivalence for any map $f:J \to J'$ in $\I$ such that $\Cyl^o(f)
\in \I$.
\end{lemma}

\proof{} In the first claim, the ``if'' part is clear: $\nu_J$ is
$\nu^\perp_f$ for $f=\id:J \to J$. For the ``only if'' part, assume
that $\C$ is separated. By Lemma~\ref{refl.J.le}, $g^*:\C_J \to
\C_{J'}$ is fully faithful and left-adjoint to $f^*$. Consider the
left-reflexive full embedding $\wt{g}:\Cyl^o(f) \cong \Cyl(g) \to J'
\times [1]$ dual to \eqref{refl.dia}, with the adjoint map
$\wt{g}$. Again by Lemma~\ref{refl.J.le}, $\wt{g}^*$ is fully
faithful and left-adjoint to $\wt{f}^*$, and we have a commutative
square
\begin{equation}\label{nu.fJ.sq}
\begin{CD}
\C_{\Cyl^o(f)} @>{\wt{g}^*}>> \C_{J' \times [1]}\\
@V{\nu^\perp_f}VV @VV{\nu_{J'}}V\\
\C_{J'} \setminus_{f^*} \C_J @>{\id \setminus g^*}>> \C_{J'}
\setminus \C_{J'}
\end{CD}
\end{equation}
with fully faithful right-reflexive horizontal arrows. Then the
square is cartesian by Lemma~\ref{epi.sq.le}, and $\nu^\perp_f$ is a
pullback of an epivalence. For the second claim, the ``if'' part is
tautological; for the ``only if'' part, consider the right-closed
full embedding $t:J' \to \Cyl^o(f)$, and note that the
left-reflexive projection $s_\dg:\Cyl^o(f) \to J$ adjoint to the
right-reflexive full embedding $s:J \to \Cyl^o(f)$ induces a
left-reflexive projection $g:\Cyl^o(t) \to \Cyl^o(f)$, with fully
faithful adjoint map $\Cyl^o(f) \to \Cyl^o(f)$. Then $g^*$ is a
fully faithful embedding by Lemma~\ref{refl.J.le}, and $\nu^\perp_f$
is obtained by retricting the epivalence $\nu^\perp_t$ to the
essential image of $g^*$.
\endproof

Explicitly, for any category $I$, objects in the category $\I \bb I$
of Example~\ref{unf.E.exa} are $I$-augmented partially ordered sets
$\langle J,\alpha \rangle$, $J \in \I$, $\alpha:J \to I$, and objects in
$\Unf(I,\I)$ are $I^o$-augmented partially ordered sets
$\langle J,\alpha \rangle$, $J \in \I$, $\alpha:J \to I^o$. Say that
a map $f$ in $\I \bb I$ or $\Unf(I,\I)$ is {\em $I$-strict} if it is
cartesian over $\I$. Then for any left-reflexive $I$-strict map $f$
in $\I \bb I$, the adjoint $f_\dg$ is canonically also a map in $\I
\bb I$, and similarly for right-reflexive $I$-strict maps in
$\Unf(\I,I)$.

\begin{lemma}\label{refl.refl.J.le}
For any fibration $\C \to \I \bb I$, the composition $\C \to \I \bb
I$ is a co-reflexive family over $I$ iff for any $I$-strict
left-reflexive full embedding $f$ in $\I \bb I$, with the adjoint
map $f_\dg$, the isomorphism $\id \cong f^* \circ f_\dg^*$ defines
an adjunction between $f_\dg^*$ and $f^*$. Dually, for any fibration
$\C \to \Unf(I,\I)$, the composition $\C \to \Unf(I,\I) \to \I$ is a
reflexive family iff for any $I$-strict reflexive full embedding $f$
in $\Unf(I,\I)$, with the adjoint map $f_\dg$, the isomorphism $f^*
\circ f_\dg^* \cong \id$ defines an adjunction between $f^*$ and
$f_\dg^*$.
\end{lemma}

\proof{} The arguments for the two claims are dual, so we only do
the second one. The composition $\C \to \Unf(I,\I) \to \I$ is a
fibration, and for any $J \in \I$, we have the fibration $\C_J \to
J^oI \cong \Fun(J,I^o)^o$ whose fiber over some map $\alpha:J \to
I^o$ is $\C_{\langle J,\alpha \rangle}$. Now use
Corollary~\ref{refl.refl.corr} and Lemma~\ref{fib.adj.le}~\thetag{ii}.
\endproof

\begin{lemma}\label{refl.pb.le}
Assume given a functor $\gamma:I' \to I$ between categories $I$,
$I'$, and a cartesian square
$$
\begin{CD}
  \E' @>>> \E\\
  @V{\pi'}VV @VV{\pi}V\\
  \Unf(I',\I) @>{\Unf(\gamma)}>> \Unf(I,\I)
\end{CD}
$$
such that $\E$ is a reflexive family over $I$, and $\pi$ is
cartesian over $I$. Then $\E'$ with the composition functor
$\pi':\C' \to \Unf(I',\I) \to \I$ is a reflexive family over $\I$,
proper and/or separated if so is $\E$.
\end{lemma}

\proof{} If we have two cartesian squares of categories interpreted
as fibrations $\C,\C' \to [1]^2$, as in Example~\ref{sq.fib.exa},
and a functor $\gamma:\C \to \C'$ cartesian over $[1]^2$ that is
left-reflexive over $0 \times 0,0 \times 1,1 \times 0 \in [1]^2$,
then by Lemma~\ref{fib.adj.le}~\thetag{ii}, \eqref{sec.sq.eq} and
Lemma~\ref{sec.adj.le}, it is also left-reflexive over $1 \times
1$. This proves the first claim. For the second, note that
$\nu_f^\perp$ for the families $\Unf(I',\I)$ and $\Unf(I,\I)$ is an
equivalence, and then $\nu_f^\perp$ for $\E'$ is a pullback of
$\nu^\perp_f$ for $\E$.
\endproof

\subsection{Truncation and unfoldings.}

If we take $J = \ppt$, then for any reflexive family $\C \to \I$,
\eqref{nu.J.eq} provides a functor $\C_{[1]} \to
\Fun([1],\C_\ppt)$. It turns out that one can do the same not only
for $[1] \in \I$ for for all $J \in \I$, and in a compatible way.

\begin{prop}\label{trunc.prop}
Assume given a reflexive family of categories $\C \to \I$, and a
category $\C'$. Then any functor $\gamma:\C_\ppt \to \C'$ extends to
a functor $\gamma':\C \to \Unf(\C',\I)$, cartesian over $\I$, and the
extension is unique up to a unique isomorphism.
\end{prop}

\proof{} Let $\nu^\hdot:\I^\hdot \to \I$ be as in \eqref{refl.fun}. By
Definition~\ref{amp.def}, $\I$ contains all subsets $J' \subset J$
of any $J \in \I$. Therefore \eqref{nu.io.pos} induces a functor
$\nu_<:\I^\hdot \to \I$, and the tautological embeddings $j
\setminus J \to J$ together provide a functorial map $a:\nu_< \to
\nu^\hdot$, while we also have the tautological map $b:\nu_< \to
\ppt$ to the constant functor with value $\ppt \in \I$. The
corresponding transition functors then fit into a diagram
\begin{equation}\label{unf.dia}
\begin{CD}
  \nu^{\hdot *}\C @>{a^*}>> \nu_<^*\C @<{b^*}<< \C_\ppt.
\end{CD}
\end{equation}
Moreover, for any $\langle J,j \rangle \in \I_\idot$, the projection
$j \setminus J \to \ppt$ is left-reflexive, with the adjoint
right-reflexive embedding $\ppt \to j \setminus J$ onto $j$, so by
Lemma~\ref{refl.J.le}, $b^*$ is left-reflexive over $\langle J,j
\rangle$, with the left-adjoint functor \eqref{ev.J.eq}. Then $b^*$
itself is left-reflexive by Lemma~\ref{fib.adj.le}~\thetag{ii}, with some
left-adjoint $b_!$, and by Definition~\ref{fl.def} and
Lemma~\ref{fun.le}, the composition $b_!  \circ a^*:\nu^{\hdot *}\C
\to \C_\ppt$ defines a functor
\begin{equation}\label{trunc.eq}
\unf(\C):\C \to \Unf(\C_\ppt,\I)
\end{equation}
cartesian over $\I$ and identified with $\Id$ over $\ppt$. This
proves existence: take $\gamma' = \Unf(\gamma) \circ \unf(\C)$. For
uniqueness, note that if we apply this construction to
$\Unf(\C',\I)$, then $b_! \circ a^*:\nu_\idot^*\Unf(\C',\I) \to \C'
\cong \Unf(\C',\I)_\ppt$ is the evaluation functor \eqref{ev.I.eq},
so that $\unf(\Unf(\C',\I)) \cong \Id$. Then for any $\gamma'':\C
\to \Unf(\C',\I)$ cartesian over $\I$, we have $\gamma'' \cong
\unf(\Unf(\C',\I)) \circ \gamma'' \cong \Unf(\gamma) \circ
\unf(\C) = \gamma'$.
\endproof

\begin{defn}
The {\em truncation} of a reflexive family of categories $\C$ over
$\I$ is the reflexive family $\Unf(\C_\ppt,\I)$, and the {\em
  truncation functor} for the family $\C$ is the functor
\eqref{trunc.eq}.
\end{defn}

For any co-reflexive family of categories $\C$ over $\I$, the
functor \eqref{trunc.eq} for $\iota^*\C$ induces a functor
\begin{equation}\label{trunc.perp}
\unf^\dg(\C):\C \to \iota^*\Unf(\C_\ppt,\iota(\I)) \cong \I \bb \C_\ppt,
\end{equation}
cartesian over $\I$. Alternatively, $\unf^\dg(\C) \cong
\unf(\C^o_\perp)^o_\perp$ is transpose-opposite to the cartesian
functor $\unf(\C^o_\perp)$ for the transpose-opposite reflexive
family $\C^o_\perp$, where we use \eqref{unf.perp} to identify the
target.

As an application of Proposition~\ref{trunc.prop}, one can define a
version of the unfolding construction of Subsection~\ref{unf.subs}
for reflexive families of categories. Namely, for any category $\E$
and a closed class of maps $v$ in $\E$, let $\Unf_v(\E,\I) \subset
\Unf(\E,\I)$ be the full subcategory spanned by functors $J^o \to
\E$ that factor through $\E_v \subset \E$. Consider the category
$\Rel$ of biordered set of Definition~\ref{rel.def}, let $\I^\dm =
U^{-1}(\I) \subset \Rel$ be the unfolding of $\I$, and note that
since $\I$ is closed under taking subsets, the forgetful functor
$U^l:\Rel \to \Pos$, $J \mapsto J^l$ induces a functor $U^l:\I^\dm
\to \I$. Then we define the {\em unfolding} $\C^\dm$ of a reflexive
family of categories $\C$ over $\I$ by the cartesian square
\begin{equation}\label{unf.sq}
\begin{CD}
  \C^\dm @>>> U^*\C\\
  @VVV @V{r^* \circ \unf(\C)}VV\\
  U^{l*}\Unf_{\Iso}(\C_\ppt,\I) @>>> U^{l*}\Unf(\C_\ppt,\I),
\end{CD}
\end{equation}
where $r:U^l \to U$ is the tautological map. Explicitly, for any
biordered $J \in \I^\dm$, $\C^\dm_J \subset \C_{U(J)}$ is the full
subcategory spanned by objects $c \in \C_{U(J)}$ whose truncation
$\unf(\C)(c):J^o \to \C_\ppt$ inverts all maps $j \leq^l j'$ in $J^l
\subset J$.

\begin{lemma}\label{unf.refl.le}
The unfolding $\C^\dm$ of a non-degenerate reflexive family of
categories $\C$ over $\I$ is constant along all reflexive maps in
$\I^\dm$.
\end{lemma}

\proof{} By the biordered version of Lemma~\ref{refl.sat.le}, it
suffices to prove that $\C^\dm$ is constant along the projection
$R([1]) \times J \to J$ for any $J \in \I^\dm$. Since $\C$ is
reflexive, $e^*:\C_{U(J)} \to \C_{U(J) \times [1]}$ is fully
faithful, and $c \in \C_{U(J) \times [1]}$ lies in its essential
image iff the adjunction map $e^*t^*c \to c$ is an isomorphism. But
since $\C$ is non-degenerate, the truncation functor
\eqref{trunc.eq} is conservative, so the latter can be checked after
restricting to $\{j\} \times [1] \in U(J) \times [1]$ for all $j \in
J$. We are then done by the definition of $\C^\dm$.
\endproof

\begin{lemma}\label{j.ind.le}
Assume given a cocartesian square $J_{01} \cong J_0 \copr_J J_1$ in
$\I$, and assume that it is perfect and universal with respect to a
map $J_{01} \to J'$, $J' \in \I$ in the sense of
Definition~\ref{pos.perf.def}. Moreover, assume given a reflexive
family $\C$ over $\I$ that is cartesian resp.\ semicartesian over
the square
\begin{equation}\label{j.ind.sq}
\begin{CD}
J \times_{J'} \wt{J} @>>> J_0 \times_{J'} \wt{J}\\
@VVV @VVV\\
J_1 \times_{J'} \wt{J} @>>> J_{01} \times_{J'} \wt{J}
\end{CD}
\end{equation}
for some map $f:\wt{J} \to J'$ in $\I$. Then for any biorder
$\wt{J}^\dm$ on $\wt{J}$ such that $f$ is a biordered map
$\wt{J}^\dm \to L(J')$, the unfolding $\C^\dm$ is cartesian
resp.\ semicartesian over the square \eqref{j.ind.sq} with the
biorder induced from $\wt{J}$, as in
Example~\ref{rel.perf.exa}. Moreover, the same holds for any perfect
square $J_{01} = J' \cong J_0 \copr_J J_1$ and bifibration
$\wt{J}^\dm \to L(J')$ in $\I^\dm$.
\end{lemma}

\proof{} By Lemma~\ref{car.rel.le} and \eqref{unf.sq}, it suffices
to check that both families $U^{l*}\Unf_{\Iso}(\C_\ppt,\I)$ and
$U^{l*}\Unf(\C_\ppt,\I)$ are cartesian along the square
\eqref{j.ind.sq}. But Lemma~\ref{pos.uni.le}, \eqref{j.ind.sq} with
$\wt{J}$ replaced by $\wt{J}^l$ is cocartesian in $\Cat$, so we are
done by \eqref{refl.fun}.
\endproof

\begin{lemma}\label{refl.sq.le}
Assume given a semicartesian square
\begin{equation}\label{refl.semi.sq}
\begin{CD}
\C_{01} @>{\gamma_0}>> \C_0\\
@V{\gamma_1}VV @VVV\\
\C_1 @>>> \C
\end{CD}
\end{equation}
of reflexive families over $I$ and functors cartesian over $I$. Then
the induced square of unfoldings $\C^\dm$, $\C_0^\dm$, $\C_1^\dm$,
$\C_{01}^\dm$ is also semicartesian.
\end{lemma}

\proof{} Since for any reflexive family $\C$, the top arrow in
\eqref{unf.sq} is a fully faithful embedding, the functor
$\C^\dm_{01} \to \C^\dm_0 \times_{\C^\dm} \C^\dm_1$ is trivially
full and conservative, and we just need to show that it is
essentially surjective. In other words, we need to check that for
any $J^\dm \in \I^\dm$, an object $c \in \C_{01U(J)}$ with
$\gamma_l(c) \in \C^\dm_{lJ} \subset \C_{lU(J)}$, $l=0,1$ lies in
$\C^\dm_{01J}$. But by definition, this can be checked after
restricting via each map $R([1]) \to J^\dm$, so we may assume right
away that $J^\dm = R([1])$. Since all unfoldings are constant along
the map $R([1]) \to \ppt$, the claim is then obvious.
\endproof

Explicitly, for any reflexive family $\C \to \I$, the truncation
functor \eqref{trunc.eq} provides a compatible system of functors
\begin{equation}\label{tr.J}
  \C_J \to J^o\C_\ppt
\end{equation}
for any $J \in \I$, and we also have the functors \eqref{nu.J.eq}
and the functors $\nu_f$ of \eqref{cyl.J.eq}. One can actually
define an even more general collection of functors that includes all
of these. It is convenient to start with a co-reflexive family $\C$
over $\I$. For any such family $\C \to \I$, the induced functor $\I
\bb \C \to \I \bb \I$ is a fibration; its fiber over an object in
$\I \bb \I$ represented by an $\I$-augmented set $\langle J,\alpha
\rangle = \langle J,\alpha_J \rangle$ is given by
\begin{equation}\label{unf.sec}
(\I \bb \C)_{\langle J,\alpha\rangle} \cong
\Sec(J,\alpha^*\C).
\end{equation}
Moreover, consider the tautological cofibration $\I_\idot \to \I$ of
\eqref{refl.fun}, and let $\I \bb^\bB \I \subset \I \bb \I$ be the
full subcategory spanned by $\langle J,\alpha \rangle$ such that
$\alpha^*\I_\idot \in \I$. Then as in Remark~\ref{pos.aug.rem}, by
the Grothendieck construction, sending $\langle J,\alpha \rangle$ to
$J_\idot = \alpha^*\I_\idot \to J$ identifies $\I \bb^\bB \I$ with
the subcategory $\Arc(\I) \subset \Ar(\I)$ of cofibrations $f:J_\idot \to
J$, with such maps from $J_0^\hdot \to J_0$ to $J_1^\hdot \to J$ that the
induced map $J_0^\hdot \to J_1^\hdot \times_{J_1} J_0$ is a cocartesian
functor over $J_0$. Under this identification, the fibration $\I
\bb^\bB \I \to \I$ of \eqref{fun.bb} becomes the projection
$\tau:\Arc(\I) \to \I$ that happens to be a fibration by
Example~\ref{comma.exa}. However, we also have the projection
$\sigma:\Arc(\I) \to \I$.

\begin{lemma}\label{C.refl.le}
For any co-reflexive family $\C \to \I$, $\sigma^*\C$ equipped with
the fibration
\begin{equation}\label{C.fib.dia}
\begin{CD}
  \sigma^*\C @>>> \Arc(\I) @>{\tau}>> \I
\end{CD}
\end{equation}
is a co-reflexive family over $\I$, proper if so is $\C$.
\end{lemma}

\proof{} For any cofibration $J_\idot \to J$ in $\I$, and any
left-reflexive full embedding $f:J' \to J$, with an adjoint
$f^\dg$, the full embedding $f_\idot:J'_\idot = f^*J_\idot \to
J_\idot$ is left-reflexive by Corollary~\ref{adm.cof.le}, with some
adjoint $f^\dg_\idot$. Then $\tau:\Arc(I) \subset \I \bb \I$ is a
co-reflexive family by Example~\ref{sub.refl.exa}, and
$\sigma^*\C$ is co-reflexive by the same argument as in
Lemma~\ref{refl.refl.J.le}. To see that it is proper as soon as so
is $\C$, note that for any right-closed embedding $f:J' \to J$,
$f_\idot:J'_\idot \to J_\idot$ is right-closed, and $\nu^\perp_f$
for $\sigma^*\C$ is then given by all the $\nu^\perp_{f_\idot}$ for
the family $\C$.
\endproof

Now, the truncation functor \eqref{trunc.perp} for the family
$\sigma^*\C$ of Lemma~\ref{C.refl.le} is a functor $\sigma^*\C \to
\I \bb \C$ over $\I \bb \I$, and if let $\I \bb^\bB \C \subset \I
\bb \C$ be the restriction of the fibration $\I \bb \C \to \I \bb
\I$ to $\I \bb^\bB \I \subset \I \bb \I$, then we obtain a functor
\begin{equation}\label{tr.big.co}
  \sigma^*\C \to \I \bb^\bB \C
\end{equation}
cartesian over $\Arc(\I) \cong \I \bb^\bB \I$. If $\C \to \I$ is a
reflexive family, then we can construct the cartesian functor
\eqref{tr.big.co} for the co-reflexive family $\C^o_\perp \to \I$,
and take the transpose-opposite functor over $\I \bb^\bB \I$ to
obtain a functor
\begin{equation}\label{tr.big}
  \sigma^*\C \to (\I \bb^\bB \C^o_\perp)^o_\perp.
\end{equation}
Explicitly, over any $\langle J,\alpha\rangle \in \I \bb \I$ with
$J_\idot = \alpha^*\I_\idot$, \eqref{tr.big} is a functor
\begin{equation}\label{tr.big.loc}
  \C_{J_\idot} \to \Sec(J^o,\alpha^{o*}\C_\perp),
\end{equation}
where $\C_\perp \to \I^o$ is the transpose cofibration, and we use
the identification \eqref{unf.sec}. If $J=\ppt$, then
\eqref{tr.big.loc} is \eqref{tr.J} for $J_\idot$, and if $J = [1]$
and $\alpha:[1] \to \I$ is the constant functor with some value
$J' \in \I$, then \eqref{tr.big.loc} is \eqref{nu.J.eq} for this
$J'$. More generally, if $J' = \Cyl(f)$ for some map $f$ in $\I$,
then \eqref{tr.big.loc} is $\nu_f$ of \eqref{cyl.J.eq}.

For any $\langle J,\alpha \rangle$, the target of the functor
\eqref{tr.big.loc} contains the full subcategory spanned by
cocartesian sections. To describe its preimage, one can use the
unfolding $\C^\dm \to \I^\dm$ of the reflexive family $\C$. Namely,
promote $\sigma:\Arc(\I) \to \I$ to a functor $\sigma_\dm:\Arc(\I)
\to \I^\dm$ by equipping $J_\idot$ with the biorder of
Example~\ref{cofib.rel.exa}, and denote the corresponding biordered
set by $J^\dm$. Then we have the fully faithful embedding
$\sigma_\dm^*:\C^\dm \to \sigma^*\C$, with fibers $\C^\dm_{J^\dm} \to
\C_{J^\hdot}$, and \eqref{tr.big.loc} fits into a cartesian square
\begin{equation}\label{tr.big.sq}
\begin{CD}
  \C^\dm_{J^\dm} @>>> \C_{J_\idot}\\
  @VVV @VVV\\
  \Sec_\Tot(J^o,\alpha^{o*}\C_\perp) @>>>
  \Sec(J^o,\alpha^{o*}\C_\perp)
\end{CD}
\end{equation}
with fully faithful horizontal functors.

By definition, the co-reflexive family $\sigma^*\C$ of
Lemma~\ref{C.refl.le} comes equip\-ped with a fibration $\sigma^*\C
\to \Arc(\I)$, and restricting it to sections $\I \to \Arc(\I)$ of
the fibration $\Arc(\I) \to \I$ produces other co-reflexive families
over $\I$. Here is one useful application of this construction.

\begin{corr}\label{refl.ev.corr}
For any co-reflexive family $\C \to I$ and $J \in \I$, there exists
a co-reflexive family $\C^{hJ} \to \I$ and a functor $\ev:\I \bb J
\times_{\I} \C^{hJ} \to \C$, cartesian over $\I$, such that for any
co-reflexive family $\C' \to \I$, any functor $\gamma:\I \bb J
\times_{\I} \C' \to \C$ cartesian over $\I$ factors as
\begin{equation}\label{refl.ev.dia}
\begin{CD}
\I \bb J \times_{\I} \C' @>{\id \times \gamma'}>> \I \bb J
\times_{\I} \C^{hJ} @>{\ev}>> \C
\end{CD}
\end{equation}
for some $\gamma':\C' \to \C^{hJ}$ cartesian over $\I$, and
$\gamma'$ is unique up to a unique isomorphism.
\end{corr}

\proof{} To construct the evaluation functor $\ev$, we first
consider the whole family $\sigma^*C \to \Arc(\I)$. We have a
functor $\Arc(\I) \to \Pos$ sending a cofibration $f:J_\idot \to J$
to the set of sections $\Sec(J,J_\idot) \in \Pos$; let
$\nu_\wr:\Arc(\I)_\idot \to \Arc(\I)$ be the corresponding
cofibration. Moreover, for any section $s:J \to J_\idot$ of a
cofibration $f:J_\idot \to J$, let $J(s) \subset J_\idot$ be the
subset of all elements $j \in J_\idot$ such that $j \leq s(f(j))$
(alternatively, $J(s) = J /^J_s J_\idot$ is the relative left
comma-set of the map $s$). Then sending $s$ to $J(s)$ defines a
functor $\nu_{\wr>}:\Arc(\I)_\idot to \Pos$, and as in
\eqref{unf.dia}, we have a map $a:\nu_{\wr>} \to \sigma \circ
\nu_\wr$, and another map $b:\nu_{\wr>} \to \tau \circ \nu_\wr$
induced by the embeddings $J(s) \subset J_\idot$ and the projections
$J(s) \to J$. These projections are right-reflexive, so again as in
\eqref{unf.dia}, we have a diagram
\begin{equation}\label{unf.wr.dia}
\begin{CD}
  \nu_\wr^*\sigma^*\C @>{a^*}>> \nu_{\wr >}^*\C @<{b^*}<< \tau^*\C,
\end{CD}
\end{equation}
and $b^*$ in \eqref{unf.wr.dia} has a left-adjoint functor
$b_!$. Now consider the functor $\eps_J:\I \to \Arc(\I)$ sending $J'
\in \I$ to the projection $J \times J' \to J'$, and let $\C^{hJ} =
\eps_J^*\sigma^*\C$. Then since $\eps_J(J' \times [1]) \cong
\eps_J(J') \times [1]$ for any $J' \in \I$, the family $\C^{hJ}$ is
co-reflexive by virtue of Lemma~\ref{C.refl.le}. Moreover,
$\eps_J^*\nu_\wr^*\sigma^*\C \cong \I \bb J \times_{\I} \C^{hJ}$ and
$\eps_J^*\tau^*\C \cong \C$, so that $\ev = b_! \circ a^*$ induces
the required evaluation functor $\ev:\I \bb J \times_{\I} \C^{hJ}
\to \C$. To check the universal property, note that the map $\sigma
\circ \eps_J \to \tau \circ \eps_J \cong \id$ induces a functor $\C
\to \C^{hJ}$ cartesian over $\I$, and we have $\I \bb J)^{hJ} \cong
\I \bb \Fun(J,J)$, so the embedding $\eps(\id):\ppt \to \Fun(J,J)$
induces a functor $\I \to (\I \bb J)^{hJ}$, again cartesian over
$\I$. The product of these two functors provides a co-evaluation
functor
$$
\ev^\dg:\C \to (I \bb J \times_{\I} \C)^{hJ},
$$
again cartesian over $\I$, and as in Subsection~\ref{cat.fun.subs}, both
compositions
$$
\begin{CD}
\I \bb J \times_{\I} \C @>{\id \times \ev^\dg}>> \I \bb J
\times_{\I} (I \bb J \times_{\I} \C)^{hJ} @>{\ev}>> I \bb J
\times_{\I} \C,\\
\C^{hJ} @>{\ev^\dg}>> (I \bb J \times_{\I} \C^{hJ})^{hJ} @>{\ev^{hJ}}>>
\C^{hJ}
\end{CD}
$$
are canonically isomorphic to the corresponding identity functors.
\endproof

In particular, in the situation of Corollary~\ref{refl.ev.corr}, we
have an identification $\Fun^\Tot_{\I}(\I \bb J,\C) \cong
\C^{hJ}_\ppt \cong \C_J$. Explicitly, we have the tautological
object $j \in (\I \bb J)_J$ represented by $\langle J,\id \rangle$,
and evaluation at $j$ provides a functor
\begin{equation}\label{ev.j.I}
  \Fun^\Tot_{\I}(\I \bb J,\C) \to \C_J.
\end{equation}
By Corollary~\ref{refl.ev.corr}, this functor is an equivalence (and
in particular, the category $\Fun^\Tot_{\I}(\I \bb J,\C)$ is
well-defined even if $\I$ is not small)

\subsection{Augmentations.}

Let us now define a version of reflexive families $\C \to \I$
augmented by a partially ordered set $I$, in the following sense.

\begin{defn}\label{refl.aug.def}
For any ample $\I \subset \Pos$ and partially ordered set $I$, an
{\em $I$-augmented reflexive family of categories} $\C$ over $\I$ is
a reflexive family of categories $\C$ over $\I$ equipped with
functor $\pi:\C \to \Unf(I,\I)$ over $\I$ that is a fibration.
\end{defn}

\begin{exa}\label{unf.E.J.exa}
For any fibration $\pi:\E \to I$, $\Unf(\pi):\Unf(\E,\I) \to
\Unf(I,\I)$ is a fibration, so that $\Unf(\E,\I)$ of
Example~\ref{unf.E.exa} is an $I$-augmented reflexive family in the
sense of Definition~\ref{refl.aug.def}, with fibers
\begin{equation}\label{unf.J.fib}
  \Unf(\E,\I)_{\langle J,\alpha\rangle} \cong \Sec(J^o,\alpha^{o*}\E)
\end{equation}
for any $\langle J,\alpha \rangle \in \Unf(I,\I)$.
\end{exa}

In particular, for any $I$-augmented reflexive family $\C$, we have
a fibration $\pi_{\ppt}:\C_{\ppt} \to I$, and then \eqref{trunc.eq}
is a cartesian functor over $\Unf(I,\I)$. Let $\flat$ be the class
of maps in $\C_\ppt$ cartesian over $I$, and define the {\em
  relative unfolding} $(\C|I)^\dm$ by upgrading \eqref{unf.sq} to a
cartesian square
\begin{equation}\label{unf.J.sq}
\begin{CD}
  (\C|I)_\dm @>>> U^*\C\\
  @VVV @V{r^* \circ \unf(\C)}VV\\
  U^{l*}\Unf_\flat(\C_\ppt,\I)  @>>> U^{l*}\Unf(\C_\ppt,\I),
\end{CD}
\end{equation}
where $\Unf_\flat$ has the same meaning as in \eqref{unf.sq}. Then
\eqref{unf.J.sq} is a cartesian square of fibrations and cartesian
functors over $\Unf^\dm(I,\I) = U^*\Unf(I,\I)$. The embedding $L:\I
\to \I^\dm$ of \eqref{rel.pos} tautologically extends to an
embedding
\begin{equation}\label{L.unf}
  L:\Unf(I,\I) \to \Unf^\dm(I,\I),
\end{equation}
and we trivally have $L^*(\C|I)^\dm \cong \C$. Explicitly, just as
for $\Unf(I,\I)$, objects in $\Unf^\dm(I,\I)$ are pairs $\langle
J,\alpha \rangle$ of a biordered set $J \in \I^\dm$ and an
augmentation map $\alpha:J \to I^o$, while morphisms from $\langle
J,\alpha \rangle$ to $\langle J',\alpha'\rangle$ are biordered maps
$f:J \to J'$ such that that $\alpha' \circ f \leq \alpha$. Say
that such a map $f$ is {\em $I$-strict} if $\alpha' \circ f =
\alpha$; then $I$-strict maps form a closed class, and the dense
embedding $b$ of \eqref{pos.tri.dia} induces an embedding
\begin{equation}\label{b.unf}
  b:\I^\dm / R(I^o) \to \Unf^\dm(I,\I)
\end{equation}
identifying $I^\dm / R(I^o)$ with the corresponding dense
subcategory in $\Unf^\dm(I,\I)$. As in Lemma~\ref{refl.refl.J.le},
if an $I$-strict map $f$ is right-reflexive in the sense of
Definition~\ref{rel.refl.def}, then the adjoint map $f_\dg:J' \to J$
is also a morphism in $\Unf^\dm(I,\I)$ (although in general not
$I$-strict). We then have the following augmented version of
Lemma~\ref{unf.refl.le}.

\begin{lemma}\label{unf.refl.J.le}
For any non-degenerate $I$-augmented reflexive family $\C$ over $\I$
and pair of maps $f$, $f_\dg$ in $\Unf^\dm(I,\I)$ such that $f$ is
$I$-strict, $f_\dg \circ f = \id$ and $f \circ f_\dg \leq \id$, the
isomorphism $f_\dg^* \circ f_* \cong \id$ defines an adjunction
between the corresponding transition functors $f_*$ and $f_\dg^*$ of
the relative unfolding $(\C|I)^\dm$ of \eqref{unf.J.sq}. Moreover,
if $f \circ f_\dg \leq^l \id$, so that $f$ is right-reflexive in the
sense of Definition~\ref{rel.refl.def}, then $f^*$ and $f_\dg^*$ are
adjoint equivalences.
\end{lemma}

\proof{}\.The first claim immediately follows from
Lemma~\ref{refl.refl.J.le} and \eqref{unf.J.sq}. For the second,
since $f_\dg \circ f = \id$, it suffices to check that for any
object $c$ in the full subcategory $(\C|I)^\dm_{\langle
  J,\alpha\rangle} \subset \C_{\langle U(J),\alpha\rangle}$, the
adjunction map $f_\dg^*(f^*(c)) \to c$ provided by
Lemma~\ref{refl.J.le} is an isomorphism. Since $\C$ is
non-degenerate, this can be checked after applying the truncation
functor $\unf(\C)$. Thus we may assume right away that
$\C=\Unf(\E,\I)$ for some fibration $\pi:\E \to I$, as in
Example~\ref{unf.E.J.exa}. But then in terms of \eqref{unf.J.fib},
objects in $(\C|I)^\dm$ correspond to sections cartesian over
$J^{lo} \subset J^o$, so we are done by Lemma~\ref{2adm.le}.
\endproof

Now assume that $I \in \iota(\I)$. Then since $\I$, hence also
$\iota(\I)$ is ample, $I \setminus_f I' \in \iota(\I)$ for any $I'
\in \iota(\I)$ equipped with a map $f:I' \to I$, so that the
extended Yoneda embedding of \eqref{yo.cat.eq} induces a fully
faithful embedding $\Y:\Unf(I,\I) = \I \bb^\iota J \to
I^o\I$. Moreover, if we interpret $I^o\I$ as the fiber $(\I \bb
\I)_{I^o}$ of the fibration $\I \bb \I \to \I$ over $I^o$, then $\Y$
takes values in the full subcategory $(\I \bb^\flat \I)_{I^o}
\subset (\I \bb \I)_{I^o}$, so that altogether, we obtain a fully
faithful embedding
\begin{equation}\label{J.yo.eq}
  \Y:\Unf(I,\I) \to (\I \bb^\flat \I)_{I^o} \cong \Arc(\I)_{I^o}
  = \I \mc I^o.
\end{equation}
Equivalently, \eqref{J.yo.eq} is the embedding $\Y$ of
\eqref{pos.tri.dia} composed with the equivalence $\iota:\iota(\I)
\mm I \cong \I \mc I^o$; explicitly, it sends $\langle J,\alpha
\rangle \in \Unf(I,\I)$ to $\tau:J /_\alpha I^o \to I^o$. The
projection $\sigma_\dm:\Arc(\I) \to \I^\dm$ then restricts to a
functor
\begin{equation}\label{sigma.dm.J}
\sigma_\dm:\Arc(\I)_{I^o} = \I \mc I^o \to \I^\dm / R(I^o),
\end{equation}
and this is exactly the fully faithful embedding
\eqref{cofib.rel.eq} of Example~\ref{cofib.rel.exa}. Moreover,
$\eta:J \to J /_\alpha I^o$ defines a functorial map $\eta:L
\to b \circ \sigma_\dm \circ \Y$, where $L$ resp.\ $b$ are as in
\eqref{L.unf} resp.\ \eqref{b.unf}, and this provides a functor
\begin{equation}\label{J.yo.dm}
\eta^*:\Y^*\sigma_\dm^*b^*(\C|I)^\dm \to L^*(\C|I)^\dm \cong \C
\end{equation}
over $\Unf(I,\I)$. The map $\eta$ is pointwise an $I$-strict
right-reflexive full embedding, so by Lemma~\ref{unf.refl.J.le}, for
a non-degenerate $\C$, \eqref{J.yo.dm} is an equivalence.

\begin{prop}\label{Y.prop}
For any $I \in \iota(\I)$ and non-degenerate $I$-augmented reflexive
family $\C$ over $\I$, the functor
\begin{equation}\label{L.yo.dm}
  \sigma_\dm^*b^*(\C|I)^\dm \to \Y_*\C
\end{equation}
adjoint to the equivalence \eqref{J.yo.dm} is itself an equivalence.
\end{prop}

\proof{} Note that the functors $L$, $\Y$ and $\sigma_\dm$ fit into
a commutative diagram
\begin{equation}\label{L.Y.sq}
\begin{CD}
  \Unf(I,\I) @>{\Y}>> \I \mc I^o\\
  @V{L}VV @VV{\sigma_\dm}V\\
  \Unf^\dm(I,\I) @>{\Y^\dm}>> \I^\dm / R(I^o),
\end{CD}
\end{equation}
where $\Y^\dm$ sends $\langle J,\alpha \rangle$ to $J /^\dm_\alpha
R(I^o)$ of Example~\ref{bifib.lc.exa}. Moreover, by adjunction
between $b$ and $a \circ \Y$ of \eqref{pos.tri.dia}, $\Y^\dm$ is
actually left-adjoint to the embedding \eqref{b.unf}, so that
$\Y^\dm_* \cong b^*$, and \eqref{L.yo.dm} is the base change functor
induced by the commutative square \eqref{L.Y.sq}. We need to check
that this is an equivalence. Were we to replace $(\C|I)^\dm$ with
$U^*\C \cong L_*\C$, the statement would hold tautologically --- we
have $\sigma_\dm^* \circ \Y^\dm_* \circ L_* \cong \sigma_\dm^* \circ
\sigma_{\dm*} \circ \Y_*$, and $\sigma_\dm^* \circ \sigma_{\dm*}
\cong \id$, $L^* \circ L_* \cong \id$ since $\sigma_\dm$ and $L$ are
fully faithful. As it happens, we have the fully faithful embedding
$(\C|I)^\dm \to U^*\C$ of \eqref{unf.J.sq}, and we need to check
that it becomes an equivalence once we apply $\sigma_\dm^* \circ
\Y^\dm_*$. It stays fully faithful, so we only need to check that it
becomes essentially surjective. By \eqref{2kan.eq}, this amounts to
checking that for any cofibration $\pi:J \to I^o$ representing an
object $\langle J,\pi \rangle \in \I \mc I^o$, any global section
\begin{equation}\label{c.J.eq}
c \in \Sec^\Tot(\Unf^\dm(I,\I) /_{\Y^\dm} \sigma_\dm(\langle
J,\pi\rangle),U^*\C)
\end{equation}
actually lies in $(\C|I)^\dm \subset U^*\C$ after evaluation at each
object in the left comma-fiber $\Unf^\dm(I,\I) /_{\Y^\dm}
\sigma_\dm(\langle J,\pi \rangle)$. Explicitly, such an object is
given by a triple $\langle J',\alpha',f\rangle$ of a biordered set
$J' \in \I^\dm$, a map $\alpha':J' \to I^o$, and a map $f:J' \to J$
such that $\pi \circ f = \alpha'$ and the order relation $f(j) \leq
f(j')$ is cocartesian over $I^o$ for any $j \leq^l j'$ in
$J'$. Since $\alpha' = \pi \circ f$, we might as well forget it and
only remember $J'$ and $f$. The conditions on the evaluation
$c_{J',f} \in \C_{U(J')}$ of the section \eqref{c.J.eq} that we need
to check are separate for each order relation $j \leq^l j'$ in $J'$,
and to check each of them, we may restrict to the image of the
corresponding map $R([1]) \to J'$. In other words, we may assume
right away that $J'=R([1])$. But then we can consider the
right-reflexive embedding $s:\ppt \to R([1])$ onto $0 \in [1]$,
$\langle \ppt,f \circ s\rangle$ is another object in the
comma-fiber, and the adjoint map $s_\dg = e:R([1]) \to \ppt$ defines
a map in the comma-fiber. Since the section $c$ in \eqref{c.J.eq} is
cartesian along all such maps, $c_{R([1]),f}$ must lie in the
essential image of the embedding $e^*$, and then we are done by
Lemma~\ref{unf.refl.J.le}.
\endproof

While Proposition~\ref{Y.prop} is quite technical, it has a very
useful corollary. Namely, any map $g:I' \to I$ between partially
ordered set induces a tautological functor $\Unf(g):\Unf(I',\I) \to
\Unf(I,\I)$, and by Lemma~\ref{refl.pb.le}, for any $I$-augmented
reflexive family $\C$, $\C' = \Unf(g)^*\C$ is an $I'$-augmented
reflexive family, proper and/or non-degenerate if so is $\C$. As it
happens, for maps in $\iota(\I)$, we can also go the other way.

\begin{corr}\label{pf.corr}
For any map $g:I' \to I$ in $\iota(\I)$, and any $I'$-augmented
reflexive family $\C$, $\C' = \Unf(f)_*\C$ is an $I$-augmented
reflexive family, proper and/or non-degenerate if so were
$\C$. Moreover, if $g$ is a cofibration, then for any map $f:I_0 \to
I$, with the induced cofibration $g_0:I'_0 = g^*I' \to I_0$ and map
$f':I_0' \to I$, the base change functor $a:\Unf(f)^*\Unf(g)_*\C \to
\Unf(g_0)_*\Unf(f')^*\C$ is an equivalence.
\end{corr}

\proof{} While $\Unf(g)$ itself does not admit a right-adjoint, it
acquires a ``relative adjoint'' in the sense of
Definition~\ref{adj.rel.def} once we compose it with the full
embedding \eqref{J.yo.eq}. Namely, let $\Unf(g)_\dg:\I \mc I^o \to
\I \mc {I'}^o$ be a functor sending some $\pi:J \to I^o$ to $g^{o*}J
\to {I'}^o$, and note that for any $\langle J',\alpha' \rangle \in
\Unf(I',\I)$ and $\langle J,\pi \rangle \in \I \mc I^o$, we have a
functorial identification
$$
\Hom(\Y \circ \Unf(g)(\langle J',\alpha' \rangle),\langle J,\pi
\rangle) \cong \Hom(\Y(\langle J',\alpha'
\rangle),\Unf(g)_\dg(\langle J,\pi \rangle)).
$$
Being functorial, this identification provides a cartesian square
\begin{equation}\label{Y.g.sq}
\begin{CD}
  \Unf(I',\I) /_{\Y \circ \Unf(g)} (\I \mc I^o) @>{\tau}>> \I \mc I^o\\
  @V{\id / \Unf(g)_\dg}VV @VV{\Unf(g)_\dg}V\\
  \Unf(I',\I) /_{\Y} (\I \mc {I'}^o) @>{\tau}>> \I \mc {I'}^o.
\end{CD}
\end{equation}
Since $\tau$ is a cofibration, \eqref{Y.g.sq} provides a base change
equivalence, and then if we compute $\Y_*\C'$ and $\Y_*\C = (\Y
\circ \Unf(g))_*\C'$ by the left-comma-categories version of
\eqref{rc.dec}, and note that $\sigma \circ (\id / \Unf(g)_\dg)
\cong \sigma$, we obtain an equivalence
$$
\Y_*\C \cong \tau_*\sigma^*\C' \cong \tau_*(\id \times
\Unf(g)_\dg)^*\sigma^*\C' \cong \Unf(g)_\dg^*\tau_*\sigma^*\C' \cong
\Unf(g)_\dg^*\Y_*\C'
$$
that restricts to an equivalence
\begin{equation}\label{Unf.g.C}
\C \cong \Y^*\Y_*\C \cong \Y^*\Unf(g)_\dg^*\Y_*\C'
\end{equation}
Then for any $I$-strict right-reflexive full embedding $f$ in
$\Unf(I,\I)$, with the adjoint $f_\dg$, $f \circ f_\dg \leq \id$
implies that $\Unf(g)_\dg(\Y(f)) \circ \Unf(g)_\dg(\Y(f_\dg)) \leq
\id$, so $\C$ is a reflexive family by Lemma~\ref{refl.refl.J.le},
Proposition~\ref{Y.prop} and Lemma~\ref{unf.refl.J.le}. It is
obviously non-degenerate if so is $\C'$, and since both $\Y$ and
$\Unf(g)_\dg$ commute with taking dual cylinders, and $\Unf(g)_\dg
\circ \Y$ sends right-closed full embeddings to right-closed full
embeddings, $\C$ is proper as soon as so is $\C'$.

Finally, for the last claim, since $a$ is cartesian over
$\Unf(I,\I)$, it suffices to check that it becomes an equivalence
after restricting with respect to the dense embedding
$\lambda:\oUnf(I,\I) = \iota(\I) / I \to \Unf(I,\I)$. But then the
induced functor $\oUnf(g):\oUnf(I',\I) \to \oUnf(I,\I)$ has a
genuine right-adjoint $\oUnf(g)_\dg$ sending $\langle J,\alpha
\rangle \in \oUnf(I,\I)$ to ${I'}^o \times_{I^o} J$. Moreover, the
map $\eta$ of \eqref{J.yo.dm} provides a functorial map $L \circ
\oUnf(g)_\dg \to \sigma_\dm \circ \Unf(g)_\dg \circ \Y \circ
\lambda$, and if $g$ is a cofibration, this map is right-reflexive,
so that \eqref{Unf.g.C} provides an equivalence
\begin{equation}\label{oUnf.g}
\lambda^*\Unf(g)_*\C' \cong \oUnf(g)_\dg^*\lambda^*\C',
\end{equation}
and similarly for $g_0$. The base change map $\oUnf(f') \circ
\oUnf(g_0)_\dg \to \oUnf(g)_\dg \circ \oUnf(f)$ is an isomorphism,
so this proves the claim.
\endproof

\subsection{Complete families.}

One way to insure that a family of categories $\C \to \I$ is
reflexive in the sense of Definition~\ref{cat.refl.def} is to
require that the transition functor $f^*$ admits a well-behaved
right-adjoint $f_*$ for {\em any} map $f$ in $\I$. The formal
definition is as follows.

\begin{defn}\label{compl.fam.def}
A reflexive family of categories $\C \to \I$ over an ample $\I
\subset \Pos$ is {\em complete} if for any map $f:J' \to J$ in $\I$,
the transition functor $f^*:\C_J \to \C_{J'}$ admits a right-adjoint
$f_*:\C_{J'} \to \C_J$, and for any commutative square
\begin{equation}\label{posf.sq}
\begin{CD}
  J'_0 @>{f'}>> J'_1\\
  @V{g_0}VV @VV{g_1}V\\
  J_0 @>{f}>> J_1
\end{CD}
\end{equation}
in $\I$ such that \thetag{i} $f$ is a cofibration and the square is
cartesian, or \thetag{ii} $J_1=J_1'=J_0'=\ppt$, $J_0=[1]$, and
$g_0=t:\ppt = [0] \to [1]$, the base change map \eqref{bc.eq}
provides an isomorphism
\begin{equation}\label{refl.bc}
  g_1^* \circ f_* \cong f'_* \circ g_0^*
\end{equation}
between the corresponding functors $\C_{J_0} \to \C_{J_1'}$.
\end{defn}

\begin{exa}
For any complete category $\E$, the family $\Unf(\E,\I)$ is complete
in the sense of Definition~\ref{compl.fam.def}, with $f_*$ given by
right Kan extensions, and \eqref{refl.bc} being an isomorphism by
(the dual version of) Lemma~\ref{bc.le}.
\end{exa}

\begin{remark}
There is an obvious clash of terminology between
Definition~\ref{compl.fam.def} and Definition~\ref{pos.compl.def},
but it is not clear what can be done about it since both usages are
very well-established. As a mitigating circumstance, we note that
the only family of groupoids over $\I$ complete in the sense of
Definition~\ref{compl.fam.def} is the trivial family $\I \to \I$.
\end{remark}

For any complete family $\C \to \I$ that is reflexive,
Definition~\ref{compl.fam.def}~\thetag{ii} is automatic, and
conversely, we have the following result.

\begin{lemma}\label{compl.refl.le}
A non-degenerate complete family of categories $\C$ over $\I$ is
reflexive.
\end{lemma}

\proof{} By Corollary~\ref{refl.refl.corr}, it suffices to shows
that for any $J \in \I$, with the projection $e:J \times [1] \to J$
and its section $t:J \to J \times [1]$, the map $t^* \to e_*$
adjoint to the isomorphism $t^* \circ e^* \cong \id$ is itself an
isomorphism. Since $\C$ is non-degenerate, this can be checked after
taking evaluations at all elements $j \in J$, and by
Definition~\ref{compl.fam.def}~\thetag{i}, this reduces to the case
$J = \ppt$, where the claim is
Definition~\ref{compl.fam.def}~\thetag{ii}.
\endproof

For any complete reflexive family $\C \to \I$ and $J \in
\I$, we will denote by $\lim^h_J = f_*$ the adjoint functor
corresponding to the projection $f:J \to \ppt$. By
Lemma~\ref{refl.J.le}, for any right-reflexive resp.\ left-reflexive
full embedding $f:J' \to J$ with adjoint map $f_\dg:J' \to J$, we
have $f_* \cong f_\dg^*$ resp $f_{\dg*} \cong f^*$. In particular,
in the left-reflexive case, we have a natural isomorphism
\begin{equation}\label{adm.refl}
\lim^h_J \cong \lim^h_{J'} \circ f_{\dg*} \cong \lim^h_{J'} \circ
f^*,
\end{equation}
a reflexive family counterpart of (the dual version of)
\eqref{adm.eq}, and for any map $g:J' \to J$ with decomposition
\eqref{lc.pos}, we have
\begin{equation}\label{lc.refl}
g_* \cong \tau_* \circ \eta_* \cong \tau_* \circ \sigma^*,
\end{equation}
so that \eqref{refl.bc} provides an identification
\begin{equation}\label{kan.refl}
\ev_j \circ g_* \cong \lim^h_{J/_g j} \circ \sigma(j)^*, \qquad j
\in J,
\end{equation}
a counterpart of (the dual version of) \eqref{kan.eq}. We note that
since $\I$ is ample, decompositions \eqref{lc.pos} for maps in $\I$
lie in $\I$.

\begin{lemma}\label{bc.J.le}
Assume given a complete family of categories $\C$ over $\I$. Then
$f_*$ is fully faithful for any fully faithful map $f:J \to J'$ in
$\I$. Moreover, in the situation of Lemma~\ref{bc.le}, assume that
$I$, $I_0$, $I_1$, $I'$ are partially ordered sets in $\I$. Then the
base change map $\nu_0^* \circ \gamma_* \to \gamma'_* \circ \nu_1^*$
is an isomorphism.
\end{lemma}

\proof{} For the first claim, note that \eqref{lc.refl} and
\eqref{refl.bc} provide an isomorphism $f^* \circ f_* \cong \tau_*
\circ \sigma^*$, where $\sigma,\tau:\Ar(J) \to J$ are standard
projections, and we have $\tau_* \circ \sigma^* \cong \tau_* \circ
\eta_* \cong \id_* \cong \id$. For the second claim, use the same
argument as in Lemma~\ref{bc.le} to reduce to the isomorphism
\eqref{refl.bc}.
\endproof

\begin{defn}\label{cat.semi.def}
A family of categories $\C$ over an ample $\I$ is {\em additive} if
for any $J \in \I$ equipped with a map $J \to S$ to a discrete set
$S$, the functor
$$
\C_J \to \prod_{s \in S}\C_{J_s}
$$
is an equivalence. A family $\C$ over $\I$ is {\em semiexact} if for
any $J \in \I$ with a map $J \to \V = \{0,1\}^<$, the functor
\begin{equation}\label{cat.semi.eq}
\C_J \to \C_{J/0} \times_{\C_{J_o}} \C_{J/1}
\end{equation}
is an epivalence.
\end{defn}

\begin{lemma}\label{refl.semi.le}
A non-degenerate additive complete proper reflexive family $\C \to
\I$ is semiexact.
\end{lemma}

\proof{} Note that for any $J \in \I$, $J \times \V \in \I$ is the
dual cylinder of the codiagonal map $f:J \copr J \to J$, and since
$\C$ is additive, we have $\C_{J \copr J} \cong \C_J \times \C_J$
and $\C_{J \copr J} \setminus_{f^*} \C_J \cong \V^o\C_J$. Assume
given a map $J \to \V$, and consider the full embedding $\eps =
\sigma \times \tau:J / \V \to J \times \V$. Then $J \times \V \to
\V$ and $\tau: J/\V \to \V$ are cofibrations, thus correspond to
some functors $X,Y:\V \to \I$ under the identification $\I \mc \V
\cong \Fun(\V,\I)$, and $\eps$ is cocartesian over $\V$, thus
defines a map $\eps^\hdot:X \to Y$. Then the functors
\eqref{tr.big.loc} for $X$ and $Y$ fit into a commutative diagram
\begin{equation}\label{semi.dia}
\begin{CD}
  \C^\dm_{J/^\dm R(\V)} @>{\eta_*}>> \C_{J/\V} @>{\eps_*}>> \C_{J
    \times \V}\\
  @V{a}VV @VVV @VV{b}V\\
  \Sec_\Tot(\V^o,X^{o*}\C_\perp) @>>> \Sec(\V^o,X^{o*}\C_\perp)
  @>{\eps^\hdot_*}>> \Sec(\V^o,Y^{o*}\C_\perp),
\end{CD}
\end{equation}
where the square on the left is the cartesian square
\eqref{tr.big.sq}. The rows in \eqref{semi.dia} are fully faithful
by Lemma~\ref{bc.J.le}, and $\Sec(\V^o,Y^{o*}\C_\perp) \cong
\V^o\C_J$, while $b$ is $\nu^\perp_f$ for the codiagonal map $f:J
\copr J \to J$. Since $\C$ is proper, $b$ is an epivalence, and then
the square on the right is cartesian by
Lemma~\ref{epi.sq.le}. Therefore $a$ is also an epivalence, and
since $\Sec_\Tot(\V^o,X^{o*}\C_\perp) \cong \Sec^\Tot(\V,X^*\C)$ is
exactly the target of the functor \eqref{cat.semi.eq}, it remains to
show that $\sigma^*$ induces an equivalence $\C_J \cong \C^\dm_{J
  /^\dm R(\V)}$. But the fully faithful embedding $\eta:L(J) \to J
/^\dm R(\V)$ is right-reflexive in $\I^\dm$, so $\eta_* \cong
\sigma^*$ is an equivalence by Lemma~\ref{unf.refl.le}.
\endproof

It is useful to note that for complete families, \eqref{cyl.J.eq}
can be considerably generalized: we have a version of a gluing
construction for left-closed embeddings. Namely, assume given a
complete family $\C$ over $\I$ and a left-closed embedding
$l:J_0 \to J$ in $\I$. Let $r:J_1 \to J$ be the right-closed
embedding of the complement $J_1 = J \ssetminus J_0 \subset J$, and
denote $\theta = r^* \circ l_*:\C_{J_0} \to \C_{J_1}$. Then for any
$c \in \C_J$, we have a functorial map $\alpha(c):r^*c \to
\theta(l^*c) = r^*l_*l^*c$ induced by the adjunction map $c \to
l_*l^*c$, and this defines a functor
\begin{equation}\label{glue.eq}
\C_J \to \C_{J_1} \setminus_\theta \C_{J_0}, \qquad c \mapsto
\langle r^*c,l^*c,\alpha(c) \rangle.
\end{equation}
If $J = \Cyl^o(f)$ is the dual cylinder of a map $f:J_1 \to J_0$,
and $l:J_0 \to \Cyl^o(f)$ is the embedding $s$ of
\eqref{cyl.o.deco}, then as soon as $\C$ is non-degenerate, it is
reflexive by Lemma~\ref{compl.refl.le}, and \eqref{glue.eq} is the
functor $\nu^\perp_f$ of \eqref{cyl.J.eq}.

\begin{defn}\label{refl.tight.def}
A complete family of categories $\C \to \I$ is {\em tight} if
\eqref{glue.eq} is an epivalence for any left-closed embedding $J_0
\subset J$.
\end{defn}

\begin{defn}\label{refl.cont.def}
A family of categories $\C \to \I$ is {\em weakly semicontinuous} if
for any $J \in \I$ equipped with a map $J \to \N$, with the
corresponding functor $J_\idot:\N \to \Posf^+$, $n \mapsto J/n$, the
functor
\begin{equation}\label{refl.cont.eq}
\nu(\C,J_\idot):\C_J \to \Sec(\N,J_\idot^*\C)
\end{equation}
is an epivalence.
\end{defn}

\begin{lemma}\label{tight.le}
An additive non-degenerate complete family $\C \to \I$ is tight if
and only if it is proper. If $\I \subset \Posf^+$, then an additive
tight complete family $\C \to \I$ is non-degenerate over $\I \cap
\Posf \subset \I$, and if it is weakly semicontinuous, it is
non-degenerate.
\end{lemma}

\proof{} The second claim is clear: $\C_J \to J^o\C_\ppt$ of
\eqref{tr.J} is conservative if $\dim J = 0$ by additivity, then if
\eqref{glue.eq} is an epivalence, it is conservative for any $J \in
\I \cap \Posf$ by induction on the skeleton filtration \eqref{sk.J},
and then for any $J \in \I$, we can apply \eqref{refl.cont.eq} for
the height function $\hht:J \to \N$.

For the first claim, one direction is tautological --- $\C$ is
reflexive by Lemma~\ref{compl.refl.le}, and since $\nu^\perp_f$ is
then a particular case of \eqref{glue.eq}, being proper is part of
being tight --- so assume given an additive $\C$ non-degenerate
complete proper reflexive family $\C \to \I$, and a left-closed
embedding $l:J_0 \to J$ in $\I$, with the complement $r:J_1 \to
J$. Let $J \to [1]$ be its characteristic map, and let $J' = [1]
\setminus J$, with the fully faithful embedding $\eta:J \to
J'$. Alternatively, we have $J' = \Cyl^o(r)$, and the functor
\eqref{glue.eq} fits into a commutative square
\begin{equation}\label{tgh.sq}
\begin{CD}
  \C_J @>{\eta_*}>> \C_{J'}\\
  @VVV @VV{\nu^\perp_r}V\\
  \C_{J_1} \setminus_\theta \C_{J_0} @>{\id \setminus l_*}>>
  \C_{J_1} \setminus_{r^*} \C_J
\end{CD}
\end{equation}
with fully faithful rows. Since $\C$ is proper, $\nu^\perp_r$ in
\eqref{tgh.sq} is an epivalence, and we are done by
Lemma~\ref{epi.sq.le}.
\endproof

\subsection{Expansion.}\label{expo.subs}

Now for any fibration $\pi:\C \to \I$ over an ample full subcategory
$\I \subset \Pos$, let $\flat = \pi^\flat(\Tot)$ be the class of
maps in $\C$ cartesian over $\I$. Then $\C_\flat$ is a family of
groupoids over $\I$. It turns out that this establishes a one-to-one
correspondence between proper non-degenerate reflexive families of
categories over $\I$, and families of groupoids satisfying certain
conditions. To state the conditions, recall that $\I$ contains
$\Delta \subset \Posff \subset \Pos$ by Definition~\ref{amp.def}, so
we have the tautological functor
\begin{equation}\label{mu.pos.eq}
\mu:\I \times \Delta \to \I, \qquad J \times [n] \mapsto J
\times [n].
\end{equation}
Let $\eps:\I \to \Delta \times \I$ be its section $J \mapsto J
\times [0]$. For any family of groupoids $\C$ over $\I$, $\mu^*\C$
is a family of groupoids over $\I \times \Delta$, with some
reduction $(\mu^*\C)^{red}$ of \eqref{seg.red}. Then for any map
$f:J \to J'$ in $\I$ with $\Cyl^o(f) \in \I$, the dual cylinder
$\Cyl^o(f)$ fits into a cocartesian square \eqref{cyl.o.sq}, so we
have a functor
\begin{equation}\label{cyl.C.dia}
\begin{CD}
\C_{\Cyl^o(f)} @>>> \C_{J'} \times_{\C_J} \mu^*\C_{J \times [1]}
@>{a}>> \C_{J'} \times_{\C_J} (\mu^*\C)^{red}_{J \times [1]},
\end{CD}
\end{equation}
where $a$ is the reduction epivalence of \eqref{grp.eq}.

\begin{defn}\label{ws.def}
A family of groupoids $\C$ over $\I$ is {\em proper} if for any
right-closed full embedding $f:J' \to J$ in $\I$ \eqref{cyl.C.dia}
is an epivalence, and {\em separated} if this holds when $f$ is also
left-reflexive. A family of groupoids $\C$ over $\I$ is {\em
  bounded} if so is $\mu^*\C$.
\end{defn}

Explicitly, $\C$ is bounded iff for any $J \in \I$, $s^* \times
t^*:\C_{J \times [1]} \to \C_J \times \C_J$ is small. Moreover, the
square \eqref{seg.del.sq} for $l = n-1$ is the square
\eqref{cyl.o.sq} for the left-reflexive full embedding $f=t:[0] \to
      [n-1]$, and its product with $J$ is \eqref{cyl.o.sq} for $\id
      \times t:J \to J \times [n-1]$, so if $\C$ is separated, the
      functor
\begin{equation}\label{ws.eq}
s^* \times t^*: \C_{J \times [n]} \to \C_{J \times [n-1]}
\times_{\C_J} (\mu^*\C)^{red}_{J \times [1]}
\end{equation}
is an epivalence. Thus if $\C$ is separated, and we let
$\iota:\Delta \to \Delta$ be the involution $[n] \mapsto [n]^o$,
then $(\id \times \iota)^*\mu^*\C^* \to \I \times \Delta$ is a weak
$2$-family of groupoids over $\I$ in the sense of
Definition~\ref{w.2.def}. Therefore by Lemma~\ref{w.2.le}, the
reduction $(\id \times \iota)^*(\mu^*\C)^{red}$ is a $2$-family of
groupoids over $\I$ in the sense of
Definition~\ref{exa.seg.def}. Moreover, it is reduced by definition,
bounded in the sense of Definition~\ref{exa.bnd.def} if $\C$ is
bounded in the sense of Definition~\ref{ws.def}, reflexive in the
sense of Definition~\ref{exa.refl.def} if $\C$ is reflexive in the
sense of Definition~\ref{pos.refle.def}, and then also complete in
the sense of Definition~\ref{exa.compl.def} if $\C$ is complete in
the sense of Definition~\ref{pos.compl.def}. If $\C$ has all these
properties, then Proposition~\ref{exa.seg.prop}~\thetag{ii} provides
an equivalence
\begin{equation}\label{del.exp}
(\id \times \iota)^*(\mu^*\C)^{red} \cong \Delta^\flat(\Exp(\C)|\I)
\end{equation}
for a certain family of categories $\Exp(\C) = h^+((\id \times
\iota)^*(\mu^*\C)^{red})$ over $\I$ that we call the {\em expansion}
of the family $\C$. Moreover, we also have
\begin{equation}\label{del.exp.2}
  \Exp(\C) \cong h^+((\id \times \iota)^*\mu^*\C)
\end{equation}
by virtue of Lemma~\ref{w.2.le}.

\begin{prop}\label{exp.prop}
\begin{enumerate}
\item For any reflexive family of categories $\E$ over $\I$ that is
  proper resp.\ separated in the sense of
  Definition~\ref{cat.prp.def} and non-degenerate in the sense of
  Definition~\ref{cat.nondege.def}, the family of groupoids $\C =
  \E_\flat$ is reflexive in the sense of
  Definition~\ref{pos.refle.def}, non-degenerate in the sense of
  Definition~\ref{pos.nondege.def}, complete in the sense of
  Definition~\ref{pos.compl.def}, and bounded and proper
  resp.\ separated in the sense of
  Definition~\ref{ws.def}. Moreover, $\E \cong \Exp(\C)$ over $\I$.
\item Conversely, for any reflexive non-degenerate complete bounded
  proper resp.\ separated family of groupoids $\C$ over $\I$, the
  expansion $\Exp(\C)$ is a proper resp.\ separated non-degenerate
  reflexive family of categories over $\I$, and we have a natural
  equivalence $\C \cong \Exp(\C)_\flat$ over $\I$.
\end{enumerate}
\end{prop}

\proof{} For \thetag{i}, $\C$ is reflexive and non-degenerate by
Lemma~\ref{unf.refl.le}. For any $J \in \I$, $\sigma \times
\tau:\Ar(\E_J) \to \E_J \times \E_J$ is a small discrete fibration,
and since \eqref{nu.J.eq} is an epivalence, $\C$ is bounded, and
$(\mu^*\C)^{red}_{J \times [1]} \cong \Ar(\E_J)_{\Iso}$. Therefore
the source and the target of a functor \eqref{cyl.C.dia} are the
isomorphism groupoids
\begin{equation}\label{C.iso.eq}
\C_{\Cyl^o(f)} \cong (\E_{\Cyl^o(f)})_{\Iso}, \quad \C_{J'}
\times_{\C_J} (\mu^*\C)^{red}_{J \times [1]} \cong (\E_{J'}
\setminus_{f^*} \E_J)_{\Iso},
\end{equation}
and \eqref{cyl.C.dia} is induced by $\nu^\perp_f$, so that $\C$ is
proper if so is $\E$, and separated if so is $\E$ by
Lemma~\ref{cyl.J.le}. It remains to observe that $\C$ is complete by
Proposition~\ref{exa.seg.prop}~\thetag{i} and Lemma~\ref{w.2.le},
and the equivalence $\E \cong \Exp(\C)$ is also provided by
Proposition~\ref{exa.seg.prop}.

For \thetag{ii}, the equivalence $\Exp(\C)_\flat \cong \C$ is again
provided by Proposition~\ref{exa.seg.prop}, so it suffices to prove
that $\Exp(\C)$ is reflexive, proper and non-degenerate. For
reflexivity, with the notation of Definition~\ref{cat.refl.def}, we
need to construct the second adjunction map $\id \to e^* \circ s^*$
for any $J \in \I$. To do it for all $J$ at once, let $S:\I \to \I$
be the functor $J \mapsto J \times [1]$, with the maps $e:S \to
\id$, $s:\id \to S$ that induce functors $s^*:S^*\Exp(\C) \to
\Exp(\C)$, $e^*:\Exp(\C) \to S^*\Exp(\C)$. Then what we need is a
functor
\begin{equation}\label{a.st.eq}
  a_s:S^*\Exp(\C) \times [1] \to S^*\Exp(\C)
\end{equation}
that restricts to $\id$ resp.\ to $e^* \circ s^*$ on $0 \in [1]$
resp.\ $1 \in [1]$, and we have $a_s \circ (s^* \times \id) \cong
s^* \circ p$, where $p:\Exp(\C) \times [1] \to \Exp(\C)$ is the
projection. As in the proof of Proposition~\ref{reg.prop}, to
produce the functor \eqref{a.st.eq}, we use the simplicial
replacement $\Delta [1] = \Delta /[1]$ of the category $[1]$. If
$\pi:\Delta [1] \to \Delta$ is the tautological discrete fibration,
then for any complete reduced bounded $2$-family of groupoids $\C_0
\to \Delta$ over $\ppt$, $\pi_!\pi^*\C_0 \cong \C_0 \times_\Delta
\Delta [1]$ is also a complete reduced bounded $2$-family of
groupoids, and the corresponding localizations of
Proposition~\ref{exa.seg.prop}~\thetag{ii} are related by
$h^+(\pi_!\pi^*\C_0) \cong h^+(\C_0) \times [1]$. Objects in $\Delta
[1]$ are maps $f:[n] \to [1]$. For any such map, we can form a map
$f_s:[n] \times [1] \to [n] \times [1]$ equal to $\id$ on $[n]
\times \{1\}$ and to $\id \times f$ on $[n] \times \{0\}$. Taken
together, these maps define an endomorphism $\alpha_s$ of the
functor
$$
S \circ \mu \circ (\id \times \pi) \cong \mu \circ (S \times \pi):\I
\times \Delta [1] \to \I
$$
sending $J \times \langle [n],f \rangle$ to $J \times [n] \times
[1]$. The endomorphism $\alpha_s$ then defines a functor
$\alpha_s^*:(S \times \pi)^*\mu^*\C \to (S \times \pi)^*\mu^*\C$
over $\I \times \Delta [1]$, and by \eqref{2.adj.bis}, it induces
a functor
\begin{equation}\label{pi.a}
\alpha^\dg_s:(\id \times \pi)_!(\id \times \pi)^*(S \times
\id)^*(\id \times \iota)^*\mu^*\C \to (S \times \id)^*(\id \times
\iota)^*\mu^*\C
\end{equation}
between weak $2$-families of groupoids over $\I$. Taking the
localization $h^+(-)$ then gives \eqref{a.st.eq}, and all the
required properties are immediate from the construction of the
endomorphism $\alpha_s$.

It remains to prove that $\E=\Exp(\C)$ is non-degenerate, and proper
resp.\ separated if so is $\C$. Non-degeneracy is immediate: by
definition, a map $g$ in $\E_J$ is represented by an object $c$ in
$\C_{J \times [1]}$, $g$ is invertible iff $c$ lies in the essential
image of $e^*:\C_J \to \C_{J \times [1]}$, and we are done by
Definition~\ref{pos.nondege.def}. Moreover, by \eqref{C.iso.eq},
$\nu^\perp_f$ for any resp.\ left-reflexive right-closed full
embedding $f$ is an epivalence on the isomorphism groupoids, thus
essentially surjective. It is conservative since so is
\eqref{trunc.eq}, so to finish the proof, it suffices to check that
it is full. By \eqref{C.iso.eq}, we need to show that for any
objects $\wt{c}_0,\wt{c}_1 \in \E_{\Cyl^o(f)}$, with
$\nu^\perp_f(\wt{c}_l) = \langle c'_l,c_l,\alpha_l \rangle \in
\E_{J'} \setminus_{f^*} \E_J$, $l=0,1$, a map $\langle g',g
\rangle:\langle c'_0,c_0,\alpha_0 \rangle \to \langle c'_1,c_1,
\alpha_1 \rangle$ comes from an object $\wt{c} \in \E_{\Cyl^o(f)
  \times [1]}$ with $(\id \times s)^*\wt{c} \cong \wt{c}_0$, $(\id
\times t)^*\wt{c} \cong \wt{c}_1$. Since $\C$ is proper
resp.\ separated, \eqref{cyl.C.dia} is also an epivalence for
$\Cyl^o(f) \times [1] \cong \Cyl^o(f \times \id)$ and the map $f
\times \id:J \times [1] \to J' \times [1]$, so it suffices to lift
$\langle g',g \rangle$ to an object
\begin{equation}\label{c.c.J}
\langle c',c,\alpha \rangle \in \E_{J' \times [1]} \setminus_{(f
  \times \id)^*} \E_{J \times [1]}.
\end{equation}
To do this, first lift the map $g'$ to an object $c' \in \E_{J'
  \times [1]}$, and then use \eqref{cyl.C.dia} for the
left-reflexive full embedding $t:J \to J \times [1]$ to lift
$\langle (f \times \id)^*c',c_1,\alpha_1 \rangle$ to an object $c''
\in \E_{J \times [2]}$. We then have embeddings $s,t:J \times [1]
\to J \times [2]$, and $s$ is right-reflexive, so
Lemma~\ref{refl.J.le} provides a map $\alpha'':s_\dg^*(f \times
\id)^*c' \to c''$ adjoint to the isomorphism $s^*c'' \cong (f \times
\id)^*c'$. Now to obtain \eqref{c.c.J}, take $c = t^*c''$ and
$\alpha = t^*(\alpha''):t^*s_\dg^*(f \times \id)^*c' \cong (f \times
\id)^*c' \to c = t^*c''$.
\endproof

\begin{exa}\label{exp.triv.exa}
Assume that $\E = \Unf(I,\I)$ for some category $I$. Then we have
$\E_\flat \cong \iota^*(\I \bbi I) \subset \iota^*(\I \bb I) = \E$,
where $\I \bbi I \to \I$ is the family of groupoids of
\eqref{pos.I.pi}. Taking the expansion $\Exp(\E_\flat)$ recovers
$\E$ --- indeed, we have
\begin{equation}\label{I.I.loc}
  (\id \times \iota)^*\mu^*\E_\flat \cong \Delta^\flat (\E|\I),
\end{equation}
and taking the localization $h^+(-)$, we get exactly $\E$ by
Proposition~\ref{exa.seg.prop}.
\end{exa}

To complement Proposition~\ref{exp.prop}, say that a reflexive
family of categories $\C$ over $\I$ or the corresponding family of
groupoids $\C_\flat \to \I$ is {\em stably constant} along a map $f$
in $\I$ iff it is constant along $f \times \id_J$ for any finite $J
\in \Posff \subset \I$, and say that $\C$ or $\C_\flat$ is {\em
  stably semicartesian} along a commutative square $\phi:[1]^2 \to
\I$ in $\I$ if it is semicartesian along $\phi \times J$ for any
finite $J$. Moreover, \eqref{unf.sq} defines the
unfolding $\C^\dm \to \I^\dm$ of the family $\C$, and if we again
let $\flat$ be the class of maps cartesian over $\I^\dm$, then
$(\C^\dm)_\flat \cong (\C_\flat)^\dm$ is the unfolding of $\C_\flat$
in the sense of Subsection~\ref{unf.subs}. Say that $\C^\dm$ or
$\C^\dm_\flat = (\C^\dm)_\flat \cong (\C_\flat)^\dm$ is {\em stably
  constant} along a map $f$ in $\I^\dm$ iff it is constant along $f
\times \id_{L(J)}$ for any finite $J$.

\begin{lemma}\label{exp.le}
A separated non-degenerate reflexive family of categories $\C$ over
$\I$ is weakly semicontinuous in the sense of
Definition~\ref{refl.cont.def} or stably semicartesian over a
commutative square in $\I$ if and only if so is $\C_\flat$, and its
unfolding $\C^\dm$ is stably constant along a map in $\I^\dm$ if and
only if so is $\C^\dm_\flat$.
\end{lemma}

\proof{} In all claims, the ``only if'' part is obvious, and the
``if'' part for squares immediately follows from
Lemma~\ref{w.2.le}. For semicontinuity, if $\nu(\C_\flat,J_\idot)$
of \eqref{refl.cont.eq} is an epivalence, then $\nu(\C,J_\idot)$ is
essentially surjective, and if $\nu(\C_\flat,J_\idot \times [1])$ is
also an epivalence, then $\nu(\C,J_\idot)$ is also conservative and
full. For unfoldings, promote \eqref{mu.pos.eq} to a
functor $\mu^\dm:\I^\dm \times \Delta \to \I^\dm$, $J \times [n]
\mapsto J \times L([n])$, and note that for any $J \in \I^\dm$, we
have the $2$-fully faithful embedding $\mu^{\dm*}\C^\dm_\flat \to (U
\times \id)^*\mu^*\C_\flat$, so that by Lemma~\ref{exa.loc.le},
\eqref{del.exp} yields an idenfication $\C^\dm \cong
h^+(\mu^{\dm*}\C^\dm_\flat)$.
\endproof

Another additional result concerns augmented reflexive families of
Definition~\ref{refl.aug.def}. For any partially ordered set $I$,
and any reflexive family $\E$ over $\I$ equipped with an
$I$-augmentation $\E \to \Unf(I,\I)$, let $\hash \supset \flat$ be
the class of maps in $\E$ cartesian over $\Unf(I,\I)$. Then
$\E_\hash$ is a family of groupoids over $\Unf(I,\I)$. Moreover,
\eqref{pos.I.eq} induces a dense embedding $\lambda:\Unf(I,\I)_\flat
\to \Unf(I,\I)$ and the fibrations $\pi:\Unf(I,\I) \to \I$, $\bpi =
\pi \circ \lambda:\Unf(I,\I)_\flat \to \I$, where we identify
$\Unf(I,\I)_\flat \cong \iota^*(\I \bbi I)$ as in
Example~\ref{exp.triv.exa}, and then $\E_\flat \cong
\bpi_!\lambda^*\E_\hash$. Say that a family of groupoids $\C$ over
$\Unf(I,\I)$ is {\em separated} resp. {\em bounded} if the family
$\bpi_!\lambda^*\C \to \I$ is separated rep.\ bounded in the sense
of Definition~\ref{ws.def}.

\begin{lemma}\label{exp.I.le}
For any partially ordered set $I$ and separated non-degene\-rate
$I$-augmented reflexive family of categories $\E$ over $\I$, the
family of group\-oids $\iota^*(\E_\hash) \to \I \bb I \cong
\iota^*\Unf(I,\I)$ is non-degenerate, reflexive in the sense of
Definition~\ref{rel.refle.I.def}, complete in the sense of
Definition~\ref{pos.compl.I.def} and coherent in the sense of
Definition~\ref{cohe.def}. Conversely, for any non-degenerate
complete coherent reflexive family of groupoids $\pi:\C \to \I \bb
I$ such that $\C^\iota = \iota^*\C$ is separated and bounded, the
functor
\begin{equation}\label{E.I.loc}
\E = \Exp(\bpi_!\lambda^*\C^\iota) \to
\Exp(\bpi_!\lambda^*\iota^*(\I \bb I)) \cong \Unf(I,\I)
\end{equation}
is a fibration, thus an augmentation for $\E$, and we have $\C \cong
\iota^*\E_\hash$.
\end{lemma}

\proof{} For the first claim, non-degeneracy, reflexivity and
completeness only depend on $\lambda^*\iota^*(\E_\hash)$, thus hold
by Proposition~\ref{exp.prop}~\thetag{i}. For coherence, note that
in the situation of Definition~\ref{cohe.def}, since $\E$ is a
reflexive family, both $e^*$ and $e_0^*$ are fully faithful, so that
\eqref{a.cyl.eq} is tautologically an equivalence, and then
\eqref{a.red.sq} being semicartesian is equivalent to
\eqref{nu.J.eq} being an epivalence which holds since $\E$ is
separated.

Conversely, since the family of groupoids $\C$ is complete,
non-degenerate, bounded, and reflexive, so is
$\bpi_!\lambda^*\C^\iota$, so that $\E =
\Exp(\bpi_!\lambda^*\C^\iota)$ is a well-defined separated reflexive
family of categories over $\I$. Then \eqref{I.I.loc} fits into a
commutative square
$$
\begin{CD}
(\id \times \iota)^*\mu^*\Unf(I,\I)_\flat @>{\sim}>> \Delta^\flat
  (\Unf(I,\I)|\I)\\
@V{\lambda}VV @VV{\lambda}V\\
(\id \times \iota)^*\mu^*\Unf(I,\I) @>{\sim}>> \Delta^\Tot
(\Unf(I,\I)|\I),
\end{CD}
$$
where the horizontal arrows are equivalences. If we then denote the
projection $\Delta^\Tot(\Unf(I,\I)|I) \cong (\id \times
\iota)^*\mu^*\Unf(I,\I) \to \Unf(I,\I)$ by $\mu_I$, we have an
identification $(\id \times \iota)^*\mu^*\bpi_!\lambda^*\C^\iota
\cong \bpi_!\lambda^*\mu_I^*\C^\iota$, so that $\E \cong
h^+(\bpi_!\lambda^*\mu_I^*\C^\iota)$. Denote by $\eps:\Unf(I,\I) \cong
\Delta^\Tot(\Unf(I,\I)|\I)_{[0]} \to \Delta^\Tot(\Unf(I,\I)|\I)$
the embedding, and define the reduction $(\mu_I^*\C^\iota)^{red}$ by
taking the decomposition
$$
\begin{CD}
\mu_I^*\C^\iota @>>> (\mu_I^*\C^\iota)^{red} @>>>
\eps_*\eps_*\mu_I\C^\iota
\end{CD}
$$
of \eqref{grp.eq} for the functor $\mu_I^*\C^\iota \to
\eps_*\eps^*\mu_I^*\C^\iota$. Then we tautologically have
$\bpi_!\lambda^*(\mu_I^*\C^\iota)^{red} \cong
(\bpi_!\lambda^*\mu_I^*\C^\iota)^{red}$, so that $\E \cong
h^+(\bpi_!\lambda^*(\mu_I^*\C^\iota)^{red})$. But then
\eqref{be.eq} refines $\eps$ to a fully faithful functor
$\beta:\Delta \times \Unf(I,\I) \to \Delta^\Tot(\Unf(I,\I)|I)$, and
the fact that $\C$ is coherent immediately implies that the functor
$$
\bpi_!\lambda^*(\mu_I^*\C^\iota)^{red} \to
\bpi_!\lambda^*\beta_*\beta^*(\mu_I^*\C^\iota)^{red}
$$
is an equivalence of $2$-families of groupoids over $\I$. Therefore
we have $\E \cong
h^+(\bpi_!\lambda^*\beta_*\beta^*(\mu_I^*\C^\iota)^{red})$, and by
\eqref{pd.disc}, we can drop $\bpi_!$ before taking localizations,
while by Lemma~\ref{beta.le} and Lemma~\ref{seg.tw.le}, we
can also drop $\lambda^*\beta_*$, so that $\E \cong
h^+(\beta^*(\mu_I^*\C^\iota)^{red})$. However,
$\beta^*(\mu_I^*\C^\iota)^{red}$ is the reduction of a weak
$2$-family of groupoids $\beta^*\mu_I^*\C^\iota$ over $\Unf(I,\I)$,
so that at the end of the day, $\E \cong
h^+(\beta^*\mu_I^*\C^\iota)$, and then $\E \to \Unf(I,\I)$ is a
fibration by Lemma~\ref{rel.loc.le}, while the equivalence $\C^\iota
\cong \E_\hash$ is provided by Proposition~\ref{exa.seg.prop}.
\endproof

\section{Enhanced categories and functors.}\label{enh.sec}

\subsection{Categories.}

We are now ready to give our main definitions. We start with the
following version of Definition~\ref{pos.def} and
Definition~\ref{refl.exci.def}.

\begin{defn}\label{cat.exci.def}
Assume given a reflexive family of categories $\C$ over an ample
full subcategory $\I \subset \Pos$.
\begin{enumerate}
\item The family $\C$ {\em satisfies excision} if it is
  semicartesian along the square \eqref{J.S.sq} for any set $S$ such
  that $S^<,S^\hush \in \I$, and any $J \in \I$ equipped with a map
  $J \to S^<$.
\item A family $\C$ {\em satisfies the cylinder axiom} if it is
  semicartesian along the square \eqref{J.cyl.sq} for any $J \in \I$
  equipped with a map $J \to [1]$.
\item A family $\C$ is {\em semicontinuous} if $\ZZ_\infty,\N \in
  \I$, and $\C$ is semicartesian along the square \eqref{J.pl.sq}
  for any $J \in \I$ equipped with a map $J \to \N$.
\end{enumerate}
\end{defn}

\begin{lemma}\label{cat.unf.le}
A non-degenerate reflexive family of categories $\C$ over an ample
full subcategory $\I \subset \Pos$ satisfies excision resp.\ the
cylinder axiom resp.\ is semicontinuous iff for any $J \in \I$
equipped with a map $J \to S^<$ such that $S^<,S^\hush \in \I$
resp.\ a map $J \to [1]$ resp.\ a map $J \to \N$, the unfolding
$\C^\dm$ is constant along the corresponding map \eqref{J.S} with
biorder \eqref{B.S.J.bi} resp.\ the map $\zeta_3^{\dm*}L(J/[1]) \to
L(J/[1])$ resp.\ the map $\zeta^{\dm*}L(J/\N) \to L(J/\N)$.
\end{lemma}

\proof{} The argument is the same as in
Proposition~\ref{restr.seg.prop}, but with Lem\-ma~\ref{E.1.bis.le}
replaced by Lemma~\ref{unf.refl.le}.
\endproof

\begin{defn}\label{sm.enh.cat.def}
A {\em restricted enhanced category} is a small additive semi\-exact
separated non-degenerate reflexive family of categories $\C \to
\Posf$ that satisfies excision and the cylinder axiom in the sense
of Definition~\ref{cat.exci.def}.
\end{defn}

\begin{defn}\label{enh.cat.def}
An {\em enhanced category} is an additive semi\-exact separated
non-degenerate reflexive family of categories $\pi:\C \to \Posf^+$
that is semicontinuous, and satisfies excision and the cylinder
axiom in the sense of Definition~\ref{cat.exci.def}. An enhanced
category is {\em small} if so is $\pi$, and {\em $\kappa$-bounded},
for a regular cardinal $\kappa$, if $\|\C_J\| < \kappa$ for any finite
$J \in \Posff \subset \Posf^+$.
\end{defn}

For any two enhanced categories or restricted enhanced categories
$\C$, $\C'$, the coproduct $\C \copr \C'$ is obivously an enhanced
category resp.\ restricted enhanced category, and so is the product
over $\Posf^+$ resp.\ $\Posf$ denoted by
\begin{equation}\label{enh.prod}
  \C \times^h \C' = \C \times_{\Posf^+} \C'.
\end{equation}
More generally, the same holds for coproducts and products of a
family $\C_s$ of enhanced categories resp.\ restricted enhanced
categories indexed by elements $s \in S$ in some set $S$. For any
family of categories $\pi:\C \to \Posf$, we denote by $\flat =
\pi^\flat(\Iso)$ the class of maps in $\C$ that are cartesian over
$\Posf$ and similarly for families over $\Posf^+$.

\begin{lemma}\label{cat.restr.le}
A small additive semiexact separated non-degenerate reflexive family
of categories $\C$ over $\Posf$ resp.\ $\Posf^+$ is a restricted
enhanced category resp.\ a small enhanced category if and only if
the underlying family of groupoids $\C_\flat$ is a restricted Segal
family resp.\ a semicanonical extension of such in the sense of
Definition~\ref{red.seg.def}.
\end{lemma}

\proof{} Combine Proposition~\ref{exp.prop}, Lemma~\ref{exp.le},
and Lemma~\ref{cat.unf.le}.
\endproof

\begin{corr}
A restricted enhanced category $\C \to \Posf$ is cartesian along all
squares \eqref{J.S.sq}, \eqref{J.cyl.sq}, and an enhanced category
is cartesian along all squares \eqref{J.S.sq}, \eqref{J.cyl.sq},
\eqref{J.pl.sq}.
\end{corr}

\proof{} Combine Corollary~\ref{seg.exa.corr} and
Lemma~\ref{exp.le}.
\endproof

\begin{lemma}\label{refl.biano.le}
An additive semiexact separated non-degenerate reflexive family $\C$
over $\Posf^+$ is an enhanced category if and only if its unfolding
$\C^\dm \to \Relf^+$ is constant along all $+$-bianodyne maps of
Definition~\ref{biano.def}, and a small additive semiexact separated
non-degenerate reflexive family $\C$ over $\Posf$ is a restricted
enhanced category if and only if its unfolding $\C^\dm$ is constant
along all bianodyne maps of Definition~\ref{biano.def}
\end{lemma}

\proof{} The ``if'' part immediately follows from
Lemma~\ref{cat.unf.le}. For the ``only if'' part, note that
$\C^\dm_\flat \to \Relf^+$ is a semicanonical extension of a Segal
family by Lemma~\ref{cat.unf.le} and
Proposition~\ref{restr.seg.prop}, and apply Lemma~\ref{bi.del.le}
and Lemma~\ref{exp.le}.
\endproof

\begin{defn}\label{cat.exci.bis.def}
A family of categories $\C \to \Posf^+$ is {\em strongly
  semicontinuous} if for any cofibration $J \to \N$ in $\Posf^+$, it
is cartesian along the square
\begin{equation}\label{J.pl.bis.sq}
\begin{CD}
[1] \times \coprod_{n \geq 1} J_n @>>> \coprod_{n \geq 1} J_n\\
@VVV @VVV\\
\zeta^*J @>>> J
\end{CD}
\end{equation}
induced by \eqref{Z.N.sq}, and it {\em satisfies the strong cylinder
  axiom} if for any map $\chi:J \to [2]$ in $\Posf^+$ such that the
composition $s_\dg \circ \chi:J \to [2] \to [1]$ is a cofibration,
it is cartesian along the square
\begin{equation}\label{J.cyl.bis.sq}
\begin{CD}
J_1 \times [1] @>>> J_1\\
@VVV @VVV\\
\zeta_3^*J @>>> J
\end{CD}
\end{equation}
induced by \eqref{Z3.sq}.
\end{defn}

\begin{corr}\label{cat.exci.bis.corr}
Any enhanced category $\C \to \Posf^+$ is strongly semicontinuous
and satisfies the strong cylinder axiom.
\end{corr}

\proof{} By Lemma~\ref{refl.biano.le}, the unfolding $\C^\dm$ is
constant along the $+$-biano\-dyne maps $\zeta^*J \to J$, $\zeta_3^*J
\to J$ of Lemma~\ref{Z.cyl.le}. If we now promote
\eqref{J.pl.bis.sq} and \eqref{J.cyl.bis.sq} to cartesian squares in
$\Relf^+$, then the top arrows are reflexive, so we
can use the same argument as in Lemma~\ref{cat.unf.le}.
\endproof

\begin{lemma}\label{enh.grp.le}
For any small family of categories $\C$ over $\Posf$ resp.\ over
$\Posf^+$, the following conditions are equivalent.
\begin{enumerate}
\item $\C$ is a restricted enhanced category resp.\ an enhanced
  category, and $\C_\ppt$ is a groupoid.
\item $\C$ is a small enhanced groupoid resp.\ a semicanonical extension
  of such in the sense of Definition~\ref{grp.enh.def}.
\end{enumerate}
\end{lemma}

\proof{} For \thetag{i}$\Rightarrow$\thetag{ii}, note that since
$\C$ is in particular non-degenerate, $\C_J$ is a groupoid for any
$J$, so that $\C$ is a family of groupoids. It is then trivially additive
and semiexact, and since adjoint
functors between groupoids are equivalences, it is also
$[1]$-invariant, thus constant along reflexive maps by
Lemma~\ref{refl.sat.le}. Then excision of
Definition~\ref{cat.exci.def} immediately implies excision of
Definition~\ref{pos.def}.

Conversely, for \thetag{ii}$\Rightarrow$\thetag{i}, additivity and
semiexactness hold by Definition~\ref{pos.def}, reflexivity and
separatedness follows from $[1]$-invariance, non-degeneracy is
trivial, and by virtue of $[1]$-invariance, excision and
semicontinuity of Definition~\ref{pos.def} and
Definition~\ref{pos.pl.def} imply excision and semicontinuity of
Definition~\ref{cat.exci.def}. The cylinder axiom is then trivial,
since all the maps in the corresponding square \eqref{J.cyl.sq} are
reflexive.
\endproof

\begin{lemma}\label{restr.prop.le}
Any restricted Segal family of groupoids $\C \to \Posf$ and any
semicanonical extension $\C^+$ of such a family is proper in the sense
of Definition~\ref{ws.def}. If $\C \to \Posf$ is small, then it is
also bounded, and so is any semicanonical extension $\C^+$ of the family
$\C$.
\end{lemma}

\proof{} For the first claim, by Lemma~\ref{bi.del.le}, $\C^\dm$
resp.\ $\C^{+\dm}$ is constant along bianodyne resp.\ $+$-bianodyne
maps. In particular, both are constant along the map $\xi:B^\dm(J)
\to J$ for any $J \in \Posf$ resp.\ $\Posf^+$. Then for any
right-closed full embedding $f:J \to J'$ in $\Posf$, with the
characteristic map $\chi:{J'}^o \to [1]$, we have $\Cyl^o(f)^o \cong
({J'}^o / [1])^o$, and Lemma~\ref{bary.cyl.le} together with
semiexactness of $\C^\dm$ show that already the first functor in
\eqref{cyl.C.dia} is an epivalence. For $\C^+$, consider the height
function $\hht:J' \to \N$, and apply
Definition~\ref{biano.def}~\thetag{ii} over $\N$ and
Lemma~\ref{bary.cyl.le} to deduce that the embeddings
\eqref{B.J.1.dm} for $\chi:{J'}^o \to [1]$ are $+$-bianodyne. For
the second claim, $\C$ is bounded since it has small fibers, and
then $\C^+$ is bounded by semicontinuity and Lemma~\ref{tel.I.le}.
\endproof

\begin{corr}\label{cat.prop.corr}
For any small complete restricted Segal family of grou\-p\-oids $\C$
over $\Posf$, the expansion $\Exp(\C)$ of Proposition~\ref{exp.prop}
is a restricted enhanced category, and for any semicanonical extension
$\C^+ \to \Posf^+$, the expansion $\Exp(\C^+)$ is an enhanced
category. Any enhanced category $\C^+ \to \Posf^+$ is proper in the
sense of Definition~\ref{cat.prp.def}.
\end{corr}

\proof{} Immediately follows from Lemma~\ref{cat.unf.le} and
Lemma~\ref{exp.le}.
\endproof

\begin{corr}\label{cat.iota.corr}
For any restricted enhanced category $\C \to \Posf$, the fibration
$\C^\iota = \iota^*\C^o_\perp \to \Posf$ is a restricted enhanced
category. For any enhanced category $\C \to \Posf^+$, the fibration
$\C^\iota = (B_\dm^*\C^\dm)^o_\perp \to \Posf^+$ is an enhanced
category, and $\xi_\perp:B_\dm \to \iota$ provides an equivalence
$\C^\iota |_{\Posf} \cong \iota^*\C^o_\perp$.
\end{corr}

\proof{} The fibration $\C^\iota$ is clearly a non-degenerate
reflexive family; it is then a restricted enhanced category
resp.\ enhanced category by Lemma~\ref{cat.restr.le} and
Corollary~\ref{iota.corr} resp.\ Proposition~\ref{iota.prop}.
\endproof

\begin{defn}
The enhanced category resp.\ restricted enhanced category $\C^\iota$ of
Corollary~\ref{cat.iota.corr} is the {\em enhanced opposite} to the
enhanced category resp.\ restricted enhanced category $\C$.
\end{defn}

\begin{corr}\label{cat.epi.corr}
For any enhanced category or restricted enhaced category $\C$, and any
map $f:J \to J'$ in $\Posf$, the functors \eqref{cyl.J.eq} are
epivalences.
\end{corr}

\proof{} For $\nu^\perp_f$, this is Lemma~\ref{refl.J.le} combined
with Corollary~\ref{cat.prop.corr}. For $\nu_f$, note that $\nu_f$
for $\C$ is opposite to $\nu^\perp_{f^o}$ for $\C^\iota$.
\endproof

\begin{exa}\label{cat.triv.exa}
For any category $\E$, let us simplify notation and denote by
$\Unf(\E) = \Unf(\E,\Posf^+) \to \Posf^+$ the reflexive family of
\eqref{unf.eq}. Then $\Unf(\E)$ is trivially proper and
non-degenerate, and it is an enhanced category in the sense of
Definition~\ref{enh.cat.def}: indeed, all the squares
\eqref{J.S.sq}, \eqref{J.cyl.sq}, \eqref{J.pl.sq} are cocartesian
squares in $\Cat$, and moreover, so is any standard pushout square
\eqref{pos.st.sq} (so that the family $\Unf(\E) \to \Posf^+$ is not
only semiexact but actually exact). We have $\Unf(\E)^\iota \cong
\Unf(\E^o)$, and if $\E$ is small, then $\Unf(\E)$ restricted to
$\Posf \subset \Posf^+$ is a restricted enhanced category in the
sense of Definition~\ref{sm.enh.cat.def}.
\end{exa}

\begin{exa}\label{ppt.h.exa}
If $\E = \ppt$, then $\Unf(\ppt) = \Posf^+$ is the terminal enhanced
category; we denote it by $\ppt^h$.
\end{exa}

\subsection{Functors.}

Having defined enhanced categories, we now define enhanced functors.

\begin{defn}\label{enh.fun.defn}
An {\em enhanced functor} between restricted enhanced categories
resp.\ enhanced categories is a functor cartesian over $\Posf$
resp.\ $\Posf^+$, and an {\em enhanced morphism} between two such
functors is a morphism over $\Posf$ resp.\ $\Posf^+$.
\end{defn}

\begin{exa}\label{fun.triv.exa}
For any enhanced category or restricted enhanced category $\C$, the
truncation functor \eqref{trunc.eq} is an enhanced functor. For any
categories $\E$, $\E'$ and functor $\gamma:\E \to \E'$, the functor
$\Unf(\gamma):\Unf(\E) \to \Unf(\E')$ between the corresponding
enhanced categories of Example~\ref{cat.triv.exa} is an enhanced
functor.
\end{exa}

\begin{defn}\label{enh.obj.def}
An {\em enhanced object} $c$ in an enhanced category $\C$ is an
enhanced functor $\eps^h(c):\ppt^h \to \C$, and an {\em enhanced
  morphism} betweem two enhanced objects $c$, $c'$ is an enhanced
morphism $\eps^h(c) \to \eps^h(c')$ in the sense of
Definition~\ref{enh.fun.def}.
\end{defn}

Since $\ppt \in \Posf^+ = \ppt^h$ is the terminal object, enhanced
functors $\ppt^h \to \C$ are uniquely determined by their values at
$\ppt$, so that effectively, an enhanced object in an enhanced
category $\C$ is the same thing as an object of $\C_\ppt$. An
enhanced morphism between two such objects $c,c' \in \C_\ppt$ is a
morphism $c \to c'$ in $\C_\ppt$.

\begin{lemma}\label{sq.cat.le}
Assume given a functor $\gamma:\E \to \E'$, and a cartesian square
\begin{equation}\label{enh.pb.sq}
\begin{CD}
\C @>>> \C'\\
@VVV @VV{\pi}V\\
\Unf(\E) @>{\Unf(\gamma)}>> \Unf(\E'),
\end{CD}
\end{equation}
where $\C'$ is an enhanced category, and $\pi$ is an enhanced
functor. Then $\C$ is an enhanced category.
\end{lemma}

\proof{} By Lemma~\ref{refl.pb.le}, $\C \to \Posf^+$ is a separated
reflexive family. It is obviously additive, and to check the rest of
the properties, apply the criterion of Lemma~\ref{car.rel.le}.
\endproof

\begin{exa}
Say that a {\em closed class of enhanced morphisms} in an enhanced
category $\C$ is a closed class $v$ of morphisms in $\C_\ppt$.  For
any such class $v$, we can define a category $\C_{hv}$ by the
cartesian square
\begin{equation}\label{h.v.sq}
\begin{CD}
\C_{hv} @>>> \C\\
@VVV @VV{\unf(\C)}V\\
\Unf((\C_{\ppt})_v) @>>> \Unf(\C_\ppt),
\end{CD}
\end{equation}
where the vertical arrow on the right is the truncation functor
\eqref{trunc.eq}, and $(\C_{\ppt})_v \subset \C_\ppt$ is the dense
subcategory defined by the class $v$. Then $\C_{hv}$ is an enhanced
category by Lemma~\ref{sq.cat.le}. In particular, if $v = \Iso$ is
the class of all invertible maps, we obtain the {\em enhanced
  isomorphism groupoid} $\C_{h\Iso}$ given by the cartesian square
\begin{equation}\label{h.iso.sq}
\begin{CD}
\C_{h\Iso} @>>> \C\\
@VVV @VV{\unf(\C)}V\\
\Unf((\C_{\ppt})_{\Iso}) @>>> \Unf(\C_\ppt).
\end{CD}
\end{equation}
It is an enhanced groupoid by Lemma~\ref{enh.grp.le}, and is
universal: any enhanced functor $\C' \to \C$ from an enhanced
groupoid $\C'$ factors through $\C_{h\Iso} \to \C$, uniquely up to a
unique isomorphism.
\end{exa}

An enhanced functor $\gamma:\C_0 \to \C_1$ between restricted
enhanced categories resp.\ enhanced categories trivially extends to
a functor $\gamma^\dm:\C_0^\dm \to \C_1^\dm$ between their
unfoldings, cartesian over $\Relf$ resp.\ $\Relf^+$, and then in
particular, we have a well-defined {\em opposite enhanced functor}
$\gamma^\iota:\C_0^\iota \to \C_1^\iota$.

\begin{defn}\label{cat.aug.def}
For any category $I$, an {\em $I$-augmentation} of an enhanced
category $\C$ is an enhanced functor $\pi:\C \to \Unf(I)$ that is a
fibration, and an {\em $I$-coaugmentation} of $\C$ is an enhanced
functor $\pi:\C \to \Unf(I)$ such that $\pi^\iota$ is an
$I$-augmentation. An enhanced functor $\gamma:\C' \to \C$ between
$I$-augmented enhanced categories $\C$, $\C'$ is {\em $I$-augmented}
if it is cartesian over $\Unf(I)$, and if $\C$, $\C'$ are
$I$-coaugmented, then $\gamma$ is {\em $I$-coaugmented} if
$\gamma^\iota$ is $I$-augmented.
\end{defn}

By definition, for any $I$-augmented enhanced category $\C$,
$\C^\iota$ is $I^o$-co\-aug\-mented, and vice versa. If $I$ is a
partially ordered set, then we also have an enhanced version of the
transpose-opposite fibration construction: for any $I$-augmented
enhanced category $\C$, we can consider the category
\begin{equation}\label{iota.perp.eq}
  \C^\iota_{h\perp} = B_{I\dm}^*(\C|I)^{\dm o}_\perp,
\end{equation}
where $B_{I\dm}$ is the barycentric subdivision of
Example~\ref{bary.rel.I.exa}, and by Lem\-ma~\ref{iota.I.le},
Lemma~\ref{seg.I.le}, Lemma~\ref{exp.I.le} and
Lemma~\ref{cat.restr.le}, $\C^\iota_{h\perp}$ is also an
$I$-augmented enhanced category (while $\C_{h\perp} =
(\C^\iota_{h\perp})^\iota$ is then $I^o$-coaugmented). For any
category $\C$ equipped with an enhanced functor $\pi:\C \to
\Unf(I)$, and any $i \in I$, we define the {\em enhanced
  fiber} $\C_i$ by
\begin{equation}\label{enh.aug.fib}
  \C_i = \Unf(\eps(i))^*\C,
\end{equation}
where $\eps(i):\ppt \to I$ is as in \eqref{eps.c}. By
Lemma~\ref{sq.cat.le}, $\C_i$ is an enhanced category; if $\pi$ is
an $I$-augmentation, then $(\C_{h\perp})_i \cong \C_i$ and
$(\C^\iota_{h\perp})_i \cong \C_i^\iota$.

\begin{remark}\label{I.rig.rem}
Note that for any partially ordered set $I$, the fibers of the
fibration $\Unf(I) \to \Posf^+$ are rigid. Therefore if we have
enhanced categories $\C_0$, $\C_1$ equipped with enhanced functors
$\pi_l:\C_l \to \Unf(I)$, $l=0,1$, then for any enhanced functor
$\gamma:\C_0 \to \C_1$, being a functor over $\Unf(I)$ is a
condition and not a structure. If $\pi_0$ and $\pi_1$ are
augmentations resp.\ coaugmentations, then the same of course holds
for being augmented resp.\ coaugmented.
\end{remark}

If $I$ is a partially ordered set, then $I$-augmentations of
Definition~\ref{cat.aug.def} are the same thing as $I$-augmentations
of Definition~\ref{refl.aug.def} over $\I = \Posf^+$. If $\C =
\Unf(I')$ for some $I'$, then any enhanced functor $\pi:\C \to
\Unf(I)$ is of the form $\pi \cong \Unf(\gamma)$ for some functor
$\gamma:I' \to I$, and $\pi$ is an $I$-augmentation resp.\ an
$I$-coaugmentation iff $\gamma$ is a fibration resp.\ a
cofibration. For any functor $\gamma:I' \to I$ and an $I$-augmented
resp.\ $I$-coaugmented enhanced category $\C$, the pullback enhanced
category $\C' = \Unf(\gamma)^*\C$ of Lemma~\ref{sq.cat.le} is
$I'$-augmented resp.\ $I'$-coaugmented. From now on, by abuse of
notation, we will denote
\begin{equation}\label{f.h.pb}
  f^{h*}\C = \Unf(f)^*\C.
\end{equation}
One can also construct a pushforward operation $f^h_*$ for augmented
categories, but in general, this requires some preparations (see
Lemma~\ref{pf.cat.le} below). One case when it is easy is when we
consider a right-closed embedding $f:I' \to I$ of finite-dimensional
partially ordered sets. In this case, $\Unf(f):\Unf(I',\I) \to
\Unf(I,\I)$ admits a right-adjoint $\Unf(f)_\dg$, $\langle J,\alpha
\rangle \mapsto J \times_I I'$, and while $\Unf(f)_\dg$ is not an
enhanced functor --- it is not even a functor over $\Posf^+$ --- we
still have the following.

\begin{lemma}\label{enh.cl.le}
For any right-closed embedding $f:I' \to I$ in $\Posf$, and
any $I'$-augmented enhanced category $\C'$, $\C=\Unf(f)_*\C'$ is an
$I$-augmented enhanced category.
\end{lemma}

\proof{} By Corollary~\ref{pf.corr}, $\C$ is a proper reflexive
family, and $\C \cong \Unf(f)_\dg^*\C'$. But
$\Unf(f)_\dg$ preserves coproducts and sends cocartesian squares
\eqref{pos.st.sq}, \eqref{J.S.sq}, \eqref{J.pl.sq}, \eqref{J.cyl.sq}
to squares of the same form.
\endproof

Combining Lemma~\ref{cat.restr.le}, Corollary~\ref{cat.prop.corr}
and Proposition~\ref{exp.prop}, we see that an enhanced category $\C
\to \Posf^+$ is completely determined by the underlying family of
groupoids $\C_\flat \to \Posf^+$. It is useful to observe that this
also works for enhanced functors, in the following form.

\begin{lemma}\label{enh.exp.le}
Assume given enhanced categories $\C$, $\C'$, with underlying
families of groupoids $\C_\flat,\C'_\flat \to \Posf^+$. Then any
functor $\gamma_\flat:\C'_\flat \to \C'_\flat$ extends to an
enhanced functor $\gamma:\C' \to \C$, and for any two enhanced
functors $\gamma,\gamma':\C' \to \C$, a morphism between their
restrictions to $\C'_\flat$ uniquely extends to a enhanced morphism
$\gamma \to \gamma'$.
\end{lemma}

\proof{} Immediately follows from Proposition~\ref{exp.prop}.
\endproof

\begin{lemma}\label{enh.epi.le}
An enhanced functor $\gamma:\C' \to \C$ that is an epivalence is an
equivalence.
\end{lemma}

\proof{} Combine Proposition~\ref{exp.prop} and
Corollary~\ref{semi.epi.corr}.
\endproof

\begin{defn}\label{enh.refl.def}
An enhanced functor $\gamma:\C \to \C'$ between restricted enhanced
categories resp.\ enhanced categories is {\em left} resp.\ {\em
  right-reflexive} if it is left resp.\ right-reflexive over $\Posf$
resp.\ $\Posf^+$, and the adjoint functor $\gamma^\dg$ is also an
enhanced functor.
\end{defn}

\begin{exa}\label{unf.adj.exa}
For any functor $\gamma:\E \to \E'$ between categories $\E$, $\E'$,
the enhanced functor $\Unf(\gamma):\Unf(\E) \to \Unf(\E')$ is left
resp.\ right-reflexive if and only if so is $\gamma$, and if
$\gamma_\dg$ is left resp.\ right-adjoint to $\gamma$, then
$\Unf(\gamma_\dg)$ is left resp.\ right-adjoint to
$\Unf(\gamma_\dg)$. In one direction, this is obvious from the
identification $\Unf(\E)_\ppt \cong \E$, and in the other direction,
the claim immediately follows from Example~\ref{fun.adj.exa}.
\end{exa}

\begin{lemma}
An enhanced functor $\gamma:\C_0 \to \C_1$ is right
resp. left-reflexive if and only if $\gamma^\iota:\C_0^\iota \to
\C_1^\iota$ is left resp.\ right-reflexive. An enhanced functor
$\gamma$ between enhanced categories is left or right-reflexive if
and only if so is its restriction to $\Posf$. A left or
right-reflexive enhanced functor $\gamma:\C \to \C'$ is an
equivalence iff this holds for the fiber $\gamma(\ppt):\C_\ppt \to
\C'_\ppt$.
\end{lemma}

\proof{} The first claim directly follows from
Lemma~\ref{fib.adj.le}. For the second one, use semicontinuity and
the same argument as in Lemma~\ref{refl.pb.le}. For the last claim,
let $\gamma^\dg$ be the adjoint, and let $a$, $a^\dg$ be the
adjunction maps; then $\gamma$ is an equivalence iff $a$ and $a^\dg$
are invertible, and since $\C$ and $\C'$ are non-degenerate, this
holds iff it holds over $\ppt \in \Posf^+$.
\endproof

\begin{defn}\label{enh.o.def}
An enhanced object $c$ in an enhanced category $\C$ is {\em initial}
resp.\ {\em terminal} if the corresponding enhanced functor
$\eps^h(c)$ is right resp.\ left-reflexive. An enhanced functor
$\gamma:\C \to \C'$ between enhanced categories with initial
resp.\ terminal objects $o$, $o'$ {\em preserves the initial}
resp.\ {\em terminal object} if the base change map $\eps(o') \to
\gamma \circ \eps(o)$ resp.\ $\gamma \circ \eps(o) \to \eps(o')$ is
invertible.
\end{defn}

\begin{exa}\label{enh.term.exa}
An initial or terminal enhanced object $o$ in an enhanced category
$\C$ is also the initial resp.\ terminal object in $\C_\ppt$, but
the converse certainly need not be true. For instance, an enhanced
groupoid has an initial or terminal enhanced object only if it is
trivial, but there are plenty of non-trivial enhanced groupods $\C$
with $\C_\ppt = \ppt$ -- take any simply connected but non-trivial
homotopy type.
\end{exa}

\begin{exa}\label{enh.gt.exa}
If $\C = \Unf(\E)$ for some category $\E$, then $e \in \C_\ppt = \E$
defines an initial resp.\ terminal enhanced object in $\C$ iff $e$
is initial resp.\ terminal in $\E$. For example, if $\E = [1]$, then
$\Unf([1])$ has an initial and a terminal enhanced object
$\eps^h(0),\eps^h(1):\Posf^+ \to \Unf([1])$, and no other enhanced
objects. For a general enhanced category $\C$, we can define
families $\C^{h>},\C^{h<} \to \Posf^+$ by
\begin{equation}\label{gt.lt.eq}
  \C^{h<} = \eps^h(1)_*\C, \qquad \C^{h>} = (\C^{\iota h<})^\iota.
\end{equation}
Then $\C^{h<}$ is a $[1]$-augmented enhanced category by
Lemma~\ref{enh.cl.le}, and dually, $\C^{h<}$ is a $[1]$-coaugmented
enhanced category. We have $\C^{h<}_\ppt \cong \C_\ppt^<$,
$\C^{h>}_\ppt \cong \C_\ppt^>$, and the enhanced object defined by
$o \in \C^<_\ppt$ resp.\ $o \in \C^>_\ppt$ is initial
resp.\ terminal. If $\C = \Unf(\E)$, then $\C^{h<} \cong \Unf(\E^<)$
and $\C^{h>} \cong \Unf(\E^>)$.
\end{exa}

An enhanced functor $\gamma:\C \to \C'$ is {\em fully faithful} if
it is fully faithful in the usual sense, and for such functors,
there is the following useful result.

\begin{lemma}\label{ff.cat.le}
A fully faithful enhanced functor $\gamma:\C' \to \C$ between
enhanced categories or restricted enhanced categories is an
equivalence if and only if it is an equivalence over $\ppt$.
\end{lemma}

\proof{} The ``only if'' part is tautological. For the ``if'' part,
since $\gamma$ is fully faithful, we just have to check that it is
essentially surjective over any $J \in \Posf^+$ resp.\ $\Posf$. Say
that $J$ is {\em good} if this happens. Then by
Definition~\ref{cat.exci.def}~\thetag{i}, being good is stable under
standard pushouts, and by \thetag{ii}, it is stable under taking
coproducts. Then as in the proof of Lemma~\ref{E.3.le}, excision
shows that if $J$ admits a map $J \to S^<$ with good comma-fibers,
then $J$ itself is good, and, for enhanced categories,
semicontinuity shows that if $J$ admits a map $J \to \N$ with good
comma-fibers, it is again good. Thus it suffices to show that all $J
\in \Posf$ are good, and induction on dimension reduces this to
showing that $J^>$ is good for a good $J \in \Posf$. This follows
from Corollary~\ref{cat.epi.corr} applied to the map $f:J \to \ppt$.
\endproof

\begin{corr}\label{enh.ff.corr}
For any enhanced category $\C$ equipped with a full subcategory $\C'_\ppt
\subset \C_\ppt$, the category $\C'$ defined by the cartesian square
\begin{equation}\label{ff.cat.sq}
\begin{CD}
  \C' @>>> \C\\
  @VVV @VVV\\
  \Unf(\C'_\ppt) @>>> \Unf(\C_\ppt)
\end{CD}
\end{equation}
is an enhanced category, the enhanced functor $\C' \to \C$ is fully
faithful, and all fully faithful enhanced functors arise in this
way.
\end{corr}

\proof{} The fact that $\C'$ is an enhanced category is
Lemma~\ref{sq.cat.le}. Then $\C' \to \C$ is the pullback of a fully
faithful functor, thus fully faithful, and every fully faithful
functor is of this form by Lemma~\ref{ff.cat.le}.
\endproof

In terms of Definition~\ref{enh.obj.def},
Corollary~\ref{enh.ff.corr} says that a collection of enhanced
objects in an enhanced category $\C$ uniquely defines a full
enhanced subcategory $\C' \subset \C$, and every full enhanced
subcategory appears in this way (we will say that $\C'$ is {\em
  spanned} by the collection of enhanced objects $c \in \C'_\ppt
\subset \C_\ppt$). If we understand $\C$ as an ``enhancement'' for
$\C_\ppt$, then the claim is that a full subcategory $\C'_\ppt
\subset \C_\ppt$ inherits a unique enhancement. For any enhanced
functor $\gamma:\C' \to \C$, we define the {\em enhanced essential
  image} of $\gamma$ as the full enhanced subcategory in $\C$
spanned by the essential image of $\gamma_\ppt:\C'_\ppt \to
\C_\ppt$; any $\gamma:\C' \to \C$ uniquely factors through its
enhanced essential image. As in the unenhanced setting, we will say
that an enhanced full subcategory $\C' \subset \C$ is {\em left}
resp.\ {\em right-admissible} if the embedding enhanced functor $\C'
\to \C$ is left resp.\ right-reflexive in the sense of
Definition~\ref{enh.refl.def}.

\begin{lemma}\label{enh.filt.ff.le}
Let $I$ be a filtered category, and assume given a fibration $\C \to
\Posf^+ \times I^o$ such that $\C_i \to \Posf^+$ is an enhanced
category for any $i \in I$, and $f^{o*}:\C_i \to \C_{i'}$ is fully
faithful for any map $f:i \to i'$ in $I$. Let $\C_\perp \to I$ be
the cofibration tranpose to $\C \to \Posf^+ \times I^o \to
I^o$. Then the $2$-colimit $\colim_I\C_\perp$ of
Definition~\ref{filt.ff.def} is an enhanced category, and for any $i
\in I$, the enhanced functor $\C_i \to \colim_I\C_\perp$ is fully
faithful.
\end{lemma}

\proof{} All the claims that we need to check to insure that a
family of categories $\E \to \Posf^+$ is an enhanced category
involve only a finite number of objects and morphisms in $\E$ (for
reflexivity, this is Lemma~\ref{a.b.adj.le}). Since
by Lemma~\ref{filt.ff.le}, $\colim_I\C_\perp$ is the union of its
full subcategories $\C_i$, and $I$ is filtered, everything can
checked at each individual $\C_i$ where it holds by assumption.
\endproof

\subsection{Barycentric extensions.}

So far, our only examples of enhanced categories are the trivial
ones of Example~\ref{cat.triv.exa}. To obtain more interesting
examples, one can use complete families of
Definition~\ref{compl.fam.def}.

\begin{lemma}\label{comp.C.le}
Assume given an additive complete family of categories $\C$ over
$\Posf^+$ that is tight in the sense of
Definition~\ref{refl.tight.def}, and assume in addition that $\C$ is
either non-degenerate, or weakly semicontinuous in the sense of
Definition~\ref{refl.cont.def}. Then $\C$ is reflexive, proper,
weakly semicontinuous and non-degenerate, and it is an enhanced
category in the sense of Definition~\ref{enh.cat.def}.
\end{lemma}

\proof{} If $\C$ is weakly semicontinuous, then it is non-degenerate
and proper by Lemma~\ref{tight.le}. If $\C$ is non-degenerate, then
it is proper by Lemma~\ref{tight.le}. In either case, $\C$ is then
reflexive by Lemma~\ref{compl.refl.le} and then also semiexact by
Lemma~\ref{refl.semi.le}. Since by Lemma~\ref{exp.le} and
Lemma~\ref{tel.I.le}, any separated reflexive family $\C \to
\Posf^+$ that is semicontinuous in the sense of
Definition~\ref{cat.exci.def} is weakly semicontinuous, it remains
to prove that $\C$ is en enhanced category, and by
Lemma~\ref{refl.biano.le}, it suffices to prove that the unfolding
$\C^\dm$ is constant along all $+$-bianodyne maps. Indeed, let $W$
be the class of maps $f:J \to J'$ in $\Relf^+$ such that the adjoint
pair of functors $U(f)^*$, $U(f)_*$ restricts to a pair of mutually
inverse equivalences between $\C^\dm_J \subset \C_{U(J)}$ and
$\C^\dm_{J'} \subset \C_{U(J')}$ --- that it, for any $c \in
\C^\dm_J$, $c' \in \C^\dm_{J'}$, the adjunction maps $U(f)^*U(f)_*c
\to c$, $c' \to U(f)_*U(f)^*c'$ are isomorphisms. Then since $\C$ is
reflexive, $W$ satisfies the condition \thetag{i} of
Definition~\ref{biano.def}, and by Lemma~\ref{bc.J.le} it also
satisfies \thetag{ii}, thus contains all $+$-bianodyne maps.
\endproof

One problem with Lemma~\ref{comp.C.le} is that $\C$ is required to
be complete (in particular, $\C \to \Posf^+$ cannot have small
fibers). To circumvent this, one can use the barycentric subdivision
functor $B^\dm:\Pos \to \Relf^\pm$, as in
Subsection~\ref{bary.dua.subs}. Say thay an ample full subcategory
$\I \subset \Pos$ is {\em bar-invariant} if the barycentric
subdivision functor $B$ sends $\I$ into itself.

\begin{defn}
A reflexive family of categories $\C$ over a bar-invariant ample
subcategory $\I \subset \Pos$ is {\em bar-invariant} if for any $J
\in \I$, the unfolding $\C^\dm$ is constant along the map
$\xi:B^\dm(J) \to J$.
\end{defn}

\begin{exa}\label{bar.exa}
The category $\Posf^+$ is bar-invariant, and any enhanced category
$\C \to \Posf^+$ is bar-invariant by Corollary~\ref{B.biano.corr}
and Lemma~\ref{refl.biano.le}.
\end{exa}

\begin{defn}\label{lf.good.def}
A reflexive family of categories $\C$ over $\Posf^\pm$ is {\em good}
if its unfolding $\C^\dm$ is constant along all the $\pm$-bianodyne
maps of Definition~\ref{pm.biano.def}.
\end{defn}

\begin{defn}\label{exci.mod.def}
A family of categories $\C \to \Posf^\pm$ {\em satisfies modified
  excision} if for any set $S$ and $J \in \Posf^\pm$ equipped with a
map $J \to S^>$, $\C$ is semicartesian along the square
\begin{equation}\label{J.S.mod.sq}
\begin{CD}
J_o \times B(S^>) @>>> J_o\\
@VVV @VVV\\  
J^\hhush @>>> J,
\end{CD}
\end{equation}
where $J^\hhush$ is as in \eqref{B.S.J.bis.bi}, and the top
horizontal arrow is the map \eqref{J.S.hh}.
\end{defn}
  
\begin{lemma}\label{lf.good.le}
An additive semiexact reflexive family of categories $\C$ over
$\Posf^\pm$ is good in the sense of Definition~\ref{lf.good.def} iff
it satisfies modified excision in the sense of
Definition~\ref{exci.mod.def}, and is semicontinuous and satisfies
the clynder axiom in the sense of Definition~\ref{cat.exci.def}.
\end{lemma}

\proof{} As in Lemma~\ref{refl.biano.le}, the ``only if'' part holds
by the same argument as in Lemma~\ref{cat.unf.le}, and for the
``if'' part, one needs to also invoke
Proposition~\ref{pm.biano.prop}.
\endproof

\begin{prop}\label{lf.prop}
For any bar-invariant ample subcategory $\I \subset \Pos$ containing
$\Posf^\pm$, a reflexive family $\C \to \Posf^\pm$ good in the sense
of Definition~\ref{lf.good.def} extends to a proper bar-invariant
reflexive family of categories $\C^\Tot$ over $\I$. Moreover, a
functor $\gamma:\C_0 \to \C_1$ between two good families with
bar-invariant extensions $\C_0^\Tot$, $\C_1^\Tot$ extends to a
functor $\gamma^+:\C_0^\Tot \to \C_1^\Tot$, and a morphism $f_0 \to
f_1$ between two such functors extends to a unique morphism
$f_0^\Tot \to f_1^\Tot$, so that the bar-invariant extension
$\C^\Tot$ of a good family $\C$ is defined uniquely up to an
equivalence unique up to a unique isomorphism. In addition to this,
the bar-invariant extension $\C^\Tot$ of a good family $\C \to
\Posf^\pm$ is semicontinuous and satisfies excision and the cylinder
axiom in the sense of Definition~\ref{cat.exci.def}, and it is
additive, semiexact, non-degenerate, separated or proper if so is
$\C$, so that if $\C$ has all these properties, and $\I = \Posf^+$,
then $\C^\Tot$ is an enhanced category. If $I \supset \Posf^+$, then
$\C^\Tot$ is also strongly semicontinuous and satisfies the strong
cylinder axiom in the sense of Definition~\ref{cat.exci.bis.def}.
\end{prop}

\proof{} To construct the extension $\C^\Tot$, consider the barycentric
subdivision functor $B^\dm:\I \to \Relf^\pm \subset \I^\dm$, and let
$\C^\Tot = B^{\dm*}\C^\dm$. For any $J \in \Posf^+$, we have a
$\pm$-bianodyne map $B^\dm(J \times [1]) \to B^\dm(J) \times L([1])$
of Example~\ref{B.i.exa}, so that $\C^\Tot$ is a reflexive family of
categories, non-degenerate and/or separated if so is $\C$. More
generally, for any right-closed embedding $f:J' \to J$ with
characteristic map $\chi:J \to [1]$, we have $\Cyl^o(f) \cong J
\times_{[1]} [2]$, so that Lemma~\ref{pm.B.le} provides a
$\pm$-bianodyne map $B^\dm(\Cyl^o(f)) \to \Cyl^o(B^\dm(f))$, and
therefore $\C^\Tot$ is proper if so is $\C$. Moreover, $\C^\Tot$ is
bar-invariant by Lemma~\ref{B.B.le}, and for any left-finite $J \in
\Posf^\pm \subset \I$, the functorial map $\xi:B^\dm(J) \to L(J)$ is
$\pm$-bianodyne and provides an equivalence $\C^\Tot|_{\Posf^\pm} \cong
\C$. For any other bar-invariant extension $\C'$, the same map
provides a canonical equivalence $\C' \cong \C^\Tot$, and the
extensions for functors and morphisms are then also obtained by
applying $B^\dm$. Since $B^\dm$ preserves coproducts and standard
pushout squares, $\C^\Tot$ is additive and semiexact if so is
$\C$. Moreover, being good, $\C$ is semicontinuous and satisfies the
cylinder axiom by the same argument as in
Lemma~\ref{cat.unf.le}. Then for any $J \in \Posf^+$ equipped with a
map $J \to \N$, Lemma~\ref{pm.B.le} provides a $\pm$-bianodyne map
$B^\dm(J / \N) \to B^\dm(J) / L(\N)$, and similarly for all the
other terms in the cartesian square \eqref{J.pl.sq}, so that
semicontinuity for $\C$ and the map $B(J) \to \N$ implies
semicontinuity for $\C^\Tot$ and $J \to \N$. Analogously, the cylinder
axiom for $\C$ implies the cylinder axiom for $\C^\Tot$. For excision,
assume given $J \in \I$ equipped with a map $J \to S^<$. Let
$B(J)^\hhush$ be the underlying partially ordered set of the
biordered set $B^\dm(J)^\hhush$ of \eqref{B.S.J.bis.bi} (and
$B^\dm(J)^\fhush$ of \eqref{B.S.J.iii.bi}). Then the map
\eqref{J.S.hh} for $B(J)$ fits into a cartesian cocartesian square
\begin{equation}\label{J.S.hh.sq}
\begin{CD}
  B(J)^\hhush @>>> B(J)\\
  @AAA @AAA\\
  B(J_o) \times B(S^>) @>>> B(J_o).
\end{CD}
\end{equation}
By the same argument as in Lemma~\ref{cat.unf.le},
Lemma~\ref{J.hh.le} implies that $\C$ is cartesian along the square
\eqref{J.S.hh.sq}, and then the unfolding $\C^\dm$ is cartesian
along the corresponding square
$$
\begin{CD}
  B^\dm(J)^\fhush @>>> B^\dm(J)\\
  @AAA @AAA\\
  B^\dm(J_o) \times B^\dm(S^>) @>>> B^\dm(J_o).
\end{CD}
$$
in $\Posf^\pm$. By Corollary~\ref{pm.B.corr}, this implies excision
for $\C^\Tot$. Finally, assume that $\I \supset \Posf^+$, and define
the class of {\em $\I$-bianodyne maps} in $\I^\dm$ as the smallest
saturated closed class that is closed under standard pushouts,
contains all reflexive maps, and satisfies the following version of
Definition~\ref{biano.def}~\thetag{ii}:
\begin{itemize}
\item for any bicofibrations $J',J'' \to J$ in $\I^\dm$ such that $J
  \in \Posf^+$, a map $f:J' \to J''$ cocartesian over $J$ and such
  that the fibers $f:J'_j \to J''_j$ are $\I$-bianodyne for any $j
  \in J$ is itself $\I$-bianodyne.
\end{itemize}
Then exactly the same induction as in Proposition~\ref{biano.prop}
shows that the unfolding $\C^{\Tot\dm}$ of the family $\C^\Tot$ is
constant along all the $\I$-bianodyne maps, and then
Lemma~\ref{Z.cyl.le} also holds for $\I$-bianodyne maps with the
same proof, so that $\C^\Tot$ is strongly semicontinuous and satifies
the strong cylinder axiom by the same argument as in
Corollary~\ref{cat.exci.bis.corr}.
\endproof

\begin{defn}\label{fin.compl.def}
A family of categories $\C$ over $\Posf^\pm$ is {\em finitely
  complete} if for any left-finite $f:J' \to J$ in $\Posf^\pm$,
the transition functor $f^*:\C_J \to \C_{J'}$ admits a right-adjoint
$f_*:\C_{J'} \to \C_J$, and for any square \eqref{posf.sq} in
$\Posf^\pm$ satisfying \thetag{i} or \thetag{ii} of
Definition~\ref{compl.fam.def}, the base change map \eqref{bc.eq}
is an isomorphism \eqref{refl.bc} between the corresponding
functors $\C_{J_0} \to \C_{J_1'}$.
\end{defn}

\begin{remark}
Since any full embedding in $\Posf^\pm$ is left-finite,
Definition~\ref{refl.tight.def} makes sense for finitely complete
families of categories over $\Posf^\pm$.
\end{remark}

\begin{lemma}\label{lf.le}
An additive non-degenerate tight finitely complete family of
categories $\C$ over $\Posf^\pm$ is good in the sense of
Definition~\ref{lf.good.def}, so that it extends to an enhanced
category $\C^\Tot \to \Posf^+$, with the same functoriality properties
as in Proposition~\ref{lf.prop}.
\end{lemma}

\proof{} Same as Lemma~\ref{comp.C.le}.
\endproof

We note that even if we start with an enhanced category $\C \to
\Posf^+$ that is already defined over $\Posf^+$,
Proposition~\ref{lf.prop} is still useful. For example, one can take
$\I = \Pos$, and extend $\C$ to a family of categories $\C^\Tot \to
\Pos$ over all partially ordered sets. We will call $\C^\Tot$ the
    {\em canonical extension} of the enhanced category $\C$; here
    are some of its properties.

\begin{lemma}\label{bar.enh.le}
Let $\C^\Tot = B^{\dm*}\C^\dm \to \Pos$ be the bar-invariant extension
of an enhanced category $\C \to \Posf^+$ provided by
Proposition~\ref{lf.prop}.
\begin{enumerate}
\item The fibration $\iota^*(\C^\Tot)^o_\perp \to \Pos$ is the
  bar-invariant extension of the opposite enhanced category
  $\C^\iota$.
\item For any left-closed embedding $l:J' \to J$ in $\Pos$, $\C^\Tot$
  is semicartesian along the corresponding cylinder square
  \eqref{cyl.o.sq}.
\item $\C^\Tot$ is weakly semicontinuous in the sense of
  Definition~\ref{refl.cont.def}.
\end{enumerate}
\end{lemma}

\proof{} For \thetag{i}, $\iota^*(\C^\Tot)^o_\perp$ is obviously
bar-invariant, and then its restriction to $\Posf^+$ is identified
with $\C^\iota = B_\dm^*\C^\dm$ simply by definition: the
constructions of Proposition~\ref{lf.prop} and
Proposition~\ref{iota.prop} are literally the same. For \thetag{ii},
apply extended cylinder axiom and semiexactness, as in
Corollary~\ref{cyl.2.corr}. Finally, for \thetag{iii}, for any map
$J \to \N$ in $\Pos$, semicontinuity for $\C^\Tot$ immediately implies
that the unfolding $\C^{\Tot\dm}$ is constant along the corresponding
map $L(J)^+ \to L(J)$ of \eqref{J.pl.bi.eq}. Then if we compose the
projection $U(L(J)^+) \to \ZZ_\infty$ with the map \eqref{x.z.eq},
and rewrite the corresponding standard pushout square in the
alternative form \eqref{pos.2.st.sq}, as in \eqref{Z.V.N.sq}, then
$\C^{\Tot\dm}$ is still semicartesian along this square by
semiexactness. If we now denote $\C_\infty = \C^\Tot_J$,
$\C_n=\C^\Tot_{J/n}$, $n \geq 0$, this means that the square
\begin{equation}\label{C.wc.sq}
\begin{CD}
  \C^\dm_{L(J)^+} \cong \C_\infty @>>> \prod_{n \geq 0}\C_{(2n+1)}\\
  @VVV @VVV\\
  \prod_{n \geq 0}\C_{2(n+1)} @>>> \prod_{n \geq 0}\C_n
\end{CD}
\end{equation}
is semicartesian. But $\Sec^\Tot(\N,\C_\idot) \cong
\Sec^\Tot(\ZZ_\infty,\zeta^*\C_\idot)$ by
Lemma~\ref{Z.le}. Moreover, if we evaluate this category by
\eqref{Z.V.N.sq}, we obtain \eqref{C.wc.sq} with $\C_\infty$
replaced by $\Sec^\Tot(\N,\C_\idot)$, but the square becomes
cartesian. Thus the functor $\C_\infty \to \Sec^\Tot(\N,\C_\idot)$
is an epivalence, and this is what we had to prove.
\endproof

\subsection{Examples and applications.}\label{enh.exa.subs}

We can now use Proposition~\ref{lf.prop} and Lemma~\ref{lf.le} to
construct two series of examples of enhanced categories in the sense
of Definition~\ref{enh.cat.def}. First, assume given a model
category $\C = \langle \C,C,F,W \rangle$, as in
Section~\ref{mod.sec}. Consider the enhanced category $\Unf(\C) =
\Unf(\C,\Posf^+)$ of Example~\ref{cat.triv.exa} and its restriction
$\Unf(\C,\Posf^\pm)$ to $\Posf^\pm \subset \Posf^+$. Since $\C$ is
by definition finitely complete, $\Unf(\C,\Posf^\pm)$ is a finitely
complete family over $\Posf^\pm$ in the sense of
Definition~\ref{fin.compl.def}, and if $\C$ is in fact complete,
then $\pi:\Unf(\C) \to \Posf^+$ is a complete family in the sense of
Definition~\ref{compl.fam.def}. For any $J \in \Posf^\pm$, the fiber
$\Unf(\C,\Posf^\pm)_J \cong J^o\C$ carries a Reedy model structure,
and if $\C$ is complete, the same holds for $J \in \Posf^+$ and
$\Unf(\C)_J$. In particular, extend $W$ to a closed class of maps in
$\Unf(\C)$ by defining a dense subcategory $\Unf(\C)_W$ by the
cartesian square
\begin{equation}\label{unf.mod.W.sq}
\begin{CD}
  \Unf(\C)_W @>>> \Unf(\C)\\
  @VVV @VVV\\
  \eps_*\C_W @>>> \eps_*\C,
\end{CD}
\end{equation}
where $\eps = \eps(\ppt):\ppt \to \Posf^\hdot$ is the embedding onto
$\ppt \in \Posf^+$. Then the relative category $\langle
\Unf(\C,\Posf^\pm),W \rangle$ is localizable by
Lemma~\ref{rel.loc.le}, and we have the fibration
$h^W(\Unf(\C,\Posf^\pm)) \to \Posf^\pm$. If $\C$ is complete, then
the same holds for $\langle \Unf(\C),W \rangle$.

\begin{lemma}\label{enh.mod.le}
For any model category $\langle \C,C,F,W \rangle$, the family of
categories $h^W(\Unf(\C,\Posf^\pm))$ over $\Posf^\pm$ is additive,
tight, non-degenerate and finitely complete, so that its extension
$\Hh^W(\C) = h^W(\Unf(\C,\Posf^\pm))^+$ provided by Lemma~\ref{lf.le}
is an enhanced category. Moreover, $\Hh^W(\C) \cong h^W(\Unf(\C))$
if $\C$ is complete.  A Quillen-adjoint pairs of functors $\lambda:\C \to
\C'$, $\rho:\C' \to \C$ between model categories $\C$, $\C'$ induces
an adjoint pair of enhanced functors $\Hh(\lambda):\Hh^W(\C) \to
\Hh^W(\C')$, $\Hh(\rho):\Hh^W(\C') \to \Hh^W(\C)$.
\end{lemma}

\proof{} The family $h^W(\Unf(\C,\Posf^+))$ is additive since
localization of model categories commutes with products, and
non-degenerate since the functors \eqref{mod.ev} are
conservative. For any left-finite map $f:J_0 \to J_1$, the adjoint
$f_*:h^W(J_0^o\C) \to h^W(J_1^o\C)$ exists by the Quilled Adjunction
Theorem~\ref{qui.adj.thm}, and to check that a map \eqref{refl.bc}
is an isomorphism, it suffices to take fibrant replacements and
apply \eqref{refl.bc} for the finitely complete family
$\Unf(\C,\Posf^\pm)$. Therefore $h^W(\Unf(\C,\Posf^\pm))$ is
finitely complete, and it is tight by
Proposition~\ref{glue.prop}. If $\C$ is complete, then
$h^W(\Unf(\C))$ is additive, complete, tight and non-degenerate by
the same argument, thus bar-invariant by Example~\ref{bar.exa}, and
then $h^W(\Unf(\C)) \cong \Hh^W(\C)$ by the uniqueness part of
Proposition~\ref{lf.prop}. Finally, for any $J \in \Posf^+$, a
Quillen-adjoint pairs of functors between some $\C$ and $\C'$
induces a Quillen-adjoint pair of functors between $J^o\C$ and
$J^o\C'$, so the last claim also immediately follows from
Theorem~\ref{qui.adj.thm}.
\endproof

Secondly, assume given an additive category $\A$, and consider the
category $C_\idot(\A)$ of chain complexes in $\A$, with the
saturated CW-structure of Example~\ref{CW.add.exa}. Again, consider
the enhanced category $\Unf(C_\idot(\A),\Posf^\pm)$ of
Example~\ref{cat.triv.exa}, and extend $W$ to a closed class of maps
in $\Unf(C_\idot(\A),\Posf^\pm)$ by
\begin{equation}\label{unf.add.W.sq}
\begin{CD}
  \Unf(\C)_W @>>> \Unf(\C)\\
  @VVV @VVV\\
  \eps_*\C_W @>>> \eps_*\C,
\end{CD}
\end{equation}
where $\eps$ is as in \eqref{unf.mod.W.sq}.

\begin{lemma}\label{enh.ho.le}
The relative category $\langle \Unf(C_\idot(\A)),\Posf^\pm,W\rangle$
is localizable for any additive category $\A$, and the localization
$h^W(\Unf(C_\idot(\A)))$ is an additive tight non-degenerate and
finitely complete family of categories over $\Posf^\pm$, so that its
extension $\Hho(\A)=h^W(\Unf(C_\idot(\A)))^+$ of Lemma~\ref{lf.le}
is an enhanced category.
\end{lemma}

\proof{} Note that the opposite category $\A^o$ is also
additive. Then the relative category $\langle
\Unf(C_\idot(\A),\Posf^\pm),W\rangle$ is localizable by
Lemma~\ref{rel.loc.le} and Corollary~\ref{add.A.corr} applied to
$\A^o$, and $h^W(\Unf(C_\idot(\A),\Posf^\pm)$ is finitely complete
by Corollary~\ref{add.A.corr} and \eqref{pos.kan.eq}. Moreover,
$h^W(\Unf(C_\idot(\A),\Posf^\pm)$ is obviously additive, it is
non-degenerate by Lemma~\ref{add.loc.le}, and tight by
\eqref{ho.A.J.eq} and Lemma~\ref{add.glue.le} applied to
$\Cof(C_\idot(\A^o)^J)^o$.
\endproof

As another application of Proposition~\ref{lf.prop}, we can
construct a pushforward operation for augmented enhanced categories
of Definition~\ref{cat.aug.def}. Namely, assume given a map $f:I'
\to I$ between partially ordered sets $I,I' \in \Posf^+$. Then for
any $I$-augmented enhanced category $\C$, we have the $I'$-augmented
category $f^{h*}\C$ of \eqref{f.h.pb}.

\begin{lemma}\label{pf.cat.le}
For any map $f:I \to I'$ in $\Posf^+$ and $I$-augmented enhanced
category $\C$, there exists an $I$-augmented enhanced category
$f^h_*\C$ and an $I$-augmented enhanced functor
$\alpha:f^{h*}f^h_*\C \to \C$ such that for any $I'$-aug\-mented
enhanced category $\C'$, any $I$-augmented enhanced functor
$f^{h*}\C' \to \C$ factors as
\begin{equation}\label{ab.enh.dia}
\begin{CD}
f^{h*}\C' @>{f^{h*}\beta}>>
f^{h*}f^h_*\C @>{\alpha}>> \C,
\end{CD}
\end{equation}
for a certain $I'$-augmented enhanced functor $\beta:\C' \to
f^h_*\C$, unique up to a unique isomorphism. Moreover, for any
commutative square \eqref{posf.sq} in $\Posf^+$ satisfying
\thetag{i} or \thetag{ii} of Definition~\ref{compl.fam.def}, and any
$J_0$-augmented enhanced category $\C$, the base change functor
\begin{equation}\label{bc.enh.eq}
g_1^{h*}f^h_*\C \to {f'}^h_*g_0^{h*}\C
\end{equation}
induced by \eqref{ab.enh.dia} is an equivalence.
\end{lemma}

\proof{} Let $\I = \Pos$, and denote $\Unf^\Tot(I) = \Unf(I,\I)$,
$\Unf^\Tot(I') = \Unf(I',\I)$. By Proposition~\ref{lf.prop}, any
enhanced category $\C$ admits a unique bar-invariant extension
$\C^\Tot \to \I$, and an augmentation $\pi:\C \to \Unf(I')$ extends
to an augmentation $\pi^\Tot:\C^\Tot \to \Unf^\Tot(I')$ in the sense
of Definition~\ref{refl.aug.def}. By Lemma~\ref{pm.B.le},
$\Unf^\Tot(f)^*\C^\Tot$ is then bar-invariant, thus canonically
identified with the extension $(f^{h*}\C)^\Tot$. Now, for any
$I$-augmented enhanced category $\C$ with the canonical
bar-invariant extension $\C^\Tot \to \I$, we can let $f^h_*\C =
\Unf^\Tot(f)_*\C^\Tot|_{\Unf(I')}$, and then the decompositions
\eqref{ab.enh.dia} are provided by Corollary~\ref{kan.corr}, so to
prove the first claim, it suffices to check that $f^h_*\C$ is an
$I'$-augmented enhanced category. Since $\C^\Tot$ is a proper
non-degenerate $I$-augmented reflexive family by
Proposition~\ref{lf.prop}, $f^h_*\C$ is a proper non-degenerate
$I'$-augmented reflexive family by Corollary~\ref{pf.corr} that is
moreover given by \eqref{Unf.g.C}. Then both $\Y$ and $\Unf(f)_\dg$
preserve coproducts and standard pushout squares, so that $f^h_*\C$
is additive and semiexact. Moreover, $\Unf(f)_\dg \circ \Y$ sends
cocartesian squares of the form \eqref{J.S.sq}, \eqref{J.cyl.sq},
\eqref{J.pl.sq} to the corresponding squares of Lemma~\ref{j.ind.le}
in $\Relf^+$, so that $f^h_*\C$ is semicontinuous and satisfies
excision and the cylinder axiom. The second claim now
follows from \eqref{2kan.eq}, by the same argument as in
Corollary~\ref{pf.corr}.
\endproof

Yet another simple but useful corollary of Proposition~\ref{lf.prop}
is the following observation: if we think of an enhanced category
$\C$ as providing an enhancement for $\C_\ppt$, then $\C$ considered
simply as a category also inherits a canonical enhancement. Namely,
we have the following.

\begin{lemma}\label{C.wr.le}
For any enhanced category $\C$, with the bar-invariant extension
$\C^\Tot \to \Pos$ and the co-reflexive family $\sigma^*\C^{\Tot
  o}_\perp \to \Pos$ provided by Lemma~\ref{C.refl.le}, the restriction
$\C_\wr \to \Posf^+$ of the transpose-opposite fibration
$(\sigma^*\C^{\Tot o}_\perp)^o_\perp$ to $\Posf^+ \subset \Pos$ is
an enhanced category, and $(\C_\wr)_\ppt \cong \C^\Tot$.
\end{lemma}

\proof{} The family $\C_\wr \to \Posf^+$ is reflexive and proper by
Lemma~\ref{C.refl.le} and Corollary~\ref{cat.prop.corr}. It is
obviously additive, so we need to check that it is semiexact,
semicontinuous and satisfies excision and the cylinder axiom. Each
of these checks amount to showing that $\C_\wr$ is cartesian or
semicartesian along a certain cocartesian square \eqref{j.ind.sq} in
$\Posf^+$. Moreover, \eqref{C.fib.dia} induces a fibration $\C_\wr
\to \Unf(\Posf^+)$, and $\Unf(\Posf^+)$ is cartesian along all these
squares, so we can apply the criterion of
Lemma~\ref{car.rel.le}. For semiexactness or excision, we start with
$J \in \Posf^+$ equipped with a map $J \to \V$ resp.\ $J \to S^<$,
and then choosing a cartesian section $[1]^2 \to \Unf(\Posf^+)$ of
the fibration $\Unf(\Posf^+) \to \Posf^+$ over the corresponding
square in $\Posf^+$ amounts to choosing a $\Posf^+$-augmentation
$\alpha$ for $J$. Then if we let $J_\idot = \alpha^*\Pos_\idot \to
J$ be the corresponding cofibration, we are done by semiexactness
resp.\ excision for $\C^\Tot$ along the square corresponding to
$J_\idot \to \V$ resp.\ $J_\idot \to S^<$. For semicontinuity, what
we choose is a $\Posf^+$-augmentation for $J/\N$, not $J$; however,
then the composition $J_\idot \to J/\N \to \N$ is a cofibration, and
by Proposition~\ref{lf.prop}, $\C^\Tot$ is strongly
semicontinuous. Analogously, for the cylinder axiom, $J_\idot \to
J/[1]$ is a cofibration, so the composition map $J_\idot \to J/[1]
\to [1]$ is a cofibration as well, and $\C^\Tot$ satisfies the
strong cylinder axiom, so again, we are done.
\endproof

\begin{defn}\label{enh.fun.def}
For any enhanced categories $\E$, $\E'$, the {\em enhanced functor
  category} $\fFun^h(\E,\E')$ is an enhanced category equipped with
an enhanced functor
\begin{equation}\label{enh.ev.eq}
\ev:\E \times^h \fFun^h(\E,\E') \to \E'
\end{equation}
such that for any enhanced category $\C$, an enhanced functor
$\gamma:\E \times^h \C \to \E'$ factors as
\begin{equation}\label{enh.ev.dia}
\begin{CD}
\E \times^h \C @>{\id \times \gamma'}>> \E \times^h \fFun^h(\E,\E')
@>{\ev}>> \E'
\end{CD}
\end{equation}
for some $\gamma':\C \to \fFun^h(\E,\E')$, uniquely up to a unique
isomorphism.
\end{defn}

\begin{corr}\label{J.C.corr}
For any enhanced category $\C$ and $J \in \Pos$, there exists an
enhanced functor category $J^o_h\C = \fFun^h(\Unf(J^o),\C)$.
\end{corr}

\proof{} Let $\eps_J:\Posf^+ \to \Unf(\Pos)$ be the unique
cartesian section of the fibration $\Unf(\Pos) \to \Posf^+$
sending $\ppt$ to $J \in \Unf(\Pos)_\ppt \cong \Pos$ ---
explicitly, $\eps_J$ sends $J' \in \Posf^+$ to $J' \times J$ --- and
let $J^o_h\C \cong \eps_J^*\C_\wr$. Then $J^o_h\C$ is an enhanced
category by Lemma~\ref{C.wr.le} and Lemma~\ref{sq.cat.le}, and the
transpose-opposite co-reflexive family $(J^o_h\C)^o_\perp$ is
$(\C^o_\perp)^{hJ}$, so it has the required universal property by
Corollary~\ref{refl.ev.corr}.
\endproof

If $J=[1]$, we denote $\Ar_h(\C)= [1]^o_h\C \cong
\fFun(\Unf([1]),\C)$ and call it the {\em enhanced arrow category}
of the enhanced category $\C$; we denote by
\begin{equation}\label{enh.ar.eq}
\sigma,\tau:\Ar_h(\C) \to \C, \qquad \eta:\C \to \Ar_h(\C)
\end{equation}
the enhanced functors induced by $s,t:\ppt \to [1]$ and $e:[1] \to
\ppt$. Reflexivity of $\C$ immediately implies that just as in the
unenhanced situation, $\sigma$ resp.\ $\tau$ is right
resp.\ left-adjoint to $\eta$. The functors \eqref{nu.J.eq} taken
together define an epivalence
\begin{equation}\label{ar.epi}
  \nu:\Ar^h(\C) \to \Ar(\C|\Posf^+)
\end{equation}
whose target is the relative arrow category of \eqref{rel.ar.sq},
and this is an equivalence if $\C = \Unf(\E)$ for some category
$\E$. If the set $J$ is arbitrary but $\C$ is equipped with an
enhanced functor $\pi:\C \to \Unf(J^o)$ -- for example, if it is
$J^o$-augmented or $J^o$-coaugmented in the sense of
Definition~\ref{cat.aug.def} -- then we can also define the enhanced
section category $\sSec^h(J^o,\C)$ by the cartesian square
\begin{equation}\label{ssec.enh.I}
\begin{CD}
  \sSec^h(J^o,\C) @>>> J^o_h\C\\
  @VVV @VVV\\
  \ppt^h @>{\eps^h(\id)}>> J^o_h\Unf(J^o),
\end{CD}
\end{equation}
where $\id$ in the bottom arrow is
$\id \in \Fun^h(\Unf(J^o),\Unf(J^o))_\ppt \cong \Fun(J^o,J^o)$. We let
$\Sec^h(J^o,\C) = \sSec(J^o,\C)_\ppt$, and if we are given some
other $J' \in \Pos$ equipped with a map $\alpha:J' \to J$, we
will simplify notation by writing $\Sec^h({J'}^o,\C) =
\Sec^h({J'}^o,\Unf(\alpha^o)^*\C)$ and $\sSec^h({J'}^o,\C) =
\sSec^h({J'}^o,\Unf(\alpha^o)^*\C)$. Equivalently, we can define
$\sSec^h({J'}^o,\C)$ by the cartesian square
\begin{equation}\label{ssec.enh.II}
\begin{CD}
  \sSec^h({J'}^o,\C) @>>> {J'}^o_h\C\\
  @VVV @VVV\\
  \ppt^h @>{\eps^h(\alpha^o)}>> {J'}^o_h\Unf(J^o).
\end{CD}
\end{equation}
Explicitly, we have an equivalence
\begin{equation}\label{sec.h.eq}
\Sec^h({J'}^o,\C) \cong \C_{\langle J',\alpha^o \rangle}
\end{equation}
between $\Sec^h({J'}^o,\C)$ and the fiber of $\pi:\C \to \Unf(J^o)$
over the $J^o$-augmented set $\langle {J'},\alpha^o \rangle$, and if
we have two such sets $J'_0$, $J'_1$, a map $f:J'_0 \to J'_1$ over
$J$ induces an enhanced functor $f^{o*}:\sSec^h({J'}^o_1,\C) \to
\sSec^h({J'}^o_0,\C)$. In particular, for any $j \in J$, we have the
evaluation enhanced functor $\sSec^h(J^o,\C) \to \C_i$ to the
enhanced fiber \eqref{enh.aug.fib}. If $\pi$ is a $J^o$-augmentation
resp.\ coaugmentation, there is more functoriality: $f^{o*}$ is
defined when $f$ is only lax resp.\ co-lax over $J$. We also have
the following enhanced version of Lemma~\ref{2adm.le}.

\begin{lemma}\label{enh.2adm.le}
Assume given a partially ordered set $J$, a right-reflexive full
embedding $f:J' \to J$, with adjoint map $f_\dg:J \to J'$, and a
$J^o$-aug\-mented enhanced category $\C$. Then the restriction
functors
$$
f^{o*}:\Sec^h(J^o,\C) \to \Sec^h({J'}^o,\C), \qquad
f^*:\Sec^h(J,\C_{h\perp}) \to \Sec^h(J',\C_{h\perp})
$$
are right resp.\ left-reflexive, with fully faithful adjoints given
by $f_\dg^{o*}$ resp.\ $f_\dg^*$.
\end{lemma}

\proof{} For $\C$, apply Lemma~\ref{unf.refl.J.le} to $L(f):\langle
L(J'),f^o \rangle \to \langle L(J),\id \rangle$. For $\C_{h\perp}$,
apply the same argument to $\C^\iota_{h\perp}$.
\endproof

\begin{lemma}\label{sec.sq.le}
For any $J \in \Posf^+$ and $J$-augmented enhanced category $\C$,
the evaluation functor
\begin{equation}\label{nondeg.sec}
  \sSec^h(J,\C) \to \prod_{j \in J} \C_j
\end{equation}
is conservative. Moreover, assume given a commutative square
\begin{equation}\label{J.sq}
\begin{CD}
  J_{01} @>>> J_0\\
  @VVV @VVV\\
  J_1 @>>> J
\end{CD}
\end{equation}
in $\Posf^+$ that is either \thetag{i} a standard pushout square, or
\thetag{ii} one of the squares \eqref{J.S.sq}, \eqref{J.pl.sq},
\eqref{J.cyl.sq}, \eqref{J.pl.bis.sq}, \eqref{J.cyl.bis.sq}. Then
for any $J$-augmented enhanced category $\C$, the induced square
\begin{equation}\label{sec.J.sq}
\begin{CD}
\sSec^h(J,\C) @>>> \sSec^h(J_0,\C)\\
@VVV @VVV\\
\sSec^h(J_1,\C) @>>> \sSec^h(J_{01},\C)
\end{CD}
\end{equation}
of enhanced categories and functors is \thetag{i} semicartesian
resp.\ \thetag{ii} cartesian.
\end{lemma}

\proof{} Immediately follows from \eqref{sec.h.eq} and the
corresponding properties of the enhanced category $\C$.
\endproof

\section{Representability and its consequences.}\label{enh.repr.sec}

\subsection{Representability for enhanced categories.}\label{enh.repr.subs}

Let us now prove a version of representability theorems of
Section~\ref{repr.sec} for small enhanced categories. Assume given a
complete fibrant Segal space $X \in \Delta^o\Delta^o\Sets$ in the
sense of Definition~\ref{bi.compl.def}, and consider the
corresponding discrete fibration $\Delta\Delta X \to \Delta \times
\Delta$. As in Section~\ref{seg.sp.sec}, we refer to the two factors
in the product $\Delta \times \Delta$ as to the {\em Reedy} and {\em
  simplicial} factor. Let
\begin{equation}\label{phi.re.eq}
  \phi = ((F \circ \iota) \times \id) \circ \delta:\Delta \to
  \Posf^+ \times \Delta,
\end{equation}
where $\delta:\Delta \to \Delta \times \Delta$ is the diagonal
embedding, $F:\Delta \to \Posf^+$ is the standard full embedding,
and $\iota:\Delta \to \Delta$ is the involution $[n] \mapsto
[n]^o$. Then we have the discrete fibration
\begin{equation}\label{pi.re.eq}
\pi:\E(X) = (\phi \times \id)_*\Delta\Delta X \to \Posf^+ \times \Delta
\times \Delta,
\end{equation}
where $\phi \times \id:\Delta \times \Delta \to \Posf^+ \times
\Delta \times \Delta$ is induced by \eqref{phi.re.eq}, and the
category $\E(X)$ in \eqref{pi.re.eq} comes equipped with the closed
class of maps
\begin{equation}\label{E.re.eq}
  E = \pi^*(\Iso \times + \times \operatorname{\Tot}).
\end{equation}
Explicitly, $f \in E$ iff $\pi(f) = f_0 \times f_1 \times f_2$,
where $f$ is a invertible map in $\Posf^+$, $f_1$ is a special map
in the Reedy factor $\Delta$, and $f_2$ is any map in the simplicial
factor $\Delta$. If $X=\ppt$ is the constant Segal space with value
$\ppt$, so that \eqref{pi.re.eq} is an equivalence, then $\E(X)
\cong \Posf^+ \times \Delta \times \Delta$ is localizable with
respect to the class $E$ of \eqref{E.re.eq}, and $h^E(\Posf^+ \times
\Delta \times \Delta) \cong \Posf^+$.

\begin{lemma}\label{cat.repr.le}
For any complete fibrant Segal space $X \in \Delta^o\Delta^o\Sets$,
the category $\E(X)$ of \eqref{pi.re.eq} is localizable with respect
to the class $E$ of \eqref{E.re.eq}, and the localization
\begin{equation}\label{unf.re.eq}
  \Unf(X) = h^E(\E(X))
\end{equation}
together with the projection $\Unf(X) \to \Posf^+$ induced by
\eqref{pi.re.eq} is an enhanced category in the sense of
Definition~\ref{enh.cat.def}.
\end{lemma}

\proof{} The projection $\E(X) \to \Posf^+ \times \Delta \times
\Delta \to \Posf^+$ is a fibration with small fibers, so by
Lemma~\ref{rel.loc.le}, the relative category $\langle \E(X),E
\rangle$ is localizable, and $\Unf(X) \to \Posf^+$ is a
fibration. Moreover, we can first localize with respect to a smaller
class $\pi^*(\Iso \times \Iso \times \Tot)$, and this amounts to
taking $\theta_!$ with respect to the projection $\theta:\Posf^+
\times \Delta \times \Delta \to \Posf^+ \times \Delta$ onto the
Reedy factor, so that we have
\begin{equation}\label{unf.psi.eq}
\Unf(X) \cong h^+(\theta_!\E(X)),
\end{equation}
where $+$ is the class of special maps with respect to the Segal
fibration $\theta_!\E(X) \to \Posf^+ \times \Delta \to \Delta$. Since
$F \circ \iota:\Delta \to \Posf^+$ is fully faithful, for any
fibration $\E \to \Delta$, we have
\begin{equation}\label{phi.mu.eq}
\phi_*\E \cong (\id \times (F \circ \iota))^*\delta_*(F \circ
\iota)_*\E \cong (\id \times \iota)^*\mu^*(F \circ \iota)_*\E,
\end{equation}
where $\delta:\Posf^+ \to \Posf^+ \times \Posf^+$ is the diagonal
embedding, and $\mu$ is the product functor \eqref{mu.pos.eq}. The
same identification then holds for any fibration $\E \to \Delta
\times \Delta$ and $(\phi \times \id)_*\E$, and $(\id \times
\iota)^*\mu^*$ commutes with $\theta_!$ by \eqref{bc.fib}, so
that if we let
\begin{equation}\label{C.psi.eq}
  \C = \theta_!((F \circ \iota) \times \id)_*\Delta\Delta X,
\end{equation}
then \eqref{unf.psi.eq} and \eqref{phi.mu.eq} provide an equivalence
$\Unf(X) \cong h^+((\id \times \iota)^*\mu^*\C)$, and this is in
turn identified with the expansion $\Exp(\C)$ by
\eqref{del.exp.2}. Then by Corollary~\ref{cat.prop.corr}, it
suffices to check that $\C \to \Posf^+$ is a complete restricted
Segal family of groupoids. However, we have the fully faithful nerve
functor $\Nn:\Posf^+ \to \Delta^o\Sets$, and the composition $\Nn
\circ F:\Delta \to \Delta^o\Sets$ is the Yoneda embedding $\Y$, so
that we have
\begin{equation}\label{C.N.eq}
\C \cong \theta_!(\Nn \times \id)^*(\Y \times \id)_*(\iota \times
\id)_*\Delta\Delta X \cong \Nn^*\theta_!(\Y \times
\id)_*\Delta\Delta X^\iota,
\end{equation}
where we denote $X^\iota = (\iota \times \id)^{o*}X \cong (\iota
\times \id)^o_!X$. But $X^\iota$ is also a complete fibrant Segal
space, it represents a family of groupoids $\Hh(X^\iota) \to
\Delta^o\Delta^o\Sets$ described in Lemma~\ref{I.cosimp.le}, and as
in Remark~\ref{I.cosimp.rem}, \eqref{2kan.eq} together with this
description provide an equivalence $\theta_!(\Y \times
\id)_*\Delta\Delta X^\iota \cong L^*\Hh(X^\iota)$. We are then done
by Proposition~\ref{css.prop}.
\endproof

\begin{exa}
Assume given a small category $I$, and let $X = L(N(I))$ be its
nerve considered as a bisimplicial set constant in the simplicial
direction. Then $X$ is a fibrant Segal space, complete if $I$ is
rigid, and in this case, $\Unf(X) \cong \Unf(I)$ is the enhanced
category corresponding to $I$, as in Example~\ref{cat.triv.exa}.
\end{exa}

\begin{theorem}\label{enh.thm}
Let $i:\Posf \to \Posf^+$ be the natural full embedding. Then for
any restricted enhanced category $\E \to \Posf$ in the sense of
Definition~\ref{sm.enh.cat.def}, there exists a complete fibrant
Segal space $X$ and an equivalence $\E \cong i^*\Unf(X)$ over
$\Posf$, and for any such $X$ and any enhanced category $\E^+ \to
\Posf^+$ in the sense of Definition~\ref{enh.cat.def}, an
equivalence $i^*\E^+ \cong \E$ over $\Posf$ comes from an
equivalence $\E^+ \cong \Unf(X)$ over $\Posf^+$. If $\E$ is
$\kappa$-bounded by an uncountable regular cardinal $\kappa$, then
one choose $X \in \Delta^o\Delta^o\Sets_\kappa$. Moreover, for any
two complete fibrant Segal spaces $X$, $X'$ and enhanced functor
$\gamma:i^*\Unf(X) \to i^*\Unf(X')$, there exists a map $f:X \to X'$
such that $i^*\Unf(f) \cong \gamma$, and for any enhanced functor
$\gamma^+:\Unf(X) \to \Unf(X')$, there exists a map $f:X \to X'$
such that $\Unf(f) \cong \gamma^+$. Finally, two maps $f,f':X \to
X'$ are homotopic if and only if $\Unf(f) \cong \Unf(f')$, and this
happens if and only if $i^*\Unf(f) \cong i^*\Unf(f')$.
\end{theorem}

\proof{} By Proposition~\ref{exp.prop}, for any restricted enhanced
category $\E \to \Posf$, we have $\E \cong \Exp(\C)$ for a complete
restricted Segal family of groupoids $\C \to \Posf^+$, and any
enhanced $\E^+ \to \Posf^+$ equipped with an equivalence $i^*\E^+
\cong \E$ is of the form $\Exp(\C^+)$ for a semicanonical extension
$\C^+ \to \Posf^+$ of the family $\C$. Moreover, for any two
complete restricted Segal families $\C_0,\C_1 \to \Posf$, an
enhanced functor $\gamma:\Exp(\C_0) \to \Exp(\C_1)$ is of the form
$\Exp(\gamma)$ for a functor $\gamma:\C_0 \to \C_1$ over $\Posf$,
the same holds for their semicanonical extensions, and $\Exp(\gamma)
\cong \Exp(\gamma')$ if and only if $\gamma \cong \gamma'$. Now for
first two claims, it remains to observe that Lemma~\ref{cat.repr.le}
actually provides an equivalence $\Unf(X) \cong
\Exp(\Nn^*L^*\Hh(X))$, and apply Proposition~\ref{css.prop},
Proposition~\ref{seg.prop} and Theorem~\ref{coho.bisi.thm}. For the
last claim, note that $f,f':X \to X'$ are homotopic iff $\Hh(f)
\cong \Hh(f')$, and again apply Proposition~\ref{seg.prop}.
\endproof

\begin{corr}\label{enh.bnd.corr}
For any restricted enhanced categories $\E$, $\E'$, isomorphism
classes of enhanced functors $\E \to \E'$ form a set.
\end{corr}

\proof{} Clear. \endproof

\begin{corr}\label{enh.ext.corr}
For any restricted enhanced category $\E$, there exists an enhanced
category $\E^+$ equipped with an equivalence $i^*\E^+ \cong
\E$. Moreover, any such extension $\E^+ \to \Posf^+$ of the
restricted enhanced category $\E \to \Posf$ is small, and for any
two restricted enhanced categories $\E_0$, $\E_1$ with extensions
$\E_0^+$, $\E_1^+$, an enhanced functor $\gamma:\E_0 \to \E_1$
extends to an enhanced functor $\gamma^+:\E_0^+ \to \E_1^+$, unique
up to an isomorphism.
\end{corr}

\proof{} Clear. \endproof

\begin{corr}\label{fin.enh.equi.corr}
An enhanced functor $\gamma:\C \to \C'$ between small enhanced
categories is an equivalence resp.\ fully faithful if its is an
epivalence resp.\ fully faithful over $\Posff \subset \Posf^+$.
\end{corr}

\proof{} For the first claim, by Corollary~\ref{semi.epi.corr}, the
unfolding $\gamma^\dm$ is an equivalence over the category of finite
biordered sets. Thus if we represent $\gamma$ by a map $f:X \to X'$
of fibrant complete Segal spaces, and let $f_n:X_n \to X'_n$ be its
restriction to $[n] \times \Delta \subset \Delta \times \Delta$, for
any $n$, then $f_n$ is an isomorphism on all the homotopy groups,
thus a weak equivalence, and then so is $f$. For the second claim,
replace $\C'$ with the essential image of $\gamma$.
\endproof

\begin{corr}\label{fin.noenh.corr}
An enhanced category $\C$ is $\kappa$-bounded by the countable
cardinal $\kappa$ if and only if $\|\C_\ppt\| < \kappa$ and $\C \cong
\Unf(\C_\ppt)$.
\end{corr}

\proof{} The ``if'' part is clear. For the ``only if'', $\|\C_\ppt\| <
\kappa$ by definition, and also $\|\C_J\| < \kappa$ for any $J \in
\Posff$. Therefore all these categories $\C_J$ are rigid, and then
as in Lemma~\ref{fin.enh.grp.le}, \eqref{nu.J.eq} is an equivalence
for any $J \in \Posff$, and $\C$ is cartesian over standard pushout
squares in $\Posff$, so that \eqref{trunc.eq} is an equivalence over
$\Posff$. Then it is an equivalence by
Corollary~\ref{fin.enh.equi.corr}.
\endproof

\begin{corr}\label{zorn.corr}
Assume given a map $f:J' \to J$ in $\Pos$, with Noetherian $J$, and
a small enhanced category $\C$ such that for any element $j \in J$,
any enhanced functor $\Unf(f^{-1}(J /' j)) \to \C$ extends to an
enhanced functor $\Unf(f^{-1}(J/j)) \to \C$. Then for any
left-closed subset $J_0 \subset J$, any enhanced functor
$\gamma_0:\Unf(f^{-1}(J_0)) \to \C$ extends to an enhanced functor
$\gamma:\Unf(J') \to \C$.
\end{corr}

\proof{} Choose a fibrant complete Segal space $X$ such that
$\Unf(X) \cong \C$ provided by Theorem~\ref{enh.thm}. Represent the
enhanced functor $\gamma_0$ by a map $g_0:\Nn_\dm(f^{-1}(J_0)) \to
X$, and let $U$ be the set of pairs $\langle J_1,g_1 \rangle$ of a
left-closed subset $J_1 \subset J$ containing $J_0$, and a map
$g_1:\Nn_\dm(f^{-1}(J_1)) \to X$ extending $g_0$. Order $U$ by
saying that $\langle J_1,g_1 \rangle \leq \langle J_1',g_1' \rangle$
iff $J_1 \subset J_1'$ and $g_1 =
g_1'|_{\Nn_\dm(f^{-1}(J_1))}$. Then by the Zorn Lemma, $U$ contains
a maximal element $\langle J_1,g_1\rangle$. If $J_1 \subset J$ is
not the whole $J$, then since $J$ is Noetherian, so is $J \ssetminus
J_1$, and we can find a minimal element $j \in J \ssetminus
J_1$. Then $J_1' = J_1 \cup \{j\} \subset J$ fits into a standard
pushout square $J_1' = J_1 \copr_{J/'j} (J/j)$ of
Example~\ref{elem.exa} that induces a corresponding decomposition of
its preimage $f^{-1}(J_1')$, and then by semiexactness, our
assumption implies that the enhanced functor
$\gamma_1:\Unf(f^{-1}(J_1)) \to \C$ represented by $g_1$ extends to
an enhanced functor $\gamma_1':\Unf(f^{-1}(J_1')) \to \C$. By
Lemma~\ref{lift.le} and Proposition~\ref{lift.prop}, $g_1$ then
extends to a map $g_1':\Nn_\dm(f^{-1}(J_1')) \to X$; this
contradicts maximality.  \endproof

\subsection{Categories of enhanced categories.}

By Corollary~\ref{enh.bnd.corr}, restricted enhanced categories and
isomorphism classes of enhanced functors between them form a
well-defined category that we denote by $\Cat^h$. By
Corollary~\ref{enh.ext.corr}, the same category can be defined as
the category of small enhanced categories $\E \to \Posf^+$ and
again, isomorphism classes of enhanced functors between
them. Theorem~\ref{enh.thm} then provides an equivalence
\begin{equation}\label{cat.h.eq}
h^W_{css}(\Delta^o\Delta^o\Sets) \cong \Cat^h, \qquad X \mapsto
\Unf(X),
\end{equation}
where $h^W_{css}(\Delta^o\Delta^o\Sets) \subset
h^W(\Delta^o\Delta^o\Sets)$ is the full subcategory spanned by
complete fibrant Segal spaces. For any small enhanced category $\E$,
we have a well-defined category $\E \bbd^h \Cat^h$ of small
categories $\C$ equipped with an enhanced functor $\alpha:\E \to
\C$, with morphisms given by isomorphism classes of lax functors
$\langle \phi,a \rangle:\langle \C,\alpha \rangle \to \langle
\C',\alpha' \rangle$ under $\E$ such that $\phi:\C \to \C'$ is an
enhanced functor, and $a: \phi \circ \alpha \to \alpha'$ is a
morphism over $\Posf^+$. We also have a category $\Cat^h \bb^h \E$
of small enhanced categories $\C$ equipped with an enhanced functor
$\alpha:\C \to \E$, with morphisms given by isomorphism classes of
lax functors $\langle \phi,a \rangle:\langle \C,\alpha \rangle \to
\langle \C',\alpha' \rangle$ over $\E$ such that $\phi:\C \to \C'$
is an enhanced functor, and $a:\alpha \to \alpha' \circ \phi$ is a
morphism over $\Posf^+$. Inside $\Cat^h \bb^h \E$, we have a dense
subcategory $\Cat^h \bbi^h \E$ defined by the closed class of functors
over $\E$ --- that is, lax functors $\langle \phi,a \rangle$ such
that $a$ is an isomorphism. We have the involution $\iota:\Cat^h \to
\Cat^h$ sending $\C$ to the enhanced-opposite category $\C^\iota$,
we can consider the twisted version $\Cat^h \bb^h_\iota \E =
\iota^*(\Cat^h \bb^h \E)$ of the category $\Cat^h \bb^h \E$, and as
in Example~\ref{iota.bb.exa}, we have $\Cat^h \bbi^h_\iota \E \cong
\Cat^h \bbi^h \E^\iota$. Enhanced groupoids of
Definition~\ref{grp.enh.def} form a full subcategory $\Sets^h
\subset \Cat^h$, and \eqref{cat.h.eq} identifies it with the full
subcategory $h^W(\Delta^o\Sets) \subset h^W(\Delta^o\Delta^o\Sets)$
spanned by objects constant in the Reedy direction; this is the
equivalence $h^W(\Delta^o\Sets) \cong \Sets^h$ of
Corollary~\ref{grp.corr}. Sending a small category $I \in \Cat$ to
$\Unf(I)$ defines a functor
\begin{equation}\label{cat.h.cat}
  \Cat^0 \to \Cat^h,
\end{equation}
and by Proposition~\ref{trunc.prop}, the functor \eqref{cat.h.cat}
is right-adjoint to the forgetful functor $\Cat^h \to \Cat^0$, $\C
\mapsto \C_\ppt$, with the adjuntion map provided by the truncation
functor \eqref{trunc.eq}. In particular, \eqref{cat.h.cat} is a
fully faithful embedding. Note that its source is $\Cat^0$ and not
$\Cat$ -- from the point of view of enhanced category theory, $\Cat$
is simply a wrong category to consider. For any regular cardinal
$\kappa$, we have a full subcategory $\Cat^h_\kappa \subset \Cat^h$
spanned by $\kappa$-bounded enhanced categories, and if $\kappa$ is
uncountable, \eqref{cat.h.eq} restricts to an equivalence
\begin{equation}\label{cat.h.kappa.eq}
h^W_{css}(\Delta^o\Delta^o\Sets_\kappa) \cong \Cat^h_\kappa,
\end{equation}
while for the countable $\kappa$, we have $\Cat^h_\kappa \cong
\Cat_\kappa \cong \Cat^0_\kappa$ by Corollary~\ref{fin.noenh.corr}.

We also have an augmented version of Theorem~\ref{enh.thm}. As in
Section~\ref{aug.sec}, let $I$ be a $\Hom$-finite directed cellular
Reedy category (for example, $I$ can be a left-bounded partially
ordered set). Then for any complete fibrant $I$-Segal space $X$ in
the sense of Definition~\ref{I.bi.seg.def}, its twisted simplicial
expansion $N_\sigma(X|I)$ of \eqref{tw.seg.eq} is a complete fibrant
Segal space by Lemma~\ref{tw.seg.le}, and if we let $\Unf(X|I) =
\Unf(N_\sigma(X^\iota|I)^\iota)$, then we have an enhanced functor
\begin{equation}\label{enh.N.eq}
  \Unf(X|I) \to \Unf(I) \cong \Unf(\ppt|I)
\end{equation}
induced by the map $X \to \ppt$ to the constant functor $\ppt:I^o
\times \Delta^o \times \Delta^o \to \Sets$ with value $\ppt$.

\begin{theorem}\label{enh.I.thm}
Assume given a directed $\Hom$-finite cellular Reedy category $I$,
and let $i:\Unf(I,\Posf) \to \Unf(I) = \Unf(I,\Posf^+)$ be the
natural full embedding. Then for any complete fibrant $I$-Segal
space $X$, the functor \eqref{enh.N.eq} is an $I$-augmentation for
the enhanced category $\Unf(X|I)$. For any $I$-augmented restricted
enhanced category $\E$, there exists a complete fibrant $I$-Segal
space $X$ and an equivalence $\E \cong i^*\Unf(X|I)$ over
$\Unf(I,\Posf)$, and for any such $X$ and any $I$-augmented enhanced
category $\E^+$, an equivalence $i^*\E^+ \cong \E$ over
$\Unf(I,\Posf)$ comes from an equivalence $\E^+ \cong \Unf(X|I)$
over $\Unf(I)$. If $\E$ is $\kappa$-bounded by an uncountable
regular cardinal $\kappa > |I|$, then one choose $X \in
I^o\Delta^o\Delta^o\Sets_\kappa$. Moreover, for any two complete
fibrant $I$-Segal spaces $X$, $X'$ and $I$-augmented enhanced
functor $\gamma:i^*\Unf(X|I) \to i^*\Unf(X'|I)$, there exists a map
$f:X \to X'$ such that $i^*\Unf(f) \cong \gamma$, and for any
$I$-augmented enhanced functor $\gamma^+:\Unf(X|I) \to \Unf(X'|I)$,
there exists a map $f:X \to X'$ such that $\Unf(f) \cong
\gamma^+$. Finally, two maps $f,f':X \to X'$ are homotopic if and
only if $\Unf(f) \cong \Unf(f')$, and this happens if and only if
$i^*\Unf(f) \cong i^*\Unf(f')$.
\end{theorem}

\proof{} The first claim is Proposition~\ref{seg.I.prop}~\thetag{i}
together with Lemma~\ref{seg.bc.le} and Lemma~\ref{exp.I.le}. For
the second claim, use Proposition~\ref{seg.I.prop}~\thetag{ii},
again combined with Lemma~\ref{seg.bc.le} and
Lemma~\ref{exp.I.le}. The other claims then follow from
Lemma~\ref{exp.I.le} and Proposition~\ref{seg.I.prop}.
\endproof

Just as Theorem~\ref{enh.thm}, Theorem~\ref{enh.I.thm} has
corollaries, namely, the obvious $I$-augmented counterparts of
Corollary~\ref{enh.bnd.corr} and Corollary~\ref{enh.ext.corr}. As a
consequence, $I$-augmented restricted enhanced categories and
isomorphism classes of $I$-augmented functors between them form a
well-defined category $\Cat^h(I)$, the same category can be
described as the category of small $I$-augmented enhanced categories
$\E \to \Unf(I)$ and isomorphism classes of $I$-augmented enhanced
functors, and we have an equivalence
\begin{equation}\label{cat.h.I.eq}
h^W_{css}(I^o\Delta^o\Delta^o\Sets) \cong \Cat^h(I), \qquad X
\mapsto \Unf(X|I),
\end{equation}
where as in \eqref{cat.h.eq}, the left-hand side stands for the full
subcategory in $h^W(I^o\Delta^o\Delta^o\Sets)$ formed by complete
fibrant $I$-Segal spaces. Inside $\Cat^h(I)$, we have the full
subcategory $\Sets^h(I)$ spanned by $I$-augmented enhanced
categories $\C \to \Unf(I)$ that are families of groupoids, and
\eqref{cat.h.I.eq} identifies it with $h^W(I^o\Delta^o\Sets)$. For
any map $f:I' \to I$ in $\Posf^+$, we have an adjoint pair of
functors
\begin{equation}\label{cat.h.pb.eq}
f^{h*}:\Cat^h(I) \to \Cat^h(I'), \quad f^h_*:\Cat^h(I') \to
\Cat^h(I),
\end{equation}
where $f^{h*}$ sends an $I$-augmented small enhanced category
$\C \to \Unf(I)$ to the $I'$-augmented small enhanced category
$\Unf(f)^*\C \to \Unf(I')$, and its right-adjoint $f^h_*$ is
provided by Lemma~\ref{pf.cat.le}. For any $I$-augmented small enhanced
category $\C$, the augmentation $\pi:\C \to \Unf(I)$ is in
particular a morphism in $\Cat^h$, and treating it as such gives a
functor
\begin{equation}\label{aug.comma.eq}
\Cat^h(I) \to \Cat^h / [\Unf(I)],
\end{equation}
where as per our general convention, $[\Unf(I)]$ in the right-hand
side stands for the object in $\Cat^h$ given by $\Unf(I)$. By
Remark~\ref{I.rig.rem}, \eqref{aug.comma.eq} is a faithful functor.

\begin{corr}\label{zorn.I.corr}
Assume given a partially ordered set $J'$ equipped with a fibration
$J' \to I$ to some $I \in \Posf$, and a map $f:J' \to J$ to a
Noetherian partially ordered set $J$ that inverts maps cartesian
over $I$. Moreover, assume given a small $I$-augmented enhanced
category $\C$ such that for any $j \in J$, any $I$-augmented
enhanced functor $\Unf(f^{-1}(J/'j)) \to \C$ extends to an
$I$-augmented enhanced functor $\Unf(f^{-1}(J/j)) \to \C$. Then for
any left-closed subset $J_0 \subset J$, any $I$-augmented enhanced
functor $\gamma_0:\Unf(f^{-1}(J_0)) \to \C$ extends to an
$I$-augmented enhanced functor $\gamma:\Unf(J') \to \C$.
\end{corr}

\proof{} Same as Corollary~\ref{zorn.corr}, with
Theorem~\ref{enh.thm} replaced by Theorem~\ref{enh.I.thm}.
\endproof

\subsection{Universal objects.}\label{univ.subs}

While to define small enhanced categories, it suffices to consider
families of categories over $\Posf$, it is quite helpful to consider
their semicanonical extensions to $\Posf^+$ provided by
Corollary~\ref{enh.ext.corr}. The reason for this is that, roughly
speaking, it is these extensions that have the correct
$2$-categorical structure. This is formalized as follows.

\begin{defn}\label{enh.univ.def}
For any enhanced category $\pi:\E \to \Posf^+$ with unfolding
$\pi:\E^\dm \to \Relf^+$, a {\em universal object} for $\E$ is an
object $e \in \E^\dm$ such that for any enhanced category $\E'$ and
an object $e' \in {\E'}^\dm_{\pi^\dm(e)}$, there exists an enhanced
functor $\gamma:\E \to \E'$ equipped with an isomorphism $\alpha:e'
\cong \gamma^\dm(e)$, and for any two enhanced functors
$\gamma_0,\gamma_1:\E \to \E'$, a map $\gamma_0^\dm(e) \to
\gamma_1^\dm(e)$ uniquely extends to a map $\gamma_0 \to \gamma_1$.
\end{defn}

By definition, for any complete fibrant Segal space $X$, with the
corresponding small enhanced category $\pi:\E = \Unf(X) \to
\Posf^+$, an object $e \in \E^\dm$ is given by a pair $\langle J,a
\rangle$ of a biordered set $I = \pi^\dm(e) \in \Relf^+$ and a map
$a:\Nn_\dm(I) \to X$. We then have the following criterion for
universality.

\begin{lemma}\label{univ.enh.le}
For any complete fibrant Segal space $X$, with the corresponding
small enhanced category $\E = \Unf(X)$, an object $e \in \E^\dm$ is
universal in the sense of Definition~\ref{enh.univ.def} if and only
if the corresponding map $a:\Nn_\dm(I) \to X$ is a Segal equivalence
in the sense of Example~\ref{seg.eq.exa}.
\end{lemma}

\proof{} The ``only if'' part is clear: since the pair $\langle
\gamma,a \rangle$ in Definition~\ref{enh.univ.def} is unique up to a
unique isomorphism, universality for small enhanced categories $\E'
= \Unf(X')$ together with Theorem~\ref{enh.thm} immediately implies
that $a$ is a Segal equivalence. Conversely, assume that $a$ is a
Segal equivalence. Note that since $\E$ in
Definition~\ref{enh.univ.def} is small, it suffices to check the
condition for small $\E'$ --- indeed, for any object $e' \in
{\E'}^\dm_{\pi(e)}$, we may replace $\E'$ with the full subcategory
spanned by $j^*e' \in \E'_{\ppt}$ for all maps $j:\ppt \to
\pi(e)$. Then a pair $\langle \gamma,\alpha \rangle$ for any $\E' =
\Unf(X')$ is provided by Theorem~\ref{enh.thm}. Moreover, a map
$f:\gamma_0^\dm(e) \to \gamma_1^\dm(e)$ defines an object $\wt{f}
\in \Ar_h(\E')^\dm_I$, where $\Ar_h(\E')$ is the enhanced arrow
category provided by Corollary~\ref{J.C.corr}, and the same argument
applied to $\Ar_h(\E')$ instead of $\E'$ then shows that $f$ extends
to a map of functors $f:\gamma_0^\dm \to \gamma_1^\dm$. To finish
the proof, it remains to show uniqueness: we assume given enhanced
functors $\gamma_0,\gamma_1:\E \to \E'$ and two maps
$f,f':\gamma_0^\dm \to \gamma_1^\dm$ such that $f(e) = f'(e)$, and
we need to show that $f(e')=f'(e')$ for any object $e' \in
\E^\dm$. However, any such $e'$ is represented by a map $\langle
I',a' \rangle$, $I' \in \Relf^+$, $a':\Nn_\dm(I') \to X$, so that we
have a diagram
\begin{equation}\label{ee.dia}
\begin{CD}
I' \longleftarrow \QQ_\dm(\Nn_\dm(I')) @>{\QQ_\dm(a')}>> \QQ_\dm(X)
@<{\QQ_\dm(a)}<< \QQ_\dm(\Nn_\dm(I)) \longrightarrow I,
\end{CD}
\end{equation}
in $\Relf^+$, and a cartesian section of the fibration $\E^\dm \to
\Relf^+$ over \eqref{ee.dia} that evaluates to $e'$ resp.\ $e$ at
$I'$ resp.\ $I$. Since both $\gamma_0^\dm$ and $\gamma_1^\dm$ send
this section to a cartesian section of the fibration ${\E'}^\dm \to
\Relf^+$, and this fibration is constant along all maps in
\eqref{ee.dia} except for possibly $\QQ_\dm(a')$, we are done by
Lemma~\ref{map.le}.
\endproof

\begin{exa}
For any partially ordered set $I$, an object $e \in \Unf(I)^\dm$ is
the same thing as a biordered set $J \in \Relf^+$ equipped with an
augmentation $\alpha:J \to L(I^o)$. Then in particular, the functor
$\xi^o_\perp:B_\dm(I)^o \to I$ defines an object $e \in
\Unf(I)^\dm_{B_\dm(I)}$, and by Proposition~\ref{lf.prop} and
\eqref{ev.j.I}, this object is universal in the sense of
Definition~\ref{enh.univ.def}. If $I \in \Posf^-$, then $J=L(I^o)$
with the tautological map $\id:L(I^o) \to L(I^o)$ is also universal.
\end{exa}

\begin{prop}\label{univ.prop}
Assume given a partially ordered set $I \in \Posf^+$ and a small
$I$-augmented enhanced category $\pi:\C \to \Unf(I)$, with the
opposite enhanced category $\C^\iota$. Then there exists a biordered
set $J \in \Relf^+$ and an object $c \in \C^{\iota\dm}_J \subset
\C^{\iota\dm}$ in the unfolding $\C^{\iota\dm}$ such that \thetag{i}
the map $J \to L(I)$ corresponding to $\pi^\iota(c) \in \Unf(I^o)_J$
is a biordered fibration, \thetag{ii} the functor $\Unf(U(J)) \to
\C$ corresponding to $c \in \C^{\iota\dm}_J \subset \C^\iota_{U(J)}$
by \eqref{ev.j.I} is $I$-augmented, and \thetag{iii} for any map
$f:I' \to I$ in $\Posf^+$, the object $\Unf(f^o)^*c \in
\Unf(f^o)^*\C^\iota$ is universal in the sense of
Definition~\ref{enh.univ.def}. Moreover, if $\C$ is $\kappa$-bounded
by a regular cardinal $\kappa > |I|$, then one can choose $J \in
\Relf^+_\kappa$, and if $I = [1]$ and the fiber $\C_1$ is
$\kappa$-bounded, then one can choose $J$ such that $|J_1| <
\kappa$.
\end{prop}

\proof{} By Lemma~\ref{univ.enh.le}, for any $J^\dm \in \Relf^+$ and
complete fibrant Segal space $X$, a map $a:\Nn_\dm(J^\dm) \to X$
defines an object in $\Unf(X)^{\iota\dm} \cong \Unf(X^\iota)^\dm$,
and this object is universal if and only if $a$ is a Segal
equivalence. Take a fibrant $I$-Segal space $X$ and an equivalence
$\C \cong \Unf(X|I) = \Unf(N_\sigma(X^\iota|I)^\iota)$ provided by
Theorem~\ref{enh.I.thm}. Consider the $I$-augmented biordered set
$J=\QQ_\dm(I)(X^\iota)$ of \eqref{bi.Q.N.I}, and let $J^\hdot = I
\setminus J$ with the biorder of Example~\ref{facto.fib.exa}, so
that $J^\hdot \to L(I)$ is a biordered fibration corresponding to a
functor $I^o \to \bRelf^+$. Then Proposition~\ref{seg.I.prop}
provides a weak equivalence $w:\Nn_\dm(I)(J) \to X^\iota$ of
$I$-Segal spaces, and $w$ induces a map
$a=N_\sigma(w):\Nn_\dm(J^\hdot) \to N_\sigma(X^\iota|I)$ of
\eqref{N.w} that is a Segal equivalence by
Proposition~\ref{tw.seg.prop}. Moreover, by the same Proposition,
$f^*(a)$ is also a Segal equivalence for any map $f:I' \to I$ in
$\Posf^+$. To finish the proof, it remains to convert $J^\hdot$ to a
left-bounded biordered set. To do this, use the same procedure as in
\eqref{Q.pl}: filter $X^\iota$ by the skeleton filtration $\sk_nX$,
with the induced exhaustive filtration on $J_\idot$, and let
$J^{\hdot +} \in \Relf^+$ be the corresponding incidence subset of
Example~\ref{N.ano.exa}. Then the map $\Nn_\dm(J^{\hdot +}) \to
\Nn_\dm(J^\hdot)$ is a weak equivalence, so that its composition
$a^+:\Nn_\dm(J^{\hdot +}) \to N_\sigma(X^\iota|I)$ with $a$ is a Segal
equivalence, and the same goes for $f^*(a^+)$ for any $f:I' \to I$.

In the $\kappa$-bounded case, if $\kappa$ is uncountable, then one
can choose $X$ to lie in $I^o\Delta^o\Delta^o\Sets_\kappa$ by
Theorem~\ref{enh.I.thm}, and then $\QQ_\dm(I)(X^\iota)$ is in
$\Relf^+_\kappa$, and so is $J^{\hdot +}$. If $\kappa$ is the
smallest infinite cardinal, then $I$ is finite, and $\C \cong
\Unf(\C_\ppt)$, while $\C_\ppt$ is finite by
Corollary~\ref{fin.noenh.corr}. We can then just take $X =
L(N(\C_\ppt))$ constant in the simplicial direction; it is a finite
$I$-bisimplicial set, $\QQ_\dm(I)(X^\iota)$ is finite, and there it
no need to pass to $J^{\hdot +}$.

Finally, if $I=[1]$ and only $\C_1$ is $\kappa$-bounded, we can
still choose a fibrant $[1]$-Segal space $X$ with an equivalence $\C
\cong \Unf(X|[1])$, and then choose a fibrant Segal space $X_1' \in
\Delta^o\Delta^o\Sets_\kappa$ with an equivalence $\C_1 \cong
\Unf(X_1')$, or take $X_1'=L(N(\C_{1,\ppt}))$ if $\kappa$ is the
smallest infinite cardinal. Then since $t:[0] \to [1]$ is
left-reflexive, $X_1 = t^*X$ is also a fibrant Segal space, and
since $\Unf(X_1) \cong \Unf(X_1') \cong \C_1$, there exists a weak
equivalence $w:X_1' \to X_1$. We can then replace $X:[1]^o \to
\Delta^o\Delta^o\Sets$ with $X':[1]^o \to \Delta^o\Delta^o\Sets$
corresponding to the map $f:X_1' \to X_1 \to X_0$, and $w$ extends
to a weak equivalence $X' \to X$. Now take
$J=\QQ_\dm([1])({X'}^\iota)$, and let $J^\hdot = [1] \setminus J$
with the biorder of Example~\ref{R.groth.exa}. Then we still have a
weak equivalence $\Nn_\dm([1])(J) \to X' \to X$ of $[1]$-Segal
spaces, so Proposition~\ref{tw.seg.prop} produces a Segal
equivalence $\Nn_\dm(J^\hdot) \to N_\sigma(X^\iota|[1])$. On the
other hand, Lemma~\ref{reedy.copr.le} applies to the projection $[1]
\times \Delta \times \Delta \to [1]$ shows that $1 \setminus J_\idot
\cong J_1 = \QQ_\dm({X'}^\iota_1)$, so if $\kappa$ is uncountable,
it remains to convert $J^\hdot$ to $J^{\hdot +} \in \Relf^+$ by taking
the incidence set for the skeleton filtration $\sk_n X'$. If
$\kappa$ is countable, so that $X_1'$ is finite, it is enough to do
this for the fiber over $0 \in [1]$, by filtering $X_0$ by subsets
$\sk_n(X_0) \cup f(X'_1) \subset X_0$.
\endproof

\begin{corr}\label{fun.h.corr}
For any enhanced categories $\E$, $\E'$ such that $\E$ is small,
morphisms between any two enhanced functors $\E \to \E'$ form a set,
so that these enhanced functors form a well-defined category
$\Fun^h(\E,\E')$. Moreover, there exist an enhanced functor category
$\fFun^h(\E,\E')$ in the sense of Definition~\ref{enh.fun.def}, and
$\fFun^h(\E,\E')_\ppt \cong \Fun^h(\E,\E')$, so that
$\fFun^h(\E,\E')$ is an enhancement for $\Fun^h(\E,\E')$.
\end{corr}

\proof{} For the first claim, choose a universal object $e \in
\E^\dm_J \subset \E^\dm$, $J \in \Relf^+$ provided by
Proposition~\ref{univ.prop}, so that an enhanced functor $\gamma:\E
\to \E'$ is uniquely determined by the object $\gamma^\dm(e) \in
\E'_J$ and morphisms $\gamma \to \gamma'$ correspond bijectively to
morphisms $\gamma^\dm(e) \to {\gamma'}^\dm(e)$, and note that these
of course form a set. Thus we indeed have a well-defined category
$\Fun^h(\E,\E')$, and moreover, evaluation at $e$ provides an
equivalence
\begin{equation}\label{ev.uni}
\ev_e:\Fun^h(\E,\E') \cong {\E'}^\dm_J.
\end{equation}
Now we can just let $\fFun^h(\E,\E') \subset U(J)^o_h\E'$ be the
full enhanced subcategory corresponding to ${\E'}^\dm_J \subset
\E'_{U(J)} \cong (U(J)^o_h\E')_\ppt$; the evaluation functor and its
universal property are provided by Corollary~\ref{J.C.corr}, and
\eqref{ev.uni} identifies $\fFun^h(\E,\E')_\ppt \cong \Fun^h(\E,\E')$.
\endproof

As in the unenhanced case, we will sometimes simplify notation and
write $\E^\iota\E' = \fFun^h(\E^\iota,\E')$, and
$I^o_h\E=\Unf(I)^\iota\E$ for an essentially small category $I$ (if
$I=J$ is a partially ordered set, this is consistent with
Corollary~\ref{J.C.corr}). If $\E$ and $\E'$ are equipped with
enhanced functors $\E,\E' \to \Unf(I)$ for some $I \in \Posf^+$ ---
for example, if both are $I$-augmented or $I$-coaugmented --- we
will denote by $\Fun^h_I(\E,\E') \subset \Fun^h(\E,\E')$ the full
subcategory spanned by functors over $\Unf(I)$ (we recall that since
$I$ is rigid, being a functor over $\Unf(I)$ is a condition and not
a structure). As a useful application of Corollary~\ref{fun.h.corr},
we note that if we have two enhanced functors $\gamma,\gamma':\E \to
\E'$ with $\E$ small, then a morphism $f:\gamma \to \gamma'$ gives
rise to an enhanced functor
\begin{equation}\label{ar.enh.eq}
  f^h:\E \times^h \Unf([1]) \to \E'
\end{equation}
equipped with isomorphisms $f^h \circ s \cong \gamma$, $f^h \circ t
\cong \gamma'$. Indeed, $\gamma$ and $\gamma'$ define enhanced
objects in $\fFun^h(\E,\E')$, and to obtain \eqref{ar.enh.eq}, it
suffices to lift $f$ to an enhanced morphism between these two
enhanced objects.

\begin{corr}\label{sq.h.corr}
Assume given a semicartesian square \eqref{refl.semi.sq} of enhanced
categories and enhanced functors, and another commutative square
\begin{equation}\label{enh.h.sq}
\begin{CD}
\C'_{01} @>{\gamma'_0}>> \C_0\\
@V{\gamma'_1}VV @VVV\\
\C_1 @>>> \C
\end{CD}
\end{equation}
of enhanced categories and enhanced functors such that $\C'_{01}$ is
small. Then there exists an enhanced functor $\alpha:\C'_{01} \to
\C_{01}$ equipped with isomorphisms $a_l:\gamma_l \circ \alpha \cong
\gamma'_l$, $l=0,1$, and the triple $\langle \alpha,a_0,a_1 \rangle$
is unique up to a (non-unique) isomorphism.
\end{corr}

\proof{} By Lemma~\ref{refl.sq.le}, the unfolding of the square
\eqref{refl.semi.sq} is also semicartesian. Evaluate this unfolding
at $\pi(c)$ for a universal $c \in {\C'}^\dm_{01}$.
\endproof

In particular, Corollary~\ref{sq.h.corr} implies that any
semicartesian square \eqref{enh.h.sq} with small $\C_{01}$ is unique
up to a unique equivalence, if it exists; if it does, we will call
$\C_{01}$ the {\em semicartesian product} of $\C_0$ and $\C_1$ over
$\C$, and denote it by $\C_{01} = \C_0 \times^h_{\C} \C_1$. If
$\C_0$, $\C_1$ and $\C$ are themselves small --- or at least, if so
are $\C_0$ and the functor $\C_1 \to \C$ --- existence immediately
follows from Theorem~\ref{enh.thm}.

\begin{lemma}\label{sq.h.le}
For any enhanced categories $\C$, $\C_0$, $\C_1$ and enhanced
functors $\C_0,\C_1 \to \C$ such that $\C_0$ and $\C_1 \to \C$ are
small, there exists a small semicartesian product $\C_0
\times^h_{\C} \C_1$, and if $\C$, $\C_0$, $\C_1$ are $I$-augmented
for some $I \in \Posf^+$, then so is $\C_0 \times^h_{\C} \C_1$.
\end{lemma}

\proof{} If the functor $\C_0 \to \C$ is fully faithful, then
already the cartesian product $\C_0 \times_{\C} \C_1$ is an enhanced
category by Corollary~\ref{enh.ff.corr} and
Lemma~\ref{sq.cat.le}. In the general case, we may replace $\C$ with
a small full enhanced subcategory that contains the essential image
of the functor $\C_0 \to \C$, so it suffice to consider the case
when $\C$, hence also $\C_1$ is small. Choose complete Segal spaces
$X$, $X_0$ and $X_1$ equipped with equivalences $\C \cong \Unf(X)$,
$\C_0 \cong \Unf(X_1)$, $\C_1 \cong \Unf(X_1)$ provided by
Theorem~\ref{enh.thm}, represent functors $\C_0,\C_1 \to \C$ by maps
$X_0,X_1 \to X$, take replacements if necessary to insure that both
maps are in $F$, and let $\C_{01} = \Unf(X_0 \times_X X_1)$. By
Corollary~\ref{semi.corr} and Lemma~\ref{w.2.le}, this gives a
semicartesian square. In the augmented case, use
Theorem~\ref{enh.I.thm} instead of Theorem~\ref{enh.thm}.
\endproof

\begin{remark}
More generally, in the situation of Lemma~\ref{sq.h.le}, the
semicartesian product exists if $\C_1 \to \C$ is small, and $\C_0 =
\copr_s\C_s$ is a (possibly large) disjoint union of small enhanced
categories $\C_s$. Indeed, it is given by $\C_0 \times^h_{\C} \C_1 =
\copr_s \C_s \times^h_{\C} \C_1$. In particular, the product always
exists if $\C_0=\Unf(\E)$ for some discrete category $\E$ (that can
be large).
\end{remark}

\subsection{Immediate applications.}\label{univ.app.subs}

As the first application of semicartesian products provided by
Lemma~\ref{sq.h.le}, we can define enhanced categories of sections
generalizing the enhanced categories $\sSec^h(I^o,\C)$ of
\eqref{ssec.enh.I}. Namely, by the universal property of the
enhanced functor categories of Corollary~\ref{fun.h.corr}, for any
enhanced categories $\C$, $\E$ with small $\C$, an enhanced functor
$\gamma:\C' \to \C$ from a small $\C'$ induces a pullback enhanced
functor
\begin{equation}\label{pb.enh}
  \gamma^*:\fFun^h(\C,\E) \to \fFun(\C',\E),
\end{equation}
and an enhanced functor $\E \to \E'$ to some $\E'$ induces a
postcomposition enhanced functor $\fFun^h(\C,\E) \to
\fFun^h(\C,\E')$. Then for any enhanced functor $\E \to \C$ between
small enhanced categories, we can define an enhanced category
$\sSec^h(\C,\E)$ by the semicartesian square
\begin{equation}\label{ssec.enh.C}
\begin{CD}
  \sSec^h(\C,\E) @>>> \fFun^h(\C,\E)\\
  @VVV @VVV\\
  \ppt^h @>>> \fFun^h(\C,\C),
\end{CD}
\end{equation}
where as in \eqref{ssec.enh.I}, the bottom arrow corresponds to
$\id:\C \to \C$. More generally, if we have another small enhanced
category $\E'$ equipped with an enhanced functor $\E' \to \C$, we
can define the enhanced category $\fFun^h_{\C}(\E',\E)$ by the
semicartesian square
\begin{equation}\label{fun.C.enh.sq}
\begin{CD}
\fFun^h_{\C}(\E',\E) @>>> \fFun^h(\E',\E)\\
@VVV @VVV\\
\ppt^h @>>> \fFun^h(\E',\C),
\end{CD}
\end{equation}
where the bottom arrow corresponds to the given enhanced functor
$\E' \to \C$. Note that we actually have $\fFun^h_{\C}(\E',\E) \cong
\sSec^h(\E',\E' \times^h_{\C} \E)$, and by abuse of notation, we
will sometimes write $\sSec^h(\E',\E) = \sSec^h(\E',\E'
\times^h_{\C} \E)$. Just as $\fFun^h(-,-)$, the enhanced category
$\fFun^h_{\C}(-,-)$ can be characterized by a universal property:
for any small enhanced category $\C'$, we have a natural isomorphism
\begin{equation}\label{ffun.C.univ}
\Cat^h(\C',\fFun^h_{\C}(\E',\E)) \cong (\Cat^h \bbi^h \C)(\C' \times^h
\E',\E),
\end{equation}
where in the right-hand side, we have the set of morphisms in the
category $\Cat^h \bbi^h \C$ (that is, the set of isomorphism classes
of functors over $\C$). Moreover, as in
Subsection~\ref{rel.fun.subs}, postcomposition with an arbitrary
enhanced functor $\gamma:\C' \to \C$ from some small enhanced
category $\C'$ provides a functor $\gamma_\trr:\Cat^h \bbi^h \C' \to
\Cat^h \bbi^h \C$, and Corollary~\ref{sq.h.corr} together with
Lemma~\ref{sq.h.le} show that $\gamma_\trr$ admits a right-adjoint
functor $\gamma^*:\Cat^h \bbi^h \C \to \Cat^h \bbi^h \C'$ sending
$\E$ to $\C' \times^h_{\C} \E$. Then by the universal property
\eqref{ffun.C.univ}, this lifts to an enhancded functor
\begin{equation}\label{gamma.st.fun}
\gamma^*:\fFun^h_{\C}(\E',\E) \to
\fFun^h_{\C'}(\gamma^*\E',\gamma^*\E),
\end{equation}
where $\E$, $\E'$ are arbitrary small enhanced categories equipped
with enhanced functors to $\C$.

\begin{remark}\label{fun.C.rem}
While up to an isomorphism, enhanced objects of the enhanced
category $\fFun^h_{\C}(\E',\E)$ are enhanced functors $\E' \to \E$
over $\C$, already for morphisms this is usually not true. Namely,
define $\Fun^h_{\C}(\E',\E)$ as $\fFun^h_{\C}(\E',\E)_\ppt$; then
there is a natural functor $\Fun^h_{\C}(\E',\E) \to
\Fun_{\C}(\E',\E)$, but it is usually {\em not} faithful. What
happens is, a morphism in $\Fun_{\C}(\E',\E)$ is a morphism between
functors that satifies a certain commutativity condition, and this
condition becomes a structure in the enhanced setting.
\end{remark}

We can also use semicartesian products to define enhanced fibers of
small enhanced functors. Namely, for any small enhanced functor
$\pi:\C \to \E$, and any enhanced object $e \in \E_\ppt$, the {\em
  enhanced fiber} $\C_e$ is given by $\C_e = \ppt^h \times^h_{\E}
\C$, where $\ppt^h \to \E$ is the enhanced functor $\eps^h(e)$ (if
$\E = \Unf(I)$ for some $I \in \Posf^+$, then this is consistent
with our earlier notation \eqref{enh.aug.fib}). More generally, for
any $e \in \E_J \subset \E$, $J \in \Posf^+$, we define the {\em
  generalized enhanced fiber} $\C_e$ by the semicartesian square
\begin{equation}\label{enh.fib.sq}
\begin{CD}
\C_e @>>> \C\\
@VVV @VVV\\
\Unf(J^o) @>>> \E,
\end{CD}
\end{equation}
where the bottom arrow corresponds to $e$. By
Corollary~\ref{sq.h.corr}, an enhanced functor $\gamma:\C' \to \C$
over $\E$ then defines an enhanced functor $\gamma_e:\C'_e \to
\C_e$.

\begin{lemma}\label{enh.ff.loc.le}
For enhanced categories $\C$, $\C'$ equipped with small enhanced
functors to an enhanced category $\E$, an enhanced functor
$\gamma:\C' \to \C$ over $\E$ is fully faithful resp.\ an
equivalence if and only if for any $e \in \E$, this holds for the
enhanced functor $\gamma_e:\C'_e \to \C_e$.
\end{lemma}

\proof{} The ``only if'' part immediately follows from
Corollary~\ref{sq.h.corr}. For the ``if'' part, we need to check
that $\gamma$ induces an equivalence between $\C'$ and the enhanced
essential image $\C'' \subset \C$ of $\gamma$. Since for any $e \in
\E$, the generalized fiber $\C''_e$ is the essential image of
$\gamma_e$, we may replace $\C$ with $\C''$ and assume that
$\gamma_e$ is an equivalence for any $e$. Then by
\eqref{enh.fib.sq}, we have a commutative diagram
\begin{equation}\label{enh.epi.sq}
\begin{CD}
\C'_e @>{\gamma_e}>> \C_e\\
@VVV @VVV\\
\Unf(J^o) \times_{\E} \C' @>{\id \times \gamma}>> \Unf(J^o)
\times_{\E} \C,
\end{CD}
\end{equation}
and the top arrow is an equivalence, while both vertical arrows are
epivalences. Therefore the bottm arrow is an epivalence. Since $e$
was arbitrary, Lemma~\ref{epi.loc.le} then implies that $\gamma$ is
an epivalence, and being an enhanced functor, it is an equivalence
by Lemma~\ref{enh.epi.le}.
\endproof

Another observation is that Corollary~\ref{sq.h.corr} actually shows
that in the enhanced setting, semicartesian squares of enhanced
categories and enhanced functors can be characterized by a universal
property analogous to the universal property of cartesian squares in
the usual category theory. If one calls squares satisfying this
property ``enhanced-semicartesian'', then the claim is that
semicartesian squares are enhanced-semicartesian. It is useful
to also introduce the following dual version of this universal
property.

\begin{defn}\label{enh.cocart.def}
A commutative square
\begin{equation}\label{enh.cocart.sq}
\begin{CD}
  \C @>>> \C_0\\
  @VVV @VV{\gamma_0}V\\
  \C_1 @>{\gamma_1}>> \C_{01}
\end{CD}
\end{equation}
of enhanced categories and enhanced functors is {\em
  enhanced-semicocartesian} resp.\ {\em enhanced-cocartesian} if for
any other such square
$$
\begin{CD}
  \C @>>> \C_0\\
  @VVV @VV{\gamma'_0}V\\
  \C_1 @>{\gamma'_1}>> \C'_{01},
\end{CD}
$$
there exists an enhanced functor $\gamma:\C_{01} \to \C'_{01}$
equipped with isomorphisms $\alpha_l:\gamma \circ \gamma_l \cong
\gamma'_l$, $l=0,1$, and the triple $\langle
\gamma,\alpha_0,\alpha_1 \rangle$ is unique up to an isomorphism
resp.\ up to a unique isomorphism.
\end{defn}

\begin{lemma}\label{enh.cocart.le}
A commutative square \eqref{enh.cocart.sq} of small enhanced
categories and enhanced functors is enhanced-semicocartesian
resp.\ enhanced-co\-cart\-e\-sian if and only if for any enhanced
category $\E$, the corresponding commutative square
\begin{equation}\label{enh.fun.sq}
\begin{CD}
  \fFun^h(\C_{01},\E) @>>> \fFun^h(\C_0,\E)\\
  @VVV @VVV\\
  \fFun^h(\C_1,\E) @>>> \fFun^h(\C,\E)
\end{CD}
\end{equation}
is semicartesian resp.\ cartesian. Moreover, it suffices to consider
small enhanced categories $\E$, and if \eqref{enh.fun.sq} is
semicartesian, then the square
\begin{equation}\label{enh.sec.sq}
\begin{CD}
  \sSec^h(\C_{01},\E) @>>> \sSec^h(\C_0,\E)\\
  @VVV @VVV\\
  \sSec^h(\C_1,\E) @>>> \sSec^h(\C,\E)
\end{CD}
\end{equation}
is semicartesian for any small enhanced category $\E$ equipped with
an enhanced functor $\E \to \C_{01}$.
\end{lemma}

\proof{} Immediately follows from Corollary~\ref{fun.h.corr}.
\endproof

\begin{exa}\label{enh.univ.exa}
Assume given a small enhanced category $\E$, and a universal object
$e \in \E^\dm$ represented by a pair $\langle I^\dm,a \rangle$. Then
Definition~\ref{univ.def} immediately implies that $\E$ fits into a
enhanced-cocartesian square
\begin{equation}\label{enh.univ.sq}
\begin{CD}
  W \times \Unf([1]) @>>> \Unf(I^o)\\
  @VVV @VV{a}V\\
  W \times \ppt^h @>>> \E,
\end{CD}
\end{equation}
where $I = U(I^\dm)$ is the underlying partially ordered set of the
biordered set $I^\dm$, and $W$ is the set of maps $R([1]) \to
I^\dm$.
\end{exa}

\begin{exa}\label{cyl.coc.exa}
For any map $f:J_0 \to J_1$ between partially ordered sets, the
cylinder $\Cyl(f)$ and the dual cylinder $\Cyl^o(f)$ fit into
cocartesian squares \eqref{cyl.o.sq}, and then $\Unf(-)$ sends both
into enhanced-semicocartesian squares of small enhanced categories:
indeed, for any enhanced category $\E$, the corresponding square
\eqref{enh.fun.sq} is semicartesian by
Lemma~\ref{bar.enh.le}~\thetag{ii}.
\end{exa}

\subsection{Functoriality and gluing.}\label{univ.funct.subs}

Another fact that follows immediately from Corollary~\ref{sq.h.corr}
is that small semicartesian products provided by Lemma~\ref{sq.h.le}
are functorial. We will need a slightly more general statement. For
any cofibration $f:I \to I'$ in $\Posf^+$, define a full subcategory
$\Unf(I|I') \subset \Unf(I)$ by the cartesian square
\begin{equation}\label{unf.rel.sq}
\begin{CD}
\Unf(I|I') @>>> \Unf(I)\\
@VVV @VV{\Unf(f)}V\\
I' \times \Posf^+ @>>> \Unf(I'),
\end{CD}
\end{equation}
where the bottom arrow is the functor of Example~\ref{cyl.fib.exa}
for the terminal object $\ppt \in \Posf^+$ and the fiber
$\Unf(I')_\ppt \cong I'$. Note that the left vertical arrow in
\eqref{unf.rel.sq} is a cofibration with fibers $\Unf(I|I')_i \cong
\Unf(I_i)$, $i \in I'$. If we have a commutative diagram
\begin{equation}\label{I.pos.sq}
\begin{CD}
I_0 @>{f}>> I_1\\
@VVV @VVV\\
I_0' @>{f'}>> I_1'
\end{CD}
\end{equation}
in $\Posf^+$ such that the vertical arrows are cofibrations, then
$\Unf(f)$ and $\Unf(f')$ give rise to a functor
$\Unf(f|f'):\Unf(I_0|I_0') \to \Unf(I_1|I_1')$, cartesian over
$\Posf^+$.

Now denote by $\Cat^h(I|I')$ the category of small fibrations $\E
\to \Unf(I|I')$ whose restrictions $\E_i \to \Unf(I_i)$, $i \in I'$
are enhanced categories, with morphisms given by isomorphism classes
of functors over $\Unf(I|I')$. If $I' = \ppt$, then $\Cat^h(I|\ppt)
\cong \Cat^h(I)$; if $I = I'$, with $\id:I \to I$, then
$\Cat^h(I|I)$ is the category of families of small enhanced
categories and enhanced functors parametrized by $I$. A commutative
square \eqref{I.pos.sq} induces a functor $(f|f')^{h*}:\Cat^h(I_1|I_1')
\to \Cat^h(I_0|I_0')$ sending $\C$ to $\Unf(f|f')^*\C$, and if $f'$
is an isomorphism, the same argument as in Lemma~\ref{pf.cat.le}
shows that $(f|f')^{h*}$ has a right-adjoint functor $(f|f')^h_*$.

\begin{lemma}\label{ccat.adj.le}
For any $I,I' \in \Posf^+$ such that $\dim I' \leq 1$, with the
projection $p:I' \to \ppt$, the functor $(\id \times
p|p)^{h*}:\Cat^h(I|\ppt) \to \Cat^h(I \times I'|I')$ admits a
right-adjoint $(\id \times p|p)^h_*$.
\end{lemma}

\proof{} Assume first that $I= \ppt$.  Consider the dual barycentric
subdivision $B(I')^o$, with the map $\xi^o_\perp:B(I')^o \to
I'$. Then by \eqref{2kan.eq}, for any fibration $\C \to I'$, the
functor $\C \to \xi^o_{\perp *}\xi^{o*}_\perp\C$ is an
equivalence. Therefore $\xi^{o*}_\perp:\Cat^h(I'|I') \to
\Cat^h(B(I')^o|B(I')^o)$ is fully faithful, and by
Example~\ref{cl.adj.exa}, it suffices to prove the statement for
$B^o(I')$ instead of $I'$. Now consider the characteristic map
$\chi:I' \to [1]$ of the $0$-skeleton $\sk_0I'$, and let $q =
B(\chi)^o:B^o(I') \to B([1])^o \cong \V$ be the corresponding
cofibration. Then again by \eqref{2kan.eq}, for any $\E \in
\Cat^h(B(I')^o|B(I')^o)$, the fibers of the fibration $q_*\C$ are
products of enhanced categories, thus enhanced categories
themselves. Therefore $(q|q)^{h*}:\Cat^h(\V|\V) \to
\Cat^h(B(I')^o|B(I')^o)$ admits a right-adjoint $(q|q)^h_*$ induced
by $q_*$, and we are reduced to the case $I' = \V$. This is
Corollary~\ref{sq.h.corr}.

For a general $I$, we have the faithful embedding
\eqref{aug.comma.eq}, so it suffices to check that $(p|p)^h_*$ sends
augmentationg to augmentations and augmented maps to augmented map;
this immediately follows from the augmented part of
Lemma~\ref{sq.h.le}.
\endproof

One situation where semicartesian products come in handy is when one
wants to assemble an $I$-augmented small enhanced category $\C$ out
of smaller pieces. Namely, assume given a right-closed embedding
$r:I_1 \to I$ in $\Posf^+$, with the characteristic map $I \to [1]$,
and let $I' = [1] \setminus I \cong \Cyl^o(r)$, with the
corresponding cocartesian square \eqref{cyl.o.sq} interpreted as a
cofibration $I^\hdot \to [1]^2 \cong \V^>$ with fibers $I^\hdot_o =
I_1$, $I^\hdot_0 = I$, $I^\hdot_1 = I_1 \times [1]$ over $o,0,1 \in
\V \subset [1]$. Let $\eps:I^\hdot \to I' \times [1]^2$ be the full
embedding cocartesian over $[1]^2$.

\begin{lemma}\label{co.I.le}
For any $I'$-augmented enhanced category $\C$, the square
\begin{equation}\label{co.I.sq}
\begin{CD}
\C @>>> \eps_{0h*}\eps^{h*}_0\C\\
@VVV @VVV\\
\eps_{1h*}\eps^{h*}_1\C @>>> \eps_{oh*}\eps^{h*}_o\C
\end{CD}
\end{equation}
of $I'$-augmented enhanced categories is semicartesian.
\end{lemma}

\proof{} For any map $\eps:I' \to I$ in $\Posf^+$, $I$-augmented
enhanced category $\C$, and an object $\langle J,\alpha \rangle \in
\Unf(I)$, with the bicofibration $J_\idot = J /^\dm R(I^o) \to
R(I^o)$, \eqref{Unf.g.C} provides an equivalence
\begin{equation}\label{pf.pb.eq}
(\eps^h_*\eps^{h*}\C)_{\langle J,\alpha \rangle} \cong
  (\C^\Tot|I)^\dm_{\langle R(\eps^o)^*J_\idot,\alpha \circ \eps^o
    \rangle},
\end{equation}
where $\C^\Tot \to \Pos$ is the canonical bar-invariant extension of
Proposition~\ref{lf.prop}. Therefore in our situation, we need to
check that for any $\langle J,\alpha \rangle \in \Unf(I')$, with the
corresponding bicofibration $J_\idot \to R({I'}^o)$, the relative
unfolding $(\C^\Tot|I)^\dm$ is semicartesian along the square
\begin{equation}\label{J.pb.sq}
\begin{CD}
  \eps_o^{o*}J_\idot @>>> \eps_0^{o*}J_\idot\\
  @VVV @VVV\\
  \eps_1^{o*}J_\idot @>>> J_\idot
\end{CD}
\end{equation}
in $\Rel$. But the square in $\Pos$ underlying \eqref{J.pb.sq} is
the cylinder square \eqref{cyl.o.sq} for the embedding
$\eps_o^{o*}U(J_\idot) \to \eps_0^{o*}U(J_\idot)$, so $\C^\Tot$ is
semicartesian along this square by Lemma~\ref{bar.enh.le}~\thetag{ii},
and then $(\C^+|I)^\dm$ is semicartesian along \eqref{J.pb.sq} by
the same argument as in Lemma~\ref{j.ind.le}.
\endproof

\begin{corr}\label{ccat.dim1.corr}
Assume given a cofibration $\pi:I \to I'$ in $\Posf^+$ such that
$\dim I' \leq 1$, and as in Lemma~\ref{ccat.adj.le}, let $p:I' \to
\ppt$ be the projection. Then the functor $(\id|p)^{h*}:\Cat^h(I) =
\Cat^h(I|\ppt) \to \Cat^h(I|I')$ is an equivalence.
\end{corr}

\proof{} Since $\dim I' \leq 1$, $B^o(I') = \Tw(I')$ is its twisted
arrow category, and $I / I' \in \Posf^+$ can be expressed as the
cylindrical colimit of Lemma~\ref{pos.phi.le}. Let $\pi':J \to
B^o(I')$ be the corresponding cofibration, with the universal cone
map $g:J \to I / I'$, and let $g':B^o(I') \to \ppt$ be the
tautological projection.  Then $(g|g')^{h*}:\Cat^h(I/I') \to
\Cat^h(J|B^o(I')$ admits a right-adjoint $(g|g')^h_* \cong (\id
\times g'|g')^h_* \circ (g \times \pi'|\id)^h_*$, where since $\dim
B^o(I') \leq I'$, the adjoint $(\id \times g'|g')^h_*:\Cat^h((I/I')
\times B^o(I')|B^o(I')) \to \Cat^h(I/I')$ is provided by
Lemma~\ref{ccat.adj.le}. Recall that explicitly, $(\id \times
g'|g')^h_*$ is computed by considering the characteristic map
$\chi:I' \to [1]$ of the $0$-skeleton $\sk_0I'$, taking the
puhsforward $(\id \times q|q)^h_*$ with respect to the corresponding
cofibration $q:B^o(I') \to B^o([1]) = \V$, and then taking the
semicartesian product provided by Lemma~\ref{sq.h.le}. The product
in question is then exactly \eqref{co.I.sq} for $\chi \circ \pi:I
\to [1]$, so by Lemma~\ref{co.I.le}, the adjunction map $\id \to
(g|g')^h_* \circ (g|g')^{h*}$ is an isomorphism. We have the full
embedding $\eta:I \to I / I'$ with its left and right-adjoint maps
$\sigma,\eta^\dg:I / I' \to I$, and $\pi \circ \sigma \circ g =
\xi_\perp \circ \pi'$, where $\xi_\perp:B^o(I') \to I'$ is the map
\eqref{B.xi}. Therefore we have a well-defined functor $(\id|p)^h_*
= \eta^{\dg h}_* \circ (g|g')^h_* \circ (\sigma \circ
g|\xi_\perp)^{h*}:\Cat^h(I|I') \to \Cat^h(I)$, and it comes equipped
with an isomorphism
$$
\begin{aligned}
(\id|p)^h_* \circ (\id|p)^{h*} &\cong \eta^{\dg h}_* \circ
(g|g')^h_* \circ (\sigma \circ g|\xi_\perp)^{h*} \circ (\id|p)^{h*} \cong\\
&\cong \eta^{\dg h}_* \circ (g|g')^h_* \circ (g|g')^{h*} \circ
  \sigma^{h*} \cong 
\eta^{\dg h}_* \circ \sigma^{h*} \cong \id,
\end{aligned}
$$
while the adjunction map $(g|g')^{h*} \circ (g|g')^h_* \to \id$ induces a
map
\begin{equation}\label{ccat.dim1.eq}
(\sigma \circ g|\xi_\perp)^{h*} \circ (\id| p)^{h*} \circ (\id | p)^h_* \to
(\sigma \circ g|\xi_\perp)^{h*}.
\end{equation}
Moreover, if we let $B^o(I|I') = \xi_\perp^*I$, with the induced map
$\xi'_\perp:B^o(I|I') \to I$, then $\sigma \circ g = \xi'_\perp
\circ \sigma'$ for some map $\sigma':J \to B^o(I|I')$ over
$B^o(I')$, and then $\sigma'$ is right-reflexive, with a fully
faithful right-adjoint, so that $(\sigma|\id)^{h*}$ is fully faithful,
while $(\xi'_\perp|\xi_\perp)^{h*}$ is fully faithful for the same
reason as in the proof of Lemma~\ref{ccat.adj.le}. Therefore
$(\sigma \circ g|\xi_\perp)^{h*}$ is also fully faithful, and
\eqref{ccat.dim1.eq} comes from a map $(\id| p)^{h*} \circ (\id |
p)^h_*\to \id$. To finish the proof, it suffices to check that this
map is an isomorphism; this can be checked after restricting to each
fiber of the cofibration $J \to B^o(J)$, and then it immediately
follows from the construction of $(\id|p)^h_*$ via the semicartesian
square \eqref{co.I.sq}.
\endproof

\begin{corr}\label{co.I.corr}
For any right-closed embedding $r:I_1 \to I$ in $\Posf^+$, with the
complementary left-closed embedding $l:I \ssetminus I_1 \to I$, and
any $I$-aug\-mented enhanced category $\C$, the square
\begin{equation}\label{co.I.bis.sq}
\begin{CD}
\C @>>> r^h_*r^{h*}\C\\
@VVV @VVV\\
l^h_*l^{h*}\C @>>> r^h_*r^{h*}l^h_*l^{h*}\C
\end{CD}
\end{equation}
is semicartesian.
\end{corr}

\proof{} As in Lemma~\ref{co.I.le}, let $I' = \Cyl^o(r) \cong [1]
\setminus I$, consider the full embedding $\eta:I \to I'$, construct
the semicartesian square \eqref{co.I.sq} for the $I'$-augmented
enhanced category $\C' = \eta^h_*\C$, restrict it back to $I$ via
$\eta^{h*}$, and use \eqref{bc.enh.eq} to identify the restriction
with \eqref{co.I.bis.sq}.
\endproof

\subsection{Enhancement for $\Cat^h$.}\label{cath.subs}

As a more advanced application of Proposition~\ref{univ.prop} and
Lemma~\ref{sq.h.le}, let us now construct a natural enhancement for
the category $\Cat^h$ of small enhanced categories. Namely, for any
$I \in \Posf^+$, we have the category $\Cat^h(I)$ of
\eqref{cat.h.I.eq}, and for any map $f:I' \to I$ in $\Posf^+$, we
have the functor $f^{h*}:\Cat^h(I) \to \Cat^h(I')$, so that the
categories $\Cat^h(I)$ form a family of categories $\cCat^h \to
\Posf^+$ with fibers $\cCat^h_I = \Cat^h(I)$. Explicitly, objects in
$\cCat^h$ are pairs of an object $I \in \Posf^+$ and an
$I$-aug\-mented small enhanced category $\pi:\C \to \Unf(I)$, with
maps from $\langle I',\C' \rangle$ to $\langle I,\C \rangle$ given
by pairs of a map $f:I' \to I$ and a commutative square
\begin{equation}\label{ccat.sq}
\begin{CD}
\C' @>{\phi}>> \C\\
@V{\pi'}VV @VV{\pi}V\\
\Unf(I') @>{\Unf(f)}>> \Unf(I),
\end{CD}
\end{equation}
with $\phi$ cartesian over $\Unf(f)$, as in \eqref{rel.cart},
considered up to an isomorphism over $\Posf^+$. Inside $\cCat^h$, we
have a full subcategory $\sSets^h \subset \cCat^h$ spanned by pairs
$\langle I,\C \rangle$ such that $\pi:\C \to \Unf(I)$ is a family of
groupoids.

\begin{prop}\label{ccat.prop}
  The family $\cCat^h \to \Posf^+$ is an enhanced category, and
  $\sSets^h \subset \cCat^h$ is a full enhanced subcategory.
\end{prop}

In principle, Proposition~\ref{ccat.prop} follows from the results
we have already proved. Namely, the equivalences \eqref{cat.h.I.eq}
commute with pullbacks with respect to maps $f:I' \to I$, so that
taken together, they define an equivalence
\begin{equation}\label{ccat.h.eq}
\cCat^h \cong \Hh^W_{css}(\Delta^o\Delta^o\Sets),
\end{equation}
where $\Hh^W(\Delta^o\Delta^o\Sets)$ is the enhanced localization of
$\Delta^o\Delta^o\Sets$ provided by Lemma~\ref{enh.mod.le}, and
$\Hh^W_{css}(\Delta^o\Delta^o\Sets) \subset
\Hh^W(\Delta^o\Delta^o\Sets)$ is the enhanced full subcategory of
Corollary~\ref{enh.ff.corr} corresponding to the full subcategory
$h^W_{css}(\Delta^o\Delta^o\Sets) \subset h^W(\Delta^o\Delta^o\Sets)
= \Hh^W(\Delta^o\Delta^o\Sets)_{\ppt}$. The equivalence
\eqref{ccat.h.eq} then induces an equivalence
\begin{equation}\label{sset.h.eq}
  \sSets^h \cong \Hh^W(\Delta^o\Sets),
\end{equation}
where $\Delta^o\Sets \subset \Delta^o\Delta^o\Sets$ is embedded as
the full subcategory spanned by bisimplicial sets constant in the
Reedy direction. However, we can also give a direct proof using
Lemma~\ref{co.I.le} and its corollaries of
Subsection~\ref{univ.funct.subs}.

\proof[Proof of Proposition~\ref{ccat.prop}.] The family $\cCat^h
\to \Posf^+$ is obviously additive, it is complete by
Lemma~\ref{pf.cat.le}, and the subfamily $\sSets^h \subset \cCat^h$
is then also additive and preseved by the pushforward operations
$f^h_*$ of Lemma~\ref{pf.cat.le}, thus complete. By
Lemma~\ref{comp.C.le}, it then suffices to prove that $\cCat^h \to
\Posf^+$ is tight and weakly semicontinuous. To see that it is
tight, assume that as in Corollary~\ref{co.I.corr}, we are given a
right-closed embedding $r:I_1 \to I$ in $\Posf^+$, with the
complement $l:I_0 = I \ssetminus I_1 \to I$ and the corresponding
gluing functor $\theta = r^{h*} \circ l^h_*:\Cat^h(I_0) \to
\Cat^h(I_1)$, and define a category $\Cat^h(I)'$ by the cartesian
square
\begin{equation}\label{ccat.tight}
\begin{CD}
  \Cat^h(I)' @>>> \Cat^h(I \times \V|\V)\\
  @VVV @VV{b}V\\
  \Cat^h(I_1) \setminus_{\theta} \Cat^h(I_0) @>{a}>> \V^o\Cat^h(I),
\end{CD}
\end{equation}
where $b$ is the functor \eqref{cat.ev}, and $a$ is a full embedding
sending a triple $\langle \C_0,\C_1,\alpha \rangle$ of categories
$\C_l \in \Cat^h(I_l)$, $l = 0,1$ and a functor $\alpha:\C_1 \to
\theta(\C_0)$ to a functor $\V^o \to \Cat^h(I)$ given by
$$
\begin{CD}
  r^h_*\C_1 @>{r^h_*(\alpha)}>> r^h_*r^{h*}l^h_*\C_0 @<<< l^h_*\C_0,
\end{CD}
$$
where the map on the right is the tautological map. Then the adjoint
functor $\Cat^h(I \times \V|\V) \to \Cat^h(I)$ of
Lemma~\ref{ccat.adj.le} restricts to a functor $q:\Cat^h(I)' \to
\Cat^h(I)$, and since $\dim \V=1$, by the same argument as in
Corollary~\ref{ccat.dim1.corr}, Corollary~\ref{co.I.corr} implies
that $q$ is an equivalence. But then, again since $\dim \V=1$, $b$
in \eqref{ccat.tight} is an epivalence by Lemma~\ref{cat.epi.le}, so
the vertical arrow on the left is an epivalence as well.

Analogously, to show that $\cCat^h$ is weakly semicontinuous, assume
given $I \in \Posf^+$ equipped with a map $I \to \N$, and let
$l_n:I/n \to I$, $n \in \N$ be the embeddings of its
comma-fibers. Then for any small $I$-augmented enhanced category
$\C$, we have family of small $I$-augmented enhanced categories
$\C_n = l^h_{n*}l^{h*}_n\C$ over $\N$ that define a category
$\C_\idot \in \Cat^h(I \times \N|\N)$, hence also a category $(\id
\times \zeta|\zeta)^{h*}\C_\idot \in \Cat^h(I \times
\ZZ_\infty|\ZZ_\infty)$. On the other hand, the categories
$\Cat^h(I)_n = \Cat^h(I/n)$ themselves form a family
$\Cat^h(I)_\idot \to \N$, and for any $n \geq 0$,
$l^h_{n*}:\Cat^h(I)_n \to \Cat^h(I)$ is a fully faithful embedding,
so we have a fully faithful embedding $\Sec^\Tot(\N,\Cat^h(I)_\idot)
\to \N^o\Cat^h(I)$ and a fully faithful embedding
$a:\Sec^\Tot(\ZZ_\infty,\zeta^*\Cat^h(I)_\idot) \to
\ZZ_\infty^o\Cat^h(I)$. We can now define a category $\Cat^h(I)'$ by
the cartesian square
\begin{equation}\label{ccat.wc}
\begin{CD}
  \Cat^h(I)' @>>> \Cat^h(I \times \ZZ_\infty|\ZZ_\infty)\\
  @VVV @VV{b}V\\
  \Sec^\Tot(\ZZ_\infty,\zeta^*\Cat^h(I)_\idot) @>{a}>>
  \ZZ_\infty^o\Cat^h(I),
\end{CD}
\end{equation}
where again, $b$ is the functor \eqref{cat.ev}, thus an epivalence
since $\dim \ZZ_\infty = 1$. Again, the adjoint functor $\Cat^h(I
\times \ZZ_\infty|\ZZ_\infty) \to \Cat^h(I)$ of
Lemma~\ref{ccat.adj.le} restricts to a functor $q:\Cat^h(I)' \to
\Cat^h(I)$, and since by Lemma~\ref{Z.le}, we have
$\Sec^\Tot(\ZZ_\infty,\zeta^*\Cat^h(I)_\idot) \cong
\Sec^\Tot(\N,\Cat^h(I)_\idot)$, it suffices to show that $q$ is an
equivalence. We again have a functor $q':\Cat^h(I) \to \Cat^h(I)'$
sending $\C \in \Cat^h(I)$ to the family $(\id \times \zeta|
\zeta)^{h*}\C_\idot$, and $q' \circ q \cong \id$ by
\eqref{bc.enh.eq}, so it further suffices to shows that the
adjunction map $\id \to q \circ q'$ is an isomorphism. But the
adjoint functor of Lemma~\ref{ccat.adj.le} is constructed by first
taking the pushforward $c_*$ with respect to a cofibration
$c:\ZZ_\infty \to \V$, and then the semicartesian product,
and if we take the cartesian product instead, we obtain
$\Sec^\Tot(\ZZ_\infty,(\id \times \zeta|\zeta)^{h*}\C_\idot) \cong
\Sec^\Tot(\N,\C_\idot)$. Thus at the end of the day, we need to
check that for any $\C \in \Cat^h(I)$, the functor $\C \to
\Sec^{\Tot}(\N,\C_\idot)$ is an epivalence. As in
Lemma~\ref{co.I.le}, this can be checked pointwise over $\Unf(I)$,
and then reduces to Lemma~\ref{bar.enh.le}~\thetag{iii} by
\eqref{pf.pb.eq}.
\endproof

We note that the involution $\iota:\Cat^h \to \Cat^h$ sending a
small enhanced category $\C$ to the opposite small enhanced category
$\C^\iota$ has an obvious enhancement
\begin{equation}\label{iota.enh.eq}
\iota:\cCat^h \to \cCat^h, \qquad \langle I,\C \rangle \mapsto
\langle I,\C^\iota_{h\perp} \rangle
\end{equation}
in term of the enhancement $\cCat^h$ provided by
Proposition~\ref{ccat.prop}. Another useful application of the
enhancement $\cCat^h$ is its truncation functor
\eqref{trunc.eq}. Namely, for any $i \in I$ and $I$-augmented small
enhanced category $\C$, the enhanced fiber $\C_i$ is given by
\eqref{enh.aug.fib} --- that is, by the cartesian square
\begin{equation}\label{enh.I.sq}
\begin{CD}
\C_i @>>> \C\\
@VVV @VVV\\
\ppt^h @>{\Unf(\eps(i))}>> \Unf(I),
\end{CD}
\end{equation}
where the bottom arrow is the enhanced object in $\Unf(I)$
corresponding to $i$ --- and then for any map $f$ in $I$
corresponding to an order relation $i \leq i'$, \eqref{trunc.eq}
provides an enhanced functor
\begin{equation}\label{trans.I.eq}
  f^*:\C_{i'} \to \C_i.
\end{equation}
If $g$ represents another order relation $i' \leq i''$, we have
$f^* \circ g^* \cong (g \circ f)^*$. If we have another
$I$-augmented small enhanced category $\C'$ and an $I$-augmented
enhanced functor $\gamma:\C \to \C'$, then for any $i \in I$, we
have an enhanced fiber $\gamma_i:\C_i \to \C'_i$, and for any $f$,
we have an isomorphism
\begin{equation}\label{enh.fu.fib}
\gamma_f:\gamma_i \circ f^* \to f^* \circ \gamma_{i'},
\end{equation}
an enhanced version of the isomorphisms \eqref{fu.fib}. Note that
the isomorphisms \eqref{enh.fu.fib} are not unique --- in fact,
already $\gamma_i$ is only defined up to a non-unique isomorphism --
and it is certainly {\em not} true that an $I$-augmented enhanced
functor $\gamma$ is defined by ``the'' enhanced functors $\gamma_i$ and
``the'' isomorphisms \eqref{enh.fu.fib}. What we have, at the end of
the day, is a well-defined commutative square
\begin{equation}\label{enh.fu.fib.sq}
\begin{CD}
\C_{i'} @>{\gamma_i}>> \C'_{i'}\\
@V{f^*}VV @VV{f^*}V\\
\C_i @>{\gamma_i}>> \C_{i'}
\end{CD}
\end{equation}
in the category $\Cat^h$ that does {\em not} come equipped with a
preferred lifting to a commutative square \eqref{cat.dia}.

To complement Proposition~\ref{ccat.prop}, we note that the
definition of the categories $\Cat^h(I)$ actually makes sense for
any $I \in \Pos$. Therefore the family of categories $\cCat^h \to
\Posf^+$ extends to a family $\cCat^{h\Tot} \to \Pos$, with the same
definition. Let us show that this is actually the canonical
bar-invariant extension of the enhanced category $\cCat^h$ provided
by Proposition~\ref{lf.prop}.

\begin{lemma}\label{ccat.le}
The family $\cCat^{h\Tot} \to \Pos$ is bar-invariant.
\end{lemma}

\proof{} It suffices to show that for any $I \in \Pos$, with the
barycentric subdivison $BI$ and the map $\xi:BI \to I$ of
\eqref{B.xi}, and for any $\C \in \Cat^h(I)$ with canonical
extension $\C^\hash \to \Unf(I)^\hash$, the functor $\C^\hash \to
\Unf(\xi)_*\Unf(\xi)^*\C^\hash$ of \eqref{2.adj} is an equivalence,
while for any $\C \in \cCat^{h\dm}_{B_\dm(I)}$, the same holds for
the functor $\Unf(\xi)^*\Unf(\xi)_*\C^\hash \to \C^\hash$. By
\eqref{pf.pb.eq}, the first claim reduces to observing that for any
$\langle J,\alpha \rangle \in \Unf(I)^\hash$, with $J_\idot = J
/^\dm_{\alpha} R(I^o)$, the family $\C^{\hash\dm}$ is constant along
the map $R(\xi^o)^*J_\idot \to J_\idot$, and this follows from
Lemma~\ref{pm.B.le} and the bar-invariance of $\C^\hash$. Since we
already know that $\cCat^h$ is non-degenerate, the second claim can
be checked pointwise over $BI$, and since we are working with
categories in $\cCat^{h\dm}_{B_\dm(I)}$, it further suffices to
consider one-element subsets $S = \{i\} \in BI$. Then by
\eqref{bc.enh.eq}, we can replace $I$ with $I/i$, so that $I =
I_0^>$, $I_0 = I \ssetminus \{j\}$, and then as in Lemma~\ref{B.le},
the embedding $\eps(o):\ppt \to B(I^>)$ is the composition of the
right-reflexive embedding $\ppt \to B(I_0)^<$ and the left-reflexive
embedding $B(I_0)^< \to B(I_0^<)$. Since by
Example~\ref{unf.adj.exa}, the correspondence $I \mapsto \Unf(I)$
sends adjoint pairs of functors to adjoint pairs, this immediately
yields the claim.
\endproof

As one immediate application of Proposition~\ref{ccat.prop}, we
obtain a natural enhancement $\cCat$ for the category $\Cat^0$: this
is simply the full enhanced subcategory $\cCat \subset \cCat^h$ of
Corollary~\ref{enh.ff.corr} corresponding to the full embedding
\eqref{cat.h.cat}. Explicitly, objects in $\cCat$ are pairs $\langle
I,\C \rangle$ of $I \in \Posf^+$ and a small category $\C$ equipped
with a fibration $\C \to I$, and morphisms from $\langle I,\C
\rangle$ to $\langle I',\C' \rangle$ are pairs of a map $f:I \to I'$
and a commutative square
\begin{equation}\label{ccat.0.sq}
\begin{CD}
  \C @>{\phi}>> \C'\\
  @VVV @VVV\\
  I @>{f}>> I'
\end{CD}
\end{equation}
such that $\phi$ is cartesian over $f$, considered up to an
isomorphism. By Lemma~\ref{ccat.le}, the same description gives the
canonical extension $\cCat^\Tot$. We can further compose
\eqref{cat.h.cat} with the full embedding $\Pos \to \Cat^0$, and
this gives an enhancement $\pPos$ for $\Pos$ and enhanced full
embeddings
\begin{equation}\label{ppos.eq}
  \pPos \to \cCat \to \cCat^h,
\end{equation}
while the involution $\iota:\Pos \to \Pos$ extends to an equivalence
between $\pPos$ and the reflexive family $\Arc(\Pos) \to \Pos$ of
\eqref{C.fib.dia}.

\subsection{Enhancing lax functor categories.}\label{lax.subs}

Let us now upgrade Proposition~\ref{ccat.prop} by constructing
enhancements for the categories $\E \bbd^h \Cat^h$ and $\Cat^h \bb^h
\E$ of Subsection~\ref{enh.repr.subs}. This uses the full force of
Proposition~\ref{univ.prop}. We start with constructing an
enhancement $\E \bbd^h \cCat^h$ for $\E \bbd^h \Cat^h$. By
definition, its objects are triples $\langle I,\C,\alpha \rangle$,
$\langle I,\C \rangle \in \cCat^h$, $\alpha:\E \times^h \Unf(I^o)
\to \C_{h\perp}$ an enhanced functor over $\Unf(I^o)$, where
$\C_{h\perp} \to \Unf(I^o)$ is the transpose $I^o$-coaugmented
enhanced category of \eqref{iota.perp.eq}. Maps from $\langle
I',\C',\alpha' \rangle$ to $\langle I,\C,\alpha \rangle$ are
represented by triples $\langle f,\phi,a \rangle$ of a map $f:I' \to
I$, an enhanced functor $\phi:\C' \to \C$ that fits into a
commutative square \eqref{ccat.sq} and is cartesian over $\Unf(f)$,
and a map $a:\phi_{h\perp} \circ \alpha' \to \alpha \circ (\id
\times \Unf(f^o))$. Two triples $\langle f,\phi,a \rangle$, $\langle
f',\phi',a' \rangle$ define the same map iff $f=f'$, and there
exists an isomorphism $b:\phi' \cong \phi$ over $\Unf(I)$ such that
$a' = a \circ (b_{h\perp} \circ \alpha')$. The forgetful functor $\E
\bbd^h \cCat^h \to \Posf^+$ sending a triple $\langle \C,I,\phi
\rangle$ to $I$ is a fibration, and we have $(\E \bbd^h
\cCat^h)_\ppt \cong \E \bbd^h \Cat^h$. We also have the forgetful
functor
\begin{equation}\label{no.E.d}
  \E \bbd^h \cCat^h \to \cCat^h, \qquad \langle I,\C,\alpha \rangle
  \mapsto \langle I,\C \rangle,
\end{equation}
cartesian over $\Posf^+$, but even describing its fibers is somewhat
delicate. Namely, for any object $\langle I,\C \rangle \in \cCat^h$,
let $\F(I,\C,\E) = \Fun^h_{I^o}(\E \times^h \Unf(I^o),\C_{h\perp})$
be the enhanced category of enhanced functors $\E \times^h \Unf(I^o)
\to \C_{h\perp}$ over $\Unf(I^o)$. Then we have a functor
$\nu(I,\C):\F(I,\C,\E) \to (\E \bbd^h \cCat^h)_{\langle
  I,\C\rangle}$ sending $\alpha$ to $\langle I,\C,\alpha
\rangle$. However, in general, $\nu(I,\C)$ is an epivalence but {\em
  not} an equivalence: different maps $a,a':\alpha \to \alpha'$ in
the enhanced functor category can given triples $\langle \id,\id,a
\rangle$, $\langle \id,\id,a' \rangle$ that represent the same map
in the fiber $(\E \bbd^h \cCat^h)_{\langle I,\C \rangle}$ (this
happens when $a' = a \circ (b_{h\perp} \circ \alpha)$ for some
automorphism $b$ of the identity enhanced functor $\id:\C \to
\C$). The best we can get is a natural presentation
\begin{equation}\label{lax.pres.dia}
\begin{CD}
\F'(I,\C,\E) \quad\rightrightarrows\quad \F(I,\C,\E)
@>{\nu(I,\C)}>> (\E \bbd^h \cCat^h)_{\langle I,\C \rangle}
\end{CD}
\end{equation}
of the epivalence $\nu(I,\C)$ in the sense of
Definition~\ref{pres.def}. To obtain \eqref{lax.pres.dia}, one takes
as $\F'(I,\C,\E)$ the category of quadruples $\langle
\alpha,\alpha',\phi,a \rangle$ of enhanced functors
$\alpha,\alpha':\E \times^h \Unf(I^o) \to \C){h\perp}$, an
$I$-augmented equivalence $\phi:\C \to \C$, and an isomorphism
$a:\alpha' \cong \phi_{h\perp} \circ \alpha$, with morphisms given
by isomorphism classes of functors. The two projections $\pi$, $\pi'$ in
\eqref{lax.pres.dia} send $\langle \alpha,\alpha',\phi,a \rangle$ to
$\alpha$ resp.\ $\alpha'$, and $\langle \id,\phi,a \rangle$ provies
an isomorphism $\nu(I,\C,\E) \circ \pi' \cong \nu(I,\C,\E) \circ
\pi$. Forgetting $\alpha'$ and $a$ but keeping $\phi$ provides an
equivalence
\begin{equation}\label{lax.pres.eq}
\F'(I,\C,\E) \cong \F(I,\C,\E) \times \Fun^h_I(\C,\C)^{\id}_{\Iso},
\end{equation}
where the second factor is the full subcategory in the isomorphism
groupoid $\Fun^h_I(\C,\C)_{\Iso}$ spanned by the identity functor
$\id$.

\begin{prop}\label{ccat.J.prop}
The family of categories $\E \bbd^h \cCat^h \to \Posf^+$ is an
enhanced category, and \eqref{no.E.d} is an enhanced functor.
\end{prop}

\proof{} Choose a universal object for $\E$ provided by
Proposition~\ref{univ.prop}, and consider the corresponding
semicocartesian square \eqref{enh.univ.sq} of
Example~\ref{enh.univ.exa}. Then the enhanced functor $a$ induces a
fully faithful embedding $\E \bbd^h \cCat^h \to \Unf(I^o) \bbd^h
\cCat^h$, and once we know that its target is an enhanced category,
its source is an enhanced category by
Corollary~\ref{enh.ff.corr}. Therefore we may assume right away that
$\E \cong \Unf(I^o)$ for some $I \in \Posf^+$. To simplify notation,
let us denote $I^o \bbd^h \cCat^h = \Unf(I^o) \bbd^h \cCat^h$.

To check that $I^o \bbd^h \cCat^h \to \Posf^+$ is a reflexive family
of categories, take some $J \in \Posf^+$, with the maps $e:J \times
[1] \to J$, $s:J \to J \times [1]$. By Lemma~\ref{a.b.adj.le}, we
need to show that for any $\langle J \times [1],\C,\alpha \rangle
\in I^o \bbd^h \cCat^h$, there exists a map $\langle \id,\phi,a
\rangle:\langle J \times [1],\C,\alpha \rangle \to e^*s^*\langle J
\times [1],\C,\alpha\rangle$ satisfying the required universal
property. Since $\cCat^h$ is an enhanced category by
Proposition~\ref{ccat.prop}, it is in particular a reflexive family,
and this provides the enhanced functor $\phi:\C \to e^*s^*\C$. But
the required map $a$ is then a map in the enhanced sections category
$\sSec(I^o \times J^o \times [1]^o,(s^*\C)_{h\perp})$ of
\eqref{ssec.enh.I}, and its existence and the universal property
immediately follow from the reflexivity of the enhanced category
$\C_{h\perp} \to \Posf^+$ over $I \times J$.

Analogously, the fact that $I^o \bbd^h \cCat^h$ is non-degenerate
and separated immediately follows from the corresponding properties
of $\cCat^h$ and the enhanced categories $\C_{h\perp}$. Since $I^o
\bbd^h \cCat^h$ is obviously additive, it remains to check that it
is semiexact, semicontinuous, and satisfies excision and the
cylinder axiom. In other words, we have to check that $I^o \bbd^h
\cCat^h$ is semicartesian over any square \eqref{J.sq} of
Lemma~\ref{sec.sq.le}. Since we already know this for $\cCat^h$, by
Lemma~\ref{car.rel.le}, it suffices to check that for any small
$J$-augmented enhanced category $\C$, the fiber of \eqref{no.E.d}
over the corresponding cartesian section $[1]^2 \to \cCat^h$ is
semicartesian. Then by Lemma~\ref{semi.epi.le}, it suffices to check
that the squares formed by categories $\F(-,\C,\Unf(I^o))$,
$\F'(-,\C,\Unf(I^o))$ of \eqref{lax.pres.dia} are semicartesian. For
$\F$, we have $\F(J,\C,\Unf(I^o)) \cong \Sec^h(I^o \times
J^o,\C_{h\perp})$, so the claim immediately follows from
Lemma~\ref{sec.sq.le}. For $\F'$, \eqref{lax.pres.eq} and
Lemma~\ref{car.epi.le} further reduce us to showing that the square
formed by $\Fun^h(\C,\C)^{\id}_{\Iso}$ is semicartesian. However,
for any map $f:J' \to J$, $\Fun^h(f^*\C,f^*\C)^{\id}_{\Iso}$ is the
full subcategory in $\Fun^h(f^*\C,\C)_{\Iso}$ spanned by the
enhanced functor $f^*\C \to \C$, or equivalently, the subcategory in
$\Fun^h(f^*\C_{h\perp},\C_{h\perp})$ spanned by $f^*\C_{h\perp} \to
\C_{h\perp}$. Then to finish the proof, it suffices to apply the
following general result to $\E = \C_{h\perp}$.

\begin{lemma}\label{sec.C.sq.le}
In the situation of Lemma~\ref{sec.sq.le}, fix a small $J$-augment\-ed
enhanced category $\C$, and let $\C_0$, $\C_1$, $\C_{01}$ be its
restrictions to $J_0$, $J_1$, $J_{01}$. Then for any enhanced
category $\E$, the square
\begin{equation}\label{C.E.J.sq}
\begin{CD}
\fFun^h(\C_{h\perp},\E) @>>> \fFun^h(\C_{0h\perp},\E)\\
@VVV @VVV\\
\fFun^h(\C_{1h\perp},\E) @>>> \fFun^h(\C_{01h\perp},\E)
\end{CD}
\end{equation}
induced by \eqref{J.sq} is semicartesian.
\end{lemma}

\proof{} Choose a biordered set $I^\dm \in \Relf^+$ and a universal
object $c \in \C^\dm_{h\perp}$ over $I^\dm$ provided by
Proposition~\ref{univ.prop}. Then in particular, the corresponding
map $I = U(I^\dm) \to J^o$ is a fibration, so that $I^o \to J$ is a
cofibration, and the enhanced functor $\Unf(I^o) \to \C_{h\perp}$
induces a full embedding $\fFun^h(\C_{h\perp},\E) \to
\sSec^h(I^o,\E)$, and similarly for $\C_0$, $\C_1$ and
$\C_{01}$. This reduces us to the case $\C \cong \Unf(I)$, and in
this case, the claim is Lemma~\ref{sec.sq.le}.
\endproof

\begin{exa}\label{enh.idot.exa}
If $I=\ppt$, then $\cCat^h_\idot = \ppt \bbd^h \cCat^h$ is an
enhancement for the category $\Cat^h_\idot$ of pairs of a small
enhanced category $\C$ and an enhanced object $c \in \C_\ppt$. This
is an enhanced version of Example~\ref{cat.exa}. Explicitly, objects
in $\cCat^h_\idot$ are triples $\langle I,\C,c \rangle$, $\langle
I,\C \rangle \in \cCat^h$, $c \in \C_{h\perp I}$, with morphisms
$\langle I',\C',c' \rangle \to \langle I,\C,c \rangle$ represented
by isomorphism classes of triples $\langle f,\phi,g \rangle$, $f:I'
\to I$, $\phi:\C' \to \Unf(f)^*\C$ an $I$-augmented enhanced
functor, $g:\phi_{h\perp}(c') \to f^{o*}c$ a map in $\C_{h\perp
  I'}$.
\end{exa}

For $\Cat^h \bb^h \E$, the definition is similar but the proof is
more difficult. Assume given an enhanced category $\E$, and note
that by Corollary~\ref{fun.h.corr}, even if $\E$ is not small,
$\Cat^h \bb^h \E$ is still well-defined, with the same definition as
in Subsection~\ref{enh.repr.subs}. Let us construct an enhancement
$\cCat^h \bb^h \E$ for $\Cat^h \bb^h \E$. By definition, its objects
are triples $\langle I,\C,\alpha \rangle$, $\langle I,\C \rangle \in
\cCat^h$, $\alpha:\C_{h\perp} \to \E$ an enhanced functor, where
$\C_{h\perp} \to \Unf(I^o)$ is the transpose $I^o$-coaugmented
enhanced category of \eqref{iota.perp.eq}. Morphism from $\langle
I',\C',\alpha'\rangle$ to $\langle I,\C,\alpha \rangle$ are
represented by triples $\langle f,\phi,a \rangle$ of a map $f:I' \to
I$, an enhanced functor $\phi:\C' \to \C$ that fits into a
commutative square \eqref{ccat.sq} and is cartesian over $\Unf(f)$,
and a map $a:\alpha' \to \alpha \circ \phi_{h\perp}$. Two triples
$\langle f,\phi,a \rangle$, $\langle f',\phi',a' \rangle$ define the
same morphism iff $f=f'$, and there exists an isomorphism $b:\phi'
\cong \phi$ over $\Unf(I)$ such that $a = a' \circ
\alpha(b_{h\perp})$. The forgetful functor $\cCat^h \bb^h \E \to
\Posf^+$ sending a triple $\langle \C,I,\phi \rangle$ to $I$ is a
fibration, and we have $(\cCat^h \bb^h \E)_\ppt \cong \Cat^h \bb^h
\E$. We also have the forgetful functor
\begin{equation}\label{no.E}
  \cCat^h \bb^h \E \to \cCat^h, \qquad \langle I,\C,\alpha \rangle
  \mapsto \langle I,\C \rangle,
\end{equation}
cartesian over $\Posf^+$, and as for \eqref{no.E.d}, its fiber over
some $\langle I,\C \rangle \in \cCat^h$ is equipped with an
epivalence $\nu^\perp(I,\C):\Fun^h(\C_{h\perp},\E) \to (\cCat^h \bb^h
\E)_{\langle I,\C \rangle}$ that admits a presentation
$$
\begin{CD}
\Fun^h(\C_{h\perp},\E) \times \Fun^h(\C,\C)^{\id}_{\Iso}
\rightrightarrows \Fun^h(\C_{h\perp},\E)
@>{\nu^\perp(I,\C)}>> (\cCat^h \bb^h \E)_{\langle I,\C \rangle}.
\end{CD}
$$
Since partially ordered sets are rigid, $\nu^\perp(I,\C)$ is an
equivalence as soon as $\C \cong \Unf(J)$, $J \in \Posf^+$. The
simplest case is $I=\ppt$, $\C=\ppt^h$ that gives a canonical
embedding
\begin{equation}\label{y.eq}
\y:\E \cong \fFun^h(\ppt^h,\E) \cong (\cCat^h \bb^h \E)_{\langle
  \ppt,\ppt^h \rangle} \to \cCat^h \bb^h \E.
\end{equation}
More generally, if we restrict \eqref{no.E} with respect to the full
embedding $\pPos \to \cCat^h$ of \eqref{ppos.eq}, it becomes a
fibration: when one defines a category $\pPos \bb^h \E$ by the
cartesian square
\begin{equation}\label{ppos.E}
\begin{CD}
  \pPos \bb^h \E @>>> \cCat^h \bb^h \E\\
  @VVV @VVV\\
  \pPos @>>> \cCat^h,
\end{CD}
\end{equation}
the vertical arrow on the left is a fibration, and $\pPos \bb^h \E$
is an enhanced category. Explicitly, $\Pos \bb^h \E = (\pPos \bb^h
\E)_\ppt$ is identified with $\iota^*\E^\Tot$, where $\E^\Tot \to
\Pos$ is the canonical bar-invariant extension of the family $\E \to
\Posf^+$, and the equivalence $\iota:\Pos \bb^h \E \cong
\iota^*\E^\Tot \to \E^\Tot$ extends to an equivalence $\pPos \bb^h
\E \cong \E_\wr$, where $\E_\wr$ is the enhanced category of
Lemma~\ref{C.wr.le}. Even more explicitly, for any fibration $J \to
I$ in $\Pos$ with the induced $I$-augmentation $\Unf(J) \to
\Unf(I)$, we have $\Unf(J)_{h\perp} \cong \Unf(J_\perp)$, where
$J_\perp \to I^o$ is the transpose cofibration. Then by
\eqref{ev.j.I}, objects in $\pPos \bb^h \E$ are triples $\langle
I,J,e \rangle$ of $I \in \Posf^+$, a fibration $J \to I$ in $\Pos$,
and an object $e \in \E_{J^o_\perp}$. Moreover, for any $\langle
I,\C,\alpha \rangle \in \cCat^h \bb^h \E$, we have
\begin{equation}\label{ccat.hom}
(\cCat^h \bb^h \E)(\langle I,J,e \rangle,\langle I,\C,\alpha
  \rangle) \cong \pi_0(e \setminus_\alpha
  (\C_{h\perp})^\dg_{J^o_\perp})_{\Iso},
\end{equation}
where $(\C_{h\perp})^\dg_{J^o_\perp} \subset
(\C_{h\perp})_{J^o_\perp} \cong \Fun^h(\Unf(J_\perp),\C_{h\perp})$
denotes the full subcategory spanned by functors co-augmented over
$I$.

\begin{prop}\label{ccat.E.prop}
  The family of categories $\cCat^h \bb^h \E \to \Posf^+$ is an
  enhanced category, and \eqref{no.E} is an enhanced functor.
\end{prop}

\proof{} The family $\cCat^h \bb^h \E \to \Posf^+$ is obviously
additive. To see that it is reflexive, assume given some $I \in
\Posf^+$ and an object $\langle I \times [1],\C,\alpha \rangle$ in
$\cCat^h \bb^h \E$. As in Definition~\ref{cat.refl.def}, let $s:I
\to I \times [1]$ be the embedding, and let $e:I \times [1] \to I$
be the projection. Choose $J^\dm \in \Relf^+$ and a universal object
$c \in (\C^\dm_{h\perp})_{J^\dm} \subset \C_{h\perp}^\dm$ provided
by Proposition~\ref{univ.prop}. Then in particular, the map $J^\dm
\to L(I \times [1]) \to L([1])$ is a biordered fibration, with some
fibers $J^\dm_0,J^\dm_1 \in \Relf^+$, and its transition functor
$g_J:J^\dm_1 \to J^\dm_0$ provides a map $a:J^\dm \to J^\dm_0 \times
L([1])$ equal to $\id$ resp.\ $g_J$ on $J^\dm_0 \subset J^\dm$
resp.\ $J^\dm_1 \subset J^\dm$. Moreover, since $s^{o*}c$ is
universal for $(s^{h*}\C)_\perp \cong \Unf(s^o)^*\C_\perp$, $a$
induces a map $\langle a,\phi \rangle:\langle I \times [1],\C,\alpha
\to \langle I,s^{h*}\C,\alpha_s \rangle$ in $\cCat^h \bb^h \E$,
where $\alpha_s$ is the restriction of $\alpha$ to
$(s^{h*}\C)_{h\perp} \subset \C_{h\perp}$. Then by
Lemma~\ref{a.b.adj.le}, it suffices to show that for any $\langle
I,\C',\alpha' \rangle \in \cCat^h \bb^h \E$ and map $\langle b,\phi'
\rangle:\langle I \times [1],\C,\alpha_s \rangle \to e^*\langle I
\times [1],\C',\alpha' \rangle$, we have
\begin{equation}\label{a.b.phi}
\langle b,\phi' \rangle = e^*s^*(\langle b,\phi' \rangle) \circ
\langle a,\phi \rangle
\end{equation}
as morphisms in $\cCat^h \bb^h \E$. Moreover, consider the fibration
$U(J^\dm) \to I$, and let $J = U(J^\dm)^o_\perp$, with the fibration
$J \to I$. Then by \eqref{ev.j.I}, $c$ defines an enhanced functor
$\wt{c}:\Unf(J) \to \C$, and Proposition~\ref{univ.prop} insures
that $\wt{c}$ is $I$-augmented, so that we obtain a morphism $\langle
I \times [1],\Unf(J),\alpha \circ \wt{c} \rangle \to \langle I
\times [1],\C,\alpha \rangle$ in $\cCat^h \bb^h \E$. By
\eqref{ev.uni} and \eqref{ev.j.I}, the induced map
$$
\cCat^h \bb^h \E(\langle I \times [1],\C,\alpha \rangle,-) \to
\cCat^h \bb^h \E(\langle I \times [1],\Unf(J),\wt{c} \circ \alpha
\rangle,-)
$$
is injective. Therefore when checking \eqref{a.b.phi}, we may
replace $\C$ with $\Unf(J)$. Then by \eqref{ccat.hom}, the two maps
in \eqref{a.b.phi} correspond to two objects $c_0$, $c_1$ in $e
\setminus_\alpha (\C'_{h\perp})^\dg_{J^o_\perp}$, and we need to
check that these two objects are isomorphic. Moreover, since
$\C'_{h\perp}$ is reflexive, we have the adjoint pair of functors
$$
s^*:e \setminus_\alpha (\C'_{h\perp})^\dg_{J^o_\perp} \to e
\setminus_\alpha (\C'_{h\perp})^\dg_{s^*J^o_\perp},\quad\
s_\dg^*:e
\setminus_\alpha (\C'_{h\perp})^\dg_{s^*J^o_\perp} \to e \setminus_\alpha
(\C'_{h\perp})^\dg_{J^o_\perp}
$$
induced by $s:s^*J^o_\perp \to J^o_\perp$ and its adjoint
$s_\dg:J^o_\perp \to s^*J^o_\perp$, and we are already given an
isomorphism $s^*(c_0) \cong s^*(c_1)$, so it suffices to check that
both $c_0$ and $c_1$ lie in the essential image of the adjoint
fully faithful embedding $s_\dg^*$. Thus at the end of the day, it
suffices to check that for any $[1]$-augmented enhanced category
$\C$, fibration $J \to [1]$ and $[1]$-augmented functor
$\gamma:\Unf(J) \to \C$ corresponding to an object $c \in
\C_{J^o}$, the adjunction map $c \to e^*s^*c$ is an
isomorphism. Since $\C$ is non-degenerate, this can be checked
pointwise over $J$, and this further reduces us to the case
$J=[1]$. In this case, by Lemma~\ref{unf.refl.J.le}, the essential
image of the embedding $e^*:\C_\ppt \to \C_{[1]}$ is
$\C^\dm_{R([1])} \subset \C_{[1]}$, and the claim immediately
follows from the definition of the unfolding $\C^\dm$.

To prove that the reflexive family $\cCat^h \bb^h \E$ is
non-degenerate, recall that $\cCat^h$ is non-degenerate by
Proposition~\ref{ccat.prop}. For any $\langle I,\C \rangle \in
\cCat^h$, one can again choose a universal object for $\C_{h\perp}$
and use \eqref{ccat.hom} to reduce non-degeneracy for $\cCat^h \bb^h
\E$ to non-degeneracy of $\E$. Then since $\cCat^h \bb^h \E$ is
non-degenerate, the corresponding functors \eqref{nu.J.eq} are
automatically conservative, so to prove that $\cCat^h \bb^h \E$ is
separated, it suffices to check that they are essentially surjective
and full. For essential surjectivity, assume given some $I \in
\Posf^+$ and a morphism $\langle \id,\phi, a \rangle:\langle
I,\C_0,\alpha_0 \rangle \to \langle I,\C_1,\alpha_1 \rangle$ in
$\cCat^h \bb^h \E$ representing an object $c \in \Ar((\cCat^h \bb^h
\E)_I)$. Then since we know that $\cCat^h$ is separated, we can find
$\langle I \times [1],\C \rangle \in \cCat^h$ such that
$\nu_I(\langle I \times [1],\C \rangle)$ is identified with the
image of $c$ under the projection \eqref{no.E}, and we just need to
construct an enhanced functor $\alpha:\C_{h\perp} \to \E$ such that
$\nu_I(\langle I \times [1],\C,\alpha \rangle \cong c$. To do this,
choose a universal object for $\C_{h\perp}$ over some $J^\dm \in
\Relf^+$ equipped with a biordered fibration $J^\dm \to L(I \times
     [1])$, and use \eqref{ccat.hom} and
     Corollary~\ref{cat.epi.corr} for $\E$. To check that the
     essentially surjective functor $\nu_I$ for $\cCat^h \bb^h \E$
     is full, assume given another object $\langle I \times
     [1],\C',\alpha' \rangle$ in $\cCat^h \bb^h \E$, with $c' =
     \nu_I(\langle I \times [1],\C',\alpha' \rangle$, and take a map
     $c \to c'$ that we need to lift to a morphism in $\cCat^h \bb^h
     \E$. Then to do this, first lift our map to a morphism $\langle
     \id,\phi \rangle:\langle I \times [1],\C \rangle \to \langle
     I,\C' \rangle$ in $\cCat^h$, and then again use
     \eqref{ccat.hom} for $\langle I \times [1],\C',\alpha' \rangle$
     and Corollary~\ref{cat.epi.corr} for $\E$ to further lift
     $\langle \id,\phi \rangle$ to a morphism $\langle \id,\phi,a
     \rangle:\langle I \times [1],\C,\alpha \rangle \to \langle I
     \times [1],\C',\alpha' \rangle$ in $\cCat^h \bb^h \E$.

To finish the proof, it remains to check that the separated
non-degenerate reflexive family $\cCat^h \bb^h \E \to \Posf^+$ is
semiexact, semicontinuous, and satisfies excision and the cylinder
axiom. As in Proposition~\ref{ccat.J.prop},
Lemma~\ref{car.epi.le} and Lemma~\ref{semi.epi.le} immediately
reduce this to Lemma~\ref{sec.C.sq.le}.
\endproof

Finally, let us observe that we can combine
Proposition~\ref{ccat.E.prop} and Proposition~\ref{ccat.J.prop},
although we will only need the special case of the latter given in
Example~\ref{enh.idot.exa}. Namely, assume given an enhanced
category $\E$, and let $\cCat^h_\idot \bb \E$ be the category of
quadruples $\langle I,\C,\alpha,c\rangle$ of $I \in \Posf^+$, an
$I$-augmented enhanced category $\C$, an enhanced functor
$\alpha:\C_{h\perp} \to \E$, and a object $c \in \C_{h\perp
  I}$. Morphisms $\langle I',\C',\alpha',c' \rangle \to \langle
I,\C,\alpha,c \rangle$ are isomorphism classes of quadruples
$\langle f,\phi,a,g \rangle$, $f:I' \to I$, $\phi:\C' \to
\Unf(f)^*\C$ an $I$-augmented enhanced functor, $g:\phi_{h\perp}(c')
\to f^{o*}c$ a map in $\C_{h\perp I'}$. Then sending $\langle
I,\C,\alpha,c \rangle$ to $I$ defines a fibration $\cCat^h_\idot
\bb^h \E \to \Posf^+$, and the forgetful functors sending it to
$\langle I,\C,\alpha \rangle$ resp.\ $\langle I,\C,c \rangle$ fit
into a commutative square
\begin{equation}\label{ccat.dot.sq}
\begin{CD}
\cCat^h_\idot \bb^h \E @>{\pi}>> \cCat^h \bb^h \E\\
@V{\pi'}VV @VVV\\
\cCat^h_\idot @>>> \cCat^h
\end{CD}
\end{equation}
of fibrations over $\Posf^+$ and functors cartesian over
$\Posf^+$. Moreover, by its very definition, the square
\eqref{ccat.dot.sq} is semicartesian, and we have a functor
\begin{equation}\label{ccdot.pair}
\cCat^h_\idot \bb^h \E \to \E, \qquad \langle I,\C,\alpha,c \rangle
\mapsto \alpha(c),
\end{equation}
also cartesian over $\Posf^+$. In fact, \eqref{ccdot.pair} gives an
enhancement of Example~\ref{cat.bb.exa}.

\begin{prop}\label{ccat.dot.prop}
For any enhanced category $\E$, the family of categories
$\cCat^h_\idot \bb^h \E \to \Posf^+$ is an enhanced category.
\end{prop}

\proof{} We already know that $\cCat^h \bb^h \E$ in
\eqref{ccat.dot.sq} is an enhanced category by
Proposition~\ref{ccat.E.prop}. Then as in
Proposition~\ref{ccat.J.prop}, the fibers of the projection $\pi$
admits presentations \eqref{lax.pres.dia}, and one can apply exactly
the same argument as in Proposition~\ref{ccat.J.prop}.
\endproof

\section{The Yoneda package.}\label{yo.sec}

\subsection{Enhanced cylinders and comma-categories.}\label{enh.co.subs}

For any enhanced category $\C$, the simplest of the epivalences
\eqref{nu.J.eq} is $\nu_\ppt:\C_{[1]} \to \Ar(\C_\ppt)$. When $\C =
\cCat^h$ is the enhanced category of small enhanced categories
provided by Proposition~\ref{ccat.prop}, the existence of the
epivalence $\nu_\ppt$ means that a pair $\C_0$, $\C_1$ of small
enhanced categories equipped with an enhanced functor $\gamma:\C_0
\to \C_1$ defines a $[1]$-augmented enhanced category $\C$, uniquely
up to an equivalence. If one spells out how the proof of
Proposition~\ref{ccat.prop} works in this case, then the main fact
used is Corollary~\ref{co.I.corr} for $I=[1]$ and the embedding
$r=t:\ppt = [0] \to [1]$, and this reduces to the following enhanced
version of the cylinder construction.

\begin{defn}\label{enh.cyl.def}
The {\em enhanced cylinder} $\Cyl_h(\gamma)$ and the {\em dual
  enhanced cylinder} $\Cyl^\iota_h(\gamma)$ of an enhanced functor
$\gamma:\C_0 \to \C_1$ between small enhanced categories $\C_0$,
$\C_1$ are given by semicartesian squares
\begin{equation}\label{enh.cyl.sq}
\begin{CD}
\Cyl_h^\iota(\gamma) @>>> \C_1 \times^h \Unf([1])\\
@VVV @VVV\\
\C_0^{h<} @>{\gamma_{h<}}>> \C_1^{h<},
\end{CD}
\qquad
\begin{CD}
\Cyl_h(\gamma) @>>> \C_1 \times^h \Unf([1])\\
@VVV @VVV\\
\C_0^{h>} @>{\gamma_{h>}}>> \C_1^{h>},
\end{CD}
\end{equation}
where $\C_1^{h<} = t_*\C_1$ and $\C_1^{h>} = (\C_1^{h<})_{h\perp}$ are
the enhanced categories of \eqref{gt.lt.eq}, and the vertical arrows
on the right are the adjoint to the equivalence $t^*(C_1 \times^h
\Unf([1])) \cong \C_1$ and the transpose-opposite functor.
\end{defn}

By construction, the enhanced cylinder $\Cyl_h(\gamma)$ resp.\ the
dual enhanced cylinder $\Cyl^\iota_h(\gamma)$ is $[1]$-coaugmented
resp.\ $[1]$-augmented in the sense of Definition~\ref{cat.aug.def},
and we have $\Cyl_h(\gamma)^\iota \cong
\Cyl^\iota_h(\gamma^\iota)$. Composing the projection $\C_1 \times^h
\Unf([1]) \to \C_1$ with the top arrows in \eqref{enh.cyl.sq}
provides enhanced functors $t_\dg:\Cyl_h(\gamma) \to \C_1$,
$s_\dg:\Cyl^\iota_h(\gamma) \to \C_1$. By Corollary~\ref{sq.h.corr},
the embeddings $s,t:\C_1 \to \C_1 \times^h \Unf([1])$ provide
enhanced functors $s:\C_0 \to \Cyl_h(\gamma)$, $t:\C_1 \to
\Cyl_h(\gamma)$, $s:\C_1 \to \Cyl^\iota_h(\gamma)$, $t:\C_0 \to
\Cyl^\iota_h(\gamma)$ that come equipped with isomorphisms $t_\dg
\circ t \cong \id$, $s_\dg \circ s \cong \id$, and fit into
factorizations
\begin{equation}\label{enh.cyl.facto}
\begin{CD}
\C_0 @>{s}>> \Cyl_h(\gamma) @>{t_\dg}>> \C_1,
\end{CD}
\qquad
\begin{CD}
\C_0 @>{t}>> \Cyl^\iota_h(\gamma) @>{s_\dg}>> \C_1
\end{CD}
\end{equation}
of the enhanced functor $\gamma$. In particular, both
$\Cyl_h(\gamma)$ and $\Cyl^\iota_h(\gamma)$ are categories under
$\C_0 \copr \C_1$. If $\C_0$, $\C_1$ are equipped with enhanced
functors $\pi_l:\C_l \to \E$, $l=0,1$ to some enhanced category
$\E$, and $\gamma$ is a functor over $\E$, we will treat
$\Cyl_h(\gamma)$ resp.\ $\Cyl^\iota_h(\gamma)$ as categories over
$\E$ via the functors $\pi_1 \circ t_\dg$ resp.\ $\pi_1 \circ
s_\dg$, and then the embeddings $s$, $t$ and the factorizations
\eqref{enh.cyl.facto} are all over $\E$. Any $[1]$-augmented small
enhanced category $\C$ is of the form $\C \cong
\Cyl^\iota_h(\gamma)$ for an enhanced functor $\gamma:\C_1 = t^*\C
\to \C_0 = s^*\C$, unique up to an isomorphism.

\begin{lemma}\label{enh.cyl.le}
For any enhanced functor $\gamma:\C_0 \to \C_1$ between small
enhanced categories, the isomorphisms $s_\dg \circ s \cong \id$ resp.\
$\id \cong t_\dg \circ t$ define an adjunction between $s$ and
$s_\dg$ resp.\ $t_\dg$ and $t$.
\end{lemma}

\proof{} Since $\Cyl_h(\gamma) \cong \Cyl_h^\iota(\gamma^\iota)$,
and this equivalence interchanges $s$, $s_\dg$ and $t$, $t_\dg$, it
suffices to consider the dual enhanced cylinder
$\Cyl^\iota_h(\gamma)$. Then this is a $[1]$-augmented enhanced
category, and $\eps^h(0):\ppt^h \to \Unf([1])$ is right-admissible,
so $s:\C_1 \to \Cyl^\iota_h(\gamma)$ is right-admissible by
Lemma~\ref{adm.cof.le}. Moreover, the top arrow in
\eqref{enh.cyl.sq} is $[1]$-augmented, and this identifies its
right-adjoint $s_\dg$ with the enhanced functor \eqref{enh.cyl.facto}.
\endproof

For any enhanced functor $\gamma:\C_0 \to \C_1$ between small
enhanced categories, we also have the enhanced functor $\gamma
\times^h \id:\C_0 \times^h \Unf([1]) \to \C_1 \times^h \Unf([1])$,
and again by Corollary~\ref{sq.h.corr}, it induces enhanced functors
from the product $\C_0 \times^h \Unf([1])$ to $\Cyl_h(\gamma)$ and
$\Cyl^\iota_h(\gamma)$ that fit into commutative squares
\begin{equation}\label{enh.cyl.o.sq}
\begin{CD}
  \C_0 @>{\gamma}>> \C_1\\
  @V{\id \times^h t}VV @VV{t}V\\
  \C_0 \times^h \Unf([1]) @>>> \Cyl_h(\gamma),
\end{CD}
\qquad
\begin{CD}
  \C_0 @>{\gamma}>> \C_1\\
  @V{\id \times^h s}VV @VV{s}V\\
  \C_0 \times^h \Unf([1]) @>>> \Cyl^\iota_h(\gamma),
\end{CD}
\end{equation}
an enhanced version of the squares \eqref{cyl.o.sq}.

\begin{lemma}\label{enh.cyl.o.le}
The commutative squares \eqref{enh.cyl.o.sq} are
enhanced-semi\-co\-car\-te\-sian in the sense of
Definition~\ref{enh.cocart.def}.
\end{lemma}

\proof{} As in Lemma~\ref{enh.cyl.le}, it suffices to consider
enhanced dual cylinders. Then $\Cyl^\iota_h(\gamma)$ is
$[1]$-augmented, and Proposition~\ref{univ.prop} provides a
universal object in some fiber
$\Cyl^\iota_h(\gamma)^{\iota\dm}_{I^\dm}$, $I^\dm \in \Relf^+$ such
that the underlying partially ordered set $I = U(I^\dm) \in \Posf^+$
is fibered over $[1]$. By universality, it then suffices to prove
the claim for $\Cyl^\iota_h(\gamma) \cong \Unf(I)$, and in this
case, the claim is Example~\ref{cyl.coc.exa}.
\endproof

As in the unenhanced situation, the dual notion to cylinders and
dual cylinders is the notion of a left and right comma-category. For
any enhanced functor $\gamma:\C_0 \to \C_1$ between small enhanced
categories, the {\em left} resp.\ {\em right enhanced
  comma-category} are defined by semicartesian squares
\begin{equation}\label{enh.rl.comma}
\begin{CD}
\C_1 \setminus^h_\gamma \C_0 @>{\tau}>> \C_0\\
@VVV @VV{\gamma}V\\
\Ar_h(\C_1) @>{\tau}>> \C_1,
\end{CD}
\qquad
\begin{CD}
\C_0 /^h_\gamma \C_1 @>{\sigma}>> \C_0\\
@VVV @VV{\gamma}V\\
\Ar_h(\C_1) @>{\sigma}>> \C_1,
\end{CD}
\end{equation}
where $\Ar_h(\C_1)$ and $\sigma$, $\tau$ are as in
\eqref{enh.ar.eq}. We have $\C_0^\iota /^h_{\gamma^\iota} \C_1^\iota
\cong (\C_1 \setminus^h_\gamma \C_0)^\iota$ and $\C_1^\iota
\setminus^h_{\gamma^\iota} \C_0^\iota \cong (\C_0 /^h_\gamma
\C_1)^\iota$. As in the unenhanced case, we drop the subscript
$\gamma$ when it is clear from the context. By
Corollary~\ref{sq.h.corr}, the enhanced functor $\eta:\C_1 \to
\Ar_h(\C_1)$ of \eqref{enh.ar.eq} gives rise to factorizations
\begin{equation}\label{enh.comma.facto}
\begin{CD}
\C_0 @>{\eta}>> \C_1 \setminus^h_\gamma \C_0 @>{\sigma}>> \C_1,
\end{CD}
\qquad
\begin{CD}
\C_0 @>{\eta}>> \C_0 /^h_\gamma \C_1 @>{\tau}>> \C_1
\end{CD}
\end{equation}
of the enhanced functor $\gamma$, where $\sigma$ and $\tau$ are
obtained by composing the projections to $\Ar_h(\C_1)$ with
$\sigma,\tau:\Ar_h(\C_1) \to \C_1$. For any enhanced object $c$ in
the enhanced category $\C_1$, the {\em left} resp.\ {\em right
  enhanced comma-fibers} of the enhanced functor $\gamma$ over $c$
are the enhanced fibers $\C_0 /^h_\gamma c = (\C_0 /^h_\gamma
\C_1)_c$ resp.\ $c \setminus^h_\gamma \C_0 = (\C_1
\setminus^h_\gamma \C_1)_c$ of the projections $\tau$
resp.\ $\sigma$ of \eqref{enh.comma.facto} in the sense of
\eqref{enh.fib.sq}, and we denote by
\begin{equation}\label{si.ta.c.eq}
\sigma(c):\C_0 /^h_\gamma c \to \C_0, \qquad \tau(c):c
\setminus^h_\gamma \C_0 \to \C_0
\end{equation}
the enhanced functors induced by \eqref{enh.rl.comma}. More
generally, for any enhanced functor $\C \to \E$ between small
enhanced categories, we define the {\em relative enhanced arrow
  category} by the cartesian square
\begin{equation}\label{enh.ar.E.sq}
\begin{CD}
  \Ar^h(\C|\E) @>>> \Ar^h(\C)\\
  @VVV @VVV\\
  \E @>{\eta}>> \Ar^h(\E),
\end{CD}
\end{equation}
where $\eta$ is the fully faithful enhanced functor of
\eqref{enh.ar.eq}, and if small enhanced categories $\C_0$, $\C_1$
are equipped with enhanced functors $\pi_l:\C_l \to \E$, $l=0,1$, we
define the {\em enhanced relative comma-categories} of an enhanced
functor $\gamma:\C_0 \to \C_1$ over $\E$ by cartesian squares
\begin{equation}\label{enh.rel.lc.sq}
\begin{CD}
  \C_1 \setminus_{\gamma,\E} \C_0 @>>> \C_1 \setminus_\gamma \C_0\\
  @VVV @VV{\pi_1/\pi_0}V\\
  \E @>{\eta}>> \Ar^h(\E),
\end{CD}
\qquad\qquad
\begin{CD}
  \C_0 /_{\gamma,\E} \C_1 @>>> \C_0 /_\gamma \C_1\\
  @VVV @VV{\pi_0/\pi_1}V\\
  \E @>{\eta}>> \Ar^h(\E).
\end{CD}
\end{equation}
Then we have $\C_0^\iota /^h_{\gamma^\iota,\E^\iota} \C_1^\iota
\cong (\C_1 \setminus^h_{\gamma,\E} \C_0)^\iota$ and $\C_1^\iota
\setminus^h_{\gamma^\iota,\E^\iota} \C_0^\iota \cong (\C_0
/^h_{\gamma,\E} \C_1)^\iota$, and \eqref{enh.rl.comma} induce
semicartesian squares
\begin{equation}\label{enh.rl.E.comma}
\begin{CD}
\C_1 \setminus^h_{\gamma,\E} \C_0 @>{\tau}>> \C_0\\
@VVV @VV{\gamma}V\\
\Ar_h(\C_1|\E) @>{\tau}>> \C_1,
\end{CD}
\qquad
\begin{CD}
\C_0 /^h_{\gamma,\E} \C_1 @>{\sigma}>> \C_0\\
@VVV @VV{\gamma}V\\
\Ar_h(\C_1|\E) @>{\sigma}>> \C_1,
\end{CD}
\end{equation}
where $\Ar_h(\C_1|\E)$ is the enhanced relative arrow category of
\eqref{enh.ar.E.sq} (in particular, for any $\C \to \E$, we have $\C
\setminus_{\id,\E} \C \cong \C /_{\id,\C} \cong \Ar_h(\C|\E)$). The
factorizations \eqref{enh.comma.facto} induce factorizations
\begin{equation}\label{enh.comma.E.facto}
\begin{CD}
\C_0 @>{\eta}>> \C_1 \setminus^h_{\gamma,\E} \C_0 @>{\sigma}>> \C_1,
\end{CD}
\quad
\begin{CD}
\C_0 @>{\eta}>> \C_0 /^h_{\gamma,\E} \C_1 @>{\tau}>> \C_1
\end{CD}
\end{equation}
of the enhanced functor $\gamma$.

By Corollary~\ref{sq.h.corr}, the comma-category construction is
functorial in the following sense: for any commutative square
\begin{equation}\label{C.E.sq}
\begin{CD}
\C' @>{\gamma}>> \C\\
@V{\pi'}VV @VV{\pi}V\\
\E' @>{\phi}>> \E
\end{CD}
\end{equation}
of enhanced categories and enhanced functors, we have natural
functors
\begin{equation}\label{fun.comma}
\gamma \setminus^h \phi:\E' \setminus^h_{\pi'} \C' \to E
\setminus^h_\pi \C, \quad \gamma /^h \phi:\C' /^h_{\pi'} \E' \to \C
/^h_\pi \E,
\end{equation}
unique up to an isomorphism. The duality between cylinders and
comma-categories of \eqref{cyl.com} carries over to the
enhanced setting: for any enhanced functor $\gamma:\C_0 \to \C_1$
between small enhanced categories, and any small enhanced category
$\E$, Lemma~\ref{enh.cocart.le} and Lemma~\ref{enh.cyl.o.le} provide
identifications
\begin{equation}\label{ffun.cyl.co}
\begin{aligned}
\fFun^h(\Cyl_h(\gamma),\E) &\cong \fFun^h(\C_0,\E) \setminus^h_{\gamma^*}
\fFun^h(\C_1,\E),\\
\fFun^h(\Cyl^\iota_h(\gamma),\E) &\cong \fFun^h(\C_1,\E) /^h_{\gamma^*}
\fFun^h(\C_0,\E).
\end{aligned}
\end{equation}
Moreover, we also have en enhanced version of \eqref{cyl.co.sec}.

\begin{lemma}\label{cyl.comma.le}
For any enhanced functor $\gamma:\C_0 \to \C_1$ between small
enhanced categories, we have
\begin{equation}\label{enh.cyl.co.sec}
\C_0 /^h_\gamma \C_1 \cong \sSec^h([1],\Cyl_h(\gamma)), \qquad
\C_1 \setminus^h_\gamma \C_0 \cong
\sSec^h([1],\Cyl^\iota_h(\gamma)).
\end{equation}
\end{lemma}

\proof{} If $\C_1 = \ppt^h$, so that $\Cyl_h(\gamma) \cong
\C_0^{h>}$ and $\Cyl^\iota_h(\gamma) \cong \C_0^{h<}$, the claim
immediately follows from \eqref{kan.2adj} applied to
$\Unf(s),\Unf(t):\ppt^h = \Unf(\ppt) \to \Unf([1])$. In the general
case, the semicartesian squares \eqref{enh.cyl.sq} are
$[1]$-augmented, and for any $J \in \Posf^+$, $\sSec^h(J,-)$ sends
$J$-augmented semicartesian squares to semicartesian squares.
\endproof

In terms of \eqref{enh.cyl.co.sec}, the projections $\sigma$ and
$\tau$ in \eqref{enh.rl.comma} and \eqref{enh.comma.facto} are
induced by the embeddings $s,t:\ppt = [0] \to [1]$, and
Lemma~\ref{enh.2adm.le} then immediately shows that $\eta:\C_0 \to
\C_0 /^h \C_1$ in \eqref{enh.comma.facto} is left-adjoint to $\sigma$ of
\eqref{enh.rl.comma}, and dually, $\eta:\C_0 \to \C_1 \setminus^h \C_0$ is
right-adjoint to $\tau$. In the relative situation, since the top
arrows in \eqref{enh.rel.lc.sq} are fully faithful embeddings, these
adjunctions induce the corresponding adjunction for $\eta$, $\tau$
and $\sigma$ in \eqref{enh.rl.E.comma} and \eqref{enh.comma.E.facto}.

\begin{lemma}\label{adj.eq.le}
Assume given enhanced functors $\lambda:\C_0 \to \C_1$, $\rho:\C_1
\to \C_0$ between small enhanced categories. Then \thetag{i}
$\lambda$ is left-adjoint to $\rho$ iff \thetag{ii} there exists an
equivalence $\eps:\Cyl_h(\lambda) \cong \Cyl^\iota_h(\rho)$ under
$\C_0 \copr \C_1$ iff \thetag{iii} there exists an equivalence $\C_0
/^h_\lambda \C_1 \cong \C_0 \setminus^h_\rho \C_1$ over $\C_0
\times^h \C_1$. Moreover, if $\C_0$, $\C_1$ are equipped with
enhanced functors $\pi_l:\C_l \to \E$, $l=0,1$ to a small enhanced
category $\E$, and $\rho$, $\lambda$ are functors over $\E$ adjoint
over $\E$, then in \thetag{ii}, one can choose an equivalence $\eps$
that fits into a commutative diagram
\begin{equation}\label{cyl.adj.dia}
\begin{CD}
\C_0 \copr \C_1 @>{s \copr t}>> \Cyl_h(\lambda) @>{t_\dg \circ
  \pi_1}>> \E\\
@| @V{\eps}VV @|\\
\C_0 \copr \C_1 @>{t \copr s}>> \Cyl^\iota_h(\rho) @>{s_\dg \circ
  \pi_1}>> \E,
\end{CD}
\end{equation}
and in \thetag{iii}, one has an equivalence $\C_0 /^h_{\lambda,\E}
\C_1 \cong \C_0 \setminus^h_{\rho,\E} \C_1$ over $\C_0 \times^h_{\E}
\C_1$.
\end{lemma}

\proof{} Both \thetag{ii} and \thetag{iii} immediately imply
\thetag{i}, by \eqref{enh.cyl.facto}
resp.\ \eqref{enh.comma.facto}. In the absolute case,
\thetag{ii} implies \thetag{iii} by \eqref{enh.cyl.co.sec}, and if
we are over $\E$, then \eqref{cyl.adj.dia} insures that
$\sSec^h([1],\eps)$ is an equivalence over
$$
(\C_0 \times \C_1) \times^h_{\E \times \E} \sSec^h([1],[1] \times^h
\E) \cong (\C_0 \times \C_1) \times^h_{\E \times \E} \Ar^h(\E)
$$
which then restricts to an equivalence $\C_0 /^h_{\lambda,\E} \C_1
\cong \C_0 \setminus^h_{\rho,\E} \C_1$ over $\C_0 \times^h_{\E}
\C_1$ via the full embedding $\eta:\E \to \Ar^h(\E)$. To finish the
proof, it remains to check that \thetag{i} implies \thetag{ii}.
Decompose $\rho \cong s_\dg \circ t$ by \eqref{enh.cyl.facto}, let
$a:\id \to \rho \circ \lambda$ be the adjunction map, and let $a':s
\to t \circ \lambda$ be the adjoint map of enhanced functors $\C_0
\to \Cyl^\iota_h(\rho)$ that exists by virtue of the adjunction of
Lemma~\ref{enh.cyl.le}. Then by \eqref{ar.enh.eq}, $a'$ defines an
enhanced functor $\alpha:\C_0 \times^h \Unf([1]) \to
\Cyl^\iota_h(\rho)$ equipped with isomorphisms $\alpha \circ s \cong
s$, $\alpha \circ t \cong \lambda \circ t$, and by
Lemma~\ref{enh.cyl.o.le}, this in turn gives rise to a functor
$\eps:\Cyl_h(\lambda) \to \Cyl^\iota_h(\rho)$ and
isomorphisms $\eps \circ s \cong s$, $\eps \circ t \cong
t$. Moreover, if we are over $\E$, then by the same semicocartesian
property of Lemma~\ref{enh.cyl.o.le}, we obtain \eqref{cyl.adj.dia}.
Dually, the adjunction map $a_\dg:\lambda \circ \rho \to \id$
defines a functor $\eps_\dg:\Cyl^\iota_h(\rho) \to \Cyl_h(\lambda)$,
and $a$ and $a_\dg$ define an adjunction between $\lambda$ and
$\rho$ if and only if $\eps \circ \eps_\dg \cong \id$ and $\eps_\dg
\circ \eps \cong \id$.
\endproof

\begin{corr}\label{enh.adj.corr}
Assume given a commutative diagram
\begin{equation}\label{l.r.enh.dia}
\begin{CD}
  \C_0' @>{\lambda'}>> \C_1' @>{\rho'}>> \C_0' @>>> \E'\\
  @VVV @VVV @VVV @VVV\\
  \C_0 @>{\lambda}>> \C_1 @>{\rho}>> \C_0 @>>> \E
\end{CD}
\end{equation}
of small enhanced categories and enhanced functors, and assume that
all squares in \eqref{l.r.enh.dia} are semicartesian, and $\lambda$
is left-adjoint to $\rho$ over $\E$. Then $\lambda'$ is left-adjoint
to $\rho'$ over $\E'$.
\end{corr}

\proof{} Use Corollary~\ref{sq.h.corr} and the criterion of
Lemma~\ref{adj.eq.le}~\thetag{iii}.
\endproof

\begin{corr}\label{enh.adj.fun.corr}
Assume given enhanced categories $\C$, $\E$ with $\C$ small. Then
for any right-reflexive enhanced functor $\lambda:\E \to \E'$ with
right-adjoint $\rho:\E' \to \E$, the enhanced functors between
$\fFun^h(\C,\E)$ and $\fFun^h(\C,\E')$ given by postcompostion with
$\lambda$ and $\rho$ are adjoint, and for any left-reflexive
enhanced functor $\rho:\C' \to \C$ from a small enhanced category
$\C$, with left-adjoint $\lambda:\C \to \C'$, the pullback
$\rho^*:\fFun^h(\C,\E) \to \fFun^h(\C',\E)$ is right-adjoint to the
pullback $\lambda^*:\fFun^h(\C',\E) \to \fFun(\C,\E)$.
\end{corr}

\proof{} The first claim immediately follows from
Lemma~\ref{adj.eq.le}~\thetag{iii} applied to $\E$, $\E'$ and the
functor categories. For the second, use
Lemma~\ref{adj.eq.le}~\thetag{ii} for $\C$, $\C'$, and then apply
\eqref{ffun.cyl.co} and Lemma~\ref{adj.eq.le}~\thetag{iii}.
\endproof

\begin{lemma}\label{enh.adj.co.le}
Assume given an enhanced functor $\lambda:\C_0 \to \C_1$ between
small enhanced categories equipped with an enhanced section
$\rho:\C_1 \to \C_0$. Then the isomorphism $\lambda \circ \rho \cong
\id$ defines an adjunction between $\lambda$ and $\rho$ if and only
if $\sigma:\C_0 \setminus^h_{\rho,\C_1} \C_1 \to \C_0$ is an
equivalence.
\end{lemma}

\proof{} Same as Lemma~\ref{adj.co.le}, with enhanced
comma-categories instead of the usual ones, and
Lemma~\ref{adj.eq.le}~\thetag{iii} in the ``only if'' part.
\endproof

\begin{corr}\label{enh.adm.sq.corr}
\begin{enumerate}
\item Assume given a semicartesian square \eqref{cat.dia} of small
  enhanced categories and enhanced functors, and left
  resp.\ right-admissible enhanced full subcategories $\C' \subset
  \C$, $\C'_l \subset \C_l$, $l=0,1$, with adjoint enhanced functors
  $\lambda:\C \to \C'$, $\lambda_l:\C_l \to \C'_l$, $l=0,1$, such
  that $\gamma_l(\C'_l) \subset \C'$, $l = 0,1$, and \eqref{bc.eq}
  is an isomorphism $\gamma_l \circ \lambda_l \cong \lambda \circ
  \gamma_l$, again for $l=0,1$. Then the full subcategory $\C'_{01}
  = \C'_0 \times^h_{\C'} \C'_l \subset \C_{01}$ is left
  resp.\ right-admissible, with some adjoint $\lambda_{01}:\C_{01}
  \to \C'_{01}$, and we have$\gamma^l_{01} \circ \lambda_{01} \cong
  \lambda_l \circ \gamma^l_{01}$ for $l=0,1$.
\item Assume given a commutative square \eqref{cat.dia} of small
  enhanced categories and enhanced functors, and left
  resp.\ right-admissible enhanced full subcategories $\C' \subset
  \C$, $\C'_l \subset \C_l$, $l=0,1,01$, with adjoint enhanced
  functors $\lambda:\C \to \C'$, $\lambda_l:\C_l \to \C'_l$,
  $l=0,1,01$, such that $\gamma_l(\C'_l) \subset \C'$,
  $\gamma^l_{01}(\C'_{01}) \subset \C'_l$, $l = 0,1$, and
  \eqref{bc.eq} provides isomorphisms $\gamma_l \circ \lambda_l
  \cong \lambda \circ \gamma_l$, $\gamma^l_{01} \circ \lambda_{01}
  \cong \lambda_l \circ \gamma^l_{01}$, $l=0,1$. Then if the square
  is semicartesian, cartesian, enhanced-semicocartesian or
  enhanced-cocartesian, so is the square of categories $\C'$,
  $\C'_l$, $l=0,1,01$.
\end{enumerate}
\end{corr}

\proof{} The proof for \thetag{i} is the same as
Corollary~\ref{adm.sq.corr}, with Lemma~\ref{enh.adj.co.le} instead
Lemma~\ref{adj.co.le}, and semicartesian products instead of
cartesian products. For \thetag{ii}, in the semicartesian and
cartesian case, since the adjoint functr $\lambda_l:\C_l \to
\C_l'$$, l = 0,1$ are essentially surjective, the adjoint
$\lambda_{01} \cong (\lambda_0 \circ \gamma^0_{01}) \times^h
(\lambda_1 \circ \gamma^1_{01}):\C_{01} \to \C_{01}' \subset
\C_{01}$ must be essentially surjective onto $\C_0' \times^h_{\C'}
\C_1' \supset \C'_{01}$, so the two full subcategories in $\C_{01}$
coincide. In the semicocartesian and cocartesian case, combine
Lemma~\ref{enh.cocart.le}, Corollary~\ref{enh.adj.fun.corr} and
\thetag{ii} in the semicartesian resp.\ cartesian case.
\endproof

In the situation of Corollary~\ref{enh.adj.fun.corr}, if $\lambda:\C
\to \C'$ is fully faithful, so that $\rho \circ \lambda \cong \id$,
then $\lambda^* \circ \rho^* \cong \id$ as well, so that
$\rho^*:\fFun^h(\C,\E) \to \fFun^h(\C',\E)$ is also fully
faithful. For applications, it is useful to have a description of
its essential image. To obtain such a description, note that the
evaluation pairing \eqref{enh.ev.eq} restricts to a pairing over
$\ppt \in \Posf^+$, and this gives rise to a comparison functor
\begin{equation}\label{fun.pt.eq}
  \fFun^h(\C,\E)_{\ppt} = \Fun^h(\C,\E) \to \Fun(\C_{\ppt},\E_{\ppt}).
\end{equation}
In general, the functor \eqref{fun.pt.eq} is very far from being an
equivalence. However, in the setting of
Corollary~\ref{enh.adj.fun.corr}, it gives rise to a commutative
square
\begin{equation}\label{fun.pt.adj.sq}
\begin{CD}
\fFun^h(\C,\E) @>{\rho^*}>> \fFun^h(\C',\E)\\
@VVV @VVV\\
\Unf(\Fun(\C_{\ppt},\E_{\ppt})) @>{\Unf(\rho(\ppt)^*)}>>
\Unf(\Fun(\C'_{\ppt},\E_{\ppt})),
\end{CD}
\end{equation}
where the vertical arrows are induced by the comparison functors
\eqref{fun.pt.eq} via Proposition~\ref{trunc.prop}. We then have the
following result.

\begin{lemma}
For a left-reflexive enhanced functor $\rho:\C' \to \C$ between
small enhanced categories, with fully faithful left-adjoint
$\lambda:\C \to \C'$, and for any enhanced category $\E$, the
commutative square \eqref{fun.pt.adj.sq} is cartesian.
\end{lemma}

\proof{} Both horizontal arrows in \eqref{fun.pt.adj.sq} are fully
faithful, so by Corollary~\ref{enh.ff.corr}, it suffices to show
that the comparison functor \eqref{c.01.f} for \eqref{fun.pt.adj.sq}
is essentially surjective over $\ppt \in \Posf^+$. In other words,
we need to show that for any enhanced functor $E:\C' \to \E$ such
that $E(\ppt):\C'_\ppt \to \E_\ppt$ factors through $\C_\ppt$, $E$
itself factors through $\C$. If $\C' = \C \times^h \Unf([1])$ and
$\lambda = s:\C \to \C'$ is induced by $s:\ppt \to [1]$, then
choosing a universal object for $\C$ reduces us to the case $\C =
\Unf(J)$, $J \in \Pos$, and in this case, the claim holds since $\E$
is a separated reflexive family. In general, consider the enhanced
functor $a:\C' \times^h \Unf([1]) \to \C'$ obtained as the
composition
$$
\begin{CD}
\C' \times^h \Unf([1]) @>>> \Cyl_h^\iota(\rho) \cong \Cyl_h(\lambda)
@>{t_\dg}>> \C',
\end{CD}
$$
and note that as we have proved already, $a^*(E)$ factors
through the projection $e:\C' \times^h \Unf([1]) \to \C'$. Therefore
the adjunction map $a^*E \to e^*s^*a^*E$ is an isomorphism, and $E
\cong t^*a^*E \cong t^*e^*s^*a^*E \cong \rho^*\lambda^*E$.
\endproof

\subsection{Enhanced fibrations and cofibrations.}\label{enh.fib.subs}

Just as in Lemma~\ref{enh.adj.co.le} and
Corollary~\ref{enh.adm.sq.corr}~\thetag{i}, if one replaces
cartesian products with semicartesian products, and uses enhanced
cylinders and enhanced comma-categories of
Subsection~\ref{enh.co.subs}, then most of the material of
Subsection~\ref{fib.subs} and Subsection~\ref{cart.subs} extends to
the enhanced setting with basically the same definitions and
proofs. However, it is actually more convenient to use a slightly
more general definition that makes sense for enhanced categories
that are not small.

\begin{defn}\label{enh.fib.def}
An enhanced functor $\pi:\C \to \E$ between enhanced categories is
an {\em enhanced fibration} if the enhanced arrow category
$\Ar^h(\C)$ admits a left-admissible enhanced full subcategory $\E
\setminus^h_\pi \C \subset \Ar^h(\C)$, with some embedding functor
$\nu:\E \setminus^h \C \to \Ar^h(\C)$ such that the commutative
square
\begin{equation}\label{enh.fib.co.sq}
\begin{CD}
\E \setminus^h_\pi \C @>{\tau \circ \nu}>> \C\\
@V{\Ar^h(\pi) \circ \nu}VV @VV{\pi}V\\
\Ar^h(\E) @>{\tau}>> \E
\end{CD}
\end{equation}
is semicartesian. An enhanced functor $\pi$ is an {\em enhanced
  cofibration} if the enhanced-opposite enhanced functor $\pi^\iota$
is an enhanced fibration, and an {\em enhanced bifibration} if it is
both an enhanced fibration and an enhanced cofibration.
\end{defn}

If the enhanced categories $\E$, $\C$ in
Definition~\ref{enh.fib.def} are small, then by
Corollary~\ref{sq.h.corr}, $\E \setminus^h_\pi \C$ must be the
comma-category of \eqref{enh.rl.comma}; this explains our
notation. For enhanced cofibrations, we let $\C /^h_\pi \E =
(\E^\iota \setminus^h_{\pi^\iota} \C^\iota)^\iota$. In general, for
any enhanced fibration $\pi:\C \to \E$, say that a full enhanced
subcategory $\C' \subset \C$ is a {\em enhanced subfibration} if
$\tau^{-1}(\C') \cap (\E \setminus^h \C) \subset \Ar^h(\C')$; in
this case, the induced enhanced functor $\pi':\C' \to \E$ is
trivially an enhanced fibration. For any full enhanced subcategory
$\E' \subset \E$, the induced enhanced functor $\pi':\C' =
\pi^{-1}(\E') \to \E'$ is also trivially an enhanced fibration, with
$\E' \setminus^h_{\pi'} \C' = \Ar^h(\C') \cap (\E \setminus^h_\pi
\C)$.

\begin{lemma}\label{fib.uni.le}
For any enhanced fibration $\pi:\C \to \E$ with small $\E$, $\C$ is
a filtered $2$-colimit of small enhanced subfibrations $\C' \subset
\C$. Conversely, assume given enhanced functor $\pi:\C \to \E$ such
that $\C$ is a filtered $2$-colimit of full subcategories $\C'
\subset \C$ such that the induced functor $\pi':\C' \to \E$ is a
small enhanced fibration. Then $\pi$ is an enhanced fibration, and
all the full subcategories $\C' \subset \C$ are subfibrations.
\end{lemma}

\proof{} For the first claim, for any small enhanced full
subcategory $\C' \subset \C$, let $T(\C') \subset \C$ be the full
enhanced subcategory spanned by $\C'$ and the essential image of the
projection $\sigma:\tau^{-1}(\C') \cap (\E \setminus^h \C) \to \C$,
and define by induction $\C'_0=\C'$, $\C'_{n+1} = T(\C'_n)$, $n \geq
0$. Then $T(\C')$ is still small, and so are all the full
subcategories $\C'_n$ and their union $\C'_\infty = \cup_n\C'_n$
that is obvously an enhanced subfibration. For the second claim,
note that by the uniqueness claim of Corollary~\ref{sq.h.corr}, for
any two subcategories $\C' \subset \C''$ in our fltered $2$-colimit,
with the embedding functor $\eps:\C' \to \C''$, $\Ar^h(\eps)$ sends
$\E \setminus_{\pi'}^h \C'$ into $\E \setminus^h_{\pi''} \C''$, and
we in fact have $\E \setminus^h_{\pi'} \C = \Ar^h(\C') \cap (\E
\setminus_{\pi''} \C''$, so we can define $\E\setminus^h_\pi \C
\subset \Ar^h(\C)$ as the union of $\E \setminus^h_{\pi'} \C'
\subset \Ar^h(\C') \subset \Ar^2(\C)$.
\endproof

\begin{corr}\label{fib.uni.corr}
For any enhanced fibration $\pi:\C \to \E$, the enhanced full
subcategory $\E \setminus^h_\pi \C$ of Definition~\ref{enh.fib.def}
is unique.
\end{corr}

\proof{} If $\E$ and $\C$ are small, this is
Corollary~\ref{sq.h.corr}. In general, $\E$ is a filtered
$2$-colimit of small full enhanced subcategories $\E' \subset \E$,
and then $\C$ is the filtered $2$-colimit of the corresponding
enhanced full subcategories $\C' = \pi^{-1}(\E')$, and $\E
\setminus^h_\pi \C$ is the filtered $2$-colimit of $\E'
\setminus^h_{\pi'} \C'$. Thus it suffices to consider the case when
$\E$ is small, and this reduces to the case when $\C$ is also small
by the first claim of Lemma~\ref{fib.uni.le}.
\endproof

\begin{defn}
For any enhanced fibration $\pi:\C \to \E$, an enhanced morphism in
$\C$ is {\em cartesian over $\E$} if the corresponding enhanced
object in $\Ar^h(\C)$ lies in the unique enhanced full subcategory
$\E \setminus^h_\pi \C \subset \Ar^h(\C)$ of
Corollary~\ref{fib.uni.corr}. For a commutative square
\eqref{C.E.sq} of enhanced categories and functors such that $\pi$,
$\pi'$ are enhanced fibrations, $\gamma$ is {\em cartesian over}
$\phi$ if $\Ar^h(\gamma):\Ar^h(\C') \to \Ar^h(\C)$ sends $\E'
\setminus^h_{\pi'} \C' \subset \Ar^h(\C')$ into $\E
\setminus^h_{\pi'} \C \subset \Ar^h(\C)$ (that is, $\gamma$ sends
enhanced morphisms in $\C'$ cartesian over $\E'$ to enhanced
morphisms in $\C$ cartesian over $\E$). Dually, for an enhanced
cofbration $\pi:\C \to \E$, an enhanced morphism in $\C$ is {\em
  cocartesian over $\E$} iff the opposite enhanced morphism in
$\C^\iota$ is cartesian over $\E^\iota$, and for a commutative
square \eqref{C.E.sq} with $\pi$, $\pi'$ enhanced cofibrations,
$\gamma$ is {\em cocartesian over $\phi$} iff $\gamma^\iota$ is
cartesian over $\phi^\iota$.
\end{defn}

\begin{exa}\label{enh.grp.fib.exa}
For any small enhanced groupoid $\E$, any enhanced functor $\pi:\C
\to \E$ from a small enhanced category $\C$ is an enhanced fibration
and an enhanced cofibration --- indeed, we have $\Ar^h(\E) \cong
\E$, so that we can take $\C /^h_\pi \E \cong \E \setminus^h_\pi \C
\cong \E$.
\end{exa}

\begin{exa}\label{enh.grp.exa}
Say that an enhanced functor $\pi:\C \to \E$ is an {\em enhanced
  family of groupoids} if the commutative square
\begin{equation}\label{enh.fib.grp.sq}
\begin{CD}
\Ar^h(\C) @>{\tau \circ \nu}>> \C\\
@V{\Ar^h(\pi) \circ \nu}VV @VV{\pi}V\\
\Ar^h(\E) @>{\tau}>> \E
\end{CD}
\end{equation}
is semicartesian. Then an enhanced family of groupoids is trivially
an enhanced fibration.
\end{exa}

As in the unenhanced case, an enhanced fibration $\pi:\C \to \E$ is
an enhanced family of groupoids if and only if all enhanced maps in
$\C$ are cartesian over $\E$. In general, while it is not true that
$\pi_\ppt:\C_\ppt \to \E_\ppt$ is a fibration, enhanced maps
cartesian over $\E$ still form a closed class $\flat$ of maps in
$\E_\ppt$, and we can define an enhanced category $\C_{h\flat}$ by
the cartesian square
\begin{equation}\label{h.flat.sq}
\begin{CD}
\C_{h\flat} @>{\nu}>> \C\\
@VVV @VVV\\
\Unf((\C_\ppt)_\flat) @>>> \Unf(\C_\ppt),
\end{CD}
\end{equation}
where as in \eqref{h.iso.sq}, the vertical arrow on the right is the
truncation functor \eqref{trunc.eq}. Then $\C_{h\flat}$ is an
enhanced category by Lemma~\ref{sq.cat.le}, the induced functor
$\pi:\C_{h\flat} \to \E$ is an enhanced family of groupoids, and
$\nu$ in \eqref{h.flat.sq} is cartesian over $\E$. If we have a
commutative square \eqref{C.E.sq} such that $\C' \to \E'$ is an
enhanced family of groupoids, and $\gamma$ is cartesian over $\phi$,
then $\gamma$ uniquely factors through $\nu$. If $\E = \ppt^h$,
then $\flat=\Iso$, and \eqref{h.flat.sq} recovers the enhanced
isomorphism groupoid \eqref{h.iso.sq} and its universal
property. Dually, for any enhanced cofibration $\pi:\C \to \E$
between small enhanced categories, we let
\begin{equation}\label{h.shrp.eq}
  \C_{h\hash} = (\C^\iota_{h\flat})^\iota,
\end{equation}
with the induced enhanced cofibration $\C_{h\sharp} \to \E$ and
functor $\nu:\C_{h\sharp} \to \C$ cocartesian over $\E$. If $\E =
\ppt^h$, then again, $\C_{h\sharp} \cong \C_{h\Iso}$.

\begin{exa}\label{noenh.exa}
For any category $I$, the epivalence \eqref{ar.epi} for the enhanced
category $\C = \Unf(I)$ is an equivalence, and then for any
$\gamma:I' \to I$, we have
\begin{equation}\label{comma.unenh}
\Unf(I) \setminus^h_{\Unf(\gamma)} \Unf(I') \cong \Unf(I)
\setminus^{\Posf^+}_\gamma \Unf(I') \cong \Unf(I \setminus_h I'),
\end{equation}
where the semicartesian product in \eqref{enh.rl.comma} is actually
cartesian by virtue of Lemma~\ref{sq.cat.le}.
Lemma~\ref{cart.fib.le} then shows that $\gamma$ is an fibration if and
only if $\Unf(\gamma)$ is an enhanced fibration. Moreover, if we have two
fibrations $I',I'' \to I$, then a functor $\gamma:I' \to I''$ is
cartesian over $I$ if and only if $\Unf(\gamma)$ is cartesian in the
sense of Definition~\ref{enh.fib.def}.
\end{exa}

As in the unenhanced setting, given enhanced fibrations $\pi:\C \to
\E$, $\pi':\C' \to \E$, we say that an enhanced functor $\gamma:\C'
\to \E$ over $\E$ is {\em cartesian over $\E$} if it is cartesian
over $\id:\E \to \E$, and dually for cofibrations and cocartesian
functors. For any enhanced fibration $\pi:\C \to \E$, with the
left-admisible embedding $\nu:\E \setminus^h_\pi \C$, the
left-adjoint $\nu^\dg$ satisfies $\nu^\dg \circ \nu \cong \id$ since
$\nu$ is full. Then by \eqref{enh.fib.co.sq}, if $\C$ and $\E$ are
small, we actually have $\nu^\dg \cong (\pi \setminus^h \id)$, so
that $\pi$ ia a fibration iff $\pi \setminus^h \id:\Ar^h(\C) \to \E
\setminus^h_\pi \C$ admits a fully faithful enhanced right-adjoint
$(\pi /^h \id)_\dg$. For a square \eqref{C.E.sq} of small enhanced
categories with $\pi$, $\pi'$ enhanced fibrations, $\gamma$ is
cartesian over $\phi$ iff the base change map $\Ar^h(\gamma) \circ
(\pi' \setminus^h \id)_dg \to (\pi \setminus^h \id)_\dg \circ (\phi
\setminus \gamma)$ is an isomorphism. In particular, by
Corollary~\ref{enh.adj.corr}, this immediately implies that the
composition $\phi \circ \gamma:\C' \to \E$ of small enhanced
fibrations $\gamma:\C' \to \C$, $\pi:\C \to \E$ is a fibration, and
$\gamma$ is cartesian over $I$.

\begin{lemma}\label{pi.id.le}
An enhanced functor $\pi:\C \to \E$ between small enhanced
categories is an enhanced fibration if and only if the fully
faithful embedding $\eta:\C \to \E \setminus^h_\pi \C$ of
\eqref{enh.comma.facto} admits a right-adjoint enhanced functor
$\eta_\dg:\E \setminus^h_\pi \C \to \C$ over $\E$. Moreover, for a
commutative square \eqref{C.E.sq} of small enhanced categories and
enhanced functors, with $\pi$, $\pi'$ fibrations, $\gamma$ is
cartesian over $\phi$ iff the base change map $\gamma \circ
\eta'_\dg \to \eta_\dg \circ (\phi \setminus^h \gamma)$ is an
isomorphism.
\end{lemma}

\proof{} For the first claim, $\eta_\dg = \sigma \circ (\pi
\setminus^h \id)_\dg$ is right-adjoint to $\eta \cong (\pi
\setminus^h \id) \circ \eta$, and conversely, if $\pi$ is an
enhanced fibration, then the adjoint $(\pi \setminus^h \id)_\dg$ is
given by the composition
\begin{equation}\label{pi.id.dia}
\begin{CD}
\E \setminus^h_\pi \C @>{\tau^{-1}}>> \C /^h_{\eta,\C} (\E
\setminus^h_\pi \C) @>{\id /^h \tau}>> \C /^h \C \cong \Ar^2(\C),
\end{CD}
\end{equation}
where $\tau:\C /^h_{\eta,\C} (\E \setminus^h_\pi, \C) \to \E
\setminus^h_\pi \C$ is an equivalence by (the dual version of)
Lemma~\ref{enh.adj.co.le}, and $\tau^{-1}$ is the inverse
equivalence. For the second claim, since $\sigma \times
\tau:\Ar^2(\C) \to \C \times^h \C$ is conservative for any $\C$, the
base change map for $(\pi \setminus^h \id)_\dg$ is an isomorphism
iff the same hold for its compositions with $\sigma$ and $\tau$, and
the latter is automatic since $\tau \circ (\pi \setminus^h \id)_\dg
\cong \tau$ by \eqref{pi.id.dia}.
\endproof

As a further slight repackaging of the data described in
Lemma~\ref{pi.id.le}, we note that for any enhanced fibration
$\pi:\C \to \E$ between small enhanced categories, the adjoint
$\eta_\dg: \E \setminus^h_\pi \C \to \E$ of Lemma~\ref{pi.id.le}
defines an enhanced functor
\begin{equation}\label{la.ro.eq}
\rho = \eta_\dg \times \tau:\E \setminus^h_\pi \C \cong \tau^*\C \to
\sigma^*\C \cong \C /^h_\pi \E
\end{equation}
cartesian over $\Ar^h(\E)$, where  $\tau:\tau^*\C \to \Ar^h(\E)$ is
the projection.

\begin{corr}\label{eta.la.corr}
Assume given an enhanced fibration $\pi:\C \to \E$ of small enhanced
categories such that the corresponding enhanced functor
\eqref{la.ro.eq} is left-reflexive over $\Ar^h(\E)$. Then $\pi$ is
an enhanced bifibration.
\end{corr}

\proof{} The embedding $\eta:\C \to \C /^h_\pi \E$ factors as
\begin{equation}\label{eta.la}
\begin{CD}
\C @>{\eta}>> \tau^*\C @>{\rho}>> \sigma^*\C,
\end{CD}
\end{equation}
and $\eta$ in \eqref{eta.la} is right-adjoint to the projection
$\tau:\tau^*\C \to \C$. Thus if $\rho$ has some left-adjoint
$\lambda$ over $\Ar^h(\E)$, then $\tau \circ \lambda$ is
left-adjoint to $\eta:\C \to \C /^h_\pi \E$ over $\E$, and
$\pi^\iota$ is then an enhanced fibration by Lemma~\ref{pi.id.le}.
\endproof

\begin{lemma}\label{cl.fib.le}
For any right-closed embedding $f:J' \to J$ of partially ordered
sets, and any $J$-coaugmented enhanced category $\C$, the pullback
functor $f^*:\sSec^h(J,\C) \to \sSec^h(J',\C)$ is an enhanced
cofibration.
\end{lemma}

\proof{} Let $\chi:J \to [1]$ be the characteristic map of the
right-closed embedding $f$. Then the full embedding $[1]
\setminus_\chi J \subset [1] \times J$ is right-reflexive, with some
adjoint $g:[1] \times J \to [1] \setminus_\chi J$, and if we let
$$
\sSec^h(J,\C) /^h_{f^*} \sSec^h(J',\C) \cong \sSec^h([1]
\setminus_\chi J,\C),
$$
then $g^*:\sSec^h(J,\C) /^h_{f^*} \sSec^h(J',\C) \to
\Ar^h(\sSec^h(J,\C))$ is a right-admissible full embedding by
Lemma~\ref{refl.J.le}, and the corresponding square
\eqref{enh.fib.co.sq} is semicartesian by Lemma~\ref{enh.cocart.le}
and Example~\ref{cyl.coc.exa}.
\endproof

\begin{corr}\label{enh.si.ta.corr}
For any enhanced functor $\pi:\C \to \E$ between small enhanced
categories $\C$, $\E$, the enhanced functor $\sigma$ resp.\ $\tau$
in \eqref{enh.comma.facto} is an enhanced fibration
resp.\ cofibration.
\end{corr}

\proof{} For $\tau$, consider the $[1]$-coaugmented enhanced
cylinder $\Cyl_h(\pi)$, and combine Lemma~\ref{cl.fib.le} for
the embedding $t:\ppt \to [1]$ and Lemma~\ref{cyl.comma.le}. For
$\sigma$, pass to the enhanced-opposite categories.
\endproof

\begin{exa}\label{enh.si.grp.exa}
If $\C$ in Corollary~\ref{enh.si.ta.corr} is an enhanced groupoid, then
$\sigma$ is an enhanced family of groupoids in the sense of
Example~\ref{enh.grp.exa} --- indeed, the adjoints of
Lemma~\ref{enh.2adm.le} are then cartesian functors between families
of groupoids, so both are equivalences.
\end{exa}

\begin{exa}
If we apply Corollary~\ref{enh.si.ta.corr} to $\id:\C \to \C$, for
some enhanced category $\C$, then we obtain a fully faithful
left-admissible embedding $\C \setminus^h_\sigma \Ar^h(\C) \subset
\Ar^(\Ar^h(\C))$ identifying $\C \setminus^h_\sigma \Ar^h(\C) \cong
\Ar^h(\Ar^h(\C)|\C) \cong [2]^o_h\C$. The square \eqref{segal.sq}
for $n=2$, $l=1$ then provides a semicartesian square
\begin{equation}\label{enh.ar.ar.sq}
\begin{CD}
\Ar^h(\Ar^h(\C)|\C) @>{\tau}>> \Ar^h(\C)\\
@V{\beta}VV @VV{\sigma}V\\
\Ar^h(\C) @>{\tau}>> \C,
\end{CD}
\end{equation}
and enhanced version of the cartesian square \eqref{ar.ar.sq}.
\end{exa}

\begin{lemma}\label{enh.fib.ar.le}
For any enhanced fibration $\pi:\C \to \E$ between small enhanced
categories, the relative enhanced arrow category $\Ar^h(\C|\E)
\subset \Ar^h(\C)$ is right-admissible, and we have a diagram with
semicartesian squares
\begin{equation}\label{enh.fib.facto.sq}
\begin{CD}
\Ar^h(\C|\E) @>{\nu}>> \Ar(\C) @>{\nu_\dg}>> \Ar(\C|\E)\\
@V{\tau}VV @V{\pi \setminus \id}VV @VV{\tau}V\\
\C @>{\eta}>> \E \setminus^h_\pi \C @>{\eta_\dg}>> \C,
\end{CD}
\end{equation}
where $\nu$, $\nu_\dg$ are the embedding and its right-adjoint.
\end{lemma}

\proof{} Same as Lemma~\ref{fib.ar.le}, with \eqref{ar.ar.sq}
replaced by \eqref{enh.ar.ar.sq}.
\endproof

\begin{lemma}\label{enh.fib.sq.le}
Assume given a semicartesian square \eqref{cat.dia} of small
enhanced categories and functors, and an enhanced fibration $\pi:\C
\to \E$ such that $\pi_l = \pi \circ \gamma_l$ is an enhanced
fibration, and $\gamma_l$ is cartesian for $l=0,1$. Then $\C_{01}
\to \E$ is an enhanced fibration, and $\gamma^0_{01}$,
$\gamma^1_{01}$ are cartesian over $\E$.
\end{lemma}

\proof{} Same as Lemma~\ref{fib.sq.le}, with
Corollary~\ref{adm.sq.corr} replaced by
Corollary~\ref{enh.adm.sq.corr}~\thetag{i}.
\endproof

\begin{lemma}\label{enh.right.fib.le}
Assume given an enhanced fibration $\C \to \E$ and a
right-admissible enhanced full subcategory $\C' \subset \C$ such
that the enhanced functor $\gamma:\C \to \C'$ right-adjoint to the
embedding is right-adjoint over $\E$. Then $\C' \to \E$ is an
enhanced fibration, and $\gamma$ is cartesian over $\E$.
\end{lemma}

\proof{} Same as Lemma~\ref{right.fib.le}.
\endproof

\begin{lemma}\label{enh.cart.fib.le}
Assume given enhanced fibrations $\pi_l:\C_l \to \E$, $l=0,1$, and
an enhanced functor $\gamma:\C_0 \to \C_1$ cartesian over $\E$. Then
$\sigma$ in \eqref{rel.rlc.facto} is an enhanced fibration, and the
functor $\gamma$ is an enhanced fibration if and only if $\gamma
\setminus^h \id:\Ar^h(\C_0|\E) \to \C_1 \setminus_{\gamma,\E} \C_0$
admits a fully faithful right-adjoint.
\end{lemma}

\proof{} Same as Lemma~\ref{cart.fib.le}, replacing Lemma~\ref{adj.pb.le}
resp.\ Lemma~\ref{right.fib.le} resp.\ \eqref{fib.facto.sq}
by Corollary~\ref{enh.adj.corr} resp.\ Lemma~\ref{right.fib.le}
resp.\ \eqref{enh.fib.facto.sq}.
\endproof

\begin{exa}\label{enh.fam.grp.exa}
If $\C_1 \to \E$ in Lemma~\ref{enh.cart.fib.le} is an enhanced
family of groupoids in the sense of Example~\ref{enh.grp.exa}, then
{\em any} $\gamma:\C_0 \to \C_1$ is an enhanced fibration.
\end{exa}

\begin{lemma}\label{enh.adm.cof.le}
Assume given an enhanced fibration $\pi:\C \to \E$ of small enhanced
categories, and a right-admissible enhanced full subcategory $\E'
\subset \E$. Then the enhanced full subcategory $\C' = \pi^{-1}(\E')
\subset \C$ is right-admissible.
\end{lemma}

\proof{} Let $\lambda:\E' \to \E$ we the embedding functor, with the
right-adjoint enhanced functor $\rho:\E \to \E'$. Then we have $\E
\cong \E' /^h_{\lambda,\E'} \E$ by (the dual version of)
Lemma~\ref{enh.adj.co.le}, and we can consider the composition
\begin{equation}\label{proj.pb}
\begin{CD}
\C \cong (\E' /^h_{\lambda,\E'} \E) \times^h_{\E} \C @>{\nu}>> \E
\setminus^h_\pi \C @>{(\pi \setminus^h \id)_\dg}>> \Ar^h(\C),
\end{CD}
\end{equation}
where $\nu$ is the full embedding induced by $\E' /^h_{\lambda,\E'}
\E \subset \Ar^h(\E)$, and $(\pi \setminus^h \id)_\dg$ is the
adjoint functor provided by Lemma~\ref{pi.id.le}. If we denote the
enhanced functor \eqref{proj.pb} by $\phi$, then $\tau \circ \phi
\cong \id$, and $p = \sigma \circ \phi:\C \to \C$ factors through
$\C' \subset \C$. Moreover, by \eqref{ar.epi}, $\phi$ defines an
enhanced map $a:\id \to p$, and it is clear from the construction
this map is invertible on $\C' \subset \C$, so $\langle p,a \rangle$
is a left projector in the sense of Definition~\ref{proj.def}.
\endproof

\begin{lemma}\label{enh.fib.le}
Assume given a semicartesian square \eqref{C.E.sq}, and assume that
$\pi$ is small. If $\pi$ is an enhanced fibration resp.\ an enhanced
family of groupoids resp.\ an enhanced cofibration, then so is
$\pi'$. Moreover, if we have enhanced fibrations resp.\ cofibrations
$\pi_l:\C_l \to \E$, and an enhanced functor $\gamma:\C_0 \to \C_1$
cartesian resp.\ cocartesian over $\E$, then $\gamma':\C_0
\times^h_{\E} \E' \to \C_1 \times^h_{\E} \E'$ is cartesian
resp.\ cocartesian over $\E'$.
\end{lemma}

\proof{} By replacing $\E'$ with an arbitrary small full enhanced
subcategory $\E'_0$, and then replacing $\E$ with a small full
enhanced subcategory containing the essential image of $\phi$, we
may assume right away that all the enhanced categories in
\eqref{C.E.sq} are small. But then the commutative squares
$$
\begin{CD}
\Ar_h(\C') @>{\pi' \setminus^h \id}>> \E' \setminus^h \C'\\
@V{\gamma \setminus^h \gamma}VV @VV{\phi \setminus^h \gamma}V\\
\Ar_h(\C) @>{\pi \setminus^h \id}>> \E \setminus^h \C,
\end{CD}
\qquad
\begin{CD}
\Ar_h(\C') @>{\id /^h \pi'}>> \C' /^h \E'\\
@V{\gamma /^h \gamma}VV @VV{\gamma /^h \phi}V\\
\Ar_h(\C) @>{\id /^h \pi}>> \C /^h \E
\end{CD}
$$
are semicartesian, and to prove all the claims, it suffices to
combine Corollary~\ref{enh.adj.corr} and Lemma~\ref{pi.id.le}.
\endproof

\begin{exa}\label{ffun.fib.exa}
For any enhanced functor $\pi:\C \to \E$ between small enhanced
categories, and any small enhanced category $\C'$, define an
enhanced category $\fFun^h(\C',\C|\E)$ by the semicartesian square
\begin{equation}\label{fun.enh.rel.sq}
\begin{CD}
\fFun^h(\C',\C|\E) @>>> \fFun^h(\C',\C)\\
@V{\nu}VV @VV{\phi}V\\
\E @>{\theta^*}>> \Fun(\C',\E),
\end{CD}
\end{equation}
where $\phi$ is postcomposition with $\pi$, and $\theta:\C' \to
\ppt^h$ is the tautological projection. Then as soon as $\pi$ is
an enhanced fibration or cofibration, so is $\phi$, by
Corollary~\ref{enh.adj.fun.corr}, and then also $\nu$, by
Lemma~\ref{enh.fib.le}.
\end{exa}

\subsection{Enhanced fibers.}\label{enh.fiber.subs}

As we saw in Example~\ref{noenh.exa}, an enhanced functor $\pi$
between small enhanced categories of the form $\Unf(I)$ is an
enhanced fibration if and only if it is a fibration in the usual
sense. It turns out that this already holds if only the target of
the functor $\pi$ is of the form $\Unf(I)$.

\begin{lemma}\label{noenh.fib.le}
For any enhanced category $\C$ and category $I$, a small enhanced
functor $\pi:\C \to \Unf(I)$ is an enhanced fibration iff it is a
fibration in the usual unenhanced sense, and for two such
fibrations, a functor $\gamma:\C' \to \C$ over $\Unf(I)$ is
cartesian in the sense of Definition~\ref{enh.fib.def} iff it is
cartesian in the usual sense.
\end{lemma}

\proof{} The enhanced comma-category $\Unf(I) \setminus^h_\pi \C$ is
given by \eqref{comma.unenh}, and the functor $\pi \setminus^h \id$
factors as
\begin{equation}\label{ar.epi.dia}
\begin{CD}
\Ar_h(\C) @>{\nu}>> \Ar(\C|\Posf^+) @>{\pi \setminus \id}>> \Unf(I)
\setminus_{\pi}^{\Posf^+} \C,
\end{CD}
\end{equation}
where $\nu$ is the epivalence \eqref{ar.epi}. Therefore if $\pi$ is
an enhanced fibration, it is a usual fibration by
Lemma~\ref{cart.fib.le} and
Corollary~\ref{a.b.adj.corr}. Conversely, assume that $\pi$ is a
fibration. Then both $\Ar_h(\C)$ and $\Unf(I) \setminus^h_{\pi} \C$
are fibered over $\Ar_h(\Unf(I)) \cong \Unf(\Ar(I))$, and $\pi
\setminus^h \id$ is cartesian over $\Unf(\Ar(I))$, so by
Lemma~\ref{fib.adj.le}~\thetag{ii}, it suffices to find a functor
$(\pi \setminus^h \id)_\dg:\Unf \setminus^h_{\pi} \C \to \Ar_h(\C)$
cartesian over $\Unf(\Ar(I))$, and an isomorphism $(\pi \setminus^h
\id) \circ (\pi \setminus^h \id)_\dg \cong \id$ that defines an
adjunction between $\pi \setminus^h \id$ and $(\pi \setminus^h
\id)_\dg$ restricted to the fiber over any object in
$\Unf(\Ar(I))$. Explicitly, such an object is given by a pair
$\langle J,\alpha \rangle$, $J \in \Posf^+$, $\alpha:J^o \times [1]
\to I$, and $\pi \setminus^h \id$ restricted to the
corresponding fiber is the functor
\begin{equation}\label{sec.t.eq}
t^*:\Sec^h(J^o \times [1],\Unf(\alpha)^*\C) \to \Sec^h(J^o,\Unf(t
\circ \alpha)^*\C).
\end{equation}
Then the fully faithful right-adjoint to \eqref{sec.t.eq} is given
by Lemma~\ref{enh.2adm.le}, and it is functorial with respect
to $\langle I,\alpha \rangle$. This proves the first claim; for the
second, recall again that $\nu$ is \eqref{ar.epi.dia} is an
epivalence, thus conservative, and apply the last claim of
Lemma~\ref{cart.fib.le}.
\endproof

In particular, Lemma~\ref{noenh.fib.le} shows that for any small
enhanced category $\C$ and $I \in \Posf^+$, an enhanced fibration
$\C \to \Unf(I)$ is the same thing as an $I$-augmentation in the
sense of Definition~\ref{cat.aug.def}. Recall that for any small
enhanced functor $\C \to \E$ and object $e \in \E_{\ppt}$, we have
the enhanced fiber $\C_e = \ppt^h \times^h_{\E} \C$ of
\eqref{enh.fib.sq}; then by Lemma~\ref{enh.fib.le}, this means for
any small enhanced fibration $\C \to \E$ and object $e \in \E_I
\subset \E$, $I \in \Posf^+$, the generalized enhanced fiber $\C_e$
is an $I$-augmented enhanced category. This can be used to define
transition functors for small enhanced fibrations and cofibrations.

Namely, if we have two objects $e,e' \in \E_\ppt$, then a morphism
$f:e \to e'$ lifts to an enhanced morphism $f^h:\Unf([1]) \to \E$,
and then if $\C \to \E$ is a small enhanced fibration, $\C_f =
\Unf([1]) \times^h_{\E} \C$ is $[1]$-augmented by
Lemma~\ref{enh.fib.le}, so that $\C_f \cong \Cyl^o_h(f^*)$ for an
{\em enhanced transition functor} $f^*:\C_{e'} \to \C_{e}$, unique
up to an isomorphism, as in \eqref{trans.I.eq}. If $\pi$ is a small
enhanced cofibration, then analogously, $\C_f \cong \Cyl_h(f_!)$ for
a unique enhanced transition functor $f_!:\C_e \to \C_{e'}$. If we
have a composable pair of maps $f$, $g$, then $(f \circ g)_! \cong
f_! \circ g_!$ and $(f \circ g)^* \cong g^* \circ f^*$. Moreover, if
we have two small enhanced fibrations $\C,\C' \to \E$ and a functor
$\gamma:\C \to \C'$ over $\E$, then for any map $f:e \to e'$ in
$\E_\ppt$, the induced functor $\gamma_f:\C_f \to \C'_f$ comes
equipped with isomorphism $\gamma_f \circ s \cong s$, $\gamma_f
\circ t \cong t$, and then \eqref{enh.cyl.facto} induces a base
change map
\begin{equation}\label{enh.fu.fib.E}
  \gamma_e \circ f^* \to f^* \circ \gamma_{e'},
\end{equation}
an analog of \eqref{fu.fib} and \eqref{enh.fu.fib}. By definition,
$\gamma$ is cartesian over $\E$ iff all the functors $\gamma_f$ are
$[1]$-augmented, and this happens iff all the maps
\eqref{enh.fu.fib.E} are invertible, so that we have commutative
squares
\begin{equation}\label{enh.E.fu.fib.sq}
\begin{CD}
\C_{e'} @>{\gamma_i}>> \C'_{e'}\\
@V{f^*}VV @VV{f^*}V\\
\C_e @>{\gamma_i}>> \C_{e'}
\end{CD}
\end{equation}
in $\Cat^h$, a generalization of \eqref{enh.E.fu.fib.sq}. The story
for cofibrations is dual.

\begin{lemma}\label{enh.ff.fib.le}
Assume given a small enhanced functor $\pi:\C \to \E$, and a full
enhanced subcategory $\C' \subset \C$ such that $\pi$ is an enhanced
fibration resp.\ cofibration, and for any morphism $f:e \to e'$ in
$\E_\ppt$, the enhanced transition functor $f^*:\C_{e'} \to \C_e$
resp.\ $f_!:\C_e \to \C_{e'}$ sends $\C'_{e'} \subset \C_{e'}$ into
$\C'_e \subset \C_e$ resp.\ $\C'_e$ into $\C'_{e'}$. Then the
induced functor $\pi':\C' \to \E$ is an enhanced fibration
resp.\ cofibration, and the embedding functor $\C' \to \C$ is
cartesian resp.\ cocartesian over $\E$. Conversely, if we have
enhanced fibrations resp.\ cofibrations $\C,\C' \to \E$ and an
enhanced functor $\nu:\C' \to \C$ cartesian resp.\ cocartesian over
$\E$, then $\nu$ is fully faithful resp.\ an equivalence if and only
if so are its fibers $\gamma_e:\C'_e \to \C_e$ for all enhanced
objects $e \in \E_\ppt$.
\end{lemma}

\proof{} The claims for fibrations and cofibrations are equivalent,
by passing to enhanced-opposite categories, so we only consider the
fibrations.  For the first claim, it suffices to see that the
adjoint $(\pi \setminus \id)_\dg:\E \setminus^h \C \to \Ar_h(\C)$
sends the full enhanced subcategory $\E \setminus^h \C' \subset \E
\setminus^h \C$ into the full enhanced subcategory $\Ar_h(\C')
\subset \Ar_h(\C)$. Since $(\pi \setminus \id)_\dg$ is fully
faithful, Corollary~\ref{enh.ff.corr} shows that it suffices to
check this over $\ppt \in \Posf^+$. But then up to an isomorphism,
objects in $(\E \setminus^h \C')_\ppt$ are triples $\langle
e,c,\alpha \rangle$, $e \in \E_\ppt$, $c \in \C'_\ppt$, $\alpha:e
\to \pi(c')$, and $(\pi \setminus \id)_\dg$ sends such a triple to
the arrow $\alpha^*(c) \to c$ that lies in $\Ar_h(\C') \subset
\Ar_h(\C)$ by assumption. For the second claim, the ``only if'' part
is trivial, and Lemma~\ref{enh.ff.loc.le} reduces the ``if'' part to
the case $\E = \Unf(I)$, $I \in \Posf^+$, so that $\C$, $\C'$ and
$\gamma$ are $J$-augmented. Moreover, the enhanced essential image
of $\gamma$ is then also $J$-augmented by the first claim, so we may
replace $\C$ with this essential image and assume that $\gamma_i$ is
an equivalence for any $i \in I$. Then since $\cCat^h$ is an
enhanced category, it is a non-degenerate family of categories over
$\Posf^+$, and this finishes the proof.
\endproof

\begin{lemma}\label{enh.fib.adj.le}
Assume given two enhanced fibrations $\C,\C' \to \E$ between small
enhanced categories, and an enhanced functor $\gamma:\C \to \C'$
cartesian over $\E$. Then if for any $e \in \E_\ppt$, $\gamma_e$ is
left-reflexive, $\gamma$ itself is left-reflexive over $\E$, and if
for any map $f:e \to e'$ in $\E_\ppt$, the square
\eqref{enh.E.fu.fib.sq} is left-reflexive, then the adjoint enhanced
functor $\gamma_\dg$ is cartesian over $\E$.
\end{lemma}

\proof{} By Lemma~\ref{enh.cart.fib.le},
Lemma~\ref{adj.eq.le}~\thetag{iii} and \eqref{enh.comma.E.facto}, it
suffices to prove the claim after replacing $\C$ with $\C'
\setminus^h_{\gamma,\E} \C$ and $\E$ with $\C$; in other words, we
may assume right away that $\C' = \E$ and $\gamma$ is an enhanced
fibration. We note that $\C_e$ has an initial enhanced object for
any $e \in \E_\ppt$, and we need to show that $\gamma$ admits a
fully faithful left-adjoint $\gamma_\dg$, cartesian if all the
squares \eqref{enh.E.fu.fib.sq} are left-reflexive. Let $\C' \subset
\C$ be the full enhanced subcategory spanned by the inital objects
$o \in \C_{e,\ppt}$, $e \in \E_\ppt$. Then to construct the adjoint
$\gamma_\dg$, it suffices to show that firstly, the induced functor
$\C_0 \to \E$ is an equivalence, so that the inverse equivalence
defines a fully faithful enhanced section $\gamma_\dg:\E \to \C$ of
the functor $\gamma$, and secondly, that $\tau:\E
\setminus^h_{\gamma_\dg,\E} \C \to \C$ is an equivalence, so that
$\gamma_\dg$ is adjoint to $\gamma$ by Lemma~\ref{enh.adj.co.le}. By
Lemma~\ref{enh.ff.loc.le}, both statements can be checked after
passing to generalized enhanced fibers of the fibration $\gamma:\C
\to \E$. By definition, checking whether the adjoint $\gamma_\dg$ is
cartesian can also be done after passing to enhanced fibers, so
altogether, we are reduced to the case when $\E = \Unf(J)$, $J \in
\Posf^+$, and $\C$ is a $J$-augmented small enhanced category. Then
$\gamma$ is a fibration in the usual sense, so by
Lemma~\ref{fib.adj.le}~\thetag{i}, it suffices to check that all its
fibers have an inital object. By \eqref{sec.h.eq}, this amounts to
proving the following:
\begin{itemize}
\item for any $J \in \Posf^+$ and $J$-augmented small enhanced
  category $\C$ such that $\C_j$ has an initial enhanced object
  for any $j \in J$, the category $\sSec^h(J,\C)$ also has an
  initial enhanced object preserved by the evaluation functor
  \eqref{nondeg.sec}.
\end{itemize}
Now as in Lemma~\ref{ff.cat.le}, say that $J \in \Posf^+$ is {\em
  good} if it satisfies \thetag{$\bullet$} for any $\C$. Then $J =
\ppt$ is trivially good. Moreover, since \eqref{nondeg.sec} is
conservative by Lemma~\ref{sec.sq.le}, for any map $f:J_0 \to J_1$
between good $J_0,J_1 \in \Posf^+$, the pullback functor
$f^*:\sSec^h(J_1,-) \to \sSec^h(J_0,-)$ preserves the initial
enhanced objects required in \thetag{$\bullet$}. Then applying
Corollary~\ref{enh.adm.sq.corr}~\thetag{i} to semicartesian squares
\eqref{sec.J.sq}, we see that for any square \eqref{J.sq} of
Lemma~\ref{sec.sq.le} with good $J_0$, $J_1$, $J_{01}$, the set $J$
is also good. Then the same induction as in Lemma~\ref{ff.cat.le}
reduces us to showing that $J^>$ is good for a good $J \in \Posf
\subset \Posf^+$. More generally, assume given good $J_0,J_1 \in
\Posf$ a map $f:J_0 \to J_1$, and let $J = \Cyl(f)$. Then by
Lemma~\ref{cl.fib.le}, $s^*:\sSec^h(J,\C)
\to \sSec^h(J_0,\C)$ is an enhanced fibration for any small
$J$-augmented enhanced category $\C$, and then by
Lemma~\ref{enh.adm.cof.le}, an enhanced initial object in
$\sSec^h(J_0,\C)$ induces a left-admissible full embedding
$\sSec^h(J_1,\C) \to \sSec^h(J,\C)$. An enhanced initial object in
$\sSec^h(J_1,\C)$ then gives an enhanced initial object in
$\sSec^h(J,\C)$, so that $J$ is good.
\endproof

\begin{corr}\label{adj.comma.corr}
An enhanced functor $\gamma:\C \to \E$ between small enhanced
categories is left resp.\ right-reflexive iff for any enhanced
object $e$ in $\E$, the right resp.\ left enhanced comma-fiber $e
\setminus^h_\gamma \C$ resp.\ $\C /^h_\gamma e$ has an initial
resp.\ terminal enhanced object.
\end{corr}

\proof{} For the first claim, apply Lemma~\ref{enh.fib.adj.le} to
the enhanced fibration $\sigma:\E \setminus^h_\gamma \C \to \E$ of
Corollary~\ref{enh.si.ta.corr}, and recall that by
Lemma~\ref{adj.eq.le}~\thetag{iii}, $\gamma$ is left-reflexive if
and only if so is $\sigma$. For the second claim, apply the first
claim to the enhanced-opposite functor $\gamma^\iota$.
\endproof

\begin{corr}\label{enh.bifib.corr}
An enhanced fibration $\pi:\C \to \E$ between small enhanced
categories is an enhanced bifibration if and only if for any enhanced
morphism $f:e \to e'$ in $\E_\ppt$, the transition functor
$f^*:\C_{e'} \to \C_e$ is left-reflexive.
\end{corr}

\proof{} For the ``only if'' part, it suffices to consider the case
$\E = \Unf([1])$ where the claim holds more-or-less by the
definition of the enhanced transition functors. For the ``if'' part,
Lemma~\ref{enh.fib.adj.le} over $\Ar^h(\E)$ shows that $\pi$
satisfies the assumptions of Corollary~\ref{eta.la.corr}.
\endproof

We can now give example of enhanced fibrations and cofibrations
between enhanced categories that are not small. Namely, assume given
a small enhanced category $\E$, and consider the enhanced categories
$\E \bbd^h \cCat^h$, $\cCat^h \bb^h \E$ of
Subsection~\ref{lax.subs}, with the functors \eqref{no.E.d} and
\eqref{no.E}.

\begin{prop}\label{no.E.prop}
For any small enhanced category $\E$, the enhanced functor
\eqref{no.E.d} is an enhanced cofibration, with enhanced transition
functors $\phi_!:\fFun(\E,\C) \to \fFun^h(\C')$ corresponsing to
$\phi:\C \to \C'$ given by postcomposition with $\phi$, and the
enhanced functor \eqref{no.E} is an enhanced fibration, with
enhanced transition functors $\phi^*:\fFun^h(\C',\E) \to
\fFun^h(\C,\E)$.
\end{prop}

\proof{} For the first claim, as in the proof of
Proposition~\ref{ccat.J.prop}, choosing a universal object for the
small enhanced category $\E$ provides a fully faithful embedding $\E
\bbd^h \cCat^h \to I^o \bbd^h \cCat^h$ for some $I \in \Posf^+$, and
then Lemma~\ref{enh.ff.fib.le} reduces us to the case $\E =
\Unf(I^o)$, $\E \bbd^h \cCat^h = I^o \bbd^h \cCat^h$. To simplify
notation, let $\Kk = I^o \bbd^h \cCat^h$, and let $\pi:\Kk \to
\cCat^h$ be the enhanced functor \eqref{no.E.d}. Its fibers are
given by $\Kk_{J,\C} = \Kk_{\langle J,\C \rangle} \cong \Sec^h((I
\times J)^o,\C_{h\perp})$, where $\langle J,\C \rangle \in \cCat^h$
is a pair of some $J \in \Posf^+$ and a $J$-augmented small enhanced
category $\C$. Objects in the enhanced arrow category
$\Ar_h(\cCat^h)$ are pairs $\langle J,\C \rangle$, $J \in \Posf^+$,
$\C$ a $(J \times [1])$-augmented category, and the fibers of
$\Ar_h(\pi):\Ar_h(\Kk) \to \Ar_h(\cCat^h)$ are then given by
$\Ar_h(\Kk)_{J,\C} = \Kk_{J \times [1],\C}$. The embedding $t:J \to
J \times [1]$ gives a restriction functor $\rho =
t^*:\Ar_h(\Kk)_{J,\C} \to \Kk_{J,t^*\C}$, and
Lemma~\ref{enh.2adm.le} provides a fully faithful functor
\begin{equation}\label{a.a0.eq}
\lambda:\Kk_{J,t^*\C} \to \Ar_h(\Kk)_{J,\C}
\end{equation}
left-adjoint to $\rho$. We now let $\Ar_h(\Kk)^0 \subset \Ar_h(\Kk)$
be the full subcategory spanned by the essential images of all the
functors \eqref{a.a0.eq}, and we observe that the claim reduces to
the following two statements:
\begin{enumerate}
\item $\Ar_h(\Kk)^0 \subset \Ar_h(\Kk)$ is a full enhanced
  subcategory that admits an enhanced functor $\Ar_h(\Kk) \to
  \Ar_h(\Kk)^0$ right-adjoint to the embedding, and
\item the projection $\id / \pi:\Ar_h(\Kk)^0 \subset \Ar_h(\Kk) \to
  \Kk /^h \cCat^h$ is an equivalence.
\end{enumerate}
For \thetag{i}, note that for any object $\alpha \in
\Ar_h(\Kk)_{J,\C}$, we have the adjunction map $a(\alpha):\alpha \to
\lambda(\rho(\alpha))$, and $\alpha \in \Ar_h(\Kk)^0$ iff
$a(\alpha)$ is an isomorphism. Since $\rho$, $\lambda$ and $a$
commute with pullbacks with respect to maps $f:J' \to J$, and
$\Ar_h(\Kk)$ is an enhanced category, thus non-degenerate, the
condition can be checked after restriction with respect to all maps
$\ppt \to J$. This means that $\Ar_h(\Kk)^0 \subset \Ar_h(\Kk)$ is
indeed the full enhanced subcategory corresponding to the full
subcategory $\Ar_h(\Kk)^0_\ppt \subset \Ar_h(\Kk)$ by
Corollary~\ref{enh.ff.corr} (and in particular, $\Ar_h(\Kk)^0$ is an
enhanced category). Moreover, for any $\langle J,\C' \rangle$ in
$\Ar_h(\cCat^h)$ and $(J \times [1])$-aug\-mented enhanced functor
$\phi:\C \to \C'$, we have a base change isomorphism
$b(\phi,\alpha):\lambda(\rho(\phi_{h\perp} \circ \alpha)) \cong
\phi_{h\perp} \circ \lambda(\rho(\alpha))$, functorial with respect
to $\phi$ and $J$. Therefore sending $\langle J,\C,\alpha \rangle$
to $\langle J,\C,\lambda(\rho(\alpha)) \rangle$ and $\langle
\id,\phi,a \rangle:\langle J,\C,\alpha \rangle \to \langle
J,\C',\alpha' \rangle$ to $\langle \id,\phi,\lambda(\rho(a)) \circ
b(\phi,\alpha)\rangle$ provides a well-defined functor $\lambda \circ
\rho:\Ar_h(\Kk)_J \to \Ar^h(\Kk)^0_J \subset \Ar_h(\Kk)_J$, and the
maps $\langle \id,\id,a(\alpha) \rangle:\langle J,\C,\alpha \rangle
\to \langle J,\C,\lambda(\rho(\alpha))\rangle$ give a functorial map
$a:\id \to \lambda \circ \rho$ with invertible $\lambda \circ
\rho(a)$. Thus $\lambda \circ \rho$ is a projector on $\Ar_h(\Kk)_J$
in the sense of Definition~\ref{proj.def}, and since everything is
functorial with respect to $J$, we obtain a projector $\Ar_h(\Kk)
\to \Ar_h(\Kk)^0 \subset \Ar_h(\Kk)$.

For \thetag{ii}, since we already know that $\Ar_h(\Kk)^0$ is an
enhanced category, the claim amounts to checking that
\begin{equation}\label{epi.ga.k}
\gamma = \id / \pi:\Ar_h(\Kk)^0 \to \Kk /^{\Posf^+} \cCat^h
\end{equation}
is an epivalence. By Lemma~\ref{epi.le}, since $\gamma$ is cartesian
over $\Posf^+$, it suffices to check this over each $J \in
\Posf^+$. An object in the target of the functor $\gamma$ is given
by $\langle J \times [1],\C \rangle \in \cCat^h$, $\langle
J,\C',\alpha\rangle \in \Kk$ and a $J$-augmented equivalence
$\eps:\C' \cong t^*\C$, and then $\langle J,\C,\lambda(\alpha \circ
\eps)\rangle$ is an object in $\Kk^0$, so that \eqref{epi.ga.k} is
essentially surjective. For any two objects $\langle J,\C,\alpha
\rangle$, $\langle J,\C',\alpha' \rangle$ in $\Ar_h(\Kk)^0_J$, a
morphism between their images under $\gamma$ is represented by a
pair of a $(J \times [1])$-augmented enhanced functor $\phi:\C \to \C'$ and a
morphism $\langle \id,\phi',a'\rangle:\langle J,t^*\C,t^*\alpha
\rangle \to \langle J,t^*\C',t^*\alpha'\rangle$ such that $\phi'$ is
isomorphic to $t^*\phi$, and then choosing such a isomorphism and
replacing $\phi'$ with $t^*\phi$ defines a morphism $\langle
\id,\phi,\lambda(a') \rangle:\langle J,\C,\alpha \rangle \to \langle
J,\C',\alpha' \rangle$ in $\Ar_h(\Kk)^0$, so that $\gamma$ is
full. Finally, since $\lambda$ is fully faithful, $a$ is invertible
if and only if so is $\lambda(a)$, and this shows that $\gamma$ is
conservative.

This finishes the proof of the first claim of the Proposition. The
argument for the second claim is similar, so we only give a
sketch. Objects in the enhanced arrow category $\Ar_h(\cCat^h \bb^h
\E)$ are triples $\langle I,\C,\alpha \rangle$ of some $I \in
\Posf^+$, an $(I \times [1])$-augmented small enhanced category
$\C$, and an enhanced functor $\alpha:\C_{h\perp} \to \E$. The
embedding $t:I \to I \times [1]$ gives rise to a restriction functor
$\lambda=t^*:\Fun^h(\C_{h\perp},\E) \to
\Fun^h(t^*\C_{h\perp},\E)$. As in the proof of
Proposition~\ref{ccat.E.prop}, choosing a universal $(I \times
[1])$-augmented object $J$ for $\C$ provides a fully faithful
embedding $\Fun^h(\C_{h\perp},\E) \subset J^o_h\E$, and then
Lemma~\ref{enh.2adm.le} provides a fully faithful functor
$\rho:\Fun^h(t^*\C_{h\perp},\E) \to \Fun^h(\C_{h\perp},\E)$
right-adjoint to $\lambda$. Then again, let $\Ar_h(\cCat^h \bb^h
\E)^0 \subset \Ar_h(\cCat^h \bb^h \E)$ be the full subcategory
spanned by essential images of the functors $\rho$ for all $\langle
J,\C \rangle \in \Ar_h(\cCat^h)$, and observe that this is a full
enhaced subcategory, and $\rho \circ \lambda$ induces an enhanced
functor $\Ar_h(\cCat^h \bb^h \E) \to \Ar_h(\cCat^h \bb^h \E)^0$
left-adjoint to the full embedding $\Ar_h(\cCat^h \bb^h \E)^0 \to
\Ar_h(\cCat^h \bb^h \E)$. Finally, by the same argument as in the
first claim, the projection $\Ar_h(\cCat^h \bb^h \E) \to \cCat^h
\setminus^h (\cCat^h \bb^h \E)$ restricts to an equivalence between
its target and $\Ar_h(\cCat^h \bb^h \E)^0$.  \endproof

\begin{remark}
Even if an enhanced category $\E$ is not small, it is a filtered
$2$-colimit of its small enhanced full subcategories $\E' \subset
\E$, and then $\cCat^h \bb^h \E$ is a filtered $2$-colimit of
enhanced full subcategories $\cCat^h \bb^h \E'$. Thus by the second
claim of Lemma~\ref{fib.uni.le}, Proposition~\ref{no.E.prop} implies
that \eqref{no.E} is an enhanced fibration for any enhanced category
$\E$.
\end{remark}

\begin{exa}\label{ccat.I.fib.exa}
By \eqref{lax.pres.dia}, the enhanced fiber of the enhanced
cofibration $\E \bbd^h \cCat^h \to \cCat^h$ over some $\C \in
\Cat^h$ is the enhanced category $\fFun^h(\E,\C)$. For example, if
$\E = \ppt^h$, so that $\E \bbd^h \cCat^h = \cCat^h_\idot$ is
the enhanced category of Example~\ref{enh.idot.exa}, this fiber is
$\C$ itself. More generally, for any $\langle I,\C \rangle \in
\cCat^h$, we have a semicartesian square
\begin{equation}\label{ccat.fib.sq}
\begin{CD}
\C @>>> \cCat^h_\idot\\
@VVV @VVV\\
\Unf(I) @>>> \cCat^h,
\end{CD}
\end{equation}
where the bottom arrow corresponds to $\langle I,\C \rangle$, and
the top arrow sends an object in $\C$ represented by a pair $\langle
I',s \rangle$, $s \in \Sec^h(I',\C)$, $f:I' \to I$ a map in $\Posf^+$
to $\langle I',f^*\C,s \rangle \in \cCat^h_\idot$. Thus the
generalized enhanced fiber \eqref{enh.fib.sq} of $\cCat^h_\idot$
over $\langle I,\C \rangle \in \cCat^h$ is again $\C$.
\end{exa}

\subsection{Enhanced Grothendieck construction.}\label{enh.groth.subs}

Recall that for any enhanced category $\E$, we have a dense
subcategory $\Cat^h \bbi^h \E \subset \Cat^h \bb^h \E$ whose
morphisms are functors over $\E$. If $\E$ is small, we can then
consider the subcategory $\Cat^h \mmh \E \subset \Cat^h \bbi^h \E$
whose objects are small enhanced categories $\C$ equipped with an
enhanced fibration $\C \to \E$, and whose morphisms are isomorphism
classes of functors cartesian over $\E$. The enhancement $\cCat^h
\bb^h \E$ for the category $\Cat^h \bb^h \E$ induces enhancements
for $\Cat^h \mmh \E \subset \Cat^h \bbi^h \E \subset \Cat^h \bb^h
\E$, by cartesian squares
\begin{equation}\label{3.h.enh.dia}
\begin{CD}
\cCat^h \mmh \E @>{a}>> \cCat^h \bbi^h \E @>{b}>> \cCat^h \bb^h \E\\
@VVV @VVV @VVV\\
\Unf(\Cat^h \mmh \E) @>>> \Unf(\Cat^h \bbi^h \E) @>>> \Unf(\Cat^h \bb^h
\E),
\end{CD}
\end{equation}
where the vertical arrow on the right is the truncation functor
\eqref{trunc.eq}. We can also consider the twisted version $\Cat^h
\bb^h_\iota \E$ of the category $\Cat^h \bb^h \E$; then it has an
enhancement $\cCat^h \bb^h_\iota \E \cong \iota^*(\cCat^h \bb \E)$,
where $\iota$ stands for the enhancement \eqref{iota.enh.eq} of the
involution $\C \mapsto \C^\iota$, and then all the other categories
in \eqref{3.h.enh.dia} can be similarly twisted by $\iota$; as in
Example~\ref{iota.bb.exa}, we have $\cCat^h \bbi^h_\iota \E \cong
\cCat^h \bbi^h \E^\iota$.

\begin{lemma}\label{al.be.le}
For any enhanced category $\E$, with enhanced isomorphism groupoid
$\E_{h\Iso}$ and embedding $\nu:\E_{h\Iso} \to \E$, the enhanced
category $\cCat^h \bbi^h \E$ of \eqref{3.h.enh.dia} is equivalent to
the category of triples $\langle I,\C,\beta \rangle$, $I \in
\Posf^+$, $\C$ a small enhanced category, $\beta:\C \to \Unf(I)
\times^h \E$ an enhanced functor such that $(\id \times \nu)^*\C \to
\Unf(I) \times^h \E_{h\Iso}$ is an enhanced fibration. Morphisms
from $\langle I',\C',\beta' \rangle$ to $\langle I,\C,\beta \rangle$
are represented by pairs $\langle f,\phi \rangle$ of a map $I' \to
I$ and a commutative square
\begin{equation}\label{ccat.E.sq}
\begin{CD}
  \C' @>{\phi}>> \C\\
  @V{\beta'}VV @VV{\beta}V\\
  \Unf(I') \times^h \E @>{\Unf(f) \times \id}>> \Unf(I) \times^h \E
\end{CD}
\end{equation}
such that $\phi$ is cartesian over $\Unf(f)$; pairs $\langle f,\phi
\rangle$ are considered modulo isomorphism over $\E$. If
$\E$ is small, $\cCat^h \mmh \E \subset \cCat^h \bbi^h \E$ is then
equivalent to the subcategory of triples $\langle I,\C,\beta
\rangle$ such that $\beta:\C \to \Unf(I) \times^h \E$ is an enhanced
fibration, with morphisms given by squares \eqref{ccat.E.sq} such
that $\phi$ is cartesian over $\Unf(f) \times^h \id$.
\end{lemma}

\proof{} For the first claim, note that for any $\langle I,\C,\alpha
\rangle \in \cCat^h \bb^h \E$, with $I$-augmentation $\pi:\C \to
\Unf(I)$, the enhanced functor $\pi \times \alpha:\C_{h\perp} \to
\Unf(I) \times^h \E$ is coaugmented iff all the corresponding maps
\eqref{enh.fu.fib} are isomorphisms. Spelling out the construction
of the functor \eqref{nu.J.eq} for the reflexive family $\cCat^h
\bb^h \E \to \Posf^+$ given in the proof of
Proposition~\ref{ccat.E.prop}, we see that this happens if and only
if $\langle I,\C,\alpha \rangle$ lies in $\cCat^h \bbi^h \E \subset
\cCat^h \bb^h \E$. We can then take $\beta = (\pi \times
\alpha)_{h\perp}$, and the claim for morphisms is also clear. For
the second claim, note that by Lemma~\ref{enh.cart.fib.le},
$\beta:\C \to \Unf(I) \times^h \E$ is an enhanced fibration iff
$\Ar^h(\C|\Unf(I)) \to \E \setminus^h \C$ has a fully faithful
$I$-augmented right-adjoint, and then applying
Lemma~\ref{enh.fib.adj.le} to the transpose-opposite $I$-augmented
enhanced categories, we see that this happens if $\C_i \to \E$ is an
enhanced fibration for any $i \in I$, and all the transition
functors \eqref{trans.I.eq} are cartesian over $\E$. The claim for
morphisms is again clear.
\endproof

For any enhanced category $\E$, Lemma~\ref{sq.h.le} and
Corollary~\ref{sq.h.corr} show that the category $\Cat^h \bbi^h \E$ has
products (given by $- \times^h_{\E} -$). However, there is more
functoriality. Assume that $\E$ is small, and we have a small
enhanced functor $\pi:\C' \to \C$. Then for any enhanced functor
$\gamma:\E \to \C$, we can form a semicartesian square
\begin{equation}\label{E.C.s.sq}
\begin{CD}
\E'(\gamma) @>>> \C'\\
@V{\beta(\gamma)}VV @VV{\pi}V\\
\E @>{\gamma}>> \C.
\end{CD}
\end{equation}
If we have another enhanced functor $\gamma':\E \to \C$ and an
isomorphism $b:\gamma \to \gamma'$, then replacing $\gamma$ in
\eqref{E.C.s.sq} with $\gamma'$ and composing the connecting
isomorphism with $a \circ \beta$ defines another commutative square,
and Corollary~\ref{sq.h.corr} then provides an enhanced functor
$\phi(b):\E'(\gamma) \to \E'(\gamma')$ equipped with an isomorphism
$a(b):\beta(\gamma) \to \beta(\gamma') \circ \phi$ (note that even
if $\gamma'=\gamma$ but $b$ is non-trivial, $\phi(b)$ might well be
a non-trivial endofunctor of $\E'(\gamma)=\E'(\gamma')$). The pair
$\langle \phi(b),a(b) \rangle$ is unique up to an isomorphism, thus
gives a well-defined map in $\Cat^h \bbi^h \E$, and sending $\gamma$ to
$\langle \E'(\gamma),\beta(\gamma) \rangle$ and $b$ to $\langle
\phi(b),a(b) \rangle$ defines a functor from the isomorphism
groupoid $\Fun^h(\E,\C)_{\Iso}$ to $\Cat^h \bbi^h \E$. Passing to the
enhanced-opposite categories, we get a functor
$\Fun^h(\E^\iota,\C)_{\Iso} \to \Cat^h \bbi^h \E$, $\gamma \mapsto
\E'(\gamma)^\iota$, and if $\pi:\C' \to \C$ is an enhanced
cofibration, we can refine this construction to an enhanced functor
$\E^\iota\C = \fFun^h(\E^\iota,\C) \to \cCat^h \mmh \E$. Namely,
explicitly, objects in $\E^\iota\C$ are pairs $\langle I,\gamma
\rangle$, $I \in \Posf^+$, $\gamma:\E^\iota \times^h \Unf(I^o) \to
\C$ an enhanced functor, and then for any such $\gamma$,
$\beta(\gamma)^\iota:\E'(\gamma)^\iota \to \E \times^h \Unf(I)$ is
an enhanced fibration by Lemma~\ref{enh.fib.le}, thus defines an
object in $\cCat^h \mmh \E$ by Lemma~\ref{al.be.le}. Then again,
Corollary~\ref{sq.h.corr} shows that this is functorial with respect
to maps in $\E^\iota\C$ cartesian over $\Posf^+$, and this then
extends to an enhanced functor
\begin{equation}\label{enh.pb.fun.eq}
\E^\iota\C \to \cCat^h \mmh \E, \qquad \langle I,\gamma \rangle
\mapsto \langle I,\E'(\gamma)^\iota,\beta(\gamma)^\iota \rangle
\end{equation}
by Lemma~\ref{enh.exp.le}. In particular, we can take the enhanced
cofibration \eqref{no.E.d} of Proposition~\ref{no.E.prop} for the
enhanced category $\cCat^h_\idot = \ppt^h \bbd^h \cCat^h$ of
Example~\ref{enh.idot.exa}. We then obtain an enhanced functor
\begin{equation}\label{enh.groth.eq}
\Gr:\E^\iota\cCat^h \to \cCat^h \mmh \E, \quad \langle I,\gamma
\rangle \mapsto \langle I,\E'(\iota \circ \gamma)^\iota,\beta(\iota
\circ \gamma)^\iota \rangle
\end{equation}
by composing \eqref{enh.pb.fun.eq} with the involution
$\iota:\E^\iota\cCat^h \to \E^\iota\cCat^h$, $\gamma \mapsto \iota
\circ \gamma$ induced by \eqref{iota.enh.eq}.

\begin{prop}\label{groth.prop}
The enhanced functor \eqref{enh.groth.eq} is an equivalence for any
small enhanced category $\E$.
\end{prop}

\proof{} To construct an enhanced functor inverse to
\eqref{enh.groth.eq}, consider the product $\E^\iota \times^h
(\cCat^h \mmh \E)$. By Lemma~\ref{al.be.le}, its objects are
quadruples $\langle I,e,\C,\beta\rangle$, $I \in \Posf^+$,
$e:\Unf(I) \to \E$ an enhanced functor, $\C$ a small enhanced
category, $\beta:\C \to \Unf(I) \times^h \E$ an enhanced
fibration. For any such quadruple, we can form a semicartesian
square
\begin{equation}\label{groth.1.sq}
\begin{CD}
\C(I,e,\C,\beta) @>>> \C\\
@V{\pi(I,e,\C,\beta)}VV @VV{\beta}V\\
\Unf(I) @>{\id \times e}>> \Unf(I) \times^h \E,
\end{CD}
\end{equation}
and by Lemma~\ref{enh.fib.le} together with
Lemma~\ref{noenh.fib.le}, $\pi(I,e,\C,\beta)$ in \eqref{groth.1.sq}
is an $I$-aug\-ment\-ation. Again, by Corollary~\ref{sq.h.corr} and
Lemma~\ref{enh.fib.le}, the square \eqref{groth.1.sq} is functorial
with respect to morphisms cartesian over $\Posf^+$, so by
Lemma~\ref{enh.exp.le}, sending $\langle I,e,\C,\beta \rangle$ to
the corresponding pair $\langle \C(I,e,\C,\beta),\pi(I,e,\C,\beta)
\rangle$ provides an enhanced functor
\begin{equation}\label{groth.1.eq}
\E^\iota \times^h (\cCat^h \mmh \E) \to \cCat^h.
\end{equation}
By the universal property of the functor categories,
\eqref{groth.1.eq} then induces an enhanced functor $\Gr':\cCat^h
\mmh \E \to \E^\iota\cCat$, and the semicartesian squares
\eqref{ccat.E.sq} provide an isomorphism $\Gr' \circ \Gr \cong
\id$. Define enhanced categories $\cCat^h_\idot \mmh_\iota \E$ and
$\cCat^h_\idot \bbi^h_\iota \E \cong \cCat^h_\idot \bbi^h \E^\iota$ by
semicartesian squares
\begin{equation}\label{ccat.I.E.sq}
\begin{CD}
\cCat^h_\idot \mmh_\iota \E @>{a}>> \cCat^h_\idot \bbi^h_\iota \E
@>{\sim}>> \cCat^h_\idot \bbi^h \E^\iota @>>> \cCat^h_\idot\\
@V{\pi}VV @VVV @VVV @VVV\\
\cCat^h \mmh_\iota \E @>{a}>> \cCat^h \bbi^h_\iota \E @>{\sim}>> \cCat^h \bbi^h
\E^\iota @>>> \cCat^h,
\end{CD}
\end{equation}
and note that \eqref{ccdot.pair} then induces an enhanced functor
\begin{equation}\label{iota.cc.pair}
\begin{CD}
  \cCat^h_\idot \mmh_\iota \E @>{a}>> \cCat^h_\idot \bbi^h \E^\iota
  @>{b}>> \cCat^h_\idot \bb^h \E^\iota @>>> \E^\iota,
\end{CD}
\end{equation}
and its product with $\pi$ of \eqref{ccat.I.E.sq} fits into a
commutative square
\begin{equation}\label{ccat.I.E.d.sq}
\begin{CD}
\cCat^h_\idot \mmh_\iota \E @>>> \cCat^h_\idot\\
@V{\pi'' \times^h \pi}VV @VVV\\
\E^\iota \times^h (\cCat^h \mmh \E) @>>> \cCat^h,
\end{CD}
\end{equation}
where the bottom arrow is the pairing \eqref{groth.1.eq}. Spelling
out the definitions, we see that the square \eqref{ccat.I.E.d.sq} is
semicartesian as well.

Now to finish the proof, take some $\langle I,\C,\beta \rangle \in
\cCat^h \mmh \E$, with the corresponding enhanced functor
$\gamma:\Unf(I^o) \to \cCat^h \mmh \E$, and note that since
\eqref{ccat.I.E.sq} and \eqref{ccat.E.sq} are semicartesian, we have
a semicartesian square
\begin{equation}\label{ccat.CE.sq}
\begin{CD}
\C^\iota @>>> \cCat^h_\idot \mmh \E\\
@V{\beta^\iota}VV @VV{\pi'' \times^h \pi}V\\
\E^\iota \times \Unf(I^o) @>{\id \times \gamma}>> \E^\iota \times^h
(\cCat^h \mmh \E).
\end{CD}
\end{equation}
The square \eqref{ccat.CE.sq} is then functorial with respect to
$\C$, and combining it with \eqref{ccat.I.E.d.sq} defines an isomorphism
$\Gr \circ \Gr' \cong \id$.
\endproof

Proposition~\ref{groth.prop} is an enhanced version of the
Grothendieck construction, and we see that it actually works better
than the unenhanced version. Namely, while by Example~\ref{IX.exa},
a discrete small fibration over some $I$ is the fibration $IX \to I$
for a unique functor $X:I^o \to \Sets$, one has to use
pseudofunctors for fibrations that are not discrete. In the enhanced
setting, the ``pseudo'' part is already incorporated into
enhancement for functors, so enhanced fibrations $\E' \to \E$ over a
small enhanced base $\E$ simply correspond to enhanced functors
$\E^\iota \to \cCat^h$. It is convenient to keep the notation of
Example~\ref{IX.exa}, and denote the enhanced fibration
corresponding to some $X:\E^\iota \to \cCat^h$ by $\E X \to \E$. The
role of discrete fibrations is played by enhanced families of
groupoids, and by Lemma~\ref{enh.fib.le}, $\E X \to \E$ is an
enhanced family of groupoids iff $X$ factors through $\sSets^h
\subset \cCat^h$.

As an immediate application, Proposition~\ref{groth.prop} allows one
to define transpose enhanced cofibrations and transpose-opposite
enhanced fibrations. Namely, for any enhanced fibration $\C \to \E$
between small enhanced categories, we can compose the corresponding
enhanced functor $\E^\iota \to \cCat^h$ with the involution $\iota$,
and obtain an enhanced fibration $\C^\iota_{h\perp} \to \E$ whose
enhanced fibers are enhanced-opposite to those of $\E$. The
enhanced-opposite category $\C_{\perp} = (\C^\iota_{h\perp})^\iota$
then comes equipped with an enhanced cofibration $\C_{h\perp} \to
\E^\iota$. If we have another enhanced fibration $\C' \to \E$ with
small $\C'$, then an enhanced functor $\gamma:\C \to \C'$ cartesian
over $\E$ induces a transpose functor $\gamma_{h\perp}:\C_{h\perp}
\to \C_{h\perp}'$ cocartesian over $\E^\iota$ and the
transpose-opposite functor $\gamma^\iota_{h\perp}:\C^\iota_{h\perp}
\to {\C'}^\iota_{h\perp}$ cartesian over $\E$. This provides an
enhanced equivalence
\begin{equation}\label{perp.eq}
\fFun^{\Toth}_{\E}(\C,\C') \cong
\fFun^h_{\Tot\E}(\C_{h\perp},\C'_{h\perp}),
\end{equation}
where $\fFun^h_{\E}(-,-)$ stands for the enhanced functor category
over $\E$, as in \eqref{fun.C.enh.sq}, and $\fFun^{\Toth}_{\E}$
resp.\ $\fFun^h_{\Tot\E}$ are the full enhanced subcategories
spanned by enhanced functors cartesian resp.\ cocartesian over
$\E$. If we have $\E = \Unf(I)$ for some $I \in \Posf^+$, this
recovers our earlier construction for $I$-aug\-mented enhanced
categories. In general, we have a covariant version of the
Grothendieck construction: any enhanced cofibration $\E' \to \E$
with a small base is of the form $\E' \cong X^*\cCat^h_\idot$ for a
unique enhanced functor $X:\E \to \cCat^h$. Extending
Example~\ref{IX.exa}, we denote $\E^\perp X = X^*\cCat^h_\idot$.

\begin{exa}\label{tw.h.exa}
For any small enhanced category $\E$, with the enhanced arrow
category $\Ar^h(\E)$, one can define the {\em enhanced twisted arrow
  category} $\Tw^h(\E)$ by $\Tw^h(\E)^\iota \cong
\Ar^h(\E)_{h\perp}$, where $\sigma:\Ar^h(\E)_{h\perp} \to \E^\iota$
is transpose to the enhanced fibration $\sigma:\Ar^h(\E) \to
\E$. The enhanced cofibration $\tau:\Ar^h(\E) \to \E$ then induces
an enhanced cofibration $\tau:\Tw^h(\E)^\iota \to \E$ by
\eqref{perp.eq}, and the resulting enhanced functor $\Tw^h(\E) \to
\E \times^h \E^\iota$ is an enhanced family of groupoids.
\end{exa}

\subsection{Enhanced relative functor categories.}\label{enh.rel.subs}

As another direct application of Proposition~\ref{groth.prop}, for
any enhanced fibration $\C \to \E$ with small $\C$, we obtain a
universal characterization of the relative functor category
$\fFun^h(\C',\C|\E)$ of \eqref{fun.enh.rel.sq}. Indeed, in terms of
the equivalence \eqref{enh.groth.eq}, passing from $\C$ to
$\Fun(\C',\C|\E)$ corresponds to composing an enhanced functor
$\E^\iota \to \cCat^h$ with the functor $\fFun^h(\C',-):\cCat^h \to
\cCat^h$, and then by Corollary~\ref{enh.adj.fun.corr}, this is
right-adjoint to the functor $\C' \times^h -$, so that we have a
natural identification
\begin{equation}\label{fun.fun.1}
(\Cat^h \mmh \E)(\E' \times^h \C',\C) \cong (\Cat^h \mmh
  \E)(\E',\fFun^h(\C',\C|\E))
\end{equation}
for any enhanced fibration $\E' \to \E$ with small $\E'$. We note
that the left-hand side in \eqref{fun.fun.1} can be also canonically
expressed via \eqref{ffun.C.univ}: we have
\begin{equation}\label{fun.fun.2}
\Cat^h(\C',\fFun^{\Toth}_{\C}(\E',\C)) \cong (\Cat^h \mmh \E)(\C'
\times^h \E',\C),
\end{equation}
where $\fFun^{\Toth}_{\C}(\E',\C) \subset \fFun^h_{\C}(\E',\E)$ is
the full enhanced subcategory spanned by enhanced functors cartesian
over $\E$.

Now assume given an enhanced functor $\gamma:\E' \to \E$ between
small enhanced categories, and recall that by \eqref{cat.phi},
postcomposition with $\gamma$ provides a functor $\gamma_\trr:\Cat^h
\bb^h \E' \to \Cat^h \bb^h \E$. The functor $\gamma_\trr$ has an
obvious enhancement $\langle I,\C,\alpha \rangle \mapsto \langle
I,\C,\gamma \circ \alpha \rangle$, so $\cCat^h \bb^h -$ can be
thought of as a $2$-functor sending small enhanced categories $\E$
to enhanced categories $\cCat^h \bb^h \E$, and enhanced functors
$\gamma:\E' \to \E$ to enhanced functors $\gamma_\trr:\cCat^h \bb^h
\E' \to \cCat^h \bb^h \E$. By virtue of Proposition~\ref{no.E.prop},
$\cCat^h \bb^h -$ inherits the nice properties of the enhanced
functor categories $\fFun^h(\C,-)$; in particular, it sends
semicartesian squares to semicartesian squares, and enhanced
fibrations to enhanced fibrations (combine
Corollary~\ref{enh.adj.fun.corr} and Lemma~\ref{enh.cart.fib.le}).

If we restrict our attention to the subcategory $\Cat^h \bbi^h \E'
\subset \Cat^h \bb^h \E'$, then $\gamma_\trr$ sends it into $\Cat^h
\bbi^h \E \subset \Cat^h \bb^h \E$, and by Corollary~\ref{sq.h.corr}
and Lemma~\ref{sq.h.le}, the resulting functor $\gamma_\trr:\Cat^h
\bbi^h \E' \to \Cat^h \bbi^h \E$ admits a right-adjoint $\gamma^*$ sending
$\C$ to $\E' \times^h_{\E} \C$. As it happens, this adjoint also has
a natural enhancement. Namely, $\langle I,\C,\alpha \rangle \in
\cCat^h \bb^h \E$ represents an object in $\cCat^h \bbi^h \E$ iff the
product $\alpha \times \pi:\C_{h\perp} \to \Unf(I^o) \times^h \E$ of
the enhanced functor $\alpha$ and the coaugmentation
$\pi:\C_{h\perp} \to \Unf(I^o)$ is $I^o$-coaugmented, and in this
case, $\gamma^*\C_{h\perp} \to \Unf(I^o)$ is an $I^o$-coaugmentation
by Lemma~\ref{enh.fib.le}. Then $\gamma^*\C_{h\perp} \cong
\C'_{h\perp}$ for a unique $I$-augmented small enhanced category
$\C' = (\gamma^*\C_{h\perp})^{h\perp}$, and we can let
$\gamma^*(\langle I,\C,\alpha \rangle) = \langle
I,\C',\gamma^*(\alpha)\rangle$. Moreover, if $\gamma$ is an enhanced
cofibration, then $\gamma^*\C_{h\perp} \to \C_{h\perp}$ is an
enhanced cofibration as well, and then $\gamma^*\C_{h\perp} \to
\Unf(I^o)$ is an $I^o$-coaugmentation for any $\langle I,\C,\alpha
\rangle \in \cCat^h \bb^h \E$. Thus for an enhancd cofibration
$\gamma$, we have a well-defined enhanced functor
\begin{equation}\label{ga.bb.eq}
\gamma^*:\cCat^h \bb^h \E \to \cCat^h \bb^h \E', \quad \langle
I,\C,\alpha \rangle \mapsto \langle
I,(\gamma^*\C_{h\perp})^{h\perp},\gamma^*(\alpha) \rangle,
\end{equation}
enhanced right-adjoint to $\gamma_\trr$. Let us now show that we
also have an enhanced version of Lemma~\ref{fun.le} and the relative
functor categories \eqref{rel.fun.eq}.

\begin{lemma}\label{enh.fun.le}
For any enhanced cofibration $\gamma:\E' \to \E$ between small
enhanced categories, the enhanced functor $\gamma^*:\cCat^h \bbi^h \E
\to \cCat^h \bbi^h \E'$ admits an enhanced right-adjoint
$\gamma_\trl:\cCat^h \bbi^h \E' \to \cCat^h \bbi^h \E$. Moreover, let
$\E'_{h\sharp}$ be as in \eqref{h.shrp.eq}, with the cocartesian
functor $\nu:\E'_{h\sharp} \to \E$. Then for any small $\C \to \E'$ such
that $\nu^*\C \to \E'_{h\sharp}$ is an enhanced fibration,
$\gamma_\trl\C \to \E$ is an enhanced fibration, and for any
small enhanced fibration $\C' \to \E$ and a functor
$\phi:\gamma^*\C' \to \C$, the adjoint functor $\phi':\C' \to
\gamma_\trl\C$ is cartesian over $\E$ iff $\nu^*(\phi)$ is cartesian
voer $\E'_{h\sharp}$.
\end{lemma}

\proof{} Let $E:\E \to \cCat^h$ be the enhanced functor
corresponding to $\E' \to \E$, so that $\E' \cong E^*\cCat^h_\idot$,
and for any small enhanced category $\C$, define an enhanced
fibration $\fFun^h(\E'|\E,\C) \to \E$ by the semicartesian square
\begin{equation}\label{rel.ffun.sq}
\begin{CD}
\fFun^h(\E'|\E,\C) @>>> \cCat^h \bb^h \C\\
@VVV @VV{\pi}V\\
\E @>{E}>> \cCat^h,
\end{CD}
\end{equation}
where $\pi$ is the enhanced fibration \eqref{no.E} of
Proposition~\ref{no.E.prop}. Then $\fFun^h(\E'|\E,\C)$ is functorial
with respect to $\C$ via the postcomposition functors
\eqref{cat.phi}, and $\fFun^h(\E'|\E,-)$ sends enhanced fibrations
to enhanced fibrations and semicartesian squares to semicartesian
squares. Moreover, by \eqref{ccat.dot.sq}, we have
$\gamma^*\fFun^h(\E'|\E,\C) \cong E^*(\cCat^h_\idot \bb^h \C)$, and
\eqref{ccdot.pair} induces an enhanced functor
\begin{equation}\label{ev.E.C.eq}
\ev:\gamma^*\fFun^h(\E'|\E,\C) \to \C.
\end{equation}
On the other hand, composing the tautological embedding \eqref{y.eq}
with \eqref{ga.bb.eq} provides an enhanced functor $\gamma^* \circ
\y:\E \to \cCat^h \bb^h \E'$, and we tautologicaly have $\pi \circ
\gamma^* \circ \y \cong E$. Altogether, we obtain an enhanced
section
\begin{equation}\label{E.s.eq}
s:\E \to \fFun(\E'|\E,\E')
\end{equation}
of the enhanced fibration $\fFun^h(\E'|\E,\E') \to \E$, and
moreover, this section factors through
$\nu_\trr:\fFun^h(\E'|\E,\E'_{h\sharp}) \to
\fFun^h(\E'|\E,\E')$. More generally, for any small enhanced
category $\C$ and enhanced functor $\alpha:\C \to \E$, the same
construction applies to the enhanced cofibration $\C' = \alpha^*\E'
\cong \gamma^*\C \to \C$, and provides an enhanced functor
\begin{equation}\label{C.E.s.sq}
s:\C \to \fFun^h(\C'|\C,\C') \cong
\alpha^*\fFun^h(\E'|\E,\gamma^*\C)
\to \fFun^h(\E'|\E,\gamma^*\C)
\end{equation}
over $\E$. Now for any small enhanced category $\C$ equipped with an
enhanced functor $\pi:\C \to \E'$, we can define
$\gamma_\trl(\pi):\gamma_\trl\C \to \E$ by the semicartesian square
\begin{equation}\label{enh.trl.sq}
\begin{CD}
\gamma_\trl\C @>>> \fFun^h(\E'|\E,\C)\\
@V{\gamma_\trl(\pi)}VV @VV{\pi_\trr}V\\
\E @>{s}>> \fFun^h(\E'|\E,\E),
\end{CD}
\end{equation}
where $s$ is the section \eqref{E.s.eq}. Then $s$ factors through
$\fFun^h(\E'|\E,\E'_{h\sharp})$, so $\gamma_\trl(\pi)$ is an
enhanced fibration as soon as so is $\nu^*(\pi):\nu^*\C \to
\E'_{h\sharp}$. Moreover, for any small enhanced category $\E''$ and
enhanced functor $\phi:\E'' \to \E$, we have $\phi^*\gamma_\trl\C
\cong \phi^*(\gamma)_\trl\phi^*\C$.  In particular, by
Lemma~\ref{al.be.le}, objects in $\cCat^h \bbi^h \E'$ correspond to
triples $\langle I,\C,\beta \rangle$, where $\beta:\C \to \Unf(I)
\times^h \E'$ becomes a fibration after restricting to $\Unf(I)
\times^h \E_{h\Iso}$, and then for any such triple, $\langle I,(\id
\times \gamma)_\trl\C,(\id \times \gamma)_\trl(\beta)\rangle$ is a
well-defined object in $\cCat^h \bbi^h \E$. This is functorial with
respect to maps cartesian over $\Posf^+$, thus defines an enhanced
functor
\begin{equation}\label{enh.trl.eq}
\gamma_\trl:\cCat^h \bbi^h \E' \to \cCat^h \bbi^h \E
\end{equation}
by Lemma~\ref{enh.exp.le}, while \eqref{ev.E.C.eq} and
\eqref{C.E.s.sq} extend to functorial maps
$$
\gamma^* \circ \gamma_\trl \to \id, \qquad \id \to \gamma_\trl \circ
\gamma^*
$$
that define adjunction between $\gamma^*$ and $\gamma_\trl$.
\endproof

Just as in the unenhanced setting of Subsection~\ref{rel.fun.subs},
Lemma~\ref{enh.fun.le} implies a universal property for the relative
enhanced functor category \eqref{rel.ffun.sq} and the evaluation
functor \eqref{ev.E.C.eq}: for any small enhanced category $\C'$
over $\E$, an enhanced functor $\phi:\gamma^*\C' \to \C$ factors as
\begin{equation}\label{enh.ev.facto}
\begin{CD}
\gamma^*\C' @>{\gamma^*\wt{\phi}}>> \gamma^*\fFun^h(\E'|\E,\C)
@>{\ev}>> \C
\end{CD}
\end{equation}
for an enhanced functor $\wt{\phi}:\C' \to \fFun^h(\E'|\E,\C)$ over
$\E$, uniquely up to an isomorphism. As in
Subsection~\ref{rel.fun.subs}, this property characterizes
$\fFun^h(\E'|\E,\C)$ uniquely up to a unique equivalence,
Lemma~\ref{enh.fun.le} proves existence when $\gamma:\E' \to \E$ is
an enhanced cofibration, but $\fFun^h(\E'|\E,-)$ actually exists for
more general enhanced functors $\gamma$. In particular, if
$\gamma:\E' \to \E$ is an enhanced fibration, we can define
\begin{equation}\label{rel.fun.perp}
\fFun^h(\E'|\E,\C) = \fFun^h({\E'}^\iota|\E^\iota,\C^\iota)^\iota
\cong \fFun^h(\E'_{h\perp}|\E^\iota,\C)_{h\perp},
\end{equation}
where $\gamma^\iota:{\E'}^\iota \to \E^\iota$ is the
enhanced-opposite cofibration, and then \eqref{ev.E.C.eq} provides
an enhanced functor
\begin{equation}\label{ev.E.C.perp.eq}
\ev^\iota:\gamma^{\iota *}\fFun^h({\E'}^\iota|\E^\iota,\C^\iota) \to \C^\iota
\end{equation}
that gives \eqref{ev.E.C.eq} after again passing to the
enhanced-opposite categories. The resulting pairing has exactly the
same universal property as in the cofibration case. If
$\gamma=\Unf(g)$ for some functor $g:I' \to I$ between essentially
small categories that is either a fibration or a cofibration, we
will simplify notation and write $\fFun^h(I'|I,-) =
\fFun^h(\Unf(I')|\Unf(I),-)$.

\subsection{Enhanced Yoneda embedding.}\label{enh.yo.subs}

Let us now give our main application of
Proposition~\ref{groth.prop}, namely, an enhanced version of the
extended Yoneda embedding \eqref{yo.cat.eq}. To construct it,
observe that for any small enhanced functor $\alpha:\C \to \E$, the
enhanced comma-category $\E \setminus^h_\alpha \C$ comes equipped
with the fibration $\sigma:\E \setminus^h_\alpha \C \to \E$ of
Corollary~\ref{enh.si.ta.corr}. More generally, for any object
$\langle I,\C,\alpha \rangle \in \cCat^h \bb^h \E$, the enhanced
comma-category $\E \setminus^h_\alpha \C_{h\perp}$ is
$I$-coaugmented via $\tau:\E \setminus^h \C_{h\perp} \to
\C_{h\perp}$, and the projection $\sigma:\E \setminus^h_\alpha
\C_{h\perp} \to \E$ is an enhanced fibration. Moreover, for any $i
\in I$, the induced projection $(\E \setminus^h_\alpha \C)_i \cong
\E \setminus^h \C_i \to \E$ is also an enhanced fibration, and for
any $i \leq i'$, the transition functor $\E \setminus^h \E_i \to \E
\setminus^h \C_{i'}$ is cartesian over $\E$. Since the
correspondence $\C \mapsto (\E \setminus^h \C_{h\perp})^{h\perp}$ is
functorial with respect to $I$ by Corollary~\ref{sq.h.corr},
Lemma~\ref{enh.exp.le} provides an enhanced functor
\begin{equation}\label{yo.enh.eq}
\Y:\cCat^h \bb^h \E \to \cCat^h \mmh \E \cong \E^\iota\cCat^h
\end{equation}
sending $\langle I,\C,\alpha \rangle$ to $\langle I,(\E
\setminus^h_\alpha \C_{h\perp})^{h\perp},\sigma \rangle$. Moreover,
for any $\alpha:\C \to \E$, $\tau:\E \setminus^h_\alpha \C \to \C$
is a lax functor over $\E$, its right-adjoint $\eta:\C \to \E
\setminus^h_\alpha \C$ is a functor over $\E$, and if $\alpha$ is an
enhanced fibration, then the second right-adjoint $\eta_\dg:\E
\setminus^h_\alpha \C \to \C$ to $\eta$ is cartesian over $\E$. More
generally, for any object $\langle I,\C,\alpha \rangle$ in $\cCat^h
\bb^h \E$, the lax functor $\tau:\E \setminus^h_\alpha \C_{h\perp}
\to \C_{h\perp}$ over $\E$ is $I$-coaugmented by definition, thus
defines a morphism in $\cCat^h \bb^h \E$, and its right-adjoint
$\eta:\C_{\perp} \to \E \setminus^h_\alpha \C_{h\perp}$, a functor
over $\E$, is also $I$-coaugmented. Moreover, if the triple $\langle
I,\C,\alpha \rangle$ is in the subcategory $\cCat^h \mmh \E \subset
\cCat^h \bb^h \E$, so that $\pi \times \alpha = \beta^{h\perp}$ for
an enhanced fibration $\beta:\C \to \Unf(I) \times^h \E$, as in
Lemma~\ref{al.be.le}, then $\E \setminus^h_\alpha \C_{h\perp} \cong
(\E \setminus^h_{\pi' \circ \beta} \C)_{h\perp}$, where
$\pi':\Unf(I) \times^h \E \to \E$ is the projection, and the
right-adjoint $\eta_\dg$ defines a morphism in $\cCat^h \mmh \E$. By
Lemma~\ref{enh.exp.le}, we then have enhanced morphisms
\begin{equation}\label{a.b.Y.eq}
\tau:b \circ a \circ \Y \to \id, \qquad \eta:\id \to a \circ \Y
\circ b, \qquad \eta_\dg:\Y \circ b \circ a \to \id,
\end{equation}
where $\Y$ is \eqref{yo.enh.eq}, and $a$, $b$ are the enhanced
functors of \eqref{3.h.enh.dia}. We also have adjunction
isomorphisms $\eta_\dg \circ \eta \cong \id$, $\tau \circ \eta \cong
\id$.

\begin{prop}\label{yo.prop}
For any small enhanced category $\E$, the morphisms \eqref{a.b.Y.eq}
define an adjunction between enhanced functors $a$ and $\Y \circ b$,
and between $a \circ \Y$ and $b$, and the enhanced functor
\eqref{yo.enh.eq} is fully faithful.
\end{prop}

\proof{} For the first claim, it suffices to check the identities of
Lemma~\ref{3.adj.le}, and all four immediately follow from the
existence of the adjunction isomorphisms $\eta_\dg \circ \eta \cong
\id$, $\tau \circ \eta \cong \id$. By Lemma~\ref{3.adj.le}, this
also proves that $\Y$ is fully fiathful over $\ppt \in
\Posf^+$. Over a more general $I \in \Posf^+$, this is not
automatic, since $b$ is no longer essentially surjective. However,
let $A$ be the class of morphisms in $\cCat^h \bb^h \E$ represented
by triples $\langle f,\phi,a \rangle$ such that $a$ is an
isomorphism. Then the dense subcategory $(\cCat^h \bb^h \E)_A
\subset \cCat^h \bb^h \E$ is no longer an enhanced category, but the
projection $(\cCat^h \bb^h \E)_A \to \Posf^+$ is still a fibration,
the embedding functor $b':(\cCat^h \bb^h \E)_A \to \cCat^h \bb^h \E$
is cartesian over $\Posf^+$, and we have a diagram
$$
\begin{CD}
\cCat^h \mmh \E @>{a'}>> (\cCat^h \bb^h \E)_A @>{b'}>> \cCat^h \bb^h
\E @>{\Y}>> \cCat^h \mmh \E,
\end{CD}
$$
where $a'$ is cartesian over $\Posf^+$, and $\Y$ is the same
enhanced functor \eqref{yo.enh.eq}. Moreover, $\eta$ in
\eqref{a.b.Y.eq} is a well-defined functor on the whole category
$(\cCat^h \bb^h \E)_A$, and Lemma~\ref{3.adj.le} still
applies. Since $b'$ is essentially surjective by definition, this
finishes the proof.
\endproof

The fully faithful embedding \eqref{yo.enh.eq} of
Proposition~\ref{yo.prop} is the enhanced version of the extended
Yoneda embedding \eqref{yo.cat.eq}. To deduce the usual version, one
can first restrict the enhanced fibration \eqref{no.E} to the full
enhanced subcategory $\sSets^h \subset \cCat^h$ of enhanced
groupoids. The resulting category $\sSets^h \bb^h \E$ is then the
enhanced category of enhanced groupoids $\C$ equipped with an
enhanced functor $\alpha:\C \to \E$. For any such $\C$, $\Y(\C) \to
\E$ is an enhanced family of groupoids by
Example~\ref{enh.si.grp.exa}, and then by Lemma~\ref{enh.fib.le},
\eqref{enh.groth.eq} sends $\Y(\C)$ to an enhanced functor that
factors through $\sSets^h \subset \cCat^h$, so that
\eqref{yo.enh.eq} induces a fully faithful enhanced functor
$\sSets^h \bb^h \E \to \E^\iota\sSets^h$. Restricting even further,
we can compose $\Y$ with the full embedding \eqref{y.eq}, or in
other words, restrict our attention to the full enhanced subcategory
in $\sSets^h \bb^h \E$ spanned by $\C = \ppt$. This enhanced
subcategory is tautologically identified with $\E$ itself, so that
\eqref{yo.enh.eq} provides a fully faithful embedding
\begin{equation}\label{yo.yo.enh.eq}
  \Y(\E):\E \to \E^\iota\sSets^h.
\end{equation}
This is the enhanced version of \eqref{yo.yo.eq}, and here is its
essential image.

\begin{lemma}\label{yo.im.le}
An enhanced family of small groupoids $\E' \to \E$ corresponding to
an enhanced functor $E:\E^\iota \to \sSets^h$ lies in the essential
image of the Yoneda embedding \eqref{yo.yo.enh.eq} iff $\E'$ admits
a terminal enhanced object.
\end{lemma}

\proof{} By definition, for any enhanced object $e \in \E$, $\Y(\E)(e)$
corresponds to the enhanced family of groupoids $\E /^h e \to \E$,
and this has the enhanced-terminal object given by the embedding
$\eta:\ppt^h \to \E /^h e \cong \E \setminus^h
\ppt^h$. Conversely, for any enhanced family of groupoids
$\pi:\E' \to \E$ and enhanced object $e \in \E'_\ppt$, the
equivalence $\pi \setminus^h \id:\Ar^h(\E') \to \E \setminus^h_\pi
\E'$ induces an equivalence $\E' /^h e \to \E /^h \pi(e)$, and if
$e$ is terminal, we have $\E' \cong \E' /^h e$.
\endproof

\begin{remark}
By the universal property of the enhanced functor category,
\eqref{yo.yo.enh.eq} defines and is defined by the enhanced {\em
  Yoneda pairing}
\begin{equation}\label{yo.pair.eq}
  \E^\iota \times^h \E \to \sSets^h.
\end{equation}
In words, \eqref{yo.pair.eq} shows that for any enhanced objects
$e$, $e'$ in an enhanced category $\E$, the groupoids of enhanced
morphisms from $e$ to $e'$ comes equipped with a canonical
enhancement. The enhanced family of groupoids corresponding to
\eqref{yo.pair.eq} by the Grothendieck construction is $\Tw^h(\E)
\to \E \times^h\E^\iota$, where $\Tw^h(\E)$ is the enhanced twisted
arrow category of Example~\ref{tw.h.exa}.
\end{remark}

Going back to the extended Yoneda embedding of
Proposition~\ref{yo.prop}, we note that as in Remark~\ref{ind.rem},
we can express the enhanced left-adjoint $\Y \circ b$ to $a:\cCat^h
\mmh \E \to \cCat^h \bbi^h \E$ via the enhanced arrow category
$\Ar^h(\E)$ and enhanced functors $\sigma,\tau:\Ar^h(\E)$ ---
namely, we have $\Y \circ b \cong \sigma_\trr \circ
\tau^*$. Explictly, for any enhanced fibrations $\C,\C' \to \E$, we
have an equivalence
\begin{equation}\label{Y.co.eq}
  \fFun^{\Toth}_\E(\E \setminus^h \C,\C') \cong \fFun^h_\E(\C,\C'),
\end{equation}
where $\fFun^{\Toth}_\E$ has the same meaning as in \eqref{perp.eq}.
Let us show that as in Corollary~\ref{ind.corr}, a dual procedure
provides a right-adjoint.

\begin{corr}\label{enh.ind.corr}
For any small enhanced category $\E$, the enhanced functor
$a:\cCat^h \mmh \E \to \cCat^h \bbi^h \E$ is right-reflexive.
\end{corr}

\proof{} As in Corollary~\ref{ind.corr}, Lemma~\ref{enh.fun.le}
shows that for any small enhanced category $\C$ equipped with an
enhanced functor $\C \to \E$, the enhanced category
$\tau_\trl\sigma^*\C$ exists since $\tau$ is an enhanced
cofibration, and moreover, the functor $\tau_\trl\sigma^*\C \to \E$
is an enhanced fibration. Moreover, for any other enhanced fibration
$\C' \to \E$, a functor $\phi:\C' \to \tau_\trl\sigma^*\C$ over $\E$
is cartesian over $\E$ iff the adjoint functor $\phi':\tau^*\C' \to
\sigma^*\C$ becomes cartesian after restriction via $\nu:\E_{h\Iso}
/^h \E \to \Ar^h(\E)$. By Corollary~\ref{sq.h.corr}, isomorphism
classes of functors $\phi':\tau^*\C' \to \sigma^*\C$ correspond
bijectively to isomorphism classes of functors $\phi'':\tau^*\C'
\cong \E \setminus^h \C' \to \C$ over $\E$, and since $\C' \to \E$
is an enhanced fibration, we have an adjoint pair of enhanced
functors $\eta:\C' \to \E \setminus^h \C'$, $\eta_\dg:\E \setminus^h
\C' \to \C'$. Then $\phi'$ becomes cartesian over $\E_{h\Iso} /^h
\E$ iff and only if the adjunction map $\phi'' \to \phi'' \circ \eta
\circ \eta_\dg$ is an isomorphism, so at the end of the day, sending
$\phi:\C' \to \tau_\trl\sigma^*\C$ to the corresponding enhanced
functor $\eta \circ \phi'':\C' \to \C$ establises a bijection
between the isomorphism classes of enhanced functors $\C' \to \C$
over $\E$ and isomorphism classes of enhanced functors $\C' \to
\tau_\trl\sigma^*\C'$ cartesian over $\E$.
\endproof

\section{Limits and colimits.}\label{enh.lim.sec}

\subsection{Complete and cocomplete enhanced categories.}

Unfortuna\-tely, representability theorems of
Subsection~\ref{enh.repr.subs} and all their corollaries discussed
in Section~\ref{enh.repr.sec} only work for small enhanced
categories. The theory of large enhanced categories is rather
deficient, even more so than in the unenhanced setting, and it seems
that to get a reasonable theory, one has to require the existence of
some sort of limits or colimits. As a starting point, let us
require them all.

\begin{defn}\label{enh.compl.def}
An enhanced category $\E$ is {\em complete} if it is a complete
family over $\Posf^+$ in the sense of
Definition~\ref{compl.fam.def}, and {\em cocomplete} if the
enhanced-opposite category $\E^\iota$ is complete.
\end{defn}

We recall that by definition, $\E$ is complete iff for any map $f:I
\to I'$ in $\Posf^+$, the pullback functor $f^*:\E_{I'} \to \E_I$
has a right-adjoint functor $f_*$, and for any cartesian square
\eqref{posf.sq} in $\Posf^+$ such that $f$ is a cofibration, we have
the base change isomorphism \eqref{refl.bc}.

\begin{lemma}\label{compl.bar.le}
Assume given an enhanced category $\E$ with the canonical
bar-invariant extension $\E^\hash \to \Pos$ of
Proposition~\ref{lf.prop}. Then $\E$ is complete resp.\ cocomplete
iff for any map $f:J' \to J$ in $\Pos$, the transition functor
$f^*:\E^\Tot_{J} \to \E^\Tot_{J'}$ has a right-adjoint $f_*$
resp.\ a left-adjoint $f_!$, and for any $j \in J$, with the induced
maps $f/j:J'/j \to J/j$, $j \setminus f:j \setminus J' \to j
\setminus J'$, the base change map $\sigma(j)^* \circ f_* \to
(f/j)_* \circ \sigma(j)^*$ resp.\ $(j \setminus f)_! \circ \tau(j)^*
\to \tau(j)^* \circ f_!$ is an isomorphism.
\end{lemma}

\proof{} In the cocomplete case, recall that by definition,
$\E^\Tot$ is the transpose-opposite fibration to
$\iota^*\E^{\iota\Tot}$, so it suffices to prove the claim for the
complete enhanced category $\E^\iota$. In the complete case, for the
``if'' part, note that since $\E \to \Posf^+$ is non-degenerate, it
suffices to check that \eqref{refl.bc} is an isomorphism when
$J_0'=\ppt$. In this case, it reduces to observing that for any
cofibration $f:J' \to J$ in $\Posf^+$, and any element $j \in J$,
the embedding $\eps:J'_j \to J' / j$ is left-admissible, with some
left-adjoint $\eps^\dg:J' / j \to J'_j$, and then since $\E$ is
reflexive, $\eps^\dg_* \cong \eps^*$. For the ``only if'' part,
assume given a map $f:J' \to J$ in $\Pos$, and consider the induced
map $f \circ \xi:B(J') \to J$. Then since $\xi^*:\E^\Tot_{J'} \to
\E^\Tot_{B(J')}$ is fully faithful by bar-invariance, it suffices to
construct the adjoint $(f \circ \xi)_*$, and since $\E^\Tot$ is
reflexive, it further suffices to construct $\tau_*$ for $\tau:B(J')
/_{f \circ \xi} J \to J$. This is a cofibration with fibers $B(J') /
j \cong B(J'/j)$, and if we consider the corresponding cartesian
square
\begin{equation}\label{B.compl.sq}
\begin{CD}
B^\dm(J) \times_{L(J)} L(B(J')/J) @>>> L(B(J')/J)\\
@V{g}VV @VV{\tau}V\\
B^\dm(J) @>{\xi}>> L(J),
\end{CD}
\end{equation}
then $\E^{\Tot\dm}$ is constant along the horizontal arrows by
Lemma~\ref{pm.B.le} and bar-invariance, while $U(g)$ is a
cofibration in $\Posf^\pm$. Then $\tau_*$ is induced by $g_*$, and
the base change isomorphisms are provided by \eqref{refl.bc}.
\endproof

\begin{defn}\label{enh.lim.exa.def}
An enhanced functor $\gamma:\E \to \E'$ between complete enhanced
categories is {\em enhanced-left-exact} if for any map $f:I \to I'$
in $\Posf^+$, the base change map $\gamma_{I'} \circ f_* \to f_*
\circ \gamma_I$ is an isomorphism. An enhanced functor $\gamma$
between cocomplete enhanced categories is {\em enhanced-right-exact}
if $\gamma^\iota$ is enhanced-left-exact.
\end{defn}

\begin{exa}
For any category $\E$, $\Unf(\E)$ is complete resp.\ cocomplete in
the sense of Definition~\ref{enh.compl.def} iff $\E$ is complete
resp.\ cocomplete in the usual sense, and a functor $\gamma:\E
\to \E'$ between complete resp.\ cocomplete categories is left
resp.\ right-exact iff $\Unf(\gamma)$ is enhanced-left-exact
resp.\ enhanced-right-exact.
\end{exa}

\begin{exa}\label{enh.adj.exa.exa}
Any enhanced functor $\gamma:\E \to \E'$ between
complete enhanced categories that is left-reflexive in the enhanced
sense is left-exact (the base change maps of
Definition~\ref{enh.lim.exa.def} are adjoint to the isomorphisms
\eqref{fu.fib} for the left-adjoint enhanced functor $\gamma^\dg:\E'
\to \E$). Dually, a right-reflexive enhanced functor between cocomplete
enhanced categories is right-exact.
\end{exa}

\begin{exa}\label{compl.ccat.exa}
The enhanced category $\cCat^h$ of Proposition~\ref{ccat.prop} is
complete, and so is the full enhanced subcategory $\sSets^h \subset
\cCat^h$. The embedding functor $\sSets^h \to \cCat^h$ is
enhanced-left-exact.
\end{exa}

\begin{exa}\label{enh.adm.compl.exa}
A right-admissible full enhanced subcategory $\E' \subset \E$ in a
complete enhanced category $\E$ is complete, and the adjoint
enhanced functor $\nu^\dg:\E \to \E'$ to the embedding $\nu:\E' \to
\E$ is enhanced-left-exact (for any morphism $f:I \to I'$, the
required right-adjoint is given by $\nu^\dg_{I'} \circ f_*$). In the
situation of Example~\ref{compl.ccat.exa}, $\sSets^h \subset
\cCat^h$ is right-admissible, and the right-adjoint enhanced functor
$\cCat^h \to \sSets^h$ sends a small enhanced category $\C \in
\Cat^h$ to its enhanced isomorphism groupoid $\C_{h\Iso}$ of
\eqref{h.iso.sq}.
\end{exa}

\begin{exa}\label{enh.cmpl.mod.exa}
For any complete model category $\C$, the enhanced localization
$\Hh^W(\C)$ of Lemma~\ref{enh.mod.le} is complete.
\end{exa}

\begin{prop}\label{enh.compl.prop}
For any enhanced functor $\gamma:\C \to \C'$ between small enhanced
categories, and any complete enhanced category $\E$, the pullback
functor $\gamma^*:\fFun^h(\C',\E) \to \fFun^h(\C,\E)$ has a
right-adjoint $\gamma_*$, and for any complete enhanced category
$\E'$, enhanced functor $F:\E \to \E'$ that preseves limits, and
functor $E:\C' \to \E$, the base change map $F \circ f_*(E) \to
f_*(F \circ E)$ is an isomorphism. Moreover, for any semicartesian
square
\begin{equation}\label{enh.bc.sq}
\begin{CD}
\C_1 @>{\nu}>> \C_0\\
@V{\gamma_1}VV @VV{\gamma_0}V\\
\C_1' @>{\nu'}>> \C_0'
\end{CD}
\end{equation}
of small enhanced categories and enhanced functors such that
$\gamma_0$, hence also $\gamma_1$ is an enhanced fibration, the
base change map ${\nu'}^* \circ \gamma_{0*} \to \gamma_{1*} \circ
\nu^*$ is an isomorphism as well.
\end{prop}

\proof{} If $\gamma$ is enhanced right-reflexive, with adjoint
enhanced functor $\gamma_\dg$, then the adjoint $\gamma_* \cong
\gamma_\dg^*$ is provided by Corollary~\ref{enh.adj.fun.corr}, and
all the claims are clear. By \eqref{enh.comma.facto} and
Corollary~\ref{enh.si.ta.corr}, we may then assume that $\gamma$ is
an enhanced fibration. If $\C=\Unf(I^o)$, $\C'=\Unf({I'}^o)$ for
some $I,I' \in \Posf^+$, so that $\gamma=\Unf(f^o)$ for some
cofibration $f:I \to I'$, all claims hold by definition. Next,
assume that $\C'=\Unf(I^o)$ but $\C$ is arbitrary, so that $\gamma$
is an $I^o$-coaugmentation. Then Proposition~\ref{univ.prop} applied
to the augmentation $\C^\iota \to \Unf(I)$ provides a fibration $g:J
\to I$ in $\Pos$, and an $I^o$-coaugmented enhanced functor
$\nu:\Unf(J^o) \to \C$ such that the pullback enhanced functor
$\nu^*:\fFun^h(\C,\E) \to \fFun^h(\Unf(J^o),\E) \cong J^o_h\E$ is
the fully faithful embedding of Example~\ref{enh.univ.exa}. Then
$\gamma_* = \Unf(g^o)_* \circ \nu^*$ is right-adjoint to $\gamma^*$
by Lemma~\ref{adj.ff.le}, and all the claims for $\gamma$ follow
from the corresponding claims for $\Unf(g^o)$.

Finally, assume that $\C'$ is arbitrary, again choose a universal
object for $\C'$ provided by Proposition~\ref{univ.prop}, with the
corresponding enhanced functor $\nu:\Unf(J) \to \C'$ of
\eqref{enh.univ.sq}, and consider the semicartesian square
\begin{equation}\label{pb.bc.sq}
\begin{CD}
\wt{\C} @>{\eta}>> \C\\
@V{\phi}VV @VV{\gamma}V\\
\Unf(J) @>{\nu}>> \C'.
\end{CD}
\end{equation}
Then we already know all the claims for $\phi$, and to prove them
for $\gamma$, it suffices to check that $\phi_* \circ \eta^*$
factors through the enhanced essential image of the full embedding
$\nu^*$. However, since the square \eqref{enh.univ.sq} is
cocartesian, this immediately follows from the base change
isomorphisms for $\phi_*$.
\endproof

\begin{corr}\label{ffun.compl.corr}
For any complete enhanced category $\E$ and small enhanced category
$\C$, the enhanced functor category $\C^\iota\E$ is complete.
\end{corr}

\proof{} Clear. \endproof

\begin{corr}\label{enh.bc.corr}
Assume given a commutative diagram 
\begin{equation}\label{enh.bc.dia}
\begin{CD}
\C_1' @>{\gamma'}>> \C'_0 @>{\pi'}>> \C'\\
@V{\nu_1}VV @VV{\nu_0}V @VV{\nu}V\\
\C_1 @>{\gamma}>> \C_0 @>{\pi}>> \C
\end{CD}
\end{equation}
of small enhanced categories and enhanced functors, such that both
squares are semicartesian, and assume that $\pi$ and $\pi \circ
\gamma$ are enhanced fibrations, and $\gamma$ is cartesian over
$\C$. Then for any complete enhanced category $\E$, the base change
map
\begin{equation}\label{enh.bc.eq}
  \gamma^* \circ \nu_{0*} \to \nu_{1*} \circ {\gamma'}^*
\end{equation}
is an isomorphism.
\end{corr}

\proof{} Use the same argument as in Lemma~\ref{bc.le}, with
Lemma~\ref{cart.fib.le} replaced by Lemma~\ref{enh.cart.fib.le}.
\endproof

Taking enhanced-opposite categories gives a dual version of
Proposition~\ref{enh.compl.prop} for cocomplete $\E$ and
left-adjoint functors; we denote $\gamma_! =
(\gamma^\iota_*)^\iota$. As in the unenhanced setting, if $\gamma:\C
\to \ppt^h$ is the projection to the point, we denote $\gamma_!  =
\colimh_{\C}$ and $\gamma_* = \limh_{\C}$. For any $I \in \Posf^+$,
we will write $\colimh_I = \colimh_{\Unf(I)}$,
$\limh_I=\limh_{\Unf(I)}$. The simplest situation when this is
non-trivial is $\colimh_{\V}$ resp.\ $\limh_{\V^o}$,

\begin{exa}\label{V.holim.exa}
It immediately follows from semiexactness and properness that for
any enhanced category $\E$, the functor $(\V^\iota\E)_\ppt \cong
\E_{\V} \to \V^o\E_\ppt$ induced by \eqref{trunc.eq} is an
epivalence. Therefore for any commutative square $E:[1]^2 \to
\E_\ppt$ in $\E_\ppt$, its restriction $E_\idot:\V^o \to \E_\ppt$ to
$\V^o \subset (\V^>)^o \cong [1]^2$ lifts to an enhanced functor
$E^h_\idot:\Unf(\V^o) \to \E$, uniquely up to an isomorphism, and
then if we denote by $p:\V^o \to \ppt$ the projection, the map
$p^*E(o) \to E_\idot$ lifts to an enhanced map $\Unf(p)^*E^h(o) \to
E_\idot^h$, where $E^h(o):\ppt^h \to \E$ is the enhanced object
corresponding to $E(o) \in \E$. If $\E$ is complete, we obtain an
adjoint map $a:E^h(o) \to \limh_{\V^o}E^h_\idot$. The square $E$ is
{\em homotopy cartesian} if $a$ is an isomorphism. Note that this
property only depends on $E$, not on the liftings, and moreover, if
$\E$ is complete, any functor $\V^o \to \E_\ppt$ extends to a
homotopy cartesian square $[1]^2 \to \E_\ppt$ (lift the functor to
an enhanced functor $E^h_\idot:\Unf(V^o) \to \E$, and take
$\Unf(\eps^o)_*E^h_\idot$, where $\eps:\V \to \V^> \cong [1]^2$ is
the embedding). The extension is unique up to an isomorphism (but
the isomorphism is {\em not} unique). Dually, if $\E$ is cocomplete,
then a square $E$ is {\em homotopy cocartesian} if $E^\iota:[1]^2
\to \E^o_\ppt$ is homotopy cartesian, and again, for a cocomplete
$\E$, any functor $\V \to \E_\ppt$ extends to a homotopy cocartesian
square, uniquely up to a non-unique isomorphism.
\end{exa}

\begin{defn}\label{cc.def}
A complete enhanced category $\E$ is {\em cartesian-closed} if for
any enhanced object $e \in \E_\ppt$, the enhanced functor $e \times
-:\E \to \E$ is enhanced-right-reflexive.
\end{defn}

\begin{exa}
The categories $\sSets^h$ and $\cCat^h$ are cartesian-closed, with
the adjoint to $\C \times^h -$ given by $\fFun^h(\C,-)$.
\end{exa}

\begin{lemma}\label{cc.le}
For any enhanced object $e \in \E_\ppt$ in a complete cocomplete
cartesian-closed enhanced category $\E$, the enhanced functor $e
\times -:\E \to \E$ is enhanced-right-exact.
\end{lemma}

\proof{} Example~\ref{enh.adj.exa.exa}. \endproof

\subsection{Limits and colimits in $\cCat^h$.}

The first example of a complete enhanced category is $\Cat^h$ of
Example~\ref{compl.ccat.exa}, with its right-admissible complete
full enhanced subcategory $\sSets^h \subset \cCat^h$. Both are
actually also cocomplete.

\begin{prop}\label{ccat.coco.prop}
The enhanced category $\cCat^h$ of Proposition~\ref{ccat.prop} is
cocomplete, and the full enhanced subcategory $\sSets^h \subset
\cCat^h$ is enhanced left-admissible and cocomplete.
\end{prop}

\proof{} We have an enhanced full embedding $\cCat^h \subset
\Hh^W(\Delta^o\Delta^o\Sets)$ provided by \eqref{ccat.h.eq}. Since
the model category $\Delta^o\Delta^o\Sets$ is both complete and
cocomplete, so is its enhanced localization, by
Example~\ref{enh.cmpl.mod.exa} applied to $\Delta^o\Delta^o\Sets$
and the opposite model category $\Delta^o\Delta^o\Sets^o$. By
Remark~\ref{I.seg.eq.rem} and Lemma~\ref{bous.le}~\thetag{ii},
$\cCat^h(I) \subset h^W(I^o\Delta^o\Delta^o\Sets)$ is
left-admissible for any $I \in \Posf^+$, so by
Lemma~\ref{fib.adj.le}, $\cCat^h \subset
\Hh^W(\Delta^o\Delta^o\Sets)$ is left-admissible as a category, with
some functor $\lambda:\Hh^W(\Delta^o\Delta^o\Sets) \to \cCat^h$ over
$\Posf^+$ left-adjoint to the embedding. Moreover, by
Lemma~\ref{T.pb.le} and \eqref{refl.bc}, the embedding enhanced
functor $\cCat^h \to \Hh^W(\Delta^o\Delta^o\Sets)$ is
enhanced-left-exact. Then as in Example~\ref{enh.adj.exa.exa},
$\lambda$ is cartesian over $\Posf^+$, so that $\cCat^h \subset
\Hh^W(\Delta^o\Delta^o\Sets)$ is left-admissible in the enhanced
sense, and then it is cocomplete by (the dual version of)
Example~\ref{enh.adm.compl.exa}.

For the second claim, again, $\Delta^o\Sets$ is complete and
cocomplete, so that by \eqref{sset.h.eq}, so is its enhanced
localization $\sSets^h$. Then it suffices to check that $\sSets^h
\cong \Hh^W(\Delta^o\Sets)$ is enhanced left-admissible in
$\Hh^W(\Delta^o\Delta^o\Sets)$. This is again clear, since
$h^W(I^o\Delta^o\Sets) \subset h^W(I^o\Delta^o\Delta^o\Sets)$ is
left-admis\-sible for any $I$, and the embedding functor
$\Hh^W(\Delta^o\Sets) \to \Hh^W(\Delta^o\Delta^o\Sets)$ is
enhanced-left-exact.
\endproof

\begin{exa}\label{groth.l.e.exa}
By Corollary~\ref{ffun.compl.corr}, the enhanced functor category
$\E^\iota\cCat^h$ is cocomplete for any small enhanced category
$\E$. The embedding $a:\E^\iota\cCat^h \cong \cCat^h \mmh \E \to
\cCat^h \bbi^h \E$ of \eqref{3.h.enh.dia} is enhanced right-reflexive
by Corollary~\ref{enh.ind.corr}. Moreover, the forgetful fibration
$\cCat^h \bbi^h \E \to \cCat^h$ of \eqref{no.E} is actually given by
$\gamma_\trr$, where $\gamma:\E \to \ppt^h$ is the tautological
projection, so as in Subsection~\ref{enh.rel.subs}, it is
right-reflexive, with right-adjoint $\gamma^*$. Combining the two
observations, we see that the enhanced functor $\E^\iota\cCat^h \to
\cCat^h$ that sends $X:\E^\iota \to \cCat^h$ to $\E X$ is enhanced
right-reflexive, thus enhanced-right-exact by
Example~\ref{enh.adj.exa.exa}.
\end{exa}

\begin{defn}\label{enh.tot.loc.def}
The {\em total localization functor} $\Hh^{\Tot}:\cCat^h \to
\sSets^h$ is the enhanced functor left-adjoint to the embedding
$\sSets^h \subset \cCat^h$. An enhanced functor $\gamma:\C \to \C'$
between small enhanced categories is {\em $h$-invertible}
resp.\ {\em $h$-trivial} if $\Hh^\Tot(\gamma)$ is an equivalence
resp.\ factors through $\ppt^h$, and a small enhanced category $\C$
is {\em hyperconnected} if $\C \to \ppt^h$ is $h$-invertible.
\end{defn}

By adjunction, for any $\C \in \Cat^h$, its total localization
$\Hh^{\Tot}(\C)$ comes equipped with a $h$-invertible enhanced
functor $\lambda:\C \to \Hh^{\Tot}(\C)$, and then for any enhanced
cataegory $\E$, it induces an enhanced functor
\begin{equation}\label{loc.fun.eq}
\lambda^*:\fFun^h(\Hh^{\Tot}(\C),\E) \to \fFun^h(\C,\E).
\end{equation}
For purely formal reasons --- namely, $\sSets^h$ and $\cCat^h$ are
cartesian-closed, with internal $\Hom$ given by $\fFun^h(-,-)$ ---
the functor \eqref{loc.fun.eq} is an equivalence as soon as $\E$ is
an enhanced groupoid. In the general case, for the same purely
formal reasons, it is an enhanced fully faithful embedding onto
$\fFun^h(\C,\E_{h\Iso}) \subset \fFun^h(\C,\E)$.

\begin{exa}
In terms of the representing Segal spaces, $\Hh^{\Tot}$ is given by
the homotopy colimit functor $R^h_!$ of \eqref{colim.del.eq}. Thus
by Lemma~\ref{loc.N.le}, if a small enhanced category $\C$ has a
universal object $c \in \C^{\dm}_{J^\dm}$ for some $J^\dm \in
\Relf^+$, and $J = U(J^\dm)$, with the adjunction map $a:J^\dm \to
R(J)$, then $\Hh^{\Tot}(\C)^\dm$ is constant over $a$, and
$\lambda(c) = a^*c'$ for a unique universal object $c' \in
\Hh^{\Tot}(\C)^\dm_{R(J)}$. In words, if $J^\dm$ represents
$\C^\iota$, then $R(U(J^\dm))$ represents $\Hh^{\Tot}(\C)$. In
particular, the enhanced functor $\Unf(J^o) \to \C$ of
\eqref{enh.univ.sq} is $h$-invertible.
\end{exa}

\begin{exa}
For any anodyne map $f:J' \to J$ in $\Posf^+$, the enhanced functor
$\Unf(f):\Unf(J') \to \Unf(J)$ is $h$-invertible --- again, this
immediately follows from Lemma~\ref{loc.N.le} and Lemma~\ref{NQ.le}.
\end{exa}

\begin{exa}
Since total localization commutes with products, the product $\id
\times^h \gamma:\E \times^h \C' \to \E \times^h \C$ of the identity
endofunctor of a small enhanced category $\E$ and an $h$-invertible
enhanced functor $\gamma:\C' \to \C$ between small enhanced
categories $\C$, $\C'$ is $h$-invertible. In particular, for any
small enhanced $\C$, the projection $\C \times^h \Unf([1]) \to \C$
is $h$-invertible. This immediately implies that any left or
right-reflexive enhanced functor $\gamma:\C' \to \C$ between small
enhanced categories is $h$-invertible.
\end{exa}

\begin{exa}\label{iota.exa.exa}
The enhanced involution $\iota:\cCat^h \to \cCat^h$ of
\eqref{iota.enh.eq} is left and right-exact, and commutes with the
total localization functor: we have $\Hh^\Tot(\C) \cong
\Hh^\Tot(\C^\iota)$ for any $\C \in \cCat^h$.
\end{exa}

We note that while the proof of Proposition~\ref{ccat.coco.prop} is
almost tautological, it fully deserves to be called a proposition,
not a lemma. Firstly, it is crucially important for many
applications, and secondly, it really uses the equivalences
\eqref{ccat.h.eq}, \eqref{sset.h.eq} in an essential way. This is
the only place in the enhanced category theory where
representability is used directly, not through the universal objects
of Proposition~\ref{univ.prop}. Both the existence and the
properties of the total localization functor $\Hh^\Tot$ depend on
the existence and properties of the model structure on
$\Delta^o\Sets$. In particular, we have the following corollary of
Lemma~\ref{kappa.X.le}.

\begin{lemma}\label{kappa.E.le}
Assume given a small enhanced category $\E \in \Cat^h_\kappa$, for a
regular cardinal $\kappa$, and an enhanced functor $\gamma:\E \to
\C$ to a hyperconnected small enhanced category $\C$. Then there
exists a commutative square
\begin{equation}\label{hyper.sq}
\begin{CD}
\E_0 @>{\beta}>> \E_1\\
@V{\eta}VV @VVV\\
\E @>{\gamma}>> \C
\end{CD}
\end{equation}
of enhanced categories and enhanced functors such that $\E_0,\E_1
\in \Cat^h_\kappa$, $\eta$ is $h$-invertible, and $\beta$ is
$h$-trivial.
\end{lemma}

\proof{} The enhanced dual cone $\Cyl^\iota_h(\gamma)$ is
$[1]$-augmented, and $\Cyl^\iota(\gamma)_1 \cong \E$ is
$\kappa$-bounded. Thus by Proposition~\ref{univ.prop}, we have a
universal object for $\Cyl^\iota_h(\gamma)^\iota$ that defines a
fibration $J \to [1]$ with $J \in \Posf^+$ and $J_1 \in
\Posf^+_\kappa$, and a $[1]$-augmented enhanced functor $\Unf(J) \to
\Cyl^\iota_h(\gamma)$ $h$-invertible over both $0,1 \in [1]$. Then
since $\C$ is hyperconnected, so is $\Unf(J_0)$. Now take $\E_0 =
\Unf(J_1)$, and construct $\E_1 = \Unf(J_1')$ by applying
Lemma~\ref{kappa.X.le} to the transition map $J_1 \to J_0$ of the
fibration $J \to [1]$.
\endproof

Since $\cCat^h$ is complete and cocomplete, for any small enhanced
category $\C$, we have the functors
$\limh_\C,\colimh_\C:\fFun^h(\C,\cCat^h) \to \cCat^h$, and similarly
for $\sSets^h$. To compute them explicitly, it helps to use the
enhanced Grothendieck construction. Namely, assume given an enhanced
functor $\gamma:\C' \to \C$ between small enhanced categories, and
assume given enhanced functors $E:\C^\iota \to \cCat^h$,
$E':{\C'}^\iota \to \cCat^h$ with the corresponding enhanced
fibrations $\E \to \C$, $\E' \to \C'$. By abuse of notation, let
$\gamma_!\E',\gamma_*\E' \to \C$ be the enhanced fibrations
corresponding to $\gamma^\iota_!E'$, $\gamma^\iota_*E'$. Then the
adjunction maps $\gamma^*\gamma_*E' \to E' \to \gamma^*\gamma_!E'$
correspond to enhanced functors $\gamma^*\gamma_*\E' \to \E' \to
\gamma^*\gamma_!\E'$ cartesian over $\C'$, and together with
\eqref{gamma.st.fun}, these produce enhanced functors
\begin{equation}\label{ffun.cat.adj}
\begin{aligned}
\fFun^{\Toth}_{\C}(\gamma_!\E',\E) &\to
\fFun^{\Toth}_{\C'}(\gamma^*\gamma_!\E',\gamma^*\E) \to
\fFun^{\Toth}_{\C'}(\E',\gamma^*\E)\\
\fFun^{\Toth}_{\C}(\E,\gamma_*\E') &\to
\fFun^{\Toth}_{\C'}(\gamma^*\E,\gamma^*\gamma_*\E') \to
\fFun^{\Toth}_{\C'}(\gamma^*\E,\E').
\end{aligned}
\end{equation}
Now, \eqref{fun.fun.1} immediately implies that for any small
enhanced category $\C''$, we have $\gamma_!\E' \times^h \C'' \cong
\gamma_!(\E' \times^h \C'')$, and then \eqref{fun.fun.2} shows that
both enhanced functors \eqref{ffun.cat.adj} are actually
equivalences, so that the adjunction between $\gamma_!$, $\gamma^*$
and $\gamma_*$ also induces an adjunction on the level of enhanced
categories of enhanced functors.

\begin{lemma}\label{cath.lim.le}
Let $\C$ be a small enhanced category, and let $E:\C^\iota \to
\cCat^h$ be an enhanced functor, with the enhanced fibration $\E =
\C E \to \C$. Then $\limh_{\C^\iota}E \cong \sSec^{\Toth}(\C,\E)
\subset \sSec^h(\C,\E)$ is the full enhanced subcategory spanned by
cartesian sections, and $\colimh_{\C^\iota}E$ fits into an
enhanced-cocartesian square
\begin{equation}\label{enh.colim.sq}
\begin{CD}
\E_{h\flat} @>>> \E\\
@VVV @VVV\\
\Hh^\Tot(\E_{h\flat}) @>>> \colimh_{\C^\iota}E,
\end{CD}
\end{equation}
where $\E_{h\flat} \to \C$ is the enhanced family of groupoids of
\eqref{h.flat.sq}, and $\Hh^\Tot$ is the total localization functor
of Definition~\ref{enh.tot.loc.def},
\end{lemma}

\proof{} Since we have $\E' \cong \fFun^h(\ppt^h,\E')$ for any
enhanced category $\E'$, in particular for $\E' = \limh_{\C^\iota}E$,
the first claim immediately follows from the second equivalence in
\eqref{ffun.cat.adj}. For the second claim, if $\theta:\C^\iota \to
\ppt^h$ is the tautological projection, then the adjunction map
$E \to \theta^*\theta_!E = \theta^*\colimh_{\C^\iota}E$ provides an
enhanced functor $\lambda:\E \to \colimh_{\C^\iota}E$, functorially
with respect to $E$. If $E$ factors through $\sSets^h \subset
\cCat^h$, so that $\E$ is an enhanced family of groupoids, then
$\colimh_{\C^\iota}E$ is also an enhanced groupoid, and then
$\lambda$ factors through an enhanced functor
$\lambda':\Hh^{\Tot}(\E) \to \colimh_{\C^\iota}E$; to see that
$\lambda'$ is an equivalence, note that by \eqref{loc.fun.eq} and
the first equivalence in \eqref{ffun.cat.adj}, $\lambda'$ becomes an
equivalence after applying $\fFun^h(-,\E')$ for any $\E' \in
\Sets^h$. In the general case, since $\lambda$ is functorial, we
have a commutative square \eqref{enh.colim.sq}. By
Lemma~\ref{enh.cocart.le}, to check that it is enhanced-cocartesian,
it suffices to check that it becomes cartesian after applying
$\fFun(-,\E')$ for any $\E' \in \Cat^h$. Then again, by the first
equivalence in \eqref{ffun.cat.adj} and the first claim, this
reduces to observing that an enhanced functor $\E \to \E'$ over $\C$
is cartesian over $\C$ iff it becomes cartesian after restricting to
$\E_{h\flat} \subset \E$.
\endproof

\begin{remark}
For small enhanced categories $\Unf(J)$, $J \in \Pos$, the square
\eqref{enh.colim.sq} is an upgraded version of the square
\eqref{cone.sq} of Example~\ref{cone.sq.exa}.
\end{remark}

\begin{remark}\label{cocat.rem}
Note that since the involution $\iota$ of Example~\ref{iota.exa.exa}
preserves all limits and colimits, there exists a version of
Lemma~\ref{cath.lim.le} that uses a covariant version of the
Grothendieck construction and encodes enhanced functor to $\cCat^h$
by enhanced cofibrations rather than fibrations. Namely, for any
small enhanced category $\C$ and enhanced functor $E:\C \to
\cCat^h$, we have the corresponding enhanced cofibration $\C^\perp E
\to \C$. Then by Lemma~\ref{cath.lim.le},
\begin{equation}\label{cocath.fib.eq}
\limh_{\C}E \cong (\limh_{\C} (\iota \circ E))^\iota \cong
\sSec^{\Toth}(\C^\iota,(\C^\perp E)^\iota)^\iota  \cong
\sSec^h_\Tot(\C,\C^\perp E),
\end{equation}
where $\sSec^h_\Tot$ stands for the full enhanced subcategory in
$\sSec^h$ spanned by cocartesian sections. For colimit, one can also
write down an analogous covariant version of \eqref{enh.colim.sq};
we leave it to the reader.
\end{remark}

\begin{exa}\label{yo.colim.exa}
For any small enhanced category $\C$, the enhanced functor
$a:\C^\iota\cCat^h \cong \cCat^h \mmh \C \to \cCat^h \bbi^h \C$ has a
left-adjoint $\Y \circ b$ of Proposition~\ref{yo.prop}, and for any
enhanced functor $\gamma:\C \to \C'$ with small $\C$, we have the
adjoint pair of enhanced functors $\gamma^*$, $\gamma_\trr$ between
$\cCat^h \bbi^h \C'$ and $\cCat^h \bbi^h \C$. Then $a$ commutes with
$\gamma^*$, and by adjunction, if we take $\C'=\ppt^h$, we see that
$\colimh_{\C^\iota} \circ (\Y \circ b):\cCat^h \bbi^h \C \to \cCat^h$ is
the forgetful functor \eqref{no.E}. Explicitly, this means that for
any small enhanced category $\E$ and enhanced functor $\E \to \C$,
we have $\E \cong \colimh_{\C^\iota} c \setminus^h \E$, where $c
\setminus^h \E:\C^\iota \to \cCat^h$ is understood as the enhanced
functor corresponding to the enhanced fibration $\C \setminus^h \E
\to \C$ by \eqref{enh.groth.eq}. In terms of \eqref{enh.colim.sq},
this can be seen by observing that $\eta:\E_{h\Iso} \to \C
\setminus^h \E_{h\Iso} \cong (\C \setminus^h \E)_{h\flat}$, being
left-reflexive, is $h$-invertible.
\end{exa}

\begin{exa}\label{discr.pf.exa}
If an enhanced functor $\gamma:\E' = \E E \to \E$ is an enhanced
family of groupoids corresponding to some $E:\E^\iota \to \sSets^h
\subset \cCat^h$, then the enhanced functor
$\gamma^\iota_!:{\E'}^\iota\sSets^h \to \E^\iota\sSets^h$ is
particularly simple. Indeed, since all enhanced morphisms in a
family of groupoids are automatically cartesian, the embedding
$a:\E^\iota\sSets^h \to \cCat^h \bbi^h \E$ is fully faithful, and
similarly for $\E'$. Then since $\gamma$ is an enhanced family of
groupoids, the postcomposition functor $\gamma_\trr:\cCat^h \bbi^h \E'
\to \cCat^h \bbi^h \E$ sends ${\E'}^\iota\sSets^h$ into
$\E^\iota\sSets^h$, and we then simply have $\gamma^\iota_!\C \cong
\gamma_\trr\C$, for any enhanced family of groupoids $\C \to
\E'$. In particular, the constant family with fiber $\ppt^h$
corresponds to $\id:\E' \to \E'$, we have $\gamma^\iota_!(\ppt^h)
\cong E$, and then for any $\C \to \E'$, we have a semicartesian
square
\begin{equation}\label{discr.sq}
\begin{CD}
\C @>>> \gamma^{\iota *}\gamma^\iota_!\C\\
@VVV @VVV\\
\E' @>>> \gamma^{\iota *}\gamma^\iota_!\E',
\end{CD}
\end{equation}
where the horizontal arrows represent adjunction maps. For any
enhanced object $e \in \E'_\ppt$, the enhanced family of groupoids
$\sigma(e):\E' /^h e \to \E'$ corresponds to $\Y(\E)(e)$, and then
$\gamma^\iota_!\Y(\E')(e) \cong \Y(\E)(\gamma(e))$ since $\E' /^h e
\cong \E /^h \gamma(e)$, so that $\gamma^\iota_!\Y(\E') \cong
\gamma^{\iota *}\Y(\E)$. Moreover, we have $\ppt^h \cong
\colimh_{\E'}\Y(\E')$ by Example~\ref{yo.colim.exa}, so that
\begin{equation}\label{E.coli.eq}
E \cong \gamma^\iota_! \ppt^h \cong
\colimh_{\E'}\gamma^\iota_!\Y(\E') \cong \colimh_{\E'}\gamma^{\iota
  *}\Y(\E).
\end{equation}
Thus just as in the usual category theory, every enhanced functor
$\E^\iota \to \sSets^h$ is a colimit of representable functors.
\end{exa}

The simplest examples of limits and colimits in $\cCat^h$ are given
by the homotopy cocartesian and cocartesian squares of
Example~\ref{V.holim.exa}, and these actually reduce to
(enhanced-)semicartesian and enhanced-semicocartesian squares of
Subsection~\ref{univ.app.subs}. Namely, by Example~\ref{sq.fib.exa},
a commutative square of small enhanced categories $\C_{00}$,
$\C_{01}$, $\C_{10}$, $\C_{11}$ and enhanced functors between them
gives a fibration $\C \to \Posf^+ \times [1]^2$, thus an object in
the category $\Cat^h([1]^2|[1]^2)$. We have $[1]^2 \cong \V^> \cong
\V^{o<}$, this gives embeddings $\lambda:\V^o \to [1]^2$, $\rho:\V
\to [1]^2$, and Example~\ref{cyl.fib.exa} provides comparison
functors $r:\C_{11} \times \V \to \rho^*\C$ resp.\ $l:\lambda^*\C
\to \C_{00} \times \V^o$ cartesian over $\V$ resp. $\V^o$. These
define maps in $\Cat^h(\V|\V)$ resp.\ $\Cat^h(\V^o|\V^o)$. However,
since $\dim \V = \dim \V^o = 1$, we can identify $\Cat^h(\V|\V)$
resp.\ $\Cat^h(\V^o|\V^o)$ with $\Cat^h(\V)$ resp.\ $\Cat^h(\V^o)$
by Corollary~\ref{ccat.dim1.corr}, so $l$ resp.\ $l$ actually give
maps in $\Cat^h(\V)$ resp.\ $\Cat^h(\V^o)$. Then if we let
$\C_\idot:\Unf([1]^2)^\iota \to \Cat^h$ be the enhanced functor
corresponding to $\C$ by \eqref{enh.groth.eq}, we have natural
comparison functors
\begin{equation}\label{cocart.cmp}
\colimh_{\V} \lambda^{o*}\C_\idot \to \C_{00}, \qquad \C_{11} \to
\limh_{\V^o}\rho^{o*}\C^\hdot
\end{equation}
adjoint to the maps $l$ resp.\ $r$.

\begin{lemma}\label{ccat.sq.le}
A commutative square $\C \to \Posf^+ \times [1]^2$ of small enhanced
categories and enhanced functors is enhanced-semicocartesian
resp.\ semicartesian if and only if the first resp.\ second of the
comparison functors \eqref{cocart.cmp} is an equivalence.
\end{lemma}

\proof{} By Lemma~\ref{cath.lim.le}, the first claim amounts to an
enhanced version of \eqref{sec.sq.eq}, and this is immediately
provided by Lemma~\ref{sec.sq.le} applied to the tautological
standard pushout square $\V \to \V$. For the second claim, apply
\eqref{ffun.cat.adj} and Lemma~\ref{enh.cocart.le}.
\endproof

\begin{corr}
For any small enhanced categories $\C$, $\C_0$, $\C_1$ and enhanced
functors $\gamma_l:\C \to \C_l$, $l=0,1$, there exists an
enhanced-semi\-cocar\-tesian square \eqref{enh.cocart.sq}.
\end{corr}

\proof{} Clear. \endproof

\subsection{Limits and Kan extensions.}

Informally, Proposition~\ref{enh.compl.prop} says that enhanced
functors to a complete enhanced category admit right Kan extensions
that satisfy base change. To make this into a formal statement, we
need to define those Kan extensions. Fix an enhanced functor
$\gamma:\C \to \C'$ between small enhanced categories, with the
enhanced cylinder $\Cyl_h(\gamma)$ and the dual enhanced cylinder
$\Cyl^\iota_h(\gamma)$. For any enhanced functor $E:\C \to \E$ to
some enhanced category $\E$, define a {\em relative cone}
resp.\ {\em relative dual cone} for the functor $E$ as functor
$E_>:\Cyl_h(\gamma) \to \E$ resp.\ $E_<:\Cyl^\iota_h(\gamma) \to \E$
equipped with an isomorphism $s^*E_> \cong E$ resp.\ $t^*E_< \cong
E$. If $\E$ is small, then by \eqref{ffun.cyl.co}, relative cones
resp.\ relative dual cones form an enhanced category
$\cCone(E,\gamma) = E \setminus^h_{\gamma^*} \fFun(\C',\E)$
resp.\ $\cCone^\iota(E,\gamma) = \fFun(\C',\E) /^h_{\gamma^*} E$,
where $E$ is understood as an enhanced object in $\fFun(\C,\E)$. We
say that a relative cone resp.\ a relative dual cone is {\em
  universal} if it is the initial resp.\ terminal enhanced object in
the corresponding enhanced category. We note that since
semicartesian products preserve fully faithful embedding, a
universal relative cone $E_>$ that factors through a full enhanced
subcategory $\E' \subset \E$ is also universal as a relative cone
for a functor to $\E'$, and similarly for a universal relative dual
cone. If the target enhanced category $\E$ is not small, we say that
a relative cone resp.\ dual cone is {\em universal} if it is
universal as a relative cone resp.\ relative dual cone for functors
to any small full enhanced subcategory $\E' \subset \E$ through
which it factors.

\begin{defn}\label{enh.kan.def}
For any enhanced functor $\gamma:\C \to \C'$ between small enhanced
categories, an enhanced functor $E:\C \to \E$ to some enhanced
category $\E$ {\em admits a left} resp.\ {\em right enhanced Kan
  extension} with respect to $\gamma$ if it admits a universal
relatie cone resp.\ relative dual cone. If this happens, then an
enhanced functor $F:\E \to \E'$ {\em preserves} the Kan extension if
it sends the universal relative cone resp.\ relative dual cone to a
universal relative cone resp.\ relative dual cone.
\end{defn}

If an enhanced functor $E:\C \to \E$ admits a right Kan extension
with respect to some enhanced functor $\gamma:\C \to \C'$ in the
sense of Definition~\ref{enh.kan.def}, then this Kan extension is
defined as $\gamma_*E = s^*E_<:\C' \to \E$, where $E_<$ is the
universal relative dual cone, and it comes equipped with an enhanced
morphism $a:\gamma^*\gamma_*E \to E$. Dually, if $E$ admits a left
Kan extension, then this extension $\gamma_!E = t^*E_>:\C' \to \E$
comes equipped with an enhanced morphism $a:E \to
\gamma^*E\gamma_!E$. If $\C' = \ppt^h$, so that $\gamma:\C \to
\ppt^h$ is the tautological projection, then $\Cyl_h(\gamma)
\cong \C^{h>}$, $\Cyl^\iota_h(\gamma) \cong \C^{h<}$, relative cones
resp.\ dual cones are called simply cones resp.\ dual cones, we
write $\cCone(E)=\cCone(E,\gamma)$,
$\cCone^\iota(E)=\cCone^\iota(E,\gamma)$, and $\gamma_!E =
\colimh_{\C}E$, $\gamma_*E = \limh_{\C}E$ are called the {\em enhanced
  colimit} resp.\ {\em enhanced limit} of the enhanced functor $E$.

\begin{exa}
For any functor $\gamma:I \to I'$ between essentially small
categories, a functor $E:I \to \E$ to a category $\E$ admits a left
resp.\ right Kan extension with respect to $\gamma$ iff the enhanced
functor $\Unf(E)$ admits a left resp.\ right Kan extension with
respect to $\Unf(\gamma)$, and $\Unf(\gamma)_!\Unf(E) \cong
\Unf(\gamma_!E)$ resp.\ $\Unf(\gamma)_*\Unf(E) \cong \Unf(f_*E)$. A
functor $F:\E \to \E'$ preserves these Kan extensions iff so does
$\Unf(F)$.
\end{exa}

\begin{exa}\label{pos.enh.sq.exa}
The embedding $\Unf(\Posf^+) \to \cCat^h$ preserves any of the
standard cocartesian squares \eqref{J.sq} of Lemma~\ref{sec.sq.le}
(to see this, combine Lemma~\ref{sec.sq.le},
Lemma~\ref{enh.cocart.le} and Lemma~\ref{ccat.sq.le}).
\end{exa}

\begin{lemma}\label{kan.compo.le}
Assume given enhanced functors $\gamma;\C \to \C'$,
$\gamma':\C' \to \C''$ between small enhanced categories, and an
enhanced functor $E:\C \to \E$ that admits a right Kan extension
$\gamma_*\E$ with respect to $\gamma$. Then $\gamma_*E$ admits a
right Kan extension with respect to $\gamma'$ iff $E$ admits a right
Kan extension with respect to $\gamma''=\gamma' \circ \gamma$, and
if this happens, we have $\gamma'_*\gamma_*E \cong \gamma''_*E$.
\end{lemma}

\proof{} By definition, it suffices to consider the case when $\E$
is small. Then we have a semicartesian square
$$
\begin{CD}
\fFun^h(\C'',\E) /^h_{{\gamma'}^*} \fFun^h(\C',\E) /^h_{\gamma^*}
E @>>> \fFun^h(\C',\E) /^h_{\gamma^*} E\\
@VVV @VVV\\
\fFun^h(\C'',\E) /^h_{{\gamma'}^*} \fFun^h(\C',\E) @>>>
\fFun^h(\C',\E)
\end{CD}
$$
induced by \eqref{segal.sq} for $n=2$, $l=1$, and since $\gamma_*E$
exists, $\gamma'_*(\gamma_*E)$ exists iff $\fFun^h(\C'',\E)
/^h_{{\gamma'}^*} \fFun^h(\C',\E) /^h_{\gamma^*} E$ has a terminal
enhanced object. But then the left-reflexive map $s_\dg:[2] \to
[1]$, $0 \mapsto 0$, $1,2 \mapsto 1$ provides a left-admissible
fully faithful embedding
$$
\fFun(\C'',\E) /^h_{{\gamma''}^*} E \to \fFun^h(\C'',\E)
/^h_{{\gamma'}^*} \fFun^h(\C',\E) /^h_{\gamma^*} E,
$$
so this happens iff $\gamma''_*E$ exists, and then $\gamma''_*E
\cong \gamma'_*\gamma_*E$ by uniqueness.
\endproof

\begin{lemma}\label{kan.univ.le}
Assume given an enhanced functor $\gamma:\C \to \C'$ between small
enhanced categories, and an enhanced functor $E:\C \to \E$ that
admits a right Kan extension $\gamma_*E$, with the corresponding map
$a:f^*f_*E \to E$. Then for any functor $E':\C' \to \E$ equipped
with a map $f:\gamma^*E' \to E$, there exists a unique map $f':E'
\to \gamma_*E$ such that $f=a \circ \gamma^*(f)$.
\end{lemma}

\proof{} As in Lemma~\ref{kan.compo.le}, it suffices to consider the
case when $\E$ is small. Then by the definition of enhanced
comma-fibers, we have an epivalence
\begin{equation}\label{fun.co.fib}
\cCone(E,\gamma)_\ppt=(\fFun^h(\C',\E) /^h_{\gamma^*} E)_\ppt \to
\Fun^h(\C',\E) /_{\gamma^*} E,
\end{equation}
and if its source has a terminal object, then so does its target.
\endproof

Both Lemma~\ref{kan.compo.le} and Lemma~\ref{kan.univ.le} have
obvious dual version for left Kan extensions obtained by taking the
enhanced opposite categories; we leave the precise formulations to
the reader. Note that the universal property of
Lemma~\ref{kan.univ.le} does {\em not} insure that an enhaned
functor is the right Kan extension $f_*E$. Firstly,
\eqref{fun.co.fib} is only an epivalence, so even if its target has
a terminal object, this need not hold for its source; secondly, even
if the terminal object lifts to the source, it need not be
enhanced-ter\-mi\-nal (see Example~\ref{enh.term.exa}). However, if we
do know that the right Kan extension $f_*E$ exists,
Lemma~\ref{kan.univ.le} completely characterizes it. Moreover, it is
functorial: if Kan extensions $\gamma_*E$, $\gamma_*E'$ exist for
two enhanced functors $E,E':\C \to \E$, then an enhanced map $f:\E
\to \E'$ gives rise to a unique enhanced map $\gamma_*(f):\gamma_*E
\to \gamma_*E'$. There is even more functoriality if $f_*E$ exists
for all $E$.

\begin{lemma}\label{enh.kan.com.le}
For any enhanced functor $\gamma:\C \to \C'$ between small enhanced
categories $\C$, $\C'$, and any enhanced category $\E$, the
following conditions are equivalent.
\begin{enumerate}
\item The pullback enhanced functor $\gamma^*:\fFun^h(\C',\E) \to
  \fFun^h(\C,\E)$ admits a left resp.\ right-adjoint enhanced
  functor $\gamma_!$ resp.\ $\gamma_*$.
\item Any enhanced functor $E:\C \to \E$ admits a left resp.\ right
  Kan extension with respect to $\gamma$.
\end{enumerate}
\end{lemma}

\proof{} The two statements are equivalent by passing to the
enhanced-opposite categories, so we only consider the right Kan
extensions and right-adjoints. If $\E$ is small, the claim
immediately follows from Corollary~\ref{adj.comma.corr}. In the
general case, note that $f_*E$ is well-defined for any $E:\C \to \E$
under either of the assumptions \thetag{i}, \thetag{ii}, and say
that a full enhanced subcategory $\E' \subset \E$ is {\em
  $\gamma$-closed} if for any enhanced functor $E:\C \to \E$ that
factors through $\E$', so does $f_*E:\C' \to \E'$. Then for any
small enhanced full subcategory $\E' \subset \E$, the smallest
$\gamma$-closed enhanced full subcategory containing $\E$ is still
small. Therefore $\E$ is covered by $\gamma$-closed small enhanced
full subcategories $\E' \subset \E$, and since $\C$, $\C'$ are
small, any enhanced functor $\C,\C' \to \E$ factors through such a
$\E'$. Then \thetag{i} for $\E$ implies \thetag{i} for each $\E'$,
and this yields \thetag{ii}. Conversely, \thetag{ii} implies
\thetag{i} for any $\gamma$-closed small $\E' \subset \E$, and this
implies \thetag{i} for $\E$ by the uniqueness of adjoints.
\endproof

By virtue of Lemma~\ref{enh.kan.com.le},
Proposition~\ref{enh.compl.prop} indeed says that any enhanced
functor to a complete resp.\ cocomplete enhanced category $\E$
admits a right resp.\ left Kan extension, and our notation is
consistent. Moreover, these Kan extensions are preserved by any
enhanced-left-exact resp.\ enhanced-right-exact enhanced functor
$F:\E \to \E'$. In particular, $\sSets^h$ is complete and
cocomplete, and by Corollary~\ref{ffun.compl.corr}, so is the
functor category $\E^\iota\sSets^h$ for any small enhanced category
$\E$.

\begin{defn}\label{univ.kan.def}
For any enhanced functor $\gamma:\C \to \C'$ between small enhanced
categories $\C$, $\C'$, and an enhanced functor $E:\C \to \E$ to
some enhanced category $\E$ that admits a right Kan extension
$\gamma_*E$, the Kan extension is {\em universal} if for any small full
enhanced subcategory $\E_0 \subset \E$ such that $\gamma_*E$ factors
through $\E_0$, it is preserved by the Yoneda embedding $\Y(\E_0)$
of \eqref{yo.yo.enh.eq}.
\end{defn}

\begin{lemma}\label{lim.univ.le}
All enhanced limits are universal in the sense of
Definition~\ref{univ.kan.def}.
\end{lemma}

\proof{} For any small enhanced $\C$, $\E$ and an enhanced functor
$E:\C \to \E$, we can compute $\limh_{\C} (\Y(\E) \circ E)$ by
Lemma~\ref{cath.lim.le}, and this shows that the corresponding
enhanced family of groupoids over $\E$ is exactly the enhanced cone
category $\cCone(E)$. If $\limh_{\C}E$ exists, then by definition,
this category has a terminal enhanced object, and then $\limh_{\C}
(\Y(\E) \circ E)$ lies in the image of $\Y(\E)$ by
Lemma~\ref{yo.im.le}.
\endproof

\begin{corr}\label{enh.kan.corr}
Assume given an enhaced functor $\gamma:\C \to \C'$ between small
enhanced categories, and an enhanced functor $E:\C \to \E$ to some
enhanced category $\E$. Moreover, assume that for any enhanced
object $c \in \C'_\ppt$, with the corresponding comma-fiber $c
\setminus^h_\gamma \C$ and the projection $\tau(c)$ of
\eqref{si.ta.c.eq}, $\limh_{c \setminus^h \C}\tau(c)^*E$ exists. Then
$E$ admits a universal right Kan extension $\gamma_*E$, and we have
a base change isomorphism
\begin{equation}\label{enh.kan.eq}
  \gamma_*E(c) \cong \limh_{c \setminus^h \C}\tau(c)^*E
\end{equation}
for any $c \in \C'_\ppt$.
\end{corr}

\proof{} For any small full enhanced subcategory $\E_0 \subset \E$
that contains $\gamma(c)$ for any $c \in \C_\ppt$ and $\limh_{c
  \setminus^h \C}\tau(c)^*E$ for any $c \in \C'_\ppt$, with the
Yoneda embedding $\Y(\E_0):\E_0 \to \E_0^\iota\sSets^h$, the Kan
extension $\gamma_*(\Y(\E_0) \circ E)$ satisfies \eqref{enh.kan.eq}
by Proposition~\ref{enh.compl.prop}, so that the corresponding
universal relative cone factors through $\E_0 \subset
\E_0^\iota\sSets^h$ by Lemma~\ref{lim.univ.le}.
\endproof

\begin{remark}
The converse to Corollary~\ref{enh.kan.corr} is wrong --- already in
the unenhanced setting, there are exotic functors $E$ such that the
Kan extension $\gamma_*E$ exists, but some of the limits in the
right-hand side of \eqref{enh.kan.eq} do not. For example, take the
discrete category $\E = \{0,1\}$. Then $\Fun(I,\E) \cong \E$ for any
connected small category $I$, so for any functor $\gamma:I' \to I$
between connected small categories and any $E:I' \to E$, $\gamma_*E$
trivially exists, but if a comma-fiber $i \setminus_\gamma I'$ is
empty, $\lim_{i \setminus I'}\sigma(i)^*E$ does not. Such Kan
extensions of course cannot be universal.
\end{remark}

\begin{corr}
An enhanced category $\E$ is complete resp.\ cocomplete if and only
if for any enhanced functor $E:\C \to \E$ from a small enhanced
cateogory $\C$, there exists $\limh_\C E$ resp.\ $\colimh_\C E$.
\end{corr}

\proof{} The ``only if'' part is
Proposition~\ref{enh.compl.prop}. Conversely, to prove that $\E$ or
$\E^\iota$ is complete, it suffices to check the condition of
Lemma~\ref{compl.bar.le}, and by virtue of
Lemma~\ref{enh.kan.com.le}, all these conditions reduce to checking
something for enhanced functors to $\E$ that factor through a small
full enhanced subcategory $\E_0 \subset \E$. Then by
Lemma~\ref{lim.univ.le}, everything can be checked after
composing our enhanced functors with the Yoneda embedding
$\Y(\E_0)$, and the category $\E_0^\iota\sSets^h$ is complete
resp.\ cocomplete by Corollary~\ref{ffun.compl.corr}.
\endproof

\begin{exa}\label{lim.colim.exa}
Assume given small enhanced categories $\C_0$, $\C_1$, and an
enhanced functor $X:\C_0 \times^h \C_1 \to \E$ to some complete
cocomplete enhanced category $\E$. By abuse of notation, we will
denote $\colimh_{\C_0}X = \pi_{0!}X:\C_1 \to \E$, $\limh_{\C_1}X =
\pi_{1*}X:\C_0 \to \E$, where $\pi_l:\C_0 \times^h \C_1 \to \C_l$,
$l=0,1$ are the projections. By the universal properties of limits
and colimits, we then have a natural map $\colimh_{\C_0}\limh_{\C_1}X
\to \limh_{\C_1}\colimh_{\C_0}X$ that may or may not be an
isomorphism. To check whether it is, it is useful to observe that if
we let $Z = \colimh_{\C_0}X$, $Y=\limh_{\C_1}X$, then the map
$\colimh_{\C_0}Y \to \limh_{\C_1}Z$ factors
as
\begin{equation}\label{cofun.facto}
\begin{CD}
\colimh_{\C_0} Y @>{\alpha}>> \colimh_{\fFun^h(\C_1,\C_0)}Y'
@>{\beta}>> \limh_{\C_1}Z,
\end{CD}
\end{equation}
where $\alpha$ is induced by $\gamma^*:\C_0 \to \fFun^h(\C_1,\C_0)$,
with $\gamma:\C_1 \to \ppt^h$ being the tautological projection, and
we denote $Y' = \limh_{\C_1}(\ev \times \id)^*X$, where $\ev:\C_1
\times^h \fFun^h(\C_1,\C_0) \to \C_0$ is the evaluation pairing of
\eqref{enh.ev.eq}.

If $\E = \sSets^h$, then both $\alpha$ and $\beta$ in
\eqref{cofun.facto} admit a convenient description in terms of the
Grothendieck construction. Namely, let $\Zz = \C_1^\perp Z \to
\C_1$, $\Xx = (\C_0 \times^h \C_1)^\perp X \to \C_0 \times^h \C_1$
be the cofibrations corresponding to $Z$ and $X$. Then the
adjunction map $X \to \pi_0^*Z$ induces an enhanced functor $\pi:\Xx
\to \Zz$ over $\C_1$, and by Example~\ref{enh.fam.grp.exa}, $\pi$ is
an enhanced cofibration. Moreover, by Example~\ref{discr.pf.exa}, we
have $\pi_!\ppt^h \cong \ppt^h$, so that $\pi$ has hyperconnected
fibers. Then if we denote $\Yy' = \fFun^h(\C_1,\C_0)^\perp Y'$,
Lemma~\ref{cath.lim.le} provides equivalences $\Yy' \cong
\sSec^h(\C_1,\Xx)$, $\limh_{\C_1}Z \cong \sSec^{\Toth}(\C_1,\Zz)$,
and $\beta$ in \eqref{cofun.facto} is the total localization of
$\sSec^h(\C_1,\pi)$ (where $\sSec^h(\C_1,\Zz) \cong
\sSec^{\Toth}(\C_1,\Zz)$ since the fibers of the cofibration $\Zz
\to \C_1$ are enhanced groupoids). On the other hand, if we let $\Yy
= \C_0^\perp Y$, then $\Yy \cong \sSec^{\Toth}(\C_1,\Xx)$, and
$\alpha$ in \eqref{cofun.facto} is the total localization of the
embedding $\sSec^{\Toth}(\C_1,\Xx) \to \sSec^h(\C_1,\Xx)$.
\end{exa}

\subsection{Kan extensions and functoriality.}

By passing to the enhanced-oppo\-site categories,
Corollary~\ref{enh.kan.corr} has a dual version for left Kan
extensions: we say that a left Kan extension $\gamma_!E$ is
universal if so is $(\gamma_!E)^\iota = \gamma^\iota_*E^\iota$, and
then for an enhanced functor $\gamma:\C \to \C'$ between small
enhanced categories, and any enhanced functor $E:\C \to \E$ such
that $\colimh_{\C /^h c}\sigma(c)^*E$ exists for any $c \in
\C'_\ppt$, a universal left Kan extension $\gamma_!E$ exists, and we
have
\begin{equation}\label{enh.cokan.eq}
  \gamma_!E(c) \cong \colimh_{\C /^h c}\sigma(c)^*E, \qquad c \in
  \C'_\ppt.
\end{equation}
This is the enhanced counterpart of \eqref{kan.eq}. In a nutshell,
Corollary~\ref{enh.kan.corr} and its dual show that if we restrict
our attention to universal Kan extensions, then they inherit all the
properties of the adjoint functors of
Proposition~\ref{enh.compl.prop}. In particular, we have the base
change isomorphisms of Corollary~\ref{enh.bc.corr} and its dual ---
these follow immediately from \eqref{enh.kan.eq} and
\eqref{enh.cokan.eq}. Another corollary is that if $\gamma:\C \to
\C'$ is fully faithful, then for any $E:\C \to \E$ that admits a
universal left resp.\ right Kan extension, the adjunction map $E \to
\gamma^*\gamma_!E$ resp.\ $\gamma^*\gamma_*E \to E$ is an
isomorphism. Thus for a fully faithful $\gamma$,
Lemma~\ref{kan.univ.le} provides a factorization
\begin{equation}\label{enh.kan.facto}
\begin{CD}
\C @>{\gamma}>> \C' @>{E'}>> \E
\end{CD}
\end{equation}
of the enhanced functor $E$ with $E' = \gamma_*E$, and for any other
factorization \eqref{enh.kan.facto}, there is a unique map $E' \to
\gamma_*E$ that makes the diagram commute. For the left Kan
extensions, the claim is dual: we have a factorization
\eqref{enh.kan.facto} with $E' = \gamma_!E$, and a unique map
$\gamma_!E \to E'$ for any other such factorization. This is an
enhanced version of Example~\ref{kan.facto.exa}.

As an application of the factorizations \eqref{enh.kan.facto}, we
can define Kan extensions along enhanced functors $\gamma:\C \to
\C'$ whose target is not small. This uses the filtered $2$-colimits
of Lemma~\ref{enh.filt.ff.le}. Assume given an enhanced functor
$\gamma:\C \to \C'$ such that $\C$ is small, and some enhanced
functor $E:\C \to \E$. Say that $E$ {\em admid a right Kan
  extension} with respect to $\gamma$ if for any small full enhanced
subcategory $\C'_0 \subset \C$ such that $\gamma$ factors through
$\C'_0$, $E$ admits a universal right Kan extension with respect to
the corresponding functor $\gamma_0:\C \to \C'_0$. Then in such a
situation, for any two such small full enhanced subcategories $\C'_0
\subset \C'_1 \subset \C'$, $\gamma_{0*}E$ admits a right Kan
extension with respect to the fully faithful embedding functor
$\nu:\C'_0 \to \C'_1$ by Lemma~\ref{kan.compo.le}, we have
$\nu_*\gamma_{0*}E \cong \gamma_{1*}E$, and since $\gamma_{l*}E$,
$l=0,1$ are universal, so is $\nu_*\gamma_{0*}$. In particular,
since $\nu$ is fully faithful we have the adjunction isomorphism
$\nu^*\gamma_{1*}E \cong \gamma_{0*}E$.  Now as in
Example~\ref{filt.ff.exa}, we can represent $\C' \cong
\colim^2_I\C'_\idot$, where $I$ is formed by all small full enhanced
subcategories in $\C'$ through which $\gamma$ factors, and then for
any $i \in I$, we have an enhanced functor $\gamma_{i*}E:\C'_i \to
\E$, and for any map $i \to i'$ in $I$, we have the adjunction
isomorphism $\nu^*\gamma_{i'}E \cong \gamma_{i*}E$. Taken together,
these define a functor $\C'_\idot \to E$ that inverts maps
cocartesian over $I$, thus descends to an enhanced functor
$\gamma_*E:\C' \to \E$. This is our right Kan extension. Note that
with this definition, Corollary~\ref{enh.kan.corr} still holds (and
all the enhanced comma-fibers in the right-hand side of
\eqref{enh.kan.eq} are small). If $\gamma$ is fully faithful, we
still have a universal factorization \eqref{enh.kan.facto}. Passing
to the enhanced-opposite categories, we obtain the dual story for
left Kan extensions including \eqref{enh.cokan.eq}.

\begin{exa}\label{com.col.exa}
When computing Kan extensions along an enhanced functor $\gamma:\C
\to \C'$ whose target is not small, one has to be careful:
technically speaking, to define the comma-fibers in the right-hand
side of \eqref{enh.kan.eq} or \eqref{enh.cokan.eq}, we have to
factor $\gamma$ through a small enhanced full subcategory $\C'_0$
that contains the enhanced object $c$. Enlarging $\C'_0$ changes the
comma-fibers by a canonical equivalence over $\C$, and since all
small full enhanced subcategories containg $\gamma(\C) \subset \C'$
and $c$ form a filtered category, $\C /^h_\gamma c$ is well-defined
as an object in $\Cat^h \mmh \C$, and similarly for $c
\setminus^h_\gamma \C$. To see that the correspondence $c \mapsto
\langle \C /^h_\gamma c,\sigma(c) \rangle$ is functorial with
respect to $c$, one can look at the universal situation: take the
Yoneda embedding $\Y(\C):\C \to \C^\iota\sSets^h$ of
\eqref{yo.yo.enh.eq}, and consider the Kan extension
$\gamma_!Y(\C):\C' \to \C^\iota\sSets^h$ that exists since its
target is cocomplete. Then since for any $c$, $\sigma(c):\C
/_\gamma^h c \to \C$ is an enhanced family of groupoids,
$\sigma(c)^*\Y(\C)$ in \eqref{enh.cokan.eq} is the enhanced functor
$\Y(\C,\C /^h_\gamma c)$ of Example~\ref{yo.colim.exa}, and we have
\begin{equation}\label{com.Y.col}
\gamma_!\Y(\C)(c) \cong \colimh_{\C /^h c}\sigma(c)^*\Y(\C) \cong \colimh_{\C
  /^h c}\Y(\C,\C /^h c) \cong \C /^h c,
\end{equation}
with the enhanced fibration $\sigma(c):\C /^h c \to \C$. If
$\C'=\C^\iota\sSets^h$, $\gamma=\Y(\C)$, we tautologically have
$\Y(\C)_!\Y(\C) \cong \Y(\C)$ by \eqref{enh.kan.facto}, and
\eqref{com.Y.col} simply reads as
\begin{equation}\label{Y.comma.eq}
  \C E \cong \C /^h_{\Y(\C)} E,
\end{equation}
for any enhanced functor $E:\C^\iota \to \sSets^h$.
\end{exa}

\begin{exa}
For any enhanced object $e \in \E_\ppt$ in an enhanced category
$\E$, we have the embedding $\eps^h(e):\ppt^h \to \E$, and if we
consider the constant enhanced functor $\ppt^h:\ppt^h \to \sSets^h$
with value $\ppt^h$, then by \eqref{enh.cokan.eq}, it admits a left
Kan extension
\begin{equation}\label{yo.e.eq}
  \Y(\E,e) = \eps^h(e)_!\ppt^h:\E \to \sSets^h
\end{equation}
with respect to $\eps^h(e)$, with $\Y(\E,e)(e') \cong (e \setminus^h
\E)_{e'}$ for any enhanced object $e' \in \E_\ppt$. If $\E$ is
small, $\Y(\E,e)$ corresponds to the enhanced cofibration $\tau:e
\setminus^h \E \to \E$ by the covariant version of the Grothendieck
construction, and $\Y(\E,e)$ is simply the restriction of the Yoneda
pairing \eqref{yo.pair.eq} with respect to the embedding
$\eps^h(e)^\iota \times^h \id:\E \cong \ppt^h \times^h \E \to
\E^\iota \times^h \E$. But even if $\E$ is not small, the enhanced
functor \eqref{yo.e.eq} is still well-defined.
\end{exa}

\begin{defn}\label{enh.gen.def}
A small full enhanced subcategory $\C \subset \C'$ in an enhanced
category $\C'$, with the embedding functor $\gamma:\C \to \C'$, is
{\em generating} if the enhanced functor $\gamma_!\Y(\C):\C' \to
\C^\iota\sSets^h$ of Example~\ref{com.col.exa} is fully faithful.
\end{defn}

\begin{exa}
Any small enhanced category $\C$, with the Yoneda embedding,
is trivially generating in $\C^\iota\sSets^h$ (since $\Y(\C)_!\Y(\C)
\cong \id$). In general, $\C \subset \C'$ is generating iff $\C'
\subset \C^\iota\sSets^h$ is a full enhanced subcategory containing
the Yoneda image of $\C$.
\end{exa}

\begin{exa}\label{gen.CC.exa}
For any full embedding $\gamma:\C \subset \C'$ of small enhanced
categories, $\C$ is generating in $\C'$ iff $\Y(\C') \cong
\gamma^\iota_!  \circ \Y(\C)$. Thus for any generating small full
enhanced subcategory $\C \subset \E$ in some $\E$, with the
embedding functor $\gamma:\C \to \E$, and any small full enhanced
subcategory $\C' \subset \E$, with the embedding functor
$\gamma':\C' \to \E$, we have $\gamma'_!\Y(\C') \cong \gamma^\iota_!
\circ \gamma_!\Y(\C)$, so that $\C' \subset \E$ is also generating.
\end{exa}

We note that in the situation of Definition~\ref{enh.gen.def}, the
adjunction map $\gamma_!\Y(\C) \circ \gamma = \gamma^*\gamma_!\Y(\C)
\to \Y(\C)$ induces an adjunction between $\gamma_!\Y(\C)$ and
$\gamma$ over $\Y(\C)$ in the sense of
Definition~\ref{adj.rel.def}. If the target enhanced category $\C'$
is also cocomplete, then we have a Kan extension
$\Y(\C)_!\gamma:\C^\iota\sSets^h \to \C'$ that is left-adjoint to
$\gamma_!\Y(\C)$ in the usual sense. Thus if $\C \subset \C'$ is
generating, $\Y(\C)_!\gamma:\C^\iota\sSets^h \to \C'$ is essentially
surjective, so that by \eqref{E.coli.eq}, any enhanced object $c'
\in \C'$ is a colimit of objects in $\C \subset \C'$. This explains
our terminology.

As another application of Kan extensions, we can describe
functoriality of the Yoneda embeddings \eqref{yo.enh.eq}. The
simplest statement here immediately follows from
\eqref{enh.cokan.eq}: for any enhanced functor $\gamma:\E_0 \to
\E_1$ between small enhanced categories $\E_0$, $\E_1$, we have a
commutative square
\begin{equation}\label{Y.ga.sq}
\begin{CD}
\E_0 @>{\Y(\E_0)}>> \E_0^\iota\sSets^h\\
@V{\gamma}VV @VV{\gamma^\iota_!}V\\
\E_1 @>{\Y(\E_1)}>> \E_1^\iota\sSets^h.
\end{CD}
\end{equation}
More generally, assume given an enhanced fibration $\E \to \C$ with
a small base, with enhanced fibers $\E_c$. Then we cannot say that
the functor categories $\E_c^o\sSets^h$ form an enhanced fibration
over $\C^\iota$ with transition functors \eqref{Y.ga.sq} since these
enhanced categories are not small. However, the category of sections
of this non-existent enhanced fibration is well-defined: this is
simply the enhanced functor category $\E^\iota\sSets^h$. If $\C =
\Unf(I)$ for a partially ordered set $I$, so that $\E$ is an
$I$-augmented small enhanced category, then we can also consider a
family $(\E|I)^o\sSets^h \to \Unf(I|I) = \Posf^+ \times I$ of enhanced
categories in the sense of Subsection~\ref{univ.funct.subs}, with
fibers $\E_i^\iota\sSets^h$, $i \in I$, and transition functors
$\gamma^o_!$ of \eqref{Y.ga.sq}. We then have a comparison functor
\begin{equation}\label{E.I.Y.eq}
\E^o\sSets^h \to \Sec(I^o,(\E|I)^o\sSets^h_\perp)
\end{equation}
obtained by taking a universal object for $\E$ of
Proposition~\ref{univ.prop}, with the corresponding fibration $J \to
I$, and then taking the functor \eqref{tr.big.loc}. For any small
enhanced category $\C$, we can consider the transpose enhanced
cofibration $\E_{h\perp} \to \C^\iota$, and then the Yoneda pairings
\eqref{yo.pair.eq} induce a relative Yoneda pairing $\E_{h\perp}
\times^h_{\C^\iota} \E^\iota \to \sSets^h$, and we also have the
evaluation functor $\sSec^h(\C_0,E_{h\perp}) \times^h \C_0 \to
\E_{h\perp}$. Taken together, these define a pairing
\begin{equation}\label{yo.sec.pair}
\begin{aligned}
\sSec^h(\C^\iota,\E_{h\perp}) \times^h \E^\iota &\cong
(\sSec^h(\C^\iota,\E_{h\perp}) \times^h \C^\iota)
\times^h_{\C^\iota} \E^\iota \to\\ &\to \E_{h\perp}
\times^h_{\C^\iota} \E^\iota \to \sSets^h
\end{aligned}
\end{equation}
that produces a relative version
\begin{equation}\label{yo.sec.eq}
\Y(\E|\C):\sSec^h(\C^\iota,\E_{h\perp}) \to \E^\iota\sSets^h
\end{equation}
of the Yoneda embedding \eqref{yo.enh.eq}. This is functorial with
respect to $\C$, and reduces to \eqref{yo.enh.eq} when $\C =
\ppt^h$. Explicitly, \eqref{perp.eq} and \eqref{Y.co.eq} provide an
equivalence
\begin{equation}\label{ssec.tw.eq}
\sSec(\C^o,\E_{h\perp}) \cong
\fFun^h_{\Tot\C^o}(\Ar^h(\C^o),\E_{h\perp}) \cong
\fFun^{\Toth}_{\C}(\Tw^h(\C),\E),
\end{equation}
and if some enhanced object $E \in \sSec^h(\C^o,\E_{h\perp})_\ppt$
corresponds to an enhanced functor $\Tw(E):\Tw^h(\C) \to \E$ under
\eqref{ssec.tw.eq}, then we have
\begin{equation}\label{Y.tw.eq}
\Y(\E|\C)(E) \cong \Tw(E)^o_!\ppt^h,
\end{equation}
a relative version of \eqref{yo.e.eq}. For any enhanced functor
$X:\E^\iota \to \sSets^h$, with the corresponding enhanced family of
groupoids $\E X \to \E$, \eqref{Y.tw.eq} combined with
\eqref{ssec.tw.eq} provides an equivalence
\begin{equation}\label{tw.Y.col}
\sSec^h(\C^o,\E_{h\perp}) /h_{\Y(\E|\C)} X \cong \sSec^h(\C^o,(\E
X)_{h\perp}),
\end{equation}
where $(\E X)_{h\perp} \to \C^o$ is transpose to the enhanced
fibration $\E X \to \E \to \C$. This is the relative version of
\eqref{Y.comma.eq}.

One can also consider contravariant enhanced functors $\E^\iota \to
\C$ from a small enhanced category $\E$ to a cocomplete enhanced
category $\C$. For any enhanced objects $e \in \E_\ppt$ and $c \in
\C_\ppt$, we define the representable functor $\Y(\E,e,c) =
\eps(e)^o_!c \in \E^\iota\C$. To construct all these functors at
once, observe that any cocomplete enhanced $\C$ admits a natural
action $\sSets^h \times^h \C \to \C$ of $\sSets^h$ obtained as
$(\eps(\ppt^h) \times \id)_!\id$, where $\eps(\ppt^h):\ppt^h \to
\sSets^h$ corresponds to the enhanced object $\ppt^h \in
\Sets^h$. Combining it with the Yoneda pairing \eqref{yo.pair.eq}
for $\E$, we obtain a pairing $\E^\iota \times^h \E \times^h \C \to
\sSets^h \times^h \C \to \C$ that gives rise to a $\C$-valued
version
\begin{equation}\label{yo.E.C.eq}
\Y(\E,\C):\E \times^h \C \to \E^\iota\C
\end{equation}
of the Yoneda embedding. On objects, we have $\Y(\E,\C)(e \times c)
= \Y(\E,e,c)$. As an application of \eqref{yo.E.C.eq}, we can
construct a version of \eqref{E.coli.eq} that shows that in fact for
any cocomplete $\C$, every enhanced functor $X:\E^\iota \to \C$ is a
colimit of representable ones; explicitly, we have
\begin{equation}\label{Y.C.X}
X \cong \colimh_{\Tw^h(\E)}\Y(\E,\C) \circ \pi(X),
\end{equation}
where $\pi(X) = \sigma \times^h (X \circ \tau):\Tw^h(\E) \to \E
\times^h \C$ is the product of the projection $\sigma:\Tw^h(\E) \to
\E$ and the composition $X \circ \tau:\Tw^h(\E) \to \E \to \C$, and
$\Y(\E,C)$ is the enhanced functor \eqref{yo.E.C.eq}.

\begin{exa}
Take $\C=\cCat^h$. Then $\E^\iota\C \cong \cCat^h \mmh \E$, and
\eqref{yo.E.C.eq} says that for any small enhanced fibration $\C \to
\E$, we have
\begin{equation}\label{Y.cat.eq}
\C \cong \colimh_{\Tw^h(\E)} (\E /^h e) \times^h \C_{e'},
\end{equation}
where the colimit is over all enhanced objects in $\Tw^h(\E)$
represented by enhanced morphisms $f:e \to e'$ in $\E$. This is a
generalization of the cylindrical colimits of
Subsection~\ref{pos.cyl.subs}.
\end{exa}

\subsection{Localizations and cofinality.}

To compute left and right Kan extensions by \eqref{enh.kan.eq} and
\eqref{enh.cokan.eq}, it helps to know when particular limits or
colimits exists. To control this, the following notion is very
useful.

\begin{lemma}\label{cofi.le}
Assume given an enhanced functor $\gamma:\C \to \C'$ between small
enhanced categories. Then the following conditions are equivalent.
\begin{enumerate}
\item For any enhanced object $c \in \C'_\ppt$, the enhanced
  comma-fiber $\C /^h_\gamma c$ is hyperconnected.
\item We have $\gamma_!\ppt^h \cong \ppt^h$.
\item For any enhanced category $\E$ and enhanced functor $E:\C' \to
  \E$, $\limh_{\C'}E$ exists if and only if so does
  $\limh_{\C}\gamma^*E$, and if both do, the adjunction map
  $\limh_{\C'} E \to \limh_{\C}\gamma^*\E$ is an isomorphism.
\end{enumerate}
\end{lemma}

\proof{} To see that \thetag{iii} implies \thetag{i}, take $\E =
(\sSets^h)^\iota$ and the enhanced functor $E = \Y(\C',c)^\iota:\C'
\to \E$. Then since the total localization of a small enhanced
category with an initial or terminal enhanced object is trivial, the
limit $\limh_{\C'} E \cong \colimh_{{\C'}^\iota}\Y(\C',c) \cong
\ppt^h)$ is trivial by Lemma~\ref{yo.im.le} and
Lemma~\ref{cath.lim.le}, and then $\limh_{\C}\gamma^*E \cong
\colimh_{\C^o}\gamma^*\Y(\C'c) \cong \Hh^{\Tot}(\C /^h_\gamma c)$.

Conversely, \thetag{i} implies \thetag{ii} by \eqref{enh.kan.eq},
and to see that this implies \thetag{ii}, it suffices to check it
for a small $\E$. Then by Lemma~\ref{lim.univ.le}, it suffices to
check that the map $\limh_{\C'}(\Y(\E) \circ E) \to \limh_{\C}(\Y(\E)
\circ E \circ \gamma)$ is an isomorphism in $\E^\iota\sSets^h$. This
can be checked after evaluating at each enhanced object $e \in \E$,
so at the end of the day, it suffices to check \thetag{iii} for $\E
= \sSets^h$. This immediately follows from \thetag{ii} by
adjunction.
\endproof

\begin{defn}\label{cofi.def}
An enhanced functor $\gamma:\C \to \C'$ between small enhanced
categories $\C$, $\C'$ is {\em final} if it satisfies the equivalent
condtions of Lemma~\ref{cofi.le}, and {\em cofinal} if its
enhanced-opposite $\gamma^\iota$ is final.
\end{defn}

\begin{exa}\label{refl.cof.exa}
A left resp.\ right-reflexive enhanced functor $\gamma:\C \to \C'$
between small $\C$, $\C'$ is cofinal resp.\ final (since the total
localization of a small enhanced category with an initial or
terminal enhanced object is trivial, Lemma~\ref{cofi.le}~\thetag{i}
immediately follows from Corollary~\ref{adj.comma.corr}).
\end{exa}

\begin{exa}\label{fib.cof.exa}
In particular, for any enhanced fibration $\gamma:\C \to \C'$
between small enhanced categories, and any enhanced object $c \in
\C'_\ppt$, the right-reflexive embedding $\eta:\C_c \subset c
\setminus^h_\gamma \C$ is final. Therefore one can compute right Kan
extensions by replacing the enhanced comma-fibers in
\eqref{enh.kan.eq} with the usual enhanced fibers. Dually, for an
enhanced cofibration and left Kan extensions, one can also use the
usual enhanced fibers.
\end{exa}

\begin{exa}
The composition of final resp.\ cofinal enhanced functors is final
resp.\ cofinal. For any enhanced fibrations $\C_0,\C_1 \to \C$
between small enhanced categories, an enhanced functor
$\gamma:\C_0 \to \C_1$ cartesian over $\C$ with cofinal enhanced
fibers $\gamma_c$, $c \in \C_{\ppt}$ is cofinal, and dually for
enhanced cofibrations and final functors (apply
Corollary~\ref{enh.bc.corr} and its dual).
\end{exa}

\begin{exa}\label{EX.exa}
Assume given an enhanced functor $\gamma:\E' \to \E$ between small
enhanced categories $\E$, $\E'$, and enhanced functors $X:\E^\iota
\to \sSets^h$, $X':{\E'}^\iota \to \sSets^h$, with the corresponding
enhanced families of groupoids $\pi':\E'X' \to \E'$, $\pi:\E X \to
\E$. Then a map $a:X' \to \gamma^{\iota*}X$ gives rise to a functor
$\alpha:\E'X' \to \E X$, and $\alpha$ is cofinal iff the adjoint map
$a_\dg:\gamma^\iota_!X' \to X$ is an isomorphism. Indeed, $a_\dg$ is
obtained by applying $\pi^\iota_!$ to the tautological map
$\alpha^\iota_!\ppt^h \to \ppt^h$, and $\pi^\iota_!$ is conservative
by Example~\ref{discr.pf.exa}.
\end{exa}

It is useful to also consider conditions stronger than finality or
cofinality. To state one such condition, note that for any enhanced
functor $\gamma:\C \to \C'$ between small enhanced categories, we
have the commutative square
\begin{equation}\label{enh.hyper.sq}
\begin{CD}
\C_{h\gamma^*{\Iso}} @>{\lambda}>> \Hh^\Tot(\C_{h\gamma^*\Iso})\\
@VVV @VVV\\
\C @>>> \C',
\end{CD}
\end{equation}
where $\C_{h\gamma^*\Iso}$ is as in \eqref{h.v.sq}, and $\lambda$
and $\Hh^\Tot$ are as in \eqref{loc.fun.eq}.

\begin{defn}\label{loc.enh.def}
An enhanced functor $\gamma:\C \to \C'$ between small enhanced
categories is a {\em localization} if $\gamma^*:\fFun^h(\C',\E) \to
\fFun^h(\C,\E)$ is fully faithful for any enhanced category $\E$,
and {\em hyperconnected} if the commutative square
\eqref{enh.hyper.sq} is enhanced-semicocartesian.
\end{defn}

An enhanced functor $\C \to \C'$ between small enhanced categories
whose target $\C'$ is an enhanced groupoid is hyperconnected iff it
is $h$-invertible in the sense of
Definition~\ref{enh.tot.loc.def}. In the general case, we have the
following.
  
\begin{lemma}\label{loc.enh.le}
Any localization $\gamma$ in the sense of
Definition~\ref{loc.enh.def} is both final and cofinal. If we have a
enhanced-semicocartesian square \eqref{enh.cocart.sq} of small
enhanced categories and enhanced functors such that the projection
$\C \to \C_0$ is a localization, then so is $\gamma_0:\C_1 \to
\C_{01}$, and the square is in fact enhanced-cocartesian. An
enhanced functor $\gamma:\C \to \C'$ between small enhanced
categories that is an enhanced fibration or an enhanced cofibration
is a localization iff all its enhanced fibers $\C_c$, $c \in
\C'_\ppt$ are hyperconnected.
\end{lemma}

\proof{} For the first claim, take $\E=\sSets^h$, note that if
$\ppt$ is the constant enhanced functor with value $\ppt$, then
$\gamma_!\ppt \cong \gamma_!\gamma^*\ppt \cong \ppt \cong
\gamma_*\gamma^*\ppt \cong \gamma_*\ppt$ since $\gamma^*$ is fully
faithful, and compute $\gamma_*\ppt$ resp.\ $\gamma_!\ppt$ by
\eqref{enh.kan.eq} resp.\ \eqref{enh.cokan.eq}. For the second
claim, apply Lemma~\ref{enh.cocart.le}, and note that fully faithful
embeddings are stable under semicartesian products, and the
semicartesian squares are in fact cartesian. For the third claim,
the ``only if'' part immediately follows from the first claim, since
the embeddings $\C_c \subset c \setminus^h \C$ resp.\ $\C_c \subset
\C /^h c$ are final resp.\ cofinal when $\gamma$ is an enhanced
fibration resp.\ enhanced cofibration. Conversely, in either of the
cases, for any $\E$, we have the relative functor category
$\fFun^h(\C|\C',\E)$ and the equivalence
\begin{equation}\label{fun.rel.equi}
  \fFun^h(\C,\E) \cong \sSec^h(\C,\fFun^h(\C|\C',\E)),
\end{equation}
and since $\sSec^h(\C,-)$ sends semicartesian squares to
semicartesian squares, we are done by Lemma~\ref{enh.cocart.le}.
\endproof

\begin{corr}\label{loc.hyper.corr}
Any enhanced functor hyperconnected in the sense of
Definition~\ref{loc.enh.def} is a localization, and the square
\eqref{enh.hyper.sq} is enhanced-cocartesian. For any
enhanced-semicocartesian square \eqref{enh.cocart.sq} of small
enhanced categories and enhanced functors such that the projection
$\C \to \C_0$ is a localization, so is $\gamma_0:\C_1 \to
\C_{01}$. An enhanced fibration or cofibration $\gamma:\C \to \C'$
between small enhanced categories is hyperconnected iff all its
enhanced fibers $\C_c$, $c \in \C'_\ppt$ are hyperconnected.
\end{corr}

\proof{} For the first claim, note that for any small enhanced
category $\C$ with total localization $\Hh^\Tot(\C)$, the enhanced
functor $\lambda:\C \to \Hh^\Tot(\C)$ is a localization by the third
claim of Lemma~\ref{loc.enh.le}, and apply its second claim to the
square \eqref{enh.hyper.sq}. Then for the second claim, it suffices
to consider the case when $\C_0 \cong \Hh^\Tot(\C)$, and the
enhanced functor $\C \to \C_0$ is the enhanced functor
$\lambda$. Then $\C \to \C_1$ factors through
$\C_1'=\C_{1h\gamma^*(\Iso)}$, any for any small enhanced category
$\E$, we have fully faithful enhanced functors
\begin{equation}\label{hyper.dia}
\begin{CD}
  \fFun^h(\C_{01},\E)\\ @V{\alpha}VV\\ \fFun^h(\C_1,E)
  \times^h_{\fFun^h(\C'_1,\E)} \fFun^h(\Hh^{\Tot}(\C_1'),\E)\\
  @V{\beta}VV\\ \fFun^h(\C_1,E) \times^h_{\fFun^h(\C,\E)}
  \fFun^h(\C_0,\E)\\ @VVV\\ \fFun^h(\C_1,\E)
\end{CD}
\end{equation}
such that $\beta \circ \alpha$ is an equivalence by
Lemma~\ref{enh.cocart.le}. Therefore $\alpha$ is also an
equivalence, and \eqref{enh.hyper.sq} for $\gamma$ is
enhanced-cocartesian by Lemma~\ref{enh.cocart.le}. For the third
claim, the ``only if'' part follows from the first claim and
Lemma~\ref{loc.enh.le}. For the ``if'' part, if $\gamma:\C \to \C'$
is an enhanced fibration resp.\ cofibration, then the left vertical
arrow in \eqref{enh.hyper.sq} is cartesian resp.\ cocartesian over
the right vertical arrow, and as in Lemma~\ref{loc.enh.le}, the
claim follows from Lemma~\ref{enh.cocart.le} and
\eqref{fun.rel.equi}.
\endproof

\begin{remark}
It is {\em not} true that all localizations are hyperconnected, see
Example~\ref{loc.P.exa}.
\end{remark}

\begin{exa}\label{loc.cof.exa}
By Lemma~\ref{loc.enh.le}, the total localization functor $\lambda$
of \eqref{loc.fun.eq} is both final and cofinal. For any small
enhanced category $\E$ with a universal object, the corresponding
enhanced functor $\Unf(I^o) \to \E$ of \eqref{enh.univ.sq} is a
localization, and in the situation of Lemma~\ref{cath.lim.le}, so is
the enhanced functor $\E \to \colimh_{\C^\iota}E$ of
\eqref{enh.colim.sq}. Moreover, as in Remark~\ref{cocat.rem}, since
$\iota:\cCat^h \to \cCat^h$ preserves colimits and sends
localizations to localizations, we also have a localization
$\E_{h\perp} \to \colimh_{\C^\iota}E$, where $\E_{h\perp} \to
\C^\iota$ is the transpose enhanced cofibration.
\end{exa}

\begin{exa}\label{2.cof.exa}
The $2$-connected functors of Subsection~\ref{V.I.subs} ---
specifically, those of Lemma~\ref{Z.le}, Lemma~\ref{Z.pl.le} and
Lemma~\ref{II.le} --- are also hyperconnected in the sense of
Definition~\ref{loc.enh.def}; the same proofs work.
\end{exa}

\subsection{Replacements and expansions.}\label{enh.repl.subs}

Looking at Example~\ref{2.cof.exa}, one may wonder whe\-ther other
$2$-connected functors constructed in Section~\ref{loc.sec} and
Section~\ref{seg.loc.sec} are hyperconnected in the sense of
Definition~\ref{loc.enh.def}. The most useful one is that of
Lemma~\ref{del.prop}. This indeed works, and to prove it, we need to
find an enhanced version of simplicial replacements of
Subsection~\ref{del.subs}. Consider again the category $\Delta$,
with its embedding $\phi:\Delta \to \Pos \to \Cat$, and upgrade it
to an enhanced embedding $\Phi:\Unf(\Delta) \to \Unf(\Pos) \cong
\pPos \to \cCat^h$, where the second embedding is
\eqref{ppos.eq}. Note that $\Phi$ then induces an enhanced functor
$\Phi_!\Y:\cCat^h \to \Delta^o+h\cCat^h$ that restricts to an
enhanced functor
\begin{equation}\label{enh.del}
  \rho \circ \Phi_!\Y:\cCat^h \to \Delta^o_h\sSets^h,
\end{equation}
where $\rho:\cCat^h \to \sSets^h$ is right-adjoint to the embedding
$\sSets^h \subset \cCat^h$. If we let $\delta:\Delta \to \Delta
\times \Delta$ be the diagonal embedding, then as in
\eqref{bi.compl.bis.sq}, for any $[n],[m] \in \Delta$ and enhanced
family of groupoids $\C \to \Unf(\Delta)$ corresponding to some
$X:\Delta^\iota_h \to \sSets^h$, we denote $\C_{[n],[m]} =
\Unf(\delta)^\iota_*X([n] \times [m]$.

\begin{defn}\label{enh.seg.def}
An {\em enhanced Segal family} is an enhanced family of groupoids
$\C \to \Unf(\Delta)$ that is semicartesian over each square
\eqref{segal.sq}, and an enhanced Segal family $\C$ is {\em
  complete} if for any $[n] \in \Delta$, the square
\begin{equation}\label{C.compl.sq}
\begin{CD}
\C_{[n],[0]}  @>>> \C_{[n],[0]} \times \C_{[n],[0]}\\
@VVV @VVV\\
\C_{[n],[3]} @>>> \C_{[n],[1]} \times \C_{[n],[1]}
\end{CD}
\end{equation}
induced by \eqref{bi.compl.sq} is also semicartesian.
\end{defn}

\begin{prop}\label{enh.del.prop}
The enhanced functor \eqref{enh.del} induces an equivalence between
$\cCat^h$ and the full enhanced subcategory in
$\Delta^o_h\sSets^h$ spanned by small complete enhanced
Segal families in the sense of Definition~\ref{enh.seg.def}.
\end{prop}

\begin{lemma}\label{enh.mod.loc.le}
For any complete cocomplete model category $\C$, with the enhanced
localization $\Hh^W(\C)$ provided by Lemma~\ref{enh.mod.le}, and for
any Reedy category $I$, the enhanced functor
\begin{equation}\label{enh.mod.loc}
\Hh^W(I^o\C) \to I^o_h\Hh^W(\C)
\end{equation}
induced by the projection $\Unf(I^o\C) \cong I^o_h\Unf(\C)
\to I^o_h\Hh^W(\C)$ is an equivalence.
\end{lemma}

\proof{} If $I=J \in \Posf^+$ is a left-bounded partially ordered
set, with the Reedy structure of Example~\ref{reedy.exa}, the claim
holds tautologically. In general, choose a universal object of
Proposition~\ref{univ.prop} for $\Unf(I^o)$, with the corresponding
$J \in \Posf^+$ and functor $a:J^o \to I^o$ such that $\Unf(a)$ is a
localization that fits into an enhanced-cocartesian square
\eqref{enh.univ.sq}. Then on one hand,
$a^{h*}:I^o_h\Hh^W(\C) \to J^o_h\Hh^W(\C) \cong
\Hh^W(J^o\C)$ is fully faithful, and on the other hand, so is
$a^*:I^o\C \to J^*\C$. Moreover, since the Reedy structure on $J$ is
directed, Lemma~\ref{J.proj.le} immediately shows that $a^*$ and
$a_*$ define a Quillen adjunction between $I^o\C$ and $J^o\C$, so by
Lemma~\ref{enh.mod.le}, $\Hh^W(a^*):\Hh^W(I^o\C) \to \Hh^W(J^o\C)$
is a right-admissible full embedding. Then so is
\eqref{enh.mod.loc}, and to finish the proof, it remains to check
that it is essentially surjective. However, its essential image
trivially contains representable functors $\Y(\Unf(I),i,c)$, $i \in
I$, $c \in \C$, and being right-admissible, \eqref{enh.mod.loc} is
enhanced-right-exact, so its essential image then contains any
$X:\Unf(I)^\iota \to \Hh^W(\C)$ by \eqref{Y.C.X}.
\endproof

\begin{remark}
In particular, Lemma~\ref{enh.mod.loc.le} shows that any enhanced
functor $\Unf(I)^\iota \to \Hh^W(\C)$ lifts to a functor $I^o \to
\C$. This is actually true for any small category $I$, and can be
shown by applying Lemma~\ref{enh.mod.loc.le} to the simplicial
replacement $\Delta I$, along the lines of \cite{DKHS}. We will not
need this.
\end{remark}

\proof[Proof of Proposition~\ref{enh.del.prop}.] Apply
Lemma~\ref{enh.mod.loc.le} to $\C = \Delta^o\Sets$ and $I=\Delta$,
and combine the resulting equivalence \eqref{enh.mod.loc} with
\eqref{ccat.h.eq}.
\endproof

Explicitly, the enhanced functor \eqref{enh.del} can be computed by
\eqref{enh.cokan.eq}: for any small enhanced category $\C$, we have
$\rho(\Phi_!\Y(\C)) \cong \Delta_h\C$, where the {\em enhanced
  simplicial expansion} $\Delta^{\Tot}_h\C$ and the {\em enhanced
  simplicial replacement} $\Delta_h\C \subset \Delta^{\Tot}_h\C$ are
defined by
\begin{equation}\label{enh.del.wr}
\begin{aligned}
\Delta^{\Tot}_h\C &= \fFun^h(\Delta_\idot | \Delta,\C) \cong
\Unf(\Delta) \bb^h \C,\\
\Delta_h\C &= \fFun^h(\Delta_\idot | \Delta,\C)_{h\flat} \cong
\Unf(\Delta) \bbi^h \C.
\end{aligned}
\end{equation}
Here as in \eqref{del.wr}, $\nu_\idot:\Delta_\idot \to \Delta$ is
the tautological cofibration \eqref{nu.eq}, with its right-adjoint
$\nu^\dg:\Delta \to \Delta_\idot$, and
$\fFun^h(\Delta_\idot|\Delta,\C)$ is the relative enhanced functor
category \eqref{rel.ffun.sq}. The evaluation pairing provides an
enhanced functor
\begin{equation}\label{enh.xi.tot.d}
 \xi^{\Tot}_\idot:\Delta^{\Tot}_{h\idot}\C =
 \Unf(\nu_\idot)^*\Delta^{\Tot}_h\C \to \C
\end{equation}
that restricts to an enhanced functor
\begin{equation}\label{enh.xi.d}
 \xi_\idot:\Delta_{h\idot}\C = \Unf(\nu_\idot)^*\Delta_h\C \to \C.
\end{equation}
Restricting further with respect to $\nu_\idot$, we obtain
enhanced functors
\begin{equation}\label{enh.xi}
\xi^\Tot = \xi^\Tot_\idot \circ \Unf(\nu^\dg):\Delta^{\Tot}_h\C \to
\C, \quad \xi = \xi_\idot \circ \Unf(\nu^\dg):\Delta_h\C \to \C,
\end{equation}
an enhanced version of \eqref{xi.eq}. Moreover, one can also
consider the fibration $\nu^\hdot:\Delta^\hdot \to \Delta^o$
transpose to $\nu$, and define the {\em enhanced cosimplicial
  expansion} $\Delta^\Tot_{h\perp}\C$ by
\begin{equation}\label{enh.codel.wr}
\Delta^\Tot_{h\perp}\C = \fFun^h(\Delta^\hdot|\Delta^o,\C) =
\fFun^h(\Delta_\idot|\Delta,\C)_{h\perp},
\end{equation}
where the right-hand side is as in \eqref{rel.fun.perp}. Then
$\Delta^\Tot_{h\perp}\C$ comes equipped with an enhanced functor
\begin{equation}\label{enh.xi.perp.d}
\xi^\idot:\Unf(\nu^\hdot)^*\Delta^\Tot_{h\perp}\C \to \C,
\end{equation}
again provided by evaluation.

\begin{lemma}\label{enh.del.loc.le}
For any small enhanced category $\C$, the enhanced functors
$\xi_\idot$ resp.\ $\xi$ of \eqref{enh.xi.d} resp.\ \eqref{enh.xi}
are hyperconnected in the sense of Definition~\ref{loc.enh.def}.
\end{lemma}

\proof{} For any enhanced family of groupoids $\E \to \Unf(\Delta)$,
define its {\em contraction} $\Con(\E) \in \Cat^h$ by the
enhanced-cocartesian square
\begin{equation}\label{enh.xi.d.sq}
\begin{CD}
\Unf(\overline{\nu}_\idot)^*\E @>{\lambda}>>
\Hh^{\Tot}(\Unf(\overline{\nu}_\idot)^*\E)\\
@VVV @VVV\\
\Unf(\nu_\idot)^*\E @>>> \Con(\E),
\end{CD}
\end{equation}
where $\Delta_{\idot\Tot} \cong [0] \setminus \Delta \subset
\Delta_\idot$ is the dense subcategory spanned by maps cocartesian
over $\Delta$, and $\overline{\nu}_\idot:\Delta_{\idot\Tot} \to
\Delta$ is induced by $\nu_\idot$. Note that since
$\Delta_{\idot\Tot}$ has an initial object $[0]$,
Lemma~\ref{enh.adm.cof.le} provides a right-admissible full
embedding $\E_{[0]} \to \Unf(\overline{\nu}_\idot)^*\Delta_h\E$, so
we can identify $\Hh{\Tot}(\Unf(\overline{\nu}_\idot)^*\E)$ in
\eqref{enh.xi.d.sq} with the enhanced fiber $\E_{[0]}$. Then
\eqref{enh.xi.d.sq} is functorial with respect to $\E$, and
\eqref{enh.del.wr} together with the adjunction of
Lemma~\ref{enh.fun.le} provide a functorial identification
\begin{equation}\label{eexp.adj}
\fFun^h_{\Delta}(\E,\Delta_h\C) \cong \fFun^h(\Con(\E),\C)
\end{equation}
for any $\C \in \Cat^h$, so that $\Con:\Delta^\iota_h\sSets^h \to
\cCat^h$ is a left-adjoint enhanced functor to
\eqref{enh.del}. Since \eqref{enh.del} is fully faithful by
Proposition~\ref{enh.del.prop}, we have $\C \cong \Con(\Delta_h\C)$
for any small enhanced category $\C$, and \eqref{enh.xi.d} then
appears as the bottom arrow in the corresponding square
\eqref{enh.xi.d.sq}, so it is hyperconnected by
Corollary~\ref{loc.hyper.corr}. For \eqref{enh.xi}, recall that the
right-adjoint $\nu^\dg$ to $\nu_\idot$ itself has a right-adjoint
$\nu_>:\Delta_\idot \to \Delta$ such that $\xi_\idot \cong
\nu_>^*(\xi)$; then \eqref{enh.xi.d.sq} for $\E = \Delta_h\C$ has an
right-admissible full subsquare that looks like
\begin{equation}\label{enh.xi.sq}
\begin{CD}
\Delta_{+h}\C @>{\lambda}>> \Hh^{\Tot}(\Delta_{h+}\C)\\
@VVV @VVV\\
\Delta_h\C @>{\xi}>> \C,
\end{CD}
\end{equation}
where $\Delta_{+h}\C$ is the restriction of $\Delta_h\C$ to
$\Delta_+ \subset \Delta$. Being right-admissible,
\eqref{enh.xi.sq} is itself enhanced-cocartesian by
Corollary~\ref{enh.adm.sq.corr}~\thetag{ii}.
\endproof

As an application of the cocartesian squares \eqref{enh.xi.d.sq} of
Lemma~\ref{enh.del.loc.le}, we can give a characterization of the
transpose enhanced cofibrations of Subsection~\ref{enh.groth.subs}
that does not use the enhanced Grothendieck construction. Namely,
assume given an enhanced fibration $\gamma:\C \to \E$ between small
enhanced categories, with the enhanced functor
$\xi^\Tot:\Delta_h^\Tot\E \to \E$ of \eqref{enh.xi}, and note that
the adjunction map $a:\nu^\dg \circ \nu_\idot \to \id$ induces an
enhanced functor $\alpha:\Delta_h^\Tot\C \to \xi^{\Tot *}\C$ over
$\Unf(\Delta)$. Then by Lemma~\ref{enh.fib.adj.le} and
Lemma~\ref{enh.2adm.le}, $\alpha$ admits a fully faithful
right-adjoint
\begin{equation}\label{al.dg}
\alpha^\dg:\xi^{\Tot *}\C \to \Delta^\Tot_h\C
\end{equation}
cartesian over $\Unf(\Delta)$. Explicitly, enhanced objects in the
enhanced simplicial expansion $\Delta^\Tot_h\C$ are triples $\langle
[n],\phi,s \rangle$, $[n] \in \Delta$, $\phi:\Unf([n]) \to \E$ an
enhanced functor, $s \in \Sec^h([n],\phi^*\C)$ a section, and then
\eqref{al.dg} identifies $\xi^{\Tot *}\C$ with the enhanced full
subcategory in $\Delta^\Tot_h\C$ spanned by triples such that $s$ is
cartesian. Since \eqref{al.dg} is cartesian over $\Unf(\Delta)$, we
have the transpose enhanced functor
$\alpha^\dg_{h\perp}:(\xi^{\Tot *}\C)_{h\perp} \to \Delta^\Tot_{h\perp}\C$,
and then restricting it to $(\xi^*\C)_{h\perp}$ and composing with
with the evaluation functor
\eqref{enh.xi.perp.d}, we obtain an enhanced functor
$\Unf(\nu^\hdot)^*(\xi^*\C)_{h\perp} \to \C$ that
fits into a commutative square
\begin{equation}\label{xi.CE.sq}
\begin{CD} 
\Unf(\overline{\nu}^\hdot)^*(\xi^*\C)_{h\perp} @>>> \C_{h\gamma^*(\Iso)}\\
@VVV @VVV\\  
\Unf(\nu^\hdot)^*(\xi^*\C)_{h\perp} @>>> \C,
\end{CD}
\end{equation}
where $\Delta^\hdot_\flat \subset \Delta^\hdot$ is the dense
subcategory spanned by maps cartesian over $\Delta^o$,
$\overline{\nu}^\hdot:\Delta^\hdot_\flat \to \Delta^o$ is induced by
$\nu^\hdot$, and $\C_{h\gamma^*(\Iso)} = \C \times^h_{\E}
\E_{h\Iso}$ is as in \eqref{h.v.sq}.

\begin{lemma}\label{enh.perp.le}
For any enhanced fibration $\C \to \E$ between small enhanced
categories, the square \eqref{xi.CE.sq} is enhanced-semicocartesian.
\end{lemma}

\proof{} If $\C \to \E$ is an enhanced family of groupoids, then the
claim reduces to Lemma~\ref{enh.del.loc.le} applied to
$\C^\iota$. More generally, since $\Cat^h$ is cartesian-closed, the
cartesian product of an enhanced-semicocartesian square with a fixed
small enhanced category $\C'$ is enhanced-semicocartesian. Therefore
the claim also holds when $\C = \C_0 \times^h \C_1$ for some
enhanced family of groupoids $\C_0 \to \E$ and small enhanced
category $\C_1$. In the general case, it is easy to see that both
squares are actually functorial with respect to $\C$ --- that is,
define enhanced functors $\E^\iota\cCat^h \to [1]^2_h\cCat^h$ ---
and we need to show that these functors factor through the full
enhanced subcategory spanned by enhanced-semicocartesian
squares. Since this subcategory is closed under enhanced colimits,
and both functors are obviously enhanced-right-exact, this
immediately follows from \eqref{Y.cat.eq}.
\endproof

\begin{corr}\label{enh.perp.corr}
For any small enhanced category $\E$ and enhanced fibration
$\gamma:\C \to \E$, with transpose enhanced cofibration $\gamma':\C'
= \C_{h\perp} \to \E^\iota$, there exist enhanced-semicocartesian
squares
\begin{equation}\label{perp.sq}
\begin{CD}
\Unf(\iota \circ \overline{\nu}_\idot)^*\xi^*\C @>>>
\C_{h\gamma^*(\Iso)}\\
@VVV @VVV\\
\Unf(\iota \circ \nu_\idot)^*\xi^*\C @>>> \C',
\end{CD}
\qquad
\begin{CD}
\Unf(\overline{\nu}_\idot)^*\xi^{\iota *}\C' @>>>
\C'_{h{\gamma'}^*(\Iso)}\\
@VVV @VVV\\
\Unf(\nu_\idot)^*\xi^{\iota *}\C' @>>> \C,
\end{CD}
\end{equation}
where as in \eqref{xi.CE.sq}, $\C_{h\gamma^*(\Iso)} = \C \times^h_{\E}
\E_{h\Iso} \cong \C' \times^h_{\E^\iota} \E^\iota_{h\Iso} =
\C'_{h{\gamma'}(\Iso)}$.
\end{corr}

\proof{} The second square is \eqref{xi.CE.sq} on the nose, and the
first one is enhanced-opposite to the square \eqref{xi.CE.sq} for the
transpose-opposite enhanced fibration $\C_{h\perp}^\iota \to \E$.
\endproof

For a slightly different interpretation of \eqref{enh.xi.d.sq},
note that Proposition~\ref{enh.del.prop} in particular means that
$\Unf(\Delta) \subset \cCat^h$, with the embedding $\Phi$, is a
generating subcategory in the sense of
Definition~\ref{enh.gen.def}. Then by \eqref{E.coli.eq}, we have
\begin{equation}\label{C.coli.eq}
\C \cong \colimh_{\Delta_h\C}\Phi \circ \pi,
\end{equation}
where $\pi:\Delta_h\C \to \Unf(\Delta)$ is the projection, and
\eqref{enh.xi.d.sq} is the corresponding square
\eqref{enh.colim.sq}. One can also consider larger generating
subcategories --- for example, the category $\Posff$ of all finite
partially ordered sets. Then if one denotes by $\Psi:\Unf(\Posff)
\to \cCat^h$ the embedding functor, we have an enhanced functor
\begin{equation}\label{enh.pos}
  \rho \circ \Psi_!\Y:\cCat^h \to \Posff^o_h\sSets^h
\end{equation}
that reduces to \eqref{enh.del} after further restricting with
respect to the embedding $\eps:\Delta \subset \Posff$.

\begin{defn}\label{Psi.def}
A cone $E_>$ of a functor $E:J \to \Posff$ from a finite partially
ordered set $J$ is {\em $\Psi$-universal} if $\Psi \circ
\Unf(E_>):\Unf(J)^{h>} \to \cCat^h$ is universal. An enhanced
functor $X:\Unf(\Posff)^\iota \to \sSets^h$ is {\em $\Psi$-exact} if
the enhanced-opposite $X^\iota:\Unf(\Posff) \to \sSets^{h\iota}$
sends $\Psi$-universal cones to universal cones.
\end{defn}

\begin{lemma}\label{enh.del.le}
For any $\Psi$-exact enhanced functor $X:\Unf(\Posff)^\iota \to
\sSets$, $\eps^{h*}\Unf\Posff)X \to \Unf(\Delta)$ is a
complete enhanced Segal family in the sense of
Definition~\ref{enh.seg.def}, and the map $X \to
\Unf(\eps^o)_*\Unf(\eps^o)^*X$ is an isomorphism.
\end{lemma}

\proof{} Since $\Psi$ commutes with finite coproducts, any exact
enhanced functor $X:\Unf(\Posff)^\iota \to \sSets^h$ is
additive. Moreover, say that a commutative square $[1]^2 \cong \V^>
\to \Posff$ is $\Psi$-cocartesian if it is $\Psi$-universal in the
sense of Definition~\ref{Psi.def}. Then all squares \eqref{segal.sq}
are $\Psi$-universal, and so is both the square \eqref{compl.sq} and
its product with any $[n] \in \Delta$. Therefore
$\eps^{h*}\Unf(\Posff)X \to \Unf(\Delta)$ is indeed a complete
enhanced Segal family for any $\Psi$-exact $X$. Moreover, for any $J
\in \Posff$, we have the perfect colimit \eqref{BJ.coli.eq}, and by
Remark~\ref{BJ.N.rem}, it is preserved by the nerve functor
$N$. Therefore firstly, it is $\Psi$-universal, and secondly, by
\eqref{enh.kan.eq}, it is preserved by
$(\Unf(\eps)_*Y)^\iota:\Unf(\Posff) \to \sSets^{h\idot}$ for any
$Y:\Delta^\iota_h \to \sSets$. Thus the map $X \to
\Unf(\eps^o)_*\Unf(\eps^o)^*X$, being an equivalence over
$\Unf(\Delta) \subset \Unf(\Posff)$, is an equivalence over the
whole $\Unf(\Posff)$.
\endproof

\begin{corr}\label{enh.del.corr}
The enhanced functor \eqref{enh.pos} induces an equivalence between
$\cCat^h$ and the full subcategory in $\Posff^o_h\sSets^h$
spanned by exact enhanced functors.
\end{corr}

\proof{} Lemma~\ref{enh.del.prop} implies that
$\eps^{oh*}:\Posff^o_h\sSets^h \to
\Delta^\iota_h\sSets^h$ induces a fully faithful embedding from
the category of $\Psi$-exact enhanced functors to the category of
complete enhanced Segal families. Since the equivalence of
Proposition~\ref{enh.del.prop} factors through this full embedding,
it must be essentially surjective.
\endproof

\subsection{Adding colimits.}\label{env.subs}

Let us now consider the case when an enhanced category $\E$ has some
limits or colimits but not all of them. To formalize the situation,
we introduce the following.

\begin{defn}\label{ka.co.def}
For any full enhanced subcategory $\I \subset \cCat^h$, an enhanced
category $\E$ is {\em $\I$-complete} if the limit $\limh_{\C} E$
exists for any $\C \in \I$ and enhanced functor $E:\C \to \E$, and
an enhanced functor $F:\E \to \E'$ between two $\I$-complete
enhanced categories is {\em $\I$-left-exact} if it preserves all
those limits. An enhanced category $\E$ is {\em $\I$-cocomplete} if
$\E^\iota$ is $\I$-complete, and an enhanced functor $F$ is {\em
  $\I$-right-exact} if $F^\iota$ is $\I$-left-exact. An enhanced
object $e \in \E_\ppt$ in an $\I$-cocomplete enhanced category $\E$
is {\em $\I$-projective} if the enhanced functor $\Y(\E,e)$ of
\eqref{yo.e.eq} is $\I$-right-exact.
\end{defn}

\begin{exa}\label{empty.exa}
Consider the full enhanced subcategory $\{\emptyset\} \subset
\cCat^h$ span\-ned by the empty category $\emptyset$. Then an
enhanced category is $\{\emptyset\}$-cocomplete iff it has an
initial enhanced object $o$, and an enhanced functor $\gamma$ is
$\{\emptyset\}$-right-exact if it preserves it.
\end{exa}

\begin{exa}
For any regular cardinal $\kappa$, we have the full subcategory
$\Cat^h_{\kappa} \subset \Cat^h$ of \eqref{cat.h.kappa.eq}, with the
induced enhancement $\cCat^h_\kappa$; to simplify terminology, we
say that an enhanced category is {\em $\kappa$-complete} resp.\ {\em
  $\kappa$-cocomplete} if it is $\Cat^h_\kappa$-complete
resp.\ $\cCat^h_\kappa$-cocomplete in the sense of
Definition~\ref{ka.co.def}, and an enhanced functor is {\em
  $\kappa$-left-exact} resp.\ {\em $\kappa$-right-exact} if it is
$\cCat^h_\kappa$-left-exact resp.\ $\cCat^h_\kappa$-right-exact
\end{exa}

\begin{exa}
For any regular cardinal $\kappa$, Lemma~\ref{cath.lim.le}
immediately shows that $\cCat^h_\kappa$ is $\kappa$-cocomplete.
\end{exa}

\begin{exa}
In the particular case when $\kappa$ is the countable cardinal,
$\kappa$-complete resp.\ $\kappa$-cocomplete is also called
``finitely complete'' resp.\ ``finitely cocomplete''. This is
consistent with our earlier terminology; in particular, the same
argument as in Proposition~\ref{enh.compl.prop} shows that an
enhanced category $\E$ is finitely complete iff the family
$\E|_{\Posf^\pm} \to \Posf^\pm$ is finitely complete in the sense of
Definition~\ref{fin.compl.def}.
\end{exa}

\begin{exa}
We have the full embedding $\sSets = \Unf(\Sets) \subset \Cat^h$
sending a set to the corresponding discrete small enhanced category,
with the induced embedding $\sSets_\kappa \subset \Cat^h_\kappa$ of
sets $S$ with $|S| < \kappa$, for any regular $\kappa$. Then an
enhanced category {\em has products} resp.\ {\em coproducts} if it
is $\sSets$-complete resp.\ $\sSets$-cocomplete, and has {\em
  $\kappa$-bounded products} resp.\ {\em coproducts} if it is
$\sSets_\kappa$-complete resp.\ cocomplete. If $\kappa$ is the
countable cardinal, $\kappa$-products resp.\ $\kappa$-coproducts are
finite products resp.\ coproducts.
\end{exa}

\begin{lemma}\label{I.fibb.le}
Assume given an enhanced family of groupoids $\pi:\C' \to \C$ such
that $\C$ is $I$-cocomplete, and for any enhanced functor $E:I \to
\C$ from some $I \in \I$, with the universal cone $E_>:I^{h>} \to
\C$, the functor $I^{h>} \to \sSets^{h\iota}$ corresponding to
$E_>^*\C' \to I^{h>}$ is a universal cone. Then $\C'$ is
$\I$-cocomplete, and $\pi$ is $\I$-right-exact.
\end{lemma}

\proof{} The condition actually insures that for any $I \in \I$ and
$E:I \to \C'$, the projection $\cCone(E) \to \cCone(\pi \circ E)$
induced by $\pi$ is an equivalence.
\endproof

\begin{exa}\label{cat.E.coc}
Lemma~\ref{I.fibb.le} immediately shows that for any small enhanced
category $\E$, the enhanced category $\cCat^h \bbi^h \E$ is
cocomplete, and the enhanced family of groupoids $\cCat^h \bbi^h \E
\to \cCat^h$ induced by \eqref{no.E} is enhanced-right-exact. For
any regular cardinal $\kappa$, the enhanced category $\cCat^h_\kappa
\bbi^h \E$ is $\kappa$-cocomplete, and $\cCat^h_\kappa \bbi^h \E \to
\cCat^h_\kappa$ is $\kappa$-right-exact.
\end{exa}

\begin{lemma}\label{I.fib.le}
Assume given an enhanced cofibration $\pi:\C' \to \C$ between small
enhanced categories, and assume that $\C$ is $\I$-complete for some
full subcategory $\I \subset \cCat^h$, the enhanced fiber $\C'_c$ is
$\I$-cocomplete for any enhanced object $c \in \C_\ppt$, and for any
enhanced morphism $f:c \to c'$ in $\C_\ppt$, the transition functor
$f_!:\C'_c \to \C'_{c'}$ is $\I$-right-exact. Then $\C'$ is
$\I$-cocomplete, and $\pi$ is $\I$-right-exact.
\end{lemma}

\proof{} For any enhanced functor $E:\E \to \C'$ from a small $\E$,
postcomposition with $\pi$ provides an enhanced functor
$\pi(E):\cCone(\E) \to \cCone(\pi \circ E)$. Since $\pi$ is an
enhanced cofibration, $\pi(E)$ is also an enhanced cofibration by
Example~\ref{ffun.fib.exa}. To see what its enhanced fibers are,
recall that $\E^{h>}$ is the enhanced cylinder of the tautological
projection $\E \to \ppt^h$, thus $[1]$-coaugmented. Then for any
object in $\cCone(\pi \circ E)$ represenseted by an enhanced functor
$\gamma = (\pi \circ E)_>:\E^{h>} \to \C$, $\gamma^*\C'$ is also
$[1]$-coaugmented, so $\gamma^*\C' \cong \Cyl_h(\nu(\gamma)')$ for
some enhanced functor $\nu(\gamma)':(\pi \circ E)^*\C' \to \C'_c$,
where $c = \gamma(o)$ is the vertex of the cone $\gamma$. Then
composing $\nu(\gamma)'$ with the taulological section $\E \to (\pi
\circ E)^*\C'$ give an enhanced functor $\nu(\gamma):\E \to \C'_c$,
and Lemma~\ref{enh.cyl.o.le} together with Lemma~\ref{enh.cocart.le}
provide an equivalence
\begin{equation}\label{cone.seq}
\cCone(E)_\gamma \cong \cCone(\nu(\gamma)).
\end{equation}
If $\E \in \I$, then all these fibers have initial enhanced
objects preserved by the transition functors, so by
Lemma~\ref{enh.fib.adj.le} applied to the transpose enhanced
fibrations, $\pi(E)$ admits a left-adjoint section
$\cCone(\pi \circ E) \to \cCone(E)$. But $\cCone(\pi \circ E)$ has
an initial enhanced object.
\endproof

Note that different full enhanced subcategories $\I \subset \cCat^h$
can lead to the same class of $\I$-complete enhanced categories and
$\I$-right-exact enhanced functors -- at the very least, both do not
change if we enlarge $\I$ by adding the targets $\C'$ of all
enhanced functors $\gamma:\C \to \C'$, $\C \in \I$ that are final in
the sense of Definition~\ref{cofi.def}. The same goes for
$\I$-cocomplete categories and cofinal functors. Again, let us
formalize the situation (we concentrate on colimits).

\begin{defn}\label{colim.cl.def}
A full enhanced subcategory $\I \subset \I'$ in a full enhanced
subcategory $\I' \subset \cCat^h$ is {\em $\colim$-dense} if any
$\I$-cocomplete enhanced category is $\I'$-cocomplete, and any
$\I$-right-exact enhanced functor is $\I'$-right-exact. The {\em
  $\colim$-closure} of a full enhanced subcategory $\I \subset
\cCat^h$ is the maximal full enhanced subcategory $\I' \subset
\cCat^h$ containing $\I$ such that $\I \subset \I'$ is
$\colim$-dense, and a full enhanced subcategory $\I \subset \cCat^h$
is {\em $\colim$-closed} if it coincides with its $\colim$-closure.
\end{defn}

\begin{lemma}\label{colim.cl.le}
Let $\I \subset \cCat^h$ be a $\colim$-closed full enhanced
subcategory. Then \thetag{i} for any cofinal enhanced functor
$\gamma:\C \to \C'$ in $\cCat^h$ with source $\C \in \I$, the target
$\C'$ is also in $\I$, \thetag{ii} $\ppt^h \in \I$, \thetag{iii}
for any enhanced cofibrations $\gamma:\C' \to \C$ such that $\C \in
\I$ and $\C'_c \in \I$ for any $c \in \C_\ppt$, $\C' \in \I$, and
\thetag{iv} $\I$ is $\I$-cocomplete.
\end{lemma}

\proof{} The first claim is tautological. For \thetag{ii}, note that
$\colimh_{\ppt}E \cong E$ by Lemma~\ref{enh.kan.com.le}, and it is
trivially preserved by any enhanced functor. For \thetag{iii}, use
\eqref{enh.cokan.eq} in the form of Example~\ref{fib.cof.exa} to
shows that the Kan extensions $\gamma_!$ exist in any
$\I$-cocomplete enhanced category and are preserved by
$\I$-right-exact functors, and then identify $\colimh_{\C'}E \cong
\colimh_{\C}\gamma_!E$. For \thetag{iv}, note that for any enhanced
functor $C:\C \to \I \subset \cCat^h$, with the corresponding
enhanced cofibration $\C' \to \C$, we have $\C' \in \I$ by
\thetag{iii}, and then $\colimh_{\C}C \in \I$ by \thetag{i} and
Example~\ref{loc.cof.exa}.
\endproof

Of course, by Corollary~\ref{enh.ff.corr}, a full enhanced
subcategory $\I \subset \cCat^h$ is completely determined by the
full subcategory $\I_\ppt \subset \Cat^h$. We use the whole $\I$
in Definition~\ref{colim.cl.def} to allow for
Lemma~\ref{colim.cl.le}~\thetag{iv}, but we will use the same
terminology as in Definition~\ref{colim.cl.def} for full
subcategories $\I_\ppt \subset \Cat^h$. Here are some very standard
examples of $\colim$-dense subcategories in $\Cat$ that are also
$\colim$-dense in $\Cat^h$.

\begin{lemma}\label{coli.kr.le}
The subcategory $\Posf^+ \subset \Cat^h$ is $\colim$-dense, and so is
the subcategory $\Posf^+_\kappa \subset \Cat^h_\kappa$ for any
regular cardinal $\kappa$. Moreover, denote by $\rho:\Sets \to \Posf^+$
the embedding sending a set $S$ to the partially ordered set $S^<$;
then $\rho(\Sets) \cup \{\emptyset\} \subset \Posf^+$,
$\rho(\Sets_\kappa) \cup \{\emptyset\} \subset \Posf^+_\kappa$ are
also $\colim$-dense. Finally, let $\{\Unf(\II)\} \in \Cat^h$ be the
full subcategory spanned by the Kronecker category $\II$; then
$\Sets \copr \{\Unf(\II)\} \subset \Cat^h$, $\Sets_\kappa \copr
\{\Unf(\II)\} \subset \Cat^h_\kappa$ are also $\colim$-dense.
\end{lemma}

In words, Lemma~\ref{coli.kr.le} says that en enhanced category has
colimits iff it has the initial enhanced object and cofibered
coproducts iff it has coproducts and coequalizers, and similarly for
$\kappa$-bounded colimits.

\proof{} The first claim immediately follows from
Proposition~\ref{univ.prop} and Example~\ref{loc.cof.exa}. For the
second claim, let $I \subset \Cat^h$ be the $\colim$-closure of
$\rho(\Sets) \cup \{\emptyset\}$; we need to show that $I$ contains
$\Posf^+$. But then by Lemma~\ref{colim.cl.le}~\thetag{ii} and
Example~\ref{refl.cof.exa}, for any map $J' \to J$ in $\Posf^+$ such
that $j \in I$ and $J' / j \in I$ for all $j \in J$, $J'$ is also in
$I$. Moreover, by Lemma~\ref{colim.cl.le}~\thetag{iii}, $I \cap
\Posf^+$ is closed under colimits of cofibrant functors $S^< \to \I$
--- indeed, by Lemma~\ref{sec.sq.le}, these are preserved by the
embedding $\Unf(\Posf^+) \to \cCat^h$. Then the standard induction
on skeleta reduces us to showing that $J^> \in I$ for any $J \in
\Posf$, and this is again Example~\ref{refl.cof.exa}. Finally, for
the third claim, let $I \subset \Cat^h$ be the $\colim$-closure of
$\Sets \cup \{\II\}$, and note that it is closed under coproducts
and coequalizers by Lemma~\ref{colim.cl.le}~\thetag{iii}, and again,
contains $\Unf(J^>)$ for any $J \in \Posf$. Then for any $S \in
\Sets_\kappa$, consider $S^\hush \in \Posf_\kappa$ of \eqref{B.S},
with the corresponding projection $S^\hush \to \V$, and compose it
with the map $\V \to \II$, $o \mapsto 0$, $0,1 \mapsto
1$. Then by Lemma~\ref{cath.lim.le}, we have $\Unf(S^\hush) \cong
\colimh_{i \in \II}(\Unf(S^\hush / i))$, and the comma-fibers are
$\Unf(S^\hush / 0) \cong S \times \ppt^h$, $\Unf(S^\hush / 1)
\cong (S \times \Unf([1])) \copr \Unf(S^>)$. Thus all these
comma-fibers are in $I$, so that $\Unf(S^\hush) \in I$, and since
the excision map $\Unf(S^\hush) \to \Unf(S^<)$ is cofinal,
$\Unf(S^<) \in I$ by Lemma~\ref{colim.cl.le}~\thetag{i}. Thus we are
done by the second claim.
\endproof

\begin{defn}\label{env.def}
For any full enhanced subcategory $\I \subset \cCat^h$, the {\em
  $\I$-envelope} of an enhanced category $\E$ is an $\I$-cocomplete
enhanced category $\Env(\E,\I)$ with a fully faithful enhanced
functor $\Y(\E,\I):\E \to \Env(\E,\I)$ such that for any
$\I$-cocomplete enhanced category $\E'$, any enhanced functor $E:\E
\to \E'$ factors as
\begin{equation}\label{env.facto}
\begin{CD}
\E @>{\Y(\E,\I)}>> \Env(\E,I) @>{\Env(E,\I)}>> \E',
\end{CD}
\end{equation}
with an $\I$-right-exact enhanced functor $\Env(E,\I)$, and such a
factorization \eqref{env.facto} is unique up to a unique
isomorphism.
\end{defn}

\begin{lemma}\label{env.le}
For any full enhanced subcategory $\I \subset \cCat^h$, any enhanced
category $\E$ admits an $\I$-envelope $\Y(\I):\E \to \Env(\E,\I)$.
\end{lemma}

\proof{} Replacing $\I$ by its $\colim$-closure, we may assume that
$\I$ is $\colim$-closed. Assume first that $\E$ is small, and let
$\Env(\E,\I) \subset \E^\iota\sSets^h$ be the full enhanced
subcategory spanned by enhanced families of groupoids $\E' \to \E$
such that $\E' \in \I$. Then by Lemma~\ref{yo.im.le},
Lemma~\ref{colim.cl.le}~\thetag{ii} and Example~\ref{refl.cof.exa},
$\Env(\E,\I)$ contains the essential image of the Yoneda embedding
$\Y(\E)$ of \eqref{yo.yo.enh.eq}, so it factors through a fully
faithful embedding $\Y(\E,\I):\E \to \Env(\E,\I)$. By
Lemma~\ref{colim.cl.le}~\thetag{iv} and Example~\ref{groth.l.e.exa},
$\Env(\E,\I)$ is $\I$-cocomplete, and then for any enhanced functor
$\E \to \E'$ to an $\I$-cocomplete $\E'$, the left Kan extension
$\Y(\E,\I)_!E$ exists by \eqref{enh.cokan.eq} and
\eqref{Y.comma.eq}, and a canonical decomposition \eqref{env.facto}
is provided by \eqref{enh.kan.facto}.

In the general case, represent $\E \cong \colim^2_I \E_i$ as the
filtered $2$-colimit of Lemma~\ref{enh.filt.ff.le} of its small full
enhanced subcategories $\E_i \subset \E$, and take $\Env(\E,\I) =
\colim^2_I \Env(\E_i,\I)$.
\endproof

\begin{corr}\label{Y.I.dg.corr}
An enhanced category $\C$ is $\I$-cocomplete if and only if
$\Y(\I):\C \to \Env(\C,\I)$ is left-reflexive.
\end{corr}

\proof{} Clear. \endproof

\begin{corr}\label{colim.cl.corr}
A full enhanced subcategory $\I \subset \cCat^h$ is $\colim$-closed
if and only if it satisfies the conditions \thetag{i}-\thetag{iii}
of Lemmca~\ref{colim.cl.le}.
\end{corr}

\proof{} The ``only if'' part is Lemma~\ref{colim.cl.le}. For the
``if'' part, note that the construction of the envelope
$\Env(\E,\I)$ given in Lemma~\ref{env.le} only uses the fact that
$\I$ satisfies the conditions \thetag{i}, \thetag{ii} of
Lemma~\ref{colim.cl.le}, together with \thetag{iv} that is a formal
corollary of \thetag{i} and \thetag{iii}. Therefore for any small
enhanced category $\E$, $\Env(\E,\I) \subset \E^\iota\sSets^h \cong
\sSets^h \mmh \E$ consists of families of enhanced groupoids $\E'
\to \E$ such that $\E' \in \I$. Then if $\E \in \cCat^h$ lies in the
$\colim$-closure $\I' \supset \I$ of $\I$, $\Env(\E,\I')$ contains
the constant enhanced family of groupoids $\E \to \E$, and $\E \in
\I$ since $\Env(\E,\I')=\Env(\E,\I)$.
\endproof

\begin{lemma}\label{env.adj.le}
For any small enhanced categories $\C$, $\C'$ and adjoint pair of
enhanced functors $\lambda:\C \to \C'$, $\rho:\C' \to \C$, and any
enhanced full subcategory $\I \subset \cCat^h$, $\Env(\lambda,\I)$
is adjoint to $\Env(\rho,\I)$.
\end{lemma}

\proof{} Immediately follows from the uniqueness of the
$\I$-right-exact extensions $\Env(\lambda,\I)$ and $\Env(\rho,\I)$.
\endproof

\begin{exa}\label{empty.env.exa}
Take $\I = \{\emptyset\}$, as in Example~\ref{empty.exa}. Then for
any enhanced category $\E$, $\Env(\E,\I) = \E^{h>}$ is the category
of Example~\ref{enh.gt.exa}.
\end{exa}

\begin{exa}
Take $\I = \cCat^h$ and $\E = \ppt^h$. Then $\Env(\E,\I)=\sSets^h$,
and Lemma~\ref{env.le} says that for any enhanced object $e$ in a
cocomplete enhanced category $\E$, there exists a unique right-exact
enhanced functor $\sSets^h \to \E$ sending $\ppt^h$ to $e$.
\end{exa}

\begin{exa}\label{env.cat.exa}
More generally, for $\I = \cCat^h$ and any small enhanced category
$\E$, we have $\Env(\E,\I) \cong \E^\iota\sSets^h$, and $\Y(\I) =
\Y(\E)$ is the Yoneda embedding \eqref{yo.yo.enh.eq}. For any
enhanced functor $\gamma:\E_0 \to \E_1$ between small enhanced
categories, $\Env(\gamma,\I):\Env(\E_0,\I) \to \Env(\E_1,\I)$ is the
left Kan extension functor $\gamma_!^o$, and in particular, it has a
right-adjoint enhanced functor $\Env(\gamma,\I)^\dg$ given by
$\gamma^{o*}$. If $\E$ is not small, the envelope $\Env(\E,\I)$ is
still well-defined, and its enhanced objects correspond to enhanced
functors $E:\E^o \to \sSets^h$, but not all such functors ---
indeed, all functors do not form a well-defined enhanced category.
Instead, we only take functors of the form $\eps^o_!E$ for some
small full enhanced subcategory $\E_0 \subset \E$, with the
embedding $\eps:\E_0 \to \E$, and some enhanced functor
$E:\E_0^\iota \to \Sets$. For any enhanced functor $\gamma:\E_0 \to
\E_1$ such that $\E_0$ is small, the enhanced functor
$\Env(\gamma,\I)$ still admits a right-adjoint
$\Env(\gamma,\I)^\dg$; for any $\E_0' \subset \E_1$ containing the
essential image of $\gamma$, with the embedding functor $\eps:\E'_0
\to \E_1$, and any $E:{\E'}_0^\iota \to \sSets$, we have
$\Env(\gamma,\I)^\dg(\eps^o_!E) \cong \gamma^{o*}\E$.
\end{exa}

\begin{lemma}\label{ff.env.le}
Assume given a $\colim$-closed full subcategory $\I \subset \cCat^h$
and an $\I$-cocomplete enhanced category $\C$, and let $\C_c \subset
\C$ be a full subcategory spanned by $\I$-projective enhanced
objects. Then the $\I$-right-exact extension $\eps':\Env(\C_c,\I)
\to \C$ of the fully faithful embedding $\eps:\C_c \to \C$ is fully
faithful.
\end{lemma}

\proof{} It suffices to consider the case when $\C_c$ is small. Let
$\I'=\cCat^h$, and note that we have a commutative diagram
\begin{equation}\label{env.dia}
\begin{CD}
\C_c @>{\Y(\I)}>> \Env(\C_c,\I) @>{\nu_c}>> \Env(\C_c,\I')\\
@V{\eps}VV @VV{\Env(\eps,\I)}V @VV{\Env(\eps,\I')}V\\
\C @>{\Y(\I)}>> \Env(\C,\I) @>{\nu}>> \Env(\C,\I'),
\end{CD}
\end{equation}
where $\nu$ is a fully faithful embedding such that $\nu \circ
\Y(\I) \cong \Y(\I')$, and similarly for $\nu_c$. Since $\C$ is
$\I$-cocomplete, $\Y(\I)$ in the bottom row has a left-adjoint
enhanced functor $\Y(\I)^\dg:\Env(\C,\I) \to \C$ provided by
Corollary~\ref{Y.I.dg.corr}, and we have $\eps' \cong \Y(\I)^\dg
\circ \End(\eps,\I)$. On the other hand, by
Example~\ref{env.cat.exa}, $\Env(\C_c,\I')$ in the top row is
identified with $\C_c^\iota\sSets^h$, and $\nu \circ \Y(\I)$ is then
the Yoneda embedding $\Y:\C_c \to \C_c^\iota\sSets$, while
$\Env(\eps,\I')$ has a right-adjoint enhanced functor
$\Env(\eps,\I')^\dg$. Then our two adjunctions provide a map
\begin{equation}\label{env.adj.eq}
\begin{aligned}
\nu_c &\to \Env(\eps,\I')^\dg \circ
\Env(\eps,\I') \circ \nu_c \cong
\Env(\eps,\I')^\dg \circ \nu \circ \Env(\eps,\I)\\
&\to \Env(\eps,\I')^\dg \circ \nu \circ \Y(\I) \circ \Y(\I)^\dg \circ
\Env(\eps,\I) \cong \eps^\dg \circ \eps',
\end{aligned}
\end{equation}
where we denote $\eps^\dg = \Env(\eps,\I')^\dg \circ \Y(\I')$, and
the map \eqref{env.adj.eq} defines an adjunction between $\eps^\dg$
and $\eps'$ over $\nu_c$ in the sense of
Definition~\ref{adj.rel.def}. Since $\nu_c$ is fully faithful,
\eqref{adj.rel.iso} reduces us to proving that \eqref{env.adj.eq} is
an isomorphism. However, since $\C_c \subset \C$ is spanned by
$\I$-projective enhanced objects, $\eps^\dg$ is $\I$-right-exact,
and then so is $\eps^\dg \circ \eps'$. Then by the uniqueness part
of Definition~\ref{env.def}, it suffices to check that
\eqref{env.adj.eq} becomes an isomorphism after restricting to $\C_c
\subset \Env(\C_c,\I)$, and this holds since $\eps$ is fully
faithful.
\endproof

\section{Large categories.}\label{enh.large.sec}

\subsection{Karoubi envelopes.}\label{enh.kar.sus}

We note that if an enhanced category $\E$ is $\I$-cocomplete in the
sense of Definition~\ref{ka.co.def}, for some $\I \subset \cCat^h$,
then the universal property of Definition~\ref{env.def} provides an
$\I$-right-exact enhanced functor $\Env(\E,\I) \to \E$ left-adjoint
to the embedding $\Y(\I)$, but $\Y(\I)$ itself is usually not an
equivalence and not $\I$-right-exact. For instance, in
Example~\ref{empty.env.exa}, even if $\E$ already has a terminal
enhanced object, passing to $\E^{h>}$ just ignores this fact and
formally adds a new one. As a consequence of this, the envelope
construction is not idempotent, with one important exception: a
subcategory $\pP \subset \cCat^h$ that controls idempotents.

\begin{defn}\label{enh.kar.def}
An enhanced category $\E$ is {\em enhanced-Karoubi-complete} if it
is $\pP$-cocomplete for the subcategory $\pP = \{\Unf(P)\} \subset
\cCat^h$, where $P$ is as in Example~\ref{P.exa}. The {\em enhanced
  Karoubi envelope} of an enhanced category $\E$ is its envelope
$\Env(\E,\pP)$.
\end{defn}

\begin{lemma}\label{kar.le}
Assume given an enhanced category $\E$.
\begin{enumerate}
\item A cone $E_>:\Unf(P^>) \to \E$ of an enhanced functor
  $E:\Unf(P) \to \E$ is universal if and only if it further extends
  to an enhanced functor $E_=:\Unf(P^=) \to \E$, where $P^= \supset
  P^>$ is as in Example~\ref{kar.lim.exa}.
\item Any limits or colimits of an enhanced functor $E:\Unf(P) \to
  \E$ are preserved by any enhanced functor $F:\E \to \E'$.
\item An enhanced category $\E$ is enhanced-Karoubi-complete iff so
  is the enhanced-opposite category $\E^\iota$, and this happens if
  and only if the embedding $\Y(\pP):\E \to \Env(\E,\pP)$ is an
  equivalence.
\item For any enhanced category $\E$, $\E_\ppt \subset
  \Env(\E,\pP)_\ppt$ is Karoubi-dense, so that if $\E$ is small,
  then so is $\Env(\E,\pP)$.
\end{enumerate}
\end{lemma}

\proof{} In \thetag{i}, note that $P \cong P^o$ is self-opposite,
and so is $P^=$, so by passing to the enhanced-opposite categories
and functors, we reduce the statement to a dual one for dual
cones. Then by Lemma~\ref{lim.univ.le}, it further suffices to
consider the case when $\E$ is complete, and in this case, it
suffices to compute the right Kan extension $\Unf(\eps)_*E$ with
respect to the embedding $\eps:P \to P^=$ by \eqref{enh.kan.eq}, and
use the uniqueness of retracts. Then \thetag{ii} immediately follows
from \thetag{i}, and since $P^=$ is self-opposite, so does the first
part of \thetag{iii}. For the second part, the ``if'' part is
tautological, and if $\E$ is enhanced-Karoubi-complete, then by
\thetag{ii}, any functor $\E \to \E'$ to an
enhanced-Karoubi-complete $\E'$ factors uniquely through $y =
\Y(\pP)$. Thus we have an enhanced functor $y':\Env(\E,\pP) \to \E$
such that $y' \circ y \cong \id$, and then by uniqueness, $y \circ
y' \cong \id$ as well. Finally, for \thetag{iv}, let $\Env(\E,\pP)'
\subset \Env(\E,\pP)$ be the full enhanced subcategory corresponding
to the Karoubi closure of $\E_\ppt$ in $\Env(\E,\pP)_\ppt$; then
$\Env(\E,\pP)'$ is enhanced-Karoubi-complete by \thetag{i}, so
again, $\Env(\E,\pP)'=\Env(\E,\pP)$ by \thetag{ii} and the unique
extension property of envelopes.
\endproof

\begin{exa}\label{loc.P.exa}
The taulotogical projection $\nu:\Unf(P) \to \ppt^h$ of course
factors through the embedding $\Unf(P) \to \Unf(P^=)$, and then
Lemma~\ref{kar.le}~\thetag{i} and \eqref{enh.kan.eq} immediately
show that for any enhanced category $\E$ and enhanced functor
$E:\ppt^h \to \E$, the Kan extension $\nu_*\nu^*E$ exists, and the
adjunction map $E \to \nu_*\nu^*E$ is an isomorphism. Therefore
$\nu$ is a localization in the sense of
Definition~\ref{loc.enh.def}, so that $\Unf(P)$ is
hyperconnected. In particular, this means that any enhanced groupoid
is Karoubi-closed. By Lemma~\ref{kar.le}~\thetag{ii} and
Lemma~\ref{I.fib.le}, this implies that any enhanced family of
groupoids $\E \to \C$ over a Karoubi-closed base $\C$ is itself
Karoubi-closed. In fact, the embedding $\Unf(P) \to \Unf(P^=)$
itself is a localization, but it is not hyperconnected (a fully
faithful enhanced functor is hyperconnected only if it is an
equivalence).
\end{exa}

\begin{exa}\label{noenh.nokar.exa}
Lemma~\ref{kar.le}~\thetag{iv} does {\em not} mean that
$\Env(\E,\pP)_\ppt$ is the unenhanced Karoubi closure of $\E_\ppt$;
in fact, $\Env(\E,\pP)_\ppt$ need not be Karoubi-complete. This is
because not all functors $P \to \E_\ppt$ lift to enhanced functors
from $\Unf(P)$ to $\E$. This happens already for $\E = \sSets^h$,
and in fact, already for the full subcategory in $\sSets^h$ spanned
by the usual unenhanced groupoids. A counterexample is a group $G$
equipped with an endomorphism $p:G \to G$ such that $p^2$ is
conjugate to $p$, but $p$ is not conjugate to a projector.
\end{exa}

As Example~\ref{noenh.nokar.exa} shows, constructing enhanced
functors $\Unf(P) \to \E$ to some enhanced category $\E$ is
non-trivial. In fact, the difficulty is essential, since $\Unf(P)$
does not admits a finite universal object --- thus lifting a
projector in $\E_\ppt$ to an enhanced projector in $\E$ requires
providing an infinite amount of additional data. However, the
situation is much better for projectors that have an image. Namely,
consider again the category $P^=$ of Example~\ref{kar.lim.exa}, with
its objects $\ppt,o \in P^+$, and note that the comma-fiber $P^= /
\ppt$ can be identified with the ordinal $[2]$ --- the object $a:o
\to \ppt$ corresponds to $1 \in [2]$, and $0,2 \in [2]$ are $f:\ppt
\to \ppt$ for $f=\id$ and $f=p$. Then the fibration $\sigma:[2]
\cong P^= / \ppt \to P^=$ fits into a cartesian square
\begin{equation}\label{P.sq}
\begin{CD}
[1] @>{m}>> [2]\\
@V{e}VV @VV{\sigma}V\\
\ppt @>{\eps(o)}>> P^=,
\end{CD}
\end{equation}
where $m:[1] \to [2]$ is the embedding onto $\{0,2\} \subset [2]$,
and $\eps(o):\ppt \to P^=$ is the embedding onto $o \in P^=$.

\begin{lemma}\label{P.eq.le}
The commutative square \eqref{P.sq} induces an enhanced-co\-cart\-esian
square of small enhanced categories.
\end{lemma}

\proof{} The fibers of the fibration $\sigma:[2] \to P^=$ are given
by $([2])_\ppt \cong \ppt$ and $([2])_o \cong [1]$. Therefore by
\eqref{enh.kan.eq}, for any enhanced category $\E$, the Kan
extension $\Unf(\sigma)_*\Unf E'$ exists for any enhanced functor
$E':\Unf([2]) \to \E$, and $E \cong \Unf(\sigma)_*\Unf(\sigma)^*E$
for any enhanced functor $E:\Unf(P^=) \to \E$, so that
$\Unf(\sigma)$ is a localization. Moreover, an enhanced functor
$E':\Unf([2]) \to \E$ lies in the essential image of the fully
faithful embedding $\Unf(\sigma)^*$ iff the adjunction map
$\Unf(\sigma)^*\Unf(\sigma)_*E' \to E'$ is an isomorphism, this
always holds over $\ppt \in P^=$, and over $o \in P^=$, this holds
iff $\Unf(m)^*E'$ lies in the essential image of $\Unf(e)^*$ by
virtue of \eqref{enh.bc.eq}. We are then done by
Lemma~\ref{enh.cocart.le}.
\endproof

By virtue of Lemma~\ref{P.eq.le}, for any enhanced category $\E$,
any retract of an object $e \in \E_\ppt$ gives rise to an enhanced
functor $\Unf(P^=) \to \E$, uniquely up to an isomorphism. Here is
one curious application of this fact.

\begin{lemma}\label{P.o.le}
Assume given a Karoubi-closed small enhanced category $\C$ such that
the identity functor $\id:\C \to \C$ admits a cone. Then $\C$ has a
terminal enhanced object.
\end{lemma}

\proof{} Let $\alpha:\C \to \C^{h>}$ be the embedding, and let
$\beta:\C^{h>} \to \C$ be the cone of $\id:\C \to \C$, so that we
have an isomorphism $\id \cong \beta \circ \alpha$.  Consider the
enhanced functor $\Unf(P^=) \to \cCat^h$ provided by
Lemma~\ref{P.eq.le}, sending $o$ resp.\ $\ppt$ to $\C$
resp.\ $\C^{h>}$ and $a$, $b$ to $\alpha$, $\beta$, and let
$\pi^=:\C^+ \to \Unf(P^=)$ be the corresponding enhanced cofibration
with fibers $\C$ and $\C^{h>}$. Note that $\C^{h>}$ is
$[1]$-coaugmented, with fibers $\C$, $\ppt^h$, so it is
Karoubi-closed by Lemma~\ref{I.fib.le}, and then so is $\C^=$, by
the same Lemma~\ref{I.fib.le}. Now let $\eps:P \to P^=$ be the
embedding. Then the only enhanced fiber of the induced cofibration
$\pi:\Unf(\eps)^*\C^= \to \Unf(P)$ is $\C^{h>}$, and it has an
enhanced terminal object, so by Lemma~\ref{enh.fib.adj.le}, $\pi$
admits a right-adjoint enhanced section $s:\Unf(P) \to
\Unf(\eps)^*\C^=$. Since the enhanced category $\C^=$ is
Karoubi-closed, $\Unf(\eps) \circ s:\Unf(P) \to \Unf(\eps)^*\C^= \to
\C^=$ extends to a section $s^=:\Unf(P^=) \to \C^=$ of the
cofibration $\pi^=$ by Lemma~\ref{kar.le}~\thetag{i}. But since $s$
is the terminal enhanced object in
$\sSec^h(\Unf(P),\Unf(\eps)^*\C^=)$, its extension $s^=$ is the
terminal enhanced object in $\sSec^h(\Unf(P^=),\C^=)$ by the same
Lemma~\ref{kar.le}~\thetag{i}, so that the section $s^=$ is
right-adjoint to the projection $\pi^=$. Then $s^=(o)$ is the
terminal enhanced object in $\C = \C^+_o$.
\endproof

\begin{corr}\label{I.proj.corr}
For any Karoubi-closed enhanced category $\C$ and full subcategory
$\I \subset \Cat^h$, an enhanced object in the envelope
$\Env(\C,\I)$ of Lemma~\ref{env.le} is $\I$-projective iff it lies
in the image of the full embedding $\Y(\I):\C \to \Env(\C,\I)$.
\end{corr}

\proof{} The ``if'' part is clear: for any $E = \Y(\I)(c)$, $c \in
\C_\ppt$, the enhanced functor $\Y(\C^\iota\sSets^h,E):\Env(\C,\I)
\subset \C^\iota\sSets^h \to \sSets^h$ is evaluation at
$c$. Conversely, assume given some $\I$-projective $E \in
\Env(\C,\I) \subset \C^\iota\sSets^h$, with the corresponding
enhanced family of groupoids $\pi:\E \to \C$, and replace $\I$ with
its $\colim$-closure so that $\E \in \I$ by Lemma~\ref{env.le}. Then
$E \cong \pi_!\ppt^h$ by Example~\ref{discr.pf.exa}, and since
$\ppt^h \cong \colimh_{\E}\Y(\E)$ tautologically, we have $E \cong
\colimh_{\E}\pi_! \circ \Y(\E)$. Since $\E \in \I$,
$\Y(\C^\iota\sSets^h,E)$ preserves this colimit. Enhanced objects in
the enhanced groupoid $\Y(\C^\iota\sSets^h,E)(E)$ are enhanced
morphisms $E \to E$, so by Lemma~\ref{cath.lim.le}, the identity
morphism $\id:E \to E$ factors through a map $f:E \cong \pi_!\ppt^h
\to \pi_!\Y(\E,e)$ for some $e \in \E_\ppt$, and
then by \eqref{discr.sq}, $f \cong \pi_!f'$ for an enhanced map
$f':\ppt^h \to \Y(\E,e)$. But then $\ppt^h$ is the vertex of a
universal cone for $\Y(\E):\E \to \E^\iota\sSets^h$, and $f'$ gives
rise to another cone for $\Y(\E)$, with vertex $\Y(\E,e)$. This cone
then takes values inside the image of the full embedding $\Y(\E)$,
so effectively, we obtain a cone $\E^{h>} \to \E$ of the identity
functor $\id:\E \to \E$. But since $\C$ is Karoubi-closed, $\E$ is
Karoubi-closed by Example~\ref{loc.P.exa}, and then $E$ lies in the
essential image of $\Y(\E)$ by Lemma~\ref{P.o.le} and
Lemma~~\ref{yo.im.le}.
\endproof

Lemma~\ref{kar.le}~\thetag{iv} shows that the Karoubi envelope
$\Env(-,\pP)$ provides a well-defined functor $\Cat^h \to
\Cat^h$. This functor is not especially well-behaved; in particular,
it does not always preserve semicartesian squares (for a trivial
example already in the unenhanced setting, note that $\ppt
\times_{P^=} P \cong \emptyset$, where $\ppt \to P^=$ is the
embedding $\eps(o)$). However, this does not happen when one of the
enhanced functors in the square is an enhanced fibration or
cofibration.

\begin{lemma}\label{kar.C.E.le}
Assume given a semicartesian square \eqref{C.E.sq} such that $\pi$
is an enhanced fibration or cofibration in the sense of
Definition~\ref{enh.fib.def}. Then the corresponding commutative
square of Karoubi envelopes $\Env(-,\pP)$ is semicartesian.
\end{lemma}

\proof{} Since $\Env(\C,\pP)^\iota \cong \Env(\C^\iota,\pP)$ for any
$\C$, it suffices to consider the case of enhanced
fibrations. Furthermore, since a functor is an epivalence iff it is
an epivalence over any small subcategory in its target, we may
assume that $\E'$ and $\E$ are small, and then by
Lemma~\ref{fib.uni.le}, we can assume the same for $\C$, hence also
for $\C'$. Then \eqref{C.E.sq} induces a fully faithful embedding
\begin{equation}\label{kar.C.E}
\Env(\C',\pP) \to \Env(\E',\pP) \times^h_{\Env(\E,\pP)} \Env(\C,\pP),
\end{equation}
and it suffices to check that it is essentially
surjective. Explicitly, an enhanced object in the target of the
enhanced functor \eqref{kar.C.E} is represented by a triple $\langle
e,c,\alpha \rangle$ of enhanced functors $e:\Unf(P) \to \E'$,
$c:\Unf(P) \to \C$ and an isomorphism $\alpha:\colimh_P(\phi \circ
e) \to \limh_P (\pi \circ c)$. Such an isomorphism lifts to a map
$\beta:\phi \circ e \to \pi \circ c$ such that
$\colimh_P(\beta)=\id$, and then since $\fFun^h(\Unf(P),-)$ sends
enhanced fibrations to enhanced fibrations and semicartesian
products to semicartesian products, $\langle e,\beta^*c,\id \rangle$
lifts to an enhanced object $c' \in \fFun^h(\Unf(P),\C)'_\ppt$, thus
an enhanced functor $c':\Unf(P) \to \C'$ such that $\colimh_P c' \in
\Env(\C',\pP)_\ppt$ is sent to $\langle c,e,\alpha\rangle$ under
\eqref{kar.C.E}.
\endproof

\subsection{Filtered colimits.}

While there are plenty of full subcategories in $\Cat^h$ that are
$\colim$-dense in the sense of Definition~\ref{colim.cl.def}, they
are usually obtained as the $\colim$-closure of some explicit
subcategory, and do not admit an explicit description themselves
(note that even $\Cat^h_\kappa \subset \Cat^h$ is not
$\colim$-closed since by Example~\ref{refl.cof.exa}, it does not
satisfy Lemma~\ref{colim.cl.le}~\thetag{i}). Remarkably, there is one
exception, or in fact, one exception for any regular cardinal
$\kappa$.

\begin{defn}
For any regular cardinal $\kappa$, an enhanced category $\C$ is {\em
  $\kappa$-filtered} if any enhanced functor $E:\E \to \C$ from some
$\E \in \Cat^h_\kappa$ admits a cone. An enhanced category $\C$ is
{\em filtered} if it is $\kappa$-filtered for the smallest regular
cardinal $\kappa$.
\end{defn}

\begin{exa}
For any category $I$, an enhanced functor $\E \to \Unf(I)$ factors
through the truncation $\Unf(\E_\ppt)$ by
Proposition~\ref{trunc.prop}, so $\Unf(I)$ is filtered iff so is
$I$. For a partially ordered set $I$, $\Unf(I)$ is $\kappa$-filtered
iff $I$ is $\kappa$-directed in the sense of
Definition~\ref{pos.filt.def}.
\end{exa}

\begin{exa}
Let $P$ be the universal idempotent category of Example~\ref{P.exa}.
Then $\id:P \to P$ admits a cone $P^> \to P$. Thus {\em any}
enhanced functor $\E \to \Unf(P)$ admits a cone, and $\Unf(P)$ is
$\kappa$-filtered for any $\kappa$.
\end{exa}

\begin{lemma}\label{filt.fib.le}
Assume given an enhanced cofibration $\pi:\C' \to \C$ such that for some
regular cardinal $\kappa$, $\C$ is $\kappa$-filtered, and the
enhanced fiber $\C'_c$ is $\kappa$-filtered for any enhanced object
$c \in \C_\ppt$. Then $\C'$ is $\kappa$-filtered.
\end{lemma}

\proof{} Use \eqref{cone.seq} and exactly the same argument as in
Lemma~\ref{I.fib.le} --- the only difference is that we need any
enhanced object in the cone category, not necessarily an initial
one, so we do not need to require that the transition functors are
right-exact.
\endproof

\begin{lemma}\label{filt.sec.le}
An enhanced category $\C$ is $\kappa$-filtered iff for any $I \in
\Posf^+_\kappa$, any enhanced functor $\Unf(I) \to \C$ admits a
cone. A $\kappa$-filtered category is non-empty. For any small
$\kappa$-filtered enhanced category $\C$ and enhanced functor $E:\E
\to \C$, $\E \in \Cat^h_\kappa$, the enhanced category $\cCone(E)$
is $\kappa$-filtered.
\end{lemma}

\proof{} For the first claim, note that if an enhanced functor
$\gamma:\E' \to \E$ is a localization in the sense of
Example~\ref{loc.cof.exa}, then $\gamma^*:\cCone(E) \to
\cCone(\gamma^*E)$ is an equivalence for any $E:\E \to \C$, and
apply Proposition~\ref{univ.prop}. The second claim is obvious (take
$\E = \emptyset$). For the third claim, the product of two
$\kappa$-bounded enhanced categories is $\kappa$-bounded; then to
obtain a cone for an enhanced functor $E':\E' \to \cCone(E)$, it
suffices to take a cone for the corresponding enhanced functor $\E'
\times^h \E^{h>} \to \C$.
\endproof

\begin{lemma}\label{filt.hyper.le}
Assume given a small $\kappa$-filtered enhanced category $\C$ and an
enhanced functor $X:\C \to \sSets$, with the corresponding enhanced
cofibration $\pi:\C^\perp X \to \C$. Then $\C' = \C^\perp X$ is
$\kappa$-filtered if and only if it is hyperconnected (or
equivalently, $\colimh_\C X \cong \ppt^h$).
\end{lemma}

\proof{} For the ``only if'' part, we may replace $\C'$ with
$\C$. The the total localization functor $\lambda:\C \to
\Hh^\Tot(\C)$ of \eqref{loc.fun.eq} is an enhanced cofibration by
Example~\ref{enh.grp.fib.exa}, and for any $J \in \Posf_\kappa$ and
$\gamma:\Unf(J) \to \Hh^\Tot(\C)$, $\gamma^*\C$ is a $J$-coaugmented
enhanced category with $\kappa$-filtered fibers. So it has a section
by Lemma~\ref{filt.sec.le}, and then since $\C$ is
$\kappa$-filtered, this section has a cone that projects to a cone
of $\gamma$. Thus in particular, for any finite $J \in \Posff
\subset \Posf_\kappa$ with the projection $p:J \to \ppt$, the
transition map $p^*:\pi_0(\Hh^\Tot(\C)_\ppt) \to
\pi_0(\Hh^\Tot(\C)_J)$ is surjective, and then by additivity,
$\Hh^\Tot(\C_J)$ is connected for any such $J$, so that
$\Hh^\Tot(\C) \cong \ppt^h$ by
Corollary~\ref{fin.enh.equi.corr}.

For the ``if'' part, assume given an enhanced functor $E:\E \to \C'$
from some $\E \in \Cat^h_\kappa$. Then if $\C'$ is hyperconnected,
Lemma~\ref{kappa.E.le} provides a square \eqref{hyper.sq}, and if we
let $\wt{\E}$ be the $\V$-coaugmented category obtained by removing
$\C$ from this square, then $\wt{\E}$ is still $\kappa$-bounded. We
then have the enhanced functor $\wt{\gamma}:\wt{\E} \to \C'$, and
the composition $\pi \circ \wt{\gamma}:\wt{\E} \to \C$ has a cone
with some vertex $c$. This can be restricted to a cone $\gamma=(\pi
\circ E)_>$ for $\pi \circ E$, and since $\C'_c$ is an enhanced
groupoid, the corresponding functor $\nu(\gamma):\E \to \C'_c$ of
\eqref{cone.seq} factors through the total localization
$\Hh^\Tot(\E)$. But since $\eta$ in \eqref{hyper.sq} is
$h$-invertible, while $\beta$ is $h$-trivial, this means that
$\nu(\gamma):\E \to \C'_c$ factors through $\ppt^h$. Therefore it
admits a cone, and by \eqref{cone.seq}, so does $E$.  \endproof

\begin{exa}
Take $\C = \ppt^h$. Then Lemma~\ref{filt.hyper.le} says that a small
enhanced groupoid $\C$ is $\kappa$-filtered iff $\C \cong \ppt^h$.
\end{exa}

\begin{corr}\label{cof.filt.corr}
An enhanced functor $\gamma:\C \to \C'$ between $\kappa$-filtered
small enhanced categories is cofinal if and only if the comma-fiber
$c \setminus^h_\gamma \C$ is $\kappa$-filtered for any enhanced
object $c \in \C'_\ppt$.
\end{corr}

\proof{} For any $c \in \C_\ppt$, we have $c \setminus^h \C \cong
\C^\perp \gamma^*\Y(\C',c)$, where the enhanced functor
$\Y(\C',c):\C' \to \sSets$ is as in \eqref{yo.e.eq}.
\endproof

Next, consider the following situation. Assume given a partially
ordered set $J \in \Posf^+$ and $J$-coaugmented small enhanced
categories $\C$, $\E$. For any enhanced functor $E:\C \to \E$ over
$\Unf(J)$, we can consider relative enhanced cones $E_>:\Cyl(\gamma)
\to \E$ over $\Unf(J)$ with respect to the coaugmentation enhanced
functor $\gamma:\C \to \Unf(J)$, and these form a small enhanced
category $\cCone(E,J) = E \setminus^h_{\gamma^*} \sSec(J,\E)$. For
any $J' \in \Posf_\kappa$ equipped with a map $f:J' \to J$, we then
have the functor
\begin{equation}\label{cone.rel.f.eq}
\Unf(f)^*:\cCone(E,J) \to \cCone(\Unf(f)^*E,J')
\end{equation}
given by restriction with respect to $\Unf(f):\Unf(J') \to \Unf(J)$.

\begin{lemma}\label{sec.filt.le}
Assume given a left-$\kappa$-bounded $J \in \Posf^+$ and
$J$-coaug\-mented small enhanced categories $\C$, $\E$ such that for
any $j \in J$, $\C_j$ is $\kappa$-bounded and $\E_j$ is
$\kappa$-filtered. Then any enhanced functor $E:\C \to \E$ over
$\Unf(J)$ admits a relative enhanced cone $E_>:\Cyl(\gamma) \to \E$
over $\Unf(J)$ with respect to the coaugmentation $\gamma:\C \to
\Unf(J)$, and for any left-closed embedding $f:J' \to J$, the
functor \eqref{cone.rel.f.eq} has non-empty enhanced fibers.
\end{lemma}

\proof{} The second claim trivially implies the first one (take $J'
= \emptyset$). Moreover, for the second claim, Lemma~\ref{sec.sq.le}
and the standard induction on skeleta reduce us to the case when
$f$ is the embedding $f:J \to J^>$ for some $J \in
\Posf_\kappa$. Then $\Cyl_h(\gamma)$ is $[1]^2$-coaugmented, with
the second coaugmentation $\Cyl_h(\gamma) \to \Unf([1])$ induced by
the cofibration $J^> \to [1]$, and if we identify $[1]^2 \cong \V^>$
via an embedding $\eps:\V \to [1]^2$ and let $\Cyl_h(\gamma)_+ =
\eps^*\Cyl_h(\gamma)$, $J_+ = \eps^*([1] \times J^>)$, then
$\Cyl_h(\gamma) \cong \Cyl_h(\gamma)_+^{h>}$ and $J^> \times [1] =
J_+^>$. We are then given an enhanced functor $\Cyl_h(\gamma)_+ \to
\Unf([1]) \times^h \E$ over $\Unf(J_+)$, and we need to show that it
admits a cone over $\Unf(J_+^>)$. Since $\Cyl_h(\gamma)_+$ is
$\kappa$-bounded, this immediately follows from \eqref{cone.seq}.
\endproof

\begin{corr}\label{sec.filt.corr}
For any $J \in \Posf_\kappa$ and $J$-coaugmented small enhanced
category $\E$ such that $\E_j$ is $\kappa$-filtered for any $j \in
J$, the enhanced category $\sSec^h(J,\E)$ is $\kappa$-filtered, and
for any map $f:J' \to J$ im $\Posf_\kappa$, the restriction functor
$f^*:\sSec^h(J,\E) \to \sSec(J',\E)$ is cofinal.
\end{corr}

\proof{} The first claim is clear. For the second claim, for any
enhanced section $e:\Unf(J') \to \E$, let $J'' = \Cyl(f)$ and $\E''
= \Cyl_h(e)$. Then $\E''$ is $J''$-co\-augmented, satisfies the
assumptions, and $e \setminus^h_{f^*} \sSec^h(J) \cong
\sSec^h(J'',\E'')$, so we are done by Corollary~\ref{cof.filt.corr}.
\endproof

\begin{exa}\label{comma.filt.exa}
If $J=[1]$, then Corollary~\ref{sec.filt.corr} says that for any
enhanced functor $\C_0 \to \C_1$ between $\kappa$-filtered small
enhanced categories, the comma-category $\C_0 /^h \C_1$ is
$\kappa$-filtered.
\end{exa}

\begin{corr}\label{J.filt.corr}
For any $J \in \Posf_\kappa$ and $\kappa$-filtered small enhanced
category $\C$, the enhanced functor $\gamma^*:\C \to J^o_h\C$
induced by $\gamma:J^o \to \ppt$ is cofinal.
\end{corr}

\proof{} The enhanced category $J^o_h\C \cong \sSec(J^o,\Unf(J^o)
\times^h \C)$ is $\kappa$-filtered by Corollary~\ref{sec.filt.corr},
and then the comma-fibers $E \setminus^h_{\gamma^*} \C \cong
\cCone(E)$ are $\kappa$-filtered by Lemma~\ref{filt.sec.le}, so that
$\gamma^*$ is cofinal by Corollary~\ref{cof.filt.corr}.
\endproof

\begin{corr}\label{prod.filt.corr}
Assume given a filtered small enhanced category $\C$ and enhanced
functors $X_0,X_1:\C \to \cCat^h$, with the corresponding enhanced
cofibrations $\C_l = \C^\perp X_l \to \C$, $l=0,1$. Then if $\C_0$
and $\C_1$ are non-empty, $\C_0 \times^h_\C \C_1$ is non-empty as
well.
\end{corr}

\proof{} We have $\C_0 \times^h_\C \C_1 \cong \C^\perp \delta^*(X_0
\times^h X_1)$, where $\delta:\C \to \C \times^h \C$ is the diagonal
embedding. If $\C_0$, $\C_1$ are not empty, then $\colimh_{\C \times^h
  \C} X_0 \times^h X_1$ is not empty either, and $\delta$ is cofinal
by Corollary~\ref{J.filt.corr} (for $J = \{0,1\}$).
\endproof
  
\begin{lemma}\label{ff.filt.le}
Let $\gamma:\C' \to \C$ be a fully faithful enhanced functor such
that $\C$ is a small $\kappa$-filtered enhanced category, and the
comma-fiber $c \setminus^h_\gamma \C'$ is non-empty for any $c \in
\C_\ppt$. Then $\C'$ is $\kappa$-filtered, and $\gamma$ is cofinal.
\end{lemma}

\proof{} For the first claim, any enhanced functor $E:\E \to \C'$
from some $\E \in \Cat^h_\kappa$ admits a cone $E_>:\E \to \C$, with
some vertex $c \in \C$. Then since $c \setminus^h_\gamma \C'$ is
non-empty, we can find $c' \in \C'$ equipped with an enhanced
morphism $f:c \to \gamma(c')$, and since the projection $\cCone(E)
\to \C$ sending a cone to its vertex is a cofibration, we can
replace $E_>$ with $f_!(E_>)$ and insure that it factors through
$\C' \subset \C$. For the second claim, the functor $c
\setminus^h \C' \to c \setminus^h \C$ induced by $\gamma$ is also
fully faithful and has non-empty right comma-fibers, and then $c
\setminus^h \C \cong \cCone(\eps^h(c))$ is $\kappa$-filtered by
Lemma~\ref{filt.sec.le}, so $c \setminus^h \C'$ is $\kappa$-filtered
by the first claim. We are done by
Corollary~\ref{cof.filt.corr}.
\endproof

\begin{corr}\label{comma.filt.corr}
For any cofinal enhanced functor $\gamma:\C' \to \C$ between
$\kappa$-filtered small enhanced categories, the embedding $\eta:\C'
\to \C' /^h_\gamma \C$ is cofinal.
\end{corr}

\proof{} The tagret of the fully faithful embedding $\eta$ is
$\kappa$-filtered by Example~\ref{comma.filt.exa}, so by
Lemma~\ref{ff.filt.le}, it suffices to check that $\eta$ has
non-empty enhanced right comma-fibers. Take some enhanced object
$\wt{c}$ in $\C /^h_\gamma \C$, with projections $c' =
\sigma(\wt{c}) \in \C'_\ppt$, $c = \tau(\wt{c}) \in \C$, and note
that we have a semicartesian square
\begin{equation}\label{co.fi.sq}
\begin{CD}
\wt{c} \setminus^h_\eta \C' @>>> c' \setminus^h \C'\\
@VVV @VV{\rho}V\\
c \setminus^h_\gamma \C' @>{\lambda}>> \gamma(c') \setminus^h_\gamma
\C',
\end{CD}
\end{equation}
where $\lambda$ is induced by the enhanced morphism $\gamma(c') \to
c$, and $\rho$ is induced by the enhanced functor $\gamma$. Now, $c'
\setminus^h \C'$ in \eqref{co.fi.sq} is trivially non-empty --- it
has an initial enhanced object --- and by
Corollary~\ref{cof.filt.corr}, $c \setminus^h_\gamma \C'$ is
filtered, thus non-empty since $\gamma$ is cofinal. But $\rho$ and
$\lambda$ are enhanced cofibrations, and $\gamma(c')
\setminus^h_\gamma \C'$ is filtered as well, so we are done by
Corollary~\ref{prod.filt.corr}.
\endproof

\begin{prop}\label{enh.filt.prop}
For any regular cardinal $\kappa$, a small enhanced category $\C$ is
$\kappa$-filtered iff $\colimh_\C:\fFun(\C,\sSets^h) \to \sSets^h$
preserves $\kappa$-bounded limits.
\end{prop}

\proof{} For the ``if'' part, assume that $\colimh_\C$ preserves
$\kappa$-bounded limits. Then it suffices to check that $\cCone(E)$
is hyperconnected, hence non-empty for any enhanced functor $E:\E
\to \C$ from a $\kappa$-bounded $\E$. But if we let $E':\E^\iota
\times^h \C \to \sSets^h$ be the composition of $E^\iota \times^h
\id$ with the Yoneda pairing \eqref{yo.pair.eq} for $\C$, then
$\Hh^\Tot(\cCone(E)) \cong \colimh_\C\limh_{\E^\iota}E'$, whereas
$\limh_{\E^\iota}\colimh_{\C}E' \cong \limh_{\E^\iota}\Hh^\Tot(E(e)
\setminus^h \C) \cong \limh_{\E^\iota}\ppt^h \cong \ppt^h$,
since all the comma-fibers $c \setminus^h \C$, $c \in \C_\ppt$ are
hyperconnected.

Conversely, assume that $\C$ is $\kappa$-filtered. Then by
Lemma~\ref{coli.kr.le}, it suffices to check that $\colimh_\C$
preserves limits $\limh_J$ in $\sSets^h$ over all $J \in
\Posf_\kappa$. Fix such a $J$, assume given an enhanced functor
$X:\C \times \Unf(J) \to \sSets^h$, and consider the corresponding
map \eqref{cofun.facto} of Example~\ref{lim.colim.exa}, with $\C_0 =
\C$, $\C_1 = \Unf(J)$. Then $\alpha$ is invertible by
Corollary~\ref{J.filt.corr}, and $\beta$ is induced by
$\sSec^h(J,\pi)$. However, $\pi$ is an enhanced cofibration with
hyperconnected fibers that are then $\kappa$-filtered by
Lemma~\ref{filt.hyper.le}, and therefore $\sSec^h(J,\pi)$ has
$\kappa$-filtered fibers by Corollary~\ref{sec.filt.corr}. Thus
again by Lemma~\ref{filt.hyper.le}, $\beta=\Hh^\Tot(\sSec^h(J,\pi))$
is invertible as well.
\endproof

\begin{corr}\label{filt.dense.corr}
For any regular cardinal $\kappa$, the full enhanced subcategory
$\F_k \subset \cCat^h$ spanned by small $\kappa$-filtered enhanced
categories is $\colim$-closed in the sense of
Definition~\ref{colim.cl.def}.
\end{corr}

\proof{} By Corollary~\ref{colim.cl.corr}, it suffices to check that
$\F_\kappa$ satisfies the conditions \thetag{i}-\thetag{iii} of
Lemma~\ref{colim.cl.le}. Indeed it does: \thetag{iii} is
Lemma~\ref{filt.fib.le}, \thetag{ii} is Lemma~\ref{filt.sec.le}, and
\thetag{i} immediately follows from Proposition~\ref{enh.filt.prop}.
\endproof

\subsection{Inductive completions and compact objects.}

Let us now discuss what happens when one adds filtered colimits to
an enhanced category using the envelope construction of
Lemma~\ref{env.le}.

\begin{defn}
For any regular cardinal $\kappa$, an enhanced category $\E$ is {\em
  $\kappa$-filtered-cocomplete} if it is $\F_\kappa$-cocomplete with
respect to the full enhanced subcategory $\F_\kappa \subset \cCat^h$
of Corollary~\ref{filt.dense.corr}, and an enhanced functor $\E \to
\E'$ between $\kappa$-filtered-cocomplete enhanced categories $\E$,
$\E'$ is {\em $\kappa$-filtered-right-exact} if it is
$\F_\kappa$-right-exact. An enhanced object $e \in \E_\ppt$ in a
$\kappa$-filtered-complete enhanced category $\E$ is {\em
  $\kappa$-compact} if it is $\F_\kappa$-projective, and
$\Comp_\kappa(\E) \subset \E$ is the full enhanced subcategory
spanned by $\kappa$-compact objects. An enhanced category is {\em
  filtered-cocomplete} if it is $\kappa$-filtered-cocomplete for the
countable cardinal $\kappa$, and then a $\kappa$-compact object in
such a category is simply {\em compact}, and a
$\kappa$-filtered-right-exact functor between such categories is {\em
  filtered-right-exact}. The {\em $\kappa$-inductive completion} of
an enhanced category $\C$ is the enhanced category $\iInd_\kappa(\C)
= \Env(\C,\F_\kappa)$, and $\iInd(\C)$ is $\iInd(\C)_\kappa$ for the
countable cardinal $\kappa$.
\end{defn}

\begin{lemma}\label{ind.filt.le}
For any small enhanced category $\C$ and regular cardinal $\kappa$,
an enhanced functor $X:\C^\iota \to \sSets$ with the corresponding
enhanced family of groupoids $\C X \to \C$ lies in $\iInd_\kappa(\C)
\subset \C^\iota\sSets^h$ iff $\C X$ is $\kappa$-filtered.
\end{lemma}

\proof{} Combine Lemma~\ref{env.le} and
Corollary~\ref{filt.dense.corr}.
\endproof

\begin{exa}
A cocomplete enhanced category is trivially
$\kappa$-filtered-cocomplete for any $\kappa$; in particular, this
applies to the functor category $\E^\iota\sSets^h_\kappa$ for any
small enhaced category $\E$.
\end{exa}

For any $\kappa$-filtered-cocomplete enhanced category $\E$, the
full embedding $\Comp_\kappa(\E) \to \E$ extends to an enhanced
functor
\begin{equation}\label{ind.comp.eq}
\iInd_\kappa(\Comp_\kappa(\E)) \to \E,
\end{equation}
and this functor is also fully faithful by Lemma~\ref{ff.env.le}.

\begin{lemma}\label{enh.comp.le}
For any small enhanced category $\C$ and regular cardinal $\kappa$,
say that an enhanced functor $X:\C^\iota \to \sSets^h$ is {\em
  $\kappa$-generated} if we have $X \cong \gamma^o_!\ppt^h$ for some
$J \in \Posf^+_\kappa$ and enhanced functor $\gamma:\Unf(J) \to
\C$. Then $\Comp_\kappa(\C^\iota\sSets^h)$ is the Karoubi closure of
the full subcategory in $\C^\iota\sSets^h$ spanned by
$\kappa$-generated enhanced functors, and the fully faithful
embedding \eqref{ind.comp.eq} for $\E = \C^\iota\sSets^h$ is an
equivalence.
\end{lemma}

\proof{} Let $\Comp_\kappa' \subset \C^\iota\sSets^h$ be the Karoubi
closure in question. Proposition~\ref{enh.filt.prop} immediately
implies that any $\kappa$-generated $X$ is $\kappa$-compact in
$\C^\iota\sSets^h$, and that $\Comp_\kappa(\C^\iota\sSets)$ is
closed under $\kappa$-bounded limits, in particular
Karoubi-closed. Thus $\Comp_\kappa' \subset
\Comp_\kappa(\C^\iota,\sSets^h)$, and then by
Corollary~\ref{I.proj.corr}, it suffices to prove that the fully
faithful embedding $\iInd_\kappa(\Comp_\kappa') \to
\C^\iota\sSets^h$ of Lemma~\ref{ff.env.le} is essentially
surjective. In other words, we have to show that any $X:\C^\iota \to
\sSets^h$ is a $\kappa$-filtered colimit of $\kappa$-generated
functors. This is clear: by Proposition~\ref{univ.prop}, we have $X
\cong \gamma_!^o\ppt^h$ for some $J \in \Posf^+$ and $\gamma:\Unf(J)
\to \C$, and then $J$ is the $\kappa$-filtered colimit of its
$\kappa$-bounded left-closed subsets.
\endproof

\begin{exa}
Take $\C = \ppt^h$. Then Lemma~\ref{enh.comp.le} says that $S \in
\sSets^h$ is $\kappa$-compact iff it is $\kappa$-bounded.
\end{exa}

\begin{lemma}\label{kappa.f.le}
For any regular cardinals $\kappa < \kappa'$, the enhanced category
$\sSets^h_\kappa$ is $\kappa'$-filtered-cocomplete
\end{lemma}

\proof{} Note that if $\C$ is an essentially small groupoid such
that $\|\C\| \geq \kappa$, then there exists a partially ordered set
$I$ of cardinality $|I|=\kappa$ and a functor $I \to \C$ that does
not factor through an essentially small category $\C'$ with $\|\C'\|
< \kappa$. Indeed, if $|\pi_0(\C)| \geq \kappa$, it suffices to let
$I$ be a discrete set $S$ with $|S|=\kappa$ numbering mutually
non-isomorphic objects in $\C$. Otherwise, $|\pi_1(\C,c)| \geq
\kappa$ for some $c \in \C$, we can choose an injective map $S \to
\pi_1(\C,c)$, $|S|=\kappa$ that is not conjugate to any map that
factors through a set $S'$ with $|S'| < \kappa$, and then we can
take as $I$ the coproduct of $S$ copies of the partially ordered set
$\VV = \V^o \copr_{\{0,1\}} \V^o$, glued together along say $0 \in
\VV$. Note also that in both cases, $\dim I \leq 1$. Then for any
small enhanced groupoid $\C$ that is not $\kappa$-bounded, we have
$\|\C_J\| \geq \kappa$ for some finite partially ordered set $J$. We
then have $I \in \Posf_\kappa$, $|I|=\kappa$, $\dim I \leq 1$, and a
functor $\gamma:I \to \C_J \cong (J^o_h\C)_\ppt$ that does not
factor through a smaller category, and since $\dim I \leq 1$,
$\gamma$ lifts to an enhanced functor $\gamma':\Unf(I) \to J^o_h\C$,
or equivalently, $\gamma'':\Unf(I \times J^o) \to \C$. Now if $\C
\cong \colimh_{\E}C$ for some $\kappa'$-filtered $\E$ and enhanced
functor $C:\E \to \sSets^h_\kappa$, with the corresponding enhanced
cofibration $\C' \to \E$, then $\gamma''$ factors through $\C'$ by
Lemma~\ref{filt.sec.le}, and then since $\E$ is $\kappa'$-filtered,
it factors through the enhanced fiber $\C'_e$ for some $e \in
\E_\ppt$. Since $\C'_e \in \sSets^h_\kappa$, this is a
contradiction. We conclude that $\sSets^h_\kappa \subset \sSets^h$
is closed under $\kappa'$-filtered colimits, thus $\kappa'$-filtered
cocomplete.  \endproof

\begin{corr}\label{sm.acc.corr}
For any regular cardinals $\kappa < \kappa'$, a Karoubi-closed
$\kappa$-bounded enhanced category $\C$ is
$\kappa'$-filtered-cocomplete, and moreover, we have
$\Comp_{\kappa'}(\C)=\iInd_{\kappa'}(\C)=\C$.
\end{corr}

\proof{} It suffices to check that any $E \in \iInd_{\kappa'}(\C)
\subset \C^\iota\sSets^h$ actually lies in $\C \subset
\iInd_{\kappa'}(\C)$. Indeed, let $\pi:\E \to \C$ be the corresponding
enhanced family of groupoids. Then $\E$ is $\kappa'$-filtered by
Corollary~\ref{filt.dense.corr} and Karoubi-closed by
Example~\ref{loc.P.exa}, so by Lemma~\ref{yo.im.le} and
Lemma~\ref{P.o.le}, it suffices to show that $\E$ is
$\kappa$-bounded. But then as in Corollary~\ref{I.proj.corr}, $E
\cong \colimh_{\E}\pi_! \circ \Y(\E)$, and since $\C$ is
$\kappa$-bounded, $\pi_! \circ \Y(\E) \cong \Y(\C) \circ \pi:\E \to
\C^\iota\sSets^h$ factors through $\C^\iota\sSets^h_\kappa$. It
remains to apply Lemma~\ref{kappa.f.le}.
\endproof

\begin{lemma}\label{J.acc.le}
For any regular cardinal $\kappa$ and $\kappa$-bounded $J \in
\Pos_\kappa$, $\C=\Unf(J)$ is $\kappa$-filtered-cocomplete, and
$\Comp_\kappa(\C)=\iInd_\kappa(\C)=\C$. Moreover, for any
$\kappa$-filtered small enhanced category $\E$ and enhanced functor
$E:\E \to \C$, with colimit $j = \colimh E \in J$, there exists a
cofinal full enhanced subcategory $\E' \subset \E$ such that
$E|_{\E'}$ factors through $\eps^h(j):\ppt^h \to \C$.
\end{lemma}

\proof{} Since $J \in \Pos$, the Yoneda embedding
$\C \to \C^\iota\sSets^h$ factors through
$\C^\iota\Unf([1])$, where the enhanced full
subcategory $\Unf([1]) \subset \sSets^h$ is spanned by $\emptyset$
and $\ppt^h$. This full subcategory is obviously closed under
$\kappa$-filtered colimits, so that $\iInd_\kappa(\C) \subset
\C^\iota\sSets^h$ lies inside $\C^\iota\Unf([1]) \subset
\C^\iota\sSets^h$. Then as in Corollary~\ref{sm.acc.corr}, for any
$X \in \iInd_\kappa(\C)$, we have $\C X \cong \Unf(J')$ for some
left-closed subset $J' \subset J$, and this $J'$ must be
$\kappa$-directed, so that $J' = J/j$ for some $j \in J$. This
proves the first claim. For the second, take $X = E^o_!\ppt^h \in
\Ind_\kappa(\C)$, identify $\C X = \Unf(J/j)$, and note that the
induced functor $E':\E \to \Unf(J/j)$ is then cofinal by
Example~\ref{EX.exa}, so that $\E' = \E_j \cong j \setminus^h \E
\subset \E$ is cofinal by Lemma~\ref{filt.hyper.le}.
\endproof

Next, we want to show that $\iInd_\kappa$ commutes with some natural
operations such that taking the comma-categories. This has to be
stated carefully since even for a small enhanced category $\C$,
$\iInd_\kappa(\C)$ is often not small, so that comma-categories are
not well-defined. However, for any $J \in \Pos$ and small
$J$-augmented enhanced category $\C$, with transpose $J$-coaugmented
enhanced category $\C_{h\perp}$, fibers $\C_j$, $j \in J$, and
embeddings $e(j):\C_j \to \C$, we can consider the full enhanced
subcategory
\begin{equation}\label{indsec.eq}
\Indsec^h_\kappa(J^o,\C_{h\perp})  \subset \C^\iota\sSets^h
\end{equation}
spanned by enhanced functors $X:\C^\iota \to \sSets^h$ such that for
any $j \in J$, $e(j)^*X \in \iInd_\kappa(\C_j) \subset
\C_j^\iota\sSets$. Note that $\Indsec^h_\kappa(J^o,\C_{h\perp})$ is
trivially $\kappa$-filtered-cocomplete.

\begin{lemma}\label{sec.ind.le}
For any $J \in \Posf^+_\kappa$ and $J^o$-augmented small enhanced
category $\C$, with the transpose $J$-coaugmented enhanced category
$\C_{h\perp}$, we have $\iInd_\kappa(\sSec^h(J,\C_{h\perp})) \cong
\Indsec^h_\kappa(J,\C_{h\perp})$.
\end{lemma}

\proof{} The relative Yoneda embedding \eqref{yo.sec.eq} provides an
enhanced functor
\begin{equation}\label{ind.sec.eq}
\Y(\C|J^o):\sSec^h(J,\C_{h\perp}) \to \C^\iota\sSets^h,
\end{equation}
and since it is functorial with respect to $J$ and reduces to
\eqref{yo.enh.eq} for $J=\ppt$, it factors through
$\Indsec^h_\kappa(J,\C_{h\perp}) \subset \C^\iota\sSets^h$. Moreover,
explicitly, \eqref{ind.sec.eq} is given by \eqref{Y.tw.eq}, so since
$\Tw(J)$ is $\kappa$-bounded, \eqref{ind.sec.eq} sends all enhanced
objects to $\kappa$-compact objects by
Proposition~\ref{enh.filt.prop}. Then by Lemma~\ref{ff.env.le},
\eqref{ind.sec.eq} extends to a fully faithful enhanced functor
\begin{equation}\label{ind.sec.ind.eq}
  \iInd_\kappa(\sSec^h(J,\C_{h\perp})) \to \Indsec^h_\kappa(J,\C_{h\perp}),
\end{equation}
and it suffices to check that it is essentially surjective. Indeed,
for any $X \in \Indsec^h_\kappa(J,\C_{h\perp})$ and $j \in J$, the
forgetful functor $\C X \to C \to \Unf(J)$ is a $J$-augmentation,
and for any $j \in J$, its enhanced fiber $(\C X)_j$ is
$\kappa$-filtered by Lemma~\ref{ind.filt.le}. Then by
\eqref{tw.Y.col} and Corollary~\ref{sec.filt.corr}, so is the
comma-category $\sSec^h(J,\C_{h\perp}) /^h_{\Y(\C|J)} X \cong
\sSec^h(J,(\C X)_{h\perp})$, so to see that $X$ lies in the
essential image of \eqref{ind.sec.ind.eq}, it suffices to see that
$\Y(\C X|J^o)_!\ppt^h \cong \ppt^h$. This can be checked after
composing with $e(j)^*$ for any $j$, and since $\Y(-|J)$ is
functorial with respect to $J$, this amounts to checking that the
evaluation enhanced functor $\ev_j:\sSec^h(J,(\C X)_{h\perp}) \to
(\C X)_j$ is cofinal. This is again Corollary~\ref{sec.filt.corr}.
\endproof

\begin{corr}\label{fun.ind.corr}
For any $J \in \Posf^+_\kappa$ and small enhanced category $\C$, we
have $\iInd_\kappa(J^o_h\C) \cong J^o_h\iInd_\kappa(\C)$.
\end{corr}

\proof{} If we consider the constant augmentation $\Unf(J) \times^h
\C \to \Unf(\C)$, then $\Indsec^h_\kappa(J^o,(\Unf(J) \times^h
\C)_{h\perp}) \cong J^o_h\iInd_\kappa(\C)$.
\endproof

\begin{exa}\label{sec.ind.exa}
Take $J = [1]$. For any enhanced functor $\gamma:\C_0 \to \C_1$
between small enhanced $\C_0$, $\C_1$, we can apply
Lemma~\ref{sec.ind.le} to $\C = \Cyl_h^\iota(\gamma)$. Then
$\sSec^h([1],\C_{h\perp}) \cong \C_0 /^h_\gamma \C_1$ by
Lemma~\ref{cyl.comma.le}, and if $\iInd_\kappa(\C_l)$, $l=0,1$ are
small, then $\Indsec^h_\kappa([1],\C) \cong \iInd_\kappa(\C_0)
/^h_{\iInd_\kappa(\gamma)} \iInd_\kappa(\C_1)$ by
\eqref{ffun.cyl.co}.
\end{exa}

Lemma~\ref{sec.ind.le} does not apply to right comma-categories
directly (indeed, an analog of Corollary~\ref{sec.filt.corr} for
$J$-augmented categories is very wrong). However, note that since for
any $\kappa$-bounded partially ordered set $J$, we have
$\iInd_\kappa(\Unf(J))\cong \Unf(J)$ by Lemma~\ref{J.acc.le}, we
have an enhanced functor
\begin{equation}\label{ind.J}
\iInd_\kappa(\phi):\iInd(\C) \to \Unf(J) = \iInd_\kappa(\Unf(J))
\end{equation}
for any enhanced category $\C$ and enhanced functor $\phi:\C \to
\Unf(J)$. If we have another $\kappa$-bounded partially ordered set
$J'$ equipped with a map $f:J' \to J$, then \eqref{ind.J} fits into a
commutative square
\begin{equation}\label{ind.J.sq}
\begin{CD}
\iInd_\kappa(\C') @>>> \iInd_\kappa(\C)\\
@V{\iInd_\kappa(\phi')}VV @VV{\iInd_\kappa(\phi)}V\\
\Unf(J') @>{\Unf(f)}>> \Unf(J),
\end{CD}
\end{equation}
where $\C' = \Unf(f)^*\C$, with the induced projection $\phi':\C' \to
\Unf(J')$.

\begin{lemma}\label{sec.ind.bis.le}
For any $\kappa$-bounded partially ordered set $J$ and small
enhanced category $\C$ equipped with an enhanced functor $\phi:\C
\to \Unf(J)$, and for any other $\kappa$-bounded partially ordered
set $J'$ and a map $f:J' \to J$, the square \eqref{ind.J.sq} is
cartesian, and
\begin{equation}\label{sec.ind.bis.eq}
  \iInd_\kappa(\sSec^h(J',\C)) \cong \sSec^h(J',\iInd_\kappa(\C)).
\end{equation}
Moreover, if $\phi$ is an augmentation or a coaugmentation, then so
is the enhanced functor \eqref{ind.J}.
\end{lemma}

\proof{} If $J'=\ppt$, so that $f:\ppt \to J$ is the embedding onto
an element $j \in J$, then the corresponding embedding functor
$e(j):\C_j \to \C$ is fully faithful, thus induces a fully faithful
functor
\begin{equation}\label{ind.J.j}
\iInd_\kappa(\C_j) \to \iInd_\kappa(\C)_j.
\end{equation}
By Lemma~\ref{J.acc.le}, for any enhanced functor $E:\E \to \Unf(J)$
from a $\kappa$-filtered $\E$ such that $\colimh_\E E = j$, we have
$E=j$ on a cofinal subcategory $\E' \subset \E$. Therefore
\eqref{ind.J.j} is essentially surjective, thus an equivalence. This
proves the first claim for $J'=\ppt$. By
Corollary~\ref{fun.ind.corr}, we have
$\iInd_\kappa(\fFun^h(\Unf(J),\C)) \cong
\fFun^h(J,\iInd_\kappa(\C))$, so applying the first claim to
\eqref{ssec.enh.II} proves \eqref{sec.ind.bis.eq} for any $f:J' \to
J$. This in turns proves that for any $f:J' \to J$, \eqref{ind.J.sq}
is cartesian over $\Posf^+_\kappa \subset \Posf^+$, so by
Corollary~\ref{fin.enh.equi.corr}, the first claim holds in full
generality as soon as $\C$ is small. Since any enhanced category
$\C$ is a filtered $2$-colimit of its small full enhanced
subcategories, it then holds for any $\C$. Finally, if $\phi$ is an
augmenation or coaugmentation, then to show the same for
$\iInd_\kappa(\phi)$, use Lemma~\ref{env.adj.le} and the first claim
for $J' = \Ar(J)$ and $f = \tau$ resp.\ $f=\sigma$.
\endproof

\begin{remark}\label{ind.sec.rem}
{\em A posteriori}, we see that in Lemma~\ref{sec.ind.le}, the
enhanced category \eqref{indsec.eq} is actually given by
\begin{equation}\label{indsec.bis.eq}
\Indsec^h_\kappa(J^o,\C_{h\perp}) \cong
\sSec^h(J^o,\iInd_\kappa(\C_{h\perp})),
\end{equation}
where $\iInd_\kappa(\C_{h\perp}) \to \Unf(J^o)$ is a
$J^o$-coagumentation by Lemma~\ref{sec.ind.bis.le}.
\end{remark}

\begin{exa}\label{noprod.exa}
An analog of Lemma~\ref{sec.ind.le} for $\sSec^h_{\Tot}$ instead of
$\sSec^h$ is {\em not} true already for $J=\V$. Indeed, consider $\N
\in \Pos$, let $q_l:\N \to \N$, $l=0,1$ be the map $n \mapsto 2n +
l$, and let $\N_\idot \to \V$ be the fibration with fibers $\N$ and
transition maps $q_0$, $q_1$. Then $\iInd(\Unf(\N)) \cong
\Unf(\N^>)$, with the new object $o \in \N^>$ corresponding to the
filtered colimit of $\id:\Unf(\N) \to \Unf(\N)$, and since both
$q_0$ and $q_1$ are cofinal, $q_0(o) = q_1(o) = o$. But the full
subcategory $\sSec^h_\Tot(\V^o,\iInd(\Unf(\N))_\idot) \subset
\iInd(\Unf(\N_\idot)|V) \cong
\sSec^h_{\Tot}(\V^o,\iInd(\Unf(\N))_\idot)$ is then the product
$\Unf(\N^> \times_{\N^>} \N^>)$, and contains a non-trivial object
$o \times o$, while $\sSec^h_\Tot(\V^o,\Unf(\N_\idot)_\perp)$ is
$\Unf(\N \times_{\N} \N) \cong \emptyset$.

More generally, for any regular cardinal $\kappa$, take a
$\kappa$-directed well-ordered set $Q_\kappa$, $|Q_\kappa|=\kappa$
of Example~\ref{ord.exa}, and let $Q = Q_\kappa \times \N$ with
lexicographical order, with the cofinal maps $\id \times q_l:Q \to
Q$, $l=0,1$. Then again, we have $\iInd_\kappa(\Unf(Q)) \cong
\Unf(Q^>)$, and $\id \times q_l$, $l=0,1$ send the new object $o \in
Q^>$ to itself, so that $Q^> \times_{Q^>} Q^>$ is not empty, while
$Q \times_Q Q \cong \emptyset$.
\end{exa}

\subsection{Semicartesian products.}

In words, Example~\ref{noprod.exa} says that $\iInd_\kappa$ does not
compute with semicartesian products. It turns out that this can be
circumvented if one increases $\kappa$. In general, for any regular
cardinals $\kappa \leq \mu$ and enhanced category $\C$, denote
$\iInd_\kappa^\mu(\C) = \Comp_\mu(\iInd_\kappa(\C))$, and note that
$\iInd_\kappa^\kappa(\C)$ is the Karoubi closure of $\C$ by
Corollary~\ref{I.proj.corr}.

\begin{prop}\label{ka.mu.prop}
Assume given regular cardinals $\kappa < \mu$ such that $\kappa \trl
\mu$ in the sense of Definition~\ref{sharp.def}. Then for any
$\kappa$-filtered small enhanced category $\C$,
$\iInd^\mu_\kappa(\C)$ is $\mu$-filtered. Moreover, for any small
enhanced category $\C$, the fully faithful enhanced functor
\begin{equation}\label{ind.ka.mu}
  \iInd_\mu(\iInd^\mu_\kappa(\C)) \to \iInd_\kappa(\C)
\end{equation}
of \eqref{ind.comp.eq} is an equivalence, and $\iInd^\mu_\kappa(\C)
= \Comp_\mu(\C^\iota\sSets^h) \cap \iInd_\kappa(\C)$.
\end{prop}

\proof{} Denote $\iInd^\mu_\kappa(\C)' = \Comp_\mu(\C^\iota\sSets^h)
\cap \iInd_\kappa(\C)$. Since the embedding $\iInd_\kappa(\C) \to
\C^\iota\sSets^h$ is $\kappa$-filtered-right-exact, we have
$\iInd^\mu_\kappa(\C)' \subset \iInd^\mu_\kappa(\C)$.  For the first
claim, if $\C$ is $\kappa$-filtered, then $\ppt^h \in
\C^\iota\sSets^h$ is in $\iInd_\kappa(\C)$ by
Lemma~\ref{ind.filt.le}, and by the same Lemma~\ref{ind.filt.le}, it
suffices to show that it lies in the essential image of the full
embedding \eqref{ind.ka.mu}. We thus have to show that $\ppt^h$ is a
$\mu$-filtered colimit of $\mu$-compact objects in
$\iInd_\kappa(\C)$. By Lemma~\ref{enh.comp.le},
$\Comp_\mu(\C^\iota\sSets^h)$ is $\mu$-filtered, and we have
\begin{equation}\label{comp.mu}
  \ppt^h \cong \colimh_{\Comp_\mu(\C^\iota\sSets)}\nu,
\end{equation}
where $\nu:\Comp_\mu(\C^\iota\sSets^h) \to \C^\iota\sSets^h$ denotes
the embedding functor. By Lemma~\ref{ff.filt.le}, it suffices to
show that $\iInd^\mu_\kappa(\C)' \subset
\Comp_\mu(\C^\iota,\sSets^h)$ is cofinal --- then \eqref{comp.mu}
reduces to a colimit over $\iInd_\kappa^\mu(\C)'$ that is moreover
$\mu$-filtered. In addition to this, Lemma~\ref{ff.filt.le} says
that it suffices to check that for any $X \in
\Comp_\mu(\C^\iota\sSets^h)$, there exists a map $X \to Y$ to some
$Y \in \iInd_\kappa^\mu(\C)'$. By Lemma~\ref{enh.comp.le}, we may
further assume that $X \cong \gamma^o_!\ppt^h$ for some $J \in
\Posf^+_\mu$ and enhanced functor $\gamma:\Unf(J) \to \C$. But then,
we can take a left-closed embedding $J \to J'$ provided by
Corollary~\ref{sharp.corr}, with $J'$ Noetherian, $\kappa$-filtered
and $\mu$-bounded, and by Corollary~\ref{zorn.corr} for the
identifity map $f = \id:J' \to J'$, the enhanced functor $\gamma$
extends to an enhanced functor $\gamma':\Unf(J') \to \C$. Now take
$Y = {\gamma'}^o_!\ppt^h$.

For the second claim, since $\iInd^\mu_\kappa(\C)'$ is
Karoubi-closed, it actually suffices to show that the full embedding
$\iInd_\mu(\iInd^\mu_\kappa(\C)') \to \iInd_\kappa(\C)$ induced by
\eqref{ind.ka.mu} is essentially surjective. In other words, we
need to check that any $X \in \iInd_\kappa(\C) \subset
\C^\iota\sSets^h$ is a $\mu$-filtered colimit of objects in
$\iInd_\kappa^\mu(\C)'$. Indeed, let $\gamma:\C X \to \C$ be the
enhanced family of groupoids corresponding to $X$. Then $\gamma^o_!$
sends $\iInd^\mu_\kappa(\C X)'$ into $\iInd^\mu_\kappa(\C)'$. But
$\C X$ is $\kappa$-filtered by Lemma~\ref{env.le}, so by the first
claim, $\iInd^\mu_\kappa(\C X)'$ is $\mu$-filtered and cofinal in
$\Comp_\mu((\C X)^\iota\sSets^h)$.  Then the required colimit
presentation for $X=\gamma_!^o\ppt^h$ is induced by \eqref{comp.mu}.
\endproof

\begin{corr}\label{ka.mu.corr}
For any regular cardinals $\kappa \leq \mu$ and small enhanced
category $\C$, the enhanced category $\iInd^\mu_\kappa(\C)$ is
small.
\end{corr}

\proof{} For any $\mu' \geq \mu$, we have $\iInd^{\mu'}_\kappa(\C)
\supset \iInd^\mu_\kappa(\C)$, so by Example~\ref{P.ka.exa}, we may
assume that $\kappa \trl \mu$. Then $\Comp_\mu(\C^\iota\sSets^h)$ is
small by Lemma~\ref{enh.comp.le}, hence so is $\iInd^\mu_\kappa(\C)
= \Comp_\mu(\C^\iota\sSets^h) \cap \iInd^\mu_\kappa(\C)$.
\endproof

\begin{corr}\label{ka.mu.X.corr}
Assume given regular cardinals $\kappa \trl \mu$ and a small
enhanced category $\C$. Then for any enhanced functor $X:\C^\iota
\to \sSets^h$, we have $\iInd^\mu_\kappa(\C X) \cong
\iInd^\mu_\kappa(\C) /^h X$, where the comma-fiber is taken with
respect to the full embedding $\iInd^\mu_\kappa(\C) \to
\Ind_\kappa(\C) \to \C^\iota\sSets^h$.
\end{corr}

\proof{} Let $\gamma:\C X \to \C$ be the enhanced family of
groupoids corresponding to $X$. Then by adjunction, $\gamma^o_!$
sends $\mu$-compact objects to $\mu$-compact objects, thus induces
an enhanced functor $\iInd^\mu_\kappa(\C X) \to \iInd^\mu_\kappa(\C)
/^h X$, and it suffices to prove that the corresponding fully
faithful embedding
\begin{equation}\label{ka.mu.X}
\iInd_\kappa(\C X) \cong \iInd_\mu(\iInd^\mu_\kappa(\C X)) \to
\iInd_\mu(\iInd^\mu_\kappa(\C) /^h X)
\end{equation}
provided by Lemma~\ref{ff.env.le} is essentially surjective. Indeed,
up to an isomorphism, an enhanced object in its target is given by
an enhanced object $Y$ in $\iInd_\mu(\iInd^\kappa_\mu(\C)) \cong
\iInd_\kappa(\C)$ and a morphism $Y \to X$, and any such
pair manifestly lies in the essential image of the embedding
\eqref{ka.mu.X}.
\endproof

Corollary~\ref{ka.mu.corr} shows that for any regular cardinals
$\kappa \leq \mu$, we have a well-defined enhanced functor
\begin{equation}\label{ind.ka.mu.eq}
\iInd_\kappa^\mu:\cCat^h \to \cCat^h.
\end{equation}
For any $\kappa$-bounded $J \in \Posf^+_\kappa$ and $X:\Unf(J) \to
\cCat^h$, with the corresponding $J$-coaugmented $\C = \Unf(J)^\perp
X$, the $\kappa$-inductive completion $\iInd_\kappa(\C)$ is
$J$-coaugmented by \eqref{ind.J} by Lemma~\ref{sec.ind.bis.le}, and
$\iInd^\mu_\kappa(\C) \subset \iInd_\kappa(\C)$ is a $J$-coaugmented
full subcategory corresponding to $\iInd^\mu_\kappa \circ X:\Unf(J)
\to \cCat^h$.

\begin{lemma}\label{ind.ka.mu.le}
For any regular cardinals $\kappa \trl \mu$, $J \in \Posf_\kappa$
and $J$-coaug\-men\-t\-ed small enhanced category $\C$, we have
\begin{equation}\label{ind.sec.ka.mu}
\iInd^\mu_\kappa(\sSec^h(J,\C)) \cong
\sSec^h(J,\iInd^\mu_\kappa(\C)).
\end{equation}
\end{lemma}

\proof{} By Lemma~\ref{sec.ind.le} and \eqref{indsec.bis.eq}, it
suffices to check that an enhanced object $X$ in $\Indsec^h(J,\C)$ is
$\mu$-compact if and only if so is $e(j)^*X \in \iInd_\kappa(\C_j)$
for any $j \in J$. For the ``if'' part, $X$ is a $\kappa$-bounded
colimit of objects $e(j)_!e(j)^*X$, $j \in J$, and since $e(j)$ is
trivially $\mu$-filtered-right-exact, all these objects are
$\mu$-compact by adjunction, and then $X$ is $\mu$-compact by
Proposition~\ref{enh.filt.prop}. For the ``only if'' part, again by
Proposition~\ref{enh.filt.prop} and \eqref{enh.kan.eq}, $e(j)_*$ is
$\mu$-filtered-right-exact for any $j \in J$, so by adjunction,
$e(j)^*$ sends $\mu$-compact objects to $\mu$-compact objects.
\endproof

\begin{remark}\label{ind.ka.mu.rem}
Note that even for $\kappa=\mu$, Lemma~\ref{ind.ka.mu.le} is
non-trivial: it says that $\sSec^h(J,-)$ commutes with Karoubi
envelopes.
\end{remark}

In particular, by Example~\ref{ord.exa}, we always have $\kappa \trl
\kappa^+$ for any regular cardinal $\kappa$ with successor cardinal
$\kappa^+$. For any small enhanced category $\C$, let $\C^+ =
\iInd^{\kappa^+}_\kappa(\C)$. Then by Proposition~\ref{ka.mu.prop},
$\iInd_{\kappa^+}(\C^+) \cong \iInd_\kappa(\C)$, and if $\C$ is
$\kappa$-filtered, then $\C^+$ is $\kappa^+$-filtered. By
Corollary~\ref{ka.mu.corr}, it is also small, and by
Corollary~\ref{ka.mu.X.corr}, $(\C X)^+ \cong \C^+ /^h X$ for any
$X:\C^\iota \to \sSets^h$.

\begin{lemma}\label{ind.prod.le}
Assume given $\kappa$-filtered small enhanced categories $\C$,
$C_0$, $\C_1$ and cofinal enhanced functors $\gamma_l:\C_l \to \C$,
$l=0,1$. Then $\C_0^+ \times^h_{\C^+} \C_1^+$ is
$\kappa^+$-filtered, and the projections $\C_0^+ \times^h_{\C^+}
\C_1^+ \to \C_0^+,\C_1^+$ are cofinal.
\end{lemma}

\proof{} Let $\C^+_\idot$ be the $\V^o$-coaugmented enhanced
category with fibers $\C^+_0$, $\C^+_1$, $\C^+$ and transition
functors $\gamma_0^+$, $\gamma_1^+$ provided by
Corollary~\ref{ccat.dim1.corr}. Then by Lemma~\ref{sec.ind.le},
Example~\ref{sec.ind.exa} and Corollary~\ref{sec.filt.corr}, it
suffices to show that the fully faithful embedding $\C_0^+
\times^h_{\C^+} \C_1^+ \cong \sSec^h_\Tot(\V^o,\C_\idot^+) \to
\sSec^h(\V^o,\C^+_\idot)$ is cofinal, and by Lemma~\ref{ff.filt.le},
it further suffices to check that it has non-empty right
comma-fibers. Let $X$ be an enhanced object in
$\sSec^h(\V^o,\C^+_\idot)$ corresponding to a $V^o$-coaugmented
enhanced functor $\gamma:\Unf(J) \to \C^+_\idot$, where $J =
\Ar(\V)^o$ is $\V^o$-coaugmented via the projection
$\sigma^o:\Ar(\V)^o \to \V^o$. Let $Q = Q_\kappa \times \N$ be the
well-ordered set of Example~\ref{noprod.exa}, with embeddings $q_l:Q
\to Q$, $l=0,1$, and let $Q_\idot \to \V^o$ be the cofibration with
fibers $Q$ and transition functors $q_0$, $q_1$. Let $J' = J
\copr_\lambda Q_\idot$, where the gluing map $\lambda:Q \to L(J)$
sends any $q \in (Q_\idot)_v$, $v \in \V^o$ to the whole fiber
$J_v$. Then $\pi:J' \to \V^o$ is a cofibration with cofinal
transition functors, and it suffices to check that $\gamma$ extends
to a $\V^o$-coaugmented enhanced functor $\gamma':\Unf(J') \to
\C^+_\idot$ (then $Y=\Unf(\pi)_!\gamma':\V^o \to \C^+_\idot$ exists
since $|Q| = \kappa < \kappa^+$, and gives the desired enhanced
object $Y \in \sSec^h_\Tot(\V^o,\C^+_\idot)$ with a map $X \to Y$).

To construct $\gamma'$, we use Corollary~\ref{zorn.I.corr}, or
rather, its obvious dual version with augmented enhanced categories
and functors replaced with coaugmented ones. Take $I = \V^o$ and $J
= J_o$, with the map $f:J' \to J$ adjoint to the embedding $J_o \to
J'$, and take our original $J \subset J'$ as $J_0$. We then have to
check that for any $q \in Q = J \ssetminus J_0$, a
$\V^o$-coaugmented enhanced functor $\Unf(f^{-1}(J /' q)) \to
\C^+_\idot$ extends to a $V^o$-coaugmented enhanced functor
$\Unf(f^{-1}(J/q)) \to \C^+_\idot$. But we have $Q \cong q_0(Q)
\copr q_1(Q)$ as sets, so the complement $f^{-1}(q) = f^{-1}(J/q)
\ssetminus f^{-1}(J/'q)$ is $[1]$-coaugmented via one of the two
embeddings $\eps(0)^<,\eps(1)^<:[1] = \ppt^< \to \V^o =
\{0,1\}^<$. Let $\eps$ be this embedding, and note that by
Lemma~\ref{sec.sq.le}, we can restrict our attention to what happens
over $\eps([1]) \subset \V^o$. We thus have a $[1]$-coaugmented
small enhanced category $\eps^*\C^+_\idot$, a cofibration $\eps^*J'
\to [1]$, a subset $\eps^*f^{-1}(J/'q) \subset \eps^*J'$, and a
$[1]$-coaugmented enhanced functor $\Unf(\eps^*f^{-1}(J'/q)) \to
\C^+_\idot$, and we need a $[1]$-coaugmented relative cone
for this enhanced functor. Then a relative cone over $[1]$ is
provided by Lemma~\ref{sec.filt.le}, and to turn it into a
$[1]$-coaugmented relative cone, use
Corollary~\ref{comma.filt.corr} and Example~\ref{comma.filt.exa}.
\endproof

\begin{corr}\label{ind.prod.corr}
For any small enhanced categories $\C$, $\C_0$, $\C_1$ and enhanced
functors $\gamma_l:\C_l \to \C$, $l=0,1$, and any regular cardinal
$\kappa$, the square
\begin{equation}\label{ind.prod.sq}
\begin{CD}
\iInd_{\kappa^+}(\C_0^+ \times^h_{\C^+} \C_1^+) @>>>
\iInd_{\kappa^+}(\C_1^+)\\
@VVV @VVV\\
\iInd_{\kappa^+}(\C_0^+) @>>> \iInd_{\kappa^+}(\C^+)
\end{CD}
\end{equation}
is semicartesian.
\end{corr}

\proof{} For any regular cardinal $\mu$, $I \in \Pos$, and
$I$-augmented small enhanced category $\C$, let
$\Indsec_\mu^{\Toth}(I^o,\C_{h\perp}) \subset
\Indsec^h_\mu(I^o,\C_{h\perp})$ be the full enhanced subcategory
defined by the cartesian square
\begin{equation}\label{ind.rel.sq}
\begin{CD}
\Indsec^{\Toth}_\mu(I^o,\C_{h\perp}) @>>>
\Indsec^h_\mu(I^o,\C_{h\perp})\\
@VVV @VVV\\
\Sec_\Tot(I^o,\iInd_\mu(\C)_\idot) @>>>
\Sec(I^o,\iInd_\mu(\C)_\idot),
\end{CD}
\end{equation}
where $\iInd_\mu(\C)_\idot \to I^o$ is the cofibration with fibers
$\iInd_\mu(\C_i)$, $i \in I$, and transition functors $\gamma^o_!$
of \eqref{Y.ga.sq} corresponding to the transition functors $\gamma$
of the $I$-augmented enhanced category $\C$, while the vertical
arrow on the right in \eqref{ind.rel.sq} is induced by the
comparison functor \eqref{E.I.Y.eq}. Explicitly, by virtue of
Example~\ref{EX.exa}, $\Indsec^{\Toth}_\mu(I^o,\C_{h\perp}) \subset
\C^\iota\sSets^h$ is the full enhanced subcategory spanned by
enhanced functors $X:\C^\iota \to \sSets^h$ such that for any $i \in
I$, with the embedding functor $\nu_i:\C_i \to \C$, $X_i =
\nu^{o*}_iX \in \iInd_\mu(\C_i)$, and for any $i \leq i'$, the
corresponding transition functor $\C_iX_i \to \C_{i'}X_{i'}$ of the
$I$-augmented enhanced category $\C X$ is cofinal. Then if $I$ has a
largest element $o \in I$, it immediately follows from
\eqref{enh.cokan.eq} and Example~\ref{EX.exa} that
$\nu^o_!:\iInd_\mu(\C_o) \to \Indsec^h_\mu(I^o,\C_{h\perp})$ induces
an equivalence $\iInd_\mu(\C_o) \cong
\Indsec^{\Toth}_\mu(I^o,\C_{h\perp})$. Explicitly, for any $X \in
\iInd_\mu(\C_o)$ and $i \in I$, we have $(\nu^o_!X)_i \cong
\gamma^o_!X$, where $\gamma:\C_o \to \C_i$ is the transition functor
of the $I$-augmented category $\C$.

Now let $\C^+_\idot \in \Cat^h(\V) \cong \Cat^h(\V|\V)$ be the
$\V$-augmented small enhanced category with fibers $\C_0^+$,
$\C_1^+$, $\C^+$ and transition functors $\gamma_0^+$, $\gamma_1^+$,
identify $[1]^2 = \V^>$, with the embedding $\eps:\V \to [1]^2$, and
extend $\C^+_\idot$ to a $[1]^2$-augmented small enhanced category
$\eps^h_*\C^+_\idot$ with enhanced fiber $\C_{01}^+ = \C^+_0
\times^h_{\C^+} \C_1^+$ over the largest element $o \in [1]^2$. Then
on one hand, we have
$\Indsec_{\kappa^+}^{\Toth}([1]^2,(\eps^h_*\C^+_\idot)_{h\perp})
\cong \iInd_{\kappa^+}(\C_{01}^+)$, and on the other hand,
Corollary~\ref{ccat.dim1.corr} and Lemma~\ref{cat.epi.le}
immediately imply that the comparison functor
\begin{equation}\label{sec.V.eq}
\Indsec_{\kappa^+}^{\Toth}(\V^o,(\C^+_\idot)_{h \perp}) \to
\iInd_{\kappa^+}(\C_0^+) \times_{\iInd_{\kappa^+}(\C)}
\iInd_{\kappa^+}(\C_1^+)
\end{equation}
induced by \eqref{E.I.Y.eq} is an epivalence. Thus to prove the
claim, it suffices to check that
$\eps^{o*}:\Indsec_{\kappa^+}^{\Toth}([1]^2,(\eps^h_*\C^+_\idot)_{h\perp})
\to \Indsec^{\Toth}_{\kappa^+}(\V^o,(\C^+_\idot)_{h\perp})$ is an
equivalence. This reduces to checking that for any $X^+ \in
\iInd_{\kappa^+}(\C^+)$ and $X^+_l \in \iInd_{\kappa^+}(\C^+_l)$,
$l=0,1,01$, any commutative square
\begin{equation}\label{C.p.sq}
\begin{CD}
\C_{01}^+X_{01}^+ @>{\gamma_1'}>> \C_1^+X_1^+\\
@V{\gamma_0'}VV @VV{\gamma_1}V\\
\C_0^+X_0^+ @>{\gamma_0}>> \C^+X^+
\end{CD}
\end{equation}
over $\eps^h_*\C^+_\idot$ with cofinal $\gamma_0$, $\gamma_1$ is
semicartesian if and only if $\gamma_0'$ and $\gamma_1'$ are
cofinal. Indeed, $\C^+X^+ \cong (\C X)^+$ for $X = X^+ \in
\iInd_\kappa(\C) \cong \iInd_{\kappa^+}(\C^+)$ by
Corollary~\ref{ka.mu.X.corr}, and similarly for $\C_l^+X_l^+$,
$l=0,1$, so that the ``only if'' part follows from
Lemma~\ref{ind.prod.le}. For the ``if'' part, we may replace $\C$,
$\C_0$, $\C_1$ with $\C X$, $\C_0 X_0$, $\C_1 X_1$, so that it
suffices to consider the case when $X=\ppt^h$, $X_l= \ppt^h$,
$l=0,1$. Then for any enhanced object $c_{01}$ in $\C_{01}^+$, with
$c_l=\gamma'_l(c_{01})$, $l=0,1$, and $c = \gamma_0(c_0) \cong
\gamma_1(c_1)$, we have
$$
c_{01} \setminus^h \C_{01}^+X_{01}^+ \cong (c_0 \setminus^h \C_0^+)
\times^h_{(c \setminus^h \C^+)} (c_1 \setminus^h \C_1^+),
$$
and this comma-fiber is $\kappa^+$-filtered by
Lemma~\ref{ind.prod.le}. Since the embedding $X(c_{01}) \subset
c_{01} \setminus^h \C_{01}^+X_{01}^+$ is right-admissible, the
enhanced groupoid $X^+(c_{01})$ is then hyperconnected for any
$c_{01}$, so that $X_{01}^+ \cong \ppt^h$.
\endproof

\subsection{Accessible categories.}

We can now use filtered colimits and inductive completions to
characterize a class of large enhanced categories that admit a
reasonable theory in the spirit of Section~\ref{yo.sec}.

\begin{defn}\label{acc.def}
For any regular cardinal $\kappa$, an enhanced category $\E$ is {\em
  $\kappa$-accessible} if it is $\kappa$-filtered-cocomplete,
$\Comp_\kappa(\E) \subset \E$ is small, and the fully faithful
functor \eqref{ind.comp.eq} is essentially surjective. An enhanced
functor $\E \to \E'$ between $\kappa$-accessible categories is {\em
  $\kappa$-accessible} if it is $\kappa$-filtered-right-exact. An
enhanced category $\E$ is {\em accessible} if it is accessible for
some regular cardinal $\kappa$, and an enhanced functor $\E \to \E'$
between accessible categories is {\em accessible} if it is
$\kappa$-accessible for some $\kappa$ such that $\E$, $\E'$ are
$\kappa$-accessible.
\end{defn}

Equivalently, an enhanced category $\E$ is $\kappa$-accessible iff
$\E \cong \iInd_\kappa(\C)$ for a small enhanced category $\C$, and
it actually suffices to require that for some small full enhanced
subcategory $\C \subset \Comp_\kappa(\E)$, the fully faithful
embedding $\iInd_\kappa(\C) \to \E$ of Lemma~\ref{ff.env.le} is
essentially surjective. Such a $\C \subset \Comp_\kappa(\E)$ is
called {\em a category of $\kappa$-compact generators} of $\E$, it
is generating in the sense of Definition~\ref{enh.gen.def}, and
$\Comp_\kappa(\E) = \iInd^\kappa_\kappa(\C)$ is then the Karoubi
closure of $\C$. For any regular cardinals $\kappa \trl \mu$, a
$\kappa$-accessible category $\E$ is $\mu$-accessible by
Proposition~\ref{ka.mu.prop}, and by Example~\ref{P.ka.exa}, for any
collection of accessible enhanced categories and accessible enhanced
functors between them, we can choose a single regular cardinal
$\kappa$ such that all of them are $\kappa$-accessible. Moreover,
$\kappa$ can be chosen to be arbitrarily large (that is, larger than
any given $\kappa'$). By Lemma~\ref{enh.comp.le}, any enhanced
object $X$ in an accessible enhanced category $\E$ is $\mu$-compact
for a large enough $\mu$, so that we have
\begin{equation}\label{mu.cup}
\E = \bigcup_\mu\Comp_\mu(\E),
\end{equation}
and by Corollary~\ref{ka.mu.prop}, all the categories in the
right-hand side are small.

\begin{exa}\label{acc.pt.exa}
For any regular cardinal $\kappa$ and $\kappa$-accessible enhanced
category $\C$, the enhanced categories $\C^{h>}$ and $\C^{h<}$ are
$\kappa$-accessible, and $\Comp_\kappa(\C^{h>}) \cong
\Comp_\kappa(\C)^{h>}$, $\Comp_\kappa(\C^{h<}) \cong
\Comp_\kappa(\C)^{h<}$.
\end{exa}

\begin{defn}
An enhanced functor $\gamma:\E \to \E'$ between $\kappa$-accessible
enhanced categories $\E$, $\E'$ is {\em $\kappa$-defined} if $\gamma
\cong \iInd_\kappa(\gamma_\kappa)$ for an enhanced functor
$\gamma_\kappa:\Comp_\kappa(\E) \to \Comp_\kappa(\E')$.
\end{defn}

\begin{lemma}\label{acc.def.le}
For any set of accessible categories and functors between them, and
any regular cardinal $\kappa$, there exists a regular cardinal $\mu
> \kappa$ such that all the categories are $\mu$-accessible, and all
the functors are $\mu$-defined.
\end{lemma}

\proof{} By induction, it suffices to consider one accessible
functor $\gamma:\E \to \E'$ between accessible categories. Choose a
regular cardinal $\kappa' > \kappa$ such that $\E$ and $\E'$ are
$\kappa'$-accessible, and then choose $\mu$ such that $\kappa' \trl
\mu$ and $\gamma$ sends $\Comp_{\kappa'}(\E)$ into
$\Comp_\mu(\E')$. Then by Proposition~\ref{ka.mu.prop}, any
$\mu$-compact object in $\E$ is a $\mu$-bounded
$\kappa'$-filtered colimit of $\kappa'$-compact object, and since
$\Comp_\mu(\E') \subset \E'$ is closed under $\mu$-bounded colimits
by Proposition~\ref{enh.filt.prop}, $\gamma$ is $\mu$-defined.
\endproof

\begin{corr}\label{acc.adj.corr}
In any adjoint pair $\lambda:\E \to \E'$, $\rho:\E' \to \E$ of
enhanced functors between accessible enhanced categories, both
$\lambda$ and $\rho$ are accessible.
\end{corr}

\proof{} The left-adjoint $\lambda$ is trivially accessible by
Example~\ref{enh.adj.exa.exa}. For $\rho$, choose $\mu$ such that
$\E$, $\E'$ are $\mu$-accessible and $\lambda$ is $\mu$-defined,
with the Yoneda embeddings $\Y_\mu(\E):\E \to
\Comp_\mu(\E)^\iota\sSets^h$, $\Y_\mu(\E'):\E' \to
\Comp_\mu(\E')^\iota\sSets^h$, and note that $\Y_\mu(\E) \circ \rho
\cong \lambda^* \circ \Y_\mu(\E')$ is then $\mu$-accessible, so that
$\rho$ is also $\mu$-accessible.
\endproof

\begin{defn}
For any regular cardinal $\kappa$, an enhanced category $\E$ is {\em
  $\kappa$-presentable} if it is $\kappa$-accessible and
$\kappa$-cocomplete, and an enhanced category $\E$ is {\em
  presentable} if it if $\kappa$-presentable for some regular
cardinal $\kappa$.
\end{defn}

\begin{lemma}\label{enh.pres.le}
For any regular cardinal $\kappa$, a $\kappa$-presentable enhanced
category $\E$ is cocomplete, and conversely, a cocomplete enhanced
category $\E$ is $\kappa$-presentable if and only if it admits a
generating small enhanced subcategory $\C \subset \Comp_\kappa(\E)$
in the sense of Definition~\ref{enh.gen.def}. A $\kappa$-presentable
category is complete, and the full embedding $\E \to
\Comp_\kappa(\E)^\iota\sSets^h$ identifies $\E$ with the enhanced
full subcategory $\Comp_\kappa(\E)^\iota_{ex}\sSets^h \subset
\Comp_\kappa(\E)^\iota\sSets^h$ spanned by enhanced functors
$X:\Comp_\kappa^\iota \to \sSets^h$ such that $X^\iota$ is
$\kappa$-right-exact.
\end{lemma}

\proof{} For the first claim, by Lemma~\ref{coli.kr.le}, an enhanced
category $\E$ is cocomplete iff it admits colimits of enhanced
functors $\Unf(J) \to \E$, $J \in \Posf^+$, and any such $J$ is a
$\kappa$-filtered colimit of its $\kappa$-bounded
subsets. Conversely, let $\C_\kappa \subset \E$ be the smallest
$\kappa$-cocomplete full enhanced subcategory containing $\C$, with
the embedding functor $\gamma:\C_\kappa \to \E$. Then $\C_\kappa$ is
still small, it is generating by Example~\ref{gen.CC.exa}, and by
Proposition~\ref{enh.filt.prop}, $\gamma$ still factors through
$\Comp_\kappa(\E)$, so by adjunction, $\gamma_!\Y(\C_\kappa):\E \to
\C_\kappa^\iota\sSets^h$ is $\kappa$-filtered-cocomplete, thus
factors through $\iInd_\kappa(\C_\kappa)$. We then have an adjoint
pair of fully faithful enhanced functors $\iInd_\kappa(\C_\kappa)
\to \E$, $\E \to \iInd_\kappa(\C_\kappa)$, so both are equivalences.

For the second claim, $\Comp_\kappa(\E)^\iota_{ex}\sSets^h \subset
\Comp_\kappa(\E)^\iota\sSets^h$ is manifestly closed under limits,
thus complete, and the fully faithful embedding
$\gamma_!\Y(\Comp_\kappa):\E \to \Comp_\kappa(\E)^\iota\sSets^h$
factors through $\Comp_\kappa(\E)^\iota_{ex}\sSets^h$, so it
suffices to check that $\gamma_!\Y(\Comp_\kappa):\E \to
\Comp_\kappa(\E)^\iota_{ex}\sSets^h$ is essentially surjective.  By
\eqref{com.Y.col} and Lemma~\ref{env.le}, this amounts to
checking that $\Comp_\kappa(\E) /^h e$ is $\kappa$-filtered for any
$e \in \E_\ppt$; however, by Proposition~\ref{enh.filt.prop} and
Lemma~\ref{I.fibb.le}, it is actually $\kappa$-cocomplete.
\endproof

\begin{corr}\label{pres.corr}
The enhanced categories $\sSets^h$ and $\cCat^h$ are
$\kappa$-present\-able for any regular cardinal $\kappa$, and so is
the enhanced category $\cCat^h \bbi^h \E$ for any small enhanced
category $\E$.
\end{corr}

\proof{} For $\sSets^h$, $\{\ppt^h\} \subset \sSets^h$ is generating,
and $\ppt^h \in \sSets^h$ is $\kappa$-compact by
Lemma~\ref{enh.comp.le}. For $\cCat^h$, Corollary~\ref{enh.del.corr}
implies that $\Unf(\Posff)$, with the full embedding
$\Psi:\Unf(\Posff) \subset \cCat^h$, is generating. Moreover, by
Proposition~\ref{enh.filt.prop}, the property of being $\Psi$-exact
in the sense of Definition~\ref{Psi.def} is closed under
$\kappa$-filtered colimits, so $\kappa$-filtered colimits in
$\Posff^o_h\sSets^h$ induce $\kappa$-filtered colimits in
the full subcategory $\cCat^h \subset
\Posff^o_h\sSets^h$. This means that $\Unf(J)$ is
$\kappa$-compact for any $J \in \Posff$, so that $\Psi(\Unf(\Posff))
\subset \Comp_\kappa(\cCat^h)$. For $\cCat^h \bbi^h \E$, it is
cocomplete by Example~\ref{cat.E.coc}, and then the full embedding
of Corollary~\ref{enh.del.corr} immediately extends to a
$\kappa$-filtered-cocomplete full embedding $\cCat^h \bbi^h \E \subset
(\Unf(\Posff) \bbi^h \E)^\iota\sSets^h$.
\endproof

By Lemma~\ref{env.le}, for any $\kappa$-accessible enhanced
categories $\E$, $\E'$, an enhanced functor $\Comp_\kappa(\E) \to
\E'$ extends to a $\kappa$-accessible enhanced functor $\E \to \E'$,
uniquely up to a unique isomorphism, and any $\kappa$-accessible
enhanced functor is of this form. Therefore these enhanced functors
form a well-defined enhanced category $\fFun^h_\kappa(\E,\E') =
\fFun^h(\Comp_\kappa(\E),\E')$. For any $\kappa \leq \mu$ such that
both $\E$ and $\E'$ are both $\kappa$-accessible and
$\mu$-accessible, a $\kappa$-accessible enhanced functor $\E \to
\E'$ is trivially $\mu$-accessible, so we have a natural full
embedding $\fFun^h_\kappa(\E,\E') \subset \fFun^h_\mu(\E,\E')$;
taking the filtered $2$-colimit with respect to these embeddings, we
obtain an enhanced category
\begin{equation}\label{ffun.acc}
\fFun^h_{acc}(\E,\E') = \bigcup_\kappa \fFun^h_\kappa(\E,\E'),
\end{equation}
where the union is over all regular cardinals $\kappa$ such that
$\E$ and $\E'$ are $\kappa$-accessible.

\begin{prop}\label{acc.prop}
For any accessible enhanced categories $\C$, $\C_0$, $\C_1$ and
accessible enhanced functors $\C_0,\C_1 \to \C$, there exists a
semicartesian square \eqref{refl.semi.sq} with accessible $\C_{01}$
and $\gamma_l$, $l=0,1$, and for any cardinal $\kappa'$, there
exists a regular cardinal $\kappa > \kappa'$ such $\Comp_\kappa(-)$
sends this square to a semicartesian square. Moreover, for any such
semicartesian square, and any other commutative square
\eqref{enh.h.sq} with accessible $\C_{01}'$ and $\gamma_l'$,
$l=0,1$, there exists an accessible enhanced functor
$\alpha:\C'_{01} \to \C_{01}$ equipped with isomorphisms
$a_l:\gamma_l \circ \alpha \cong \gamma'_l$, $l=0,1$, and the triple
$\langle \alpha,a_0,a_1 \rangle$ is unique up to a (non-unique)
isomorphism.
\end{prop}

\proof{} To construct a semicartesian square \eqref{refl.semi.sq},
choose $\kappa$ such that all enhanced categories and functors are
$\kappa$-accessible, and apply Corollary~\ref{ind.prod.corr}. To
check the universal property, increase $\kappa$ if needed so that
$\C_{01}'$, $\gamma_0'$ and $\gamma_1'$ are also
$\kappa$-accessible, replace $\C_{01}'$ with
$\Comp_\kappa(\C_{01}')$, and apply Corollary~\ref{sq.h.corr}.
\endproof

\begin{corr}\label{acc.corr}
For any regular cardinal $\kappa$ and $\kappa$-accessible enhanced
categories $\E$, $\E'$, the enhanced category
$\fFun^h_\kappa(\E,\E')$ of \eqref{ffun.acc} is accessible.
\end{corr}

\proof{} Replacing $\E$ with $\Comp_\kappa(\E)$, we may assume that
$\E$ is small, and then by Lemma~\ref{enh.cocart.le},
Example~\ref{enh.univ.exa} and Proposition~\ref{acc.prop}, we may
further assume that $\E \cong \Unf(J^o)$ for some $J \in
\Posf^+$. Then choose $\mu$ such that $\E'$ is $\mu$-accessible and
$J$ is $\mu$-bounded, and note that by Corollary~\ref{fun.ind.corr},
we have $J^o_h\E' \cong J^o_h\iInd_\mu(\Comp_\mu(\E')) \cong
\iInd_\mu(J^o_h\Comp_\mu(\E'))$.
\endproof

\begin{lemma}\label{acc.aug.le}
For any $J \in \Posf^+$, an enhanced category $\C$ equipped with an
enhanced functor $\C \to \Unf(J)$ that is either an augmentation or
a coaugmentation is accessible if and only if all the fibers $\C_j$,
$j \in J$ and all the transition functors $\C_{j'} \to \C_j$ are
accessible. A $J$-augmented enhanced category $\C$ is accessible iff
so is the transpose $J$-caougmented enhanced category $\C_{h\perp}$,
and if this holds, for any cardinal $\kappa'$, there exists a
regular cardinal $\kappa > \kappa'$ such that $\Comp_\kappa(\C)$ is
$J$-augmented, $\C \cong \iInd_\kappa\Comp_\kappa(\C))$, and
$\C_{h\perp} \cong \iInd_\kappa(\Comp_\kappa(\C)_{h\perp})$.
\end{lemma}

\proof{} If $\C$ is accessible, then $\C_j$ are accessible by
Proposition~\ref{acc.prop}, and the transition functors are
accessible by Proposition~\ref{acc.prop} and
Corollary~\ref{acc.adj.corr}. Conversely, if $\C_j$ and the
transition functors are accessible, then by Lemma~\ref{acc.def.le},
we can choose a single regular cardinal $\kappa$ such that $\C_j$,
$j \in J$ are $\kappa$-accessible, the transition functors are
$\kappa$-defined, and $J$ is $\kappa$-bounded. Then by
Lemma~\ref{sec.ind.bis.le}, all the embeddings $e(j):\C_j \to \C$
are $\kappa$-filtered-right-exact, so that $\C$ is
$\kappa$-filtered-cocomplete. Moreover, an enhanced object $c$ in
some $\C_j \subset \C$ is $\kappa$-compact in $\C$ iff it is
$\kappa$-compact in $\C_j$, so that $\C \cong
\iInd_\kappa(\Comp_\kappa(\C))$ is accessible, and $\Comp_\kappa(\C)
\to \Unf(J)$ is an augmentation resp.\ coaugmentation with fibers
$\Comp_\kappa(\C)_j \cong \Comp_\kappa(\C_j)$, and the full
embedding $\Comp_\kappa(\C) \to \C$ is augmented resp.\ coaugmented.
\endproof

\begin{corr}\label{acc.prod.corr}
In the situation of Proposition~\ref{acc.prop}, if $\C$ and $\C_l$,
$\gamma_l$, $l=0,1$ are $J$-augmented or coaugmented for some $J \in
\Posf^+$, then so is $\C_{01}$ and $\gamma'_l$, $l=0,1$.
\end{corr}

\proof{} By Lemma~\ref{acc.aug.le}, Lemma~\ref{acc.def.le} and
Proposition~\ref{acc.prop}, we can choose $\kappa$ so that all the
categories in the semicartesian square \eqref{refl.semi.sq} are
$\kappa$-accessible, the functors are $\kappa$-defined,
$\Comp_\kappa(-)$ sends the square to a semicartesuan square, and
all the categories $\Comp_\kappa(\C)$, $\Comp_\kappa(\C_l)$, $l=0,1$
and enhanced functors $\gamma_l$, $l=0,1$ between them are
$J$-augmented resp,\ $J$-coaugmented. It then remains to apply the
augmented part of Lemma~\ref{sq.h.le} and
Lemma~\ref{sec.ind.bis.le}.
\endproof

With Proposition~\ref{acc.prop} and its corollaries, most of the
results of Subsection~\ref{enh.co.subs} and
Subsection~\ref{enh.fib.subs} generalize to accessible categories
with the same proofs. In particular, for any accessible enhanced
functor $\gamma:\C_0 \to \C_1$ between accessible categories, we
have the accessible cylinder $\Cyl_h(\gamma)$ and dual cylinder
$\Cyl^\iota_h(\gamma)$ defined by the semicartesian squares
\eqref{enh.cyl.sq}, with the corresponding factorizations
\eqref{enh.cyl.facto}, and by Corollary~\ref{acc.prod.corr},
$\Cyl_h(\gamma)$ resp.\ $\Cyl^\iota_h(\gamma)$ is $[1]$-coaugmented
resp.\ $[1]$-augmented, so that Lemma~\ref{enh.cyl.le} holds with
the same proof. Moreover, by virtue of Lemma~\ref{acc.aug.le}, we
have $\Comp_\kappa(\Cyl^h(\gamma)) \cong
\Cyl^h(\Comp_\kappa(\gamma))$ for a cofinal family of regular
cardinals $\kappa$, and similarly for $\Cyl^\iota_h(\gamma)$, and
then Lemma~\ref{enh.cyl.o.le} immediately implies that the
corresponding squares for $\gamma$ are semicocartesian with respect
to accessible enhanced functors. Then one defines the accessible
enhanced comma-categories $\C_0 /^h_\gamma \C_1$, $\C_1
\setminus^h_\gamma \C_0$ by semicartesian squares
\eqref{enh.rl.comma}, with the corresponding decompositions
\eqref{enh.comma.facto}, and Lemma~\ref{cyl.comma.le} holds with the
same proof. The same goes for Lemma~\ref{adj.eq.le} and its
Corollary~\ref{enh.adj.corr}. For accessible functors that are
enhanced fibrations, we have Lemma~\ref{pi.id.le}, and its
counterpart for accessible cofibrations. We then have
Corollary~\ref{eta.la.corr}, Lemma~\ref{cl.fib.le} and its
Corollary~\ref{enh.si.ta.corr}, Lemma~\ref{enh.fib.sq.le},
Lemma~\ref{enh.right.fib.le}, Lemma~\ref{enh.adm.cof.le},
Lemma~\ref{enh.fib.le}, and their counterparts for cofibrations
obtained by literally dualizing the proofs. Enhanced fibers of
accessible enhanced functors are accessible, and so are the enhanced
transition functors for enhanced fibrations and cofibrations. The
material of Subsection~\ref{enh.fiber.subs} also carries over
literally to the accessible setting; in particular, this includes
Lemma~\ref{enh.ff.fib.le}, Lemma~\ref{enh.fib.adj.le},
Corollary~\ref{adj.comma.corr} and Corollary~\ref{enh.bifib.corr}.

\begin{remark}\label{acc.rem}
As a bit of a warning, one operation on small enhanced categories
that definitely does {\em not} extend to accessible categories is
the involution $\iota$ sending $\C$ to its enhanced-opposite
$\C^\iota$. In fact, already for the unenhanced category $\Sets$,
the opposite category $\Sets^o$ is not accessible. Thus in the
accessible world, passing to the enhanced-opposite category is not
always possible, and has to be carefully justified. As a consequence
of this, one cannot expect a meaningful generalization of the Yoneda
embedding. Another consequence is that facts about enhanced
fibrations and cofibration no longer simply imply each other, and
have to be proved separately.
\end{remark}

\subsection{Tame fibrations and cofibrations.}

By Corollary~\ref{sm.acc.corr}, any small Karoubi-closed enhanced
category is accessible, or in other words, any small enhanced
category is accessible ``up to Karoubi closure''. To handle the
remaining discrepancy, one can introduce the following.

\begin{defn}
An enhanced category $\C$ is {\em tame} if its Karoubi envelope
$\Env(\C,\pP)$ is accessible, and an enhanced functor $\gamma:\C \to
\C'$ between tame enhanced categories is {\em tame} if
$\Env(\gamma,\pP):\Env(\C,\pP) \to \Env(\C',\pP)$ is accessible.
\end{defn}

Unfortunately, since $\Env(-,\pP)$ is badly behaved, there is a
price to pay for this generality: Proposition~\ref{acc.prop} does
not hold anymore. Namely, since the embedding $\C \to \Env(\C,\pP)$
is enhanced fully faithful for any $\C$, we can apply
Proposition~\ref{acc.prop} to construct a semicartesian product
$\C_{01} = \C_0 \times^h_{\C} \C_1$ for any tame enhanced categories
$\C$, $\C_0$, $\C_1$ and tame enhanced functors $\gamma_l:\C_l \to
\C$, $l=0,1$, but it is not true in general that $\C_{01}$ is
tame. By Lemma~\ref{kar.C.E.le}, for enhanced fibrations and
cofibrations, this problem does not occur.

\begin{lemma}\label{tame.prod.le}
For any tame enhanced functor $\pi:\C \to \E$ between tame enhanced
categories that is an enhanced fibration or cofibration in the sense
of Definition~\ref{enh.fib.def}, the Karoubi envelope
$\Env(\pi,\pP)$ is an enhanced fibration resp.\ cofibration. For any
semicartesian square \eqref{C.E.sq} of enhanced categories such that
$\C$, $\E$, $\E'$, $\pi$ and $\phi$ are tame, and $\pi$ is an
enhanced fibration or cofibration, $\C'$, $\pi'$ and $\gamma$ are
also tame, and $\pi'$ is also an enhanced fibration resp.\ enhanced
cofibration.
\end{lemma}

\proof{} For the first claim, since $\Env(-,\pP)$ commutes with
$2$-filtered colimits of full embeddings, we may assume that $\E$ is
small, and then by Lemma~\ref{fib.uni.le}, we may assume the same
for $\C$. Then $\Env(-,\pP)$ commutes with taking enhanced right
comma-categories by by Lemma~\ref{ind.ka.mu.le} and
Remark~\ref{ind.ka.mu.rem}, and since it also commuted with
$\iota:\Cat^h \to \Cat^h$, $\C \mapsto \C^\iota$, the same holds for
left comma-categories. The claim then immediately follows from
Lemma~\ref{env.adj.le}. For the second claim, $\C'$ is tame by
Lemma~\ref{kar.C.E.le}, and then again by Lemma~\ref{fib.uni.le}, it
suffices to check that $\pi'$ is a fibration when $\C$ and $\E$ are
small. This is Lemma~\ref{enh.fib.le}.
\endproof

\begin{exa}\label{tame.comma.exa}
For any enhanced category $\E$, with the enhanced arrow category
$\Ar^h(\E)$, the projections $\sigma,\tau:\Ar^h(\E) \to \E$ are an
enhanced fibration resp.\ cofibration (apply Lemma~\ref{cl.fib.le}
to the $[1]$-coaugmented enhanced category $\E \times^h \Unf([1])$
and its enhanced-opposite). Therefore Lemma~\ref{tame.prod.le}
applies to semicartesian squares \eqref{enh.rl.comma}, and for any
tame enhanced functor $\gamma:\C \to \E$ between tame enhanced
categories, we have tame comma-categories $\C /^h_\gamma \E$, $\E
\setminus^h_\gamma \C$. Moreover, Lemma~\ref{tame.prod.le} also
applies to the squares \eqref{enh.cyl.sq}, so that the enhanced
cylinder $\Cyl_h(\gamma)$ and the enhanced dual cylinder
$\Cyl^\iota_h(\gamma)$ are also tame.
\end{exa}

To construct tame enhanced fibrations and cofibrations, note that
for any regular cardinals $\kappa' < \kappa$ and $\kappa'$-bounded
Karoubi-closed small enhanced category $\E \in \Cat^h_{\kappa'}$, we
have $\iInd_\kappa(\E) = \E$ by Corollary~\ref{sm.acc.corr}, so as
in \eqref{ind.J}, any enhanced functor $\phi:\C \to \E$ gives rise
to an enhanced functor
\begin{equation}\label{ind.E}
\iInd_\kappa(\phi):\iInd_\kappa(\C) \to \E = \iInd_\kappa(\E).
\end{equation}
We then have the following weaker version of
Lemma~\ref{sec.ind.bis.le}.

\begin{lemma}\label{ind.E.le}
For any regular cardinals $\kappa' < \kappa$, $\kappa'$-bounded
Karoubi-closed small enhanced category $\E$, and enhanced functor
$\phi:\C \to \E$ from a small enhanced category $\C$ that is an
enhanced fibration or an enhanced cofibration, \eqref{ind.E} is an
enhanced fibration resp.\ cofibration. Morever, for any
Karoubi-closed enhanced category $\E' \in \Cat^h_{\kappa'}$ and
enhanced functor $\gamma:\E' \to \E$, the commutative square
\begin{equation}\label{ind.E.sq}
\begin{CD}
\iInd_\kappa(\gamma^*\C) @>>> \iInd_\kappa(\C)\\
@VVV @VVV\\
\E' @>{\gamma}>> \E
\end{CD}
\end{equation}
is semicartesian, and $\sSec^h(\E',\Ind_\kappa(\C)) \cong
\iInd_\kappa(\sSec^h(\E',\C))$.
\end{lemma}

\proof{} For the first claim, $\iInd_\kappa(-)$ commutes with taking
right and left enhanced comma-categories by
Lemma~\ref{sec.ind.bis.le}, so the claim immediately follows from
Lemma~\ref{env.adj.le}. For the second claim, as in
Lemma~\ref{sec.ind.bis.le}, it suffices to show that
\eqref{ind.E.sq} is semicartesian in the case when $\E' = \ppt$, so
that $\gamma = \eps^h(e)$ is the embedding onto an enhanced object
$e \in \E_\ppt$. Then if we let $\C' \to \Unf(\V^o)$ be the
$\V^o$-coaugmented enhanced category with fibers $\ppt^h$, $\C$,
$\E$ and transition functors $\eps^h(e):\ppt^h \to \E$, $\phi:\C \to
\E$, Lemma~\ref{sec.ind.le} shows that
$\iInd_\kappa(\sSec^h(\V^o,\C')) \cong
\sSec^h(\V^o,\iInd_\kappa(\C'))$, the full embedding $\C_e \to
\sSec(\V^o,\C')$ induces a full embedding
\begin{equation}\label{ind.E.eq}
\iInd_\kappa(\C_e) \to \iInd_\kappa(\C)_e,
\end{equation}
and we need to check that \eqref{ind.E.eq} is essentially
surjective. Indeed, by Lemma~\ref{ind.filt.le}, enhanced objects in
$\iInd_\kappa(\C)_e$ are given by enhanced functors $X:\C^\iota \to
\sSets^h$ such that $\C X$ is $\kappa$-filtered and $\phi^o_!X \cong
\Y(e)$. Take such a functor $X$, with the cofinal enhanced functor
$\C X \to \E /^h e$ of Example~\ref{EX.exa}. If $\phi$ is an
enhanced cofibration, let $\C' = \C X /^h (\E /^h e)$. Then the
projection $\pi:\C X \to \C$ extends to an enhanced functor
$\pi':\C' \to \C /^h e$ cocartesian over $\E /^h e$, we have the
full embedding $\eta:\C X \to \C'$ that is cofinal by
Corollary~\ref{comma.filt.corr}, and since $e$ is the terminal
enhanced object in $E /^h e$, we also have $\C X \cong \C'_e$, so
that we obtain another full embedding $\psi:\C X \cong \C'_e \to
\C'$ that is left-admissible, hence also cofinal. We then have
\begin{equation}\label{X.Y.p.eq}
X \cong \colimh_{\C X}\Y(\C) \circ \pi \cong \colimh_{\C'}
\Y(\C) \circ \pi' \cong \colimh_{\C'_e}\Y(\C) \circ \pi' \circ \psi,
\end{equation}
and the right-hand side manifestly lies in the essential image of
\eqref{ind.E.eq}. Similarly, if $\phi$ is an enhanced fibration, let
$\C' = (\E /^h e) \setminus^h \C X$, and extend $\pi$ to an enhanced
functor $\pi':\C' \to \C /^h e$ cartesian over $\E /^h e$. Then we
again have a full embedding $\eta:\C X \to \C'$ and another full embedding
$\psi:\C'_e \to \C'$, $\eta$ is left-admissible, hence cofinal, so
to conclude that \eqref{X.Y.p.eq} still holds, we need to check that
$\psi$ is cofinal as well. The argument for that is essentially the
same as in Corollary~\ref{comma.filt.corr}; we leave it to the
reader.
\endproof

\begin{lemma}\label{tame.fib.le}
Assume given an accessible enhanced functor $\phi:\C \to \E$ to a
small enhanced category $\E$, and assume that $\phi$ is either an
enhanced fibration or an enhanced cofibration. Then for any cardinal
$\kappa'$, there exists a regular cardinal $\kappa > \kappa'$ such
that $\Comp_\kappa(\C) \to \E$ is an enhanced fibration
resp.\ cofibration, and an enhanced object $c \in \C_\ppt$ is
$\kappa$-compact iff it is $\kappa$-compact as an object in the
enhanced fiber $\C_{\phi(c)}$.
\end{lemma}

\proof{} By Lemma~\ref{acc.def.le} and Corollary~\ref{sm.acc.corr},
we may choose a regular cardinal $\kappa > \kappa'$ such that the
enhanced fibers $\C_e$, $e \in \E_\ppt$ are $\kappa$-compact, the
transition functors are $\kappa$-defined, and
$\E=\iInd_\kappa(\E)$. Then as in Lemma~\ref{ind.E.le}, all
$\kappa$-filtered colimits in $\C$ can be replaced with cofinal
$\kappa$-filtered colimits that factor through an enhanced fiber
$\C_e$, $c \in \E_\ppt$, and by adjunction, this immediately implies
that $c \in \C_\ppt$ is $\kappa$-compact in $\C$ iff it is
$\kappa$-compact in $\C_{\phi(c)}$. Since all the transition
functors are $\kappa$-defined, $\Comp_\kappa(\C) \subset \C$ is then
a subfibration resp.\ subcofibration by Lemma~\ref{fib.uni.le} and
Lemma~\ref{enh.ff.fib.le}.
\endproof

A natural next question is whether tame enhanced fibrations and
cofibrations admits some form of the Grothendieck construction, but
this question is delicate and presents obvious problems. Firstly,
for fibrations, we cannot pass to the opposite categories, as
explained in Remark~\ref{acc.rem}. Secondly, tame categories do not
themselves form an enhanced category. Indeed, even if we fix a
regular cardinal $\kappa$ and only consider $\kappa$-accessible
enhanced categories and enhanced functors, the functors categories
$\fFun^h_\kappa(\C,\C')$ are usually large. We thus leave the
general Grothendieck construction aside, and limit ourselves to one
of its applications --- namely, transpose enhanced cofibrations and
fibrations.

\begin{prop}\label{tame.perp.prop}
Assume given a small enhanced category $\E$. Then for any tame
enhanced fibration $\C \to \E$, there exists a unique tame enhanced
cofibration $\C' = \C_{h\perp} \to \E^\iota$, and dually, for any
tame enhanced cofibration $\C' \to \E^\iota$, there exists a unique
tame enhanced fibration $\C = {\C'}^{h\perp} \to \E$ such that the
squares \eqref{perp.sq} are enhanced-cocartesian with respect to
tame enhanced functors. Moreover, for any enhanced functor
$\gamma:\E' \to \E$ from a small enhanced category $\E'$, we have
$(\gamma^*\C)_{h\perp} \cong \gamma^{\iota *}\C_{h\perp}$.
\end{prop}

\proof{} Assume first that we are given a tame fibration $\C \to
\E$, and $\C$ and $\E$ are Karoubi-closed, thus accessible. Then by
Lemma~\ref{tame.fib.le}, we have $\C \cong
\iInd_\kappa(\Comp_\kappa(\C))$ for some $\kappa$ such that
$\E=\iInd_\kappa(\E)$, $\E^\iota = \iInd_\kappa(\E^\iota)$ and
$\Comp_\kappa(\C) \to \E$ is an enhanced fibration. Then by
Lemma~\ref{enh.perp.le}, the enhanced cofibration
$\Comp_\kappa(\E)_{h\perp} \to \E^\iota$ fits into an
enhanced-cocartesian square \eqref{perp.sq}, and by
Lemma~\ref{ind.E.le}, it it induces an analogous commutative square
for the enhanced category $\C'_\kappa =
\iInd_\kappa(\Comp_\kappa(\C)_{h\perp})$ that is
enhanced-cocartesian with respect to $\kappa$-accessible enhanced
functors. Moreover, the category $\C'_\kappa$ and the square are
functorial with respect to enhanced functors $\gamma:\E' \to \E$ as
required. Then for any $\mu$ such that $\kappa \trl \mu$,
$\Comp_\mu(\C)_{h\perp}$ also fits into an enhanced-cocartesian
square \eqref{perp.sq}, and by universality, we have a comparison
functor $\Comp_\mu(\C)_{h\perp} \to \Comp_\mu(\C'_\kappa)$,
cocartesian over $\E^\iota$. By functoriality with respect to
$\gamma$, this functor is an equivalence over any enhanced object $e
\in \E_\ppt$, thus an equivalence. We conclude that the square
\eqref{perp.sq} for $\C'_\kappa$ is also enhanced-cocartesian with
respect to $\mu$-accessible functors, and since $\mu$ is arbitrary,
it is enhanced-cocartesian with respect to all accessible
functors. This implies uniqueness.

If we are given $\C'$ rather than $\C$, the argument is exactly the
same. Moreover, in either case, by functoriality with respect to
$\gamma$, $\C$ and $\C'$ have exactly the same enhanced
objects. Then to finish the proof, it remains to consider the case
when $\C$ or $\C'$ and $\E$ are only tame, not accessible; in this
case, it suffices to take $\Env(\C,\pP)_{h\perp}$
resp.\ $\Env(\C',\pP)^{h\perp}$, and consider the full subcategory
spanned by enhanced objects in $\C$ resp.\ $\C'$.
\endproof

Finally, let us consider the case when the base of a tame enhanced
fibration or cofibration is not small. Here all we have is one
non-example and one half-example.

\begin{exa}
The enhanced cofibration $\cCat^h_\idot \to \cCat^h$ described in
Example~\ref{enh.idot.exa} is {\em not} tame. Indeed, were it tame,
the enhanced fibers of its Karoubi closure would have been
$\kappa$-accessible for some fixed regular cardinal
$\kappa$. However, these enhanced fibers are all small enhanced
categories, so it is not possible to fix a single $\kappa$ for all
of them.

What one can do is fix a regular cardinal $\kappa$, and construct an
enhanced cofibration $\cCat^h_{\kappa\idot} \to \cCat^h$ with
enhanced fibers $\iInd_\kappa(\C)$, $\C \in \Cat^h$. It is possible
to do this, and the resulting category $\cCat^h_{\kappa\idot}$ is
$\kappa$-accessible; however, it is not clear whether so is the
projection $\cCat^h_{\kappa\idot} \to \cCat^h$.
\end{exa}

\begin{prop}\label{pres.prop}
For any regular cardinal $\kappa$ and $\kappa$-filtered-cocomplete
enhanced category $\E$, the enhanced category $\cCat^h \bb^h \E$ is
is $\kappa$-filtered-cocom\-ple\-te, and the enhanced fibration
\eqref{no.E} is $\kappa$-right-exact. If $\E$ is
$\kappa$-present\-able, then $\cCat^h \bb^h \E$ is
$\kappa$-presentable, and \eqref{no.E} is enhanced-right-exact. If
$\E$ is accessible, then for any small Karoubi-closed full enhanced
subcategory $\I \subset \cCat^h$, the enhanced full subcategory $\I
\bb^h \E \subset \cCat^h \bb^h \E$ is accessible.
\end{prop}

\proof{} The enhanced category $\cCat^h$ is $\kappa$-presentable,
hence $\kappa$-accessible by Corollary~\ref{pres.corr}. By
definition, for any small enhanced category $\I$, an enhanced
functor $E:\I \to \cCat^h \bb^h \E$ is determined by the pair
$\langle E,EX \rangle$ of its composition $X:\I \to \cCat^h$ with
\eqref{no.E}, and an enhanced functor $EX:\I^\perp X \to \E$. If we
let $\C = \colimh_{\I} X$, then \eqref{enh.colim.sq} provides an
enhanced functor $\gamma:\I^\perp X \to \C$, and then $\colimh_{\I}E
= \langle \C,\gamma_!EX \rangle$, provided the left Kan extension
exists. Then for the first claim, we need to check that $\gamma_!EX$
exists if $\E$ is $\kappa$-filtered-cocomplete and $\I$ is
$\kappa$-filtered. Indeed, composing $X$ with the embedding
\eqref{enh.del} gives an enhanced functor $X':\Unf(\Delta^o)
\times^h \I \to \sSets^h$ and an enhanced functor
$EX':(\Unf(\Delta^o) \times^h \I)^\perp X \to \E$; then $\gamma$
induces an enhanced functor $\gamma':(\Unf(\Delta^o) \times^h
\I)^\perp X' \to (\Delta_h\C)^\iota$ cocartesian over
$\Unf(\Delta^o)$, and by Lemma~\ref{enh.del.loc.le}, it suffices to
check that $\gamma'_!EX'$ exists. But then, \eqref{enh.del}
preserves filtered colimits, and if one computes $\gamma'_!$ by
\eqref{enh.cokan.eq}, all the comma-fibers are $\kappa$-filtered by
Lemma~\ref{filt.hyper.le}.

For the second claim, $\gamma_!$ always exists since presentable
enhanced categories are enhanced-cocomplete, so by
Lemma~\ref{enh.pres.le}, it suffices to find a small generating full
subcategory in $\cCat^h \bb^h \E$ spanned by $\kappa$-compact
objects. By Corollary~\ref{enh.del.corr}, $\Posff \bb^h
\Comp_\kappa(\E) \subset \cCat^h \bb^h \E$ does the job.

Finally, for the third claim, by Corollary~\ref{acc.corr}, we may
choose a regular cardinal $\mu > \kappa$ that is large enough so
that $\I$ is $\mu$-bounded and $\mu$-accessible, $\fFun^h(\C,\E)$ is
$\mu$-accessible for any enhanced object $\C \in \I_{\ppt}$, and the
pullback functor $\gamma^*:\fFun^h(\C',\E) \to \fFun^h(\C,\E)$ is
$\mu$-accessible for any enhanced morphism $\gamma:\C \to \C'$ in
$\I$. Then $\I \bb^h \E$ is $\mu$-filtered-cocomplete by the first
claim, and $\Comp_\mu(\I \bb^h \E) \subset \I \bb^h \E$ is an
enhanced subfibration over $\I$ with enhanced fibers
$\Comp_\mu(\fFun^h(\C,\E))$, $\C \in \I_{\ppt}$, while the fully
faithful embedding $\iInd_\mu(\Comp_\mu(\I \bb^h \E)) \to \I \bb^h
\E$ is cartesian over $\I$, and essentially surjective on each
enhanced fiber. Therefore it is essentially surjective.
\endproof

It is probably {\em not} true that $\cCat^h \bb^h \E$ is accessible
for an arbitrary accessible enhanced category $\E$ (roughly
speaking, were it $\kappa$-accesible for some regular cardinal
$\kappa$, $\fFun^h(\C,\E)$ would have been $\kappa$-accessible for
any $\kappa$-bounded small enhanced category $\E$, and by
Corollary~\ref{acc.corr}, it is only $\kappa'$-accessible for some
$\kappa' > \kappa$). However, even the restricted statement given in
Proposition~\ref{pres.prop} has its uses. In particular, for any
enhanced cofibration $\C' \to \C$ between small enhanced categories,
and any accessible enhanced category $\E$, we can still define the
relative functor category $\fFun^h(\C'|\C,\E)$ by the semicartesian
square 
\begin{equation}\label{rel.acc.ffun.sq}
\begin{CD}
\fFun^h(\C'|\C,\E) @>>> \cCat^h \bb^h \E\\
@VVV @VV{\pi}V\\
\C @>{X}>> \cCat^h,
\end{CD}
\end{equation}
where $X:\C \to \cCat^h$ corresponds to $\C' \to \C$ by the
covariant enhanced Grothendieck construction, and the square
\eqref{rel.acc.ffun.sq} exists since $X$ factors through a small
enhanced subcategory $\I \subset \cCat^h$. The relative functor
category $\fFun^h(\C'|\C,\E)$ is accessible by
Proposition~\ref{acc.prop}, and has the same universal property with
respect to accessible enhanced functors as the relative functor
category \eqref{rel.ffun.sq} of
Subsection~\ref{enh.rel.subs}. Moreover, the forgetful functor
$\fFun^h(\C'|\C,\E) \to \C$ is an accessible enhanced fibration, and
we have an enhanced fibration $\C' \to \C$ between small enhanced
categories, we can define
\begin{equation}\label{rel.acc.fun.perp}
\fFun^h(\C'|\C,\E) \cong \fFun^h(\C'_{h\perp}|\C^\iota,\E)_{h\perp},
\end{equation}
where $\gamma^\iota:{\E'}^\iota \to \E^\iota$ is the
enhanced-opposite cofibration, and the right-hand side is the
transpose tame enhanced cofibration provided by
Proposition~\ref{tame.perp.prop}. This again has the same universal
property as the relative functor category \eqref{rel.fun.perp} of
Subsection~\ref{enh.rel.subs}.


\bigskip

{\noindent
Affiliations:
\begin{enumerate}
\renewcommand{\labelenumi}{\arabic{enumi}.}
\item Steklov Mathematics Institute (main affiliation).
\item National Research University Higher School of Economics.
\end{enumerate}}

{\small\noindent
{\em E-mail address\/}: {\tt kaledin@mi-ras.ru}
}

\end{document}